%% file: JvN-Functional_Analysis-2nd-ed.tex
\documentclass{my-cambridge7A}

\usepackage{amssymb}
\usepackage{amsmath}
\usepackage{amsthm}
\usepackage{amsfonts}
\usepackage{mathptmx}
\usepackage{natbib}
\usepackage[figuresright]{rotating}
\usepackage{floatpag}
  \rotfloatpagestyle{empty}

\newcommand{\comment}[1]{}

\usepackage[colorlinks, 
pagebackref, citecolor={blue}]{hyperref}


\usepackage{physics}
\usepackage{mathrsfs}
  \renewcommand{\thesubsection}{\arabic{chapter}.\arabic{section}.\alph{subsection}}
\usepackage[active]{srcltx}
\usepackage[shortlabels]{enumitem}
\usepackage[fleqn,tbtags]{mathtools}

\makeatletter
\DeclareRobustCommand\widecheck[1]{{\mathpalette\@widecheck{#1}}}
\def\@widecheck#1#2{%
    \setbox\z@\hbox{\m@th$#1#2$}%
    \setbox\tw@\hbox{\m@th$#1%
       \widehat{%
          \vrule\@width\z@\@height\ht\z@
          \vrule\@height\z@\@width\wd\z@}$}%
    \dp\tw@-\ht\z@
    \@tempdima\ht\z@ \advance\@tempdima2\ht\tw@ \divide\@tempdima\thr@@
    \setbox\tw@\hbox{%
       \raise\@tempdima\hbox{\scalebox{1}[-1]{\lower\@tempdima\box\tw@}}}%
    {\ooalign{\box\tw@ \cr \box\z@}}}
\makeatother

\input{macros.tex}
\allowdisplaybreaks

\newcommand\blfootnote[1]{%
  \begingroup
  \renewcommand\thefootnote{}\footnote{#1}%
  \addtocounter{footnote}{-1}%
  \endgroup
}

\usepackage{pgfplots} 

\pgfplotsset{compat=newest}
\usepgfplotslibrary{fillbetween}

\usepackage{float, tikz, tikz-cd}
\usepackage[verbose]{wrapfig}

\usepackage{needspace}

\tikzcdset{scale cd/.style={every label/.append style={scale=#1},
    cells={nodes={scale=#1}}}}

\usepackage{graphicx}
\usepackage{makeidx}
 \makeindex

  \theoremstyle{plain}
  \newtheorem{theorem}{Theorem}[chapter]
  \newtheorem{lemma}[theorem]{Lemma}
  \newtheorem{proposition}[theorem]{Proposition}
  \newtheorem{corollary}[theorem]{Corollary}

  \newtheorem*{theorem*}{Theorem}
  \newtheorem*{lemma*}{Lemma}
  \newtheorem*{proposition*}{Proposition}
  \newtheorem*{corollary*}{Corollary}
  \newtheorem*{conjecture*}{Conjecture}

  \theoremstyle{definition}
  \newtheorem{definition}[theorem]{Definition}
  \newtheorem{example}[theorem]{Example}
  \newtheorem{remark}[theorem]{Remark}

  \newtheorem*{definition*}{Definition}
  \newtheorem*{example*}{Example}
  \newtheorem*{prob*}{Problem}
  \newtheorem*{remark*}{Remark}
  \newtheorem*{notation*}{Notation}
  \newtheorem*{exer*}{Exercise}

\newenvironment{pictureproof}[1][\proofname]
  {\par\addvspace{\topsep}%
   \pushQED{\qed}\normalfont
   \noindent\textit{#1}\enskip\ignorespaces}
  {\popQED\par\addvspace{\topsep}\ignorespacesafterend}

  \def\makeRRlabeldot#1{\hss\llap{#1}}
  \renewcommand\theenumi{{\rm (\roman{enumi})}}
  \renewcommand\theenumii{{\rm (\alph{enumii})}}
  \renewcommand\theenumiii{{\rm (\arabic{enumiii})}}
  \renewcommand\theenumiv{{\rm (\Alph{enumiv})}}

\begin{document}

\title{Functional Analysis}

\subtitle{Second Edition}

\author{Jan van Neerven \vskip11cm
{\rm \small \begin{quote}This book has been published by Cambridge University Press in the series `Cambridge Studies in Advanced Mathematics'. The present version is free to view and download for personal use only. Not for re-distribution, re-sale or use in derivative works.\newline {\copyright} Jan van Neerven\end{quote}}}

\frontmatter
\maketitle

\setcounter{tocdepth}{1}

\setcounter{page}{5}

\addtocontents{toc}{\protect\setcounter{tocdepth}{2}}
\tableofcontents

\cleardoublepage 
\include{preface-2nd-edition}

\cleardoublepage 
\include{preface-1st-edition}

\cleardoublepage 
\include{how-to-use-this-book}

\cleardoublepage 
\include{notation}

\mainmatter
\cleardoublepage 
\include{ch01-BanachSpaces}

\cleardoublepage 
\include{ch02-ClassicalBanachSpaces}

 \cleardoublepage 
 \include{ch03-HilbertSpaces}

 \cleardoublepage 
 \include{ch04-Duality}

 \cleardoublepage  
  \include{ch05-BoundedOperators}

  \cleardoublepage 
  \include{ch06-Spectrum}

  \cleardoublepage 
  \include{ch07-CompactOperators}

  \cleardoublepage 
  \include{ch08-HilbertOperators}

  \cleardoublepage 
  \include{ch09-SpectralTheorem-bdd}

  \cleardoublepage 
  \include{ch10-UnboundedOperators}

   \cleardoublepage 
   \include{ch11-BoundaryValueProblems}

   \cleardoublepage 
   \include{ch12-LaplaceOperator}

   \cleardoublepage 
   \include{ch13-Semigroups}

   \cleardoublepage  
   \include{ch14-HS-TC}

   \cleardoublepage 
   \include{ch15-QM}

   \cleardoublepage 
   \include{ch16-SymmQuant}

\backmatter

\appendix

 \addtocontents{toc}{\vspace{\baselineskip}}

\include{appendix}

 \endappendix

 \addtocontents{toc}{\vspace{\baselineskip}}

 \renewcommand{\refname}{References}

\cleardoublepage 

 \bookreferences

\bibliography{literature}\label{refs}
\bibliographystyle{cambridgeauthordate}

\cleardoublepage
\phantomsection
\printindex
\end{document}

%% file: macros.tex
\allowdisplaybreaks

\usepackage{newsymbol}

\let\mathcal \undefined
\def\mathcal{\mathscr}
\let\emptyset \undefined
\newsymbol\emptyset    203F
\newsymbol\notle       230A
\newsymbol\notge       230B

\renewcommand*{\ge}{\geqslant}
\renewcommand*{\le}{\leqslant}
\renewcommand*{\geq}{\geqslant}
\renewcommand*{\leq}{\leqslant}

\newcommand\R{{\mathbb R}}
\newcommand\N{{\mathbb N}}
\newcommand\C{{\mathbb C}}
\newcommand\Q{{\mathbb Q}}
\newcommand\Z{{\mathbb Z}}
\newcommand\K{{\mathbb K}}
\newcommand\T{{\mathbb T}}

\newcommand{\al}{\alpha}
\newcommand{\be}{\beta}
\newcommand{\la}{\lambda}
\newcommand{\om}{\omega}
\newcommand{\Om}{\Omega}

\renewcommand{\P}{{\mathbb P}}
\newcommand{\E}{{\mathbb E}}
\renewcommand{\H}{{\mathbb H}}
\newcommand{\HH}{{\mathscr H}}
\newcommand{\GG}{{\mathscr G}}

\newcommand{\calA}{{\mathscr A}}
\newcommand{\calB}{{\mathscr B}}
\newcommand{\calC}{{\mathscr C}}
\newcommand{\calD}{{\mathscr D}}
\newcommand{\calE}{{\mathscr E}}
\newcommand{\calF}{{\mathscr F}}
\newcommand{\calG}{{\mathscr G}}

\newcommand{\calI}{{\mathscr I}}

\newcommand{\calL}{{\mathscr L}}
\newcommand{\calM}{{\mathscr M}}
\newcommand{\calP}{{\mathscr P}}
\newcommand{\calS}{{\mathscr S}}

\newcommand{\calX}{{\mathscr X}}

\newcommand{\dps}{\displaystyle}

\newcommand{\lb}{\langle}
\newcommand{\rb}{\rangle}

\newcommand{\wh}{\widehat}
\newcommand{\wt}{\widetilde}
\newcommand{\ran}{\mathsf{R}}
\renewcommand{\ker}{\mathsf{N}}

\newcommand{\nn}{|\!|\!|}

\newcommand{\limn}{\lim_{n\to\infty}}
\newcommand{\sumn}{\sum_{n\ge 1}}

\newcommand{\limk}{\lim_{k\to\infty}}

\newcommand{\limj}{\lim_{j\to\infty}}

\newcommand{\limN}{\lim_{N\to\infty}}

\newcommand{\iprod}[2]{( #1|#2 )}
\newcommand{\bigiprod}[2]{\bigl( #1\big|#2 \bigr)}

\newcommand{\one}{{\bf 1}}
\DeclareMathOperator*{\esssup}{ess\,sup}
\DeclareMathOperator*{\codim}{codim}
\DeclareMathOperator*{\ind}{ind}
\DeclareMathOperator*{\sign}{sign}
\DeclareMathOperator*{\supp}{supp}

\newcommand{\X}{\mathscr{X}}

\newcommand{\ud}{\,{\rm d}}

\renewcommand{\aa}{\mathfrak{a}}
\newcommand{\bb}{\mathfrak{b}}

\newcommand{\Gr}{\mathsf{G}}
\newcommand{\Dom}{\mathsf{D}}
\newcommand{\Ran}{\mathsf{R}}

\newcommand{\Ker}{\mathsf{N}}

\newcommand{\e}{\varepsilon}
\newcommand{\n}{\Vert}
\newcommand{\s}{^*}
\newcommand{\ov}{\overline}
\newcommand{\ot}{\otimes}
\newcommand{\B}{{\mathscr{B}}}
\newcommand{\M}{{\mathscr{M}}}
\newcommand{\PP}{{\mathscr{P}}}
\renewcommand{\tr}{{\rm{tr}}}
\renewcommand{\var}{{\rm var}}
\newcommand{\F}{{\mathscr F}}

\newcommand{\eps}{\varepsilon}
\newcommand{\rh}{\varrho}
\newcommand{\si}{\sigma}
\newcommand{\inv}{^{-1}}
\newcommand{\dtt}{\frac{{\rm d}t}{t}}

\newcommand{\bfm}{{\bf m}}
\newcommand{\bfn}{{\bf n}}
\newcommand{\hot}{\otimes}
\newcommand{\co}{{\rm co}}

\newcommand{\Op}{\mathrm{Op}}

%% file: preface-2nd-edition.tex
\chapter*{Preface to the Second Edition}

In preparing this second edition, I have taken the opportunity to improve the presentation in several places, to test and revise the problem sets, and to add new material.
The additions include treatments of conditional expectations (Chapter 2), marginal problems and optimal transport (Chapter 4), the Rayleigh--Ritz principle (Chapter 8), the ideal of compact operators (Chapter 9), Sobolev inequalities and the direct variational method (Chapter 11), symmetric forms, Kato's square root theorem, and the landscape function (Chapter 12), Strichartz's estimate (Chapter 13), Markov kernels, quantum correlations, and Grothendieck's inequality (Chapter 15), and the time--energy pair and Weyl quantisation (Chapter 16). Some of these sections benefitted from valuable suggestions by Bas Janssens and Jan Maas.
Much of the new material was written during several retreats in the peaceful village of Herbeumont. It is a pleasure to thank Sonja Cox and Lenny Taelman for making these retreats possible.

Recent advances in artificial intelligence have made it possible to subject a manuscript of this length to systematic searches for typographical errors, inconsistencies, gaps in
exposition, mathematical oversights, and editorial issues.
Successive generations of OpenAI's ChatGPT and Codex, as well as Claude's Sonnet and Fable, assisted with these checks and suggested clarifications, additional detail, occasional simplifications
of proofs, and editorial improvements. I thank Emiel Lorist and Shreyas Kalvankar for setting up the workflows using Codex and the API. Every suggestion was critically examined, and many were substantially revised before being incorporated. Responsibility for the final text rests entirely with me.

I am deeply grateful to Roger Astley and his team at Cambridge University Press for the generosity and support that made this second edition possible.

Last but not least, I would like to thank Silvania for her patience, love, and support throughout the preparation of this edition.
\bigskip\bigskip

\hfill Delft, October 2026

%% file: preface-1st-edition.tex
\chapter*{Preface to the First Edition}

This book is based on notes compiled during the many years I taught the course `Applied Functional Analysis' in the first year of the master's programme at Delft University of Technology, for students with prior exposure to the basics of Real Analysis and the theory of Lebesgue integration.
Starting with the basic results of the subject covered in a typical Functional Analysis course, the text progresses towards a treatment of several
advanced topics, including Fredholm theory, boundary value problems, form methods, semigroup theory, trace formulas, and some mathematical aspects of Quantum Mechanics. With a few exceptions in the later chapters, complete and detailed proofs are given throughout.
This makes the text ideally suited for students wishing to enter the field.

Great care has been taken to present the various topics in a connected and integrated way, and to illustrate abstract results with concrete (and sometimes nontrivial) applications. For example, after introducing Banach spaces and discussing some of their abstract properties, a substantial chapter is devoted to the study of the classical Banach spaces $C(K)$, $L^p(\Om)$, $M(\Om)$, with some emphasis on compactness, density, and approximation techniques. The abstract material in the chapter on duality is complemented by a number of nontrivial applications, such as a characterisation of translation-invariant operators on $L^1(\R^d)$ and Prokhorov's theorem about weak convergence of probability measures. The chapter on bounded operators contains a discussion of the Fourier transform and the Hilbert transform, and includes proofs of the Riesz--Thorin and Marcinkiewicz interpolation theorems. After the introduction of the Laplace operator as a closable operator in $L^p$\!, its closure $\Delta$ is revisited in later chapters from different points of view: as the operator arising from a suitable sesquilinear form, as the operator $-\nabla^\star \nabla$ with its natural domain, and as the generator of the heat semigroup. In parallel, the theory of its Gaussian analogue, the Ornstein--Uhlenbeck operator, is developed and the connection with orthogonal polynomials and the quantum harmonic oscillator is established. The chapter on semigroup theory, besides developing the general theory, includes a detailed treatment of some important examples such as the heat semigroup, the Poisson semigroup, the Schr\"odinger group, and the wave group. By presenting the material in this integrated manner, it is hoped that the reader will appreciate Functional Analysis as a subject that, besides having its own depth and beauty, is deeply connected with other areas of Mathematics and Mathematical Physics.

In order to contain this already lengthy text within reasonable bounds, some choices had to be made. Relatively abstract subjects such as topological vector spaces, Banach algebras, and $C^\star$-algebras are not covered. Weak topologies are introduced {\em ad hoc}, the use of distributions in the treatment of weak derivatives is avoided, and the theory of Sobolev spaces is developed only to the extent needed for the treatment of boundary value problems, form methods, and semigroups. The chapter on states and observables in Quantum Mechanics is phrased in the language of Hilbert space operators.

A work like this makes no claim to originality and most of the results presented here
belong to the core of the subject. Not just the statements, but often
their proofs too, are part of the established canon. Most are taken from, or represent minor variations of, proofs in the many excellent Functional Analysis textbooks in print.

Special thanks go to my students, to whom I dedicate this work. Teaching them has always been a great source of inspiration. Arjan Cornelissen,
Bart van Gisbergen, Sigur Gouwens, Tom van Groeningen, Sean Harris, Sasha
Ivlev, Rik Ledoux, Yuchen Liao, Eva Maquelin, Garazi Muguruza, Christopher Reichling, Floris Roodenburg, Max Sauerbrey,
Cynthia Slotboom, Joop Vermeulen, Matthijs Vernooij, Anouk Wisse, and Timo Wortelboer
pointed out many misprints and more serious errors in earlier versions of this manuscript. The
responsibility for any remaining ones is of course with me. A list of errata will be maintained on my personal webpage.
 I thank Emiel Lorist, Lukas Mias\-kiw\-skyi, and Ivan Yaroslavtsev
for suggesting some interesting problems, Jock Annelle and Jay Kangel for typographical comments,
and Francesca Arici, Martijn Caspers,
Tom ter Elst, Markus Haase, Bas Janssens, Kristin Kirchner, Klaas Landsman, Ben de Pagter, Pierre Portal, Fedor Sukochev, Walter van Suijlekom, and Mark Veraar for helpful discussions and valuable suggestions.

A significant portion of this book was written in the extraordinary circumstances of the global pandemic.
The sudden decrease in overhead and the opportunity to work from home created the time and serenity needed for this project. Paraphrasing the epilogue of W. F. Hermans's novel {\em Onder Professoren} (Among Professors), the book was written entirely in the hours otherwise spent on departmental meetings, committee meetings, evaluations, accreditations, visitations, midterms, reviews, previews, etcetera, and so forth. All that precious time has been spent in a very useful way by the author.

\bigskip\bigskip

\hfill Delft, April 2022

\newpage
\noindent
In the present corrected version we have fixed numerous small misprints, a few misformulations and editing errors, as well as a small number of mathematical oversights. In some proofs, additional details have been written out, and some arguments have been streamlined. I thank Jan Maas for some valuable suggestions in this direction.

\bigskip\bigskip
\hfill Delft, May 2023

\bigskip\bigskip
\noindent
In the present version we have fixed further misprints, most of which were kindly pointed out by Quinten Donker, Niels Goedegebure, Norman Goldstein, Robert Spek\-snijder, and Chris van Vliet. Some paragraphs underwent minor polishing, a few proofs have been simplified, and several new problems have been added.

\bigskip\bigskip
\hfill Delft, October 2024

%% file: how-to-use-this-book.tex
\chapter*{How to Use This Book}

\noindent
The chapters progress from the foundations to more advanced
theory and applications. The core material comprises
Sections~\ref{sec:BanachSpaces}--\ref{sec:compactness},
\ref{sec:sequencespaces}--\ref{sec:measures},
\ref{sec:HilbertSpaces}--\ref{subsec:OS},
\ref{sec:duals-classical}--\ref{sec:HB-separation},
\ref{sec:UBT}--\ref{sec:FT},
\ref{sec:spectrum}--\ref{sec:DC},
\ref{sec:compact}--\ref{sec:spectrum-compact},
\ref{sec:sa-un-no}--\ref{sec:contFC},
and~\ref{sec:spectral-thm-normal-compact}--%
\ref{sec:spectral-thm-normal-bdd}.
The remaining sections of
Chapters~\ref{ch:Banach}--\ref{chap:spectral-theorem},
together with
Chapters~\ref{ch:unbdd}--\ref{ch:SymmQuant},
provide a broader panorama of advanced topics and connect
the core material with the modern functional analytic
approach to partial differential equations, spectral theory, and the foundations of quantum
mechanics.

Several routes through the book are therefore possible
without reading every chapter in full. The following diagram summarises the principal ones.

\begin{center}
\begin{tikzpicture}[
  box/.style={
    draw,
    rounded corners=1pt,
    align=center,
    text width=4.7cm,
    inner sep=5pt,
    font=\small
  },
  >=stealth
]
\node[box] (foundation) at (0,0)
  {Basic functional analysis\\
   Chapters~\ref{ch:Banach}--\ref{ch:operators},
   selected sections};

\node[box] (spectral) at (0,-1.7)
  {Spectral theory and compact operators\\
   Chapters~\ref{ch:spectral}--\ref{ch:unbdd},
   selected sections};

\node[box] (pde) at (-3,-3.7)
  {Sobolev spaces and forms\\
   Chapters~\ref{chap:boundaryvalueproblems}--%
   \ref{chap:forms}};

\node[box] (trace) at (3,-3.7)
  {Trace ideals\\
   Sections~\ref{sec:Hilbert-Schmidt}--%
   \ref{subsec:partialtrace}};

\node[box] (evolution) at (-3,-5.8)
  {Semigroup theory\\
   Chapter~\ref{chap:semigroups}, selected sections};

\node[box] (quantum) at (3,-5.8)
  {States and observables\\
   Sections~\ref{sec:states-observables-classical}--\ref{sec:POVM}};

\node[box] (symmetry) at (3,-7.8)
  {Symmetry and covariance\\
   Sections~\ref{sec:symmetries}--%
   \ref{subsec:StonevonNeumann}};

\draw[->] (foundation) -- (spectral);
\draw[->] (spectral) -- (pde);
\draw[->] (spectral) -- (trace);
\draw[->,dashed] (pde) -- (evolution);
\draw[->] (trace) -- (quantum);
\draw[->] (quantum) -- (symmetry);

\draw[->]
  (spectral.south)
  -- (0,-2.5)
  -- (0,-5.8)
  -- (evolution.east);

\draw[->,dashed]
  (evolution.south)
  -- node[above,sloped,font=\small]
     {\phantom{xxxxx}\ref{sec:Stone}, \ref{subsec:OU}}
  (symmetry.west);
\end{tikzpicture}
\end{center}

\Needspace{6\baselineskip}
An arrow indicates that material in the upper box supplies
tools for the lower box; it does not mean that every section
in the upper box is required.
The dashed arrows indicate the use of Sobolev spaces and
forms in the PDE applications of semigroup theory, and the roles of unitary groups and the Ornstein--Uhlenbeck semigroup in quantum mechanics.
The direct arrow from spectral and operator theory to
semigroups records that the abstract semigroup theory
can be studied before those applications.

The division into core material and remaining sections is not meant to imply a strict division into elementary and advanced material: some advanced sections prepare tools used later.
This is exemplified by the following links, which are useful when making a selection.

\begin{center}
\small
\begin{tabular}{@{}p{.30\textwidth}p{.63\textwidth}@{}}
\hline
Earlier material & Later use\\
\hline
Subsection~\ref{sec:appl-DV}:
fixed point arguments
&
Semilinear Cauchy problems in
Subsection~\ref{subsec:semilinACP}.
\\[4pt]
Subsections~\ref{subsec:mollification}--%
\ref{subsec:Frechet-Kolmogorov}:
approximation and compactness
&
The development of Sobolev spaces in Section~\ref{sec:Sobolev} and compact embeddings
in Section~\ref{sec:Poisson-problem}.
\\[4pt]
Subsection~\ref{subsec:Hermite}:
Hermite polynomials
&
The Ornstein--Uhlenbeck semigroup in
Subsection~\ref{subsec:OU} and second quantisation
in Section~\ref{sec:SQ}.
\\[4pt]
Section~\ref{sec:Sz-Nagy}:
unitary dilation
&
Von Neumann's inequality in
Section~\ref{sec:spectral-thm-normal-bdd}
and the POVM construction in Section~\ref{sec:POVM}.
\\[4pt]
Section~\ref{sec:bicommutant}:
bicommutants
&
The representation of commutative POVMs by Markov
kernels in Subsection~\ref{subsec:Markov-kernels-POVMs}.
\\[4pt]
Section~\ref{sec:Weyl}:
Weyl's theorem
&
The Dirichlet heat trace formula in
Subsection~\ref{subsec:Dirichlet-heat-trace}.
\\[4pt]
Section~\ref{sec:Stone}:
Stone's theorem
&
The treatment of time evolution in
Subsection~\ref{subsec:unitary-time-evolution}.
\\
\hline
\end{tabular}
\end{center}

The applications in Section~\ref{sec:duality-applications}
and topics such as Fredholm theory, the landscape function,
the trace formulas, and classical and quantum correlations
can serve as endpoints for suitable courses.
They provide opportunities to pursue an application in
depth without requiring the reader to continue through
all the subsequent chapters.

\subsection*{Course Use}

The detailed and self-contained presentation, with careful explanations and complete proofs, makes the book suitable for a variety of graduate-level courses.
The problems at the end of every chapter range from exercises consolidating definitions
and theorems to substantial developments of complementary
material. Some ask for examples and counterexamples,
others give alternative proofs, others establish relations
to material developed elsewhere in the book, and longer problems
introduce additional material and can be used for further study
or as individual projects.

%% file: notation.tex
\chapter*{Notation and Conventions}

We write $\N=\{0,1,2,\ldots\}$\index{$Nn$@$\N$} for the set of natural numbers, and $\Z$, $\Q$, $\R$, and $\C$ for the sets of
integer, rational, real, and complex numbers. We set $\R_+ = (0,\infty)$. In statements valid over the real and complex scalars,
we use the symbol $\K$\index{$Kk$@$\mathbb{K}$} to denote either $\R$ or $\C$. Given a complex number $z = a+bi$ with $a,b\in\R$, we denote by $\ov z = a-bi$\index{$Z$@$\ov z$} its complex conjugate\index{complex conjugate} and by $\Re z = a$\index{$R$@$\Re z$}\index{real part!of a complex number} and $\Im z = b$\index{$Im$@$\Im z$}\index{imaginary part!of a complex number} its real and imaginary parts.
We write $|x|$\index{$X$@$\abs{x}$}
for the absolute value of a real number $x\in\R$, the modulus\index{modulus!of a complex number} of a complex number $x\in \C$, and the Euclidean norm of an element $x = (x_1,\dots,x_d)\in \K^d$\!. When dealing with functions $f$ defined on some domain $D$, we write $f\equiv c$ on $S\subseteq D$ if $f(x) = c$ for all $x\in S$.
We use the symbols $\mathbb{D}$\index{$Dd$@$\mathbb{D}$} and $\mathbb{T}$\index{$Tt$@$\mathbb{T}$} for the open unit disc and the unit circle in the complex plane, respectively.
The indicator function\index{indicator function} of a set $A$ is denoted by $\one_A$.\index{$1A$@$\one_A$}

We write $S'\subseteq S$ to express that $S'$ is a subset of $S$. The complement of a set $S$ is denoted by $\complement S$\index{$C$@$\complement{}$}\index{complement} when the larger ambient set, of which $S$ is a subset, is understood. We write $S\setminus S':=\{s\in S:\, s\not\in S'\}$ and define the {\em symmetric difference}\index{symmetric difference}\index{$D$@$A\Delta B$} of $S_1$ and $S_2$ as $S_1\Delta S_2:= (S_1\setminus S_2)\cup (S_2\setminus S_1)$. The interior and closure of a subset $S$ of a topological space are denoted by $S^\circ$\index{$Sa$@$S^\circ$}\index{interior} and $\ov{S}$\index{$Sa$@$\ov S$}\index{closure}, respectively. In the context of metric and normed spaces, $B(x;r)$ denotes the open ball with radius $r$ centred at $x$.

With regard to summation, we write $\sum_{n\ge 1}$ to express a sum that extends over a finite or countably infinite set labelled by the index $n$. When we wish to emphasise either alternative we write $\sum_{n=1}^N$ and $\sum_{n=1}^\infty$, respectively.

Unless explicitly stated otherwise, the symbols $X$ and $Y$ denote Banach spaces and $H$ and $K$  Hilbert spaces. In order to avoid frequent repetitions in the statements of results, these spaces are always thought of as being given and fixed.

The null space\index{null space} and range\index{range} of a linear operator $A$ are denoted by $\Ker(A)$\index{$N$@$\ker(A)$} and $\Ran(A)$,\index{$R$@$\ran(A)$} respectively. When $A$ is unbounded, its domain is denoted by $\Dom(A)$.\index{$DZ$@$\Dom(A)$} A comprehensive list of symbols is contained in the index.

Conventions regarding the choice of scalar field, if applicable, are usually stated at the beginning of a chapter or, in some cases, at the beginning of a section. As a rule of thumb, in Chapters \ref{ch:Banach}--\ref{ch:operators} the scalar field $\K$ can be either $\R$ or $\C$, with a small number of exceptions where this is explicitly stated, such as in our treatment of the Hahn--Banach theorem, the Fourier transform, and the Hilbert transform. From Chapter \ref{ch:spectral} onwards, spectral theory and Fourier transforms are used extensively and the default choice of scalar field is $\C$. In many cases, however, statements not explicitly involving complex numbers or constructions involving them admit counterparts over the real scalars which can be obtained by simple complexification arguments. We leave it to the reader to check this in particular instances.

%% file: ch01-BanachSpaces.tex
\chapter{Banach Spaces}\label{ch:Banach}

\blfootnote{\rm \small \begin{quote}This book has been published by Cambridge University Press in the series `Cambridge Studies in Advanced Mathematics'. The present version is free to view and download for personal use only. Not for re-distribution, re-sale or use in derivative works.\newline {\copyright} Jan van Neerven\end{quote}}

\noindent
The foundations of modern analysis were laid in the early decades of the twentieth century, through the work of Maurice Fr\'echet, Ivar Fredholm,  David Hilbert, Henri Lebesgue, Frigyes Riesz,
and many others. These authors realised that it is fruitful to study
linear operations in a setting of abstract spaces endowed with further structure to accommodate the notions of convergence
and continuity. This led to the introduction of abstract topological and metric spaces and, when combined with linearity,
of topological vector spaces, Hilbert spaces, and Banach spaces.
Since then, these spaces have played a prominent role in all branches of analysis.

The main impetus came from the study of ordinary and partial differential equations where linearity
is an essential ingredient, as evidenced by the linearity of the main operations involved: point
evaluations, integrals, and derivatives. It was discovered that
 many theorems known at the time, such as existence and uniqueness results for ordinary differential
equations and the Fredholm alternative for integral equations, can be conveniently abstracted into general theorems about linear operators in infinite-dimensional spaces of functions.

A second source of inspiration was the discovery, in the 1920s by John von Neumann, that the -- then brand-new -- theory of quantum mechanics can be put on a solid mathematical foundation by means of the spectral theory of selfadjoint operators on Hilbert spaces.
It was not until the 1930s that these two lines of
mathematical thinking were brought together in the theory of Banach spaces,
named after its creator Stefan Banach (although this class of spaces was also discovered, independently and at about the same time, by Norbert Wiener). This theory provides a unified perspective on Hilbert spaces and the various spaces of functions encountered in analysis, including the spaces $C(K)$ of continuous functions and the Lebesgue spaces $L^p(\Om)$.

\section{Banach Spaces}\label{sec:BanachSpaces}

The aim of the present chapter is to introduce the class of Banach spaces and discuss some elementary properties of these spaces. The main classical examples are only briefly mentioned here; a more detailed treatment is deferred to the next two chapters.
Much of the general theory applies to both the real and complex scalar fields. Whenever this applies, the symbol $\mathbb{K}$ is used to denote the scalar field, which is $\R$ in the case of real vector spaces and $\C$ in the case of complex vector spaces.

\subsection{\hbox{Definition and General Properties}}\label{subsec:def-BS}

\begin{definition}[Norms]\label{def:norm} A {\em normed space}\index{normed space} is a pair $(X,\n \cdot \n)$, where $X$ is a vector space
over $\K$ and $\n \cdot\n: X\to [0,\infty)$ is a {\em norm},\index{norm} that is, a mapping with the following properties:
 \begin{enumerate}[leftmargin=*, label=(\roman*)]
  \item\label{it:norm3} $\n x\n = 0$ implies $ x= 0$;
  \item\label{it:norm1} $\n cx\n = |c| \, \n x\n$ for all $c\in \K$ and $x\in X$;
  \item\label{it:norm2} $\n x+x'\n \le \n x\n+\n x'\n$ for all $x,x'\in X$.
 \end{enumerate}
\end{definition}

When the norm $\n \cdot\n$ is understood we simply write $X$ instead of $(X,\n \cdot \n)$. If we wish to emphasise the role of $X$ we write $\n \cdot\n_X$ instead of $\n \cdot\n$.

The properties \ref{it:norm1} and \ref{it:norm2} are referred to as {\em scalar homogeneity} and the {\em triangle inequality}.

The defining properties of a norm imply that every normed space is a metric space, with distance function $$ d(x,y):= \n x-y\n.$$
This observation allows us to introduce
notions such as openness, closedness, compactness, denseness,
limits, convergence, completeness, and continuity in the context of normed spaces by carrying them over from the theory of metric spaces. For instance, a sequence $(x_n)_{n\ge 1}$ in $X$ is said to {\em converge}\index{convergent} if there exists an element $x\in X$ such that $\limn \n x_n-x\n=0$. This element, if it exists, is unique and is called the {\em limit}\index{limit} of the sequence $(x_n)_{n\ge 1}$. We then write $\limn x_n = x$ or simply `$x_n\to x$ as $n\to\infty$'.

The triangle inequality \ref{it:norm2} implies both
$ \n x\n - \n x'\n \le \n x-x'\n$
and
$\n x'\n - \n x\n \le \n x'-x\n.$
Since $\n x'-x\n = \n (-1)\cdot(x-x')\n = \n x-x'\n$ by scalar homogeneity,
we obtain the {\em reverse triangle inequality}\index{triangle inequality!reverse}
$$|\,\n x\n - \n x'\n\,| \le \n x-x'\n.$$
It shows that taking norms $x\mapsto \n x\n$ is a continuous operation.

If $\limn x_n = x$ and $\limn x_n' = x'$ in $X$ and $c_n\to c$ in $\K$, then
$\n c_n x_n - cx\n = \n c_n (x_n-x) + (c_n - c)x\n \le |c_n|\n x_n-x\n + |c_n-c|\n x\n$ implies
$$\limn\n c_n x_n - cx\n =0.$$
Likewise,
$\n (x_n+x_n') - (x+x')\n = \n (x_n-x) + (x_n'-x')\n \le \n x_n-x\n + \n x_n'-x'\n$
implies
$$\limn \n(x_n+x_n') - (x+x')\n = 0.$$
This proves sequential continuity, hence continuity, of the vector space operations.

Throughout this work we use the notation
$$B(x_0;r) := \{x\in X:\, \n x-x_0\n < r\}$$ for the {\em open ball} centred at $x_0\in X$ of radius $r>0$, and
$$\ov B(x_0;r) := \{x\in X:\, \n x-x_0\n \le r\}$$ for the corresponding {\em closed ball}. The {\em open unit ball}\index{unit ball} and {\em closed unit ball} are the balls
\index{$Baa$@$B_X$}\index{$Bab$@$\ov B_X$}
$$ B_X:= B(0;1)=\{x\in X: \, \n x\n < 1\}, \qquad \ov B_X := \ov B(0;1) =\{x\in X: \, \n x\n \le 1\}.$$

\begin{definition}[Banach spaces]
 A {\em Banach space}\index{Banach!space} is a complete normed space.
\end{definition}

Thus a Banach space is a normed space $X$ in which every Cauchy sequence is convergent, that is,
$\lim_{m,n\to\infty} \n x_n-x_m\n = 0$ implies the existence of an $x\in X$ such that
$\limn \n x_n-x\n=0$.

The following proposition gives a necessary and sufficient condition for a normed space to be a Banach space. We need the following terminology. Given a sequence $(x_n)_{n\ge 1}$ in a normed space $X$, the sum $\sum_{n\ge 1}x_n$ is said to be {\em convergent}
if there exists $x\in X$ such that $$ \lim_{N\to\infty} \Big\n x - \sum_{n=1}^N x_n\Big\n =0.$$
The sum $\sum_{n\ge 1}x_n$ is said to be
{\em absolutely convergent}\index{absolutely!convergent} if $\sum_{n\ge 1} \n x_n\n < \infty$.

\begin{proposition}\label{prop:competeness-series}
 A normed space $X$ is a Banach space if and only if every absolutely convergent sum in $X$ converges in $X$.
\end{proposition}
\begin{proof}
 `Only if': \ Suppose that $X$ is complete, and let $\sum_{n\ge 1}x_n$ be absolutely convergent. Then the sequence of
partial sums $(\sum_{j=1}^n x_j)_{n\ge 1}$ is a Cauchy sequence,
 for if $n > m$ the triangle inequality implies
 $$ \Big\n \sum_{j=1}^n x_j - \sum_{j=1}^m x_j\Big\n = \Big\n \sum_{j=m+1}^n x_j \Big\n
 \le \sum_{j=m+1}^n \n x_j\n,$$
 which tends to $0$ as $m,n\to\infty$. Hence, by completeness, the sum  $\sum_{n\ge 1}x_n$ converges.

 \smallskip
 `If': \ Suppose that every absolutely convergent sum in $X$ converges in $X$,
and let $(x_n)_{n\ge 1}$ be a Cauchy sequence in $X$. We must prove that $(x_n)_{n\ge 1}$ converges in $X$.

Choose indices $n_1< n_2< \dots$ in such a way that
 $\n x_i-x_j\n < \frac1{2^k}$ for all $i,j\ge n_k$, $k=1,2,\dots$. The sum $x_{n_1} + \sum_{k\ge 1} (x_{n_{k+1}}-x_{n_k})$
is absolutely convergent since $$ \sum_{k\ge 1} \n x_{n_{k+1}}-x_{n_k}\n \le  \sum_{k\ge 1} \frac1{2^k} <\infty.$$
By assumption it converges to some $x\in X$. Then, by cancellation,
 $$ x = \lim_{m\to\infty} \Big(x_{n_1} + \sum_{k= 1}^m (x_{n_{k+1}}-x_{n_k})\Big) = \lim_{m\to\infty} x_{n_{m+1}},$$
 and therefore the subsequence $(x_{n_m})_{m\ge 1}$ is convergent, with limit $x$. To see that $(x_{n})_{n\ge 1}$ converges to $x$,
we note that
 $$ \n x_m - x\n \le \n x_m - x_{n_m}\n + \n x_{n_m} - x\n \to 0$$
 as $m\to \infty$ (the first term since we started from a Cauchy sequence and the second term by what we just proved).
\end{proof}

The next theorem asserts that every normed space can be completed to a Banach space.
For the rigorous formulation of this result we need the following
terminology.

\begin{definition}[Isometries] A linear mapping $T$ from a normed space $X$ into a normed space $Y$ is said to be an {\em isometry}\index{isometry} if it preserves norms. A normed space $X$ is {\em isometrically contained} in a normed space $Y$
if there exists an isometry from  $X$ into $Y$.
\end{definition}

\begin{theorem}[Completion]\label{thm:completion}\index{completion}
Let $X$ be a normed space. Then:
\begin{enumerate}[label={\rm(\arabic*)}, leftmargin=*]
 \item\label{it:completion1}
there exists a Banach space $\ov X$ containing $X$ isometrically as a dense subspace;
 \item\label{it:completion2}
the space $\ov X$ is unique up to isometry in the following sense: if $X$ is isometrically contained as a dense subspace
 in the Banach spaces $\ov X$ and $\ov{\ov X}$, then the identity mapping on $X$ has a unique extension to an
 isometry from $\ov X$ onto $\ov{\ov X}$.
\end{enumerate}
\end{theorem}
\begin{proof}
As a metric space, $X = (X,d)$ has a completion
$\ov X = (\ov X,\ov d)$ by Theorem \ref{thm:completion-metric}, defined as the set of all equivalence classes of Cauchy sequences in $X$, obtained by declaring Cauchy sequences $(x_n)_{n\ge 1}$ and $(x_n')_{n\ge 1}$ to be equivalent if $\limn d(x_n,x_n')=0$, and endowed with the metric
$\ov d( [(x_n)_{n\ge 1}], [(x_n')_{n\ge 1}]) = \limn d(x_n,x_n')$;
the brackets denote equivalence classes.

The space $\ov X$ is a vector space under the scalar multiplication and addition
$$ c[(x_n)_{n\ge 1}] := [c(x_n)_{n\ge 1}], \qquad [(x_n)_{n\ge 1}]+[(x_n')_{n\ge 1}] := [(x_n+x_n')_{n\ge 1}].$$ The proof that these operations are well defined is routine, and so is the proof that the inclusion mapping $X \subseteq \ov X$ is linear.

We now prove that
$\ov X$ is a Banach space in a natural way, with a norm $\n \cdot\n_{\ov X}$ such that $\ov d(x,x') = \n x-x'\n_{\ov X}$. Properties \ref{it:completion1} and \ref{it:completion2}  follow from the corresponding assertions for metric spaces; in \ref{it:completion2}, the distance-preserving extension is linear because it agrees with the identity on the dense subspace $X$ and the vector space operations are continuous.

If $(x_n)_{n\ge 1}$ is a Cauchy sequence in $X$, the reverse triangle inequality implies that the nonnegative sequence $(\n x_n\n)_{n\ge 1}$ is Cauchy, and hence convergent by the completeness of the real numbers. We now define a norm on $\ov X$ by
\begin{align*} \n [(x_n)_{n\ge 1}]\n_{\ov X} := \limn \n x_n\n.
\end{align*}
It is easy to check that this norm is well defined.
Denoting by $\ov d$ the metric on $\ov X$ given by $ \ov d(x,x'):= \limn d(x_n,x_n')$, where $x=[(x_n)_{n\ge 1}]$ and $x'=[(x_n')_{n\ge 1}]$, it is clear that $\ov d(x,x') = \n x-x'\n_{\ov X}$. This gives the completeness of $\ov X$ as a normed space.
\end{proof}

Most Banach spaces of interest in analysis are {\em infinite-dimensional}\index{infinite-dimensional} in the sense that they do not have a finite spanning set. In this context the following definition is useful.

\begin{definition}[Separability]
A metric space is called {\em separable}\index{separable!normed space} if it contains a countable dense set.
\end{definition}

\begin{proposition}\label{prop:countable} A normed space $X$ is separable if and only if $X$ contains a countable set whose linear span is dense.
\end{proposition}
\begin{proof}
 The `only if' part is trivial. To prove the `if' part, let $(x_n)_{n\ge 1}$ have dense span in $X$.
 Let $Q$ be a countable dense set in $\K$ (for example, one could take $Q = \Q$ if $\K=\R$ and $Q = \Q+i\Q$ if $\K=\C$). Then
 the set of all $Q$-linear combinations of the $x_n$, that is, all linear combinations  involving coefficients from $Q$, is dense in $X$.
\end{proof}

\subsection{Subspaces, Quotients, and Direct Sums}\label{subsec:sub-quot-dirsum}

Several abstract constructions enable us to create new Banach spaces from given ones. We take a brief look at the three most basic constructions, namely, passing to closed subspaces and quotients and taking direct sums.

\paragraph{Subspaces}\index{closed!subspace}
A subspace $Y$ of a normed space $X$ is a normed space with respect to the norm inherited from $X$.
A subspace $Y$ of a Banach space $X$ is a Banach space with respect to the norm inherited from $X$ if and only if
$Y$ is closed in $X$.

To prove the `if' part,
suppose that $(y_n)_{n\ge 1}$ is a Cauchy sequence
in the closed subspace $Y$ of a Banach space $X$. Then it has a limit in $X$, by the completeness of $X$,
and this limit belongs to $Y$, by the closedness of $Y$. The proof of the `only if' part is equally simple and does not require $X$ to be complete. If $(y_n)_{n\ge 1}$ is a sequence in the complete subspace $Y$ such that $y_n\to x$ in $X$, then $(y_n)_{n\ge 1}$ is a Cauchy sequence in $X$, hence also in $Y$, and therefore it has a limit $y$ in $Y$, by the completeness of $Y$. Since $(y_n)_{n\ge 1}$ also converges to $y$ in $X$, it follows that $y= x$ and therefore $x\in  Y$.

\paragraph{Quotients}\index{quotient!Banach spaces} If $Y$ is a closed subspace of a Banach space $X$, the quotient space $X/Y$ can be endowed with a norm by
$$ \n [x]\n := \inf_{y\in Y} \n x-y\n,$$
where for brevity we write $[x]:= x+Y$ for the equivalence class of $x$ modulo $Y$.
Let us check that this indeed defines a norm. If $\n [x]\n = 0$, then there is a sequence $(y_n)_{n\ge 1}$ in $Y$
such that
$\limn \n x-y_n\n = 0$, so $x\in Y$ by the closedness of $Y$. This implies that $[x] = [0]$, the zero element of $X/Y$.
The identity $\n c[x]\n = |c|\n [x]\n$ is trivially verified, and so is the triangle inequality. No completeness was needed up to this point.

To see that the normed space $X/Y$ is complete we use the completeness of $X$ and Proposition \ref{prop:competeness-series}.
If $\sumn \n[x_n]\n <\infty$ and the $y_n\in Y$ are such that $\n x_n - y_n\n \le \n[x_n]\n+\frac{1}{n^2}$,
the proposition implies that $ \sumn (x_n - y_n)$ converges in $X$, say to $x$.
Then, for all $N\ge 1$,
\begin{align*} \Bigl\n [x] - \sum_{n=1}^N [x_n]\Bigr \n
=  \Bigl\n \Bigl[x - \sum_{n=1}^N x_n\Bigr] \Bigr\n
\le  \Bigl\n  x -  \sum_{n=1}^N x_n + \sum_{n=1}^N  y_{n}\Bigr\n
  = \Bigl\n  x  - \Bigl(\sum_{n=1}^N (x_n - y_n)\Bigr)\Bigr\n.
\end{align*}
As $N\to \infty$, the right-hand side tends to $0$ and therefore $\lim_{N\to\infty} \sum_{n=1}^N [x_n] = [x]$ in $X/Y$.

\paragraph{Direct Sums}
A {\em product norm}\index{product!norm} on a finite
product $X = X_1\times\cdots\times X_N$ of normed spaces is a
norm $\n\cdot\n$ satisfying
$$\bigl\n (0,\dots,0,\underbrace{x_n}_{n{\rm th}},0,\dots,0)\bigr\n =  \n x_n\n  \le \n (x_1,\dots,x_N)\n$$
for all $x=(x_1,\dots,x_N)\in X$ and $n=1,\dots,N$.
For instance, every norm $|\cdot|$ on $\K^N$ which is monotone in the sense that $|x_n| \le |x_n'|$ for $n=1,\dots,N$ implies $|x|\le |x'|$ and assigns norm one to the standard unit vectors induces a product norm on $X$ by the formula
\begin{align*}
\n (x_1,\dots,x_N)\n := \bigl|  (\n x_1\n,\dots, \n x_N\n)\bigr|.
\end{align*}
As a normed space endowed with a product norm, the product will be denoted\index{$X_0\oplus X_1$}
$$ X = X_1 \oplus \cdots\oplus X_N$$
and called a {\em direct sum} of $X_1,\dots,X_N$. If every $X_n$ is a Banach space, then the normed space $X$ is a Banach space. Indeed, from
\begin{align}\label{eq:prod-norm}
\n x\n = \Bigl\n \sum_{n=1}^N (0,\dots,0,x_n,0,\dots,0)\Bigr\n \le \sum_{n=1}^N \n x_n \n \le N \n x\n
\end{align}
we see that a sequence $(x^{(k)})_{k\ge 1}$ in $X$ is Cauchy if and only if all its coordinate sequences $(x_n^{(k)})_{k\ge 1}$ are Cauchy. If the spaces $X_n$ are complete, these coordinate sequences have limits $x_n$ in $X_n$, and these limits serve as the coordinates of an element $x=(x_1,\dots,x_N)$ in $X$ which is the limit of the sequence $(x^{(k)})_{k\ge 1}$.

\subsection{First Examples}\label{subsec:examples-BS}

The purpose of this brief section is to present a first catalogue of Banach spaces. The presentation is not self-contained; the examples will be revisited in more detail in the next two chapters, where the relevant terminology is introduced and proofs are given.

\begin{example}[Euclidean spaces]\label{ex:eucliden}
On $\K^d$ we may consider the Euclidean norm\index{Euclidean norm}
$$ \n a\n_2  := \Bigl(\sum_{j=1}^d |a_j|^2\Bigr)^{1/2}\!,$$
and more generally the $p$-norms
$$ \n a\n_p  := \Bigl(\sum_{j=1}^d |a_j|^p\Bigr)^{1/p}\!, \quad 1\le p<\infty,$$
as well as the supremum norm
$$\n a\n_\infty  := \sup_{1\le j\le d} |a_j|.$$
\begin{figure}
\begin{center}
 \begin{tikzpicture}[scale=1]
\draw (1,-1.3) -- (1,1.3);
\draw (-0.3,0) -- (2.3,0);

\draw (5,-1.3) -- (5,1.3);
\draw (3.7,0) -- (6.3,0);

\draw (9,-1.3) -- (9,1.3);
\draw (7.7,0) -- (10.3,0);

\draw (2,-0.1) -- (2,0.1) node[above right] {$1$};
\draw (0.9,1) -- (1.1,1) node[above right] {$1$};

\draw (6,-0.1) -- (6,0.1) node[above right] {$1$};
\draw (4.9,1) -- (5.1,1) node[above right] {$1$};

\draw (10,-0.1) -- (10,0.1) node[above right] {$1$};
\draw (8.9,1) -- (9.1,1) node[above right] {$1$};

\draw [dotted]  (0,0) -- (1,1) -- (2,0) -- (1,-1) -- (0,0);
\draw [dotted] (5,0) circle (28pt);
\draw [dotted] (10,-1) -- (10,1) -- (8,1) -- (8,-1) -- (10,-1);
\filldraw [nearly transparent]  (0,0) -- (1,1) -- (2,0) -- (1,-1) -- (0,0);
\filldraw [nearly transparent](5,0) circle (28pt);
\filldraw [nearly transparent]  (10,-1) -- (10,1) -- (8,1) -- (8,-1) -- (10,-1);
 \end{tikzpicture}
 \caption{The open unit balls of $\R^2$ with respect to the norms $\n\cdot\n_1$,  $\n\cdot\n_2$,  $\n\cdot\n_\infty$.}
\end{center}
\end{figure}
It is not immediately obvious that the $p$-norms are indeed norms; the triangle inequality for these norms will be proved in the next chapter.
It is an easy matter to check that the above norms are all {\em equivalent} in the sense defined in
Section \ref{sec:eqnorm}.
{\em In what follows the Euclidean norm of an element $x\in \K^d$ is denoted by $|x|$.}
\end{example}

\begin{example}[Sequence spaces] The preceding example may be generalised to infinite sequences as follows.
For $1\le p<\infty$ the space $\ell^p$ is defined as the space of all scalar sequences $a = (a_k)_{k\ge 1}$ satisfying
$$ \n a\n_p:= \Bigl(\sum_{k\ge 1}|a_k|^p \Bigr)^{1/p}<\infty.$$ The mapping $a\mapsto \n a\n_p$ is a norm which turns $\ell^p$ into a Banach space. The space $\ell^\infty$ of all bounded scalar sequences $a = (a_k)_{k\ge 1}$ is a Banach space with respect to the norm
$$\n a\n_\infty := \sup_{k\ge 1} |a_k|<\infty.$$
The space $c_0$ consisting of all bounded scalar sequences $a = (a_k)_{k\ge 1}$ satisfying $$\limk a_k = 0$$ is a closed subspace of $\ell^\infty$\!. As such it is a Banach space in its own right.
\end{example}

\begin{example}[Spaces of continuous functions]
 Let $K$ be a nonempty compact topological space.
The space $C(K)$ of all continuous functions $f:K\to \K$ is a Banach space with respect to the supremum norm
$$ \n f\n_{\infty} := \sup_{x\in K} |f(x)|.$$
This norm captures the notion of uniform convergence\index{uniform convergence!of functions}:
for functions in $C(K)$ we have
$\limn\n f_n - f\n_\infty = 0$ if and only if $\limn f_n = f$ uniformly.
If $K$ is a compact metric space, then $C(K)$ is separable (Proposition \ref{prop:CK-separable}).
\begin{figure}
\begin{center}

\begin{tikzpicture}
\begin{axis}[
    axis lines = middle,
    xlabel = {$x$},
    ylabel = {$y$},
    xtick={0,0.5,1},
    xmin=0, xmax=1,
    ymin=0, ymax=5],
\addplot [dotted, name path = A,
    domain = 0:1,
    samples = 1000] {-x^2+2-2*x+3*x^3};

\addplot [name path = B,
    domain = 0:1] {-x^2+3-2*x+3*x^3}
    node [left=95pt, below=10pt] {$f(x)$};

\addplot [dotted, name path = C,
    domain = 0:1] {-x^2+4-2*x+3*x^3};
\addplot [nearly transparent] fill between [of = A and C, soft clip={domain=0:1}];
\end{axis}
\end{tikzpicture}
\caption{The open ball $B(f;1)$ in $C[0,1]$ consists of all functions in $C[0,1]$ whose graph lies inside the shaded area.}
\end{center}
\end{figure}
\end{example}

\begin{example}[Spaces of integrable functions] Let $(\Om,\calF\!,\mu)$ be a measure space.
For $1\le p<\infty$, the space $L^p(\Om)$ consisting of
all measurable functions $f:\Om\to \K$ such that
$$ \n f\n_p := \Bigl(\int_\Om |f|^p\ud \mu\Bigr)^{1/p}< \infty,$$
identifying functions that are equal $\mu$-almost everywhere,
is a Banach space with respect to the norm $\n \cdot\n_p$.
The space $L^\infty(\Om)$ consisting of
all measurable and $\mu$-essentially bounded functions $f:\Om\to \K$,
identifying functions that are equal $\mu$-almost everywhere,
is a Banach space with respect to the norm given by the $\mu$-essential supremum
$$\n f\n_\infty := \hbox{$\mu$-}\esssup_{\om\in\Om} |f(\om)|:= \inf\bigl\{r>0: \, |f|\le r \ \hbox{$\mu$-almost everywhere}\bigr\}.$$
If $1\le p<\infty$, the measure $\mu$ is $\sigma$-finite, and the $\sigma$-algebra $\calF$ is countably generated, then $L^p(\Om,\mu)$ is separable (Proposition \ref{prop:Lp-separable}).
\end{example}

\begin{example}[Spaces of measures]
 Let $(\Om,\calF)$ be a measurable space.
 The space $M(\Om)$ consisting of all $\K$-valued measures on $(\Om,\calF)$
 is a Banach space with respect to the variation norm
 $$ \n \mu\n := |\mu|(\Om) :=
 \sup_{\mathscr{A}\in \mathbb{F}} \sum_{A\in \mathscr{A}}|\mu(A)|,$$
where $\mathbb{F}$ denotes the set of all
finite collections of
pairwise disjoint sets in $\calF$.
\end{example}

\begin{example}[Hilbert spaces]
A {\em Hilbert space} is an inner product space $(H,\iprod{\cdot}{\cdot})$ that is complete with respect to the norm
 $$ \n h\n := \iprod{h}{h}^{1/2}\!.$$
Examples include the spaces $\K^d$ with the Euclidean norm, $\ell^2$\!, and the spaces $L^2(\Om)$. Precise definitions and further examples will be given in later chapters.
\end{example}

\section{Bounded Operators}\label{sec:bdd-operators}

Having introduced normed spaces and Banach spaces, we now introduce a class of linear operators acting between them which interact with the norm in a meaningful way.

\subsection{Definition and General Properties}

Let $X$ and $Y$ be normed spaces.

\begin{definition}[Bounded operators]\label{def:bdd-op} A linear operator $T: X\to Y$ is {\em bounded}\index{bounded!operator}
if there exists a finite constant $C\ge 0$ such that
$$\n Tx\n \le C \n x\n, \quad x\in X.$$
Here, and in the rest of this work, we write $Tx$ instead of the more cumbersome $T(x)$.
A {\em bounded operator} is a linear operator that is bounded.
\end{definition}

The infimum $C_T$ of all admissible constants $C$ in Definition \ref{def:bdd-op} is itself admissible. Thus $C_T$ is the least admissible constant. We claim that it equals the number
\begin{align*} \n T\n := \sup_{\n x\n\le 1} \n Tx\n.
\end{align*}
To see this, let $C$ be an admissible constant in Definition \ref{def:bdd-op}, that is, we assume that $\n Tx\n \le C \n x\n$ for all $x\in X$. Then $\n T\n = \sup_{\n x\n\le 1}\n Tx\n \le C$. Since this holds for all admissible constants $C$, it follows that
$\n T\n \le C_T$. The opposite inequality $C_T\le \n T\n$ follows by observing that for all $x\in X$ we have
$$ \n Tx\n \le \n T\n \n x\n,$$
which means that $\n T\n$ is an admissible constant.
This inequality is trivial for $x=0$, and for $x\not=0$ it follows from scalar homogeneity, the linearity of $T$ and the definition of the number $\n T\n$:
$$ \n Tx\n = \Bigl\n \frac{1}{\n x\n}Tx\Bigr\n {\n x\n}
= \Bigl\n T\frac{x}{\n x\n}\Bigr\n {\n x\n} \le \n T\n \n x\n.$$

\begin{proposition}\label{prop:bdd-cont} For a linear operator $T: X\to Y$ the following assertions are equivalent:
 \begin{enumerate}[label={\rm(\arabic*)}, leftmargin=*]
  \item\label{it:bdd-cont1} $T$ is bounded;
  \item\label{it:bdd-cont2} $T$ is continuous;
  \item\label{it:bdd-cont3} $T$ is continuous at some point $x_0\in X$.
 \end{enumerate}
\end{proposition}
\begin{proof}
 The implication \ref{it:bdd-cont1}$\Rightarrow$\ref{it:bdd-cont2} follows from $$\n Tx - Tx'\n = \n T(x-x')\n \le \n T\n \n x-x'\n$$ and
the implication \ref{it:bdd-cont2}$\Rightarrow$\ref{it:bdd-cont3} is trivial. To prove the implication \ref{it:bdd-cont3}$\Rightarrow$\ref{it:bdd-cont1}, suppose that $T$ is continuous at $x_0$. Then there exists a $\delta>0$ such that $\n x_0-y\n<\delta$ implies $\n Tx_0 - Ty\n < 1.$
Since every $x\in X$ with $\n x\n<\delta$ is of the form $x = x_0-y$ with $\n x_0-y\n<\delta$ (take $y = x_0-x$), it follows that $\n x\n<\delta$ implies $\n Tx\n < 1$.
By scalar homogeneity and the linearity of $T$ we may scale both sides by a factor $\delta$, and obtain that $\n x\n<1$ implies $\n Tx\n<1/\delta$.
Applying this to $rx$ with $0<r<1$, letting $r\uparrow1$, and using continuity of scalar multiplication shows that $T$ is bounded and $\n T\n\le1/\delta$.
\end{proof}

Easy manipulations involving the properties of norms and linear operators, such as those used in the above proofs, will henceforth be omitted.

The set of all bounded operators from $X$ to $Y$ is
a vector space in a natural way with respect to pointwise scalar multiplication and addition by putting
$$ (cT)x:= c(Tx), \qquad (T+T')x:= Tx +T'x.$$
This vector space will be denoted by $\calL(X,Y)$.\index{$L$@$\calL(X,Y)$}
We further write $\calL(X):=\calL(X,X)$.\index{$L$@$\calL(X)$}

For all $T,T'\in \calL(X,Y)$ and $c\in\K$ we have $$\n cT\n = |c|\n T\n, \qquad \n T+T'\n \le \n T\n + \n T'\n.$$
Let us prove the second assertion; the proof of the first is similar. For all $x\in X$, the triangle inequality gives
$$ \n (T+T')x \n \le \n Tx\n + \n T'x\n \le (\n T\n + \n T'\n)\n x\n,$$
and the result follows by taking the supremum over all $x\in X$ with $\n x\n\le 1$.

Since $\n T\n = 0$ implies $T=0$, it follows that $T\mapsto \n T\n$ is a norm on $\calL(X,Y)$. Endowed with this norm, $\calL(X,Y)$ is a normed space.
If $T:X\to Y$ and $S:Y\to Z$ are bounded, then so is their composition $ST$ and we have
$$\n ST\n \le \n S\n \n T\n.$$ Indeed, for all $x\in X$ we have
$$ \n STx \n \le \n S\n \n Tx\n \le \n S\n \n T\n \n x\n$$
and the result follows by taking the supremum over all $x\in X$ with $\n x\n\le 1$.

\begin{proposition}
If $Y$ is complete, then $\calL(X,Y)$ is complete.
\end{proposition}
\begin{proof}
 Let $(T_n)_{n\ge 1}$ be a Cauchy sequence in $\calL(X,Y)$. From $\n T_n x - T_m x\n \le \n T_n - T_m\n \n x\n$ we see that
 $(T_n x)_{n\ge 1}$ is a Cauchy sequence in $Y$ for every $x\in X$. Let $Tx$ denote its limit. The linearity of each of the operators $T_n$ implies that the mapping $T: x\mapsto Tx$
is linear and we have
$ \n Tx \n = \limn \n T_n x\n \le M \n x\n,$
 where $M := \sup_{n\ge 1} \n T_n\n$ is finite since Cauchy sequences in normed spaces are bounded.
This shows that the linear operator $T$ is bounded, so it is an element of $\calL(X,Y)$. To prove that $\limn \n T_n - T\n =0$,
fix $\e>0$ and let $N\ge 1$ be so large that $\n T_n - T_m\n< \e$ for all $m,n\ge N$. Then, for $m,n\ge N$, from
$ \n T_n x - T_m x \n \le \e \n x\n$ we obtain, upon letting $m\to\infty$, that
$$ \n T_n x - T x \n \le \e \n x\n.$$ Since this holds for all $x\in X$ and $n\ge N$,
it follows that $\n T_n  -T\n \le \e$ for all $n\ge N$.
\end{proof}

The important special case $Y = \K$ leads to the following definition.

\begin{definition}\label{def:dual} The {\em dual space} of a normed space $X$ is the Banach space\index{$X^*$} $$X\s := \calL(X,\K).$$
\end{definition}

For $x\in X$ and $x^*\in X^*$ one usually writes $ \lb x,x^*\rb := x^*(x).$
The elements of the dual space $X^*$ are often referred to as {\em bounded functionals} or simply {\em functionals}.\index{functional} Duality is a subject in its own right which will be taken up in Chapter \ref{ch:duality}. In that chapter,
explicit representations of duals of several classical Banach spaces are given. For Hilbert spaces, this duality takes a particularly simple form, described by the Riesz representation theorem, to be proved in Chapter \ref{ch:Hilbert-spaces}.

It often happens that a linear operator can be shown to be well defined and bounded on a dense subspace.
In such cases, a {\em density argument}\index{density!argument} can be used to extend the operator to the whole space.

\begin{proposition}[Density argument -- extending operators]\label{prop:extendT}
Let $X$ be a normed space and $Y$ be a Banach space, and let $X_0$ be a dense subspace of $X$.
 If $T_0: X_0\to Y$ is a bounded operator, there exists a unique bounded operator $T:X\to Y$
 extending $T_0$. The norm of this extension satisfies $\n T\n = \n T_0\n$.
 \end{proposition}
\begin{proof}
 Fix $x\in X$, and let $\limn x_n = x$ with $x_n\in X_0$ for all $n\ge 1$.
 The boundedness of $T_0$ implies that $\n T_0 x_n - T_0 x_m\n \le \n T_0\n \n x_n - x_m \n \to 0$
 as $m,n\to \infty$, so $(T_0 x_n)_{n\ge 1}$ is a Cauchy sequence in $Y$. Since $Y$ is complete,
 we have $T_0 x_n \to y$ for some $y \in Y$.

 If $x_n'\to x$ as well, the same argument shows that $T_0 x_n'\to y'$ for some (possibly different) $y'\in Y$.
 From $$ \n T_0 x_n'-T_0 x_n\n \le \n T_0\n \n x_n'-x_n\n \le \n T_0\n (\n x_n'-x\n + \n x - x_n\n)$$
 it follows that
 $  \n y'-y\n =\limn \n T_0 x_n'-T_0 x_n\n =0$
 and therefore $y'=y$.

 Denoting the common limit $y=y'$ by
 $Tx$, we thus obtain a well-defined mapping $x\mapsto Tx$. It is evident that this mapping extends $T_0$,
 for if $x\in X_0$ we may take $x_n = x$ and then $Tx = \limn T_0 x_n = T_0 x$.

 It is easily checked that $T$ is linear. To show that it is bounded, with $\n T\n \le \n T_0\n$,
 we just note that
 $$ \n T x\n = \limn \n T_0 x_n\n \le \n T_0\n\limn \n x_n\n = \n T_0\n \n x\n.$$
 The converse inequality $\n T\n \ge \n T_0\n$ trivially holds since $T$ extends $T_0$.

Finally, if the bounded operators $T$ and $T'$ both extend $T_0$, then the bounded operator $T-T'$ equals $0$ on the dense subspace $X_0$ and hence, by continuity, on all of $X$.
\end{proof}

Under a uniform boundedness assumption, a similar density argument can be used to extend the existence
of limits from a dense subspace to the whole space.

\begin{proposition}[Density argument -- extending convergence of operators]\label{prop:approxTn} Let $X$ be a normed space and $Y$ a Banach space, and let $X_0$ be a dense subspace of $X$. Let
$(T_n)_{n\ge 1}$ be a sequence of operators in $\calL(X,Y)$ satisfying $\sup_{n\ge 1}\n T_n\n<\infty$.
If the limit $\limn T_n x_0$ exists in $Y$ for all $x_0\in X_0$, then
the limit $Tx:= \limn T_n x$ exists in $Y$ for all $x\in X$. Moreover,
the operator $T: x\mapsto Tx$ is linear and bounded from $X$ to $Y$, and
$$ \n T\n \le \liminf_{n\to\infty}\n T_n\n.$$
\end{proposition}

\begin{proof} Let $M:= \sup_{n\ge 1}\n T_n\n$.
We will show that the sequence $(T_nx)_{n\ge 1}$ is Cauchy for every $x\in X$. Fix arbitrary $x\in X$ and $\e>0$ and choose $x_0\in X_0$
in such a way that $\n x-x_0\n< \e/(M+1)$.
Since $(T_nx_0)_{n\ge 1}$ is a Cauchy sequence, there is an $N\ge 1$ such that $\n T_n x_0 -T_m x_0\n<\e$
for all $m,n\ge N$. Then, for all $m,n\ge N$,
\begin{align*} \n T_n x-T_m x\n & \le  \n T_n x-T_n x_0\n + \n T_n x_0-T_m x_0\n + \n T_mx_0 -T_m x\n
\\ & \le M \n x- x_0\n + \e + M\n x_0 - x\n < 3\eps.
\end{align*}
The sequence $(T_nx)_{n\ge 1}$ is thus Cauchy. Since $Y$ is complete, this sequence has a limit, which we denote by $Tx$.
Linearity of $T: x\mapsto Tx$ is clear, and boundedness, along with the estimate for the norm, follows from
\[\n T x\n =\limn \n T_n x\n =\liminf_{n\to\infty}\n T_n x\n \le \liminf_{n\to\infty}\n T_n \n \n x\n.
\qedhere \]
\end{proof}

This proposition should be compared with Proposition \ref{prop:strong-limits}, which provides the following partial
converse: if $X$ is a Banach space, $Y$ is a normed space, and $(T_n)_{n\ge 1}$ is a sequence in $\calL(X,Y)$ such that $Tx:= \limn T_n x$ exists in $Y$
for all $x\in X$, then $\sup_{n\ge 1} \n T_n\n <\infty$.

\begin{definition}[Null space and range]
The {\em null space}\index{null space} of a bounded operator $T\in \calL(X,Y)$ is the subspace
$$ \ker(T) := \{x\in X:\, Tx = 0\}.$$
The {\em range}\index{range} of $T$ is the subspace
$$\ran(T):= \{Tx:\, x\in X\}.$$
\end{definition}

By linearity, both the null space $\ker(T)$ and the range $\Ran(T)$ are subspaces. By continuity, the null space of a bounded operator is closed. The following result gives a useful sufficient criterion for the range of a bounded operator to be closed.

\begin{proposition}\label{prop:closed-range} Let $X$ be a Banach space and $Y$ be a normed space.
 If $T\in\calL(X,Y)$ satisfies $\n Tx \n \ge C \n x\n$ for some $C>0$ and all $x\in X$, then
 $T$ is injective and has closed range.
\end{proposition}
\begin{proof}
Injectivity is clear. Suppose that $Tx_n \to y$ in $Y$; we must prove that $y\in \ran(T)$. From $\n x_n -x_m\n \le C^{-1} \n Tx_n - Tx_m\n$
 it follows that $(x_n)_{n\ge 1}$ is a Cauchy sequence in $X$ and therefore converges to some $x\in X$.
 Then $y = \limn T x_n = Tx$.
\end{proof}

We conclude by introducing some terminology that will be used throughout this work. In the next four definitions, $X$ and $Y$ are normed spaces.

\begin{definition}[Isomorphisms]
An {\em isomorphism}\index{isomorphism} is a bijective operator $T\in\calL(X,Y)$ whose inverse is bounded as well.
An {\em isometric isomorphism}\index{isometric isomorphism} is an isomorphism that is also isometric.
The spaces $X$ and $Y$ are called {\em (isometrically) isomorphic} if there exists an (isometric) isomorphism from $X$ to $Y$.
\end{definition}

\begin{definition}[Contractions]\label{def:contraction}
A {\em contraction}\index{contraction} is an operator $T\in\calL(X,Y)$ satisfying $\n T\n \le 1$.
\end{definition}

\begin{definition}[Uniform boundedness] A subset $\mathscr{T}$ of $\calL(X,Y)$ is said to be
{\em uniformly bounded}\index{bounded!uniformly} if it is a bounded subset of $\calL(X,Y)$, that is, if
$ \sup_{T\in\mathscr{T}}\n T\n <\infty.$
\end{definition}

\begin{definition}[Uniform, strong, and weak convergence of operators]
A sequence $(T_n)_{n\ge 1}$ in $\calL(X,Y)$ is said to
\begin{enumerate}[label={\rm(\arabic*)}, leftmargin=*]
 \item {\em converge uniformly}\index{uniform convergence!of operators}
 to an operator $T\in \calL(X,Y)$ if $\limn \n T_n-T\n=0$;
 \item {\em converge strongly}\index{strong!convergence, of operators}
 to an operator $T\in \calL(X,Y)$ if $\limn \n T_n x-Tx\n=0$ for all $ x\in X$;
  \item {\em converge weakly}\index{weak!convergence, of operators}
 to an operator $T\in \calL(X,Y)$ if $\limn \lb T_n x-Tx,y^*\rb=0$ for all $x\in X, \ y^*\in Y^*$, where $Y^*$ is the dual of $Y$
and $\lb y,y^*\rb := y^*(y)$ for $y\in Y$.
\end{enumerate}
In these situations we call $T$ the {\em uniform limit}, the {\em strong limit}, and the {\em weak limit}, respectively, of the sequence $(T_n)_{n\ge 1}$.
Uniqueness of weak limits is assured by the Hahn--Banach theorem (see Corollary \ref{cor:HB-norming}).
\end{definition}

Uniform convergence implies strong convergence and strong convergence implies weak convergence, but the converses generally fail. For instance,  the projections onto the first $n$ coordinates in $\ell^p$\!, $1\le p<\infty$, converge strongly to the identity operator, but not uniformly; and the operators $T^n$, where $T$ is the right shift in $\ell^p$\!, $1< p<\infty$, converge weakly to the zero operator but not strongly (for the case $p=1$ see Problem \ref{prob:Schur}).

\subsection{Subspaces, Quotients, and Direct Sums}\label{subsec:sub-quot-dirsum-op}

\paragraph{Restrictions} If $T$ is a bounded operator from a normed space $X$ into a normed space $Y$, then the restriction of $T$ to a subspace $X_0$ of $X$ defines a bounded operator $T|_{X_0}$ from $X_0$ into $Y$ of norm $\n T|_{X_0}\n \le \n T\n$.

\paragraph{Quotients} Let $Y$ be a closed subspace of a Banach space $X$. By the definition of the quotient norm, the {\em quotient map} $q:x\mapsto x+Y$ is bounded from $X$ to $X/Y$ of norm $\n q\n\le 1$.

Let $Z$ be a normed space and let $T\in\calL(X,Z)$ be a bounded operator with the property that $Y$ is contained in the null space $\Ker(T)$.
We claim that  $$T_{/Y}(x+Y):= Tx, \quad x\in X,$$
is well defined and defines a bounded {\em quotient operator} $T_{/Y}:X/Y\to Z$ of norm $\n T_{/Y}\n = \n T\n$.
Well-defined\-ness of $T_{/Y}$ is clear, and for all $x\in X$ and $y\in Y$
we have $\n Tx\n = \n T(x-y) \n \le \n T\n \n x-y\n.$ Taking the infimum over all $y\in Y$
gives the bound $$\n T_{/Y}(x+Y)\n = \n Tx\n \le \n T\n \inf_{y\in Y} \n x-y\n = \n T\n \n x+Y\n.$$
Hence $T_{/Y}$ is bounded and $\n T_{/Y}\n \le \n T\n$. For the converse inequality we note that
$$ \n Tx \n = \n T_{/Y}(x+Y)\n \le \n T_{/Y}\n \n x+Y\n = \n T_{/Y}\n \inf_{y\in Y} \n x-y\n  \le \n T_{/Y}\n \n x\n.$$

\paragraph{Direct Sums} If $X_n$ is a normed space and $T_n\in \calL(X_n)$ for $n=1,\dots,N$, then the {\em direct sum operator}
 $$T = \bigoplus_{n=1}^N T_n : (x_1,\dots,x_N)\mapsto (T_1 x_1,\dots,T_Nx_N)$$
 is bounded on $X = \bigoplus_{n=1}^N X_n$ with respect to any product norm;
 this follows from
 \eqref{eq:prod-norm}.

\subsection{First Examples}

We revisit the examples of Section \ref{subsec:examples-BS}
and discuss how various natural operations used in analysis give rise to bounded operators.

\begin{example}[Matrices]\label{ex:matrices-bdd}
 Every $m\times n$ matrix $A=(a_{ij})_{1\le i\le m,\,1\le j\le n}$ defines a bounded operator in $\calL(\K^n\!,\K^m)$
 and its norm satisfies
\begin{align}\label{eq:norm-A-upper} \n A\n^2 = \sup_{|x|\le 1} |Ax|^2 = \sup_{|x|\le 1} \sum_{i=1}^m \Big|\sum_{j=1}^n a_{ij} x_j\Big|^2
  \le \sum_{i=1}^{m} \sum_{j=1}^{n}|a_{ij}|^2\!,
\end{align}
where the last step follows from the Cauchy--Schwarz inequality.
 More generally, every linear operator from a finite-dimensional normed space $X$ into a normed space $Y$ is bounded; this will be shown in Corollary \ref{cor:equiv-norms-fd}.
\end{example}

The upper bound \eqref{eq:norm-A-upper} for the norm of a matrix $A$ need not be an equality. An explicit method to determine the operator norm of a matrix is described in Problem \ref{prob:hermit}.

\begin{example}[Point evaluations]\label{ex:bdd-op-evaluation}
Let $K$ be a compact topological space. For each $x_0\in K$ the point evaluation\index{point!evaluation} $E_{x_0}: f \mapsto f(x_0)$
is bounded as an operator from $C(K)$ into $\K$ with norm $\n E_{x_0}\n = 1$. Boundedness with norm $\n E_{x_0}\n \le 1$ follows
from
$$ |E_{x_0}f| = | f(x_0)| \le \sup_{x\in K} |f(x)|  = \n f\n_\infty.$$
By taking $f=\one$, the constant-one function on $K$, we see that $\n E_{x_0}\n=1$.
\end{example}

\begin{example}[Integration]\label{ex:bdd-op-integration}
Let $(\Omega,\calF\!,\mu)$ be a measure space. The mapping
$I_\mu: f\mapsto \int_\Omega f\ud \mu$
is bounded from $L^1(\Om)$ to $\K$, with norm $\n I_\mu\n =1$ if there exists a set $F\in\calF$ with $0<\mu(F)<\infty$, and $\n I_\mu\n =0$ otherwise.
Indeed, for all $f\in L^1(\Om)$ we have
$$ |I_\mu f|= \Bigl|\int_\Omega f\ud \mu\Bigr| \le \int_\Omega |f|\ud \mu = \n f\n_1.$$
If $0<\mu(F)<\infty$, taking $f=\one_F$ shows that $\n I_\mu\n=1$. If no such set exists, then every $f\in L^1(\Om)$ vanishes almost everywhere, since each set $\{|f|>1/n\}$ has finite measure; hence $\n I_\mu\n=0$.
\end{example}

\begin{example}[Pointwise multipliers]\label{ex:pointwise-multiplier}\index{multiplier}
 Let $(\Omega,\calF\!,\mu)$ be a measure space and fix $1\le p\le \infty$.
For any $m\in L^\infty(\Om)$, the pointwise multiplier
$T_m: f\mapsto mf$
defines a bounded operator on $L^p(\Om)$  with norm
$\n T_m\n \le \n m\n_\infty$, with equality if $p=\infty$; for $1\le p<\infty$, equality holds if $\mu$ is $\sigma$-finite. Indeed, for $\mu$-almost all
$\om\in\Om$ we have
$$ |(mf)(\om) |=  |m(\om)||f(\om)| \le \n m\n_\infty |f(\om)|.$$
For $1\le p<\infty$, upon integration we obtain
$$ \n T_m f\n_p^p  = \int_\Om |mf|^p\ud \mu \le \n m\n_\infty^p  \int_\Om |f|^p\ud \mu =  \n m\n_\infty^p \n f\n_p^p,$$
so $T_m$ is bounded on $L^p(\Om)$ and $\n T_m \n \le \n m\n_\infty$.
For $p=\infty$ the analogous bound follows by taking essential suprema.

To prove equality $\n T_m\n = \n m\n_\infty$, first note that there is
nothing to prove if $\n m\n_\infty=0$. Otherwise, let $0<\eps<1$.
By the definition of the essential supremum, the set
$$
  F_\eps := \{|m|>(1-\eps)\n m\n_\infty\}
$$
has strictly positive measure. By $\sigma$-finiteness it contains a
measurable subset $G_\eps$ of strictly positive finite measure. Then
$\one_{G_\eps}\in L^p(\Om)$ and
$$
  \n T_m \one_{G_\eps}\n_p
  \ge (1-\eps)\n m\n_\infty \n \one_{G_\eps}\n_p .
$$
Thus $\n T_m\n\ge (1-\eps)\n m\n_\infty$. Since $\eps>0$ was arbitrary,
this implies $\n T_m\n\ge \n m\n_\infty$. For $p=\infty$ this argument works
without assuming $\sigma$-finiteness.
\end{example}

\begin{example}[Integral operators]\label{ex:kernel} Let $\mu$ be a finite Borel measure on a compact metric space $K$.
With respect to the {\em product metric}
$d((s,t),(s',t')) := d(s,s') +d(t,t')$,
$K\times K$ is a compact metric space (see Proposition \ref{prop:Tychonov-metric}).
Let $k\in C(K\times K)$ and define, for $f\in C(K)$, the function $Tf: K\to\K$ by
$$ Tf(s):= \int_K k(s,t)f(t)\ud \mu(t), \quad s\in K.$$
By the uniform continuity of $k$ (see Theorem \ref{thm:compact-UC}), $Tf$ is a continuous function. Indeed, given $\eps>0$, choose $\delta>0$ so small that
$d((s,t),(s',t'))<\delta$ implies $|k(s,t)-k(s',t')|<\eps$. Then $d(s,s')<\delta$ implies $$|Tf(s) -Tf(s')|
\le \eps \int_K |f(t)|\ud\mu(t) \le \eps \mu(K)\n f\n_\infty .$$
As a result, $T$ acts as a linear operator
on $C(K)$. To prove boundedness, we estimate
$$ |Tf(s)| \le \int_K |k(s,t)||f(t)|\ud \mu(t) \le \mu(K)\n k\n_\infty \n f\n_\infty .$$
Taking the supremum over $s\in K$ gives
$$ \n Tf \n_\infty \le \mu(K)\n k\n_\infty\n f\n_\infty.$$
It follows that $T$ is bounded and $\n T\n \le \mu(K)
\n k\n_\infty$.

For kernels $k\in L^\infty(K\times K,\mu\times\mu)$,
after a measurable representative of $k$ has been chosen, the same formula,
defined for almost all $s\in K$, determines a bounded operator on
$L^\infty(K,\mu)$. The resulting equivalence class is independent of
the chosen representative. If $k\in L^2(K\times K,\mu\times\mu)$,
the formula similarly determines a bounded operator on
$L^2(K,\mu)$ satisfying
$\n T\n \le \n k\n_2$.
Indeed, by the Cauchy--Schwarz inequality (its abstract version for Hilbert spaces will be proved in Chapter \ref{ch:Hilbert-spaces}) and Fubini's theorem we obtain
\begin{align*}
\ & \int_K \Bigl|  \int_K k(s,t)f(t)\ud \mu(t)\Bigr|^2\ud \mu(s)
\\ & \qquad \le \int_K\Bigl( \int_K |k(s,t)|^2\ud \mu(t)\Bigr)\Bigl(  \int_K |f(t)|^2\ud \mu(t)\Bigr)\ud\mu(s)
=  \n k\n_2^2 \n f\n_2^2
\end{align*}
and the claim follows. This inequality generalises the estimate of Example \ref{ex:matrices-bdd}.
\end{example}

\begin{example}[Volterra operator]\label{ex:Volterra}
 For all $f\in L^2(0,1)$, the Cauchy--Schwarz inequality implies that the indefinite integral
$$ Tf(s):= \int_0^s f(t)\ud t, \quad s\in [0,1],$$
is well defined and that $|Tf(s) - Tf(s')| \le |s-s'|^{1/2} \n f\n_2$ for all $s,s'\in [0,1].$ From this we infer that $Tf\in C[0,1]$ and, by taking $s'=0$, that $\n Tf\n_\infty \le \n f\n_2.$
This implies that
$T$ is bounded from $L^2(0,1)$ into $C[0,1]$ with norm $\n T\n \le 1$.

By composing $T$ with the natural inclusion mapping from $C[0,1]$ into $L^2(0,1)$, we may view the indefinite integral as a bounded operator on $L^2(0,1)$ of norm at most $1$. A sharper bound is obtained
by applying the last part of the preceding example (with $k(s,t) = \one_{(0,s)}(t)$).
This shows that $T$ is bounded as an operator on $L^2(0,1)$ with norm  $$\n T\n \le \n k\n_2 = 1/\sqrt 2 \approx 0.7071\dots.$$
Interestingly, this norm bound is not sharp; it can be shown that
$$\n T\n = 2/\pi \approx 0.6366\dots.$$
This will be proved using the spectral theory of selfadjoint operators
in Chapter \ref{chap:Hilbert-operators}.
\end{example}

As demonstrated by this brief list of examples, operators that naturally occur in analysis tend to be bounded. This raises the natural question whether linear operators acting between Banach spaces
$X$ and $Y$ are always bounded.
If one is willing to accept the Axiom of Choice, the answer is
negative, even for separable Hilbert spaces $X$ and $Y=\K$
(see Problem \ref{prob:nonbdd-operator}). In Zermelo--Fraenkel set theory without the Axiom of Choice, it is consistent that every linear operator acting between Banach spaces is bounded. The reader is referred to the Notes for a further discussion of this topic.

\section{Finite-Dimensional Spaces}\label{sec:eqnorm}

The aim of this section is to prove that every finite-dimensional normed space is a Banach space. This will be deduced as an easy consequence of the fact that any two norms on a finite-dimensional normed space are equivalent,
in the sense made precise in the next definition.

\begin{definition}[Equivalent norms]
 Two norms $\n \cdot \n$ and $\nn\cdot\nn$ on a vector space $X$ are {\em equivalent}\index{equivalent norm}
 if there exist constants $0<c\le C<\infty$ such that for all $x\in X$ we have
 $$ c \n x\n \le \nn x\nn \le C\n x\n.$$
\end{definition}

\begin{example}\label{ex:prod-norm-eq}
Any two product norms on the product  $X = X_1\times\cdots\times X_N$ of normed spaces are equivalent.
Indeed, \eqref{eq:prod-norm} shows that every product norm on $X$ is equivalent to the product norm $\n x\n_1:= \sum_{n=1}^N \n x_n\n$ on $X$.
\end{example}

In the above situation we have the inclusions of open balls $$B_{\n \cdot\n}(x;r/C) \subseteq B_{\nn \cdot\nn}(x;r) \subseteq
B_{\n \cdot\n}(x;r/c).$$ Hence if two norms on a given vector space are equivalent, the resulting normed spaces
have the same open sets. This implies that topological notions such as openness, closedness, compactness, and convergence are preserved when passing to an equivalent norm.

\begin{theorem}[Equivalence of norms in finite dimensions]\label{thm:equiv-norms-fd}
Any two norms on a finite-dimensional vector space are equivalent.
\end{theorem}
\begin{proof}
Let $(X,\n\cdot\n)$ be a finite-dimensional normed space,
say of dimension $d$. If $d=0$, there is nothing to prove, so we may assume $d\ge1$. Let $(x_j)_{j=1}^d$ be a basis for $X$.
Relative to this basis, every $x\in X$ admits a unique representation $x = \sum_{j=1}^d c_j x_j$.
We may use this to define a norm $\n \cdot\n_2$ on $X$ by
\begin{align*}
\Big\n \sum_{j=1}^d c_jx_j\Big\n_2 := \Bigl(\sum_{j=1}^d |c_j|^2\Bigr)^{1/2}\!.
\end{align*}
The theorem follows once we have shown that the norms $\n \cdot\n$ and $\n \cdot\n_2$ are equivalent.

Let $M := \max_{1\le j\le d}\n x_j\n$.
Applying the triangle and Cauchy--Schwarz inequalities, we find that for any
$x= \sum_{j=1}^d c_j x_j$ we have
\begin{align}\label{eq:equiv-norms-fd} \n x\n \le \sum_{j=1}^d | c_j| \n x_j\n \le M \sum_{j=1}^d |c_j| \le Md^{1/2}
\Bigl(\sum_{j=1}^d |c_j|^2\Bigr)^{1/2} = Md^{1/2}\n x\n_2.
\end{align}
This gives one of the two inequalities in the definition of equivalence of norms.

To prove that a similar inequality holds in the opposite direction,
let $S_2$ denote the unit sphere in $(X,\n \cdot\n_2)$.
Since $(c_1,\dots,c_d)\mapsto \sum_{j=1}^d c_j x_j$ maps the unit sphere of $\K^d$
isometrically (hence continuously)
onto $S_2$, $S_2$ is compact. Consider the identity mapping
$I: x\mapsto x$, viewed as a mapping from $(X,\n\cdot\n_2)$ to $(X,\n\cdot\n)$. The inequality \eqref{eq:equiv-norms-fd}
implies that $I$ is bounded and therefore continuous. Since taking norms is continuous as well and
$S_2$ is compact, the mapping $x\mapsto \n Ix\n$ is continuous from $S_2$ to $[0,\infty)$ and attains a minimum at some
point $x_0\in S_2$.

Denoting this minimum by $m$, we claim that $m>0$. It is clear that $m\ge 0$.
Reasoning by contradiction, if we had $m=\n Ix_0\n = 0$, then
$Ix_0= 0$ in $X$, hence $x_0 = 0$ in $X$. Then $\n x_0\n_2 = 0$, while at the same time
$\n x_0 \n_2 =1$ because  $x_0\in S_2$, a contradiction.

For any nonzero $x\in X$ we have $\frac{x}{\n x\n_2}\in S_2$ and therefore
$\n I\frac{x}{\n x\n_2}\n \ge m.$
This gives $$m\n x\n_2 \le \n Ix\n = \n x\n$$
for nonzero $x\in X$; for trivial reasons it also holds for $x=0$.
\end{proof}

\begin{corollary}\label{cor:fdns-Banach} Every $d$-dimensional normed space is isomorphic to $\K^d$\!. In particular, every finite-dimensional normed space is a Banach space.
\end{corollary}
\begin{proof}
The first assertion has been proved in the course of the proof of Theorem \ref{thm:equiv-norms-fd}, and the second assertion follows from it since $\K^d$ is complete.
\end{proof}

\begin{corollary}\label{cor:findim-closed} Every finite-dimensional subspace of a normed space is closed.
\end{corollary}
\begin{proof}
By Corollary \ref{cor:fdns-Banach}, every finite-dimensional subspace of a normed space is complete, and it has been shown in
the first paragraph of Section \ref{subsec:sub-quot-dirsum} that every complete subspace of a normed space is closed.
\end{proof}

\begin{corollary}\label{cor:equiv-norms-fd} Every linear operator from a finite-dimensional normed space $X$ into a normed space $Y$ is bounded.
\end{corollary}
\begin{proof}
We may assume that $\dim(X)\ge 1$. Let $(x_j)_{j=1}^d$ be a basis for $X$. If $T:X\to Y$ is linear, for $x = \sum_{j=1}^d c_j x_j$ we
obtain, by the Cauchy--Schwarz inequality,
$$ \n T x\n = \Big\n \sum_{j=1}^d c_j Tx_j\Big\n \le  \sum_{j=1}^d |c_j| \n Tx_j\n
\le Md^{1/2}\n x\n_2,$$
where $\n x\n_2 := (\sum_{j=1}^d |c_j|^2)^{1/2}$ as in Theorem \ref{thm:equiv-norms-fd} and $M:= \max_{1\le n\le d} \n Tx_n\n$.
By Theorem \ref{thm:equiv-norms-fd} there exists a constant $K\ge 0$ such that $\n x\n_2 \le K \n x\n$ for all $x\in X$. Combining this with the preceding estimate we obtain
$$ \n Tx\n \le Md^{1/2} \n x\n_2 \le KMd^{1/2} \n x\n.$$
This means that $T$ is bounded with norm at most $KMd^{1/2}$.
\end{proof}

Every bounded subset of a finite-dimensional normed space $X$ is relatively compact; this follows from the corresponding result for $\K^d$ and the fact that $X$ is isomorphic to $\K^d$ for some $d\in \N$ by Corollary \ref{cor:fdns-Banach}. Conversely, a normed space with the property that every bounded subset is relatively compact is finite-dimensional:

\begin{theorem}[Finite-dimensional spaces]\label{thm:FD-compact}
The closed unit ball of a normed space $X$ is compact if and only if $X$ is finite-dimensional.
\end{theorem}

The proof depends on the following lemma:

\begin{lemma}[Riesz]\label{lem:Riesz}\index{lemma!Riesz}
 If $Y$ is a proper closed subspace of a normed space $X$, then for every $\e>0$ there exists a norm-one vector $x\in X$ with
 $d(x,Y) \ge 1-\e$.
\end{lemma}
Here, $d(x,Y) = \inf_{y\in Y} \n x-y\n$ is the distance from $x$ to $Y$.
\begin{proof}
 Fix any $x_0\in X\setminus Y$.
Note that $d(x_0,Y) >0$: otherwise, we could select elements $y_n\in Y$ such that $\limn y_n = x_0$;  the closedness of $Y$ would then imply $x_0\in Y$.
 Fix $\e>0$ and choose $y_0\in Y$ such that $\n x_0-y_0 \n \le (1+\e)d(x_0,Y)$.
 The vector $(x_0-y_0)/\n x_0-y_0\n$ has norm one, and for all $y\in Y$ we have
 $$ \Bigl\n \frac{x_0-y_0}{\n x_0-y_0\n } - y\Bigr\n = \Bigl\n \frac{x_0 - y_0 - \n x_0-y_0\n y}{\n x_0-y_0\n}\Bigr\n
 \ge \frac{d(x_0,Y)}{(1+\e)d(x_0,Y)} = \frac1{1+\e} \ge 1-\e.$$
 It follows that
 $$ d\Bigl(\frac{x_0-y_0}{\n x_0-y_0\n }, Y\Bigr) \ge \frac1{1+\e}.$$
 Since $\e>0$ was arbitrary, this completes the proof.
\end{proof}

\begin{proof}[Proof of Theorem \ref{thm:FD-compact}]
It remains to prove the `only if' part. Suppose that $X$ is infinite-dimensional and pick an arbitrary norm-one vector $x_1\in X$. Proceeding by induction, suppose that norm-one vectors $x_1,\dots, x_n\in X$ have been chosen such that
$\n x_k - x_j\n  \ge \frac12$ for all $1\le j,k\le n$ with $j\not=k$. Choose a norm-one vector $x_{n+1}\in X$ by applying Riesz's lemma
to the proper closed subspace $Y_n = \operatorname{span}\{x_1,\dots,x_n\}$ and $\e = \frac12$
(that $Y_n$ is closed follows from Corollary \ref{cor:findim-closed}). Then
$\n x_{n+1}-x_j\n  \ge \frac12$ for all $1\le j\le n$.
The resulting sequence $(x_n)_{n\ge 1}$ is contained in the closed unit ball of $X$ and satisfies $\n x_j-x_k\n\ge \frac12$ for all $j,k\ge1$ with $j\not=k$, so $(x_n)_{n\ge 1}$ has no convergent subsequence. It follows that the closed unit ball of $X$ is not compact.
\end{proof}

\section{Compactness}\label{sec:compactness}

Let $X$ be a normed space. By Theorem \ref{thm:FD-compact}, the collections of bounded subsets of $X$ and relatively compact subsets of $X$ coincide if and only if $X$ is finite-dimensional. Thus, in infinite-dimensional spaces, relative compactness is a stronger property than boundedness. The purpose of the present section is to record some easy but useful general results on compactness that will be frequently used. Compactness in the spaces $C(K)$ and $L^p(\R^d)$ will be studied in the next chapter, and {\em compact operators}, which map bounded sets into relatively compact sets, are studied in Chapter \ref{ch:CompactOperators}.

By a general result in the theory of metric spaces (Theorem \ref{thm:totally-bounded}), every relatively compact set in a normed space is totally bounded, and the converse holds in Banach spaces. This fact is used in the proof of the following necessary and sufficient condition for relative compactness. For sets $A$ and $B$ in a vector space $V$ we write $$A+B := \{u+v:\, u\in A,\, v\in B\}.$$

\begin{proposition}\label{prop:compact-totbdd}
 A subset $S$ of a Banach space $X$ is relatively compact if and only if for all $\eps>0$
 there exists a relatively compact set $K_\eps\subseteq X$ such that $S\subseteq K_\eps + B(0;\eps)$.
\end{proposition}
\begin{proof}
 `If': \ The existence of the sets $K_\eps$ implies that $S$ is totally bounded and hence relatively compact, for if the balls $B(x_{1,\eps};\eps),\dots,B(x_{n_\eps,\eps};\eps)$ cover $K_\eps$, then the balls $B(x_{1,\eps};2\eps),\dots,B(x_{n_\eps,\eps};2\eps)$ cover $S$.

 \smallskip
 `Only if': This is trivial (if $S$ is relatively compact, take $K_\eps = S$ for all $\eps>0$).
\end{proof}

The {\em convex hull}\index{convex!hull} of a subset $F$ of a vector space $V$ is the smallest convex set in $V$ containing $F$.
This set is denoted by $\co(F)$.\index{$C00$@${\rm co}(F)$} When $F$ is a subset of a normed space, the closure of $\co(F)$ is denoted by $\ov{\co}(F)$ and is referred to as the {\em closed convex hull}\index{closed!convex hull} of $F$.

As a first application of
Proposition \ref{prop:compact-totbdd} we have the following result.

\begin{proposition}
The closed convex hull of a compact set in a Banach space is compact.\index{compactness!of the closed convex hull}
\end{proposition}

\begin{proof}
Let $K$ be a compact subset of the Banach space $X$. If $K=\emptyset$, there is nothing to prove, so we assume that $K$ is nonempty. Let $\eps>0$ be arbitrary, and choose $\xi_1,\ldots,\xi_M\in K$ so that the open balls $B(\xi_1;\eps),\ldots,B(\xi_M;\eps)$ cover $K$, and put $F:=\{\xi_1,\ldots,\xi_M\}$. The set $\co(F)$ is compact, being the image of the compact set $\Delta_M:=\{(\la_1,\ldots,\la_M)\in[0,1]^M:\,\sum_{m=1}^M\la_m=1\}$ under the continuous mapping $(\la_1,\ldots,\la_M)\mapsto\sum_{m=1}^M\la_m\xi_m$.

We claim that $\co(K)\subseteq\co(F)+B(0;\eps)$. Indeed, if $x\in\co(K)$, say $x=\sum_{j=1}^k\la_jx_j$ with $x_j\in K$, $0\leq\la_j\leq 1$, and $\sum_{j=1}^k\la_j=1$, choose $m_j\in\{1,\ldots,M\}$ such that $x_j\in B(\xi_{m_j};\eps)$. Then $\sum_{j=1}^k\la_j\xi_{m_j}\in\co(F)$, and
\begin{align*}
\Bigl\n x-\sum_{j=1}^k\la_j\xi_{m_j}\Bigr\n
&\leq \sum_{j=1}^k\la_j\n x_j-\xi_{m_j}\n < \sum_{j=1}^k\la_j\eps = \eps.
\end{align*}
This proves the claim. By Proposition \ref{prop:compact-totbdd}, $\co(K)$ is relatively compact. Therefore $\ov{\co}(K)$ is compact.
\end{proof}

The second result asserts that strong convergence implies uniform convergence on relatively compact sets.

 \begin{proposition}\label{prop:uniform-limits-on-K}  Let $X$ and $Y$ be normed spaces, let the operators $T_n\in \calL(X,Y)$, $n\ge 1$, be uniformly bounded,
 and let $T\in \calL(X,Y)$. If $\limn T_n = T$ strongly, then for all relatively compact subsets $K$ of $X$ we have
   $$\limn \sup_{x\in K} \n T_n x-Tx\n = 0.$$
 \end{proposition}

 It will be shown in Proposition \ref{prop:strong-limits} that if $X$ is a Banach space, strong convergence $T_n\to T$ already implies uniform boundedness of the operators $T_n$.

 \begin{proof}
Let $K$ be a relatively compact subset of $X$, let $\eps>0$ be arbitrary, and select finitely many open balls $B(x_1;\eps),\dots,B(x_k;\eps)$ covering $K$. Choose $N\ge 1$ so large that $\n T_n x_j- Tx_j\n <\eps$ for all $n\ge N$ and $j=1,\dots,k$. Let $M:= \sup_{n\ge 1} \n T_n\n$; this number is finite by assumption, and from $\n Tx\n =\limn \n T_n x\n \le M\n x\n$ we infer that $\n T\n \le M$. Fixing an arbitrary $x\in K$, choose $1\le j_0\le k$ such that $\n x-x_{j_0}\n<\eps$. Then, for $n\ge N$,
\begin{align*} \n T_n x- Tx\n &\le \n T_n x- T_n x_{j_0}\n +\n T_n x_{j_0}- Tx_{j_0}\n + \n T x_{j_0}- Tx\n
 \\ & \le M\eps + \eps +M\eps = (2M+1)\eps.
\end{align*}
Taking the supremum over $x\in K$, we find that if $n\ge N$, then
$$ \sup_{x\in K}  \n T_n x- Tx\n \le (2M+1)\eps. $$
Since $\eps>0$ was arbitrary, this proves the final assertion.
\end{proof}

\section{Integration in Banach Spaces}

In a variety of circumstances, some of which will be encountered in later chapters, one wishes to integrate $X$-valued functions, where $X$ is a Banach space. In order to have the tools available when they are needed, we insert a brief discussion of the $X$-valued counterparts of the Riemann and Lebesgue integrals.

\subsection{The Riemann Integral}\label{subsec:integral-Riemann}

Let $K$ be a compact metric space and let $\mu$ be a finite Borel measure on $K$.
We will set up the Riemann integral with respect to $\mu$ for continuous functions $f:K\to X$.
To this end we need the following terminology. A {\em partition} of $K$ is a finite collection of nonempty pairwise disjoint Borel
subsets of $K$ whose union equals $K$. The {\em mesh} of a partition is the maximum of the diameters of its subsets.

\begin{proposition}[Riemann integral]
Let $\mu$ be a finite Borel measure on a nonempty compact metric space $K$, let $X$ be a Banach space, and let $f:K\to X$
be a continuous function. There exists a unique element in $X$, denoted by $\int_K f\ud \mu$, with the following property: for every $\eps>0$ there exists a $\delta>0$ such that whenever $(K_n)_{n=1}^N$ is a partition of $K$ of mesh less than $\delta$
and $(t_n)_{n=1}^N$ is a collection of points in $K$ with $t_n\in K_n$ for all $n=1,\dots,N$, then
$$ \Big\n \int_K f\ud \mu - \sum_{n=1}^N \mu(K_n) f(t_n)\Big\n < \eps.$$
\end{proposition}

The proof of this proposition follows the undergraduate construction of the Riemann integral for continuous
functions $f:[0,1]\to \K$ step by step and is therefore omitted.
The element $\int_K f\ud \mu$ is called the {\em Riemann integral of $f$ with respect to $\mu$}.\index{Riemann integral} Whenever this is convenient we use the more elaborate notation $\int_K f(t)\ud \mu(t)$.

\begin{proposition}\label{prop:Riem-int-norm}
Let $\mu$ be a finite Borel measure on a nonempty compact metric space $K$, let $X$ be a Banach space, and let $f:K\to X$
be a continuous function.
Then
$$ \Bigl\n \int_K f\ud \mu\Bigr\n \le \int_K \n f\n\ud \mu.$$
\end{proposition}
\begin{proof}
For any partition $(K_n)_{n=1}^N$ of $K$ and any collection of points  $(t_n)_{n=1}^N$ in $K$ with $t_n\in K_n$ for all $n=1,\dots,N$ we have
$$
 \Bigl\n \sum_{n=1}^N \mu(K_n)f(t_n)\Bigr\n
 \le \sum_{n=1}^N \mu(K_n)\n f(t_n)\n$$
 by the triangle inequality. The result follows by taking the limit along any sequence of
 partitions whose meshes tend to zero.
\end{proof}

In the special case where $K = [0,1]$ and $\mu$ is the Lebesgue measure,
the usual calculus rules apply (defining differentiability of an $X$-valued function in the obvious way):

\begin{proposition}\label{prop:der-f-zero}
 Let $X$ be a Banach space and let $f: [0,1]\to X$ be a function. Then:
 \begin{enumerate}[label={\rm(\arabic*)}, leftmargin=*]
  \item\label{it:der-f-zero0} if $f$ is differentiable at the point $t_0\in [0,1]$, then $f$ is continuous at $t_0$;
  \item\label{it:der-f-zero1} if $f$ is differentiable on $(0,1)$ and $f'\equiv 0$ on $(0,1)$, then $f$ is constant on $(0,1)$;
  \item\label{it:der-f-zero2} if $f$ is continuously differentiable on $[0,1]$, then
  $$ \int_0^1 f'(t)\ud t = f(1)-f(0).$$
 \end{enumerate}
\end{proposition}
\begin{proof}
\smallskip
\ref{it:der-f-zero0}: \
Fix an arbitrary $\eps>0$. The assumption implies that there exists
$\delta>0$ such that, for all $t\in[0,1]$ with
$0<|t-t_0|<\delta$,
$$
 \Big\n \frac{f(t)-f(t_0)}{t-t_0}-f'(t_0)\Big\n <\eps .
$$
For such $t$ we have
$$
  \n f(t)-f(t_0)\n
  <(\eps+\n f'(t_0)\n)|t-t_0|.
$$
For $t=t_0$ the same continuity conclusion is trivial. Hence $f$ is
continuous at $t_0$.

\smallskip\ref{it:der-f-zero1}: \ The usual calculus proof via Rolle's theorem does not extend to the present setting, as it uses
the order structure of the real numbers.

Fix $\eps>0$, and let $[a,b]\subseteq(0,1)$ be a closed
subinterval. For each $t\in[a,b]$, the assumption $f'(t)=0$ gives a number
$\delta_t>0$ such that $[t-\delta_t,t+\delta_t]\subseteq(0,1)$ and
$$\|f(s)-f(t)\|\le\eps|s-t|,\quad |s-t|\le\delta_t.
$$
The intervals $I_t:=(t-\delta_t,t+\delta_t)$, $t\in[a,b]$,
form an open cover of $[a,b]$. By compactness, this cover
admits a finite subcover. Choose one with the fewest members,
and list its intervals as $I_{t_1},\ldots,I_{t_k}$ with
$t_1<\cdots<t_k$.
No interval in this subcover contains another, since otherwise
the smaller interval could be removed. Hence the left
endpoints, and likewise the right endpoints, are strictly
increasing. Consequently, $a\in I_{t_1}$, $b\in I_{t_k}$,
and consecutive intervals overlap: otherwise a point between
them would remain uncovered.

For $j=1,\ldots,k-1$, the intersection
$I_{t_j}\cap I_{t_{j+1}}$ intersects $[t_j,t_{j+1}]$.
Indeed, if a point of the intersection lies to the left of
$t_j$, then $t_j$ belongs to both intervals; if it lies to
the right of $t_{j+1}$, then $t_{j+1}$ belongs to both.
Choose $s_j\in I_{t_j}\cap I_{t_{j+1}}\cap[t_j,t_{j+1}]$,
and set $s_0:=a$ and $s_k:=b$. Then
$s_{j-1}\le t_j\le s_j$ for $j=1,\ldots,k$.
In particular we have $a=s_0\le s_1\le\cdots\le s_k=b$.
Since both endpoints belong to $I_{t_j}$, we also have
$[s_{j-1},s_j]\subseteq I_{t_j}$.

For each $j=1,\ldots,k$, since $s_{j-1},s_j\in I_{t_j}$, we have
\begin{align*}
 \|f(s_j)-f(s_{j-1})\|
 &\le \|f(s_j)-f(t_j)\|+\|f(t_j)-f(s_{j-1})\|  \\
 &\le \eps |s_j-t_j|+\eps|t_j-s_{j-1}|
  = \eps(s_j-s_{j-1}),
\end{align*}
where in the last step we used $t_j\in[s_{j-1},s_j]$. Summing over
$j$ gives
\begin{align*}
 \|f(b)-f(a)\|
 \le \sum_{j=1}^k \|f(s_j)-f(s_{j-1})\|
 \le \eps\sum_{j=1}^k(s_j-s_{j-1})
  = \eps(b-a).
\end{align*}
Since $\eps>0$ was arbitrary, this proves $f(a)=f(b)$. Since $[a,b]\subseteq(0,1)$ was arbitrary, $f$ is
constant on $(0,1)$.

\smallskip
\ref{it:der-f-zero2}: \ For the function $g:[0,1]\to X$,  $g(t) := f(t) - \int_0^t f'(s)\ud s$,
we have
$$ \lim_{h\to 0} \frac1h(g(t+h)-g(t)) = f'(t) -\lim_{h\to 0} \frac1h\int_t^{t+h} f'(s)\ud s = 0$$
by continuity, and therefore $g$ is continuously differentiable on $[0,1]$ with de\-riv\-a\-tive $g' = 0$.
It follows from \ref{it:der-f-zero1} that $g$ is constant on $(0,1)$, hence on $[0,1]$ by continuity, and then $g(0)=f(0)$ implies
$$f(t)-  \int_0^t f'(s)\ud s = g(t) = g(0) = f(0), \quad t\in [0,1].$$
Taking $t=1$ gives the result.
\end{proof}

In Chapter \ref{ch:duality} we will sketch a different proof using duality.

\subsection{The Bochner Integral}\label{subsec:Bochner}

We turn next to the more delicate problem of generalising the Lebesgue integral to functions taking values
in a Banach space $X$.

In what follows we fix a measurable space $(\Om,\calF)$.
Experience shows that if one attempts to define the measurability of a function  $f:\Omega\to X$ by imposing that $f^{-1}(B)$ be in $\calF$ for all Borel (equivalently, for all open) subsets $B$ of $X$, one arrives at a notion of measurability that is not very practical, the problem being that it does not connect well with approximation theorems such as the dominated convergence theorem. In fact it is better to start from the following necessary and sufficient condition for measurability in the scalar-valued setting: {\em A scalar-valued function is measurable if and only if it is the pointwise limit of a sequence of simple functions.}

For a function $f: \Omega\to \K$ and $x\in X$ we define $f\otimes x:\Omega\to X$ by\index{$F$@$f\otimes x$}
\begin{align}\label{eq:fotx} (f\otimes x)(\omega):= f(\omega)x.
\end{align}

\begin{definition}[Simple functions, strong measurability]\label{def:strongmeas}
  A function $f:\Om\to X$ is called {\em simple}\index{simple function!vector-valued} if it is a finite linear combination of functions of the form $\one_F\otimes x$ with $F\in \calF$ and $x\in X$, and
{\em strongly measurable}\index{strongly measurable} if it is the
  pointwise limit of a sequence of simple functions.
\end{definition}

A scalar-valued  function is strongly measurable if and only if it is measurable, and for such functions we omit the adjective `strongly'\!.

\begin{theorem}[Pettis measurability theorem, first version]\label{thm:Pettis}\index{theorem!Pettis measurability} Let $(\Om,\calF)$ be a measurable space.
 A function $f:\Om\to X$ is strongly measurable if and only if $f$ takes its values in a separable closed subspace $X_0$ of $X$ and the nonnegative functions $\n f(\cdot)-x_0\n$
 are measurable for all $x_0\in X_0$. In this case, $\n f(\cdot)-x\n$ is measurable for every $x\in X$.
\end{theorem}

A second version of this theorem will be proved in Chapter \ref{ch:duality} (see Theorem \ref{thm:Pettis-secondversion}).

\begin{proof}
`If': \
If $X_0=\{0\}$, the assertion is trivial. Otherwise, let
$(x_n)_{n\ge1}$ be a dense sequence in $X_0$ with all terms distinct, and
define the functions $\phi_n:X_0\to\{x_1,\dots,x_n\}$ as follows.
For each $y\in X_0$ let $k(n,y)$ be the least integer $1\le k\le n$
such that
$$\| y - x_{k}\| = \min_{1\le j\le n}\| y - x_j\|,$$
and put $\phi_n(y) := x_{k(n,y)}.$
Since $(x_n)_{n\ge 1}$ is dense in $X_0$ we have
$$\limn \| \phi_n(y)-y\| = 0, \quad y\in X_0.$$
Now define $\psi_n(\om) := \phi_n(f(\om))$.
We have
\begin{align*}
  \{\om\in\Om:\,\psi_n(\om)=x_1\}
  =\Big\{\om\in\Om:\, \|f(\om)-x_1\|=\min_{1\leq j\leq n}\|f(\om)-x_j\|\Big\}
\end{align*}
and, for $2\le k\le n$,
\begin{align*}
 \ & \{\om\in\Om:\,\psi_n(\om)=x_k\}
 \\ & \quad =\Big\{\om\in\Om:\, \|f(\om)-x_k\|=\min_{1\leq j\leq n}\|f(\om)-x_j\|<\min_{1\leq j\le k-1}\|f(\om)-x_j\|\Big\}.
\end{align*}
In both identities, the set on the right-hand side is in $\calF$\!.
Hence each $\psi_n$ is simple, takes values in $X_0$, and for all $\om\in\Om$ we have
$$ \limn \| \psi_n(\om)-f(\om)\| =  \limn \| \phi_n(f(\om))- f(\om)\| = 0.
$$

`Only if': \ Let $f_n\to f$ pointwise with each $f_n$ simple.
Let $X_0$ be the closed linear span of the ranges of the
functions $f_n$. Then $X_0$ is separable and $f$ takes its
values in $X_0$. For every $x\in X$, each scalar-valued function
$\om\mapsto\n f_n(\om)-x\n$ is measurable, since $f_n$ is simple.
Therefore
$$
\n f(\om)-x\n=\limn\n f_n(\om)-x\n,\quad \om\in\Om,
$$
is measurable as well. This proves both the required condition
for $x\in X_0$ and the final assertion.
\end{proof}

\begin{corollary}\label{cor:ptw-lim-strmeas}
If $\limn f_n = f$ pointwise, with each $f_n$ strongly measurable, then $f$ is strongly measurable.
\end{corollary}
\begin{proof} We check the conditions of the Pettis measurability theorem.
Every function $f_n:\Om\to X$ is the pointwise limit of a sequence of simple functions $f_{nm}:\Om\to X$, and every $f_{nm}$ takes only finitely many values. It follows that $f$ takes its values in the closed linear span of these countably many finite sets, which is a separable subspace of $X$. For every $x_0$ in this separable subspace, the final assertion
of Theorem \ref{thm:Pettis} shows that $\n f_n-x_0\n$ is
measurable for every $n$. Since
$\n f-x_0\n=\limn\n f_n-x_0\n$, this function is measurable
as well. The theorem therefore implies that $f$ is strongly
measurable.
\end{proof}

From now on we let $(\Om,\calF\!,\mu)$ be a fixed measure space.

\begin{definition}[$\mu$-Simple functions]\label{def:mu-simple} A simple function
$ f$ is called {\em $\mu$-simple}\index{simple function!$\mu$-} if it admits a representation $ f = \sum_{n=1}^N \one_{F_n}\otimes x_n$ with $\mu(F_n)<\infty$ whenever $x_n\ne0$. After omitting zero terms, we define
$$ \int_\Om f\ud \mu:= \sum_{n=1}^N \mu(F_n)x_n.$$
\end{definition}

We leave it as a simple exercise to
verify that $\int_\Om f\ud \mu$ is well defined in the sense that it does not depend on the representation
of $f$ as a linear combination of functions $\one_{F_n}\otimes x_n$ with $\mu(F_n)<\infty$.
If $f$ is $\mu$-simple, then by passing to a finite measurable disjoint
refinement of the sets in a representation of $f$, we may write
$f=\sum_{n=1}^N \one_{F_n}x_n$
with the sets $F_1,\ldots,F_N$ pairwise disjoint and of finite measure.
For such a representation the triangle inequality implies
\begin{align}\label{eq:triangle-Bochner}
 \Big\n  \int_\Om f\ud \mu\Big\n
 = \big\|\sum_{n=1}^N \mu(F_n) x_n\big\|\le \sum_{n=1}^N \mu(F_n)\|x_n\|
 = \int_\Om \n f\n\ud \mu.
\end{align}

\begin{definition}[Bochner integral]
 A strongly measurable function $f:\Om\to X$ is said to be {\em Bochner integrable}\index{Bochner!integrable} with respect to $\mu$
 if there is a sequence of $\mu$-simple functions $f_n:\Om\to X$ such that
 \begin{align}\label{eq:def-BochnerInt1}
 \limn \int_\Om \n f-f_n\n\ud \mu = 0.
 \end{align}
 In that case we define the {\em Bochner integral}\index{Bochner!integral} of $f$ by
 \begin{align}\label{eq:def-BochnerInt2}
 \int_\Om f\ud \mu:= \limn \int_\Om f_n \ud \mu.
 \end{align}
\end{definition}

On a set where $f_n$ takes the value $x\in X$,
the function $\n f-f_n\n$ equals $\n f-x\n$, which is
measurable by the final assertion of Theorem \ref{thm:Pettis}.
Hence $\n f-f_n\n$ is measurable, and the integral in
\eqref{eq:def-BochnerInt1} is well defined.

The limit in \eqref{eq:def-BochnerInt2} exists since the assumption together with \eqref{eq:triangle-Bochner} (applied to $f_n-f_m$) implies
that $(\int_\Om f_n \ud \mu)_{n\ge 1}$ is a Cauchy sequence in $X$. We leave it as another simple exercise to
verify that $\int_\Om f\ud \mu$ is well defined in the sense that it does not depend on the sequence
of approximating functions $f_n$. It is equally elementary to verify that if $\Om = K$ is a compact metric space and $\calF$ is its Borel $\sigma$-algebra, then every continuous function $f:K\to X$ is Bochner integrable with respect to every finite Borel measure $\mu$ on $K$, and the Riemann and Bochner integrals coincide.

\begin{proposition}\label{prop:Bochner-integrable} A strongly measurable function $f:\Om\to X$ is Bochner integrable
with respect to $\mu$ if and only if $$ \int_\Om \n f\n \ud \mu < \infty.$$
In this situation we have
\begin{align*}
 \Big\n  \int_\Om f\ud \mu\Big\n \le  \int_\Om \n f\n\ud \mu.
\end{align*}

\end{proposition}
\begin{proof}
`If': \
Let $f$ be a strongly measurable function satisfying
$\int_\Om \n f\n \ud\mu < \infty.$
Let $g_n$ be simple functions such that $\limn g_n = f$  pointwise
and define $$f_n := \one_{\{\n g_n\n\le 2\n f\n\}}g_n.$$
Since $\{\|g_n\|\le 2\|f\|\}$ is measurable,
each $f_n$ is simple. Moreover, $\limn f_n = f$ pointwise. Since $\n f_n\n\le 2\n f\n$ pointwise and $\int_\Om \n f\n\ud \mu<\infty$,
each $f_n$ is $\mu$-simple and by dominated convergence we obtain
$$\limn \int_\Om \n f-f_n\n\ud\mu=0.$$

`Only if': If $f$ is Bochner integrable and the $\mu$-simple function $g:\Om\to X$ is such that
 $ \int_\Om \n f-g \n\ud \mu \le 1,$
then
$$ \int_\Om \n f\n\ud \mu \le 1 + \int_\Om \n g\n\ud \mu <\infty.$$

The final assertion follows from \eqref{eq:triangle-Bochner} by approximation.
\end{proof}

\begin{problems}

\item\label{prob:openball}
Show that in any normed space $X$, for all $x_0\in X$ and $r>0$ the
following hold:
\begin{enumerate}[\rm(a), leftmargin=*]
 \item $B(x_0;r) = \{x\in X: \, \n x-x_0\n < r\}$ is an open set.
 \item $\ov B(x_0;r) = \{x\in X: \, \n x-x_0\n \le r\}$ is a closed set.
 \item $\ov B(x_0;r) = \ov{B(x_0;r)}$, that is, $\ov B(x_0;r)$ is the closure of $B(x_0;r)$.
\end{enumerate}

\item
 Let $X$ be a normed space.
\begin{enumerate}[\rm(a), leftmargin=*]
\item\label{prob:1-2} Show that if
$x,y\in X$ satisfy $\n x-y\n < \eps$ with $0<\eps<\n x\n$,
then $y\not=0$ and $$\Big\n x - \frac{\n x\n}{\n y\n}y\Big\n < 2\eps.$$
\item Conclude that $x\mapsto x/\n x\n$ maps $X\setminus\{0\}$ continuously onto the unit sphere $S_X$.
\item Show that the constant $2$ in part \ref{prob:1-2}
is the best possible.
\end{enumerate}

\item
Show that a norm $\n\cdot\n$ on the product $X = X_1\times\cdots\times X_N$ of normed spaces is a product norm if and only if $\n x\n_\infty\le \n x\n\le \n x\n_1$ for all $x = (x_1,\dots,x_N) \in X$, where
$$ \n x\n_\infty := \max_{1\le n\le N}\n x_n\n, \qquad  \n x\n_1 :=  \sum_{n=1}^N \n x_n\n.$$

\item
Show that if $X = X_1\oplus\cdots\oplus X_N$ is a direct sum of normed spaces, then each summand $X_n$ is closed as a subspace of $X$.

\item
Determine whether the following sets are open and/or closed in $C[0,1]$:
\begin{enumerate}[\rm(a), leftmargin=*]
  \item $\{f\in C[0,1]:\, f(t)\ge 0$ for all $t\in [0,1]\}$;
  \item $\{f\in C[0,1]:\, f(t)>0$ for all $t\in [0,1]\}$.
  \end{enumerate}

\item
Determine whether the following sets are open and/or closed in $\ell^1$:
\begin{enumerate}[\rm(a), leftmargin=*]
  \item $\{a\in \ell^1:\, a_n\ge 0$ for all $n\ge 1\}$;
  \item $\{a\in \ell^1:\, a_n > 0$ for all $n\ge 1\}$.
\end{enumerate}

\item
For $1\le p\le \infty$ and integers $n_0\ge 1$, show that the coordinate mapping $E_{n_0}:\ell^p\to \K$ given by $E_{n_0}(a) := a_{n_0}$
is bounded, and find its norm.

\item
Let $c_{00}$ be the space of finitely supported scalar sequences.
Show that its completion in the supremum norm is $c_0$, whereas,
for $1\le p<\infty$, its completion in the $\ell^p$-norm is $\ell^p$.

\item
Let $X$ be a normed space.
\begin{enumerate}[\rm(a), leftmargin=*]
 \item Show that if $X$ is separable, then the completion of $X$ is separable.
 \item Show that if $X$ is a Banach space and $Y$ is a closed subspace of $X$, then $X$ is separable if and only if both $Y$ and $X/Y$ are separable.
\end{enumerate}

\item
Prove that if $X\not=\{0\}$ and  $T\in \calL(X,Y)$ is bounded, then
$$\n T\n =\sup_{\n x\n=1} \n Tx\n = \sup_{\n x\n< 1} \n Tx\n.$$

\item
Let  $1\le p< \infty$. Show that $\ell^p$ is a dense subspace of $c_0$, and that the inclusion mapping of $\ell^p$ into $c_0$ is bounded. What is its norm?

\item
This problem gives an example of a bounded operator that does not attain its norm.
 Let $X$ be the space of continuous functions $f:[0,1]\to \K$ satisfying $f(0)=0$.
\begin{enumerate}[\rm(a), leftmargin=*]
  \item Show that $X$ is a closed subspace of $C[0,1]$.
\end{enumerate}
Thus, with the norm inherited from $C[0,1]$, $X$ is a Banach space.
\begin{enumerate}[\rm(a), leftmargin=*]\setcounter{enumii}{1}
  \item Show that the operator $T: X\to\K$,
  $$ Tf := \int_0^1 f(t)\ud t,$$
  is bounded  and has norm $\n T\n =1$.
  \item Prove that $|Tf| <1$ for all $f\in X$ with $\n f\n_\infty \le 1$.
\end{enumerate}

\item\label{prob:3cex1}
Let $X_0 := C_{\rm c}^1(0,1)$ be the vector space of all $C^1$-functions $f:(0,1)\to \K$ with
compact support in $(0,1)$.
\begin{enumerate}[\rm(a), leftmargin=*]
  \item Show that $X :=  \{f\in C[0,1]:\, f(0)=f(1)=0\}$ is a Banach space and that $X_0$ can be naturally identified with a dense subspace of $X$.
  \item Show that for each $f\in X_0$ the limit $\limn T_n f$ exists in the norm of $C[0,1]$ and equals $f'$\!, where $T_n:X\to C[0,1]$ is defined by $$T_n f(t) = \frac{f(t+1/n) - f(t)}{1/n},$$ using zero extensions outside $[0,1]$.
  \item Show that for general $f\in X$ the limit $\limn T_n f$ need not exist in $C[0,1]$.
\end{enumerate}
This example shows that the uniform boundedness assumption cannot be omitted in Proposition \ref{prop:approxTn}.

\item
Show that for any normed space $X$ the restriction mapping $\ov x^* \mapsto \ov x^*|_X$ is an isometric isomorphism from $\ov X^*$ onto $X^*$.

\item
This problem gives an example of a bounded operator whose range is not closed. Consider the linear operator $T$ on $C[0,1]$ given by the indefinite integral
$$ Tf(t) = \int_0^t f(s)\ud s, \quad t\in [0,1].$$
\begin{enumerate}[\rm(a), leftmargin=*]
  \item Show that $T$ is bounded, and find its norm.
  \item Show that
  $\Ran(T)$ is not closed in $C[0,1]$.
\end{enumerate}

\item
Let $V$ be a vector space, $X$ a normed space, and $T: V\to X$ be injective and linear.
\begin{enumerate}[\rm(a), leftmargin=*]
 \item Show that $\n v\n_V := \n Tv\n$ defines a norm on $V$.
 \item Show that $T: (V, \n \cdot\n_V)\to X$ is an isometry.
\end{enumerate}

\item\label{prob:sokal-lemma}
Let $X$ and $Y$ be normed spaces
and let $T\in \calL(X,Y)$.
\begin{enumerate}[\rm(a), leftmargin=*]
\item Show that for all $x\in X$ and $r>0$ we have
$\sup_{y\in B(x;r)} \n Ty\n \ge r\n T\n.$
\item Show that if $X$ is nonzero and $T$ has a bounded inverse, then $\inf_{\n x\n=1}\n Tx\n=1/\n T^{-1}\n$.
\end{enumerate}

\item
 Let $X$ be a Banach space and $Y$ be a normed space.
Show that if $T: X\to Y$ is a bounded operator satisfying $\n Tx \n \ge C \n x\n$ for some $C>0$ and all $x\in X$, then
its range $\Ran(T)$ is complete and $T$ is an isomorphism from $X$ to $\Ran(T)$.

\item\label{prob:quotient-universal}
Let $X$ be a Banach space, $Z$ a normed space, and $T\in\mathcal L(X,Z)$. Show that $T$ induces a bounded bijection $\widetilde T:X/\ker(T)\to\Ran(T)$ given by
$$\widetilde T(x+\ker (T))=Tx, \quad x\in X,$$
and that its norm is given by $\|\widetilde T\|=\|T\|$. Show that $\widetilde T$ is an isomorphism if and only if there exists a constant $c>0$ such that
$$ \|Tx\|\ge c\,d(x,\ker(T)),\quad x\in X, $$
and show that in this case $\Ran(T)$ is closed.

\item\label{prob:BanachMazur}
Let $X$ be a nonzero separable Banach space and let $S$ be a countable dense symmetric
subset of the unit sphere of $X$. Prove that every $x\in X$ with $\n x\n<1$
can be expressed as a convergent sum
$x=\sum_{n\ge1}2^{-n}x_n$
with $(x_n)_{n\ge1}\subseteq S$.

\noindent{\em Hint:}\ Given $x\in X$ with $\|x\|<1$, choose $x_1\in S$
such that $\|2x-x_1\|<1$ if $x\neq0$. If $x=0$, choose
$x_1\in S$ and then take $x_n=-x_1$ for all $n\ge2$.
For $x\neq0$, set $r_1:=x-\tfrac12x_1$. Continue inductively.

\item\label{prob:quotient-l1}
Prove that every separable Banach space $X$ is isometric to a quotient of $\ell^1$.
\index{quotient!of $\ell^1$}

\noindent{\em Hint:}\ The case $X=\{0\}$ is immediate. Otherwise, let $S=\{s_k:\,k\ge1\}$ be a countable dense subset of the
unit sphere of $X$ which satisfies $S=-S$, and define $T:\ell^1\to X$ by
$Ta:=\sum_{k\ge1}a_ks_k$.
Use Problem~\ref{prob:BanachMazur} to show that
$B_X\subseteq T(\ov B_{\ell^1})$.
Deduce that $T$ is surjective and that the quotient operator
$T_{/\ker (T)}:\ell^1/\ker (T)\to X$
is an isometry.

\item\label{prob:Ulam-Mazur}
Let $X$ and $Y$ be normed spaces. A mapping $\phi:X\to Y$ is said to be {\em distance-preserving} if for all $x_1,x_2\in X$ we have $$\n \phi (x_1)-\phi(x_2)\n = \n x_1-x_2\n,$$ and {\em affine}\index{affine} if
for all $x_1,x_2\in X$ and all real numbers $0<\la<1$ we have $$\phi ((1-\la)x_1+\la x_2) = (1-\la)\phi (x_1) + \la \phi (x_2).$$
The aim of this problem is to prove the
{\em Mazur--Ulam theorem}:\index{theorem!Mazur--Ulam} Every bijective distance-preserving mapping $\phi :X\to Y$ between normed spaces is affine.

For any mapping $\phi:X\to Y$, define the `affine defect' relative to the pair $(x_1,x_2)\in X\times X$ by
$$
\text{def}_{(x_1,x_2)}(\phi) := \Bigl\|\phi\Bigl(\frac{x_1 + x_2}{2}\Bigr) - \frac{\phi(x_1) + \phi(x_2)}{2}\Bigr\|.
$$

We now fix a bijective distance-preserving mapping $\phi:X\to Y$.

\begin{enumerate}[\rm(a), leftmargin=*]
 \item\label{it:Ulam-Mazur1} Show that for all $x_1,x_2\in X$ we have $$\text{def}_{(x_1,x_2)}(\phi) \leq \frac{1}{2}\|x_1 - x_2\|.$$
\end{enumerate}
Let $\rho(y) = \phi(x_1) + \phi(x_2) - y$ be the reflection with respect to $\frac12(\phi(x_1) + \phi(x_2))$.
\begin{enumerate}[\rm(a), leftmargin=*]\setcounter{enumii}{1}
\item\label{it:Ulam-Mazur2}  Show that $\psi := \phi^{-1}\rho\phi: X\to X$ is a bijective distance-preserving mapping which satisfies
$$ \text{def}_{(x_1,x_2)}(\psi) =2 \,\text{def}_{(x_1,x_2)}(\phi).$$
\item\label{it:Ulam-Mazur3}  By iterating part \ref{it:Ulam-Mazur2}, conclude from part \ref{it:Ulam-Mazur1} that
$\text{def}_{(x_1,x_2)}(\phi)= 0$.
\item\label{it:Ulam-Mazur4} Deduce from part \ref{it:Ulam-Mazur3} that $\phi$ is affine.
\end{enumerate}

\item\label{prob:complexification}
Let $X$ be a real vector space. The product $X\times X$ can be given the
structure of a complex vector space by introducing a complex scalar
multiplication as follows:
$$(a+ib)(x,y) := (ax-by, bx+ay).$$
The idea is to think of the pair $(x,y)\in X\times X$ as `$x+iy$'.
\begin{enumerate}[\rm(a), leftmargin=*]
  \item\label{it:complexification1} Check that this formula for the scalar multiplication does indeed turn $X\times X$ into a complex vector space.
\end{enumerate}
The resulting complex vector space is denoted by $X_\C$.

Suppose now that $X$ is a real normed space.
\begin{enumerate}[\rm(a), leftmargin=*]\setcounter{enumii}{1}
  \item Prove that
  $$ \n (x,y) \n := \sup_{\theta\in [0,2\pi]} \n (\cos \theta) x + (\sin\theta)y\n $$
  defines a norm on $X_\C$ which turns $X_\C$ into a complex normed space. Show that $X_\C$ is a Banach space if and only if $X$ is a Banach space.
  \item\label{it:complexification3} Show that this norm on $X_\C$
  extends the norm of $X$ in the sense that
  $\n (x,0)\n = \n (0,x)\n = \n x\n$ for all $x\in X.$
  \item\label{it:complexification4} Show that $\n (x,y)\n = \n (x,-y)\n$ for all $x,y\in X.$
\item Show that if $T$ is a (real-)linear bounded operator on $X$, then $T$ extends to a bounded
(complex-)linear operator
$T_\C$ on $X_\C$ by putting $ T_\C(x,y) := (Tx,Ty)$.
\item Show that $\n T_\C\n = \n T\n$.
\end{enumerate}

\item\label{prob:nonbdd-operator}
Using Zorn's lemma, choose a basis $(h_i)_{i\in I}$ of $\ell^2$ (in the sense of linear algebra: every $h\in \ell^2$ can be represented uniquely as a finite linear combination of basis vectors) containing the standard unit vectors $(e_n)_{n\ge1}$.
Pick $i_0\in I$ with $h_{i_0}\notin\{e_n:\,n\ge 1\}$.
Let $Y:=\operatorname{span}\{h_i:\, i\in I\setminus\{i_0\}\}$.

\begin{enumerate}[\rm(a),leftmargin=*]
\item Show that $Y$ is dense in $(\ell^2,\|\cdot\|)$ and $\ell^2 = \{y+ch_{i_0}:\, y\in Y,\, c\in \K\}$.
\item Show that the mapping $\phi:\ell^2\to\K$ defined by $\phi(y+ch_{i_0}):=c$ is well defined and linear, but not bounded.
\end{enumerate}

\item
Show that two finite-dimensional normed spaces are isomorphic if and only if they have the same dimension.

\item
For each pair of integers $m,n\ge 1$, find an isomorphism from $\calL(\K^n\!,\K^m)$ to $\K^{mn}$.

\item
Let $\n\cdot\n$ and $\n\cdot\n'$ be two norms on a vector space $X$. Show that the following assertions are equivalent:
\begin{enumerate}[label={\rm(\arabic*)}, leftmargin=*]
 \item there exists a constant $C\ge 0$ such that $\n x\n \le C\n x\n'$ for all $x\in X$;
 \item every open set in $(X,\n\cdot\n)$ is open in $(X,\n\cdot\n')$;
 \item convergent sequences in $(X,\n\cdot\n')$ converge in
       $(X,\n\cdot\n)$ to the same limit;
 \item Cauchy sequences in $(X,\n\cdot\n')$ are Cauchy in $(X,\n\cdot\n)$.
\end{enumerate}

\item
Show that if two norms $\n\cdot\n$ and $\n \cdot\n'$ on a normed space $X$ are equivalent, then the norms of the
completions of $(X,\n\cdot\n)$ and $(X,\n\cdot\n')$ are equivalent (under the natural identification of these completions).

\item\label{prob:norm-dense-subspace}
Let $X$ be a Banach space with respect to the norms $\n\cdot\n$ and $\n \cdot\n'$. Suppose that $\n \cdot\n$ and $\n\cdot \n'$ agree on a subspace $Y$ that is dense in $X$ with respect to both norms.
We ask whether the norms agree on all of $X$.
\begin{enumerate}[\rm(a), leftmargin=*]
\item\label{it:norm-dense-subspace1}
Comment on the following argument:
Apply Proposition \ref{prop:extendT} to the identity mapping on $Y$, viewed as a mapping from
the normed space $(Y, \n \cdot\n)$ to the normed space $(Y,\n\cdot \n')$ and as a mapping in the opposite direction.
\item
Comment on the following argument:
Let $x\in X$ be fixed and choose a sequence $x_n\to x$ with $x_n\in Y$ for all $n\ge 1$. Then $$\n x\n = \limn \n x_n\n = \limn \n x_n\n' = \n x\n'.$$
\end{enumerate}

\item\label{prob:isol1l2linfty}
For $1\le p\le \infty$ and integers $d\ge 1$ let $\ell_d^p := (\K^d, \n \cdot\n_p)$ as in Example \ref{ex:eucliden}.
\begin{enumerate}[\rm(a), leftmargin=*]
 \item\label{it:isol1l2linfty-1} Show that if $T\in\calL(\ell_d^1, \ell_d^2)$ is an isomorphism, then $\n T\n \n T^{-1} \n\ge \sqrt{d}$.
 \item\label{it:isol1l2linfty-2} Show that if $T\in\calL(\ell_d^2, \ell_d^\infty)$ is an isomorphism, then $\n T\n \n T^{-1} \n\ge \sqrt{d}$.
 \item\label{it:isol1l2linfty-3} Show that the identity map
 $I:\ell_d^1\to \ell_d^\infty$ is an isomorphism with $\|I\|\,\|I^{-1}\|=d$.
\end{enumerate}

\item
Show that $\ell^1$ and $\ell^2$ are not isomorphic.

\item
Let $X$ and $Y$ be finite-dimensional normed spaces. Prove that if $T_n,T\in \calL(X,Y)$, then the following assertions are equivalent:
\begin{enumerate}[label={\rm(\arabic*)}, leftmargin=*]
  \item $\limn T_n = T$ uniformly;
  \item $\limn T_n = T$ strongly;
  \item $\limn T_n = T$ weakly.
\end{enumerate}

\item
Show that if $K$ is a compact subset of a Banach space $X$, then $K$ is contained in a separable closed subspace of $X$.

\item
Show that if $K$ and $K'$ are compact subsets of a Banach space $X$, then the set $K+K':=\{x+x':\, x\in K,\,x'\in K'\}$ is compact.

\item
Let $K$ and $F$ be nonempty disjoint subsets of a Banach space $X$, with
$K$ compact and $F$ closed. Show that $d(K,F):= \inf\{\n x-y\n:\, x\in K, \, y\in F\}>0.$

\item
In each of the Banach spaces $\ell^1$, $\ell^\infty$, $L^1(0,1)$, and $L^\infty(0,1)$, find a sequence $(x_n)_{n\ge 1}$ of norm-one vectors satisfying
$$ \n x_n - x_m \n = 2,\quad m\not=n.$$

\item\label{prob:compact-nullseq}
Show that a subset $K$ of a Banach space $X$ is relatively compact if and only if $K$ is contained in the closed convex hull of a
sequence $(x_n)_{n\ge 1}$ in $X$ with $\limn x_n = 0$.

\noindent{\em Hint:}\ For the `only if' part, for each $n\in\N$ cover $K$ with finitely many open balls of radius $2^{-n}$ and centres in $K$. Assign each centre
at level $n+1$ a centre at level $n$ within distance $2^{-n}$.
Choose positive numbers $a_n$ such that
$\sum_{n\ge0}a_n\le1$ and $2^{-n}/a_n\to0$.
Express each centre as its initial centre plus the increments
along its chain of predecessors, and absorb the coefficients
$a_n$ into these vectors. Include $0$ to account for unused weights.

\item
Using Riesz's lemma, show that there exists a real number $\delta\in (0,1)$ with the property that the open unit ball of any infinite-dimensional normed space contains infinitely many disjoint open balls of radius $\delta$.

\item\label{prob:tr-inv-Borel}
Using the preceding problem, show that if $X$ is a normed space supporting a translation-invariant Borel measure $\mu$ such that every nonempty open ball has positive measure and some open ball has finite measure, then $X$ is finite-dimensional.

\item
Let $(\Om,\calF)$ be a measurable space. Adapting the proof of Theorem \ref{thm:Pettis}, show that if $f:\Omega\to X$ is strongly measurable, there are simple functions $f_n:\Omega\to X$ such that $f_n\to f$ and $\n f_n\n \le \n f\n$ pointwise.

\item
Let $\mu$ be a finite Borel measure on a nonempty compact metric space $K$, and let $X$ be a Banach space. Prove that every continuous function $f:K\to X$ is Bochner integrable with respect to $\mu$ and that its Bochner and Riemann integrals are equal.

\item
Let $(\Om,\calF\!,\mu)$ be a measure space.
Show that if $T:X\to Y$ is a bounded operator and $f:\Om\to X$ is Bochner integrable with respect to $\mu$, then
 $Tf: \Om\to Y$ is Bochner integrable with respect to $\mu$ and $$T\int_\Om f\ud\mu = \int_\Om Tf\ud\mu.$$

\item\label{prob:Bochner-int-subspace}
Let $(\Om,\calF\!,\mu)$ be a measure space and let $X_0$ be a closed subspace of the Banach space $X$.
Let $f:\Om\to X$ satisfy $f(\om)\in X_0$ for all $\om\in\Om$.

\begin{enumerate}[\rm(a), leftmargin=*]
 \item Show that if $f$ is strongly measurable as an $X$-valued function, then $f$ is strongly measurable as an $X_0$-valued function.
\item Show that if $f$ is Bochner integrable as an $X$-valued function, then $f$ is Bochner integrable as an
$X_0$-valued function.
\end{enumerate}

\item
Let $(\Om,\calF\!,\mu)$ be a measure space and let $(\Om'\!,\calF')$ be a measurable space.
Let $\phi:\Om\to \Om'$ be measurable
and let $f:\Om'\to X$ be strongly measurable. Let $\nu = \mu\circ \phi^{-1}$ be the image measure of $\mu$ under $\phi$.
\begin{enumerate}[\rm(a), leftmargin=*]
 \item Show that $f\circ \phi$ is strongly measurable.
 \item Show that $f\circ \phi$ is Bochner integrable with respect to $\mu$ if and only if $f$ is Bochner integrable with respect to $\nu$, and that in this situation we have
$$\int_\Om f\circ \phi \ud \mu = \int_{\Om'} f \ud \nu.$$
\end{enumerate}

\item
Let $(\Om,\calF\!,\mu)$ be a probability space and let $f:\Om\to X$ be Bochner integrable.

\begin{enumerate}[\rm(a), leftmargin=*]
\item Show that  $\int_{\Om} f\ud \mu$ is contained in the closed convex hull of $\{f(\om): \, \om\in \Om\}.$

\item Show that the set $\{\int_F f \ud \mu: \, F\in \calF\}$ is relatively compact in $X$.

\noindent{\em Hint:} \ Approximate $f$ by simple functions and use Proposition \ref{prop:compact-totbdd}.

\end{enumerate}

\end{problems}

%% file: ch02-ClassicalBanachSpaces.tex
\chapter{The Classical Banach Spaces}\label{ch:ClassicalBanach}

\blfootnote{\rm \small \begin{quote}This book has been published by Cambridge University Press in the series `Cambridge Studies in Advanced Mathematics'. The present version is free to view and download for personal use only. Not for re-distribution, re-sale or use in derivative works.\newline {\copyright} Jan van Neerven\end{quote}}

\noindent
In this chapter we undertake a closer examination of the classical Banach spaces introduced in Chapter \ref{ch:Banach}. These spaces will provide both a testing ground for the general theory and a recurring source of examples throughout the book.

\section{Sequence Spaces}\label{sec:sequencespaces}

Besides the finite-dimensional spaces $\K^d\!$, perhaps the simplest examples of Banach spaces are provided by
the class of sequence spaces. By definition, these are spaces of sequences which, endowed with a suitable norm, turn into Banach spaces. Here we introduce the most important sequence spaces, namely, $c_0$ and
$\ell^p$\!, $1\le p\le \infty$.

\paragraph{The Spaces \texorpdfstring{$c_0$}{c\textunderscore 0} and \texorpdfstring{$\ell^\infty$}{ell-infinity}}\index{$L$@$\ell^\infty$}\index{$C$@$c_0$}

The space $c_0$ consisting of all scalar sequences $a = (a_k)_{k\ge 1}$ satisfying $\limk a_k = 0$ is a Banach space with respect to the
supremum norm $$\n a\n_\infty := \sup_{k\ge 1} |a_k|.$$
A justification of this notation is given in the next paragraph.
That this is indeed a norm is left as an exercise; the proof of completeness runs as follows.
Suppose $(a^{(n)})_{n\ge 1}$ is a Cauchy sequence in $c_0$. Then each coordinate sequence $(a_k^{(n)})_{n\ge 1}$
is Cauchy in $\K$ and therefore has a limit which we denote by $a_k$. We wish to prove that the sequence $a:= (a_k)_{k\ge 1}$
belongs to $c_0$ and that $\limn \n a^{(n)} - a\n_\infty = 0.$

Fix $\e>0$ and choose $N$ so that $\n a^{(n)} - a^{(m)}\n_\infty<\e$ for all $m,n\ge N$.
Choose $N'$ so large that $|a_k^{(N)}| < \e$ for all $k\ge N'$\!.
Then, for $k\ge N'$\!,
$$ |a_k| \le |a_k - a_k^{(N)}| + |a_k^{(N)}|  = \lim_{m\to\infty}|a_k^{(m)} - a_k^{(N)}| + |a_k^{(N)}| \le  \e + \e = 2\e.$$
It follows that $\limk a_k = 0$, so $a\in c_0$.

Finally, for all $k\ge 1$ and $m,n\ge N$ we have
$ |a_k^{(n)} - a_k^{(m)}| < \e.$
Letting $m\to\infty$ while keeping $n$ fixed, for all $k\ge 1$ we obtain
$$ |a_k^{(n)} - a_k| \le \e.$$
Taking the supremum over $k\ge 1$ we infer that
$\n a^{(n)}- a\n_\infty \le \e$ for all $n\ge N$, and the convergence $a^{(n)}\to a$ in $c_0$ follows.

Similarly, one proves that the space $\ell^\infty$ consisting of all bounded scalar sequences $a = (a_k)_{k\ge 1}$
is a Banach space with respect to the supremum norm. This space contains $c_0$ isometrically as a closed subspace.

\paragraph{The Spaces \texorpdfstring{$\ell^p$}{ell-p}}\index{$L$@$\ell^p$}

For $1\le p<\infty$, the space  $\ell^p$ of
scalar sequences $a = (a_k)_{k\ge 1}$ for which
$$ \n a\n_p:= \Bigl(\sum_{k\ge 1}|a_k|^p \Bigr)^{1/p}$$ is finite is a Banach space
with respect to the norm $\n\cdot\n_p$. That this is indeed a norm on $\ell^p$ is nontrivial; the validity of the triangle inequality
$\n a+b \n_p \le \n a\n_p + \n b \n_p$ can be proved by following the line of proof of Proposition \ref{prop:Minkowski}. Completeness of $\ell^p$ can be proved as in Theorem \ref{thm:Lp-complete}. Alternatively, these facts can be deduced as special cases of Proposition \ref{prop:Minkowski} and Theorem \ref{thm:Lp-complete} by taking $\Omega = \{1,2,3,\dots\}$ with the counting measure, that is, the measure which assigns mass $1$ to every element of $\Om$.

It is easy to see (see Problem \ref{prob:ellp}) that $1\le p\le q\le \infty$ implies $\ell^p\subseteq \ell^q$
and $$\n a\n_q\le \n a\n_p$$
for all $a\in \ell^p$\!,
and that if $a\in \ell^{p}$ for some $1\le p<\infty$, then
$$\lim_{\substack{q\to\infty\\ q\ge p}} \n a\n_q = \n a\n_\infty.$$ This justifies the notation $\n \cdot\n_\infty$ for the supremum norm.

\begin{remark}\label{rem:ellpI}
In some applications it is useful to use countable index sets $I$ other than the positive integers.
We then define
$$ \n a\n_{\ell^p(I)} := \Bigl( \sum_{n\ge 1} |a_{i_n}|^p\Bigr)^{1/p}\!,$$
where $(i_n)_{n\ge 1}$ is an enumeration of $I$, with the obvious adjustment if $p=\infty$. This definition is independent of the choice of the enumeration,
and the space $\ell^p(I)$ of all mappings $a:I\to\K$ for which this expression is finite is again a Banach space.
\end{remark}

\section{Spaces of Continuous Functions}\label{sec:CK}

In this section we study some properties of the space $C(K)$ of continuous functions defined on a compact topological space $K$.

\subsection{Completeness}\label{subsec:complete}

It is a standard result in any introductory course in analysis that the uniform limit of a sequence of continuous functions is continuous. The following theorem recasts this result as a completeness result.

\begin{theorem}[Completeness]\label{thm:CK} Let $K$ be a nonempty compact topological space.
The space $C(K)$\index{$C1a$@$C(K)$} is a Banach space with respect to the {\em supremum norm}\index{supremum norm}
$$ \n f\n_{\infty} := \sup_{x\in K} |f(x)|.$$
\end{theorem}

The elementary verification that this is indeed a norm is left to the reader.
The above supremum is finite (and actually a maximum) since $K$ is compact.

\begin{proof}
Suppose that $(f_n)_{n\ge 1}$ is a Cauchy sequence in $C(K)$. Then for each $x\in K$, $(f_n(x))_{n\ge 1}$ is a Cauchy sequence in $\K$
and therefore convergent to some limit in $\K$ which we denote by $f(x)$.
We will prove that the function $f$ thus defined is continuous and that $\limn \n  f_n-f\n_\infty = 0$.

Fix $\e>0$ and choose $N\ge 1$ so large that $\n f_n-f_m\n_\infty< \e$ for all $m,n\ge N$.
Then in particular for all $m,n\ge N$ and all $x\in K$ we have
$ |f_n(x) - f_m(x) | < \e.$
Passing to the limit $m\to \infty$ for each fixed $n\ge N$, we obtain
\begin{align}\label{eq:lim-cont} |f_n(x) - f(x) | \le \e, \quad x\in K, \ n\ge N.\end{align}
Now fix $x\in K$ arbitrary and let $U\subseteq K$ be an open set containing $x$ such that $|f_N(x) - f_N(x')| < \e$ whenever $x'\in U$.
Then, for $x'\in U$,
\begin{align*} |f(x) - f(x') |  \le  |f(x) - f_N(x) | + |f_N(x) - f_N(x') | + |f_N(x') - f(x') |
\le \e + \e +\e = 3\e,
\end{align*}
where we applied \eqref{eq:lim-cont} to $n=N$ and the points $x$ and $x'$\!.
An argument of this type is called a {\em $3\e$-argument}.
This proves the continuity of $f$ at the point $x$. Since $x\in K$ was arbitrary, $f$ is continuous and therefore
belongs to $C(K)$.
Finally, since
\eqref{eq:lim-cont}
holds for all $x\in K$, it follows that
$$\n f_n-f\n_\infty = \sup_{x\in K}|f_n(x) - f(x) | \le \e, \quad n\ge N.$$
This proves that $\limn \n f_n- f\n_\infty = 0$.
\end{proof}

We give three more examples of spaces of functions that are Banach spaces with respect to the supremum norm.
The proofs that these spaces are complete are similar to the ones for $c_0$, $\ell^\infty$, and $C(K)$, and are left as an exercise.

\begin{itemize}
 \item The space $B_{\rm b}(X)$ of bounded Borel measurable functions on a topological space $X$.

 \item The space $C_{\rm b}(X)$\index{$C_{\rm b}(X)$} of bounded continuous functions on a topological space $X$.

 \item The space $C_0(X)$\index{$C_0(X)$} of continuous functions on a locally compact topological space $X$ which vanish at infinity (the precise definitions
 are given in Section \ref{subsec:dual-CK}).
\end{itemize}

\subsection{The Stone--Weierstrass Approximation Theorem}

The Stone--Weierstrass theorem provides a useful density criterion for the spaces $C(K)$.
We begin with the more elementary Weierstrass approximation theorem for $K = [a,b]$.

\begin{theorem}[Weierstrass approximation theorem]\label{thm:Weierstrass}\index{theorem!Weierstrass approximation}
The polynomials with coefficients in $\K$ are dense in $C[a,b]$.
\end{theorem}

\begin{proof}
By translation and scaling it suffices to prove the theorem for the space $C[0,1]$.
Our proof is constructive in that it produces an actual sequence of polynomials approximating a given function.
Let $f\in C[0,1]$ be arbitrary and fixed and define the {\em Bernstein polynomials}\index{Bernstein polynomials}
associated with $f$ by
$$ B_n^{(f)}(x):= \sum_{k=0}^n \binom{n}{k} f(\frac{k}{n}) x^k (1-x)^{n-k}\!, \quad x\in [0,1], \
n\ge 1.
$$
We will show that $\limn \n B_n^{(f)} -f\n_\infty = 0.$
To begin with, the binomial identity
\begin{align}\label{eq:binom} \sum_{k=0}^n \binom{n}{k} x^k (1-x)^{n-k} = [x + (1-x)]^n = 1
\end{align}
implies
\begin{align*}
 B_n^{(f)}(x) -f(x) = \sum_{k=0}^n \binom{n}{k} x^k (1-x)^{n-k} \bigl(f(\frac{k}{n}) - f(x) \bigr).
\end{align*}
Fix an arbitrary $\e>0$. Since $f$ is uniformly continuous, there is a real number $0<\delta<1$ such that
$|f(x) - f(x')|<\e$ whenever $x,x'\in [0,1]$ satisfy $|x-x'|<\delta$.
Fix $x\in [0,1]$ and set
$I := \{0\le k\le n:\, |\frac{k}{n}-x| < \delta\}$ and $I' = \{0\le k\le n:\, k\not\in I\}$.
The sum over the indices $k\in I$ can be estimated by
\begin{align*}
\sum_{k\in I} \binom{n}{k} x^k (1-x)^{n-k} \bigl|f(\frac{k}{n}) - f(x) \bigr|
 \le \e \sum_{k=0}^n \binom{n}{k} x^k (1-x)^{n-k} = \e,
\end{align*}
while for $k\in I'$ we have $\delta^2 \le ( \frac{k}{n}-x)^2$ and therefore
\begin{align*}
\delta^2\sum_{k\in I'} \binom{n}{k} x^k (1-x)^{n-k} \bigl|f(\frac{k}{n}) - f(x) \bigr|
 & \le \sum_{k\in I'} (\frac{k}{n}-x)^2\binom{n}{k} x^k (1-x)^{n-k} \bigl|f(\frac{k}{n}) - f(x) \bigr|
\\ & \le 2\n f\n_\infty \sum_{k=0}^n (\frac{k}{n}-x)^2\binom{n}{k}  x^k (1-x)^{n-k}
\\ & \stackrel{(*)}{=} 2\n f\n_\infty\frac{ x(1-x)}{n} \le  \frac{2}{n}\n f\n_\infty,
\end{align*}
where the identity in $(*)$ follows from the following computation.
Differentiating the polynomial identity
$$
\sum_{k=0}^n\binom{n}{k}x^k(1-x)^{n-k}t^k
=(1-x+xt)^n
$$
once and twice with respect to $t$, and then setting $t=1$, gives
\begin{align*}
\sum_{k=0}^n k\binom{n}{k}x^k(1-x)^{n-k}
=nx,\qquad
\sum_{k=0}^n k(k-1)\binom{n}{k}x^k(1-x)^{n-k}
=n(n-1)x^2.
\end{align*}
Using $k^2=k(k-1)+k$ and \eqref{eq:binom}, we therefore obtain
\begin{align*}
\sum_{k=0}^n\Bigl(\frac{k}{n}-x\Bigr)^2
\binom{n}{k}x^k(1-x)^{n-k}
=\frac{n(n-1)x^2+nx}{n^2}-2x^2+x^2=\frac{x(1-x)}{n}.
\end{align*}
Combining things, we obtain
$$ | B_n^{(f)}(x) -f(x)| \le \e +  \frac{2}{\delta^2 n}\n f\n_\infty.$$
Taking the supremum over $x\in [0,1]$ and letting $n\to\infty$, we obtain
$$\limsup_{n\to\infty} \n B_n^{(f)} -f\n_\infty \le \eps.$$
This shows that $f$ can be approximated arbitrarily well by polynomials.
\end{proof}

\begin{remark}\label{rem:Bernstein-pointwise}
The same argument shows that if $f:[0,1]\to\K$ is {\em any}
bounded
function which is continuous at a point $x_0\in [0,1]$,
then $\limn B_n^{(f)}(x_0) = f(x_0)$.
\end{remark}

The proof of Theorem \ref{thm:Weierstrass} has an interesting connection with the law of large numbers.
Suppose $\xi_1,\xi_2,\xi_3,\dots$ are independent and identically distributed random variables taking the values $0$ and $1$
with probabilities $1-p$ and $p$, respectively.
Let $$S_n := \frac1n \sum_{k=1}^n \xi_k.$$
Suppose that $f:[0,1]\to\K$ is bounded and continuous at $p$. Denoting expectation\index{expected value} and probability by $\E$\index{$Ee$@$\E$} and $\mathbb{P}$,\index{$Pp$@$\mathbb{P}$} respectively, we have
$$\E f(S_n) = \sum_{k=0}^n f(\frac{k}{n}) \cdot \mathbb{P}\Big(S_n = \frac{k}{n}\Big).$$
Since $\mathbb{P}(S_n = \frac{k}{n})
= \binom{n}{k} p^k (1-p)^{n-k}\!,$ from Remark \ref{rem:Bernstein-pointwise} we therefore obtain
$$
\lim_{n\to\infty}\E f(S_n)= \lim_{n\to\infty}\sum_{k=0}^n  \binom{n}{k} f(\frac{k}{n}) p^k (1-p)^{n-k} =
\limn B^{(f)}_n(p)  = f(p).$$
In particular, taking $f=\one_{\{|\,\cdot-p|\ge\varepsilon\}}$ we recover the {\em weak law of large numbers}
$$\limn  \P(|S_n - p|\ge \eps) = 0.$$

The proof of the Weierstrass theorem using Bernstein polynomials offers little room for generalisation,
but the theorem itself does admit a far-reaching generalisation:

\begin{theorem}[Stone--Weierstrass theorem, algebra version]\label{thm:Stone-Weierstrass}\index{theorem!Stone--Weierstrass}
 Let $K$ be a compact Hausdorff space and suppose that $Y$ is a subspace of $C(K)$ with the following properties:
 \begin{enumerate}[label={\rm(\roman*)}, leftmargin=*]
  \item\label{it:SW1} $\one \in Y$;
  \item\label{it:SW2} $g\in Y$ implies $\ov g \in Y$;
  \item\label{it:SW3} $g\in Y$ and $h\in Y$ implies $gh\in Y$;
  \item\label{it:SW4} $Y$ separates the points of $K$.
 \end{enumerate}
Then $Y$ is dense in $C(K)$.
\end{theorem}
By definition, condition \ref{it:SW4} means that for any two distinct points $x,y\in K$ there exists a function $g\in Y$ such that $g(x)\neq g(y)$.

As a preliminary observation we note that it suffices to prove the theorem for real-valued functions. Indeed, the complex version of the theorem
follows from the real version as follows. If $g\in Y$, then the real-valued functions
$\Re g = \frac12(g + \ov g)$ and $\Im g = \frac1{2i}(g -\ov g)$ belong to $Y$. From this it is easy to see that the real-linear space $Y_\R$
of all real-valued functions contained in $Y$ satisfies \ref{it:SW1}--\ref{it:SW4} again.
Now if $f = u+iv\in C(K)$, we may use the real version of the theorem, with $Y$ replaced by $Y_\R$, to approximate $u$ and $v$ by functions
$u_n,v_n\in Y_\R$. Then the functions $u_n+iv_n$ approximate $f$.

The real version of the theorem will be deduced from its companion, in which condition \ref{it:SW3}
is replaced by closedness under taking pointwise absolute values:

\begin{theorem}[Stone--Weierstrass theorem, lattice version]\label{thm:Stone-Weierstrass-lattice}\index{theorem!Stone--Weierstrass}
 Let $K$ be a compact Hausdorff space and suppose that $Y$ is a subspace of $C(K)$ with the following properties:
 \begin{enumerate}[label={\rm(\roman*)}, leftmargin=*]
  \item\label{it:SWlatt1} $\one \in Y$;
  \item\label{it:SWlatt2} $g\in Y$ implies $\ov g\in Y$;
  \item\label{it:SWlatt3} $g\in Y$ implies $|g|\in Y$;
  \item\label{it:SWlatt4} $Y$ separates the points of $K$.
 \end{enumerate}
Then $Y$ is dense in $C(K)$.
\end{theorem}
\begin{proof}
Reasoning as before, it suffices to prove the theorem over the real scalars.

For the minimum $a\wedge b:= \min\{a,b\}$\index{$A$@$a\wedge b$} and maximum $a\vee b:=\max\{a,b\}$\index{$A$@$a\vee b$} we have the formulas
$$ a\wedge b = \frac12\bigl((a+b) - |a-b|\bigr), \qquad a\vee b = \frac12\bigl((a+b) + |a-b|\bigr).$$
They imply that $Y$ is closed under taking pointwise maxima and minima.

Fix $f\in C(K)$ and $\e>0$.

\smallskip
{\em Step 1} -- We prove that for each $x\in K$ there exists a function  $g_x\in Y$ such that $g_x(x) = f(x)$ and $g_x < f+\e$ pointwise.

Since $Y$ is a linear subspace containing the constant functions and separating the points of $K$,
for all $y\in K$ there exists a function $g_{xy}\in Y$ such that $g_{xy}(x) = f(x)$ and $g_{xy}(y) = f(y)$.
Indeed, if $x=y$, we may take $g_{xx}:=f(x)\one\in Y$. If $x\neq y$,
by \ref{it:SWlatt4} there exists a function $h\in Y$ with $h(x)\neq h(y)$.
Set
$$
a:=\frac{f(x)-f(y)}{h(x)-h(y)},
\qquad
b:= f(x)-a h(x),
$$
and define $g_{xy}:=a h+b\one$.
Since $Y$ is a linear subspace and $\one\in Y$, we have $g_{xy}\in Y$.
Moreover, $g_{xy}(x)=a h(x)+b=f(x)$ and $g_{xy}(y)=a h(y)+b=f(y)$.

The set $U_{xy} = \{z\in K:\,g_{xy}(z)<f(z)+\e\}$ is open and contains both $x$ and $y$. Since $K$ is compact, the open cover $\{U_{xy}:\, y\in K\}$ has a finite subcover,
say $\{U_{xy_n}: \, n=1,\dots,N_x\}$. The function $g_x := g_{xy_1}\wedge \cdots\wedge g_{xy_{N_x}}$ has the required properties.

\smallskip
{\em Step 2} -- We prove that there exists a function  $g\in Y$ such that $f-\eps < g < f+\e$; this implies $\n f-g\n_\infty\le \eps$ and concludes the proof.

For each $x\in K$ the set $U_x = \{z\in K: \, f(z)-\eps < g_x(z)\}$ is open and contains $x$. Since $K$ is compact, the open cover $\{U_{x}:\, x\in K\}$ has a finite subcover,
say $\{U_{x_n}: \, n=1,\dots,N\}$. The function $g := g_{x_1}\vee \cdots\vee g_{x_{N}}$ has the required properties.
\end{proof}

\begin{proof}[Proof of Theorem \ref{thm:Stone-Weierstrass}]
As already noted, it suffices to prove the theorem over the real scalar field.
Let $Y$ be a subspace of $C(K)$ with the properties \ref{it:SW1}--\ref{it:SW4}  stated in Theorem \ref{thm:Stone-Weierstrass}.
If we can approximate any $f\in C(K)$ with
functions from the closure $\ov Y$, we can also approximate with functions from $Y$. Since $\ov Y$ also satisfies the properties \ref{it:SW1}--\ref{it:SW4} of Theorem \ref{thm:Stone-Weierstrass},
we may assume that $Y$ is closed.
The strategy of the proof is then to show, under this additional closedness assumption, that $Y$ satisfies the assumptions of
Theorem \ref{thm:Stone-Weierstrass-lattice}. For this we need to show that if $f\in Y$, then $|f|\in Y$ as well.

Fix a function $g\in Y$ and let $\e>0$.
Since $K$ is compact, the range of $g$ is contained in some compact interval $[a,b]$. By Theorem \ref{thm:Weierstrass}
there exists a polynomial $p:[a,b]\to \R$ such that $\n p-q \n_\infty < \e$, where $q(t) := |t|$ is the absolute value function.
Since $Y$ is an algebra containing the constant functions, $p\circ g$ belongs to $Y$ and satisfies $\n p\circ g - |g| \n_{C(K)} < \eps$.
Since $\e>0$ was arbitrary and $Y$ is closed, it follows that $|g|\in \ov Y = Y$.
\end{proof}

\begin{remark}
Theorems \ref{thm:Stone-Weierstrass} and \ref{thm:Stone-Weierstrass-lattice} are stated for compact Hausdorff spaces, but the Hausdorff property was not used in the proofs. Note, however, that the separation-of-points assumptions in these theorems already imply the Hausdorff property.
\end{remark}

Theorem \ref{thm:Stone-Weierstrass} implies the following characterisation of separability.

\begin{proposition}\label{prop:CK-separable}
Let $K$ be a compact Hausdorff space. Then $C(K)$ is separable
if and only if $K$ is {\em metrisable},\index{metrisable} that is, $K$ admits a metric which generates the topology of $K$.
\end{proposition}

\begin{proof}
The assertion is immediate if $K$ is finite, so we assume that
$K$ is infinite.

Suppose first that $K$ is metrisable. Choose a compatible
metric $d$ and a countable dense set $(x_n)_{n\ge1}$ in $K$;
such a set exists by compactness. Put $d_n(x):=d(x,x_n)$.
These functions separate points: if $x\ne y$, choose $n$
such that $d(x,x_n)<\frac12d(x,y)$. Then
$$
d(y,x_n)\ge d(x,y)-d(x,x_n)>\frac12d(x,y)>d(x,x_n).
$$
The algebra generated by $\one$ and the functions $d_n$
contains the constants, separates points, and is closed under
complex conjugation. By Theorem \ref{thm:Stone-Weierstrass},
it is dense in $C(K)$.
Finite linear combinations of finite products of these generators
with coefficients in $\Q$ in the real case, or $\Q+i\Q$ in the
complex case, form a countable dense subset of $C(K)$.

Conversely, suppose that $C(K)$ is separable, and choose a
countable dense family $(f_n)_{n\ge1}$.
By Urysohn's lemma, continuous functions separate the points
of $K$, and density implies that the family $(f_n)_{n\ge1}$
also separates points. Define
$$
d(x,y):=\sum_{n=1}^\infty
2^{-n}\frac{|f_n(x)-f_n(y)|}{1+\n f_n\n_\infty},
\quad x,y\in K.
$$
This is a metric on $K$, and $x\mapsto d(x,y)$ is continuous in the original
topology for every $y\in K$. Hence the identity mapping
from $K$ to $(K,d)$ is continuous.
Being a continuous bijection from a compact space onto a
Hausdorff space, it is a homeomorphism by Proposition \ref{prop:compact-homeomorphism}. Thus $d$ induces the original topology of $K$.
\end{proof}

The next two examples give further illustrations of the Stone--Weierstrass theorem.

\begin{example}\label{ex:trigon-polCK-dense}
Over $\K=\C$, the {\em trigonometric polynomials},\index{trigonometric!polynomials} that is, linear combinations of the functions
$ e_n(\theta):= \exp(in\theta)$, $n\in\Z$,
are dense in $C(\mathbb{T})$, where $\mathbb{T}$ denotes the unit circle, which we think of as parametrised by $[-\pi,\pi]$. Indeed, they satisfy the requirements of Theorem \ref{thm:Stone-Weierstrass}.
An explicit approximation of functions in $C(\T)$ with trigonometric polynomials is described in Section \ref{subsec:Fourier-basis}.
\end{example}

\begin{example}\label{ex:tensorCK-dense}
Let $K_1,\dots,K_k$ be compact Hausdorff spaces. The linear combinations of functions of the form
$$ f(x) = f_1(x_1)\cdots f_k(x_k), \quad x = (x_1,\hdots,x_k)\in K_1\times\cdots\times K_k,$$
with $f_j\in C(K_j)$ for all $j=1,\dots,k$,
are dense in $C(K_1\times\cdots\times K_k)$. Indeed, they satisfy the requirements of Theorem \ref{thm:Stone-Weierstrass}.
\end{example}

\subsection{The Arzel\`a--Ascoli Compactness Theorem}
The next theorem gives a necessary and sufficient condition for
relative compactness in $C(K)$. We need the following terminology.
A subset $S\subseteq C(K)$ is said to be {\em equicontinuous}\index{equicontinuous}
at the point $x\in K$ if for all $\e>0$ there exists an open set $U$ in $K$ containing $x$ such that for all $x'\in U$ and $f\in S$ we have $|f(x)-f(x')|<\e$, and it is said to be {\em equicontinuous} if it is equicontinuous at every point of $K$.
The set $S$ is said to be {\em pointwise bounded}\index{pointwise!bounded} if for all $x\in K$ we have
$\sup_{f\in S}|f(x)|<\infty$.

\begin{theorem}[Arzel\`a--Ascoli]\label{thm:Arzela-Ascoli}\index{theorem!Arzel\`a--Ascoli}\index{compactness!in $C(K)$}
Let $K$ be a compact topological space.
For any subset $S$ of $C(K)$, the following assertions are equivalent:
\begin{enumerate}[label={\rm(\arabic*)}, leftmargin=*]
\item\label{it:Arzela-Ascoli1} $S$ is relatively compact;
\item\label{it:Arzela-Ascoli2} $S$ is bounded and equicontinuous;
\item\label{it:Arzela-Ascoli3} $S$ is pointwise bounded and equicontinuous.
\end{enumerate}
\end{theorem}

An equivalent way of formulating the theorem is that a
subset of $C(K)$ is compact if and only if it is closed, (pointwise) bounded, and equicontinuous.

\begin{proof}
\ref{it:Arzela-Ascoli1}$\Rightarrow$\ref{it:Arzela-Ascoli2}: \ Let $S\subseteq C(K)$ be relatively compact. Then $S$ is bounded, so all we need to do is to prove that $S$ is equicontinuous. Assuming without loss of generality that $K$ is nonempty, pick $x_0\in K$ and let $\eps>0$ be arbitrary and fixed. We can cover the compact set $\ov S$ with finitely many (say, $n$) open balls of radius $\eps$. Let $f_1,\dots,f_n$ be their centres. Using the continuity of these (finitely many) functions we can find an open set $U$ containing $x_0$ such that $|f_j(x)-f_j(x_0)| < \eps$ for all $x\in U$ and $j=1,\dots, n$. Now consider an arbitrary $f\in S$.
Choose $j_0\in \{1,\dots,n\}$ such that $\n f-f_{j_0}\n_\infty < \eps$; this is possible by the choice of the functions $f_1,\dots,f_n$. Then for all $x\in U$ we have
\begin{align*}
 |f(x)-f(x_0)| & \le |f(x)-f_{j_0}(x)| +|f_{j_0}(x)-f_{j_0}(x_0)|+ |f_{j_0}(x_0)-f(x_0)| <  3\eps.
\end{align*}
This verifies the equicontinuity condition.
\smallskip

\ref{it:Arzela-Ascoli2}$\Rightarrow$\ref{it:Arzela-Ascoli3}: \ This implication is trivial.
\smallskip

\ref{it:Arzela-Ascoli3}$\Rightarrow$\ref{it:Arzela-Ascoli1}: \ Let $S\subseteq C(K)$ be pointwise bounded and equicontinuous, and fix $\eps>0$.
By equicontinuity, for every $x\in K$ there is an open set $U_x$ in $K$ such that $|f(x)-f(x')|<\e$ for all $x'\in U_x$ and $f\in S$. By compactness, finitely many of these open sets cover $K$, say $U_{x_1},\dots,U_{x_k}$. By pointwise boundedness, for each $j=1,\dots, k$ the set $\{f(x_j):\, f\in S\}$ is bounded. It follows that we can find
$c_1,\dots,c_N \in \K$ such that for all $f\in S$ and $j=1,\dots,k$ we have $\min_{1\le n\le N} |f(x_j)-c_n| < \eps$.
Let $\mathscr{N} = \{n=(n_1,\dots,n_k):\, 1\le n_j\le N$ for all $j=1,\dots,k\}$. For $n\in \mathscr{N}$ let
$$ B_n = \{f\in S:\, |f(x_j)-c_{n_j}| < \eps \ \hbox{for all} \ j=1,\dots,k\}.$$
By what we just observed,
$$ S = \bigcup_{n\in\mathscr{N}}B_n.$$
Suppose that $f,g\in B_n$ and let $x\in K$ be arbitrary. Then $x$ belongs to at least one of the sets $U_{x_j}$.
Then,
\begin{align*} |f(x) - g(x)| & \le  |f(x) - f(x_j)| + |f(x_j) - g(x_j)| + |g(x_j) - g(x)|
\\ & \le \eps + |f(x_j) -c_{n_j}| +| c_{n_j}- g(x_j)| + \eps < 4\eps,
\end{align*}
where the last inequality holds uniformly with respect to $x\in K$. It follows that
$\n f-g\n_\infty<4\eps$.
If, for each $n\in \mathscr{N}$ for which $B_n$ is nonempty, we pick a function $f_n\in B_n$ and consider the open balls $B(f_n;4\eps)$, we obtain a finite cover of $S$ by $4\eps$-balls. Since $\eps>0$ was arbitrary, this means that $S$ is totally bounded and hence relatively compact, by Theorem \ref{thm:totally-bounded}.
\end{proof}

\subsection{Applications to Differential Equations}\label{sec:appl-DV}

As an interlude to the main development of the theory,
in this section we apply the completeness result of Theorem \ref{thm:CK} and the compactness result
of Theorem \ref{thm:Arzela-Ascoli} to study the following initial value problem:
\begin{equation*} \label{DV}\tag{IVP}
\left\{
\begin{aligned}
u'(t) &= f(t,u(t)), \quad t\in [0,T], \\
 u(0) &= u_0,
\end{aligned}
\right.
\end{equation*}
where $T>0$, the function $f:[0,T]\times \K^d\to \K^d$ is continuous, and $u_0\in\K^d$ is given.

\paragraph{Global Existence and Uniqueness}

A {\em global solution} is a continuously differentiable function $u:[0,T]\to \K^d$
satisfying $u(0)=u_0$ and $u'(t) = f(t,u(t))$ for all $t\in [0,T]$.

\begin{theorem}[Existence \& uniqueness, Picard--Lindel\"of]\label{thm:DV-Lip}\index{theorem!Picard--Lindel\"of}
\index{existence and uniqueness!under a Lipschitz assumption}
If $f:[0,T]\times \K^d\to \K^d$ is continuous and there is a constant $L\ge 0$ such that for all
$t\in [0,T]$ and $x,x'\in \K^d$ we have $$ |f(t,x) - f(t,x')| \le L|x-x'|,$$
then \eqref{DV} admits a unique global solution.
\end{theorem}

The condition on $f$ is often summarised by saying that $f$ is {\em Lipschitz continuous  in its second variable, uniformly with respect to its first variable.}\index{Lipschitz continuous}

The proof of Theorem \ref{thm:DV-Lip}
is based on the following abstract fixed point theorem.

\begin{theorem}[Banach fixed point theorem]\label{thm:fixed-point-banach}\index{theorem!Banach fixed point}
Let $X$ be a
nonempty
complete metric space and let
$f:X\to X$ be {\em uniformly contractive,}\index{uniformly!contractive} that is,
there exists a constant
$0\le c<1$ such that  $$d(f(x),f(x')) \le c\, d(x,x'), \quad  x,x'\in X.$$
Then $f$ has a unique {\em fixed point}, that is, there exists a unique element
$x\in X$ with the property that $f(x)=x$.
\end{theorem}
\begin{proof}
 If $x$ and $x'$ are both fixed points, then
$d(x,x') = d(f(x),f(x')) \le c\, d(x,x'),$
which is only possible if $d(x,x')=0$, that is, if
$x=x'$\!. It follows that a fixed point, if it exists, is unique.

To prove that a fixed point exists, choose an arbitrary $x_0\in X$ and define
the sequence $(x_n)_{n\ge 0}$ by
$x_{n+1} := f(x_n)$ for $n\ge 0.$
 We claim that this is a Cauchy sequence. Indeed, for all $n\ge 1$ we have
$$d(x_{n+1},x_n) = d(f(x_n), f(x_{n-1})) \le c\,d(x_n, x_{n-1}),$$
and therefore by induction one sees that
$d(x_{n+1},x_n) \le c^{n-1}\, d(x_2,x_1)$ for all $n\ge 1.$
For all  $m\ge n\ge N$ we have
$$\begin{aligned}
 d(x_{m},x_n)
 & \le d(x_m, x_{m-1}) + \cdots + d(x_{n+1},x_n)
 \\ &  \le (c^{m-2} + \cdots + c^{n-1}) \cdot d(x_2,x_1) \le \Bigl( \sum_{k=N-1}^\infty c^k\Bigr)\cdot  d(x_2,x_1)
   =\frac{c^{N-1}}{1-c} \cdot d(x_2,x_1),
 \end{aligned}
 $$
and the right-hand side can be made small by taking $N$ large.
This proves the claim.

Since $X$ was assumed to be complete, the sequence $(x_n)_{n\ge 0}$ converges in $X$. Let $x$ be its limit.
Then the continuity of $f$ implies
$f(x) = \lim_{n\to\infty} f(x_n) = \lim_{n\to\infty} x_{n+1} = x,$
which shows that $x$ is a fixed point for $f$.
\end{proof}

It is possible to replace uniform contractivity by a compactness assumption; this is
the {\em Schauder fixed point theorem} discussed in Problem \ref{prob:Schauder}.

By $C([0,T];\K^d)$\index{$C1$@$C([0,T];\K^d)$} we denote the space of all continuous functions $f:[0,T]\to \K^d$\!. Endowed with the supremum norm, this space is a Banach space. Indeed, suppose that $(f^{(n)})_{n\ge 1}$ is a Cauchy sequence in $C([0,T];\K^d)$. Then the $d$ sequences of coordinate functions $(f_j^{(n)})_{n\ge 1}$ are Cauchy in $C[0,T]$ and therefore converge to limits $f_j$ in $C[0,T]$. This easily implies that the sequence $(f^{(n)})_{n\ge 1}$
converges in $C([0,T];\K^d)$ to the function $f$ with coordinate functions $f_j$.

We will use the Banach fixed point theorem to prove that the (nonlinear) mapping $I_T:C([0,T];\K^d)\to C([0,T];\K^d)$
defined by
$$ (I_T u)(t):= u_0 + \int^t_0 f(s,u(s))\ud s, \quad t\in [0,T],$$
where the integral is interpreted as a $\K^d$-valued Riemann integral, has a fixed point.
This will prove the theorem in view of the next lemma.

\begin{lemma}\label{lem:ODE-fixedpoint}\index{fixed point!argument}
A function $u\in C([0,T];\K^d)$ satisfies \eqref{DV} for all $t\in [0,T]$
if and only if $u$ is a fixed point of $I_T$.
\end{lemma}
\begin{proof}
Indeed, $u$ is a fixed point of $I_T$ if and only if
$u(t)  = u_0 + \int^t_0 f(s,u(s))\ud s$ for all $t\in [0,T]$. By integration, this identity holds if $u$ is a solution, and conversely, if the identity holds, then $u$ is continuously differentiable (since the right-hand side is) and by differentiation we find that $u$ is a solution.
\end{proof}

\begin{proof}[Proof of Theorem \ref{thm:DV-Lip}]
Let us start with a preliminary estimate that will be refined shortly.
For all $u,v\in C([0,T];\K^d)$ and all $t\in [0,T]$ we have
\begin{align*}
|(I_T(u))(t) - (I_T(v))(t)| & = \Bigl|\int^t_0 f(s,u(s)) - f(s,v(s))\ud s\Bigr|
\\ & \le \int^t_0 L |u(s) - v(s)|\ud s
 \le \int^t_0 L \n u-v\n \ud s
 \le LT\n u-v\n.
\end{align*}
Taking the supremum over $t\in [0,T]$, we find that
$$ \n I_T(u) - I_T(v)\n \le LT\n u-v\n.$$
If $LT <1$, then $I_T$ is uniformly contractive and the Banach fixed point theorem
guarantees the existence of a unique fixed point. This proves Theorem \ref{thm:DV-Lip}
in the special case in which the smallness condition $LT< 1$ is satisfied.

To circumvent this condition
we modify the norm of $C([0,T];\K^d)$. Fix a real number $\lambda>0$ (in a moment
we will see that we need $\lambda>L$), and
define $$\n g\n_\lambda := \sup_{t\in [0,T]} e^{-\lambda t} |g(t)|.$$
It is clear that this defines a norm on $C([0,T];\K^d)$ and we have
$$ e^{-\lambda T}\n g\n \le \n g\n_\lambda \le \n g\n.
$$
This implies that a sequence in $C([0,T];\K^d)$ is Cauchy with respect to the norm $\n \cdot\n_\lambda$
if and only if it is Cauchy with respect to the norm $\n \cdot\n$, and since $C([0,T];\K^d)$ is complete with respect to the latter, we conclude that $C([0,T];\K^d)$ is a Banach space with respect to the norm $\n \cdot\n_\lambda$.
Using this norm, we redo the above computations and find
\begin{align*}
\n I_T(u) - I_T(v)\n_\lambda
& = \sup_{t\in [0,T]} e^{-\lambda t} |(I_T(u))(t) - (I_T(v))(t)|
\\ & \le \sup_{t\in [0,T]} e^{-\lambda t} \int^t_0 L e^{\lambda s}e^{-\lambda s}|u(s) - v(s)|\ud s
\\ & \le \sup_{t\in [0,T]} e^{-\lambda t} \n u-v\n_\lambda\int^t_0 L e^{\lambda s}\ud s
\\ & = \sup_{t\in [0,T]} e^{-\lambda t} \n u-v\n_\lambda \cdot \frac{L}\lambda(e^{\lambda t}-1)
 \le \frac{L}\lambda \n u-v\n_\lambda.
\end{align*}
Hence if we choose $\lambda > L$, then $I_T$ is uniformly contractive on $C([0,T];\K^d)$
with respect to the norm $\n \cdot\n_\lambda$. Now an application of the Banach fixed point
theorem produces a unique fixed point for $I_T$.
\end{proof}

\begin{remark}
Theorem \ref{thm:DV-Lip} remains true if we replace the interval $[0,T]$ by $[0,\infty)$
and assume that $f: [0,\infty)\times \K^d\to\K^d$ is continuous and satisfies  $$ |f(t,x) - f(t,x')| \le L|x-x'|$$ for all
$t\in [0,\infty)$ and $x,x'\in \K^d\!$.
Indeed, the preceding argument produces a solution $u_T$ on every interval $[0,T]$. We may now define $u:[0,\infty)\to \K^d$
by setting $u:= u_T$ on the interval $[0,T]$.
Since by uniqueness we have $u_T = u_S$ on $[0,S\wedge T]$, this is well defined. The resulting function is continuously differentiable and satisfies
\eqref{DV} on every interval $[0,T]$, hence on all of $[0,\infty)$, and is therefore a global solution on $[0,\infty)$.
\end{remark}

\begin{remark}
 All that has been said extends to the case where $\K^d$ is replaced by a Banach space $X$.
\end{remark}

\paragraph{Local Existence}

As an application of Theorem \ref {thm:Arzela-Ascoli} we present a local existence result for differential equations with continuous
right-hand side. In contrast to the situation in Theorem \ref{thm:DV-Lip}, where a Lipschitz continuity
assumption was made, we do not get uniqueness of the solution.

The problem \eqref{DV} {\em admits a local solution}
if there exists $0<\delta\le T$ and a continuously differentiable function $u:[0,\delta]\to \K^d$
satisfying $u(0)=u_0$ and $u'(t) = f(t,u(t))$ for all $t\in [0,\delta]$.

\begin{theorem}[Local existence, Peano]\label{thm:DV-cont}\index{existence!under a continuity assumption}\index{theorem!Peano}
If $f:[0,T]\times \K^d\to \K^d$ is continuous, then the problem \eqref{DV} admits a local solution.
\end{theorem}

 A global solution need not always exist. Consider the problem
\begin{equation*}
\left\{
\begin{aligned}
u'(t) &= (u(t))^2,\quad t\in[0,1],\\
u(0) &= 1.
\end{aligned}
\right.
\end{equation*}
The function $u(t)=1/(1-t)$ is a solution on every interval
$[0,\delta]$ with $0<\delta<1$.
To see why no solution exists on all of $[0,1]$, let $u$ be
any solution on $[0,\delta]$, where $0<\delta\le1$, and put
$$
v(t):=\exp\Bigl(-\int_0^t u(s)\ud s\Bigr),
\quad t\in[0,\delta].
$$
Differentiation gives $v'=-uv$, and therefore
$$
(uv)'=u'v+uv'=u^2v-u^2v=0.
$$
Since $u(0)=v(0)=1$, it follows that $uv=1$ on $[0,\delta]$.
Consequently, $v'=-uv=-1$, and hence $v(t)=1-t$.
But the exponential defining $v$ never vanishes, so
$\delta\ge 1$ is impossible.
For $\delta<1$, the identity $uv=1$ also gives
$u(t)=1/(1-t)$, proving uniqueness on $[0,\delta]$.

The right-hand side $z\mapsto z^2$ in this example is locally
Lipschitz continuous. For a merely continuous right-hand side,
however, local solutions need not be unique. Indeed,
both $u_1(t) =0$ and $u_2(t) = t^{3/2}$
are solutions of the problem
\begin{equation*}
\left\{
\begin{aligned}
u'(t) &= \frac32 |u(t)|^{1/3}\!, \quad t\in [0,1], \\
 u(0) &= 0.
\end{aligned}
\right.
\end{equation*}
The function $\phi(t,x):= \frac32 |x|^{1/3}$ fails to be
locally Lipschitz continuous in the variable $x$ in a neighbourhood of the initial value $0$.

As a first step towards the proof of Theorem \ref{thm:DV-cont} we show
that the problem \eqref{DV} is equivalent to an
integrated version of it.
For the remainder of this section we assume that $f$ is continuous.

\begin{lemma}\label{lem:opl} A function $u\in C([0,\delta];\K^d)$ is a local solution of \eqref{DV}
if and only if for all $t\in [0,\delta]$ we have
\begin{equation} \label{IV}
u(t) = u_0 + \int^t_0 f(s,u(s))\ud s.
\end{equation}
\end{lemma}

The function $s\mapsto f(s,u(s))$
is continuous on $[0,\delta]$, so the integral is well defined as a Riemann integral with values in $\K^d\!$.

\begin{proof}
If $u\in C([0,\delta];\K^d)$ is a local solution, then for all $t\in [0,\delta]$ we have
$$ u(t) - u_0 = u(t) - u(0) = \int^t_0 u'(s)\ud s = \int^t_0 f(s,u(s))\ud s$$
and \eqref{IV} holds for all $t\in [0,\delta]$.
Conversely, if $u\in C([0,\delta];\K^d)$ satisfies \eqref{IV} for all
$t\in [0,\delta]$, then $u$ is continuously differentiable on $[0,\delta]$. Using $u(0) = u_0$ and differentiating, we find
$$ u'(t) = \frac{{\rm d}}{{\rm d}t} \int^t_0 f(s,u(s))\ud s = f(t,u(t))$$
for all $t\in [0,\delta]$, that is, \eqref{DV} holds.
\end{proof}

Now we are ready for the proof of Theorem \ref{thm:DV-cont}. It relies on a compactness argument. The idea is to construct, for small enough
$\delta\in (0,T]$, a
sequence of approximate solutions $(u_n)_{n\ge 1}$ in $C([0,\delta];\K^d)$ and show its equicontinuity (the definition of which
extends to vector-valued functions in the obvious way). An appeal to the Arzel\`a--Ascoli theorem (which extends to the vector-valued
case as well, without change in the proof) then produces a subsequence $(u_{n_k})_{k\ge 1}$ that converges in $C([0,\delta];\K^d)$.
The limit $u\in C([0,\delta];\K^d)$ will be shown to solve \eqref{DV} on the interval $[0,\delta]$.

\begin{proof}[Proof of Theorem \ref{thm:DV-cont}]
By continuity of $f$ and compactness, the quantity
$$
M:=\sup_{(t,x)\in[0,T]\times\ov B(u_0;1)}|f(t,x)|
$$
is finite. Put $\delta:=\min\{T,1/M\}$, taking $\delta:=T$
if $M=0$.

For $n\ge1$, divide $[0,\delta]$ into $2^n$ equal intervals,
with mesh size $h_n:=2^{-n}\delta$ and partition points
$t_{j,n}:=jh_n$.
Set $u_n(0):=u_0$ and define the remaining values at the
partition points by
\begin{align*}
u_n(t_{j+1,n})
=u_n(t_{j,n})+h_n f(t_{j,n},u_n(t_{j,n})),
\quad j=0,\dots,2^n-1.
\end{align*}
Between successive partition points, define $u_n$ by linear
interpolation.

An induction shows that
$$
|u_n(t_{j,n})-u_0|\le Mt_{j,n}\le1,
\quad j=0,\dots,2^n.
$$
Indeed, the assertion holds for $j=0$. If it holds for $j$,
then $u_n(t_{j,n})\in\ov B(u_0;1)$, so
$|f(t_{j,n},u_n(t_{j,n}))|\le M$.
The recursion therefore gives
$$
|u_n(t_{j+1,n})-u_0|
\le Mt_{j,n}+Mh_n=Mt_{j+1,n}.
$$
Since the closed ball is convex, the interpolating line segments
also lie in $\ov B(u_0;1)$.
Their slopes have modulus at most $M$, and consequently
$$
|u_n(t)-u_n(s)|\le M|t-s|,
\quad s,t\in[0,\delta].
$$
Thus the sequence $(u_n)_{n\ge1}$ is uniformly bounded and
equicontinuous. By the Arzel\`a--Ascoli theorem, some subsequence
$(u_{n_k})_{k\ge1}$ converges uniformly to a function
$u\in C([0,\delta];\K^d)$.
The limit also takes its values in $\ov B(u_0;1)$.

To identify the limit, let $\kappa_n(s):=t_{j,n}$ when
$t_{j,n}\le s<t_{j+1,n}$, and put $\kappa_n(\delta):=\delta$.
By integrating the constant slope on each partition interval, we obtain
$$
u_n(t)=u_0+\int_0^t
f\bigl(\kappa_n(s),u_n(\kappa_n(s))\bigr)\ud s,
\quad t\in[0,\delta].
$$
We have $|\kappa_n(s)-s|\le h_n$, and the Lipschitz estimate
for $u_n$ gives
$$
|u_n(\kappa_n(s))-u_n(s)|\le Mh_n,
\quad s\in[0,\delta].
$$
Hence
$$
\sup_{s\in[0,\delta]}
|u_{n_k}(\kappa_{n_k}(s))-u(s)|
\le Mh_{n_k}+\n u_{n_k}-u\n_\infty
\to 0.
$$
Thus both arguments of $f$ in the integrand converge uniformly
to $s$ and $u(s)$, respectively.
Since $f$ is uniformly continuous on
$[0,T]\times\ov B(u_0;1)$, the integrands converge uniformly
to $s\mapsto f(s,u(s))$.
Passing to the limit in the integral identity, we obtain
$$
u(t)=u_0+\int_0^t f(s,u(s))\ud s,
\quad t\in[0,\delta].
$$
By Lemma \ref{lem:opl}, the function $u$ is a local solution
of \eqref{DV}.
\end{proof}

An alternative proof, based on the Schauder fixed point theorem, is discussed in Problem \ref{prob:local-existence}.

\section{Spaces of Integrable Functions}\label{sec:Lp}

Let $(\Om,\calF\!,\mu)$ be a measure space and
fix $1\le p<\infty$. We define $\calL^p(\Om)$ as the set of
all measurable functions $f:\Om\to \K$ such that
$$ \int_\Om |f|^p\ud \mu< \infty.$$
For such functions we set
$$\n f\n_p := \Big(\int_\Om |f|^p\ud \mu \Big)^{1/p}\!.$$
For $p=\infty$ we define $\calL^\infty(\Om)$ as the set of
all measurable functions $f:\Om\to \K$ that are {\em $\mu$-essentially bounded},\index{essentially!bounded} meaning that
there exists a set $N\in\calF$ of $\mu$-measure $0$ such that $f$ is bounded on $\complement N$. For such functions we define
$\n f\n_\infty$ as the {\em $\mu$-essential supremum} of $|f|$,
$$\n f\n_\infty :=  \hbox{$\mu$-}\esssup_{\om\in\Om} |f(\om)|  \,:=\, \inf\bigl\{r>0: \, |f|\le r \ \hbox{$\mu$-almost everywhere}\bigr\}.$$
When there is no risk of confusion, the measure $\mu$ is omitted from this notation.

The spaces $\calL^p(\Om)$ are vector spaces:

\begin{proposition}[Minkowski inequality]\label{prop:Minkowski}\index{inequality!Minkowski}
Let $1\le p\le \infty$. For all functions $f,g\in \calL^p(\Om)$ we have $f+g\in \calL^p(\Om)$ and
$$\n f+g\n_p \le \n f\n_p + \n g\n_p.$$
\end{proposition}

\begin{proof}
The result is trivial for $p=\infty$, so we only consider the case $1\le p<\infty$.

An elementary calculus argument shows that for all nonnegative real numbers $a$ and $b$ one has
$$ (a+b)^p = \inf_{t\in (0,1)} t^{1-p} a^p + (1-t)^{1-p} b^p\!.$$
Applying this identity to $|f(\om)|$ and $|g(\om)|$ with $\om\in\Om$ and integrating with respect to
$\mu$,
for each fixed $t\in(0,1)$ we obtain
$$
\int_\Om |f+g|^p \ud \mu
\le  \int_\Om (|f|+|g|)^p \ud \mu
 \le t^{1-p}  \int_\Om |f|^p\ud \mu + (1-t)^{1-p}\int_\Om  |g|^p\ud \mu.
$$
Stated differently, this says that
$$ \n f+g\n_p^p \le t^{1-p} \n f\n_p^p + (1-t)^{1-p}\n g\n_p^p.$$
Taking the infimum over all $t\in (0,1)$ gives the result.
\end{proof}

In spite of this result, $\n \cdot\n_p$ is not a norm on $\calL^p(\Om)$, because
$\n f\n_p = 0$ only implies that $f=0$ $\mu$-almost everywhere. To circumvent this imperfection, we define an equivalence relation $\sim$ on $\calL^p(\Om)$ by
$$ f\sim g \Leftrightarrow f = g \hbox{ $\mu$-almost everywhere}.$$
The equivalence
class of a function $f$ modulo $\sim$ is denoted by $[f]$.
On the quotient space\index{$L^p(\Om)$}
$$ L^p(\Om):= \calL^p(\Om)/\sim,$$
whose elements are the equivalence classes $[f]$ of functions $f\in \calL^p(\Om)$,
we define a scalar multiplication and addition in the natural way:
$$ c[f] := [cf], \qquad [f] + [g] := [f+g].$$
We leave it as an exercise to check that both operations are well defined.
With these operations, $L^p(\Om)$
is a normed vector space with respect to the norm
$$ \n [f] \n_p := \n f\n_p.$$
When we explicitly wish to express the dependence on $\calF$ or $\mu$ we write
$L^p(\Om,\calF)$ or $L^p(\Om,\mu)$.
{\em Following common practice we make no
distinction between functions in $\calL^p(\Om)$ and their equivalence classes in $L^p(\Om)$}, and call the latter `functions' as well. In the same vein, we will not hesitate to talk
about the `sets' $$\{\om\in\Om:\, f(\om)\in B\}$$
when $B\subseteq\K$ is a Borel set and $f$ is an element of $L^p(\Om)$. The rigorous interpretation is that
$\{\om\in\Om:\, f(\om)\in B\}$
defines an equivalence class
of sets in $\calF$, the representatives of which are obtained
by selecting pointwise defined measurable representatives for $f$.

\subsection{Completeness}\label{subsec:Lp-complete}

For the remainder of this section we fix a measure space $(\Om,\calF\!,\mu)$.
The main result of this section is the following completeness result for the spaces $L^p(\Om)$.

\begin{theorem}[Completeness]\label{thm:Lp-complete}
For all $1\le p\le \infty$ the normed space $L^p(\Om)$ is complete.
\end{theorem}
\begin{proof}
First let $1\le p<\infty$, and let $(f_n)_{n\ge 1}$ be a Cauchy sequence with respect to the norm $\n
\cdot\n_p$ of $L^p(\Om)$. Passing to a subsequence, we may assume that
$$ \n f_{n+1} - f_{n}\n_p \le \frac1{2^n}, \quad n\ge 1.$$
Define the nonnegative measurable functions
$$ g_N := \sum_{n=0}^{N-1} |f_{n+1} - f_{n}|, \qquad
g := \sum_{n=0}^\infty |f_{n+1} -
f_{n}|,$$ with the convention that $f_0 = 0$.
By the monotone convergence theorem,
$$ \int_\Om g^p \ud \mu = \lim_{N\to\infty} \int_\Om g_N^p\ud \mu.$$
Taking $p$th roots and using Minkowski's inequality we obtain
$$ \n g\n_p = \lim_{N\to\infty} \n g_N\n_p \le \lim_{N\to\infty}
\sum_{n=0}^{N-1} \n
f_{n+1} - f_{n}\n_p  =  \sum_{n=0}^\infty \n
f_{n+1} - f_{n}\n_p  \le 1 + \n f_1\n_p.$$
It follows that $g$ is finitely valued $\mu$-almost everywhere, which means that
the sum defining $g$ converges absolutely  $\mu$-almost everywhere.
As a result,  the sum
$$ f:=  \sum_{n=0}^{\infty} (f_{n+1} - f_{n})$$
converges on the set $\{g<\infty\}$.
On this set we have
$$ f=   \lim_{N\to\infty} \sum_{n=0}^{N-1} (f_{n+1} - f_{n})
= \lim_{N\to\infty} f_N.$$
Defining $f$ to be identically zero on the null set
$\{g=\infty\}$, the resulting function $f$
is measurable.

Since $|f_N|\le g$ and $f_N\to f$ almost everywhere, we have
$|f|\le g$ almost everywhere. As $g\in L^p(\Om)$, this proves
that $f\in L^p(\Om)$.
Moreover, the series defining $f$ gives, $\mu$-almost everywhere,
$$
|f-f_N|
\le\sum_{n=N}^\infty|f_{n+1}-f_n|
=g-g_N\le g.
$$
The left-hand side tends to zero almost everywhere, and its
$p$th power is bounded by the integrable function $g^p$.
The dominated convergence theorem therefore gives
$$
\n f-f_N\n_p^p
=\int_\Om|f-f_N|^p\ud\mu
\to 0.
$$

We have proved that a subsequence of the original Cauchy sequence converges to
$f$ in $L^p(\Om)$. As is easily verified, this implies that the original Cauchy
sequence converges to $f$ as well. This completes the proof for exponents $1\le p<\infty$.

It remains to establish the result for $p=\infty$. Let $(f_n)_{n\ge 1}$ be a Cauchy sequence with respect to the norm $\n\cdot\n_\infty$ of $L^\infty(\Om)$. Passing to a subsequence, we may assume that
$$ \n f_{n+1} - f_{n}\n_\infty \le \frac1{2^n}, \quad n\ge 1.$$
Choose $\mu$-null sets $F_n$, $n\in\N$, such that $f_1$ is bounded on $\complement F_0$ and $|f_{n+1}(\om) - f_{n}(\om)| \le \frac1{2^n}$ for all $n\ge 1$ and $\om\in\complement F_n$.
Outside the null set $F:=F_0\cup\bigcup_{n\ge1}F_n$, the tails satisfy
$$
\sum_{n=N}^{M-1}|f_{n+1}(\omega)-f_n(\omega)|
\le \sum_{n=N}^{M-1} \frac1{2^n},
\quad \omega\in \Omega\setminus F.
$$
Hence the telescoping sums $f_N$ converge uniformly on
$\Omega\setminus F$. Taking $f$ to be their limit on $\Omega\setminus F$, and setting $f$
equal to $0$ on $F$, we obtain $f\in L^\infty(\Omega)$ and
$f_N\to f$ in $L^\infty(\Omega)$.
\end{proof}

In the course of the proof we obtained the following result:

\begin{corollary} \label{cor:Lp-ae-subseq}
Every convergent sequence $(f_n)_{n\ge 1}$ in $L^p(\Om)$, with $1\le p\le \infty$, has a $\mu$-almost everywhere
convergent subsequence $(f_{n_k})_{k\ge 1}$, and this subsequence may be chosen to satisfy
$|f_{n_k}|\le g$ almost everywhere for some fixed $0\le g\in L^p(\Om)$.
\end{corollary}

In the majority of applications the first part of this corollary suffices, but the second part is sometimes helpful in setting the stage for an application of the dominated convergence theorem.

\begin{remark}
 Except when $p=\infty$, in the setting of Corollary \ref{cor:Lp-ae-subseq} it need not be the case that the sequence $(f_n)_{n\ge 1}$ itself is $\mu$-almost everywhere convergent to its $L^p(\Om)$-limit $f$ (see Problems \ref{prob:Lp-subseq1} and \ref{prob:Lp-subseq2}).
\end{remark}

Recall from Definition \ref{def:mu-simple} that a {\em $\mu$-simple function} is a simple function supported on sets of finite $\mu$-measure.

\begin{proposition}[Approximation by $\mu$-simple functions]
\label{prop:approx-simple}
Let $(\Om,\calF\!,\mu)$ be a measure space. For every
$1\le p<\infty$, the $\mu$-simple functions are dense in
$L^p(\Om)$. The same result holds for $L^\infty(\Om)$ if $\mu$ is
finite.
\end{proposition}

\begin{proof}
Let $1\le p<\infty$ and fix $f\in L^p(\Om)$.
By Proposition \ref{prop:approx}, there exist simple functions
$s_n$ such that $s_n\to f$ pointwise and $|s_n|\le |f|$.
Write $s_n = \sum_{j=1}^{N_n} c_{j,n} \one_{B_{j,n}}$, where for each $n\ge 1$ the sets $B_{j,n}$ are disjoint and the coefficients $c_{j,n}$ are nonzero. From
$$|c_{j,n}|^p \mu(B_{j,n}) \le \int_\Om |s_n|^p\ud \mu \le \int_\Om |f|^p\ud \mu < \infty$$
we then see that each $s_n$ is $\mu$-simple.
Since $|s_n-f|^p\to0$ pointwise and
$|s_n-f|^p\le 2^p|f|^p$, dominated convergence therefore gives
$ \n s_n-f\n_p^p\to 0$.

If $\mu(\Om)<\infty$ and $f\in L^\infty(\Om)$, we choose a bounded measurable representative of $f$.
By Proposition \ref{prop:approx}, this representative can be
approximated uniformly by simple functions.
Since $\mu$ is finite, every simple function is $\mu$-simple,
and uniform convergence implies convergence in $L^\infty(\Om)$.
\end{proof}

The next inequality is known as {\em H\"older's inequality}.\index{inequality!H\"older}
For $p=q=2$ and $r=1$ it reduces to a special case of the Cauchy--Schwarz inequality (see Proposition \ref{prop:CS}).

\begin{proposition}[H\"older's inequality]\label{prop:Holder} Let $1\le p,q,r\le \infty$ satisfy $\frac1p+ \frac1q=\frac1r$.
If $f\in L^p(\Om)$ and $g\in L^q(\Om)$, then $fg\in L^r(\Om)$ and
$$ \n fg\n_r \le \n f\n_p \n g\n_q.$$
\end{proposition}

For $r=1$ the condition on $p$ and $q$ reads $\frac1p+\frac1q=1$; we call such $p$ and $q$ {\em conjugate exponents}.\index{conjugate exponents}

\begin{proof}
If $r=\infty$, then $p=q=\infty$ and the result is trivial. To prove the result for
$1\le r<\infty$ it suffices to consider the case $r=1$; the general case follows by applying this special case to the functions $|f|^r$ and $|g|^r$.

For $p=1$, $q=\infty$ and for $p=\infty$, $q=1$, the inequality follows by a direct
estimate. Thus we may assume from now on that $1< p,q<\infty$.
The inequality is then proved in the same way as Minkowski's inequality,
this time using the identity
\[ ab = \inf_{t>0} \Big(\frac{t^p a^p}{p} + \frac{b^q}{qt^q}\Big).\qedhere\]
\end{proof}

\begin{remark}\label{rem:Holder}
Let $1\le p_1,\dots,p_N,r\le \infty$ satisfy $\frac1{p_1}+\cdots+ \frac1{p_N}=\frac1r$.
If $f_n\in L^{p_n}(\Om)$ for $n=1,\dots,N$, then $\prod_{n=1}^N f_n \in L^r(\Om)$ and
$$ \Big\n \prod_{n=1}^N f_n\Big\n_r \le \prod_{n=1}^N \n f_n\n_{p_n}.$$
This more general version of H\"older's inequality follows from Proposition \ref{prop:Holder} by an easy induction argument.
\end{remark}

As an immediate corollary of H\"older's inequality we have the following result.

\begin{corollary}\label{cor:Holder} Let $1\le p,q,r\le \infty$ satisfy $\frac1p+ \frac1q=\frac1r$.
Then the mapping
$$ (f,g) \mapsto fg$$
is jointly continuous from $L^p(\Om)\times L^q(\Om)$ into $L^r(\Om)$.
In particular, if $1\le p,q\le \infty$ satisfy $\frac1p+ \frac1q=1$ and $f_n\to f$ in $L^p(\Om)$, then for all $g\in L^q(\Om)$ we have
$$ \limn \int_\Om f_n g\ud\mu = \int_\Om fg\ud \mu.$$
\end{corollary}
\begin{proof}
If $f_n\to f$ in $L^p(\Om)$ and $g_n\to g$ in $L^q(\Om)$,  then
H\"older's inequality implies that $f_ng_n$, $f_ng$, $fg_n$, $fg$ belong to $L^r(\Om)$ and
\begin{align*}
 \n f_n g_n - fg\n_r &\le \n (f_n  - f)g\n_r+ \n f_n (g_n - g)\n_r
 \\ & \le \n f_n  - f\n_p \n g\n_q+ M\n g_n - g\n_q \to 0 \quad\hbox{as $n\to\infty$},
\end{align*}
where $M:= \sup_{n\ge 1}\n f_n\n_p$.
By Proposition \ref{prop:sequentially-cont}, this proves the asserted continuity.
\end{proof}

A useful special case of H\"older's inequality concerns finite measures. If $\mu(\Omega)<\infty$ and $1\le r\le  p\le \infty$, then H\"older's inequality implies that if $f\in L^p(\Om)$, then $f\in L^r(\Om)$ and
$$ \n f\n_r  \le \mu(\Om)^{\frac1r - \frac1p} \n f\n_p.$$
In the case of a probability measure $\mu$ this takes the simpler form
$ \n f\n_r \le \n f\n_p$.

The following result provides a converse to H\"older's inequality.  We formulate it for exponents $\frac1p+\frac1q =1$;
as in the proof of H\"older's inequality, this implies a more general version for
exponents $\frac1p+\frac1q =\frac1r$. A further variation will be given in Proposition \ref{prop:Lp-via-Lq-2}.

 \begin{proposition}\label{prop:Lp-via-Lq-1}\index{inequality!H\"older, converse to}
Let $1\le p, q\le \infty$ satisfy $\frac1p+\frac1q =1$. Let $(\Omega,\mu)$ be a
$\sigma$-finite measure space.
A measurable function $f$ belongs to $L^p(\Om)$ if and only if $$fg\in L^1(\Om) \ \hbox{ and } \ \n fg\n_1 \le M \n g\n_q$$ for some constant $M\ge 0$ and all
$g\in Y$, where $Y$ is a dense subspace of $L^q(\Om)$.
If $p=1$, the result also holds if $Y$ is taken to be the space of bounded simple functions on $\Om$ supported on sets of finite measure.
In each of these cases we have $\n f\n_p \le M.$
 \end{proposition}

\begin{proof}
The `only if' part follows from H\"older's inequality.
For the converse, choose measurable sets $E_n\uparrow\Om$
with $\mu(E_n)<\infty$; this is possible by $\sigma$-finiteness.

\smallskip
{\em Step 1} --
First suppose that $Y$ is dense in $L^q(\Om)$.
We show that the assumed estimate holds for every
$g\in L^q(\Om)$.
Choose $g_n\in Y$ such that $g_n\to g$ in $L^q(\Om)$.
By Corollary \ref{cor:Lp-ae-subseq}, after passing to a
subsequence we may also assume that $g_n\to g$
$\mu$-almost everywhere.
Fatou's lemma then gives
$$
\int_\Om |fg|\ud\mu
\le\liminf_{n\to\infty}\int_\Om |fg_n|\ud\mu
\le M\lim_{n\to\infty}\n g_n\n_q
=M\n g\n_q.
$$
Thus $fg\in L^1(\Om)$ and
\begin{align}\label{eq:fg-L1}
\n fg\n_1\le M\n g\n_q,
\quad g\in L^q(\Om).
\end{align}

{\em Step 2} --
We next prove the result for $1\le p<\infty$.
If $p=1$, apply \eqref{eq:fg-L1} with $g=\one_{E_n}$ to obtain
$$
\int_{E_n}|f|\ud\mu
\le M\n\one_{E_n}\n_\infty\le M.
$$
Monotone convergence gives
$\int_\Om|f|\ud\mu\le M$.
This argument also proves the additional assertion for $p=1$:
each $\one_{E_n}$ is a bounded simple function supported on a
set of finite measure, so in that case the displayed estimate
follows directly from the assumption.

Now let $1<p<\infty$. Each function $u_n:=(|f|\wedge n)\one_{E_n}$
is bounded and supported on a set of finite
measure. Hence $u_n^{p-1}\in L^q(\Om)$, and
$(p-1)q=p$ gives $\n u_n^{p-1}\n_q=\n u_n\n_p^{p-1}$.
Since $0\le u_n\le |f|$, applying \eqref{eq:fg-L1} with
$g=u_n^{p-1}$ gives
$$
\n u_n\n_p^p
=\int_\Om u_n^p\ud\mu
\le\int_\Om |f|u_n^{p-1}\ud\mu
\le M\n u_n\n_p^{p-1}.
$$
If $\n u_n\n_p>0$, this gives
$\n u_n\n_p\le M$; the same inequality holds trivially if
$\n u_n\n_p=0$.
Since $u_n\uparrow |f|$, by the monotone convergence theorem we obtain
\begin{align*}
\int_\Om |f|^p\ud\mu
=\lim_{n\to\infty}\int_\Om u_n^p\ud\mu
\le M^p.
\end{align*}
This proves that $f\in L^p(\Om)$ and $\n f\n_p\le M$.

\smallskip
{\em Step 3} --
Finally, let $p=\infty$, so that $q=1$.
Suppose, for a contradiction, that $|f|\le M$ does not hold
$\mu$-almost everywhere. Then for some $\eps>0$ the set
$$
A:=\{x\in\Om:\, |f(x)|>M+\eps\}
$$
has positive measure.
Since $E_n\uparrow\Om$, for some $n$ the set $E:=A\cap E_n$
has positive and finite measure.
Applying \eqref{eq:fg-L1} with $g=\one_E$, we obtain
$$
(M+\eps)\mu(E)
\le\int_E|f|\ud\mu
\le M\n\one_E\n_1
=M\mu(E),
$$
which is impossible.
Thus $f\in L^\infty(\Om)$ and $\n f\n_\infty\le M$.
\end{proof}

\begin{remark}\label{rem:Lp-via-Lq-1}
The argument of Step 3 proves, more generally, that if $(\Omega,\calF\!,\mu)$ is a $\sigma$-finite measure space and $1\le p\le \infty$, then a measurable function $f$ belongs to $L^\infty(\Om)$ if and only if $fg\in L^p(\Om)$ and $\n fg\n_p \le M \n g\n_p$ for some constant $M\ge 0$ and all $g\in Y$, where $Y$ is a dense subspace of $L^p(\Om)$; in that case we have $\n f\n_\infty \le M.$
\end{remark}

The final result of this subsection characterises separability.

\begin{proposition}[Separability of $L^p$-spaces]
\label{prop:Lp-separable}
Let $(\Om,\F,\mu)$ be a $\sigma$-finite measure space and let
$1\le p<\infty$. The following assertions are equivalent:
\begin{enumerate}[label={\rm(\arabic*)}, leftmargin=*]
\item\label{it:Lp-separable-sigma1} $L^p(\Om,\mu)$ is separable;
\item\label{it:Lp-separable-sigma2} there exists a countably generated sub-$\sigma$-algebra
$\mathscr G\subseteq\F$ such that for every $F\in\F$ there is
a $G\in\mathscr G$ with $\mu(F\triangle G)=0$.
\end{enumerate}
In particular, separability does not depend on the choice of
$p\in[1,\infty)$.
\end{proposition}

\begin{proof}
Choose measurable sets $E_n\uparrow\Om$ with
$\mu(E_n)<\infty$.

\smallskip
\ref{it:Lp-separable-sigma1}$\Rightarrow$\ref{it:Lp-separable-sigma2}: \
Let $(f_j)_{j\ge1}$ be a dense sequence in $L^p(\Om,\mu)$,
choose a measurable representative for each $f_j$, and put
$\mathscr G:=\sigma(f_j:\,j\ge 1)$, the $\sigma$-algebra generated by these representatives.
Since the Borel $\sigma$-algebra of $\K$ is countably generated,
so is $\mathscr G$.
Every function in $L^p(\Om,\mu)$ is an $L^p$-limit of a
sequence of these representatives. Passing to an almost
everywhere convergent subsequence shows that it has a
$\mathscr G$-measurable representative.
Applying this to $\one_{F\cap E_n}$, with $F\in\F$, we obtain
sets $G_n\in\mathscr G$ such that
$\mu((F\cap E_n)\triangle G_n)=0$.
Then $G:=\bigcup_{n\ge1}G_n$ belongs to $\mathscr G$ and
satisfies $\mu(F\triangle G)=0$.

\smallskip
\ref{it:Lp-separable-sigma2}$\Rightarrow$\ref{it:Lp-separable-sigma1}: \
Enlarging $\mathscr G$ if necessary, we may assume that it
contains all sets $E_n$. Choose a countable algebra
$\mathscr A$ generating $\mathscr G$ and containing these sets.
Let $V$ be the closed linear span in $L^p(\Om,\mu)$ of $\{\one_{A\cap E_n}:\,A\in\mathscr A,\ n\ge1\}$. For each fixed $n$, dominated convergence shows that
$$
\mathscr D_n:=
\{G\in\mathscr G:\,\one_{G\cap E_n}\in V\}
$$
is a Dynkin system. It contains $\mathscr A$, so Lemma \ref{lem:pi-lambda} gives $\mathscr D_n=\mathscr G$. It follows that $\one_{F\cap E_n}\in V$ for every $F\in\F$.
If $\mu(F)<\infty$, these indicators converge to $\one_F$
in $L^p(\Om,\mu)$. Density of the $\mu$-simple functions now gives
$V=L^p(\Om,\mu)$.
The rational linear combinations of the displayed countable
family (using coefficients in $\Q+i\Q$ in the complex case)
form a countable dense subset of $L^p(\Om,\mu)$.
\end{proof}

\subsection{Approximation by Mollification}\label{subsec:mollification}

It is generally difficult to handle $L^p$-functions directly. There are two ways of dealing with this: approximation reduces questions to simpler functions, while interpolation reduces them to more tractable exponents. The present section is devoted to approximation techniques; interpolation is treated in Section \ref{sec:interpolation}.

We start by proving that if $D$ is an open subset of $\R^d\!$, the vector space $C_{\rm c}^\infty(D)$\index{$Cac$@$C_{\rm c}^\infty(D)$} of all compactly supported smooth functions $f: D\to \K$ is dense in $L^p(D)$.
Here, and in what follows, the {\em support}\index{support!of a continuous function}
${\rm supp}(f)$ of a continuous function $f:D\to\K$ is defined as the
closure in $D$ of the set $\{x\in D:\,f(x)\not=0\}$.

\begin{proposition}[Approximation by compactly supported continuous functions]
\label{prop:Cc-dense}
Let $1\le p<\infty$ and let $D\subseteq\R^d$ be open.
Then $C_{\rm c}^\infty(D)$ is dense in $L^p(D)$.
\end{proposition}

\begin{proof}
We start by proving the density of the space $C_{\rm c}(D)$\index{$CC$@$C_{\rm c}(D)$} of all compactly supported continuous functions on $D$.
By Proposition \ref{prop:approx-simple}, for this purpose it suffices to
approximate indicator functions $\one_B$, where $B\subseteq D$
is measurable and has finite Lebesgue measure, with functions in $C_{\rm c}(D)$.

Fix $\eps>0$. By regularity of Lebesgue measure, there exist
a compact set $K$ and an open set $O$ such that
$$
K\subseteq B\subseteq O,
\qquad |O\setminus K|<\eps.
$$
If $K=\emptyset$, take $\phi=0$.
Otherwise, since $K$ is compact and contained in the open set
$D\cap O$, choose $\delta>0$ so small that
$$
K_\delta := \{x\in\R^d:\, d(x,K)\le\delta\}\subseteq D\cap O.
$$
Define
\begin{align}\label{eq:Ury-explicit}
\phi(x):=\max\Bigl\{0,1-\frac{d(x,K)}{\delta}\Bigr\},
\quad x\in D.
\end{align}
This function is continuous, equals one on $K$, and its support is
contained in the compact set $K_\delta$.
Thus, in either case, $\phi\in C_{\rm c}(D)$, $0\le\phi\le1$,
$\phi=1$ on $K$, and $\phi=0$ outside $O$.
The functions $\phi$ and $\one_B$ agree on $K$ and outside $O$,
and their difference has modulus at most one. Consequently,
$$
\n\phi-\one_B\n_p^p
=\int_D|\phi-\one_B|^p\ud x
\le |O\setminus K|<\eps.
$$
This proves the density
of $C_{\rm c}(D)$ in $L^p(D)$.

The proof that $C_{\rm c}^\infty(D)$ is dense will be given after
Theorem \ref{thm:approx-identity}.
\end{proof}

The corresponding result for $p=\infty$ is wrong: if $D$ is nonempty and open, then $C_{\rm c}(D)$, the vector space of all compactly supported continuous functions $f: D\to \K$,
fails to be dense in $L^\infty(D)$. Indeed, if $D'$ is a nonempty open set properly contained in $D$, then
$$\n f - \one_{D'}\n_{L^\infty(D)}\ge \frac12$$ for every $f\in C_{\rm c}(D)$.
On the positive side, it is possible to approximate functions in $L^\infty(D)$ {\em almost everywhere} by functions in $C_{\rm c}(D)$; see Theorem \ref{thm:Lusin}.

\begin{remark}\label{rem:CKdenseLpK}
The same argument shows that when $\mu$ is a finite Radon measure on a locally compact Hausdorff space $X$, then $C_{\rm c}(X)$ is dense in $L^p(X,\mu)$ for every $1\le p<\infty$; instead of \eqref{eq:Ury-explicit} we use Proposition \ref{prop:Urysohn2}.
The corresponding assertion for metric spaces, even when they are separable and complete, is wrong (see Problem \ref{prob:no-Cc-completemetric}).
\end{remark}

As an application of Proposition \ref{prop:Cc-dense} we prove an $L^p$-continuity result for translation.

\begin{proposition}[Continuity of translation]\label{prop:Lp-transl}\index{translation!continuity of, in $L^p$} Let  $f\in L^p(\R^d)$ with $1\le p<\infty$.
Then $$\lim_{|h|\to 0} \n f(\cdot+h) - f(\cdot)\n_p = 0.$$
\end{proposition}

\begin{proof}
 Define, for $h\in \R^d$ and $x\in\R^d\!$, $(\tau_h f)(x) := f(x+h)$, that is, $\tau_h f$
 is the translate of $f$ by $h$. Clearly, $\n \tau_h f\n_p = \n f\n_p$.

 First consider a function $f\in C_{\rm c}(\R^d)$. Such a function is uniformly continuous,
 so given an $\e>0$, we may choose $\delta>0$ such that $|x-x'|<\delta$ implies $|f(x)- f(x')|<\eps$. Hence if $|h|<\delta$, then
 for all $x\in \R^d$ we have $|\tau_h f(x) - f(x)| < \e$. Choose $r>0$ large enough such that
 the compact support of $f$ is contained in the rectangle $(-r,r)^d\!$. If $|h|<\delta$ is so small that the  support of
 $\tau_h f$ is also contained in $(-r,r)^d\!$, then
 $$\n \tau_h f - f\n_p^p = \int_{(-r,r)^d} |f(x+h)-f(x)|^p\ud x \le \e^p \int_{(-r,r)^d}\ud x = \e^p(2r)^{d}\!.$$
 This proves that $\lim_{|h|\to 0} \n \tau_h f-f\n_p = 0$
 for all $f\in C_{\rm c}(\R^d)$.

 Now let $f\in L^p(\R^d)$ be arbitrary.
 Since $C_{\rm c}(\R^d)$ is dense in $L^p(\R^d)$ by Proposition \ref{prop:Cc-dense}, we can find
 $g\in C_{\rm c}(\R^d)$ such that $\n f-g\n_p < \e$. Choose $\eta>0$ so small that $|h|<\eta$ implies
 $\n \tau_h g - g\n_p <\eps$; this is possible by what we just proved. Then, for $|h|<\eta$,
 $$ \n \tau_h f-f\n _p \le \n \tau_h f-\tau_h g\n_p + \n \tau_h g-g\n_p + \n g-f\n _p < \e+\e+\e = 3\e,$$
 noting that $\n \tau_h f- \tau_h g\n_p=\n \tau_h (f-g)\n_p = \n f-g\n_p <\e$.
\end{proof}

We now turn to an approximation technique based on convolution.
It relies on the following fundamental inequality.

\begin{proposition}[Young's inequality]\label{prop:Young}\index{inequality!Young}
Let $1\le p,q,r\le \infty$ satisfy $\frac1p+\frac1q = 1+\frac1r$, and let $f\in L^p(\R^d)$ and $g\in L^q(\R^d)$.
Then:
\begin{enumerate}[label={\rm(\arabic*)}, leftmargin=*]
 \item \label{it:Young1}
for almost all $x\in \R^d$ the function $y\mapsto f(x-y)g(y)$ is integrable;
\item \label{it:Young2} the function $f*g:\R^d\to \K$, defined for almost all $x\in \R^d$ by
$$ (f* g)(x) := \int_{\R^d} f(x-y)g(y)\ud y,$$
belongs to $L^r(\R^d)$ and we have
$$ \n f*g \n_r \le \n f\n_p \n g\n_q.$$
\end{enumerate}
Moreover, we have $f*g = g*f$ in $L^r(\R^d)$; and if  $1\le s,t,u\le \infty$ satisfy $\frac1r+\frac1s = 1+ \frac1t$ and
$\frac1q+\frac1s = 1+ \frac1u$,
then $\frac1p+\frac1u=1+\frac1t$ and for all $h\in L^s(\R^d)$ we have $(f*g)*h = f*(g*h)$ in $L^t(\R^d)$.
\end{proposition}

The most important special case corresponds to the choices $q=1$ and $r=p$, for which the proof below simplifies considerably.

\begin{proof} We prove the result for $1\le r<\infty$; the case $r=\infty$ is an immediate consequence of H\"older's inequality.

The identity $\frac1p+\frac1q = 1+\frac1r$ implies
$ \frac1r + \frac{r-p}{pr} + \frac{r-q}{qr}=1$ with $r\ge p,q$. Hence,
by elementary rewriting and
H\"older's inequality (for three functions, see  Remark \ref{rem:Holder}), for almost all $x\in \R^d$ we have
\begin{align*}
\ & \int_{\R^d} |f(x-y)g(y)|\ud y \\ &
\qquad = \int_{\R^d} (|f(x-y)|^p|g(y)|^q)^{1/r}\cdot |f(x-y)|^{(r-p)/r} \cdot |g(y)|^{(r-q)/r}\ud y
\\ & \qquad \le \Big\n \bigl(|f(x-\cdot)|^p|g(\cdot)|^q\bigr)^{1/r}\Big\n_r
\Big\n  |f(x-\cdot)|^{(r-p)/r}  \Big\n_{\frac{pr}{r-p}}
\Big\n  |g(\cdot)|^{(r-q)/r}  \Big\n_{\frac{qr}{r-q}}
\\ & \qquad = (I)\cdot(II)\cdot(III).
\end{align*}
Now
\begin{align*}
(I) & = \Bigl(\int_{\R^d}|f(x-y)|^p |g(y)|^q \ud y \Bigr)^{1/r}\!, \\
(II) & = \Bigl(\int_{\R^d}|f(x-y)|^p \ud y \Bigr)^{(r-p)/pr} = \n f\n_p^{(r-p)/r}\!, \\
(III) & =  \Bigl(\int_{\R^d}|g(y)|^q  \ud y \Bigr)^{(r-q)/qr} = \n g\n_q^{(r-q)/r}\!.
\end{align*}
Putting things together and using Fubini's theorem, we obtain
\begin{align*}
\int_{\R^d} | (f*g)(x)|^r\ud x
& \le \n f\n_p^{r-p}\n g\n_q^{r-q}\int_{\R^d}\int_{\R^d}|f(x-y)|^p |g(y)|^q \ud y\ud x
\\ & = \n f\n_p^{r-p}\n g\n_q^{r-q}\int_{\R^d}|g(y)|^q\int_{\R^d} |f(x-y)|^p\ud x\ud y
 = \n f\n_p^r\n g\n_q^r.
\end{align*}
This implies the first assertion as well as the second.

The identity $f*g=g*f$ follows by the change of variables
$y\mapsto x-y$.
For associativity, first note that the assumptions on the
exponents give
$$
\frac1p+\frac1u
=\frac1p+\frac1q+\frac1s-1
=\frac1r+\frac1s
=1+\frac1t.
$$
Applying Young's inequality twice shows that
$(|f|*|g|)*|h|\in L^t(\R^d)$.
Consequently, by Fubini's theorem, for almost all $x\in\R^d$
we have
$$
\int_{\R^d}\int_{\R^d}
|f(x-y-z)g(y)h(z)|\ud y\ud z
=\bigl((|f|*|g|)*|h|\bigr)(x)<\infty.
$$
Fubini's theorem and the substitution $w=y+z$ therefore give
\begin{align*}
((f*g)*h)(x)
&=\int_{\R^d}\int_{\R^d}
f(x-y-z)g(y)h(z)\ud y\ud z\\
&=\int_{\R^d}f(x-w)
\Bigl(\int_{\R^d}g(w-z)h(z)\ud z\Bigr)\ud w\\
&=(f*(g*h))(x).
\end{align*}
Both sides belong to $L^t(\R^d)$ by Young's inequality,
so the asserted identity holds in this space.
\end{proof}

The function $f*g$ is called the {\em convolution}\index{convolution}
of $f$ and $g$.

\begin{theorem}[Approximation by mollification]\label{thm:approx-identity}\index{mollification}
Let  $f\in L^p(\R^d)$ with $1\le p<\infty$. Let $\phi\in L^1(\R^d)$ satisfy $\int_{\R^d}\phi(x)\ud x=1.$
For $\eps>0$ define $$\phi_\eps(x):=\eps^{-d}\phi(\eps^{-1}x), \quad x\in\R^d\!.$$
Then $$\lim_{\eps\downarrow 0}\n \phi_\eps*f - f\n_p =0.$$
\end{theorem}
\begin{proof}
By Young's inequality, $\phi_\eps*f$ is well defined almost
everywhere and belongs to $L^p(\R^d)$.
Since $\int_{\R^d}\phi(y)\ud y=1$, a change of variables gives
\begin{align}\label{eq:millif-change-var}
\phi_\eps*f(x)-f(x)
=\int_{\R^d}\phi(y)[f(x-\eps y)-f(x)]\ud y
\end{align}
for almost every $x\in\R^d$.

Suppose first that $p>1$.
By H\"older's inequality with exponents $p/(p-1)$ and $p$, applied
to the product
$$
|\phi(y)|^{1-1/p}
\bigl(|\phi(y)|^{1/p}|f(x-\eps y)-f(x)|\bigr),
$$
together with the triangle inequality for the integral we obtain
\begin{align*}
|\phi_\eps*f(x)-f(x)|^p
&\le
\Bigl(\int_{\R^d}|\phi(y)|
      |f(x-\eps y)-f(x)|\ud y\Bigr)^p\\
&\le
\n\phi\n_1^{p-1}
\int_{\R^d}|\phi(y)|
      |f(x-\eps y)-f(x)|^p\ud y.
\end{align*}
For $p=1$, the same estimate follows directly from the triangle
inequality.

Integrating with respect to $x$ and applying Fubini's theorem,
we obtain
\begin{align*}
\n\phi_\eps*f-f\n_p^p
&\le
\n\phi\n_1^{p-1}
\int_{\R^d}|\phi(y)|
\int_{\R^d}|f(x-\eps y)-f(x)|^p\ud x\ud y\\
&=
\n\phi\n_1^{p-1}
\int_{\R^d}|\phi(y)|
\n f(\cdot-\eps y)-f\n_p^p\ud y.
\end{align*}
For each fixed $y\in\R^d$, Proposition \ref{prop:Lp-transl}
gives
$\n f(\cdot-\eps y)-f\n_p\to 0$ as $\eps\downarrow0$.
Moreover, since translation preserves the $L^p$-norm,
$$
\n f(\cdot-\eps y)-f\n_p^p\le 2^p\n f\n_p^p.
$$
The integrand in the preceding estimate is therefore bounded
by the function
$y\mapsto 2^p\n f\n_p^p|\phi(y)|$, which is integrable.
The dominated convergence theorem now implies
$\n\phi_\eps*f-f\n_p\to0$ as $\eps\downarrow0$.
\end{proof}

\begin{proof}[Completion of the proof of Proposition \ref{prop:Cc-dense}]
It remains to prove that $C_{\rm c}^\infty(D)$ is dense in $L^p(D)$.
Let $f\in C_{\rm c}(D)$ be given, and let $\phi\in C_{\rm c}^\infty(\R^d)$ be such that
$\int_{\R^d}\phi(x)\ud x = 1$, and for $\eps>0$ set $\phi_\eps(x):= \eps^{-d}\phi(\eps^{-1}x)$. Viewing $f$ as a function in $C_{\rm c}(\R^d)$ by zero extension, by standard calculus arguments we see that the convolution $f*\phi_\eps$ belongs to $C_{\rm c}^\infty(\R^d)$. Moreover, it is supported in $D$ for small enough $\eps>0$, and by the change of variable identity \eqref{eq:millif-change-var} and uniform continuity we have $f*\phi_\eps \to f$ uniformly. Since for small enough $\eps>0$, say $0<\eps<\eps_0$, the supports of $f*\phi_\eps$ are contained in a compact subset $K$ of $D$, this implies $f*\phi_\eps \to f$ in $L^p(D)$.
\end{proof}

\subsection{The Fr\'echet--Kolmogorov Compactness Theorem}\label{subsec:Frechet-Kolmogorov}

In this section we prove a characterisation of relatively compact sets in $L^p(\R^d)$. It will be used in Chapter \ref{chap:boundaryvalueproblems} to prove the Rellich--Kondrachov theorem on compact embeddings of Sobolev spaces.

Recall the notation $\tau_h f$ for the translate of a function $f$ by $h$,  $$\tau_hf(x) = f(x+h).$$

\begin{theorem}[Fr\'echet--Kolmogorov]\label{thm:Frechet-Kolmogorov}\index{theorem!Fr\'echet--Kolmogorov}\index{compactness!in $L^p(\Om)$} Let $1\le p<\infty$.
A subset $S$ of $L^p(\R^d)$ is relatively compact if and only if it satisfies the following two conditions:
\begin{enumerate}[label={\rm(\roman*)}, leftmargin=*]
 \item\label{it:FK1} $\displaystyle\lim_{|h|\to 0} \sup_{f\in S} \n \tau_h f-f\n_p = 0;$
 \item\label{it:FK2} $\displaystyle\lim_{\rho\to\infty}  \sup_{f\in S}\int_{\complement B(0;\rho)} |f(x)|^p\ud x = 0$.
\end{enumerate}
\end{theorem}

\begin{pictureproof}
`If': \
Let us begin by proving that \ref{it:FK1} and \ref{it:FK2} together imply that the set $S$ is bounded in $L^p(\R^d)$.
Choose $r>0$ such that $ \sup_{f\in S} \n \tau_h f-f\n_p \le 1$ for all $h\in \R^d$ with $|h|\le r$,
and choose $R>0$ such that $\sup_{f\in S}\int_{\complement B(0;R)} |f(x)|^p\ud x \le 1$.
Fix $h\in \R^d$ with $|h|= r$. For all $f\in S$ and $x\in \R^d$ we have
\begin{align*} \n \one_{B(x;R)}f\n_p
& \le  \n \one_{B(x;R)}(f-\tau_h f )\n_p + \n  \one_{B(x;R)} \tau_h f\n_p
\\ & = \n \one_{B(x;R)}(f-\tau_h f )\n_p + \n  \one_{B(x+h;R)} f\n_p
 \le 1 + \n  \one_{B(x+h;R)} f\n_p.
\end{align*}
Hence, by induction,
$$ \n \one_{B(0;R)}f\n_p \le N + \n \one_{B(Nh;R)}f\n_p.$$
Choose $N\ge 1$ such that $Nr = N|h| > 2R$. Then $B(Nh;R)\subseteq \complement B(0;R)$ and
\begin{align*}
\n f\n_p &
\le
\n \one_{B(0;R)}f\n_p + \n \one_{\complement B(0;R)}f\n_p
 \le  N + \n \one_{B(Nh;R)}f\n_p + \n \one_{\complement B(0;R)}f\n_p
\le N+2.
\end{align*}
This proves that $S$ is bounded in $L^p(\R^d)$.

Let us now prove that if $S$ satisfies \ref{it:FK1} and \ref{it:FK2}, then it is relatively compact. Let $M\ge0$ be a common
bound for the norms of the functions in $S$.
Fix $\eps>0$. By \ref{it:FK1} and \ref{it:FK2}, choose
$\delta>0$ and $R>0$ such that, writing $B:=B(0;R)$, for all
$f\in S$ we have
$\n\tau_hf-f\n_p\le\eps$ for $|h|\le\delta$, and $\n\one_{\complement B}f\n_p<\eps$.
Choose a nonnegative function
$\phi\in C_{\rm c}(B(0;\delta))$ with
$\int_{\R^d}\phi(y)\ud y=1$.
The estimate in the proof of Theorem
\ref{thm:approx-identity} gives
$$
\n\phi*f-f\n_p^p
\le\int_{\R^d}\phi(y)\n\tau_{-y}f-f\n_p^p\ud y
\le\eps^p,
\quad f\in S.
$$
Consequently,
\begin{align*}
\n f-\one_B(\phi*f)\n_p
&\le
\n\one_B(f-\phi*f)\n_p+\n\one_{\complement B}f\n_p
\le
\n f-\phi*f\n_p+\n\one_{\complement B}f\n_p
<2\eps.
\end{align*}
Thus every element of $S$ is within distance $2\eps$ of the set
$$S_\eps:=\{\one_B(\phi*f):\, f\in S\}.$$

We will prove that $S_\eps$ is relatively compact.

Let $q$ be conjugate to $p$.
By H\"older's inequality, for every $f\in S$ the integral
$$
(\phi*f)(x)=\int_{\R^d}\phi(x-y)f(y)\ud y
$$
is absolutely convergent for every $x\in\R^d$, and
$$
|(\phi*f)(x)|\le\n\phi\n_q\n f\n_p\le M\n\phi\n_q.
$$
Moreover, for $x,h\in\R^d$,
\begin{align*}
|(\phi*f)(x+h)-(\phi*f)(x)|
&\le
\int_{\R^d}|\phi(x+h-y)-\phi(x-y)|\,|f(y)|\ud y\le M\n\tau_h\phi-\phi\n_q.
\end{align*}
If $q<\infty$, Proposition \ref{prop:Lp-transl} shows that
the last expression tends to zero as $h\to0$.
If $q=\infty$, the same conclusion follows from the uniform
continuity of $\phi$.
These estimates show that the functions $\phi*f$, $f\in S$,
are continuous, uniformly bounded, and equicontinuous.

By the Arzel\`a--Ascoli theorem, their restrictions to $\ov B$
form a relatively compact subset of $C(\ov B)$.
The mapping from $C(\ov B)$ to $L^p(\R^d)$ which assigns to
$g$ the function equal to $g$ on $B$ and to zero outside $B$
is bounded; this follows from the bound
$\n\one_Bg\n_p\le |B|^{1/p}\n g\n_\infty$ for $g\in C(\ov B)$.
The image of the preceding family under this mapping is precisely $S_\eps$,
so $S_\eps$ is relatively compact in $L^p(\R^d)$.
Since $\eps>0$ was arbitrary, the approximation estimate
and Proposition \ref{prop:compact-totbdd} imply that $S$
is relatively compact.

\smallskip
`Only if': \ If $S$ is relatively compact in $L^p(\R^d)$, then \ref{it:FK1} and \ref{it:FK2} follow from Proposition \ref{prop:uniform-limits-on-K} applied to the operators $f\mapsto \tau_h f$ for $|h|\downarrow 0$ and
$f\mapsto \one_{\complement B(0;\rho)}f$ for $\rho\to\infty$, respectively.
\end{pictureproof}

We have the following immediate corollary for bounded domains.

\begin{corollary}\label{cor:RF}
Let $1\le p<\infty$ and let $D$ be a bounded open subset of $\R^d\!$.
A subset $S$ of $L^p(D)$ is relatively compact if and only if
$$\lim_{|h|\to 0} \sup_{f\in S} \n \tau_h f-f\n_p = 0.$$
\end{corollary}

 Here we
identify functions in $L^p(D)$ with their zero extensions in $L^p(\R^d)$.

\subsection{The Lebesgue Differentiation Theorem}\label{subsec:Lebesgue-diff}

By $L^1_{\rm loc}(\R^d)$\index{$L^1_{\rm loc}(\R^d)$} we denote the vector
space of functions $f:\R^d\to \K$ that are {\em locally integrable},\index{locally!integrable}
that is, integrable on every compact subset of $\R^d\!$, identifying two such functions when they are equal almost everywhere.
The aim of this section is to prove the Lebesgue differentiation theorem,
which says that if $f\in L^1_{\rm loc}(\R^d)$, then at almost every
point $x\in \R^d$ one has
$$\lim_{\substack{B\owns x \\ |B|\to 0}} \frac1{|B|}\int_{B} |f(y)-f(x)|\ud y = 0,$$
where the limit is taken as $B$ ranges over the balls in $\R^d$ containing $x$.
Here, $|B|$ denotes the Lebesgue measure of a measurable set $B$.
The proof of this theorem is based on the following lemma. For balls $B = B(x;r)$ in $\R^d$ and
real numbers $\la>0$ we set $\la B := B(x;\la r)$.

\begin{lemma}[Vitali covering lemma]\label{lem:Vitali}\index{lemma!Vitali covering}
Every finite collection $\mathscr{B}$ of open balls in $\R^d$ has a
subcollection $\mathscr{B}_0$ of pairwise disjoint balls such that each ball $B\in\mathscr{B}$ is
contained in $3B_0$ for some ball $B_0\in \mathscr{B}_0$.
\end{lemma}

\begin{proof}
The assertion is trivial if $\mathscr B$ is empty.
Otherwise, select a ball of largest radius and remove it,
together with all balls intersecting it, from the collection.
Repeat this procedure with the remaining balls until none remain,
and denote the collection of selected balls by $\mathscr B_0$.
Since $\mathscr B$ is finite, the procedure terminates.
The selected balls are pairwise disjoint.

Every ball $B=B(x;r)$ in the original collection intersects a
selected ball $B_0=B(x_0;r_0)$ whose radius satisfies $r\le r_0$:
either $B$ was selected itself, or it was removed when $B_0$
was selected.
For every $z\in B$ we therefore have
$$
|z-x_0|
\le|z-x|+|x-x_0|
<r+(r+r_0)\le3r_0.
$$
Thus $B\subseteq3B_0$, as required.
\end{proof}

For $f\in L^1_{\rm loc}(\R^d)$ define the {\em Hardy--Littlewood maximal function}\index{Hardy--Littlewood maximal function} $Mf:\R^d\to [0,\infty]$\index{$M2$@$Mf$} by
\begin{align*}
  Mf(x):=\sup_{B\owns
x}\frac{1}{|B|}\int_B|f(y)|\ud y,
\end{align*}
where the supremum is taken over all open balls $B$ containing $x$.
For every $t\ge0$ the set $\{Mf>t\}$ is the union of all open balls $B$ satisfying
$$|B|^{-1}\int_B|f(y)|\ud y>t.$$
Thus the sets $\{Mf>t\}$ are open. In particular, $Mf$ is
lower semicontinuous and hence measurable.

In the proof of the next theorem we need the following simple lemma.

\begin{lemma}\label{lem:layers}
Let $(\Omega,\mu)$ be a $\sigma$-finite measure space and let $g:\Omega\to\K$
be a measurable function. Then for all $0<p<\infty$ we have
$$
\int_\Omega |g|^{p}\ud\mu
= p\int_{0}^{\infty} t^{p-1}\,\mu(\{|g|>t\})\ud t,
$$
where the integrals are allowed to take the value $\infty$.
\end{lemma}

\begin{proof}
For all $\om\in\Om$ we have
$$
|g(\omega)|^{p}=\int_{0}^{|g(\omega)|^{p}} \one \ud s
=\int_{0}^{\infty}\one_{\{s<|g(\omega)|^{p}\}}\ud s.
$$
Integrating over $\Omega$ and using Fubini's theorem, we obtain
$$
\int_\Omega |g|^{p}\ud\mu
=\int_{0}^{\infty}\mu(\{|g|^{p}>s\})\ud s.
$$
The desired identity follows by substituting $s=t^{p}$ on the right-hand side.
\end{proof}

\begin{theorem}[Hardy--Littlewood maximal theorem]\label{thm:HL-maximal}\index{theorem!Hardy--Littlewood maximal}
For all  $f\in L^1(\R^d)$ and $t>0$ we have the {\em weak $L^1$-bound}\index{weak!$L^1$-bound}
$$t |\{Mf >t\}|\le 3^d \|f\|_1.$$
Moreover, for all $1<p\le \infty$ there exists a constant $C_{d,p}\ge 0$ such that
for all $f\in L^p(\R^d)$ we have $Mf\in L^p(\R^d)$ and
$$\n Mf\n_p\le C_{d,p}\n f\n_p.
$$
\end{theorem}

\begin{proof}
We begin with the proof of the first assertion. By the definition of $Mf$, for every $x\in\{Mf>t\}$ there exists a ball $B$ containing $x$ such that
$\frac1{|B|}\int_B|f(y)|\ud y>t$. If $K\subseteq\{Mf>t\}$ is a compact subset, it can be covered by a
finite collection $\mathscr{B}$ of such balls. Let $\mathscr{B}_0$ be a disjoint
subcollection of this cover provided by the Vitali covering lemma. Then,
\begin{align*}
|K| \leq \Big|\bigcup_{B\in\mathscr{B}}B\Big|
    \leq \Big|\bigcup_{B\in\mathscr{B}_0}3B\Big|
\le
\sum_{B\in\mathscr{B}_0} 3^d|B|
    \leq \sum_{B\in\mathscr{B}_0} \frac{3^d}{t}\int_{B} |f(y)|\ud y
    \leq \frac{3^d}{t}\n f\n_{1}.
\end{align*}
Since this holds for all compact sets $K$ contained in the open set $\{Mf>t\}$, the first assertion follows.

For $1<p<\infty$, the second assertion follows from the first
as follows.
For any $f\in L^p(\R^d)$ and $t>0$ the function $f_t(x) := \one_{\{|f|\ge t/2\}}f$ belongs to $L^1(\R^d)\cap L^p(\R^d)$, with $\n f_t\n_1\le (2/t)^{p-1}\n f\n_p^p$, and satisfies the pointwise bound
$$ Mf \le  \sup_{B\owns x}\frac{1}{|B|}\int_B|f_t(y)|\ud y + \sup_{B\owns
x}\frac{1}{|B|}\int_B|\one_{\{|f|< t/2\}}f(y)|\ud y\le Mf_t + t/2,$$ which implies
$$ \{Mf > t\}\subseteq \{Mf_t > t/2\}.$$ Hence, by the first part of the theorem,
\begin{align}\label{eq:ibp-p-2} |\{Mf > t\}| \le |\{Mf_t > t/2\}| \le \frac{2\cdot 3^d}{t} \|f_t\|_1 = \frac{2\cdot 3^d}{t}\int_{\{|f|\ge t/2\}} |f(x)|\ud x.
\end{align}
By Lemma \ref{lem:layers},
\eqref{eq:ibp-p-2}, and Fubini's theorem,
\begin{align*}
 \int_{\R^d} |Mf(x)|^p \ud x
 & \le  p\int_0^\infty t^{p-1} \Bigl(\frac{2\cdot 3^d}{t}\int_{\{|f|\ge t/2\}} |f(x)|\ud x\Bigr)\ud t
 \\ & = 3^d\cdot 2p \int_{\R^d} |f(x)| \int_0^{2|f(x)|} t^{p-2} \ud t \ud x
 \\ & = \frac{3^d\cdot 2p}{p-1} \int_{\R^d} |f(x)| (2|f(x)|)^{p-1}\ud x = C_{d,p}^p  \int_{\R^d} |f(x)|^p \ud x,
\end{align*}
where $C_{d,p} = 2(\frac{3^d p}{p-1})^{1/p}$\!.

For $p=\infty$ the second assertion follows trivially from the pointwise inequality $Mf \le \n f\n _\infty$,
with constant $C_{d,\infty} =1$.
\end{proof}

Inspection of this proof reveals that the derivation of the $L^p$-bound
uses only the subadditivity of the maximal operator, its weak $L^1$-bound,
and its $L^\infty$-bound. This observation forms the basis of the Marcinkiewicz interpolation theorem in Chapter \ref{ch:operators}
(see Theorem \ref{thm:Marcinkiewicz}).

As a corollary of Theorem \ref{thm:HL-maximal}, we have the following fundamental result.

\begin{theorem}[Lebesgue differentiation theorem]\label{thm:Leb-diff}\index{theorem!Lebesgue differentiation}
If $f\in L_{\rm loc}^1(\R^d)$, then for almost all $x\in\R^d$ we have
\begin{align}\label{eq:Leb-pt} \lim_{\substack{B\owns x\\ |B|\to 0}} \frac{1}{|B|}\int_B  |f(y)-f(x)| \ud  y = 0.\end{align}
\end{theorem}

The correct way of interpreting this theorem is as follows. For every pointwise defined locally integrable function $\wt f$ on $\R^d$\!, the limit in \eqref{eq:Leb-pt} (with $f$ replaced by $\wt f$) exists for almost all $x\in \R^d$, say on a Borel set $\Om\subseteq\R^d$ such that $|\R^d\setminus \Om| = 0$. If both $\wt f_1$ and $\wt f_2$ are pointwise representatives, the symmetric difference of corresponding sets $\Om_1$ and $\Om_2$ has measure $0$.

The set of all points $x\in \R^d$ for which \eqref{eq:Leb-pt} holds is called the set of {\em  Lebesgue points}\index{Lebesgue!point} of $f$. As just explained, this set is uniquely determined only up to a set of measure $0$
(see also Problem \ref{prob:Lebesgue-set}).

\begin{proof} A point $x\in \R^d$ is a Lebesgue point of $f$ if and only if it is a Lebesgue point of $\one_U f$, for any bounded open set $U$ containing $x$. Hence, upon replacing $f$ by $\one_U f$ if necessary, and keeping in mind that $\R^d$ can be covered by countably many bounded open sets,
it suffices to prove the theorem under the stronger assumption $f\in L^1(\R^d)$.

Fixing a pointwise defined representative of $f$, for all $x\in \R^d$ let
\begin{equation*}
  N f(x):=\limsup_{\substack{B\owns x \\ |B|\to 0}}\frac{1}{|B|}\int_B  | f(y)-f(x)| \ud  y.
\end{equation*} As for $Mf$, the limsup may be taken over rational balls, so $Nf$ is measurable. We wish to prove that $N f(x)=0$ for almost all $x\in\R^d\!$.
For this purpose it suffices to show that $|\{N f>\eps\}|=0$ for any fixed $\eps>0$.

For any fixed $\delta>0$, Proposition \ref{prop:Cc-dense} provides us with a function
$g\in C_{\rm c}(\R^d)$ such that $\n{f-g}\n_{1}<\delta$. Then,
$$ N f \le N(f-g)+N g \le M(f-g)+|f-g|+0, $$
using the pointwise inequality
$N h(x)\leq Mh(x)+|h(x)|$ and the continuity of $g$, which implies $N g(x)=0$.
Therefore, by Theorem \ref{thm:HL-maximal},
\begin{align*}
|{\{N f>\eps\}}| & \le |\{M(f-g)>\eps/2\}| + |\{|f-g|>\eps/2\}|
\\ & \le \frac{2\cdot 3^d}{\eps}\n f-g\n_{1}+\frac{2}{\eps}\n f-g\n_{1} \le \frac{2}{\eps}(3^d+1)\delta.
\end{align*}
Since $\delta>0$ was arbitrary, it follows that $|\{N f>\eps\}|=0$.
\end{proof}

\begin{corollary}[Lebesgue density points]\label{cor:density-points}
Let $E$ be a Borel subset of $\R^d$. Then almost every point of $E$ is a
{\em Lebesgue density point}\index{density!point} of $E$, that is, for almost all $x\in E$ we have
$$
\lim_{r\downarrow 0}\frac{|E\cap B(x;r)|}{|B(x;r)|}=1.
$$
\end{corollary}

\begin{proof}
Applying Theorem \ref{thm:Leb-diff} to the locally integrable function
$\one_E$, for almost every $x\in\R^d$ we obtain
$$
\lim_{r\downarrow 0}\frac1{|B(x;r)|}
       \int_{B(x;r)} \one_E(y)\ud y
 = \one_E(x).
$$
Since the average on the left is precisely $|E\cap B(x;r)|/|B(x;r)|$,
the result follows.
\end{proof}

\section{Spaces of Measures}\label{sec:measures}

In this section we introduce the space $M(\Om)$\index{$M3$@$M(\Om)$} of $\K$-valued measures on a given measurable space $(\Om,\F)$ and discuss some of its properties. From the functional-analytic point of view, the importance of this space resides in the fact that $M(\Om)$ is a vector space in a natural way by setting
$$(c\mu)(F):= c\mu(F), \qquad (\mu+\nu)(F):= \mu(F)+\nu(F), \quad F\in \calF\!,$$
and that it is a Banach space with respect to the variation norm introduced below in Definition \ref{def:variation}
(see Theorem \ref{thm:completeness-M}).

\subsection{The Banach Space \texorpdfstring{$M(\Om)$}{M(Omega)}}\label{subsec:MOmega}

In what follows we fix a measurable space $(\Om,\F)$.

\begin{definition}[$\K$-valued measures] A {\em $\K$-valued measure}\index{measure!real}\index{measure!complex} on $(\Om,\F)$ is a mapping
$ \mu: \F \to \K$ with the following properties:

\begin{enumerate}[label={\rm(\roman*)}, leftmargin=*]
\item $\mu(\emptyset) = 0$;
\item for all disjoint $\calF$-measurable sets $F_1, F_2,\dots$ we have
      $$\mu\big(\bigcup_{n\ge 1} F_n\big) = \sum_{n\ge 1} \mu(F_n).$$
\end{enumerate}
\end{definition}

\begin{remark}
An ordinary measure (in the sense of Definition \ref{def:measure})
is a real-valued measure if and only if it is finite.
\end{remark}

The terms `real-valued measure' and `complex-valued measure' are often abbreviated to `real measure' and `complex measure'.

\begin{definition}[Variation]\label{def:variation} Let $\mu$ be a $\K$-valued measure on $(\Om,\F)$.
 The {\em variation of $\mu$ on the set  $F\in \F$}\index{variation} is defined by
 $$ |\mu|(F) := \sup_{\mathscr{A}\in \mathbb{F}_F} \sum_{A\in \mathscr{A}}|\mu(A)|,$$
where $\mathbb{F}_F$ denotes the set of all
finite collections of pairwise disjoint $\calF$-measurable subsets of $F$.
\end{definition}

It is immediate that $|\mu(F)|\le |\mu|(F)$ and $|\mu|(F) \le |\mu|(F')$ whenever $F,F'\in \F$ and $F\subseteq F'$. Also, $|c\mu| = |c||\mu|$ for all $c\in\K$ and $|\mu+\nu|(F) \le |\mu|(F)+|\nu|(F)$ for all $F\in\calF$\!. If $\mu$ takes values in $[0,\infty)$, then $|\mu| = \mu$.

We next show that $\K$-valued measures have finite variation.

\begin{proposition}\label{prop:ansmu} If $\mu$ is a $\K$-valued measure, then $|\mu|$ is a finite measure.
\end{proposition}

\begin{proof} We proceed in two steps.

\smallskip
 {\em Step 1} -- We prove that $|\mu|$ is a measure. It is clear that $|\mu|(\emptyset) = 0$. Let $(F_n)_{n\ge 1}$ be a sequence
 of pairwise disjoint measurable sets and let $F$ be their union.
 We must prove that $|\mu|(F) = \sum_{n\ge 1}|\mu|(F_n).$

 If $\mathscr{A}_n\in \mathbb{F}_{F_n}$ is a finite collection of pairwise disjoint
 measurable subsets of $F_n$, then for every $N\ge 1$ the union $\bigcup_{n=1}^N \mathscr{A}_n$ is a finite collection of pairwise disjoint
 measurable subsets of $F$ and therefore
 $$\sum_{n=1}^N \sum_{A\in \mathscr{A}_n} |\mu(A)| \le |\mu|(F).$$
 Taking the supremum over all $\mathscr{A}_n\in \mathbb{F}_{F_n}$, we find that
 $\sum_{n=1}^N |\mu|(F_n) \le |\mu|(F).$ Since this holds for all $N\ge 1$, we conclude that
 $$\sum_{n\ge 1}|\mu|(F_n) \le |\mu|(F).$$
 In the converse direction, suppose that the measurable subsets $A_1,\dots, A_k$ of $F$ are disjoint.
 With $F_{j,n}:= A_j\cap F_n$ we have
 $$ \sum_{j=1}^k |\mu(A_j)| \le \sum_{j=1}^k \sum_{n\ge 1} |\mu(F_{j,n})|
 = \sum_{n\ge 1}  \sum_{j=1}^k|\mu(F_{j,n})| \le \sum_{n\ge 1}|\mu|(F_n).$$
 Taking the supremum over all finite families of pairwise disjoint measurable subsets of $F$,
 we find that
 $$ |\mu|(F) \le \sum_{n\ge 1}|\mu|(F_n).$$

 {\em Step 2} -- We prove that the measure $|\mu|$ is finite, that is,  $|\mu|(\Om)<\infty$.
 By considering real and imaginary parts, it suffices to do this in the real-valued case.

  If a finite set  $\{r_1,\dots,r_N\}$ of real numbers is given, then either the positive or the negative  numbers in this set (or both) contribute at least half to the sum $\sum_{n=1}^N |r_n|$.
  Enumerating one such contributing subset as $r_{n_1},\dots, r_{n_k}$, we thus have
  \begin{align}\label{eq:supzn} \Big|\sum_{j=1}^k r_{n_j}\Big| \ge \frac12\sum_{n=1}^N|r_n|.
  \end{align}

  Suppose, for a contradiction, that some measurable set $F$ satisfies $|\mu|(F) = \infty$.
  Choose disjoint measurable subsets $F_1,\dots,F_N$ of $F$ such that
  $$\sum_{n=1}^N |\mu(F_n)| \ge 2 (1+|\mu(F)|).$$
  By \eqref{eq:supzn} the union $F'$ of a suitable subcollection of the $F_j$ satisfies
  $$ |\mu(F')| \ge \frac{1}{2} \sum_{n=1}^N |\mu(F_n)| \ge
  1+|\mu(F)|.$$
  For $F'':= F\setminus F'$ we then have
  $$ |\mu(F'')| \ge \bigl||\mu(F)| - |\mu(F')|\bigr| \ge 1.$$
  Thus if a set $F\in \calF$ satisfies $|\mu|(F)=\infty$, there is a disjoint decomposition $F = F'\cup F''$
  with $|\mu(F')|\ge 1$ and $|\mu(F'')|\ge 1$. Since $|\mu|$ is a measure, at least one of the numbers
  $|\mu|(F')$ and $|\mu|(F'')$ equals $\infty$, and therefore we can continue this procedure inductively. For each $k\ge 1$ let $G_k$ be the set discarded at the $k$th step. This produces a sequence of pairwise disjoint measurable sets $G_1,G_2,\dots$, each of which satisfies
  $|\mu(G_k)|\ge 1$.
Let $G$ be their union. Since $\mu$ is a $\K$-valued measure
  we have $\mu(G) = \sum_{k\ge 1} \mu(G_k)$. This sum cannot converge since its terms fail to converge to $0$.
  This is the required contradiction.
\end{proof}

\begin{theorem}[Completeness]\label{thm:completeness-M}
Endowed with the variation norm $$\n \mu\n := |\mu|(\Om),$$
the space $M(\Om)$ of $\K$-valued measures on $(\Om,\calF)$ is a Banach space.
\end{theorem}

\begin{proof}
We leave it as an exercise to prove that $|\mu|(\Om)$
defines a norm.
By Proposition \ref{prop:competeness-series}, to prove
completeness it suffices to show that every series
$\sum_{n\ge1}\mu_n$ with
$\sum_{n\ge1}\n\mu_n\n<\infty$ converges in this norm.

Let such a series be given. For every $F\in\calF$ we have
$$
\sum_{n\ge1}|\mu_n(F)|
\le\sum_{n\ge1}\n\mu_n\n<\infty.
$$
We may therefore define
$$
\mu(F):=\sum_{n\ge1}\mu_n(F),
\quad F\in\calF.
$$
Clearly $\mu(\emptyset)=0$.
To prove countable additivity, let $(F_j)_{j\ge1}$ be a
sequence of pairwise disjoint measurable sets.
Since each $|\mu_n|$ is a finite measure, we have
$$
\sum_{n\ge1}\sum_{j\ge1}|\mu_n(F_j)|
\le\sum_{n\ge1}|\mu_n|(\Om)
=\sum_{n\ge1}\n\mu_n\n<\infty.
$$
Thus the double series is absolutely convergent. Using the countable additivity
of each $\mu_n$, we obtain
\begin{align*}
\mu\Bigl(\bigcup_{j\ge1}F_j\Bigr)
=\sum_{n\ge1}\mu_n\Bigl(\bigcup_{j\ge1}F_j\Bigr)
=\sum_{n\ge1}\sum_{j\ge1}\mu_n(F_j)
=\sum_{j\ge1}\sum_{n\ge1}\mu_n(F_j)
=\sum_{j\ge1}\mu(F_j).
\end{align*}
Hence $\mu\in M(\Om)$.

It remains to prove convergence in the variation norm.
For any finite family $F_1,\dots,F_k$ of pairwise disjoint
measurable sets,
\begin{align*}
\sum_{j=1}^k
\Bigl|\Bigl(\mu-\sum_{n=1}^N\mu_n\Bigr)(F_j)\Bigr|
\le
\sum_{j=1}^k\sum_{n>N}|\mu_n(F_j)|
=
\sum_{n>N}\sum_{j=1}^k|\mu_n(F_j)|
\le\sum_{n>N}\n\mu_n\n.
\end{align*}
Taking the supremum over all such finite families gives
$$
\Big\n\mu-\sum_{n=1}^N\mu_n\Big\n
\le\sum_{n>N}\n\mu_n\n\to 0.
$$
Thus the series converges in $M(\Om)$, proving completeness.
\end{proof}

If $\mu$ is a complex measure on $(\Om,\F)$, then
  \begin{align*}
  (\Re \mu)(F)  := \Re (\mu(F)), \qquad  (\Im \mu)(F) := \Im (\mu(F)),
  \end{align*}
define real measures on $(\Om,\F)$ and we have $\mu = \Re\mu + i\Im \mu.$
The next result shows that real measures allow decompositions into positive and negative parts:

\begin{theorem}[Hahn--Jordan decomposition]\label{thm:Jordan}\index{theorem!Hahn--Jordan}
If $\mu$ is a real measure on a measurable space $(\Om,\F)$, then\index{$Ma$@$\mu^{\pm}$}
  \begin{align*}
  \mu^+(F) &:= \sup\big\{\mu(A):\, A\in \calF\!, \ A\subseteq F\big\}, \\
     \mu^-(F)& := -\inf\big\{\mu(A):\, A\in \calF\!, \ A\subseteq F\big\}
  \end{align*}
  are finite measures on $(\Om,\F)$ and $$\mu = \mu^+ - \mu^-, \qquad |\mu| = \mu^+ + \mu^-\!.$$
The measures $\mu^+$ and $\mu^-$ are supported on disjoint sets, in the sense that
there exists a disjoint decomposition $\Om = \Om^+\cup\Om^-$ with $\Om^\pm\in\calF$ such that for all $F\in \calF$ we have
$$ \mu^+(F) = \mu(F\cap \Om^+), \qquad \mu^-(F) = -\mu(F\cap \Om^-).$$
If $\nu_1$ and $\nu_2$ are finite measures on $(\Om,\calF)$ such that
 $\mu = \nu_1-\nu_2$, then for all $F\in \calF$ we have
 $$\nu_1(F)\ge \mu^+(F), \qquad \nu_2(F)\ge \mu^-(F).$$
\end{theorem}

The decomposition $\mu = \mu^+-\mu^-$ for real measures $\mu$ is called the {\em Jordan decomposition}\index{Jordan!decomposition}; the existence of a corresponding decomposition  $\Om = \Om^+\cup\Om^-$ for their supports is often referred to as the {\em Hahn decomposition theorem}.\index{Hahn decomposition}

\begin{proof}
By Proposition \ref{prop:ansmu}, the variation $|\mu|$ is
a finite measure. In particular, $\mu$ is continuous under
increasing and decreasing limits of measurable sets.
For example, if $F_n\uparrow F$, then
$$
|\mu(F)-\mu(F_n)|
\le|\mu|(F\setminus F_n)\to 0.
$$
The argument for $F_n\downarrow F$ is the same, using
$F_n\setminus F$.

\smallskip
{\em Step 1} --
Put
$$
a:=\inf\{\mu(F):\, F\in\calF\}.
$$
Since $|\mu(F)|\le|\mu|(\Om)$, this infimum is finite.
We first show that there exists a measurable set $\Om^-$
such that $\mu(\Om^-)=a$.

Choose measurable sets $E_n$ such that
$\mu(E_n)<a+2^{-n}$.
For any measurable sets $F$ and $G$, finite additivity gives
$$
\mu(F\cup G)
=\mu(F)+\mu(G)-\mu(F\cap G)
\le\mu(F)+\mu(G)-a.
$$
Applying this inequality repeatedly, for $m\ge n$ we obtain
\begin{align*}
a
\le\mu\Bigl(\bigcup_{k=n}^m E_k\Bigr)
\le\sum_{k=n}^m\mu(E_k)-(m-n)a\le a+\sum_{k=n}^m2^{-k}.
\end{align*}
Set
$F_n:=\bigcup_{k\ge n}E_k$ and $\Om^-:=\bigcap_{n\ge1}F_n$.
Letting $m\to\infty$ in the preceding estimate and using
continuity under increasing limits gives
$$
a\le\mu(F_n)\le a+\sum_{k\ge n}2^{-k}.
$$
Since $F_n\downarrow\Om^-$, continuity under decreasing
limits now gives $\mu(\Om^-)=a$.

Put $\Om^+:=\Om\setminus\Om^-$.
If $A\in\calF$ and $A\subseteq\Om^-$, then
$$
a\le\mu(\Om^-\setminus A)=a-\mu(A),
$$
so $\mu(A)\le0$.
If instead $A\subseteq\Om^+$, then
$$
a\le\mu(\Om^-\cup A)=a+\mu(A),
$$
so $\mu(A)\ge0$.
Thus $\Om^-$ is a negative set and $\Om^+$ is a positive
set: every measurable subset of $\Om^-$ has nonpositive
measure, and every measurable subset of $\Om^+$ has
nonnegative measure.

\smallskip
{\em Step 2} --
We now identify the measures $\mu^+$ and $\mu^-$.
If $F\in\calF$ and $A\in\calF$ satisfies $A\subseteq F$,
the sign properties just proved imply
$$
\mu(A)
=\mu(A\cap\Om^+)+\mu(A\cap\Om^-)
\le\mu(A\cap\Om^+)
\le\mu(F\cap\Om^+).
$$
Equality is attained by taking $A=F\cap\Om^+$.
Taking the supremum over $A$, we therefore obtain
$$
\mu^+(F)=\mu(F\cap\Om^+).
$$
Similarly, the infimum of $\mu(A)$ over measurable subsets
$A$ of $F$ is attained at $A=F\cap\Om^-$, so
$$
\mu^-(F)=-\mu(F\cap\Om^-).
$$
These identities show that $\mu^+$ and $\mu^-$ are finite
measures supported on $\Om^+$ and $\Om^-$, respectively,
and that $\mu=\mu^+-\mu^-$.

The triangle inequality for variation gives
$$
|\mu|(F)
=|\mu^+-\mu^-|(F)
\le\mu^+(F)+\mu^-(F).
$$
For the reverse inequality, use the disjoint measurable
subsets $F\cap\Om^+$ and $F\cap\Om^-$ in the definition
of variation:
\begin{align*}
|\mu|(F)\ge|\mu(F\cap\Om^+)|+|\mu(F\cap\Om^-)|=\mu^+(F)+\mu^-(F).
\end{align*}
This proves $|\mu|=\mu^++\mu^-$.

Finally, suppose that $\mu=\nu_1-\nu_2$, where $\nu_1$
and $\nu_2$ are finite measures.
For every measurable subset $A$ of $F$ we have
$$
\mu(A)\le\nu_1(A)\le\nu_1(F),
\qquad
-\mu(A)\le\nu_2(A)\le\nu_2(F).
$$
Taking suprema over $A$ gives
$\mu^+(F)\le\nu_1(F)$ and $\mu^-(F)\le\nu_2(F)$,
as required.
\end{proof}

\subsection{The Radon--Nikod\'ym Theorem}\label{subsec:RN}

If $f:\Om\to \K$ is integrable with respect to the measure $\mu$,
by dominated convergence the formula
 $$ \nu(F):= \int_F f \ud \mu, \quad F\in \calF\!,$$
 defines a $\K$-valued measure $\nu$. This measure
is {\em absolutely continuous} with respect to $\mu$,\index{absolutely!continuous, measure} that is, $\mu(F)=0$ implies $\nu(F) = 0$. The following theorem provides a converse under a $\sigma$-finiteness assumption.

\begin{theorem}[Radon--Nikod\'ym]\label{thm:RN}\index{theorem!Radon--Nikod\'ym}
Let $(\Om,\calF\!,\mu)$ be a $\sigma$-finite measure space.
If the $\K$-valued measure $\nu:\calF\to \K$ is absolutely continuous
with respect to $\mu$, then there exists a unique $g\in L^1(\Om,\mu)$
such that
$$
\nu(F) = \int_F g\ud \mu, \quad F\in \calF\!.
$$
\end{theorem}

\begin{proof} Uniqueness is clear, so it remains to prove existence.
By considering real and imaginary parts separately, we may restrict attention to the case of real scalars.

If $\mu(F)=0$, then every measurable $A\subseteq F$ has $\mu(A)=0$, hence $\nu(A)=0$, so $\nu^\pm(F)=0$. Thus $\nu^+$ and $\nu^-$ are absolutely continuous with respect to $\mu$.
We may therefore restrict attention to the case where $\nu$ is a finite
nonnegative measure.

Consider the set
$$ S:= \Bigl\{f\in L^1(\Om,\mu):\, f\ge 0,\, \int_F f\ud\mu \le \nu(F)\ \hbox{ for all $F\in \calF$}\Bigr\}.$$
Then $0\in S$, so $S$ is nonempty. Let $$M:= \sup_{f\in S}\, \int_\Om f\ud \mu.$$
For all $f\in S$ we have
$\int_{\Om} f\ud\mu  \le  \nu(\Om)$ and therefore $M \le \nu(\Om)<\infty$.

\smallskip
{\em Step 1} -- In this step we prove that there exists a function $g\in S$ for which the supremum in the definition of $M$ is attained. Let $(f_n)_{n\ge 1}$ be a sequence in $S$ with the property that $\limn \int_\Om f_n\ud \mu= M$.
Set $g_n:= f_1\vee \cdots\vee f_n$. Any set $F\in \calF$ can be written as a disjoint union of
sets $F_1^{(n)}\!, \dots,F_n^{(n)}\in \calF$ such that $g_n = f_j$ on $F_j^{(n)}$
and therefore
$$ \int_F g_n\ud\mu = \sum_{j=1}^n \int_{F_j^{(n)}}f_j\ud\mu \le  \sum_{j=1}^n \nu(F_j^{(n)})
= \nu(F).$$
It follows that $g_n\in S$. The sequence $(g_n)_{n\ge 1}$ is nondecreasing and therefore
its pointwise limit $g:= \limn g_n$ is well defined as a $[0,\infty]$-valued function. By the monotone convergence theorem, for all $F\in \calF$ we have
\begin{align}\label{eq:RN-intFg} \int_F g\ud \mu = \limn \int_F g_n\ud \mu \le \nu(F)
\end{align}
and therefore $g$ takes finite values $\mu$-almost everywhere and belongs to $S$. Moreover,
$$ M = \limn \int_\Om f_n\ud \mu \le \limn\int_\Om g_n\ud \mu = \int_\Om g\ud \mu \le M$$
and therefore equality holds throughout. This proves that the supremum in the definition of $M$ is attained by the function $g\in S$.

\smallskip
{\em Step 2} --
We show that the maximising function $g$ obtained in Step 1
represents $\nu$. By \eqref{eq:RN-intFg}, the formula
$$
\eta(F):=\nu(F)-\int_F g\ud\mu,\quad F\in\calF,
$$
defines a finite nonnegative measure. Moreover, $\eta$ is
absolutely continuous with respect to $\mu$.

Suppose, for a contradiction, that $\eta(\Om)>0$.
By $\sigma$-finiteness, there exists a measurable set $E$ such
that $\mu(E)<\infty$ and $\eta(E)>0$.
Choose $c>0$ so small that $\eta(E)>c\mu(E)$, and apply the
Hahn--Jordan decomposition to the finite real-valued measure
$\eta|_E-c\mu|_E$. Let $P\subseteq E$ be its positive set.
Since the measure is nonpositive on $E\setminus P$, we have
$$
\eta(P)-c\mu(P)
\ge\eta(E)-c\mu(E)>0.
$$
In particular, $\eta(P)>0$, and absolute continuity implies
$\mu(P)>0$.

By the defining property of $P$, for every $F\in\calF$ we have
$$
c\mu(F\cap P)\le\eta(F\cap P)\le\eta(F).
$$
Consequently, the function $h:=g+c\one_P$ is nonnegative,
belongs to $L^1(\Om,\mu)$, and satisfies
$$
\int_F h\ud\mu
=\int_F g\ud\mu+c\mu(F\cap P)
\le\int_F g\ud\mu+\eta(F)
=\nu(F).
$$
Thus $h\in S$. But
$$
\int_\Om h\ud\mu
=M+c\mu(P)>M,
$$
contradicting the definition of $M$.
It follows that $\eta=0$, so $g$ has the required properties.
\end{proof}

An alternative proof of the Radon--Nikod\'ym theorem, based on Hilbert space methods, is outlined in Problem \ref{prob:RN}.

\begin{example}\label{ex:densityfunction}
If $f:\Om\to \K$ is integrable with respect to $\mu$, then
 $$ \nu(F):= \int_F f \ud \mu, \quad F\in \calF\!,$$
 defines a $\K$-valued measure and for all $F\in\F$ we have
 $$|\nu|(F) = \int_F |f| \ud \mu.$$
If $f$ is real-valued, then $\nu$ is a real measure and
 $$ \nu^\pm(F) = \int_F f^\pm \ud \mu.$$

 To prove the first assertion, let $A_1,\dots,A_n\in \calF$ be disjoint subsets of $F$.
Then
$$ \sum_{j=1}^n |\nu(A_j)| \le \sum_{j=1}^n\int_{A_j} |f| \ud \mu \le \int_F |f| \ud \mu,
$$
which gives the upper bound `$\le$'\!.

To prove the lower bound `$\ge$'\!, define the measurable function
$$
  \theta :=
  \begin{cases}
  \overline f/|f|, & f\ne0,\\
  0, & f=0.
  \end{cases}
$$
Then $|\theta|\le1$ and $\theta f=|f|$. Choose simple functions
$\theta_n$ supported in $F$, with $|\theta_n|\le1$, such that
$\theta_n\to \one_F\theta$ pointwise. Writing
$\theta_n=\sum_{j=1}^{N_n} c_j^{(n)}\one_{A_j^{(n)}}$, with the
sets $A_j^{(n)}\in\calF$ pairwise disjoint subsets of $F$ and
$|c_j^{(n)}|\le1$, we obtain
$$
  \Big|\int_\Omega \theta_n f\ud\mu\Big|
  =
  \Big|\sum_{j=1}^{N_n} c_j^{(n)}\nu(A_j^{(n)})\Big|
  \le
  \sum_{j=1}^{N_n} |\nu(A_j^{(n)})|
  \le |\nu|(F).
$$
By dominated convergence,
$$
  \int_F |f|\ud\mu
  =
  \int_\Omega \one_F\theta f\ud\mu
  =
  \lim_{n\to\infty}\int_\Omega \theta_n f\ud\mu
  \le |\nu|(F).
$$

To prove the second assertion we note that
the sets $\Om^+ = \{f\ge 0\}$ and
$\Om^- = \{f<0\}$
satisfy the requirements of the second part of Theorem \ref{thm:Jordan}, and the second part of the proof of the theorem shows that the decomposition $\nu = \nu^+-\nu^-$ is obtained from any such decomposition of $\Omega$. For real scalars this also gives a second proof of the first assertion:
$$ |\nu|(F) = \nu^+(F)+\nu^-(F) =  \int_F f^++f^- \ud \mu =  \int_F |f| \ud \mu.$$
\end{example}

\begin{example}\label{ex:K-val-meas} If $\mu$ is a $\K$-valued measure, then $\mu$ is absolutely continuous with respect to its variation $|\mu|$.
By the Radon--Nikod\'ym theorem, there exists $h\in L^1(\Om,|\mu|)$ such that
$$\mu(F) = \int_F h \ud|\mu|, \quad F\in \calF\!.$$
By the result of Example \ref{ex:densityfunction},
$$|\mu|(F) = \int_F |h|\ud|\mu|, \quad F\in \calF\!,$$
so $|h| = 1$ $|\mu|$-almost everywhere.
\end{example}

A measurable function $f$ is said to be {\em integrable} with respect to a $\K$-valued measure $\mu$ if it is integrable
with respect to $|\mu|$. The function $f$ is integrable with respect to a real measure $\mu$ if and only if it is integrable with respect to the measures $\mu^+$ and $\mu^-$\!, where $\mu = \mu^+-\mu^-$ is the Jordan decomposition, and $f$ is integrable with respect to a complex measure $\mu$ if and only if $f$ is integrable with respect to the real and imaginary parts of $\mu$.

The integral of an integrable function $f$ with respect to a real measure $\mu$ is defined by
$$ \int_\Om f\ud \mu := \int_\Om f\ud \mu^+ - \int_\Om f\ud \mu^-\!,$$
and the integral of an integrable function $f$ with respect to a complex measure $\mu$ by
$$ \int_\Om f\ud \mu := \int_\Om f\ud \Re\mu + i \int_\Om f\ud \Im\mu.$$

\begin{proposition}\label{prop:int-compl-meas-ineq}
 If $f$ is integrable with respect to a $\K$-valued measure $\mu$, then
 $$ \Big|\int_\Om f\ud \mu\Big| \le \int_\Om |f|\ud |\mu|.$$
\end{proposition}
\begin{proof}
 First let $f = \sum_{n=1}^N c_n \one_{F_n}$ be a simple function, with the sets $F_n\in\calF$ pairwise disjoint. Then
 $$ \Big|\int_\Om f\ud \mu\Big| = \Big|\sum_{n=1}^N c_n \mu(F_n)\Big|\le \sum_{n=1}^N |c_n| |\mu(F_n)|
\le \sum_{n=1}^N |c_n| |\mu|(F_n) =  \int_\Om |f|\ud |\mu|.$$
The general case follows from this by observing that the simple functions are dense in $L^1(\Om,|\mu|)$
and that $f_n\to f$ in $L^1(\Om,|\mu|)$ implies $\int_\Om |f_n - f|\ud\nu\to 0$ for each of the measures
$\nu\in \{(\Re\mu)^\pm,(\Im\mu)^\pm\}$.
\end{proof}

A more elegant, but less elementary, alternative definition of the integral $\int_\Om f\ud \mu$ can be given with the help of the Radon--Nikod\'ym theorem. Indeed, defining $\int_\Om f\ud \mu$ as above,
by the result of Example \ref{ex:K-val-meas} for functions $f\in L^1(\Om,|\mu|)$
we have the identity
$$ \int_\Om f\ud \mu = \int_\Om fh\ud |\mu|,$$
where $\ud \mu = h\ud|\mu|$ as in the example (note that $fh\in L^1(\Om,|\mu|)$ since $|h| = 1$ $|\mu|$-almost everywhere). This identity could be taken as an alternative definition for the integral $\int_\Om f\ud \mu$.

\subsection{Conditional Expectations}
\label{subsec:conditional-expectations}
Let
$(\Omega,\calF\!,\mu)$ be a finite measure space, and let
$\mathscr G\subseteq \calF$ be a sub-$\sigma$-algebra. If
$f\in L^1(\Omega,\calF\!,\mu)$, then the formula
$$
   \nu(G):=\int_G f\ud\mu,\quad G\in\mathscr G,
$$
defines a $\K$-valued measure $\nu$ on $(\Omega,\mathscr G)$ which is absolutely continuous with respect to the restricted measure $\mu|_{\mathscr G}$ on $(\Omega,\mathscr G)$. By the Radon--Nikod\'ym theorem there exists a unique
$g\in L^1(\Omega,\mathscr G,\mu)$ such that
$$
   \nu(G)=\int_G g\ud\mu,\quad G\in\mathscr G.
$$
This function is called the {\em conditional expectation}\index{conditional expectation}
of $f$ with respect
to $\mathscr G$ and we write
$$
   g=\E(f|\mathscr G).
$$
The following theorem characterises $\E(f|\mathscr G)$ as the unique
$\mathscr G$-measurable function having the same integrals as $f$
over all sets in $\mathscr G$.

\begin{theorem}[Conditional expectations]
\label{thm:existence-conditional-expectations}
Let $(\Omega,\calF\!,\mu)$ be a finite measure space, and let
$\mathscr G\subseteq\calF$ be a sub-$\sigma$-algebra. For every
$f\in L^1(\Omega,\calF\!,\mu)$ there exists a unique element
$\E(f|\mathscr G)\in L^1(\Omega,\mathscr G,\mu)$ satisfying
$$
   \int_G \E(f|\mathscr G)\ud\mu
   =
   \int_G f\ud\mu,
   \quad G\in\mathscr G.
$$
The mapping
$$
   \E(\,\cdot\,|\mathscr G):L^1(\Omega,\calF\!,\mu)
   \to
   L^1(\Omega,\mathscr G,\mu)
$$
is linear, positive, and contractive.
\end{theorem}

\begin{proof} As explained above, existence and uniqueness follow from the Radon--Nikod\'ym theorem. Linearity follows from the defining identity
and uniqueness, and positivity follows because a positive measure has
a nonnegative Radon--Nikod\'ym derivative.

Let $f\in L^1(\Omega,\calF,\mu)$ and put
$g:=\E(f|\mathscr G)$.
By linearity, the defining identity implies
\begin{align}\label{eq:CE-bounded-test}
\int_\Omega \varphi g\ud\mu
=\int_\Omega \varphi f\ud\mu
\end{align}
for every simple $\mathscr G$-measurable function $\varphi$.
By uniform approximation, this identity extends to every bounded
$\mathscr G$-measurable function $\varphi$.

Define the $\mathscr G$-measurable function $\theta$ by
$$
\theta(\omega):=
\begin{cases}
\overline{g(\omega)}/|g(\omega)|, & g(\omega)\ne0,\\
0, & g(\omega)=0.
\end{cases}
$$
Then $|\theta|\le1$ and $\theta g=|g|$.
Applying \eqref{eq:CE-bounded-test} with
$\varphi=\one_G\theta$, where $G\in\mathscr G$, gives
\begin{align*}
\int_G|g|\ud\mu
&=\Bigl|\int_G\theta f\ud\mu\Bigr|
\le\int_G|f|\ud\mu
=\int_G\E(|f|\,|\,\mathscr G)\ud\mu.
\end{align*}
Both $|g|$ and $\E(|f|\,|\,\mathscr G)$ are
$\mathscr G$-measurable. Since the inequality holds for every
$G\in\mathscr G$, it follows that
\begin{align*}
|\E(f|\mathscr G)|
\le\E(|f|\,|\,\mathscr G)
\quad\hbox{almost everywhere}.
\end{align*}
Taking $G=\Omega$ in the preceding integral inequality gives
$\n\E(f|\mathscr G)\n_1\le\n f\n_1$, proving contractivity.
\end{proof}

\begin{proposition}[Properties of the conditional expectation]
\label{prop:basic-properties-conditional-expectations}
Let $(\Omega,\calF\!,\mu)$ be a finite measure space, and let
$\mathscr G\subseteq\calF$ be a sub-$\sigma$-algebra. Then:
\begin{enumerate}[label={\rm(\arabic*)}, leftmargin=*]
\item\label{it:CE3} if $f\in L^1(\Omega,\calF\!,\mu)$ and
$h\in L^\infty(\Omega,\mathscr G,\mu)$, then
$$\E(hf|\mathscr G)=h\,\E(f|\mathscr G);$$

\item\label{it:CE4} if $\mathscr H\subseteq\mathscr G\subseteq\calF$ are
sub-$\sigma$-algebras, then for all
$f\in L^1(\Omega,\calF\!,\mu)$ we have
$$\E(\E(f|\mathscr G)|\mathscr H)=\E(f|\mathscr H).$$
\end{enumerate}
\end{proposition}

\begin{proof}
\ref{it:CE3}:\
The function $h\,\E(f|\mathscr G)$ is integrable and
$\mathscr G$-measurable.
For every $G\in\mathscr G$, applying
\eqref{eq:CE-bounded-test} with $\varphi=\one_Gh$ gives
$$
\int_G h\,\E(f|\mathscr G)\ud\mu
=\int_G hf\ud\mu.
$$
Thus $h\,\E(f|\mathscr G)$ satisfies the defining identity for
$\E(hf|\mathscr G)$, and the assertion follows by uniqueness.

\smallskip
\ref{it:CE4}:\
For every $G\in\mathscr H$, we have $G\in\mathscr G$ and therefore
$$
\int_G\E(\E(f|\mathscr G)|\mathscr H)\ud\mu
=\int_G\E(f|\mathscr G)\ud\mu
=\int_Gf\ud\mu.
$$
Since $\E(\E(f|\mathscr G)|\mathscr H)$ is integrable and
$\mathscr H$-measurable, uniqueness gives the asserted identity.
\end{proof}

The conditional expectation restricts to the $L^p$-setting:

\begin{proposition}[Conditional expectations on $L^p$]
\label{prop:conditional-expectations-Lp}
Let $(\Omega,\calF\!,\mu)$ be a finite measure space, and let
$\mathscr G\subseteq\calF$ be a sub-$\sigma$-algebra. For all
$1\leq p\leq\infty$, the conditional expectation restricts to a
positive contraction
$$
   \E(\,\cdot\,|\mathscr G):
   L^p(\Omega,\calF\!,\mu)
   \to
   L^p(\Omega,\mathscr G,\mu).
$$
Moreover, this operator is a projection onto
$L^p(\Omega,\mathscr G,\mu)$.
\end{proposition}

\begin{proof}
The case $p=1$ has already been proved above. For $p=\infty$, let
$M:=\|f\|_\infty$ and put $g:=\E(f|\mathscr G)$. Define
$$
  \theta(\omega):=
  \begin{cases}
    \overline{g(\omega)}/|g(\omega)|, & g(\omega)\ne0,\\
    0, & g(\omega)=0.
  \end{cases}
$$
Then $|\theta|\le1$ and $\theta g=|g|$. As in the proof of Theorem \ref{thm:existence-conditional-expectations}, for all $G\in\mathscr G$ we have
$$
  \int_G |g|\ud\mu
  =
  \Big|\int_G \theta g\ud\mu\Big|
  =
  \Big|\int_G \theta f\ud\mu\Big|
  \le
  M\mu(G).
$$
It follows that $|g|\le M$ almost everywhere, and hence
$  \|\E(f|\mathscr G)\|_\infty\le \|f\|_\infty$.

Now let $1<p<\infty$ and $\frac1p+\frac1q =1$. Put
$u:=\E(f|\mathscr G)$. We prove that $u\in L^p(\Omega,\mathscr G,\mu)$ by using
Proposition \ref{prop:Lp-via-Lq-1} on the finite measure space
$(\Omega,\mathscr G,\mu)$.

Let $h$ be a bounded $\mathscr G$-measurable function. Define the
$\mathscr G$-measurable function $\theta$ by
$$
   \theta(\omega) :=
   \overline{u(\omega)h(\omega)}/|u(\omega)h(\omega)|
$$
if $u(\omega)h(\omega)\ne0$,
and $\theta(\omega):=0$ if $u(\omega)h(\omega)=0$.
Then $\theta h$ is bounded and $\mathscr G$-measur\-able. As in the proof of Theorem \ref{thm:existence-conditional-expectations},
$$
   \int_\Omega |uh|\ud\mu
   =
   \int_\Omega \theta h u\ud\mu
   =
   \int_\Omega \theta h f\ud\mu.
$$
Therefore H\"older's inequality gives
$$
   \|uh\|_1
   =\Big|\int_\Omega \theta h f\ud\mu\Big|
   \le
   \int_\Omega |h||f|\ud\mu
   \le
   \|h\|_q\|f\|_p .
$$
Since the bounded $\mathscr G$-measurable functions are dense in
$L^q(\Omega,\mathscr G,\mu)$, Proposition
\ref{prop:Lp-via-Lq-1} implies that $u\in L^p(\Omega,\mathscr G,\mu)$
and $\|u\|_p\le \|f\|_p$.

The final assertion follows from the second part of Proposition \ref{prop:basic-properties-conditional-expectations} and the observation that every $h\in L^p(\Omega,\mathscr G,\mu)$ satisfies $\E(h|\mathscr G)=h$ by uniqueness.
\end{proof}

In the Hilbert space case $p=2$, conditional expectation has a
particularly simple interpretation. The space
$L^2(\Omega,\mathscr G,\mu)$ is a closed subspace of
$L^2(\Omega,\calF\!,\mu)$,
and, borrowing some terminology from the next chapter, we may say that
$$
   \E(\,\cdot\,|\mathscr G):
   L^2(\Omega,\calF\!,\mu)
   \to
   L^2(\Omega,\mathscr G,\mu)
$$
is the orthogonal projection onto this subspace. Indeed, if
$f\in L^2(\Omega,\calF\!,\mu)$ and
$h\in L^2(\Omega,\mathscr G,\mu)$, then, first for bounded $h$ and
then by approximation,
$$
   \int_\Omega (f-\E(f|\mathscr G))\overline h\ud\mu=0.
$$
Thus $f-\E(f|\mathscr G)$ is orthogonal to
$L^2(\Omega,\mathscr G,\mu)$.

This $L^2$-characterisation of the conditional expectation can be taken as the point of departure for an alternative construction of the
conditional expectation which avoids the use of the Radon--Nikod\'ym theorem; see Problem \ref{prob:CE-Hilbert}.

\section{Banach Lattices}\label{sec:BL}

Over the real scalar field, the classical scalar-valued Banach spaces
introduced so far in this chapter are examples of Banach lattices, a
class of Banach spaces that will be briefly discussed in the present
section. The main result, Theorem \ref{thm:latticenorm}, shows that any complete norm on a Banach lattice $X$ which is monotone with respect to the partial order of $X$ is equivalent to the given norm of $X$.

Let $(S,\le)$ be a partially ordered set and let $S'$ be a subset of $S$.
An element $x\in S$ is said to be a {\em lower bound}
for $S'$ if we have $x\le x'$ for all $x'\in S'$\!.
 Such an element is
called a {\em greatest lower bound}\index{lower!bound, greatest} for
$S'$ if $y\le x$ holds for every lower bound $y$ for $S'$\!.
Similarly, an element $x\in S$ is said to be an {\em upper bound}
 for $S'$ if we have $x'\le x$ for all $x'\in S'$\!, and such an element is called a {\em least upper bound}\index{upper!bound, least} for $S'$ if $x\le y$ holds for every upper bound $y$ for $S'$\!.
Greatest lower bounds and least upper bounds, if they exist, are unique.

\begin{definition}[Lattices]\label{def:lattice}
A partially ordered set $(S,\le)$ is called a {\em lattice}\index{lattice} if every pair of elements has a greatest lower bound and a least upper bound.
\end{definition}

The greatest lower bound and the least upper bound of the pair $\{x,y\}\subseteq S$ in a partially ordered set $S$
will be denoted by $x\wedge y$ and $x\vee y$, respectively.

\begin{definition}[Vector lattices]\label{def:vectorlattice}
A {\em vector lattice}\index{vector!lattice} is a partially ordered real vector space $(V,\le)$ with the following properties:
\begin{enumerate}[label={\rm(\roman*)}, leftmargin=*]
 \item\label{it:vectorlattice1} $(V,\le)$ is a lattice;
 \item\label{it:vectorlattice2} for all $0\le c \in\R$ and $u,v\in V$ we have $u\le v\Rightarrow cu\le cv$;
 \item\label{it:vectorlattice3} for all $u,v,w\in V$ we have $u\le v\Rightarrow u+w\le v+w$.
\end{enumerate}
\end{definition}

Let $(V,\le)$ be a vector lattice. If $u,u'\!,v,v'\in V$ satisfy $u\le v$ and $u'\le v'$\!, then $u+u' \le v+u'$ and $v+u' \le v+v'$\!. Thus,
$$[u\le v \ \hbox{ and } \ u'\le v'] \Rightarrow u+u' \le v+v'$$ by transitivity.
Also, if $u\le v$, then
$-v = u+ (-u-v) \le v+ (-u-v) = -u$, and the converse implication is obtained similarly. Thus,
\begin{align}\label{eq:orderedVS0} u\le v \Leftrightarrow (-v)\le (-u).
\end{align}

For $v\in V$ we define its {\em positive part},\index{positive!part}  {\em negative part},\index{negative part} and {\em modulus}\index{modulus!in a vector lattice}
$$ v^+:= v\vee 0, \qquad v^-:= (-v)\vee 0, \qquad |v| := v \vee (- v).$$
Furthermore, from $\pm(u+v) \le |u|+|v|$ it follows that
\begin{align}\label{eq:orderedVS5} |u+v| \le |u|+|v|.
\end{align}

The next two propositions list some slightly less trivial identities.

\begin{proposition}\label{prop:negpart-hom}
Let $(V,\le)$ be a vector lattice. Then for all $u,v\in V$ and $c\ge 0$ we have
$$c (u\wedge v) = c u\wedge c v, \qquad c (u\vee v) = c u\vee c v.$$
In particular, $
(cv)^+ = cv^+$ and $(c v)^- = c v^-$.
\end{proposition}

\begin{proof}
Fix $u,v\in V$ and $c\ge 0$. Then $u\wedge v\le u$ and $u\wedge v\le v$ imply
$c(u\wedge v)\le cu$ and $c(u\wedge v)\le cv$, and therefore
$c(u\wedge v)\le cu\wedge cv$.

Conversely, if $c=0$, there is nothing to prove. Now let $c>0$, and let $w\in V$ be such that
$w\le cu$ and $w\le cv$. Then $\frac1c\,w\le u$ and $\frac1c\,w\le v$, so
$\frac1c\,w\le u\wedge v$, and hence $w\le c(u\wedge v)$. This shows that every
lower bound of $cu$ and $cv$ is dominated by $c(u\wedge v)$, and therefore
$cu\wedge cv\le c(u\wedge v)$.

Combining both inequalities we obtain the first identity. The second identity is proved similarly.
The concluding statement is obtained by taking $u=0$ (and replacing $v$ by $-v$ in the second identity).
\end{proof}

\begin{proposition}\label{prop:orderedVS} Let $(V,\le)$ be a vector lattice. Then for all $u,v,w\in V$ we have:
 \begin{enumerate}[label={\rm(\arabic*)}, leftmargin=*]
\item\label{it:ordered-VS1} $(-u)\wedge (-v) = -(u\vee v)$;
\item\label{it:ordered-VS2} $u+ (v\vee w) = (u+v)\vee (u+w)$;
\item\label{it:ordered-VS3} $u+v = u\wedge v + u\vee v$;
\item\label{it:ordered-VS4} $v = v^+-v^-$;
\item\label{it:ordered-VS5} $|v| = v^++v^-$\!.
 \end{enumerate}
\end{proposition}

The representation of $v$ as the difference of two positive elements in \ref{it:ordered-VS4}
is minimal in a sense explained in Problem  \ref{prob:vplusvminus}.

\begin{proof}
Let $u,v,w\in V$.

\smallskip
\ref{it:ordered-VS1}:\
 We have $u \le u\vee v$, so $ -(u\vee v)\le -u$ by \eqref{eq:orderedVS0}. The inequality $ -(u\vee v)\le -v$ is proved similarly.
 It follows that $-(u\vee v)$ is a lower bound for
 $\{-u,-v\}$. To prove that it is the greatest lower bound, we must show that if $w\le -u$ and $w\le -v$, then $w \le -(u\vee v)$.
 This follows by noting that $u\le -w$ and $v\le -w$, so $-w$ is an upper bound for $\{u,v\}$ and therefore $u\vee v\le -w$.
 By \eqref{eq:orderedVS0} this implies $w \le -(u\vee v)$ as required.

\smallskip
\ref{it:ordered-VS2}:\ We have $u+v \le u+ (v\vee w)$ and $u+w \le u+ (v\vee w)$, so
$u+ (v\vee w)$ is an upper bound for $\{u+v,u+w\}$. To prove that it is the least upper bound we must show that if
$u+v\le x$ and $u+w \le x$, then $u+ (v\vee w)\le x$. But $v\le x-u$ and $w\le x-u$ imply
$v\vee w \le x-u$ and therefore $u+ (v\vee w)\le x$ as desired.

\smallskip
\ref{it:ordered-VS3}:\
By \ref{it:ordered-VS2} and \ref{it:ordered-VS1} we have
$ (u\vee v)-u-v =(-v)\vee(-u) =-(u\wedge v).$
Rearranging gives the desired identity.

\smallskip
\ref{it:ordered-VS4}:\ Taking $u=0$ in \ref{it:ordered-VS3} and using \ref{it:ordered-VS1}
we obtain $v = 0\wedge v + 0\vee v =
-(0\vee (-v)) + v^+ = v^+ - v^-$\!.

\smallskip
\ref{it:ordered-VS5}:\ By \ref{it:ordered-VS2} and Proposition \ref{prop:negpart-hom}, $|v| =  v\vee (-v)$ implies $|v|-v = 0\vee (-2v) = (2v)^-
= 2 v^-$\!. It follows that $|v| = v + 2v^- = v^+- v^- +2v^- = v^++v^-$\!.
\end{proof}

\begin{definition}[Normed vector lattices]\label{def:normedVL}
A {\em normed vector lattice} is a triple $(X,\n\cdot\n,\le)$ with the following properties:
\begin{enumerate}[label={\rm(\roman*)}, leftmargin=*]
 \item\label{it:normedVL1} the pair $(X,\n\cdot\n)$ is a real
 normed space;
 \item\label{it:normedVL2} the pair $(X,\le)$ is a vector lattice;
 \item\label{it:normedVL3} for all $x,y\in X$ we have $|x|\le |y| \Rightarrow \n x\n \le \n y\n$.
\end{enumerate}
\end{definition}

In any normed vector lattice,
as an immediate consequence of \ref{it:normedVL3} we have
\begin{align}\label{eq:pos-mod} \bigl\n \, |x|\,\bigr\n = \n x\n.
\end{align}
Moreover, the lattice operations $x\mapsto x^+$\!, $x\mapsto x^-$\!,
$x\mapsto |x|$ are continuous. Indeed, by \eqref{eq:orderedVS5} we have
$|x|-|y|\le |x-y|$ and $|y|-|x|\le |x-y|$, hence $|\,|x|-|y|\,|\le |x-y|$.
By \ref{it:normedVL3} and \eqref{eq:pos-mod},
$$
  \bigl\n\, |x|-|y|\,\bigr\n
  \le
  \bigl\n\, |x-y|\,\bigr\n
  =
  \n x-y\n .
$$
This gives the continuity of $x\mapsto |x|$. The other two assertions follow from this by noting that $x^+ = \frac12(x+|x|)$ and
$x^- = x^+-x$. As a consequence we note that the {\em positive cone}\index{positive!cone}
$$X^+:= \{x\in X:\, x\ge 0\}$$
is closed.

\begin{definition}[Banach lattices]
 A {\em Banach lattice}\index{Banach!lattice} is a complete normed vector lattice.
\end{definition}

The spaces $c_0$, $\ell^p$\!, $C(K)$, and $L^p(\Om)$ with
$1\le p\le \infty$, all over the real scalar field, are Banach lattices
with respect to their natural pointwise ordering. The real Banach space of
real-valued measures on $(\Om,\calF)$, endowed with the variation norm,
is a Banach lattice with respect to the ordering given by declaring
$\mu\le \nu$ if and only if the measure $\nu-\mu$ is nonnegative; the
greatest lower bound and least upper bound of $\mu$ and $\nu$ are given by
\begin{align*}
 \mu \wedge \nu & = \mu - (\mu-\nu)^+ = \nu - (\nu-\mu)^+, \\
 \mu\vee \nu & = \mu+(\nu-\mu)^+ = \nu+(\mu-\nu)^+,
\end{align*}
respectively (cf. Problem \ref{prob:lub}); here, $(\nu-\mu)^+$ is defined
as in Theorem \ref{thm:Jordan}. Alternatively, one may use the analogues for
real-valued measures of the formulas appearing in Theorem
\ref{thm:dual-BL}.

\begin{definition} Let $V$ and $W$ be vector lattices.
A linear operator $T:V\to W$ is said to be {\em positivity preserving}\index{positivity preserving} if $v\ge 0$ implies $Tv\ge 0$.
\end{definition}

If $T:V\to W$ is positivity preserving, then
\begin{align}\label{eq:Tpos-mod}
|Tv| \le T|v|, \quad v\in V.
\end{align}
Indeed, from $-|v|\le v\le |v|$ and positivity we obtain $-T|v|\le Tv\le T|v|$, and therefore $|Tv|\le T|v|$.

\begin{theorem}\label{thm:pos-op}
 Let $X$ and $Y$ be Banach lattices. Every positivity-preserving linear operator $T:X\to Y$ is bounded.
\end{theorem}
\begin{proof}
Arguing by contradiction, suppose that $T$ is not bounded.
Then for all $n\ge 1$ there is a norm-one vector $x_n\in X$ such that $\n Tx_n\n\ge n^3$\!.
By \eqref{eq:pos-mod} and \eqref{eq:Tpos-mod}, upon replacing $x_n$ by $|x_n|$ we may assume that $x_n\ge 0$ for all $n\ge 1$.
In view of $ \sum_{n\ge 1} {\n x_n\n}/{n^2} < \infty$ the sum $\sum_{n\ge 1} {x_n}/{n^2}$ converges in $X$.

Clearly, for all $n\ge1$ and $N\ge n$ we have
$x_n/n^2\le\sum_{m=1}^N x_m/m^2$.
By the closedness of the positive cone $X^+$, upon letting $N\to \infty$ we obtain
${x_n}/{n^2} \le \sum_{m\ge 1} x_m/m^2$ for all $n\ge 1.$
Hence, for all $n\ge 1$,
$$ n \le\frac1{n^2} \n T x_n \n \le \Big\n T  \sum_{m\ge 1} \frac{x_m}{m^2}\Big\n. $$
This contradiction completes the proof.
\end{proof}

This theorem has an interesting consequence:

\begin{theorem}\label{thm:latticenorm}
Any two norms which turn a vector lattice into a Banach lattice are equivalent.
\end{theorem}
\begin{proof}
Suppose the vector lattice $(X,\le)$ is a Banach lattice with respect to the norms $\n\cdot\n$ and $\n\cdot\n'$\!. Then by Theorem \ref{thm:pos-op} the identity mapping from $(X,\n\cdot\n)$ to $(X,\n\cdot\n')$ and its inverse, both of which are positivity preserving, are bounded.
\end{proof}

\begin{problems}

\item\label{prob:ellp}
Let $1\le p <\infty$.
Show that if $a\in \ell^p$\!, then $a\in \ell^q$ for all $p\le q\le \infty$ and
$$\n a\n_\infty\le \n a\n_q\le \n a\n_p, \qquad \lim_{q\to\infty} \n a\n_q = \n a\n_\infty.$$
\noindent{\em Hint:}\ First show that
it suffices to consider sequences $a = (a_n)_{n\ge 1}$ satisfying $|a_n|\le 1$ for all $n\ge 1$.

\item
Show that $c_0$ and $\ell^p$ with $1\le p<\infty$ are separable, but $\ell^\infty$ is not.

\noindent {\em Hint:}\  $\ell^\infty$ contains an uncountable family $(a^{(i)})_{i\in I}$ such that $\n a^{(i)} - a^{(i')}\n = 1$ for all $i\not=i'$\!. Prove that a normed space $X$ containing such a family is nonseparable.

\item
Let $K$ be a compact topological space and $X$ be a Banach space.
Prove that the space $C(K;X)$\index{$C1b$@$C(K;X)$} of all continuous functions $f:K\to X$ is a Banach space with respect to the supremum norm $\n f\n_\infty := \sup_{x\in K}\n f(x)\n$.

\item\label{prob:Ck}
Let $D$ be a bounded open subset of $\R^d\!$.
By $C^k(\overline D)$ we denote the space of functions $f:\overline D\to \K$ that are
$k$ times continuously differentiable on $D$ and all of whose partial derivatives
$\partial^\alpha f$ extend continuously to $\overline D$ for all multi-indices  $\alpha = (\alpha_1,\dots,\alpha_d)\in \N^d$ satisfying $|\alpha|:= \alpha_1+\cdots+\alpha_d\le k$. Here, $\partial^{\alpha}  := \partial_1^{\alpha_1}\circ\cdots\circ \partial_d^{\alpha_d}$, where $\partial_j$ is the partial derivative in the $j$th direction.
Prove that $C^k(\ov{D})$ is a Banach space with respect to the norm
$$ \n f\n_{C^k(\overline D)} := \max_{|\alpha|\le k} \n \partial^\alpha f\n_\infty.$$

\item\label{prob:Dini}
Let $K$ be a compact metric space and let $(f_n)_{n\ge 1}$ be a
sequence of real-valued functions in $C(K)$ satisfying
$$
f_1(x)\ge f_2(x)\ge\cdots\ge0
$$
for all $x\in K$. Our aim is to prove that if $\limn f_n(x) =0$ for all $x\in K$, then $\limn f_n = 0$ uniformly on $K$.
This result, known as {\em Dini's theorem}\index{theorem!Dini},
provides one of the rare instances where pointwise convergence implies uniform convergence.

\begin{enumerate}[\rm(a), leftmargin=*]
  \item Reasoning by contradiction, show that if the result is false, then there exists an $\eps>0$, a sequence $(x_n)_{n\ge 1}$ in $K$, and an $x\in K$, such that $\limn x_n = x$ and $f_{n}(x_{n})\ge \frac12\e$ for all $n\ge 1$.
  \item
  Using also the facts that $f_n(x)\downarrow0$ as $n\to\infty$ and
$f_n(x_n)\le f_m(x_n)$ when $n\ge m$, show that this leads to a contradiction.
\end{enumerate}

\item
Find a sequence $(f_n)_{n\ge 1}$ of real-valued functions in
$C_{\rm b}(0,1)$ such that we have $$0 \leq f_{n+1}(x) \leq f_n(x)\le 1$$ for all $x\in (0,1)$ and $n\ge 1$ and $f_n(x)\downarrow 0$ for all $x\in (0,1)$, but $\n f_n\n_{\infty} \not \rightarrow 0$. Compare this with Dini's theorem (Problem \ref{prob:Dini}).

\item\label{prob:norm-dense-subspace-2}
Consider the vector space $P[0,1]$ of all polynomials on $[0,1]$.
\begin{enumerate}[\rm(a), leftmargin=*]
  \item Show that
  $$ \Bigl\n t\mapsto \sum_{n=0}^N c_n t^n \Bigr\n := \max\Bigl\{ \Bigl\n t\mapsto\sum_{k=0}^{\lfloor\frac{N}{2}\rfloor} c_{2k} t^{2k}\Bigr\n_\infty,    2\Bigl\n t\mapsto\sum_{k=1}^{\lceil\frac{N}{2}\rceil} c_{2k-1} t^{2k-1}\Bigr\n_\infty  \Bigr\}$$
  defines a norm on $P[0,1]$. Here, $\lfloor y \rfloor$ is the greatest integer $n\le y$ and $\lceil y \rceil$ is the least integer $n\ge y$.
  \item Show that the functions $\one, t^2\!, t^4\!, \dots$ span a subspace of $P[0,1]$ which is dense with respect to the supremum norm.
  \item Conclude that two different norms on a vector space may agree on a subspace which is dense with respect to one of these norms.
\end{enumerate}

\item\label{prob:no-fixedpoint}
Let $\K = \R$. Consider the set
$$F:= \big\{f\in C[0,1]:\, 0\le f\le \one, \ f(0)=0, \ f(1)=1\big\}$$
and the linear mapping $T: C[0,1]\to C[0,1]$ defined by
$$(Tf)(t) = tf(t), \quad t\in [0,1].$$
\begin{enumerate}[\rm(a), leftmargin=*]
  \item Show that $F$ is bounded, convex, and closed in $C[0,1]$.
  \item Show that $T$ maps $F$ into $F$ and satisfies $\|Tf - Tg\|_\infty < \|f-g\|_\infty$ for all $f,g\in F$, $f\not=g$.
  \item Show that $T$ has no fixed point in $F$.
  \item Compare this result with the Banach fixed point theorem.
\end{enumerate}

\item\label{prob:Schauder}
Prove the {\em Schauder fixed point theorem}:\index{theorem!Schauder fixed point}
Let $C$ be a nonempty closed convex subset of a Banach space $X$ and
let $T:C\to C$ be continuous. If $T(C)$ is relatively compact, then
$T$ has a fixed point.

\noindent{\em Hint:}\
Put $K:=\ov{T(C)}$. Note that $K$ is compact and $K\subseteq C$.
Given $\eps>0$, choose
$y_1,\dots,y_N\in K$ such that the balls
$B(y_1;\eps),\dots,B(y_N;\eps)$ cover $K$. For $y\in K$ put
$$
\alpha_j(y):=(\eps-\n y-y_j\n)^+,
\qquad
\beta_j(y):=
\frac{\alpha_j(y)}{\sum_{k=1}^N\alpha_k(y)}.
$$
Show that $$F_\eps(y):=\sum_{j=1}^N\beta_j(y)y_j, \quad y\in K,$$
defines a continuous mapping
$F_\eps:K\to\co\{y_1,\dots,y_N\}$ satisfying
$\n F_\eps(y)-y\n<\eps$ for all $y\in K$.

For $y\in K$ and $\la\in \Delta_{N-1}$ define
$$
Q_\eps(y):=(\beta_1(y),\dots,\beta_N(y)), \qquad
A(\la_1,\dots,\la_N):=\sum_{j=1}^N\la_jy_j.
$$
Applying the Brouwer fixed point theorem (Theorem \ref{thm:Brouwer}) to
the continuous mapping
$Q_\eps\circ T\circ A:\Delta_{N-1}\to\Delta_{N-1}$, deduce that there exists an $x_\eps\in C$ such that
$x_\eps=F_\eps(Tx_\eps)$, and therefore
$$\n Tx_\eps - x_\eps\n<\eps.$$
Finally, let $\eps\downarrow0$ and use the compactness of $K$.

\item\label{prob:local-existence}
Give an alternative proof of the Peano existence theorem
(Theorem \ref{thm:DV-cont}) using the Schauder fixed point theorem.

\noindent{\em Hint:}\
Put
$M:=\sup_{(t,x)\in[0,T]\times\ov B(u_0;1)}|f(t,x)|$ and $\delta:=\min\{T,1/M\}$,
setting $\delta:=T$ if $M=0$. On the closed convex subset
$$
C:=\{u\in C([0,\delta];\K^d):\,u(0)=u_0,\
\n u-u_0\n_\infty\le1\}
$$
consider
$$
(\Phi u)(t):=u_0+\int_0^t f(s,u(s))\ud s.
$$
Show that $\Phi:C\to C$ is continuous and that $\Phi(C)$ is relatively
compact by the Arzel\`a--Ascoli theorem. Then apply Problem \ref{prob:Schauder}.

\item\label{prob:Perron-Frobenius}
Let $K$ be a nonempty compact metric space, let $\mu$ be a Borel probability
measure on $K$ with full support, in the sense that $\mu(O)>0$ for every nonempty open subset $O\subseteq K$, and let $k\in C(K\times K)$ be
nonnegative. Define
$$
Tf(x):=\int_K k(x,y)f(y)\ud\mu(y),\quad f\in C(K).
$$

\begin{enumerate}[\rm(a), leftmargin=*]
\item\label{it:Perron-Frobenius1} Suppose that $k(x,y)>0$ for all $x,y\in K$. Prove that
$Tf = \la f$ for some $\la>0$ and a strictly positive function $f\in C(K)$.

\noindent{\em Hint:}\
Apply the Schauder fixed point theorem to
$$
C:=\Bigl\{f\in C(K):\,f\ge0,\ \int_Kf\ud\mu=1\Bigr\}
$$
and the mapping
$$
Ff:=\frac{Tf}{\int_KTf\ud\mu}.
$$
Use the Arzel\`a--Ascoli theorem to prove that $F(C)$ is relatively
compact.

\item\label{it:Perron-Frobenius2} Under the assumption of part \ref{it:Perron-Frobenius1}, show that if
$f,g\in C(K)$ are nonzero nonnegative eigenfunctions of $T$,
with eigenvalues $\la$ and $\la'$, then $f=cg$ for some
$c>0$ and $\la=\la'$.

\noindent{\em Hint:}\
For two positive eigenfunctions $f$ and $g$, consider
$\max_{x\in K}f(x)/g(x)$.
%
\end{enumerate}

\item
Let $(M,d)$ be a nonempty metric space and fix $x_0\in M$. For $x\in M$
define the {\em Kuratowski embedding}\index{Kuratowski embedding} $Jx\in C_{\rm b}(M)$ by
$$
(Jx)(y):=d(x,y)-d(x_0,y),\quad y\in M.
$$
Show that
$$
\n Jx-Jx'\n_\infty=d(x,x'),\quad x,x'\in M.
$$
Thus every metric space is isometric to a subset of a Banach space.

\item\label{prob:no-Cc-completemetric}
Show that every compactly supported continuous function on the separable complete metric space $\ell^2$ is zero.

\item
  Show that for $1\le p<\infty$ the space $L^p(0,1)$ is the completion of $C[0,1]$ with respect to the norm
  $$\n f\n_p = \Bigl(\int_0^1 |f(t)|^p \ud t \Bigr)^{1/p}\!,$$
  and that $C[0,1]$ can be identified with a proper closed subspace of $L^\infty(0,1)$.

\item\label{prob:Lp-subseq1}
Let $(\Om,\calF\!,\mu)$ be a measure space and let the functions $f_n:\Om\to \K$ and $f:\Om\to \K$ be bounded and measurable.
Show that if $f_n\to f$ in $L^\infty(\Om)$, then there is a $\mu$-null set $N$ such that
$\limn \sup_{\om\in\complement N}|f_n(\om)- f(\om)| = 0$.

\item\label{prob:Lp-subseq2}
Show that passing to a subsequence is necessary for Corollary \ref{cor:Lp-ae-subseq} to be true.

\noindent{\em Hint:}\ Consider $\mathbb{T} = \{z\in \C:\, |z|=1\}$ with its Lebesgue measure. Let $t_n:= \sum_{m=1}^n \frac1m$ and consider the indicator functions of the arcs $I_n := \{e^{2\pi i t}:\, t\in (t_n,t_{n+1})\}$.

\item
Let $\K=\R$ and consider the set
$$
S:=\{f\in L^1(0,1):\, f(t)\ge 0
\text{ for almost all }t\in(0,1)\}.
$$
\begin{enumerate}[\rm(a), leftmargin=*]
  \item Determine whether $S$ is a closed subset of $L^1(0,1)$.
  \item Characterise the functions belonging to the interior of $S$.
\end{enumerate}
Consider the set
$$
S':=\{f\in L^1(0,1):\, f(t)>0
\text{ for almost all }t\in(0,1)\}.
$$
\begin{enumerate}[\rm(a), leftmargin=*, resume]
  \item Determine whether $S'$ is an open subset of $L^1(0,1)$.
  \item Characterise the functions belonging to the closure of $S'$.
\end{enumerate}

\item
Let $(\Om,\calF\!,\mu)$ be a finite measure space and let $1\le p\le \infty$. Prove that
$$ \nn f \nn_p:= \Big| \int_\Om f\ud \mu\Big| + \Big\n f - \Bigl(\int_\Om f\ud \mu\Bigr)\one\Big\n_p$$
defines an equivalent norm on $L^p(\Om)$. Here, $\n \cdot\n_p$ is the usual norm on $L^p(\Om)$.

\item\label{prob:intersect-sum}
Let  $(\Om,\calF\!,\mu)$ be a measure space and fix $1\le p,q\le \infty$.
\begin{enumerate}[\rm(a), leftmargin=*]
  \item Prove that  $L^p(\Om)\cap L^q(\Om)$ is a Banach space with respect to the norm
  $$\n f\n_{L^p(\Om)\cap L^q(\Om)}:= \max\bigl\{\n f\n_{L^p(\Om)},\, \n f\n_{L^q(\Om)}\bigr\}.$$
  \item Prove that  $L^p(\Om) + L^q(\Om) = \{g +h: \, g\in L^p(\Om),\, h\in L^q(\Om)\}$ is a Banach space with respect to the norm
  \begin{align*}
  \, & \n f\n_{L^p(\Om) + L^q(\Om)}
  \\ & \qquad := \inf\bigl\{\n g\n_{L^p(\Om)} + \n h\n_{L^q(\Om)}:\, f = g + h, \ g\in L^p(\Om),\, h\in L^q(\Om)\bigr\}.
  \end{align*}

  \noindent{\em Hint:}\ For the proof of completeness use Proposition \ref{prop:competeness-series}.
  \item Prove that if $1\le p \le r \le q \le \infty$, then
  $$L^p(\Om)\cap L^q(\Om) \subseteq L^r(\Om) \subseteq L^p(\Om)+ L^q(\Om)$$
  and that the inclusion mappings are continuous.

  \noindent{\em Hint:}\ Write $f = \one_{\{|f|\le 1\}}f + \one_{\{|f|>1\}}f$.

  \item Prove that if $1\le p \le r \le q \le \infty$ and $0\le \theta\le 1$ are such that
  $\frac1{r} = \frac{1-\theta}{p}+
  \frac\theta{q}$, then for all $f\in L^{p}(\Omega,\mu)\cap L^{q}(\Omega,\mu)$ we have
  $$ \| f \|_{r} \le \| f\|_{p}^{1-\theta}\| f\|_{q}^\theta.
  $$
  \noindent {\em Hint:}\ Use H\"older's inequality with suitable exponents.
\end{enumerate}

\item
Let $(\Om,\calF\!,\mu)$ be a measure space, let $1\le p\le \infty$, and suppose that $(f_n)_{n\ge 1}$ is a bounded sequence in $L^p(\Om)$ converging to a measurable function $f$ $\mu$-almost everywhere.
\begin{enumerate}[\rm(a), leftmargin=*]
  \item\label{it:LpFatou1} Using Fatou's lemma for $1\le p<\infty$ (and the essential supremum for $p=\infty$), show that $f\in L^p(\Om)$ and $\n f\n_p \le \liminf_{n\to\infty} \n f_n\n_p$.
\end{enumerate}
In addition to the above assumptions, assume now that $\mu(\Om)<\infty$ and $1<p\le \infty$.
\begin{enumerate}[\rm(a), leftmargin=*]\setcounter{enumii}{1}
  \item\label{it:LpFatou2} Show that $\limn f_n = f$ in $L^1(\Om)$.

 \noindent{\em Hint:}\ First show that for every $\eps>0$ there exists an $r\ge 0$ such that $$\sup_{n\ge 1} \int_{\{|f_n| >r\}}|f_n|\ud\mu <\eps.$$

  \item\label{it:LpFatou3} Do we also have $\limn f_n = f$ in $L^p(\Om)$?
\end{enumerate}

\item\label{prob:Jensen}
Let $(\Om,\calF\!,\mu)$ be a probability space and let
$\phi:\R\to\R$ be {\em convex}\index{convex!function}, that is,
$$
\phi((1-\la)x+\la y)
\le(1-\la)\phi(x)+\la\phi(y),
\quad x,y\in\R,\ 0\le\la\le1.
$$
Prove {\em Jensen's inequality}\/:\index{inequality!Jensen}
if a real-valued function $f\in L^1(\Om)$ has the property that
$\phi\circ f$ is integrable, then
$$
 \phi\Big(\int_{\Om} f\ud \mu\Big)\le \int_{\Om} \phi\circ f\ud \mu.
$$

\noindent{\em Hint:}\ A real-valued convex function $\phi$ is the pointwise supremum of its affine minorants $\psi(x)=ax+b$.

\item
Using Young's inequality, show that if $f\in L^p(\R^d)$, $g\in L^q(\R^d)$, $h\in L^r(\R^d)$ with $1\le p,q,r\le \infty$ such that $\frac1p+\frac1q+\frac1r = 2$, then
$$ \int_{\R^d}\int_{\R^d}|f(x)g(x-y)h(y)|\ud x\ud y \le \n f\n_p \n g\n_q\n h \n_r.$$

\item\label{prob:LpG-closed}
Let $(\Om,\calF\!,\mu)$ be a measure space and let $\mathscr{G}\subseteq\calF$ be a sub-$\sigma$-algebra.
Show that for all $1\le p\le \infty$ the set $L^p(\Om,\mathscr{G})$ consisting of all $f\in L^p(\Om)$ that are equal $\mu$-almost everywhere to a $\mathscr{G}$-measurable function
is a closed subspace of $L^p(\Om)$.

\item\label{prob:CI-unitinterval}
For a finite measurable partition $\pi$ of the unit interval $(0,1)$ and functions $f\in L^1(0,1)$, compute $\E\iprod{f}{\calF_\pi}$, where $\calF_\pi$ is the $\sigma$-algebra generated by $\pi$.

\item\label{prob:Young-mult}
On $\R_+ = (0, \infty)$ we consider the Borel measure $\mu$ given by
\begin{equation*}
\mu(B) := \int_B \frac{1}{t} \ud t, \quad B \in \mathscr{B}(\R_+).
\end{equation*}
\begin{enumerate}[\rm(a), leftmargin=*]
  \item Show that for all $B\in \mathscr{B}(\R_+)$ and $s \in \R_+$ we have
  \begin{align*}
  \mu(B) = \mu(sB), \qquad \mu(B) &= \mu(B^{-1}),
  \end{align*}
  where $sB :=\{st: \,t \in B\}$ and $B^{-1}:= \{t^{-1}:\,t \in B\}$.
  \item For $h \in L^1(\R_+,\dtt)$ and $s \in \R_+$ show that
  \begin{equation*}
  \int_0^\infty h(st)\dtt = \int_0^\infty h(t)\dtt = \int_0^\infty h(t^{-1})\dtt.
  \end{equation*}
\end{enumerate}
Fix $1\le p< \infty$ and assume that $f \in L^p(\R_+,\dtt)$ and $g \in L^1(\R_+,\dtt)$. Define the {\em multiplicative convolution}\index{convolution!multiplicative}
\begin{equation*}
  f\diamond g(t):= \int_{0}^\infty f({t}/{s})g(s) \frac{\!\ud s}{s}, \quad t \in \R_+.
\end{equation*}
\begin{enumerate}[(a)]\setcounter{enumii}{2}
\item Show that the multiplicative convolution is well defined for almost all $t \in \R_+$ and that the following analogue of Young's inequality holds:
  \begin{equation*}
  \n{f\diamond g}\n_{L^p(\R_+,\dtt)} \leq \n{f}\n_{L^p(\R_+,\dtt)} \n{g}\n_{{L^1(\R_+,\dtt)}}.
  \end{equation*}
  \item Show that the preceding result also holds if
  $f\in L^1(\R_+,\dtt)$ and $g \in L^p(\R_+,\dtt)$ and that, with the same notation, $f \diamond g = g \diamond f$.
\end{enumerate}

\item\label{prob:Schur-estimate}
Let $k: (0,1)\times (0,1) \to \K$ be measurable and suppose that
\begin{align*}
 A :=  \esssup_{y \in (0,1)}\int_0^1 |{k(x,y)}|\ud x <\infty,\qquad
 B :=  \esssup_{x \in (0,1)}\int_0^1 |{k(x,y)}|\ud y <\infty.
\end{align*}
Let $1\le p\le \infty$ and define, for $f \in L^p(0,1)$,
\begin{equation*}
  T_kf(x) := \int_0^1 k(x,y)f(y)\ud y, \quad x \in (0,1).
\end{equation*}
Use H\"older's inequality to show that $T_k: L^p(0,1) \to L^p(0,1)$ is a well-defined linear operator which satisfies
the so-called {\em Schur estimate}\index{Schur!estimate}
\begin{equation*}
 \|T_k f\|_p \le A^{1/p}B^{1-1/p}\|f\|_p, \quad  f\in L^p(0,1).
\end{equation*}

\item\label{prob:LpX}
Let $(\Om,\calF\!,\mu)$ be a measure space and let $X$ be a Banach space.
 For $1\le p \le \infty$ we denote by $L^p(\Om;X)$\index{$L^p(\Om;X)$} the space of all
(equivalence classes of) strongly measurable functions $f:\Om\to X$ for which $\om\mapsto \n f(\om)\n$ belongs to $L^p(\Om)$.
\begin{enumerate}[\rm(a), leftmargin=*]
  \item Prove that $L^p(\Om;X)$ is a Banach space with respect to the norm
   $$ \n f\n_p := \big\n \omega\mapsto\n f(\omega)\n\big\n_{L^p(\Om)}.$$
\end{enumerate}
By $L^p(\Om)\otimes X$\index{$L^p(\Om)\otimes X$} we denote the vector space obtained as the linear span in $L^p(\Om;X)$ of the set of all functions of the form $f\otimes x$ (cf. \eqref{eq:fotx}) with $f\in L^p(\Om)$ and $x\in X$.
\begin{enumerate}[\rm(a), leftmargin=*]\setcounter{enumii}{1}
  \item Show that if $1\le p<\infty$, then $L^p(\Om)\otimes X$ is dense in $L^p(\Om;X)$.
\end{enumerate}

\item\label{prob:LpX-pos}\index{extension!vector-valued}
Let $(\Om,\calF\!,\mu)$ be a measure space and $1\le p<\infty$.
Let $T$ be a bounded operator on $L^p(\Om)$ and let $X$ be a nonzero Banach space.
Consider the linear mapping
from $L^p(\Om)\otimes X$ into itself defined by $$(T\otimes I) (f\otimes x) :=  (Tf)\otimes x.$$
\begin{enumerate}[\rm(a), leftmargin=*]
  \item Show that the operator $T\otimes I$ is well defined.
  \item Prove that if $T$ is a positivity-preserving operator,
  then $T \otimes I$ admits a unique extension to a bounded operator on $L^p(\Om;X)$ of norm $\n T\otimes I\n = \n T\n$.
\end{enumerate}

\item\label{prob:cont-Mink}
Let $(\Om,\calF\!,\mu)$ and $(\Om'\!,\calF'\!,\mu')$ be measure spaces and let $1\le p\le q< \infty$.
\begin{enumerate}[\rm(a), leftmargin=*]
  \item
  Show that the identity mapping on the linear span of the functions
  $(\one_A\otimes \one_B)(\om,\om'):=  \one_{A}(\om)\one_{B}(\om')$ with $\mu(A)<\infty$ and $\mu'(B)<\infty$ extends uniquely to a contraction operator from $L^p(\Om;L^q(\Om'))$ into $L^q(\Om';L^p(\Om)).$

  \noindent{\em Hint:}\ Use \eqref{eq:triangle-Bochner}.

  \item Deduce that if $(\Om,\calF\!,\mu)$ and $(\Om'\!,\calF'\!,\mu')$ are $\sigma$-finite and $f: \Om\times\Om' \to \K$ is a measurable function, then the {\em continuous Minkowski inequality} holds:\index{inequality!Minkowski, continuous}
  \begin{align*} \ &  \Bigl(\int_{\Om'} \Bigl(\int_{\Om} \!|f(\om,\om')|^p\ud \mu(\om)\Bigr)^{q/p}\ud\mu'(\om')\Bigr)^{\!1/q}
  \\ & \qquad \le \Bigl(\int_{\Om} \Bigl(\int_{\Om'} \!|f(\om,\om')|^q\ud \mu'(\om')\Bigr)^{p/q}\ud\mu(\om)\Bigr)^{\!1/p}\!\!.
  \end{align*}
\end{enumerate}

\item\label{prob:Lebesgue-set}
Let $f\in L_{\rm loc}^1(\R^d)$ be given.
\begin{enumerate}[\rm(a), leftmargin=*]
\item  Show that for all $c\in \K$ there exists a Borel null set $N_c\subseteq \R^d$ such that for all $x\in\complement N_c$ we have (with $B$ ranging over all balls containing $x$)
\begin{align}\label{eq:Leb-set} \lim_{\substack{B\owns x\\ |B|\to 0}} \frac{1}{|B|}\int_B  |f(y)-c| \ud  y = |f(x)-c|.\end{align}
\item\label{it:Leb-set} Show that
there exists a Borel set $L\subseteq\R^d$ with $|\complement L|=0$ such that \eqref{eq:Leb-set} holds for all $x\in L$ and $c\in \K$.

\noindent {\em Hint:}\ Consider $\bigcup_{n\ge 1} N_{c_n}$ with $(c_n)_{n\ge 1}$ a dense sequence in $\K$.
\end{enumerate}
A set $L$ with the properties of part \ref{it:Leb-set} is called a {\em Lebesgue set}\index{Lebesgue!set} for $f$.

\item
Prove the following one-sided version of the Lebesgue differentiation theorem for $d=1$:
 If $x$ is a Lebesgue  point of $f\in L^1_{\rm loc}(\R)$, then
 $$\lim_{h\downarrow 0} \frac{1}{h}\int_x^{x+h} |f(y) - f(x)| \ud y =
 \lim_{h\downarrow 0} \frac{1}{h}\int_{x-h}^{x} |f(y) - f(x)| \ud y =0.$$

\item\label{prob:compact-c0}
Show that a subset $K$ of $c_0$ is relatively compact if and only if there is a $y\in c_0$ such
that for all $x\in K$ and $n\ge 1$ we have $ |x_n| \le |y_n|.$

\item
Let $1\le p<\infty$. Show that a bounded subset $S$ of $\ell^p$ is relatively compact if and only if for every $\eps>0$ there exists
an $N\ge 1$ such that $$\sup_{x\in S} \sum_{n\ge N} |x_n|^p \le \eps^p\!.$$

\item
Show that a $\K$-valued measure $\nu$ is absolutely continuous with respect to a measure $\mu$ if and only if for every $\eps>0$ there exists a $\delta>0$ such that $|\nu(F)|<\eps$ whenever $F\in\calF$ satisfies $\mu(F)<\delta$.

\item\label{prob:Radon-closed}
A $\K$-valued measure $\mu$ on a Hausdorff topological space $X$ is said to be
{\em regular} if its variation $|\mu|$ is regular (see Definition \ref{def:regular}).
Prove that the set $M_{\rm r}(X)$ of all regular $\K$-valued Borel measures on $X$ is a closed subspace of $M(X)$.

\item
Let ${\rm Lip}[0,1]$\index{$Lip$@${\rm Lip}[0,1]$} be the vector space of all functions $f:[0,1]\to \K$ for which
$$ \n f \n_{{\rm Lip}[0,1]}:= |f(0)| + \sup_{\substack{0\le x,y\le 1\\ x\not = y}} \Big| \frac{f(x) - f(y)}{x-y}\Bigr| $$
is finite. Show that ${\rm Lip}[0,1]$ is a Banach space with respect to the norm $\n\cdot \n_{{\rm Lip}[0,1]}$.

\item\label{prob:abscont}
A function $f:[0,1]\to \K$ is said to be
{\em absolutely continuous}\index{absolutely!continuous, function} if
for every $\eps>0$ there exists a $\delta>0$ such that, whenever
\([a_1,b_1],\ldots,[a_N,b_N]\) are finitely many pairwise disjoint
intervals in \([0,1]\) with \(a_n\le b_n\) and
$\sum_{n=1}^N (b_n-a_n)<\delta$,
one has
\[
  \sum_{n=1}^N |f(b_n)-f(a_n)|<\eps .
\]
Show that a function $f:[0,1]\to\K$ is absolutely continuous and satisfies $f(0)=0$
if and only if there exists a function $g\in L^1(0,1)$ such that $$ f(t) = \int_0^t g(s)\ud s, \quad t\in [0,1],$$
and that this function $g$, if it exists, is unique and satisfies
$$ \n g \n_{1} = \sup_{N,\pi} \sum_{n=1}^{N}|f(t_n)-f(t_{n-1})|,$$
where the supremum is taken over all $N\ge 1$ and all partitions
$\pi=\{0=t_0<t_1<\cdots<t_N=1\}$ of $[0,1]$.

\noindent{\em Hint:}\
For the `only if' part, first show that $f$ determines a
$\K$-valued Borel measure $\nu$ on $[0,1]$ satisfying
$\nu((a,b])=f(b)-f(a)$ for all $0\le a<b\le1.$
Use the absolute continuity of $f$ to prove that $\nu$ is absolutely
continuous with respect to Lebesgue measure, and then apply the
Radon--Nikod\'ym theorem.

\item\index{$BMO$@$BMO(\R^d)$}
A function $f\in L_{\rm loc}^1(\R^d)$ is said to have {\em bounded mean oscillation}\index{bounded!mean oscillation} if $$|f|_{BMO(\R^d)} := \sup_B \frac1{|B|}\int_B |f(x) - {\rm av}_B (f)|\ud x$$
is finite, where the supremum is taken over all balls $B$ in $\R^d\!$, $|B|$ is the Le\-besgue measure of $B$, and
${\rm av}_B (f) := \frac1{|B|}\int_B f(y)\ud y$ is the average of $f$ on $B$.
\begin{enumerate}[\rm(a), leftmargin=*]
  \item Show that if $f$ and $g$ have bounded mean oscillation, then:
  \begin{enumerate}[label={\rm(\roman*)}, leftmargin=*]
    \item $cf$ has bounded mean oscillation and $$|c f|_{BMO(\R^d)} = |c||f|_{BMO(\R^d)};$$
    \item $f+g$ has bounded mean oscillation and $$|f+g|_{BMO(\R^d)} \le |f|_{BMO(\R^d)} + |g|_{BMO(\R^d)}.$$
  \end{enumerate}
  \item Show that $|f|_{BMO(\R^d)} = 0$ if and only if $f$ is almost everywhere constant.
  \item Show that every $f\in L^\infty(\R^d)$ has bounded mean oscillation and $$|f|_{BMO(\R^d)} \le 2 \n f\n_\infty.$$
  \item Show that the unbounded function $x\mapsto \log |x|$ has bounded mean oscillation.
  \item
  Show that the quotient, denoted $BMO(\R^d)$, of the vector space of functions with bounded mean oscillation by the constant functions is a Banach space in a natural way.
\end{enumerate}

\item\label{prob:lub}
Show that in a vector lattice $(V,\le)$, the greatest lower bound and the least upper bound of two elements $u,v\in V$ satisfy $u \wedge v = u - (u-v)^+ = v - (v-u)^+ $ and $u\vee v = u+(v-u)^+ = v+(u-v)^+$\!.

\item\label{prob:vplusvminus}
Let $V$ be a vector lattice.
\begin{enumerate}[\rm(a), leftmargin=*]
 \item Show that for all $v\in V$ we have $v^+\wedge v^-=0$ and $v^+\vee v^- = |v|$.
 \item Show that if $v,w,w'\in V$ satisfy $v = w-w'$ with $w\ge 0$ and $w'\ge 0$, then $w \ge v^+$ and $w' \ge v^-$\!.
 \item Show that if $v,w,w'\in V$ satisfy $v = w-w'$ with $w\ge 0$, $w'\ge 0$, and $w\wedge w'=0$, then $w = v^+$ and $w'= v^-$\!.
\end{enumerate}

\item
Prove that if $X$ is a normed vector lattice, the lattice operations $(x,y)\mapsto x\wedge y$ and $(x,y)\mapsto x\vee y$
are continuous from $X\times X$ to $X$.

\item
Provide the missing details of the proof, outlined in
Section \ref{sec:BL}, that the real Banach space of real-valued measures
in $M(\Om)$ is a Banach lattice.

\end{problems}

%% file: ch03-HilbertSpaces.tex
\chapter{Hilbert Spaces}\label{ch:Hilbert-spaces}

\blfootnote{\rm \small \begin{quote}This book has been published by Cambridge University Press in the series `Cambridge Studies in Advanced Mathematics'. The present version is free to view and download for personal use only. Not for re-distribution, re-sale or use in derivative works.\newline {\copyright} Jan van Neerven\end{quote}}

\noindent
Arguably the most important class of Banach spaces is the class of Hilbert spaces. These spaces play a central role in the theory and in various areas of application, some of which will be discussed in later chapters. The present introductory chapter develops the basic geometric properties of Hilbert spaces arising from the presence of an inner product generating the norm, such as the orthogonal complementation of closed subspaces, the existence of orthonormal bases, and the self-duality of Hilbert spaces embodied by the Riesz representation theorem.

\section{Hilbert Spaces}\label{sec:HilbertSpaces}

Let $V$ be a vector space.
A mapping $\phi: V\times V\to \K$ is called {\em sesquilinear}\index{sesquilinear} if it is linear in the first
variable and conjugate-linear\index{conjugate-linear} in the second variable, that is,
\begin{align}
 \phi(v+v'\!,w) = \phi(v,w) + \phi(v'\!,w), & \qquad \phi(cv,w) = c\phi(v,w), \label{eq:sesqui1}\\
 \phi(v,w+w') = \phi(v,w) + \phi(v,w'), & \qquad \phi(v,cw) = \ov c \phi(v,w),\label{eq:sesqui2}
\end{align}
 for all $c\in \K$ and $v,v'\!,w,w'\in V$.
The complex conjugation in \eqref{eq:sesqui2} is of course redundant when the scalar field is real and sesquilinearity reduces to bilinearity in that case.

\begin{definition}[Inner products]\label{def:ip} An {\em inner product space}\index{inner product!space} is a pair $(H, \iprod{\cdot}{\cdot})$, where $H$ is a vector space and
 $\iprod{\cdot}{\cdot}$\index{$X$@$\iprod{x}{y}$} is an {\em inner product}\index{inner product} on $H\times H$, that is, a sesquilinear mapping from $H\times H$ to $\K$
 with the following properties:
 \begin{enumerate}[label={\rm(\roman*)}, leftmargin=*]
 \item\label{it:ip1} $\iprod{x}{x}\ge 0$ for all $x\in H$ and $\iprod{x}{x} = 0 \Rightarrow  x=0$;
 \item\label{it:ip2} $\iprod{x}{y} = \overline{\iprod{y}{x}}$ for all $x,y\in H$.
 \end{enumerate}
\end{definition}

The conjugation bar in \ref{it:ip2} is again redundant when the scalar field is real.
If \ref{it:ip2} holds, then \eqref{eq:sesqui1} implies \eqref{eq:sesqui2}.

It will be used frequently without further comment that
 $$ \hbox{ if} \ \iprod{x}{y} = 0 \hbox{ for all }  y\in H,  \hbox{ then }  x=0.$$
Indeed, the hypothesis implies that $\iprod{x}{x} = 0$, and then $x=0$ by the definition of an inner product.

When the inner product $\iprod{\cdot}{\cdot}$ is understood we simply write $H$ instead of $(H, \iprod{\cdot}{\cdot})$.

\begin{example}\label{ex:innerprod}
Here are some examples of inner products:
 \begin{enumerate}[label={\rm(\roman*)}, leftmargin=*]
  \item on $\K^d$ an inner product is given by
  $ \iprod{x}{y} = \sum_{n=1}^d x_n \ov{y_n}$;
  \item on $\ell^2$ an inner product is given by
  $ \iprod{a}{b} = \sum_{n\ge 1} a_n \ov {b_n}$;
  \item on $L^2(\Om,\mu)$ an inner product is given by
  $ \iprod{f}{g} = \int_\Om f\ov{g}\ud \mu.$
 \end{enumerate}
\end{example}

In order to turn inner product spaces into normed vector spaces we need the following inequality. Its finite-dimensional version has already been used in various places in Chapters \ref{ch:Banach} and \ref{ch:ClassicalBanach}.

\begin{proposition}[Cauchy--Schwarz inequality]\label{prop:CS}\index{inequality!Cauchy--Schwarz}
Let $H$ be a vector space and consider a sesquilinear mapping $\iprod{\cdot}{\cdot}: H\times H \to \K$
 with the following properties:
 \begin{enumerate}[label={\rm(\roman*)}, leftmargin=*]
 \item\label{it:CS1} $\iprod{x}{x}\ge 0$ for all $x\in H$;
 \item\label{it:CS2} $\iprod{x}{y} = \overline{\iprod{y}{x}}$ for all $x,y\in H$.
 \end{enumerate}
 Then for all $x,y\in H$ we have
 $$ |\iprod{x}{y}|^2 \le \iprod{x}{x}\iprod{y}{y}.$$
\end{proposition}
\begin{proof}
For every $c\in\K$ we have
\begin{align*}
0\le \iprod{x-cy}{x-cy}
  & = \iprod{x}{x}- \iprod{x}{cy} - \iprod{cy}{x} + \iprod{cy}{cy}
 \\ & = \iprod{x}{x}-2\Re(\overline c\,\iprod{x}{y})  + |c|^2 \iprod{y}{y}.
\end{align*}
If $\iprod{y}{y}=0$, the choice $c=t\,\iprod{x}{y}$ with $t>0$ gives
$0\le \iprod{x}{x}-2t|\iprod{x}{y}|^2$.
Letting $t\to\infty$, we see that this can only hold if $\iprod{x}{y}=0$.
If $\iprod{y}{y}\neq 0$, the choice $c =
{\iprod{x}{y}}/\iprod{y}{y}$ gives
$$
0\le \iprod{x}{x} - 2\Re\frac{|\iprod{x}{y}|^2}{\iprod{y}{y}} +\frac{|\iprod{x}{y}|^2}{\iprod{y}{y}} = \iprod{x}{x} - \frac{|\iprod{x}{y}|^2}{\iprod{y}{y}},
$$
and multiplication by $\iprod{y}{y}$ gives the desired result.
\end{proof}

\begin{proposition}
Every inner product space $H$ can be made into a normed space by defining
\begin{align*} \n x\n := \iprod{x}{x}^{1/2}\!, \quad x\in H.\end{align*}
\end{proposition}
\begin{proof} We must check that $\n \cdot\n$ defines a norm on $H$.
It is immediate that $\n x\n = 0$ implies $x=0$ and that $\n c x\n = |c|\n x\n$.
The triangle inequality follows from the Cauchy--Schwarz inequality:
\[
\n x+y\n^2
 = \n x\n^2 + 2 \Re\iprod{x}{y}  + \n y\n^2
 \le \n x\n^2 + 2\n x\n\n y\n + \n y\n^2
= (\n x\n + \n y\n)^2\!.
\qedhere\]
\end{proof}

Henceforth it is understood that inner product spaces are endowed with this norm.

When $H$ and $H'$ are inner product spaces, then
$\iprod{(x,y)}{(x',y')}:=\iprod{x}{x'} + \iprod{y}{y'}$
defines an inner product on $H\times H'$, and the associated norm is given by
$\n (x,y)\n^2 = \n x\n^2+\n y\n^2$. If $H$ and $H'$ are Hilbert spaces, then so is $H\times H'$.

 \begin{corollary}\label{cor:ip-joint}
Every inner product $\iprod{\cdot}{\cdot}$ is continuous as a mapping from $H\times H$ to $\K$.
\end{corollary}
\begin{proof}
It suffices to show that if $x_n\to x$ and $y_n\to y$, then
$\iprod{x_n}{y_n}\to \iprod{x}{y}$. We have
\begin{align*}
 |\iprod{x_n}{y_n} -  \iprod{x}{y}| & \le |\iprod{x_n}{y_n} -  \iprod{x_n}{y}| + |\iprod{x_n}{y} -  \iprod{x}{y}|
 \\ & =  |\iprod{x_n}{y_n-y}| + |\iprod{x_n-x}{y}|
  \le \n x_n\n \n y_n-y\n  + \n x_n-x\n \n y\n.
\end{align*}
Since convergent sequences are bounded, the number $M:= \sup_{n\ge 1}\n x_n\n$ is finite, and we find
$$ |\iprod{x_n}{y_n} -  \iprod{x}{y}| \le M \n y_n-y\n  + \n x_n-x\n \n y\n.$$
Both terms on the right-hand side tend to $0$ as $n\to \infty$.
\end{proof}

\begin{proposition}[Parallelogram identity]\label{prop:parallelogram}
In every inner product space $H$ the {\em parallelogram identity}\index{parallelogram identity} holds: for all $x, y\in H$
we have $$
  2\|x\|^2 + 2\|y\|^2 = \|x+y\|^2 + \|x-y\|^2\!.
 $$
Conversely, if $X$ is a normed space with the property that the parallelogram identity holds for all $x, y\in X$,
then there exists an inner product on $X$ generating the norm of $X$.
\end{proposition}

In what follows we only need the first part of the proposition. Its converse is included for reasons of completeness and can be safely skipped upon first reading.
The parallelogram identity has $L^p$-analogues, known as Clarkson's inequalities (see Problem \ref{prob:Clarkson}).

\begin{proof}
 The proof of the first part is routine:
\begin{align*} \|x+y\|^2 + \|x-y\|^2
 = \iprod{x+y}{x+y} + \iprod{x-y}{x-y}
  =   2\iprod{x}{x} +  2\iprod{y}{y} =  2\|x\|^2 + 2\|y\|^2\!.
\end{align*}

The proof of the second assertion is quite involved and relies on
finding a formula for the inner product in terms of the norm of $X$.
To get an idea of what this formula might look like we first {\em assume} there is such an inner product and denote it by $\iprod{\cdot}{\cdot}$.
Then
\begin{align*} \n x \pm y\n^2
= \n x\n^2 \pm 2\Re \iprod{x}{y} + \n y\n^2\!.
\end{align*}
From this we obtain
\begin{align}\label{eq:ip-real}  \Re\iprod{x}{y} = \frac14\Bigl(\n x+y\n^2 - \n x- y\n^2\Bigr).
\end{align}
If the scalar field $\K$ is real, then $\Re\iprod{x}{y} = \iprod{x}{y}$ and the above identity expresses the inner product in terms of the norm.
If the scalar field $\K$ is complex, by the previous identity we obtain
\begin{align*} \Im\iprod{x}{y} =  \Re\iprod{x}{iy} = \frac14\Bigl(\n x+iy\n^2 - \n x- iy\n^2\Bigr).
\end{align*}
This leads to the identity
\begin{equation}\label{eq:ip-compl}
\begin{aligned}
\iprod{x}{y}
& = \Re\iprod{x}{y}+i\Im\iprod{x}{y}
\\ & = \frac14\Bigl( \n x+y\n^2 - \n x- y\n^2 +i \n x+iy\n^2 - i \n x - iy\n^2 \Bigr) .
\end{aligned}
\end{equation}

In an arbitrary normed space we could now try to {\em define} an inner product by \eqref{eq:ip-real} if $\K=\R$ and by \eqref{eq:ip-compl} if $\K=\C$, but this does not always give an inner product. If, however, the parallelogram identity holds, then it does. Let us check this for the case $\K = \R$. First,
$$ \iprod{x}{x}= \frac14\Bigl(\n x+x\n^2 - \n x- x\n^2\Bigr) = \n x\n^2 \ge 0$$
and $\iprod{x}{x} = 0$ implies $\n x\n = 0$ and hence $x=0$. Second, the identity $\iprod{x}{y} = \iprod{y}{x}$ is immediate. Third,
\begin{align*}\iprod{x+x'}{y} & = \frac14\Bigl(\n (x+x')+y\n^2 - \n (x+x')- y\n^2\Bigr)
\\ & = \frac14\Bigl(2\n x\n^2 +2\n x'+y\n^2 - \n x-(x'+y)\n^2\Bigr)
\\ & \qquad\qquad - \frac14\Bigl(2\n x\n^2 + 2\n x'- y\n^2 - \n  x-(x'- y)\n^2\Bigr)
\\ & = \frac12\Bigl(\n x'+y\n^2 - \n x'- y\n^2\Bigr)
\\ & \qquad\qquad - \frac14\Bigl(  \n x-(x'+y)\n^2- \n  x-(x'- y)\n^2\Bigr)
\\ & = 2\iprod{x'}{y} + \iprod{x-x'}{y}.
\end{align*}
This proves that
\begin{align}\label{eq:Hilb-parallel}\iprod{x+x'}{y} - \iprod{x-x'}{y} = 2\iprod{x'}{y}.
\end{align}
Taking $x = x'$ in \eqref{eq:Hilb-parallel}, we obtain $\iprod{2x}{y} = 2\iprod{x}{y}$.
Also, $$\iprod
{-x}{y}=\frac 14(\|y-x\|^2-\|y+x\|^2)=
-\frac14(\|x+y\|^2-\|x-y\|^2)=-\iprod {x}{y}.$$
Now let $x_0,x_1\in X$ be arbitrary. Applying \eqref{eq:Hilb-parallel} to $x = \frac12(x_0-x_1)$ and $x' = \frac12(x_0+x_1)$
gives
$$\iprod{x_0}{y} + \iprod{x_1}{y}= \iprod{x_0}{y} - \iprod{-x_1}{y} =  2\Bigl(\frac12(x_0+x_1)\Big|y\Bigr) = \iprod{x_0+x_1}{y}.$$
This gives additivity in the first coordinate. Additivity in the second coordinate is proved similarly.
Using this inductively, for positive integers $k$ we obtain
\begin{align*}\iprod{kx}{y} & = \iprod{(k-1)x+x}{y} = \iprod{(k-1)x}{y} + \iprod{x}{y} \\ & =  \iprod{(k-2)x+x}{y} + \iprod{x}{y} = \iprod{(k-2)x}{y} + 2 \iprod{x}{y}
\\ & =\cdots = \iprod{x}{y}+(k-1)\iprod{x}{y} = k\iprod{x}{y}.
\end{align*}
Applying this to $k^{-1}x$ we also find $k^{-1}\iprod{x}{y} =\iprod{k^{-1}x}{y}$.
Next let $q = m/n$ with $m,n$ positive integers.
By what we just proved,
$$ \iprod{qx}{y} = m \iprod{n^{-1}x}{y} = mn^{-1}\iprod{x}{y} = q\iprod{x}{y}.$$
This proves homogeneity with respect to multiplication by the positive rationals. For such rationals we also have $\iprod{-qx}{y} = -\iprod{qx}{y} = -q\iprod{x}{y}$,
and therefore homogeneity holds for all rationals.
Finally, by the continuity of the norm $\n\cdot\n$ and the continuity of scalar multiplication, the mappings $q \mapsto \iprod{qx}{y}$ and $q\mapsto q\iprod{x}{y}$ are continuous.
Therefore, the identity $\iprod{cx}{y} = c\iprod{x}{y}$ for arbitrary $c\in \R$ follows by approximation by rationals.

This completes the proof that for $\K = \R$, the formula for $\iprod{\cdot}{\cdot}$ in \eqref{eq:ip-real} indeed defines an inner product. We have already seen that $\iprod{x}{x} = \n x\n^2$\!, so this
inner product generates the norm of $X$.

For $\K=\C$, define $\iprod{x}{y}$ by \eqref{eq:ip-compl}. Its real part is given by \eqref{eq:ip-real}, and its imaginary part is obtained by applying \eqref{eq:ip-real} to $x$ and $iy$. The additivity and real homogeneity identities therefore follow from the real case, applied once to $y$ and once to $iy$. Moreover, the identities $\iprod{ix}{y}=i\iprod{x}{y}$ and $\iprod{x}{iy}=-i\iprod{x}{y}$ follow directly from \eqref{eq:ip-compl}; hence the form is complex linear in the first variable and conjugate-linear in the second. Finally, $\iprod{y}{x}=\ov{\iprod{x}{y}}$ and $\iprod{x}{x}=\n x\n^2$ are immediate from \eqref{eq:ip-compl}.
\end{proof}

\begin{definition}[Hilbert spaces]
A {\em Hilbert space}\index{Hilbert space} is a complete inner product space.
\end{definition}

Thus, by definition, every Hilbert space is a Banach space.

\begin{example}
By the completeness results in the preceding chapter, all three spaces $\K^d\!$, $\ell^2$\!, $L^2(\Om,\mu)$ appearing in Example \ref{ex:innerprod} are Hilbert spaces.
\end{example}

Further examples of Hilbert spaces will be given in the problems section.

\begin{proposition}[Completion]\label{prop:completion-H}\index{completion}
Let $H$ be an inner product space. On its completion $\ov H$ as a normed space, a well-defined inner product is obtained by setting
 $$\iprod{\ov x}{\ov y} := \limn \iprod{x_n}{y_n}, \quad \ov x,\ov y\in \ov H,$$
 whenever $x_n,y_n\in H$ satisfy $\limn x_n = \ov x$ and $\limn y_n = \ov y$ in $\ov H$.
The norm associated with this inner product agrees with the norm of $\ov H$ obtained by completing $H$ as a normed space.
\end{proposition}
\begin{proof}
The proof relies on a few routine verifications, some of which are left as an exercise to the reader.

Arguing as in Corollary \ref{cor:ip-joint}, we see that the limit $\limn \iprod{x_n}{y_n}$ exists and is independent of the approximating sequences. Clearly it agrees with the original inner product in $H$ when $\ov x,\ov y\in H$.
To see that $\iprod{\cdot}{\cdot}$ is indeed an inner product on $\ov H$, suppose that $\iprod{\ov x}{\ov x}=0$ and let $\limn x_n = \ov x$ with each $x_n\in H$. Then $\limn\iprod{x_n}{x_n} =0$ and therefore $\limn x_n = 0$, which, by the construction of the completion, means that the Cauchy sequence $(x_n)_{n\ge 1}$ defines the zero element $\ov 0$ of $\ov H$. It follows that $\ov x = \limn x_n= \ov 0$ in $\ov H$.
The remaining properties of an inner product follow trivially by limiting arguments.

Finally, the norm of $\ov H$ agrees with the norm generated by its inner product. This follows again by approximation: if $\limn x_n = \ov x$ with each $x_n\in H$,
then for the norm $\n\cdot\n$ obtained via the completion procedure in the proof of Theorem \ref{thm:completion} we have $$\n \ov x\n^2 = \limn \n x_n\n^2 = \limn \iprod{x_n}{x_n} = \iprod{\ov x}{\ov x},$$ the first step being true by
the definition of this norm, the second because the norm of $H$ is generated by its inner product, and the third by the definition of the inner product on $\ov H$.
\end{proof}

\section{Orthogonal Complements}

In this section we introduce the notion of orthogonality and prove that every closed subspace has an orthogonal complement.
Throughout, we fix a Hilbert space $H$.

\begin{definition}[Orthogonality]\label{def:orthogonality}
The elements $x,x'\in H$ are said to be {\em orthogonal}\index{orthogonal}
if $\iprod{x}{x'} = 0$; in this case we write $$ x\perp x'\!.$$
Two subsets $A$ and $B$ of $H$ are called {\em orthogonal} if $a\perp b$ for all $a\in A$ and $b\in B$.
\end{definition}

Orthogonal elements $x\perp x'$ satisfy the Pythagorean identity
$$ \n x+x'\n^2 = \n x\n^2 + \n x'\n^2\!,$$
as is seen by expanding the squared norms in terms of inner products.

\begin{definition}[Orthogonal complement]\label{def:orthogonal-complement}
The {\em orthogonal complement}\index{orthogonal!complement} of a subset $A$ of $H$ is the set\index{$Aa$@$A^\perp$}
$$ A^\perp:= \{x\in H:\, x\perp a\ \hbox{ for all } \ a\in A\}.$$
\end{definition}

The orthogonal complement $A^\perp$ of a subset $A$ is a closed subspace of $H$. Indeed, it is immediate that $A^\perp$ is a vector space. To prove its closedness, let $x_n\to x$ in $H$
with $x_n\in A^\perp$\!. Then, by the continuity of the inner product, for all $a\in A$ we obtain $\iprod{x}{a} = \limn \iprod{x_n}{a} = 0$.

The most important result on orthogonality is certainly the fact that every closed subspace $Y$ of a Hilbert space is orthogonally complemented by $Y^\perp$\!. This is the content
of Theorem \ref{thm:ONC} below. For its proof we need the following approximation theorem for convex closed sets in Hilbert space. Recall
that a subset $C$ of a vector space is called {\em convex}\index{convex!set} if for all $x_0,x_1 \in C$  we have $(1-\la)x_0+\la x_1\in C$ for all $0\le \la\le 1$.

\begin{theorem}[Best approximation]\label{thm:convex-best-approx}
 Let $C$ be a nonempty convex closed subset of $H$.
 Then for all $x\in H$ there exists a unique $c\in C$ that minimises\index{minimiser} the distance from $x$ to the points of $C$:
 $$ \n x-c\n = \min_{y\in C} \n x-y\n.$$
\end{theorem}
\begin{proof}
 Let $(y_n)_{n\ge 1}$ be a sequence in $C$ such that $$\limn \n x-y_n\n = \inf_{y\in C} \n x-y\n =: D.$$
 We claim that this sequence is Cauchy. By the parallelogram identity of Proposition \ref{prop:parallelogram}, applied to the vectors
 $x-y_m$ and $x-y_n$,
 $$ \n y_n-y_m\n^2 + \n 2x-(y_n+y_m)\n^2 = 2  \n x-y_m\n^2 + 2\n x-y_n\n^2\!.$$
 As $m,n\to\infty$, the right-hand side tends to $2D^2+2D^2 = 4D^2$\!, whereas from $\frac12(y_m+y_n)\in C$ (by convexity)
 it follows that
 $$ \n 2x-(y_n+y_m)\n^2 = 4\n x-\tfrac12(y_n+y_m)\n^2 \ge 4D^2\!.$$
 It follows that
 $$ \limsup_{m,n\to\infty} \n y_n-y_m\n^2 \le 4D^2 - 4D^2 = 0.$$
The limit superior is also nonnegative, and therefore it equals $0$.
This proves the claim.

Since $H$ is complete, we have $\limn y_n =c$ for some $c\in H$, and since $C$ is closed we have $c\in C$.
Now $\n x-c\n = \limn \n x-y_n\n = D$, so $c$ minimises the distance to $x$. This proves existence.

To prove uniqueness, suppose that $c_1,c_2\in C$ are minimisers, that is,
$\|x-c_1\|=\|x-c_2\|=D$.
By convexity $m:=\frac12(c_1+c_2)$ belongs to $C$, so $\|x-m\|\ge D$.
Applying the parallelogram identity to $x-c_1$ and $x-c_2$, we obtain
$$
\|c_1-c_2\|^2+\|2x-(c_1+c_2)\|^2
=2\|x-c_1\|^2+2\|x-c_2\|^2
=4D^2\!.
$$
But $\|2x-(c_1+c_2)\|^2=4\|x-m\|^2\ge 4D^2$,
and therefore $\|c_1-c_2\|^2 \le 4D^2-4D^2=0$, that is,  $c_1=c_2$.
\end{proof}

Both the existence and uniqueness parts of this theorem fail for general Banach spaces; see Problems \ref{prob:minimiser-C1} and \ref{prob:minimiser-C2}.

\begin{theorem}[Closed subspaces are orthogonally complemented]\label{thm:ONC}\index{$H_1 \oplus H_2$}
 If $Y$ is a closed linear subspace of $H$,
then we have an orthogonal direct sum decomposition\index{orthogonal!direct sum}\index{orthogonal!decomposition} $$ H = Y\oplus Y^\perp\!,$$
that is, we have $Y\cap Y^\perp =\{0\}$, $Y+Y^\perp = H$, and $Y\perp Y^\perp$\!.
\end{theorem}

\begin{proof}
We have already seen that $Y^\perp$ is a closed subspace. If $y\in Y\cap Y^\perp$\!, then
$y\perp y$, so $\iprod{y}{y}=0$ and $y=0$.
It remains to show that $Y+Y^\perp = H$.

Let $x\in H$ be arbitrary but fixed. We must show that $x\in Y+Y^\perp$\!. Let  $\pi_{\,Y}: H \to Y$ denote the mapping arising from Theorem \ref{thm:convex-best-approx}, that is, $\pi_{\,Y} x$ is the unique element of $Y$ minimising the distance to $x$:
$$ \n x-\pi_{\,Y} x  \n  = \min_{y\in Y} \n x- y \n.$$
Set $y_0:= \pi_{\,Y} x$ and $y_1:= x-y_0$.
Then $y_0\in Y$,  and for all $y\in Y$ we have
\begin{align}\label{eq:Yperp} \n y_1\n = \n x-y_0\n \le \n x-\underbrace{(y_0-y)}_{\in Y}\n = \n y+(x-y_0)\n = \n y+y_1\n.
\end{align}
We claim that \eqref{eq:Yperp} implies $y_1\in  Y^\perp$\!. To see this, fix a nonzero $z\in Y$.
For any $c\in \K$ we have, by \eqref{eq:Yperp},
$$ \n y_1\n^2 \le \n c z + y_1\n^2 = |c|^2\n z\n^2 + 2\Re\iprod{c z}{y_1} + \n y_1\n^2\!.$$
Taking $c = -\ov{\iprod{z}{y_1}}/\n z\n^2$ gives
$$0\le \frac{|\iprod{z}{y_1}|^2}{\n z\n ^2} - 2\frac{|\iprod{z}{y_1}|^2}{\n z\n^2} ,$$
which is only possible if $\iprod{z}{y_1} =0$. Since $0\neq z\in Y$ was arbitrary, this shows that $y_1\in Y^\perp$\!.
This proves the claim. It follows that $x = y_0+y_1$ belongs to $Y+Y^\perp$\!.
\end{proof}

\begin{definition}[Orthogonal projection onto a closed subspace]\label{def:orth-proj-closed-subspace}
The projection $\pi_{\,Y}$ onto $Y$ along $Y^\perp$ given by
$$ \pi_{\,Y} (y+y^\perp):= y, \quad y\in Y, \ y^\perp\in Y^\perp\!,$$
is called the {\em orthogonal projection}\index{orthogonal!projection} onto $Y$.
\end{definition}

The Pythagorean identity implies that $\n \pi_{\,Y}\n \le 1$ (with equality if $Y\neq\{0\}$).
As was shown in the course of the proof of Theorem \ref{thm:ONC}, the projection $\pi_{\,Y}$
coincides with the mapping arising from Theorem \ref{thm:convex-best-approx}. For general closed convex sets $C$, the distance-minimising mapping of Theorem \ref{thm:convex-best-approx} is generally nonlinear.

As a corollary of Theorem \ref{thm:ONC}, for closed subspaces $Y$ we find that $(Y^{\perp})^\perp = Y$,
and more generally for any subspace $Y$ we obtain (since $Y^\perp = (\ov Y)^\perp$)
$$(Y^{\perp})^\perp = \overline Y.$$

By way of example, if $(h_n)_{n=1}^N$ is an {\em orthonormal sequence}, that is, $\iprod{h_j}{h_k} = \delta_{jk}$ for all $1\le j,k\le N$, the orthogonal projection $\pi_{\,N}$ in $H$
onto the span of $h_1,\dots,h_N$ is given by
$\pi_{\,N} x = \sum_{n=1}^N \iprod{x}{h_n}h_n, $
as the reader will have no difficulty checking.

\section{The Riesz Representation Theorem}\label{sec:Riesz-representation}

 Let $H$ be a Hilbert space.
As a first application of Theorem \ref{thm:ONC} we prove the  {\em Riesz representation theorem}\index{theorem!Riesz representation, for Hilbert spaces}, which sets up
a conjugate-linear identification of a Hilbert space $H$ and its dual $H^* = \calL(H,\K)$.

By the Cauchy--Schwarz inequality, every $h\in H$ defines a bounded functional
$\psi_h: H \to \K$ by taking inner products:
$$ \psi_h (x):= \iprod{x}{h}, \quad x\in H.$$
Boundedness is evident from
$ |\psi_h(x)| = |\iprod{x}{h}| \le \n x\n \n h\n,$
which shows that $\n \psi_h\n \le \n h\n$. From
$ \psi_h(h) = \iprod{h}{h} = \n h\n^2$ we see that also $\n \psi_h\n \ge \n h\n$, so that
$ \n \psi_h \n = \n h\n.$

All bounded functionals $\phi:H\to \K$ arise in this way:

\begin{theorem}[Riesz representation theorem]  \label{thm:Riesz}
If $\phi:H\to \K$ is a bounded functional, there exists a unique element $h\in H$ such that $\phi = \psi_{h}$, that is,
$$ \phi(x) = \iprod{x}{h}, \quad x\in H.$$
\end{theorem}

\begin{proof}
If $\phi(x) = 0$ for all $x\in H$, we take $h = 0$. If $\phi\neq0$,
the null space $\ker(\phi)$ is a proper closed subspace of $H$. Using Theorem \ref{thm:ONC}, choose a norm-one vector $y_0\in(\ker(\phi))^\perp$\!. Since $y_0\notin\ker(\phi)$, we have $\phi(y_0)\neq0$.
Fix an arbitrary $x\in H$ and choose the scalar $c\in \K$ so that $\phi(x) = c\phi(y_0)$. Then $x-c y_0 \in\ker(\phi)$,
and from $y_0\perp \ker(\phi)$ we infer $\iprod{x-c y_0}{y_0} = 0$, or equivalently
$$\iprod{x}{y_0} =c\iprod{y_0}{y_0} = c, $$ recalling that $\iprod{y_0}{y_0} = \n y_0\n^2 = 1$. It follows that
\begin{align*} \phi(x) = c\phi(y_0) = \iprod{x}{y_0}\phi(y_0)= \iprod{x}{\ov{\phi(y_0)}y_0}.
\end{align*}
This proves that $\phi = \psi_{\ov{\phi(y_0)}y_0}$.

To prove uniqueness, suppose that $\phi = \psi_h = \psi_{h'}$ for $h,h'\in H$. Then $\n h-h'\n^2 = \iprod{h-h'}{h-h'} = \psi_h(h-h') - \psi_{h'}(h-h') = 0$ and therefore $h'=h$.
\end{proof}

\begin{proposition}\label{prop:Hilbert-weakconvergence}
For every bounded sequence $(x_n)_{n\ge 1}$ in $H$
there exist a subsequence $(x_{n_k})_{k\ge 1}$ and an
 $x\in H$ such that
 $$ \limk \iprod{h} {x_{n_k}} = \iprod{h} {x}, \quad h\in H.$$
\end{proposition}

\begin{proof} Let $Y$ denote the linear span of $(x_n)_{n\ge 1}$, and let
$H_0$ denote its closure.
The usual diagonal argument, using boundedness of the scalar sequences $(\iprod{x_j}{x_n})_{n\ge 1}$, gives a subsequence $(x_{n_k})_{k\ge 1}$ such that the limit $\limk\iprod{x_j}{x_{n_k}}$ exists for every $j\ge 1$. By linearity, the limit $\limk\iprod{y}{x_{n_k}}$ exists for all $y\in Y$.

For $k\ge 1$ define $T_k\in\calL(H_0,\K)$ by $T_kh:=\iprod{h}{x_{n_k}}$. With $M:=\sup_{n\ge 1}\n x_n\n$ we have $\n T_k\n\leq M$ for all $k\ge 1$. Since the sequence $(T_ky)_{k\ge 1}$ converges for all $y\in Y$, Proposition \ref{prop:approxTn} implies that the limit $\phi(h):=\limk T_kh$ exists for all $h\in H_0$, and that $\phi\in H_0^*$. By the Riesz representation theorem, applied to the Hilbert space $H_0$, there exists an $x\in H_0$ such that $\phi(h)=\iprod{h}{x}$ for all $h\in H_0$. Thus $\limk\iprod{h}{x_{n_k}}=\iprod{h}{x}$ for all $h\in H_0$.

It remains to pass from $H_0$ to $H$. Let $h\in H$, and write $h=h_0+h_0^\perp$ with respect to the orthogonal decomposition $H=H_0\oplus H_0^\perp$. Since $x_{n_k}\in H_0$ for all $k\ge 1$, and since $x\in H_0$, we have
$$
\iprod{h}{x_{n_k}}=\iprod{h_0}{x_{n_k}}, \qquad \iprod{h}{x}=\iprod{h_0}{x}.
$$
The asserted convergence now follows from the result already proved on $H_0$.
\end{proof}

\section{Orthonormal Systems}\label{subsec:OS}
We now turn from orthogonal decompositions to orthonormal systems. These
play the role of coordinate systems in Hilbert spaces. We first record the
basic convergence criterion for orthogonal series and then pass to Parseval's identity and orthonormal bases.

\begin{proposition}\label{prop:orth-sum}
Let $(x_n)_{n\ge 1}$ be a sequence in $H$ with $x_m\perp x_n$ for all $m\neq n$.
The following assertions are equivalent:
\begin{enumerate}[label={\rm(\arabic*)}, leftmargin=*]
 \item\label{it:orth-sum1} $\sum_{n\ge 1} x_n$ converges in $H$;
 \item\label{it:orth-sum2} $\sum_{n\ge 1} \n x_n\n^2 < \infty.$
\end{enumerate}
In this situation we have
$ \n \sum_{n\ge 1} x_n \n^2 =\sum_{n\ge 1} \n x_n\n^2\!.$
\end{proposition}
\begin{proof}
Let us first note that if $I$ is any finite set of positive integers, then
\begin{align*}\Big\n \sum_{n\in I}x_n \Big\n^2
 &= \Big(\sum_{m\in I} x_m \Big| \sum_{n\in I} x_n\Big) = \sum_{m\in I}\sum_{n\in I} \iprod{x_m}{x_n}
 =\sum_{n\in I} \n x_n\n^2\!,
 \end{align*}
since $\iprod{x_m}{x_n} = 0$ if $m\neq n$.

\smallskip
\ref{it:orth-sum1}$\Rightarrow$\ref{it:orth-sum2}: \  If $\sum_{n\ge 1} x_n$ converges in $H$, say to $x$,
then $x = \lim_{N\to\infty}
 \sum_{n=1}^N x_n$ in $H$ and therefore
 \begin{align*}\n x\n^2 & = \lim_{N\to\infty} \Big\n \sum_{n= 1}^N x_n \Big\n^2
  = \lim_{N\to\infty} \sum_{n=1}^N \n x_n\n^2\!.
 \end{align*}
It follows that $\sum_{n\ge 1} \n x_n\n^2 = \n x\n^2$\!. This also proves the final identity
in the statement of the proposition,
since by definition $\sum_{n\ge 1} x_n = x$.

\smallskip
\ref{it:orth-sum2}$\Rightarrow$\ref{it:orth-sum1}: \  Suppose, conversely, that $\sum_{n\ge 1} \n x_n\n^2 < \infty.$ Then
\begin{align*}\lim_{\substack{M,N\to \infty \\ N>M}} \Big\n \sum_{n=1}^N x_n - \sum_{n=1}^M x_n\Big\n^2
 = \lim_{\substack{M,N\to \infty \\ N>M}} \Big\n \sum_{n=M+1}^N x_n \Big\n^2
= \lim_{\substack{M,N\to \infty \\ N>M}} \sum_{n=M+1}^N \n x_n\n^2 =0.
\end{align*}
It follows that $(\sum_{n=1}^N x_n)_{N\ge 1}$ is Cauchy, and hence convergent.
\end{proof}

As a special case of Proposition \ref{prop:orth-sum} we record:

\begin{proposition}[Parseval]\label{prop:Parseval}\index{theorem!Parseval}
Let $(h_n)_{n\ge 1}$ be an orthonormal sequence in $H$. For a scalar sequence $(c_n)_{n\ge 1}$, the following assertions are equivalent:
\begin{enumerate}[label={\rm(\arabic*)}, leftmargin=*]
 \item $\sum_{n\ge 1} c_n h_n$ converges in $H$;
 \item $\sum_{n\ge 1} |c_n|^2 < \infty.$
\end{enumerate}
In this situation the {\em Parseval identity}\index{Parseval identity} $\n\sum_{n\ge 1} c_n h_n \n^2 =\sum_{n\ge 1} |c_n|^2$ holds.
\end{proposition}

\begin{definition}[Orthonormal system, orthonormal basis]\label{def:orthogonal-system} Let $I$ be a set.
A family $(h_i)_{i\in I}$ in $H$ is called an {\em orthonormal system}\index{orthonormal!system} if for all $ i,j\in I$ we have \index{$D$@$\delta_{ij}$}  $$\iprod{h_i}{h_j} = \delta_{ij}:=\begin{cases}1 & \hbox{if} \ \ i=j \\ 0 &\hbox{otherwise.}\end{cases}$$
In the case of a finite or countably infinite set $I$, an orthonormal system $(h_i)_{i\in I}$
is called an {\em orthonormal basis}\index{orthonormal!basis} if every $x\in H$ can be represented as a finite or convergent series
\begin{align}\label{eq:basis-hi}
x = \sum_{i\in I } c_i h_i
\end{align}
for suitable coefficients $c_i\in \K$.
\end{definition}

The sum \eqref{eq:basis-hi} is understood as a finite sum if $I$ is finite. If $I$ is countably infinite, we pick an enumeration $(i_n)_{n\ge 1}$ of $I$ and ask for convergence of $\sum_{n\ge 1} c_{i_n} h_{i_n}$. By Parseval's theorem,
this sum converges if and only if $\sum_{n\ge 1} |c_{i_n}|^2<\infty$. Thus convergence and the value of the sum are independent of the enumeration.
Indeed, if $(j_m)_{m\ge 1}$ is another enumeration and
$x:= \sum_{n\ge 1} c_{i_n} h_{i_n}$, $y:=\sum_{m\ge 1} c_{j_m} h_{j_m}$,
then for any $i\in I$, say $i = i_{n} = j_{m}$,
$$\iprod{x}{h_i}= \iprod{x}{h_{i_n}} = c_{i_n} = c_i = c_{j_m} = \iprod{y}{h_{j_m}} = \iprod{y}{h_i}.$$
Since both $x$ and $y$ belong to the closed linear span of the family $(h_i)_{i\in I}$, this implies $x=y$.
This argument also shows that
the coefficients $c_i$ in \eqref{eq:basis-hi} are uniquely determined by $x$ and given by $c_i = \iprod{x}{h_i}$.

\begin{example}
 The standard unit vectors $(0,\dots,0,1,0,\dots)$, $n\ge 1$ (with the $1$ in the $n$th place),
form an orthonormal basis for $\ell^2$\!.

In Section \ref{sec:ONB-examples} we prove that the trigonometric system is an orthonormal basis for $L^2(0,1)$
and the (suitably normalised) Hermite polynomials form an orthonormal basis for $L^2(\R,\gamma)$, where $\gamma$ is the standard Gaussian measure on $\R$.
\end{example}

\begin{theorem}[Orthonormal bases]\label{thm:ONB-dense} Let $(h_n)_{n\ge 1}$ be an orthonormal sequence in $H$.
The following assertions are equivalent:
\begin{enumerate}[label={\rm(\arabic*)}, leftmargin=*]
 \item\label{it:ONB-dense1} $(h_n)_{n\ge 1}$ is an orthonormal basis;
 \item\label{it:ONB-dense2} $(h_n)_{n\ge 1}$ has dense linear span;
 \item\label{it:ONB-dense3} if $x\in H$ satisfies $\iprod{x}{h_n}=0$ for all $n\ge 1$, then $x=0$.
\end{enumerate}
\end{theorem}

\begin{proof}
The equivalence \ref{it:ONB-dense2}$\Leftrightarrow$\ref{it:ONB-dense3} is immediate from the fact that a subspace is dense if and only if its orthogonal complement is trivial.

\smallskip
 \ref{it:ONB-dense1}$\Rightarrow$\ref{it:ONB-dense2}: \ This implication is trivial, because by assumption every $x\in H$ can be approximated by the partial sums of
 a series representation of the form $\sum_{n\ge 1} c_nh_n$, and these partial sums belong to the linear span
 of $(h_n)_{n\ge 1}$.

\smallskip
 \ref{it:ONB-dense2}$\Rightarrow$\ref{it:ONB-dense1}: \ Suppose the linear span of $(h_n)_{n\ge 1}$ is dense in $H$ and fix an arbitrary $x\in H.$ We must prove that $x$ admits a representation as a convergent sum $\sum_{n\ge 1} c_n h_n$.

For each $N\ge 1$ the mapping
$$P_N x:= \sum_{n=1}^N \iprod{x}{h_n}h_n$$ is a projection that maps $H$ onto the span $H_N$ of $(h_n)_{n=1}^N$.
If $x\perp H_N$, then $\iprod{x}{h_n}=0$ for $n=1,\dots,N$ and therefore $P_Nx =0$. It follows that the projection $P_N$ is orthogonal. This implies that $P_N$ is contractive. Therefore, with $c_n:= \iprod{x}{h_n}$,
$$ \sum_{n=1}^N |c_n|^2 = \n P_N x \n^2 \le \n x\n^2\!.$$
Since this holds for all $N\ge 1$, it follows that $ \sum_{n\ge 1} |c_n|^2 < \infty$ and therefore
the sum $y:= \sum_{n\ge 1} c_n h_n$ is convergent in $H$ by Proposition \ref{prop:Parseval}.
For all $n\ge 1$ we have
$ \iprod{y}{h_n} = c_n = \iprod{x}{h_n}$, and since the span of $(h_n)_{n\ge 1}$ is dense in $H$ this implies $x= y = \sum_{n\ge 1} c_n h_n$.
\end{proof}

 When $(x_n)_{n\ge 1}$ is a (finite or infinite) linearly independent sequence in $H$, we may construct a (finite or infinite) orthonormal sequence $(h_n)_{n\ge 1}$ with the property that $${\rm span}\{x_1,\dots,x_k\} = {\rm span}\{h_1,\dots,h_k\}, \quad k\ge 1,$$ as follows. Set $h_1:= x_1/\n x_1\n$. Suppose the orthonormal vectors $h_1,\dots,h_k$ have been chosen subject to the condition that $H_j:= {\rm span}\{x_1,\dots,x_j\}$ equals ${\rm span}\{h_1,\dots,h_j\}$ for all $j=1,\dots,k$. By linear independence, the subspace $H_{k+1}:= {\rm span}\{x_1,\dots,x_{k+1}\}$ has  dimension $k+1$. The orthogonal complement in $H_{k+1}$ of the $k$-dimensional subspace $H_k$ has dimension $1$, and therefore we may select a norm-one vector $h_{k+1}\in H_{k+1}$ orthogonal to $H_k$. Then $H_{k+1} = {\rm span}\{h_1,\dots,h_{k+1}\}$ as desired. This procedure is called {\em Gram--Schmidt orthogonalisation}.\index{Gram--Schmidt orthogonalisation}

\begin{theorem}[Orthonormal bases and separability]\label{thm:orth-bases-sep}
A Hilbert space has an orthonormal basis if and only if it is separable.
\end{theorem}
\begin{proof}
 `If': \ If $H$ is finite-dimensional we apply Gram--Schmidt to a basis of $H$. Otherwise choose a sequence with dense span, pass to a linearly independent subsequence with the same span, and apply Gram--Schmidt. The resulting orthonormal sequence has dense span, hence is an orthonormal basis by Theorem \ref{thm:ONB-dense}.

\smallskip
`Only if': \ If $(h_i)_{i\in I}$ is an orthonormal basis of $H$, the set of finite linear combinations of the $h_i$ with coefficients in $\Q$ if $\K=\R$ and in $\Q+i\Q$ if $\K=\C$ is countable and dense.
\end{proof}

\begin{corollary} Any two infinite-dimensional separable Hilbert spaces are isometrically isomorphic.
\end{corollary}
\begin{proof}
Suppose that the Hilbert spaces $H_1$ and $H_2$ are separable and pick orthonormal bases $(h_n^{(1)})_{n\ge 1}$ and $(h_n^{(2)})_{n\ge 1}$. The linear operator $U$, defined on the span of $(h_n^{(1)})_{n\ge 1}$ by sending $h_n^{(1)}$ to $h_n^{(2)}$ for each $n\ge 1$, has dense range. It is isometric by the Parseval identity, and therefore extends uniquely to an isometry $U:H_1\to H_2$. Its range is dense and closed, and therefore $U$ is surjective.
\end{proof}

\begin{definition}[Maximal orthonormal systems]\label{def:maximal-orth-system}
A {\em maximal orthonormal system}\index{orthonormal!system, maximal} is a family $(h_i)_{i\in I}$,
where $I$ is a nonempty set and the following conditions hold:
\begin{enumerate}[label={\rm(\roman*)}, leftmargin=*]
 \item\label{it:ONS1} $\iprod{h_i}{h_j} =\delta_{ij}$ for all $i,j\in I$;
 \item\label{it:ONS2} if $h\perp h_i$ for all $i\in I$, then $h=0$.
\end{enumerate}
\end{definition}

In a separable Hilbert space, every maximal orthonormal system is countable and can therefore be relabelled as an orthonormal basis.

\begin{theorem}[Maximal orthonormal systems]\label{thm:maximal-orth-system}
 Every nonzero Hilbert space has a maximal orthonormal system.
\end{theorem}
\begin{proof}
Partially order the set of all orthonormal systems in the nonzero Hilbert space $H$ by set inclusion. Since the union of every chain of orthonormal systems is again an orthonormal system and is an upper bound for the chain, Zorn's lemma (Theorem \ref{lem:Zorn}) applies and gives that
this set has a maximal element, say $(h_i)_{i\in I}$, where $I$ is some index set. It is clear that
condition \ref{it:ONS1} in the above definition holds. If there were a nonzero $h\in H$ perpendicular to each $h_i$, after normalising $h$ to unit length, we would obtain a new orthonormal system properly containing $(h_i)_{i\in I}$, contradicting the maximality of $(h_i)_{i\in I}$. Therefore also condition \ref{it:ONS2} holds.
\end{proof}

\section{Examples}\label{sec:ONB-examples}

In this final section we present several important orthonormal bases. The
examples illustrate how the abstract Hilbert space results of the preceding
section lead to concrete expansions: Fourier series on the circle and on
intervals, tensor product bases, and the Hermite basis for $L^2(\R,\gamma)$.

\subsection{The Trigonometric System}\label{subsec:Fourier-basis}

We will work over the complex scalar field.
In this example $\T$ denotes the unit circle in the complex plane, parametrised by the interval $[-\pi,\pi]$ and equipped with the normalised Lebesgue measure ${\rm d}\theta/2\pi$. We will prove that the functions
$$ e_n(\theta):= \exp(in\theta),\quad  \theta\in [-\pi,\pi], \ n\in\Z,$$
form an orthonormal basis for $L^2(\T)$.

That $(e_n)_{n\in \Z}$ is an orthonormal sequence in $L^2(\T)$ is evident from
\begin{align*} \iprod{e_j}{e_k} & = \frac1{2\pi}\int_{-\pi}^\pi \exp(ij\theta) \ov{\exp(ik\theta)}\ud \theta
= \frac1{2\pi}\int_{-\pi}^\pi{\exp(i(j-k)\theta)}\ud \theta = \delta_{jk}.
\end{align*}
By Theorem \ref{thm:ONB-dense}, to prove that $(e_n)_{n\in\Z}$ is an orthonormal basis, it remains to prove that the functions of the form
$ \sum_{n=-N}^N c_n e_n$ are dense in $L^2(\T)$.
This can be deduced from the Stone--Weierstrass theorem,
but we prefer the following argument from Fourier analysis which gives explicit approximants and some error bounds.

\begin{definition}[Fourier coefficients]\label{def:Fourier-coefficients} The {\em Fourier coefficients}\index{Fourier!coefficient}
of a function $f\in L^1(\T)$ are defined as
$$ \wh f(n) := \frac1{2\pi}\int_{-\pi}^\pi f(\theta)\exp(-in\theta)\ud \theta, \quad n\in\Z.$$
\end{definition}

\begin{theorem}\label{thm:Cesaro-convergence} For all
$f\in C(\T)$
we have
$$  \lim_{N\to\infty} \Bigl\n f - \frac1{N} \sum_{n=0}^{N-1}\sum_{k=-n}^n \wh f(k) e_k \Bigr\n_\infty = 0.$$
\end{theorem}
\begin{proof}
By scaling, it suffices to prove this for $f\in C(\T)$ with $\n f\n_\infty = 1$.  Fixing such a function, we have
$$ \wh f(n) \exp(in\theta) = \frac1{2\pi}\int_{-\pi}^\pi \exp(in(\theta-\si))f(\si)\ud \si =\frac1{2\pi}\int_{-\pi}^\pi \exp(in\si)f(\theta-\si)\ud \si$$
and therefore
\begin{equation}\label{eq:fNKN}f_N(\theta):= \frac1{N} \sum_{n=0}^{N-1}\sum_{k=-n}^n \wh f(k) \exp(ik\theta) = \frac1{2\pi}\int_{-\pi}^\pi K_N(\si)f(\theta-\si)\ud \si,
\end{equation}
where the {\em Fej\'er kernel}\index{Fej\'er kernel} $K_N$ is defined by
$K_N(0)=N$ and
\begin{equation}\label{eq:Fejer} K_N(\theta):= \frac1{N} \sum_{n=0}^{N-1}\sum_{k=-n}^n \exp(ik\theta) = \frac1N \frac{\sin^2(\frac12 N \theta)}{\sin^2(\frac12 \theta)}, \quad \theta\not=0;
\end{equation}
the identity on the right-hand side is readily deduced from the geometric series. In particular, $K_N\ge 0$ pointwise. Note that \eqref{eq:fNKN} expresses the averaged partial sums as a normalised convolution $K_N*f$ on $\mathbb T$. In view of the fact that $\frac1{2\pi}\int_{-\pi}^\pi K_N(\theta)\ud \theta = 1$, we have
\begin{equation}\label{eq:Fejer1}
\begin{aligned}
|f_N(\theta) - f(\theta)| & = \Bigl|\frac1{2\pi}\int_{-\pi}^\pi K_N(\si)(f(\theta-\si) - f(\theta))\ud \si\Bigr|
\\ & \le \frac1{2\pi}\int_{-\pi}^\pi K_N(\si)|f(\theta-\si) - f(\theta)|\ud \si.
\end{aligned}
\end{equation}
Fix $\eps>0$ and choose $0<\delta<\pi$ so small that $\n f(\cdot-\si) - f\n_\infty<\eps$ for all
$|\si|<\delta$; this is possible since $f$ is uniformly continuous.
Then
\begin{align}\label{eq:Fejer2}
\frac1{2\pi}\int_{-\delta}^\delta K_N(\si) |f(\theta-\si) - f(\theta)|\ud \si
\le \frac{\eps}{2\pi}\int_{-\pi}^\pi K_N(\si)\ud \si = \eps
\end{align}
and, by \eqref{eq:Fejer} and the normalisation $\n f\n_\infty =1$,
\begin{align}\label{eq:Fejer3}
\frac1{2\pi}\int_{\complement (-\delta,\delta)} K_N(\si)|f(\theta-\si) - f(\theta)| \ud \si
\le  \frac{2}{N} \frac1{\sin^2(\frac12 \delta) }.
\end{align}
Combining \eqref{eq:Fejer1} with \eqref{eq:Fejer2} and \eqref{eq:Fejer3} we obtain
$$ \n f_N - f\n_\infty \le \eps +  \frac{1}{N} \frac2{\sin^2(\frac12 \delta) } ,$$
so $\limsup_{N\to\infty} \n f_N - f\n_\infty \le \eps$. Since $\eps>0$ was arbitrary,
this completes the proof.
\end{proof}

This theorem implies that the trigonometric polynomials are dense in $C(\T)$.
Since this space is dense in $L^2(\T)$ by Remark \ref{rem:CKdenseLpK},
it follows that the trigonometric polynomials are dense in $L^2(\T)$. Therefore, Theorem \ref{thm:ONB-dense}
implies:

\begin{theorem} The trigonometric system\index{trigonometric!system} forms an orthonormal basis for $L^2(\T)$.
\end{theorem}

The theory of orthonormal bases can now be applied. It entails that every function $f\in L^2(\T)$
has a unique series representation of the form $f = \sum_{n\in \Z} c_n e_n$, with convergence in $L^2(\T)$ and coefficients given by
$c_n = \iprod{f}{e_n} = \wh f(n)$. The resulting expansion
$$ f = \sum_{n\in \Z} \wh f(n) e_n$$
is called the {\em Fourier series}\index{Fourier!series} of $f$. By Parseval's identity we have
$$ \frac1{2\pi}\int_{-\pi}^\pi |f(\theta)|^2\ud \theta= \sum_{n\in \Z} |\wh f(n)|^2\!,$$
and the mapping $f \mapsto (\wh f(n))_{n\in\Z}$ is an isometry from $L^2(\T)$ onto $\ell^2(\Z)$.

\medskip
By translation and scaling, the functions $$\wt e_n(\,\cdot):= \exp(2\pi i n\,\cdot),  \quad n\in\Z,$$ form an orthonormal basis for $L^2(0,1)$. This can be used to prove:

\begin{corollary}[Euler]\index{Euler's identity $\sum_{n=1}^\infty \frac1{n^2} = \frac{\pi^2}{6}$}
$\displaystyle \sum_{n=1}^\infty \frac1{n^2} = \frac{\pi^2}{6}.$
\end{corollary}
\begin{proof}
For the function $f(t) = t $ in $L^2(0,1)$ we have
$\iprod{f}{\wt e_0} = \frac12$ and $\iprod{f}{\wt e_n} = -\frac1{2\pi i n}$ for $n\neq 0$, so by Parseval's identity
\begin{align*}
\frac13 = \int_{0}^1 t ^2\ud t  = \n f\n^2
 = \sum_{n\in\Z} |\iprod{f}{\wt e_n}|^2
 = \bigl(\frac12\bigr)^2 + 2 \sum_{n=1}^\infty \bigl(\frac1{2\pi  n}\bigr)^2
= \frac14 + \frac1{2\pi^2} \sum_{n=1}^\infty \frac1{n^2},
\end{align*}
and the result follows.
\end{proof}

Another proof will be given in Section \ref{subsec:Euler}.

\begin{figure}
\begin{center}
\begin{minipage}{0.47\textwidth}
\begin{center}
\begin{tikzpicture}[scale=2]
 \draw[very thin,color=gray] (-1.15,-1.15) grid (1.15,1.15);
 \draw[->] (-1.15,0) -- (1.2,0) node[right] {$t$};
 \draw[->] (0,-1.15) -- (0,1.2);
 \draw[very thick, domain=0:1] plot (\x, {(\x+0.5)/2});
 \draw[very thick, domain=-1:0] plot (\x, -{((-\x)+0.5)/2});
 \draw (-1,-0.08) node[fill=white] {$-1$};
 \draw (0,-0.08)  node[fill=white] {$0$};
 \draw (1,-0.08)  node[fill=white] {$1$};
 \node at (0,-1.32) {odd extension};
\end{tikzpicture}
\end{center}
\end{minipage}
\hfill
\begin{minipage}{0.47\textwidth}
\begin{center}
\begin{tikzpicture}[scale=2]
 \draw[very thin,color=gray] (-1.15,-1.15) grid (1.15,1.15);
 \draw[->] (-1.15,0) -- (1.2,0) node[right] {$t$};
 \draw[->] (0,-1.15) -- (0,1.2);
 \draw[very thick, domain=0:1] plot (\x, {(\x+0.5)/2});
 \draw[very thick, domain=-1:0] plot (\x, {((-\x)+0.5)/2});
 \draw (-1,-0.08) node[fill=white] {$-1$};
 \draw (0,-0.08)  node[fill=white] {$0$};
 \draw (1,-0.08)  node[fill=white] {$1$};
 \node at (0,-1.32) {even extension};
\end{tikzpicture}
\end{center}
\end{minipage}
\caption{Left: the odd extension $f_{\rm odd}$ of a function $f$ from
$(0,1)$ to $(-1,1)$, satisfying
$f_{\rm odd}(-t)=-f_{\rm odd}(t)$.
Right: the even extension $f_{\rm ev}$, satisfying
$f_{\rm ev}(-t)=f_{\rm ev}(t)$.
The odd extension is orthogonal to the even functions
$\one$ and $\cos(\pi nt)$, while the even extension is orthogonal
to the odd functions $\sin(\pi nt)$.
For two functions of the same parity, their inner product over
$(-1,1)$ is twice the corresponding inner product over $(0,1)$.}
\label{fig:odd-even-extension}
\end{center}
\end{figure}

\smallskip
The system of functions
\begin{align*}
\one, \sqrt{2} \sin(2\pi n\,\cdot), \  \sqrt{2}\cos(2\pi n\,\cdot): \ \  n\ge 1,
\end{align*}
is orthonormal in $L^2(0,1)$ and the trigonometric functions $\exp(2\pi i n \,\cdot )$, $n\in\Z$, are contained in its linear span. Hence this system forms an orthonormal basis by Theorem \ref{thm:ONB-dense}. Interestingly, from this we can deduce:

\begin{theorem}\label{thm:L2bases-sin-cos} Each one of the two systems
\begin{align}\label{eq:trig-sin}  \sqrt{2}\sin(\pi n\,\cdot): \ \    n\ge 1,
\end{align}
and
\begin{align}\label{eq:trig-cos} \one, \sqrt{2} \cos(\pi n\,\cdot): \ \  n\ge 1,
\end{align}
forms an orthonormal basis for $L^2(0,1)$.
\end{theorem}

These bases arise naturally as the eigenvector bases of the Dirichlet and Neumann Laplacians in $L^2(0,1)$, respectively (see Example \ref{ex:ev-Dir}).

\begin{proof}
The two systems are orthonormal by direct computation. By Theorem \ref{thm:ONB-dense}, it remains to prove that in both cases the only function orthogonal to all elements of the system is the zero function.

Translating and rescaling the preceding trigonometric basis on $L^2(0,1)$ gives the orthonormal basis
$2^{-1/2}\one, \sin(\pi n\,\cdot), \cos(\pi n\,\cdot)$, $n\ge 1$, for $L^2(-1,1)$. Suppose first that $f\in L^2(0,1)$ is orthogonal to all functions $\sin(\pi n\,\cdot)$, $n\geq 1$. Let $f_{\rm odd}\in L^2(-1,1)$ be the odd extension of $f$. Then $f_{\rm odd}$ is orthogonal to $\one$ and to all functions $\cos(\pi n\,\cdot)$, $n\geq 1$. Moreover, for all $n\geq 1$,
$$
\iprod{f_{\rm odd}}{\sin(\pi n\,\cdot)}_{L^2(-1,1)}
=2\iprod{f}{\sin(\pi n\,\cdot)}_{L^2(0,1)}=0.
$$
Thus $f_{\rm odd}$ is orthogonal to an orthonormal basis of $L^2(-1,1)$, and therefore $f_{\rm odd}=0$. It follows that $f=0$. This proves that the system \eqref{eq:trig-sin} is an orthonormal basis for $L^2(0,1)$.

Next suppose that $f\in L^2(0,1)$ is orthogonal to $\one$ and to all functions $\cos(\pi n\,\cdot)$, $n\geq 1$. Let $f_{\rm ev}\in L^2(-1,1)$ be the even extension of $f$. Then $f_{\rm ev}$ is orthogonal to all functions $\sin(\pi n\,\cdot)$, $n\geq 1$. Furthermore, its inner products with $\one$ and with the functions $\cos(\pi n\,\cdot)$ are twice the corresponding inner products on $(0,1)$, and hence vanish. Thus $f_{\rm ev}$ is orthogonal to the same orthonormal basis of $L^2(-1,1)$. Hence $f_{\rm ev}=0$, and therefore $f=0$. This proves that the system \eqref{eq:trig-cos} is an orthonormal basis for $L^2(0,1)$.
\end{proof}

\subsection{The Hermite Polynomials}\label{subsec:Hermite}

In this section we prove that the (suitably normalised) Hermite polynomials form an orthonormal basis for $L^2(\R,\gamma)$, where
$\gamma$\index{$Gamma1$@$\gamma$} is the standard Gaussian measure\index{Gaussian!measure, standard} on $\R$. This is the Borel probability measure on $\R$ which is given, for Borel sets $B\subseteq \R$, by
$$\gamma(B) = \frac1{\sqrt{2\pi}} \int_B \exp\Bigl(-\frac12 x^2\Bigr)\ud x.$$
The Hermite polynomials will resurface in Chapters \ref{chap:spectral-theorem}, \ref{chap:semigroups}, and \ref{ch:SymmQuant} in connection with the spectral theorem, the Ornstein--Uhlenbeck semigroup, and second quantisation, respectively.

\begin{definition}
For $n\in\N$ the {\em Hermite polynomials}\index{Hermite!polynomials} $H_n:\R\to\R$ are defined by the identity
\begin{align}\label{eq:Hermite-gen-fc} H(t,x):= \exp\Bigl(tx-\frac12 t^2\Bigr) = \sum_{n=0}^\infty \frac{t^n}{n!}H_n(x),\quad t,x\in\R.
\end{align}
\end{definition}

The first five Hermite polynomials are given by
\begin{align*}
H_0(x) &= 1,\\
H_1(x) &= x,\\
H_2(x) &= x^2-1,\\
H_3(x) &= x^3 - 3 x, \\
H_4(x) &= x^4 - 6x^2 + 3.
\end{align*}
By induction one shows that
$$H_n(x)=\sum_{k=0}^{\lfloor n/2\rfloor } \frac{(-1)^k}{2^k} \frac{n!}{k!(n-2k)!}x^{n-2k}\!, \quad n\in\N.$$

\begin{proposition}\label{prop:Hermite-pols}
The Hermite polynomials have the following properties:
\begin{enumerate}[label={\rm(\roman*)}, leftmargin=*]
\item\label{it:prop:Hermite-pols-1} ${H_n}(-x)=(-1)^n H_n(x)$;
\item\label{it:prop:Hermite-pols-2} $H_{n+2}(x) = x H_{n+1}(x)-(n+1)H_n(x)$;
\item\label{it:prop:Hermite-pols-3} $H_{n+1}'(x) = (n+1)H_n(x)$;
\item\label{it:prop:Hermite-pols-4} $H_n$ is a monic polynomial of degree $n$.
\end{enumerate}
\end{proposition}
\begin{proof}
Property \ref{it:prop:Hermite-pols-1} follows from the identity $H(t,-x) = H(-t,x)$,
\ref{it:prop:Hermite-pols-2}
from the identity $\frac{\partial H}{\partial t}(t,x) =
(x-t)H(t,x),$ and \ref{it:prop:Hermite-pols-3}
from $\frac{\partial H}{\partial x}(t,x) =
t H(t,x).$ Assertion \ref{it:prop:Hermite-pols-4} follows from \ref{it:prop:Hermite-pols-2} and the
fact that $H_0(x)=1$ and $H_1(x)=x$.
\end{proof}

\begin{theorem}\label{thm:Hermite-ONB} The sequence
  $(\frac1{\sqrt{n!}}{H_n})_{n\in\N}$ forms an orthonormal basis for
$L^2(\R,\gamma)$.
\end{theorem}
\begin{proof}
For all $s,t\in\R$ we have
\begin{equation}\label{eq:Hermite-c}
\begin{aligned}
\int_{-\infty}^\infty H(s,x)H(t,x) \ud \gamma(x)
 & =  \frac{1}{\sqrt{2\pi}}
\int_{-\infty}^\infty  \exp\Bigl(-\frac12({s^2}+{t^2}) + (s+t)x-\frac12 x^2\Bigr)\ud x
\\ & =  \frac{1}{\sqrt{2\pi}} \exp(st)
\int_{-\infty}^\infty  \exp\Bigl(-\frac12 (s+t-x)^2\Bigr)\ud x
\\ & =  \frac{1}{\sqrt{2\pi}} \exp(st)
\int_{-\infty}^\infty  \exp\Bigl(-\frac12y^2\Bigr)\ud y= \exp(st).
\end{aligned}
\end{equation}
Taking derivatives
$\frac{\partial^{m+n}}{\partial s^m \partial t^n}$
at $s=t=0$
on both sides of the identity in \eqref{eq:Hermite-c}
gives
$$\int_{-\infty}^\infty H_m(x)H_n(x)\ud \gamma(x)= m!\delta_{mn}.$$
Since $m!\delta_{mn} = \sqrt{m!}\sqrt{n!}\delta_{mn}$, this shows that the sequence $(\frac1{\sqrt{n!}}{H_n})_{n\in\N}$ is orthonormal in $L^2(\R,\gamma)$.

It remains to show that the span of the Hermite polynomials is dense in $L^2(\R,\gamma)$.
A quick proof is obtained by making use of the injectivity of the Fourier transform (which is an immediate consequence of Theorem \ref{thm:FT-inversion}).
If $f\in L^2(\R,\gamma)$ is orthogonal to every Hermite polynomial, then
it is orthogonal to every polynomial. From this it follows that for all $z\in \C$,
$$ F(z) := \frac1{\sqrt{2\pi}}\int_{-\infty}^\infty  f(x) e^{zx-\frac12 x^2} \ud x
= \sum_{k=0}^\infty \frac{z^k}{k!}
\int_{-\infty}^\infty  f(x) x^k \ud \gamma(x) = 0,$$
where we use dominated convergence to justify interchanging sum and integration.
In particular, we have $F(-it)=0$. From this we infer that the Fourier transform of $x\mapsto f(x)e^{-\frac12 x^2}$ vanishes identically.
By the injectivity of the Fourier transform, this implies that $f(x)e^{-\frac12 x^2} = 0$ for almost all $x\in\R$, and
therefore $f(x) = 0$ for almost all $x\in\R$.
\end{proof}

The identity  $H_{n+2}(x) = x H_{n+1}(x)-(n+1)H_n(x)$ of Proposition \ref{prop:Hermite-pols}
is an example of a so-called {\em three-term recurrence relation}. As we will see in Section \ref{sec:orthpol},
orthogonal polynomials in $L^2(\R,\mu)$, with $\mu$ a finite Borel measure on $\R$, always satisfy a three-term recurrence relation, and conversely if a sequence of polynomials on $\R$ satisfies such a relation, then under a mild additional assumption these polynomials are orthogonal in $L^2(\R,\mu)$ for a suitable finite Borel measure $\mu$ on $\R$.

Theorem \ref{thm:Hermite-ONB} admits an extension to higher dimensions; see
Section \ref{sec:SQ}.

\subsection{Tensor Bases}\label{subsec:tensorbases}\index{tensor!basis}

Let $\mu_j$ be a finite Borel measure on a compact metric space $K_j$ for each $j=1,\dots,k$, and let $K = K_1\times\cdots\times K_k$ and $\mu = \mu_1\times\cdots\times \mu_k$ be their products. If $(f_n^{(j)})_{n\ge 1}$ is an orthonormal basis for $L^2(K_j,\mu_j)$ for each $j=1,\dots,k$, then the functions $$ f_{n}(x):= f_{n_1}^{(1)}(x_1)\cdots f_{n_k}^{(k)}(x_k), \quad n\in \{1,2,\dots\}^k\!,$$ form an orthonormal
basis for $L^2(K,\mu)$. Orthonormality being clear, in view of Theorem \ref{thm:ONB-dense} it remains to check that
the span of the functions $f_n$ is dense. This follows from the fact that $C(K)$ is dense in $L^2(K,\mu)$ by the observation in Remark \ref{rem:CKdenseLpK} and the fact that the linear span of functions of the form
$g(x):= g^{(1)}(x_1)\cdots g^{(k)}(x_k)$
with $g^{(j)}\in C(K_j)$ for all $j=1,\dots, k$ is dense in $C(K)$ by Example \ref{ex:tensorCK-dense}. Since each of the functions $g^{(j)}$ can be approximated in $L^2(K_j,\mu_j)$ by linear combinations of the functions $f_n^{(j)}$\!, $g$ can be
approximated in $L^2(K,\mu)$ by linear combinations of the functions $f_n$.

This is a special case of a more general construction involving tensor products of Hilbert spaces (see Chapters \ref{chap-HilbertSchmidt-TraceClass} and \ref{chap:QM}, in particular \eqref{eq:L2-tensor}): if $(h_n^{(j)})_{n\ge 1}$ is an orthonormal basis for the Hilbert space $H_j$, $j=1,\dots,k$, then the vectors $$h_n := h_{n_1}^{(1)}\ot\cdots\ot h_{n_k}^{(k)}, \quad n = (n_1,\dots,n_k)\in \{1,2,\dots\}^k,$$ form an orthonormal basis for the Hilbert space tensor product $H = H_1\otimes\cdots\otimes H_k$.

\begin{problems}
\item
Show that in a Hilbert space, equality $|\iprod{x}{y}| = \n x\n \n y\n$ holds in the Cauchy--Schwarz inequality
if and only if $x$ and $y$ are linearly dependent.

\item
Determine when two elements in a Hilbert space satisfy $\n x + y\n = \n x \n + \n y\n$.

\item
Let $(x_n)_{n\geq 1}$ be a sequence in a Hilbert space $H$. Suppose that there exists $x\in H$ such that:
\begin{enumerate}[label={\rm(\roman*)}, leftmargin=*]
  \item  $\iprod{x_n}{y} \to \iprod{x}{y}$ for all $y\in H$;
  \item  $\|x_n\|\to\|x\|$.
\end{enumerate}
Show that $x_n\to x$ in $H$.

\item\label{prob:Hilbert-strictly-convex}
A Banach space $X$ is called {\it strictly convex}\index{strictly convex} if for all norm-one vectors
$x_0,x_1\in X$ with $x_0\neq x_1$ and $0<\la<1$ we have $\n (1-\la) x_0 + \la x_1\n <1$.
Prove that every Hilbert space is strictly convex.

\item\label{prob:minimiser-C1}
Give an example of a nonempty compact convex set $C$ in a two-dimen\-sional Banach space $X$ along with a vector $x\in X$ such that the set
 $$ \Bigl\{c\in C:\, \n x - c\n = \min_{y\in C}\n x-y\n\Bigr\}$$
consists of more than one element.

\item\label{prob:minimiser-C2}
Let
$$ C:= \Big\{f\in C[0,1]:\, f \ \hbox{is real-valued}, \ f(0)=0, \ \int_0^1 f(t)\ud t = 0\Big\}.$$
\begin{enumerate}[\rm(a), leftmargin=*]
  \item Check that $C$ is a closed and convex subset of $C[0,1]$.
\end{enumerate}
Let $g\in C[0,1]$ be the function defined by $g(t):= t$.
\begin{enumerate}[\rm(a), leftmargin=*]\setcounter{enumii}{1}
  \item Show that for any $f\in C$ we have $\n f - g\n > \frac12$.
  \item Show that $\inf_{f\in C} \n f - g\n = \frac12.$
\end{enumerate}
This shows that $C$ contains no point minimising the distance $d(g,C)$.

\item
Show that a projection $P\in \calL(H)$ is an orthogonal projection if and only if $$\n Ph\n \leq \n h\n, \quad h \in H.$$

\noindent {\em Hint:}\ The latter condition implies
that $\n P(g+ch)\n ^2 \leq \n g+ch\n ^2$ for all $c\in \mathbb{K}$ and $g,h \in H$. Now consider
$g \in \ran(P)$ and $h \in \ker(P)$, and vary $c$.

\item
Show that two orthogonal projections $P$ and $Q$ commute if and only if $PQ$ is an orthogonal projection.

\item
Let $P$ and $Q$ be orthogonal projections on $H$. Show that
$\Ran(P)\subseteq\Ran(Q)$ if and only if $QP=P$ and $PQ=P$.

\item
Let $P$ and $Q$ be orthogonal projections on $H$. Show that $P+Q$
is an orthogonal projection if and only if
$\Ran(P)\perp\Ran(Q).$
In this case, identify the range of $P+Q$.

\item
In this problem we determine some orthogonal complements.
\begin{enumerate}[\rm(a), leftmargin=*]
  \item Let
  $$ Y := \bigl\{f\in L^2(0,1):\, f(t)=0  \ \text{ for almost all }\ t\in \bigl(0,\tfrac 12\bigr)\bigr\}.$$
  Show that $Y$ is a closed subspace of $L^2(0,1)$ and find $Y^{\perp}$\!.
  \item
  Let
  $$ Y := \Bigl\{f\in L^2(0,1):\, \int_0^1f(t)\ud t=0\Bigr\}. $$
  Show that $Y$ is a closed subspace of $L^2(0,1)$ and find $Y^{\perp}$\!.
\end{enumerate}

\item\label{prob.orthXY}
For any two closed subspaces $X$ and $Y$ of a Hilbert space $H$, show:

\begin{enumerate}[\rm(a), leftmargin=*]
\item\label{it.orthXY1} $(X+ Y)^\perp = X^\perp \cap Y^\perp$.
\item\label{it.orthXY2} $(X\cap Y)^\perp=\ov{X^\perp+Y^\perp}$.
\end{enumerate}

\item
Let $Y_1\subseteq Y_2\subseteq\cdots$ be closed subspaces of $H$, let
$Y:=\ov{\bigcup_{n\ge1}Y_n}$,
and denote by $P_n$ and $P$ the orthogonal projections onto $Y_n$
and $Y$, respectively.

\begin{enumerate}[label={\rm(\roman*)}, leftmargin=*]
 \item Show that $P_nx\to Px$ for all $x\in H$.
 \item Deduce that if $(Y_n)_{n\ge1}$ is a decreasing sequence of closed subspaces, then
$P_{n}x\to P_{\cap}x$ for all $x\in H$, where $P_{\cap}$ is the orthogonal projection onto $\bigcap_{n\ge1}Y_n$.
\end{enumerate}

\item
Let $Y$ be a closed subspace of $H$. Show that for every $x\in H$,
$$
d(x,Y)
=
\sup_{\substack{h\in Y^\perp\\ \n h\n\le1}}
|\iprod{x}{h}|.
$$

\item
Let $0\ne\phi\in H^*$ and let $h\in H$ be its Riesz representative.
\begin{enumerate}[label={\rm(\roman*)}, leftmargin=*]
\item Show that $\Ker(\phi)=\{h\}^\perp$.
\item Show that, for every $x\in H$,
$$
d(x,\Ker(\phi))
=
\frac{|\phi(x)|}{\n\phi\n}.
$$
\end{enumerate}

\item
Let $Y$ be a closed subspace of $H$. Show that $\codim(Y)=1$ if and
only if $Y=\{h\}^\perp$ for some nonzero $h\in H$.

\item
Let $Y$ be a closed subspace of a Hilbert space $H$ and let
$\phi\in Y^*$. Show that there exists a unique $\wt\phi\in H^*$ such
that $\wt\phi|_Y=\phi$ and $\n\wt\phi\n=\n\phi\n$.
This provides a special Hilbert space instance of the extension problem
considered more generally in Chapter~\ref{ch:duality}.

\item
Show that if $(h_n)_{n\ge 1}$ is an infinite orthonormal sequence in $H$, then $\iprod{h_n}{y}\to 0$ for every $y\in H$.

\item
Let $(e_j)_{j\ge 1}$ be the sequence of standard unit vectors in $\ell^2$\!, and let
$F=\overline{\text{span}}\{e_{2n-1}:\, n\ge 1 \}$ and let $G=\overline{\text{span}}\{e_{2n-1} + \frac{1}{n} e_{2n}:\, n\ge 1 \}$.

\begin{enumerate}[\rm(a), leftmargin=*]
  \item Give explicit expressions for the orthogonal projections $P_F$ and $P_G$.
  \item Show that $F \cap G = \{0 \}$.
  \item Show that $F+G$ is dense, but not closed in $\ell^2$\!.
\end{enumerate}

\item
Show that if $H$ is a separable Hilbert space supporting a finite Borel measure\index{Borel!measure} $\mu$ that is invariant under every isometric isomorphism of $H$ and satisfies $\mu(B)>0$ for some open ball $B$ in $H\setminus\{0\}$, then $H$ is finite-dimensional.

\noindent{\em Hint:}\ Modify the solution to Problem \ref{prob:tr-inv-Borel}.

\item
Prove the identity \eqref{eq:Fejer}.

\item
Let $(h_n)_{n\ge1}$ be an orthonormal sequence in $H$ and put
$$
H_N:=\operatorname{span}\{h_1,\ldots,h_N\},
\qquad
P_Nx:=\sum_{n=1}^N\iprod{x}{h_n}h_n.
$$
Show that $P_Nx$ is the unique best approximation to $x$ from $H_N$
and that
$$
\n x-P_Nx\n^2
=
\n x\n^2-\sum_{n=1}^N|\iprod{x}{h_n}|^2.
$$
Deduce that the $N$th Fourier partial sum of
$f\in L^2(\T)$ is the best $L^2$-approximation to $f$ among
trigonometric polynomials of degree at most $N$.

\item\label{prob:unif-sep-H}
Let $(h_n)_{n\ge1}$ be a sequence of norm-one vectors in a Hilbert
space $H$, and suppose that for some $r>0$,
\begin{align}\label{eq:sep-r}
\n h_n-h_m\n\ge r,\quad n\neq m.
\end{align}
\begin{enumerate}[\rm(a), leftmargin=*]
\item\label{it:unif-sep-H1} Show that the standard unit vectors in $\ell^2$ have pairwise
distance $\sqrt2$.

\item\label{it:unif-sep-H2} Show that necessarily $r\le\sqrt2$. Thus the constant $\sqrt2$
in part \ref{it:unif-sep-H1} is optimal.

\noindent{\em Hint:}\ For $n\neq m$, use the assumption to estimate
$\Re\iprod{h_n}{h_m}$ from above. Then expand
$0\le \n h_1+\cdots+h_N\n^2
$
and let $N\to\infty$.
\end{enumerate}

\item
Prove the following binomial identity for the Hermite polynomials: for all $n\in\N$ and $x,y\in\R$,
$$ H_n(x+y) = \sum_{k=0}^n \binom{n}{k} x^{n-k}H_k(y).$$

\item
Prove the following formula for the Hermite polynomials: For all $n\in\N$ and $x\in\R$,
$$ H_n(x) = \frac1{\sqrt{2\pi}} \int_{-\infty}^\infty(x+iy)^n \exp\Bigl(-\frac12 y^2\Bigr)\ud y.$$

\item\label{prob:Laguerre}
Define the polynomials $L_n$, $n\in\N$, by the generating function expansion
$$ \frac{\exp(tx/(1+t))}{1+t} = \sum_{n\in\N} \frac{t^n}{n!} L_n(x).$$
These are the {\em Laguerre polynomials}\index{Laguerre polynomials} normalised so as to become monic.
\begin{enumerate}[\rm(a), leftmargin=*]
  \item Compute the polynomials $L_n$ for $n=0,1,2,3$.
  \item Show that the polynomials $L_n$ are monic, have degree $n$, and satisfy
  the recurrence relation $$L_{n+2}(x) = (x-2n-3) L_{n+1}(x) - (n+1)^2L_{n}(x), \quad n\in\N.$$
  \item Prove that the sequence $(\frac{1}{n!}L_n)_{n\ge 0}$ is an orthonormal basis of $L^2(\R_+, e^{-x}\ud x)$.

  \noindent{\em Hint:}\
For orthogonality, integrate the product of two generating functions.
For completeness, imitate the density argument in the proof of
Theorem~\ref{thm:Hermite-ONB}. If $f$ is orthogonal to all
polynomials, consider
$$
F(z):=\int_0^\infty f(x)e^{zx}e^{-x}\ud x,
\quad \Re z<\frac12,
$$
and use the vanishing of its derivatives at zero.
\end{enumerate}

\item\label{prob:disc-algebra2} \index{$A(\mathbb{D})$}
The {\em disc algebra} $A(\mathbb{D})$ is the closed subspace of the Banach space $C(\ov{\mathbb D})$ consisting of those functions that are holomorphic on $\mathbb{D}$.
\begin{enumerate}[\rm(a), leftmargin=*]
\item Explain why
$\n f\n_\infty = \sup_{\theta\in [-\pi,\pi]} |f(e^{i\theta})|.$
  \item Show that for all $f\in A({\mathbb D})$ and $z_0\in \mathbb{D}$ we have
  $$ f(z_0) = \frac1{2\pi i}\int_{\mathbb {T}} \frac{f(z)}{z-z_0}\ud z.$$
  \item Show that a function $\varphi\in C(\mathbb T)$ is the restriction of a function $f\in A(\mathbb D)$ if and only if $\wh\varphi(n)=0$ for all $n<0$, where
  $$ \wh\varphi(n) := \frac1{2\pi}\int_{-\pi}^\pi \varphi(e^{i\theta})e^{-in\theta}\ud \theta, \quad n\in\Z.$$
  \noindent{\em Hint:}\ For $0<r<1$ consider the Poisson kernel
$$P_r(\theta):=\frac{1-r^2}{1-2r\cos\theta+r^2}.$$
Show that its normalised convolution with $\varphi$ tends uniformly
to $\varphi$ as $r\uparrow1$; first check trigonometric polynomials
and then use their density in $C(\mathbb T)$.
If the negative Fourier coefficients of $\varphi$ vanish, show
that each such convolution is the boundary restriction of a
function holomorphic on a neighbourhood of $\ov{\mathbb D}$.
Use the maximum principle to pass to the limit.
\end{enumerate}

\item\label{prob:H2}
The {\em Hardy space} $H^2(\mathbb{D})$ is the vector space of all holomorphic functions on $\mathbb{D}$ of the form $\sum_{n\in\N} c_n z^n$ with $\sum_{n\in\N} |c_n|^2 < \infty$.
\begin{enumerate}[\rm(a), leftmargin=*]
  \item\label{it:H2-1} Prove that $H^2(\mathbb D)$ is a Hilbert space with respect to the norm
  $$ \n f\n_{H^2(\mathbb D)} := \Bigl(\sum_{n\in\N} |c_n|^2\Bigr)^{1/2}\!.$$
\end{enumerate}
Let the functions $e_n\in
L^2(\T)$ be defined by $e_n(\theta):= \exp(in\theta)$, $n\in\N$, $\theta\in [-\pi,\pi]$.
\begin{enumerate}[\rm(a), leftmargin=*]\setcounter{enumii}{1}
  \item\label{it:H2-2} Show that for all $f= \sum_{n\in\N} c_n z^n\in H^{2}(\mathbb{D})$ the sum $f|_{\T} := \sum_{n\in\N} c_n e_n$ converges in $L^2(\T)$ and that the mapping
  $$ f\mapsto  f|_{\T}$$
  defines an isometry from $H^2(\mathbb{D})$ onto the closed subspace of
$L^2(\mathbb{T})$ consisting of all functions whose negative Fourier coefficients
vanish.

  \item\label{it:H2-3} For holomorphic functions $f:\mathbb{D}\to\C$ with power series expansion $f(z) = \sum_{n\in\N} c_n z^n$ and  $0<r<1$ define
  $$ f_r:= \sum_{n\in\N} c_n r^n e_n.$$
  Show that for $f\in H^2(\mathbb{D})$ and $0<r<1$ we have $f_r\in L^2(\T)$ and $\n f_r\n_{ L^2(\T)} \le\n f|_{\T}\n_{L^2(\T)}$, and show that $$\lim_{r\uparrow 1} f_r = f|_{\T}\ \hbox{in  $L^2(\T)$}.$$

  \item\label{it:H2-4} Show that if $f: \mathbb D\to \C$ is holomorphic, then $f\in H^2(\mathbb{D})$ if and only if
  $\sup_{0<r<1} \n f_r \n _{L^2(\T)} < \infty,$ and in this case we have
  $$ \n f\n_{H^2(\mathbb{D})} = \n f|_{\T}\n_{L^2(\T)}  = \sup_{0<r<1} \n f_r\n _{L^2(\T)} .$$

  \item\label{it:H2-5} Show that for all $f\in H^2(\mathbb{D})$, $0<r<1$, and $\theta\in [-\pi,\pi]$ we have
  $$ f(re^{i\theta}) = \frac1{2\pi}\int_{-\pi}^\pi  f|_{\T}(\eta)P_r(\theta-\eta)\ud\eta,$$
  where the {\em Poisson kernel}\index{Poisson!kernel, for the disc} is given by
  $$ P_r(\theta) = \frac{1-r^2}{1 - 2 r \cos(\theta) + r^2}.$$
  \noindent{\em Hint:}\ Begin by showing that $P_r(\theta) =\sum_{n\in\Z} r^{|n|}e^{in\theta}$.
\end{enumerate}

\item\label{prob:H2-continued}
We continue our study of the space $H^2(\mathbb{D})$.
\begin{enumerate}[\rm(a), leftmargin=*]
  \item\label{it:H2-continued1} Show that for all $z_0\in \mathbb{D}$ the function
  $$ k_{z_0}(z):= \frac1{1-z\ov{z_0}}$$ belongs to $H^2(\mathbb{D})$ and
  $$ \n k_{z_0}\n_{H^2(\mathbb{D})} = \frac1{(1-|z_0|^2)^{1/2}}.$$
  \item\label{it:H2-continued2} Show that for all $f\in H^2(\mathbb{D})$ and $z_0\in \mathbb{D}$
  we have
  $$ f(z_0) = \iprod{f}{k_{z_0}}.$$
  \item Use parts \ref{it:H2-continued1} and \ref{it:H2-continued2}
  to show that if $f_n\to f$ in $H^2(\mathbb{D})$,
  then $f_n\to f$ uniformly on every compact subset of $\mathbb{D}$.
\end{enumerate}

\item\label{prob:A2}
On $\C$ we consider the measure $$\ud \gamma_{\,\C}(z) = \frac1\pi e^{-|z|^2}\ud z,$$
where $\!\ud z$ is the Lebesgue measure on $\C$.
\begin{enumerate}[\rm(a), leftmargin=*]
  \item Show that $\gamma_{\,\C}$ is a probability measure satisfying
  $\int_\C |z|^2\ud \gamma_{\,\C}(z) = 1$.
\end{enumerate}
Let $A^2(\C)$\index{$A^2(\C)$} denote the complex vector space of all entire functions $f:\C \to \C$ that are
square-integrable with respect to $\gamma_{\,\C}$,
$$ \int_\C |f(z)|^2\ud \gamma_{\,\C}(z)<\infty.$$
\begin{enumerate}[\rm(a), leftmargin=*]\setcounter{enumii}{1}
  \item Show that $A^2(\C)$ is a Hilbert space with respect to the inner product
  $$ \iprod{f}{g} = \int_\C f(z)\overline{g(z)}\ud \gamma_{\,\C}(z).$$
 \item Show that the functions $$e_n(z):= \frac{z^n}{\sqrt{n!}}, \quad n\in\N,$$
 form an orthonormal basis for $A^2(\C)$.
\end{enumerate}

\item\label{prob:A2-ctd}
In this problem we continue our study of the Hilbert space $A^2(\C)$.

\begin{enumerate}[\rm(a), leftmargin=*]
\item  Using the mean value property, show that for all $w\in \C$ the mapping $f\mapsto f(w)$ is continuous from $A^2(\C)$ to $\C$. Deduce that for all $w\in \C$ there exists a unique function $k_w\in A^2(\C)$ such that
  $$ f(w) = \int_\C f(z) \ov{k_w(z)}\ud\gamma_{\,\C}(z), \quad f\in A^2(\C), \ w\in \C.$$
  \item  Show that this function is given by $$k_w(z) = \exp(z\ov w),\quad w,z\in \C.$$
  \item  Show that the orthogonal projection $P$ in $L^2(\C,\gamma_{\,\C})$ onto $A^2(\C)$ is given by $$
  Pf(w) = \int_\C f(z)\ov{k_w(z)}\ud \gamma_{\,\C}(z), \quad  f \in L^2(\C,\gamma_{\,\C}), \ w\in\C.$$
\end{enumerate}

\item\label{prob:RN}\index{theorem!Radon--Nikod\'ym}
In this problem we outline an alternative proof of the existence part of the Radon--Nikod\'ym theorem (Theorem \ref{thm:RN}) based on Hilbert space methods.
Let $(\Om,\calF\!,\mu)$ be a $\sigma$-finite measure space and let the $\K$-valued measure $\nu$ on $(\Om,\calF)$ be absolutely continuous with respect to $\mu$.

We first assume that $\nu$ is a finite nonnegative measure.
\begin{enumerate}[\rm(a), leftmargin=*]
  \item\label{it:RN1} Show that there exists a measurable function $w \in L^1(\Om,\mu)$ such that $w(\om)>0$ for $\mu$-almost all $\om\in\Om$.
  \item\label{it:RN2}  Show that the mapping $f\mapsto \int_\Om f\ud \nu$ is bounded on $L^2(\Om,\la)$, where $\la$ is the finite measure  on $(\Om,\calF)$ given by $\la(F):= \nu(F)+ \int_F w\ud \mu$. Conclude that there exists a unique $h\in L^2(\Om,\la)$ such that
  $$ \int_\Om f\ud \nu = \int_\Om fh\ud \la, \quad f\in L^2(\Om,\la).$$
  \item\label{it:RN3}  Show that $0\le h\le 1$ for $\la$-almost all $\om\in\Om$.

  \noindent {\em Hint:}\ Apply the identity in part \ref{it:RN2}
  to $f = \one_F$ with $F\in \calF$\!.
  \item\label{it:RN4}  Show that $\mu(B)=0$, where $B = \{\om\in\Om:\, h(\om)=1\}$.

  \noindent {\em Hint:}\ Apply the identity in part \ref{it:RN2}
  to $f = \one_B$.
  \item\label{it:RN5} Show that there exists a nonnegative function $g\in L^1(\Om,\mu)$ such that
  $$\nu(F) = \int_F g\ud \mu, \quad F\in \calF\!.$$
  \noindent {\em Hint:}\ Apply the identity in part \ref{it:RN2}
  to the function $f = \one_F(1+h+\cdots+h^n)$ and use parts \ref{it:RN3}, \ref{it:RN4},
  and the monotone convergence theorem to show that the limit $g:= \limn (1+h+\cdots+h^n)hw$ exists $\mu$-almost everywhere and belongs to $L^1(\Om,\mu)$.
\end{enumerate}
This proves the Radon--Nikod\'ym theorem for finite nonnegative measures $\nu$.
\begin{enumerate}[\rm(a), leftmargin=*]\setcounter{enumii}{5}
  \item\label{it:RN6} Deduce from this the general case.
\end{enumerate}

\item\label{prob:CE-Hilbert}
This problem provides a Hilbert space construction of the conditional expectation, thus avoiding the use of the Radon--Nikod\'ym theorem.
Let $(\Om,\calF\!,\mu)$ be a probability space and let $\calG$ be a sub-$\sigma$-algebra of $\calF$\!.
Let $1\le p\le \infty$ and denote by $L^p(\Om,\calG)$ the subspace of $L^p(\Om)$ consisting of all $f\in L^p(\Om)$ having a
$\calG$-measurable pointwise defined representative. From Problem \ref{prob:LpG-closed} we recall that $L^p(\Om,\calG)$ is a closed subspace of $L^p(\Om)$.
\begin{enumerate}[\rm(a), leftmargin=*]
  \item\label{it:CEb} Let $P_\calG$ denote the orthogonal projection of $L^2(\Om)$ onto $L^2(\Om,\calG)$ and
  fix a function $f\in L^2(\Om)$. Prove that $P_\calG f$ is the unique element of $L^2(\Om,\calG)$
  such that the following identity holds for all $G\in \calG$:
  $$ \int_G f\ud\mu =  \int_G P_\calG f\ud\mu.$$
  \noindent {\em Hint:}\ $f-P_\calG f \perp {\bf 1}_G$ in $L^2(\Om)$.
  \item Prove that if $f\in L^2(\Om)$ satisfies $f\ge 0$ $\mu$-almost everywhere, then  $P_\calG f \ge 0$ $\mu$-almost everywhere.
  \item Prove that if $f\in L^2(\Om)$ satisfies  $0\le f\le {\bf 1}$ $\mu$-almost everywhere, then  $0\le P_\calG f \le {\bf 1}$ $\mu$-almost everywhere.
  \item\label{it:CEe} Prove that $|P_\calG f|\le P_\calG |f|$ $\mu$-almost everywhere.
  \item\label{it:CEf} Prove that $P_\calG$ restricts to a contractive projection of $L^\infty(\Om)$ onto $L^\infty(\Om,\calG)$
  and extends to a contractive projection of $L^1(\Om)$ onto $L^1(\Om,\calG)$, and that the properties
  described in parts \ref{it:CEb}--\ref{it:CEe}
  extend to functions in these spaces. Here, a {\em projection} is understood to be a bounded operator $P$ satisfying $P^2 = P$.
  \item Discuss the relation of this problem with Problem \ref{prob:CI-unitinterval}.
  \item Give an explicit expression for $P_\calG$ in each of the following two cases:
  \begin{enumerate}[label={\rm(\roman*)}, leftmargin=*]
     \item $\Omega = (0,1)$, $\calF = \mathscr{B}(\Omega)$, $\mu$ the Lebesgue measure, and
     $\calG = \{\emptyset, \Omega\}$;
     \item $\Omega = (-\frac12,\frac12)$, $\calF= \mathscr{B}(\Omega)$, $\mu$ the Lebesgue measure, and
     $\calG = \{B\in \calF:\, B = -B\}$.
  \end{enumerate}
  \item Use the Radon--Nikod\'ym theorem  (Theorem \ref{thm:RN}) to give an alternative proof of the existence of the projection $P_\calG$ in $L^1(\Om)$ of part \ref{it:CEf}.
\end{enumerate}

\item\label{prob:Ito}
This problem assumes familiarity with the language of probability theory.
Let $(B_t)_{t\in [0,1]}$ be a standard Brownian motion on a probability space $(\Om,\calF\!,\P)$. For each $t\in [0,1]$, let $\calF_t$ denote the $\sigma$-algebra generated by the family
$(B_s)_{s\in [0,t]}$.
\begin{enumerate}[\rm(a), leftmargin=*]
  \item Show that if $0\le s<t\le 1$, the increment $B_{t} - B_{s}$ is independent of $\calF_s$.
\end{enumerate}
Let $H_0^2((0,1)\times\Om)$ denote the subspace of $L^2((0,1)\times\Om)$ consisting of all stochastic processes $\xi = (\xi_t)_{t\in [0,1]}$ of the form
\begin{align*}
\xi_t(\om) = \sum_{n=0}^{N-1} a_n(\om) \one_{(t_{n}, t_{n+1}]}(t), \quad t\in [0,1],\ \om\in\Om,
\end{align*}
where $0=t_0<t_1<\dots<t_N = 1$ and each $a_n\in L^2(\Om)$ is $\calF_{t_n}$-measurable.
For such processes we define
$$ \int_0^1 \xi_t \ud B_t := \sum_{n=0}^{N-1} a_n (B_{t_{n+1}} - B_{t_{n}}).$$
\begin{enumerate}[\rm(a), leftmargin=*, resume]
\item Show that the above definition is independent of the chosen representation, and that for all $\xi \in H_0^2((0,1)\times\Om)$ we have $\int_0^1 \xi_t \ud B_t \in L^2(\Om)$ and
  $$ \Bigl\n\int_0^1 \xi_t \ud B_t\Bigr\n_{L^2(\Om)}^2 = \n \xi\n_{L^2((0,1)\times\Om)}^2\!.$$
\end{enumerate}
Let $H^2((0,1)\times\Om)$ denote the closure of $H_0^2((0,1)\times \Om)$ in $L^2((0,1)\times \Om)$.
\begin{enumerate}[\rm(a), leftmargin=*, resume]  \item Deduce that the mapping $\xi\mapsto \int_0^1 \xi_t \ud B_t$ admits a unique extension to an isometry
  from $H^2((0,1)\times \Om)$ into $L^2(\Om)$, the so-called {\em It\^o isometry}.\index{It\^o isometry}
\end{enumerate}
The random variable $\int_0^1 \xi_t\ud B_t$ is called the {\em It\^o stochastic integral}\index{stochastic integral} of $\xi$ with respect to the Brownian motion $(B_t)_{t\in [0,1]}$.

\end{problems}

%% file: ch04-Duality.tex
\chapter{Duality}\label{ch:duality}\index{duality}

\blfootnote{\rm \small \begin{quote}This book has been published by Cambridge University Press in the series `Cambridge Studies in Advanced Mathematics'. The present version is free to view and download for personal use only. Not for re-distribution, re-sale or use in derivative works.\newline {\copyright} Jan van Neerven\end{quote}}

\noindent
Continuous linear functionals allow us to study a Banach space through scalar quantities. For the classical spaces, these functionals often admit concrete descriptions; more generally, the Hahn--Banach theorems ensure that there are enough of them to distinguish points and to separate a point from a closed convex set. In this way duality connects the geometry of Banach spaces with adjoint operators and with forms of convergence and compactness weaker than those defined by the norm.

\section{Duals of the Classical Banach Spaces}\label{sec:duals-classical}

Recall that the {\em dual} of a Banach space $X$ is the Banach space
$X^*:= \calL(X,\K)$.
For $x\in X$ and $x^*\in X^*$, the scalar $x^*(x)\in \K$ is denoted by $\lb x,x^*\rb$, that is, we write
$$ x^*(x)=:\lb x,x^*\rb.$$
The Hahn--Banach theorems guarantee an abundance of nontrivial functionals in the dual of any Banach space. In many concrete situations, however, it is possible to completely describe the dual space. It will be our first task to do this for some classical Banach spaces discussed in Chapter \ref{ch:ClassicalBanach}.

\subsection{Finite-Dimensional Spaces}\label{subsec:dual-fd-spaces}

It is instructive to start with duality of finite-dimensional spaces. As we have seen, every nonzero finite-dimensional Banach space is isomorphic to $\K^d$ for some integer $d\ge 1$. The dual of $\K^d$ is determined as follows.

Every $\xi\in \K^d$ determines an element $\phi_\xi\in (\K^d)\s$ by the prescription
$$ \phi_{\xi}(x):= x\cdot \xi = \sum_{n=1}^d x_n \xi_n, \quad x\in \K^d\!.$$ Indeed, the Cauchy--Schwarz inequality implies
$|\phi_{\xi}(x)| \le \n x\n \n \xi\n$, from which it follows that $\phi_{\xi}$ is bounded and $\n \phi_{\xi}\n \le \n \xi\n$.
Conversely, every $\phi\in (\K^d)\s$ is of this form. To see this, let $e_1,\dots,e_d$ be the standard unit vectors of $\K^d$ and set
$\xi_n:= \phi(e_n)$. Then $\xi := (\xi_1,\dots,\xi_d)\in \K^d$ and, for all $x = (x_1,\dots,x_d) = \sum_{n=1}^d x_n e_n$,
$$
\phi(x)= \phi\Big(\sum_{n=1}^d x_n e_n\Big) =   \sum_{n=1}^d x_n \phi(e_n)  =
\sum_{n=1}^d x_n \xi_n = x\cdot \xi = \phi_{\xi}(x) .$$
It follows that $\phi=\phi_{\xi}$.
Moreover, $\n \xi\n^2 = \phi_\xi(\ov\xi)\le  \n \phi_\xi\n\n\xi\n$. Together with the inequality $\n \phi_{\xi}\n \le \n \xi\n$ it follows that $\n \phi_\xi \n = \n\xi\n$.

In summary, the correspondence $\phi_\xi \leftrightarrow \xi$ establishes an isometric isomorphism
$$ (\K^d)\s \simeq \K^d\!.$$

\subsection{Sequence Spaces}\label{subsec:dual-seqspaces}

The above proof scheme can easily be extended to identify the duals of the infinite-dimensional sequence spaces $c_0$ and $\ell^p$\!.
We begin by proving that the dual of $c_0$ can be identified with $\ell^1$\!. Every $\xi\in \ell^1$ determines an element $\phi_\xi\in (c_0)\s$ by the prescription
$$ \phi_{\xi}(x):= \sum_{n\ge 1} x_n \xi_n, \quad x\in c_0.$$ Indeed,
$$|\phi_{\xi}(x)| \le (\sup_{n\ge 1}|x_n|) \sum_{n\ge 1} |\xi_n| = \n x\n_{\infty}\n \xi\n_{1},$$
 so $\phi_{\xi}$ is bounded and $\n \phi_{\xi}\n \le \n \xi\n_1$.
Conversely, every $\phi\in (c_0)\s$ is of this form. To see this, let $(e_n)_{n\ge 1}$ be the sequence of standard unit vectors of $c_0$ and set
$\xi_n:= \phi(e_n)$. We claim that $\sum_{n\ge 1} |\xi_n|<\infty$. To see this, choose scalars $c_n\in \K$ of modulus one
such that $c_n \xi_n = |\xi_n|$. The sequence $(c_1,\dots,c_N, 0,0,\dots) = \sum_{n=1}^N c_n e_n$ belongs to $c_0$ and has norm one,
and $$\sum_{n=1}^N |\xi_n| =\sum_{n=1}^N  c_n \xi_n = \phi\Bigl(\sum_{n=1}^N c_n e_n \Bigr)\le \n \phi\n.$$
Since $N\ge 1$ was arbitrary, this establishes the claim, with bound $\n \xi\n_1\le \n \phi\n$.
It follows that $\xi = (\xi_1,\xi_2,\dots)$ belongs to $\ell^1$ and for all $x\in c_0$ we have
$$
\phi(x) = \phi\Big(\sum_{n\ge 1} x_n e_n\Big) = \sum_{n\ge 1} x_n \phi(e_n)
 = \sum_{n\ge 1} x_n \xi_n = \phi_{\xi}(x)  .$$
It follows that $\phi=\phi_{\xi}$, and the preceding bounds combine to give the norm equality $\n \phi_\xi\n = \n \xi\n_1$. In summary, the correspondence $\phi_\xi \leftrightarrow \xi$ establishes an isometric isomorphism
$$ (c_0)\s \simeq \ell^1\!.$$

\medskip
In much the same way one proves that the dual of $\ell^p$\!, $1\le p<\infty$, can be represented as $\ell^q$,
where $\frac1p+\frac1q = 1$. More precisely, every
element $\xi\in \ell^q$ defines a bounded functional $\phi_\xi\in (\ell^p)^*$ of norm $\n\phi_\xi\n\le \n\xi\n_q$ by the same formula as before, this time using H\"older's inequality
$$|\phi_\xi(x)| = \Big|\sum_{n\ge 1} x_n \xi_n\Big| \le \n x\n_p\n \xi\n_q.$$
Conversely, every bounded functional is of this form. To see this let $(e_n)_{n\ge 1}$ be the sequence of standard unit vectors of $\ell^p$ and set
$\xi_n:= \phi(e_n)$. We claim that $(\xi_n)_{n\ge 1}$ belongs to $\ell^q$. The case $p=1$ and $q=\infty$ is trivial,
for $|\xi_n| \le \n \phi\n \n e_n\n = \n \phi\n$, $n\ge 1$, so $\n \xi\n_\infty\le \n \phi\n$. Therefore we only consider the case $1<p<\infty$, in which case also $1<q<\infty$.
To prove that
$\sum_{n\ge 1} |\xi_n|^q<\infty$ it obviously suffices to show that
\begin{align}\label{eq:ell-q} \sum_{n=1}^N |\xi_n|^q \le \n \phi\n^q\!, \quad  N\ge 1.
\end{align}
Fix $N\ge 1$ and put $x^{(N)}:= (c_1|\xi_1|^{q/p}\!, \dots, c_N|\xi_N|^{q/p}\!, 0,0,\dots)$, where the scalars $c_n\in\K$ are chosen in such a way that $c_n \xi_n = |\xi_n|$.
This sequence belongs to $\ell^p$\!, with norm
$$ \n x^{(N)}\n_{p}^p = \sum_{n=1}^N |\xi_n|^{q}\!.$$
Since $\frac1p+\frac1q = 1$ implies $ \frac{q}{p} +1 =q$,
\begin{align*} \sum_{n=1}^N |\xi_n|^q
 = \Big|\sum_{n=1}^N c_n|\xi_n|^{{q/p}}\xi_n \Big|  = |\phi(x^{(N)})| \le \n x^{(N)}\n_p\n \phi\n = \Big(\sum_{n=1}^N |\xi_n|^{q}\Big)^{1/p}\n \phi\n
\end{align*}
and therefore
$(\sum_{n=1}^N |\xi_n|^{q})^{1/q} \le \n \phi\n$, using once more that $\frac1p+\frac1q = 1$.
This proves \eqref{eq:ell-q}.
Since $N\ge 1$ was arbitrary,
it follows that $\xi = (\xi_1,\xi_2,\dots)$ belongs to $\ell^q$ with norm $\n \xi\n_q\le  \n \phi\n$,
and for all $x\in \ell^p$ we have
$$
\phi(x) = \phi\Big(\sum_{n\ge 1} x_n e_n\Big) = \sum_{n\ge 1} x_n \phi(e_n)
 = \sum_{n\ge 1} x_n \xi_n = \phi_{\xi}(x) .$$
It follows that $\phi=\phi_{\xi}$, and the preceding bounds combine to give the norm equality $\n \phi_\xi\n = \n \xi\n_q$.
In summary, the correspondence $\phi_\xi \leftrightarrow \xi$ establishes an isometric isomorphism
$$ (\ell^p)\s \simeq \ell^q, \quad 1\le p<\infty, \ \frac1p+\frac1q=1.$$
At the end of Section \ref{sec:HB-extension} we show that this result does not extend to $p=\infty$.

\subsection{Spaces of Continuous Functions}\label{subsec:dual-CK}

For locally compact Hausdorff spaces $X$, we denote by $C_0(X)$\index{$C_0(X)$} the space of continuous functions $f:X\to\K$ vanishing at infinity, that is, for every $\eps>0$ there exists a compact set $K\subseteq  X$ such that $|f(x)| <\eps$ for all $x\in \complement K$. With respect to the supremum norm, $C_0(X)$ is a Banach space. Note that $C_0(X)=C(X)$ if $X$ is compact.

The main result of this section is the Riesz representation theorem, which identifies the dual of $C_0(X)$, with $X$ locally compact Hausdorff, as the space $M_{\rm R}(X)$ of $\K$-valued Radon measures on $X$.
A Borel measure $\mu\in M(X)$ is said to be {\em Radon}
if its variation $|\mu|$ is Radon
(see  Definition \ref{def:Radon}), that is, if for every Borel subset $B$ of $X$ and all $\eps>0$
there is a compact set $K\subseteq  X$ and an open set $U\subseteq  X$ such that $K\subseteq  B\subseteq  U$ and $|\mu|(U\setminus K)<\eps$; these properties are referred to as {\em inner regularity with respect to compact sets} and {\em outer regularity}, respectively.

As was shown in Section \ref{sec:measures}, the space $M(X)$ of $\K$-valued Borel measures on $X$ is a Banach space with respect to the variation norm $\n \mu\n = |\mu|(X)$.

\begin{proposition} $M_{\rm R}(X)$ is a closed subspace of $M(X)$.
\end{proposition}
\begin{proof}
For $\mu,\nu\in M_{\rm R}(X)$ and $a,b\in\K$, we have
$|a\mu+b\nu|\le |a||\mu|+|b||\nu|$.
Taking unions of the compact approximating sets and intersections
of the open approximating sets shows that $a\mu+b\nu$ is Radon.
Thus $M_{\rm R}(X)$ is a linear subspace.

To prove closedness, suppose that $\mu_n\in M_{\rm R}(X)$ and
$\n\mu_n-\mu\n\to0$ in $M(X)$.
Fix a Borel set $B\subseteq X$ and $\eps>0$.
Choose $n$ such that $\n\mu-\mu_n\n<\eps/2$.
Since $\mu_n$ is Radon, there exist a compact set $K$ and
an open set $U$ such that
$K\subseteq B\subseteq U$ and $|\mu_n|(U\setminus K)<\eps/2$.
Consequently,
$$
|\mu|(U\setminus K)
\le |\mu-\mu_n|(U\setminus K)+|\mu_n|(U\setminus K)
<\eps.
$$
Hence $\mu$ is Radon.
\end{proof}

Every $\mu\in M(X)$ determines a bounded functional $\phi_\mu\in (C_0(X))\s$
given by $$\phi_\mu(f) := \int_X f\ud \mu.$$
By Proposition \ref{prop:int-compl-meas-ineq} it satisfies $$ |\phi_\mu(f)| \le \int_X |f|\ud |\mu| \le \n f\n_\infty \n \mu\n$$
and therefore
\begin{align*}
\n \phi_\mu\n \le \n \mu\n.
\end{align*}

\begin{theorem}[Riesz representation theorem]\label{thm:CK-dual}\index{theorem!Riesz representation, for $C_0(X)$}
For every $\phi \in (C_0(X))\s$ there exists a unique Radon measure $\mu\in M_{\rm R}(X)$ such that $\phi = \phi_\mu$, that is,
$$ \lb f,\phi \rb = \int_X f \ud\mu, \quad f\in C_0(X).$$
This measure satisfies $\n \mu\n = \n \phi\n$.
The correspondence $\phi\leftrightarrow \mu$ establishes an isometric isomorphism
$$ (C_0(X))\s \simeq M_{\rm R}(X).$$
The representing measure $\mu$ is nonnegative if and only if $\phi$ is positivity preserving.
\end{theorem}

In the special case of a compact metric space we have $M_{\rm R}(X) = M(X)$ by Proposition \ref{prop:Radon} and we obtain an isometric isomorphism
$$ (C_0(X))\s \simeq M(X).$$

For the proof of Theorem \ref{thm:CK-dual} we use Propositions \ref{prop:Urysohn2} and
\ref{prop:Cc-dense-in-C0} from the appendix, together with the following approximation theorem.

\begin{theorem}[Lusin approximation]\label{thm:Lusin}\index{theorem!Lusin}
Let $\mu$ be a finite Radon
measure on a locally compact Hausdorff space $X$. For every $f\in L^\infty(X,\mu)$ there exists a sequence
$(f_n)_{n\ge 1}$ in $C_{\rm c}(X)$ such that  $\limn f_n = f$ $\mu$-almost everywhere on $X$ and  $\sup_{n\ge 1} \n f_n\n_\infty \le \n f\n_\infty$.
\end{theorem}

\begin{proof}
The space $C_{\rm c}(X)$ is dense in $L^1(X,\mu)$.
This is proved in the same way as in Proposition \ref{prop:Cc-dense},
using the Radon property to obtain inner regularity with compact sets, and using the Urysohn functions of Proposition \ref{prop:Urysohn2} instead of those given by \eqref{eq:Ury-explicit}.

Now let $f\in L^\infty(X,\mu)$ and put $M:=\n f\n_\infty$.
Since $\mu$ is finite, we have $f\in L^1(X,\mu)$.
By the preceding paragraph, there exist $g_n\in C_{\rm c}(X)$
such that $g_n\to f$ in $L^1(X,\mu)$.
By Corollary \ref{cor:Lp-ae-subseq}, after passing to a subsequence
we may assume that $g_n\to f$ $\mu$-almost everywhere.

Define the continuous function $\rho:\K\to\K$ by
$$
\rho(z):=
\begin{cases}
z, & |z|\le M,\\
Mz/|z|, & |z|>M,
\end{cases}
$$
and put $f_n:=\rho\circ g_n$. Since $\rho(0)=0$, the support of
$f_n$ is contained in the support of $g_n$, and therefore
$f_n\in C_{\rm c}(X)$. Moreover, $\n f_n\n_\infty\le M$.
Since $|f|\le M$ almost everywhere, we have $\rho(f)=f$ almost
everywhere. The continuity of $\rho$ now gives
$f_n\to f$ $\mu$-almost everywhere, as required.
\end{proof}

We are now prepared for the proof of Theorem \ref{thm:CK-dual}.

\begin{proof}[Proof of Theorem \ref{thm:CK-dual}]
We first prove existence. The norm identity and uniqueness will be
established in Step 4.

\smallskip
{\em Step 1} --
We begin with the case of positivity-preserving functionals $\phi$. In this step we prove the existence of a nonnegative representing Radon measure $\mu\in M(X)$ for such functionals.
In what follows we freely use the notation and terminology introduced in Appendix \ref{sec:app-measures}.
Let $\mathscr U$ denote the collection of open subsets of $X$. For
$U\in\mathscr U$ and $f\in C_{\rm c}(X)$ we write $f\prec U$ if
$0\le f\le\one$ and ${\rm supp}(f)$ is compact and contained in $U$.
Define
$$
        \mu(U):=\sup\{\lb f,\phi\rb:\, f\in C_{\rm c}(X),\ f\prec U\},
        \quad U\in\mathscr U,
$$
with the convention $\mu(\emptyset)=0$. Then $0\le\mu(U)\le\|\phi\|$ for
all $U\in\mathscr U$, and $\mu$ is monotone.

We claim that $\mu$ is countably subadditive on $\mathscr U$. Indeed, let
$U_1,U_2,\ldots\in\mathscr U$ and suppose that
$f\prec\bigcup_{j\ge1}U_j$. Since ${\rm supp}(f)$ is compact, it is
contained in $\bigcup_{j=1}^kU_j$ for some $k$.
By the locally compact case of Theorem \ref{thm:partition-unity},
there exist functions $g_1,\ldots,g_k\in C_{\rm c}(X)$ such that
$g_j\prec U_j$ and $g_1+\cdots+g_k=1$ on ${\rm supp}(f)$. Hence
$fg_j\prec U_j$ and, by positivity of $\phi$,
$$
        \lb f,\phi\rb
        =\sum_{j=1}^k\lb fg_j,\phi\rb
        \le \sum_{j=1}^k\mu(U_j)
        \le \sum_{j\ge1}\mu(U_j).
$$
Taking the supremum over all $f\prec\bigcup_{j\ge1}U_j$ shows that
$$
        \mu\Bigl(\bigcup_{j\ge1}U_j\Bigr)
        \le \sum_{j\ge1}\mu(U_j).
$$

We next verify the closed-set approximation condition in Theorem~\ref{thm:open-Caratheodory}. Let
$U,V\in\mathscr U$ and let $\varepsilon>0$. Choose
$f\prec U\cap V$ such that
$\mu(U\cap V)<\lb f,\phi\rb+\varepsilon$.
Put $K:={\rm supp}(f)$. Then $K$ is closed and $K\subseteq U$. Choose
$g\prec V\setminus K$ such that
$\mu(V\setminus K)<\lb g,\phi\rb+\varepsilon$.
Since $f+g\prec V$, we obtain
\begin{align*}
        \mu(U\cap V)+\mu(V\setminus K)
        < \lb f,\phi\rb+\lb g,\phi\rb+2\varepsilon
        = \lb f+g,\phi\rb+2\varepsilon
         \le \mu(V)+2\varepsilon.
\end{align*}
As $\varepsilon>0$ was arbitrary, Theorem~\ref{thm:open-Caratheodory} applies.
Thus the formula
$$
        \mu^*(A):=\inf\{\mu(U):\, A\subseteq U,\ U\in\mathscr U\},
        \quad A\subseteq X,
$$
defines an outer measure whose restriction to $\mathscr B(X)$ is a finite
Borel measure satisfying $\mu(X)\le \n\phi\n$, which extends the original set function $\mu$ on open sets. We denote this Borel measure again by $\mu$.

We next show that $\mu$ is Radon. Outer regularity follows
from its construction. If $g\prec X$ and $L:={\rm supp}(g)$,
then $g\prec U$ for every open set $U$ containing $L$.
Consequently,
$$
\mu(L)=\inf_{\substack{U\supseteq L \\ U\ {\rm open}}}\mu(U)
\ge\lb g,\phi\rb.
$$
By the definition of $\mu(X)$, for every $\eps>0$ we can
therefore find a compact set $L$ such that
$\mu(X\setminus L)<\eps/2$.

Now let $B\subseteq X$ be Borel. By outer regularity, choose
an open set $U\supseteq X\setminus B$ such that
$\mu(U)<\mu(X\setminus B)+\eps/2$.
Then $\mu(U\cap B)<\eps/2$, and $K:=L\setminus U$ is a
compact subset of $B$ satisfying
$$
\mu(B\setminus K)
\le\mu(X\setminus L)+\mu(B\cap U)<\eps.
$$
Thus $\mu$ is also inner regular with respect to compact sets.

\smallskip
{\em Step 2} --
We prove that $\mu$ represents $\phi$. First let
$0\le f\in C_{\rm c}(X)$. Fix $\eps>0$ and choose $N\ge1$
such that $N\eps\ge\n f\n_\infty$. Put
$$
U_k:=\{f>k\eps\},\quad k=1,\ldots,N.
$$
For any functions $g_k\prec U_k$, we have $\eps\sum_{k=1}^N g_k\le f$.
Indeed, at a point $x$, only indices with $k\eps<f(x)$
can contribute, and each $g_k(x)$ is at most $1$.
By positivity of $\phi$,
$$
\lb f,\phi\rb
\ge\eps\sum_{k=1}^N\lb g_k,\phi\rb.
$$
Taking the supremum independently over each $g_k\prec U_k$
gives
$$
\lb f,\phi\rb\ge\eps\sum_{k=1}^N\mu(U_k).
$$
The simple function $\eps\sum_{k=1}^N\one_{U_k}$ differs
from $f$ by a function taking values in $[0,\eps]$.
Hence
$$
\int_X f\ud\mu
\le\eps\sum_{k=1}^N\mu(U_k)+\eps\mu(X)
\le\lb f,\phi\rb+\eps\mu(X).
$$
Letting $\eps\downarrow0$, we obtain
$\int_X f\ud\mu\le\lb f,\phi\rb$.

To prove the reverse inequality, put $C:=\n f\n_\infty$.
Given $\delta>0$, choose a compact set $K$ such that
$\mu(X\setminus K)<\delta$. By Proposition
\ref{prop:Urysohn2}, there exists $g\prec X$ such that
$g=1$ on $K\cup{\rm supp}(f)$.
Then $Cg-f$ is nonnegative and belongs to $C_{\rm c}(X)$.
Applying the inequality just proved to this function gives
\begin{align*}
\lb f,\phi\rb-\int_X f\ud\mu
\le C\Bigl(\lb g,\phi\rb-\int_X g\ud\mu\Bigr)
\le C\bigl(\mu(X)-\mu(K)\bigr)
<C\delta,
\end{align*}
where we used $\lb g,\phi\rb\le\mu(X)$ and
$\int_X g\ud\mu\ge\mu(K)$.
Since $\delta>0$ was arbitrary, equality follows.

By linearity, the representing identity holds for every
$f\in C_{\rm c}(X)$. Both functionals are bounded on $C_0(X)$,
so Proposition \ref{prop:Cc-dense-in-C0} extends the identity
to every $f\in C_0(X)$.

\smallskip
{\em Step 3} --
We first consider an arbitrary bounded real-linear functional
$\phi$ on $C_0(X;\R)$. For $0\le f\in C_0(X;\R)$ define
\begin{align}\label{eq:phi-pos}
\phi^+(f):=
\sup\bigl\{\lb g,\phi\rb:\, g\in C_0(X;\R),\ 0\le g\le f\bigr\}.
\end{align}
Then
$0\le\phi^+(f)\le\n\phi\n\n f\n_\infty$,
and $\phi^+(cf)=c\phi^+(f)$ for all scalars $c\ge0$.

We claim that $\phi^+$ is additive on the cone of nonnegative
functions. If $0\le g_1\le f_1$ and $0\le g_2\le f_2$, then
$0\le g_1+g_2\le f_1+f_2$, so
$$
\phi^+(f_1+f_2)
\ge\lb g_1+g_2,\phi\rb
=\lb g_1,\phi\rb+\lb g_2,\phi\rb.
$$
Taking the supremum over the admissible $g_1$ and $g_2$ gives
$\phi^+(f_1+f_2)\ge\phi^+(f_1)+\phi^+(f_2)$.
For the reverse inequality, let $0\le g\le f_1+f_2$ and put
$g_1:=f_1\wedge g$ and $g_2:=g-g_1$.
Then $0\le g_1\le f_1$ and $0\le g_2\le f_2$, and therefore
$$
\lb g,\phi\rb
=\lb g_1,\phi\rb+\lb g_2,\phi\rb
\le\phi^+(f_1)+\phi^+(f_2).
$$
Taking the supremum over $g$ proves the claim.

For arbitrary $f\in C_0(X;\R)$ we now define
$$
\phi^+(f):=\phi^+(f^+)-\phi^+(f^-).
$$
If $f=u-v$ with $u,v\ge0$, then $u+f^-=v+f^+$.
By additivity on the positive cone, it follows that
$$
\phi^+(u)-\phi^+(v)
=\phi^+(f^+)-\phi^+(f^-)=\phi^+(f).
$$
Thus the difference on the left depends only on $u-v$.
Together with additivity and positive homogeneity on the positive
cone, this proves that the extension is real-linear.
It is positive, and hence
$$
|\phi^+(f)|
\le\phi^+(|f|)
\le\n\phi\n\n f\n_\infty.
$$
In particular, $\phi^+$ is bounded.

The functional $\phi^-:=\phi^+-\phi$ is also bounded and
real-linear. For $f\ge0$, taking $g=f$ in \eqref{eq:phi-pos}
gives $\phi^+(f)\ge\lb f,\phi\rb$, so $\phi^-$ is positive.
We have therefore written $\phi=\phi^+-\phi^-$ as a difference
of two bounded positive functionals.
The preceding two steps apply over the real scalar field.
Taking the difference of the representing measures of
$\phi^+$ and $\phi^-$ gives a real-valued Radon measure
representing $\phi$.

Finally, suppose that $\phi$ is a bounded complex-linear
functional on $C_0(X;\C)$. On $C_0(X;\R)$ define
$$
a(u):=\Re\lb u,\phi\rb,\qquad
b(u):=\Im\lb u,\phi\rb.
$$
These are bounded real-linear functionals. By what we have
just proved, they have real-valued Radon representing measures
$\mu_a$ and $\mu_b$. Put $\mu:=\mu_a+i\mu_b$.
For every real-valued $u\in C_0(X)$ we then have
$$
\lb u,\phi\rb
=a(u)+ib(u)=\int_Xu\ud\mu.
$$
By complex linearity, this identity holds for every
$f\in C_0(X;\C)$, completing the existence proof.

\smallskip {\em Step 4} --
It remains to prove the norm identity and uniqueness. Let $\mu\in M_{\rm R}(X)$
represent $\phi$. By Example \ref{ex:K-val-meas}, there exists
$h\in L^1(X,|\mu|)$ such that
$$
        \mu(B)=\int_B h\ud|\mu|,
        \quad B\in\mathscr B(X),
$$
and $|h|=1$ $|\mu|$-almost everywhere.

Applying Theorem \ref{thm:Lusin} to $\ov h\in L^\infty(X,|\mu|)$, we obtain functions $f_n\in C_{\rm c}(X)$ such that $\n f_n\n_\infty\le 1$ and $f_n\to\ov h$ $|\mu|$-almost everywhere. By dominated convergence,
$$
\lb f_n,\phi\rb
  =\int_X f_n\ud\mu =\int_X f_nh\ud|\mu|
  \to \int_X |h|^2\ud|\mu| =|\mu|(X).
$$
It follows that $\n\phi\n\ge |\mu|(X)=\n\mu\n$.
The reverse inequality $\n\phi\n\le\n\mu\n$ was established before the
statement of the theorem, and hence $\n\phi\n=\n\mu\n$.

Finally, if $\mu_1,\mu_2\in M_{\rm R}(X)$ both represent $\phi$, then
$\mu_1-\mu_2$ represents the zero functional. Applying the norm identity
to this measure, we obtain $\n\mu_1-\mu_2\n=0$,
and therefore $\mu_1=\mu_2$.
\end{proof}

The duality between spaces of continuous functions and spaces of Borel measures shows how
elements of measure theory emerge naturally from considerations involving only linearity and topology (namely, from the problem of finding the continuous linear functionals of a space of continuous functions).

\subsection{Spaces of Integrable Functions}\label{subsec:LpLq-duality}

 Let $(\Omega, \calF\!,\mu)$ be a measure space and let $1\le p,q\le \infty$ satisfy $\frac1p+\frac1q=1$.
 By H\"older's inequality, every function $g\in L^q(\Om)$ defines a functional
 $\phi_g\in (L^p(\Om))\s$ by setting
 $$ \phi_g(f) := \int_\Omega f g\ud\mu, \quad f\in L^p(\Om),$$
 and we have $\n \phi_g\n\le \n g\n_q$.
 If $1\le p<\infty$ and the measure space is $\sigma$-finite, every functional arises in this way:

 \begin{theorem}[Dual of $L^p(\Om)$]\label{thm:Lp-dual}
 Let $(\Omega, \calF\!,\mu)$ be a $\sigma$-finite measure space and let $1\le p<\infty$ and $\frac1p+\frac1q = 1$.
 For every $\phi \in (L^p(\Om))\s$ there exists a unique $g\in L^q(\Om)$ such that $\phi = \phi_g$, that is,
 $$ \lb f,\phi \rb = \int_\Omega f g\ud\mu, \quad f\in L^p(\Om),$$
 and it satisfies $\n g\n_q = \n \phi\n$.
The correspondence $\phi_g \leftrightarrow g$ establishes an isometric isomorphism
$$ (L^p(\Om))\s \simeq L^q(\Om), \quad 1\le p<\infty, \ \frac1p+\frac1q=1.$$
\end{theorem}

\begin{proof}
Let $\phi\in(L^p(\Om))^*$.
By $\sigma$-finiteness, we can write $\Om$ as a disjoint union
of measurable sets $E_n$ of finite measure.
Put $F_N:=\bigcup_{n=1}^N E_n$.

\smallskip
{\em Step 1} --
For each $n\ge1$ define
$$
\nu_n(A):=\lb\one_A,\phi\rb,
\quad A\in\calF,\ A\subseteq E_n.
$$
This is a $\K$-valued measure on $E_n$. Indeed, if the measurable
sets $A_j\subseteq E_n$ are disjoint, dominated convergence gives
$\sum_{j=1}^m\one_{A_j} \to \one_{\bigcup_{j\ge1}A_j}$ in $L^p(\Om)$.
The countable additivity of $\nu_n$ follows from this by applying $\phi$.
Moreover, $\nu_n$ is absolutely continuous with respect to the
restriction of $\mu$ to $E_n$.

By the Radon--Nikod\'ym theorem (Theorem \ref{thm:RN}), there exists
$g_n\in L^1(E_n)$ such that
$$
\nu_n(A)=\int_Ag_n\ud\mu,
\quad A\in\calF,\ A\subseteq E_n.
$$
Choose measurable representatives and define $g$ on $\Om$ by
$g=g_n$ on $E_n$. Since the sets $E_n$ are disjoint, this defines
a measurable function, and $g\in L^1(F_N)$ for every $N$.

Let $Y$ denote the space of bounded measurable functions
vanishing outside some $F_N$. By the definition of the measures
$\nu_n$ and linearity, we have
\begin{align}\label{eq:Lp-duality-prelim3}
\lb f,\phi\rb=\int_\Om fg\ud\mu
\end{align}
for every simple function in $Y$.
If $f\in Y$ vanishes outside $F_N$, we may approximate it uniformly
by simple functions vanishing outside $F_N$.
Since $\mu(F_N)<\infty$ and $g\in L^1(F_N)$, we can pass to
the limit on both sides. Thus \eqref{eq:Lp-duality-prelim3}
holds for every $f\in Y$.

The space $Y$ is dense in $L^p(\Om)$.
Indeed, if $f\in L^p(\Om)$, the functions
$f_N:=\one_{F_N}\one_{\{|f|\le N\}}f$
belong to $Y$ and converge to $f$ in $L^p(\Om)$ by dominated
convergence.

\smallskip
{\em Step 2} --
We next prove that $g\in L^q(\Om)$ and $\n g\n_q\le\n\phi\n$.
Let $h\in Y$. Since $h$ is bounded and vanishes outside some
$F_N$, we have $hg\in L^1(\Om)$.
Choose a measurable function $u$ with $|u|\le1$ such that
$uhg=|hg|$ almost everywhere.
Then $uh\in Y$, and \eqref{eq:Lp-duality-prelim3} gives
$$
\n hg\n_1
=\int_\Om uhg\ud\mu
=|\lb uh,\phi\rb|
\le\n\phi\n\n uh\n_p
\le\n\phi\n\n h\n_p.
$$
Since $Y$ is dense in $L^p(\Om)$, Proposition
\ref{prop:Lp-via-Lq-1}, with the exponents $q$ and $p$
interchanged, implies that $g\in L^q(\Om)$ and
$\n g\n_q\le\n\phi\n$.

By H\"older's inequality, integration against $g$ defines
a bounded functional on $L^p(\Om)$. Since it agrees with $\phi$
on the dense subspace $Y$, it agrees with $\phi$ on all of
$L^p(\Om)$. H\"older's inequality also gives
$\n\phi\n\le\n g\n_q$, and therefore $\n g\n_q=\n\phi\n$.

Finally, suppose that $h\in L^q(\Om)$ also represents $\phi$.
For each $n$, H\"older's inequality implies that
$h|_{E_n}\in L^1(E_n)$, and for every measurable $A\subseteq E_n$
we have
$$
\int_A h\ud\mu=\lb\one_A,\phi\rb=\nu_n(A).
$$
By uniqueness in the Radon--Nikod\'ym theorem, $h=g_n$
almost everywhere on $E_n$. Since the sets $E_n$ cover $\Om$,
we conclude that $h=g$ almost everywhere.
\end{proof}

\begin{remark}\label{rem:no-sigma-finite}
If $1<p<\infty$, the $\sigma$-finiteness assumption in the preceding
theorem can be omitted. Indeed, let $\phi\in (L^p(\Om))^*$, choose a sequence
$(f_n)_{n\ge1}$ in the closed unit ball of $L^p(\Om)$ such that
$ |\lb f_n,\phi\rb|\to \n\phi\n$, and put
$$
 S:=\bigcup_{n,m\ge1}\{|f_n|>1/m\}.
$$
Then $S$ is $\sigma$-finite and all $f_n$ vanish outside $S$.

We claim that $\phi(g)=0$ for every $g\in L^p(\Om)$ vanishing on $S$.
After normalising we may assume that $\n g\n_p=1$. Multiplying $f_n$
and $g$ by suitable scalars of modulus one, for every $t>0$ we obtain
$$
 |\lb f_n,\phi\rb|+t|\lb g,\phi\rb|
 \le \n\phi\n\,\n f_n+tg\n_p
 \le \n\phi\n(1+t^p)^{1/p},
$$
since $f_n$ and $g$ have disjoint supports. Letting $n\to\infty$ gives
$$
 t|\lb g,\phi\rb|
 \le \n\phi\n\big((1+t^p)^{1/p}-1\big).
$$
Dividing by $t$ and letting $t\downarrow0$, using $p>1$, gives
$\lb g,\phi\rb=0$.

The restriction of $\phi$ to $L^p(S)$ therefore has norm $\n\phi\n$.
Since $S$ is $\sigma$-finite, the preceding theorem provides
$h\in L^q(S)$ such that $\n h\n_q=\n\phi\n$ and
$$ \lb f,\phi\rb=\int_S fh\ud\mu,\quad f\in L^p(S).$$
Extending $h$ by zero outside $S$ gives
$$
 \lb f,\phi\rb=\int_\Om fh\ud\mu,\quad f\in L^p(\Om).
$$
Thus $(L^p(\Om))^*=L^q(\Om)$ isometrically.
The corresponding assertion for $p=1$ may fail; see
Problem~\ref{prob:no-sigma-finite}.
\end{remark}

\subsection{Hilbert Spaces}

The Riesz representation theorem (Theorem \ref{thm:Riesz}) states that
every bounded functional on a Hilbert space $H$ is of the form $\psi_h$ for some unique $h\in H$, where
$$ \psi_h(g) = \iprod{g}{h}, \quad g\in H.$$
Moreover, we have equality of norms $\n h\n = \n \psi_h\n$.
The identification $\psi_h\leftrightarrow  h$
therefore provides a bijective and norm-preserving correspondence $$H^*\longleftrightarrow  H.$$
It is important to observe that this correspondence is linear if $\K=\R$, but conjugate-linear if $\K=\C$. This is a consequence of the conjugate-linearity of inner products with respect to their second variable. Indeed, from
$ \psi_{ch}(x) = \iprod{x}{c h} = \ov c \iprod{x}{h} = \ov c \psi_h(x)$
it follows that $$\psi_{ch} = \ov c \psi_h.$$
In contrast, the correspondence $\phi_x\leftrightarrow x$ in each of the Sections \ref{subsec:dual-fd-spaces}--\ref{subsec:LpLq-duality} is linear both when $\K =\R$ and $\K=\C$.

\subsection{Banach Lattices}

In this subsection we work over the real scalar field. The objective is to prove that the dual of a Banach lattice is a Banach lattice in a natural way. The proof relies on several lemmas.

If $V$ is a vector lattice, then for $0\le v\in V$ we write $[0,v]:= \{u\in V:\, 0\le u\le v\}$.

\begin{lemma}[Riesz decomposition]\label{lem:Riesz-dec}\index{Riesz decomposition}
  If $V$ is a vector lattice, then for all $0\le v,v'\in V$ we have
  $$ [0,v]+[0,v'] = [0,v+v'].$$
\end{lemma}
\begin{proof}
The inclusion $[0,v]+[0,v']\subseteq[0,v+v']$ is immediate.
For the reverse inclusion, let $0\le w\le v+v'$ and put
$u:=v\wedge w$ and $u':=w-u$.
Then $0\le u\le v$ and $u+u'=w$.
By the lattice identities in Proposition \ref{prop:orderedVS},
$$
u'=w-(v\wedge w)=(w-v)\vee0.
$$
Both $w-v$ and $0$ are bounded above by $v'$, so
$0\le u'\le v'$.
Thus $u\in[0,v]$ and $u'\in[0,v']$, as required.
\end{proof}

\begin{lemma}[Extension from $X^+$]\label{lem:cone-extension}
Let $X$ be a vector lattice and let $F: X^+\to\mathbb{R}$ be additive and
positively homogeneous, that is, for all $x,y\in X^+$ and scalars $c\ge 0$ we have
$$
F(x+y)=F(x)+F(y),\qquad F(c x)=c F(x).
$$
Then there exists a unique linear functional $\widetilde F: X\to\mathbb{R}$
such that $\widetilde F|_{X^+}=F$.
\end{lemma}

\begin{proof}
Define, for $x\in X$,
$$
\widetilde F(x):=F(x^+)-F(x^-),
$$
where $x = x^+ - x^-$ with $x^+$ and $x^-$ the positive and negative parts of $x$.
To see that this indeed extends $F$, note that if $x\in X^+$ then $x^+=x$ and
$x^-=0$, hence $\widetilde F(x)=F(x)$.

To check that $\widetilde F$ is well-defined in the sense that it agrees with
$F(u)-F(v)$ for \emph{any} representation $x=u-v$ with
$u,v\in X^+$, observe that $x=u-v$ implies $u+x^- = v+x^+$, and both belong to $X^+$.
By additivity of $F$ on $X^+$, this gives
$F(u)+F(x^-)=F(v)+F(x^+)$ and hence $F(u)-F(v)=F(x^+)-F(x^-)=\widetilde F(x)$.

Linearity of $\widetilde F$ now follows from the identity
$\widetilde F(u-v)=F(u)-F(v)$ and the additivity/homogeneity of $F$ on $X^+$.

Uniqueness follows by observing that if $L: X\to\mathbb{R}$ is linear with $L|_{X^+}=F$, then
$L(x)=L(x^+)-L(x^-)=F(x^+)-F(x^-)=\widetilde F(x)$ for all $x\in X$.
\end{proof}

We can now prove the main result of this section:

\begin{theorem}\label{thm:dual-BL} With respect to the natural partial order given by
$$ x\s\le y\s \Leftrightarrow \lb x, x\s\rb\le \lb x,y\s\rb \ \hbox{ for all } \ 0\le x\in X,$$
the dual $X\s$ of a Banach lattice $X$ is a Banach lattice. Moreover, for all $0\le x\in X$ we have
\begin{align*}
 \lb x,x\s\wedge y\s\rb & = \inf\bigl\{\lb x-y, x\s\rb + \lb y,y\s\rb:\, 0\le y\le x \bigr\}, \\
 \lb x,x\s\vee y\s \rb   & = \sup\bigl\{\lb x-y, x\s\rb + \lb y,y\s\rb:\, 0\le y\le x \bigr\}.
\end{align*}
\end{theorem}

A special case of the second formula has already been encountered in \eqref{eq:phi-pos}.

\begin{proof}
The relation $\le$ defines a partial order on $X\s$ which is compatible with addition
and multiplication by nonnegative scalars; antisymmetry follows
from the fact that every element of $X$ is a difference of two
positive elements.

\smallskip
{\em Step 1} --
Fix $\phi\in X^*$. For $x\in X^+$ define
$$
F_\phi(x):=\sup\bigl\{\lb y,\phi\rb:\, 0\le y\le x\bigr\}.
$$
Since $0\le y\le x$ implies $\n y\n\le\n x\n$, we have
$$
0\le F_\phi(x)\le\n\phi\n\n x\n.
$$
The mapping $F_\phi$ is positively homogeneous.
Moreover, Lemma \ref{lem:Riesz-dec} gives, for $x,x'\in X^+$,
\begin{align*}
F_\phi(x+x')
&=\sup\bigl\{\lb y+y',\phi\rb:\,
             0\le y\le x,\ 0\le y'\le x'\bigr\}\\
&=\sup\bigl\{\lb y,\phi\rb:\,0\le y\le x\bigr\}
 +\sup\bigl\{\lb y',\phi\rb:\,0\le y'\le x'\bigr\}
=F_\phi(x)+F_\phi(x').
\end{align*}
By Lemma \ref{lem:cone-extension}, $F_\phi$ extends uniquely
to a linear functional on $X$, which we denote by $\phi^+$.
This extension is positive. Since $-|x|\le x\le|x|$, we obtain
$$
|\lb x,\phi^+\rb|
\le\lb |x|,\phi^+\rb
=F_\phi(|x|)
\le\n\phi\n\n x\n.
$$
Thus $\phi^+\in X^*$.

\smallskip
{\em Step 2} --
We show that $\phi^+$ is the least upper bound of $\phi$ and $0$.
Taking $y=x$ and $y=0$ in the definition of $F_\phi$ shows that
$\phi^+\ge\phi$ and $\phi^+\ge0$.
If $\chi\in X^*$ is another upper bound, then for $0\le y\le x$
we have
$$
\lb y,\phi\rb\le\lb y,\chi\rb\le\lb x,\chi\rb,
$$
where the second inequality uses $\chi\ge0$.
Taking the supremum over $y$ gives
$\lb x,\phi^+\rb\le\lb x,\chi\rb$ for every $x\ge0$.
Hence $\phi^+\le\chi$, proving the claim.

Subtracting $\psi$ turns the upper bounds of $\{\phi,\psi\}$
into the upper bounds of $\{\phi-\psi,0\}$.
Similarly, $\chi$ is a lower bound of $\{\phi,\psi\}$ if and
only if $\phi-\chi$ is an upper bound of $\{0,\phi-\psi\}$.
Consequently, every pair $\phi,\psi\in X^*$ has a least upper
bound and a greatest lower bound, given by
$$
\phi\vee\psi=\psi+(\phi-\psi)^+,\qquad
\phi\wedge\psi=\phi-(\phi-\psi)^+.
$$
Thus $X^*$ is a vector lattice.

For $x\ge0$, the formula for the least upper bound gives
\begin{align*}
\lb x,\phi\vee\psi\rb
&=\lb x,\psi\rb
  +\sup\bigl\{\lb y,\phi-\psi\rb:\,0\le y\le x\bigr\}\\
&=\sup\bigl\{\lb y,\phi\rb+\lb x-y,\psi\rb:\,
             0\le y\le x\bigr\}\\
&=\sup\bigl\{\lb x-y,\phi\rb+\lb y,\psi\rb:\,
             0\le y\le x\bigr\},
\end{align*}
where the last equality follows by replacing $y$ by $x-y$.
This is the asserted formula for the supremum.
The formula for the infimum follows by applying it to
$(-\phi)\vee(-\psi)$ and taking negatives.

\smallskip
{\em Step 3} --
It remains to verify the lattice property of the dual norm.
Taking $\psi=-\phi$ in the formula just proved gives, for $x\ge0$,
\begin{align*}
\lb x,|\phi|\rb
=\sup\bigl\{\lb x-2y,\phi\rb:\,0\le y\le x\bigr\}
=\sup\bigl\{|\lb z,\phi\rb|:\,|z|\le x\bigr\}.
\end{align*}
Indeed, $z=x-2y$ ranges over the interval $[-x,x]$, which
equals $\{z\in X:\,|z|\le x\}$ and is symmetric under $z\mapsto-z$.

Now suppose that $|\phi|\le|\psi|$. For every $x\in X$, the
preceding identity gives
\begin{align*}
|\lb x,\phi\rb|
&\le\lb |x|,|\phi|\rb
\le\lb |x|,|\psi|\rb
=\sup\bigl\{|\lb z,\psi\rb|:\,|z|\le|x|\bigr\}
\le\n x\n\n\psi\n.
\end{align*}
The last inequality follows because $|z|\le|x|$ implies
$\n z\n\le\n x\n$.
Taking the supremum over $\n x\n\le1$, we see that
$\n\phi\n\le\n\psi\n$.
Since $X^*$ is complete in its usual norm, it is a Banach lattice.
\end{proof}

\section{The Hahn--Banach Extension Theorem}\label{sec:HB-extension}

We now turn to one of the main pillars of functional analysis, the Hahn--Banach theorem. This is a collective name for a number of closely related results, all of which assert the existence
of certain nontrivial functionals with desirable properties. These results come in two flavours: as extension theorems asserting the extendability\index{extension!of a functional} of functionals that are given  {\em a priori}  on a subspace and as separation theorems asserting that certain disjoint subsets can be separated by means of functionals.

The present section is concerned with Hahn--Banach extension theorems.
We begin with a version for {\em real} vector spaces whose proof exploits the order structure of the real line.

 \begin{theorem}[Hahn--Banach extension theorem for real vector spaces]\label{thm:appHB1}\index{theorem!Hahn--Banach extension, for vector spaces}Let $V$ be a real vector space and let
 $p:V\to\R$ be {\em sublinear}\index{sublinear},
 that is,  for all $v,v'\in V$ and $t\geq 0$ we have
 \begin{equation*}
    p(v+v')\leq p(v)+p(v'),\qquad p(tv)=tp(v).
 \end{equation*}
 If $W\subseteq  V$ is a subspace and $\phi:W\to\R$ is a linear mapping satisfying
 $$\phi(w)\leq p(w), \quad  w\in W,$$
 then there exists a linear mapping $\Phi:V\to \R$ that extends $\phi$, and satisfies
  $$
  \Phi(v)\leq p(v),\quad v\in V.$$
  \end{theorem}
  \begin{proof}
   We may assume that $W$ is a proper subspace of $V$. Fix $w,w'\in W$ and $v\in V\setminus W$.
   From $\phi(w)+\phi(w') = \phi(w+w')\le p(w+w')$ and $p(w+w')\leq p(w-v)+p(v+w')$
   we obtain $\phi(w) - p(w-v) \le p(w'+v)-\phi(w')$. By subadditivity, $\alpha:= \sup_{w\in W} (\phi(w) - p(w-v))$ is finite and we have
   $$ \phi(w)+\alpha \le p(w+v), \qquad \phi(w') - \alpha \le p(w'-v), \quad  w,w'\in W.$$
   Let $W_1$ denote the linear span of $W$ and $v$. Then $\phi_1: W_1\to \R$, $\phi_1(w+tv):= \phi(w)+t\alpha$,
   is well defined, is
   linear, extends $\phi$, and satisfies
   $\phi_1(w_1)\le p(w_1)$ for all $w_1\in W_1$. To see this,
   note that for $w\in W$ and $t>0$,
   $$ \phi_1(w+tv)= \phi(w)+t\alpha = t(\phi(t^{-1}w)+\alpha) \le tp(t^{-1}w+v) = p(w+tv)$$
   and
   $$ \phi_1(w-tv)= \phi(w)-t\alpha = t(\phi(t^{-1}w)-\alpha) \le tp(t^{-1}w-v) = p(w-tv)$$
   while for $t=0$ we have $ \phi_1(w) = \phi(w)$. This proves that $\phi_1(w_1)\le p(w_1)$ for all $w_1\in W_1$.
   The proof can now be finished by an appeal to Zorn's lemma (Theorem \ref{lem:Zorn}),
   applied to all linear extensions $\phi'$ of $\phi$ satisfying the inequality
   $\phi'\le p$. More precisely, if we order pairs $(W',\phi')$, where
   $W\subseteq W'\subseteq V$ and $\phi'$ is such an extension on $W'$,
   by extension, the union of every chain is again an admissible extension and
   hence is an upper bound for the chain. A maximal pair must have $W'=V$,
   since otherwise the one-dimensional extension constructed above produces a
   strictly larger admissible pair.
\end{proof}

This result is used to give a second version of the theorem which is also valid over the complex scalars.

   \begin{theorem}[Hahn--Banach extension theorem for vector spaces]\label{thm:appHB2}\index{theorem!Hahn--Banach extension, for vector spaces}Let $V$ be a (real or complex) vector space and let $p:V\to [0,\infty)$
  be a {\em seminorm}\index{seminorm}, that is,  for all $v,v'\in V$ and $t\in\K$ we have
  \begin{equation*}
     p(v+v')\leq p(v)+p(v'),\qquad p(tv)=|t|p(v).
  \end{equation*}
  If $W$ is a subspace of $V$ and $\phi:W\to\K$ is a linear mapping satisfying
  $$|\phi(w)|\leq p(w), \quad  w\in W,$$
  then there exists a linear mapping $\Phi:V\to \K$ that extends $\phi$ and satisfies
   $$ |\Phi(v)|\leq p(v), \quad   v\in V.$$
  \end{theorem}
\begin{proof} First we consider the case $\K=\R$. The assumptions imply $\phi(w)\leq p(w)$ for all $w\in W$,
and therefore by Theorem \ref{thm:appHB1} the mapping $\phi:W\to\R$ admits a linear extension $\Phi:V\to \R$ satisfying $\Phi(v)\leq p(v)$ for all $v\in V.$
Also, $-\Phi(v) = \Phi(-v) \le p(-v) = p(v)$, and therefore $|\Phi(v)| \le p(v)$ for all $v\in V$.

Next we consider the case $\K = \C$.
  Let us write  $\phi = \Re\phi+i \Im\phi$, where $\Re\phi$ and $\Im\phi$ are the real and imaginary parts of $\phi$;
  these functions are real-linear. From
$$    \phi(x) = -i\phi(ix) = -i ( \Re\phi(ix)+i \Im\phi(ix))= \Im\phi(ix)-i\Re\phi(ix),$$ with $ \Im\phi(ix)$ and $\Re\phi(ix)$ real,
we infer that $\Im\phi(x) = -\Re\phi(ix).$ Hence,
$$ \phi(x) = \Re\phi(x) - i\Re\phi(ix).$$
  The real-valued function $\psi:=\Re\phi$ satisfies the assumptions of Theorem \ref{thm:appHB1} and thus
  extends to a real-linear mapping $\Psi:V\to \R$ satisfying $\Psi\le p$. We now define
  $$\Phi(v):= \Psi(v)-i\Psi(iv).$$
  Then $\Phi$ extends $\phi$, $\Phi$ is real-linear, and since also $$\Phi(iv) = \Psi(iv)-i\Psi(-v) =  \Psi(iv)+i\Psi(v) = i(\Psi(v) - i\Psi(iv)) = i\Phi(v),$$
  $\Phi$ is actually complex-linear.
  Finally, for $v\in V$, choose $t\in\C$ such that $|t|=1$
and $t\Phi(v)=|\Phi(v)|$.
Since $\Phi$ is complex-linear and $\Re\Phi=\Psi$, we obtain
\[
|\Phi(v)| =\Re\Phi(tv) =\Psi(tv) \le p(tv) =p(v).
\qedhere\]
\end{proof}

For normed spaces we obtain the following version.

\begin{theorem}[Hahn--Banach extension theorem for normed spaces]\label{thm:HB}\index{theorem!Hahn--Banach extension, for normed spaces} Let\, $X$\, be\, a\, normed space and let $Y\subseteq  X$ be a subspace. Then every functional $y\s\in Y\s$ has an extension to a
functional $x\s\in X\s$ that satisfies
$$\n x\s\n =\n y\s\n.$$
\end{theorem}

Here, of course, $\n x\s\n$ is the norm of $x^*$ as an element of $X^*$ and $\n y\s\n$  is the norm of $y^*$ as an element of $Y^*$.
Such obvious conventions will be in place throughout the text.

\begin{proof}
Given a functional $y\s\in Y\s$\!, we apply Theorem \ref{thm:appHB2} to $V = X$, $W = Y$, $\phi(y):= \lb y,y\s\rb$ for $y\in Y$, and $p(x):= \n x\n\n y\s\n $
for $x\in X$.
This gives a functional $x^*\in X^*$ extending $y^*$ and satisfying
$$
   |\lb x,x^*\rb|\le \|x\|\,\|y^*\|,\quad x\in X.
$$
Hence $\|x^*\|\le \|y^*\|$. Since $x^*|_Y=y^*$, the reverse inequality
$\|y^*\|\le \|x^*\|$ is immediate. Therefore $\|x^*\|=\|y^*\|$.
\end{proof}

\begin{remark}
The proof of Theorem \ref{thm:appHB1} depends on the Axiom of Choice through the use of Zorn's lemma. If $X$ is separable and a countable dense sequence $(x_n)_{n\ge 1}$ is given, Theorem \ref{thm:HB} can be proved without invoking Zorn's lemma.

Over the real scalars, we revisit the proof of Theorem \ref{thm:appHB1}. Starting from a functional $y^*\in Y^*$, let $Y_n$ be the span of $Y$ and $\{x_1,\dots,x_n\}$ and extend successively from $Y_{n-1}$ to $Y_n$. On the union of the spaces $Y_n$ this gives a well-defined functional of norm at most $\n y^*\n$. Since this union is dense, Proposition \ref{prop:extendT} gives a unique extension to all of $X$ with the same norm.

Over the complex scalars, apply the real construction to $\Re y^*$ on the underlying real spaces and recover the complex-linear extension from the resulting real-linear functional $\Psi$ by the formula
$$
\Phi(x):=\Psi(x)-i\Psi(ix),\quad x\in X,
$$
as in the proof of Theorem \ref{thm:appHB2}.
\end{remark}

Recall the identity
$$ \n x\s\n = \sup_{\n x\n \le 1} |\lb x,x\s\rb|,$$
which is nothing but the definition of the operator norm of $x\s$ as an element of $\calL(X,\K)$.
As a consequence of the Hahn--Banach theorem
we obtain the following dual expression for the norm of elements $x\in X$:

\begin{corollary}\label{cor:HB-norming} For all $x\in X$ we have
 $$ \n x\n = \sup_{\n x\s\n \le 1} |\lb x,x\s\rb|.$$
In particular, if $\lb x,x\s\rb = 0$ for all $x\s\in X\s$\!, then $x=0$.
\end{corollary}
\begin{proof} Fix $x\in X$. If $x=0$, the asserted identity trivially holds, so we may assume that $x\not=0$.
 Let $Y$ be the one-dimensional subspace of $X$ spanned by $x$ and define $y_0\s\in Y\s$ by $\lb tx,y_0\s\rb:= t\n x\n$.
 Then $\n y_0\s\n =1$. Let $x_0\s\in X\s$ be a Hahn--Banach extension provided by Theorem \ref{thm:HB}, that is, $x_0\s|_Y = y_0\s$
 and $\n x_0\s\n = \n y_0\s\n = 1$. Then
$$\n x\n = \lb x,y_0\s\rb = \lb x,x_0\s\rb \le  \sup_{\n x\s\n \le 1} |\lb x,x\s\rb|,$$
while trivially
\[
\sup_{\n x\s\n \le 1} |\lb x,x\s\rb|\le \sup_{\n x\s\n \le 1}\n x\n \n x\s\n = \n x\n.
\qedhere\]
\end{proof}

The next application of Theorem \ref{thm:HB} provides a condition for recognising proper closed subspaces.

\begin{corollary}\label{cor:proper-subspace}
 If $Y$ is a proper closed subspace of a Banach space $X$, then for every $x_0\in  X\setminus Y$ there exists an $x\s\in X\s$ such that
 $$ \lb x_0,x^*\rb \not=0 \  \hbox{ and }  \ \lb y,x^*\rb = 0 \ \hbox { for all } \ y\in Y.$$
\end{corollary}

\begin{proof}
Let $q:X\to X/Y$ denote the quotient mapping. Since $x_0\notin Y$, we have $qx_0\not=0$. By Corollary \ref{cor:HB-norming}, applied to the normed space $X/Y$, there exists a $\psi\in (X/Y)\s$ such that $\lb qx_0,\psi\rb\not=0$. Define $x\s\in X\s$ by $\lb x,x\s\rb:=\lb qx,\psi\rb$. Then $x\s$ is bounded, since $q$ and $\psi$ are bounded. Moreover, $\lb y,x\s\rb=0$ for all $y\in Y$, whereas $\lb x_0,x\s\rb=\lb qx_0,\psi\rb\not=0$. This proves the result.
\end{proof}

\begin{definition}[Complemented subspaces] A closed linear subspace $X_0$ of a normed space $X$ is said to be {\em complemented}\index{complemented}
if there exists a closed linear subspace $X_1$ of $X$ such that
\begin{itemize}
 \item $X_0 + X_1 = X$;
 \item $X_0\cap X_1 = \{0\}$.
\end{itemize}
Here,  $X_0 + X_1 := \{x_0+x_1:\, x_0\in X_0,\, x_1\in X_1\}$.
In this situation we write $$ X = X_0\oplus X_1.$$
\end{definition}

\begin{definition}[Projections]
A {\em projection}\index{projection} is an operator $P\in \calL(X)$ satisfying $P^2=P$.
\end{definition}

Notice that the boundedness of $P$ is taken to be part of the definition. If $P$ is a projection, then so is $I-P$ and the
range $\ran(P)$ of $P$ equals the null space $\ker (I-P)$. This implies that $\ran(P)$ is closed and we have a direct sum decomposition
$$ X = \ker(P)\oplus \ran(P).$$ Thus we have shown the following simple result:

\begin{proposition}\label{prop:proj-compl}
If a closed subspace $X_0$ of a normed space $X$ is the range of a projection on $X$, then $X_0$ is complemented.
\end{proposition}

Conversely, if $X = X_0\oplus X_1$ is a direct sum decomposition of a Banach space, then the natural projections associated with it are bounded; this will be proved in the next chapter (see Proposition \ref{prop:compl-proj}).

We have seen in Corollary \ref{cor:findim-closed} that finite-dimensional subspaces are always closed. As an application of the Hahn--Banach theorem we next prove the stronger assertion that they are always complemented.
For later use we also include an analogue for subspaces of finite codimension, which does not require the use of the Hahn--Banach theorem. A subspace $X_0$ of a Banach space $X$ is said to have {\em finite codimension} if there exists a
finite-dimensional subspace $Y$ of $X$ such that $X_0\cap Y = \{0\}$ and $X_0+ Y = X$. In this situation we define the {\em codimension}\index{codimension} of $X_0$ to be the dimension of $Y$ and denote this number by $\codim(X_0)$. The following argument shows that this number is well defined. If $Y_0$ and $Y_1$ are subspaces with these properties, then for every $y_0\in Y_0$ there are unique $x_0\in X_0$ and $y_1\in Y_1$ such that $y_0 = x_0+y_1$. The mapping $y_0\mapsto y_1$ from $Y_0$ to $Y_1$ is easily seen to be linear. By the same procedure we obtain a well-defined mapping from $Y_1$ to $Y_0$, and it is clear that these mappings are each other's inverses. Hence they are isomorphisms of finite-dimensional vector spaces and therefore $\dim Y_0 = \dim Y_1$.

Subspaces of finite dimension are closed, but
subspaces of finite codimension need not be closed,
at least when we accept the Axiom of Choice: under this assumption,
Problem \ref{prob:nonbdd-operator} constructs a dense subspace of $\ell^2$ with codimension one, and such subspaces cannot be closed.

If $X_0$ is a {\em closed} subspace of finite co\-dimension,
it is easily checked that
$$ \codim (X_0) = \dim (X/X_0).$$

\begin{proposition}\label{prop:findim-compl} Let $Y$ be a subspace
of a normed space $X$. Then the following assertions hold:
\begin{enumerate}[label={\rm(\arabic*)}, leftmargin=*]
 \item\label{it:findim-compl1} if $\dim(Y)<\infty$, then $Y$ is closed and complemented;
 \item\label{it:findim-compl2} if $\codim(Y) <\infty$ and $Y$ is closed, then $Y$ is complemented.
 \end{enumerate}
\end{proposition}

\begin{proof}
\ref{it:findim-compl1}: \  Let $Y$ be a finite-dimensional subspace of $X$. By Corollary \ref{cor:findim-closed}, $Y$ is closed. To prove that $Y$ is complemented we show that $Y$ is the range of a projection on $X$.

Let $(y_n)_{n=1}^N$ be a basis for $Y$.
 Then every $y\in Y$ admits a unique representation $y = \sum_{n=1}^N c_n(y) y_n$ with coefficients $c_n(y)\in\K$.
 The mappings $c_n: y\mapsto c_n(y)$ are linear, and since linear mappings on finite-dimensional normed spaces are
 bounded, we have $c_n\in Y\s$\!. By the Hahn--Banach theorem we may extend each $c_n$ to a functional $x_n^*\in X\s$\!.
Consider the (bounded) linear operator $P$ on $X$ defined by
$$ Px := \sum_{n=1}^N \lb x,x_n^*\rb y_n.$$
It is clear that $P$ maps $X$ into $Y$ and from $\lb y_m,x_n^*\rb = \delta_{mn}$ we see that
$$Py_m = \sum_{n=1}^N\lb y_m,x_n^*\rb y_n = y_m.$$
This shows that $P$ maps $X$ {\em onto} $Y$. The preceding identity, applied to the element $Px\in Y$,
also shows that $P^2 x = P(Px) = Px$, so $P$ is a projection.

\smallskip\ref{it:findim-compl2}: \ Since finite-dimensional subspaces are closed,
this is immediate from the definitions.
\end{proof}

We next identify the duals of closed subspaces, quotients, and direct sums. For this purpose we need the first part of the following definition. The second part is included for symmetry of presentation and will be needed later.

\begin{definition}[Annihilators and pre-annihilators]\label{def:annih} Let $X$ be a Banach space.
\begin{enumerate}[leftmargin=*, label=(\roman*)]
 \item The {\em annihilator}\index{annihilator} of a subset $A\subseteq  X$ is the set\index{$Aa$@$A^\perp$}
 $$A^\perp:= \{x^*\in X^*:\, \lb x,x^*\rb = 0, \quad x\in A\}.$$
 \item The {\em pre-annihilator}\index{pre-annihilator} of a subset $B\subseteq  X^*$ is the set\index{$Aa$@${}^\perp A$}
 $${}^\perp B:= \{x\in X:\, \lb x,x^*\rb = 0, \quad x^*\in B\}.$$
\end{enumerate}
\end{definition}

\begin{proposition}\label{prop:dualOfQuotient}\label{prop:dualOfSubspace}
Let $X$ be a Banach space and let $Y$ be a closed subspace of $X$. Then $Y^\perp$
is a closed subspace of $X^*$ and we have the following assertions:
\begin{enumerate}[label={\rm(\arabic*)}, leftmargin=*]
 \item\label{it:dualOfQuotient1} the mapping $i: X\s/Y^\perp \to Y\s$ defined by $i(x\s+Y^\perp) := x\s|_{Y}$ is well defined and
 induces an isometric isomorphism
 $$ Y\s \simeq X\s/Y^\perp.$$

\item\label{it:dualOfQuotient2} the mapping $j: Y^\perp \to (X/Y)\s$ defined by $j x\s (x+Y) := \lb x, x\s\rb$  is well defined and
 induces an isometric isomorphism
  $$(X/Y)\s \simeq Y^\perp\!.$$
\end{enumerate}
\end{proposition}

\begin{proof} The easy proof that $Y^\perp$ is a closed subspace of $X^*$ is left as an exercise.

\smallskip
\smallskip
\ref{it:dualOfQuotient1}:\
Two functionals in the same equivalence class modulo $Y^\perp$ have the
same restriction to $Y$, so $i$ is well defined and linear.
If $i(x^*+Y^\perp)=0$, then $x^*$ vanishes on $Y$, and hence
$x^*\in Y^\perp$. Thus $i$ is injective.
Surjectivity follows from the Hahn--Banach theorem.

To prove that $i$ is an isometry, fix $x^*\in X^*$.
For every $y^*\in Y^\perp$ we have
$$
\n x^*|_Y\n
=\n(x^*+y^*)|_Y\n
\le\n x^*+y^*\n.
$$
Taking the infimum over $y^*\in Y^\perp$ gives
$$
\n i(x^*+Y^\perp)\n
\le\n x^*+Y^\perp\n_{X^*/Y^\perp}.
$$
Conversely, by the Hahn--Banach theorem, the restriction
$x^*|_Y$ has an extension $z^*\in X^*$ satisfying
$\n z^*\n=\n x^*|_Y\n$.
Since $z^*-x^*\in Y^\perp$, we obtain
$$
\n x^*+Y^\perp\n_{X^*/Y^\perp}
\le\n z^*\n
=\n x^*|_Y\n
=\n i(x^*+Y^\perp)\n.
$$
Together these inequalities prove the isometry.

\smallskip
\ref{it:dualOfQuotient2}: \ Let $x\s\in Y^\perp$. If $x+Y=x'+Y$, then $x-x'\in Y$, and therefore $\lb x-x',x\s\rb=0$. This proves that $jx\s$ is well defined.

For all $x\in X$ and $y\in Y$ we have
$$
|\lb x+Y,jx\s\rb|=|\lb x,x\s\rb|=|\lb x-y,x\s\rb|\leq \n x\s\n\,\n x-y\n .
$$
Taking the infimum over all $y\in Y$, we obtain $|\lb x+Y,jx\s\rb|\leq \n x\s\n\,\n x+Y\n_{X/Y}$. Hence $\n jx\s\n_{(X/Y)\s}\leq \n x\s\n$. Conversely, since $\n x+Y\n_{X/Y}\leq \n x\n$, we have
$$
|\lb x,x\s\rb|=|\lb x+Y,jx\s\rb|\leq \n jx\s\n_{(X/Y)\s}\,\n x+Y\n_{X/Y}\leq \n jx\s\n_{(X/Y)\s}\,\n x\n .
$$
Taking the supremum over all $x\in X$ with $\n x\n\leq 1$, we find that $\n x\s\n\leq \n jx\s\n_{(X/Y)\s}$. Thus $j$ is an isometry.

It remains to prove that $j$ is onto. Let $\phi\in (X/Y)\s$ be given, and define $x_\phi\s\in X\s$ by $\lb x,x_\phi\s\rb:=\lb x+Y,\phi\rb$. This is a bounded functional on $X$, being the composition of $\phi$ with the quotient mapping $x\mapsto x+Y$. Moreover, $x_\phi\s\in Y^\perp$, since $y+Y=0$ for all $y\in Y$. Finally, for all $x\in X$,
$$
\lb x+Y,jx_\phi\s\rb=\lb x,x_\phi\s\rb=\lb x+Y,\phi\rb .
$$
Hence $jx_\phi\s=\phi$, and $j$ is surjective.
\end{proof}

The duality of direct sum decompositions is discussed in Proposition \ref{prop:dual-iso}.

The Hahn--Banach theorem, through Corollary \ref{cor:HB-norming}, offers a technique to reduce certain vector-valued questions to their scalar-valued counterparts.
By way of example we demonstrate this technique by reproving some calculus rules for the vector-valued Riemann integral of
Proposition \ref{prop:der-f-zero}. A second example is given in Problem \ref{prob:weak-Cauchy} where the Cauchy integral formula is extended
to vector-valued holomorphic functions.

\medskip\noindent
{\it Second proof of Proposition \ref{prop:der-f-zero}, parts \ref{it:der-f-zero1} and \ref{it:der-f-zero2}}. \

\ref{it:der-f-zero1}: \
If $f(t_0)\not=f(t_1)$ for some $t_0,t_1\in (0,1)$, Corollary \ref{cor:HB-norming} provides a functional $x^*\in X^*$ such that
$\lb f(t_0),x^*\rb\not=\lb f(t_1),x^*\rb$.
The scalar-valued function
$\lb f,x^*\rb(t):=\lb f(t),x^*\rb$ is differentiable on $(0,1)$ and satisfies
$\lb f,x^*\rb'=\lb f',x^*\rb=0$.
It is therefore constant on $(0,1)$ by the scalar-valued version of the proposition, contradicting the choice of $x^*$.

\smallskip
\ref{it:der-f-zero2}: \ By the scalar-valued version of the proposition,
$$ \Bigl\lb f(1)-f(0)-\int_0^1 f'(t)\ud t, x\s\Bigr\rb = \lb f(1),x\s\rb -\lb f(0),x\s\rb -\int_0^1 \lb f'(t),x\s\rb \ud t = 0$$ for all $x\s\in X\s$\!. Corollary \ref{cor:HB-norming} implies that
$  f(1)-f(0)-\int_0^1 f'(t)\ud t = 0$. \hfill$\square$

\medskip
Using duality we can give the following version of the Pettis measurability theorem (Theorem \ref{thm:Pettis}):

\begin{theorem}[Pettis measurability theorem, second version]\label{thm:Pettis-secondversion}\index{theorem!Pettis measurability}
A function $f:\Om\to X$ is strongly measurable if and only if $f$ takes its values in a separable closed subspace of $X$ and is weakly measurable,\index{weakly!measurable} that is, $\lb f,x\s\rb:\Om\to\K$ is measurable for all $x\s\in X\s$\!.
\end{theorem}
\begin{proof}
The `only if' part follows from the first version of the Pettis measurability theorem and the trivial fact that strong measurability implies weak measurability.

For the `if' part, let $X_0$ be a closed separable subspace of $X$ containing the range of $f$. If $X_0=\{0\}$, there is nothing to prove. Otherwise choose a dense sequence $(x_k)_{k\ge 1}$ of nonzero elements of $X_0$.
By the Hahn--Banach theorem, for every $k\ge 1$ there is a unit vector $x_k^*\in X^*$ such that
$|\lb x_k,x_k^*\rb|=\n x_k\n$. Then for all $k\ge 1$ we have
$\sup_{n\ge 1}|\lb x_k,x_n\s\rb| = \n x_k\n$, and by a simple approximation argument this implies that
$\sup_{n\ge 1}|\lb x,x_n\s\rb| = \n x\n$ for all $x\in X_0$. Then, for all $x_0\in X_0$,
 $$\om\mapsto \n f(\om)-x_0\n = \sup_{n\ge 1} |\lb f(\om)-x_0,x_n\s\rb|$$
 is a measurable function. Now the result follows from Theorem \ref{thm:Pettis}.
\end{proof}

Theorem \ref{thm:Pettis-secondversion} is accompanied by the following uniqueness result.

\begin{proposition}  Let $(\Om,\calF\!,\mu)$ be a measure space.
If $f:\Omega\to X$ is a strongly measurable function and for all $x\s\in X\s$ we have $\lb f,x\s\rb = 0$ $\mu$-almost everywhere, then $f = 0$ $\mu$-almost everywhere.
\end{proposition}

Similarly, one proves that if $\lb f,x\s\rb = 0$ pointwise for all $x\s\in X\s$\!, then
 $f = 0$ pointwise.
\begin{proof}
Let $(x_n^*)_{n\ge 1}$ be a sequence in $X\s$ separating
the points of a closed subspace $X_0$ in which $f$ takes its values; such a sequence exists by the argument in the proof of Theorem \ref{thm:Pettis-secondversion}.
Since $\lb f,x_n\s\rb = 0$
outside a $\mu$-null set $N_n$, we conclude that $f=0$ on the complement of the $\mu$-null set
$\bigcup_{n\ge 1} N_n$.
\end{proof}

Corollary \ref{cor:HB-norming} has the interesting consequence that every Banach space can be isometrically identified with a
closed subspace of the bi-dual $X^{**} := (X^*)^*$ in a natural way.
More specifically, given an element $x\in X$ we define a mapping $Jx: X\s \to \K$ by
$$ Jx(x\s) := \lb x,x\s\rb.$$
It is clear that this mapping is bounded and therefore it defines an element of the bi-dual $X^{**}$. By the corollary, its norm  is given by
$$ \n Jx\n = \sup_{\n x\s\n\le 1} |\lb x\s\!, Jx\rb| = \sup_{\n x\s\n\le 1} |\lb x, x\s\rb| = \n x\n,$$
and therefore the mapping $J: x\mapsto Jx$ is isometric. It is also linear, and therefore we have proved:

\begin{proposition}[Isometric embedding into the bi-dual]\label{prop:bidual} The operator $J$ is an
isometric embedding of $X$ into $X^{**}$.
\end{proposition}

The image of $X$ under $J$ is closed in $X^{**}$; this is
immediate from the fact that $J$ is isometric.
It may happen that $J(X)$ is a proper subspace of $X^{**}$. For instance, the bi-dual of $c_0$ is $\ell^\infty$\!.
Examples of Banach spaces for which we have $J(X) = X^{**}$ include all Hilbert spaces and the spaces $\ell^p$ and $L^p(\Om)$ for $1<p<\infty$. Spaces with this property are called {\em reflexive} and enjoy some pleasant properties, some of which will be discussed in Section \ref{subsec:reflexivity}.

We conclude this section by filling in a detail that was left open in our treatment of the duality of $\ell^p$\!, namely, that the duality $(\ell^p)^*\simeq \ell^q$ with $\frac1p+\frac1q=1$, which has been shown to hold for $1\le p<\infty$ in Section \ref{subsec:dual-seqspaces}, does not hold for $p=\infty$.

Consider the closed subspace $Y$ of $\ell^\infty$
consisting of all convergent sequences, and define $y\s\in Y\s$ by
$$ \lb y,y\s\rb := \limn y_n$$
for $y = (y_n)_{n\ge 1} \in Y$.
Let $x^*\in(\ell^\infty)^*$ be any Hahn--Banach extension
of $y^*$.
Suppose that some $z\in\ell^1$ represents $x^*$, that is,
$\lb x,x^*\rb=\lb z,x\rb$ for every $x\in\ell^\infty$.
Testing this identity on the unit vectors $e_n$, which are
convergent sequences with limit zero, gives
$$
z_n=\lb e_n,z\rb=\lb e_n,x^*\rb=\lb e_n,y^*\rb=0,
\quad n\ge1.
$$
Thus $z=0$. Testing on the constant sequence $\one$, which
belongs to $Y$, gives the contradiction
$$
0 = \lb z,\one\rb = \lb\one,x^*\rb=\lb\one,y^*\rb=1.
$$
This shows that no element of $\ell^1$ represents $x^*$.

\section{Adjoint Operators}\label{subsec:BSadjoint}

The Hahn--Banach theorem will now be used to show that when $X$ and $Y$ are Banach spaces and $T\in \calL(X,Y)$ is a bounded operator, there exists a unique bounded operator
$T\s\in \calL(Y\s\!,X\s)$ of norm $\n T\s\n=\n T\n$ such that
$$ \lb Tx, y\s\rb = \lb x, T\s y\s\rb, \quad x\in X, \ y\s\in Y\s\!.$$
When $X$ and $Y$ are Hilbert spaces, the Riesz representation theorem will be used to prove the existence of a unique bounded operator $T^\star \in \calL(Y,X)$ of norm $\n T^\star\n=\n T\n$ such that
$$\iprod{Tx}{y} = \iprod{x}{T^\star y}, \quad  x\in X, \ y\in Y.$$

\subsection{The Banach Space Adjoint}

Let $X$ and $Y$ be Banach spaces.

\begin{proposition}
For every bounded operator $T\in \calL(X,Y)$ there exists a unique bounded operator
$T\s\in \calL(Y\s\!,X\s)$ such that
\begin{align*}
\lb Tx, y\s\rb = \lb x, T\s y\s\rb, \quad x\in X, \ y\s\in Y\s.
\end{align*}
Furthermore,
$$\n T\s\n = \n T\n.$$
\end{proposition}
\begin{proof} Given a functional $y\s\in Y\s$, we may define a linear mapping
$ T\s y\s: X \to \K$ by
$$ (T\s y\s)x := \lb Tx, y\s\rb.$$
This mapping is bounded, of norm $\n T\s y\s\n\le \n T\n \n y\s\n$, since
$$|(T\s y\s)x| \le \n Tx\n \n y\s\n \le \n T\n \n x\n \n y\s\n.$$
Accordingly $T\s y\s$ defines an element of $X\s$\!. The resulting mapping $T\s:Y\s\to X\s$ which maps
$y\s \in Y\s$ to the element $T\s y\s\in X\s$ is linear. By the above estimate $T\s$ is bounded, of norm $\n T\s\n\le \n T\n$.
It is clear from the definitions that $ \lb Tx, y\s\rb = \lb x, T\s y\s\rb$ for all $x\in X$ and $y\s\in Y\s$\!.

Turning to uniqueness, if $S:Y\s\to X\s$ is an operator satisfying
$ \lb Tx, y\s\rb = \lb x, S y\s\rb$ for all $x\in X$ and $y\s\in Y\s$\!, then
$ \lb x, T\s y\s\rb =  \lb x, S y\s\rb$ for all $x\in X$ and $y\s\in Y\s$\!, so $T\s y\s = S y\s$ for all $y\s\in Y\s$\!,
and so $T\s = S$.

Finally,
\begin{align*}
 \n T\s\n= \sup_{\n y\s\n\le 1} \n T\s y\s\n
& = \sup_{\n y\s\n\le 1} \sup_{\n x\n\le 1} |\lb x,T\s y\s\rb|
\\ & = \sup_{\n x\n\le 1} \sup_{\n y\s\n\le 1}  |\lb Tx, y\s\rb|
 = \sup_{\n x\n\le 1}  \n Tx\n = \n T\n,
\end{align*}
using Corollary \ref{cor:HB-norming} in the penultimate step.
\end{proof}

It is clear that $I_X\s = I_{X\s} $, where $I_X$ and $I_{X\s}$ are the identity operators on $X$ and $X\s$\!, respectively. For all
$T_1, T_2\in \calL(X,Y)$ and $c_1,c_2\in\K$ we have
$$(c_1 T_1 + c_2 T_2)\s = c_1 T_1\s + c_2 T_2\s$$
and for all $T\in \calL(X,Y)$ and $S \in \calL(Y,Z)$ we have
$$ (S\circ T)\s = T\s \circ S\s\!.$$

\begin{definition}[Adjoint operator] The bounded operator $T^*$
 is called the {\em adjoint} of $T$.\index{adjoint!of a bounded operator}\index{$T\s$}
\end{definition}

\begin{example}\label{ex:adjoint-fd} Let $X = \K^n$, $Y=\K^m$,
and let $A\in \calL(\K^n\!,\K^m)$ be represented, with respect to the standard bases, by the $m\times n$ matrix with coefficients $a_{ij}\in\K$ (see Example \ref{ex:matrices-bdd}).
Using the identification $\xi\leftrightarrow\phi_\xi$ of Section
\ref{subsec:dual-fd-spaces}, we represent its adjoint $A^*$, with respect to the same
bases, by the $n\times m$ matrix with coefficients $a_{ji}$.
Stated differently, the matrix associated with $A\s$ is the transpose of the matrix associated with $A$.
\end{example}

\begin{example}\label{ex:adjoint-ko}
 The adjoint of the kernel operator $T_k$ on $L^2(0,1)$ given by
 $$ T_k f(t):= \int_0^1 k(t,s)\,f(s)\ud s,$$
 where we assume that $k\in L^2((0,1)\times (0,1))$ (see  Example \ref{ex:kernel}),
 is the kernel operator
$T_{k^*}$ on $L^2(0,1)$ given by
 $$ T_{k^*}g(t) = \int_0^1 {k^*(t,s)}\,g(s)\ud s \ \hbox{ with } \ k^*(t,s) = k(s,t).$$
\end{example}

\begin{example}\label{ex:adjoint-shift} Let $1\le p<\infty$ and $\frac1p+ \frac1q = 1$. The adjoint of the left (right) shift on $\ell^p$ is the right (left) shift on $\ell^q$.
 The adjoint of the left (right) translation on $L^p(\R)$ is the right (left) translation on $L^q(\R)$.
\end{example}

As a first application we prove some simple duality results. The first is an immediate consequence of Corollary \ref{cor:proper-subspace}:

\begin{proposition}\label{prop:dual-dense-T} Let $X$ and $Y$ be Banach spaces. For all $T\in \calL(X,Y)$ we have $$ \ov{\Ran(T)} = {}^\perp(\Ker(T^*)) .$$
In particular, $T$ has dense range if and only if $T^*$ is injective.
\end{proposition}
\begin{proof}
  $\subseteq$: \
 If $y =Tx \in \Ran(T)$, then for all $y^*\in \Ker(T^*)$ we have
  $ \lb y, y\s\rb = \lb x,T\s y\s\rb = 0,$
 and therefore $y\in {}^\perp(\Ker(T^*))$. This proves the inclusion $\Ran(T)\subseteq {}^\perp(\Ker(T^*))$.
The inclusion $\ov{\Ran(T)} \subseteq {}^\perp(\Ker(T^*))$ now follows from the fact that ${}^\perp(\Ker(T^*))$ is closed.

  \smallskip
 $\supseteq$: \
  If $y_0\not\in \ov{\Ran(T)}=: Y_0$, by
  Corollary \ref{cor:proper-subspace}
  there exists a $y\s\in Y\s$ such that $\lb y_0,y^*\rb\not=0$ and $y\s|_{Y_0} \equiv 0$. For all $x\in X$,
  $ \lb x,T\s y\s\rb =  \lb Tx, y\s\rb = 0$ and therefore $y\s \in \ker(T^*)$. Since $\lb y_0,y^*\rb\not=0$
  we have $y_0\not\in  {}^\perp(\Ker(T^*))$.
 \end{proof}

The second application considers inclusions and restrictions.

\begin{proposition}\label{prop:dual-iso} Let $X$ be a Banach space.
Then:
\begin{enumerate}[label={\rm(\arabic*)}, leftmargin=*]
 \item \label{it:dial-iso1}
if $i:X\to Y$ is an (isometric) isomorphism of $X$ onto another Banach space $Y$, then the adjoint operator $i^*:Y^*\to X^*$ is an (isometric) isomorphism of their duals.
 \item \label{it:dial-iso2}
if $Y$ is a closed subspace of $X$ and
$i:Y\to X$  is the inclusion mapping, then the adjoint operator $i\s:X\s\to Y\s$ is the restriction mapping given by $i\s x\s  = x\s|_Y$;
 \item \label{it:dial-iso3}
if $X$ admits a direct sum decomposition $$X = X_0\oplus X_1$$ with projections $\pi_0\in \calL(X)$ and $\pi_1\in\calL(X)$, then
$$X\s = X_0\s \oplus X_1\s$$ isomorphically. More precisely, we have the direct sum decomposition $X\s = \pi_0\s X\s \oplus \pi_1\s X\s$,
and for $k=0,1$, with $i_k:X_k\to X$ denoting the inclusion mapping, the mapping $\pi_k\s x\s \mapsto i_k\s x\s$ defines an isomorphism from $\pi_k\s X\s$ onto $X_k\s$
 with $$
\|i_k^* x^*\| \ \le\ \|\pi_k^* x^* \| \le \|\pi_k\| \| i_k^* x^* \|.
$$
In particular, it is an isometric isomorphism when $\|\pi_k\|=1$ for $k=0,1$.
\end{enumerate}
\end{proposition}

It will be shown in the next chapter that the boundedness of $\pi_0$ and $\pi_1$ in part \ref{it:dial-iso3} is automatic.

\begin{proof}
We leave the proofs of \ref{it:dial-iso1} and \ref{it:dial-iso2}
to the reader and prove \ref{it:dial-iso3}.

Taking adjoints of the identities satisfied by $\pi_0$ and
$\pi_1$ shows that $\pi_0^*$ and $\pi_1^*$ are bounded
projections satisfying
$\pi_0^*+\pi_1^*=I$ and $\pi_0^*\pi_1^*=\pi_1^*\pi_0^*=0$.
Their ranges are closed, every $x^*\in X^*$ is the sum
$\pi_0^*x^*+\pi_1^*x^*$, and the intersection of these ranges
is $\{0\}$. This gives the direct sum decomposition
$$
X^*=\pi_0^*X^*\oplus\pi_1^*X^*.
$$

For $k=0,1$ and $\psi\in X_k^*$ define $J_k\psi\in X^*$ by
$$
\lb x,J_k\psi\rb:=\lb\pi_kx,\psi\rb,\quad x\in X.
$$
Then $J_k$ is bounded with $\n J_k\n\le\n\pi_k\n$.
Since $\pi_k$ is the identity on $X_k$, we have
$(J_k\psi)|_{X_k}=\psi$.
Moreover, $\pi_k^2=\pi_k$ implies
$\pi_k^*(J_k\psi)=J_k\psi$, so $J_k\psi\in\pi_k^*X^*$.

Conversely, if $\xi=\pi_k^*x^*$ with $x^*\in X^*$,
then $\pi_k^2=\pi_k$ gives $\xi\circ\pi_k=\xi$.
It follows that $J_k(\xi|_{X_k})=\xi$.
Thus restriction to $X_k$ and $J_k$ are inverse mappings
between $\pi_k^*X^*$ and $X_k^*$, and for all $\psi\in X_k^*$ we have
$$
\n\psi\n
=\n(J_k\psi)|_{X_k}\n
\le\n J_k\psi\n
\le\n\pi_k\n\n\psi\n.
$$
Since $(\pi_k^*x^*)|_{X_k}=i_k^*x^*$, these are
the asserted identification and norm bounds.
In particular, the identification is isometric when
$\n\pi_k\n=1$.
\end{proof}

We conclude with a simple observation about the bi-adjoint operator $T^{**}:= (T^*)^*$\!. Identifying $X$ with a closed subspace of $X^{**}$ by means of the natural isometric embedding $J:X\to X^{**}$ (see  Proposition \ref{prop:bidual}), the restriction of
$T^{**}$ to $X$ equals $T$. Indeed, the claim follows from
$$ \lb y\s\!, T^{**}Jx\rb = \lb T\s y\s\!, Jx\rb = \lb x,T\s y\s\rb = \lb Tx,y\s\rb = \lb y\s\!, JTx\rb,$$
so that $T^{**}Jx = JTx$ as claimed.

\subsection{The Hilbert Space Adjoint}\label{subsec:HSadjoint}

Let $H$ and $K$ be Hilbert spaces.
If $T\in \calL(H, K)$ is a bounded operator,
its adjoint $T\s\in\calL(K\s\!,H\s)$ is a bounded operator acting in the reverse direction between their duals.
By the Riesz representation theorem, the duals $H\s$ and $K^*$ can be canonically identified with $H$ and $K$. Under these identifications, the adjoint of an operator $T\in\calL(H,K)$ can be reinterpreted as an operator acting from $K$ to $H$.
Although the identifications are conjugate-linear, as an operator from  $K$ to $H$ the adjoint of $T$ is nevertheless linear. This is the content of the next proposition, which, incidentally, admits a straightforward direct proof which does not call upon the Hahn--Banach theorem.

\begin{proposition}\label{prop:HS-adjoint} For every bounded operator $T\in \calL(H,K)$ there exists a unique bounded operator $T^\star\in \calL(K,H)$ such that
$$ \iprod{Tx}{y} = \iprod{x}{T^\star y}, \quad  x\in H, \ y\in K.$$
Furthermore,
\begin{align*}\n T\n = \n T^\star\n = \n T^\star T\n^{1/2}\!.
\end{align*}
\end{proposition}
\begin{proof}
 Let $y\in K$ be fixed and define a mapping $\phi = \phi_{y,T}: H\to \K$ by
 $$ \phi(x):= \iprod{Tx}{y}.$$
From $|\phi(x)|\le \n Tx\n \n y\n \le \n T\n \n x\n \n y\n$ we see that $\phi$ is bounded with norm at most
$\n T\n \n y\n$. Hence by the Riesz representation theorem there is a unique element $T^\star y\in H$
with norm $\n T^\star y\n = \n \phi\n$ such that
$$ \phi(x) = \iprod{x}{T^\star y}.$$ Combining the two identities we obtain $\iprod{Tx}{y} = \iprod{x}{T^\star y}$.

We must show that the mapping $T^\star:y \mapsto T^\star y$ is linear and bounded.
Additivity is easy and homogeneity with respect to scalar multiplication follows from
$$ \iprod{x}{T^\star (cy)} = \iprod{Tx}{cy} = \ov c \iprod{Tx}{y} = \ov c\iprod{x}{T^\star y} = \iprod{x}{cT^\star y},$$
which implies $T^\star (cy) = cT^\star y$.

Next we show that $T^\star$ is bounded. This follows from what we already know. Indeed, we have
$$ \n T^\star y\n = \n \phi\n \le \n T\n \n y\n,$$
so $T^\star$ is bounded with norm $\n T^\star \n \le \n T\n$.
Write $T^{\star\star} := (T^\star)^\star$\!. From the identity
$$ \iprod{T^{\star\star} x}{y} = \ov{\iprod{y}{T^{\star\star} x}} = \ov{\iprod{T^\star y}{x}} = \iprod{x}{T^\star y} = \iprod{Tx}{y}$$
it follows that $T^{\star\star} x=Tx$ for all $x\in H$ and therefore $T^{\star\star}=T$.
Hence, by what we just proved applied to $T^\star$\!, $\n T\n = \n T^{\star\star}\n \le \n T^\star\n$.
We conclude that equality holds: $\n T\n =\n T^\star\n$.

Next we prove the identity $\n T^\star T\n^{1/2}=\n T\n$. The inequality
$\n T^\star T\n\le \n T^\star\n\,\n T\n=\n T\n^2$
is immediate and gives the upper estimate. Conversely, for $\n x\n\le1$ we have
$$
        \n Tx\n^2
        =
        \iprod{Tx}{Tx}
        =
        \iprod{x}{T^\star Tx}
        \le
        \n T^\star T\n .
$$
Taking the supremum over all such $x$ gives the lower estimate
$\n T\n^2\le \n T^\star T\n$.

Finally, we prove uniqueness. If $U\in \calL(K,H)$ is bounded and
$\iprod{x}{Uy} = \iprod{Tx}{y}$ for all $x\in H$ and $y\in K$, then $Uy = T^\star y$ for all $y\in K$, so $U=T^\star$.
\end{proof}

It is clear that $I^\star = I$ and for all $T,U\in \calL(H)$ and $c\in \K$ we have
\begin{align*}
(T+U)^\star = T^\star + U^\star\!, \qquad (cT)^\star = \ov c T^\star\!,
\end{align*}
and, as we have seen in the proof of Proposition \ref{prop:HS-adjoint},
$$ T^{\star\star} = T.$$

\begin{definition}[Hilbert space adjoint]
The bounded operator $T^\star$\index{$T^\star$}
is called the {\em Hilbert space adjoint}\index{adjoint!of a bounded Hilbert space operator} of $T$.
\end{definition}

The relation between the Hilbert space adjoint $T^\star$ and the adjoint operator $T^*$ is given by
$$  T\s \psi_{y} = \psi_{T^\star y}$$
as elements of $H\s$\!, where $\psi_y\in K^*$ and  $\psi_{T^\star y}\in H^*$ are the functionals associated with $y$ and $T^\star y$, respectively.
Indeed, this follows from
$$  \lb x, T\s \psi_{y}\rb = \lb Tx, \psi_y\rb = \iprod{Tx}{y} = \iprod{x}{T^\star y} = \lb x, \psi_{T^\star y}\rb.$$

\begin{example}\label{ex:HSadjoints}
Here are some examples of Hilbert space adjoints. They should be compared with Examples \ref{ex:adjoint-fd}, \ref{ex:adjoint-ko}, and \ref{ex:adjoint-shift}, respectively.

\begin{enumerate}[leftmargin=*, label=(\roman*)]
 \item
 As in Example \ref{ex:adjoint-fd}, let $A\in \calL(\K^n\!,\K^m)$ be represented as an $m\times n$ matrix with coefficients $a_{ij}\in \K$.
Viewing $\K^n$ and $\K^m$ as finite-dimensional Hilbert spaces,
its Hilbert space adjoint $A^\star$ may be represented as the $n\times m$ matrix with coefficients $\ov{a_{ji}}.$
Stated differently, the matrix associated with $A^\star$ is the Hermitian transpose of the matrix associated with $A$.

 \item
The Hilbert space adjoint of the kernel operator $T_k$ on $L^2(0,1)$ of Example \ref{ex:adjoint-ko}
is the kernel operator $T_{k^\star}$ on $L^2(0,1)$ given by
 $$ T_{k^\star}g(t) = \int_0^1 {k^\star(t,s)}\,g(s)\ud s \ \hbox{ with } \ k^\star(t,s) = \ov{k(s,t)}.$$

 \item
 The adjoint of the left (right) shift on $\ell^2(\Z)$ is the right (left) shift.
 Similarly, the adjoint of the left (right) translation on $L^2(\R)$ is the right (left) translation.
\end{enumerate}
\end{example}

The following decomposition result refines Proposition \ref{prop:dual-dense-T}.

\begin{proposition}\label{prop:injective-denserange} If $T\in \calL(H,K)$ is a bounded operator,
then $H$ and $K$ admit orthogonal decompositions
$$
H = \Ker(T) \oplus \overline{\Ran(T^\star)}, \qquad
K = \Ker(T^\star) \oplus \overline{\Ran(T)}.
$$
In particular,
\begin{enumerate}[label={\rm(\arabic*)}, leftmargin=*]
 \item $T$ is injective if and only if $T^\star$ has dense range;
 \item $T$ has dense range if and only if $T^\star$ is injective.
\end{enumerate}
\end{proposition}
\begin{proof}
 If $x\perp \ov{\Ran(T^\star)}$, then
 $\iprod{Tx}{y} = \iprod{x}{T^\star y} = 0$ for all $y\in K$ and therefore $Tx = 0$, so $x\in \ker(T)$.
 Conversely, if $x\in \ker(T)$, then
 $ \iprod{x}{T^\star y} = \iprod{Tx}{y} =0$ for all $y\in K$ implies that $x\perp\ran(T^\star)$ and hence
 $x\perp \ov{\Ran(T^\star)}$.
 This proves the orthogonal decomposition for $H$.
 The decomposition for $K$ follows from it by applying it to $T^\star$ and using that $T^{\star\star} = T$.
\end{proof}

\section{The Hahn--Banach Separation Theorem}\label{sec:HB-separation}

In what follows, $X$ is a Banach space.
Corollary \ref{cor:proper-subspace} can be interpreted as a separation theorem, in that it guarantees the existence
of a functional separating a closed subspace from a given element not contained in it. The following result provides a far-reaching generalisation:

\begin{theorem}[Hahn--Banach separation theorem]\label{thm:HB-separation}\index{theorem!Hahn--Banach separation}
Let $C$ and $D$ be disjoint nonempty convex sets in $X$, with $C$ open.
Then there exist a nonzero $x\s\in X\s$ and a $c\in\R$ such that
$$
\Re\lb x,x\s\rb<c\le\Re\lb y,x\s\rb,
\quad x\in C,\ y\in D.
$$
In particular, the sets $\lb C,x\s\rb$ and $\lb D,x\s\rb$
are disjoint.
\end{theorem}
\begin{proof}
We prove the theorem in three steps.

\smallskip
{\em Step 1} -- First we prove the theorem for the real scalar field and $D = \{x_0\}$.
Replacing $C$ and $x_0$ by $C-y_0$ and $x_0-y_0$ for some fixed $y_0\in C$, we may assume without loss of generality that $0\in C$.

Define the {\em Minkowski functional}\index{Minkowski functional} of $C$
 as the mapping $\la_{\,C}: X\to [0,\infty)$ given by $$
 \la_{\,C}(x):= \inf\{t>0:\, t^{-1}x\in C\}.$$
 This functional is finite-valued: since $C$ is open and contains
 $0$, it contains a ball $B(0;r)$, and hence $t^{-1}x\in C$ whenever
 $t>\|x\|/r$.
 Using that $C$ is convex and open, and contains $0$, one checks that
 $\la_{\,C}(x) < 1$ if and only if $x\in C$.
 We claim that $\la_{\,C}$ enjoys the following two properties:
 \begin{enumerate}[leftmargin=*, label=(\roman*)]
  \item\label{it:Mink1} $\la_{\,C}(x+y) \le \la_{\,C}(x) + \la_{\,C}(y)$ for all $x,y\in X$;
  \item\label{it:Mink2} $\la_{\,C}(tx) = t\la_{\,C}(x)$ for all $t\ge 0$.
 \end{enumerate}
To prove \ref{it:Mink1}, fix $\eps>0$ and let $s,t>0$ be such that
$s^{-1}x \in C$ and $t^{-1}y\in C$, with $s< \la_{\,C}(x)+\eps$ and $t< \la_{\,C}(y)+\eps$.
Then
$$ (s+t)^{-1}(x+y) = \frac{s}{s+t} s^{-1}x +\frac{t}{s+t} t^{-1}y$$
is a convex combination of the elements $s^{-1}x, t^{-1}y\in C$ and therefore belongs to $C$.
It follows that
$$ \la_{\,C}(x+y) \le s+t \le \la_{\,C}(x) + \la_{\,C}(y)+2\eps.$$
Since $\eps>0$ was arbitrary, this establishes \ref{it:Mink1}. Assertion \ref{it:Mink2} is obvious.

We now apply Theorem \ref{thm:appHB1} to the linear span $W$ of $x_0$
and the linear mapping $\phi:W\to\R$ given by $\phi(tx_0):= t$ for $t\in\R$.
In view of $\la_{\,C}(x_0)\ge 1$, for all $t\ge 0$ it satisfies
$$\phi(tx_0) = t  \leq t\la_{\,C}(x_0) = \la_{\,C}(tx_0);$$
for $t<0$ one trivially has $\phi(tx_0)=t<0\le \la_{\,C}(tx_0)$.
Hence we may apply the theorem and obtain a linear mapping $x\s:X\to \R$ extending $\phi$
which satisfies $ \lb x,x\s\rb\leq \la_{\,C}(x)$ for all $x\in X$.
For all $x\in C$ it satisfies $\lb x,x\s\rb \leq \la_{\,C}(x) < 1$ and for all
$x\in -C$ it satisfies $-\lb x,x\s\rb = \lb -x,x\s\rb \le \la_{\,C} (-x) =\la_{-C}(x)< 1.$
It follows that $|\lb x,x\s\rb|<1$ for all $x\in C\cap -C$. Since $C\cap -C$ is open and contains $0$, this proves that $x\s\in X\s$\!. Moreover, $\lb x,x\s\rb \le \la_{\,C}(x)<1$ for all $x\in C$ and $\lb x_0,x\s\rb = \phi(x_0) = 1$.
We have thus proved that
$\lb x,x\s\rb<\lb x_0,x\s\rb$ for all $x\in C$. Therefore this functional has the required properties, with $c=1$.

\smallskip
{\em Step 2} -- In the case of complex scalars and $D=\{x_0\}$,
we regard $X$ as a real normed space. Step 1 gives a bounded
real-linear functional $x_\R\s:X\to\R$ such that
$x_\R\s(x)<x_\R\s(x_0)$ for all $x\in C$.
As in the proof of Theorem \ref{thm:appHB2}, the formula
$$
x\s(x):=x_\R\s(x)-ix_\R\s(ix)
$$
defines a bounded complex-linear functional $x^*\in X^*$. Since
$\Re\lb x,x\s\rb=x_\R\s(x)$, it satisfies
$\Re\lb x,x\s\rb<\Re\lb x_0,x\s\rb$ for all $x\in C$.

\smallskip
{\em Step 3} -- In the general case, the set $C-D$ is nonempty,
open, and convex, and $0\notin C-D$. The preceding two steps
provide a nonzero $x\s\in X\s$ such that
$$
\Re\lb z,x\s\rb<0,\quad z\in C-D.
$$
Consequently,
$$
\Re\lb x,x\s\rb<\Re\lb y,x\s\rb,
\quad x\in C,\ y\in D.
$$
Put $c:=\sup_{x\in C}\Re\lb x,x\s\rb$. This is a finite real
number, since $C$ and $D$ are nonempty and the preceding inequality
bounds the supremum from above. It also implies
$c\le\Re\lb y,x\s\rb$ for every $y\in D$.
Since $x\s\ne0$, there exists a $v\in X$ with
$\Re\lb v,x\s\rb>0$. For every $x\in C$, openness of $C$ gives
$x+tv\in C$ for all sufficiently small $t>0$. Hence
\[
\Re\lb x,x\s\rb
<\Re\lb x+tv,x\s\rb\le c.
\qedhere\]
\end{proof}

\begin{corollary}\label{cor:Mazur} Suppose $(x_n)_{n\ge 1}$ is a sequence in $X$ and suppose that there exists an $x\in X$ such that
 $$ \limn \lb x_n, x^*\rb = \lb x,x^*\rb, \quad x^*\in X^*\!.$$
Then there exists a sequence $(y_n)_{n\ge 1}$ in the convex hull of $(x_n)_{n\ge 1}$ such that
 $$ \limn y_n = x$$
with convergence in norm.
\end{corollary}

\begin{proof}
Let $D$ be the closure of the convex hull of $(x_n)_{n\ge1}$.
Suppose that $x\notin D$, and choose an open ball $C$ centred at
$x$ and disjoint from $D$. By Theorem \ref{thm:HB-separation},
there exist $x^*\in X^*$ and $c\in\R$ such that
$$
\Re\lb x,x^*\rb<c\le\Re\lb x_n,x^*\rb,\quad n\ge1.
$$
Letting $n\to\infty$, we obtain a contradiction with the assumed weak convergence. Thus $x\in D$. For each $n\ge1$, choose $y_n$ in the
convex hull of $(x_k)_{k\ge1}$ with $\n y_n-x\n<1/n$.
Then $y_n\to x$ in norm.
\end{proof}

As another application of the Hahn--Banach theorem we prove the following result.

\begin{theorem}[Helly]\label{thm:locrefl}\index{theorem!Helly}
For all $x^{**} \in X^{**}$,
$x_1\s\!,\dots,x_N\s\in X\s$\!, and $\eps>0$ there exists an $x\in X$ such that
$\n x\n < \n x^{**}\n+\eps$ and
$$ \lb x,x_n\s\rb = \lb x_n\s, x^{**}\rb, \quad n=1,\dots,N.$$
\end{theorem}
\begin{proof}
Let $E$ denote the finite-dimensional subspace of $X^*$ spanned by $x_1\s,\ldots,x_N\s$, and put $Y:={}^{\perp}E$. We first observe that $Y^\perp=E$. The inclusion $E\subseteq Y^\perp$ is immediate. Conversely, let $x^*\in Y^\perp$, and choose a basis $e_1^*,\ldots,e_d^*$ for $E$. The mapping $S:X\to\K^d$ defined by $Sx=(\lb x,e_j^*\rb)_{j=1}^d$ has range $\K^d$. Indeed, if its range were a proper subspace of $\K^d$, there would exist scalars $\la_1,\ldots,\la_d$, not all zero, such that $\sum_{j=1}^d\la_j\lb x,e_j^*\rb=0$ for all $x\in X$, contradicting the linear independence of $e_1^*,\ldots,e_d^*$. Since $\ker(S)=Y$ and $x^*\in Y^\perp$, we have $\lb x-x',x^*\rb=0$ whenever $Sx=Sx'$. Hence the formula $\ell(Sx):=\lb x,x^*\rb$ defines a linear functional $\ell$ on $\K^d$, and $x^*=\ell\circ S$. It follows that $x^*$ is a linear combination of $e_1^*,\ldots,e_d^*$, and therefore $x^*\in E$.

By Proposition \ref{prop:dualOfQuotient},
the dual of $X/Y$ is isometrically identified with $Y^\perp=E$. Hence, for $x\in X$,
$$
\n x+Y\n_{X/Y}=\sup_{\substack{x^*\in E\\ \n x^*\n\leq 1}}|\lb x,x^*\rb|.
$$
We claim that the mapping $J:X/Y\to E^*$ defined by
$$
        (J(x+Y))(x^*):=\lb x,x^*\rb,\quad x^*\in E,
$$
is an isometry.
To see that it is surjective, let $\phi\in E^*$.
Since the mapping $S$ constructed above is surjective onto $\K^d$,
there is an $x\in X$ such that
$$
\lb x,e_j^*\rb=\phi(e_j^*),\quad j=1,\dots,d.
$$
The functionals $J(x+Y)$ and $\phi$ therefore agree on a basis
of $E$, so $J(x+Y)=\phi$. We conclude that $J$ is an isometric isomorphism.

 Restricting a given $x^{**}$ to $E$, we obtain an element $x^{**}|_E$ of $E^*$ with $\n x^{**}|_E\n\leq \n x^{**}\n$. Since $J$ is surjective, there exists $z\in X/Y$ such that $Jz=x^{**}|_E$, and since $J$ is isometric we have $\n z\n_{X/Y}\leq \n x^{**}\n$. Choose $x\in X$ representing $z$ and satisfying $\n x\n<\n z\n_{X/Y}+\eps$. Then $\n x\n<\n x^{**}\n+\eps$, and we have
\[
 \lb x,x_n\s\rb=(Jz)(x_n\s)=\lb x_n\s,x^{**}\rb ,\quad n=1,\ldots,N .
\qedhere\]
\end{proof}

\section{The Krein--Milman Theorem}\label{sec:Krein--Milman}

Extreme points play an important role in many applications of functional analysis.

\begin{definition}[Extreme points]\label{def:extreme-point}
An {\em extreme point}\index{extreme point} of a convex subset $C$ of a vector space is an element $v\in C$ such that if
$v = (1-\la)v_0+ \la v_1$ with $v_0,v_1\in C$ and $0<\la<1$, then $v_0=v_1= v$.
\end{definition}

Stated differently,
extreme points are points of $C$ which cannot be realised in a nontrivial way as a convex combination of
other points of $C$.

\begin{example}\label{ex:extreme-points-K} Let $C$ denote the set of all probability measures on a given measure space.
Viewing $C$ as a closed convex subset of $M(\Om)$, the Banach space of $\K$-valued measures on $(\Om,\calF)$, we claim that a probability measure $\mu$ is an extreme point of $C$ if and only if $\mu$ is two-valued, that is, whenever $\Omega = A_0\cup A_1$ with disjoint $A_0,A_1\in\calF$\!, then $\min\{\mu(A_0),\mu(A_1)\}=0$.

To prove the claim, suppose first that $\mu\in C$ and $\mu = (1-\la)\mu_0+\la \mu_1$ with $\mu_0,\mu_1\in C$ and $0<\la<1$.
If $\mu_0\not=\mu_1$, there is a set $A\in \calF$ such that $\mu_0(A)\not=\mu_1(A)$. Interchanging $\mu_0$ and $\mu_1$ if necessary, we may assume that  $0\le \mu_0(A)< \mu_1(A)\le 1$.
Then $\mu_1(A)>0$ implies $\mu(A)>0$, and $\mu_0(\complement A)>0$ implies $\mu(\complement A)>0$,
and therefore $\mu$ is not two-valued. This proves that every two-valued measure is an extreme point of $C$.

Conversely, if $\mu\in C$ is not two-valued, then there exists a set $A\in \calF$ such that $\mu(A)>0$ and $\mu(\complement A)>0$. Consider the probability measures $\mu_0,\mu_1\in C$ given by
$$\mu_0(B) := (\mu(A))^{-1}\mu(B\cap A), \qquad \mu_1(B) := (\mu(\complement A))^{-1}\mu(B\cap \complement A),$$
for $B\in \calF$\!.
With $\la:= \mu(\complement A)$ we have $0<\la<1$ and  $\mu = (1-\la) \mu_0 + \la \mu_1$,
so $\mu$ is not extreme.

As a special case, if $K$ is a compact Hausdorff space, the extreme points
of the set of Radon probability measures on $K$ are the Dirac measures.
Indeed, by the preceding part, the extreme points are precisely the two-valued
Radon probability measures. If $\mu$ is such a measure and
$\mu(\{x\})=0$ for all $x\in K$, then by outer regularity, for every
$x\in K$ there exists an open neighbourhood $U_x$ of $x$ such that
$\mu(U_x)<1$, hence $\mu(U_x)=0$. By compactness, finitely many of these
sets cover $K$, contradicting $\mu(K)=1$. Thus $\mu(\{x\})=1$ for some
$x\in K$, and then $\mu=\delta_x$.
\end{example}

Closed convex sets need not have any extreme points:

\begin{example}
In $L^1(0,1)$, the closed convex set $$C = \{f\in L^1(0,1):\, f \ge 0, \ \n f\n_1 = 1\}$$
has no extreme points.
Indeed, let $f\in C$. The mapping $\phi(\delta):= \n\one_{(0,\delta)}f\n_1$ is continuous from $[0,1]$ to $[0,1]$,  and satisfies $\phi(0)=0$ and $\phi(1)=1$. Hence there exists $0<\delta<1$ such that $\phi(\delta)=\frac12$.
Then $f= \frac12g + \frac12h$, where $g,h\in C$ are given by $g = 2\one_{(0,\delta)}f$ and $h = 2\one_{(\delta,1)}f$,
so $f$ is not an extreme point of $C$.
\end{example}

As an application of the Hahn--Banach theorem we can prove the following result
about the existence of extreme points. Recall that the {\em (closed) convex hull}
of a subset $S$ of a Banach space is the smallest (closed) convex set containing $S$.

\begin{theorem}[Krein--Milman]\label{thm:Krein-Milman}\index{theorem!Krein--Milman}\index{compactness!and extreme points}  Every compact convex  subset of a Banach space is the closed convex hull of its extreme points.
\end{theorem}
\begin{proof} Let $K$ be a compact convex subset of the Banach space $X$. We may assume that $K$ is nonempty, for otherwise there is nothing to prove.

 We first prove the existence of extreme points of $K$, and then prove that $K$ is the closed convex hull of its extreme points.

 A {\em face}\index{face} of $K$ is a nonempty closed convex subset $F$ of $K$ whose elements can only be realised as convex combinations of elements in $F$, that is, whenever $x\in F$ satisfies
$x = (1-\la)x_0+ \la x_1$ with $x_0,x_1\in K$ and $0<\la<1$, then $x_0, x_1\in F$.
We make two useful observations:
\begin{enumerate}[leftmargin=*, label=(\roman*)]
 \item\label{it:KM1} $x\in K$ is an extreme point of $K$ if and only if $\{x\}$ is a face of $K$;
 \item\label{it:KM2} if $F$ is a face of $K$ and $F'$ is a face of $F$, then $F'$ is a face of $K$.
\end{enumerate}
Claim \ref{it:KM1} is evident. For \ref{it:KM2}, if $x\in F'$ is given as
$x = (1-\la)x_0+ \la x_1$ with $x_0,x_1\in K$ and $0<\la<1$, then $x_0,x_1\in F$ since $x\in F$ and
$F$ is a face of $K$. Then $x_0,x_1\in F'$ since $F'$ is a face of $F$.

\smallskip
{\em Step 1} -- Let $\mathscr{K}$ denote the collection of faces of $K$.
This collection is nonempty, for it contains $K$. We partially order
$\mathscr{K}$ by declaring $K_1\le K_2$ if $K_2\subseteq  K_1$. By the finite intersection property (see Appendix \ref{sec:topology}),
any totally ordered subset $\mathscr{L}$ of $\mathscr{K}$ has nonempty intersection,
and this intersection is a face, hence belongs to $\mathscr{K}$, and clearly it is an upper bound for $\mathscr{L}$.
Hence we can apply Zorn's lemma and obtain that
$\mathscr{K}$ has a maximal element, say $F$. We claim that $F$ is a singleton, say $F= \{x\}$.
By \ref{it:KM1}, this means that $x$ is an extreme point of $K$.

To prove the claim, assume the contrary and let $x_0, x_1\in F$ be two distinct points.
By the Hahn--Banach theorem there exists an $x\s\in X\s$ such that
$\lb x_0,x\s\rb \not=\lb x_1,x\s\rb.$
Multiplying $x\s$ by an appropriate scalar if necessary, we may even assume that
$\Re\lb x_0,x\s\rb \not=\Re\lb x_1,x\s\rb.$
Let $$F_0 := \bigl\{x\in  F: \, \Re\lb x,x\s\rb = \inf_{y\in F}\Re\lb y,x\s\rb\bigr\}.$$
Then $F_0$ is a proper closed subset of $F$, and $F_0$ is nonempty since the compactness of $F$ implies that
the infimum is a minimum. If an element $x\in F_0$ can be represented as
$x = (1-\la)x'+ \la x''$ with $x'\!,x''\in F$ and $0<\la<1$, then
\begin{align}\label{eq:equality-face-F0}(1-\la)\Re\lb x'\!,x\s\rb + \la \Re\lb x''\!,x\s\rb=  \inf_{y\in F}\Re\lb y,x\s\rb.
\end{align}
If $x'\not\in F_0$,  then  $\Re\lb x'\!,x\s\rb > \inf_{y\in F}\Re\lb y,x\s\rb$. Since also
$ \Re\lb x''\!,x\s\rb \ge \inf_{y\in F}\Re\lb y,x\s\rb$,
the equality \eqref{eq:equality-face-F0} cannot hold. The same contradiction is reached if $x''\not\in F_0$.
We conclude that $x'\!,x''\in F_0$ and $F_0$ is a face of $F$.

Now \ref{it:KM2} implies that $F_0$ is a face of $K$.
Since $F_0\subsetneq F$, this contradicts the maximality of $F$. This completes the proof of the claim that $F$ is a singleton. It follows that $K$ has an extreme point.

\smallskip
{\em Step 2} -- Let $L$ denote the closed convex hull of all extreme points of $K$. We wish to show that $L=K$. Arguing by contradiction, suppose that $L$ is a proper subset of $K$ and fix an element $x_0\in K\setminus L$. Choose an open ball $C$ centred at $x_0$ such that $C\cap L=\emptyset$.
By Theorem \ref{thm:HB-separation}, there exist $x^*\in X^*$ and
$c\in\R$ such that
$$
\Re\lb x_0,x^*\rb<c\le\inf_{y\in L}\Re\lb y,x^*\rb.
$$
As in Step 1 we see that the set
$$ F := \Big\{x\in K:\, \Re\lb x,x\s\rb = \inf_{y\in K} \Re\lb y,x\s\rb\Big\}$$
is a nonempty face of $K$. Moreover, Step 1 applied to $F$ shows that $F$ has an extreme point $x_1$.
By \ref{it:KM1} and \ref{it:KM2}, $x_1$ is also an extreme point of $K$. On the other hand, from
$$  \Re\lb x_1,x\s\rb = \inf_{y\in K} \Re\lb y,x\s\rb \le \Re\lb x_0,x\s\rb < \inf_{y\in L} \Re\lb y,x\s\rb$$
we infer that $x_1\not\in L$. Since $L$ contains all extreme points of $K$ we have arrived at a contradiction.
\end{proof}

\section{The Weak and Weak\texorpdfstring{$^*$}{*} Topologies}\label{sec:weak-star}

Some of the more advanced applications of the Hahn--Banach theorem can be conveniently formulated in terms of certain topologies generated by bounded functionals. The two most important ones are the weak topology of a Banach space and the weak$\s$ topology of its dual.

\begin{definition}[Weak topologies]
Let $V$ and $W$ be vector spaces and let  $\beta:V\times W \to \K$ be a
bilinear mapping.
The {\em weak topology of $V$ generated by $W$} is the smallest topology $\tau$ on $V$ with the property that the linear mapping $v\mapsto \beta(v,w)$ is continuous for all $w\in W$.
\end{definition}

This topology is obtained as the intersection of all topologies on $V$
for which all linear mappings $v\mapsto \beta(v,w)$, $w\in W$, are continuous. The family of topologies with this property is nonempty, for it always contains the power set topology of $V$.

Necessarily, the weak topology $\tau$ must contain every set of the form
$$ U_{v_0, w_0,\eps} := \{v\in V:\, | \beta(v- v_0, w_0)| < \eps\},$$
noting that this set is the inverse image under the continuous mapping $v\mapsto \beta(v,w_0)$ of the open ball $B(\beta(v_0,w_0);\eps)$ in $\K$.

We claim that $\tau$ coincides with the topology $\tau'$ generated by the sets $U_{v_0, w_0,\eps}$, where
$v_0$, $w_0$, and $\eps$ range over $V$, $W$, and $(0,\infty)$, respectively.
The observation just made implies that $\tau'\subseteq \tau$.
In the opposite direction, for every $w\in W$ the inverse image under $\beta(\cdot,w)$ of every open ball belongs to $\tau'$\!, so
every $\beta(\cdot,w)$ is continuous with respect to $\tau'$\!. Since
$\tau$ is the smallest topology with this property we have $\tau\subseteq  \tau'$\!. This establishes the claim.

It follows from the claim that a set $U\subseteq  V$ belongs to $\tau$ if and only if it can be written
as a union of finite intersections of sets of the form $U_{v_0, w_0,\eps}$. Indeed, the collection $\tau''$ of sets that can be written this way is a topology which contains every set $U_{v_0, w_0,\eps}$, and therefore we have $\tau'\subseteq  \tau''$\!. By the preceding observation, this means that $\tau\subseteq  \tau''$\!. In the converse direction, the fact that
topologies are closed under taking unions and finite intersections implies that every set in $\tau''$ belongs to $\tau$.

\begin{proposition}\label{prop:seqconv-weak} In the above setting,
a sequence $(v_n)_{n\ge 1}$ converges to $v$ with respect to $\tau$ if and only if
$$ \limn \beta(v_n-v,w) = 0, \quad w\in W.$$
\end{proposition}
\begin{proof}
The `only if' part follows from the fact that if $v_n\to v$ with respect to $\tau$,
then for all $\eps>0$ and $w\in W$ we have $v_n \in U_{v,w,\eps}$ for all large enough $n$.
For the `if' part we note that if $U\in \tau$ contains $v$, then the observation preceding the statement of the proposition
allows us to find $U_{v^{(1)}, w^{(1)},\eps^{(1)}}, \dots, U_{v^{(k)}, w^{(k)},\eps^{(k)}}$ such that
$$ v\in \bigcap_{j=1}^k U_{v^{(j)}, w^{(j)},\eps^{(j)}} \subseteq  U.$$
Since we assume that $\beta(v_n-v,w^{(j)}) \to 0$ for $j=1,\dots,k$,
for large enough $n$ we have
$v_n \in \bigcap_{j=1}^k U_{v^{(j)}, w^{(j)},\eps^{(j)}}$ and hence $v_n\in U$.
\end{proof}

The duality between a Banach space and its dual leads to two special cases of interest:

\begin{definition}[The weak and weak$^*$ topologies] \label{def:weaktop} Let $X$ be a Banach space.
 \begin{enumerate}[label={\rm(\roman*)}, leftmargin=*]
  \item\label{it:weaktop1}
The {\em weak topology of $X$} is the topology induced by $X\s$\!.\index{weak!topology}
 \item\label{it:weaktop2}
The {\em weak$^*$ topology of $X\s$} is the topology induced by $X$.\index{weak$^*$!topology}
 \end{enumerate}
\end{definition}

It is implicit that in \ref{it:weaktop1} we use the bilinear mapping from $X\times X\s$ to $\K$
given by $(x,x\s) \mapsto \lb x,x\s\rb$; in \ref{it:weaktop2} we use the bilinear mapping from $X\s\times X$ to $\K$
given by $(x\s\!,x) \mapsto\lb x,x\s\rb$. For these topologies, Proposition \ref{prop:seqconv-weak} takes the following form:

\begin{corollary}\label{cor:seqconv-weak} Let $X$ be a Banach space.  The following assertions hold:
\begin{enumerate}[label={\rm(\arabic*)}, leftmargin=*]
 \item a sequence $(x_n)_{n\ge 1}$ in $X$ converges to $x\in X$ with respect to the weak topology of $X$ if and only if
$ \limn \lb x_n-x,x\s\rb = 0$ for all $ x\s\in X\s;$
\item a sequence $(x_n\s)_{n\ge 1}$ in $X\s$ converges to $x\s\in X\s$ with respect to the weak$^*$ topology of $X\s$ if and only if
$ \limn \lb x,x_n\s-x\s\rb = 0$ for all $ x\in X.$
\end{enumerate}
\end{corollary}

Convergence with respect to the weak and weak$^*$ topologies will be referred to as {\em weak convergence}\index{weak!convergence} and {\em weak$^*$ convergence},\index{weak$^*$!convergence} respectively.

\begin{proposition}[Closed convex sets are weakly closed]\label{prop:convex-weaklyclosed}\index{weakly!closed}
Every closed convex set in a Banach space is closed in the weak topology.
\end{proposition}
\begin{proof}
The empty set is weakly closed. Let $D$ be a nonempty closed
convex set, and let $x_0\notin D$. Choose an open ball $C$ centred
at $x_0$ and disjoint from $D$. By Theorem \ref{thm:HB-separation},
there exist $x^*\in X^*$ and $c\in\R$ such that
$$
\Re\lb x_0,x^*\rb<c\le\Re\lb y,x^*\rb,\quad y\in D.
$$
The set
$$
U:=\{x\in X:\,\Re\lb x,x^*\rb<c\}
$$
is weakly open, contains $x_0$, and is disjoint from $D$.
Thus $X\setminus D$ is weakly open, so $D$ is weakly closed.
\end{proof}

The next result characterises functionals that are continuous with respect to the weak and weak$^*$ topologies.

\begin{proposition}\label{prop:weak-dual} Let $X$ be a Banach space.  The following assertions hold:
 \begin{enumerate}[label={\rm(\arabic*)}, leftmargin=*]
  \item\label{it:weak-dual1} a linear mapping $\phi: X\to \K$ is continuous with respect to the weak topology if and only if it belongs to $X^*$;
  \item\label{it:weak-dual2} a linear mapping $\phi: X^*\to \K$ is continuous with respect to the weak$^*$ topology if and only if it belongs to $X$, that is, there exists an $x\in X$ such that $\phi(x^*) = \lb x,x^*\rb$ for all $x^*\in X^*$.
 \end{enumerate}
\end{proposition}
\begin{proof}
 \ref{it:weak-dual1}: \
 We only need to prove the `only if' part. If $\phi$ is weakly continuous, then $\phi^{-1}(B_\K)$ contains a weakly open set containing the origin. Since weakly open sets are open, this set contains a ball $B(0;r)$ with $r>0$. This means that $\phi\in X^*$ and $\n \phi\n \le 1/r$.

 \smallskip
 \ref{it:weak-dual2}: \ Again we only need to prove the `only if' part. If $\phi$ is weak$^*$ continuous, then $\phi^{-1}(B_\K)$ contains a set of the form $$ U:=\{x^*\in X^*:\, |\lb x_j,x^*\rb|<\eps,\quad j=1,\dots,k\} $$
 for suitable $\eps>0$, $k\ge 1$, and $x_1,\dots,x_k\in X$.
 Let
 $$ N:=\{x^*\in X^*:\, \lb x_j,x^*\rb=0,\quad j=1,\dots,k\}. $$
 If $x^*\in N$, then $cx^*\in U$ for every $c\in\K$, and therefore $|c\phi(x^*)|<1$ for every $c\in\K$. Hence $\phi(x^*)=0$, and thus $N\subseteq\Ker(\phi)$.
 Consider the linear mapping $ A:X^*\to\K^k$,
 $$Ax^*:=(\lb x_1,x^*\rb,\dots,\lb x_k,x^*\rb). $$
 Since $\Ker(A) = N \subseteq\Ker(\phi)$, the functional $\phi$ factors through $\Ran(A)$. By finite-dimensional linear algebra there exist scalars $c_1,\dots,c_k\in\K$ such that $$ \phi(x^*)=\sum_{j=1}^k c_j\lb x_j,x^*\rb =\Big\lb\sum_{j=1}^k c_jx_j,x^*\Big\rb,\quad x^*\in X^*. $$ Thus $\phi$ is represented by the element $\sum_{j=1}^k c_jx_j$ of $X$.
 \end{proof}

We next prove a version of the Hahn--Banach separation theorem
for the weak$^*$ topology.

\begin{proposition}\label{prop:HB-sep-weakstar}
 If $F$ is a weak$\s$ closed convex subset of $X^*$ and $x_0^*\not\in F$, then there exists an element $x_0\in X$ such that $\lb x_0,x_0^*\rb \not\in \lb x_0,F\rb $.
\end{proposition}

\begin{proof}
If $F=\emptyset$, we may take $x_0=0$.
Suppose therefore that $F\ne\emptyset$.
Since $F$ is weak$^*$ closed and $x_0^*\notin F$, there exist
$x_1,\dots,x_m\in X$ and $\eps>0$ such that
$$
U:=\bigl\{x^*\in X^*:\,
|\lb x_j,x^*-x_0^*\rb|<\eps,\ j=1,\dots,m\bigr\}
$$
is disjoint from $F$.

Define the linear mapping $A:X^*\to\K^m$ by
$$
Ax^*:=(\lb x_1,x^*\rb,\dots,\lb x_m,x^*\rb),
$$
and let
$$
V:=\bigl\{z\in\K^m:\,\n z-Ax_0^*\n_\infty<\eps\bigr\}.
$$
Then $U=A^{-1}(V)$, so $V$ and $A(F)$ are disjoint.
Moreover, $V$ is nonempty, open, and convex, and $A(F)$ is
nonempty and convex.

By the Hahn--Banach separation theorem
(Theorem \ref{thm:HB-separation}), there exists a linear
functional $\ell$ on $\K^m$ such that
$\ell(V)\cap\ell(A(F))=\emptyset$.
Since $Ax_0^*\in V$, it follows that
$\ell(Ax_0^*)\notin\ell(A(F))$.

Write $\ell(z)=\sum_{j=1}^m c_jz_j$ and put
$x_0:=\sum_{j=1}^m c_jx_j$.
Then
$$
\ell(Ax^*)=\sum_{j=1}^m c_j\lb x_j,x^*\rb
=\lb x_0,x^*\rb,\quad x^*\in X^*.
$$
Consequently,
\[
\lb x_0,x_0^*\rb
=\ell(Ax_0^*)
\notin\ell(A(F))
=\lb x_0,F\rb.
\qedhere\]
\end{proof}

Using the notation introduced in Definition \ref{def:annih} we have the following characterisation of weak and weak$^*$ closures of subspaces.

\begin{proposition}\label{prop:perpperp} Let $X$ be a Banach space.  The following assertions hold:
 \begin{enumerate}[label={\rm(\arabic*)}, leftmargin=*]
  \item\label{it:perpperp1} for every subspace $Y$ of $X$ we have ${}^\perp(Y^\perp) = \ov Y^{{\rm weak}} = \ov Y$;
  \item\label{it:perpperp2} for every subspace $Y$ of $X^*$ we have $({}^\perp Y)^\perp = \ov Y^{{\rm weak}^*}$\!.
 \end{enumerate}
\end{proposition}

\begin{proof}
\ref{it:perpperp1}: \ The inclusion ${}^\perp(Y^\perp)  \supseteq \ov Y^{{\rm weak}}$
follows from the easy observation that pre-annihilators are weakly closed, and the equality $\ov Y^{{\rm weak}} = \ov Y$
is a consequence of Proposition \ref{prop:convex-weaklyclosed}. To prove the inclusion ${}^\perp(Y^\perp)  \subseteq  \ov Y$, let $x\not \in \ov Y$. By Corollary \ref{cor:proper-subspace} we can find $x_0\s\in Y^\perp$
such that $\lb x,x_0\s\rb\not=0$. This means that $x$ is not in the pre-annihilator of
$Y^\perp$\!.

\smallskip
\ref{it:perpperp2}: \ This is proved in the same way, except that now we use Proposition \ref{prop:HB-sep-weakstar}.
\end{proof}

As a second application of Proposition \ref{prop:HB-sep-weakstar} we have the following dual version of Proposition \ref{prop:dual-dense-T}.
The proof is the same, except that instead of Corollary \ref{cor:proper-subspace} we now use Proposition \ref{prop:HB-sep-weakstar}.

\begin{proposition}\label{prop:dual-dense-T-dual} Let $X$ and $Y$ be Banach spaces and let $T\in \calL(X,Y)$. The following assertions are equivalent:

\begin{enumerate}[label={\rm(\arabic*)}, leftmargin=*]
 \item \label{it:dial-dense-T-dual-1} $T$ is injective;
 \item \label{it:dial-dense-T-dual-2} $T^*$ has weak$^*$ dense range.
 \end{enumerate}
 \end{proposition}

\section{The Banach--Alaoglu Theorem} \label{sec:BA}

We have seen in Chapter \ref{ch:Banach} that the closed unit ball of a Banach space is compact if and only if the space is finite-dimensional. In this section we prove that the closed unit ball of every dual Banach space is compact with respect to the weak$^*$ topology. In the next section we use this to prove that a Banach space is reflexive if and only if its closed unit ball is compact with respect to the weak topology.

We start with an elementary sequential result, which extends Proposition \ref{prop:Hilbert-weakconvergence} to separable Banach spaces.

\begin{proposition}\label{prop:BA-seq} Let $X$ be a separable Banach space
and let $(x_n\s)_{n\ge 1}$ be a bounded sequence in the dual space $X\s$\!. Then there exists a subsequence $(x_{n_k}\s)_{k\ge 1}$ and an $x\s\in X\s$ such that
$$ \limk \lb x,x_{n_k}\s\rb = \lb x,x\s\rb, \quad  x \in X.$$
\end{proposition}

\begin{proof}
Let $(x_n\s)_{n\geq 1}$ be a bounded sequence in $X\s$, say $\n x_n\s\n\leq M$ for all $n\geq 1$. Since $X$ is separable, there exists a sequence $(y_j)_{j\geq 1}$ whose linear span $Y$ is dense in $X$. For each fixed $j\geq 1$, the scalar sequence $(\lb y_j,x_n\s\rb)_{n\geq 1}$ is bounded. By a diagonal argument, there exists a subsequence $(x_{n_k}\s)_{k\geq 1}$ such that the limit $\limk\lb y_j,x_{n_k}\s\rb$ exists for every $j\geq 1$. By linearity, the limit $\limk\lb y,x_{n_k}\s\rb$ exists for all $y\in Y$.

For $k\geq 1$ define $T_k\in\calL(X,\K)$ by $T_kx:=\lb x,x_{n_k}\s\rb$. Then $\n T_k\n=\n x_{n_k}\s\n\leq M$ for all $k\geq 1$, and $(T_ky)_{k\geq 1}$ converges for all $y\in Y$. By Proposition \ref{prop:approxTn}, the limit $x\s(x):=\limk T_kx$ exists for all $x\in X$, and $x\s:X\to\K$ is a bounded linear functional. Thus $x_{n_k}\s\to x\s$ weak$\s$ in $X\s$.
\end{proof}

In contrast to the Hilbert space case considered in Proposition \ref{prop:Hilbert-weakconvergence}, the separability assumption in Proposition \ref{prop:BA-seq} cannot be omitted:

\begin{example}
Consider the sequence $(\phi_n)_{n\ge 1}$ of coordinate functionals $x\mapsto x_n$ on $X = \ell^\infty$\!.
Given a subsequence $(\phi_{n_k})_{k\ge 1}$, let $x\in \ell^\infty$ be any element such that $x_{n_k} = (-1)^k$ for all $k\ge 1$. Then $\lb x,\phi_{n_k}\rb = (-1)^k$ fails to converge as $k\to\infty$.
\end{example}

Proposition \ref{prop:BA-seq} can be viewed as a sequential version of the following theorem.

\begin{theorem}[Banach--Alaoglu] \label{thm:BA}\index{theorem!Banach--Alaoglu}\index{weak$^*$!compactness, of $\ov B_{X^*}$}\index{compactness!weak$^*$}
 The closed unit ball of the dual of a Banach space is compact with respect to the weak$\s$ topology.\index{compact!weak$\s$}
\end{theorem}

\begin{proof}
Let $\ov B_{\K}=\{a\in\K:\, |a|\leq 1\}$ and $\ov B_{X\s}=\{x\s\in X\s:\, \n x\s\n\leq 1\}$, where $X$ is a given Banach space. In the remainder of the proof we think of $\ov B_{X\s}$ as endowed with the weak$\s$ topology inherited from $X\s$.

By Tychonoff's theorem (Theorem \ref{thm:Tychonov}), the product space $K:=\prod_{x\in X}\n x\n\cdot\ov B_{\K}$ is compact with respect to the product topology. Denoting elements of $K$ by $k=(k_x)_{x\in X}$, consider the mapping $\phi:\ov B_{X\s}\to K$ given by
\begin{align}\label{eq:BA-phi}
\phi(x\s):=(\lb x,x\s\rb)_{x\in X}.
\end{align}
Let $R$ denote its range. We prove that $R$ is closed in $K$.

For $a_0,a_1\in\K$ and $x_0,x_1\in X$, the set of all $k\in K$ satisfying
$$a_0k_{x_0}+a_1k_{x_1}=k_{a_0x_0+a_1x_1}$$ is closed in $K$, since it is the inverse image of $\{0\}$ under the continuous mapping $k\mapsto a_0k_{x_0}+a_1k_{x_1}-k_{a_0x_0+a_1x_1}$. Hence the set $R_0$ of all $k\in K$ satisfying this identity for all $a_0,a_1\in\K$ and $x_0,x_1\in X$ is closed in $K$.

Clearly $R\subseteq R_0$. Conversely, let $k\in R_0$ and define $f:X\to\K$ by $f(x):=k_x$. The defining identities of $R_0$ imply that $f$ is linear. Since $k\in K$, we have $|f(x)|=|k_x|\leq \n x\n$ for all $x\in X$. Thus $f\in\ov B_{X\s}$, and by definition $\phi(f)=k$. This proves that $R_0\subseteq R$, and hence $R=R_0$. Therefore $R$ is closed in $K$, and consequently compact.

Since the mapping $\phi$ defined by \eqref{eq:BA-phi} is injective, the inverse mapping $\phi^{-1}:R\to\ov B_{X\s}$ is well defined. We claim that this mapping is continuous. Indeed, if $x\s\in\ov B_{X\s}$, $x\in X$, and $\eps>0$, then
$$
\phi\bigl(\{y\s\in\ov B_{X\s}:\, |\lb x,y\s-x\s\rb|<\eps\}\bigr)
=\{k\in R:\, |k_x-\lb x,x\s\rb|<\eps\},
$$
and the set on the right-hand side is open in $R$. Since such sets generate the weak$\s$ topology of $\ov B_{X\s}$, the inverse mapping $\phi^{-1}$ is continuous.
The weak$\s$ compactness of $\ov B_{X\s}=\phi^{-1}(R)$ now follows from the compactness of $R$.
\end{proof}

A topology $\tau$ is said to be {\em metrisable}\index{metrisable} if the underlying set
admits a metric whose open sets are precisely the sets of $\tau$.
The next result shows that if $X$ is separable, then the weak$\s$ topology of the unit ball of $X\s$ is metrisable. As a result, Proposition \ref{prop:BA-seq} can alternatively be deduced from Theorem \ref{thm:BA} by using the fact that compactness and sequential compactness are equivalent for metric spaces.

\begin{proposition}\label{prop:metrisable}
If $X$ is a separable Banach space, then the weak$\s$ topology of the closed unit ball of $X\s$ is metrisable.
\end{proposition}
\begin{proof}
Let $(x_n)_{n\ge1}$ be a dense sequence in the closed unit ball
$\ov B_X$, and define, for $x^*,y^*\in\ov B_{X^*}$,
$$
d(x^*,y^*):=
\sum_{n\ge1}2^{-n}|\lb x_n,x^*-y^*\rb|.
$$
Since $|\lb x_n,x^*-y^*\rb|\le2$, the series converges
uniformly on $\ov B_{X^*}\times\ov B_{X^*}$.
Symmetry and the triangle inequality follow from the corresponding
properties of the absolute value.
If $d(x^*,y^*)=0$, then $x^*-y^*$ vanishes on the dense
set $\{x_n:\,n\ge1\}$, hence on $\ov B_X$ by continuity,
and therefore on all of $X$.
Thus $d$ is a metric.

For fixed $x^*\in\ov B_{X^*}$, the function
$y^*\mapsto d(x^*,y^*)$ is weak$^*$ continuous, being a
uniform limit of weak$^*$ continuous functions.
Consequently, every open ball for $d$ is weak$^*$ open.
The identity mapping
$$
(\ov B_{X^*},\hbox{weak$^*$})
\to
(\ov B_{X^*},d)
$$
is therefore continuous.

By the Banach--Alaoglu theorem, its domain is compact.
It consequently maps closed sets to compact sets, which
are closed in the metric space $(\ov B_{X^*},d)$.
Its inverse is therefore continuous as well, proving that
$d$ induces the weak$^*$ topology.
\end{proof}

As an application of the Banach--Alaoglu theorem we have the following density result.

\begin{theorem}[Goldstine]\label{thm:Goldstine}\index{theorem!Goldstine}
Let $X$ be a Banach space.  The following assertions hold:
 \begin{enumerate}[label={\rm(\arabic*)}, leftmargin=*]
  \item $X$ is weak$^*$ dense in $X^{**}$;
  \item $\ov B_X$ is weak$^*$ dense in $\ov B_{X^{**}}$.
 \end{enumerate}
\end{theorem}

\begin{proof}
It suffices to prove the second assertion; the first follows from it by scaling arbitrary elements of $X^{**}$.
Arguing by contradiction, suppose that $x_0^{**}\in \ov B_{X^{**}}$ is not contained in the weak$^*$ closure $F$ of $B_X$.
Then by Proposition \ref{prop:HB-sep-weakstar} there exists an element $x_0^*\in X^*$ such that
 $\lb x_0^*,x_0^{**}\rb \not\in \lb x_0^*,F\rb $. Multiplying by a scalar, we
 may assume that $\n x_0^*\n = 1$. Then $B_X\subseteq  F\subseteq  \ov B_{X^{**}}$ implies $B_\K\subseteq  \lb x_0^*,F\rb \subseteq  \ov B_\K$. By the Banach--Alaoglu theorem $\ov B_{X^{**}}$ is weak$^*$ compact, hence so is $F$, and therefore $\lb x_0^*,F\rb$
 is a compact subset of $\K$. We conclude that
$ \lb x_0^*,F\rb = \ov B_\K$. As a consequence we have $|\lb x_0^*,x_0^{**}\rb|>1$.
 But this contradicts the fact that $\n x_0^{**}\n \le 1$ and $\n x_0^*\n=1$.
 \end{proof}

\section{Reflexivity}\label{subsec:reflexivity}

Recall that if $X$ is a Banach space, we can use the natural isometry $J:X\to X^{**}$ given by
$\lb x^*,Jx\rb := \lb x,x^*\rb$ to identify $X$ with a closed subspace of $X^{**}$\!.

\medskip
\begin{definition}[Reflexivity]
A Banach space $X$ is called {\em reflexive}\index{reflexive}
if the mapping $J:X\to X^{**}$ is surjective.
\end{definition}

The Banach--Alaoglu theorem implies the following characterisation of reflexivity.

\begin{theorem}[Reflexivity and weak compactness of the unit ball]\label{thm:refl-weakcomp}\index{weak!compactness, of $\ov B_{X}$}\index{compactness!weak}
 A Banach space is reflexive if and only if its closed unit ball is weakly compact.
\end{theorem}
\begin{proof}
The `only if' part follows from the Banach--Alaoglu theorem, noting that the canonical embedding $J:X\to X^{**}$ maps $\ov B_X$ onto $\ov B_{X^{**}}$ and that, under the identification of
$X$ and $X^{**}$\!, the weak topology of $X$ equals the weak$\s$ topology of $X^{**}$\!.

The `if' part follows from Goldstine's theorem (Theorem \ref{thm:Goldstine}). Indeed,
if $\ov B_X$ is weakly compact, then its image under the canonical embedding $J$ is weak$^*$ compact in $X^{**}$: this follows from the observation that $J$ is continuous as a mapping from $(X,$\,weak$)$ to $(X^{**}\!,$\,weak$^*$), which in turn is a trivial consequence of the definitions of the weak and weak$^*$ topologies. It follows that $\ov B_X$ is weak$^*$ compact as a subset of the closed unit ball $\ov B_{X^{**}}$. On the other hand, Goldstine's theorem says that $\ov B_X$
is weak$^*$ dense as a subset of $\ov B_{X^{**}}$. Hence we must have $\ov B_X = \ov B_{X^{**}}$, and this implies $X = X^{**}$\!.
\end{proof}

\begin{corollary}\label{cor:refl} Let $X$ be a Banach space.  The following assertions hold:
 \begin{enumerate}[label={\rm(\arabic*)}, leftmargin=*]
  \item\label{it:cor:refl1} if $X$ is reflexive, then every closed subspace of $X$ is reflexive;
  \item\label{it:cor:refl2} $X$ is reflexive if and only if $X^*$ is reflexive;
  \item\label{it:cor:refl3} if $X$ is isomorphic to a Banach space $Y$, then $X$ is reflexive if and only if $Y$ is reflexive.
 \end{enumerate}
\end{corollary}
\begin{proof}
\ref{it:cor:refl1}: \ If $Y$ is a closed subspace of $X$, the closed unit ball $\ov B_Y$ is the intersection of the set $\ov B_X$, which is weakly compact by  Theorem \ref{thm:refl-weakcomp}, and the set $Y$, which is weakly closed by Proposition \ref{prop:convex-weaklyclosed}.
As a result, $\ov B_Y$ is compact, as a subset of $X$, in the weak topology of $X$. By the Hahn--Banach theorem, this topology restricts to the weak topology of $Y$. It follows that $\ov B_Y$ is weakly compact in $Y$, and $Y$ is reflexive by another application of Theorem \ref{thm:refl-weakcomp}.

\smallskip
\ref{it:cor:refl2}: \ If $X$ is reflexive, the weak$^*$ and weak topologies of $X^*$ coincide. As a result, $\ov B_{X^*}$ is weakly compact by the Banach--Alaoglu theorem, and therefore $X^*$ is reflexive by Theorem \ref{thm:refl-weakcomp}. If $X^*$ is reflexive, then $X^{**}$ is reflexive by what we just proved, and then $X$, viewed as a closed subspace of $X^{**}$, is reflexive by part \ref{it:cor:refl1}.

\smallskip
\ref{it:cor:refl3}: \ If $i:X\to Y$ is an isomorphism, then $i^{**}: X^{**}\to Y^{**}$ is an isomorphism by Proposition \ref{prop:dual-iso} applied twice.
Denoting the natural isometries of $X$ and $Y$ into their second duals by $J_X$ and $J_Y$, one easily checks that
$i^{**}\circ J_X = J_Y \circ i$. This identity implies that $J_X$ is surjective if and only if $J_Y$ is surjective.
\end{proof}

Part \ref{it:cor:refl3} can also be deduced from Theorem \ref{thm:refl-weakcomp};
we leave this as an easy exercise.

\begin{corollary}\label{cor:refl-subseq}
Every bounded sequence in a reflexive Banach space has a weakly convergent subsequence.
\end{corollary}

\begin{proof} Let $(x_n)_{n\ge 1}$ be a bounded sequence in the reflexive Banach space $X$, and let $Y$ be its closed span. Then $Y$ is separable and reflexive by Corollary \ref{cor:refl}.

By Proposition \ref{prop:metrisable}, the closed unit ball of $Y^*$ is weak$^*$ compact and metrisable, and therefore contains a countable weak$^*$ dense set $D$. We claim that the closed linear span of $D$ equals $Y^*$. Otherwise, the Hahn--Banach theorem gives a nonzero element $y^{**}\in Y^{**}$ vanishing on $D$. Since $Y$ is reflexive, we have $y^{**}=Jy$ for some $y\in Y$. The mapping $y^*\mapsto\lb y,y^*\rb$ is weak$^*$ continuous, so its vanishing on the weak$^*$ dense set $D$ implies that it vanishes on all of $Y^*$. Hence $y=0$, a contradiction.

Thus $Y^*$ is separable.
By Proposition \ref{prop:BA-seq}, the bounded sequence
$(Jx_n)_{n\ge1}$ in $(Y^*)^*=Y^{**}$ has a weak$^*$
convergent subsequence, with limit $y^{**}\in Y^{**}$.
Since $Y$ is reflexive, there exists $y\in Y$ such that
$y^{**}=Jy$.
Hence
$$
\lim_{k\to\infty}\lb x_{n_k},y^*\rb
=\lb y,y^*\rb,\quad y^*\in Y^*.
$$
This proves weak convergence in $Y$.
It also proves weak convergence in $X$, since every functional
in $X^*$ restricts to an element of $Y^*$.
\end{proof}

\begin{example}\label{ex:reflexive} The following are examples of reflexive Banach spaces:
 \begin{itemize}
  \item finite-dimensional Banach spaces;
  \item the spaces $\ell^p$ with $1<p<\infty$;
  \item the spaces $L^p(\Om)$ with $1<p<\infty$;
  \item Hilbert spaces.
 \end{itemize}

If $X$ is finite-dimensional, then $\dim X^{**}=\dim X$.
The canonical embedding $J:X\to X^{**}$ is injective and
therefore surjective, so $X$ is reflexive.

The second and third examples follow by applying twice the isometric identifications
$(\ell^p)^*=\ell^q$ and $(L^p(\Om))^*=L^q(\Om)$,
where $\frac1p+\frac1q=1$.
The resulting identifications with the second duals are
precisely the canonical embeddings.

That Hilbert spaces are reflexive is a consequence of the Riesz representation theorem (Theorem \ref{thm:Riesz}) which sets up a conjugate-linear identification of the dual $H^*$ of a Hilbert space with the Hilbert space $H$ itself.
Applying the theorem twice and composing the identifications, we see that the resulting identification of $H$ with $H^{**}$ is given by the natural embedding
$J:H\to H^{**}$\!.
\end{example}

\begin{example}\label{ex:nonrefl} The spaces $c_0$, $\ell^1$ and $\ell^\infty$ are nonreflexive: for $c_0$ this follows from the fact that $c_0^{**} = \ell^\infty$ and for $\ell^1 = c_0^*$ and $\ell^\infty = (\ell^1)^*$ this follows from
Corollary \ref{cor:refl}.

The spaces $C(K)$ and $L^1(\Omega)$ are nonreflexive except when they
are finite-dimensional. Indeed, if $K$ is an infinite compact Hausdorff
space, then $K$ contains pairwise disjoint nonempty open sets
$(U_n)_{n\ge1}$. By Urysohn's lemma we may choose functions
$f_n\in C(K)$ such that $\supp(f_n)\subseteq U_n$ and $\|f_n\|_\infty=1$.
Then
$$
        (a_n)_{n\ge1}\mapsto \sum_{n\ge1} a_n f_n
$$
defines an isomorphic embedding of $c_0$ into $C(K)$.

Similarly, if $L^1(\Omega)$ is infinite-dimensional, one can choose
pairwise disjoint measurable sets $E_n$ with
$0<\mu(E_n)<\infty$. The mapping
$$
        (a_n)_{n\ge1}\mapsto
        \sum_{n\ge1} a_n\,\frac{\one_{E_n}}{\mu(E_n)}
$$
is an isometric embedding of $\ell^1$ into $L^1(\Omega)$. Thus in both
cases nonreflexivity follows from Corollary \ref{cor:refl}.
\end{example}

\section{Applications}\label{sec:duality-applications}

In this section we apply the duality and compactness results developed above to trans\-lation-invariant operators, probability measures, marginal problems, and optimal transport.

\subsection{Translation-Invariant Operators on \texorpdfstring{$L^1(\R^d)$}{L\textasciicircum 1(R\textasciicircum d)}}

The first is concerned with characterising translation-invariant operators on $L^1(\R^d)$ as convolutions with a finite Borel measure. It is complemented by Theorem \ref{thm:transl-invar-L2} in the next chapter, which characterises the translation-invariant operators
on $L^2(\R^d)$ as the Fourier multiplier operators.

\begin{theorem}[Translation-invariant operators on $L^1(\R^d)$] \label{thm:transl-invar-L1}\index{translation!invariant, operators on $L^1(\R^d)$}
If $T$ is a bounded operator on $L^1(\R^d)$ commuting with every translation,
then $T$ is given by convolution with a (necessarily unique) measure
$\mu\in M(\R^d)$; more precisely, for all $f\in C_{\rm c}(\R^d)$ we have, for almost all $x\in \R^d$,
$$ Tf(x) = \int_{\R^d} f(x-y)\ud \mu(y).$$
Moreover, $\n T\n = \n\mu\n$.
\end{theorem}

\begin{proof}
Fix a nonnegative function $\eta\in L^1(\R^d)$ satisfying $\int_{\R^d} \eta(x)\ud x =1$. For $\eps>0$ the rescaled functions $\eta_\eps(x):= \eps^{-d}\eta(\eps^{-1}x)$ belong to $L^1(\R^d)$ and satisfy $\n \eta_\eps\n_1 = \n \eta\n_1$. By viewing the functions $T\eta_\eps\in L^1(\R^d)$ as densities of finite Borel measures on $\R^d$ we may identify them with finite Borel measures $\mu_\eps \in M(\R^d)$.
By Theorem \ref{thm:CK-dual} and Proposition \ref{prop:Radon}, $M(\R^d)$ can be identified isometrically with the dual of $C_0(\R^d)$.
Moreover, for all $\eps>0$,
$$
\|\mu_\varepsilon\|=\|T\eta_\varepsilon\|_1 \le \|T\|\,\|\eta_\varepsilon\|_1=\|T\|\,\|\eta\|_1\le \|T\|.
$$

Since $C_0(\R^d)$ is separable (this follows from Proposition \ref{prop:CK-separable}), by the sequential version of Banach--Alaoglu (Proposition \ref{prop:BA-seq}) some sequence $(T\eta_{\eps_n})_{n\ge 1}$ with $\eps_n\downarrow 0$ converges weak$^*$ to a measure $\mu\in M(\R^d)$.
This measure satisfies
$$ \|\mu\|\le \liminf_{n\to\infty}\|\mu_{\eps_n}\|\le \|T\|.$$

For $f\in C_{\rm c}(\R^d)$ and $g\in L^1(\R^d)$, the mapping
$y\mapsto f(y)\tau_{-y}g$ is continuous and compactly supported
as an $L^1(\R^d)$-valued function by Proposition
\ref{prop:Lp-transl}, and we have
$$
f*g=\int_{\R^d}f(y)\tau_{-y}g\ud y.
$$
Using the commutation assumption, we obtain
\begin{align*}
T(f*g)
=\int_{\R^d}f(y)T\tau_{-y}g\ud y
=\int_{\R^d}f(y)\tau_{-y}Tg\ud y
=f*Tg.
\end{align*}
In particular, Theorem \ref{thm:approx-identity} gives
$ f*T\eta_{\eps_n} =T(f*\eta_{\eps_n}) \to Tf$ in $L^1(\R^d)$.
On the other hand, for every $x\in\R^d$ the function
$y\mapsto f(x-y)$ belongs to $C_{\rm c}(\R^d)$.
Weak$^*$ convergence of the measures $\mu_{\eps_n}$ therefore gives
$$
(f*T\eta_{\eps_n})(x)
=\int_{\R^d}f(x-y)\ud\mu_{\eps_n}(y)
\to
\int_{\R^d}f(x-y)\ud\mu(y).
$$
Passing to an almost everywhere convergent subsequence, we conclude that
$$
Tf(x)=\int_{\R^d}f(x-y)\ud\mu(y)
\ \hbox{for almost all }x\in\R^d.
$$

By Fubini's theorem,
$$
\n Tf\n_1
\le\int_{\R^d}\int_{\R^d}|f(x-y)|\ud|\mu|(y)\ud x
=\n\mu\n\n f\n_1.
$$
Since $C_{\rm c}(\R^d)$ is dense in $L^1(\R^d)$, it follows
that $\n T\n\le\n\mu\n$. Combining this with the opposite
inequality proved above gives $\n T\n=\n\mu\n$.

It remains to establish the uniqueness of the measure $\mu$. Suppose that $\mu\in M(\R^d)$ induces the zero convolution operator. For every $f\in C_{\rm c}(\R^d)$, the uniform continuity of $f$ and the finiteness of $|\mu|$ imply that the function
$x\mapsto\int_{\R^d}f(x-y)\ud\mu(y)$
is continuous. Since it vanishes almost everywhere, it vanishes identically. Taking $x=0$ and replacing $f$ by the function $y\mapsto\phi(-y)$, we obtain
$$
\int_{\R^d}\phi(y)\ud\mu(y)=0,\quad \phi\in C_{\rm c}(\R^d).
$$
By Proposition \ref{prop:Cc-dense-in-C0}, the preceding identity
extends from $C_{\rm c}(\R^d)$ to $C_0(\R^d)$. Theorem
\ref{thm:CK-dual} therefore implies that $\mu=0$.
\end{proof}

\subsection{Prokhorov's Theorem}\label{subsec:Prokhorov}

The aim of this section is to prove a compactness result of fundamental importance in probability theory, known as Prokhorov's theorem. For its statement we need the following terminology.

\begin{definition}[Uniform tightness] A collection $P$ of Borel measures on a Hausdorff topological space $X$
is called {\em uniformly tight}\index{uniformly!tight}
if for every $\eps>0$ there exists a compact set
$K$ in $X$ such that $\mu(X\setminus K)<\eps$ for all $\mu\in P$.
\end{definition}

\begin{definition}[Weak convergence]
A sequence $(\mu_n)_{n\ge 1}$ of Borel probability measures on a Hausdorff topological space $X$ is said to
{\em converge weakly}\index{weak!convergence, of measures} to a Borel probability measure $\mu$ on $X$
if
$$\limn \int_X f\ud\mu_n = \int_X f \ud\mu,\quad f\in C_{\rm b}(X),$$
where $C_{\rm b}(X)$ denotes the Banach space of bounded continuous
functions on $X$.
\end{definition}

Viewing Borel probability measures on $X$
as functionals in the dual space $(C_{\rm b}(X))^*$, weak convergence in the sense
of the above definition is precisely weak$^*$ convergence in the sense discussed in the present chapter.
The terminology `weak convergence' is firmly established in the probability theory literature, however.

\begin{theorem}[Prokhorov's theorem]\label{thm:Prokhorov}\index{theorem!Prokhorov}
Let $X$ be a metric space. Then:
\begin{enumerate}[label={\rm(\arabic*)}, leftmargin=*]
\item\label{it:Prokhorov1}  if a family $P$ of Borel probability measures on $X$ is uniformly tight, then every sequence in $P$ contains a weakly convergent subsequence;
\item\label{it:Prokhorov2}  if $X$ is separable and complete,
then every weakly convergent sequence of Borel probability measures on $X$ is uniformly tight.
\end{enumerate}
\end{theorem}

For the proof of this theorem we need the following characterisation of weak convergence,
known as the {\em Portmanteau theorem}.\index{theorem!Portmanteau}

\begin{proposition}\label{prop:Portmanteau} Let $\mu_n$, $n\ge 1$, and $\mu$ be Borel probability measures on a metric space $X$. The following assertions are equivalent:
\begin{enumerate}[label={\rm(\arabic*)}, leftmargin=*]
\item\label{it:portmanteau1}
$\limn \mu_n = \mu$ weakly;
\item\label{it:portmanteau2} for all open subsets $U$ of $X$ we have
$ \mu(U) \le \liminf_{n\to\infty} \mu_n(U)$;
\item\label{it:portmanteau3} for all closed subsets $F$ of $X$ we have
$ \limsup_{n\to\infty} \mu_n(F)\le \mu(F)$.
\end{enumerate}
\end{proposition}

\begin{proof}
\ref{it:portmanteau1}$\Rightarrow$\ref{it:portmanteau2}:\
Let $U\subseteq X$ be open. If $U=X$, the assertion is immediate.
Otherwise, the functions
$$
f_k(x):=\min\{1,k\,d(x,X\setminus U)\},\quad k\ge1,
$$
belong to $C_{\rm b}(X)$ and increase pointwise to $\one_U$.
For every $k\ge1$, weak convergence gives
$$
\int_X f_k\ud\mu
=\lim_{n\to\infty}\int_X f_k\ud\mu_n
\le\liminf_{n\to\infty}\mu_n(U).
$$
Letting $k\to\infty$ and applying monotone convergence proves
\ref{it:portmanteau2}.

\smallskip
The equivalence
\ref{it:portmanteau2}$\Leftrightarrow$\ref{it:portmanteau3}
follows by taking complements.

\smallskip
\ref{it:portmanteau2}$\Rightarrow$\ref{it:portmanteau1}:\
Let $f\in C_{\rm b}(X)$ satisfy $0\le f\le M$.
For any Borel probability measure $\nu$ on $X$, Fubini's
theorem gives
$$
\int_X f\ud\nu
=\int_X\int_0^M\one_{\{t<f(x)\}}\ud t\ud\nu(x)
=\int_0^M\nu(\{f>t\})\ud t.
$$
Since the sets $\{f>t\}$ are open, assumption
\ref{it:portmanteau2} and Fatou's lemma imply
\begin{align*}
\int_X f\ud\mu
=\int_0^M\mu(\{f>t\})\ud t
&\le\int_0^M\liminf_{n\to\infty}\mu_n(\{f>t\})\ud t\\
&\le\liminf_{n\to\infty}\int_0^M\mu_n(\{f>t\})\ud t
=\liminf_{n\to\infty}\int_X f\ud\mu_n.
\end{align*}
Applying the same inequality to $M-f$, we obtain
\begin{align*}
\limsup_{n\to\infty}\int_X f\ud\mu_n
=M-\liminf_{n\to\infty}\int_X(M-f)\ud\mu_n
\le M-\int_X(M-f)\ud\mu
=\int_X f\ud\mu.
\end{align*}
Together these inequalities prove convergence of the integrals
for nonnegative $f\in C_{\rm b}(X)$.
The real-valued case follows by considering positive and negative
parts, and the complex-valued case by considering real and
imaginary parts.
\end{proof}

The proof of part \ref{it:Prokhorov1} of Prokhorov's theorem relies on the following lemma.

\begin{lemma}\label{lem:piece-together} Let $(\Om,\calF)$ be a measurable space and let $A_1\subseteq  A_2\subseteq  \dots$
be an increasing sequence of sets in $\calF$ such that $\bigcup_{j\ge 1} {A_j} = \Om$. For each $j\ge 1$ let
$\mu_j$ be a measure on $(A_j,\calF|_{A_j})$. If the measures $\mu_j$ are increasing
in the sense that $ \mu_j|_{A_i} \ge  \mu_i$ on $\calF|_{A_i}$ whenever $j\ge i,$
then
$$ \mu(B):= \lim_{j\to\infty} \mu_j(B\cap A_j),\quad B\in\calF\!,$$
defines a measure on $\Om$.
\end{lemma}
\begin{proof}
For $B\in\calF$ and $j\ge i$, the assumption gives
$$
\mu_i(B\cap A_i)
\le\mu_j(B\cap A_i)
\le\mu_j(B\cap A_j).
$$
Thus the limit defining $\mu(B)$ exists.
Clearly, $\mu(\emptyset)=0$.
If $(B_n)_{n\ge1}$ is a sequence of pairwise disjoint sets
in $\calF$, countable additivity of the measures $\mu_j$
and monotone convergence give
\begin{align*}
\mu\Bigl(\bigcup_{n\ge1}B_n\Bigr)
=\lim_{j\to\infty}\sum_{n\ge1}\mu_j(B_n\cap A_j)
=\sum_{n\ge1}\lim_{j\to\infty}\mu_j(B_n\cap A_j)
=\sum_{n\ge1}\mu(B_n).
\end{align*}
This proves that $\mu$ is a measure.
\end{proof}

\begin{proof}[Proof of Theorem \ref{thm:Prokhorov}]
\ref{it:Prokhorov1}:  \ Let $(\mu_n)_{n\ge1}$ be a sequence in the uniformly tight family
$P$. Choose compact sets $L_j\subseteq X$ such that
$\mu_n(L_j)\ge1-2^{-j}$ for all $n\ge1$.
Putting $K_j:=L_1\cup\cdots\cup L_j$, we obtain an increasing
sequence of compact sets satisfying
$$
\mu_n(K_j)\ge1-2^{-j},\quad n,j\ge1.
$$

{\em Step 1} --
For each $j$, we identify the measures $\mu_n|_{K_j}$ with
elements of the closed unit ball of $(C(K_j))^*$.
Since $K_j$ is a compact metric space, $C(K_j)$ is separable
by Proposition \ref{prop:CK-separable}.
Applying the sequential Banach--Alaoglu theorem
(Proposition \ref{prop:BA-seq}) successively for $j=1,2,\dots$
and taking a diagonal subsequence, we obtain a subsequence
$(\mu_{n_k})_{k\ge1}$ such that every sequence
$(\mu_{n_k}|_{K_j})_{k\ge1}$ converges weak$^*$.
The limits are positive functionals and therefore are represented
by nonnegative Borel measures $\nu_j$ on $K_j$. Thus
$$
\lim_{k\to\infty}\int_{K_j}f\ud\mu_{n_k}
=\int_{K_j}f\ud\nu_j,\quad f\in C(K_j).
$$
Taking $f=\one$ gives
$$
1-2^{-j}
\le\nu_j(K_j)
=\lim_{k\to\infty}\mu_{n_k}(K_j)
\le1.
$$

{\em Step 2} --
Extend each $\nu_j$ by zero outside $K_j$, so that it becomes
a Borel measure on $X$.
We claim that $\nu_i\le\nu_j$ whenever $i\le j$.
For every nonnegative $f\in C(K_j)$, Step 1 gives
\begin{align*}
\int_{K_i}f\ud\nu_i
=\lim_{k\to\infty}\int_{K_i}f\ud\mu_{n_k}
\le\lim_{k\to\infty}\int_{K_j}f\ud\mu_{n_k}
=\int_{K_j}f\ud\nu_j.
\end{align*}
By the positivity assertion in the Riesz representation theorem,
$\nu_j-\nu_i$ is a nonnegative measure on $K_j$.
Both measures vanish outside $K_j$, so $\nu_i\le\nu_j$ on $X$.

\smallskip
{\em Step 3} --
By Lemma \ref{lem:piece-together}, applied with $A_j=X$ for all
$j$, the formula
$$
\mu(B):=\lim_{j\to\infty}\nu_j(B),\quad B\in\calB(X),
$$
defines a Borel measure on $X$.
The bounds in Step 1 imply
$$
\mu(X)=\lim_{j\to\infty}\nu_j(K_j)=1,
$$
so $\mu$ is a probability measure.
Moreover, $\mu-\nu_j$ is nonnegative and
$$
(\mu-\nu_j)(X)=1-\nu_j(K_j)\le2^{-j}.
$$

We prove that $\mu_{n_k}\to\mu$ weakly.
Let $f\in C_{\rm b}(X)$ and fix $j\ge1$.
Splitting the integral against $\mu_{n_k}$ into its parts
on $K_j$ and $X\setminus K_j$, we obtain
\begin{align*}
\Bigl|\int_X f\ud\mu_{n_k}\!-\!\int_X f\ud\mu\Bigr|
&\le
\Bigl|\int_{K_j}f\ud\mu_{n_k}
      \!-\!\int_{K_j}f\ud\nu_j\Bigr|
+\n f\n_\infty\mu_{n_k}(X\setminus K_j)
      +\n f\n_\infty(\mu-\nu_j)(X).
\end{align*}
The first term tends to zero by Step 1, while each of the
last two terms is at most $2^{-j}\n f\n_\infty$.
Consequently,
$$
\limsup_{k\to\infty}
\Bigl|\int_X f\ud\mu_{n_k}-\int_X f\ud\mu\Bigr|
\le2^{1-j}\n f\n_\infty.
$$
Letting $j\to\infty$ proves the required weak convergence.

\smallskip
\ref{it:Prokhorov2}: \
Suppose that $\mu_n\to\mu$ weakly, and fix $\eps>0$.
For $j\ge1$, put
$r_j:=2^{-j}$ and $\delta_j:=\eps\,2^{-j-1}$.
We first construct, for each $j$, a finite union $U_j$ of
open balls of radius $r_j$ such that
$$
\mu_n(X\setminus U_j)<\delta_j,\quad n\ge1.
$$

Fix $j\ge1$. Since $X$ is separable, it has a countable cover
by open balls of radius $r_j$.
By continuity from below, a finite union $V_j$ of these balls
satisfies $\mu(V_j)>1-\delta_j$.
The Portmanteau theorem gives
$$
\liminf_{n\to\infty}\mu_n(V_j)
\ge\mu(V_j)>1-\delta_j.
$$
Hence there exists $N_j$ such that
$\mu_n(V_j)>1-\delta_j$ for all $n\ge N_j$.
For each of the finitely many measures
$\mu_1,\dots,\mu_{N_j-1}$, continuity from below likewise
provides a finite union of balls from the same cover having
measure greater than $1-\delta_j$.
Adding these balls to $V_j$ gives the required finite union
$U_j$.

Let $F_j$ be the union of the closed balls having the same
centres and radii as the balls making up $U_j$. The set
$K:=\bigcap_{j\ge1}F_j$ is closed, and for every $j$ it is covered
by finitely many open balls of radius $2r_j$.
Since $r_j\downarrow0$, this shows that $K$ is totally bounded.
As a closed subset of the complete space $X$, it is also
complete, and therefore compact by
Theorem \ref{thm:totally-bounded}.

Finally, since $U_j\subseteq F_j$, for every $n\ge1$ we have
$$
\mu_n(X\setminus K)
\le\sum_{j\ge1}\mu_n(X\setminus F_j)
\le\sum_{j\ge1}\delta_j
=\frac{\eps}{2}<\eps.
$$
This proves uniform tightness.
\end{proof}

\subsection{Marginal Problems}
\label{subsec:marginals}

Let $X_1,\ldots,X_n$ be nonempty compact metric spaces and put
$X:=X_1\times\cdots\times X_n$. For a nonempty subset
$C\subseteq\{1,\ldots,n\}$, write
$$
X_C:=\prod_{j\in C}X_j,\qquad
\pi_C:X\to X_C,\qquad \pi_C(x):=(x_j)_{j\in C}.
$$
We denote the set of Borel probability measures on a compact metric
space $Y$ by $\mathcal P(Y)$. If $\mu\in\mathcal P(X)$, the image
measure $(\pi_C)_\#\mu$ is called its {\em marginal} on $X_C$.
All measures under consideration are Radon, and the Borel
$\sigma$-algebra of a finite product agrees with the product of the
Borel $\sigma$-algebras of its factors.

Suppose that $\mathscr C$ is a nonempty collection of nonempty subsets
of $\{1,\ldots,n\}$, and that a measure $\mu_C\in\mathcal P(X_C)$
has been specified for every $C\in\mathscr C$. The {\em marginal
problem}\index{marginal problem} asks whether there exists a measure
$\mu\in\mathcal P(X)$ satisfying
$$
(\pi_C)_\#\mu=\mu_C,\quad C\in\mathscr C.
$$
For such a measure to exist, prescribed marginals must agree on their
common coordinates. This condition alone need not suffice.

If only the individual coordinate marginals are prescribed, a joint
measure always exists: one may take their product. When distributions
on overlapping sets of coordinates are prescribed, the existence
question is more delicate. The following theorem gives a necessary
and sufficient condition in terms of continuous functions.

\begin{theorem}[Kellerer's criterion]\label{thm:epi-kellerer}
\index{theorem!Kellerer}
A Borel probability measure on $X$ with the prescribed marginals
exists if and only if
\begin{align}\label{eq:epi-kellerer}
\sum_{C\in\mathscr C}\int_{X_C}f_C\ud\mu_C
\le \max_{x\in X}\sum_{C\in\mathscr C}f_C(\pi_Cx)
\end{align}
for every family $f_C\in C(X_C;\R)$, $C\in\mathscr C$.
\end{theorem}

\begin{proof}
If $\mu$ has the prescribed marginals, the substitution formula gives
$$
\sum_{C\in\mathscr C}\int_{X_C}f_C\ud\mu_C
=\int_X\sum_{C\in\mathscr C}f_C\circ\pi_C\ud\mu
\le\max_{x\in X}\sum_{C\in\mathscr C}f_C(\pi_Cx).
$$
This proves necessity.

\smallskip
Conversely, suppose that \eqref{eq:epi-kellerer} holds. Let $V$ be
the real vector subspace of $C(X;\R)$ consisting of all functions
of the form $g=\sum_{C\in\mathscr C}f_C\circ\pi_C$. We wish to define
a linear functional on $V$ by
$$
L\Bigl(\sum_{C\in\mathscr C}f_C\circ\pi_C\Bigr)
:=\sum_{C\in\mathscr C}\int_{X_C}f_C\ud\mu_C.
$$
To check that this is well defined, suppose that
$\sum_C f_C\circ\pi_C=0$. Applying \eqref{eq:epi-kellerer} to
$(f_C)_C$ and to $(-f_C)_C$ gives
$$
\sum_C\int_{X_C}f_C\ud\mu_C\le0,
\qquad
-\sum_C\int_{X_C}f_C\ud\mu_C\le0.
$$
Thus the sum of integrals is zero. Subtracting two representations
of the same function now shows that $L$ is well defined, and its
linearity follows from the linearity of the integral.

Applying \eqref{eq:epi-kellerer} to $g$ and $-g$ gives
$|L(g)|\le\n g\n_\infty$. Moreover, $\one\in V$ and $L(\one)=1$,
so $\n L\n=1$. By the Hahn--Banach theorem, $L$ extends to a
functional $\Lambda$ on $C(X;\R)$ with $\n\Lambda\n=1$.
This extension is positive. Indeed, if $0\le g\le\one$, then
$\n\one-g\n_\infty\le1$ and therefore
$$
\Lambda(g)=1-\Lambda(\one-g)\ge0.
$$
Scaling gives $\Lambda(g)\ge0$ for every $g\ge0$.

By the Riesz representation theorem (Theorem \ref{thm:CK-dual}),
there is a positive Borel measure $\mu$ on $X$ such that
$\Lambda(g)=\int_Xg\ud\mu$. Its total mass is $\Lambda(\one)=1$.
For $C\in\mathscr C$ and $f\in C(X_C;\R)$ we have
$$
\int_{X_C}f\ud((\pi_C)_\#\mu)
=\Lambda(f\circ\pi_C)
=L(f\circ\pi_C)
=\int_{X_C}f\ud\mu_C.
$$
Uniqueness in the Riesz representation theorem implies
$(\pi_C)_\#\mu=\mu_C$, as required.
\end{proof}

The theorem gives a useful way of disproving the existence of a joint
measure: it suffices to find continuous functions for which
\eqref{eq:epi-kellerer} fails. Such functions express an obstruction
to realising the prescribed marginals simultaneously.
We illustrate this in the special case where
$$X  =  \{-1,1\}\times \{-1,1\}\times \{-1,1\}\times \{-1,1\}.$$
We write points of $X$ as $(a_1,a_2,b_1,b_2)$. Suppose that the marginal law $\mu_{ij}$ of each pair $(a_i,b_j)$, $i,j=1,2$, has been prescribed. A common joint law is a probability measure $\mu$ on $X$ with these four pair marginals.
For $i,j=1,2$ put
$$
E_{ij}:=\int_{\{-1,1\}\times \{-1,1\}}ab\ud\mu_{ij}(a,b).
$$
A necessary condition for the existence of such a common joint law is given by
the following inequality, named after Clauser, Horne, Shimony, and Holt.

\begin{corollary}[CHSH inequality]\label{cor:epi-chsh-classical}
\index{inequality!CHSH}
If the four measures $\mu_{ij}$ defined above have a common joint law, then
\begin{align}\label{eq:epi-chsh-classical}
|E_{11}+E_{12}+E_{21}-E_{22}|\le2.
\end{align}
\end{corollary}

\begin{proof}
For every $(a_1,a_2,b_1,b_2)\in X$ we have
$$
a_1b_1+a_1b_2+a_2b_1-a_2b_2
=a_1(b_1+b_2)+a_2(b_1-b_2)\in\{-2,2\},
$$
since exactly one of $b_1+b_2$ and $b_1-b_2$ is zero and the other
equals $2$ or $-2$. Taking the four test functions $ab,ab,ab,-ab$
in Theorem \ref{thm:epi-kellerer}, and then their negatives, gives
\eqref{eq:epi-chsh-classical}.
\end{proof}

In Subsection \ref{subsec:epi-quantum-correlations} we will construct
four pair distributions arising from quantum measurements whose
one-coordinate marginals agree, but which violate
\eqref{eq:epi-chsh-classical}. Thus agreement of the one-coordinate marginals
does not suffice for the existence of a joint distribution.

\subsection{Optimal Transport and Kantorovich Duality}
\label{subsec:optimal-transport}

Having considered the existence of measures with prescribed marginals,
we now turn to an optimisation problem. If we think of a probability measure
as describing a mass distribution, two probability measures
on a metric space can be compared by the cost of transporting one
distribution of mass into the other. A {\em transport plan} describes
how this mass is redistributed. For compact metric spaces and
continuous costs, we will prove that a plan of minimum cost exists
and express this minimum in terms of continuous functions. When the
cost is the distance travelled, this construction gives a metric
on the set of probability measures.

Throughout this section, $(X,d)$ is a nonempty compact metric space.
We write $\mathcal P(X)$ and $\mathcal P(X\times X)$ for the sets of Borel probability measures on $X$ and $X\times X$,
and $C(X;\R)$ and $C(X\times X;\R)$ for the real Banach spaces of continuous real-valued functions on $X$ and $X\times X$, respectively, both with the supremum norm.
All finite Borel measures on these spaces are Radon by Proposition \ref{prop:Radon}.

Since the Borel $\sigma$-algebra of $X\times X$ agrees with the product of the
Borel $\sigma$-algebras of $X$ (the proof is the same as the one for Euclidean spaces in the appendix),
the product measure
$\mu\times\nu$ supplied by Theorem \ref{thm:product-measures} is a Borel
probability measure on $X\times X$ whenever $\mu,\nu\in\mathcal P(X)$.

For $\pi\in\mathcal P(X\times X)$, its first and second {\em marginals}
are the probability measures
$$
\pi_1(B):=\pi(B\times X),\qquad
\pi_2(B):=\pi(X\times B),\quad B\in\calB(X).
$$
A {\em coupling}\index{coupling, of probability measures} of
$\mu,\nu\in\mathcal P(X)$ is a probability measure $\pi$ with
$\pi_1=\mu$ and $\pi_2=\nu$. We denote the set of all such couplings by
$$
\Pi(\mu,\nu):=
\{\pi\in\mathcal P(X\times X):\,\pi_1=\mu,\ \pi_2=\nu\}.
$$
This set is nonempty, since it contains $\mu\times\nu$.
The marginals are the image measures under the coordinate projections,
so the substitution formula of Example \ref{ex:subst}, applied to positive
and negative parts, gives
\begin{align}\label{eq:transport-marginals}
\int_{X\times X}\bigl(\varphi(x)+\psi(y)\bigr)\ud\pi(x,y)
=\int_X\varphi\ud\mu+\int_X\psi\ud\nu
\end{align}
for all bounded real-valued Borel functions $\varphi,\psi$ on $X$ and
all $\pi\in\Pi(\mu,\nu)$.

We may think of a coupling as a transport plan, interpreting the quantity $\pi(B\times C)$ as the amount of mass transported from $B$ to $C$. If $c(x,y)$ is the cost of moving a
unit of mass from $x$ to $y$, the total cost of the plan is $\int c\ud\pi$.
The following theorem expresses the minimum cost as a supremum over
pairs of continuous functions. For $\varphi,\psi\in C(X;\R)$, we use
the notation
$$
(\varphi\oplus\psi)(x,y):=\varphi(x)+\psi(y),\quad x,y\in X.
$$

\begin{theorem}[Kantorovich duality]\label{thm:kantorovich-duality}
\index{theorem!Kantorovich duality}
Let $\mu,\nu\in\mathcal P(X)$ and let $c\in C(X\times X;\R)$.
Then
$$
\min_{\pi\in\Pi(\mu,\nu)}\int_{X\times X}c\ud\pi
=\sup_{\substack{\varphi,\psi\in C(X;\R)\\
                  \varphi\oplus\psi\le c}}
\Bigl(\int_X\varphi\ud\mu+\int_X\psi\ud\nu\Bigr).
$$
In particular, a coupling attaining the minimum exists.
\end{theorem}

\begin{proof}
We call the pair $(\varphi,\psi)$ {\em majorising} for $g$ if $g\le\varphi\oplus\psi$.
By \eqref{eq:transport-marginals}, every pair $(\varphi,\psi)$ satisfying $\varphi\oplus\psi\le c$ gives a lower bound for the cost of every coupling. It therefore suffices to construct a coupling
attaining the supremum.

On the real Banach space $E:=C(X\times X;\R)$ define
$$
p(g):=\inf\Bigl\{
\int_X\varphi\ud\mu+\int_X\psi\ud\nu:\,
\varphi,\psi\in C(X;\R),\ g\le\varphi\oplus\psi
\Bigr\}.
$$
Integration against $\mu\times\nu$ and the majorising pair
$(\|g\|_\infty,0)$ give
\begin{align}\label{eq:transport-p-bounds}
-\|g\|_\infty
\le\int_{X\times X}g\ud(\mu\times\nu)
\le p(g)\le\|g\|_\infty.
\end{align}
Adding majorising pairs and multiplying them by
nonnegative scalars shows that $p$ is sublinear.
For $g=\varphi\oplus\psi$, the pair $(\varphi,\psi)$
attains the lower bound, so
\begin{align}\label{eq:transport-p-marginals}
p(\varphi\oplus\psi)
=\int_X\varphi\ud\mu+\int_X\psi\ud\nu.
\end{align}

Since $-p(-c)\le p(c)$, the linear functional
$t(-c)\mapsto tp(-c)$ on the span of $-c$ is dominated
by $p$. The Hahn--Banach theorem
(Theorem~\ref{thm:appHB1}) extends it to a linear
functional $\Lambda:E\to\R$ satisfying
$$
\Lambda(g)\le p(g),\quad g\in E,\qquad
\Lambda(-c)=p(-c).
$$
Applying \eqref{eq:transport-p-bounds} to $g$ and $-g$
gives $|\Lambda(g)|\le\|g\|_\infty$.
Moreover, if $g\ge0$, then $p(-g)\le0$, and hence
$\Lambda(g)\ge0$. By the Riesz representation theorem,
$\Lambda$ is therefore represented by a positive Borel
measure $\pi$ on $X\times X$.

Applying the domination inequality to
$\varphi\oplus\psi$ and its negative, and using
\eqref{eq:transport-p-marginals}, gives
$$
\int_{X\times X}(\varphi\oplus\psi)\ud\pi
=\int_X\varphi\ud\mu+\int_X\psi\ud\nu.
$$
Taking $(\varphi,\psi)=(1,0)$ shows that $\pi$ is a
probability measure. Taking $\psi=0$, and then
$\varphi=0$, identifies its marginals as $\mu$ and $\nu$
by uniqueness in the Riesz representation theorem.
Thus $\pi\in\Pi(\mu,\nu)$.

Finally, replacing $\varphi,\psi$ by their negatives
in the definition of $p(-c)$ gives
\[
\int_{X\times X}c\ud\pi
=-p(-c)
=\sup_{\substack{\varphi,\psi\in C(X;\R)\\
                 \varphi\oplus\psi\le c}}
\Bigl(\int_X\varphi\ud\mu+\int_X\psi\ud\nu\Bigr).
\qedhere\]
\end{proof}

We now take the distance $d$ itself as the cost function. For
$\mu,\nu\in\mathcal P(X)$, the quantity
$$
W_1(\mu,\nu):=
\min_{\pi\in\Pi(\mu,\nu)}
\int_{X\times X}d(x,y)\ud\pi(x,y)
$$
is called the {\em Wasserstein distance of order one}.\index{Wasserstein distance}
The minimum exists by the preceding theorem.

Recall that a function $f:X\to\R$ is {\em $1$-Lipschitz} if
$|f(x)-f(y)|\le d(x,y)$ for all $x,y\in X$.

\begin{theorem}[Kantorovich--Rubinstein]\label{thm:KR-W1}
\index{theorem!Kantorovich--Rubinstein}
For all $\mu,\nu\in\mathcal P(X)$ we have
\begin{align}\label{eq:W1-duality}
W_1(\mu,\nu)
=\max\Bigl\{\int_X f\ud\mu-\int_X f\ud\nu:\,
f:X\to\R \ \hbox{ is $1$-Lipschitz}\Bigr\}.
\end{align}
Moreover, $W_1$ is a metric on $\mathcal P(X)$.
\end{theorem}

\begin{proof}
If $f$ is $1$-Lipschitz, then $(f,-f)$ is admissible in
Theorem \ref{thm:kantorovich-duality} with $c=d$. Hence the supremum
on the right-hand side of \eqref{eq:W1-duality} is at most
$W_1(\mu,\nu)$.
Conversely, let $\varphi,\psi\in C(X;\R)$ satisfy
$\varphi(x)+\psi(y)\le d(x,y)$ for all $x,y\in X$, and put
$$
f(x):=\inf_{y\in X}\bigl(d(x,y)-\psi(y)\bigr),\quad x\in X.
$$
This infimum is finite, since $d$ and $\psi$ are bounded.
For $x,x',y\in X$, the triangle inequality gives
$$
d(x,y)-\psi(y)
\le d(x,x')+d(x',y)-\psi(y).
$$
Taking infima over $y$ and interchanging $x$ and $x'$ shows that
$|f(x)-f(x')|\le d(x,x')$.
Admissibility gives $\varphi\le f$, while taking $y=x$ in the
definition of $f$ gives $f\le-\psi$. Consequently,
$$
\int_X\varphi\ud\mu+\int_X\psi\ud\nu
\le\int_X f\ud\mu-\int_X f\ud\nu.
$$
Taking the supremum over admissible pairs and applying
Theorem \ref{thm:kantorovich-duality} proves the desired equality
with a supremum in place of the maximum.

Fix $x_0\in X$. Subtracting the constant $f(x_0)$ does not change the
difference of the two integrals. We may therefore restrict the
supremum to the set
\begin{align}\label{eq:W1-normalised-Lipschitz}
\mathscr L:=\{f\in C(X;\R):\,
f(x_0)=0,\ |f(x)-f(y)|\le d(x,y)\ \hbox{for all }x,y\in X\}.
\end{align}
Every $f\in\mathscr L$ satisfies
$|f(x)|\le d(x,x_0)\le\operatorname{diam}(X)$.
Thus $\mathscr L$ is uniformly bounded and equicontinuous. It is
also closed in the supremum norm, since its defining conditions
are preserved under uniform limits. By the Arzel\`a--Ascoli theorem
(Theorem \ref{thm:Arzela-Ascoli}), $\mathscr L$ is compact in
$C(X;\R)$. The mapping
$f\mapsto\int_X f\ud\mu-\int_X f\ud\nu$ is continuous, so its
supremum on $\mathscr L$ is attained.

It remains to verify the metric properties. Nonnegativity follows
from the definition of $W_1$, and \eqref{eq:W1-duality} gives
$W_1(\mu,\mu)=0$. Replacing $f$ by $-f$ gives symmetry.
For $\lambda\in\mathcal P(X)$ and any $1$-Lipschitz $f$,
\begin{align*}
\int_X f\ud\mu-\int_X f\ud\lambda
&=\Bigl(\int_X f\ud\mu-\int_X f\ud\nu\Bigr)
 +\Bigl(\int_X f\ud\nu-\int_X f\ud\lambda\Bigr)
\\
&\le W_1(\mu,\nu)+W_1(\nu,\lambda).
\end{align*}
Taking the supremum gives the triangle inequality.
Finally, if $W_1(\mu,\nu)=0$, a minimising coupling $\pi$ satisfies
$d(x,y)=0$ for $\pi$-almost all $(x,y)$.
It is therefore concentrated on the diagonal. For every Borel set
$B\subseteq X$, the functions $\one_B(x)$ and $\one_B(y)$ agree
on this diagonal, so integration gives $\mu(B)=\nu(B)$.
Hence $\mu=\nu$.
\end{proof}

The two descriptions of $W_1$ provide a useful way to recognise an
optimal plan: a coupling gives an upper bound for $W_1$, and a
$1$-Lipschitz function gives a lower bound; if these bounds agree,
the coupling is optimal.

\begin{example}\label{ex:W1-interval}
Let $X=[0,1]$ with its usual metric, let $\mu$ be Lebesgue measure on
$[0,1]$, and let $\nu:=\frac12\delta_0+\frac12\delta_1$.
Consider the measurable map
$$
\tau(x):=
\begin{cases}
0,&0\le x\le\frac12,\\
1,&\frac12<x\le1.
\end{cases}
$$
The image of $\mu$ under $x\mapsto(x,\tau(x))$ is a coupling of
$\mu$ and $\nu$. By the substitution formula, its cost is
$$
\int_0^1|x-\tau(x)|\ud x
=\int_0^{1/2}x\ud x+\int_{1/2}^1(1-x)\ud x
=\frac14.
$$
On the other hand, the function $f(x):=\min\{x,1-x\}$ is
$1$-Lipschitz and satisfies
$$
\int_X f\ud\mu-\int_X f\ud\nu
=\int_0^1\min\{x,1-x\}\ud x
=\frac14.
$$
Hence this coupling is optimal and $W_1(\mu,\nu)=\frac14$.
\end{example}

The preceding example extends to arbitrary
$\mu,\nu\in\mathcal P([0,1])$. Define their distribution
functions by
$F_\mu(x):=\mu([0,x])$ and $F_\nu(x):=\nu([0,x]),\quad x\in[0,1]$,
and their {\em quantile functions} by
$$
G_\mu(t):=\inf\{x\in[0,1]:\,F_\mu(x)\ge t\},
\quad t\in(0,1),
$$
and similarly for $G_\nu$. These functions are nondecreasing
and satisfy
$$
G_\mu(t)\le x\ \Longleftrightarrow\ t\le F_\mu(x).
$$
Consequently, the image of Lebesgue measure on $(0,1)$ under
$G_\mu$ is $\mu$, and similarly for $\nu$. The image under
$t\mapsto(G_\mu(t),G_\nu(t))$
therefore defines a coupling $\pi$ of $\mu$ and $\nu$, the
so-called {\em monotone coupling}: since both
quantile functions are nondecreasing, $s<t$ implies
$G_\mu(s)\le G_\mu(t)$ and $G_\nu(s)\le G_\nu(t)$.

We now show that this coupling is optimal for the cost $|x-y|$.
The elementary identity
$$
|a-b|
=\int_0^1|\one_{\{a\le s\}}-\one_{\{b\le s\}}|\ud s,
\quad a,b\in[0,1],
$$
and Fubini's theorem give
\begin{align*}
\int_{[0,1]^2}|x-y|\ud\pi(x,y)
&=\int_0^1|G_\mu(t)-G_\nu(t)|\ud t\\
&=\int_0^1\int_0^1
  |\one_{\{t\le F_\mu(s)\}}-\one_{\{t\le F_\nu(s)\}}|
  \ud t\ud s =\int_0^1|F_\mu(s)-F_\nu(s)|\ud s.
\end{align*}
For the reverse bound, put
$$
f(x):=\int_0^x
\operatorname{sgn}(F_\nu(s)-F_\mu(s))\ud s,
\quad x\in[0,1],
$$
where $\operatorname{sgn}(0):=0$.
This function is $1$-Lipschitz, and Fubini's theorem gives
\begin{align*}
\int_{[0,1]}f\ud\mu-\int_{[0,1]}f\ud\nu
&=\int_0^1\operatorname{sgn}(F_\nu(s)-F_\mu(s))
  \bigl(\mu((s,1])-\nu((s,1])\bigr)\ud s\\
&=\int_0^1|F_\mu(s)-F_\nu(s)|\ud s.
\end{align*}
By Theorem~\ref{thm:KR-W1}, the coupling $\pi$ is optimal and
$$
W_1(\mu,\nu)
=\int_0^1|G_\mu(t)-G_\nu(t)|\ud t
=\int_0^1|F_\mu(x)-F_\nu(x)|\ud x.
$$
In Example~\ref{ex:W1-interval}, this
construction recovers the function $f(x)=\min\{x,1-x\}$.

We conclude by relating this metric to the weak convergence of
probability measures considered in Section \ref{subsec:Prokhorov}.

\begin{corollary}\label{cor:W1-weak-convergence}
For $\mu_n,\mu\in\mathcal P(X)$, we have
$W_1(\mu_n,\mu)\to0$ if and only if $\mu_n\to\mu$ weakly.
Moreover, $(\mathcal P(X),W_1)$ is a compact metric space.
\end{corollary}

\begin{proof}
Suppose first that $\mu_n\to\mu$ weakly, and fix $\eps>0$.
By compactness of the set $\mathscr L$ in
\eqref{eq:W1-normalised-Lipschitz}, there exist
$f_1,\dots,f_N\in\mathscr L$ such that every $f\in\mathscr L$
satisfies $\n f-f_j\n_\infty<\eps$ for some $j$.
Since $\mu_n$ and $\mu$ are probability measures,
$$
W_1(\mu_n,\mu)
\le 2\eps+
\max_{1\le j\le N}
\Bigl|\int_X f_j\ud\mu_n-\int_X f_j\ud\mu\Bigr|.
$$
The maximum tends to zero by weak convergence. Letting
$\eps\downarrow0$ proves $W_1(\mu_n,\mu)\to0$.

Conversely, suppose that $W_1(\mu_n,\mu)\to0$ and let
$g\in C(X;\R)$. Given $\eps>0$, choose $\delta>0$ such that
$|g(x)-g(y)|<\eps$ whenever $d(x,y)<\delta$.
For all $x,y\in X$ we then have
$$
|g(x)-g(y)|\le\eps+\frac{2\n g\n_\infty}{\delta}d(x,y),
$$
using that $d(x,y)\ge \delta$ implies $|g(x)-g(y)|\le
2\n g\n_\infty\le 2\n g\n_\infty \delta^{-1}d(x,y)$.
Choose a minimising coupling $\pi_n\in\Pi(\mu_n,\mu)$.
Integrating the preceding estimate gives
$$
\Bigl|\int_X g\ud\mu_n-\int_X g\ud\mu\Bigr|
\le\eps+\frac{2\n g\n_\infty}{\delta}W_1(\mu_n,\mu).
$$
Letting $n\to\infty$ and then $\eps\downarrow0$ proves weak
convergence, since it suffices to test real-valued continuous functions.

Finally, $\mathcal P(X)$ is uniformly tight because $X$ is compact.
By Prokhorov's theorem (Theorem \ref{thm:Prokhorov}), every sequence
in $\mathcal P(X)$ has a
weakly convergent subsequence, and hence a subsequence converging in
$W_1$. Sequential compactness and compactness coincide for metric
spaces, proving the last assertion.
\end{proof}

\begin{problems}

\item\label{prob:mult-fc-CK}
Let $K$ be a compact Hausdorff space and let
$\phi\in (C(K))^*$ be an element with the following two properties:
\begin{enumerate}[leftmargin=*, label=(\roman*)]
  \item $\phi(\one) = 1$;
  \item $\phi(fg) = \phi(f)\phi(g)$ for all $f,g\in C(K)$.
\end{enumerate}
Prove that there exists an $x\in K$ such that $$\phi(f)=f(x), \quad f\in C(K).$$

\noindent{\em Hint:}\ Show that $\phi$, as an element of $M(K)$, is supported on a singleton.

\item\label{prob:no-sigma-finite}
Show that the result of Remark \ref{rem:no-sigma-finite} does not extend to $p=1$.

\item
Let $X$ be a Banach lattice.

\begin{enumerate}[\rm(a), leftmargin=*]
 \item Show that $x\in X$ satisfies $x\ge 0$ if and only if $\lb x,x\s\rb \ge 0$ for all $x\s\in X\s$ satisfying $x\s\ge 0$.
\item Show that if $x\s\in X\s$ satisfies $x\s\ge 0$,
then for all $x\in X$ we have the following assertions:
\begin{enumerate}[\rm(i), leftmargin=*]
 \item $\lb x^+, x^*\rb = \sup\{\lb x,y^*\rb:\, 0\le y^*\le x^*\}$;
 \item $\lb |x|, x^*\rb = \sup\{\lb x,y^*\rb:\, |y^*|\le x^*\}$.
\end{enumerate}
\end{enumerate}

\item
In this problem we work over the real scalar field.
Under the assumptions of Theorem \ref{thm:CK-dual}, let $\mu\in M_{\rm R}(X)$ represent the functional $\phi\in (C_0(X))^*$.
Show that the measures representing $\phi^+$, $\phi^-$, and $|\phi|$ are $\mu^+$, $\mu^-$, and $|\mu|$, respectively.

\item\label{prob:FR-ops}
Let $X$ and $Y$ be Banach spaces. A bounded operator $T\in \calL(X,Y)$ is said to be of {\em finite rank} if its range is finite-dimen\-sional.
Show that every finite-rank operator $T\in \calL(X,Y)$ is of the form
$$Tx = \sum_{n=1}^N \lb x,x_n\s\rb y_n, \quad x\in X,$$ for certain $y_n\in Y$ and $x_n\s\in X\s$.

\item\label{prob:dist-formula}
Let $X_0$ be a proper closed subspace of a Banach space $X$ and let
$x\in X\setminus X_0$.
\begin{enumerate}[\rm(a), leftmargin=*]
\item\label{it:dist-formula1a}
On $Y:=X_0+\K x$ define
$$
\phi(x_0+tx):=t,\quad x_0\in X_0,\ t\in\K.
$$
Show that $\phi\in Y^*$ and
$\n\phi\n=1/{d(x,X_0)}$.
\item\label{it:dist-formula1b}
Deduce from the Hahn--Banach theorem that there exists a norm-one vector
$x^*\in X^*$ such that
$$
\lb x,x^*\rb=d(x,X_0), \quad x^*|_{X_0}=0.
$$
\item\label{it:dist-formula2}
Using part \ref{it:dist-formula1b}, show that there exists
$x^*\in(L^\infty(0,1))^*$ such that
$$
\lb f,x^*\rb=\int_0^1f(t)\ud t,\quad f\in C[0,1],
$$
but
$\lb\one_{(0,\frac12)},x^*\rb\neq\frac12$.
 \end{enumerate}

\item\label{prob:kernel-quotient}
Let $X$ be a Banach space.
\begin{enumerate}[\rm(a), leftmargin=*]
\item Show that for all nonzero $x^* \in X^*$ we have an isomorphism of Banach spaces
  $$ X/\Ker(x\s) \simeq \K,$$
  where $\Ker(x^*):= \{x\in X:\, \lb x,x^*\rb=0\}$.
\item Let $Q: X\to X/\Ker (x\s)$ denote the quotient mapping. Prove that if $\|x^*\|=1$, then for all $x\in X$ we have $$ \|Qx\|  = |\lb x,x^*\rb |.$$

\noindent{\em Hint:}\ The case $x\in\Ker(x^*)$ is immediate.
Otherwise, for $0<\eps<1$, choose a unit vector $u_\eps\in X$
such that $|\lb u_\eps,x^*\rb|>1-\eps$, and write
$u_\eps=c_\eps x+y_\eps$ with $y_\eps\in\Ker(x^*)$.
Show that
$$
|c_\eps|>\frac{1-\eps}{|\lb x,x^*\rb|},
\qquad \n Qu_\eps\n\le1.
$$
Use $Qu_\eps=c_\eps Qx$ and let $\eps\downarrow0$ to obtain
the upper bound for $\n Qx\n$. For the reverse inequality,
use $\n x^*\n=1$ and the definition of the quotient norm.
\end{enumerate}

\item
Find an example of a two-dimensional Banach space $X$ and a functional on one of its closed one-dimensional subspaces which has infinitely many extensions to a functional on $X$ of the same norm.

\item\label{prob:sepdual}
Prove that if the dual of a Banach space $X$ is separable, then $X$ is separable.

\noindent{\em Hint:}\ Let $(x_n\s)_{n\ge 1}$ be a countable dense set in the unit sphere of $X^*$\!, and let $(x_n)_{n\ge 1}$ be a sequence in the unit sphere of $X$
such that $|\lb x_n,x_n\s\rb| \ge \frac12$ for all $n\ge 1$.

\item
We consider Banach spaces containing -- and contained in -- $\ell^\infty$\!.
\begin{enumerate}[\rm(a), leftmargin=*]
  \item\label{it:containellinfty1} Let $X$ be a Banach space, $Y$ be a closed subspace of $X$, and let $T_0:Y \to \ell^{\infty}$
  be a bounded operator. Show that there exists a bounded operator $T: X \to \ell^{\infty}$ such that $T|_{Y} = T_0$
  and $\|T\| = \|T_0\|$.

  \noindent{\em Hint:}\ Apply the Hahn--Banach theorem `coordinatewise'.

  \item\label{it:containellinfty2} Using the result of part \ref{it:containellinfty1}, prove that if a closed subspace $Y$ of a Banach space $X$ is isomorphic to $\ell^{\infty}$, then $Y$ is complemented in $X$.

  \item\label{it:containellinfty3} Show that every separable Banach space $X$ is isometrically isomorphic to a closed subspace of $\ell^\infty$\!.
  More precisely, show that if $X$ is a separable Banach space, then there exists a closed subspace $Y$ of $\ell^\infty$
  and an isometric isomorphism $T$ from $X$ onto $Y$.

  \noindent{\em Hint:}\ Use the Hahn--Banach theorem in combination with the separability to make
  a clever choice of a sequence of functionals $(x_n^*)_{n\ge 1}$ in $X^*$,
  and consider the mapping $T: x\mapsto (\lb x,x_n^*\rb)_{n\ge 1}.$
\end{enumerate}

\item\label{prob:weak-Cauchy}
Consider an open set $D\subseteq \C$ and a complex Banach space $X$. A function $f:D\to X$ is said to be
{\em holomorphic}\index{holomorphic!$X$-valued function} if for all $z_0\in D$ the limit
$$ \lim_{z\to z_0} \frac{f(z)-f(z_0)}{z-z_0}$$
exists in $X$.
Use the Hahn--Banach theorem to prove that the Cauchy theorem\index{theorem!Cauchy, for Banach space-valued functions} and the Cauchy integral formula hold for holomorphic functions $f:D\to X$ defined on an open set $D$ in $\C$:
\begin{align*}
\frac1{2\pi i} \int_{\{ |z-z_0|= r\}} f(z)\ud z & = 0
\intertext{and}
\frac1{2\pi i} \int_{\{ |z-z_0|= r\}} \frac{f(z)}{z-z_0}\ud z & = f(z_0).
\end{align*}
Here it is assumed that  $z_0\in D$ and $r>0$ is so small that the closed disc $\{ |z-z_0|\le r\}$ is contained in $D$;
the contour $\{ |z-z_0|= r\}$ is oriented anticlockwise.

\item\label{prob:strictly-convex}
Recall from Problem \ref{prob:Hilbert-strictly-convex} that a Banach space $X$ is called {\it strictly convex}\index{strictly convex} if for all norm-one vectors $x_0,x_1\in X$ with $x_0\not=x_1$ and $0<\la<1$ we have $$\n (1-\la) x_0 + \la x_1\n <1.$$
Prove that if the dual $X\s$ of a Banach space is strictly convex, then every functional on a closed subspace of $X$ has a {\em unique} Hahn--Banach extension of the same norm.

\item\label{prob:Haar-property}
Let $Y$ be a closed subspace of a Banach space $X$.
\begin{enumerate}[\rm(a), leftmargin=*]
  \item\label{it:Haar-property1} Prove that for all  $x\s\in X\s$ we
  have $$d(x\s\!,Y^\perp) = \big\n x\s|_{Y}\big\n_{Y\s}.$$
\item\label{it:Haar-property2} Show that for every
$x\s\in X\s$ the infimum defining $d(x\s,Y^\perp)$ is attained.
\item\label{it:Haar-property3}
Prove that every functional $y\s\in Y\s$ has a unique extension
to a functional in $X\s$ of the same norm if and only if,
for every $x\s\in X\s$, there exists a unique $z\s\in Y^\perp$
such that $d(x\s,Y^\perp)=\n x\s-z\s\n$.
\end{enumerate}

\item
Let $(A_n)_{n\ge 1}$ be a sequence of disjoint
intervals of positive measure $|A_n|$ in the interval $[0,1]$
and define $f:[0,1]\to c_0$ by $$f(t) = \sum_{n\ge 1} \frac1{|A_n|}\one_{A_n}(t)  u_n,$$
where $(u_n)_{n\ge 1}$ is the sequence of standard unit vectors of $c_0$.
\begin{enumerate}[\rm(a), leftmargin=*]
  \item Show that $f$ is strongly measurable.
  \item Show that for all $x^*\in c_0^* = \ell^1$ the integral $\int_0^1\lb f(t),x^*\rb\ud t$ is well defined.
  \item Show that $f$ fails to be Bochner integrable.
\end{enumerate}

\item
Consider the mapping $f:(0,1)\to L^\infty(0,1)$ given by $f(t):= \one_{(0,t)}$.
\begin{enumerate}[\rm(a), leftmargin=*]
  \item Show that $\lb f,x^*\rb$ is measurable for all $x^*\in (L^\infty(0,1))^*$\!.

  \noindent{\em Hint:}\ Monotone scalar-valued functions are measurable.
  \item Show that $f$ fails to be Bochner integrable.
\end{enumerate}

\item\label{prob:dual-iso}
Let $X$ and $Y$ be Banach spaces.
Prove that if $T:X\to Y$ is an isomorphism from $X$ onto its range in $Y$, then the adjoint operator $T\s:Y\s \to X\s$ is surjective.

\item\label{prob:indef-Lp}
Let $1\leq p< \infty$.
\begin{enumerate}[\rm(a), leftmargin=*]
  \item\label{it:indef-Lp1} Let $T:L^p(0,1) \to {\mathbb K}$ be the linear mapping defined by
  $$ Tf := \int_0^1 f(s)\ud s, \quad  f\in L^p(0,1). $$
  Show that $T$ is bounded and find an expression for $T^*$\!.

  \item\label{it:indef-Lp2} Let $T:L^p(0,1) \to L^p(0,1)$ be the linear operator defined by
  $$(Tf)(t) := \int_0^t f(s)\ud s, \quad t\in (0,1), \ f\in L^p(0,1).$$
  Show that $T$ is bounded and find an expression for $T^*$\!.
  \item\label{it:indef-Lp3} Let $T$ be the linear operator of part \ref{it:indef-Lp2},
  now viewed as an operator from $L^p(0,1)$ into $C[0,1]$.
  Show that $T$ is bounded and find an expression for $T^*$\!.
\end{enumerate}

\item
Let $H_0$ be a closed subspace of a Hilbert space $H$ and let $i:H_0\to H$ be the inclusion mapping.
Show that the adjoint $i^\star: H\to H_0$ is the orthogonal projection of $H$ onto $H_0$, viewed as a mapping from $H$ to $H_0$.

\item\label{prob:ergodic}
Let $H$ be a Hilbert space and let $T\in \calL(H)$ be a contraction, that is, $\n T\n \le 1$.
\begin{enumerate}[\rm(a), leftmargin=*]
  \item\label{it:ergodic2} Show that for each $x\in H$ we have $Tx = x$ if and only if $T^\star x=x$. Conclude that
  $H$ admits an orthogonal direct sum decomposition
  $$ H = \Ker(I-T) \oplus \overline{\Ran(I-T)}. $$
  \noindent {\em Hint:}\  If $Tx=x$, show that $T^\star x-x \perp x$ and deduce
  that $T^\star x =x$.

  \item\label{it:ergodic3}
  Define $$S_{n}:=\frac 1n \sum_{k=0}^{n-1}T^k\!, \quad n\ge 1.$$
  Using Proposition \ref{prop:injective-denserange}  and the result of part \ref{it:ergodic2},
  prove that
  $$\limn S_{n} x  = \begin{cases}
                     x & \hbox{if } x\in  \Ker(I-T), \\
                     0 & \hbox{if } x\perp \Ker(I-T).
                   \end{cases}
  $$
\end{enumerate}

\item\label{prob:ch1ex1}
Let $C$ be a closed convex subset of a nonzero separable Banach space $X$.
Prove that there exists a sequence $(x_n\s)_{n\ge 1}$ of norm-one elements
in $X\s$ and a sequence $(F_n)_{n\ge 1}$ of closed sets in $\K$
such that $$C= \bigcap_{n\ge 1}
\big\{x\in X:\, \lb x,x_n\s\rb \in  F_n\big\}.$$

\noindent{\em Hint:}\ Separate $C$ from countably many open balls whose union is $\complement C$.

\item
Over the real scalar field, find the extreme points of the closed convex set $$C = \{f\in L^2(0,1):\, f \ge 0, \, \n f\n_2\le 1\}.$$

\item
\begin{enumerate}[\rm(a), leftmargin=*]
  \item Show that if $x_n\to x$ weakly in $X$, then
$\n x\n\le\liminf_{n\to\infty}\n x_n\n$.
\item Show that if $x_n^*\to x^*$ weak$^*$ in $X^*$, then
$\n x^*\n\le\liminf_{n\to\infty}\n x_n^*\n$.
\end{enumerate}

\item
Let $1<p<\infty$ and let $(x^{(n)})_{n\ge1}$ be a bounded sequence in
$\ell^p$. Show that $x^{(n)}\to x$ weakly in $\ell^p$ if and only if
$x_k^{(n)}\to x_k$ for all $k\ge 1$.

\noindent{\em Hint:}\ Approximate elements of $\ell^q$, $\frac1p+\frac1q=1$, by finitely supported sequences.

\item
Prove that the Borel $\sigma$-algebra of a separable Banach space $X$ is the smallest $\sigma$-algebra relative to which all functionals $x\s\in X\s$ are measurable.

\item\label{prob:weak-cont-bdd}
Let $X$ and $Y$ be Banach spaces. Prove the following assertions:
\begin{enumerate}[\rm(a), leftmargin=*]
  \item\label{it:weak-cont-bdd1}
  a linear operator $T:X\to Y$ is continuous with respect to the weak topologies of $X$ and $Y$
  if and only if it is bounded;
  \item\label{it:weak-cont-bdd2}
  a linear operator $S:Y^*\to X^*$ is continuous with respect to the weak$^*$ topologies of $Y^*$ and $X^*$ if and only if it is the adjoint of some $T\in\calL(X,Y)$.
\end{enumerate}
\noindent{\em Hint:}\ For the `only if' implication in part \ref{it:weak-cont-bdd1}, show that
the graph is norm closed and use the closed graph theorem
(Theorem \ref{thm:CGT}). For part (b), first construct the
pre-adjoint using Proposition \ref{prop:weak-dual}.

\item
Prove that the weak topology of a Banach space $X$ coincides with the norm topology if and only if $X$ is finite-dimensional.

\item
Let $(e_n)_{n\ge1}$ be the standard unit vectors in
$\ell^1=(c_0)^*$. Show that $e_n\to 0$ weak$^*$ in $\ell^1$,
but not weakly.

\item
Find an example of a sequence $(x_n^*)_{n\ge 1}$ in a dual Banach space $X^*$ with the following two properties:
\begin{enumerate}[leftmargin=*, label=(\roman*)]
  \item there exists an $x^*\in X^*$ such that
  $ \limn \lb x, x_n^*\rb = \lb x,x^*\rb$ for all $x\in X$;
  \item no sequence in the convex hull of $(x_n^*)_{n\ge 1}$ converges to $x^*$ strongly in $X^*$\!.
\end{enumerate}
Compare with Corollary \ref{cor:Mazur}.

\item
Prove that $c_0$, $C[0,1]$, and
$C_{\rm b}(D)$ with $D\subseteq\R^d$ open, identified in the natural way with subspaces of
$\ell^\infty$, $L^\infty[0,1]$, and $L^\infty(D)$, respectively, are
norm closed and weak$^*$ dense in these spaces.

\item\label{prob:weak-star-density-L1-in-measures}
Let $X$ be a locally compact Hausdorff space and let $\mu$ be a
positive finite Radon measure on $X$.
Prove that $L^1(X,\mu)$ is weak$^*$ dense in $M_{\rm R}(X)$
if and only if $\mu$ has full support, that is, $\mu(U)>0$ for every nonempty open set $U\subseteq X$. Here we identify functions $f\in L^1(X,\mu)$ with the measures $f\mu \in M_{\rm R}(X)$.

\noindent{\em Hint:}\
First prove that if $\mu$ has full support, then every Dirac measure
$\delta_x$, $x\in X$, belongs to the weak$^*$ closure of $L^1(X,\mu)$.

\item\label{prob:universal-CK}
Let $X$ be a Banach space and endow
$K:=\ov B_{X^*}$ with the weak$^*$ topology. By the
Banach--Alaoglu theorem, the Hausdorff space $K$ is compact. For $x\in X$ define
$Jx:K\to\K$ by
$$
(Jx)(x^*):=\lb x,x^*\rb.
$$
Show that $Jx\in C(K)$ and that $J:X\to C(K)$ is a linear isometry.
Conclude that every Banach space is isometrically isomorphic to a closed
subspace of a space $C(K)$ with $K$ compact.

\item
Prove the following converse to Proposition
\ref{prop:metrisable}:
If the weak$\s$ topology of the closed unit ball of the dual of a Banach space $X$ is metrisable, then $X$ is separable.

\noindent {\em Hint:}\ Complete the details of the following argument.
Let $d$ be a metric which induces the weak$\s$ topology of $\ov{B}_{X\s}$.
Then $(\ov{B}_{X\s},d)$ is a compact metric
space and therefore the Banach space
$C(\ov{B}_{X\s},d)$ is separable by Proposition \ref{prop:CK-separable}. Now observe that $X$ is isometrically
contained in $C(\ov{B}_{X\s},d)$ in a natural way.

\item\label{prob:BanachLimit}
As an application of the Banach--Alaoglu theorem, prove that there exist functionals $\phi\in (\ell^\infty)^*$ such that for all $x = (x_n)_{n\ge 1} \in \ell^\infty$ we have:
\begin{enumerate}[leftmargin=*, label=(\roman*)]
  \item $\lb x,\phi\rb \ge 0$ whenever $x\ge 0$;
  \item $\lb x,\phi\rb = \lb Sx,\phi\rb$, where $S:(x_n)_{n\ge 1} \mapsto (x_{n+1})_{n\ge 1}$ is the left shift;
  \item $\lb x,\phi\rb = \limn x_n$ whenever $\limn x_n$ exists.
\end{enumerate}
Functionals with these properties are called {\em Banach limits}.\index{Banach!limit}

\noindent{\em Hint:}\ Consider the functionals $\phi_n(x) := \frac1n \sum_{j=1}^n x_j$.

\item
Let $(x_n)_{n\ge 1}$ be the sequence in $\ell^\infty$ defined by
$$ x_n := (\underbrace{0,\dots,0}_{n\ {\rm{times}}}, 1,1,\dots), \quad n\ge 1.$$
\begin{enumerate}[\rm(a), leftmargin=*]
  \item
  Show that this sequence has no weakly convergent subsequence.

  \noindent{\em Hint:}\ Use the result of Problem \ref{prob:BanachLimit}.

  \item Why doesn't this contradict the Banach--Alaoglu theorem?
\end{enumerate}

\item\label{prob:Goldstine}
Deduce Theorem \ref{thm:Goldstine} from Theorem \ref{thm:locrefl}.

\item\label{prob:KM-weakstar}
This problem proves a weak$^\ast$ Krein--Milman theorem.
Let $X$ be a Banach space and let $K\subseteq X^\ast$ be nonempty, convex,
and weak$^\ast$ compact. Denote by $\operatorname{Extr}(K)$ the set of
extreme points of $K$.

\begin{enumerate}[\rm(a),leftmargin=*]
\item
Show that $\operatorname{Extr}(K)\neq\emptyset$.

\noindent{\em Hint:}\
Use Zorn's lemma on the family of nonempty weak$^\ast$ compact faces of
$K$. A minimal such face must be a singleton: otherwise expose a proper
nonempty face by maximising a nonconstant weak$^\ast$ continuous affine
functional $x^\ast\mapsto \Re x^\ast(x)$.

\item
Prove that $K=\overline{\operatorname{co}}^{\,w^\ast}\!(\operatorname{Extr}(K))$,
the closure being taken in the weak$^\ast$ topology.

\noindent{\em Hint:}\
If $C:=\overline{\operatorname{co}}^{\,w^\ast}\!(\operatorname{Extr}(K))$
is a proper subset of $K$, use Proposition \ref{prop:HB-sep-weakstar}
to separate a point of $K\setminus C$ from $C$. Then apply part (a) to
the exposed face of $K$ on which the separating functional attains its
maximum.
\end{enumerate}

\item
Prove that if $Y$ is a closed subspace of a reflexive Banach space $X$, then the quotient space $X/Y$ is reflexive.

\item
Show that the weak topology of a weakly compact subset of a separable Banach space is metrisable.

\item
Using the result of the preceding problem, show that if $K$ is a weakly compact subset of a Banach space, then every sequence $(x_n)_{n\ge 1}$ contained in $K$ has a weakly convergent subsequence.

\item\label{prob:nonweakly}
Using the result of the preceding problem, show that $C[0,1]$ and $L^1(0,1)$ are nonreflexive by checking that their closed unit balls
contain sequences that have no weakly convergent subsequence.

\item\label{prob:Schur}
Show that the Banach space $\ell^1$ has the {\em Schur property}:\index{Schur!property} if $\limn\, x_n = x$ weakly in $\ell^1$\!, then $\limn\, x_n = x$ strongly in $\ell^1$\!.

\item
Let $(\Om,\calF\!,\mu)$ be a probability space. Let $(f_n)_{n\ge 1}$ be a bounded sequence in $L^1(\Om)$ which is {\em uniformly integrable},\index{uniformly!integrable} that is,
$$ \lim_{r\to\infty} \sup_{n\ge 1} \n\one_{\{|f_n|>r\}}f_n\n_1 = 0.$$
Show that $(f_n)_{n\ge 1}$ contains a weakly convergent subsequence by completing the details of the following argument.
\begin{enumerate}[\rm(a), leftmargin=*]
  \item For each $k\ge 1$, the sequence $(f_n^{(k)})_{n\ge 1}$ defined by
  $$f_n^{(k)}:=\one_{\{|f_n|\le k\}}f_n, \quad n\ge 1,$$
  contains a subsequence that is weakly convergent in $L^2(\Om)$, and hence weakly convergent in $L^1(\Om)$.

  \item After passing to a suitable diagonal subsequence, show that one may assume that for every $k\ge 1$ the sequence $(f_n^{(k)})_{n\ge1}$ converges weakly.

\end{enumerate}
Denote by $f^{(k)}$ their weak limits in $L^1(\Om)$.

\begin{enumerate}[\rm(a), leftmargin=*]
  \item[(c)] Show that $\n f^{(k)} - f^{(\ell)}\n_1 \le \liminf_{n\to\infty} \n f_n^{(k)} - f_n^{(\ell)}\n_1$ and the latter tends to $0$ as $k,\ell\to\infty$ by uniform integrability.

  \item[(d)] Conclude that the limit $\limk f^{(k)} = f$ exists in $L^1(\Om)$ and that
  $\limn f_n = f$ weakly in $L^1(\Om)$.
\end{enumerate}

\item\label{prob:SDP}
This problem illustrates how the Hahn--Banach separation theorem leads to
a basic duality principle in {\em semidefinite programming}.\index{semidefinite programming}
Let $\mathscr S^n$ be the real vector space of symmetric $n\times n$
matrices, equipped with the inner product
$$
\iprod{A}{B}_{\rm tr}:=\tr(AB).
$$
We write $A\ge0$ if $Ax\cdot x\ge 0$ for all $x\in \R^n$
and $A>0$ if $Ax\cdot x>0$ for all $0\ne x\in\R^n$.

\begin{enumerate}[\rm(a), leftmargin=*]
\item\label{it:SDP1} Show that the cone
$$
\mathscr S^n_+:=\{A\in\mathscr S^n:\, A\ge0\}
$$
is self-dual in the sense that
$A\ge 0$ if and only if $\tr(AB)\ge0$ for all $B\ge 0$.

\noindent{\em Hint:}\ For the `if' part, take
$B=x x^\mathsf{T}$.

  \item\label{it:SDP2} Let $X$ be a linear subspace of $\mathscr S^n$. Prove that exactly
  one of the following alternatives holds:
  \begin{enumerate}[leftmargin=*, label={\rm(\roman*)}]
  \item $X$ contains a positive definite matrix;
  \item there exists a nonzero matrix $B\ge0$ such that $\tr(AB)=0$ for all $A\in X$.
  \end{enumerate}
  \noindent {\em Hint:}\ If (i) fails, separate the open convex cone
  $\{A\in\mathscr S^n:\, A>0\}$ from $X$ by the Hahn--Banach separation
  theorem. Since $X$ is a linear subspace, the separating functional must
  vanish on $X$. Use part \ref{it:SDP1} to identify its sign.

\item\label{it:SDP3} More generally, let $A_0,A_1,\dots,A_m\in\mathscr S^n$. Show that
exactly one of the following holds:
\begin{enumerate}[leftmargin=*, label={\rm(\roman*)}]
\item there exist $x_1,\dots,x_m\in\R$ such that
$$A_0+\sum_{j=1}^m x_jA_j>0;$$
\item there exists a nonzero matrix $B\ge0$ such that
$$\tr(A_jB)=0 \ \hbox{ for all } j=1,\dots,m, \qquad \tr(A_0B)\le 0.$$
\end{enumerate}
\noindent{\em Hint:}\ Apply the Hahn--Banach separation theorem to the affine subspace
$A_0+\mathrm{span}\{A_1,\dots,A_m\}$ and the open cone of positive
definite matrices.
\end{enumerate}

\end{problems}

%% file: ch05-BoundedOperators.tex
\chapter{Bounded Operators}\label{ch:operators}

\blfootnote{\rm \small \begin{quote}This book has been published by Cambridge University Press in the series `Cambridge Studies in Advanced Mathematics'. The present version is free to view and download for personal use only. Not for re-distribution, re-sale or use in derivative works.\newline {\copyright} Jan van Neerven\end{quote}}

\noindent
The bounded operators introduced in the first chapter are simple to define, but the deeper general properties of bounded operators reflect the completeness of the spaces on which they act: pointwise control of a family of operators can lead to a uniform bound, while conditions on a single operator can have strong consequences for its inverse or its range. The Baire category theorem underlies several of these principles.

For operators given by explicit formulas, such as the Fourier transform and the Hilbert transform, boundedness may itself be difficult to establish. In such cases, interpolation theorems allow bounds known at certain exponents to yield bounds at intermediate ones.

\section{The Uniform Boundedness Theorem}\label{sec:UBT}

The proof of the uniform boundedness theorem, as well as the proofs of some other results in this chapter, depends on the Baire category theorem.

\subsection{The Baire Category Theorem}

\begin{theorem}[Baire category theorem]\label{thm:Baire}\index{theorem!Baire} Let $X$ be a nonempty complete metric space.
Let $F_1,F_2,\dots$ be closed subsets of $X$ such that
$$ X = \bigcup_{n\ge 1} F_n.$$
Then at least one of the sets $F_n$ has nonempty interior.
 \end{theorem}
 \begin{proof}
Assuming that all sets $F_n$ have empty interior, we prove the existence of an $x\in X$ not contained in any one of the $F_n$'s.

Pick an $x_1\in \complement F_1$. This is possible, for otherwise
we have $F_1 = X$ and $F_1$ contains open balls.
Since $F_1$ is closed, $\complement F_1$ is open and therefore
contains an open ball ${B}(x_1;r_1)$. By shrinking the radius a bit, we may even assume that the closed ball $\ov{B}(x_1;r_1)$ is contained in $\complement F_1$ and, moreover, that $0<r_1\le 1$. The ball $B(x_1;r_1)$
is not contained in $F_2$ and consequently the open set $B(x_1;r_1)\setminus F_2$
is nonempty. By the same reasoning as before, this set
contains a closed ball $\ov{B}(x_2;r_2)$ of radius $0<r_2\le \frac12$. Continuing in this way we obtain
a decreasing sequence of closed balls $\ov{B}(x_1;r_1) \supseteq \ov{B}(x_2;r_2) \supseteq \dots$
with $0<r_n\le \frac1n$. The sequence $(x_n)_{n\ge 1}$ is a Cauchy sequence, and therefore has a limit
$x$, by the completeness of $X$. It is clear that  $x\in \bigcap_{n\ge 1} \ov{B}(x_n;r_n)$,
and therefore $x\not\in \bigcup_{n\ge 1}F_n$.
\end{proof}

\subsection{The Uniform Boundedness Theorem}

The uniform boundedness theorem allows one to infer uniform boundedness of a family of bounded operators from their pointwise boundedness.

\begin{theorem}[Uniform boundedness theorem]\label{thm:UBT}\index{theorem!uniform boundedness} Let $(T_i)_{i\in I}$ be a family of bounded operators
from a Banach space $X$ into a normed space $Y$. If
$$ \sup_{i\in I} \n T_i x\n < \infty,\quad x\in X,$$
then
$$ \sup_{i\in I} \n T_i\n < \infty.$$
\end{theorem}
\begin{proof}
We may of course assume that the index set $I$ is nonempty.
For each $i\in I$ and integer $n\ge 1$ the set $\{x\in X:\, \n T_i x\n\le n\}$ is closed by the continuity of the operator $T_i$.
Since the intersection of closed sets is closed, for each $n\ge 1$
the set $$F_n := \bigl\{x\in X:\, \sup_{i\in I} \n T_i x\n\le n\bigr\} = \bigcap_{i\in I} \bigl\{x\in X:\, \n T_i x\n\le n\bigr\} $$
is closed. Moreover, the union equals $X$. By the Baire category theorem, at least one of them, say $F_{n_0}$, has nonempty interior.
Accordingly there exist
$x_0\in X$ and $r_0>0$ such that the open ball $B(x_0;r_0)$ is contained in $F_{n_0}$.

Fix an index $i\in I$. For any $x\in X$ with norm $\n x\n < r_0$ we write $x = x_0  - (x_0-x)$ and note that both $x_0$ and $x_0-x$ belong to $B(x_0;r_0)$. As a consequence,
$$\n T_i x\n \le \n T_i x_0\n + \n T_i (x_0-x)\n \le n_0+n_0 = 2n_0.$$
Hence, for all $x\in X$ with norm $\n x\n <1$,
$$\n T_i x\n  \le 2n_0/r_0.$$
This implies that $\n T_i\n \le 2n_0/r_0$.
Since this holds for all $i\in I$, we have shown that
$\sup_{i\in I} \n T_i\n \le 2n_0/r_0.$
\end{proof}

We continue with some typical applications.

 \begin{proposition}\label{prop:strong-limits}  Let $X$ be a Banach space, $Y$ be a normed space, and suppose that $T_n: X\to Y$, $n\ge 1$, are bounded operators such that $\lim_{n\to\infty} T_n x =: Tx$ exists for all $x\in X$. Then
 the
 operators $T_n$, $n\ge 1$, are uniformly bounded, the mapping $x\mapsto Tx$ is linear and bounded, and
  $$ \n T\n \le \liminf_{n\to\infty} \n T_n\n.$$
  \end{proposition}
 \begin{proof}
 The uniform boundedness theorem implies $\sup_{n} \n T_n\n <\infty$. The remaining assertions are proved by the argument of Proposition \ref{prop:approxTn}.
 \end{proof}

\begin{proposition}[Boundedness of bilinear mappings]
Let $X,Y,Z$ be normed spaces and suppose that at least one of the spaces $X$ and $Y$ is a Banach space.
Let $B:X\times Y\to Z$ be a bilinear mapping which is bounded in each variable separately.
Then there exists a constant $C\ge 0$ such that
$$
\n B(x,y)\n\le C\n x\n\n y\n,\quad x\in X,\ y\in Y.
$$
In particular, $B$ is jointly continuous.
\end{proposition}
\begin{proof}
Assume that $X$ is a Banach space (if $Y$ is a Banach space we interchange the roles of $X$ and $Y$).
For each $y\in Y$,
 $T_y x:= B(x,y)$ defines an element of $\calL(X,Z)$ since $B$ is bounded in its first
variable. Also, for each $x\in X$ we have $\sup_{\n y\n\le 1} \n T_y x\n <\infty$ since $B$ is bounded in its second
variable. Since $X$ is a Banach space, the uniform boundedness theorem shows that
$\{T_y:\, \n y\n\le 1\}$ is uniformly bounded in $\mathscr{L}(X,Z)$.
 With $M:= \sup_{\n y\n\le 1} \n T_y\n$ we then obtain, for all $y\in Y$ with $\n y\n\le 1$,
$$ \n B(x,y)\n = \n T_y x\n \le M \n x\n .$$
By a scaling argument for the second variable, this implies the claim as stated.
\end{proof}

The same proof works if we assume that $B$ is linear in the first variable, conjugate-linear in the second variable, and bounded in both variables separately; this observation will be useful in the context of Hilbert spaces.

The following proposition and its corollary give an application of the uniform boundedness theorem to duality.

\begin{proposition}[Weakly bounded sets are bounded]\label{prop:weakly-bounded}
 A subset $S$ of a
 normed space $X$ is bounded if and only if it is {\em weakly bounded},\index{weakly!bounded} that is, the set
 $\lb S,x^*\rb := \{\lb x,x^*\rb:\, x\in S\}$ is bounded for all $x^*\in X^*$\!.
\end{proposition}
\begin{proof}
Only the `if' part needs proof. Suppose that $S$ is weakly bounded. For each $x\in S$, the mapping $T_x: x\s\mapsto \lb x,x\s\rb$ is bounded, with $$\n T_x \n = \sup_{\n x\s\n\le 1}|\lb x,x\s\rb| =  \n x\n.$$ Since for each $x\s\in X\s$ we have $\sup_{x\in S} | T_x x\s| < \infty$, the uniform boundedness theorem (which can be applied since $X^*$ is a Banach space) implies that  $\sup_{x\in S} \n x\n = \sup_{x\in S} \n T_x\n<\infty$.
\end{proof}

\begin{corollary}\label{cor:weakly-bdd-seq} Let $X$ be a Banach space.
The following assertions hold:
\begin{enumerate}[label={\rm(\arabic*)}, leftmargin=*]
 \item\label{it:weakly-bdd-seq1} if $\limn \lb x,x_n^*\rb= \lb x,x\s\rb$ for all $x\in X$, then
$(x_n^*)_{n\ge 1}$ is bounded;
 \item\label{it:weakly-bdd-seq2} if $\limn \lb x_n,x^*\rb= \lb x,x\s\rb$ for all $x^*\in X^*$, then $(x_n)_{n\ge 1}$ is bounded.
\end{enumerate}
\end{corollary}
\begin{proof} Part \ref{it:weakly-bdd-seq1} follows directly from the uniform boundedness theorem
and part \ref{it:weakly-bdd-seq2} is a special case of Proposition \ref{prop:weakly-bounded}.
\end{proof}

Observe that the completeness of $X$ is only needed for part \ref{it:weakly-bdd-seq1}.

\section{The Open Mapping Theorem}\label{sec:OMT}

The next main theorem is the open mapping theorem. Among other things, it implies that a bijective
bounded operator between Banach spaces has a bounded inverse (and hence is an isomorphism). Its proof relies on the following lemma,
in which we use subscripts to distinguish between open balls in $X$ and $Y$.

\begin{lemma}\label{lem:OMT} Let $X$ be a Banach space, $Y$ a normed space, and let $T\in \calL(X,Y)$ be
a bounded operator. If $0< r,R<\infty$ are such that
$$B_Y(0;r) \subseteq \overline{T(B_X(0;R))},$$
then $$B_Y(0;r) \subseteq {T(B_X(0;2R))}.$$
\end{lemma}

The proof shows that the constant $2$ may be replaced by $1+\eps$ for any fixed $\eps>0$.
\begin{proof}
Fix $y_0\in B_Y(0;r)$. By the assumption, we can choose
$x_0\in B_X(0;R)$ such that $y_1:=y_0-Tx_0$ satisfies $
\n y_1\n<\frac12r$.
Since $y_1 \in \frac12 B_Y(0;r)$, we can choose
$x_1\in \frac12 B_X(0;R)$ such that $y_2:=y_1-Tx_1$ satisfies $
\n y_2\n<\frac14r$.
Continuing in this way, we obtain sequences $(x_n)_{n\ge 0}$ and
$(y_n)_{n\ge 0}$ with $$ y_{n+1} = y_{n}- Tx_n, \qquad \n x_n\n < 2^{-n}R, \qquad \n y_n\n < 2^{-n}r,$$
for all $n\ge 0$.
Since $\sum_{n\ge 0}\n x_n\n < \sum_{n\ge 0}2^{-n}R = 2R$,
the completeness of $X$ implies that the sum $x:=\sum_{n\ge1}x_n$ converges in $X$, and $\n x\n<2R$.
By construction,
$$
y_N = y_{N-1} - Tx_{N-1} = \hdots = y_0-T\sum_{n=0}^{N-1} x_n.
$$
Letting $N\to\infty$, we obtain
$$Tx = \limN T\sum_{n=0}^{N-1} x_n = \limN (y_0 - y_N) = y_0,$$
where we used that $\n y_N\n<2^{-N}r$.
Thus $y_0 = Tx \in T(B_X(0;2R))$.
\end{proof}

\begin{theorem}[Open mapping theorem]\label{thm:OM}\index{theorem!open mapping} Let $X$ and $Y$ be Banach spaces.
If $T\in\calL(X,Y)$ is bounded and surjective, then $T$ maps open sets to open sets.
\end{theorem}
\begin{proof}
Set $F_n := \overline{T(B_X(0;n))}$. The surjectivity of $T$ implies that $Y = \bigcup_{n\ge 1} F_n$.
Therefore, by the Baire category theorem (Theorem \ref{thm:Baire}), some $F_{n_0}$ has nonempty interior.
This means that there exist $y_0\in Y$ and $r_0>0$ such that
$$ B_Y(y_0;r_0) \subseteq  \overline{T(B_X(0;n_0))}.$$
In view of $T(-x) = -Tx$, we then also have
$$ B_Y(-y_0;r_0) \subseteq  \overline{T(B_X(0;n_0))}.$$
Writing $y = \frac12(y_0+y) + \frac12(-y_0+y)$, we find that
$$ B_Y(0;r_0) \subseteq  \tfrac12\cdot\overline{T(B_X(0;n_0))} + \tfrac12\cdot\overline{T(B_X(0;n_0))} = \overline{T(B_X(0;n_0))}.$$
Now we can invoke Lemma \ref{lem:OMT} and find
$$ B_Y(0;r_0) \subseteq T(B_X(0;2n_0)).$$

Let $U$ be an open set in $X$; we wish to prove that $T(U)$ is open. To this end let $Tx \in T(U)$ be given, with $x\in U$; we wish to
prove that $T(U)$ contains the open ball $B_Y(Tx;\rho)$ for some $\rho>0$.

Since $U$ is open, there is an $\e>0$ such that $B_X(x;\e)\subseteq U$. Let $\delta:= \e/(2n_0)$.
Then $T(U)$ contains $Tx + T(B_X(0;\e))  = Tx + \delta T(B_X(0;2n_0))$, and the latter contains
the open ball $Tx + \delta B_Y(0;r_0) = B_Y(Tx;\delta r_0)$.
\end{proof}

\begin{corollary}\label{cor:OMT} Let $X$ and $Y$ be Banach spaces, and let $T\in \calL(X,Y)$ be bounded. If $T$ is a bijection, then $T$ is an isomorphism from $X$ onto $Y$. More generally, if $\ran(T)$ is closed, then the quotient operator $T_/$ is an isomorphism from $X/\ker(T)$ onto $\Ran(T)$.
\end{corollary}
\begin{proof}
First assume that $T$ is a bijection. The fact that $T$ maps open sets to open sets can be reformulated as saying that $(T^{-1})^{-1}(U)$ is open for every open set $U$ in $X$, so $T^{-1}$ is continuous and hence bounded. This gives the first assertion. The second follows from it by noting that the quotient operator $T_/$ is a bijection from $X/\ker(T)$ onto $\Ran(T)$.
\end{proof}

We have already encountered an example of this situation in Chapter \ref{ch:duality}. If $Y$ is a closed subspace of a Banach space $X$, the Hahn--Banach extension theorem implies that the restriction mapping
$r:X\s\to Y\s$ is surjective. Since the null space of $r$ equals the annihilator $Y^\perp$, by
Corollary \ref{cor:OMT} $r$ induces an isomorphism $r_/$ of Banach spaces
$$ X\s/Y^\perp \simeq Y\s.$$
The Hahn--Banach extension theorem shows that $r_/$ is an isometry, and therefore this isomorphism is in fact isometric. This recovers the first part of Proposition \ref{prop:dualOfQuotient}.

It was noted in Proposition \ref{prop:proj-compl} that if a subspace $X_0$ of a normed space $X$ is the range of a projection on $X$, then $X_0$ is complemented in $X$.
As an application of Corollary \ref{cor:OMT} we prove the following converse:

\begin{proposition}\label{prop:compl-proj} A closed subspace of a Banach space $X$ is complemented if and only if it is the range of a projection on $X$.
\end{proposition}
\begin{proof}
It remains to prove the `only if' part.
Let $X_1$ be a closed subspace of $X$ such that
$X=X_0\oplus X_1$. Endow $X_0\times X_1$ with the norm
$$
\n(x_0,x_1)\n:=\n x_0\n+\n x_1\n.
$$
Since $X_0$ and $X_1$ are closed subspaces of the Banach space $X$,
this product is a Banach space.

The addition mapping
$$
J:X_0\times X_1\to X,\qquad J(x_0,x_1):=x_0+x_1,
$$
is linear and contractive by the triangle inequality.
The direct sum assumption means precisely that $J$ is a bijection.
By Corollary \ref{cor:OMT}, its inverse is bounded.

For $x=x_0+x_1$ along the direct sum decomposition, define
$\pi_0x:=x_0$. Then
$$
\n\pi_0x\n
=\n x_0\n
\le\n x_0\n+\n x_1\n
=\n J^{-1}x\n
\le\n J^{-1}\n\,\n x\n.
$$
Thus $\pi_0$ is a bounded projection onto $X_0$.
\end{proof}

\section{The Closed Graph Theorem}\label{sec:CGT}

Let $X$ and $Y$ be Banach spaces.
The {\em graph}\index{graph!of a linear operator} of a mapping $T: X\to Y$ is the set\index{$Gamma2$@$\Gr(T)$}
$$ \Gr(T):= \{(x,Tx): \, x\in X\}$$
in $X\times Y$.
If $T$ is linear, then $\Gr(T)$ is a linear subspace of $X\times Y$. Endowed
with the norm
\begin{align}\label{eq:prod-norm-1}\n (x,y)\n_1 := \n x\n+\n y\n,
\end{align}
the space $X\times Y$ is a Banach space, and it is easy to check that if $T$ is bounded, then $\Gr(T)$ is closed in $X\times Y$.
Since all product norms on $X\times Y$ are equivalent (see Example \ref{ex:prod-norm-eq}), the particular choice of product norm made in \eqref{eq:prod-norm-1} is immaterial.

\begin{definition}[Closed operators]
 A linear operator $T: X\to Y$ is {\em closed}\index{closed!linear operator} if its graph is closed in $X\times Y$.
\end{definition}

Every bounded linear operator is closed. In the converse direction we have the following result.

\begin{theorem}[Closed graph theorem]\label{thm:CGT}\index{theorem!closed graph} Let $X$ and $Y$ be Banach spaces.
 If a linear operator $T: X\to Y$ is closed, then $T$ is bounded.
\end{theorem}
\begin{proof}
 Let $\pi_X: X\times Y\to X$ and $\pi_{\,Y}:X\times Y\to Y$ be given by
 $$\pi_X(x,y):= x,\qquad \pi_{\,Y}(x,y):= y.$$ Both mappings are bounded and of norm at most one.
 By assumption $\Gr(T)$ is a closed subspace of $X\times Y$, hence a Banach space with respect to the
 inherited norm.
 Consider the linear operator $S: X\to \Gr(T)$ given by $Sx := (x,Tx)$. This operator is a bijection
 whose inverse $S^{-1}$ is the bounded operator $\pi_X|_{\Gr(T)}$. By Corollary \ref{cor:OMT} the inverse $S$ of
 $S^{-1}$ is bounded. Hence also $T = \pi_{\,Y}|_{\Gr(T)} \circ S$ is bounded.
\end{proof}

\begin{figure}[ht]
 \begin{center}
\begin{tikzcd} [scale cd=1.2, sep=huge]
& \Gr(T)\arrow[dr,  "\pi_{\,Y}|_{\Gr(T)}"] & \\
X  \arrow[ur,  "S"] \arrow[rr,  "T"] &
 & Y
\end{tikzcd}
\caption{Proof of the closed graph theorem}
\end{center}
\end{figure}

As an application of the closed graph theorem we prove the following variation of Proposition \ref{prop:Lp-via-Lq-1}.

 \begin{proposition}\label{prop:Lp-via-Lq-2}\index{inequality!H\"older, converse to}
  Let $1\le p, q\le \infty$ satisfy $\frac1p+\frac1q =1$. Let $(\Omega,\mu)$ be
  a $\sigma$-finite measure space.
  A measurable function $f$ belongs to $L^p(\Om)$ if and only if $fg\in L^1(\Om)$ for all
  $g\in L^q(\Om)$. In that case we have $$\n f\n_p = \sup_{\n g\n_q\le 1} \int_\Omega |fg|\ud \mu.$$
 \end{proposition}
\begin{proof}
 The `only if' part is immediate from H\"older's inequality. To prove the `if' part we may assume that $f$ is not identically $0$. Using Corollary \ref{cor:Lp-ae-subseq}, we see that the mapping $g\mapsto fg$ from
$L^q(\Om)$ to $L^1(\Om)$ is closed: if $g_n\to g$ in $L^q(\Om)$ and $fg_n\to h$ in $L^1(\Om)$, we may pass to a subsequence such that $g_{n_k}\to g$ and $fg_{n_k}\to h$ $\mu$-almost everywhere, and therefore
$fg=\lim_{k\to\infty}fg_{n_k}=h$ $\mu$-almost everywhere.
By the closed graph theorem, the operator $g\mapsto fg$ is bounded.
It follows that the assumptions of  Proposition \ref{prop:Lp-via-Lq-1} are satisfied, with $M$ the norm of
the operator $g\mapsto fg$. This shows that $f\in L^p(\Om)$, with bound
$$\n f\n_p \le \sup_{\n g\n_q\le 1} \int_\Omega |fg|\ud \mu.$$
H\"older's inequality gives the opposite bound.
\end{proof}

As a further illustration of the use of the closed graph theorem,
we give a second proof of Proposition \ref{prop:compl-proj}.
Let $X=X_0\oplus X_1$ be a direct sum decomposition into closed
subspaces, and let $\pi_0:X\to X$ be the linear mapping given by
$\pi_0(x_0+x_1):=x_0$.

To prove that $\pi_0$ has closed graph, suppose that $x_n\to x$
and $\pi_0x_n\to y$ in $X$. Since $\pi_0x_n\in X_0$ and $X_0$
is closed, we have $y\in X_0$. Similarly,
$x_n-\pi_0x_n\to x-y$, and the closedness of $X_1$ gives $x-y\in X_1$.
Thus $x=y+(x-y)$ is the direct sum decomposition of $x$, and
therefore $\pi_0x=y$. This proves that $\pi_0$ is closed.
The closed graph theorem now implies that $\pi_0$ is bounded.

\section{The Closed Range Theorem}\label{sec:CRT}

Recall that the {\em annihilator} of a subset $A$ of $X$ is the set
$$A^\perp:= \{ x^*\in X^*:\, \lb x,x\s\rb = 0 \hbox{ for all } x\in A\}$$
and the {\em pre-annihilator} of a subset $B$ of $X^*$ is the set
$${}^\perp B:= \{ x\in X:\, \lb x,x\s\rb = 0 \hbox{ for all } x\s\in B\}.$$
It was shown in Proposition
\ref{prop:dual-dense-T} that for any operator $T\in \calL(X,Y)$, where $X$ and $Y$ are Banach spaces, we have
\begin{align}\label{eq:HB-denserange}\ov{\Ran(T)} = {}^\perp(\Ker(T^*)) .
\end{align}
In Sections \ref{sec:spectrum-compact} and \ref{sec:Fredholm-theory} we will encounter interesting classes of
operators whose ranges are closed. For such operators, the closed range theorem provides a `dual' variant of \eqref{eq:HB-denserange}
which is considerably harder to prove. We need this theorem in
our discussion of duality of Fredholm operators in Section \ref{sec:Fredholm-theory}.

\begin{theorem}[Closed range theorem]\label{thm:closed-range-dual}\index{theorem!closed range} Let $T\in\calL(X,Y)$, where $X$ and $Y$ are Banach spaces. If $\ran(T)$ is closed, then
$$ \Ran(T^*) = (\Ker(T))^\perp\!.$$
As a consequence, $\ran(T^*)$ is weak$^*$ closed.
\end{theorem}

\begin{proof}
Once we have proved the identity $\Ran(T^*)=(\Ker(T))^\perp$, the weak$\s$ closedness of $\ran(T^*)$ follows from the general observation that annihilators are weak$\s$ closed.

\smallskip
$\subseteq$: \ If $x^*=T^*y^*\in\ran(T^*)$, then for all $x\in\ker(T)$ we have $\lb x,x^*\rb=\lb Tx,y^*\rb=0$. This shows that $x^*\in(\Ker(T))^\perp$.

\smallskip
$\supseteq$: \ Suppose that $x^*\in(\Ker(T))^\perp$.
By Proposition \ref{prop:dualOfQuotient}, $x^*$ induces a bounded
functional $z^*\in(X/\ker(T))^*$ given by
$$
\lb x+\ker(T),z^*\rb:=\lb x,x^*\rb.
$$
Since $\ran(T)$ is closed, the quotient operator
$T_/:X/\ker(T)\to\ran(T)$ is an isomorphism by
Corollary \ref{cor:OMT}.
Define $\phi\in(\ran(T))^*$ by
$\phi:=z^*\circ T_/^{-1}$.
Then, for all $x\in X$,
$$
\phi(Tx)
=
\lb x+\ker(T),z^*\rb
=
\lb x,x^*\rb.
$$
By the Hahn--Banach theorem, $\phi$ extends to a functional $y^*\in Y^*$. It follows that, for all $x\in X$,
$$
\lb x,T^*y^*\rb=\lb Tx,y^*\rb=\phi(Tx)=\lb x,x^*\rb .
$$
Hence $x^*=T^*y^*$, and therefore $x^*\in\ran(T^*)$.
\end{proof}

\section{The Fourier Transform}\label{sec:FT}

In the present section and the next we study two nontrivial examples of bounded operators: the Fourier--Plancherel transform and the Hilbert transform.
It is not an exaggeration to state that, at least from the point
of view of the theory of partial differential equations, these rank among the most important bounded operators in all of analysis.

In the remainder of this chapter, unless otherwise stated, we will work over the complex scalar field.

\subsection{Definition and General Properties}

\begin{definition}[Fourier transform]
 The {\em Fourier transform} of a function $f\in L^1(\R^d)$ is the function $\wh f:\R^d\to \C$ defined by\index{$FB$@$\widehat{f}$}
 \begin{equation}\label{eq:FT} \wh f(\xi) := \frac1{(2\pi)^{d/2}} \int_{\R^d} f(x)\exp(-i x\cdot \xi)\ud x, \quad \xi\in \R^d\!,
 \end{equation}
where $x\cdot\xi:= \sum_{j=1}^d x_j\xi_j$.
\end{definition}

It is evident that $\wh f\in L^\infty(\R^d)$ and
$$\n \wh f\n_\infty\le (2\pi)^{-d/2}\n f\n_1.$$
This shows that the operator $\calF: f \mapsto \wh f$,
which will be referred to as the {\em Fourier transform},\index{Fourier!transform}\index{$FA$@$\calF$}
defines a bounded operator from $L^1(\R^d)$ to $L^\infty(\R^d)$.

\begin{remark}\label{rem:FT-normalised}
 In certain situations it is useful to absorb the constant $(2\pi)^{-d/2}$ into the measure. Denoting the resulting {\em normalised Lebesgue measure}\index{Lebesgue!measure, normalised}
by $${\rm d}m(x) = (2\pi)^{-d/2}\ud x,$$
we can interpret the Fourier transform as the operator from $L^1(\R^d\!,m)$ to $L^\infty(\R^d\!,m)$ given by
 \begin{equation*} \wh f(\xi) := \int_{\R^d} f(x)\exp(-i x\cdot \xi)\ud m(x), \quad \xi\in \R^d\!.
 \end{equation*}
The advantage of this point of view is that this operator is contractive. In many applications, however, working with the normalised Lebesgue measure is somewhat artificial, and for this reason we stick with \eqref{eq:FT} most of the time.
\end{remark}

The dominated convergence theorem implies that for all $f \in L^1(\R^d)$
 the function $\wh f$ is sequentially continuous, hence continuous.
More is true: the following lemma shows that $\wh f$ belongs to $C_0(\R^d)$, the space of continuous
functions vanishing at infinity.

\begin{theorem}[Riemann--Lebesgue lemma]\index{lemma!Riemann--Lebesgue}
 For all $f \in L^1(\R^d)$ we have $\wh f \in C_0(\R^d)$.
\end{theorem}

\begin{proof}
It has already been observed that $\wh f$ is continuous for all $f\in L^1(\R^d)$. It remains to prove that $\wh f$ vanishes at infinity.

First let $f\in C_{\rm c}^\infty(\R^d)$. For $j=1,\ldots,d$, integration by parts gives
$$\widehat{\partial_j f}(\xi)=i\xi_j\wh f(\xi).$$ Hence, if $\xi_j\not=0$,
$$
|\wh f(\xi)|\leq \frac{1}{(2\pi)^{d/2}}\frac{\n \partial_j f\n_1}{|\xi_j|}.
$$
For $\xi\not=0$, choose $j$ such that $|\xi_j|\geq |\xi|/\sqrt d$. Then
$$
|\wh f(\xi)|\leq \frac{\sqrt d}{(2\pi)^{d/2}}\Bigl(\max_{1\leq j\leq d}\n \partial_j f\n_1\Bigr) \cdot \frac1{|\xi|}.
$$
It follows that $\wh f$ vanishes at infinity.

Now let $f\in L^1(\R^d)$. Using Proposition \ref{prop:Cc-dense}, choose functions $f_n\in C_{\rm c}^\infty(\R^d)$ such that $f_n\to f$ in $L^1(\R^d)$. Since the Fourier transform is bounded from $L^1(\R^d)$ to $L^\infty(\R^d)$, we have $\wh{f_n}\to\wh f$ uniformly. Since $C_0(\R^d)$ is closed in $L^\infty(\R^d)$ and $\wh{f_n}\in C_0(\R^d)$ for all $n\geq 1$, it follows that $\wh f\in C_0(\R^d)$.
\end{proof}

We continue with an inversion theorem for the Fourier transform.
Its proof is based on a simple lemma.

\Needspace{6\baselineskip}
\begin{lemma}\label{lem:Gauss}
For $\la>0$ the Fourier transform of the function
$$g^{(\la)}(x)=\frac1{(2\pi)^{d/2}}\exp\Bigl(-\frac12\la|x|^2\Bigr)$$ equals $$\wh{g^{(\la)}}(\xi) =  \frac1{(2\pi \la)^{d/2}}\exp\Bigl(-\frac12 |\xi|^2/\la\Bigr).$$
\end{lemma}
\begin{proof}
First let $d=1$. Completing the square and using Cauchy's theorem
to shift the path of integration, we find
\begin{align*}
\ &  \frac1{(2\pi)^{1/2}}\int_{\R} g^{(\la)}(x)\exp(-i x\cdot\xi )\ud x
\\ & \qquad =  \frac1{2\pi}\int_{\R}\exp\Bigl(-\frac12\la[(x + i\xi/\la)^2 + \xi^2/\la^2]\Bigr)\ud x
\\ & \qquad = \frac1{2\pi}\exp\Bigl(-\frac12\xi^2/\la\Bigr) \int_{\R+i\xi/\la}\exp\Bigl(-\frac12\la z^2\Bigr)\ud z
\\ & \qquad = \frac1{2\pi}\exp\Bigl(-\frac12\xi^2/\la\Bigr)\int_{\R} \exp\Bigl(-\frac12\la z^2\Bigr)\ud z
 = \frac1{(2\pi \la)^{1/2}}\exp\Bigl(-\frac12\xi^2/\la\Bigr).\end{align*}
The general case follows from this by separation of variables.
\end{proof}

A different proof based on the Picard--Lindel\"of theorem is outlined in Problem \ref{prob:Gauss}.

\begin{theorem}[Fourier inversion theorem]\label{thm:FT-inversion}\index{theorem!inversion of the Fourier transform}
If $f\in L^1(\R^d)$ satisfies $\wh f\in L^1(\R^d)$, then for
almost all $x\in \R^d$ we have the identity
$$ f(x) = \frac1{(2\pi)^{d/2}}\int_{\R^d} \wh f(\xi) \exp(i x\cdot \xi)\ud \xi.$$
\end{theorem}

In particular, this result implies that the Fourier transform is injective as a mapping from $L^1(\R^d)$ to $L^\infty(\R^d)$.
A more general injectivity result will be proved in Theorem \ref{thm:uniq-FT-Rd}.

\begin{proof}
Let
$$g(x):=(2\pi)^{-d/2}\exp(-\tfrac12|x|^2),$$ and for $\la>0$ set
$g_\la(x):=\la^{-d}g(\la^{-1}x)$.
Then $\int_{\R^d}g(x)\ud x=1$.
By Lemma \ref{lem:Gauss} and a change of variables,
$$
g_\la(y)
=
\frac1{(2\pi)^{d/2}}
\int_{\R^d}g(\la\xi)\exp(iy\cdot\xi)\ud\xi.
$$
Consequently, Fubini's theorem gives, for every $x\in\R^d$,
\begin{align*}
(g_\la*f)(x)
&=
\frac1{(2\pi)^{d/2}}
\int_{\R^d}g(\la\xi)
\Bigl(\int_{\R^d}f(y)\exp(i(x-y)\cdot\xi)\ud y\Bigr)\ud\xi
\\*
&=
\int_{\R^d}g(\la\xi)\wh f(\xi)\exp(ix\cdot\xi)\ud\xi.
\end{align*}
The interchange of integrals is justified by the integrability
of $f$ and of $\xi\mapsto g(\la\xi)$.

Since $\wh f$ is integrable, $g$ is bounded, and
$g(\la\xi)\to g(0)=(2\pi)^{-d/2}$, dominated convergence gives
$$
\lim_{\la\downarrow0}(g_\la*f)(x)
=
\frac1{(2\pi)^{d/2}}
\int_{\R^d}\wh f(\xi)\exp(ix\cdot\xi)\ud\xi,
\quad x\in\R^d.
$$
On the other hand, Theorem \ref{thm:approx-identity} gives
$g_\la*f\to f$ in $L^1(\R^d)$.
By Corollary \ref{cor:Lp-ae-subseq}, along a sequence
$\la_n\downarrow0$ this convergence also holds almost everywhere.
Comparing the two limits proves the theorem.
\end{proof}

The Fourier transform of the translate $\tau_h f$ of a function $f\in L^1(\R)$ is given by
$$ \wh{\tau_h f}(\xi)
= \frac1{\sqrt{2\pi}} \int_{-\infty}^\infty f(x+h)\exp(-ix\xi)\ud x = \exp(ih\xi)\wh f(\xi).$$
 It follows that if $\wh f(\xi_0)=0$ for some $\xi_0\in\R$, then $\wh{\tau_h f}(\xi_0)=0$ for all $h\in\R$. Therefore the linear span of the set of translates of $f$ is contained in $\{g\in L^1(\R): \, \wh g(\xi_0) = 0\}$, which is a proper
closed subspace  of $L^1(\R)$. Thus, a necessary condition in order that the linear span of the set of translates of a function $f\in L^1(\R)$ be dense in $L^1(\R)$ is that $\wh f$ be zero-free.
Strikingly, this necessary condition is also sufficient. This is the content of the next theorem, which will be proved by operator-theoretic methods in Section \ref{subsec:C0groups}.

\begin{theorem}[Wiener's Tauberian theorem]\label{thm:Wiener}\index{theorem!Wiener's Tauberian}
If the Fourier transform
of a function $f\in L^1(\R)$ is zero-free, then
the span of the set of all translates of $f$
is dense in $L^1(\R)$.
\end{theorem}

\subsection{The Plancherel Theorem}\label{subsec:Plancherel}

The Fourier transform enjoys an important $L^2$-boundedness property.

\begin{theorem}[Plancherel, preliminary version]\label{thm:Plancherel}
If $f\in L^1(\R^d)\cap L^2(\R^d)$, then $\wh f \in L^2(\R^d)$
 and $$\n \wh f\n_2 = \n f\n_2.$$
\end{theorem}
\begin{proof}
Since $f\in L^1(\R^d)$,  $\wh f$ is bounded and $\xi\mapsto \exp(-\frac12\la|\xi|^2)|\wh f(\xi)|^2$ is integrable
for all $\la>0$ and, with the notation of Lemma \ref{lem:Gauss},
\begin{align*}
&  \int_{\R^d}\exp\Bigl(-\frac12\la|\xi|^2\Bigr)|\wh f(\xi)|^2\ud \xi
  \\ & \qquad =\int_{\R^d}\exp\Bigl(-\frac12\la|\xi|^2\Bigr)\wh f(\xi)\ov{\wh f(\xi)}\ud \xi
   \\ & \qquad = \frac1{(2\pi)^{d/2}}\int_{\R^d}g^{(\la)}(\xi)
\int_{\R^d}f(x)\exp(- ix\cdot\xi)\ud x \int_{\R^d}\ov{f(y)}\exp( iy\cdot\xi)\ud y\ud \xi
   \\ & \qquad = \int_{\R^d}\int_{\R^d} \Bigl(\frac1{(2\pi)^{d/2}}\int_{\R^d}  g^{(\la)}(\xi)\exp(-i(x-y)\cdot\xi)\ud \xi\Bigr) f(x)\ov {f(y)}\ud x\ud y
   \\ & \qquad = \int_{\R^d}\int_{\R^d} \frac1{(2\pi\la)^{d/2}}\exp\Bigl(-\frac12 |x-y|^2/\la\Bigr) f(x)\ov {f(y)}\ud x\ud y
   \\ & \qquad = \int_{\R^d} f*\phi_{\sqrt{\la}}(y)\ov {f(y)}\ud y,
 \end{align*}
where $\phi_\mu(x) := \mu^{-d} \phi(\mu^{-1}x)$ with $\phi(y) ={(2\pi)^{-d/2}}\exp(-\frac12|y|^2)$; the change of order of integration is justified by the absolute integrability of the integrand.
Applying Theorem \ref{thm:approx-identity} we find that
$\lim_{\la\downarrow 0} f*\phi_{\sqrt{\la}} = f$ in $L^2(\R^d)$.
Then,
\begin{align*}
 \lim_{\la\downarrow 0} \int_{\R^d} f*\phi_{\sqrt{\la}}(y)\ov {f(y)}\ud y
 = \int_{\R^d} f(y)\ov {f(y)}\ud y = \n f\n_2^2.
\end{align*}
On the other hand,
$$ \lim_{\la\downarrow 0}  \int_{\R^d}  \exp\Bigl(-\frac12\la|\xi|^2\Bigr)|\wh f(\xi)|^2\ud \xi
= \int_{\R^d}|\wh f(\xi)|^2\ud \xi =\n \wh f\n_2^2
$$
by the monotone convergence theorem. This completes the proof.
\end{proof}

Consider the vector space
$$\calF^2(\R^d):= \big\{f\in  L^1(\R^d)\cap L^2(\R^d):\,  \wh f\in L^1(\R^d)\cap L^2(\R^d)\big\}.$$
There is some redundancy in this definition, for if $f\in  L^1(\R^d)\cap L^2(\R^d)$, then $\wh f\in L^2(\R^d)$ by the Plancherel theorem. The advantage of the
above format is that it brings out the symmetry between $f$ and $\wh f$ explicitly. The interest of this space is explained by the following two observations.

\begin{lemma}\label{lem:F2} The Fourier transform maps $\calF^2(\R^d)$ bijectively onto itself.
\end{lemma}

\begin{proof}
If $f\in\calF^2(\R^d)$, Fourier inversion gives
$\widehat{\widehat f\phantom{.}}\!(x)=f(-x)$ almost everywhere; hence $\widehat{\widehat f\phantom{.}}\!\in L^1(\R^d)\cap L^2(\R^d)$, and so
$\wh f\in\calF^2(\R^d)$. Injectivity follows from the Plancherel theorem.
For surjectivity, let $f\in\calF^2(\R^d)$ and define
$g(\xi):=\widehat f(-\xi)$.
Then $g\in L^1(\R^d)\cap L^2(\R^d)$, and $\widehat g=f$ by the Fourier inversion theorem.
Moreover, since $\widehat g=f\in L^1(\R^d)\cap L^2(\R^d)$, we have
$g\in\calF^2(\R^d)$. Hence every $f\in\calF^2(\R^d)$ is the Fourier
transform of some $g\in\calF^2(\R^d)$.
\end{proof}

\begin{lemma}\label{lem:F2-dense} $\calF^2(\R^d)$ is dense in $L^2(\R^d)$.
\end{lemma}

\begin{proof}
Since $C_{\rm c}^\infty(\R^d)$ is dense in $L^2(\R^d)$ by Proposition \ref{prop:Cc-dense},
it suffices to show that every $f\in C_{\rm c}^\infty(\R^d)$ belongs to $\calF^2(\R^d)$. Integrating by parts, for all $f\in C_{\rm c}^\infty(\R^d)$ and multi-indices $\al\in\N^d$ we have
 $$\wh{\partial^\alpha f}(\xi) = i^{|\alpha|} \xi^\alpha\wh f(\xi),$$
 where
 $\partial^{\alpha}  := \partial_1^{\alpha_1}\circ\cdots\circ \partial_d^{\alpha_d}$ with $\partial_j$ the partial derivative in the $j$th direction,
 $\xi^{\alpha}  := \xi_1^{\alpha_1}\cdots\xi_d^{\alpha_d}$,
 and
 $|\al| := \al_1+\cdots+\al_d$.
 Since Fourier transforms of integrable functions are bounded, this implies that $\xi\mapsto (1+|\xi|^k)\wh f(\xi)$
is bounded for every integer $k\ge 1$. The desired result follows from this.
\end{proof}

Combining these lemmas with Proposition \ref{prop:extendT}, we obtain the following improved version of Theorem \ref{thm:Plancherel}.

\begin{theorem}[Plancherel]\label{thm:FT-main}\index{theorem!Plancherel} The restriction of the Fourier transform to $L^1(\R^d)\cap L^2(\R^d)$
 has a unique extension to an isometry from $L^2(\R^d)$ onto itself.
\end{theorem}

\begin{definition}[Fourier--Plancherel transform]
This isometry of $L^2(\R^d)$ is called
the {\em Fourier--Plancherel transform}.\index{Fourier--Plancherel transform}
\end{definition}

With a slight abuse of notation we denote the Fourier--Plancherel transform again by $\calF: f \mapsto \wh f$.
It is important to realise that $\wh f$ is no longer given by the pointwise formula
\eqref{eq:FT}. In fact, for functions $f\in L^2(\R^d)$ the integrand in \eqref{eq:FT} is not even integrable
unless $f\in L^1(\R^d)\cap L^2(\R^d)$.

\begin{remark}\label{rem:FT-normalised-L2} Theorem \ref{thm:FT-main} also holds with respect to the normalised Le\-besgue measure ${\rm d}m(x) = (2\pi)^{-d/2}\ud x$: the restriction to $L^1(\R^d\!,m)\cap L^2(\R^d\!,m)$ of the Fourier transform
as defined in Remark \ref{rem:FT-normalised} extends to an isometry from $L^2(\R^d\!,m)$ onto itself.
\end{remark}

For later use we record two further properties of the Fourier--Plancherel transform.

\begin{proposition}\label{prop:FT-fFg} For all $f,g\in L^2(\R^d)$ we have
 $$ \int_{\R^d} f(x) \wh{g}(x)\ud x = \int_{\R^d} \wh f(x) g(x)\ud x. $$
\end{proposition}
\begin{proof}
 For $f,g\in \calF^2(\R^d)$
  the identity follows by writing out the Fourier transforms and using Fubini's theorem:
  \begin{align*}
    \int_{\R^d} f(x) \wh{g}(x)\ud x
   & =   \int_{\R^d} \frac1{(2\pi)^{d/2}} \int_{\R^d}  f(x) g(\xi)\exp(-i x\cdot\xi)\ud \xi\ud x
   \\ & =  \int_{\R^d} \frac1{(2\pi)^{d/2}} \int_{\R^d}  f(x) g(\xi)\exp(-i x\cdot\xi)\ud x\ud \xi
    = \int_{\R^d} \wh f(\xi) g(\xi)\ud \xi.
  \end{align*}
 The general case follows by approximation, using the fact that $\calF^2(\R^d)$ is dense in $L^2(\R^d)$ by
 Lemma \ref{lem:F2-dense}.
\end{proof}

The appearance of the factor $(2\pi)^{d/2}$ in the next proposition is an artefact of our convention of normalising the Fourier transform with this factor, but not the convolution.

\begin{proposition}\label{prop:FT-convol}
Let $f\in L^1(\R^d)$ and let $g\in L^1(\R^d)$ or $g\in L^2(\R^d)$. For almost all $\xi\in \R^d$ we have
$$ \wh{f*g}(\xi) = (2\pi)^{d/2} \wh f(\xi) \wh g(\xi).$$
\end{proposition}

\begin{proof}
If $g\in C_{\rm c}(\R^d)$, then $f*g\in L^1(\R^d)\cap L^2(\R^d)$ by Young's inequality, and by Fubini's theorem and a change of variables
we obtain
\begin{align*}
 \wh{f*g}(\xi)&  = \frac1{(2\pi)^{d/2}}\int_{\R^d}\Bigl( \int_{\R^d}f(x-y)g(y)\ud y\Bigr) \exp(- i x\cdot\xi)\ud x
\\ &  =\int_{\R^d}\Bigl(\frac1{(2\pi)^{d/2}} \int_{\R^d}f(x-y)\exp(-i (x-y)\cdot\xi)\ud x\Bigr) \exp(- i y\cdot\xi)g(y)\ud y
\\ &  = \int_{\R^d}\Bigl(\frac1{(2\pi)^{d/2}} \int_{\R^d} f(u)\exp(-i u\cdot\xi)\ud u\Bigr) g(y)\exp(- i y\cdot\xi)\ud y
\\ &  = (2\pi)^{d/2} \wh f(\xi)\wh g (\xi).
\end{align*}
This proves the identity for $g\in C_{\rm c}(\R^d)$.
For $p\in\{1,2\}$ and $\frac1p+\frac1q=1$, Young's inequality, together with the $L^1$--$L^\infty$ bound for the Fourier transform when $p=1$ and the Plancherel theorem when $p=2$, shows that
$$
g\mapsto\wh{f*g}
$$
is bounded from $L^p(\R^d)$ to $L^q(\R^d)$. The same is true of
$$
g\mapsto(2\pi)^{d/2}\wh f\,\wh g.
$$
Since $C_{\rm c}(\R^d)$ is dense in $L^p(\R^d)$ by Proposition \ref{prop:Cc-dense}, the identity extends to all $g\in L^p(\R^d)$.
\end{proof}

From Proposition \ref{prop:ansmu} we know that if $\mu$ is a real or complex measure,
its variation $|\mu|$ is a finite measure. Accordingly, the Lebesgue integrals of bounded
Borel functions with respect to $\mu$ are well defined.
In particular we can define the {\em Fourier--Stieltjes transform}\index{Fourier--Stieltjes transform}
of a real or complex Borel measure $\mu$ on $\R^d$  by
$$ \wh \mu(\xi):= \frac1{(2\pi)^{d/2}} \int_{\R^d} \exp(-i x\cdot\xi)\ud\mu(x), \quad \xi\in \R^d\!.$$
From
$$  |\wh \mu(\xi)|\le  \frac1{(2\pi)^{d/2}}\int_{\R^d}  \ud|\mu|  = \frac1{(2\pi)^{d/2}}\n \mu\n,$$
where $\n \mu\n= |\mu|(\R^d)$ is the variation norm of $\mu$, we see that $\wh \mu$ is a bounded function, and it is continuous by dominated convergence. Thus we have $\wh\mu\in C_{\rm b}(\R^d)$ and
$\n \wh\mu\n_\infty \le (2\pi)^{-d/2} \n \mu\n$. The Riemann--Lebesgue lemma does not extend to the present setting, as is
demonstrated by the identity $\wh{\delta_0} = (2\pi)^{-d/2}\one$.

The Fourier inversion theorem (Theorem \ref{thm:FT-inversion}) implies that the Fourier transform is injective as an operator from $L^1(\R^d)$ to $L^\infty(\R^d)$.
More generally we have the following result.

\begin{theorem}[Injectivity of the Fourier--Stieltjes transform]\label{thm:uniq-FT-Rd}
 \index{theorem!injectivity, Fourier--Stieltjes transform}
If $\mu$ is a real or complex Borel measure on $\R^d$ satisfying $\widehat{\mu}(\xi) = 0$ for all $\xi\in \R^d$\!, then $\mu = 0$.
\end{theorem}

\begin{proof}
By the uniqueness part of the Riesz representation theorem
(Theorem \ref{thm:CK-dual}), it suffices to prove that
\begin{align}\label{eq:FT-mu-vanish}
\int_{\R^d}f\ud\mu=0,\quad f\in C_{\rm c}(\R^d).
\end{align}
Fix $f\in C_{\rm c}(\R^d)$, and put
$$
g(x):=(2\pi)^{-d/2}\exp(-\tfrac12|x|^2).
$$
With
$g_\eps(x):=\eps^{-d}g(\eps^{-1}x)$,
by Lemma \ref{lem:Gauss} and Fubini's theorem, for every
$x\in\R^d$ we have
\begin{align*}
(g_\eps*f)(x)
&=
\frac1{(2\pi)^{d/2}}
\int_{\R^d}g(\eps\xi)
\Bigl(\int_{\R^d}f(y)\exp(-iy\cdot\xi)\ud y\Bigr)
\exp(ix\cdot\xi)\ud\xi
\\
&=
\int_{\R^d}g(\eps\xi)\wh f(\xi)\exp(ix\cdot\xi)\ud\xi.
\end{align*}
Integrating this identity with respect to $\mu$ and applying
Fubini's theorem once more gives
\begin{align*}
\int_{\R^d}(g_\eps*f)(x)\ud\mu(x)
&=
\int_{\R^d}g(\eps\xi)\wh f(\xi)
\Bigl(\int_{\R^d}\exp(ix\cdot\xi)\ud\mu(x)\Bigr)\ud\xi
\\
&=
(2\pi)^{d/2}\int_{\R^d}
g(\eps\xi)\wh f(\xi)\wh\mu(-\xi)\ud\xi
=
0.
\end{align*}
Here the interchange of integrals is justified by the boundedness
of $\wh f$, the integrability of $\xi\mapsto g(\eps\xi)$, and the
finiteness of $|\mu|$.

Since $\int_{\R^d}g(y)\ud y=1$, we have
$$
\n g_\eps*f-f\n_\infty
\le
\int_{\R^d}g(y)\n f(\cdot-\eps y)-f\n_\infty\ud y
\to0
$$
by the uniform continuity of $f$ and dominated convergence.
Consequently,
$$
\int_{\R^d}f\ud\mu
=
\lim_{\eps\downarrow0}
\int_{\R^d}(g_\eps*f)\ud\mu
=0.
$$
This proves \eqref{eq:FT-mu-vanish}.
\end{proof}

For later reference we also record the corresponding result on the
unit circle. Its proof is a direct application of the
Stone--Weierstrass theorem (see Problem \ref{prob:uniq-FT-T}).
 The {\em Fourier--Stieltjes coefficients}
of a real or complex Borel measure $\mu$ on the unit circle $\mathbb{T}$ are defined by
$$\wh{\mu}(n) :=  \frac1{2\pi}\int_{-\pi}^\pi \exp(-i n\theta)\ud\mu(\theta), \quad n\in\Z.$$

\begin{theorem}[Injectivity of the Fourier--Stieltjes transform on the circle]\label{thm:uniq-FT-T}
 \index{theorem!injectivity, Fourier--Stieltjes transform}
If $\mu$ is a real or complex Borel measure on $\mathbb{T}$ satisfying $\widehat{\mu}(n) = 0$ for all $n\in \Z$, then $\mu = 0$.
\end{theorem}

If $\mu$ is real-valued, it suffices to have $\widehat{\mu}(n) = 0$
for all $n\in\N$, for then $\wh \mu(-n) = \ov{\wh\mu(n)}$ implies that $\wh\mu(n)= 0$ for all $n\in \Z$.

\subsection{Fourier Multiplier Operators}\label{subsec:FM}

The Plancherel theorem provides a method for constructing
nontrivial bounded operators on $L^2(\R^d)$  as follows. Given $m\in L^\infty(\R^d)$ and $f\in L^2(\R^d)$, the function
$$ m\wh f: \xi\mapsto m(\xi)\wh f(\xi)$$
belongs to $L^2(\R^d)$ and therefore the same is true for
$(m\wh f)\widecheck{\phantom{b}}.$\index{$FC$@$\widecheck{f}$} Here we write $\widecheck g:=\calF^{-1}g$ for the inverse Fourier--Plancherel transform of $g\in L^2(\R^d)$.

\begin{definition}[Fourier multiplier operators]  For  functions  $m\in L^\infty(\R^d)$, the bounded operator
on $L^2(\R^d)$ defined by
$$T_m: f \mapsto (m\wh f)\widecheck{\phantom{b}}$$
is called the {\em Fourier multiplier operator}\index{Fourier!multiplier}\index{operator!Fourier multiplier} with {\em symbol}\index{symbol!of a Fourier multiplier operator} $m$.
\end{definition}

The operator $T_m$ is bounded with norm $\n T_m\n =\n g\mapsto mg \n_{\calL(L^2(\R^d))} = \n m \n_\infty$ (cf. Example \ref{ex:pointwise-multiplier}). We have the elementary properties
$$T_{m_1+m_2} = T_{m_1}+T_{m_2}, \qquad T_{m_1m_2} = T_{m_1}\circ T_{m_2}.$$

Fourier multipliers can be
characterised by commutation properties. This fact depends on the following lemma.

\begin{lemma}\label{lem:char-mult} Let $(\Om,\calF\!,\mu)$ be a $\sigma$-finite measure space and let $1\le p\le \infty$.
 For a bounded operator $T$ on $L^p(\Om)$ the following assertions are equivalent:
 \begin{enumerate}[label={\rm(\arabic*)}, leftmargin=*]
  \item\label{it:char-mult1} $T$ is a pointwise multiplier,\index{multiplier!characterisation of} that is, there exists a function $m\in L^\infty(\Om)$ such that $Tf = mf$ for all $f\in L^p(\Om)$;
  \item\label{it:char-mult2} $T$ commutes with all pointwise multipliers, that is, $TM_\phi = M_\phi T$ for all $\phi\in L^\infty(\Om)$, where $M_\phi f = \phi f$.
 \end{enumerate}
If $\Om = \R^d$ with Lebesgue measure and $1\le p<\infty$, then {\rm\ref{it:char-mult1}} and {\rm\ref{it:char-mult2}} are equivalent to
\begin{enumerate}[label={\rm(\arabic*)}, leftmargin=*]\setcounter{enumi}{2}
 \item\label{it:char-mult3} $TM_{e_\xi} = M_{e_\xi} T$ for all $\xi\in \R^d$, where
 $e_\xi(x) = \exp( ix\cdot \xi)$ for $x\in \R^d$.
\end{enumerate}
\end{lemma}
\begin{proof}
The implication
\ref{it:char-mult1}$\Rightarrow$\ref{it:char-mult2}
is immediate, since pointwise multipliers commute.

\smallskip
\ref{it:char-mult2}$\Rightarrow$\ref{it:char-mult1}:\
Suppose that $T$ commutes with every pointwise multiplier.
If $p=\infty$, then $\one\in L^\infty(\Om)$ and
$$
Tf=T(M_f\one)=M_fT\one=fT\one,
\quad f\in L^\infty(\Om).
$$
Thus $T$ is multiplication by the bounded function $T\one$.

Now let $1\le p<\infty$. By $\sigma$-finiteness, we may choose a function $h\in L^p(\Om)$
such that $h>0$ almost everywhere.
Put $m:=Th/h$.
For every $\phi\in L^\infty(\Om)$, the commutation assumption gives
$$
T(\phi h)=TM_\phi h=M_\phi Th=m\phi h.
$$
We claim that $m\in L^\infty(\Om)$ and
$\n m\n_\infty\le\n T\n$.
To see this, put $M:=\n T\n$ and $E:=\{|m|>M\}$.
Taking $\phi=\one_E$ in the preceding identity gives
$$
\int_E |m|^p h^p\ud\mu
=
\n T(\one_Eh)\n_p^p
\le
M^p\n\one_Eh\n_p^p
=
M^p\int_E h^p\ud\mu.
$$
Since $h>0$ almost everywhere, the left-hand side would be
strictly greater than the right-hand side if $\mu(E)>0$.
Hence $\mu(E)=0$, proving the claim.

Finally, let $f\in L^p(\Om)$ and put
$$
f_n:=\one_{\{|f|\le nh\}}f,\quad n\ge1.
$$
Then $f_n/h\in L^\infty(\Om)$, so the preceding commutation
identity gives $Tf_n=mf_n$.
Moreover, $f_n\to f$ in $L^p(\Om)$ by dominated convergence.
Since both $T$ and multiplication by $m$ are bounded on
$L^p(\Om)$, passing to the limit gives $Tf=mf$.
This proves \ref{it:char-mult1}.

\smallskip
Now let $\Om=\R^d$ with Lebesgue measure and let $1\le p<\infty$.
The implication
\ref{it:char-mult2}$\Rightarrow$\ref{it:char-mult3}
is immediate. To prove the converse, suppose that
\ref{it:char-mult3} holds.

Let $q$ be the conjugate exponent of $p$.
By Theorem \ref{thm:Lp-dual}, we may identify the dual of
$L^p(\R^d)$ with $L^q(\R^d)$. Under this identification, the
adjoint $T^*$ acts on $L^q(\R^d)$ and satisfies
$$
\int_{\R^d}(Tf)g\ud x
=
\int_{\R^d}f(T^*g)\ud x,
\quad f\in L^p(\R^d),\ g\in L^q(\R^d).
$$
Fix such $f$ and $g$. For every $\xi\in\R^d$, the commutation
assumption gives
\begin{align*}
\int_{\R^d}e_\xi(Tf)g\ud x
&= \int_{\R^d}T(e_\xi f)g\ud x
= \int_{\R^d}e_\xi f(T^*g)\ud x.
\end{align*}
By H\"older's inequality, the functions $(Tf)g$ and $f(T^*g)$
belong to $L^1(\R^d)$. The preceding identities show that their
Fourier transforms agree. Injectivity of the Fourier transform,
which follows from Theorem \ref{thm:FT-inversion}, therefore gives
$$
(Tf)g=f(T^*g).
$$

Let $\phi\in L^\infty(\R^d)$. Using this identity, we obtain
$$
\int_{\R^d}T(\phi f)g\ud x
=
\int_{\R^d}\phi f(T^*g)\ud x
=
\int_{\R^d}\phi(Tf)g\ud x.
$$
Since this holds for every $g\in L^q(\R^d)$, duality implies
$T(\phi f)=\phi Tf$. Thus $T$ commutes with every pointwise
multiplier, proving \ref{it:char-mult2}.
\end{proof}

As an application we have the following characterisation of translation-invariant operators on $L^2(\R^d)$.

\begin{theorem}[Translation-invariant operators on $L^2(\R^d)$]\label{thm:transl-invar-L2}
If $T$ is a bounded operator on $L^2(\R^d)$ commuting with every translation, then $T$ is a Fourier multiplier operator, that is,
there exists a (necessarily unique) function $m\in L^\infty(\R^d)$ such that
$T f =\F^{-1} (m \F f)$ for all $f\in L^2(\R^d)$.
\end{theorem}
\begin{proof}
Using the notation of Lemma \ref{lem:char-mult} and letting
$\tau_y f(x):=f(x+y)$, we readily calculate that
$M_{e_y}\calF=\calF\tau_y$ and $\tau_y\calF^{-1}=\calF^{-1}M_{e_y}$ for all $y\in\R^d$.
These identities imply that the operator $ \wt T = \F T \F^{-1}$ has the property $M_{e_y} \wt T  = \wt T M_{e_y}$
 for all $y\in \R^d$\!, and the lemma implies that $\wt T$ is a pointwise multiplier. This means that $T$ is a Fourier multiplier.
\end{proof}

\section{The Hilbert Transform}\label{sec:HT}

A case of special interest concerns the multiplier
$$ m(\xi) = -i\,{\rm sign}(\xi), \quad \xi\in \R.$$ In order to obtain an explicit representation for
the Fourier multiplier operator $T_m$, we observe that
$m = -i(\one_{\R_+} - \one_{\R_-})$ and consider the functions
$$n_a^\pm(\xi) := \exp(-a|\xi|)\one_{\R_\pm}(\xi), \quad a>0.$$
Then
$$ \widecheck{n_a^+}(x) = \frac1{\sqrt{2\pi}}\int_0^\infty  \exp(-a\xi)\exp(ix\xi )\ud \xi
= \frac1{\sqrt{2\pi}} \frac1{a-i x}
$$
and
$$\widecheck{n_a^-}(x) =  \frac1{\sqrt{2\pi}}\int_{-\infty}^0 \exp(a\xi)\exp( ix\xi) \ud \xi
=  \frac1{\sqrt{2\pi}}\frac1{a+i x}.
$$
Formally letting $a\downarrow 0$, in view of Proposition \ref{prop:FT-convol} one expects that
\begin{align*} T_{-i\,{\rm sign}} f & = -i\lim_{a\downarrow 0} (T_{n_a^+} f - T_{n_a^-}f)
\\ & = -\frac{i}{\sqrt{2\pi}}\lim_{a\downarrow 0} (\widecheck{n_a^+}*f - \widecheck{n_a^-}*f)
\\ & = - \frac{i}{2\pi}\lim_{a\downarrow 0}\Bigl(\frac1{a-i (\cdot)} - \frac1{a+i (\cdot)}\Bigr) *f
= \frac1{\pi(\cdot)}*f.
\end{align*}
This suggests the formula
$$ T_{-i\,{\rm sign}} f(x) =  \frac{1}{\pi} \int_{-\infty}^\infty \frac{f(x-y)}{y}\ud y.$$
The above argument is nonrigorous and the convolution with the nonintegrable function $1/x$ is not even well defined
as an operator acting on $L^2(\R)$.
The next theorem turns the above formal derivation into a rigorous argument.
We need the following classical identity.

\begin{lemma}\label{lem:Dirichlet}
$$ \int_0^\infty \frac{\sin(t)}{t} \ud t = \frac{\pi}{2}.$$
\end{lemma}
\begin{proof}
Let
$$I(a):=\int_0^\infty e^{-at}\frac{\sin t}{t}\ud t,\quad a>0.
$$
Differentiation of $I$ under the integral gives
$$
I'(a)=-\int_0^\infty e^{-at}\sin t\ud t =
-\frac1{1+a^2}.
$$
Since $I(a)\to0$ as $a\to\infty$, it follows that
$I(a)=\arctan(1/a)$
and
$$\int_0^T e^{-at}\frac{\sin t}{t}\ud t
= \arctan(1/a) -  \int_T^\infty e^{-at}\frac{\sin t}{t}\ud t.
$$
The result is obtained by combining $\lim_{a\downarrow 0}\arctan(1/a) = \frac12\pi$ with the elementary estimate
\begin{align}\label{eq:Dirichlet}
\Bigl|\int_T^\infty e^{-at}\frac{\sin t}{t}\ud t\Bigr|
\le\frac2T,\quad T>0,\ 0\le a\le1.
\end{align}
This gives the result.
\end{proof}

 \begin{theorem}[Hilbert transform as a Fourier multiplier operator]\label{thm:HT} The Fourier multiplier operator $H:= T_{-i\,{\rm sign}}$ is given by
 $$ H f (\cdot) = \lim_{\eps\downarrow 0} \frac{1}{\pi} \int_{\{|y| >\eps\}} \frac{f(\cdot -y)}{y}\ud y,\quad f\in L^2(\R),$$
 the convergence being in the sense of $L^2(\R)$.
 \end{theorem}

\begin{proof}
Fix $f\in L^2(\R)$. For $R>\varepsilon>0$ set
$$
  \phi_\varepsilon(x):=\frac1{\pi x}\one_{\{|x|>\varepsilon\}}, \qquad
  \phi_{\varepsilon,R}(x)  :=\frac1{\pi x}\one_{\{\varepsilon<|x|<R\}} .
$$
Then $\phi_\varepsilon\in L^2(\R)$, and therefore
$$
  H_\varepsilon f(x):=(\phi_\varepsilon*f)(x)
  =\frac1\pi\int_{\{|y|>\varepsilon\}}\frac{f(x-y)}{y}\ud y
$$
is well defined for every $x\in\R$ by the Cauchy--Schwarz inequality, and bounded as a function of $x$. Also, $\phi_{\varepsilon,R}\in L^1(\R)\cap L^2(\R)$, so Proposition \ref{prop:FT-convol} gives
$$
  \widehat{\phi_{\varepsilon,R}*f}
  =
  \sqrt{2\pi}\,\widehat{\phi_{\varepsilon,R}}\,\widehat f.
$$
Moreover,
$$
 \sqrt{2\pi}\,\widehat{\phi_{\varepsilon,R}}(\xi)
 = \frac1\pi\int_{\varepsilon<|x|<R}\frac{e^{-ix\xi}}{x}\ud x
 = -\frac{2i}{\pi}\int_\varepsilon^R\frac{\sin(x\xi)}{x}\ud x
 =  -\frac{2i}{\pi}\bigl(S(R\xi)-S(\varepsilon\xi)\bigr),
$$
where
$$ S(t):=\int_0^t \frac{\sin s}{s}\ud s.$$
Since $\lim_{t\to\pm\infty}S(t)=\pm\frac12{\pi}$ by Lemma \ref{lem:Dirichlet}, for almost all $\xi\in\R$ we have
$$
\lim_{R\to\infty}  \sqrt{2\pi}\,\widehat{\phi_{\varepsilon,R}}(\xi) =   m_\varepsilon(\xi)
  :=
  -i\,{\rm sign}(\xi)+\frac{2i}{\pi}S(\varepsilon\xi).
$$
Also, since $S$ is bounded, the functions
$\sqrt{2\pi}\,\widehat{\phi_{\varepsilon,R}}$ are bounded, uniformly in
$R>\varepsilon$. Therefore, by dominated convergence,
$$
  \sqrt{2\pi}\,\widehat{\phi_{\varepsilon,R}}\,\widehat f
  \to
  m_\varepsilon\,\widehat f
  \ \hbox{ in }L^2(\R).
$$
By Plancherel's theorem it follows that
$$
  \phi_{\varepsilon,R}*f
  \to
  \mathscr F^{-1}(m_\varepsilon\widehat f) \ \hbox{ in }L^2(\R).
$$

On the other hand, by the Cauchy--Schwarz inequality and dominated convergence,
$$
  \|(\phi_{\varepsilon,R}-\phi_\varepsilon)*f\|_\infty
  \le
  \|\phi_{\varepsilon,R}-\phi_\varepsilon\|_2\,\|f\|_2
  \to 0 .
$$
Thus, after passing to a subsequence if necessary, the above $L^2$-limit
agrees almost everywhere with $\phi_\varepsilon*f$. Consequently,
$$
  \widehat{\phi_\varepsilon*f}=m_\varepsilon\,\widehat f \ \hbox{ in }L^2(\R).
$$

Finally, since $S$ is bounded and $S(\varepsilon\xi)\to0$ for each
$\xi\in\R$ as $\varepsilon\downarrow0$, dominated convergence gives
$$
  m_\varepsilon\,\widehat f
  \to
  -i\,{\rm sign}\,\widehat f \ \hbox{ in }L^2(\R).
$$
Another application of Plancherel's theorem therefore gives
\[
  \phi_\varepsilon*f\to T_{-i\,{\rm sign}}f=Hf
  \ \hbox{ in }L^2(\R).
\qedhere\]
\end{proof}

\begin{definition}[Hilbert transform]\label{def:HT} The operator $H= T_{-i\,{\rm sign}}$ is called the {\em Hilbert transform}\index{Hilbert transform}.
\end{definition}

The Hilbert transform has deep applications in harmonic analysis and the theory of partial differential equations. It will resurface in our treatment of the Poisson semigroup in Chapter \ref{chap:semigroups}. Its connection with harmonic functions is pointed out in Problem \ref{prob:harmonic}.

The final theorem of this section gives a characterisation of the Hilbert transform in terms of commutation properties.
A {\em dilation}\index{dilation} on $L^2(\R^d)$ is an operator of the form $$D_\delta f (x):= f(\delta x), \quad x\in \R^d\!,$$ where $\delta>0$.

\begin{theorem}[Characterisation of the Hilbert transform]\label{thm:char-HT}\index{Hilbert transform!characterisation}
If $T$ is a bounded operator on $L^2(\R)$
commuting with every translation and dilation, then $T$ is a linear combination of the identity operator and the Hilbert transform.
\end{theorem}
\begin{proof}
Theorem \ref{thm:transl-invar-L2} tells us that $T$ is a Fourier multiplier operator, say $T = T_m$ with $m\in L^\infty(\R)$.
 For $\delta>0$, simple calculations give
 $$T_m D_\delta f(x) = T_{m(\delta\cdot)} f(\delta x)\ \hbox{ and }\
  D_\delta T_m f(x) =  T_{m} f(\delta x)$$
 for almost all $x\in\R$.
It follows that $T_m D_\delta = D_\delta T_m$ if and only if
$ T_{m(\delta\cdot)} =  T_{m}$, that is, if and only if $m(\delta\xi) = m(\xi)$ for almost all $\xi\in \R$. This holds for all $\delta>0$, and we claim that $m$ is constant almost everywhere on both $\R_+$ and $\R_-$. We prove this for $\R_+$; the proof for $\R_-$ is the same.

Put $c:=\int_0^1m(s)\ud s$.
For every $t>0$, a change of variables and the dilation identity
with $\delta=t$ give
$$
\int_0^t m(s)\ud s
=t\int_0^1m(tx)\ud x
=t\int_0^1m(x)\ud x
=tc.
$$
Since $m$ is locally integrable, the left-hand side has derivative
$m(t)$ for almost every $t>0$, by the Lebesgue differentiation
theorem. Differentiating therefore gives $m(t)=c$ for almost
every $t>0$.

Applying the same argument to $s\mapsto m(-s)$ shows that $m$ is
also constant almost everywhere on $\R_-$.
Hence $m=a\one+b\sign$ for suitable $a,b\in\C$.
Since $H=T_{-i\,{\rm sign}}$, we obtain
$T=aI+ibH$, as required.
\end{proof}

\section{Interpolation}\label{sec:interpolation}

In general it can be difficult to establish $L^p$-boundedness of operators acting on spaces of measurable functions.
In such situations, interpolation theorems may be helpful. They serve to establish $L^p$-boundedness in situations
where suitable boundedness properties can be established for `endpoint' exponents $p_0$ and $p_1$ satisfying $p_0\le p\le p_1$.
In typical applications one takes $p_0\in \{1,2\}$ and $p_1\in \{2,\infty\}$ or, as in the case of the Hilbert transform, a large finite exponent $p_1=2^n$ (cf. Sections \ref{subsec:HY} and \ref{subsec:HTp}).

\subsection{The Riesz--Thorin Interpolation Theorem}\label{subsec:RT}

In the next theorem we assume that the scalar field is $\C$. When the scalar field is $\R$ we may apply the complex version of the theorem to the complexified operators at the expense of an additional constant $2$. At the end of the section we discuss some conditions under which this additional constant can be avoided.

\begin{theorem}[Riesz--Thorin interpolation theorem]\label{thm:RieszThorin}\index{theorem!Riesz--Thorin}
Let  $(\Om,\calF\!,\mu)$  and $(\Om'\!,\calF'\!,\mu')$ be measure spaces,
with $(\Om'\!,\calF'\!,\mu')$ $\sigma$-finite,
and let $1\le p_0,p_1,q_0,q_1\leq\infty$.
Let $$T_0: L^{p_0}(\Om) \to L^{q_0}(\Om'), \qquad T_1: L^{p_1}(\Om)\to L^{q_1}(\Om')$$
be bounded operators which are consistent in the sense that
$T_0 f = T_1 f$  $\mu'$-almost everywhere, for all $f\in L^{p_0}(\Om)\cap L^{p_1}(\Om)$.
Assume furthermore that
\begin{align*}
  \n{T_0f}\n_{L^{q_0}(\Om')}\leq A_0\n{f}\n_{L^{p_0}(\Om)},\quad  f\in
L^{p_0}(\Om),\\
  \n{T_1f}\n_{L^{q_1}(\Om')}\leq A_1\n{f}\n_{L^{p_1}(\Om)},\quad  f\in
L^{p_1}(\Om).
\end{align*}
Let $0<\theta<1$ and set \begin{equation*}
  \frac{1}{p_\theta}:=\frac{1-\theta}{p_0}+\frac{\theta}{p_1}, \qquad
  \frac{1}{q_\theta}:=\frac{1-\theta}{q_0}+\frac{\theta}{q_1}.
\end{equation*}
Then the common restriction $T$ of $T_0$ and $T_1$ to $L^{p_0}(\Om)\cap L^{p_1}(\Om)$ maps this space
into $L^{q_\theta}(\Om')$
and extends uniquely to a bounded operator
$$ T: L^{p_\theta}(\Om)\to L^{q_\theta}(\Om')$$
satisfying
\begin{equation*}
  \n{Tf}\n_{L^{q_\theta}(\Om')} \leq A_0^{1-\theta}A_1^\theta
\n{f}\n_{L^{p_\theta}(\Om)},\quad f\in L^{p_\theta}(\Om).
\end{equation*}
\end{theorem}

We begin with a simple lemma, which corresponds to the special case where the interpolated operator is the identity operator.

\begin{lemma}\label{lem:Lp-logconvex}
Let $(\Om,\calF\!,\mu)$ be a measure space, let
$1\le r_0,r_1\le\infty$ and $0<\theta<1$, and set
$\frac1{r_\theta}:=\frac{1-\theta}{r_0}+\frac{\theta}{r_1}$.
Then for all $f\in L^{r_0}(\Omega)\cap L^{r_1}(\Omega)$ we have
$f\in L^{r_\theta}(\Om)$ and
$$ \| f \|_{r_\theta} \le \| f\|_{r_0}^{1-\theta}\| f\|_{r_1}^\theta.
$$
\end{lemma}
\begin{proof} By symmetry, we may assume that $r_0\le r_1$.
If $1\le r_0\le r_1 < \infty$, write $|f|^{r_\theta} = |f|^{(1-\theta)r_\theta}|f|^{\theta r_\theta}$ and apply H\"older's inequality.
If $1\le r_0 < r_1 = \infty$, then
$$
\|f\|_{r_\theta}^{r_\theta}
=\int_\Om |f|^{r_0}\,|f|^{\theta r_\theta}\ud \mu
\le \Bigl(\int_\Om |f|^{r_0}\ud \mu\Bigr)\,\bigl\|\,|f|^{\theta r_\theta}\bigr\|_\infty
= \|f\|_{r_0}^{r_0}\,\|f\|_\infty^{\theta r_\theta},
$$ and taking $1/r_\theta$-powers gives the result.
Finally, if $r_0=r_1=\infty$, there is nothing to prove.
\end{proof}

Consider the open strip
$$\mathbb{S}:= \{z\in\C:\, 0<\Re z <1\}.$$

\begin{lemma}[Three lines lemma]\label{lem:three-lines}\index{lemma!three lines}
Suppose that $F: \ov{\mathbb{S}} \to \C $ is a bounded continuous function,
holomorphic on $\mathbb{S}$, and satisfying $$\sup_{v\in\R} | F(iv)| \le A_0, \qquad \sup_{v\in\R} | F(1+iv)|\le A_1.$$
Then for all $0<\theta<1$ we have
$$ \sup_{v\in \R} | F(\theta+iv)| \le A_0^{1-\theta} A_1^\theta.
$$
\end{lemma}
\begin{proof}
Let $F$ satisfy the assumptions of the lemma with constants $A_0$ and $A_1$.  Applying the argument below with $A_j+\eta$ in place of $A_j$ and letting $\eta\downarrow0$, we may assume $A_0,A_1>0$.

For each $\eps>0$ the function $$F_\eps(z) := A_0^{z-1}A_1^{-z}\exp(\eps z(z-1))F(z)$$
satisfies the assumptions of the lemma with constants  $A_{0,\eps} = A_{1,\eps} = 1$.
Moreover, we have  $\lim_{|v|\to\infty} | F_\eps(u+iv)| = 0$
uniformly with respect to $u\in [0,1]$. Hence for large enough $R$ we have $| F_\eps|\le 1$ on the boundary
of the rectangle $\Re z \in [0,1]$, $|\Im z|\le R$. The maximum modulus principle implies that $| F_\eps| \le 1$
on this rectangle,
and by letting $R\to\infty$ we find that $| F_\eps|\le 1$ on $\ov{\mathbb{S}}$. The lemma now follows by letting
$\eps\downarrow 0$.
\end{proof}

Inspection of the proof reveals that the boundedness assumption on $F$ can be relaxed. It cannot be entirely dispensed with, however: the function
$F(z)=\exp\bigl(-i\exp(\pi i z)\bigr)$
satisfies $|F|=1$ on the boundary lines $\Re z=0$ and
$\Re z=1$, while on the line $z=\frac12+iy$ we have
$|F(z)|=\exp(e^{-\pi y})$,
which tends to $\infty$ as $y\to-\infty$.

We are now ready for a proof of the Riesz--Thorin theorem.

\begin{proof}[Proof of Theorem \ref{thm:RieszThorin}]
We use the convention $1/\infty=0$.
Let $T$ be the common
restriction of $T_0$ and $T_1$ to $L^{p_0}(\Om)\cap L^{p_1}(\Om)$.
If $p_0=p_1=p$, Lemma \ref{lem:Lp-logconvex} gives
$$
  \|Tf\|_{q_\theta}
  \le \|T_0f\|_{q_0}^{1-\theta}\|T_1f\|_{q_1}^\theta
  \le A_0^{1-\theta}A_1^\theta\|f\|_p,\quad f\in L^p(\Om),
$$
which proves the theorem. We may therefore assume that $p_0\ne p_1$,
so that $1<p_\theta<\infty$.

\smallskip
{\em Step 1} --
Let $a$ be $\mu$-simple with $\|a\|_{p_\theta}=1$, and let $b\ge0$
be $\mu'$-simple with $\|b\|_{q_\theta'}\le1$.
For $z\in\ov{\mathbb S}$ define
$$
  \frac1{p_z}:=\frac{1-z}{p_0}+\frac z{p_1},\qquad
  \frac1{q_z'}:=\frac{1-z}{q_0'}+\frac z{q_1'},
$$
and put $$f_z:=\one_{\{a\ne0\}}|a|^{p_\theta/p_z}a/|a|.$$
If $q_\theta'<\infty$, put $g_z:=\one_{\{b>0\}}b^{q_\theta'/q_z'}$.
If $q_\theta'=\infty$, then $q_0=q_1=1$, and we put $g_z:=b$.
All functions are set equal to zero at the respective zeros of $a$
and $b$, and $t^w:=\exp(w\log t)$ for $t>0$.
Thus $f_\theta=a$, $g_\theta=b$, and each $f_z$ belongs to $L^{p_0}(\Om)\cap L^{p_1}(\Om)$.

Fix $j\in \{0,1\}$ and $v\in\R$. If $p_j<\infty$, then
\begin{align*}
  \|f_{j+iv}\|_{p_j}^{p_j}
  &=\int_\Om |a|^{p_\theta}\ud\mu=1;
\end{align*}
if $p_j=\infty$, then $|f_{j+iv}|=\one_{\{a\ne0\}}$ and $\n f_{j+iv}\n_\infty =1$.
Likewise, if $q_\theta'<\infty$ and $q_j'<\infty$, then
\begin{align*}\|g_{j+iv}\|_{q_j'}^{q_j'}
  &=\int_{\Om'}b^{q_\theta'}\ud\mu'\le 1;
\end{align*}
if $q_\theta'<\infty$ and $q_j'=\infty$, then
$|g_{j+iv}|=\one_{\{b>0\}}$ and $\|g_{j+iv}\|_\infty\le 1$; and if $q_\theta'=\infty$, then $g_{j+iv} = b$
and $q_j'=\infty$, and $\|g_{j+iv}\|_\infty=\|b\|_\infty\le1$.

Consequently, in all cases we have
$\|f_{j+iv}\|_{p_j}\le 1$ and $\|g_{j+iv}\|_{q_j'}\le 1$.

\smallskip
{\em Step 2} -- Put $h:=Ta$ and define
$\eta:=\overline h/|h|$ on $\{h\ne 0\}$, and $\eta:=0$ on $\{h=0\}$.
Then $|\eta|\le1$ and $h\eta=|h|$. For $z\in\ov {\mathbb S}$ define
$$
  F(z):=\int_{\Om'}(Tf_z)\eta g_z\ud\mu',
$$
so that
$F(\theta)
=\int_{\Om'}(Ta)\eta b\ud\mu'
=\int_{\Om'}|Ta|b\ud\mu'$.
We will check that $F$ satisfies the assumptions of the three lines lemma.
To this end write
$a=\sum_k a_k\one_{A_k}$ and $b=\sum_\ell b_\ell\one_{B_\ell}$, each
with nonzero coefficients and disjoint sets of finite measure.
Then $$f_z=\sum_k c_k(z)\one_{A_k\cap\{a\ne 0\}}, \qquad
g_z=\sum_\ell d_\ell(z)\one_{B_\ell\cap\{b>0\}},$$ where each coefficient is a
nonzero multiple of an exponential of the form $e^{cz}$ with $c\in\R$.
H\"older's inequality gives
$$
  \int_{\Om'}|T\one_{A_k}|\one_{B_\ell}\ud\mu'
  \le\|T_0\one_{A_k}\|_{q_0}\|\one_{B_\ell}\|_{q_0'}<\infty.
$$
Hence the integral defining $F$ is absolutely convergent and
$$
  F(z)=\sum_{k,\ell}c_k(z)d_\ell(z)
       \int_{\Om'}(T\one_{A_k})\eta \one_{B_\ell}\ud\mu'.
$$
Thus $F$ is a finite sum of exponentials $e^{cz}$ with real $c$, and hence is bounded and continuous on $\ov{\mathbb S}$ and holomorphic on $\mathbb S$.
Moreover, H\"older's inequality and the endpoint bounds give
$$
  |F(j+iv)|
  \le\|T_jf_{j+iv}\|_{q_j}\|\eta g_{j+iv}\|_{q_j'}
  \le A_j,\quad j=0,1,\  v\in\R,
$$
because $\n g_{j+iv}\n_{q_j'} \le 1$. The three lines lemma therefore implies
$$
  \int_{\Om'}|Ta|b\ud\mu'=F(\theta)\le A_0^{1-\theta}A_1^\theta.
$$
Since $Ta=T_ja\in L^{q_j}(\Om')$ for $j=0,1$,
Lemma \ref{lem:Lp-logconvex} already gives
$Ta\in L^{q_\theta}(\Om')$.
Taking the supremum in the preceding inequality over all
nonnegative $\mu'$-simple functions $b$ with
$\n b\n_{q_\theta'}\le1$, we obtain
$$
\n Ta\n_{q_\theta}\le A_0^{1-\theta}A_1^\theta.
$$
Indeed, if $q_\theta'<\infty$, these test functions suffice in
the usual dual norm formula by density.
If $q_\theta'=\infty$, then $q_\theta=1$; in this case take
$b=\one_{F_n}$, where $F_n\uparrow\Om'$ and each $F_n$ has finite
measure, and apply the monotone convergence theorem.
Here we use the $\sigma$-finiteness of $\Om'$.

By homogeneity, this proves that
\begin{align}\label{eq:RT-simple}
\n Ta\n_{q_\theta}\le A_0^{1-\theta}A_1^\theta\n a\n_{p_\theta}
\end{align}
for every $\mu$-simple function $a$.

\smallskip
{\em Step 3} --
Let $f\in L^{p_0}(\Om)\cap L^{p_1}(\Om)$.
Choose $\mu$-simple functions $f_n$ such that
$f_n\to f$ almost everywhere, $|f_n|\le |f|$ almost everywhere,
and, if one of the exponents $p_j$ equals $\infty$, such that in
addition
$\|f_n-f\|_\infty\to 0$.
For every finite endpoint exponent $p_j$, dominated convergence gives
$\|f_n-f\|_{p_j}\to0$. Hence $f_n\to f$ in both endpoint norms.

By Lemma \ref{lem:Lp-logconvex},
\begin{align*}
  \|f_n-f\|_{p_\theta}
  &\le\|f_n-f\|_{p_0}^{1-\theta}\|f_n-f\|_{p_1}^\theta\to0,
\intertext{and}
  \|Tf_n-Tf\|_{q_\theta}
  &\le\|T_0(f_n-f)\|_{q_0}^{1-\theta}
       \|T_1(f_n-f)\|_{q_1}^\theta\to0.
\end{align*}
Passing to the limit in the estimate \eqref{eq:RT-simple} for $f_n$ proves
$$\|Tf\|_{q_\theta}\le A_0^{1-\theta}A_1^\theta\|f\|_{p_\theta}, \quad f\in L^{p_0}(\Om)\cap L^{p_1}(\Om).$$
Finally, $L^{p_0}(\Om)\cap L^{p_1}(\Om)$ contains all $\mu$-simple functions and is therefore
dense in $L^{p_\theta}(\Om)$ by
Proposition \ref{prop:approx-simple}.
This gives the asserted unique extension.
\end{proof}

Up to this point we have implicitly assumed that the scalar field is complex. Suppose now that the scalar field is real. We may extend a bounded operator $T: L^p(\Om)\to L^q(\Om')$ to a bounded operator $T_\C:L^p(\Om;\C)\to L^q(\Om';\C)$ by setting
$$ T_\C (u+iv):= Tu + i Tv$$
for real-valued $u,v\in L^p(\Om)$.
The triangle inequality implies the trivial bounds $$\n T\n \le \n T_\C\n\le 2\n T\n.$$ As an immediate consequence of this, the Riesz--Thorin theorem extends to the real-valued case at the expense of an additional multiplicative constant $2$.
We will now discuss some situations where this additional constant can be avoided.

If $T$ is a positivity-preserving operator (that is, $f\ge 0$ implies $Tf\ge 0$), the identity
$$|a + ib | = \sup_{\theta\in [0,2\pi]} |a\cos\theta  + b\sin\theta|$$
(for a proof, rotate the point $(a,b)\in\R^2$ to the positive $x$-axis)
together with the inequality $|Tg|\le T|g|$ for real-valued $g\in L^p(\Om)$ (which follows from \eqref{eq:Tpos-mod}) implies the pointwise bound
\begin{align*} |T_\C f| & = |Tu + i Tv| = \sup_{\theta\in [0,2\pi]} |(Tu)\cos\theta  + (Tv)\sin\theta|
\\ & \le \!\sup_{\theta\in [0,2\pi]}\! T|u\cos\theta  + v\sin\theta|\le T|u+iv| = T|f|.
\end{align*}
This implies the norm bound $\n T_\C\n\le \n T\n$ and hence equality
\begin{align}\label{eq:TcT} \n T_\C\n = \n T\n.
\end{align}

With some additional work, the equality \eqref{eq:TcT} can be extended to arbitrary bounded operators $T$ and exponents $1\le p\le q<\infty$.
The proof is based on the observation that for all $z = a+bi \in \C$ and $1\le q<\infty$ we have
\begin{align}\label{eq:gamma-p}
|z| =  \frac{1}{\n\gamma\n_{q}}(\wt\E|a\gamma_1 + b\gamma_2|^q)^{1/q}\!,
\end{align}
where $\gamma,\gamma_1,\gamma_2$ are real-valued standard Gaussian random variables defined on some probability space $\wt\Omega$, with $\gamma_1$ and $\gamma_2$ independent, and $\wt\E$ denotes the expectation.
Indeed, if $|z|=1$, then $a \gamma_1+ b\gamma_2$ is another standard Gaussian random
variable and $$(\wt\E |a \gamma_1 +b \gamma_2|^q)^{1/q}=\n{\gamma_1}\n_{q}.$$ The general
case follows from this by scaling.

Assume that $1\le p\le q<\infty$. Combining \eqref{eq:gamma-p} with
Fubini's theorem (after replacing $\Omega$ and $\Omega'$ by $\sigma$-finite
measurable subsets supporting $u$, $v$, $Tu$, and $Tv$), we obtain
\begin{align*}
  \n{\gamma}\n_{q}^q \n {T_{\C}(u+iv)}\n_{L^q(\Om';\C)}^q
  & =  \n{\gamma}\n_{q}^q\int_{\Om'}|Tu + iTv|^q\ud \mu'
  \\ & = \int_{\Om'} \wt \E  |\gamma_1Tu + \gamma_2 Tv|^q\ud \mu'
   =\n {\gamma_1 Tu+\gamma_2 Tv}\n_{L^q(\Om';L^q(\wt\Omega))}^q
  \\ & =\n {T(\gamma_1 u+\gamma_2 v)}\n_{L^q(\wt\Omega;L^q(\Om'))}^q
  \leq\n {T}\n^q\n {\gamma_1 u+\gamma_2 v}\n_{L^q(\wt\Omega;L^p(\Om))}^q
  \\ & \le \n {T}\n^q\n {\gamma_1 u+\gamma_2 v}\n_{L^p(\Om;L^q(\wt\Omega))}^q
 =\n {\gamma}\n_{q}^q\n {T}\n^q\n{u+iv}\n_{L^p(\Om;\C)}^q.
\end{align*}
In the penultimate step we used the continuous version of Minkowski's inequality (Problem \ref{prob:cont-Mink}).
This proves \eqref{eq:TcT} for $1\le p\le q<\infty$. In the range $q<p$, \eqref{eq:TcT} is generally false, as is shown by a classical counterexample due to M. Riesz.

Suppose now that the spaces in Theorem \ref{thm:RieszThorin}
are real. If both endpoint operators are positivity preserving,
or if
$$
1\le p_j\le q_j<\infty,\quad j=0,1,
$$
then the preceding arguments show that their complexifications
$S_j:=(T_j)_\C$ satisfy the endpoint estimates with the same
constants $A_j$.
The complex Riesz--Thorin theorem therefore gives a bounded
operator
$$
S_\theta:L^{p_\theta}(\Om;\C)\to L^{q_\theta}(\Om';\C)
$$
of norm at most $A_0^{1-\theta}A_1^\theta$.
On the real space $L^{p_0}(\Om)\cap L^{p_1}(\Om)$,
this operator agrees with the common restriction of $T_0$
and $T_1$, and therefore takes real-valued functions to
real-valued functions. Since this intersection is dense in
$L^{p_\theta}(\Om)$ and the real-valued functions form a closed
subspace of $L^{q_\theta}(\Om';\C)$, the same holds for every
real-valued input. Restricting $S_\theta$ to the real spaces
proves the Riesz--Thorin theorem with the same constant
under either of the stated conditions.

\subsection{The Hausdorff--Young Theorem}\label{subsec:HY}

This brief section and the next are devoted to some applications of the Riesz--Thorin theorem.

As we have seen in Section \ref{sec:FT}, the Fourier transform is bounded from $L^1(\R^d)$ to $L^\infty(\R^d)$
and its restriction to $L^1(\R^d)\cap L^2(\R^d)$ extends to an isometry from $L^2(\R^d)$ onto itself.
The Fourier transform with respect to the normalised Lebesgue measure ${\rm d}m(x) = (2\pi)^{-d/2}\ud x$
defined in Remark \ref{rem:FT-normalised} is contractive from $L^1(\R^d\!,m)$ to $L^\infty(\R^d\!,m)$,
and its restriction to $L^1(\R^d\!,m)\cap L^2(\R^d\!,m)$ extends to an isometry from $L^2(\R^d\!,m)$ onto itself by Remark \ref{rem:FT-normalised-L2}.
Accordingly, the Riesz--Thorin theorem implies:

\begin{theorem}[Hausdorff--Young]\index{theorem!Hausdorff--Young} Let $1\le p\le 2$ and $\frac1p+\frac1q = 1$.
The restriction to $L^1(\R^d)\cap L^2(\R^d)$ of the Fourier transform has a unique extension to a bounded operator from
 $L^p(\R^d)$ to $L^q(\R^d)$. With respect to the normalised Lebesgue measure, the Fourier transform has a unique extension to a contraction from
 $L^p(\R^d\!,m)$ to $L^q(\R^d\!,m)$.
\end{theorem}

A similar result holds for the Fourier transform on the circle (see Problem \ref{prob:HY-FT-circle}).

\subsection{\texorpdfstring{$L^p$}{L\textasciicircum p}-Boundedness of the Hilbert Transform}\label{subsec:HTp}

A second application of the Riesz--Thorin theorem is the following theorem
of M. Riesz on the $L^p$-boundedness of the Hilbert transform.

\begin{theorem}[Riesz]\label{thm:HT-Lp-bdd}\index{theorem!L@$L^p$-boundedness of the Hilbert transform}
For all $1<p<\infty$ the restriction of the Hilbert transform to $L^2(\R)\cap L^p(\R)$
has a unique extension to a bounded operator on $L^p(\R)$.\index{Hilbert transform!$L^p$-boundedness}
\end{theorem}

The proof of Theorem \ref{thm:HT-Lp-bdd} is based on a couple of lemmas.

\begin{lemma}\label{lem:HC1}
If $f\in C^1_{\rm c}(\R)$, then $Hf\in  L^p(\R)$ for all $2\le p\le \infty$.
\end{lemma}
\begin{proof}
Let $I$ be a bounded interval containing the support of $f$.
The pointwise identity
\begin{align*}
H_{\eps}f(x) & = \frac{1}{\pi}\int_{\varepsilon}^\infty \frac{f(x-y) - f(x+y)}{y} \ud y
\\ & = \frac{1}{\pi}\int_{\varepsilon}^\infty \one_{(-I+x)\cup(I-x)}(y) \frac{f(x-y) - f(x+y)}{y} \ud y, \quad x\in\R,
\end{align*}
implies the bound
\begin{align}\label{eq:HC1}
|H_{\eps}f(x)| \le \frac1\pi\cdot 2|I|\cdot 2\n f'\n_\infty, \quad x\in\R.
\end{align}
As $\eps\downarrow 0$, we have $H_\eps f\to Hf$ in $L^2(\R)$
by Theorem \ref{thm:HT} and, after passing to a subsequence converging almost everywhere, \eqref{eq:HC1} implies that $Hf\in  L^\infty(\R)$. This implies $Hf\in L^p(\R)$ and, for $2\le p<\infty$, $\n Hf\n_p^p\le\n Hf\n_\infty^{p-2}\n Hf\n_2^2$.
\end{proof}

\begin{lemma}\label{lem:HplusiHf}
The Hilbert transform of a real-valued function $u\in L^2(\R)$ is the unique real-valued function $v\in L^2(\R)$ such that the Fourier--Plancherel transform of $u + iv$ vanishes on $\R_-$.
\end{lemma}
\begin{proof}
By Theorem \ref{thm:HT}, $Hu$ is the $L^2(\R)$-limit of the
real-valued functions $H_\eps u$, and is therefore real-valued
almost everywhere.
 That the Fourier transform of $u+iHu$ vanishes on $\R_-$ is immediate from Theorem \ref{thm:HT}. In the converse direction, let $u,v\in L^2(\R)$ be real-valued such that the Fourier--Plancherel transform of $u+iv$ vanishes on $\R_-$. Then for almost all $\xi>0$ we have
\begin{equation*}
  0=\wh{u}(-\xi)+i\wh{v}(-\xi)=\overline{\wh{u}(\xi)}+i\overline{\wh{v}(\xi)}
   =\overline{\wh{u}(\xi)-i\wh{v}(\xi)},
\end{equation*}
so $\wh{v}(-\xi) = i\wh{u}(-\xi)$ and $\wh{v}(\xi) = -i\wh{u}(\xi)$. Hence, for almost all $\xi\in \R$,
\begin{equation*}
  \wh{v}(\xi)=-i\sign(\xi)\wh{u}(\xi).
\end{equation*}
By Theorem \ref{thm:HT} this proves that $v = Hu$.
\end{proof}

In the next lemma we use the fact that if
$\phi\in C_{\rm c}^2(\R)$, then  $\wh{\phi''}(\xi) = -|\xi|^2 \wh\phi(\xi)$ is bounded, and therefore
$\wh\phi$ is integrable.

\begin{lemma}[Cotlar]\label{lem:Cotlar}\index{lemma!Cotlar}
Let $H$ be the Hilbert transform on $L^2(\R)$. For all real-valued $u\in C_{\rm c}^2(\R)$ we have
$$(Hu)^2=u^2+2H(u Hu).$$
\end{lemma}

\begin{proof}
Let $u,v\in C_{\rm c}^2(\R)$ be real-valued and put
$$
        U:=u+iHu,\qquad V:=v+iHv .
$$
Since $u,v\in C_{\rm c}^2(\R)$, their Fourier transforms belong to
$L^1(\R)\cap L^2(\R)$. By Theorem \ref{thm:HT} we have
$$
        \widehat{Hu}=-i\,{\rm sign}(\cdot)\widehat u(\cdot),
        \qquad
        \widehat{Hv}=-i\,{\rm sign}(\cdot)\widehat v(\cdot).
$$
Hence $\widehat U$ and $\widehat V$ belong to $L^1(\R)\cap L^2(\R)$
and vanish on $\R_-$. By Fourier inversion, $U$ and $V$ may be
represented by bounded continuous functions. In particular, $UV\in L^2(\R)$.

The product formula
$$
        \widehat{UV} = (2\pi)^{-1/2}\widehat U*\widehat V
$$
follows from Proposition \ref{prop:FT-convol} by Fourier inversion. Since
$\widehat U,\widehat V\in L^1(\R)\cap L^2(\R)$, the convolution on the
right belongs to $L^1(\R)\cap L^2(\R)$. Moreover, it vanishes on
$\R_-$. Therefore the Fourier--Plancherel transform of
$$
        -iUV = (uHv+Hu\,v)+i(Hu\,Hv-uv)
$$
vanishes on $\R_-$. Lemma \ref{lem:HplusiHf} gives
$$
        Hu\,Hv-uv = H(uHv+Hu\,v).
$$
Cotlar's identity follows by taking $u = v$.
\end{proof}

\begin{proof}[Proof of Theorem \ref{thm:HT-Lp-bdd}]
The proof consists of three steps.
First we prove the theorem for exponents $p=2^n$ with $n\ge 1$ using
Cotlar's identity, then for $2<p<\infty$ by interpolation, and finally for $1<p<2$ by duality.

 \smallskip
{\em Step 1} -- In this step we show that if $H$ is $L^p$-bounded for some $2\le p<\infty$, then $H$ is $L^{2p}$-bounded. The proof also gives a bound for $\n H\n_{2p}$ in terms of $\n H\n_p$. In what follows we set $\n H\n_p =: c_p$.

Let $u\in C_{\rm c}^2(\R)$ be real-valued. Then $Hu \in L^{2p}(\R)$ by Lemma \ref{lem:HC1}. By Cotlar's identity (Lemma \ref{lem:Cotlar}) and H\"older's inequality,
$$ \n (Hu)^2\n_p \le \n u^2\n_p + 2 \n H\n_p \n u Hu\n_p \le  \n u^2\n_p + 2c_p \n u\n_{2p}\n Hu\n_{2p}.$$
Using the identity $\n v^2\n_p = \n v\n_{2p}^2$, we obtain
$$ \n  Hu\n_{2p}^2 \le \n u\n_{2p}^2 + 2c_p \n u\n_{2p}\n Hu\n_{2p},$$
or equivalently,
$$ (\n Hu\n_{2p} - c_p\n u\n_{2p})^2 \le (1+c_p^2)\n u\n_{2p}^2.$$
It follows that
$$ \n Hu\n_{2p} - c_p\n u\n_{2p} \le \sqrt{1+c_p^2}\n u\n_{2p}$$
and hence $$ \n Hu\n_{2p} \le
\Bigl(c_p+\sqrt{1+c_p^2}\,\Bigr)\n u\n_{2p}.$$
Since the real-valued functions in $C_{\rm c}^2(\R)$ are dense in
$L^{2p}(\R;\R)$, this estimate gives a unique bounded extension
to the latter space. By \eqref{eq:TcT}, with both exponents equal
to $2p$, its complexification has the same norm.
Denote this complexified operator by $H_{2p}$.
By complex linearity, it agrees with $H$ on $C_{\rm c}^2(\R)$.

To check agreement on $L^2(\R)\cap L^{2p}(\R)$, let $f$ belong
to this intersection and choose $f_n\in C_{\rm c}^2(\R)$ such that
$f_n\to f$ in both norms. Such an approximation follows by
truncation and mollification, as in Proposition \ref{prop:Cc-dense}.
Then $Hf_n\to Hf$ in $L^2(\R)$, whereas
$Hf_n=H_{2p}f_n\to H_{2p}f$ in $L^{2p}(\R)$.
Passing to subsequences converging almost everywhere shows that
$H_{2p}f=Hf$. Thus we have obtained the required consistent
extension, with
$$
c_{2p}\le c_p+\sqrt{1+c_p^2}.
$$

{\em Step 2} -- Since $H$ is $L^2$-bounded, Step 1 implies that $H$ is $L^{2^n}$-bounded for all integers $n\ge 1$.
The Riesz--Thorin theorem then implies that $H$ is $L^p$-bounded for all $2\le p<\infty$.

To check agreement again, let $H_p$ denote the interpolated extension, and let
$2^n>p$. For $f\in L^2(\R)\cap L^p(\R)$, put
$f_k:=\one_{\{|f|\le k\}}f$. Then
$f_k\in L^2(\R)\cap L^{2^n}(\R)$ and $f_k\to f$ in both
$L^2(\R)$ and $L^p(\R)$. Hence $H_pf_k=Hf_k$ converges to
$H_pf$ in $L^p(\R)$ and to $Hf$ in $L^2(\R)$.
Passing to subsequences converging almost everywhere gives
$H_pf=Hf$.

\smallskip
{\em Step 3} --
Finally, let $1<p<2$ and let $q$ be its conjugate exponent.
Denote the bounded extension obtained in Step 2 by $H_q$ and
put $c_q:=\n H_q\n$.
For each $f\in L^p(\R)$, H\"older's inequality gives
$$
\Bigl|\int_\R fH_qg\ud x\Bigr|
\le c_q\n f\n_p\n g\n_q,\quad g\in L^q(\R).
$$
By Theorem \ref{thm:Lp-dual}, there is therefore a unique
$H_pf\in L^p(\R)$ such that
$$
\int_\R(H_pf)g\ud x
=
-\int_\R fH_qg\ud x,\quad g\in L^q(\R).
$$
This defines a bounded linear operator $H_p$ on $L^p(\R)$
satisfying $\n H_pf\n_p\le c_q\n f\n_p$.

We next check that $H_p$ agrees with the original Hilbert
transform on $L^p(\R)\cap L^2(\R)$.
The multiplier representation of $H$ and Plancherel's theorem give
$$
\int_\R(Hf)g\ud x
=
-\int_\R fHg\ud x,\quad f,g\in L^2(\R).
$$
Hence, if $f\in L^p(\R)\cap L^2(\R)$ and
$g\in L^q(\R)\cap L^2(\R)$, then
$$
\int_\R(H_pf)g\ud x
=
-\int_\R fH_qg\ud x
=
-\int_\R fHg\ud x
=
\int_\R(Hf)g\ud x.
$$
Both $H_pf$ and $Hf$ are locally integrable. Since the test
functions $g$ include all bounded measurable functions of bounded
support, this identity implies that $H_pf=Hf$ almost everywhere.

The density of $L^p(\R)\cap L^2(\R)$ in $L^p(\R)$ gives uniqueness
of the extension. Moreover, under the duality identification
$(L^q(\R))^*=L^p(\R)$, the defining identity says that
$H_p=-H_q^*$. Consequently,
$\n H_p\n=\n H_q\n=c_q$.
\end{proof}

\subsection{The Marcinkiewicz Interpolation Theorem}

In this final section we prove another $L^p$-interpolation theorem, the Marcinkiewicz interpolation theorem. It has the virtue of imposing less stringent assumptions at the endpoints
and on the operator to be interpolated, and of not even requiring the
operator to be linear. On the downside, the constant obtained from the proof is rather poor.
The theorem elaborates on the observation, made after the proof of the Hardy--Littlewood maximal theorem, that the proof of the $L^p$-bound essentially only depended on the weak $L^1$-bound.

By $(L^{p_0}+L^{p_1})(\R^d)$ we denote the vector space of all $f\in L^1_{\rm loc}(\R^d)$ that admit a decomposition
$f = f_0+ f_1$ with $f_0\in L^{p_0}(\R^d)$ and $f_1\in L^{p_1}(\R^d)$.

\begin{theorem}[Marcinkiewicz interpolation theorem]\label{thm:Marcinkiewicz}\index{theorem!Marcinkiewicz}
Let $1\le p_0<p<p_1\le\infty$. Suppose that $T$ assigns to each
$f\in (L^{p_0}+L^{p_1})(\R^d)$ a measurable function $Tf$ on $\R^d$,
modulo equality almost everywhere, and that $T$ is {\em subadditive}
\index{subadditive mapping} in the sense that
$$
|T(f+g)|\le |Tf|+|Tg| \ \hbox{ almost everywhere}
$$
for all $f,g\in (L^{p_0}+L^{p_1})(\R^d)$. Suppose furthermore that, for
$j=0,1$, there are constants $A_j\ge0$ such that
$$
\bigl|\{|Tf|>t\}\bigr| \le \Bigl(\frac{A_j\|f\|_{p_j}}{t}\Bigr)^{p_j}, \quad t>0,\ f\in L^{p_j}(\R^d);
$$
if $p_1=\infty$, we replace the assumption for $j=1$ by
$$
\|Tf\|_\infty\le A_1\|f\|_\infty, \quad f\in L^\infty(\R^d).
$$
Then $T$ maps $L^p(\R^d)$ into $L^p(\R^d)$, and there exists a
constant $C$, depending only on $p_0,p,p_1,A_0,A_1$, such that
$$
 \|Tf\|_p\le C\|f\|_p, \quad f\in L^p(\R^d).
$$
The same conclusion holds if $T$ is a {\em linear} operator defined only on
$L^{p_0}(\R^d)\cap L^{p_1}(\R^d)$,
provided that the endpoint estimates above hold for all $f$ in this intersection. In this case, $T$ has a unique bounded extension to $L^p(\R^d)$ satisfying the above $L^p$-bound.
\end{theorem}

\begin{proof}
We give the proof for $1<p_1<\infty$; the case $p_1=\infty$ proceeds along the lines of Theorem \ref{thm:HL-maximal}, requiring small changes that are left to the reader.

Let $f\in L^p(\R^d)$ and fix $t>0$.
We split $f = f_0+ f_1$ with $f_0\in L^{p_0}(\R^d)$ and $f_1\in L^{p_1}(\R^d)$ by taking $f_0 = \one_{\{|f|\ge t/2\}}f$ and $f_1 = \one_{\{|f|< t/2\}}f$.
From the subadditivity of $T$ we obtain
$$ \{|T(f)| > t\}\subseteq \{|T(f_0)| > t/2\} \cup \{|T(f_1)| > t/2\}$$
and therefore
$$ |\{|T(f)| > t\}| \le |\{|T(f_0)| > t/2\}| +  |\{|T(f_1)| > t/2\}|.$$
Using Lemma \ref{lem:layers}, the endpoint estimates, and
Fubini's theorem, we obtain
\begin{align*}
\int_{\R^d}|Tf(x)|^p\ud x
&=
p\int_0^\infty t^{p-1}|\{|Tf|>t\}|\ud t
\\
&\le
p(2A_0)^{p_0}
\int_0^\infty t^{p-p_0-1}
\int_{\{|f|\ge t/2\}}|f(x)|^{p_0}\ud x\ud t
\\
&\quad+
p(2A_1)^{p_1}
\int_0^\infty t^{p-p_1-1}
\int_{\{|f|<t/2\}}|f(x)|^{p_1}\ud x\ud t
\\
&=
p(2A_0)^{p_0}
\int_{\R^d}\int_0^{2|f(x)|}
t^{p-p_0-1}|f(x)|^{p_0}\ud t\ud x
\\
&\quad+
p(2A_1)^{p_1}
\int_{\R^d}\int_{2|f(x)|}^\infty
t^{p-p_1-1}|f(x)|^{p_1}\ud t\ud x
\\
&=
C^p\int_{\R^d}|f(x)|^p\ud x,
\end{align*}
where
$$
C^p
=
2^p p\Bigl(
\frac{A_0^{p_0}}{p-p_0}
+
\frac{A_1^{p_1}}{p_1-p}
\Bigr).
$$
The conditions $p_0<p<p_1$ ensure convergence of the two
inner integrals when $f(x)\ne0$; when $f(x)=0$, both integrands
vanish.

Suppose finally that $T$ is linear and defined only on
$L^{p_0}(\R^d)\cap L^{p_1}(\R^d)$. The preceding argument applies to every
$f$ in this intersection, since the functions
$f_0=\one_{\{|f|\ge t/2\}}f$ and $f_1=\one_{\{|f|<t/2\}}f$
belong to the same intersection. It follows that
$$
\n Tf\n_p\le C\n f\n_p,\quad
f\in L^{p_0}(\R^d)\cap L^{p_1}(\R^d).
$$
Since $C_{\rm c}^\infty(\R^d)$ is contained in this intersection and is dense in
$L^p(\R^d)$ by Proposition \ref{prop:Cc-dense}, the density argument gives a unique bounded extension of $T$ to $L^p(\R^d)$ with the same bound.
\end{proof}

\begin{problems}

\item
Let $X$ be a complete metric space and let $(U_n)_{n\ge1}$ be a sequence
of open dense subsets of $X$. Show that $\bigcap_{n\ge1}U_n$
is dense in $X$.

\noindent{\em Hint:}\ Given a nonempty open set $U\subseteq X$,
choose an open ball $B(x;r)$ whose closure
$K:=\overline{B(x;r)}$ is contained in $U$.
Apply the Baire category theorem to the complete metric space $K$.

\item
Let $X$ be a Banach space, $Y$ a normed space, and let
$(T_n)_{n\ge1}$ be a sequence in $\calL(X,Y)$ such that
$\sup_{n\ge1}\n T_n\n=\infty$. Show that
$$
\Bigl\{x\in X:\, \sup_{n\ge1}\n T_nx\n=\infty\Bigr\}
$$
is a dense subset of $X$.

\noindent{\em Hint:}\ For $m\ge1$ put
$U_m:=\bigcup_{n\ge1}\{x\in X:\,\n T_nx\n>m\}$.
Show that each $U_m$ is open and dense, and apply the preceding problem.

\item \label{prob:no-count-basis}
Let $(x_n)_{n \geq 1}$ be a sequence with dense linear span in a Banach space $X$.
Using Baire's theorem, prove that this linear span equals $X$ if and only if $\dim X < \infty$.

\item
In this problem it is shown that there is no norm on $L^1_{\rm loc}(0,1)$ which turns it
into a Banach lattice.
Reasoning by contradiction, assume that $L^1_{\rm loc}(0,1)$ carries a Banach lattice
norm $\n\cdot\n$. For $n\ge1$, put
$$
A_n:=\Big(\frac1{n+1},\frac1n\Big],
\qquad
M_nf:=\bigl(1+(n-1)\one_{A_n}\bigr)f.
$$

\begin{enumerate}[\rm(a), leftmargin=*]
\item Show that $M_n$ is a bounded operator and that $\n M_n\n=n$.

\item
Fix $f\in L^1_{\rm loc}(0,1)$ and define
$$g:=\sum_{k\ge1}k|f|\one_{A_k}.$$
Show that $g\in L^1_{\rm loc}(0,1)$ and that for all $n\ge 1$ we have
$ |M_nf|\le g$ almost everywhere.
\item Using that $\n\cdot\n$ is a lattice norm, deduce that  $\sup_{n\ge 1}\n M_n f\n < \infty$.

\item Use the uniform boundedness theorem to arrive at a contradiction.
\end{enumerate}

\item\label{prob:3cex2}
Let $X$ be the linear span of the standard basis vectors of $\ell^2$ and let $P_n\in \calL(\ell^2)$ denote the orthogonal projection of $\ell^2$ onto the $n$th coordinate subspace. Show that $nP_n x \to 0$ for all $x\in X$ and $\n n P_n \n =n \to \infty$.
Conclude that the completeness assumption cannot be omitted from the uniform boundedness theorem.

\item
Let $X$ be a Banach space and let $(x_n)_{n\ge 1}$ be a sequence in $X$ such that
$\sum_{n=1}^\infty \langle x_n,x^*\rangle$ converges for every $x^* \in X^*$.
\begin{enumerate}[\rm(a), leftmargin=*]
\item Show that the sequence $(x_n)_{n\ge 1}$ is bounded.
\item Do the assumptions imply that $\sum_{n=1}^\infty x_n$ converges in $X$?
\end{enumerate}

\item
State and prove an analogue of Proposition \ref{prop:weakly-bounded} for the weak$^*$ topology.

\item\label{prob:holoX-valued}
The aim of this problem is to prove that a weakly holomorphic function is holomorphic. Let us start with the definitions of these notions.
We fix an open set $D\subseteq\C$ and a complex Banach space $X$. A function $f:D\to X$ is said to be:
\begin{itemize}
  \item {\em holomorphic}, if for all $z_0\in D$ the limit
  $$ \lim_{z\to z_0} \frac{f(z)-f(z_0)}{z-z_0}$$
  exists in $X$ (see  Problem \ref{prob:weak-Cauchy});
  \item {\em weakly holomorphic},\index{holomorphic!weakly} if for all $z_0\in D$ and $x^*\in X^*$ the limit
  $$ \lim_{z\to z_0} \Bigl\lb\frac{f(z)-f(z_0)}{z-z_0}, x^*\Bigr\rb$$
  exists in $\C$.
\end{itemize}
Obviously every holomorphic function $f:D\to X$ is weakly holomorphic.

Now let $f:D\to X$ be weakly holomorphic, fix $z_0\in D$, and let $r>0$ be so small that the closed disc $\{z\in \C: \, |z-z_0|\le r\}$ is contained in $D$.
\begin{enumerate}[\rm(a), leftmargin=*]
  \item
  Apply Proposition \ref{prop:weakly-bounded} to the set $${}\!\!\!\!\!\!\! U = \Big\{\frac1{h-g} \Big( \frac{f(z_0+h)-f(z_0)}{h} - \frac{f(z_0+g)-f(z_0)}{g}\Big):\,  g\not=h, 0<|g|, |h| < \frac{r}2 \Big\}.$$
  After applying functionals $x^*\in X^*$ and using the Cauchy integral formula for scalar-valued holomorphic functions (see Problem \ref{prob:weak-Cauchy}),
  prove that there is a constant $M\ge 0$ such that for all $0<|g|,|h|<r/2$ we have $$  \Big\n \frac{f(z_0+h)-f(z_0)}{h} - \frac{f(z_0+g)-f(z_0)}{g}\Big\n \le M|h-g|.$$
  \item Deduce that every weakly holomorphic function $f:D\to X$  is holomorphic.
\end{enumerate}

\item\label{prob:unbbd-Four}
For $f\in L^1(-\pi,\pi)$ and $N \in\N$ define the functions
$$ s_N f (t):= \sum_{n=-N}^N \wh f(n)\exp(int), \quad t\in [-\pi,\pi],$$
where $\wh f(n)$  is the $n$th Fourier coefficient of $f$, that is,
$$ \wh f(n) :=\frac1{2\pi}\int_{-\pi}^{\pi} f(s) \exp(-ins)\ud s, \quad n\in \mathbb{Z}.$$
By the results of Section \ref{subsec:Fourier-basis}, for all $f\in L^2(-\pi,\pi)$ we have
$$ f = \lim_{N\to\infty} s_N f = \sum_{n\in\Z} \wh f(n)\exp( in(\cdot))$$
with convergence in $L^2(-\pi,\pi)$;
the series on the right-hand side is the Fourier series of $f$.
One might hope that if $f\in C[-\pi,\pi]$ is periodic, then its Fourier
series converges to $f$ in the norm of $C[-\pi,\pi]$. The aim of this
problem is to prove that this is false in a strong sense: there exists a
function $f\in C[-\pi,\pi]$ that is periodic in the sense that $f(-\pi)=f(\pi)$ and
whose Fourier series diverges at $t=0$.
\begin{enumerate}[\rm(a), leftmargin=*]
  \item\label{it:unbbd-Four2} Show that $s_N f (t) = \frac1{2\pi}\int_{-\pi}^\pi f(s) D_N(t-s)\ud s$, where
  the {\em Dirichlet kernel}\index{Dirichlet!kernel} is given by
  $$ D_N(t) := \sum_{n=-N}^N \exp( int) = \frac{\sin ((N+\frac12)t)}{\sin (\frac12 t)},\quad t\notin 2\pi\Z,$$
  and $D_N(t) = 2N+1$ for $t\in 2\pi\Z$.

  \item\label{it:unbbd-Four3} Show that $$\| \Lambda_N\| = \| D_N\|_1,$$
  where $\Lambda_N: C_{\rm per}[-\pi,\pi]\to\C$ is given by $\Lambda_N f := s_N f(0)$ and $\n\cdot\n_1$ refers to the normalised Lebesgue measure.

  \noindent{\em Hint:}\ To prove the inequality $\| \Lambda_N\| \ge \| D_N\|_1$, approximate $\sign(D_N)$ pointwise
  almost everywhere by a sequence of continuous periodic functions $f_n$ of norm $\le 1$, set $g_n(t) = f_n(-t)$,
  and use dominated convergence
  to obtain $$  \lim_{n\to\infty}\Lambda_N(g_n)  = \lim_{n\to\infty}  \frac1{2\pi}\int_{-\pi}^\pi f_n(t)D_N(t)\ud t = \frac1{2\pi}\int_{-\pi}^\pi |D_N(t)|\ud t.$$
  Fill in the missing details.

  \item\label{it:unbbd-Four4} Show that $\lim_{N\to\infty} \| D_N\|_1 = \infty$.

  \noindent{\em Hint:}\ Use the inequality $|\sin(x)| \le |x|$ for $x\in\mathbb R$, and then perform some careful estimates
  on the resulting integral.

  \item\label{it:unbbd-Four5}
  Apply the uniform boundedness theorem to prove that there exists a function $f\in C_{\rm per}[-\pi,\pi]$ such that $\sup_{N\ge1}|s_Nf(0)|=\infty$.
\end{enumerate}

\item\label{prob:complemented-quotient}
Let $X$ be a Banach space with a direct sum decomposition $X = X_0\oplus X_1$. Prove that the projections onto the summands
define isomorphisms of Banach spaces $$ X/X_0 \simeq X_1, \qquad X/X_1 \simeq X_0.$$

\item
Let $X$ be a Banach space with a direct sum decomposition $X = X_0\oplus X_1$. Show that if $T_0:X_0\to Y$ and $T_1:X_1\to Y$ are bounded operators, then the operator $T:= T_0\oplus T_1$ from $X$ to $Y$ defined by
$$ T(x_0+x_1):= T_0x_0+T_1 x_1$$
is bounded. What can be said about the norm of $T$?

\noindent{\em Hint:}\ First show that $\nn x_0+x_1\nn := \n x_0\n+\n x_1\n$ is an equivalent norm on $X$.

\item\label{prob:wrong}
Let $X$ be a vector space.
\begin{enumerate}[\rm(a), leftmargin=*]
  \item\label{it:wrong1} Suppose that $\n \cdot \n$ and $\n \cdot \n'$ are two norms on $X$, each of which turns $X$ into a Banach space.
  Show that if there exists a constant $C\ge 0$ such that $\n x\n\le C\n x\n'$ for all $x\in X$, then the two norms are equivalent.
  \item\label{it:wrong2}
  Find the mistake in the following `proof' that any two norms
  $\n \cdot \n$ and $\n \cdot \n'$ turning $X$ into a Banach space are equivalent.
  Define
  $$ \nn x \nn  := \n x\n + \n x\n'\!, \quad x\in X.$$
  This is a norm which turns $X$ into a Banach space and we have $\n x\n \le \nn x \nn $ and $\n x\n' \le \nn x \nn $. Hence part \ref{it:wrong1}
  implies that $\n x\n$ and $\nn x \nn $ are equivalent and that
  $\n x\n'$ and $\nn x \nn $ are equivalent.
  It follows that $\n \cdot\n$ and $\n \cdot \n'$ are equivalent.
\end{enumerate}

\item
Using the open mapping theorem, show that there is no complete norm $\nn \cdot\nn$ on $C[0,1]$
with the property that\index{pointwise!convergence, vs. norm convergence}  $$ \nn f_n - f\nn \to 0 \ \Leftrightarrow f_n \to f \ \hbox{ pointwise}.$$

\item
Let $X$ be a Banach space. A sequence $(x_n)_{n\ge 1}$ in $X$ with $\n x_n\n=1$ for every $n\ge 1$ is called a {\em Schauder basis}\index{Schauder basis} for $X$ if every $x\in X$ admits a unique representation as a convergent sum $x = \sumn c_n x_n$ with $c_n\in \K$ for all $n\ge 1$.

Let $(x_n)_{n\ge 1}$ be a Schauder basis for $X$, and let $Y$ be the vector space of all scalar sequences $c = (c_n)_{n\ge 1}$ such that the sum $\sumn c_n x_n$ converges in $X$.
\begin{enumerate}[\rm(a), leftmargin=*]
\item Show that $$\n c\n_Y:= \sup_{N\ge 1} \Bigl\n \sum_{n=1}^N c_n x_n\Bigr\n$$ defines a norm on $Y$ and that $Y$ is a Banach space with respect to this norm.
 \item Show that $c \mapsto \sumn c_n x_n$ is an isomorphism from $Y$ onto $X$.
 \item Conclude that the {\em coordinate functionals}
 $$ P_k: \sumn c_n x_n \mapsto c_k$$
 are bounded and that $\sup_{k\ge 1}\n P_k\n<\infty$.
\end{enumerate}

\item\label{prob:surj-open}
Let $X$ and $Y$ be Banach spaces. The aim of this problem is to prove that the set of surjective operators is open in $\calL(X,Y)$.
\begin{enumerate}[\rm(a), leftmargin=*]
  \item\label{it:surj-open1}
  Let $T\in \mathscr{L}(X,Y)$ be a surjective operator. Show that there is a constant $A>0$ such that for all $y\in Y$ there exists an $x\in X$
  such that $\n x\n \le A\n y\n$ and $Tx=y$.
  \item\label{it:surj-open2}
  Let $T\in \mathscr{L}(X,Y)$ be a bounded operator. Suppose there exist constants
  $A\ge 0$ and $0\le B< 1$ such that for all $y\in Y$ with $\n y\n\le 1$ there exists an $x\in X$ such that
  $\n x\n \le A$ and $\n Tx-y\n \le B$. Show that $T$ is surjective.

  \noindent{\em Hint:}\ Consider the proof of Lemma \ref{lem:OMT}.

  \item\label{it:surj-open3} Show that if $T\in \mathscr{L}(X,Y)$ is surjective and $S\in \mathscr{L}(X,Y)$ is a
  bounded operator satisfying $\n S\n < 1/A$, where $A$ is the constant of part \ref{it:surj-open1},
  then $T+S$ is surjective.

  \noindent{\em Hint:}\ Use part \ref{it:surj-open1} and apply part \ref{it:surj-open2} with $B = A\n S\n$.
\end{enumerate}

\item
Let $X$ and $Y$ be Banach spaces and let $T:X\to Y$ be linear.
Suppose that $y^*\circ T\in X^*$ for every $y^*\in Y^*$.
Show that $T$ is bounded.

\noindent{\em Hint:}\ Use the closed graph theorem and the fact that
$Y^*$ separates the points of $Y$.

\item
Let $(a_n)_{n\ge 1}$ be a scalar sequence with the property that
the sum $\sum_{n\ge 1} a_n b_n$ converges for all scalar sequences $(b_n)_{n\ge 1}$ satisfying
$\sum_{n\ge 1} |b_n|^2 <\infty$.

\begin{enumerate}[\rm(a), leftmargin=*]
  \item Show that $\sum_{n\ge 1} a_n b_n$ converges absolutely
  for all scalar sequences $(b_n)_{n\ge 1}$ satisfying
  $\sum_{n\ge 1} |b_n|^2 <\infty$.
  \item Show that $\sum_{n\ge 1} |a_n|^2 <\infty$.

  \noindent{\em Hint:}\ Apply the closed graph theorem to the mapping $T: \ell^2 \to \ell^1$ defined by
  $T: (b_n)_{n\ge 1}  \mapsto (a_nb_n)_{n\ge 1}$.
\end{enumerate}

\item
Let $X$ be a Banach space, $(\Omega,\calF\!,\mu)$ be a measure space, let $1\le p\le \infty$, and suppose that $f: \Omega\to X$ is a function which has the property
that the scalar-valued function $\omega\mapsto \langle f(\omega),x^*\rangle$ belongs to $L^p(\Om)$ for every $x^*\in X^*$\!.

\begin{enumerate}[\rm(a), leftmargin=*]
  \item Show that the mapping $T: X^*\to L^p(\Omega)$ defined by $x^* \mapsto \langle f(\cdot),x^*\rangle$ is closed.
  \item Deduce that there exists a constant $C\ge 0$ such that
  $$ \n \langle f(\cdot),x^*\rangle\n_{L^p(\Omega)} \le C \n x^*\n, \quad x^*\in X^*\!.$$
\end{enumerate}

\item\label{prob:Pettis}
Let $(\Om,\calF\!,\mu)$ be a finite measure space and let $X$ be a Banach space.
\begin{enumerate}[\rm(a), leftmargin=*]
 \item\label{it:Pettis1}
Let $f:\Om\to X$ be a strongly measurable function. Show that if there is an exponent $1<p\le\infty$ such that
$\lb f(\cdot),x\s\rb\in L^p(\Om)$
for all $x\s\in X\s$\!, then there exists a unique element $x_f\in X$, the {\em Pettis integral}\index{Pettis
integral} of $f$ with respect to $\mu$, such that
$$ \lb x_f, x\s\rb = \int_\Om \lb f(\om),x\s\rb \ud\mu(\om), \quad x\s\in X\s\!.$$

\noindent{\em Hint:}\ The integrals
$\int_\Om \one_{\{\n f\n\le n\}} f\ud\mu$
are well defined as Bochner integrals.

\item\label{it:Pettis2} Show that the result of part \ref{it:Pettis1} fails for $p=1$.

\noindent{\em Hint:}\ Let $(A_n)_{n\ge 1}$ be a sequence of disjoint
intervals of positive measure $|A_n|$ in the interval $(0,1)$
and consider the function $f:(0,1)\to c_0$ defined by $$f(t) := \sumn \frac1{|A_n|}\one_{A_n}(t)  e_n, \quad t\in (0,1),$$
where $(e_n)_{n\ge 1}$ is the sequence of standard unit vectors in $c_0$.
\end{enumerate}

\item
Let $T\in\calL(X,Y)$, where $X$ and $Y$ are Banach spaces. Show that
$\Ran(T)$ is closed if and only if there exists a constant $C\ge0$
such that
$$
d(x,\Ker(T))\le C\n Tx\n,\quad x\in X.
$$
Deduce that if $T$ is injective, then $\Ran(T)$ is closed if and only
if $T$ is bounded below, that is, $\n Tx\n\ge c\n x\n$ for some
$c>0$ and all $x\in X$.

\item
Let $T\in\calL(H,K)$, where $H$ and $K$ are Hilbert spaces. Show that
$\Ran(T)$ is closed if and only if $\Ran(T^\star)$ is closed.
In this case prove that
$$
\Ran(T)=(\Ker(T^\star))^\perp,
\qquad
\Ran(T^\star)=(\Ker(T))^\perp.
$$

\item\label{prob:inclusion-of-ranges-Hilbert}
Let $H_1$ and $H_2$ be Hilbert spaces and $X$ be a Banach space.
Show that for $T_1\in \calL(H_1,X)$ and $T_2\in \calL(H_2,X)$ the following assertions are equivalent:
\begin{enumerate}[label={\rm(\arabic*)}, leftmargin=*]
  \item\label{it:inclusion-of-ranges-Hilbert1} $\Ran(T_1)\subseteq  \Ran(T_2)$;
  \item\label{it:inclusion-of-ranges-Hilbert2} there is a constant $C\ge 0$ such that $\n T_1\s x\s\n \le C\n T_2\s x\s\n$ for all $x\s\in X\s$.
\end{enumerate}

\noindent{\em Hint:}\ We may assume that $T_1$ and $T_2$ are injective. For \ref{it:inclusion-of-ranges-Hilbert1}$\Rightarrow$\ref{it:inclusion-of-ranges-Hilbert2}, apply the closed graph theorem to prove boundedness of $T_2^{-1}T_1: H_1\to H_2$. For the converse implication, use the Riesz representation theorem.

\item\label{prob:RanT-Hilbert}
Let $H$ be a Hilbert space. Show that a vector $y\in H$ belongs to the range of an operator $T\in \calL(H)$ if and only if there exists
a constant $C_y\ge 0$ such that
$$ |\iprod{x}{y}| \le C_y \n T^\star x\n, \quad x\in H.$$
{\em Hint:}\ Apply the result of Problem \ref{prob:inclusion-of-ranges-Hilbert} to the orthogonal projection onto the span of $y$.

\item\label{prob:Gauss}
It was shown in Lemma \ref{lem:Gauss} that for $\la>0$ the Fourier transform of $g(x)=(2\pi)^{-d/2}\exp(-\frac12\la|x|^2)$
is given by $$\wh{g}(\xi) = (2\pi\la)^{-d/2} \exp\Bigl(-\frac12 |\xi|^2/\la\Bigr).$$ Give an alternative proof of this identity by completing the following steps:
\begin{enumerate}[\rm(a), leftmargin=*]
  \item it suffices to prove the identity for $\la=1$ and in dimension $d=1$;
  \item the function $u(x):= (2\pi)^{-1/2}\exp(-\frac12 x^2)$ solves the differential equation $$u'(x) + xu(x) = 0;$$
  \item the Fourier transform of $u$ also satisfies this equation, with $\wh u(0)=u(0)$;
  \item apply the Picard--Lindel\"of theorem (Theorem \ref{thm:DV-Lip}).
\end{enumerate}

\item\label{prob:uniq-FT-T}
Write out a proof of Theorem \ref{thm:uniq-FT-T}.

\item\label{prob:FTapprox}
For $n \ge 1$ and $f \in L^2(\R^d)$ let
\begin{align*}
\mathscr{F}_n f(\xi) :=\frac1{(2\pi)^{d/2}} \int_{[-n,n]^d} \exp(- i x \cdot \xi) f(x)\ud x.
\end{align*}
Show that $\mathscr{F}_n$ maps $L^2(\R^d)$ into itself
and defines a bounded operator on $L^2(\R^d)$, and show that for all $f\in L^2(\R^d)$ we have
the following identity for the Fou\-rier--Plancherel transform:
\begin{align*}
\mathscr{F} f = \lim_{n \to \infty} \mathscr{F}_n f,
\end{align*}
where the limit is taken in $L^2(\R^d)$.

\item\label{prob:F4=I}
Consider the Fourier--Plancherel transform $\calF: f\mapsto \wh f$ on $L^2(\R^d)$.
\begin{enumerate}[\rm(a), leftmargin=*]
 \item Show that $\calF^2 = R,$ where $Rf(x) := f(-x)$ is the reflection operator on $L^2(\R^d)$.
 \item Deduce that $\calF^4 = I$.
\end{enumerate}

\item\label{prob:tensor-extension}
We use the notation introduced in Problem \ref{prob:LpX}.
Let $\calF$ be the Fourier--Plancherel transform on $L^2(\R)$ and let $X$ be a Banach space.
 On the space $L^2(\R)\otimes X$ we define the linear operator $\calF\otimes I$ by
$$ (\calF\otimes I)(f\otimes x):= (\calF f) \otimes x, \quad f\in L^2(\R),\ x\in X.$$

\begin{enumerate}[\rm(a), leftmargin=*]
  \item Show that this operator is well defined.
  \item Show that if $X = \ell^p$ with $1\le p\le \infty$, then $\calF\otimes I$ extends to a bounded operator on $L^2(\R;\ell^p)$ if and only if $p=2$.

  \noindent{\em Hint:}\ For $1\le p<2$ consider the functions $$f_N := \sum_{n=0}^N f(\cdot+n) \otimes e_{n+1},$$
  where $(e_n)_{n\ge 1}$ is the sequence of standard unit vectors in $\ell^p$ and $0\not=f\in C_{\rm c}(\R)$ has support in the interval $(-\pi,\pi)$; for $2< p\le \infty$ use the functions $$f_N := \sum_{n=0}^N e^{-in (\cdot)} f \otimes e_{n+1}.$$
\end{enumerate}

\item
Prove the estimate \eqref{eq:Dirichlet}.

\item\label{prob:H2-I}
Prove that the Hilbert transform $H$ on $L^2(\R)$ satisfies $H^2 = -I$.

\item\label{prob:harmonic}
This problem establishes a connection between the Hilbert transform and the theory of harmonic functions.

For real-valued functions $f\in L^2(\R)$ we define $u_f: \R\times (0,\infty)\to \R$ by
$$ u_f(x,y):= p_y *f(x), \quad x\in \R, \ y>0,$$
where $$ p_y(x):= \frac1\pi\frac{y}{x^2 + y^2}, \quad x\in \R, \ y>0,$$
is the {\em Poisson kernel}.\index{Poisson!kernel}

\begin{enumerate}[\rm(a), leftmargin=*]
  \item Show that $u_f$ is {\em harmonic}\index{harmonic},
  that is, $u_f\in C^2(\R\times (0,\infty))$ and $\Delta u_f \equiv 0$.
  \item Show that $u_f + iu_{Hf}$ is holomorphic.
\end{enumerate}

\item\label{prob:Poissker}
Let $f\in L^1(\R)$ satisfy $\wh f(-\xi) =0$ for almost all $\xi\ge 0$.
\begin{enumerate}[\rm(a), leftmargin=*]
  \item\label{it:Poissker1} Show that for all $y>0$ the function $p_y*f$, where $p_y$ is the Poisson kernel introduced in the preceding problem, is integrable and its Fourier transform belongs to $L^1(\R)\cap L^2(\R)$.

  \noindent{\em Hint:}\ Compute the Fourier transform of $p_y$.
  \item\label{it:Poissker2} Using Fourier inversion, prove that the function $$g(x+iy):= p_y*f(x)$$ is holomorphic on the open half-plane $\{x+iy \in \C:\,y >0\}$.

  \item\label{it:Poissker3} Using Theorem \ref{thm:approx-identity}, show that $$\lim_{y\downarrow 0} \n g(\cdot +iy) - f(\cdot)\n_{L^1(\R)} = 0.$$
\end{enumerate}
Somewhat informally, part \ref{it:Poissker3}
states that every $L^1$-function whose Fou\-rier transform vanishes on the negative half-line is the boundary value (in the $L^1$-sense) of a holomorphic function on the upper half-plane.
\begin{enumerate}[\rm(a), leftmargin=*]\setcounter{enumii}{3}
  \item State and prove a version of part \ref{it:Poissker3}
  for the disc.
\end{enumerate}

\item
Let $1\le p\le \infty$.
Young's inequality implies that the convolution of a function $f\in L^1(\R^d)$ with a function $g\in L^p(\R^d)$ belongs to $L^p(\R^d)$ and $\n f*g\n_p \le \n f\n_1\n g\n_p$.
\begin{enumerate}[\rm(a), leftmargin=*]
  \item Write out a proof by taking the first input exponent equal to $1$ and the output exponent equal to $p$ in the proof of Proposition \ref{prop:Young}.
\end{enumerate}
The special case of Young's inequality just stated can be reformulated as saying that for every $g\in L^p(\R^d)$ the convolution operator $C_g: f\mapsto f*g$ is bounded from $L^1(\R^d)$ to $L^p(\R^d)$
with norm $$\n C_g\n_{\calL(L^1(\R^d),L^p(\R^d))}  \le \n g\n_p.$$
\begin{enumerate}[\rm(a), leftmargin=*]\setcounter{enumii}{1}
  \item Let $\frac1p+\frac1q=1$. Using H\"older's inequality, show that convolution with $g$
defines a bounded operator from $L^q(\R^d)$ to $L^\infty(\R^d)$
of norm at most $\n g\n_p$, consistent with the operator above;
for $1<p\le\infty$, this is its unique bounded extension.
  \item Use the Riesz--Thorin interpolation theorem to obtain an alternative proof of the general form of Young's inequality.
\end{enumerate}

\item
Using the Riesz--Thorin interpolation theorem, give alternative proofs of:
\begin{enumerate}[\rm(a), leftmargin=*]
 \item the $L^p$-contractivity of conditional expectations of Proposition \ref{prop:conditional-expectations-Lp};
 \item the Schur estimate of Problem \ref{prob:Schur-estimate}.
\end{enumerate}

\item
Let $s_N$ be the Fourier partial sum operator of Problem~\ref{prob:unbbd-Four},
considered with respect to normalised Lebesgue measure ${\rm d}\theta/2\pi$ on $\T$.
\begin{enumerate}[\rm(a), leftmargin=*]
\item Show that $$\n s_N\n_{\calL(L^2(\T))}=1, \qquad
\n s_N\n_{\calL(L^1(\T),L^\infty(\T))}\le 2N+1.$$
\end{enumerate}
Now let $1\le p\le2$ and $1/p+1/q=1$.
\begin{enumerate}[\rm(a), leftmargin=*]\setcounter{enumii}{1}

\item Using the Riesz--Thorin interpolation theorem, prove that, for all $f\in L^p(\T)$,
$$
\n s_Nf\n_q\le (2N+1)^{2/p-1}\n f\n_p.
$$
\item Deduce the following special case of {\em Nikolskii's inequality}:\index{inequality!Nikolskii}
every
trigonometric polynomial $P$ of degree at most $N$ satisfies
$$
\n P\n_q\le (2N+1)^{2/p-1}\n P\n_p.
$$
\end{enumerate}

\item\label{prob:Clarkson}
Let $(\Om,\calF\!,\mu)$ be a $\sigma$-finite measure space and let $1<p,q<\infty$
satisfy $\frac1p+\frac1q=1$.
The aim of this problem is to use the Riesz--Thorin interpolation theorem
to derive the {\em Clarkson inequalities}:\index{inequality!Clarkson}

\begin{enumerate}[label={\rm(\arabic*)}, leftmargin=*]
\item\label{it:Clarkson1}
  if $1<p\le2$, then for all $f,g\in L^p(\Om)$ we have
  $$ (\n {f+g}\n_p^{p} +  \n {f-g}\n_p^{p})^{1/p} \le {2^{1/p}} (\n f\n_p^p + \n g\n_p^p)^{1/p}$$
  and
  $$(\n {f+g}\n_p^{q} +  \n {f-g}\n_p^{q})^{1/q} \le {2^{1/q}} (\n f\n_p^p + \n g\n_p^p)^{1/p};$$
  \item\label{it:Clarkson2}
  if $2\le p <\infty$, then for all $f,g\in L^p(\Om)$ we have
  $$ (\n {f+g}\n_p^{p} +  \n {f-g}\n_p^{p})^{1/p} \le {2^{1/q}} (\n f\n_p^{p} + \n g\n_p^{p})^{1/p}$$
  and
  $$ (\n {f+g}\n_p^{p} +  \n {f-g}\n_p^{p})^{1/p} \le {2^{1/p}} (\n f\n_p^{q} + \n g\n_p^{q})^{1/q}\!.$$
\end{enumerate}
Let $1\le r,s\le\infty$. On the product $L^r(\Om)\times L^r(\Om)$ we consider the norm
$$ \n (f,g)\n_{r,s} := (\n f\n_r^s + \n g\n_r^s)^{1/s}\!,$$
with the obvious modification for $s=\infty$.
\begin{enumerate}[\rm(a), leftmargin=*]
  \item\label{it:probClarkson1} Show that the resulting normed space $X_{r,s}(\Om)$ is complete.
\end{enumerate}
On $X_{r,s}(\Om)$ we consider the operator
$$ T: (f,g) \mapsto (f+g,f-g).$$
Its norm will be denoted by $\n T\n_{r,s}$.
\begin{enumerate}[\rm(a), leftmargin=*]\setcounter{enumii}{1}
  \item\label{it:probClarkson2} Show that $\n T\n_{1,1}\le 2$ and
$\n T\n_{2,2}\le\sqrt 2$.
  \item\label{it:probClarkson3} Prove the first inequality in \ref{it:Clarkson1}.

  \noindent{\em Hint:}\
Identify $X_{r,r}(\Om)$ with $L^r(\Om\times\{1,2\})$, using
counting measure on $\{1,2\}$, and interpolate the bounds
in part \ref{it:probClarkson2}.

  \item\label{it:probClarkson4} Show that for all $f,g\in L^1(\Om)$ we have
  $$\n T(f,g)\n_{1,\infty}\le \n (f,g)\n_{1,1}.$$
  \item\label{it:probClarkson5} Deduce the second inequality in \ref{it:Clarkson1}.

  \noindent{\em Hint:}\
Apply the Riesz--Thorin theorem over $\C$ to the operator
$$S(a,b):=(a+b,a-b),\quad (a,b)\in\C^2.$$
The endpoint bounds from $\ell^1_2$ to $\ell^\infty_2$ and
from $\ell^2_2$ to itself follow from the arguments of
parts \ref{it:probClarkson4} and \ref{it:probClarkson2} on a one-point probability space.
This gives
$$
\|S(a,b)\|_{\ell^q_2}
\le 2^{1/q}\|(a,b)\|_{\ell^p_2}.
$$
Apply this inequality pointwise and then use Minkowski's inequality to
pass from the pointwise $\ell^q_2$-norm to the norm of
$X_{p,q}(\Om)$.

 \item\label{it:probClarkson6} Prove both inequalities in \ref{it:Clarkson2}.

\noindent{\em Hint:}\ For the second inequality, interchange
$p$ and $q$ in the scalar inequality of the preceding hint
and apply Minkowski's inequality. Deduce the first inequality
using H\"older's inequality.

\item\label{it:probClarkson7} Prove that $L^p(\Om)$ is strictly convex.
\end{enumerate}

\item\label{prob:HY-FT-circle}
State and prove an analogue of the Hausdorff--Young theorem for the circle.

\item\label{prob:HT-weak11}
This problem sketches an alternative proof of the $L^p$-boundedness of the
Hilbert transform using the Marcinkiewicz interpolation theorem.

\begin{enumerate}[\rm(a), leftmargin=*]
 \item\label{it:HT-weak11-1}
 Let $f\in L^1(\R)\cap L^2(\R)$ and let $t>0$. Prove that there exist
 pairwise disjoint bounded intervals $(I_j)_{j\ge1}$ and a decomposition
 $$f=g+\sum_{j\ge1}b_j$$ such that
 $|g|\le 2t$ almost everywhere, $\n g\n_2^2\le 2t\n f\n_1$, $\supp(b_j)\subseteq I_j$, and
$$ \int_{I_j}b_j\ud x=0,
  \qquad
  \sum_{j\ge1}|I_j|\le \frac1t\n f\n_1,
\qquad \sum_{j\ge1}\n b_j\n_1\le 2\n f\n_1.$$

 \noindent{\em Hint:}\
 Use the maximal dyadic intervals $I$ for which
 $$\frac1{|I|}\int_I|f(x)|\ud x>t.$$
 Use the Lebesgue differentiation theorem outside their union.

 \item\label{it:HT-weak11-2}
 If $I_j^*$ is the interval with the same centre as $I_j$ and
 twice its length, prove that
 $$
 \int_{\R\setminus I_j^*}|Hb_j(x)|\ud x
 \le C\n b_j\n_1
 $$
 with an absolute constant $C$.

 \noindent{\em Hint:}\
 If $c_j$ is the centre of $I_j$, use $\int_{I_j}b_j\ud x=0$ to write
 $$
 Hb_j(x)
 =
 \frac1\pi\int_{I_j}
 \Bigl(\frac1{x-y}-\frac1{x-c_j}\Bigr)b_j(y)\ud y,
 \quad x\notin I_j^*.
 $$
 Then integrate the resulting estimate with respect to $x$.

 \item\label{it:HT-weak11-3}
 Deduce that the Hilbert transform on $L^2(\R)$ is of {\em weak type}
 $(1,1)$\index{weak type $(1,1)$} on $L^1(\R)\cap L^2(\R)$, in the sense that there is an absolute constant $C\ge0$ such that
 $$
 t|\{|Hf|>t\}|\le C\n f\n_1,
 \quad t>0,\ f\in L^1(\R)\cap L^2(\R).
 $$

 \noindent{\em Hint:}\
 Use the $L^2$-boundedness of $H$ for the good part $g$, and estimate
 the bad part outside $\bigcup_{j\ge1}I_j^*$ using part \ref{it:HT-weak11-2}.

 \item\label{it:HT-weak11-4}
 Use part \ref{it:HT-weak11-3}, the $L^2$-boundedness of the Hilbert
 transform, and the Marcin\-kiewicz interpolation theorem to obtain a
 bounded extension of $H$ to $L^p(\R)$ for every $1<p<2$.

 \item\label{it:HT-weak11-5}
 Using the identity
 $$
 \int_\R Hf\,g\ud x=-\int_\R f\,Hg\ud x,
 \quad f,g\in C_{\rm c}^2(\R),
 $$
 duality, and the result of part \ref{it:HT-weak11-4}, prove that $H$ extends
 boundedly to $L^p(\R)$ for every $2<p<\infty$.
\end{enumerate}

\end{problems}

%% file: ch06-Spectrum.tex
\chapter{Spectral Theory} \label{ch:spectral}

\blfootnote{\rm \small \begin{quote}This book has been published by Cambridge University Press in the series `Cambridge Studies in Advanced Mathematics'. The present version is free to view and download for personal use only. Not for re-distribution, re-sale or use in derivative works.\newline {\copyright} Jan van Neerven\end{quote}}

\noindent
Spectral theory is the branch of operator theory that seeks to extend the theory of eigenvalues to an infinite-dimensional setting.
Much of its power derives from the observation that, away from the spectrum of a bounded operator $T$, the operator-valued function $\la\mapsto (\la I-T)^{-1}$ is holomorphic. This makes it possible to import results from the theory of functions into operator theory. For instance, the fact that bounded operators on nonzero Banach spaces have nonempty spectra is deduced from
Liouville's theorem, and the Cauchy integral formula can be used to introduce a functional calculus
for functions holomorphic in an open set containing the spectrum of $T$.

\section{Spectrum and Resolvent}\label{sec:spectrum}

In linear algebra, a complex number $\lambda$ is said to be an {\em eigenvalue} of an $n\times n$
matrix $A$ with complex coefficients if there exists a nonzero vector $x\in \C^n$ such that $Ax = \la x$. The number
$\la$ is an eigenvalue if and only if $\la I-A$ fails to be invertible, or equivalently, if and only if $\det(\lambda I - A) = 0$. Writing out the determinant, we obtain the so-called characteristic polynomial in the variable $\lambda$,
which has $n$ zeroes (counting multiplicities) by the Fundamental Theorem of Algebra. Our first task will be to investigate to what extent these results generalise to bounded operators acting on a Banach space.

Throughout the chapter, $T$ denotes a bounded operator acting on a
complex Banach space $X$.

\begin{definition}[Spectrum and resolvent]\label{def:spectrum}
 The {\em resolvent set}\index{$R$@$\varrho(T)$}\index{resolvent!set} of an operator $T\in\calL(X)$ is the set $\varrho(T)$ consisting of all $\la\in\C$ for which
 the operator $\la I - T$ is {\em bounded\-ly invertible}\index{boundedly invertible}, by which we mean that there exists a bounded operator $U$ on $X$
 such that $$(\la I-T)U = U(\la I-T) = I.$$
 The {\em spectrum}\index{spectrum}\index{$S$@$\sigma(T)$} of $T$ is the complement of the resolvent set of $T$: $$\si(T):= \C\setminus\rh(T).$$
\end{definition}

From now on we will write $\la-T$ instead of $\la I - T$. It is customary to write\index{$Ra$@$R(\la,T)$}
$$ R(\la,T) := (\la-T)^{-1}$$
for the {\em resolvent operator}\index{resolvent!operator} of $T$ at the point $\la\in \varrho(T)$.
By the open mapping theorem (Theorem \ref{thm:OM}), a complex number $\la$ belongs to $\rh(T)$ if and only if $\la - T$ is a bijection on $X$.

\begin{example} The spectrum of an $n\times n$ matrix with complex coefficients, viewed as a bounded operator acting on $\C^n$\!, equals its set of eigenvalues.
\end{example}

In the present setting, a complex number $\la$ is said to be an {\em eigenvalue}\index{eigenvalue} of the bounded operator $T\in\calL(X)$ if $Tx = \la x$ for some nonzero vector $x\in X$; such a vector is then said to be an {\em eigenvector}\index{eigenvector}. The set $\sigma_{\rm p}(T)$\index{$S$@$\sigma_{\rm p}(T)$} of all eigenvalues of $T$ is called the {\em point spectrum}\index{point!spectrum} of $T$.
If $\lambda$ is an eigenvalue of $T$, then $\la-T$ is not injective and therefore not invertible. As a result, eigenvalues belong to the spectrum.
In contrast to the finite-dimensional situation, however, points in the spectrum need not be eigenvalues:

\begin{example}\label{ex:rightshift} The right shift $T$ on $\ell^2$\!, given by
$$ T: (c_1,c_2,\dots)\mapsto (0,c_1,c_2,\dots)$$
has no eigenvalues. Indeed, the identity $Tc=\la c$ can be written out as
$$ (0,c_1,c_2,\dots) =  (\la c_1,\la c_2,\dots).$$
If $\la\not=0$, comparison of the entries of these sequences inductively gives $c_n = 0$ for all $n\ge 1$.
If $\la = 0$, the identity reads $(0,c_1,c_2,\dots) =  (0,0,\dots)$ and again we obtain $c_n=0$ for all $n\ge 1$.
In both cases we find that $Tc = \la c$ only admits the zero solution $c=0$.

The spectrum of $T$ equals the closed unit disc: $\sigma(T) = \ov{\mathbb{D}}$. This can be proved directly (see Problem \ref{prob:rightshift-spectrum}) or by the following argument based on results proved below. By Proposition \ref{prop:spectrum-dual} we have $\sigma(T) = \sigma(T^*)$. The adjoint operator $T^*$ is readily identified as the left shift $(c_1,c_2,\dots)\mapsto (c_2,c_3,\dots)$. For each $\la\in \C$ with $|\la|<1$, the element
$(1,\la,\la^2,\dots)\in \ell^2$ is an eigenvector for this operator with eigenvalue $\la$. It follows that $\mathbb{D} \subseteq\sigma(T)$. Since by Lemma \ref{lem:res-open-bdd} the spectrum of a bounded operator is closed, this forces $\ov{\mathbb{D}} \subseteq \sigma(T)$.
On the other hand, by Lemma \ref{lem:Neumann}, the fact that $T$ is a contraction implies that $\sigma(T)\subseteq \ov{\mathbb{D}}$.
\end{example}

\begin{remark} In contrast to the case of
matrices in finite dimensions, the existence of a left inverse
does not imply the existence of a right inverse and vice versa. For example, the right shift on $\ell^2$
has a left inverse, the left shift, but not a right inverse; similarly the left shift on $\ell^2$
has a right inverse, namely the right shift, but not a left inverse. This is the reason for insisting on
the existence of a two-sided inverse in Definition \ref{def:spectrum}. If both a left inverse $U_{\rm l}$ and
a right inverse $U_{\rm r}$ exist, then necessarily $U_{\rm l}=U_{\rm r}$ and this operator is a two-sided inverse.
\end{remark}

If the bounded operators $T$ and $U$ are boundedly invertible, then so is the operator
$TU$ and we have $(TU)\inv = U\inv T\inv\!$.
Invertibility of $TU$ by itself does not imply invertibility
of $T$ or $U$; for instance, if $T$ is the left shift and $U$
the right shift on $\ell^2$, then $TU=I$, although neither
operator is invertible.

Our first main result, Theorem \ref{thm:siTnonempty}, asserts that the spectrum of a bounded operator acting on a nonzero Banach space is always a {\em nonempty compact} subset of the complex plane. This will be deduced from a series of lemmas.

\begin{lemma}[Neumann series]\label{lem:Neumann}\index{Neumann!series}
 If $\n T\n < 1$, then $I-T$ is boundedly invertible and its inverse is given by the absolutely convergent series
 $$ (I - T)\inv = \sum_{n=0}^\infty T^n\!.$$
 As a consequence, the spectrum of a bounded operator $T$ is contained in the closed disc $\{z\in\C:\, |z|\le \n T\n\}$.
\end{lemma}
\begin{proof} The
absolute convergence in $\calL(X)$ of the series follows from
 $\sum_{n=0}^\infty \n T^n\n \le \sum_{n=0}^\infty \n T\n ^n < \infty$. By the completeness of $\calL(X)$,
 the series $\sum_{n= 0}^\infty T^n$ converges in $\calL(X)$.
The identity in the statement of the lemma is a consequence of the identities
$$ (I-T)\sum_{n=0}^N T^n = \sum_{n=0}^N T^n(I-T) = I - T^{N+1}\!,$$
valid for all $N\ge 1$.
Upon letting $N\to\infty$ they give
$$ (I-T)\sum_{n=0}^\infty T^n = \sum_{n=0}^\infty T^n(I-T) = I, $$
which means that the bounded operator $\sum_{n=0}^\infty T^n $ is a two-sided inverse for  $I-T$.

To prove the second assertion, let $T\in \calL(X)$ be an arbitrary bounded operator. By the first assertion,
for all $\la\in \C$ with $|\la|> \n T\n$ the operator $\la - T = \la (I - T/\la)$ is boundedly invertible.
\end{proof}

\Needspace{6\baselineskip}
As an application we prove that the spectrum is always a closed subset of $\C$.

\begin{lemma}\label{lem:res-open-bdd} The spectrum $\sigma(T)$ is a closed subset of $\C$. If $X\not=\{0\}$ and $\la\in\varrho(T)$, then
 $\varrho(T)$ contains the open ball with centre $\la$ and radius $r = 1/\n R(\la,T)\n$. If $|\la-\mu|\le \delta r$ with $0\le \delta<1$,
 then $$\n R(\mu,T)\n \le \frac{1}{1-\delta}\n R(\la,T)\n.$$
\end{lemma}
\begin{proof}
If $X=\{0\}$, then $\si(T)=\emptyset$ and there is nothing to prove, so we may assume that $X\not=\{0\}$.
If $S\in \calL(X)$ is boundedly invertible and $U\in \calL(X)$ has norm $\n U\n \le \delta r$ with $r := 1/\n S^{-1}\n$,
then $\n  S\inv U\n \le \delta<1$ and therefore, by Lemma \ref{lem:Neumann},
$ S-U = S(I-S\inv U)$ is boundedly invertible and
\begin{align*}
 \n( S-U)\inv  \n \le \n S\inv \n \n (I-S\inv U)\inv \n
  = \n S\inv\n\Big\n \sum_{n=0}^\infty ( S\inv U)^n\Big\n\le \frac1r\sum_{n=0}^\infty \delta^n
  = \frac1r
 \frac{1}{1-\delta}.
\end{align*}
Now take $U = (\lambda-\mu)I$ and $S = \lambda-T.$
\end{proof}

\begin{lemma}[Resolvent identity]\label{lem:resolvent-identity}\index{resolvent!identity}
For all $\la,\mu\in\varrho(T)$ we have $$R(\la,T) - R(\mu,T) =  (\mu-\la)R(\la,T)R(\mu,T).$$
\end{lemma}
\begin{proof}
Multiply both sides by the invertible operator $(\mu-T)(\la-T)$.
\end{proof}

\begin{definition}[Holomorphy]
 Let $\Omega$ be an open subset of $\C$. A function $f:\Omega\to X$ is {\em holomorphic} if for all $z_0\in \Omega$ the limit
$$ \lim_{z\to z_0} \frac{f(z)-f(z_0)}{z-z_0}$$
exists in $X$.
\end{definition}

\begin{lemma}\label{lem:spectr-holo-bdd} The function $\la\mapsto R(\la,T)$ is holomorphic on $\rh(T)$
and satisfies $$\lim_{|\la|\to\infty} \n R(\la,T)\n = 0.$$
\end{lemma}
\begin{proof}
Continuity of the mapping  $\la\mapsto R(\la,T)$ follows from the resolvent identity and the bound in Lemma \ref{lem:res-open-bdd}.
To prove the holomorphy of $\la\mapsto R(\la,T)$ on $\rh(T)$ we use the resolvent identity and the continuity of $\la\mapsto R(\la,T)$ to obtain
\begin{align*}\lim_{\mu\to\la} \frac{R(\la,T) - R(\mu,T)}{\la-\mu}
= - \lim_{\mu\to\la}R(\la,T)R(\mu,T)  = -(R(\la,T))^{2}\!.
\end{align*}

For  $|\la| > 2\n T\n$ the Neumann series gives
 $$\n R(\la,T)\n = |\la|\inv \n (I - \la\inv T)\inv\n \le|\la|\inv \sum_{n=0}^\infty\n\la\inv T\n^n
\le 2|\la|\inv\!.$$
This proves the second assertion.
\end{proof}

We are ready for the first main result of this chapter:

\begin{theorem}[Nonemptiness of the spectrum]\label{thm:siTnonempty} If $T$ is a bounded operator on a non\-zero Banach space, then $\sigma(T)$ is a nonempty compact subset of the closed disc $\{\la\in\C:\, |\la|\le \n T\n\}$.
\end{theorem}

\begin{proof}
Containment in $\{\la\in\C:\, |\la|\le \n T\n\}$ and closedness of the spectrum have already been proved in Lemmas \ref{lem:Neumann} and \ref{lem:res-open-bdd}, respectively. Since bounded closed subsets of $\C$ are compact, this gives the compactness of $\sigma(T)$.

 Suppose, for a contradiction, that $\sigma(T)=\emptyset$. Then the function
 $\la\mapsto R(\la,T)$ is holomorphic on $\C$. By Lemma \ref{lem:spectr-holo-bdd} it is also bounded.
 Now we are in a position to apply Liouville's theorem: for all $x\in X$ and $x\s\in X\s$ we find that
 $\la\mapsto  \lb R(\la,T)x,x\s\rb $ is constant. Its limit for $|\la|\to\infty$ is zero, and
 therefore $\lb R(\la,T)x,x\s\rb = 0$ for all $\la\in\rh(T)$ and all $x\in X$ and $x\s\in X\s$\!.
 By the Hahn--Banach theorem, $R(\la,T)x =0$ for all $\la\in\rh(T)$ and all $x\in X$.
 It follows that $R(\la,T) =0$ for all $\la\in\rh(T)$. This implies $X = \{0\}$.
\end{proof}

Instead of using duality to reduce matters to scalar-valued functions, one may note that the proof of Liouville's theorem generalises {\em mutatis mutandis} to holomorphic functions with values in a Banach space.

We have seen in Lemma \ref{lem:spectr-holo-bdd} that the resolvent $\la\mapsto R(\la,T)$ is
holomorphic on $\varrho(T)$. The next result shows that the topological boundary
$\partial\varrho(T):= \ov{\varrho(T)}\setminus \varrho(T)$ is a natural barrier to the holomorphic continuation of this function.

\begin{proposition}\label{prop:res-blowup-bdd}
 If $\la_n\to \la$ in $\C$, with each $\la_n\in \varrho(T)$ and $d(\la_n,\sigma(T))\to 0$, then
 $$\limn \n R(\la_n,T)\n = \infty.$$
\end{proposition}

\begin{proof}
By Lemma \ref{lem:res-open-bdd}, if $\mu\in \varrho(T)$, then the open ball $B(\mu;\n R(\mu,T)\n^{-1})$ is contained in $\varrho(T)$. This implies
the more precise assertion that for all $\mu\in \varrho(T)$ we have $d(\mu,\sigma(T)) \ge \n R(\mu,T)\n^{-1}$, that is,
$\n R(\mu,T)\n\ge 1/d(\mu,\si(T))$.
\end{proof}

An immediate application is the following analytic continuation result.

\begin{corollary}
If $D\subseteq\C$ is a connected open set  intersecting the resolvent set of a bounded operator $T$, and if $\la\mapsto R(\la,T)$ extends holomorphically to $D$, then
$D \subseteq \varrho(T)$.
\end{corollary}

A more substantial application of Proposition \ref{prop:res-blowup-bdd} is the following result about the spectra of isometries.
Recall that an {\em isometry} is an operator $T\in \calL(X,Y)$ such that $\n Tx\n = \n x\n$ for all $x\in X$.

\begin{corollary}\label{cor:spectrum-isometry}\index{isometry!spectrum of} If $X \not=\{0\}$, the spectrum of an isometry $T\in\calL(X)$ is either contained in the unit circle $\mathbb{T}$ or else equals the closed unit disc $\ov{\mathbb{D}}$.
\end{corollary}
\begin{proof}
Since $X\not=\{0\}$, an isometry $T$ on $X$ has norm $\n T\n =1$ and nonempty spectrum $\sigma(T)$, which is contained in $\ov{\mathbb{D}}$ by the second assertion of
Lemma \ref{lem:Neumann}.

If $\mu\in \varrho(T)$ with $0< |\mu|<1$, then, since $T$ is an isometry, for all $x\in X$ we have
\begin{align*} \n x\n  = \n(\mu-T)R(\mu,T)x\n \ge \n TR(\mu,T)x \n - \n \mu R(\mu,T)x\n  = (1-|\mu|)\n R(\mu,T) x\n
\end{align*}
and therefore
$ \n R(\mu,T) \n \le 1/(1-|\mu|).$
In view of Proposition \ref{prop:res-blowup-bdd}, this implies that for every $0<r<1$ the closed disc $\{\mu\in\C:\, |\mu|\le r\}$ contains no boundary points of $\sigma(T)$. Therefore the set $\sigma(T)\cap\mathbb D$ has empty boundary in the connected open set $\mathbb D$.
It follows that either $\sigma(T)\cap\mathbb D=\emptyset$ or
$\sigma(T)\cap\mathbb D=\mathbb D$. In the first case we have
$\sigma(T)\subseteq\mathbb T$. In the latter case, since $\sigma(T)$ is
closed and contained in $\ov{\mathbb D}$, we have $\sigma(T)=\ov{\mathbb D}$.
\end{proof}

This result is the best possible in the following sense: both the closed unit disc and
any nonempty closed subset of the unit circle can be realised as the spectrum of suitable isometries.
For instance, the right shift on $\ell^2$ provides an example of the former, and if $K\subseteq \mathbb{T}$
is a nonempty closed set, then the isometry on $C(K)$ given by $(Tf)(z) = zf(z)$ is easily verified to have spectrum equal to $K$.

We have the following continuity result for spectra:

\begin{proposition}[Upper semicontinuity of the spectrum]\label{prop:spect-cont}\index{upper!semicontinuity, of the spectrum}
 Let $\Om$ be an open set in the complex plane containing $\si(T)$.
 Then there exists a $\delta>0$ such that if the bounded operator $T'$ satisfies $\n T - T'\n < \delta$, then $\sigma(T')\subseteq \Om$.
\end{proposition}

\begin{proof}
If $X=\{0\}$ or $\Om=\C$, there is nothing to prove, so we may assume that
$X\not=\{0\}$ and $\Om\neq\C$.

In the proof of Lemma \ref{lem:spectr-holo-bdd} we have seen that
$\lim_{|\la|\to\infty} \n R(\la,T)\n = 0$.
In particular, $\sup_{|\la| > \max\{1,2\n T\n\}} \n R(\la,T)\n < \infty$. By the continuity of $\la\mapsto R(\la,T)$ we also
have $\sup_{\{|\la|\le \max\{1,2\n T\n\}\}\cap \{\la\not\in \Om\}} \n R(\la,T)\n < \infty$.
Combining these, we find that $$\sup_{\la \not\in  \Om} \n R(\la,T)\n < \infty.$$
Denote this supremum by $M$. If $\n T-T'\n<1/M$ and $\la\not\in\Om$, then from
$$ \la-T' = (\la-T)[I + R(\la,T) (T-T')]$$
we infer that $\la - T'$ is invertible, noting that $I+R(\la,T) (T-T')$ is invertible since
$\n R(\la,T) (T-T')\n < M\cdot 1/M = 1$.
\end{proof}

It has already been noted that eigenvalues belong to the spectrum, and that a bounded operator need not have any eigenvalues (see Example \ref{ex:rightshift}). We now prove a useful result that makes up for this to some extent.

\begin{definition}[Approximate point spectrum]\label{def:approx-eigenvalue}
 The number $\la\in \C$ is
 called an {\em approximate eigenvalue}\index{approximate!eigenvalue} of the operator $T$ if there exists a sequence
 $(x_n)_{n\ge 1}$ in $X$ with the following two properties:
 \begin{enumerate}[label={\rm(\roman*)}, leftmargin=*]
  \item $\n x_n\n = 1$ for all $n\ge 1$;
  \item $\n Tx_n - \la x_n\n\to 0$ as $n\to \infty.$
 \end{enumerate}
In this context the sequence $(x_n)_{n\ge 1}$ is called an {\em approximate eigensequence} for $\la$.\index{approximate!eigensequence} The set of all approximate eigenvalues is called the {\em approximate point
spectrum}.\index{approximate!point spectrum}
\end{definition}

Every eigenvalue is an approximate eigenvalue.
Approximate eigenvalues belong to the spectrum, for if $\la$ were an approximate eigenvalue belonging to the resolvent set of $T$, we would arrive at the contradiction $$ 1 = \n x_n \n = \n R(\la,T) (Tx_n - \la x_n)\n \le \n R(\la,T)\n \n Tx_n - \la x_n\n \to 0 \ \ \hbox{as} \ n\to\infty.$$

\begin{proposition}\label{prop:approx-eigenvalue-bdry} The boundary of $\sigma(T)$ consists of approximate eigenvalues.
\end{proposition}
\begin{proof}
Let $\lambda\in\partial\sigma(T)$, and choose a sequence
$(\lambda_n)_{n\ge1}$ in $\varrho(T)$ with $\lambda_n\to\lambda$.
By Proposition \ref{prop:res-blowup-bdd} we have $\|R(\lambda_n,T)\|\to\infty$, so for each $n\ge 1$ we may choose $y_n\in X$ with $\|y_n\|=1$ and
$$
        \|R(\lambda_n,T)y_n\|\ge \frac12\|R(\lambda_n,T)\|.
$$
Set
$$
        x_n:=
        \frac{R(\lambda_n,T)y_n}{\|R(\lambda_n,T)y_n\|}.
$$
Then $\|x_n\|=1$, and
\begin{align*}
\|Tx_n-\lambda x_n\|
&=
\frac{\|[(T-\lambda_n)+(\lambda_n-\lambda)]R(\lambda_n,T)y_n\|}
        {\|R(\lambda_n,T)y_n\|}
\le
\frac{\|y_n\|}{\|R(\lambda_n,T)y_n\|}
+
|\lambda_n-\lambda|,
\end{align*}
and the right-hand side tends to $0$ as $n\to\infty$. Thus $\lambda$ is an approximate eigenvalue.
\end{proof}

We have the following duality result:

\begin{proposition}\label{prop:spectrum-dual}
A bounded operator $T$ and its adjoint $T\s$ have the same spectrum,
$$\sigma(T\s) = \si(T).$$
\end{proposition}
\begin{proof}
If $\la\in\rh(T)$, then
$$ (\la-T\s)[R(\la,T)]\s = [R(\la,T) (\la-T)]\s = I_X\s = I_{X\s},$$ and
similarly $[R(\la,T)]\s (\la-T\s) = I_{X\s}$, from which it follows that $\la\in \rh(T\s)$ and $R(\la,T\s) = [R(\la,T)]\s$\!.
This proves the inclusion $\sigma(T\s)\subseteq \si(T)$.

Conversely, suppose that $\lambda\in\varrho(T^*)$. Since
$\ker(\lambda-T^*)=\{0\}$, the range of $\lambda-T$ is dense in $X$ by Proposition \ref{prop:dual-dense-T}.
Moreover, $\lambda-T$ is bounded below. Indeed, for $x\in X$,
\begin{align*}
\|x\| = \sup_{\|x^*\|\le1}|\langle x,x^*\rangle|
  & = \sup_{\|x^*\|\le1}
    |\langle x,(\lambda-T^*)R(\lambda,T^*)x^*\rangle| \\
  & = \sup_{\|x^*\|\le1}
    |\langle(\lambda-T)x,R(\lambda,T^*)x^*\rangle|
   \le \|R(\lambda,T^*)\|\,\|(\lambda-T)x\|.
\end{align*}
Thus $\lambda-T$ has closed range. Its range is therefore both dense and
closed, hence equal to $X$. Since $\lambda-T$ is also bounded below, it is
injective, and its inverse is bounded. Hence $\lambda\in\varrho(T)$.
This proves $\sigma(T)\subseteq\sigma(T^*)$.
\end{proof}

We continue with an observation about relative spectra. Let $\mathscr{A}\subseteq \calL(X)$
be a {\em unital closed subalgebra},\index{subalgebra} that is, $\calA$ is a closed subspace of $\calL(X)$
closed under composition and containing the identity operator $I$.
For an operator $T\in \calA$ we define $\varrho_\calA(T)$ as the set of all $\la\in\C$ for which
the operator $\la - T$ is boundedly invertible in $\calA$, that is, there exists an operator $U\in\calA$ such that
$(\la-T)U = U(\la-T) = I$. We further set $\si_\calA(T):= \C\setminus\rh_\calA(T).$
Repeating the proofs of Lemmas \ref{lem:Neumann} and \ref{lem:res-open-bdd}, we see that $\rh_\calA(T)$ is an open set and $\si_\calA(T)$ is a closed set
contained in the closed disc of radius $\n T\n$, and therefore $\si_\calA(T)$ is compact.

It is evident that $\rh_\calA(T)\subseteq \rh(T)$ and therefore $$\si(T)\subseteq \si_\calA(T).$$ The next result provides a partial converse.

\begin{proposition}\label{prop:siAT-bdry}
Let $\calA\subseteq \calL(X)$ be a unital closed subalgebra and let $T\in \calA$. Then
$$\partial \si_\calA(T) \subseteq \si(T).$$
\end{proposition}
\begin{proof} Let $\la\in \partial \si_\calA(T)$ and let $\la_n \to \la$ in $\C$ with each $\la_n$ in $\rh_\calA(T)$.
Repeating the proof of Proposition \ref{prop:res-blowup-bdd}, we find that $\n R(\la_n,T)\n \to \infty$ as $n\to \infty$.
Since $\rh(T)$ is open and the resolvent $\la\mapsto R(\la,T)$ is continuous with respect to the operator norm, this implies that $\la\in\si(T)$.
\end{proof}

The example of the right shift $T$ on $X = \ell^2(\Z)$ and the unital closed subalgebra $\calA$ generated by the identity and $T$ shows that Proposition \ref{prop:siAT-bdry} is the best possible: here one has $\partial\si_\calA(T) = \si(T) = {\mathbb{T}}$ and $\si_\calA(T) = \ov{\mathbb{D}}$.

\section{The Holomorphic Functional Calculus}\label{sec:DC}

If $f:\C\to \C$ is an entire function and $T$ is a bounded operator on $X$,
we may define a bounded operator $f(T)$ on $X$ as follows. Writing $f$ as a convergent power series about $z=0$,
$$ f(z) = \sum_{n=0}^\infty \frac{f^{(n)}(0)}{n!} z^n\!,$$
we define
$$ f(T):= \sum_{n=0}^\infty \frac{f^{(n)}(0)}{n!} T^n\!.$$
This series converges absolutely in $\calL(X)$ since the same is true for the power series of $f(z)$ for every $z\in\C$.
The mapping $f\mapsto f(T)$ is called the {\em entire functional calculus}\index{functional calculus!entire} of $T$ and has the following properties, the easy proofs of which we leave to the reader:

\begin{enumerate}[label={\rm(\roman*)}, leftmargin=*]
 \item if $f(z) = z^n$ with $n\in\N$, then $f(T) = T^n$;
 \item $f(T)g(T) = (fg)(T)$;
 \item $g(f(T)) = (g\circ f)(T)$.
\end{enumerate}
This calculus may be used to define operators such as
$\exp(T)$, $\sin(T)$, $\cos(T)$, and so forth. There is a beautiful way to extend the entire functional calculus to a larger class of holomorphic functions,
replacing
power series expansions by the Cauchy integral formula
\begin{align}\label{eq:Dunf-calc-scalar} f(z_0) = \frac1{2\pi i} \int_{\Gamma} \frac{f(\la)}{\la-z_0}\ud \la.
\end{align}
Here, $f$ is assumed to be a holomorphic function on an open set $\Omega$ in the complex plane containing $z_0$, and $\Gamma$ is a suitable contour in $\Om$ winding once anticlockwise around $z_0$. Formally substituting $T$ for $z_0$ and interpreting
$1/(\la-T)$ as $R(\la,T)$, one is led to conjecture that a holomorphic functional calculus may be defined by the formula
\begin{align}\label{eq:Dunf-calc}
f(T) := \frac1{2\pi i} \int_{\Gamma} f(\la)R(\la,T)\ud \la.
\end{align}
Since $\la\mapsto R(\la,T)$ is continuous on $\varrho(T)$ with respect to the operator norm of $\calL(X)$, once $\Gamma$ has been parametrised,
the integral in \eqref{eq:Dunf-calc} is well defined as a Riemann integral with values in $\calL(X)$ (see Section \ref{subsec:integral-Riemann}).

In order to formulate a set of conditions on $\Omega$, $\Gamma$\!, and $T$
under which this idea works, we first take a closer look at the precise assumptions in the Cauchy integral formula \eqref{eq:Dunf-calc-scalar}. These are that
$\Omega$ is an open set in the complex plane containing $z_0$ and $\Gamma$ is a piecewise continuously differentiable closed contour in $\Om\setminus\{z_0\}$ with the following two properties:
\begin{enumerate}[label={\rm(\roman*)}, leftmargin=*]
 \item\label{it:winding1} the winding number of $\Gamma$ about the point $z_0$ equals $1$;
 \item\label{it:winding2} the winding number of $\Gamma$ about every point in
 $\complement\Omega$ equals $0$.
\end{enumerate}
Here, the {\em winding number}\index{winding number} of $\Gamma$ about a point $z_0$ is the (integer) number
$$ w(\Gamma;z_0):= \frac{1}{2\pi i} \int_\Gamma \frac1{\la -z_0}\ud \la.$$
More generally, we can admit finite unions of such contours, as long as their union satisfies \ref{it:winding1} and \ref{it:winding2}. If $\Gamma = \Gamma_{\!1}\cup\cdots\cup\Gamma_{\!k}$ is such a union, we define
$$ w(\Gamma;z_0):= \frac{1}{2\pi i} \sum_{j=1}^k\int_{\Gamma_{\!j}} \frac1{\la -z_0}\ud \la.$$
Condition \ref{it:winding2} is satisfied if $\Om = \Om_1\cup\cdots\cup\Om_k$ is a finite union of disjoint convex sets and $\Gamma_{\!j}$ is a piecewise continuously differentiable closed contour in $\Om_j\setminus\{z_0\}$.

Turning to the discussion of \eqref{eq:Dunf-calc}, we need to fix a similar set of assumptions regarding $\Omega$ and $\Gamma$,
while letting the operator $T$ play the role of $z_0$. We require that $\Omega$ is an open set in the complex plane containing $\sigma(T)$ and $\Gamma$ is a piecewise continuously differentiable closed contour in $\Omega\setminus\sigma(T)$ with the following two properties:
\begin{enumerate}[label={\rm(\roman*)}, leftmargin=*]
 \item\label{it:windingT1} the winding number of $\Gamma$ about every point $z_0\in \sigma(T)$ equals $1$;
 \item\label{it:windingT2} the winding number of $\Gamma$ about every point in $\complement\Omega$ equals $0$.
\end{enumerate}

It is easy to see that such contours always exist. If the conditions \ref{it:windingT1} and \ref{it:windingT2} are met, we say that $\Gamma$ is an {\em admissible contour}\index{admissible contour} for $\si(T)$ in $\Omega$.
As before we admit the possibility that $\Gamma$ is a finite union of such contours; for such contours
we interpret \eqref{eq:Dunf-calc} as
\begin{align*}
f(T) := \frac1{2\pi i}\sum_{j=1}^k \int_{\Gamma_{\!j}} f(\la)R(\la,T)\ud \la.
\end{align*}
For example, if $\sigma(T)=K_1\cup K_2$ is the union of two disjoint compact sets with $K_j\subseteq \Om_j$, where $\Om:=\Om_1\cup \Om_2$ is a disjoint union of open convex sets, we may select contours $\Gamma_{\!j}$ with winding number $1$ about every point in $K_j$ and consider $\Gamma = \Gamma_{\!1}\cup\Gamma_{\!2}$ as an admissible contour for $\si(T)$ in $\Om$. This example is relevant for defining spectral projections, where one uses the holomorphic functions $f = \one_{\Om_1}$
and $f = \one_{\Om_2}$ (see Theorem \ref{thm:spec-proj}).

\begin{figure}
\vskip-.8cm
\begin{center}
\begin{tikzpicture}[scale=0.7]
\filldraw [lightgray](0,0) .. controls (3,3) and (4,1) .. (3,0);
\filldraw [lightgray](3,0) .. controls (2,-1) and (-2,-2) .. (0,0);

\filldraw[thick,  nearly transparent] (-1,0) .. controls (3.5,5.3) and (5,2.6) .. (3.5,0);
\filldraw[thick,  nearly transparent] (3.5,0) .. controls (2,-3) and (-3.5,-3) .. (-1,0);

\draw[<-] (-0.5,0.2) .. controls (3.5,4) and (4,2) .. (3.4,0.4);
\draw (3.4,0.4) .. controls (2.5,-2) and (-3,-2) .. (-0.5,0.2);

\filldraw [lightgray](-3,2) .. controls (-6,-1) and (-4,-2) .. (-3,-1);
\filldraw [lightgray](-3,-1) .. controls (-2,0) and (-1,4) .. (-3,2);

\filldraw[thick,  nearly transparent] (-6,0) .. controls (-1.5,5) and (0,3) .. (-1.5,0);
\filldraw[thick,  nearly transparent] (-1.5,0) .. controls (-3,-3) and (-8.5,-3) .. (-6,0);

\draw[<-] (-4.8,0) .. controls (-1.6,5) and (-1,2) .. (-1.8,0);
\draw (-1.8,0) .. controls (-2.5,-2) and (-6.5,-3) .. (-4.8,0);

\draw (-3.5,0) node {$K_1$};
\draw (1.5,0) node {$K_2$};
\draw (-5.5,-1.7) node {$\Gamma_{\!1}$};
\draw (-1,-1.7) node {$\Gamma_{\!2}$};
\draw (-4.5,2.3) node {$\Omega_1$};
\draw (0.7,2.3) node {$\Omega_2$};
\end{tikzpicture}
\caption{The contour $\Gamma = \Gamma_{\!1}\cup\Gamma_{\!2}$ around $\sigma(T) = K_1\cup K_{2}$ in $\Om = \Om_1\cup\Om_2$\label{fig:sigmaT}}
\end{center}
\end{figure}

For an open set $\Omega$ in the complex plane we denote by $H(\Om)$\index{$H$@$H(\Om)$} the vector space of all holomorphic functions on $\Om$.

\begin{theorem}[Holomorphic functional calculus]\index{functional calculus!holomorphic}\label{thm:Dunford}
Let $T\in \calL(X)$ be a bounded operator, and let $\Om\subseteq\C $ be open and contain $\si(T)$.
For $f\in H(\Om)$ we define
$$ f(T) := \frac1{2\pi i} \int_{\Gamma} f(\la)R(\la,T)\ud \la,$$
where $\Gamma$ is an admissible contour for $\si(T)$ in $\Omega$.
The resulting operators $f(T)$ are well defined and have the following properties:
\begin{enumerate}[label={\rm(\roman*)}, leftmargin=*]
 \item\label{it:Dunford0} the operator $f(T)$ is independent of $\Om$ and the admissible contour $\Gamma$ in $\Om$;
 \item\label{it:Dunford1} for entire functions the holomorphic and entire functional calculi agree;
 \item\label{it:Dunford2} for all $f,g\in H(\Om)$ we have $f(T)g(T) = (fg)(T);$
 \item\label{it:Dunford4} for all $f\in H(\Om)$ we have $f(T^*) = (f(T))^*$;
 \item\label{it:Dunford3} for $f_\mu(\la) = 1/(\mu-\la)$ with $\mu\in \varrho(T)$ we have $f_\mu(T) = R(\mu,T).$
\end{enumerate}
In part \ref{it:Dunford3}, it is understood that the open set $\Om$ containing $\si(T)$ is chosen so that $\mu\not\in\Om$.
\end{theorem}
\begin{proof}

\ref{it:Dunford0}: \ This follows by applying the Cauchy theorem to the scalar-valued integrals
$$ \frac1{2\pi i} \int_{\Gamma} f(\la)\lb R(\la,T)x,x^*\rb\ud \la$$
and then using the Hahn--Banach theorem.
Alternatively, one may extend the proof of the Cauchy theorem
{\em mutatis mutandis} to functions with values in a Banach space.

\smallskip
\ref{it:Dunford1}:\ Let $f$ be entire. By part
\ref{it:Dunford0}, we may take $\Gamma$ to be the circle
centred at the origin with radius $r>\n T\n$, oriented
anticlockwise. For $\la\in\Gamma$, the Neumann series gives
$R(\la,T)=\sum_{n=0}^\infty \la^{-n-1}T^n$,
with absolute convergence uniform in $\la\in\Gamma$.
We may therefore integrate term by term to obtain
\begin{align*}
f(T)
=\frac1{2\pi i}\int_\Gamma f(\la)R(\la,T)\ud\la
=\sum_{n=0}^\infty
  \Bigl(\frac1{2\pi i}\int_\Gamma f(\la)\,\frac{\!\ud\la}{\la^{n+1}}\Bigr)T^n
=\sum_{n=0}^\infty\frac{f^{(n)}(0)}{n!}T^n.
\end{align*}
The last equality follows from the Cauchy integral formula
for the derivatives of $f$. The resulting series is precisely
the definition of $f(T)$ in the entire functional calculus.

\smallskip
\ref{it:Dunford2}: \ Let $\Gamma$ be an admissible contour for $\si(T)$ in $\Om$, and let $\Gamma'$ be admissible for $\si(T)$ in $\{z\in \Om\setminus \Gamma:\, w(\Gamma;z) =1\}$. Then $\Gamma$ winds once around $\Gamma'$, while $\Gamma$ is exterior to $\Gamma'$.
By the resolvent identity and Fubini's theorem,
\begin{align*}f(T)g(T) &= \frac1{(2\pi i)^2} \int_{\Gamma}\int_{\Gamma'}f(\la)g(\mu)R(\la,T)R(\mu,T)\ud \mu\ud \la
 \\ & = \frac1{(2\pi i)^2} \int_{\Gamma}\int_{\Gamma'}f(\la)g(\mu)(\la-\mu)\inv[R(\mu,T) - R(\la,T)]\ud \mu\ud \la
 \\ & = \frac1{(2\pi i)^2} \int_{\Gamma'}\int_{\Gamma}f(\la)g(\mu)(\la-\mu)\inv R(\mu,T)\ud \la\ud \mu
 \\ & = \frac1{2\pi i} \int_{\Gamma'}f(\mu)g(\mu)R(\mu,T)\ud \mu = (fg)(T).
\end{align*}
Here we used the fact that
$$ \frac1{2\pi i}\int_{\Gamma'}g(\mu)(\la-\mu)^{-1}\ud \mu = 0 \ \ \hbox{and} \ \ \frac1{2\pi i}\int_{\Gamma}f(\la)(\la-\mu)\inv\ud \la =  f(\mu).$$

\smallskip
\ref{it:Dunford4}: \ This follows from Proposition \ref{prop:spectrum-dual} and the continuity of the mapping $S \mapsto S^*$\!, which allows one to take adjoints under the integral sign.

\smallskip
\ref{it:Dunford3}: \ As was already mentioned, we may assume that $\mu\not\in\Om$. The identity $(\mu-\la) \cdot f_\mu(\la) = f_\mu(\la) \cdot (\mu-\la) = 1$ implies, by \ref{it:Dunford1} and \ref{it:Dunford2}, that
$(\mu-T) f_\mu(T) = f_\mu(T)(\mu-T) = I$. It follows that $f_\mu(T) = R(\mu,T)$.
\end{proof}

\begin{theorem}[Spectral mapping theorem]\label{thm:SMT}\index{theorem!spectral mapping, holomorphic calculus}\index{spectral!mapping, the holomorphic calculus}
Let $\Om\subseteq\C $ be an open set containing $\si(T)$. For all $f\in H(\Om)$ we have
 $$ \si(f(T)) = f(\si(T)).$$
\end{theorem}

\begin{proof}
We begin with the proof of the inclusion
$\sigma(f(T))\subseteq f(\sigma(T))$.
Fix $\la \not\in f(\sigma(T))$; our aim is to show that $\la\not\in\si(f(T))$.
Put
$$
\Om':=\{z\in\Om:\, f(z)\ne\la\}.
$$
Since $\la\notin f(\si(T))$, this is an open set containing
$\si(T)$. Moreover, the function
$g_\la(z):=(\la-f(z))^{-1}$ is holomorphic on $\Om'$.
By the multiplicativity of the holomorphic
functional calculus applied to $H(\Omega')$,
$$({\la} - f(T))g_{{\la}}(T) = g_{{\la}}(T)({\la} - f(T)) = I,$$
from which we infer that $\la\not\in \si(f(T))$.

Turning to the converse inclusion $f(\sigma(T))\subseteq \sigma(f(T))$,
let $\la\in \si(T)$. By the theory of functions in one complex variable we have
$f(\la)-f(z) = (\la-z)h_\la(z)$ for some $h_\la\in H(\Om)$, so $f(\la) - f(T) = (\la-T)h_\la(T).$
Since $\la\in \sigma(T)$, the operator $\la-T$ is not invertible. Since $\la-T$ and $h_\la(T)$ commute, this implies that $f(\la) - f(T)$ is not invertible and therefore $f(\la)\in \sigma(f(T))$.
\end{proof}

\begin{theorem}[Composition]\label{thm:Dunford-composition}
Let $\Omega\subseteq\C$ be an open set containing $\sigma(T)$, let
$f\in H(\Omega)$, and let $\Omega'\subseteq\C$ be an open set
containing $\sigma(f(T))$ such that $f(\Omega)\subseteq\Omega'$.
Then for all $g\in H(\Omega')$,
$$
   g(f(T))=(g\circ f)(T).
$$
\end{theorem}

\begin{proof}
Choose an admissible contour $\Gamma'$ for $\sigma(f(T))$ in
$\Omega'$. Let
$$
  U:=\{u\in\Omega'\setminus \Gamma' :\,  w(\Gamma';u)=1\}.
$$
Then $U$ is an open set containing $\sigma(f(T))$. We have $f(\sigma(T))\subseteq \sigma(f(T))$ by the spectral
mapping theorem, and therefore
$\sigma(T)\subseteq f^{-1}(U)$.
Put $\Omega_0:=\Omega\cap f^{-1}(U)$.
Then $\Omega_0$ is an open set containing $\sigma(T)$.
Choose an admissible contour $\Gamma$ for $\sigma(T)$
in $\Omega_0$.

If $\mu$ is a point on $\Gamma'$\!, then $h_\mu(\lambda):= (\mu-f(\la))^{-1}$ defines a function in $H(\Om_0)$
and
$$ h_\mu(T) (\mu-f(T)) =  (\mu-f(T))h_\mu(T) =I.$$
Consequently,
$$ R(\mu,f(T)) = h_\mu(T) = \frac1{2\pi i} \int_{\Gamma} (\mu - f(\la))^{-1}R(\la,T)\ud\la.
$$
Hence, using Fubini's theorem to justify the change of order of integration,
\begin{align*}
 g(f(T)) & = \frac1{2\pi i} \int_{\Gamma'} g(\mu)R(\mu,f(T))\ud\mu
\\ & = \bigl(\frac1{2\pi i}\bigr)^2 \int_{\Gamma'}\int_{\Gamma}  g(\mu)(\mu - f(\la))^{-1}R(\la,T)\ud\la\ud\mu
\\ & =  \bigl(\frac1{2\pi i}\bigr)^2 \int_{\Gamma}\int_{\Gamma'}  g(\mu)(\mu - f(\la))^{-1}R(\la,T)\ud\mu\ud\la
\\ & \stackrel{(*)}{=} \frac1{2\pi i} \int_{\Gamma} g(f(\la)) R(\la,T)\ud\la
 = (g\circ f)(T).
\end{align*}
In $(*)$ we used the Cauchy integral formula in the form
$\frac1{2\pi i} \int_{\Gamma'}  g(\mu)(\mu - f(\la))^{-1}\ud\mu =  g(f(\la))$. The role of $U$ is to ensure that, for every $\lambda\in\Gamma$, we have
$f(\lambda)\notin\Gamma'$ and $w(\Gamma';f(\lambda))=1$, which justifies both
the definition of $h_\mu$ on $\Omega_0$ and the final application of the Cauchy integral formula.
\end{proof}

An interesting application of the holomorphic calculus arises when the spectrum is the union of two nonempty disjoint compact sets.

\begin{theorem}[Spectral projections]\label{thm:spec-proj}\index{spectral!projection}
Suppose that $\sigma(T)$ is the union of two nonempty disjoint compact sets
$K_1$ and $K_2$, and let $\Om_1$ and $\Om_2$ be disjoint open sets containing
$K_1$ and $K_2$, respectively. Let $\Gamma_{\!1}$ and $\Gamma_{\!2}$ be
admissible contours for $K_1$ and $K_2$ in $\Om_1$ and $\Om_2$, respectively. The operators
$$ P_j := \frac1{2\pi i} \int_{\Gamma_{\!j}} R(\la,T)\ud \la, \quad  j\in \{1,2\},$$
are projections, their ranges $X_1$ and $X_2$ are invariant under $T$, and we have a direct sum decomposition
$ X = X_1\oplus X_2.$
Moreover, $$\sigma(T|_{X_j}) = K_j, \quad  j\in \{1,2\}.$$
\end{theorem}

\begin{proof}
To see that $P_j$ is a projection we just note that $P_j = f_j(T)$, where $f_j:\Om_1\cup\Om_2\to \C$ is the
holomorphic function which is $1$ on $\Om_j$ and $0$ elsewhere. By the multiplicativity of the holomorphic calculus, $P_j$ is bounded and
$$ P_j^2 = (f_j(T))^2 = f_j^2(T) = f_j(T) = P_j.$$
Also, $f_1+f_2\equiv1$ implies that $P_2=I-P_1$.
Since $P_1$ is a projection, we have
$\ran(P_2)=\ran(I-P_1)=\ker(P_1)$.
Every $x\in X$ can therefore be written as $x=P_1x+P_2x$,
with the two terms in $\ran(P_1)$ and $\ker(P_1)$, respectively.
These subspaces have trivial intersection, since $P_1$ acts
as the identity on its range and as zero on its kernel.
It follows that we have a direct sum decomposition
$X = \ran(P_1)\oplus \ran(P_2).$
Since $T$ obviously commutes with $P_j$, it follows that $T$ maps $X_j = \ran(P_j)$ into itself.
It follows that $T$ restricts to a bounded operator $T_j := T|_{X_j}$ on $X_j $.

We first prove that $\si(T_j)\subseteq K_j$.
Let $\mu\in\C\setminus K_j$, and let $k\in\{1,2\}$ be the
index different from $j$. By replacing $\Om_j$ by
$\Om_j\setminus\{\mu\}$ if necessary, we may assume that
$\mu\notin\Om_j$. The independence of the holomorphic calculus
from the choice of neighbourhoods ensures that $P_j=f_j(T)$
still holds.

Define $h_{\mu,j}\in H(\Om_j\cup\Om_k)$ by
$h_{\mu,j}(\la):=(\mu-\la)^{-1}$ on $\Om_j$ and
$h_{\mu,j}(\la):=0$ on $\Om_k$, and put
$H_{\mu,j}:=h_{\mu,j}(T)$.
Then $(\mu-\cdot)h_{\mu,j}=f_j$ and $f_jh_{\mu,j}=h_{\mu,j}$. Hence
$$
(\mu-T)H_{\mu,j}=H_{\mu,j}(\mu-T)=P_j, \qquad
P_jH_{\mu,j}=H_{\mu,j}P_j=H_{\mu,j}.
$$
Thus $H_{\mu,j}$ maps $X_j$ into $X_j$. For $y\in X_j$ we obtain
$$
(\mu-T_j)H_{\mu,j}y=y,
\qquad
H_{\mu,j}(\mu-T_j)y=y.
$$
It follows that $\mu\in\rh(T_j)$. This proves that $\si(T_j)\subseteq K_j$; in particular $\si(T_1)\cap\si(T_2)=\emptyset$.

To prove the reverse inclusions, we first claim that
\begin{equation}\label{eq:spdec-dir-sum}
\si(T)=\si(T_1)\cup\si(T_2).
\end{equation}
Indeed, the projections $P_1$ and $P_2$ commute with
$T$. Hence, if $\mu\in\rh(T)$, they therefore also commute with
$R(\mu,T)$, so $R(\mu,T)$ maps each $X_j$ into itself and
its restriction to $X_j$ is a two-sided inverse to $\mu-T_j$.
Thus $\mu\in\rh(T_1)\cap\rh(T_2)$. Conversely, if $\mu\in\rh(T_1)\cap\rh(T_2)$, the bounded operator $R(\mu,T_1)P_1+R(\mu,T_2)P_2$
is a two-sided inverse to $\mu-T$. Hence $\mu\in\rh(T)$.

Combining \eqref{eq:spdec-dir-sum} with $\si(T)=K_1\cup K_2$,
$\si(T_j)\subseteq K_j$, and $K_1\cap K_2=\emptyset$, we obtain
$\si(T_j)=K_j$ for $j=1,2$.
\end{proof}

As a further application of the holomorphic calculus we prove the following exact formula for the {\em spectral radius}\index{spectral!radius}
$$r(T):= \sup\{|\la|:\, \la\in\si(T)\}$$
of a bounded operator $T\in\calL(X)$ (with the convention $\sup \emptyset = 0$ to deal with the trivial case $X = \{0\}$). Since the spectrum $\sigma(T)$ is
contained in the closed disc of radius $\n T\n$ we have
$r(T) \le \n T\n.$ More generally, by the spectral mapping theorem,
$r(T)^n = r(T^n) \le \n T^n\n$, so $r(T)\le \n T^n\n^{1/n}$\!.
The precise relation is given by the following theorem.

\begin{theorem}[Gelfand]\label{thm:rT}\index{theorem!spectral radius formula}
For every bounded operator $T\in\calL(X)$ we have
$$ r(T) = \inf_{n\ge 1} \n T^n\n^{1/n} = \limn \n T^n\n^{1/n}.$$
\end{theorem}

Existence of the limit on the right-hand side is part of the assertion.

\begin{proof}
We have already observed that $$r(T) \le  \inf_{n\ge 1} \n T^n\n^{1/n}.$$
The theorem will be proved once we show that $$\limsup_{n\to\infty} \n T^n\n^{1/n} \le r(T).$$

Fix $\e>0$ and let $\Gamma$ denote the circular contour about the origin with radius $R =  r(T) +\e$, oriented anticlockwise.
Then
  $$ T^n = \frac1{2\pi i} \int_\Gamma \la^n R(\la,T)\ud \la$$
implies
  $$ \n T^n\n \le \frac1{2\pi} \cdot 2\pi R \cdot R^n \sup_{|\la|= R}\n R(\la,T) \n,$$
the supremum being finite in view of the continuity of $\la\mapsto R(\la,T)$.
Taking $n$th roots and passing to the limit superior as $n\to\infty$, we obtain
$ \limsup_{n\to\infty} \n T^n\n^{1/n} \le R = r(T)+\e.$
Since $\e>0$ was arbitrary, this gives the result.
\end{proof}

As an application we have the following stability result.

\begin{theorem}[Lyapunov's stability theorem]\label{thm:Lyapunov} If $A\in \calL(X)$ is a bounded operator with $\sigma(A)\subseteq \{z\in\C:\, \Re z < 0\}$, then
 $$ \lim_{t\to\infty} \n e^{tA} \n =0.$$
\end{theorem}
\begin{proof}
 Since $\sigma(A)$ is compact, there is a $\delta>0$ such that $\sigma(A)\subseteq \{z\in\C: \, \Re z \le -\delta\}$.
By the spectral mapping theorem (Theorem \ref{thm:SMT}) this implies that $\sigma(e^{A})\subseteq \{z\in \C: \, |z| \le e^{-\delta}\}$. Stated differently,
we have $r(e^{A}) \le e^{-\delta} < 1$.
Choose $q$ with $r(e^A)<q<1$. By Gelfand's theorem
(Theorem \ref{thm:rT}), applied to $e^A$, we have
$\n e^{nA}\n^{1/n}\to r(e^A)$.
Consequently, $\n e^{nA}\n\le q^n$ for all sufficiently
large integers $n$.

By estimating the power series defining $e^{sA}$ we obtain
$\n e^{sA}\n\le e^{s\n A\n}$ for all $s\ge0$.
Now write $t=n+s$, where $n\ge0$ is an integer and $0\le s<1$.
For sufficiently large $t$ we obtain
$$
\n e^{tA}\n
\le \n e^{nA}\n\,\n e^{sA}\n
\le q^n e^{\n A\n}
\le q^{t-1}e^{\n A\n}.
$$
Since $0<q<1$, this proves the assertion, with exponential rate.
\end{proof}

The interpretation of this theorem is as follows. Consider the initial value problem
$$
\left\{
\begin{aligned}
 u'(t) & = Au(t), \quad t\ge 0,\\
 u(0) & = u_0,
\end{aligned}
\right.
$$
where $A \in \calL(X)$ is bounded and $u_0\in X$ is given.
By differentiating the power series defining $e^{tA}$ we see that this problem is solved by the function
$ u(t) := e^{tA}u_0.$
Lyapunov's theorem now gives the following sufficient spectral criterion for stability
of this solution: if the spectrum of $A$ is contained in the open
left half-plane, then $\lim_{t\to\infty} \n u(t)\n = 0$.

A far-reaching generalisation of Lyapunov's theorem is given by Theorem \ref{thm:strong-stability}.

\begin{problems}
\item\label{prod:commprod-inverses}
Prove that the product of two commuting bounded operators is
invertible if and only if each of the operators is invertible.

\item
Prove the following improvement to Lemma \ref{lem:Neumann}: if $T\in \calL(X)$ is such that the sum
$\sum_{n=0}^\infty T^n x$ converges for all $x\in X$, then $I-T$ is invertible. What is its inverse?

\item
Show that if $X$ is nonzero, and $T\in \calL(X)$, $\la\in\varrho(T)$, and $\mu\in \C$ satisfy $|\la-\mu|\le \delta \n R(\la,T)\n^{-1}$ with $0\le \delta<1$, then $\mu\in\varrho(T)$ and
 $$ \n R(\mu,T)- R(\la,T)\n \le \frac{\delta}{1-\delta}\n R(\la,T)\n.$$

\item\label{prob:exp}
Prove in an elementary way, by multiplying power series, that if $T\in\calL(X)$, then
$e^{w T}e^{z T} = e^{(w+z)T}$ for all complex numbers $w,z\in \C.$

\item\label{prob:rightshift-spectrum}
For $1\le p \le \infty$, consider the right shift $T\in \calL( \ell^p)$:
$$ T: (a_1,a_2,a_3,\ldots) \mapsto (0,a_1,a_2,\ldots). $$
Give a direct proof (see Example \ref{ex:rightshift}) that
$\sigma(T) = \{\la\in \C:\, |\la|\le 1\}$.

\item\label{prob:spectrum-FMO}
We compute the spectra of some multiplier operators.
\begin{enumerate}[\rm(a), leftmargin=*]
  \item Let $m\in C[0,1]$ and define $T_m\in \calL(C[0,1])$ by
  $$ (T_m f)(t) = m(t)f(t),\quad f\in C[0,1],\ t\in [0,1].$$
  Show that $\sigma(T_m)$ coincides with the range of $m$.
  \item Let $m\in L^{\infty}(0,1)$, let $1\le p\le \infty$, and define $T_m\in \calL(L^p(0,1))$ by
  $$ (T_m f)(t) = m(t)f(t),\quad f\in L^p(0,1),\ t\in (0,1).$$
  Show that $\sigma(T_m)$ coincides with the {\em essential range}\index{essential range} of $m$, that is, the set of all $\lambda\in \C$
  such that for any open set $U\subseteq \C$ with $\lambda\in U$ the set $\{t\in (0,1):\, m(t)\in U\}$ has positive measure.

  \noindent{\em Hint:}\ First show that  $\lambda$ is not contained in the essential range of $m$ if and only if $\frac{1}{m-\lambda}$ is well defined almost everywhere and essentially bounded.
\end{enumerate}

\item\label{prob:spectrum-FM}
Let $T_m$ be the Fourier multiplier on $L^2(\R^d)$ with symbol
$m\in L^\infty(\R^d)$.

\begin{enumerate}[\rm(a), leftmargin=*]
\item
Show that $\sigma(T_m)$ equals the essential range of $m$.

\noindent{\em Hint:}\ By Plancherel's theorem, $T_m$ is unitarily
equivalent to multiplication by $m$ on $L^2(\R^d)$.
Repeat the argument of part (b) of Problem
\ref{prob:spectrum-FMO} for this multiplication operator.

\item
Let $\Omega\subseteq\C$ be open with $\sigma(T_m)\subseteq\Omega$, and let
$f$ be holomorphic on $\Omega$. Show that $f\circ m\in L^\infty(\R^d)$
and that
$$
  f(T_m)=T_{f\circ m},
$$
where $f(T_m)$ is defined by the holomorphic functional calculus.
\end{enumerate}

\item
Show that every nonempty compact subset of $\C$ can be realised as the spectrum of some bounded operator.

\item
Let $T\in\calL(X)$. Show that if $S\in\calL(X)$ satisfies $\n S\n<1$, then $\Ran(T(I+S))=\Ran(T)$.

\item
Show that if $P$ and $Q$ are projections satisfying $\n P-Q\n <1$, then
$\dim(\Ran(P)) = \dim(\Ran(Q))$ (allowing the case $\infty=\infty$).

\noindent{\em Hint:}\ Use the invertibility of $I - (P-Q)$ and $I - (Q-P)$, and note that $PQ = P(I-(P-Q))$ and $QP = Q(I-(Q-P))$.

\item\label{prob:contours}
Show that if $T\in\calL(X)$ and $\Omega\subseteq \C$ is an open set containing $\si(T)$, then admissible contours for $\si(T)$ in $\Omega$ always exist.

\item\label{prob:isometry-approx-spectrum}
Show that if $T\in \calL(X)$ is an isometry, then every approximate eigenvalue of $T$ has modulus one. Use this to give an alternative proof of Corollary \ref{cor:spectrum-isometry}.

\item
Let $\varphi: [0,1]\to [0,1]$ be continuous, and for $f\in C[0,1]$ define
$$
 C_\varphi f := f\circ \varphi.
$$
\begin{enumerate}[\rm(a), leftmargin=*]
\item Prove that $C_\varphi$ is a bounded linear operator on $C[0,1]$ with $\|C_\varphi\|=1$.
\item Prove that $1\in \sigma(C_\varphi)$.
\item Prove that $0\in \sigma(C_\varphi)$ if and only if $\varphi$ is not a bijection of $[0,1]$ onto itself.
\end{enumerate}

\item
Determine the spectrum of the Hilbert transform $H$ on $L^2(\R)$.

\noindent{\em Hint:} \ Use $H^2=-I$, and show that both $i$ and $-i$ are eigenvalues by considering nonzero functions whose Fourier transforms are supported in the negative and positive half-lines, respectively.

\item\label{prob:spectrum-reduced}
Suppose that we have a direct sum decomposition $X = X_0\oplus X_1$.
Prove that if $T\in\calL(X)$ leaves both $X_0$ and $X_1$ invariant, then
$$\si(T) = \sigma(T|_{X_0})\cup \sigma(T|_{X_1}),$$
viewing
$T|_{X_0}$ and $T|_{X_1}$ as bounded operators on the Banach spaces $X_0$ and $X_1$.

\item\label{prob:sigmaSTvsTS}
Prove that for all $S,T \in \calL(X)$ we have
$$ \sigma(ST)\setminus\{0\} =  \sigma(TS)\setminus\{0\}.$$

\noindent{\em Hint:}\ Use the Neumann series to relate the resolvents of $ST$ and $TS$.

\item
Let $T\in\calL(X)$. Show that the following assertions are equivalent:
\begin{enumerate}[\rm(a), leftmargin=*]
\item $r(T)<1$;
\item there exist constants $C\ge1$ and $0<q<1$ such that
$\n T^n\n\le Cq^n,\quad n\ge 0$.
\end{enumerate}
Deduce that if $r(T)<1$, then the series $\sum_{n=0}^\infty T^n$
converges in operator norm and its sum equals
$(I-T)^{-1}$.

\item
Let $T\in \calL(X)$ and let $f$ be a nonzero holomorphic function on a connected open set $\Omega$ containing $\si(T)$. Show that if $f(T) = 0$, then $\si(T)$ is a finite set.

\item\label{prob:allan-ransford}
Let $T\in\calL(X)$ be an operator with the property that $\si(T) = \{1\}$. Assume furthermore that $T$ is {\em doubly power bounded},\index{doubly power bounded} that is, $T$ is invertible and $\sup_{k\in\Z} \n T^ k\n<\infty$.

\begin{enumerate}[\rm(a), leftmargin=*]
  \item\label{it:allan-ransford1} Using the holomorphic calculus and its properties, explain that we can define the bounded operators
  $S:= -i\log T$ and $\sin(nS)$, $n\in \N$.
  \item\label{it:allan-ransford2} Show that $\sin(nS) = \frac1{2i}(T^n - T^{-n})$, $n\in\N$.
  \item\label{it:allan-ransford3} Using the spectral mapping theorem, show that
  $$\sigma(S)=(-i\log)(\sigma(T))=\{0\}.$$ Deduce that $\sigma(nS)=\{0\}$ and, again by spectral mapping, $$\sigma(\sin(nS))=\sin(\sigma(nS))=\{0\}.$$
\end{enumerate}
We now use the fact that if $\sum_{k=0}^\infty c_kz^k$ denotes the Taylor series of
the principal branch of $\arcsin z$ at $z=0$, then
$c_k\ge 0$ for all $k\in\N$ and $\sum_{k=0}^\infty c_k =\arcsin(1)=
\frac{\pi}{2}$.
\begin{enumerate}[\rm(a), leftmargin=*]\setcounter{enumii}{3}
  \item\label{it:allan-ransford4} Show that $nS = \arcsin (\sin (nS))$, where the latter is again defined by means of the holomorphic calculus, and deduce that
  $$ \n nS\n \le \frac{\pi}{2} \sup_{k\in\Z}\n T^ k\n.$$
  Conclude that $S=0$ and $T=e^{iS}=I.$

  \noindent{\em Hint:}\
  First use part \ref{it:allan-ransford2} and the binomial theorem to show that $\n(\sin(nS))^j\n\le M:=\sup_{k\in\mathbb Z}\|T^k\|$. Then apply the power series of $\arcsin$.
\end{enumerate}

\end{problems}

%% file: ch07-CompactOperators.tex
\chapter{Compact Operators}\label{ch:CompactOperators}

\blfootnote{\rm \small \begin{quote}This book has been published by Cambridge University Press in the series `Cambridge Studies in Advanced Mathematics'. The present version is free to view and download for personal use only. Not for re-distribution, re-sale or use in derivative works.\newline {\copyright} Jan van Neerven\end{quote}}

\noindent
This chapter studies the class of compact operators. By definition, these are the operators that map bounded sets to relatively compact sets. Examples include integral operators on various Banach spaces of functions over a compact domain. For this reason, compact operators have important applications in the theory of partial differential equations and mathematical physics. After establishing some generalities, we prove the Riesz--Schauder theorem, which asserts that the nonzero part of the spectrum of a compact operator is discrete and consists of eigenvalues.

The final section of this chapter presents an introduction to the theory of Fredholm operators. These are the operators that are invertible modulo a compact operator, and the difference between their kernel and cokernel dimensions is the so-called Fredholm index. As an example, we prove the Noether--Gohberg--Krein theorem, which states that a Toeplitz operator with continuous zero-free symbol is Fredholm and its index equals the negative winding number of its symbol.

\section{Compact Operators}\label{sec:compact}

Throughout this chapter the scalar field is complex.
We fix Banach spaces $X$ and $Y$.

\begin{definition}[Compact operators]An operator $T\in\calL(X,Y)$ is {\em compact}\index{compact!operator} if it maps bounded sets to relatively compact sets.
\end{definition}

Since every bounded set in $X$
is contained in a multiple of the unit ball  $B_X = \{x\in X:\, \n x\n<1\}$, a bounded operator is compact if and only if
$T B_X$ is relatively compact.
Furthermore, since a subset of a Banach space is relatively compact if and only if it is relatively sequentially compact, a linear operator $T$ is compact
 if and only if $(Tx_n)_{n\ge 1}$ has a convergent subsequence for every bounded sequence $(x_n)_{n\ge 1}$ in $X$.

The set of all compact operators is a linear subspace of $\calL(X,Y)$.
It is clear that $cT$ is compact if $T$ is compact,
for any scalar $c\in\K$, and if $S$ and $T$ are compact, then $S+T$ is compact as well: for if $SB_X$ and $T B_X$ are contained
in the compact sets $K$ and $L$, then $(S+T)B_X$ is contained in the compact set $K+L$ (this set is
the image of the compact set $K\times L$ under the continuous mapping $(x_1,x_2)\mapsto x_1+x_2$
from $Y\times Y$ to $Y$). Moreover, if $T\in\calL(X,Y)$ is compact
and $S\in \calL(X'\!,X)$ and $U\in \calL(Y,Y')$ are bounded, then $UTS\in\calL(X'\!,Y')$ is compact.
Indeed, if $C$ is a bounded set in $X'$\!, then $S(C)$ is bounded in $X$, so $TS(C)$ is contained in a
compact set $K$ of $Y$, and then $UTS(C)$ is contained in the compact set $U(K)$.

\begin{example}\label{ex:FD-compact}
As an immediate corollary of Theorem \ref{thm:FD-compact}, the identity operator on a Banach space $X$ is compact if and only if $X$ is finite-dimensional.
\end{example}

\begin{example}\label{ex:finiterank}
A bounded operator is said to be of {\em finite rank}\index{finite!rank operator} if its range is finite-dimensional.
Since bounded sets in finite-dimensional spaces are relatively compact, every finite-rank operator is compact.
\end{example}

\begin{example}[Integral operators on \hbox{$C(K)$}]\label{ex:integral-comp}
Let $\mu$ be a finite Borel measure on a compact metric space $K$ and let $k: K\times K\to\K$ be continuous.
 Then the operator $T: C(K)\to C(K)$,
 $$ Tf(x):= \int_K k(x,y)f(y)\ud \mu(y), \quad f\in C(K),\ x\in K,$$
 is well defined and bounded by Example \ref{ex:kernel}. Let us show that $T$ is compact. Let $(f_n)_{n\ge 1}$ be a bounded sequence in $C(K)$. We claim that the bounded sequence
 $(Tf_n)_{n\ge 1}$ is equicontinuous. Once we have shown this, the
 Arzel\`a--Ascoli theorem (Theorem \ref{thm:Arzela-Ascoli}) implies that this sequence is relatively compact, hence has a convergent subsequence. This implies that $T$ is compact.

 To check the equicontinuity we first note that $K\times K$ is compact and hence $k$ is uniformly continuous,
 in the sense that given any $\e>0$ we can find $\delta>0$
 such that $d(x,x')+d(y,y')<\delta$ implies $|k(x,y)-k(x'\!,y')| <\e$. Then, if $d(x,x')<\delta$,
 $$ | Tf_n(x) - Tf_n(x')| \le \int_K |k(x,y)-k(x'\!,y)||f_n(y)|\ud \mu(y)  \le \e M \mu(K), $$
 where $M = \sup_{n\ge 1}\n f_n\n_\infty.$
 The equicontinuity follows immediately from this.
\end{example}

The next proposition shows that the compact operators form a closed subspace of $\calL(X,Y)$. This subspace will be denoted by\index{$K$@$\mathscr{K}(X,Y)$}
$\mathscr{K}(X,Y)$
and we write $\mathscr{K}(X):=\mathscr{K}(X,X)$.

\begin{proposition}\label{prop:limcomp}
 If $\lim_{n\to \infty} \n T_n-T\n = 0$ with each $T_n$ compact, then $T$ is compact. In particular, uniform limits of finite-rank operators are compact.
\end{proposition}
\begin{proof}
For any $\eps>0$ we can choose an index $n_\eps\ge 1$ such that $\n T_{n_\eps} - T\n<\eps$. Since $T_{n_\eps}B_X$ is relatively compact and $TB_X \subseteq T_{n_\eps}B_X + B(0;\eps)$, the relative compactness of $TB_X$ follows from  Proposition \ref{prop:compact-totbdd}.
\end{proof}

The following converse of the second assertion of Proposition \ref{prop:limcomp} holds for Hilbert space operators:

\begin{proposition}\label{prop:fr-dense-compH} Let $H$ and $K$ be Hilbert spaces.
 An operator $T\in\calL(H,K)$ is compact if and only if it is the uniform limit of finite-rank operators.
\end{proposition}
\begin{proof}
 It remains to prove the `only if' part. Let $T$ be compact, say $TB_H\subseteq C$ with $C\subseteq K$ compact,
 where $B_H$ is the open unit ball of $H$.
 Fix $\e>0$ and let $B(y_1;\eps), \dots, B(y_N;\eps)$ be an open cover of $C$.
 Let $Y$ denote the linear span of $\{y_1,\dots,y_N\}$ and let $P$ be the orthogonal projection of $K$ onto $Y$.
 This projection is of finite rank, and therefore $PT$ is of finite rank.
 For any $x\in H$ with $\n x\n< 1$ we have $Tx\in C$, so
 $\n Tx-y_j\n < \e$ for some $1\le j\le N$.
 Since $Py_j=y_j$ and $I-P$ is an orthogonal projection,
$$
\n Tx-PTx\n
=\n(I-P)(Tx-y_j)\n
\le\n Tx-y_j\n<\e.
$$
Taking the supremum over all $x\in H$ with $\n x\n<1$,
we obtain $\n T-PT\n\le\e$.
\end{proof}

\begin{example}[Integral operators on $L^2(K,\mu)$]\label{ex:integral-comp2}
Let $\mu$ be a finite Borel measure on a compact metric space $K$ and let $k:K\times K\to \K$ be square-integrable. By Fubini's theorem, for almost all $x\in K$ we have $k(x,\cdot)\in L^2(K,\mu)$, and for all $f\in L^2(K,\mu)$ the function
$$T f(x):= \int_K k(x,y)f(y)\ud \mu(y)$$
is well defined as an element of $L^2(K,\mu)$. The operator  $T:L^2(K,\mu)\to L^2(K,\mu)$ thus defined is bounded with norm $\n T\n \le \n k\n_2$.
Let us prove that $T$ is compact.

 First consider a continuous kernel
$\kappa\in C(K\times K)$ and define
$$
T_\kappa f(x):=\int_K\kappa(x,y)f(y)\ud\mu(y),
\quad f\in L^2(K,\mu),\ x\in K.
$$
The Cauchy--Schwarz inequality gives
$$
\n T_\kappa f\n_\infty
\le \mu(K)^{1/2}\n\kappa\n_\infty\n f\n_2,
$$
and, for $x,x'\in K$,
\begin{align*}
|T_\kappa f(x)-T_\kappa f(x')|
&\le
\Bigl(\int_K|\kappa(x,y)-\kappa(x',y)|^2\ud\mu(y)\Bigr)^{1/2}
\n f\n_2\\
&\le
\mu(K)^{1/2}\sup_{y\in K}|\kappa(x,y)-\kappa(x',y)|\,\n f\n_2.
\end{align*}
Since $\kappa$ is uniformly continuous, these estimates show that
$T_\kappa f\in C(K)$ and that the image under $T_\kappa$ of the
unit ball of $L^2(K,\mu)$ is uniformly bounded and equicontinuous.
By the Arzel\`a--Ascoli theorem, $T_\kappa$ is therefore compact
as an operator from $L^2(K,\mu)$ into $C(K)$.
The natural mapping from $C(K)$ into $L^2(K,\mu)$ is bounded,
so $T_\kappa$ is also compact as an operator on $L^2(K,\mu)$.

For a general square-integrable kernel $k$, choose
$\kappa_n\in C(K\times K)$ such that $\n\kappa_n-k\n_2\to0$;
this is possible by Remark \ref{rem:CKdenseLpK}.
The estimate for integral operators established above gives
$$
\n T_{\kappa_n}-T\n\le\n\kappa_n-k\n_2\to0.
$$
Since each $T_{\kappa_n}$ is compact, Proposition \ref{prop:limcomp}
implies that $T$ is compact.

In Example \ref{ex:kernel-HS} it will be shown, under the weaker assumption that the measure space $(K,\mu)$ be $\sigma$-finite, that the integral operator $T$ defined above is Hilbert--Schmidt. By Proposition \ref{prop:HS-compact}, this property implies compactness. The required separability of $L^2(K,\mu)$ follows from Remark \ref{rem:CKdenseLpK}.
\end{example}

We conclude this section with a duality result for compact operators.

 \begin{proposition}\label{prop:adjoint-compact}
  An operator $T\in \calL(X,Y)$ is compact if and only if its adjoint $T\s\in \calL(Y\s,X\s)$ is compact.
 \end{proposition}

\begin{proof}
 First we prove the `only if' part. Let $T$ be compact and let $K$ denote the closure of $TB_X$. By assumption, $K$ is a compact subset of $Y$.
 By restriction,
 every $y\s\in Y\s$ determines a function in $C(K)$ given by $y\s(y):= \lb y,y\s\rb$ for $y\in K$.
 Moreover, if $(y_n\s)_{n\ge 1}$ is a bounded sequence in
 $Y\s$\!, the corresponding functions are uniformly bounded and equicontinuous; the latter follows from
 $ |\lb y-y'\!, y_n\s\rb| \le M\n y-y'\n$ with $M:= \sup_{n\ge 1}\n y_n\s\n$. Hence, by the Arzel\`a--Ascoli theorem,
 there is a subsequence  $(y_{n_j}\s)_{j\ge 1}$ such that $$\lim_{j,k\to\infty} \n y_{n_k}\s - y_{n_j}\s\n_{C(K)} = 0.$$
 Since $K=\ov{TB_X}$, the definition of the adjoint gives
\begin{align*}
\n T^*y_{n_k}^*-T^*y_{n_j}^*\n
=\sup_{\n x\n<1} |\lb Tx,y_{n_k}^*-y_{n_j}^*\rb|=\n y_{n_k}^*-y_{n_j}^*\n_{C(K)}.
\end{align*}
Thus $(T^*y_{n_j}^*)_{j\ge1}$ is a Cauchy sequence in $X^*$ and therefore converges.
This proves that $T^*$ is compact.

\smallskip
To prove the `if' part, suppose that $T\s$ is compact. Then, by what we just proved, $T^{**}$ is compact as an operator from $X^{**}$
to $Y^{**}$. Identifying $X$ with a closed subspace of $X^{**}$ in the natural way, the restriction of
$T^{**}$ to $X$ maps $X$ to $Y$ and equals $T$. Since the restriction of a compact operator is compact,
the compactness of $T^{**}$ implies the compactness of $T$.
\end{proof}

\section{The Riesz--Schauder Theorem}\label{sec:spectrum-compact}

As we have seen in earlier examples, the spectrum of a bounded operator need not contain
eigenvalues. This is in sharp contrast to the situation in finite dimensions, where the
spectra of matrices consist of eigenvalues. The aim of the present section is to show that
compact operators spectrally resemble matrices to some degree.
The main result is Theorem \ref{thm:comp-spec}, which shows that the nonzero part of the spectrum of a compact operator is discrete and consists of eigenvalues.

\begin{lemma}\label{lem:kerran} If $T\in\calL(X)$ is compact, then:
\begin{enumerate}[label={\rm(\arabic*)}, leftmargin=*]
 \item\label{it:kerran1} $\ker(I-T)$ is finite-dimensional;
 \item\label{it:kerran2} $\Ran(I-T)$ is closed.
\end{enumerate}
\end{lemma}
\begin{proof}
\ref{it:kerran1}: \  Let $B_Y$ denote the unit ball of $Y:= \ker(I-T)$. We have $Ty = y$ for all $y\in Y$, so the compactness
of $T$ implies that $B_Y = TB_Y$ is relatively compact. By Theorem \ref{thm:FD-compact}, this implies that $Y$ is finite-dimensional.

\smallskip
\ref{it:kerran2}: \
Put $Y:=\ker(I-T)$. By part \ref{it:kerran1}, $Y$ is
finite-dimensional, so Proposition
\ref{prop:findim-compl}\ref{it:findim-compl1} gives a closed
subspace $Z$ of $X$ such that $X=Y\oplus Z$.
We claim that there exists a constant $c>0$ such that
$$
\n(I-T)z\n\ge c\n z\n,\quad z\in Z.
$$
If this were false, we could choose $z_n\in Z$ with
$\n z_n\n=1$ and $\n(I-T)z_n\n<1/n$.
By compactness, some subsequence $(Tz_{n_k})_{k\ge1}$ converges
to an element $z\in X$. Since
$$
z_{n_k}=(I-T)z_{n_k}+Tz_{n_k}\to z,
$$
the closedness of $Z$ implies that $z\in Z$, and we have
$\n z\n=1$. On the other hand, continuity of $I-T$ gives
$$
(I-T)z=\lim_{k\to\infty}(I-T)z_{n_k}=0.
$$
Thus $z\in Y\cap Z=\{0\}$, contradicting $\n z\n=1$.
This proves the claim.

By Proposition \ref{prop:closed-range}, the lower bound just
proved implies that $(I-T)Z$ is closed.
Since $I-T$ vanishes on $Y$ and $X=Y\oplus Z$, we have
$\Ran(I-T)=(I-T)Z$. Hence $\Ran(I-T)$ is closed.
\end{proof}

In order to describe the spectra of compact operators, we need the following lemma.

\begin{lemma}\label{lem:comp-inj-surj}
 Let $T\in\calL(X)$ be a compact operator. Then $I-T$ is injective if and only if $I-T$ is surjective.
\end{lemma}

\begin{proof}
 We begin with the proof of the `only if' part. We assume that $I-T$ is injective but not surjective and deduce a contradiction.
 By assumption, $X_1:= \ran(I-T)$ is a proper subspace of $X$ and by Lemma \ref{lem:kerran} this subspace is closed.
 It is clear that $TX_1\subseteq X_1$.
 Let $T_1$ denote the restriction of $T$ to $X_1$. This operator is compact.

 The operator $I-T_1: X_1 \to X_1$ is not surjective. Indeed, if it were, then for every $x\in X$ there would exist $y\in X_1$ such that
 $(I-T_1)y=(I-T)x$. Since $(I-T_1)y=(I-T)y$ and
 $I-T$ is injective on $X$, we would obtain $x=y$ and therefore $x\in X_1$. Thus
 it would follow that $X=X_1$, a contradiction.
 The same argument shows that $X_2:= \ran(I-T_1)$ is a proper closed subspace
 of $X_1$. It is clear that $TX_2\subseteq X_2$.
 Let $T_2$ denote the restriction of $T$ to $X_2$. Continuing as above we obtain a strictly decreasing
 sequence of closed subspaces $X_1\supsetneq X_2 \supsetneq \dots$ each of which is $T$-invariant,
 such that $$X_{n+1} = (I-T)X_n, \quad n\ge 1.$$
 By Lemma \ref{lem:Riesz} we can select vectors $x_{n} \in X_{n}\setminus X_{n+1}$ of norm one such that
 \begin{align}\label{eq:dist12}\inf_{y\in X_{n+1}} \n x_{n} - y\n \ge \frac12.
 \end{align}
 Since $T$ is compact, $(Tx_{n})_{n\ge 1}$ has a convergent subsequence $(Tx_{n_k})_{k\ge 1}$. Then, for $\ell > k$,
 $$ \n Tx_{n_k} - Tx_{n_\ell}\n = \n x_{n_k} + (T-I)x_{n_k} -  Tx_{n_\ell}\n \ge \frac12,$$
 where we use \eqref{eq:dist12} along with  $(T-I)x_{n_k} \in X_{{n_k}+1}$ and $Tx_{n_\ell} \in X_{n_\ell} \subseteq  X_{n_k+1}$.
 This contradicts the convergence of $(Tx_{n_k})_{k\ge 1}$.

 \smallskip
 Turning to the `if' part, assume that $I-T$ is surjective. If $T^* x^* = x^*$ for some $x^*\in X^*$,
 then by writing an arbitrary $x\in X$ as $(I-T)y$, we find
 $\lb x,x^*\rb = \lb y,(I-T^*)x^*\rb = 0$ for all $x\in X$, so $x^*=0$. This shows that $I-T^*$ is injective. Applying the preceding step to the operator $T^*$,
 which is compact by Proposition \ref{prop:adjoint-compact}, we conclude that
$I-T^*$ is surjective as well.
Applying the preceding argument to $T^*$, we conclude that
$I-T^{**}$ is injective, and hence $I-T$ is injective.
 \end{proof}

\begin{theorem}[Riesz--Schauder]\label{thm:comp-spec}\index{theorem!Riesz--Schauder}\index{spectrum!of a compact operator}
Let $T\in\calL(X)$ be a compact operator. Then:
\begin{enumerate}[label={\rm(\arabic*)}, leftmargin=*]
 \item\label{it:comp-spec1} every nonzero $\lambda\in\sigma(T)$ is an eigenvalue of $T$ and the eigenspace\index{eigenspace} $$E_\la := \{x\in X:\, Tx = \la x\}$$ is finite-dimensional;
 \item\label{it:comp-spec2} for every $r>0$, the number of eigenvalues satisfying $|\lambda|\ge r$ is finite;
 \item\label{it:comp-spec3} if $\dim(X) = \infty$, then $0\in \si(T)$.
\end{enumerate}
In particular, $\sigma(T)$ is at most countable and can only accumulate at $0$.
\end{theorem}

\begin{proof}
\ref{it:comp-spec1}: \ Let $0\not=\la\in \sigma(T)$ and suppose that $\la$ is not an eigenvalue of $T$. Then $I - \la^{-1}T$ is injective and hence, by the preceding lemma, surjective. It follows that $I - \la^{-1}T$ is invertible, and this implies that $\la\in \varrho(T)$.

Since $T$ acts as a multiple of the identity on the subspace $E_\la$ and $T$ is compact, the identity operator on $E_\la$ is compact. By Theorem \ref{thm:FD-compact}, this implies that $E_\la$ is finite-dimensional.

 \smallskip
\ref{it:comp-spec2}: \ Suppose there is an infinite sequence of distinct eigenvalues $\la_n$, $n\ge 1$, all of which satisfy $|\la_n|\ge r$.
Let $x_n\in X$ be eigenvectors for $\la_n$ of norm one.

Let $Y_n$ denote the linear span of $\{x_1,\dots,x_n\}$,
$n\ge 1$, and set $Y_0:=\{0\}$. By the distinctness of the eigenvalues,
the vectors $x_n$ are linearly independent (this is easily proved by induction on $n$).
Therefore $\dim(Y_n)=n$.
In particular, $Y_{n}$ is a proper subspace of $Y_{n+1}$. For $y\in Y_n$, say $y= \sum_{j=1}^n c_j x_j$,
we have
$ Ty =  \sum_{j=1}^n c_j \la_j x_j \in Y_n$
and
$$(\la_n - T)y =   \sum_{j=1}^n c_j (\la_n - \la_j) x_j =  \sum_{j=1}^{n-1} c_j (\la_n - \la_j) x_j\in Y_{n-1}.$$
Lemma \ref{lem:Riesz} shows that for every $n\ge 1$ it is possible to find a vector $y_n\in Y_n$ of
norm one such that $\n y-y_n\n \ge \frac12$ for all $y\in Y_{n-1}$.
Arguing as in the proof of Lemma \ref{lem:comp-inj-surj}, we obtain, for $n>m$, the lower
bound
\begin{align*}
\n Ty_n - Ty_m\n
 & = \n \la_n y_n + (T-\la_n)y_n  - T y_m\n
\\ & = |\la_n|\Bigl\n y_n - \frac{T y_m - (T-\la_n)y_n}{\la_n}\Bigr\n
 \ge \frac12|\la_n| \ge \frac12 r
\end{align*}
since $Ty_m\in Y_m\subseteq Y_{n-1}$ and
$(T-\la_n)y_n\in Y_{n-1}$,
and therefore $(Ty_n)_{n\ge 1}$ cannot have a convergent subsequence. This contradicts the compactness of $T$.

\smallskip
\ref{it:comp-spec3}: \ If $T$ is invertible, then $T^{-1}B_X$ is bounded and $B_X = T(T^{-1}B_X)$ is relatively compact. Therefore $X$ must be finite-dimensional.
\end{proof}

For compact normal operators on a Hilbert space, a more direct proof is sketched in Problem \ref{prob:spectrum-compact-normal}.

\begin{definition} The number $\dim(E_\la)$ is called the {\em geometric multiplicity}\index{geometric multiplicity} of $\la$.
\end{definition}

Suppose now that $T\in \calL(X)$ is compact and that $0\not=\la\in\si(T)$. Then $\la$ is an isolated point of $\si(T)$, that is, for small enough $r>0$ we have $B(\la;r)\cap \si(T) = \{\la\}$.
With $0<r'<r$ and $\Gamma = \partial B(\la;r')$ oriented anticlockwise, the spectral projection corresponding to $\la$ (see  Theorem \ref{thm:spec-proj}) is given by
$$P_{\la} := \frac1{2\pi i}\int_{\Gamma} R(\mu,T)\ud \mu.$$
By Theorem \ref{thm:spec-proj} the range $X_{\la} := P_{\la}X$ is invariant under $T$ and  $\sigma(T|_{X_{\la}}) = \{\la\}$.
In particular, $T|_{X_{\la}}$ is invertible. This operator is also compact as an operator on $X_{\la}$. This is only possible if $\dim(X_{\la})$ is finite. Thus we have proved:

\begin{corollary}\label{cor:compact-fd-spectralsubspaces} Let $T$ be a compact operator on a Banach space $X$.
For all nonzero $\la\in \sigma(T)$ the range $X_\la$ of the spectral projection $P_{\la}$ is finite-dimensional.
\end{corollary}

\begin{definition}\label{def:alg-multipl} The number $\dim(X_{\la})$ is called the {\em algebraic multiplicity}\index{algebraic!multiplicity} of $\la$.
\end{definition}

\begin{proposition} Let $T\in \calL(X)$ be a compact operator.
Then the geometric multiplicity of every nonzero $\la\in \sigma(T)$ is less than or equal to its algebraic multiplicity.
\end{proposition}
\begin{proof}
If $Tx = \la x$ and $\Gamma$ is as before, then $P_\la x = ( \frac1{2\pi i}\int_{\Gamma} (\mu-\la)^{-1}\ud \mu)x = x$. This shows that
the eigenspace $E_\la$ is contained in the range of the projection $P_\la$.
\end{proof}

\begin{example}[Jordan normal form]\label{ex:Jordan}\index{Jordan!normal form}
Consider a $k\times k$ Jordan block
$$
J_\la = \left(\begin{array}{ccccc}
\la \    & 1      & 0     & \cdots & 0      \\
 0  \    & \la    & 1     &  \ddots      & \vdots \\
\vdots \  & \ddots &\ddots &\ddots  &  0      \\
\vdots    &        & \ddots &  \ddots     &  1\\
 0   \   &  \cdots &      &   0    & \la    \\
\end{array}\right)
$$ and its resolvent
$$
R(\mu,J_\la) = \left(\begin{array}{ccccc}
(\mu-\la)^{-1} \    & (\mu-\la)^{-2}      &      & \cdots & (\mu-\la)^{-k}      \\
 0  \    & (\mu-\la)^{-1}    & (\mu-\la)^{-2}     &       & \vdots \\
 \\
  \vdots      &   \ddots & \ddots & \ddots &  \vdots      \\
  \\
 \vdots  &  &\ddots &\ddots  &  (\mu-\la)^{-2}      \\
 0   \   &      \cdots & &    0    & (\mu-\la)^{-1}    \\
\end{array}\right).
$$
If $A$ is in Jordan normal form and $\la\in\sigma(A)$, it follows that
the spectral projection corresponding to $\la$ is given by
$$P_\la =
\frac1{2\pi i}\int_\Gamma R(\mu,A)\ud\mu
= I_\la,$$
where $I_\la$ is the diagonal matrix with $1$'s on all diagonal entries
belonging to Jordan blocks with eigenvalue $\la$, and $0$'s elsewhere.
It follows that the algebraic multiplicity $\nu_\la$ of $\la$ equals
the sum of the dimensions of all Jordan blocks with $\la$ on the diagonal.
\end{example}

\section{Fredholm Theory}\label{sec:Fredholm-theory}

We now introduce the class of Fredholm operators.
Throughout this section we let $X$ and $Y$ be Banach spaces.

\subsection{The Fredholm Alternative}\label{subsec:Fredholm-alternative}

From the results in the preceding section we know that if $T$ is a compact operator on a Banach space $X$,
then every nonzero $\la\in \sigma(T)$ is an eigenvalue and the corresponding eigenspace is finite-dimensional. The next theorem (when applied to the compact operator $\la^{-1}T$) asserts that the dimension of the eigenspace is equal to the codimension of the range of $\la-T$.
This generalises the elementary result in linear algebra that for a $d\times d$ matrix $A$ we have $\dim\ker(A) + \dim\Ran(A) = d$.

\begin{theorem}[Fredholm alternative]\label{thm:Fredholm-alternative}\index{theorem!Fredholm alternative}\index{Fredholm!alternative}
If $T\in\calL(X)$ is compact, then
$$\dim\Ker(I-T) =  \codim\Ran(I-T),$$
and both numbers are finite.
\end{theorem}

This theorem contains Lemma \ref{lem:comp-inj-surj} as a special case.
The proof of the theorem is based on the following geometric lemma. Recall that if $Y\subseteq X^*$,
then $${}^\perp Y = \{x\in X:\, \lb x,x^*\rb = 0\ \hbox{ for all }\ x^*\in Y\}.$$

\begin{lemma}\label{lem:codim}
If $Y$ is a finite-dimen\-sional subspace of $X^*$, then ${}^\perp Y$ has finite codimension in $X$ and
$\codim ({}^\perp Y) = \dim Y.$
\end{lemma}
\begin{proof}
Let $x_1^*,\dots,x_d^*$ be a basis of $Y$ and define
$\psi:X\to\K^d$ by
$$
\psi(x):=(\lb x,x_1^*\rb,\dots,\lb x,x_d^*\rb).
$$
We claim that $\psi$ is surjective.
If its range were a proper subspace of $\K^d$, there would
exist scalars $c_1,\dots,c_d$, not all zero, such that
$\sum_{j=1}^d c_j\lb x,x_j^*\rb=0$ for all $x\in X$.
This would imply $\sum_{j=1}^d c_jx_j^*=0$, contrary to
the linear independence of $x_1^*,\dots,x_d^*$.
Thus $\psi$ is surjective.

By the definition of the pre-annihilator,
$\ker(\psi)={}^\perp Y$. Hence the induced mapping $x+{}^\perp Y\mapsto\psi(x)$ from $X/({}^\perp Y)$ to $\K^d$
is well defined and injective, and the surjectivity of
$\psi$ implies that it is also surjective. Therefore
$ \codim({}^\perp Y)
=\dim\bigl(X/({}^\perp Y)\bigr)
=d =\dim Y$.
\end{proof}

\begin{proof}[Proof of Theorem \ref{thm:Fredholm-alternative}]
Put $N:=\Ker(I-T)$ and $R:=\Ran(I-T)$.
By Lemma \ref{lem:kerran}, $N$ is finite-dimensional and
$R$ is closed. Since $T^*$ is compact by Proposition
\ref{prop:adjoint-compact}, the same lemma shows that
$\Ker((I-T)^*)$ is finite-dimensional. Proposition
\ref{prop:dual-dense-T} and Lemma \ref{lem:codim} therefore give
$$
R={}^\perp(\Ker(I-T^*)),
\qquad
\codim R=\dim\Ker(I-T^*)<\infty.
$$

By Proposition \ref{prop:findim-compl}, we may choose closed
subspaces $Y,Z\subseteq X$ such that
$$
X=N\oplus Y,\qquad X=R\oplus Z,
$$
with $\dim Z=\codim R$.
The restriction $(I-T)|_Y:Y\to R$ is bijective: it is injective
because $Y\cap N=\{0\}$, and surjective because $I-T$ vanishes
on $N$.

Choose a linear map $L:N\to Z$ of rank
$\min\{\dim N,\dim Z\}$, and let $\pi$ be the bounded
projection onto $N$ along $Y$.
Since $L\pi$ has finite rank, $T-L\pi$ is compact. Hence,
by Lemma \ref{lem:comp-inj-surj}, the operator
$$
(I-T)+L\pi=I-(T-L\pi)
$$
is injective if and only if it is surjective.

For $u\in N$ and $y\in Y$ we have
$$
((I-T)+L\pi)(u+y)=Lu+(I-T)y.
$$
Since $Lu\in Z$, $(I-T)y\in R$, and $(I-T)|_Y$ is bijective,
the direct sum decomposition $X=R\oplus Z$ gives
$$
\Ker((I-T)+L\pi)=\Ker(L),
\qquad
\Ran((I-T)+L\pi)=R\oplus\Ran(L).
$$
If $\dim N<\dim Z$, our choice of $L$ makes it injective
but not surjective, and the same is then true of $(I-T)+L\pi$.
If $\dim N>\dim Z$, it makes $L$ surjective but not
injective, and again the same is true of $(I-T)+L\pi$.
Both cases contradict the equivalence just established.
Therefore
\[
\dim\Ker(I-T)=\dim N=\dim Z=\codim\Ran(I-T).
\qedhere\]
\end{proof}

As an application of the foregoing theory we turn to the problem of finding a function $u\in C[0,1]$
solving inhomogeneous integral equations of the form
\begin{equation*}\label{eq:Hf}\tag{$H_f$} \lambda u(s)  = f(s) + \int_0^1 k(s,t)u(t)\ud t, \quad s\in [0,1].
\end{equation*}
Here $f\in C[0,1]$ is given, $k:[0,1]\times[0,1]\to\K$ is continuous, and $\la$ is a nonzero scalar.
By a {\em solution} of this equation we mean a function $u\in C[0,1]$ satisfying \eqref{eq:Hf} for all $s\in [0,1]$. In order to study existence of solutions, it is useful to also consider the homogeneous equation corresponding to $f=0$,
\begin{equation*}\label{eq:Hf0}\tag{$H_0$} \lambda u(s)  = \int_0^1 k(s,t)u(t)\ud t, \quad  s\in [0,1],
\end{equation*}
as well as the `dual' homogeneous problem
\begin{equation*}\label{eq:Hf0star}\tag{$H_0^*$} \lambda v(s)  = \int_0^1 k(t,s)v(t)\ud t, \quad  s\in [0,1].
\end{equation*}
Solutions to these problems are defined in the same way.

\begin{theorem}[Fredholm alternative for integral equations]\label{thm:FredholmAlternative}\index{Fredholm!alternative, for integral equations}
Let $k:[0,1]\times [0,1]\to\K$ be continuous and let $\la\not=0$ be fixed. Then:
\begin{enumerate}[label={\rm(\arabic*)}, leftmargin=*]
 \item\label{it:FredholmAlternative1} if the homogeneous problem \eqref{eq:Hf0} has no nonzero solution, then
for all $f\in C[0,1]$ the
inhomogeneous problem \eqref{eq:Hf} has a unique solution $u$ in $C[0,1]$;
 \item\label{it:FredholmAlternative2} if the homogeneous problem \eqref{eq:Hf0} has a nonzero solution, then it has at most finitely
many linearly independent nonzero solutions;
 \item\label{it:FredholmAlternative3} for a given $f\in C[0,1]$, the inhomogeneous problem \eqref{eq:Hf} admits a solution in $C[0,1]$
if and only if
$$\int_0^1 f(t) v(t)\ud t = 0$$
for all $v\in C[0,1]$ satisfying the dual homogeneous problem \eqref{eq:Hf0star}.
\end{enumerate}
\end{theorem}

\begin{proof}
By the result of Example \ref{ex:integral-comp} the operator $T: C[0,1]\to C[0,1]$,
$$  Tu(s):= \int_0^1 k(s,t)u(t)\ud t, \quad s\in [0,1],$$
is compact. Using this operator, we can formulate the problem \eqref{eq:Hf} abstractly as
$$ (\la - T) u = f.$$

\ref{it:FredholmAlternative1}: \
If the homogeneous problem $ (\la - T) u =0$ has no nonzero solution, then  $\la$ is not an eigenvalue.
Since we are assuming that $\la\not=0$, it follows that $\la\not\in\si(T)$ by Theorem \ref{thm:comp-spec}. Therefore, $\la-T$ is invertible
and the inhomogeneous problem $ (\la - T) u =f$ has the unique solution
$u=(\la-T)^{-1}f$.

\smallskip
\ref{it:FredholmAlternative2}: \
If the homogeneous problem $ (\la - T) u =0$ has a nonzero solution, then $\la$ is an eigenvalue,
in which case the space of solutions $u$ equals the eigenspace corresponding to $\la$, which is finite-dimensional.

\smallskip
\ref{it:FredholmAlternative3}: \ For a given $f\in C[0,1]$,
the inhomogeneous problem $ (\la - T) u =f$ has a solution $u\in C[0,1]$
if and only if $f\in \Ran(\la-T) = {}^\perp\Ker(\la-T^*)$; here we use Proposition \ref{prop:dual-dense-T}
along with the fact that $\la-T$ has closed range.
Stated differently, the
problem $ (\la - T) u =f$ has a solution  $u\in C[0,1]$ if and only if
$$\lb f,x\s\rb = 0, \quad x\s\in \ker(\la-T\s).$$
It remains to describe $\Ker(\la-T^*)$ explicitly.

We recall from Section \ref{subsec:dual-CK} that
the dual of $C[0,1]$ is the space of complex
Borel measures on $[0,1]$, the duality between functions $\phi$ and measures $\mu$ being given by
$\lb \phi,\mu\rb = \int_0^1 \phi\ud \mu$.
Let $\mu$ be such a measure. For $g\in C[0,1]$ we have, by Fubini's theorem,
\begin{align*} \lb g, T\s \mu\rb = \lb Tg,\mu\rb  & =  \int_0^1\int_0^1 k(s,t)g(t)\ud t\ud \mu(s)
  \\ & = \int_0^1 g(t)\int_0^1 k(s,t)\ud \mu(s)\ud t
  =: \int_0^1 g(t) h_\mu(t)\ud t.
 \end{align*}
Since $k$ is uniformly continuous on $[0,1]^2$, the function $h_\mu$ is
continuous. Thus $T^*\mu$ is the absolutely continuous measure with density $h_\mu$.

Suppose now that $\mu\in\Ker(\la-T^*)$. Then $\la\mu=T^*\mu$, and hence
$$
        \ud\mu(t)=\frac1\la h_\mu(t)\ud t .
$$
In particular, $\mu$ is absolutely continuous with respect to Lebesgue
measure, with a continuous density. Writing $\ud\mu(t)=v(t)\ud t$ with $v(t):=\frac1\la h_\mu(t)$,
we obtain
$$ \la v(t) =h_\mu(t)
            =\int_0^1 k(s,t)\ud\mu(s)
            =\int_0^1 k(s,t)v(s)\ud s,\quad t\in[0,1].
$$
Thus $v\in C[0,1]$ solves the dual homogeneous problem $\eqref{eq:Hf0star}$.

Conversely, let $v\in C[0,1]$ solve \eqref{eq:Hf0star},
and define the measure $\ud\mu_v(t):=v(t)\ud t$.
Then
$$
h_{\mu_v}(t)
=\int_0^1 k(s,t)v(s)\ud s
=\la v(t),\quad t\in[0,1].
$$
By the formula for $T^*\mu$ established above, this implies
$T^*\mu_v=\la\mu_v$, so $\mu_v\in\Ker(\la-T^*)$.

Thus the elements of $\Ker(\la-T^*)$ are precisely the
measures with a continuous density solving
\eqref{eq:Hf0star}. Combining this with the
abstract solvability criterion, we obtain
\ref{it:FredholmAlternative3}.
\end{proof}

\subsection{Fredholm Operators}\label{subsec:Fredholm-operators}

Let $X$ and $Y$ be Banach spaces.
The following definition is suggested by the Fredholm alternative (Theorem \ref{thm:Fredholm-alternative}):

\begin{definition}[Fredholm operators]
A bounded operator $T\in\calL(X,Y)$ is called a {\em Fredholm operator}\index{Fredholm!operator} if it has the following properties:
\begin{enumerate}[leftmargin=*, label=(\roman*)]
 \item $\dim\Ker(T) < \infty$;
 \item $\codim\Ran(T) < \infty$.
\end{enumerate}
The {\em index}\index{index} of such an operator is defined as
$$ \ind(T):= \dim \Ker(T) - \codim \Ran(T).$$
\end{definition}

\begin{example}\label{ex:Fredholm}
Here are some examples of Fredholm operators:
\begin{enumerate}[leftmargin=*, label=(\roman*)]
 \item If $T$ is a compact operator, then $I-T$ is Fredholm with index $\ind(I-T)=0$.
This is a restatement of the Fredholm alternative.

 \item The left and right shifts on $\ell^p$, $1\le p\le\infty$, are Fredholm with indices $1$ and $-1$, respectively. What changes if we replace $\ell^p$ by $\ell^p(\Z)$?

 \item For every zero-free $\phi\in C(\mathbb{T})$ the Toeplitz operator $T_\phi$ on the Hardy space $H^2(\mathbb{D})$ is Fredholm with index $ \ind(T_\phi) = - w(\phi)$, where $w(\phi)$ is the winding number of $\phi$. This is the content of the Noether--Gohberg--Krein theorem in Section \ref{subsec:Toeplitz}, where the relevant definitions can be found.
\end{enumerate}
\end{example}

We begin our analysis of Fredholm operators with the observation that such operators have closed range. As a result,
$\codim\Ran(T)$ equals the dimension of the quotient Banach space $Y/\Ran(T)$.

\begin{proposition}\label{prop:Fred-closedrange}
 If the range of a bounded operator $T\in \calL(X,Y)$ has finite codimension, then it is closed.
\end{proposition}
\begin{proof} Let $Y_0$ be a finite-dimensional subspace of $Y$ such that $\ran(T)\cap Y_0 = \{0\}$ and $\ran(T)+Y_0 = Y$. Then $Y_0$ is closed
and the bounded operator $S: X\times Y_0\to Y$ defined by
$ S(x,y_0):= Tx+y_0$
 is surjective. By the open mapping theorem, $S$ is open.
In particular, $S(X\times (Y_0\setminus \{0\}))$ is open.
Clearly this set is the complement of $S(X\times \{0\}) = \Ran(T)$ and therefore $\ran(T)$ is closed.
\end{proof}

Fredholm operators are closed under composition and the Fredholm index is additive:

\begin{proposition}\label{prop:Fred-multiplicative}
 If $T_1 \in \calL(X,Y)$ and $T_2\in \calL(Y,Z)$ are Fredholm, where $Z$ is another Banach space, then $T_2T_1\in \calL(X,Z)$ is Fredholm and
 $$ \ind(T_2T_1) = \ind(T_1)+\ind(T_2).$$
\end{proposition}
\begin{proof}
Set $R_1:=\Ran(T_1)$, $N_2: = \Ker(T_2)$, and $m:=\dim(R_1\cap N_2)$. Since
$T_1:\Ker(T_2T_1)\to R_1\cap N_2$
is surjective and has kernel $\Ker(T_1)$, we have
$$
\dim(\Ker(T_2T_1))=\dim(\Ker(T_1))+m.
$$
To compute the codimension of the range, consider the
surjective linear map
$$
Q:Y/R_1\to\Ran(T_2)/\Ran(T_2T_1)
$$
given by
$Q(y+R_1):=T_2y+\Ran(T_2T_1)$.
This is well defined because $T_2(R_1)=\Ran(T_2T_1)$.
We have $y+R_1\in \Ker(Q)$ if and only if
$T_2y=T_2r$ for some $r\in R_1$, or equivalently, if and only if
$y\in R_1+N_2$. Thus
$$
\Ker(Q)=(R_1+N_2)/R_1.
$$
The mapping $u\mapsto u+R_1$ from $N_2$ onto this kernel
has kernel $R_1\cap N_2$, so
$$
\dim\Ker(Q)=\dim N_2-m.
$$
Since $Y/R_1$ is finite-dimensional, the rank--nullity
formula applied to $Q$ gives
$$
\dim\bigl(\Ran(T_2)/\Ran(T_2T_1)\bigr)
=\codim R_1-\dim N_2+m.
$$
Adding the codimension of $\Ran(T_2)$ in $Z$, we obtain
\begin{align*}
\codim\Ran(T_2T_1)
&=\dim\bigl(\Ran(T_2)/\Ran(T_2T_1)\bigr)
  +\codim\Ran(T_2)\\
&=\codim R_1-\dim N_2+m+\codim\Ran(T_2).
\end{align*}
In particular, both the kernel and the cokernel of $T_2T_1$
are finite-dimensional, so $T_2T_1$ is Fredholm.
Subtracting the last identity from the formula for
$\dim\Ker(T_2T_1)$ gives
\[
\ind(T_2T_1)=\ind(T_1)+\ind(T_2).
\qedhere\]
\end{proof}

The following theorem characterises Fredholm operators as the operators which are `invertible up to compacts'.

\begin{theorem}[Atkinson]\label{thm:Atkinson}\index{theorem!Atkinson}
For a bounded operator $T\in \calL(X,Y)$ the following assertions are equivalent:
\begin{enumerate}[label={\rm(\arabic*)}, leftmargin=*]
 \item\label{it:Atk1} $T$ is Fredholm;
 \item\label{it:Atk2}  there exist a bounded operator $S\in\calL(Y,X)$ and compact operators $K\in \calL(X)$ and $L\in\calL(Y)$ such that
 $$ ST = I-K, \qquad TS = I-L.$$
\end{enumerate}
If these equivalent conditions hold, the operator $S$ is Fredholm with index $$\ind(S) = -\ind(T).$$ Moreover, $S$ can be chosen in such a way that $K = I-ST$ and $L = I-TS$ are finite-rank projections satisfying $\dim(\ran(K)) = \dim(\ker(T))$ and $\dim(\ran(L)) = \codim(\ran(T))$.
\end{theorem}

\begin{proof}
\ref{it:Atk1}$\Rightarrow$\ref{it:Atk2}: \
By Propositions \ref{prop:findim-compl} and \ref{prop:Fred-closedrange} there exist
closed subspaces $X_0\subseteq X$ and $Y_0\subseteq Y$ such that $\codim(X_0)<\infty$, $\dim(Y_0)<\infty$, and
$$ X = \Ker(T)\oplus X_0, \qquad Y = \Ran(T)\oplus Y_0.$$
Let $P$ and $Q$ denote the corresponding projections of $X$ and $Y$ onto $X_0$ and $\Ran(T)$, respectively.
The restriction $T_0:= T|_{X_0}$ is a bijection from $X_0$ onto $\Ran(T)$, injectivity and surjectivity both being clear. Since $\Ran(T)$ is closed, the open mapping theorem implies that the inverse mapping $S_0:= T_0^{-1}$
is bounded as an operator from $\Ran(T)$ onto $X_0$. Define $S\in \calL(Y,X)$ by $S:= S_0\circ Q$. Then for all $x\in X$ and $y\in Y$ we have
\begin{align*}
STx  = S_0 QTx = S_0Tx &= S_0 TPx = Px = x - Kx \quad \hbox{with }\ K = I-P
\intertext{and}
TSy =  TS_0Qy & = Qy = y - Ly \quad \hbox{with }\ L = I-Q.
\end{align*}
Since $I-P$ and $I-Q$ are the projections onto the finite-dimensional subspaces $\Ker(T)$ and $Y_0$,
these projections are of finite rank and hence compact. It also follows that $\dim(\ran(K)) = \dim(\ker(T))$ and $\dim(\ran(L)) = \codim(\ran(T))$.

\smallskip
\ref{it:Atk2}$\Rightarrow$\ref{it:Atk1}: \ We have
$\ker(T) \subseteq \Ker(ST)$ and hence $$\dim(\ker(T)) \le\dim(\Ker(ST)) = \dim(\ker(I-K))<\infty.$$
Likewise,
$\Ran(T)\supseteq \ran(TS)$ and hence
$$\codim(\ran(T)) \le\codim(\ran(TS)) = \codim(\ran(I-L))<\infty,$$
the finiteness of the codimension being a consequence of Theorem \ref{thm:Fredholm-alternative}.
This completes the proof of the equivalence
\ref{it:Atk1}$\Leftrightarrow$\ref{it:Atk2}.

It remains to prove the assertion about the index of $S$. Suppose that
$S\in\calL(Y,X)$ satisfies $ST=I-K$ and $TS=I-L$
with $K\in\mathscr K(X)$ and $L\in\mathscr K(Y)$. By the equivalence
just proved, applied now to $S$, the operator $S$ is Fredholm. Moreover,
$I-K$ is Fredholm with index $0$ by the Fredholm alternative. Hence,
by Proposition \ref{prop:Fred-multiplicative},
$$
   0=\ind(I-K)=\ind(ST)=\ind(S)+\ind(T).
$$
Thus $\ind(S)=-\ind(T)$.
\end{proof}

The equivalence \ref{it:Atk1}$\Leftrightarrow$\ref{it:Atk2} can be concisely stated by introducing the quotient space
$\calL(X,Y)/\mathscr{K}(X,Y).$
In the special case $X=Y$, the space $\calL(X)/\mathscr{K}(X)$ is called the {\em Calkin algebra}\index{Calkin algebra} of $X$.
Theorem \ref{thm:Atkinson} states that an operator $T\in \calL(X,Y)$ is Fredholm if and only if its equivalence class in $\calL(X,Y)/\mathscr{K}(X,Y)$ is invertible in the sense that there exists an operator $S\in \calL(Y,X)$ such that $ST = I\hbox{ mod }\mathscr{K}(X)$ and $TS = I\hbox{ mod }\mathscr{K}(Y)$.

The next three results show that Fredholmness is preserved under various operations.

\begin{proposition}\label{prop:Fredholm-comppert}
 If $T\in \calL(X,Y)$ is Fredholm and $K\in \mathscr{K}(X,Y)$ is compact, then
 $T+K$ is Fredholm and
 $$ \ind(T+K) = \ind(T).$$
\end{proposition}
\begin{proof} By Theorem \ref{thm:Atkinson} there exist $S\in\calL(Y,X)$ and compact operators $L_1\in \calL(X)$ and $L_2\in\calL(Y)$ such that
 $ ST = I-L_1$ and $TS = I-L_2$. Then
 $ S(T+K) = I-L_1 + SK = I-M_1$ with $M_1 = L_1 -SK$ compact, and
 $(T+K)S = I-L_2+KS = I-M_2$ with $M_2 = L_2-KS$ compact. Hence $T+K$ is Fredholm by
 another application of Theorem \ref{thm:Atkinson}.
Moreover, by Proposition \ref{prop:Fred-multiplicative},
$$ 0 = \ind(I-M_1) = \ind(S(T+K)) = \ind(S)+\ind(T+K),$$
so $\ind(T+K) = - \ind(S) = \ind(T)$ by the identity for indices in Theorem \ref{thm:Atkinson}.
\end{proof}

The set of Fredholm operators is open in $\calL(X,Y)$:

\begin{theorem}[Dieudonn\'e]\label{thm:Dieudonne}\index{theorem!Dieudonn\'e}
For any Fredholm operator $T\in \calL(X,Y)$ there exists a number $\delta>0$ such that for all $U\in \calL(X,Y)$ with $\n U\n < \delta$ the operator
$T+U$ is Fredholm and $$\ind(T+U) = \ind(T).$$
\end{theorem}
\begin{proof}
The proof is a variation on the proof of the openness of the set of invertible bounded operators.
By Theorem \ref{thm:Atkinson} there exist $S\in \calL(Y,X)$ and compact operators $K\in \calL(X)$ and $L\in\calL(Y)$
such that $ST = I-K$ and $TS = I-L$.
Then
$$ S(T+U) = I-K+SU, \qquad (T+U)S = I-L+US.$$
Choose $\delta>0$ such that $\n S\n\delta<1$. If $\n U\n<\delta$, then $I+SU$ and $I+US$ are boundedly invertible, and
\begin{align*}
(I+SU)^{-1} S(T+U) &= I-(I+SU)^{-1}K
\intertext{and}
(T+U)S(I+US)^{-1} &= I-L(I+US)^{-1}\!,
\end{align*}
where $M:= (I+SU)^{-1}K$ and $N:= L(I+US)^{-1}$ are compact.
The identity $(I+SU)S=S(I+US)$ gives
$$
(I+SU)^{-1}S=S(I+US)^{-1}.
$$
Together with the preceding two identities, this shows
by Theorem \ref{thm:Atkinson} that $T+U$ is Fredholm.

Next, by Proposition \ref{prop:Fred-multiplicative} and Theorem \ref{thm:Fredholm-alternative},
$$ \ind((I+SU)^{-1}) + \ind(S) + \ind(T+U) = \ind(I-M) = 0$$
and $\ind((I+SU)^{-1}) = 0$ by invertibility. By the identity for indices in Theorem \ref{thm:Atkinson} it follows that
$\ind(T+U) = -\ind(S) = \ind(T).$
\end{proof}

As an immediate corollary we obtain homotopy invariance of the Fredholm index: if $T:[0,1]\to\calL(X,Y)$ is continuous with respect to the operator norm and $T(t)$ is Fredholm for all $t\in[0,1]$, then $\ind(T(t))=\ind(T(0))$ for all $t\in[0,1]$.

\begin{proposition}\label{prop:Fred-dual}
 If $T\in \calL(X,Y)$ is Fredholm, then $T\s\in \calL(Y\s,X\s)$ is Fredholm and
  $$ \ind(T\s) = -\ind(T).$$
\end{proposition}
\begin{proof}
Since $T$ is Fredholm, its range is closed by Proposition
\ref{prop:Fred-closedrange}. Moreover,
$\Ker(T^*)=\Ran(T)^\perp$, and
Proposition \ref{prop:dualOfQuotient} identifies this
annihilator with $(Y/\Ran(T))^*$. Since $Y/\Ran(T)$ is
finite-dimensional, it follows that
$$
\dim\Ker(T^*)
=\dim\bigl((Y/\Ran(T))^*\bigr)
=\codim\Ran(T).
$$

By the closed range theorem
(Theorem \ref{thm:closed-range-dual}), we also have $\Ran(T^*)=(\Ker(T))^\perp$.
Proposition \ref{prop:dualOfSubspace} therefore gives
an isomorphism
$$
X^*/\Ran(T^*)\simeq(\Ker(T))^*.
$$
Since $\Ker(T)$ is finite-dimensional, we obtain
$$
\codim\Ran(T^*)
=\dim\bigl((\Ker(T))^*\bigr)
=\dim\Ker(T).
$$
Both dimensions are finite, so $T^*$ is Fredholm.
Combining the two identities gives
\[
\ind(T^*)
=\dim\Ker(T^*)-\codim\Ran(T^*)
=\codim\Ran(T)-\dim\Ker(T)
=-\ind(T).
\qedhere\]
\end{proof}

\subsection{The Noether--Gohberg--Krein Theorem}\label{subsec:Toeplitz}

Let $\mathbb{D}$ and $\mathbb{T}$ denote the open unit disc and unit circle in the complex plane, respectively. We think of $\mathbb{T}$ as  parametrised by $\theta\in [-\pi,\pi]$ and equipped with the normalised Lebesgue measure $\!\ud \theta/2\pi$.
 The {\em Hardy space}\index{Hardy space} $H^2(\mathbb{D})$\index{$H$@$H^{2}(\mathbb{D})$} is the Hilbert space of all holomorphic functions on $\mathbb{D}$ of the form $\sum_{n\in\N} c_n z^n$ with $$\sum_{n\in\N} |c_n|^2 < \infty.$$
Since every square-summable sequence $(c_n)_{n\in \N}$ defines a power series $\sum_{n\in\N} c_n z^n$ holomorphic on $\mathbb{D}$, the correspondence between a power series and its coefficient sequence sets up an isometric isomorphism between $H^2(\mathbb{D})$ and $\ell^2(\N)$. With respect to the norm
 $$ \n f\n_{H^2(\mathbb D)} := \Bigl(\sum_{n\in\N} |c_n|^2\Bigr)^{1/2}\!,$$
 $H^2(\mathbb D)$ is a Hilbert space.
For $n\in\Z$ consider the functions $e_n\in L^2(\mathbb{T})$ defined by $$e_n(\theta):= \exp(in\theta).$$ Since $(e_n)_{n\in\N}$ is an orthonormal sequence in $L^2(\mathbb{T})$, every square-summable sequence $(c_n)_{n\in\N}$ defines a convergent series $\sum_{n\in\N} c_n e_n$ in $L^2(\mathbb{T})$. Denoting the sum by $f$, its Fourier coefficients are given by $$\wh f(n) = \begin{cases}c_n, & n\ge 0, \\ 0, & n \le -1.\end{cases}$$
Conversely, if all negative Fourier coefficients of a function $f\in L^2(\mathbb{T})$ vanish, then
$f=\sum_{n\in\N}\wh f(n)e_n$ with convergence in $L^2(\mathbb T)$.
Since the Fourier coefficients of functions in $L^2(\mathbb{T})$ are square-summable, we obtain an isometric isomorphism between $H^2(\mathbb D)$ and the
closed subspace of $L^2(\mathbb{T})$ consisting of all functions whose negative Fourier coefficients vanish.
In what follows we identify $H^2(\mathbb{D})$ with this closed subspace of $L^2(\mathbb{T})$. As such, $H^2(\mathbb{D})$ is the range of the {\em Riesz projection}\index{Riesz projection}
$$ P: \sum_{n\in \Z}\wh f(n)e_n \mapsto \sum_{n\in\N}\wh f(n)e_n$$
on $L^2(\mathbb{T})$, which discards the terms in the Fourier series with negative indices.

Given a function $\phi\in L^\infty(\mathbb{T})$ we define the bounded operator $M_\phi$ on $L^2(\mathbb{T})$ by pointwise multiplication,
$$ M_\phi f := \phi f, \quad f\in L^2(\mathbb{T}).$$
When $M_\phi$ is applied to a function $f\in H^2(\mathbb{D})$, the resulting function $\phi f $ generally does not belong to $H^2(\mathbb{D})$, but the Riesz projection will take us back to $H^2(\mathbb{D})$. This motivates the following definition.

\begin{definition}[Toeplitz operators]\label{def:Toeplitz}
Given a function $\phi\in L^\infty(\mathbb{T})$, the {\em Toeplitz operator with symbol $\phi$}\index{Toeplitz operator}\index{symbol!of a Toeplitz operator}
is the operator $T_\phi$ on $H^2(\mathbb{D})$ given by
$$T_\phi f := P(\phi f), \quad f\in H^2(\mathbb{D}),$$
where $P$ is the Riesz projection.
\end{definition}

Every Toeplitz operator $T_\phi$ is bounded with norm
$\n T_\phi\n \le \n P\n \n M_\phi\n \le \n \phi\n_\infty$.
The Hilbert space adjoint is given by $T_\phi^\star = T_{\ov \phi}$; this follows from
$$ \iprod{T_\phi f}{g}_{H^2(\mathbb{D})} = \iprod{\phi f}{g}_{L^2(\mathbb{T})}
= \iprod{f}{\ov \phi g}_{L^2(\mathbb{T})} =  \iprod{ f}{T_{\ov\phi} g}_{H^2(\mathbb{D})} .
$$

The following theorem shows that a Toeplitz operator with continuous and zero-free symbol $\phi\in C(\mathbb{T})$
is Fredholm, and its index equals the negative of the winding number of the closed contour in $\C\setminus\{0\}$ parametrised by $\phi$.

The {\em winding number}\index{winding number} of a continuous
closed curve in $\C\setminus\{0\}$, parametrised by a function
$\phi:[-\pi,\pi]\to\C\setminus\{0\}$ with
$\phi(-\pi)=\phi(\pi)$, can be defined as follows.
One shows that there exists a continuous function
$g:[-\pi,\pi]\to\C$ such that
\begin{align}\label{eq:curve-a-theta}
\phi(\theta)=e^{ig(\theta)},\quad \theta\in[-\pi,\pi].
\end{align}
Since $e^{ig(\pi)}=e^{ig(-\pi)}$, we have
$g(\pi)-g(-\pi)\in2\pi\Z$.
We define
$$
w(\phi):=\frac{g(\pi)-g(-\pi)}{2\pi}\in\Z.
$$
This definition does not depend on the choice of $g$,
since any two choices differ by a constant integer
multiple of $2\pi$.

Put $n:=w(\phi)$ and $b(\theta):=g(\theta)-n\theta$.
Then $b(-\pi)=b(\pi)$, and
$$
h(t,\theta):=
\exp\bigl(i((1-t)b(\theta)+n\theta)\bigr),
\quad (t,\theta)\in[0,1]\times[-\pi,\pi],
$$
defines a homotopy through closed curves in
$\C\setminus\{0\}$ from $\phi$ to $e_n$.
Indeed, $h$ is continuous,
$h(t,-\pi)=h(t,\pi)$ for every $t\in[0,1]$, and
$h(0,\cdot)=\phi$ and $h(1,\cdot)=e_n$.

The integer-valued mapping $\phi\mapsto w(\phi)$ is locally
constant on the zero-free functions in $C(\mathbb T)$.
Indeed, if
$$
\n\psi-\phi\n_\infty
<
\min_{\theta\in\mathbb T}|\phi(\theta)|,
$$
then $\psi/\phi$ takes its values in $B(1;1)$ and has a
continuous logarithm with equal values at $-\pi$ and $\pi$.
Hence $w(\psi)=w(\phi)$.
It follows that $w$ is constant along homotopies.
Since $w(e_m)=m$ for every $m\in\Z$, the integer $n$ above
is the unique integer for which $\phi$ is homotopic to $e_n$.

For piecewise $C^1$ functions $\phi$, the winding number
is given analytically, as in Section \ref{sec:DC}, by
the contour integral
$$
w(\phi)
=\frac1{2\pi i}\int_\phi\frac{\ud z}{z}
=\frac1{2\pi i}\int_{-\pi}^{\pi}
  \frac{\phi'(t)}{\phi(t)}\ud t;
$$
see Problem \ref{prob:winding}.

\begin{theorem}[Noether--Gohberg--Krein]\label{thm:Gohberg-Krein}\index{theorem!Noether--Gohberg--Krein}
If the function $\phi\in C(\mathbb{T})$ is zero-free, then
the Toeplitz operator $T_\phi$ is Fredholm on $H^2(\mathbb{D})$ and
 $$ \ind(T_\phi) = - w(\phi).$$
\end{theorem}

This theorem is remarkable, as it computes an analytic quantity (the index) in terms of a topological one (the winding number). Problem \ref{prob:Cauchy-index} relates this theorem
to the Cauchy integral formula and the argument principle.

The main ingredient in the proof is the following lemma, which implies that the mapping $\phi\mapsto T_\phi$
from $C(\mathbb{T})$ to $\calL(H^2(\mathbb{D}))$ is multiplicative up to a compact operator.

\begin{lemma}\label{lem:Toep-comp}
 For all $\phi,\psi\in C(\mathbb{T})$ the operator $T_\phi T_\psi - T_{\phi\psi}$ is compact on $H^2(\mathbb{D})$.
\end{lemma}

\begin{proof}
By the estimate $\n T_\phi\n \le \n \phi\n_\infty$, the Stone--Weierstrass theorem (Theorem \ref{thm:Stone-Weierstrass}), and the fact that uniform limits of compact operators are compact (Proposition \ref{prop:limcomp}), it suffices to prove the lemma for trigonometric polynomials $\phi$ and $\psi$.

Let $\phi = \sum_{m=-M}^M c_m e_m$ and $\psi = \sum_{n=-N}^N d_n e_n$ be trigonometric polynomials. For $j\ge N$ we have
\begin{align*} \ & (T_\phi T_{\psi} - T_{\phi\psi})e_j
  \\ & =  \sum_{m=-M}^M  \sum_{n=-N}^N c_md_n P(e_m P (e_n e_j)) - \sum_{m=-M}^M \sum_{n=-N}^N c_m d_n P(e_me_n e_j)
  = 0,
\end{align*}
since $n+j\ge 0$ and hence $e_m P (e_n e_j) = e_m P (e_{n+j})= e_m e_{n+j} =e_me_n e_j$ in each summand.
By linearity and density, this shows that $ (T_\phi T_{\psi} - T_{\phi\psi})f = 0$
for all $f$ in the closed linear span of $\{e_j:\, j\ge N\}$. This implies that
$T_\phi T_{\psi} - T_{\phi\psi}$
is a finite-rank operator (of rank at most $N$) and hence compact.
\end{proof}

In particular, this lemma implies that the commutator $T_\phi T_\psi -T_\psi T_\phi$ is compact.
For functions $\phi,\psi\in C^2(\mathbb{T})$ a more precise result will be proved in Section \ref{subsec:Toeplitz-trace}.

\begin{proof}[Proof of Theorem \ref{thm:Gohberg-Krein}]
Since $\phi$ is continuous and zero-free, $1/\phi\in C(\mathbb T)$.
By Lemma \ref{lem:Toep-comp}, the operators
$T_{1/\phi}T_\phi-I$ and $T_\phi T_{1/\phi}-I$
are compact. Atkinson's theorem (Theorem \ref{thm:Atkinson})
therefore implies that $T_\phi$ is Fredholm.

It remains to compute the index. To this end, let $w(\phi) = n$ be the winding number of $\phi$ and let $h:[0,1]\times\mathbb T\to \C\setminus\{0\}$
be a homotopy from $\phi$ to $e_n$. By the continuity of the mapping $\phi\mapsto T_\phi$ and Dieudonn\'e's theorem (Theorem \ref{thm:Dieudonne}), the mapping
$ t\mapsto \ind(T_{h_t})$
is locally constant on $[0,1]$, and hence constant, where $h_t(s) := h(t,s)$ is zero-free. Since $h_0 = \phi$ and $h_1 = e_n$,
in particular we obtain $$ \ind(T_\phi) = \ind(T_{e_n}).$$
Moreover, from $$T_{e_n} \sum_{j\in\N} c_je_j = P\Bigl(\sum_{j\in\N} c_je_{j+n}\Bigr) =
\sum_{j=(-n)\vee 0}^\infty c_je_{j+n}
$$
we see that
$$\dim\ker(T_{e_n}) = \begin{cases}
                         0, & n\ge 0,\\
                         -n, & n\le -1,
                        \end{cases}
\qquad
\codim\ran(T_{e_n}) = \begin{cases}
                         n, & n\ge 0,\\
                         0, & n\le -1,
                        \end{cases}
$$
so that
$ \ind (T_{e_n}) = -n$.
\end{proof}

The following result clarifies why the symbol was assumed to be zero-free.

\begin{theorem}[Hartman--Wintner]\label{thm:HW}\index{theorem!Hartman--Wintner}
 If $\phi\in C(\mathbb{T})$ is such that the Toeplitz operator $T_\phi$ is Fredholm on $H^2(\mathbb{D})$, then $\phi$ is zero-free.
\end{theorem}
\begin{proof}
By Atkinson's theorem (Theorem \ref{thm:Atkinson}), there
exist a bounded operator $R$ and a finite-rank operator $K$
on $H^2(\mathbb D)$ such that
$RT_\phi=I-K$.
Fix a trigonometric polynomial
$g=\sum_{j=-N}^N c_je_j$. For $n\ge N$ put
$$
f_n:=e_ng=\sum_{j=-N}^N c_je_{n+j}\in H^2(\mathbb D).
$$
Since $|e_n|=\one$, we have $\n f_n\n_2=\n g\n_2$.
Moreover, for every $h\in H^2(\mathbb D)$,
$$
\iprod{f_n}{h}
=\sum_{j=-N}^N c_j\iprod{e_{n+j}}{h}\to0 \ \ \hbox{as $n\to\infty$},
$$
because the coefficients of $h$ with respect to the
orthonormal basis $(e_j)_{j\ge0}$ tend to $0$.
Thus $f_n\to0$ weakly. Since $K$ has finite-dimensional
range, it follows that $Kf_n\to0$ in norm.

By the identity $f_n=RT_\phi f_n+Kf_n$ and the contractivity of the Riesz projection $P$,
\begin{align*}
\n g\n_2=\n f_n\n_2
\le \n R\n\,\n T_\phi f_n\n_2+\n Kf_n\n_2
\le \n R\n\,\n\phi g\n_2+\n Kf_n\n_2.
\end{align*}
Letting $n\to\infty$ gives
$\n g\n_2\le\n R\n\,\n\phi g\n_2$.
Since trigonometric polynomials are dense in $L^2(\mathbb T)$
and multiplication by $\phi$ is bounded, we conclude that
$$
\n g\n_2\le\n R\n\,\n\phi g\n_2,
\quad g\in L^2(\mathbb T).
$$

Suppose now that $\phi$ vanishes at some point of $\mathbb T$.
Choose $\eps>0$ such that $\n R\n\eps<1$.
By continuity, there is a nonempty open arc $E\subseteq\mathbb T$
on which $|\phi|<\eps$.
Taking $g=\one_E$ in the preceding inequality gives
$$
\n\one_E\n_2
\le\n R\n\,\n\phi\one_E\n_2
\le\n R\n\eps\,\n\one_E\n_2
<\n\one_E\n_2,
$$
a contradiction, since $E$ has positive measure.
Thus $\phi$ is zero-free.
\end{proof}

The mapping $\phi\mapsto T_\phi$ is isometric:

\begin{proposition}\label{prop:normTphi} For all $\phi\in C(\mathbb{T})$,
 $$ \n T_\phi\n = \n \phi\n_\infty.$$
\end{proposition}
\begin{proof}
Let $U$ be multiplication by $z$ on $L^2(\mathbb T)$, and put
$Q_n:=U^{-n}PU^n$, where $P$ is the Riesz projection.
By approximation with trigonometric polynomials, $Q_n\to I$
strongly. Since $U$ commutes with $M_\phi$, we have strong convergence
$$
Q_nM_\phi Q_n=U^{-n}T_\phi PU^n\to M_\phi.
$$
The operators on the right have norm at most
$\n T_\phi\n$, so
\[
\n\phi\n_\infty=\n M_\phi\n
\le\n T_\phi\n\le\n\phi\n_\infty.
\qedhere\]
\end{proof}

In the next theorem we denote by $T_z$ the Toeplitz operator with symbol $\phi(z)=z$. Its adjoint is the Toeplitz operator
$T_z^\star = T_{\ov z}$
with symbol $\ov\phi(z)= \ov z$.
Under the identification of $H^2(\mathbb D)$ with $\ell^2(\N)$ that sends
the function $z\mapsto z^n$ to the $n$th unit vector $e_n$, the operators
$T_z$ and $T_{\ov z}$ correspond to the right and left shifts on $\ell^2(\N)$, respectively. With these identifications in mind, the following theorem can be interpreted as giving a precise description of the closed subalgebra of $\calL(\ell^2(\N))$ generated by the right and left shifts.

\begin{theorem}[Coburn]\label{thm:Coburn}\index{theorem!Coburn}
 Let $\mathscr{T}$ denote the closed subalgebra of $\calL(H^2(\mathbb{D}))$ generated by $T_z$ and $T_z^\star$. Denote by
$\mathscr{K}$ the space of compact operators on $H^2(\mathbb{D})$. Then:
\begin{enumerate}[label={\rm(\arabic*)}, leftmargin=*]
\item\label{it:Coburn0} for all $\phi\in C(\mathbb{T})$ and
$K\in\mathscr{K}$ we have $$\|T_\phi+K\|\ge \|T_\phi\|;$$
\item\label{it:Coburn1} $\mathscr{T} = \{T_\phi+K:\, \phi\in C(\mathbb{T}), \ K\in \mathscr{K}\}$;
\item\label{it:Coburn2} the mapping $\pi: T_\phi +K \mapsto \phi$ induces a multiplicative isometric isomorphism
$$ \mathscr{T}/\mathscr{K} \simeq C(\mathbb{T}).$$
\end{enumerate}
As a consequence, the representation of elements of $\mathscr T$ as the sum of a Toeplitz operator and a compact operator is unique, and we have the short exact sequence
$$ 0 \longrightarrow \mathscr{K} \longrightarrow \mathscr{T} \stackrel{\pi}{\longrightarrow} C(\mathbb{T})\longrightarrow 0.$$
\end{theorem}

Recall that a sequence of mappings is said to be {\em exact} if the range of every operator in the sequence equals the null space of the next operator in the sequence.

In the proof, as well as in later chapters, we need the following notation.
For elements $g,h$ of a Hilbert space $H$, we denote by $g\,\bar\otimes\, h$\index{$Gz$@$g\,\bar\otimes\, h$} the rank-one operator on $H$ defined by
$$ (g\,\bar\otimes\, h)x := \iprod{x}{h}g, \quad x\in H.$$
The bar in this notation serves to emphasise the fact that $\,\bar\otimes\,$ is not a tensor product, but rather its sesquilinear counterpart in the sense that for all $c\in \K$ and $x\in H$ we have
$$ (cg)\,\bar\otimes\, h = c(g\,\bar\otimes\, h), \qquad g\,\bar\otimes\, (c h) = \ov c(g\,\bar\otimes\, h).$$
For norm-one vectors $h\in H$, the operator $h\,\bar\otimes\, h$ is the orthogonal projection onto the one-dimensional subspace of $H$ spanned by $h$.

\begin{proof}[Proof of Theorem \ref{thm:Coburn}]
\ref{it:Coburn0}: \
Let $U$ be multiplication by $z$ on $L^2(\mathbb T)$ and
put $Q_n:=U^{-n}PU^n$. As shown in the proof of
Proposition \ref{prop:normTphi}, we have strong convergence
$$
Q_nM_\phi Q_n=U^{-n}T_\phi PU^n\to M_\phi.
$$

Fix $g\in L^2(\mathbb T)$. We claim that $PU^ng\to0$
weakly in $H^2(\mathbb D)$. Indeed,
$$
\iprod{PU^ng}{e_j}=\wh g(j-n)\to0,\quad j\ge0,
$$
and $\n PU^ng\n_2\le\n g\n_2$.
Since $(e_j)_{j\ge0}$ is an orthonormal basis, this proves
the claim.

If $K$ has finite rank, it follows that $KPU^ng\to0$ in norm.
The same conclusion holds for every compact $K$, by
finite-rank approximation (Proposition \ref{prop:fr-dense-compH})
and the boundedness of $(PU^ng)_{n\ge1}$.
Since $U$ is unitary, we therefore obtain
\begin{align*}
U^{-n}(T_\phi+K)PU^ng=Q_nM_\phi Q_ng+U^{-n}KPU^ng\to M_\phi g.
\end{align*}
The norm of the expression on the left-hand side is at most
$\n T_\phi+K\n\,\n g\n_2$. Passing to the limit, we obtain
$$
\n M_\phi g\n_2\le\n T_\phi+K\n\,\n g\n_2,
\quad g\in L^2(\mathbb T).
$$
Taking the supremum over $\n g\n_2\le1$ and using
Proposition \ref{prop:normTphi}, we conclude that
$$
\n T_\phi+K\n
\ge\n M_\phi\n
=\n\phi\n_\infty
=\n T_\phi\n.
$$

\ref{it:Coburn1}: \
Put
$$
\mathscr A:=\{T_\phi+K:\, \phi\in C(\mathbb T),\ K\in\mathscr K\}.
$$
By Lemma \ref{lem:Toep-comp}, $\mathscr A$ is a subalgebra containing
$T_z$ and $T_z^\star$. It is also closed, for if
$T_{\phi_n}+K_n$ converges in operator norm, then
\ref{it:Coburn0} and Proposition \ref{prop:normTphi} give
$$
\|\phi_n-\phi_m\|_\infty
=\|T_{\phi_n-\phi_m}\|
\le \|(T_{\phi_n}+K_n)-(T_{\phi_m}+K_m)\|,
$$
so the $\phi_n$ converge uniformly to some $\phi\in C(\mathbb T)$.
It follows that $K_n$ converges in operator norm to a compact operator
$K$, and the original sequence converges to $T_\phi+K$.
Thus $\mathscr A$ is closed, and hence $\mathscr T\subseteq\mathscr A$.

For the inclusion `$\supseteq$' we must prove that $\mathscr{T}$ contains all Toeplitz operators with continuous symbol and all compact operators. Given a function $\phi\in C(\mathbb{T})$, we use the Stone--Weierstrass theorem to
find a sequence of trigonometric polynomials $p_n \to \phi$ in $C(\mathbb{T})$. Then $\n T_{p_n}- T_\phi\n \le \n p_n - \phi\n_\infty \to 0$.
Since $T_{p_n}$ is a linear combination of powers of $T_z$ and $T_z^\star$, we have
$T_{p_n}\in \mathscr{T}$, and since $ \mathscr{T}$ is closed, this implies $T_\phi\in  \mathscr{T}$.

We next prove that $\mathscr{T}$ contains every compact operator. Let $S:=T_z$ for brevity, where $z$ is shorthand for the function $z\mapsto z$. We have $S^\star S=I$ and $I-SS^\star=P_0$, the orthogonal projection onto the constant functions. These identities show that $P_0\in\mathscr{T}$.

Let $e_n$ denote the function $z\mapsto z^n$, $n\in\N$. Since $\mathscr{T}$ is an algebra containing $S$, $S^\star$, and $P_0$, we have $S^mP_0(S^\star)^n\in\mathscr{T}$ for all $m,n\in\N$. For $f\in H^2(\mathbb{D})$,
$$
S^mP_0(S^\star)^nf
 = \iprod{(S^\star)^nf}{e_0}S^m e_0
 = \iprod{f}{S^ne_0}e_m
 = \iprod{f}{e_n}e_m.
$$
Thus $S^mP_0(S^\star)^n=e_m\,\bar\otimes\, e_n$.
It follows that all finite linear combinations of the operators $e_m\,\bar\otimes\, e_n$ belong
to $\mathscr T$.

Let $K\in\mathscr K$, and let $P_N$ denote the orthogonal
projection onto $\operatorname{span}\{e_0,\ldots,e_N\}$.
Then $P_NKP_N$ is a finite linear combination of the operators
$e_m\,\bar\otimes\,e_n$, and therefore belongs to $\mathscr T$.
The projections $P_N$ converge strongly to $I$ and have
norm at most one, so this convergence is uniform on compact
sets. Both $K$ and $K^\star$ are compact, and therefore
$(I-P_N) K\to 0$ and $(I-P_N)K^\star\to 0$ in operator norm.
It follows that
\begin{align*}
\n K-P_NKP_N\n
&\le\n(I-P_N)K\n+\n P_NK(I-P_N)\n
\le\n(I-P_N)K\n+\n(I-P_N)K^\star\n
\to0.
\end{align*}
Since $\mathscr T$ is closed, this implies $K\in\mathscr T$.
Thus we have shown that $\mathscr K\subseteq\mathscr T$.

 \smallskip
 \ref{it:Coburn2}: \
From \ref{it:Coburn0} it follows that
$ \n T_\phi + \mathscr{K}\n = \inf_{K\in \mathscr{K}}\n T_\phi +K\n
\ge \n T_\phi\n.$
Together with the trivial inequality
$ \n T_\phi\n \ge \inf_{K\in \mathscr{K}}\n T_\phi +K\n$ we conclude that
$$ \n T_\phi + \mathscr{K}\n = \n T_\phi\n = \n \phi\n_\infty.$$
Thus the map $T_\phi +  \mathscr{K}\mapsto \phi$ is well defined and isometric. Surjectivity is clear, hence it is an isometric isomorphism. Multiplicativity follows from Lemma \ref{lem:Toep-comp}.
\end{proof}

Based on elementary facts from the theory of $C^\star$-algebras, an alternative proof of Theorem \ref{thm:HW} as a corollary of Theorem \ref{thm:Coburn} is sketched in the Notes to this chapter.

\begin{remark} Under the identification of $H^2(\mathbb{D})$ with $\ell^2(\N)$ as indicated above, the short exact sequence of the theorem induces a short exact sequence
$$ 0 \longrightarrow \mathscr{K}(\ell^2(\N)) \longrightarrow \mathscr{T}(\ell^2(\N)) \stackrel{\ov\pi}{\longrightarrow} C(\mathbb{T})\longrightarrow 0,$$
where $\mathscr{K}(\ell^2(\N))$ and $\mathscr{T}(\ell^2(\N))$ denote, respectively, the compact operators acting on $\ell^2(\N)$ and the closed algebra generated by the left and right shifts on $\ell^2(\N)$, and $\ov\pi$ is the operator induced by $\pi$ under the identifications made.
\end{remark}

\begin{problems}

\item
Give an alternative proof of Proposition \ref{prop:limcomp} by using the equivalence of compactness and sequential compactness.

\item
If $T\in \calL(X,Y)$ is compact and $x_n\to x$ weakly, prove that $Tx_n\to Tx$ in norm.

\item
Let $X$ and $Y$ be Banach spaces. Prove that if $X$ is infinite-dimen\-sional and $T\in \calL(X,Y)$ is compact, then there exists a sequence $(x_n)_{n\ge 1}$ of norm-one vectors in $X$ such that $\limn Tx_n = 0$ in $Y$.

\item\label{prob:compact-mult-c0}
Let $(m_n)_{n\ge 1}$ be a bounded scalar sequence.
\begin{enumerate}[\rm(a), leftmargin=*]
  \item For $1\le p<\infty$, show that the multiplication operator
  $(c_n)_{n\ge 1} \mapsto (m_n c_n)_{n\ge 1}$ is compact on $\ell^p$ if and only if $\limn m_n = 0$.
  \item Does the same result hold for $\ell^\infty$? And for $c_0$?

  \noindent{\em Hint:}\ Compare with Problem \ref{prob:compact-c0}.
\end{enumerate}

\item\label{prob:Pitt-reverse}
Let $1\le p < q \le \infty.$
\begin{enumerate}[\rm(a), leftmargin=*]
  \item Prove that the inclusion mapping $\ell^p\subseteq \ell^q$ is not compact.
  \item Prove that the inclusion mapping $L^q(0,1)\subseteq L^p(0,1)$ is not compact.

  \noindent {\em Hint:}\ Consider the Rademacher functions $r_n (\theta)= {\rm sign}(\sin (2\pi 2^n \theta))$. Observe that if $m\not=n$, then $|\{|r_m - r_n| = 2\}| = \frac12$.
\end{enumerate}

\item
For $1\le p\le \infty$ consider the bounded operator $T_p:L^p(0,1)\to C[0,1]$,
$$ T_p f(t):= \int_0^t f(s)\ud s, \quad t\in [0,1].$$
\begin{enumerate} [\rm(a), leftmargin=*]
  \item Show that if $1<p\le \infty$, then $T_p$ is compact.
  \item Is $T_1$ compact?
\end{enumerate}

\item
For $f\in C_{\rm c}(0,\infty)$ and $t>0$ let
$$ (Tf)(t):= \frac1t \int_0^t f(s)\ud s.$$
\begin{enumerate}[\rm(a), leftmargin=*]
 \item Show that $Tf \in L^2(\R_+)$ for all $f\in C_{\rm c}(0,\infty)$, and the mapping $f\mapsto Tf$ thus defined has a unique extension to a bounded operator $T\in \calL(L^2(\R_+))$.
 \item Is this extension compact?
\end{enumerate}

\item
For any fixed $k\ge 1$, find a bounded operator $T$ acting on a Hilbert space such that $T^{k+1}$ is compact but $T^k$ is not.

\item
Let $g\in C[0,1]$ be given. Show that the multiplication operator $T_g: f\mapsto  fg$ on $C[0,1]$ is compact if and only if $g=0$.

\item\label{prob:Terzioglu}
Let $X$ and $Y$ be Banach spaces.
Show that an operator $T\in\calL(X,Y)$ is compact
if and only if there is a sequence $(x_n^*)_{n\ge 1}$ in $X^*$ such that
$ \lim_{n\to\infty} x_n^*= 0$ in $X^*$ and
$$ \| Tx\| \le \sup_{n\ge 1} |\lb x,x_n^*\rb|, \quad x\in X.$$

\noindent{\em Hint:}\
For the `if' part, regard $T$ as the composition
$x\mapsto(\lb x,x_n^*\rb)_{n\ge1}\mapsto Tx$
and argue as in Problem \ref{prob:compact-mult-c0};
for the `only if' part use the result of Problem \ref{prob:compact-nullseq} and the compactness of $T^*$\!.

\item
Let $X$ be a reflexive Banach space.
\begin{enumerate}[\rm(a), leftmargin=*]
 \item Show that every operator $T\in \calL(X,\ell^1)$ is compact.

 \noindent{\em Hint:}\ Use the result of Problem \ref{prob:Schur}.

 \item Show that every operator $T\in \calL(c_0,X)$ is compact.
\end{enumerate}

\item
Show that a compact operator has closed range if and only if it is of finite rank.

\noindent{\em Hint:}\ In one direction, use the open mapping theorem.

\item\label{prob:comp-ONS-null}
Let $H$ be an infinite-dimensional Hilbert space with orthonormal basis $(h_n)_{n\ge 1}$.
\begin{enumerate}[\rm(a), leftmargin=*]
 \item Show that for all $T\in \calL(H)$ and $x\in H$ we have $\lim_{n\to\infty} (Th_n|x) = 0$.

 \noindent {\em Hint:}\ Consider $T^\star x$.
 \item Show that if $T$ is compact, then $\limn \n Th_n\n = 0$.
\end{enumerate}

\noindent{\em Remark.} A converse to this result will be proved in Problem \ref{prob:comp-ONS-null-converse} in Chapter \ref{chap:spectral-theorem}.

\item\label{prob:contains-FR}
Let $X$ be a Banach space.
Show that if $\calI$  is a nonzero two-sided ideal in $\calL(X)$, then $\calI$ contains all finite-rank operators on $X$.

\noindent{\em Hint:}\ For given $x \in X$ and $x\s \in X^*$ define the rank-one operator $x \otimes x\s$ on $X$ by $$(x\otimes x\s)y := \lb y,x\s\rb x, \quad y\in X.$$
Show that if $x_0\in X$,  $x_0\s\in X\s$, and  $T\in\calL(X)$ are such that $\lb Tx_0,x_0^*\rb=1$, then
$$x \otimes x\s = (x \otimes x_0\s) \circ T \circ (x_0 \otimes x\s).$$
Now use the result of Problem \ref{prob:FR-ops}.

\item
Let $X$ be a Banach space and let $T\in \calL(X)$ be compact.
\begin{enumerate}[\rm(a), leftmargin=*]
  \item\label{it:comp-spec-a} Use Proposition \ref{prop:approx-eigenvalue-bdry} to show that every nonzero $\la\in \partial \sigma(T)$ is an eigenvalue.
  \item Use part \ref{it:comp-spec-a}
  and part \ref{it:comp-spec2} of Theorem \ref{thm:comp-spec} to deduce that every nonzero $\la\in\sigma(T)$ is an eigenvalue.
 \end{enumerate}

\item\label{prob:alg-mult}
Let $T$ be a compact operator on a Banach space $X$, let $0\not=\la\in\si(T)$, and let $\nu$ be its algebraic multiplicity. Let $X_{\la}:= P_\la X$ be the range of the spectral projection associated with the point $\la$.
Prove the following assertions:
\begin{enumerate}[\rm(a), leftmargin=*]
  \item a vector $x\in X$ belongs to $X_{\la}$ if and only if $(\la-T)^k x=0$ for some $k\ge 1$;
  \item for all $x\in X_{\la}$ we have $(\la-T)^{\nu}x=0$;
  \item $X_{\la} = \ker(\la-T)^{\nu}$.
\end{enumerate}

\item\label{prob:Calkin}
Let $X$ be a Banach space. In this problem we write $[T]$ for the element $T + \mathscr{K}(X)$ of the Calkin algebra $\calL(X)/\mathscr{K}(X)$.
Show that  the multiplication $[S]\circ[T]:= [ST]$ is well defined on $\calL(X)/\mathscr{K}(X)$ and satisfies
$$\n[S]\circ[T]\n_{\calL(X)/\mathscr{K}(X)}\le \n [S]\n_{\calL(X)/\mathscr{K}(X)}\n [T]\n_{\calL(X)/\mathscr{K}(X)}.$$
In the terminology introduced in the Notes to this chapter, this shows that the Calkin algebra $\calL(X)/\mathscr{K}(X)$ is a (unital) Banach algebra.

\item
Let $X$ be a Banach space and $T\in \calL(X)$ be a bounded operator. Show that if $T^k$ is compact for some integer $k\ge 1$,
then $I+T$ is Fredholm. What is its index?

\item
Show that a continuous function satisfying \eqref{eq:curve-a-theta} always exists.

\item\label{prob:winding}
Prove that the two definitions of the winding number of a piecewise $C^1$ curve discussed in Section \ref{sec:Fredholm-theory} agree.

\item
Show that an operator $T\in\mathscr{T}$ is Fredholm if and only if $T = T_\phi + K$, where $\phi\in C(\mathbb{T})$ is zero-free and $K\in \mathscr{K}(H^2(\mathbb{D}))$.

\item\label{prob:Cauchy-index}
Let $T_\phi$ denote the Toeplitz operator on $H^2(\mathbb D)$
with symbol $\phi$.

\begin{enumerate}[\rm(a), leftmargin=*]
\item\label{it:Cauchy-index1} For $a\in\mathbb D$, put $\phi_a(\zeta):=\zeta-a$.
Use the Cauchy integral formula to show that
$$
(T_{1/\phi_a}f)(z)
=\frac{f(z)-f(a)}{z-a},
\quad f\in H^2(\mathbb D),\ z\in\mathbb D\setminus\{a\}.
$$
Deduce that
$$
T_{1/\phi_a}T_{\phi_a}=I,
\qquad
T_{\phi_a}T_{1/\phi_a}=I-E_a,
$$
where $E_af:=f(a)\one$, and conclude that $T_{\phi_a}$
is Fredholm with index $-1$.

\noindent{\em Hint:}\
First let $f$ be a polynomial and use
the Cauchy integral representation
$$
(Pg)(z)=\frac1{2\pi i}\int_{\mathbb T}
\frac{g(\zeta)}{\zeta-z}\ud\zeta
$$
of the Riesz projection, identifying its range with
$H^2(\mathbb D)$. Then use density.

\item\label{it:Cauchy-index2} Let $p$ be a polynomial with no zeros on $\mathbb T$.
Deduce that $T_p$ is Fredholm and
$$
\ind(T_p)
=-N_{\mathbb D}(p)
=-\frac1{2\pi i}\int_{\mathbb T}
  \frac{p'(\zeta)}{p(\zeta)}\ud\zeta = -w(p),
$$
where $N_{\mathbb D}(p)$ counts the zeros in $\mathbb D$
with multiplicities and $\mathbb T$ is traversed
anticlockwise.

{\em Hint:}\ Factor $p$ and use the additivity of the
Fredholm index.
\end{enumerate}

\end{problems}

%% file: ch08-HilbertOperators.tex
\chapter{Bounded Operators on Hilbert Spaces}\label{chap:Hilbert-operators}

\blfootnote{\rm \small \begin{quote}This book has been published by Cambridge University Press in the series `Cambridge Studies in Advanced Mathematics'. The present version is free to view and download for personal use only. Not for re-distribution, re-sale or use in derivative works.\newline {\copyright} Jan van Neerven\end{quote}}

\noindent
The identification of a Hilbert space $H$ with its dual via the Riesz representation theorem makes it possible to consider a bounded operator $T$ and its adjoint simultaneously on $H$. This leads to the important classes of selfadjoint, unitary, and normal operators. Their spectral theory is particularly rich. Its full power comes to bear only in the next chapter, where we prove the spectral theorem for bounded normal operators. The present chapter discusses the elementary theory and, for normal operators $T$, establishes a generalisation of the holomorphic calculus to a calculus for continuous functions on the spectrum $\sigma(T)$. Using this calculus, we prove a number of nontrivial results such as the existence of a unique positive square root of a positive operator and a polar decomposition for general bounded operators. In the last section we establish the celebrated Sz.-Nagy theorem on the existence of unitary dilations for Hilbert space contractions.

\section{Selfadjoint, Unitary, and Normal Operators}\label{sec:sa-un-no}

Throughout this chapter, $H$ is a complex Hilbert space.
The following proposition is key to several proofs in this chapter.
The example of rotation by $\frac12\pi$ in $\R^2$ shows that its counterpart for real Hilbert spaces fails.

\begin{proposition}\label{prop:polarisation}
If $T\in \calL(H)$ satisfies $\iprod{Tx}{x} = 0$ for all $x\in H$, then  $T = 0$.
\end{proposition}
\begin{proof}
For all $x,y\in H$, from $\iprod{T(x+y)}{x+y} = 0$ we obtain
\begin{align}\label{eq:polarise1} \iprod{Tx}{y} + \iprod{Ty}{x} =0.
\end{align}
Replacing $y$ by $iy$, we obtain
$-i \iprod{Tx}{y} + i\iprod{Ty}{x} =0$.
Multiplying both sides by $i$ gives
\begin{align}\label{eq:polarise2} \iprod{Tx}{y} - \iprod{Ty}{x} =0.
\end{align}
Adding \eqref{eq:polarise1} and \eqref{eq:polarise2} gives $\iprod{Tx}{y} =0$ for all $x,y\in H$. This implies the result.
\end{proof}

The trick used in the proof is called {\em polarisation}.\index{polarisation}

In Proposition \ref{prop:HS-adjoint} it was shown that if $H$ and $K$ are Hilbert spaces and $T\in \calL(H,K)$ is a bounded operator, then there exists a unique bounded operator $T^\star\in \calL(K,H)$, the {\em Hilbert space adjoint}\index{adjoint!of a bounded Hilbert space operator} of $T$, such that
$$\iprod{Tx}{y} = \iprod{x}{T^\star y}, \quad x\in H, \ y\in K.$$
Furthermore,
\begin{align}\label{eq:HS-adjoint}\n T\n = \n T^\star\n = \n T^\star T\n^{1/2}\!.
\end{align}
The existence of Hilbert space adjoints permits the introduction of several interesting classes of Hilbert space operators.

\begin{definition}[Normal, unitary, selfadjoint, and positive operators]
 An operator $T\in \calL(H)$ is called:
 \begin{itemize}
  \item {\em positive},\index{positive!operator, Hilbertian} if $\iprod{Tx}{x}\ge 0$ for all $x\in H$;
  \item {\em selfadjoint},\index{selfadjoint operator} if $T = T^\star$;
  \item {\em unitary},\index{unitary!operator} if $TT^\star = T^\star T = I$;
  \item {\em normal},\index{normal!operator} if $TT^\star = T^\star T$.
 \end{itemize}
\end{definition}

Every positive operator is selfadjoint: for if $T$ is positive, then for all $x\in H$ we have $\iprod{Tx}{x}\ge0$ and therefore $\iprod{T^\star x}{x} = \ov{\iprod{x}{T^\star x}} = \ov{\iprod{Tx}{x}}= \iprod{Tx}{x}$. Proposition \ref{prop:polarisation} now implies that $T=T^\star$. As the example preceding the statement of the proposition shows, it is important here to work over the complex scalar field.
Selfadjoint operators and unitary operators are normal.

The classes of positive, selfadjoint, unitary, and normal operators can be viewed as operator analogues of the positive real numbers, the real numbers,
the complex numbers of modulus one, and the complex numbers, respectively. A number of results support this view:
\begin{itemize}
 \item[--] every selfadjoint operator is the difference of two positive operators;
 \item[--] every invertible operator is the composition of an invertible positive operator and a unitary operator;
 \item[--] a bounded operator is unitary if and only if it is the complex exponential of a self\-adjoint operator.
\end{itemize}
The first result follows from
$$
T=(T+\n T\n I)-\n T\n I,
$$
since both operators on the right-hand side are positive:
for selfadjoint $T$, the Cauchy--Schwarz inequality gives
$\iprod{Tx}{x}\ge-\n T\n\,\n x\n^2$ for all $x\in H$.
A decomposition into positive operators whose product is
zero is obtained in Problem \ref{prob:sa-diff-2pos}.
The second result and the `if' part of the third are proved
in the present chapter, and the `only if' part of the third
result is a consequence of the spectral theorem in the next
chapter (Problem \ref{prob:UT}).

The following spectral characterisations will be proved in Corollary \ref{cor:normal-spec-char}:
 \begin{itemize}
 \item[--] a normal operator is unitary if and only if
 its spectrum is contained in the unit circle;
 \item[--] a normal operator is selfadjoint if and only if
 its spectrum is contained in the real line;
 \item[--] a normal operator is positive if and only if
 its spectrum is contained in $[0,\infty)$;
 \item[--] a normal operator is an orthogonal projection if and only if
 its spectrum is contained in $\{0,1\}$.
 \end{itemize}
The last of these results indicates that orthogonal
projections can be viewed as the analogues of the `Boolean' set $\{0,1\}$.
It also implies that normal projections are orthogonal.

Let us begin by proving an operator analogue of the decomposition of a complex number into real and imaginary parts.

\begin{proposition}\label{prop:dec-T-sa} For every operator $T\in \calL(H)$
there exist unique selfadjoint operators $A,B\in \calL(H)$ such that $T = A+ iB$.
\end{proposition}
\begin{proof}
The operators $A:= \frac12(T+T^\star)$ and $B:=\frac1{2i}(T-T^\star)$ are selfadjoint and $T = A+iB$.
Suppose also that $T=A'+iB'$ with $A'$ and $B'$ selfadjoint.
Taking adjoints gives $T^\star=A'-iB'$. Adding and
subtracting these identities, we obtain
$$
A'=\frac12(T+T^\star)=A,\qquad
B'=\frac1{2i}(T-T^\star)=B.
$$
This proves uniqueness.
\end{proof}

A complex number satisfies $|z| = 1$ if and only if there is a real number $x$ such that $z = e^{ix}$.
The operator analogue of the  `if' part is contained in the next proposition.

\begin{proposition}\label{prop:sa-unitary}
 If $T\in \calL(H)$ is selfadjoint, then $e^{iT}$ is unitary.
\end{proposition}
\begin{proof}
From the expansion $e^{iT} = \sum_{n=0}^\infty \frac{i^n}{n!}T^n$ we see that
$$(e^{iT})^\star = \sum_{n=0}^\infty \frac{(-i)^n}{n!}T^n = e^{-iT}\!.$$
It is elementary to check that $e^{iT}e^{-iT} = e^{-iT}e^{iT} = I$ by writing out the defining power series and multiplying them.
Alternatively, this identity follows from the multiplicativity of the entire calculus of $T$ applied with $f(z) = \exp(iz)$ and $g(z) = \exp(-iz)$.
\end{proof}

We have the following simple characterisation of unitary operators:

\begin{proposition}\label{prop:unitary}
For an operator $U\in \calL(H)$ the following assertions are equivalent:
 \begin{enumerate}[label={\rm(\arabic*)}, leftmargin=*]
  \item\label{it:unitary1} $U$ is unitary;
  \item\label{it:unitary2} $U$ is surjective and $\n Ux\n = \n x\n $ for all $x\in H$;
  \item\label{it:unitary3} $U$ is surjective and $\iprod{Ux}{Uy} = \iprod{x}{y}$ for all $x,y\in H$.
 \end{enumerate}
\end{proposition}
\begin{proof}
\ref{it:unitary1}$\Rightarrow$\ref{it:unitary3}: \ If $U$ is unitary, then $U$ is invertible (with inverse $U\inv = U^\star$) and therefore
 $U$ is surjective. Moreover, $\iprod{Ux}{Uy} = \iprod{x}{U^\star Uy} = \iprod{x}{y}$.

 \smallskip
\ref{it:unitary3}$\Rightarrow$\ref{it:unitary2}: \ Take $x=y$.

 \smallskip
\ref{it:unitary2}$\Rightarrow$\ref{it:unitary1}: \ We have $$\iprod{U^\star Ux}{x} = \iprod{Ux}{Ux} = \n Ux\n^2 = \n x\n^2 = \iprod{x}{x}$$ for all $x\in H$ and therefore
 $U^\star U = I$. It follows that $U^\star$ is a left inverse to $U$. The assumptions further imply that
 $U$ is surjective and injective, hence invertible. The inverse must be equal to the left inverse,
 which is therefore $U^\star$\!. It follows that also $UU^\star =I$.
\end{proof}

The right shift on $\ell^2$ shows that the surjectivity assumptions cannot be omitted.

\begin{example} The left and right shifts on $\ell^2(\Z)$ are unitary. Indeed, the adjoint of the left (right) shift is the right (left) shift, so in either case the adjoint equals the inverse. Similarly, left and right translations on $L^2(\R)$ are unitary.
\end{example}

\begin{example} The Fourier--Plancherel transform on $L^2(\R^d)$ is unitary; this follows from Theorem \ref{thm:FT-main} and Proposition \ref{prop:unitary}.
\end{example}

Projections in Hilbert spaces are orthogonal if and only if they are selfadjoint:

\begin{proposition}\label{prop:orth-proj} For a projection $P\in\calL(H)$ the following assertions are equivalent:
\begin{enumerate}[label={\rm(\arabic*)}, leftmargin=*]
 \item\label{it:orth-proj1}
$P$ is {\em orthogonal},\index{orthogonal!projection}
that is, its null space and range are orthogonal;
\item\label{it:orth-proj2} $P$ is selfadjoint.
\end{enumerate}
\end{proposition}
\begin{proof}
\ref{it:orth-proj1}$\Rightarrow$\ref{it:orth-proj2}:  \ If $P$ is orthogonal, then $x-Px \perp Py$ for all $x,y\in H$, noting that
$x-Px\in \ker(P)$ (since $P(x-Px) = Px-P^2 x = Px-Px = 0$) and $Py\in \Ran(P)$.
Therefore,
 $$\iprod{x}{Py} = \iprod{Px}{Py} + \iprod{x-Px}{Py} = \iprod{Px}{Py}$$
 and similarly $$\iprod{Px}{y} = \iprod{Px}{Py}+ \iprod{Px}{y-Py} = \iprod{Px}{Py},$$
 so $\iprod{Px}{y} = \iprod{x}{Py}$ and $P$ is selfadjoint.

 \smallskip
\ref{it:orth-proj2}$\Rightarrow$\ref{it:orth-proj1}: \ If $P$ is a selfadjoint projection, then
 $$\iprod{x-Px}{Py} = \iprod{P^\star(x-Px)}{y} = \iprod{P(x-Px)}{y} = 0,$$
 since $P=P^2$\!. Since every element in $\Ker(P)$ is of the form $x-Px$, this shows that $\ker(P)\perp\Ran(P)$, that is, the projection $P$ is orthogonal.
\end{proof}

We now turn to the study of some spectral properties of Hilbert space operators.
From Proposition \ref{prop:spectrum-dual} we recall that for every bounded operator $T$ on a Banach space we have
$\sigma(T\s)= \sigma(T).$
A similar result holds for the spectrum of the Hilbert space adjoint $T^\star$:

\begin{proposition}\label{prop:spectrum-HS-dual} For all $T\in \calL(H)$ we have
$$\sigma(T^\star)  = \ov{\sigma(T)},$$
where the bar denotes complex conjugation.
\end{proposition}

\begin{proof}
The proof is similar to that of Proposition \ref{prop:spectrum-dual}, but is simpler because of the
Riesz representation theorem. The idea is to prove that $\la\in\rh(T)$ if and only if $\ov \la\in\rh(T^\star)$, and that in this case $$ (R(\la,T))^\star = R(\ov \la,T^\star).$$

First suppose that $\la\in\rh(T)$. Then
$$ (\ov\la-T^\star)(R(\la,T))^\star = (\la-T)^\star(R(\la,T))^\star =  (R(\la,T) (\la-T))^\star = I^\star= I.$$ Similarly, $(R(\la,T))^\star (\ov\la-T^\star)  = I$. It follows that $\ov\la\in \rh(T^\star)$ and
$R(\ov\la,T^\star) = (R(\la,T))^\star$\!.

If $\la\in \rh(T^\star)$, applying what we just proved to $T^\star$ gives
 $\ov\la\in \rh(T^{\star\star}) = \rh(T)$ and
$R(\ov\la,T) =R(\ov\la,T^{\star\star}) = (R(\la,T^\star))^\star$\!.
\end{proof}

For unitary operators we have the following simple result.

\begin{proposition}\label{prop:spectrum-unitary} If $U\in \calL(H)$ is unitary, then $\sigma(U)$ is contained in the unit circle.
\end{proposition}
\begin{proof}
If $|\la|>1$, then $\n\la^{-1}U\n\le|\la|^{-1}<1$,
so the Neumann series shows that
$\la-U=\la(I-\la^{-1}U)$ is invertible.
If $|\la|<1$, then
$$
\la-U=-U(I-\la U^\star).
$$
Here $U$ is invertible, and $I-\la U^\star$ is invertible
by the Neumann series since $\n\la U^\star\n\le|\la|<1$.
Thus every $\la$ with $|\la|\not=1$ belongs to $\rh(U)$,
proving the assertion.
\end{proof}

Another proof is outlined in Problem \ref{prob:unitary-spectrum}.
In the converse direction, a normal operator whose spectrum is contained in the unit circle is unitary. The proof of this fact is harder and will be given in Corollary \ref{cor:normal-spec-char}.

The next result describes the spectrum of selfadjoint operators.

\begin{theorem}[Spectrum of selfadjoint operators]\label{thm:spect-sa}
An operator $T\in \calL(H)$ is selfadjoint if and only if $\iprod{Tx}{x}\in \R$ for all $x\in H$.
If $T$ is selfadjoint and $H\not=\{0\}$, then $$\n T\n = \sup_{\n x\n\le 1} |\iprod{Tx}{x}| = \max\{|m|,|M|\}$$
and
$$ \{m,M\}\subseteq \sigma(T)\subseteq [m,M],$$  where
$ m := \inf_{\n x\n= 1} \iprod{Tx}{x}$ and $M := \sup_{\n x\n= 1} \iprod{Tx}{x}$. The infimum is attained if and only if $m$ is an eigenvalue, the supremum is attained if and only if $M$ is an eigenvalue, and in both cases the attaining norm-one vectors are precisely the normalised eigenvectors.
\end{theorem}
\begin{proof}
If $\iprod{Tx}{x}\in \R$, then $\iprod{T^\star x}{x} = \ov{\iprod{x}{T^\star x}}
= \ov{\iprod{Tx}{x}} = \iprod{Tx}{x}$. Hence, if $\iprod{Tx}{x}\in \R$ for all $x\in H$,
then $T = T^\star$ by Proposition \ref{prop:polarisation} applied to $T-T^\star$\!. Conversely, if $T = T^\star$\!,
then $ \iprod{Tx}{x} = \iprod{x}{Tx} = \ov{\iprod{Tx}{x}}$ and therefore $\iprod{Tx}{x}\in \R$.

Suppose that $T$ is selfadjoint and $H\not=\{0\}$. To prove that
$\si(T)\subseteq[m,M]$, let $\la\notin[m,M]$ and put
$$
d:=\inf_{t\in[m,M]}|\la-t|>0.
$$
For every nonzero $x\in H$, the definitions of $m$ and $M$ give
$\iprod{Tx}{x}/\n x\n^2\in[m,M]$. Hence the Cauchy--Schwarz
inequality gives
\begin{align*}
\n(\la-T)x\n\,\n x\n
&\ge |\iprod{(\la-T)x}{x}|
=\Bigl|\la-\frac{\iprod{Tx}{x}}{\n x\n^2}\Bigr|\n x\n^2
\ge d\n x\n^2.
\end{align*}
Thus $\n(\la-T)x\n\ge d\n x\n$ for all $x\in H$, and
Proposition \ref{prop:closed-range} implies that $\la-T$ is
injective and has closed range. Since $\ov\la\notin[m,M]$,
the same argument shows that $(\la-T)^\star=\ov\la-T$ is
injective. By Proposition \ref{prop:injective-denserange},
$\la-T$ therefore has dense range. Consequently, $\la-T$ is
invertible, proving that $\si(T)\subseteq[m,M]$.

We next prove the norm identities. Put
$$
\mu:=\max\{|m|,|M|\} = \sup_{\n x\n\le1}|\iprod{Tx}{x}|
=\sup_{\n x\n=1}|\iprod{Tx}{x}|.
$$
Clearly $\mu\le\n T\n$, and by homogeneity
$|\iprod{Th}{h}|\le \mu\n h\n^2$ for all $h\in H$.
To prove the converse inequality, let $\n x\n=1$ and suppose
that $Tx\not=0$. Put $y:=Tx/\n Tx\n$, so that $\n y\n=1$.
Selfadjointness and the parallelogram identity give
\begin{align*}
4\n Tx\n
=4\Re\iprod{Tx}{y}
& =\iprod{T(x+y)}{x+y}-\iprod{T(x-y)}{x-y}\\
&\le \mu\bigl(\n x+y\n^2+\n x-y\n^2\bigr)
=4\mu.
\end{align*}
Thus $\n Tx\n\le \mu$; this also holds if $Tx=0$.
Taking the supremum over all norm-one vectors gives
$\n T\n\le \mu$, and hence $\n T\n=\mu$.

It remains to prove (i) that $m,M\in\si(T)$, and (ii) that $m$ and $M$ are eigenvalues if and only if the defining infimum (for $m$) or supremum (for $M$) is attained, and that the attaining vectors are the normalised eigenvectors (for $m$ and $M$, respectively).
We prove this for $M$; the result for $m$ follows by considering $-T$.

We begin with (i) and prove that $M\in\si(T)$.
Replacing $T$ by $T-m$, we may assume that $0= m\le M$. Then $\iprod{Tx}{x}\ge 0$ for all $x\in H$ and therefore
$$M = \sup_{\n x\n= 1} \iprod{Tx}{x} = \sup_{\n x\n= 1} |\iprod{Tx}{x}| =\n T\n.$$
Choose a sequence $(x_n)_{n\ge 1}$ of norm-one vectors such that $\limn \iprod{Tx_n}{x_n} = M$.
Then,
\begin{align*}
 \n (M-T)x_n\n^2 & = \iprod{(M-T)x_n}{(M-T)x_n} \\ & = M^2 \n x_n\n^2 -2M\iprod{Tx_n}{x_n}+ \n Tx_n\n^2
 \\ & \le  M^2  -2M\iprod{Tx_n}{x_n}  + \n T\n^2
  = M^2  -2M\iprod{Tx_n}{x_n}  + M^2\!,
\end{align*}
which tends to $M^2 - 2M^2 + M^2 = 0$ as $n\to\infty$. This implies that
$M$ is an approximate eigenvalue of $T$.

Turning to (ii), suppose that $x_0\in H$ is a norm-one vector such that $\iprod{Tx_0}{x_0}=M$. Taking the constant-$x_0$ sequence in the preceding argument gives $Tx_0 = Mx_0$, that is, $M$ is an eigenvalue with normalised eigenvector $x_0$. This proves the `if' part. The `only if' part is trivial.
\end{proof}

A short alternative proof of the spectral inclusion $\sigma(T)\subseteq \R$ is obtained by combining Propositions \ref{prop:sa-unitary} and \ref{prop:spectrum-unitary} with the spectral mapping theorem: since $T$ is selfadjoint, the operator $e^{iT}$ is unitary; consequently, $\sigma(e^{iT}) = e^{i\sigma(T)}$ is contained in the unit circle and therefore $\sigma(T)$ must be real. Yet another proof, also based on Proposition \ref{prop:spectrum-unitary}, is outlined in Problem \ref{prob:selfadjoint-spectrum}.
In the converse direction, a normal operator whose spectrum is contained in the real line is selfadjoint. This will be proved in Corollary \ref{cor:normal-spec-char}.

\begin{remark}[Rayleigh--Ritz principle]\label{rem:RR} The fact that a minimiser of $\inf_{\n x\n=1} \iprod{Tx}{x}$ is an eigenvector of $T$ with eigenvalue equal to the attained minimum can be interpreted as solving the problem of minimising $\iprod{Tx}{x}$ under the constraint $\n x\n=1$. This conclusion can also be obtained by a variational argument. Set $m:= \inf_{\n x\n=1} \iprod{Tx}{x}$ and consider the `Lagrangian'
$$
L(x) := \iprod{Tx}{x} - m(\|x\|^2-1),
$$
where the term $\|x\|^2-1$ encodes the constraint that $x$ should be a unit vector.

Suppose that $x_0$ is a norm-one vector satisfying
$\iprod{Tx_0}{x_0}=m$. By the definition of $m$ and homogeneity,
$\iprod{Tx}{x}\ge m\n x\n^2$ for all $x\in H$. Consequently,
$$
L(x)\ge m=L(x_0),\quad x\in H.
$$
Thus $x_0$ minimises $L$ on all of $H$ and is therefore a
critical point of $L$.

Let $h\in H$ and consider the perturbation $x_\eps := x_0+\eps h$, where $\eps$ is a small real parameter. Then
\begin{align*}
L(x_\eps)  & = \iprod{T(x_0+\eps h)}{x_0+\eps h} - m(\|x_0+\eps h\|^2-1)
\\ & = \iprod{Tx_0}{x_0} - m \n x_0\n^2 + m+ 2\eps \Re \iprod{Tx_0}{h}
\\ & \qquad\qquad -2m \eps\Re\iprod{x_0}{h}  + \eps^2\iprod{Th}{h}- m\eps^2\n h\n^2\!.
\end{align*}
Differentiating at $\eps=0$ gives
$$
\frac{\ud}{\ud\eps}L(x_\eps)\Big|_{\eps=0}
=2\Re\iprod{Tx_0}{h}-2m\Re\iprod{x_0}{h}.
$$
If $x_0\in H$ is a {\em critical point} of $L$, that is, a point where the derivative vanishes for all directions $h$, this results in the `Euler--Lagrange equation'\index{Euler--Lagrange equation}
$$ \Re\iprod{Tx_0}{h} = m\Re\iprod{x_0}{h}.$$
Taking $h=(T-m)x_0$ in the Euler--Lagrange equation gives
$$
0=\Re\iprod{(T-m)x_0}{(T-m)x_0}
=\n(T-m)x_0\n^2.
$$
Hence $Tx_0=mx_0$. Thus every nonzero critical point of $L$
is an eigenvector of $T$ with eigenvalue $m$, the greatest lower bound of the spectrum.

In finite dimensions, this principle is known as the {\em Rayleigh--Ritz principle}.\index{Rayleigh--Ritz principle} For compact positive operators it can be considerably sharpened
(see Theorem \ref{thm:minmax}).
\end{remark}

It is an immediate consequence of Theorem \ref{thm:spect-sa}
that the norm of a selfadjoint operator equals its spectral radius. More generally, this is true for normal operators; see Proposition \ref{prop:normal-spectralradius}.
This equality of norm and spectral radius
can sometimes be used to determine the norm of an operator. We illustrate this by
determining the norm of the Volterra operator.

\begin{example}[Volterra operator]\label{ex:Volterra2}\index{Volterra operator}
From Example \ref{ex:Volterra} we recall that the {\em Volterra operator} is the operator $T\in \calL(L^2(0,1))$
 given by the indefinite integral
$$
 (T f)(s) := \int_0^s f(t)\ud t,\quad f\in L^2(0,1), \ s\in (0,1).
$$
The operator $T$ fails to be selfadjoint (it even fails to be normal, see  Problem \ref{prob:Volterra}), but we have $\n T\n = \n S\n$, where $S\in \calL(L^2(0,1))$
is defined by
 $$
 (S f)(s) := \int_0^{1-s} f(t)\ud t,\quad f\in L^2(0,1), \ s\in (0,1),
$$
as is immediate from the identity $(Sf)(s) = (Tf)(1-s)$.
For $f,g\in L^2(0,1)$, Fubini's theorem gives
\begin{align*}
\iprod{Sf}{g}
=\int_0^1\int_0^{1-s}f(t)\ov{g(s)}\ud t\ud s
=\int_0^1 f(t)
  \Bigl(\int_0^{1-t}\ov{g(s)}\ud s\Bigr)\ud t
=\iprod{f}{Sg},
\end{align*}
which shows that $S$ is selfadjoint.
By Example \ref{ex:integral-comp2} $S$ is compact, and therefore by
Theorem \ref{thm:comp-spec} every nonzero $\lambda\in\sigma(S)$ is an eigenvalue. To compute the spectral radius
$r(S)$ we therefore have to determine the set of nonzero eigenvalues of $S$.

Suppose that $\lambda\not=0$ is an eigenvalue of $S$ and let $f$ be an eigenfunction.
Then
\begin{align}\label{eq:volterra-id}
f(s) = \frac1\lambda \int_0^{1-s} f(t)\ud t
\end{align}
for almost all $s\in (0,1)$,
and the right-hand side is a continuous function of $s$. It follows that $f\in C[0,1]$ and that \eqref{eq:volterra-id}
holds for all $s\in [0,1]$.
Then the same argument shows that in fact $f\in C^1[0,1]$, and applying the same argument once more gives $f\in C^2[0,1]$. Differentiating \eqref{eq:volterra-id} twice gives
\begin{align*}
\lambda^2 f''(s) &= -f(s),\quad s\in [0,1].
\end{align*}
The solutions satisfying $f'(0)=0$ are of the form
$f(s)=C\cos(ks)$ with $k^2=1/\la^2$, and the boundary condition
$f(1)=0$ gives
$$
        k=\frac{\pi}{2}+\pi n,\quad n=0,1,2,\ldots .
$$
Substituting the corresponding functions
$f_n(s)=\cos((\tfrac{\pi}{2}+\pi n)s)$ into the original equation and setting $k_n:=\frac{\pi}{2}+\pi n$, we obtain
$$
(Sf_n)(s)=\frac{\sin(k_n(1-s))}{k_n}
=\frac{\sin k_n}{k_n}f_n(s)
=\frac{(-1)^n}{k_n}f_n(s).
$$
Thus the eigenvalues of $S$ are given by
$$\la_n=\frac{(-1)^n}{k_n}=\frac{2(-1)^n}{(2n+1)\pi},
$$
and therefore $\|T\|=\|S\| = \sup_{n\in\N} |\la_n|=2/\pi$.
\end{example}

The remainder of this section is devoted to studying some spectral properties of normal operators.
Our first aim is to prove that the spectral radius of a normal operator equals the operator norm. For selfadjoint operators this has already been observed as a consequence of Theorem \ref{thm:spect-sa}, and for unitary operators this is an immediate consequence of Proposition \ref{prop:spectrum-unitary}.

\begin{proposition}\label{prop:normal-spectralradius}
 An operator $T\in\calL(H)$ is normal if and only if
 $$ \n Tx \n = \n T^\star x\n, \quad x\in H.$$
 If $T$ is normal, then
 $ \n T \n^n = \n T^n\n$ for all $n\ge 1$,
 and therefore $r(T)=\n T\n.$
\end{proposition}
\begin{proof} If $T$ is normal, then
 $$ \n Tx \n^2 = \iprod{Tx}{Tx} = \iprod{x}{T^\star Tx} = \iprod{x}{TT^\star x} = \iprod{T^\star x}{T^\star x} = \n T^\star x\n^2\!.$$
 In the converse direction, the equality implies $\iprod{(T^\star T - TT^\star)x}{x} = 0$ for all $x\in H$ and therefore
 $T^\star T - TT^\star = 0$ by Proposition \ref{prop:polarisation}. This proves the first assertion.

If $T$ is normal, then for all norm-one vectors $x\in H$ we have
$$ \n T^\star Tx\n^2 = \iprod{(T^\star T)^2 x}{x}
= \iprod{(T^\star)^2 T^2 x}{x} = \n T^2 x\n^2$$
and therefore, since $\n T^\star T\n = \n T\n^2$ by \eqref{eq:HS-adjoint},
\begin{align*}
\n T\n^2 = \n T^\star T\n = \n T^2\n.
\end{align*}

Since $T^{2^k}$ is normal for every $k\ge0$, repeated application gives
$$
\|T^{2^k}\|=\|T\|^{2^k}, \quad k\ge 0.
$$
By the spectral radius formula,
$$
r(T) = \lim_{k\to\infty}\|T^{2^k}\|^{1/2^k} = \|T\|.
$$
Finally, $T^n$ is normal for every $n\ge1$. Applying the identity
$r(S)=\|S\|$, just proved for normal operators, to $S=T^n$, and using the
spectral mapping theorem, gives
\[
\|T^n\| = r(T^n) = r(T)^n =\|T\|^n.
\qedhere\]
\end{proof}

Recall that $\la\in \C$ is
 called an {\em approximate eigenvalue} of an operator $T$ on a Banach space $X$ if there exists a sequence
 $(x_n)_{n\ge 1}$ in $X$ such that $\n x_n\n = 1$ for all $n\ge 1$ and $\limn\n Tx_n - \la x_n\n = 0$.
By Proposition \ref{prop:approx-eigenvalue-bdry} the boundary spectrum of any bounded operator on $X$ consists of approximate eigenvalues. For normal operators on Hilbert spaces, more is true:

\begin{proposition}\label{prop:normal-approx-ev}
Every point in the spectrum of a normal operator $T\in\calL(H)$ is an approximate eigenvalue.
\end{proposition}

\begin{proof}
Suppose that $\la\in\C$ is not an approximate eigenvalue. Then there exists a constant $c>0$ such that $\n(\la-T)x\n\geq c\n x\n$ for all $x\in H$. By Proposition \ref{prop:closed-range}, the operator $S:= \la-T$ is injective and has closed range.

Since $T$ is normal, also $S$ is normal, and hence $\n Sx\n=\n S^\star x\n$ for all $x\in H$. It follows that $S^\star$ is injective. Therefore $\ran(S)$ is dense in $H$. Since $\ran(S)$ is also closed, we have $\ran(S)=H$. Thus $S=\la-T$ is bijective. Moreover, its inverse is bounded by the estimate $\n S^{-1}y\n\leq c^{-1}\n y\n$, $y\in H$. Hence $\la\in\rh(T)$.

We have proved that every $\la\in\C$ which is not an approximate eigenvalue belongs to $\rh(T)$. Therefore every point of $\si(T)$ is an approximate eigenvalue.
\end{proof}

Theorem \ref{thm:comp-spec} implies that every nonzero element of the spectrum of a compact operator
is both an isolated point and an eigenvalue.
The next theorem states that for normal operators on a Hilbert space, all isolated
points in the spectrum are eigenvalues; no compactness assumption is needed. Normality cannot be omitted:
the Volterra operator has spectrum $\{0\}$, but $0$ is not an eigenvalue (see Problem \ref{prob:Volterra}).

If $T$ is a bounded operator on $H$, for $\la\in \C$ we set
$$ E_\la := \{x\in H:\, Tx = \la x\}.$$
Thus $\la$ is an eigenvalue of $T$ if and only if $E_\la$ is nonzero.
We recall that spectral projections have been defined in Theorem \ref{thm:spec-proj}.

\begin{theorem}[Isolated points are eigenvalues]\label{thm:sa-isolated-pt}
 Let $T\in\calL(H)$ be a normal operator and let $\la$ be an isolated point in $\si(T)$.
 Then $\la$ is an eigenvalue for $T$ and the spectral projection $P^{\{\la\}}$ corresponding to $\{\la\}$
equals the orthogonal projection $P_\la$ onto $E_\la$.
\end{theorem}

The proof uses the following simple observation.

\begin{lemma}\label{lem:normal-reduce} Let $T\in \calL(H)$. If $Y$ is
a closed subspace of $H$, then:
 \begin{enumerate}[label={\rm(\arabic*)}, leftmargin=*]
  \item\label{it:normal-reduce1} if $Y$ is invariant under $T$, then $Y^\perp$ is  invariant under $T^\star$;
  \item\label{it:normal-reduce2}  if $Y$ is invariant under $T$ and $T^\star$, then $Y^\perp$ is invariant under $T$ and $T^\star$ and
  $$(T|_{Y})^\star = T^\star|_{Y}  \ \ \hbox{and} \ \ (T|_{Y^\perp})^\star = T^\star|_{Y^\perp}$$ as operators on $Y$ and $Y^\perp$, respectively.
 \end{enumerate}
In particular, if $T$ is selfadjoint (respectively, normal) and $Y$ is invariant under $T$ (respectively, under $T$ and $T^\star$), then $T|_Y$ and $T|_{Y^\perp}$ are selfadjoint (respectively, normal).
\end{lemma}
 \begin{proof}  If $Y$ is invariant under $T$,
then for all $y\in Y$ and $y^\perp\in Y^\perp$ we have $\iprod{y}{T^\star y^\perp} = \iprod{Ty}{y^\perp} = 0$. This proves \ref{it:normal-reduce1}. The first assertion of \ref{it:normal-reduce2} follows as well, and
if $Y$ is invariant under $T$ and $T^\star$, then for all $y, y'\in Y$ we have
$$ \iprod{y}{(T|_Y)^\star y'}
=  \iprod{T|_Y y}{y'} = \iprod{T y}{y'}
= \iprod{y}{T^\star y'} = \iprod{y}{(T^\star)|_Y y'}.$$
This proves the first identity of \ref{it:normal-reduce2}. The second is proved in the same way. The final assertion is an immediate consequence of \ref{it:normal-reduce2}.
\end{proof}

\begin{proof}[Proof of Theorem \ref{thm:sa-isolated-pt}]
Replacing $T$ by $T-\la$, we may assume that $\la=0$. Let $P$ denote the spectral projection corresponding to $\{0\}$.

If $\si(T)=\{0\}$, then Proposition \ref{prop:normal-spectralradius} gives $\n T\n=r(T)=0$, and hence $T=0$. In this case $E_0=H$ and $P=I$, so there is nothing left to prove.
We may therefore assume that $\si(T)\setminus\{0\}$ is nonempty.

Put $E:=\ran(P)$ and $F:=\ker(P)=\ran(I-P)$.
By Theorem \ref{thm:spec-proj}, $H=E\oplus F$, the spaces
$E$ and $F$ are invariant under $T$, and the restrictions
$T_E:=T|_E$ and $T_F:=T|_F$ satisfy
$\si(T_E)=\{0\}$ and $\si(T_F)=\si(T)\setminus\{0\}$.
In particular, $0\in\rh(T_F)$ and $E\not=\{0\}$.

Since $T$ is normal, $T^\star$ commutes with $T$, hence with
every resolvent of $T$, and therefore with $P$. Thus $E$ is
also invariant under $T^\star$. By Lemma \ref{lem:normal-reduce},
$T|_{E}$ is normal, and Proposition \ref{prop:normal-spectralradius}
gives $\n T|_{E}\n=r(T|_{E})=0$.
Consequently, $E\subseteq E_0$.

Conversely, if $y\in E_0$, then
$R(\mu,T)y=\mu^{-1}y$ for all $\mu\in\rh(T)$.
Let $\Gamma\subseteq\rh(T)$ be a positively oriented circle
enclosing $0$ and no other point of $\si(T)$. Then
$$
Py=\frac1{2\pi i}\int_\Gamma R(\mu,T)y\ud\mu
=\frac1{2\pi i}\int_\Gamma\frac{y}{\mu}\ud\mu
=y.
$$
Thus $E=E_0$, and since $E\not=\{0\}$, this proves that $0$
is an eigenvalue.

It remains to prove that $P$ is the orthogonal projection onto $E_0$. Let $y\in E_0$ and $z\in F$. Since $T$ is normal and $Ty=0$, Proposition \ref{prop:normal-spectralradius} gives $T^\star y=0$. As $0\in\rh(T_F)$, there exists a $w\in F$ such that $T_Fw=z$. Therefore
$$
\iprod{z}{y}=\iprod{Tw}{y}=\iprod{w}{T^\star y}=0.
$$
Thus $F\perp E_0$. Since $H=E_0\oplus F$, it follows that $F=E_0^\perp$, and therefore $P$ is the orthogonal projection onto $E_0$.
\end{proof}

\begin{corollary}\label{cor:copact-normal-multipicities}
The geometric and algebraic multiplicities of every nonzero element in the spectrum of a compact normal operator coincide.
\end{corollary}

This justifies using the term {\em multiplicity}\index{multiplicity}
for their common value.

We have the following commutation theorem for normal operators.

\begin{theorem}[Fuglede--Putnam--Rosenblum]\label{thm:FPR}\index{theorem!Fuglede--Putnam--Rosenblum}
If $T\in \calL(H)$ is normal and $S\in \calL(H)$ is bounded and satisfies
$ ST = TS,$
then
$ ST^\star = T^\star S.$
\end{theorem}
\begin{proof}
We prove the more general assertion that if $T_1,T_2$ are normal
and $ST_1=T_2S$, then $ST_1^\star=T_2^\star S$.
By induction, $ST_1^n=T_2^nS$ for all $n\in\N$, and expanding
the exponentials gives
$Se^{wT_1}=e^{wT_2}S$ for all $w\in\C$.

Consider the entire function
$F(z):=e^{zT_2^\star}Se^{-zT_1^\star}$.
Using the preceding identity with $w=\ov z$ and the normality
of $T_1$ and $T_2$, we obtain
\begin{align*}
F(z)
&=e^{zT_2^\star-\ov zT_2}e^{\ov zT_2}Se^{-zT_1^\star}
=e^{zT_2^\star-\ov zT_2}Se^{\ov zT_1-zT_1^\star}.
\end{align*}
Both exponents on the last line are $i$ times selfadjoint
operators, so Proposition \ref{prop:sa-unitary} shows that
their exponentials are unitary. Consequently,
$\n F(z)\n=\n S\n$ for all $z\in\C$.
Liouville's theorem, applied to the scalar functions
$z\mapsto\iprod{F(z)x}{y}$ for $x,y\in H$, implies that $F$
is constant. Differentiating at zero gives
\[
0=F'(0)=T_2^\star S-ST_1^\star.
\qedhere\]
\end{proof}

We finish  with an observation about relative spectra. Let $\calA\subseteq \calL(H)$
be a unital closed {\em $\star$-subalgebra},\index{subalgebra!$\star$-} that is, $\calA$ is a norm-closed unital subalgebra of $\calL(H)$ that is also closed under taking adjoints. For such subalgebras we have the following improvement to Proposition \ref{prop:siAT-bdry}:

\begin{proposition}\label{prop:siAT}
Let $\calA\subseteq \calL(H)$ be a unital closed $\star$-subalgebra and let $T\in \calA$. Then
$$\si_\calA(T) = \si(T).$$
\end{proposition}
\begin{proof}
The case $H=\{0\}$ is trivial, so assume that $H\not=\{0\}$.
We first prove that if $T\in\calA$ is invertible in $\calL(H)$,
then $T^{-1}\in\calA$.

Put $B:=T^\star T\in\calA$ and set
$a:=\n T^{-1}\n^{-2}>0$ and $b:=\n T\n^2>0$.
Since $a\n x\n^2\le\n Tx\n^2$, we have
$$
a\n x\n^2 \le \n Tx\n^2 =  \iprod{Bx}{x}\le b\n x\n^2,
\quad x\in H.
$$
Thus the selfadjoint operator $I-b^{-1}B$ satisfies
$$
0\le\iprod{(I-b^{-1}B)x}{x}\le1-a/b,\quad \n x\n=1.
$$
Hence, by Theorem \ref{thm:spect-sa}, $\n I-b^{-1}B\n\le1-a/b<1$, and
the Neumann series gives
$$
B^{-1}=b^{-1}\sum_{n=0}^\infty(I-b^{-1}B)^n.
$$
Since $I\in\calA$, this series converges
in $\calA$. It follows that $T^{-1}=B^{-1}T^\star\in\calA$.

Applying this result to $\la-T$, we find that $\rh(T)\subseteq\rh_\calA(T)$. The reverse inclusion is immediate, and hence $\si_\calA(T)=\si(T)$.
\end{proof}

\section{The Continuous Functional Calculus}\label{sec:contFC}

In Chapter \ref{ch:spectral} we have seen how to associate a bounded operator $f(T)$ with a bounded operator $T$
when $f$ is holomorphic in an open neighbourhood of $\sigma(T)$. Here we will prove that
for normal operators $T$ acting on a Hilbert space, the
functional calculus $f\mapsto f(T)$ can be extended to continuous functions on $\sigma(T)$.

\subsection{The Selfadjoint Continuous Functional Calculus}\label{sec:contFC-sa}

We begin with the case of selfadjoint operators.

\begin{theorem}[Continuous\, functional\, calculus\, for\, selfadjoint\, operators]\index{functional calculus!continuous}
\label{thm:sa-cont-fc}\index{theorem!continuous functional calculus}
Let\  $T\in\calL(H)$ be a selfadjoint operator. Then there exists a unique continuous linear mapping $ f\mapsto f(T)$ from $C(\si(T))$ to $\calL(H)$
 with the following properties:
 \begin{enumerate}[label={\rm(\roman*)}, leftmargin=*]
  \item\label{it:sa-cont-fc1} if $f(z) = z^n$ with $n\in\N$, then $f(T) = T^{n}$;
  \item\label{it:sa-cont-fc3} for all $f,g\in C(\si(T))$ we have $(fg)(T) = f(T)g(T)$;
  \item\label{it:sa-cont-fc2} for all $f\in C(\si(T))$ we have $\ov f(T) = (f(T))^\star$;
  \item\label{it:sa-cont-fc4} for all $f\in C(\si(T))$ we have $\n f(T)\n = \n f\n_{\infty}$.
 \end{enumerate}
The operators $f(T)$ are normal, and $f(T)$ is selfadjoint if and only if $f$ is real-valued.
\end{theorem}
\begin{proof}
 For polynomials $p(z) = \sum_{n=0}^N c_n z^n$ we define $p(T):= \sum_{n=0}^N c_n T^n $. These operators are normal and satisfy \ref{it:sa-cont-fc1}, \ref{it:sa-cont-fc3}, and \ref{it:sa-cont-fc2}.
 Moreover, by the spectral mapping theorem for the holomorphic calculus,
 \begin{align*} \n p(T)\n  = \sup\{|\la|:\, \la\in \si(p(T))\}  = \sup\{|\la|:\, \la \in p(\si(T))\}
 =  \n p\n_{\infty}
 \end{align*}
 and therefore \ref{it:sa-cont-fc4} holds.
By the Stone--Weierstrass approximation theorem, the polynomials are dense in
 $C(\sigma(T))$. Therefore, by \ref{it:sa-cont-fc4} and an approximation argument,
 the mapping $p\mapsto p(T)$ has a unique extension to an isometry
 from $C(\si(T))$ into $\calL(H)$, and \ref{it:sa-cont-fc3}--\ref{it:sa-cont-fc4} again hold.

Since normality is preserved under operator-norm limits, the operators $f(T)$ are normal.
Property \ref{it:sa-cont-fc2} implies that if $f$ is real-valued, then $f(T)$ is selfadjoint. Conversely, if $f(T)$ is selfadjoint, then
$f(T) = \ov f(T)$ by property \ref{it:sa-cont-fc2} and therefore $\n f-\ov f\n_{C(\si(T))} = 0$ by property \ref{it:sa-cont-fc4},
so $f = \ov f$ is real-valued.
\end{proof}

\subsection{The Normal Continuous Functional Calculus}\label{subsec:contFCnorma}

Every polynomial in the real variables $x$ and $y$ can be written as a polynomial in the variables $z$ and $\ov z$ by substituting $z = x+iy$, $\ov z = x-iy$. For example, $x^2+y^2 =z\ov z$. For polynomials $p(z,\ov z) = \sum_{i,j=0}^k c_{ij}z^i\ov z^j$ and normal operators $T\in\calL(H)$, we define $$p(T,T^\star):=  \sum_{i,j=0}^k c_{ij}T^i T^{\star j}.$$
The crucial result that enables us to extend the continuous functional calculus to normal operators is the following spectral mapping theorem.

\begin{proposition}\label{prop:SMT-normal}
 If $T\in\calL(H)$ is normal and $p$ is a polynomial in $z$ and $\ov z$, then
$$\si(p(T,T^\star)) = \{p(\la,\ov\la):\, \la\in \si(T)\}.$$
\end{proposition}

\begin{proof}
By Proposition \ref{prop:normal-approx-ev}, every $\la\in\si(T)$ is an approximate eigenvalue of $T$, that is, there exists a sequence $(x_n)_{n\ge 1}$ of norm-one vectors such that $\limn Tx_n - \la x_n = 0$.
Then $\limn T^\star x_n - \ov\la x_n = 0$ by Proposition \ref{prop:normal-spectralradius}. This implies
$\limn p(T,T^\star)x_n - p(\la,\ov\la)x_n = 0$, so $p(\la,\ov\la)$ is an approximate eigenvalue for $p(T,T^\star)$. In particular, $p(\la,\ov\la) \in \si(p(T,T^\star))$.
This proves the inclusion  `$\supseteq$'.

For the inclusion `$\subseteq$', fix $\mu\in\si(p(T,T^\star))$.
We wish to prove that there exists $\la\in\si(T)$ such that $p(\la,\ov\la) = \mu$.
\smallskip

{\em Step 1} -- Fixing $\eps>0$, we claim that there is a nonzero closed subspace $Y$ of $H$, invariant under both $T$ and $T^\star$, such that
\begin{align}\label{eq:pTTsY}
\big\n(p(T,T^\star)-\mu I)|_Y\big\n\le\eps^{1/2}.
\end{align}

To prove the claim, let $S:= p(T,T^\star)-\mu I$. This operator is normal and we have $0\in\si(S)$. Let $R:= S^\star S$.
Arguing as above, we find that $0\in\si(R)$.
Consider the continuous function $f:[0,\infty)\to [0,1]$ given by
$$ f(t):=
\begin{cases}
 1, & 0\le t \le \eps/2; \\
 2(1-t/\eps), & \eps/2 \le t\le \eps; \\
 0, & t\ge \eps,
\end{cases}
$$
and let $f(R)$ be the selfadjoint operator obtained from the continuous functional calculus for selfadjoint operators (Theorem \ref{thm:sa-cont-fc}).
We will show that $$Y:= \{x\in H:\, f(R)x = x\}$$
has the desired properties.

Since $T$ commutes with $R$, it commutes with $f(R)$, and therefore $Y$ is invariant under $T$.
By the same reasoning, $Y$ is invariant under $T^\star$. Moreover, by the properties of continuous functional calculus,
for all $y\in Y $ we have
$$ \n R y \n= \n R f(R)y\n \le \big\n t\mapsto t f(t)\big\n_{C(\si(R))}\n y\n
\le \big\n t\mapsto t f(t)\big\n_{C[0,\eps]} \n y\n \le \eps\n y\n.$$
This implies $$ \n Sy\n^2 = \iprod{Ry}{y} \le \n Ry\n\n y\n \le \eps \n y\n^2\!.$$
This gives \eqref{eq:pTTsY}.
The claim will be proved once we have checked that $Y$ is nonzero. If $f(2t)\not=0$ for some $t\ge 0$, then $2t\le \eps$ and therefore $f(t)=1$. By multiplicativity,
$$ \n (I - f(R))f(2R)\n = \big\n t\mapsto (1-f(t))f(2t)\big\n_{C(\si(R))} = 0,$$
where $f(2R):=g(R)$ with $g(t):= f(2t)$.
It follows that $\Ran(f(2R))\subseteq Y$. But $\ran(f(2R))$ is nonzero since
$$\n f(2R)\n = \big\n  t\mapsto |f(2t)|\big\n_{C(\si(R))} \ge |f(0)| = 1.$$

{\em Step 2} -- Fix $\eps>0$ and let $Y$ be the subspace
constructed in Step 1. Put $B:=T|_Y$. By
Lemma \ref{lem:normal-reduce}, $B$ is normal and
$B^\star=T^\star|_Y$. Since $Y\not=\{0\}$, Theorem
\ref{thm:siTnonempty} shows that $\si(B)$ is nonempty.

Choose $\la\in\si(B)$. By Proposition
\ref{prop:normal-approx-ev}, $\la$ is an approximate eigenvalue
of $B$ and hence of $T$. In particular, $\la\in\si(T)$.
The inclusion `$\supseteq$' established in the first part of
the proof can be applied to the normal operator $B$, which gives
$p(\la,\ov\la)\in\si(p(B,B^\star))$.
Since the modulus of a spectral value is bounded by the
operator norm, it follows that
\begin{align*}
|p(\la,\ov\la)-\mu|\le\n p(B,B^\star)-\mu I\n=\big\n(p(T,T^\star)-\mu I)|_Y\big\n
\le\eps^{1/2}.
\end{align*}
Thus $\mu$ has distance at most $\eps^{1/2}$ from the compact set
$\{p(\la,\ov\la):\,\la\in\si(T)\}$.
Since $\eps>0$ was arbitrary, $\mu$ belongs to this set.
This proves the inclusion `$\subseteq$'.
\end{proof}

With this result at hand we can extend the continuous functional calculus to normal operators.
Repeating the proof of Theorem \ref{thm:sa-cont-fc}, we obtain the following result.

\begin{theorem}[Continuous functional calculus for normal operators]\index{functional calculus!continuous}
\label{thm:normal-cont-fc}\index{theorem!continuous functional calculus}
 Let $T\in\calL(H)$ be a normal operator. Then there exists a unique continuous linear mapping $ f\mapsto f(T)$ from $C(\si(T))$ to $\calL(H)$
 with the following properties:
 \begin{enumerate}[label={\rm(\roman*)}, leftmargin=*]
  \item\label{it:normal-cont-fc1} if $f(z) = z^m\ov z^n$ with $m,n\in\N$, then $f(T) = T^m T^{\star n}$;
  \item\label{it:normal-cont-fc3} for all $f,g\in C(\si(T))$ we have $(fg)(T) = f(T)g(T)$;
  \item\label{it:normal-cont-fc2} for all $f\in C(\si(T))$ we have $\ov f(T) = (f(T))^\star$;
  \item\label{it:normal-cont-fc4} for all $f\in C(\si(T))$ we have $\n f(T)\n = \n f\n_{\infty}$.
 \end{enumerate}
The operators $f(T)$ are normal, and $f(T)$ is selfadjoint if and only if $f$ is real-valued.
\end{theorem}

The next theorem extends Proposition \ref{prop:SMT-normal} to continuous functions defined on $\si(T)$.

\begin{theorem}[Spectral mapping theorem]\label{thm:SMT-normal}\index{theorem!spectral mapping, bounded normal operators}
If $T\in\calL(H)$ is normal, then for all $f\in C(\si(T))$ we have
$$ \si(f(T)) = f(\si(T)).$$
\end{theorem}
\begin{proof}
The proof of the inclusion $\si(f(T))\subseteq f(\si(T))$ is similar to that in Theorem \ref{thm:SMT}.
Indeed, if $\la\not\in f(\si(T))$, the function $g_\la = 1/f_\la$ with $f_\la(z) := \la- f(z)$ is continuous on $\si(T)$, and by multiplicativity we obtain
  $$ g_\la(T)(\la-f(T)) = (\la-f(T))g_\la(T) = (f_\la g_\la)(T) = \one(T) =  I,$$
so $\la\in \rh(f(T))$ and $R(\la,f(T))  = g_\la(T)$. This gives the stated inclusion.

For the converse inclusion, put
$$
\calA:=\{g(T):\,g\in C(\si(T))\}.
$$
Note that $\calA$ is closed in operator norm since $C(\si(T))$ is complete and the calculus is isometric.

Let $\mu\in\rh(f(T))$. By Proposition \ref{prop:siAT},
$R(\mu,f(T))$ belongs to $\calA$, so there exists
$g\in C(\si(T))$ such that $R(\mu,f(T))=g(T)$.
Multiplicativity gives
$$
((\mu-f)g)(T)
=(\mu I-f(T))g(T)
=I=\one(T).
$$
Since the calculus is injective, $(\mu-f)g=\one$ on $\si(T)$.
In particular, $\mu-f$ has no zeroes on $\si(T)$, and therefore
$\mu\notin f(\si(T))$. This proves
$f(\si(T))\subseteq\si(f(T))$.
\end{proof}

\begin{corollary}\label{cor:normal-cont-fc-pos}
 Let $T\in\calL(H)$ be a normal operator and let $f\in C(\si(T))$. Then $f(T)\ge 0$ if and only if $f$ is nonnegative.
\end{corollary}
\begin{proof}
 If $f$ is nonnegative, the spectral mapping theorem gives $\si(f(T)) = f(\si(T)) \subseteq [0,\infty)$, and therefore the selfadjoint operator  $f(T)$ is positive by Theorem \ref{thm:spect-sa}.
 Conversely, if $f(T)$ is positive,  then $\si(f(T))\subseteq [0,\infty)$ by Theorem \ref{thm:spect-sa},
 and therefore $f(\si(T)) \subseteq [0,\infty)$ by the spectral mapping theorem.
\end{proof}

The next theorem extends Theorem \ref{thm:Dunford-composition} to continuous functions defined on $\si(T)$.

\begin{theorem}[Composition]\label{thm:comp-c} Let $T\in\calL(H)$ be normal.
For all $f\in C(\si(T))$ and $g\in C(f(\si(T))) = C(\si(f(T)))$ we have $g\circ f\in C(\si(T))$ and
 $$ g(f(T)) = (g\circ f)(T).$$
\end{theorem}
\begin{proof}
First let $p(z) = z^m\ov z^n$ with $m,n\in \N$. Then, by the properties of the continuous calculus,
$$p(f(T)) = (f(T))^m (\ov f(T))^n = (f^m\ov f^n)(T) = (p\circ f)(T).$$
By linearity, this identity extends
to polynomials $p$ in $z$ and $\bar z$.
If $g\in C(\si(f(T)))$ is any continuous function,
the identity follows by approximating $g$ uniformly by polynomials $p_n$ via the Stone--Weierstrass theorem to obtain
\begin{align*}
  \n p_n(f(T)) - g(f(T))\n & = \n p_n-g\n_{C(\sigma(f(T)))} \to 0 \ \ \hbox{ as $n\to\infty$}
\intertext{and}
  \n (p_n\circ f)(T) - (g\circ f)(T)\n & = \n p_n\circ f -g\circ f\n_{C(\sigma(T))} \to 0 \ \ \hbox{ as $n\to\infty$}.
  \qedhere
\end{align*}
\end{proof}

We finally check consistency with the holomorphic calculus.

\begin{theorem}\label{thm:consistent-Dunford}
Let $T\in \calL(H)$ be normal.
If $f\in H(\Om)$, where $\Om$ is an open set containing $\si(T)$, then the operator $f(T)$ defined through the continuous calculus agrees with the operator defined through the holomorphic calculus.
\end{theorem}

 \begin{proof}
Let $\Phi: C(\si(T))\to \calL(H)$ be the continuous functional calculus of $T$, so by assumption $f(T) = \Phi(f|_{\si(T)})$. Let $\Gamma$ be an admissible contour for $\si(T)$ in $\Om$.
For $z\in\Gamma$ define
$$
\rho_z(\la):=\frac{1}{z-\la},\quad \la\in\si(T).
$$
Then $\rho_z\in C(\si(T))$, and by multiplicativity of the continuous
functional calculus,
$$
(zI-T)\rho_z(T)=((z-{\rm id})\rho_z)(T)=\one(T)=I.
$$
Similarly, $\rho_z(T)(zI-T)=I$, and therefore
$\rho_z(T)=R(z,T)$.

Next define
$$
g(\la):=\frac{1}{2\pi i}\int_\Gamma \frac{f(z)}{z-\la}\ud z,
\quad \la\in\si(T).
$$
By the Cauchy integral formula, $g(\la)=f(\la)$ for all $\la\in\si(T)$, so
$g=f|_{\si(T)}$.
Since $(z,\la)\mapsto f(z)/(z-\la)$ is continuous, hence uniformly continuous, on
$\Gamma\times\si(T)$, the map $z\mapsto f(z)\rho_z$ is continuous from
$\Gamma$ into $C(\si(T))$. By linearity and continuity of $\Phi$,
$$
f(T)=g(T)
=\frac{1}{2\pi i}\int_\Gamma f(z)\rho_z(T)\ud z
=\frac{1}{2\pi i}\int_\Gamma f(z)R(z,T)\ud z,
$$
and the right-hand side is the operator defined by the holomorphic calculus.
\end{proof}

\subsection{Applications}

We now turn to some applications of the continuous functional calculus.

 \begin{proposition}[Square roots]\label{prop:pos-sq-root}
 If $T\in \calL(H)$ is positive, there exists a unique positive operator $S\in \calL(H)$ such that $S^2 = T$.
 \end{proposition}

 Henceforth, this operator $S$ will be denoted by $T^{1/2}$\!.

\begin{proof} Since $T$ is positive, $T$ is selfadjoint and its spectrum is contained in $[0,\infty)$
by Theorem \ref{thm:spect-sa}. Hence, $f(t) = \sqrt t$ is a well-defined continuous function on $\sigma(T)$. The operator $S:= f(T)$ is positive by Corollary \ref{cor:normal-cont-fc-pos}, and it satisfies $S^2 = T$ by the properties of the continuous functional calculus. It remains to prove uniqueness.
Suppose $\wt S$ is another positive operator with the property that $\wt S^2 = T$. With $f(t) = \sqrt{t}$, $g(t)=t^2$\!, and $h(t)=t$ we have, by the continuous functional calculus and Theorem \ref{thm:comp-c},
$$f(\wt S^2) = f(g(\wt S)) = (f\circ g)(\wt S) = h(\wt S) = \wt S.$$
It follows that $S = f(S^2) = f(T) = f(\wt S^2) = \wt S.$
This completes the uniqueness proof.
\end{proof}

\begin{definition}[Modulus of an operator]  The {\em modulus}\index{modulus!of an operator}
\index{$T$@$\abs{T}$} of an operator $T\in \calL(H)$ is the positive operator $ |T|:= (T^\star T)^{1/2}\!.$
\end{definition}

\begin{corollary}\label{cor:sqrt-normal} If $T\in\calL(H)$ is a normal operator, then
$|T|  = f(T)$, where $f(z) := |z|$.
\end{corollary}
\begin{proof}
Let $g(z):= \ov z z$. Then $f^2 = g$ and therefore
$ |T|^2 = T^\star T = g(T) = f^2(T) = (f(T))^2$ by the multiplicativity of the continuous functional calculus.
Since by Corollary \ref{cor:normal-cont-fc-pos} the operator $f(T)$ is positive, the result follows by
taking square roots.
\end{proof}

We continue with a polar decomposition result. In view of future applications we phrase it for bounded operators $T\in \calL(H,K)$, where $H$ and $K$ are Hilbert spaces. The {\em modulus} of such an operator is the positive operator
$ |T|:= (T^\star T)^{1/2}$ on $H$. An operator $U\in \calL(H,K)$ will be called {\em unitary}
if it is isometric and surjective.
A {\em partial isometry}\index{partial!isometry} is a bounded operator $V\in \calL(H,K)$
for which there exists an orthogonal direct sum decomposition $H = H_0\oplus H_0^\perp$
such that $V$ is isometric from $H_0$ into $K$ and zero on $H_0^\perp$.

\begin{theorem}[Polar decomposition]\label{thm:polar} Let
$T\in \calL(H,K)$.  The following assertions hold:
\begin{enumerate}[label={\rm(\arabic*)}, leftmargin=*]
 \item\label{it:polar1} $T$ admits a representation $T = U|T|$, where $U$ is the unique partial isometry from $H$ to $K$ satisfying $\ker(U)=\ker(T)$
and mapping $\overline{\ran(|T|)}$ isometrically onto
$\overline{\ran(T)}$;

 \item\label{it:polar2} if $T$ is invertible, then $T$ admits a unique representation $T = U|T|$ with $U$ unitary from $H$ onto $K$.
\end{enumerate}
\end{theorem}
\begin{proof}
\ref{it:polar1}: \
From
$$
\n Tx\n^2 = \iprod{T^\star Tx}{x} = \iprod{|T|x}{|T|x} = \n |T|x\n^2
$$
it follows that the mapping $U_0: |T|x \mapsto Tx$, viewed as a linear operator from $\ran(|T|)$ onto ${\ran(T)}$, is well defined and isometric, and by density it extends to an isometry from $\ov{\ran(|T|)}$ onto $\ov{\ran(T)}$.
Moreover, $T = U_0 |T|$.

Since $\n Tx\n=\n |T|x\n$, we have $\ker(T)=\ker(|T|)=\overline{\ran(|T|)}^\perp$.
With respect to the orthogonal decomposition $H = \ov{\ran(|T|)}\oplus (\overline{\ran(|T|)})^\perp$ we extend $U_0$ by zero on $(\overline{\ran(|T|)})^\perp$
to obtain a partial isometry $U$ with the desired properties. Uniqueness follows since $\ran(|T|)$ is dense in
$\overline{\ran(|T|)}$: the identity $T=U|T|$ determines $U$ on
$\ran(|T|)$, and the condition $\ker(U)=\ker(T)=\overline{\ran(|T|)}^\perp$
determines $U$ on the orthogonal complement.

\smallskip
\ref{it:polar2}:\
If $T$ is invertible, then $\ker(T)=\{0\}$ and $\ran(T)=K$.
Since $\ov{\ran(|T|)}=\ker(T)^\perp=H$ by part \ref{it:polar1},
the partial isometry $U$ constructed there is isometric
on all of $H$ and maps $H$ onto $K$, so it is unitary.
If $T=\wt U|T|$ is another such representation, then $U$
and $\wt U$ agree on $\ran(|T|)$. Since this subspace is dense
in $H$, we obtain $U=\wt U$.
\end{proof}

\section{The Sz.-Nagy Dilation Theorem}\label{sec:Sz-Nagy}

The last section of this chapter is devoted to a proof of the celebrated Sz.-Nagy dilation theorem, which asserts that every Hilbert space contraction has a unitary dilation.
Since this entails no additional difficulty, we take a rather general
approach, starting from an arbitrary group $G$ with identity element $e$; the Sz.-Nagy dilation theorem is obtained by considering $G = \Z$.

\begin{definition}[Positive definiteness]\label{def:posdef}
 A mapping $T: G \to \calL(H)$ is called
 {\em positive definite}\index{positive definite} if for all finite choices of
 $g_1,\dots,g_N\in G$ and $h_1,\dots,h_N\in H$ we have
$$\sum_{m,n=1}^N \iprod{T(g_m^{-1}g_n)h_n}{h_m} \ge 0.$$
\end{definition}
\begin{definition}[Representations, unitary representations]
 A mapping $U: G \to \calL(H)$ is called a {\em representation of $G$ on}\index{representation} $H$
 if $U(e)=I$ and $U(g_1)U(g_2) = U(g_1g_2)$ for all $g_1,g_2\in G$.
 A {\em unitary representation} is a representation for which every
operator $U(g)$ is unitary.
\end{definition}

The following result connects these two notions.

\begin{proposition}
Let $U:G\to\calL(\wt H)$ be a unitary representation of $G$ on a Hilbert space $\wt H$. Let $J:H\to \wt H$ be an isometric embedding of another Hilbert space $H$ into $\wt H$. Then the function $T:G\to \calL(H)$, given by
$$ T(g) :=  J^\star U(g)J,$$
is positive definite and satisfies $T(e)=I$ and $T^\star(g) = T(g^{-1})$ for all $g\in G$.
\end{proposition}
\begin{proof}
The identity $T(e)=I$ follows from $U(e)=I$ and $J^\star J=I$, and the
identity $U(g^{-1}) = (U(g))^\star$ implies $$T^\star(g) = J^\star U^\star(g)J = J^\star (U(g))^{-1}J= J^\star U(g^{-1})J  = T(g^{-1}).$$
To prove positive definiteness, let $g_1,\dots,g_N\in G$ and $h_1,\dots,h_N\in H$. Then
\begin{align*}
 \sum_{m,n=1}^N \iprod{T(g_m^{-1}g_n)h_n}{h_m} &  = \sum_{m,n=1}^N \iprod{U(g_m^{-1})U(g_n)Jh_n}{Jh_m}
  = \Bigl\n \sum_{n=1}^N U(g_n)Jh_n\Big\n^2 \ge 0.
  \qedhere
\end{align*}
\end{proof}

The next theorem establishes that, conversely, every positive definite function $T: G\to \calL(H)$
satisfying $T(e)=I$ and $T^\star(g) = T(g^{-1})$ for all $g\in G$ arises in this way.

\begin{theorem}[Unitary dilations, Sz.-Nagy]\label{thm:dilation-G}\index{theorem!Sz.-Nagy}
Let $T: G\to \calL(H)$ be a positive definite function satisfying $T(e) = I$ and $T^\star(g) = T(g^{-1})$ for all $g\in G$. There exist a Hilbert space $\wt H$, an isometric embedding $J:H\to \wt H$, and a unitary representation $U:G\to \calL(\wt H)$ such that
$$ T(g)h = J^\star U(g)J h, \quad h\in H.$$
Moreover, one may choose this dilation to be {\em minimal}\index{minimal dilation}
in the sense that
$$\wt H=\overline{\mathrm{span}}\{U(g)Jh:\,g\in G,\,h\in H\},
$$
and any two such minimal dilations are unitarily equivalent.
\end{theorem}

\begin{proof}
Let $V$ be the vector space of all functions $f:G\to H$ that
vanish outside a finite set. For $g\in G$ and $h\in H$, let
$\delta_g h$ denote the function that equals $h$ at $g$ and
zero elsewhere. Every element of $V$ is a finite sum of such
functions. We define a sesquilinear mapping from $V\times V$
to $\C$ by prescribing
\begin{equation}\label{eq:inner-posdef-1}
\iprod{\delta_g h}{\delta_{g'}h'}
:=\iprod{T((g')^{-1}g)h}{h'}
\end{equation}
and extending by sesquilinearity.

We claim that $\iprod{f_1}{f_2}=\ov{\iprod{f_2}{f_1}}$
for all $f_1,f_2\in V$ and $\iprod{f}{f}\ge0$ for all $f\in V$.
For elementary functions we have
\begin{align}\label{eq:inner-posdef-2}
\ov{\iprod{\delta_{g'}h'}{\delta_g h}}
=\iprod{h}{T(g^{-1}g')h'}.
\end{align}
Since $T^\star(g)=T(g^{-1})$ for all $g\in G$, the right-hand
sides of \eqref{eq:inner-posdef-1} and
\eqref{eq:inner-posdef-2} are equal. By sesquilinearity, this
proves the identity
$\iprod{f_1}{f_2}=\ov{\iprod{f_2}{f_1}}$ for all $f_1,f_2\in V$.
Positive definiteness implies $\iprod{f}{f}\ge0$.

The properties established in the claim imply the
Cauchy--Schwarz inequality for this positive semidefinite
sesquilinear form. It may happen, however, that
$\iprod{f}{f}=0$ for certain nonzero functions $f$ in $V$,
so this sesquilinear form may fail to be an inner product.
For this reason we consider the vector space quotient $V/N$,
where
$$
N=\{f\in V:\,\iprod{f}{f}=0\}.
$$
Let us prove that $N$ is indeed a subspace of $V$.
It is clear that $cf\in N$ for all $c\in\C$ and $f\in N$.
Furthermore, if $f,f'\in N$, then $\iprod{f}{f'}=0$ by the
Cauchy--Schwarz inequality, and from this it follows that
$f+f'\in N$.

On the quotient space $V/N$, the sesquilinear mapping
$\iprod{\cdot}{\cdot}$ induces the inner product
$$
\iprod{f+N}{g+N}:=\iprod{f}{g}.
$$
Define $\wt H$ to be the Hilbert space completion of $V/N$
with respect to this inner product. To realise $H$ as a
closed subspace of $\wt H$, we identify elements $h\in H$
with the classes modulo $N$ of the functions $f_h:G\to H$
given by $f_h=\delta_e h$. Then
$$
\iprod{f_{h_1}}{f_{h_2}}
=\iprod{\delta_e h_1}{\delta_e h_2}
=\iprod{T(e)h_1}{h_2}
=\iprod{h_1}{h_2}
$$
since $T(e)=I$. This implies that the mapping
$J:h\mapsto f_h+N$ is isometric from $H$ into $\wt H$.

The linear mapping $U_0(g):V\to V$ given by
$$
(U_0(g)f)(g'):=f(g^{-1}g'),
\quad f\in V,\ g,g'\in G,
$$
is well defined. Since $U_0(g)\delta_a h=\delta_{ga}h$,
for all $a,b\in G$ and $h,h'\in H$ we have
\begin{align*}
\iprod{U_0(g)\delta_a h}{U_0(g)\delta_b h'}
&=\iprod{T((gb)^{-1}(ga))h}{h'}
=\iprod{T(b^{-1}a)h}{h'}
=\iprod{\delta_a h}{\delta_b h'}.
\end{align*}
Thus, by sesquilinearity, $U_0(g)$ preserves the semidefinite
inner product and maps $N$ into itself. Moreover,
\begin{align}\label{eq:repres-U}
(U_0(g_1)U_0(g_2)f)(g)
=f(g_2^{-1}g_1^{-1}g)
=f((g_1g_2)^{-1}g)
=(U_0(g_1g_2)f)(g).
\end{align}
In particular, $U_0(g)^{-1}=U_0(g^{-1})$, so $U_0(g)$
maps $N$ onto $N$.

Upon passing to the quotient, we obtain a well-defined
linear mapping $\widetilde U(g)$ on $V/N$, which preserves
inner products. Therefore $\widetilde U(g)$ extends to an
isometry from $\widetilde H$ into itself, again denoted by
$\widetilde U(g)$. By passing to the quotient in
\eqref{eq:repres-U}, we see that
$$
\widetilde U(g_1)\widetilde U(g_2)=\widetilde U(g_1g_2),
$$
so the mapping $\widetilde U:G\to\calL(\widetilde H)$ is a
homomorphism. Since $U_0(g)$ maps $N$ onto $N$, the induced
isometry on $V/N$ is surjective. Its extension to
$\widetilde H$ is therefore surjective as well, and hence
$\widetilde U(g)$ is unitary.

This representation has the desired properties: for all
$g\in G$ and $h,h'\in H$ we have
$$
U_0(g)f_h=U_0(g)\delta_e h=\delta_g h
$$
and consequently, by \eqref{eq:inner-posdef-1},
\begin{align*}
\iprod{J^\star\wt U(g)Jh}{h'}
&=\iprod{\wt U(g)Jh}{Jh'}
=\iprod{\delta_g h}{\delta_e h'}
=\iprod{T(g)h}{h'}.
\end{align*}

For every $g\in G$ and $h\in H$ we have
$$
\wt U(g)Jh=(\one_{\{g\}}\otimes h)+N.
$$
These vectors span $V/N$, which is dense in $\wt H$.
Hence the dilation is minimal.

For uniqueness, let $(\wt H_i,U_i,J_i)$, $i=1,2$, be two
minimal dilations of $T$. For $g,g'\in G$ and $h,h'\in H$,
$$
\iprod{U_i(g)J_i h}{U_i(g')J_i h'}
=\iprod{T({g'}^{-1}g)h}{h'},\quad i=1,2.
$$
Consequently, the mapping
$$
W\Bigl(\sum_{k=1}^N U_1(g_k)J_1h_k\Bigr)
:=\sum_{k=1}^N U_2(g_k)J_2h_k
$$
is well defined and isometric. By minimality, it extends to
a unitary operator $W:\wt H_1\to\wt H_2$ satisfying
$WJ_1=J_2$ and $WU_1(g)=U_2(g)W$ for all $g\in G$.
These identities determine $W$ uniquely on a dense subspace.
\end{proof}

The theorem will be applied in the following situation:

\begin{lemma}\label{lem:contr-pos-def}
If $T\in \calL(H)$ is a contraction, the mapping $S: \Z\to \calL(H)$ defined by
$$ S(n):= \begin{cases}
           T^n\!, & n\ge 1, \\ I, & n=0, \\ (T^\star)^{-n}\!, & n\le -1,
          \end{cases}
$$
is positive definite and satisfies $S^\star(n) = S(-n)$ for all $n\in\Z$.
\end{lemma}
\begin{proof}
Since $T$ is a contraction,
$$
\iprod{(I-T^\star T)x}{x}
=\n x\n^2-\n Tx\n^2\ge0,
$$
so $I-T^\star T$ is positive.
As a consequence, by Proposition \ref{prop:pos-sq-root} the {\em defect operator}\index{defect operator}
$$D_T:= (I-T^\star T)^{1/2}$$ is well defined and positive, and we have
$$\n D_T x\n^2 = \iprod{(I-T^\star T)^{1/2}x}{(I-T^\star T)^{1/2}x} = \iprod{(I-T^\star T)x}{x} = \n x\n^2 - \n Tx\n^2\!.$$

Define
$$
\ell^2(H):=
\Bigl\{\textrm{h}=(h_n)_{n\ge1}:\
h_n\in H \ \hbox{ for all } n\ge1,\
\sum_{n\ge1}\n h_n\n^2<\infty
\Bigr\}.
$$
With respect to the inner product $\iprod{\textrm{g}}{\textrm{h}}:= \sum_{n\ge 1} \iprod{g_n}{h_n}$,
$\ell^2(H)$ is a
Hilbert space; completeness is proved in the same way as for $\ell^2$\!.
We define the operator $\wt T$ on $\ell^2(H)$ by
$$\wt T:\textrm{h} \mapsto (Th_1, D_Th_1, h_2,h_3, \dots).$$
Clearly
$$ \n \wt T \textrm{h}\n_{\ell^2(H)}^2 = \n Th_1\n^2 + \n D_Th_1\n^2 + \sum_{n\ge 2} \n h_n\n^2  = \n \textrm{h}\n^2\!,$$ so $\wt T$ is isometric.

The identity $S^\star(n)=S(-n)$ follows from the definition.
To prove positive definiteness, let $J:H\to\ell^2(H)$ be
the isometric embedding given by $Jh:=(h,0,0,\dots)$.
By induction, the first coordinate of $\wt T^nJh$ equals
$T^nh$ for every integer $n\ge0$.

Since $\wt T$ preserves inner products, for integers
$n\ge m\ge0$ and $h,h'\in H$ we have
\begin{align*}
\iprod{\wt T^nJh}{\wt T^mJh'}
=\iprod{\wt T^{n-m}Jh}{Jh'}
=\iprod{T^{n-m}h}{h'}
=\iprod{S(n-m)h}{h'}.
\end{align*}
For $n<m$, the same identity follows by interchanging
$(n,h)$ and $(m,h')$ and taking complex conjugates, using
$S(m-n)^\star=S(n-m)$.

Now let $k_1,\dots,k_N\in\Z$ and $h_1,\dots,h_N\in H$
be arbitrary. Choose $r\in\Z$ such that $k_j+r\ge0$ for
all $j$. Positive definiteness of $S$ now follows from
\[
\sum_{i,j=1}^N\iprod{S(k_j-k_i)h_j}{h_i}
=\sum_{i,j=1}^N
  \iprod{\wt T^{k_j+r}Jh_j}{\wt T^{k_i+r}Jh_i}
=\Bigl\n\sum_{j=1}^N\wt T^{k_j+r}Jh_j\Bigr\n^2
\ge 0.
\qedhere\]
\end{proof}

Combining the lemma with the theorem, we arrive at the following result.

\begin{theorem}[Sz.-Nagy dilation theorem]\label{thm:Nagy}\index{theorem!Sz.-Nagy}
If $T\in\calL(H)$ is a contraction, then there exist a Hilbert space $\wt H$, an isometric embedding $J:H\to\wt H$,
and a unitary operator $U\in \calL(\wt H)$ such that
$$ T^n  = J^\star U^n J, \quad  n\in\N.
$$
\end{theorem}

In this context the operator $U$ is said to be a {\em unitary dilation}\index{unitary!dilation}\index{dilation!of a contraction} of $T$. As a simple example, the left (right) shift on $\ell^2(\Z)$ is a unitary dilation of the left (right) shift on $\ell^2$\!.

\begin{problems}

\item
Prove the {\em Hellinger--Toeplitz theorem}:\index{theorem!Hellinger--Toeplitz}
if $T: H \to H$ is a linear mapping
satisfying $$\iprod{Tx}{y} = \iprod{x}{Ty}, \quad x, y \in H,$$
then $T$ is bounded.

\noindent{\em Hint:}\ Apply the uniform boundedness theorem to the family of operators
$T_y:H\to\K$, $\|y\|\le1$, given by $T_yx:=\iprod{x}{Ty}$.

\item
Let $H \not = \{0\}$ and let $U,V\in\calL(H)$ be unitary.
\begin{enumerate}[\rm(a), leftmargin=*]
 \item Show that $\frac12 (U+V)$ has norm one if and only if $1\in\sigma(U^\star V)$.

 \noindent{\em Hint:}\ Use $\n U+V\n=\n I+U^\star V\n$ and the fact that
 the norm of a normal operator equals its spectral radius.

 \item  Prove that $\frac12(U+V)$ is unitary if and only if $U=V$.

 \noindent{\em Hint:}\ For the `only if' part, use that
 $\n Ux+Vx\n=2\n x\n$ for all $x\in H$, and consider the case of equality in the triangle inequality.
\end{enumerate}

\item\label{prob:unitary-spectrum}
Deduce Proposition \ref{prop:spectrum-unitary} from Corollary  \ref{cor:spectrum-isometry}.

\item\label{prob:sec-rot}
For $\theta\in \R$, show that the rotation operator on $L^2(\mathbb{T})$ defined by $$R_\theta f(e^{ ix}) := f(e^{i(x-\theta)})$$
is unitary, and find its spectrum.

\noindent{\em Hint:}\ Distinguish the cases $\theta/2\pi\in \Q$ and
$\theta/2\pi\not\in \Q$.

\item\label{prob:selfadjoint-spectrum}
Let $T\in \calL(H)$ be selfadjoint. The aim of this problem is to deduce the inclusion $\si(T)\subseteq \R$ in an elementary way from Proposition \ref{prop:spectrum-unitary}. Fix $\la\in\si(T)$.
\begin{enumerate}[\rm(a), leftmargin=*]
  \item Show that $$e^{i\la} - e^{iT} = e^{i\la}(\la-T)\Bigl(\sum_{n=1}^\infty \frac{i^n}{n!}(T-\la)^{n-1}\Bigr)$$
  and conclude that $e^{i\la} - e^{iT}$ fails to be invertible.
  \item Combine this with Proposition \ref{prop:spectrum-unitary} to conclude that $\la\in\R$.
\end{enumerate}

\item\label{prob:Cayley}
The {\em Cayley transform}\index{Cayley transform} connects selfadjoint and unitary operators.
Let $T\in\calL(H)$ be selfadjoint.

\begin{enumerate}[\rm(a), leftmargin=*]
\item Show that $T+iI$ is invertible and that
$U:=(T-iI)(T+iI)^{-1}$
is unitary. Show moreover that $1\notin\si(U)$.

\item Conversely, let $U\in\calL(H)$ be unitary and suppose that
$1\notin\si(U)$. Show that
$T:=i(I+U)(I-U)^{-1}$ is selfadjoint and that its Cayley transform is $U$.

\item Conclude that the Cayley transform defines a bijection between
the bounded selfadjoint operators on $H$ and the unitary operators $U$
satisfying $1\notin\si(U)$.
\end{enumerate}

\item\label{prob:hermit}
Using Proposition \ref{prop:HS-adjoint}, prove that if
$A$ is an $m\times n$ matrix with complex coefficients, viewed as a bounded operator from $\C^n$ to $\C^m$\!, then
  \begin{align*}
  \|A \|^2 = \max \{ \lambda \ge 0 : \, \lambda \text{ is an eigenvalue of } A^\star A\}.
  \end{align*}

\item
Show that the norm of an operator $T\in \calL(H)$ is given by
$$\n T\n^2 = \inf\{\la \ge 0 :\, T^\star T \le \la I \},$$
where $ T^\star T \le \la I$ means that $\la - T^\star T$ is positive.

\item\index{inequality!Cauchy--Schwarz}
Show that if $T\in\calL(H)$ is a positive operator, then for all $x, y \in H$ we have the Cauchy--Schwarz-type inequality
$$|\iprod{Tx}{y}|^2 \le \iprod{Tx}{x}\iprod{Ty}{y}.$$

\item\label{prob:S1leqS2}
Show that if the operators $S_1, S_2 \in \calL(H)$ satisfy $0\le S_1\le S_2$, then for all $T\in \calL(H)$ we have $0\le T^\star S_1 T \le T^\star S_2 T$.

\item
Let $T\in\calL(H)$ be normal. Show that $\Ker(T)=\Ker(T^\star)$ and $\ov{\Ran(T)}=\ov{\Ran(T^\star)}$.

\item
Let $T\in \calL(H)$ be selfadjoint and let $\la\in\R$.
\begin{enumerate}[\rm(a), leftmargin=*]
  \item Show that if $\sigma(T) = \{\la\}$, then $T = \la I$.
  \item Show that if $\sigma(PTP) = \{\la\}$, where $P$ is an orthogonal projection with $\Ker(P)$ finite-dimensional, then $\la I-T$ is compact;
  if, in addition, $P\not=I$, then $\la=0$ and $T$ is compact.
\end{enumerate}

\item\label{prob:Volterra}\index{Volterra operator}
Consider the Volterra operator $T$ on $L^2(0,1)$ of Example \ref{ex:Volterra}:
$$ (Tf)(t) = \int_0^t f(s)\ud s, \quad f\in L^2(0,1), \ t\in [0,1].$$
We sketch two proofs that $\sigma(T) = \{0\}$.
\begin{enumerate}[\rm(a), leftmargin=*]
  \item\label{it:Volterra1} Show that for all $0\not=\la \in\C$ and $g\in L^2(0,1)$ the equation $(\la - T)f = g$ has a unique solution in $L^2(0,1)$. Deduce that $\sigma(T) = \{0\}$.
\end{enumerate}
A second proof is obtained by estimating the norm of $T^n$:
\begin{enumerate}[\rm(a), leftmargin=*]\setcounter{enumii}{1}
  \item\label{it:Volterra2} Show that $\|T^n\|\leq \frac1{n!}$ for all $n\in\N$.

  \noindent{\em Hint:}\ The case $n=0$ is trivial. For $n\ge1$, first show that
  \begin{align*}(T^n f)(t) & = \int_0^t \int_{0}^{t_{n-1}} \cdots \int_0^{t_1}f(s) \ud s \ud t_1 \cdots \ud t_{n-1}
  \\ & = \int_0^t f(s) \int_{s}^{t} \int_s^{t_{n-1}} \cdots \int_{s}^{t_2} \ud t_1 \cdots \ud t_{n-2}\ud t_{n-1}\ud s.
  \end{align*}

  \item\label{it:Volterra3} Using Theorem \ref{thm:rT}, conclude that $\sigma(T) = \{0\}$.

  \item Show that $0$ is not an eigenvalue of $T$.

  \item Deduce from \ref{it:Volterra1} or \ref{it:Volterra3}
  that $T$ is not normal.
\end{enumerate}

\item
The {\em numerical range}\index{numerical!range} of an operator $T\in \calL(H)$ with $H\not=\{0\}$ is the set
$$ W(T) := \{\iprod{Tx}{x}:\, \n x\n=1\}.$$
The {\em numerical radius}\index{numerical!radius} of $T$ is defined by
$$ w(T):= \sup\{|\la|:\, \la \in W(T)\}.$$
Prove that $$\frac12\n T\n \le w(T)\le \n T\n.$$

\noindent {\em Hint:}\ To prove the first inequality, use the identity
  \begin{align*}
  4\iprod{Tx}{y} & = \iprod{T(x+y)}{x+y}-\iprod{T(x-y)}{x-y}\\ & \qquad +i\iprod{T(x+iy)}{x+iy}-i\iprod{T(x-iy)}{x-iy}.
  \end{align*}

\item\label{prob:Toep-Haus}
This problem proves the {\em Toeplitz--Hausdorff theorem},\index{theorem!Toeplitz--Hausdorff}
which asserts that the numerical range of any operator $T\in \calL(H)$ is a convex subset of $\C$.
\begin{enumerate}[\rm(a), leftmargin=*]
  \item Show that for all $\la,\mu\in \C$ we have $$W(\la T + \mu I) = \la W(T) + \mu.$$
  Conclude that in order to prove the Toeplitz--Hausdorff theorem it suffices to establish that $\{0,1\}\subseteq W(T)$ implies $[0,1]\subseteq W(T)$.
\end{enumerate}
In what follows we fix an operator $T\in \calL(H)$ such that $\{0,1\}\subseteq W(T)$, and prove that $[0,1]\subseteq W(T)$.

Choose norm-one vectors $x,y\in H$ such that $\iprod{Tx}{x} = 0$ and  $\iprod{Ty}{y} = 1$.
\begin{enumerate}[\rm(a), leftmargin=*]\setcounter{enumii}{1}
  \item Define $g:[0,2\pi]\to \C$ by $$g(t):= e^{-it}\iprod{Tx}{y} + e^{it}\iprod{Ty}{x}.$$ Using the fact that $g(t+\pi) = -g(t)$ for $t\in [0,\pi]$, show that either there exists $t_0\in [0,\pi]$ such that $g(t_0)=0$ or else there exists $t_0\in [0,2\pi]$ such that $g(t_0)$ is a positive real number.
\end{enumerate}
Set $\wt y:= e^{it_0}y$.
\begin{enumerate}[\rm(a), leftmargin=*]\setcounter{enumii}{2}
  \item\label{it:Toep-Haus2} Show that $x$ and $\wt y$ are linearly independent.
\end{enumerate}
Define $z:[0,1]\to H$ and $f:[0,1]\to \C$ by
$$ z(t):= \frac{(1-t)x+ t\wt y}{\n(1-t)x+ t\wt y\n}, \qquad f(t):= \iprod{Tz(t)}{z(t)}.$$
These functions are well defined by part \ref{it:Toep-Haus2}.
\begin{enumerate}[\rm(a), leftmargin=*]\setcounter{enumii}{3}
  \item Show
  $f$ is continuous and real-valued and satisfies
  $f(0)=0$ and $f(1)=1$. Deduce that $[0,1]\subseteq W(T)$.
\end{enumerate}

\item
Prove that for all $T\in \calL(H)$ with $H\not=\{0\}$ we have
$$\sigma(T)\subseteq \ov{W(T)}.$$

\noindent{\em Hint:}\ First prove that approximate eigenvalues belong to $\ov{W(T)}$. Then apply the Toeplitz--Hausdorff theorem in
combination with Proposition \ref{prop:approx-eigenvalue-bdry}.

\item\label{prob:normal-numerical-range}
Prove that if $H$ is nonzero and  $T\in\calL(H)$ is normal, then $\ov{W(T)}$ is the closed
convex hull of $\si(T)$.

\noindent{\em Hint:}\ For the inclusion `$\supseteq$' use the results of the preceding two problems. For the inclusion `$\subseteq$', show that if $z\in W(T)$ does not belong
to the closed convex hull of $\si(T)$, then for a suitable $\alpha\in\C$ and $r\in\R$ we have $\Re(\alpha z)>r\ge\Re(\alpha\la)$ for all $\la\in\si(T)$.
Then apply Corollary \ref{cor:normal-cont-fc-pos} to
$f(\la):=r-\Re(\alpha\la)$.

\item\label{prob:FMO-composition}
Let $T_m$ be a Fourier multiplier operator on $L^2(\R^d)$ with symbol $m\in L^\infty(\R^d)$.
\begin{enumerate}[\rm(a), leftmargin=*]
 \item Show that $T_m$ is normal.
 \item Show that if $f$ is a continuous function
on the essential range of $m$ (see Problem \ref{prob:spectrum-FMO}), then,  after modifying $m$ on a null set if necessary, $f(T_m)$ is well defined via the continuous functional calculus and equal to the Fourier multiplier $T_{f\circ m}$ with symbol $f\circ m \in L^\infty(\R^d)$.
 \item Compare this result with Problem \ref{prob:spectrum-FM}.
\end{enumerate}

\item\label{prob:normal-invariant-subspace}
Let $T\in\calL(H)$ be a normal operator which is not a scalar
multiple of the identity.

\begin{enumerate}[\rm(a), leftmargin=*]
\item\label{it:invariant-subspace1} Show that $\si(T)$ contains two distinct points and choose nonzero
functions $f,g\in C(\si(T))$ such that $fg=0$.

\item\label{it:invariant-subspace2} Put $M:=\ker(f(T))$. Use the continuous functional calculus to show that $\{0\}\subsetneq M\subsetneq H$.
Show moreover that both $M$ and $M^\perp$ are invariant under $T$.

\noindent{\em Hint:}\ Use the identities $f(T)g(T)=(fg)(T)=0$ and
$\n h(T)\n=\n h\n_\infty$.

\item\label{it:invariant-subspace3} Use the Fuglede--Putnam--Rosenblum theorem to show that $M$ is
in fact hyperinvariant for $T$, that is, invariant under every bounded operator commuting with $T$.

\noindent{\em Hint:}\ If $ST=TS$, first show that
$S$ commutes with $T^\star$ and hence with $h(T)$ for every
$h\in C(\si(T))$.

\item\label{it:invariant-subspace4} Conclude that every normal operator on a Hilbert space of dimension greater than one has a nontrivial closed invariant subspace.
\end{enumerate}

\item
Show that if $S,T\in\calL(H)$ are positive operators, then $\si(ST)\subseteq [0,\infty)$.

\noindent{\em Hint:}\ Apply the result of Problem \ref{prob:sigmaSTvsTS}
to the operators $S^{1/2}$ and $S^{1/2}T$.

\item
Consider a nonzero operator $T \in \calL(H)$. Show that the following assertions are equivalent:
\begin{enumerate}[\rm(a), leftmargin=*]
\item $T$ is positive;
\item $T$ is selfadjoint and there exists a bounded $S$ such that $T = S^\star S$;
\item $T$ is selfadjoint and there exists a selfadjoint $S$ such that $T = S^2$;
\item $T$ is selfadjoint and $\big\|I-T/\|T\|\big\| \leq 1.$
\end{enumerate}

\item\label{prob:TandmodTequal-ranges}
Show that for all $T\in \calL(H)$ we have
$\Ran(T) = \Ran((TT^\star)^{1/2}).$

\noindent {\em Hint:}\ Apply the result of Problem \ref{prob:inclusion-of-ranges-Hilbert}.

\item
Let $T\in\calL(H)$ have polar decomposition $T=U|T|$, with $U$ the
partial isometry of Theorem \ref{thm:polar}. Show that $T$ is normal
if and only if $U$ is normal and $U|T| = |T|U$.

\noindent{\em Hint:}\ Recall that
$U^\star U$ is the orthogonal projection onto
$\ov{\Ran(|T|)}$, whereas $UU^\star$ is the orthogonal projection
onto $\ov{\Ran(T)}$. If $T$ is normal, use
$\Ker(T)=\Ker(T^\star)$ and the uniqueness of positive square roots.
For the converse, compare $T^\star T$ and $TT^\star$.

\item
Let $H\not=\{0\}$, let $T\in\calL(H)$ be invertible, with polar decomposition $T=U|T|$ with $U$  unitary. Show that
$$
\n T-U\n
=\inf\{\n T-V\n:\, V\in\calL(H)\hbox{ unitary}\}
$$
and that both are equal to $\max\{\n T\n-1,\, 1-1/\n T^{-1}\n\}$.

\noindent{\em Hint:}\ For unitary $V$, use the triangle inequality
and the identity $\inf_{\n x\n=1}\n Tx\n=1/\n T^{-1}\n$. For $V=U$, use the identity
$\n T-U\n=\n |T|-I\n$ and the continuous functional calculus.

\item\label{prob:normal-algebraic}
Let $T\in\calL(H)$ be normal and {\em algebraic}, that is,
$p(T)=0$ for some nonzero polynomial $p$. Show that there exist
distinct scalars $\la_1,\dots,\la_m$ and pairwise orthogonal
projections $P_1,\dots,P_m$ such that
$$
P_1+\cdots+P_m=I,
\qquad
T=\sum_{j=1}^m\la_jP_j.
$$
In particular, $H$ is the orthogonal direct sum of the eigenspaces
of $T$.

\noindent{\em Hint:}\ By the spectral mapping theorem,
$\si(T)$ is contained in the finite set of zeroes of $p$.
Write $\si(T)=\{\la_1,\dots,\la_m\}$ and, for each $j$, let
$f_j\in C(\si(T))$ be the indicator function of $\{\la_j\}$.
Put $P_j=f_j(T)$ and use the continuous functional calculus.

\item
Prove that if $T\in \calL(H)$ is an isometry, then there exist Hilbert spaces $G$ and $K$ such that
we have an isometric isomorphism of Hilbert spaces $$H \simeq \ell^2(G)\oplus K,$$
where $\ell^2(G)$ is the Hilbert space of all square-summable sequences ${\rm g} = (g_n)_{n\ge 1}$ in $G$ with norm $\n{\rm g}\n_{\ell^2(G)}^2 = \sumn \n g_n\n^2$
and, with respect to this decomposition,
$$T\simeq S\oplus U,$$
where $S$ is the right shift on $\ell^2(G)$, that is, $S$ maps the sequence $g_1,g_2,\hdots$ to $0,g_1,g_2,\hdots$ and $U$ is a unitary operator on $K$. This decomposition is known as the {\em Wold decomposition}\index{Wold decomposition}.

\noindent
{\em Hint:}\ For $n\in\N$ let $H_n:= \ran(T^n)$, and for $n\ge 1$
let $G_n$ denote the orthogonal complement of $H_n$ in $H_{n-1}$.
Show that the spaces $G_n$ are all isometric as Hilbert spaces and set $K:= \bigcap_{n\in\N} H_n$.

\end{problems}

%% file: ch09-SpectralTheorem-bdd.tex
\chapter{The Spectral Theorem for Bounded Normal Operators}\label{chap:spectral-theorem}

\blfootnote{\rm \small \begin{quote}This book has been published by Cambridge University Press in the series `Cambridge Studies in Advanced Mathematics'. The present version is free to view and download for personal use only. Not for re-distribution, re-sale or use in derivative works.\newline {\copyright} Jan van Neerven\end{quote}}

\noindent
In this chapter we show that normal operators admit a spectral representation as sums or integrals of orthogonal projections. We begin by showing that
every compact normal operator admits the spectral decomposition $ T = \sum_{n\ge 1}\la_n P_n,$
where $(\la_n)_{n\ge 1}$ is the sequence of eigenvalues and $(P_n)_{n\ge 1}$ is the sequence of orthogonal projections onto the corresponding eigenspaces. The main result of this chapter proves that an arbitrary bounded normal operator admits an analogous integral representation
$T = \int_{\sigma(T)} \la \ud P(\la)$ with respect to a projection-valued measure $P$ on the spectrum $\si(T)$. This, in turn, can be used to uncover various deeper properties of such operators.

\section{The Spectral Theorem for Compact Normal Operators}\label{sec:spectral-thm-normal-compact}

Throughout this chapter, $H$ is a complex Hilbert space.

From linear algebra we know that normal matrices can be orthogonally diagonalised. This result admits the following extension to compact normal operators on $H$:

\begin{theorem}[Spectral theorem for compact normal
operators]\label{thm:spect-thm-comp}\index{theorem!spectral, for compact normal operators}
 Let $T\in\calL(H)$ be a compact normal operator and let
 $(\la_n)_{n\ge 1}$ be the (finite or infinite) sequence of its distinct eigenvalues.
 Let $(E_n)_{n\ge 1}$ be the
 corresponding sequence of eigenspaces, and let $(P_n)_{n\ge 1}$ be the associated sequence of orthogonal projections.
 Then:
 \begin{enumerate}[label={\rm(\arabic*)}, leftmargin=*]
  \item\label{it:spect-thm-comp1} the spaces $E_n$ are pairwise orthogonal and have dense linear span;
  \item\label{it:spect-thm-comp2} we have $$ T = \sum_{n\ge 1}\la_n P_n$$ with convergence in the operator norm of $\calL(H)$.
 \end{enumerate}
\end{theorem}

\begin{proof}
The proof of the theorem uses the properties of spectra of compact operators on Banach spaces established in Theorem \ref{thm:comp-spec}.
In the present situation, where the compact operator acts on a Hilbert space, the proof of this theorem
can be considerably shortened; see  Problem \ref{prob:spectrum-compact-normal}.

\smallskip
\ref{it:spect-thm-comp1}: \
From Proposition \ref{prop:normal-spectralradius} (applied to $T - \la$) we see that
$Tx - \la x = 0$ if and only if $T^\star x- \ov\la x = 0$, so $\la$ is an eigenvalue of
$T$ if and only if $\ov \la$ is an eigenvalue of $T^\star$ and the eigenspaces coincide.

If $y_m\in E_m$ and $y_n \in E_n$ are nonzero vectors, then
 $$\la_m\iprod{y_m}{y_n} = \iprod{Ty_m}{y_n}= \iprod{y_m}{T^\star y_n}=
\iprod{y_m}{\ov {\la_n} y_n} = \la_n \iprod{y_m}{y_n}.$$  If $\la_m\not=\la_n$, then this
is possible only if $\iprod{y_m}{y_n}=0$. This gives $E_n \perp E_m$ for $n\not=m$.

Let $E:= \bigoplus_{n\ge 1} E_n$ denote the closed linear span of the spaces $E_n$, $n\ge 1$.
We wish to prove that $E = H$. Suppose the contrary. Then $E^\perp$ is a nonzero closed subspace of $H$.
Since $TE_n\subseteq E_n$ for all $n\ge 1$, we have $TE\subseteq E$. Furthermore, if $x\in E_n$, then $T^\star x = \ov\la_n x \in E_n$, so $T^\star E_n\subseteq E_n$.
Since this holds for all $n\ge1$, it follows that $T^\star E\subseteq E$.
Hence, by Lemma \ref{lem:normal-reduce} we have $T E^\perp \subseteq E^\perp$ and the restriction $T^\perp$ of $T$ to $E^\perp$ is normal.
Moreover, by the very construction of $E$, $T^\perp$ has
no eigenvalues. By Theorem \ref{thm:comp-spec}, this implies that $\si(T^\perp)\subseteq\{0\}$,
and since $\si(T^\perp)\not=\emptyset$, it follows that $\si(T^\perp) = \{0\}$. Now Proposition \ref{prop:normal-spectralradius}
implies that $T^\perp = 0$. This means that every nonzero element of $E^\perp$ is an eigenvector of $T^\perp$ with
eigenvalue $0$. This contradicts the observation just made that $T^\perp$ has no eigenvalues and completes the proof that $E=H$.

\smallskip
\ref{it:spect-thm-comp2}: \ Let $P_n$ denote the orthogonal projection onto $E_n$. Then for all $x\in H$ we have $x = \sum_{n\ge 1}P_n x$ with convergence in $H$. This is clear for every $x$ in the span of the spaces $E_n$, and since this span is dense in $H$ and the operators $\sum_{n=1}^N P_n$ are orthogonal projections and hence have norm one, the convergence extends to all $x\in H$ by Proposition \ref{prop:approxTn}.
It follows that the sum $\sum_{n\ge 1}TP_n x$ converges as well, with sum $Tx$.

Fix $\e>0$. The set $\Lambda_\e := \{n\ge 1:\, |\la_n| > \e\}$
is finite by Theorem \ref{thm:comp-spec}. Let $N\ge 1$ be so large that $\Lambda_\e\subseteq \{1,2,\dots,N\}$.
Fixing $x\in H$ and writing $x_n:= P_n x$, by orthogonality we have
\begin{align*}
 \Bigl\n Tx - \sum_{n=1}^N \la_n P_n x \Bigr\n^2
 & =  \Bigl\n \sum_{n\ge 1} Tx_n - \sum_{n=1}^N \la_n x_n \Bigr\n^2
  =  \Bigl\n \sum_{n\ge N+1} \la_n x_n\Bigr\n^2
 \\ & = \sum_{n\ge N+1} |\la_n|^2  \n x_n\n^2
 \le \e^2 \sum_{n\ge N+1} \n x_n\n^2
 \le \e^2 \n x\n^2\!.
 \end{align*}
Taking the supremum over all $x\in H$ with $\n x\n\le 1$, we obtain
$$ \Bigl\n T - \sum_{n=1}^N \la_n P_n \Bigr\n^2 \le \e^2\!.$$
This completes the proof.
\end{proof}

Let $H$ and $K$ be Hilbert spaces. For $h\in H$ and $k\in K$ we denote by $k\,\bar\otimes\, h$\index{$Gz$@$g\,\bar\otimes\, h$} the operator in $\calL(H,K)$ defined by
$$ (k\,\bar\otimes\, h)x := \iprod{x}{h}k, \quad x\in H.$$
If $H=K$ and $h\in H$ has norm one, then $h\,\bar\otimes\, h$ is the orthogonal projection onto the subspace spanned by $h$.
If $T\in \calL(H)$ is a compact normal operator, the eigenspaces corresponding to nonzero eigenvalues are finite-dimensional. Choosing
orthonormal bases for each of them, from
Theorem \ref{thm:spect-thm-comp} we obtain a representation
 $$ T = \sum_{n\ge 1} \la_n h_n\,\bar\otimes\,h_n$$
 with convergence in the operator norm,
 where now $(\la_n)_{n\ge 1}$ is the sequence of nonzero eigenvalues of $T$ {\em repeated according to multiplicities} and $(h_n)_{n\ge 1}$ is an associated orthonormal sequence of eigenvectors.

 Strictly speaking, Theorem \ref{thm:spect-thm-comp} gives convergence
blockwise, one eigenspace at a time. The stated rank-one expansion follows by
the same tail estimate. Indeed, if $T=\sum_{n\ge1}\lambda_n h_n\bar\otimes h_n$,
with the eigenvalues repeated according to multiplicity and ordered so
that $|\lambda_n|\to0$, then for $N<M$ and $\|x\|\le1$,
$$
\Bigl\|
   \sum_{n=N}^{M}\lambda_n
   (h_n\bar\otimes h_n)x
\Bigr\|^2
=\sum_{n=N}^{M}|\lambda_n|^2|\iprod{x}{h_n}|^2
\le
\sup_{n\ge N}|\lambda_n|^2\!.
$$
Hence the tails converge to $0$ in operator norm.
The geometric and algebraic multiplicities of the eigenvalues coincide by Corollary \ref{cor:copact-normal-multipicities}, so we can unambiguously speak of their {\em multiplicity}.

Theorem \ref{thm:spect-thm-comp} allows us to deduce the following general representation theorem for compact operators acting on a Hilbert space. It strengthens Proposition \ref{prop:fr-dense-compH}, which asserted that such operators can be approximated in operator norm by finite-rank operators.

\begin{theorem}[Singular value decomposition]\label{thm:sing-value-comp}\index{theorem!singular value decomposition}\index{singular value!decomposition}
Let $T\in \calL(H,K)$ be a compact operator, where $K$ is another Hilbert space. Then
$T$ admits a decomposition
 $$ T = \sum_{n\ge 1} \mu_n k_n\,\bar\otimes\,h_n$$
 with convergence in the operator norm,
 where $(\mu_n)_{n\ge 1}$ is the sequence of nonzero eigenvalues of the compact operator $(T^\star T)^{1/2}$ repeated according to multiplicities, and $(h_n)_{n\ge 1}$ and $(k_n)_{n\ge 1}$ are orthonormal sequences in $H$ and $K$, respectively, the former consisting of eigenvectors of $(T^\star T)^{1/2}$.
\end{theorem}

The proof depends on the following lemma.

\begin{lemma}\label{lem:absTcompact}
If $S\in \calL(H)$ is a positive compact operator, then its square root $S^{1/2}$ is compact.
\end{lemma}
\begin{proof}
Let $(x_n)_{n\ge1}$ be a bounded sequence in $H$.
Since $S$ is compact, after passing to a subsequence we may
assume that $(Sx_n)_{n\ge1}$ converges.
Then
\begin{align*}
\|S^{1/2}x_n-S^{1/2}x_m\|^2
&=\iprod{S(x_n-x_m)}{x_n-x_m}\le\|Sx_n-Sx_m\|\,\|x_n-x_m\|
\to 0
\end{align*}
as $n,m\to\infty$.
Thus $(S^{1/2}x_n)_{n\ge1}$ is a Cauchy sequence and hence
converges in $H$.
\end{proof}

\begin{proof}[Proof of Theorem \ref{thm:sing-value-comp}]
By Theorem \ref{thm:spect-thm-comp}, applied to $|T|:= (T^\star T)^{1/2}$, which is compact by Lemma \ref{lem:absTcompact}, we arrive at a representation
 $$|T| = \sum_{n\ge 1} \mu_n h_n\,\bar\otimes\, h_n,$$
with convergence in the operator norm, where $(\mu_n)_{n\ge 1}$ is the sequence of nonzero eigenvalues of $|T|$
  repeated according to multiplicities and the orthonormal sequence $(h_n)_{n\ge 1}$ consists of eigenvectors of $|T|$.
Let $T = U|T|$, with $U$ an isometry from $\ov{\Ran(|T|)}$ onto $\ov{\Ran(T)}$ as in Theorem \ref{thm:polar}.
The sequence $(k_n)_{n\ge 1}$ defined by $k_n := Uh_n$ is orthonormal in $K$ and $$ T = \sum_{n\ge 1} \mu_n k_n\,\bar\otimes\,h_n$$ with convergence in the operator norm.
\end{proof}

As a second application of Theorem \ref{thm:spect-thm-comp}, we record the following formulas for the eigenvalues of a compact positive operator.

\begin{theorem}[Min-max theorem]\label{thm:minmax}\index{theorem!min-max} Let $T\in \calL(H)$ be compact and positive, and let $\la_1\ge \la_2\ge \cdots > 0$ be the sequence of its nonzero eigenvalues repeated according to multiplicities.
Then these eigenvalues are given by
$$ \la_n = \inf_{\substack{Y\subseteq H \\ \dim(Y) = n-1}} \sup_{\substack{\n y\n=1 \\ y \perp Y}} \iprod{Ty}{y} = \inf_{\substack{Y\subseteq H \\ \dim(Y) = n-1}} \sup_{\substack{\n y\n=1 \\ y \perp Y}} \n Ty\n,$$
where the infima are taken over all subspaces $Y$ of $H$ of dimension $n-1$.
\end{theorem}
\begin{proof}
By the observation after the proof of Theorem \ref{thm:spect-thm-comp}, we may choose an orthonormal family of eigenvectors
$(h_j)_{1\le j\le N}$
if the sequence of nonzero eigenvalues is finite, say
$\lambda_1\ge\cdots\ge\lambda_N>0$, and an orthonormal sequence
$(h_j)_{j\ge 1}$ if it is infinite, such that
$ Th_j=\lambda_jh_j $ for all indices $j$
and the closed linear span of these vectors equals $(\ker(T))^\perp$.
In the finite case, all sums below over indices $j\ge n$ are understood to end at $N$.

For $n=1$, the only subspace $Y$ to be considered is $\{0\}$. Since
$T$ is positive and selfadjoint, Theorem \ref{thm:spect-sa} gives
$$
   \sup_{\|y\|=1}\iprod{Ty}{y}=\sup_{\|y\|=1}\|Ty\|=\|T\|=\lambda_1.
$$
In the remainder of the proof we may therefore assume that $n\ge 2$.

Let $Y\subseteq H$ be any subspace of dimension $n-1$ and let $H_n$ denote the linear span of the vectors $h_1,\dots,h_n$. Then $Y^\perp \cap H_n$ is a nonzero subspace of $H$, so it contains a norm-one vector $y$.
Writing $y = \sum_{j=1}^n c_j h_j$ with $\sum_{j=1}^n |c_j|^2=1$,
$$\iprod{Ty}{y} =  \sum_{j=1}^n \la_j| c_j|^2 \ge \la_n  \sum_{j=1}^n |c_j|^2 = \la_n.$$
This proves the inequality
$$ \la_n \le \inf_{\substack{Y\subseteq H \\ \dim(Y) = n-1}} \sup_{\substack{\n y\n=1 \\ y \perp Y}} \iprod{Ty}{y}.$$
The inequality
$$  \inf_{\substack{Y\subseteq H \\ \dim(Y) = n-1}} \sup_{\substack{\n y\n=1 \\ y \perp Y}} \iprod{Ty}{y} \le \inf_{\substack{Y\subseteq H \\ \dim(Y) = n-1}} \sup_{\substack{\n y\n=1 \\ y \perp Y}} \n Ty\n$$
holds trivially.
To prove the inequality
$$ \inf_{\substack{Y\subseteq H \\ \dim(Y) = n-1}} \sup_{\substack{\n y\n=1 \\ y \perp Y}} \n Ty\n \le \la_n,$$
let $y\perp H_{n-1}$ have norm one. Then $y=z+\sum_{j\ge n}\iprod{y}{h_j}h_j$ with $z\in\ker(T)$, and hence
$$ \n Ty\n^2 = \Big\n \sum_{j\ge n} \la_j\iprod{y}{h_j}h_j\Big\n^2
= \sum_{j\ge n} \la_j^2 | \iprod{y}{h_j}|^2 \le \la_n^2 \sum_{j\ge n} | \iprod{y}{h_j}|^2 \le \la_n^2,
$$
and the result follows.
\end{proof}

\begin{corollary}\label{cor:poscomp-domination}
If $S,T\in \calL(H)$ are compact operators satisfying $0\le S\le T$, and if
$\la_1\ge \la_2\ge \dots > 0$ and $\mu_1\ge \mu_2\ge \dots > 0$ are the eigenvalue sequences of $S$ and $T$, respectively, repeated according to multiplicities and extended by zeroes if necessary, then for all $n\ge 1$ we have $\la_n\le \mu_n.$
\end{corollary}

\section{Projection-Valued Measures}\label{sec:PVM}

This section and the next deal with the preliminaries needed to state and prove the spectral theorem for bounded normal operators.

Let $(\Om,\calF)$ be a measurable space.

\begin{definition}[Projection-valued measures]\index{projection-valued measure}
A {\em projection-valued measure} on a measurable space $(\Om,\calF)$
is a mapping $P: \calF\to \calL(H)$ that assigns to every set $F\in\calF$ an orthogonal projection $P_F:=P(F)\in \calL(H)$ such that the following conditions are satisfied:
\begin{enumerate}[label={\rm(\roman*)}, leftmargin=*]
 \item\label{it:PVM1} $P_\Om = I$;
 \item\label{it:PVM2} for all $x\in H$ the mapping $F\mapsto \iprod{P_F x}{x}$
 defines a measure on $(\Om,\calF)$.
\end{enumerate}
\end{definition}

For $x\in H$ the measure defined by \ref{it:PVM2} is denoted by $P_{x}$. Thus, for all $F\in \calF$\!,
$$ \iprod{P_F x}{x}  = P_{x}(F) =\int_{\Om} \one_F \ud P_{x}.$$
From $P_{x}(\Om) = \iprod{P_{\Om}x}{x} = \iprod{x}{x} = \n x\n^2$
we see that $P_x$ is a finite measure.

We make some easy observations:

\begin{itemize}
 \item $P_{\emptyset} = 0$.
\end{itemize}
Indeed, the additivity of $P_{x}$, applied to ${\Om} = {\Om}\cup\emptyset$, implies
$$ \iprod{x}{x}  = P_{x}(\Om) = P_{x}({\Om}\cup\emptyset)
= P_{x}(\Om) + P_{x}(\emptyset) = \iprod{x}{x} + P_{x}(\emptyset)
$$ and therefore $ \iprod{P_\emptyset x}{x} =  P_{x}(\emptyset) = 0$ for all $x\in H$.

\begin{itemize}
 \item {\em If $F_1,F_2\in \calF$ are disjoint, then $P_{F_1\cup F_2} = P_{F_1} + P_{F_2}$ and the ranges of $P_{F_1}$ and $P_{F_2}$ are orthogonal.}
\end{itemize}
Since each $P_F$ is an orthogonal (hence selfadjoint) projection, for all $F\in \calF$ and $x\in H$ we have $\n P_F x\n^2 = \iprod{P_F x}{P_F x} = \iprod{P^2_F x}{x} = \iprod{P_Fx}{x}$.

If $x \in \ran(P_{F_1})$, then $x = P_{F_1}x$ and therefore, since $F_1$ and $F_2$ are disjoint,
\begin{align*}
\n x\n^2 = \n P_{F_1}x\n^2 &\le \n P_{F_1}x\n^2 + \n P_{F_2}x \n^2
 = \iprod{P_{F_1}x}{x}+ \iprod{P_{F_2}x}{x}
 \\ & = P_x(F_1)+ P_x(F_2) = P_x(F_1\cup F_2) =
\iprod{P_{F_1\cup F_2}x}{x}
\le \n x\n^2
\end{align*}
and therefore each step holds with equality. In particular, the second and third terms are equal, and therefore $\n P_{F_2}x\n^2 = 0$, that is, $P_{F_2}x = 0$. This shows that $\ran(P_{F_1}) \subseteq \ker (P_{F_2}) = (\ran(P_{F_2}))^\perp$ and therefore $\ran(P_{F_1}) \perp \ran(P_{F_2})$. The identity $P_{F_1\cup F_2} = P_{F_1} + P_{F_2}$ also follows, by polarisation.

\begin{itemize}
 \item {\em For all $F_1,F_2\in \calF$ we have $P_{F_1\cap F_2} = P_{F_1} P_{F_2}  = P_{F_2}P_{F_1}$.}
\end{itemize}
In the special case of disjoint sets, this has just been proved, with all three expressions equal to $0$. From this special case it follows that
\begin{align*}
P_{F_1} P_{F_2} & = (P_{F_1\setminus F_2} + P_{F_1\cap F_2}) (P_{F_2\setminus F_1}+P_{F_1\cap F_2})
\\ & = P_{F_1\setminus F_2}P_{F_2\setminus F_1} + P_{F_1\cap F_2}P_{F_2\setminus F_1}
+P_{F_1\setminus F_2}P_{F_1\cap F_2} + P_{F_1\cap F_2}^2
 = 0+0+0+ P_{F_1\cap F_2}.
\end{align*}
Interchanging the roles of $F_1$ and $F_2$ gives the other identity.

\begin{itemize}
\item {\em If $F\in\calF$ and $x\in\Ran(P_F)$, then $P_x(\complement F) = 0$.}
\end{itemize}
Let $F\in\calF$ and pick $x\in\Ran(P_F)$, so $x=P_Fx$. For every $E\in\calF$ we have
$P_EP_F=P_{E\cap F}$, and therefore
$$
P_x(E)=\iprod{P_Ex}{x}
=\iprod{P_EP_Fx}{x}
=\iprod{P_{E\cap F}x}{x}
= P_x(E\cap F).
$$
Taking $E=\complement F$ gives $P_x(\complement F)=P_x((\complement F)\cap F)=P_x(\emptyset)=0$.

\begin{itemize}
 \item If the sets $F_1,F_2,\dots \in \calF$ are pairwise disjoint and
$F:=\bigcup_{n\ge1}F_n$, then for all $x\in H$ the series $\sum_{n\ge1}P_{F_n}x$ converges in
$H$ and $
P_Fx=\sum_{n\ge1}P_{F_n}x$.
\end{itemize}
For $m\ge 1$ set $G_m:=\bigcup_{n=1}^m F_n$.
Then $P_F=P_{G_m}+P_{F\setminus G_m}$ and, for all $x\in H$,
$$
\|P_Fx-P_{G_m}x\|^2=\|P_{F\setminus G_m}x\|^2
= \iprod{P_{F\setminus G_m}x}{x}
=P_x(F\setminus G_m).
$$
Now $F\setminus G_m=\bigcup_{n>m}F_n\downarrow\emptyset$, and therefore
$P_x(F\setminus G_m)\downarrow 0$ for all $x\in H$. Thus
$\|P_Fx-P_{G_m}x\|\to0$ as $m\to\infty$, that is, $P_{G_m}x\to P_Fx$ in $H$.
Since $P_{G_m}x=\sum_{n=1}^m P_{F_n}x$ by finite additivity, it follows that the series $\sum_{n\ge1}P_{F_n}x$ converges in $H$ and its sum equals $P_Fx$.

\begin{lemma}[Strong continuity from below and above]\label{lem:PVM-continuity}
Let $P:\calF\to\calL(H)$ be a projection-valued measure.
\begin{enumerate}[label={\rm(\roman*)}, leftmargin=*]
\item If $F_1\subseteq F_2\subseteq\cdots$ and $F=\bigcup_{n\ge1}F_n$, then for all $x\in H$ we have $P_{F_n}x\to P_Fx$.
\item If $F_1\supseteq F_2\supseteq\cdots$ and $F=\bigcap_{n\ge1}F_n$, then for all $x\in H$ we have $P_{F_n}x\to P_Fx$.
\end{enumerate}
\end{lemma}

\begin{proof}
(i) Set $G_n:=F_n\setminus F_{n-1}$ for $n\ge 1$, setting $F_0:= \emptyset$. The sets $G_n$ are pairwise disjoint and
$F=\bigcup_{n\ge1}G_n$.
By strong countable additivity, for all $x\in H$ we have
$$
 \limn P_{F_n}x = \limn \sum_{j=1}^n P_{G_j}x = \sum_{j\ge1} P_{G_j}x = P_Fx.
$$

(ii) Apply (i) to the increasing sets $\complement F_n\uparrow \complement F$ and use that $P_{F_n}=I-P_{\complement F_n}$.
\end{proof}

\begin{example} If $T\in\calL(H)$ is a compact normal operator, Theorem \ref{thm:spect-thm-comp} implies that
the mapping $\{\la\}\mapsto P_\lambda$, where $P_\lambda$ is the orthogonal
 projection onto the eigenspace of $\la$, extends to a projection-valued measure on $\sigma(T)$.
\end{example}

\section{The Bounded Functional Calculus}\label{sec:Borel-FC}

Let $({\Om},\calF)$ be a measurable space. The Banach space of all bounded measurable functions $f:\Om\to \C$, endowed with the supremum norm $\n f\n_\infty = \sup_{\om\in \Om}|f(\om)|,$
is denoted by
$B_{\rm b}(\Om)$\index{$B$@$B_{\rm b}(\Om)$}.

\begin{theorem}[Bounded functional calculus]\label{thm:Borel-FC}\index{functional calculus!bounded}
Let $P:\calF\to \calL(H)$ be a projec\-tion-valued measure. There exists a unique linear mapping $ \Phi: B_{\rm b}(\Om)\to \calL(H)$
with the following properties:
 \begin{enumerate}[label={\rm(\roman*)}, leftmargin=*]
  \item\label{it:Borel-FC2} for all $F\in \calF$ we have $\Phi(\one_F) = P_F$;
  \item\label{it:Borel-FC4} for all $f,g\in B_{\rm b}(\Om)$ we have $\Phi(fg) = \Phi(f)\Phi(g)$;
  \item\label{it:Borel-FC3} for all $f\in B_{\rm b}(\Om)$ we have $\Phi(\ov f) = (\Phi(f))^\star$;
  \item\label{it:Borel-FC5} for all $f\in B_{\rm b}(\Om)$ we have $\n \Phi(f)\n \le \n f\n_\infty$;
  \item\label{it:Borel-FC6} for all $f_n,f\in B_{\rm b}(\Om)$, if $\sup_{n\ge 1}\n f_n\n_\infty<\infty $ and $f_n\to f$ pointwise on $\Om$,
  then  for all $x\in H$ we have $\Phi(f_n)x\to \Phi(f)x$.
 \end{enumerate}
 Moreover, for all $x\in H$ and $f\in B_{\rm b}(\Om)$ we have
 \begin{align}\label{eq:intfP} \iprod{ \Phi(f)x}{x} = \int_{\Om} f\ud P_x
 \end{align}
 and
\begin{align}\label{eq:intfP2} \n \Phi(f)x\n^2 = \int_{\Om} |f|^2\ud P_x.
\end{align}
The operators $\Phi(f)$ are normal, and if $f$ is real-valued (respectively, takes values in $[0,\infty)$) they are selfadjoint (respectively, positive).
\end{theorem}

\begin{proof}
For $F\in \calF$ we set $\Phi(\one_F):= P_F$, which is \ref{it:Borel-FC2}, and extend this definition by linearity
to simple functions $f$. It is routine to verify that $\Phi(f)$ is well defined for such functions and that
 \ref{it:Borel-FC4} and \ref{it:Borel-FC3} hold.

If $f = \sum_{j=1}^k c_j\one_{F_j}$ is a simple function with disjoint supporting sets $F_j\in\calF$\!,
the orthogonality of the vectors $P_{F_j} x$ gives
$$ \n \Phi(f)x\n^2
= \sum_{j=1}^k |c_j|^2 \n P_{F_j}x\n^2 \le \max_{1\le j\le k} |c_j|^2 \n x\n^2 = \n f\n_\infty^2 \n x\n^2\!.$$
It follows that $\n\Phi(f)\n \le \n f\n_\infty.$
For general $f\in B_{\rm b}(\Om)$ we can find a sequence of simple functions $f_n$ converging to $f$
uniformly on $\Om$ and satisfying $\n f_n\n_\infty \le \n f\n_\infty$. From $\n \Phi(f_n)- \Phi(f_m)\n \le \n f_n-f_m\n_\infty$ we infer that the operators $\Phi(f_n)$ form a Cauchy sequence in $\calL(H)$. It follows that the limit
$$\Phi(f):= \limn \Phi(f_n)$$
exists, with convergence in the operator norm, and it is routine to check that the limit is independent of the choice of approximating sequence. Moreover,
$$\n \Phi(f)\n \le \limsup_{n\to\infty} \n f_n\n_\infty \le \n f\n_\infty,$$ which gives \ref{it:Borel-FC5}.
The general case of \ref{it:Borel-FC4} and \ref{it:Borel-FC3} now follows by approximation.

To prove \ref{it:Borel-FC6} we first establish \eqref{eq:intfP} and \eqref{eq:intfP2}.
If $f = \sum_{j=1}^k c_j\one_{F_j}$ is a simple function with disjoint supporting sets $F_j\in\calF$\!, then
\begin{align*}
\iprod{ \Phi(f)x}{x} = \sum_{j=1}^k c_j \iprod{P_{F_j}x}{x} = \sum_{j=1}^k c_j P_x(F_j) = \int_{\Om} f\ud P_x.
\end{align*}
This gives \eqref{eq:intfP} for simple functions. The general case follows by approximation and dominated convergence. Similarly,
\begin{align*}
 \n \Phi(f)x\n^2 & = \sum_{j=1}^k |c_j|^2 \n P_{F_j}x\n^2
 = \sum_{j=1}^k |c_j|^2 \iprod{P_{F_j}x}{x}
 = \int_{\Om} |f|^2\ud P_x.
\end{align*}
This gives \eqref{eq:intfP2} for simple functions.
Again the general case follows by approximation and dominated convergence.
Property \ref{it:Borel-FC6} follows from \eqref{eq:intfP2}, applied to the functions $f_n-f$, and dominated convergence.

To prove normality of $\Phi(f)$, note that \ref{it:Borel-FC4} and \ref{it:Borel-FC3} imply
$$
(\Phi(f))^\star \Phi(f) = \Phi(\ov f)\Phi(f) = \Phi(|f|^2) = \Phi( f)\Phi(\ov f) =
\Phi(f)(\Phi(f))^\star.
$$
Selfadjointness (respectively, positivity) for real-valued (respectively, nonnegative) $f$ is immediate from \ref{it:Borel-FC3} (respectively, \eqref{eq:intfP}).
\smallskip

It remains to prove the uniqueness assertion. Assume that $\Phi_1, \Phi_2: B_{\rm b}(\Om)\to \calL(H)$ are two linear maps satisfying properties \ref{it:Borel-FC2} and the norm bound of Theorem~\ref{thm:Borel-FC}, and
$$\Phi_1(\one_F)=\Phi_2(\one_F),\quad F\in \calF\!.$$
Let $f\in B_{\rm b}(\Om)$ and let $(f_n)_{n\ge 1}$ be a sequence of simple functions converging uniformly to $f$ with $\n f_n\n_\infty\le \n f\n_\infty$ for all $n$. Then $\Phi_1(f_n)=\Phi_2(f_n)$ for all $n$, and property \ref{it:Borel-FC5} implies
$$\n \Phi_i(f)-\Phi_i(f_n)\n \le \n f-f_n\n_\infty,\quad i=1,2.$$
It follows that $\Phi_i(f_n)$ converges in operator norm to $\Phi_i(f)$ as $n\to\infty$, and hence $\Phi_1(f)=\Phi_2(f)$. Since $f\in B_{\rm b}(\Om)$ was arbitrary, we conclude that $\Phi_1=\Phi_2$.
\end{proof}

In what follows we write
$$ \Phi(f) = \int_{\Om} f \ud P = \int_{\Om} f(\la) \ud P(\la)
$$
for functions $f \in B_{\rm b}(\Om)$; the rigorous interpretation of these integrals is through
\eqref{eq:intfP}.

\begin{proposition}[Substitution]\label{prop:projvalmeas-subst} Let $(\Om,\calF)$ and $(\Om'\!,\calF')$ be measurable spaces and let
$f:{\Om}\to \Om'$ be a measurable mapping. If $P:\calF\to \calL(H)$ is a projection-valued measure,
then the {\em image of $P$ under $f$},\index{image!of a projection-valued measure} defined by\index{$FP$@$f_\#(P)$}
$$ (f_\#(P))_{F'} := P_{f^{-1}(F')}, \quad F'\in \calF',$$
is a projection-valued measure.
Denoting by $\Phi$ and $\Psi$ the bounded functional calculi of $P$ and $f_\#(P)$, for all $g\in B_{\rm b}(\Om')$ we have
$$\Phi(g\circ f) = \Psi(g).$$
\end{proposition}

\begin{proof} Let us write $Q:= f_\#(P)$ for brevity.
The elementary verification that $Q$ is a projection-valued measure is left as an exercise. For all $F'\in \calF'$ and $x\in H$,
\begin{align*}
\int_{\Om} \one_{F'}\circ f\ud P_x = \int_{\Om}\one_{f^{-1}(F')}\ud P_x
= \iprod{P_{f^{-1}(F')}x}{x}
=   \iprod{Q_{F'} x}{x}  = \int_{\Om'}\one_{F'}\ud Q_x.
\end{align*}
By linearity and monotone convergence, it follows that for nonnegative functions $g\in B_{\rm b}(\Om')$ and $x\in H$ we have
\begin{align*} \int_{\Om} g\circ f\ud P_x = \int_{\Om'} g\ud Q_x,
\end{align*}
that is,
$\iprod{\Phi(g\circ f)x}{x} = \iprod{\Psi(g)x}{x}$.
For nonnegative $g\in B_{\rm b}(\Om')$ the result now follows from Proposition \ref{prop:polarisation}. For general
$g\in B_{\rm b}(\Om')$ the result follows by splitting into real and imaginary parts and considering their positive and negative parts.
\end{proof}

Due to the absence of a reference measure $\mu$ on $\Om$, we had to work with the Banach space $B_{\rm b}(\Om)$ rather than with
a Lebesgue space $L^\infty(\Om,\mu)$. However, the projection-valued measure $P$ can be used to define a Lebesgue-type space $L^\infty({\Om},P)$ as follows.

\begin{definition}[$P$-Essential boundedness] \label{def:LinftyP} A measurable function
$f:{\Om}\to \C$ is said to be {\em $P$-essentially bounded} if $P_{\{|f| > r\}} = 0$ for some $r\ge 0$.
We define $L^\infty({\Om},P)$\index{$L^\infty({\Om},P)$} to be the space of all equivalence classes of $P$-essentially bounded measurable functions, identifying the functions $f$ and $g$ when $P_{\{f\not=g\}} = 0$.
\end{definition}

With respect to the norm
$$ \n f\n_{L^\infty({\Om},P)}:= \inf\big\{r\ge 0:\, P_{\{|f| > r\}} = 0\big\},$$
the space $L^\infty({\Om},P)$ is easily checked to be a Banach space.

To extend $\Phi$ to $L^\infty(\Om,P)$, let $f$ be $P$-essentially bounded and choose $r\ge0$ such that
$P_{\{|f|>r\}}=0$.
The function $\widetilde f:=\one_{\{|f|\le r\}}f$
is bounded, and we define
$
\Phi(f):=\Phi(\widetilde f).
$
This definition is independent of the choice of the bounded representative. Indeed, if bounded measurable functions $f_1$ and $f_2$ satisfy
$P_{\{f_1\ne f_2\}}=0$,
then, by \eqref{eq:intfP2},
$$
\n(\Phi(f_1)-\Phi(f_2))x\n^2
=
\int_{\Om}|f_1-f_2|^2\ud P_x
=
0,\quad x\in H.
$$
For functions $f\in L^\infty({\Om},P)$ we thus obtain a well-defined bounded operator $\Phi(f)$.
The properties \ref{it:Borel-FC2} and \ref{it:Borel-FC4} hold again, and \ref{it:Borel-FC5} improves to {\em equality}:

\begin{proposition}\label{lem:LinftyP} If $P:\calF\to \calL(H)$ is a projection-valued measure, then for all $f\in L^\infty({\Om},P)$ we have
$$\n \Phi(f)\n = \n f\n_{L^\infty({\Om},P)}.$$
\end{proposition}
\begin{proof}
The upper bound `$\le$' is an immediate consequence of part \ref{it:Borel-FC5} of Theorem \ref{thm:Borel-FC}.

For the lower bound `$\ge$' we distinguish two cases. If
$\n f\n_{L^\infty(\Om,P)}=0$,
then there is nothing to prove. Assume therefore that $\n f\n_{L^\infty(\Om,P)}>0$.
Then the definition of the $P$-essential supremum implies that for all $0<\eps<1$
the projection $P_{F_\eps}$ is nonzero, where
$$
F_\eps := \{|f| > (1-\eps)\n f\n_{L^\infty(\Om,P)}\}.
$$
For all $x\in \ran(P_{F_\eps})$,
\begin{align*}
\n \Phi(f)x\n^2  = \int_{\Om} |f|^2\ud P_x
& \ge (1-\eps)^2\n f\n_{L^\infty({\Om},P)}^2\int_{\Om} \one_{F_\eps}\ud P_x
\\ & = (1-\eps)^2\n f\n_{L^\infty({\Om},P)}^2 \n P_{F_\eps}x\n^2
 =  (1-\eps)^2\n f\n_{L^\infty({\Om},P)}^2\n x\n^2\!.
\end{align*}
This shows that $\n \Phi(f)\n \ge  (1-\eps)\n f\n_{L^\infty({\Om},P)}$. Letting $\eps\to 0$, we obtain the result.
\end{proof}

We now turn to the special case of projection-valued measures defined on the Borel $\sigma$-algebra $\mathscr{B}(M)$ of a compact metric space $M$. In this setting we may define the {\em (topological) support}\index{support!of a projection-valued measure} $\operatorname{supp}(P)$ as the smallest closed subset $F\subseteq M$ with the property that $P_{\,\complement F} = 0$. Such a smallest closed set indeed exists: let $(x_n)_{n\ge 1}$ be a countable dense set in $M$ and $(r_m)_{m\ge 1}$ an enumeration of the positive rationals, and put
$$G:=\bigcup\bigl\{B(x_n;r_m):\, P_{B(x_n;r_m)}=0\bigr\}.$$ Then $P_G=0$, and $\complement G$ is contained
in every closed set $F\subseteq M$ such that $P_{\complement F}=0$.

Here, we are interested in the special case where $M = K$ is a compact subset of $\C$.

\begin{proposition}\label{prop:TP} Let $K\subseteq \C$ be compact and let $P:\mathscr{B}(K)\to \calL(H)$ be a projection-valued measure.
The bounded operator
$$ T:= \int_K \la \ud P(\la)$$
is normal, and its spectrum is given by
$$\si(T) = \operatorname{supp}(P).$$
If $B\in\calB(K)$ is such that $\sigma(T)\cap B= \emptyset$, then $P_B=0$.
\end{proposition}

\begin{proof}
Let $\Phi:B_{\rm b}(K)\to\calL(H)$ be the bounded functional calculus associated with $P$, and define ${\rm id}(\la):=\la$ for $\la\in K$. With this notation, $T=\Phi({\rm id})$.
Normality of $T$ is immediate from the properties of the bounded calculus.

If $\la_0\in\complement K$, the functions $\la\mapsto \la_0-\la$ and $\la\mapsto(\la_0-\la)^{-1}$ are bounded on $K$ and multiply to $\one$. In view of $\Phi(\one)=I$, we have $\la_0-T=\Phi(\la_0-{\rm id})$, and property \ref{it:Borel-FC4} of the bounded calculus shows that $\Phi((\la_0-{\rm id})^{-1})$ is a two-sided inverse of $\la_0-T$. It follows that $\la_0\in\varrho(T)$. Thus $\si(T)\subseteq K$.

Next we prove that $\si(T)\subseteq\operatorname{supp}(P)$. Let $\la_0\in K$, and suppose that there exists a relatively open neighbourhood $U\subseteq K$ of $\la_0$ such that $P_U=0$. If $K\setminus U=\emptyset$, then $I=P_K=P_U=0$ and there is nothing to prove. Otherwise, put $\delta:=d(\la_0,K\setminus U)>0$. Define the bounded Borel function $h$ on $K$ by
$$
h(\lambda):=
\begin{cases}
(\lambda_0-\lambda)^{-1}, & \lambda\in K\setminus U,\\
0, & \lambda\in U.
\end{cases}
$$
Then $(\la_0-{\rm id})h=\one_{K\setminus U}=\one-\one_U$. Applying the bounded functional calculus and using $\Phi(\one_U)=P_U=0$, we obtain $$(\la_0 I-T)\Phi(h)=\Phi(h)(\la_0 I-T)=I - P_U = I.$$ Hence $\la_0\in\varrho(T)$. Therefore, if $\la_0\in\si(T)$, every relatively open neighbourhood $U$ of $\la_0$ satisfies $P_U\neq0$. This proves that $\si(T)\subseteq\operatorname{supp}(P)$.

Conversely, let $\la_0\in\operatorname{supp}(P)$. Suppose, for a contradiction, that $\la_0\in\varrho(T)$. Choose $r>0$ so small that
$r\le 1/(2\n R(\la_0,T)\n)$. Since $\la_0\in\operatorname{supp}(P)$, we have $P_{K\cap B(\la_0;r)}\neq0$. Pick $x\in\Ran(P_{K\cap B(\la_0;r)})$ with $\n x\n=1$. Then $P_x$ is supported on $K\cap B(\la_0;r)$, and therefore
\begin{align*}
\n(\la_0-T)x\n^2
 =\int_K|\la_0-\la|^2\ud P_x(\la)
=\int_{K\cap B(\la_0;r)}|\la_0-\la|^2\ud P_x(\la)
 \leq r^2\!.
\end{align*}
It follows that
$$
1=\n x\n=\n R(\la_0,T)(\la_0-T)x\n\leq \n R(\la_0,T)\n\,\n(\la_0-T)x\n \le \frac12,
$$
a contradiction. Thus $\la_0\notin\varrho(T)$, and hence $\la_0\in\si(T)$. This proves $\operatorname{supp}(P)\subseteq\si(T)$, and therefore $\si(T)=\operatorname{supp}(P)$.

Finally, if $B\in\calB(K)$ satisfies $\si(T)\cap B=\emptyset$, then $B\subseteq K\setminus\operatorname{supp}(P)$. Since $P$ is supported on $\operatorname{supp}(P)$, we obtain $P_{K\setminus \operatorname{supp}(P)}=0$ and hence $P_B=0$.
\end{proof}

It follows from the proposition that $P$ restricts to a projection-valued measure on $\sigma(T)$ in a natural way. Accordingly, we have
$$ T = \int_{\sigma(T)} \la \ud P(\la).$$
The spectral theorem for bounded normal operators, which will be proved in Section \ref{sec:spectral-thm-normal-bdd},
asserts that, conversely, for every normal operator $T\in \calL(H)$ there exists a unique projection-valued measure $P$
on $\sigma(T)$ such that
$$T=  \int_{\sigma(T)} \la \ud P(\la).$$

We have the following uniqueness result:

\begin{proposition}[Uniqueness]\label{prop:spres-uniq} Let $P$ and $\widetilde P$ be projection-valued measures on a compact set $K\subseteq\C$
and define
$$ T=  \int_{K} \la \ud P(\la), \qquad \wt T=  \int_{K} \la \ud \wt P(\la).$$
 If $T = \wt T$, then $P = \wt P$.
\end{proposition}
\begin{proof}
We have $T = \Phi({\rm id})  =\wt\Phi({\rm id})$ and $T^\star = \Phi(\ov{{\rm id}})  =\wt\Phi(\ov{{\rm id}})$, where $\Phi$ and $\wt \Phi$ are the bounded calculi associated with $P$ and $\wt P$, respectively,
and where ${\rm id}(\la) = \la$ as before. By the multiplicativity of the calculi,
$$T^mT^{\star n} = \Phi({\rm id}^m\ov{{\rm id}}^n) = \Phi({\rm id})^m\Phi(\ov{{\rm id}})^n = \wt\Phi({\rm id})^m\wt\Phi(\ov{{\rm id}})^n= \wt\Phi({\rm id}^m\ov{{\rm id}}^n).$$
It follows that $p(T) = \Phi(p) = \wt\Phi(p)$ for all functions $p(z) = q(z,\ov z)$ with $q$ a polynomial in two variables, and then
$f(T) = \Phi(f) = \wt\Phi(f)$ for all $f\in C(K)$ by approximation using the Stone--Weierstrass
theorem (Theorem \ref{thm:Stone-Weierstrass}).
This means that
$$ \int_{K} f\ud P = \int_{K} f\ud \wt P, \quad f\in C(K).$$
By the uniqueness part of the Riesz representation theorem
(Theorem \ref{thm:CK-dual}), the finite Borel measures $P_x$
and $\wt P_x$ agree for every $x\in H$. Thus, for every
$B\in\calB(K)$,
$$
\iprod{P_Bx}{x}
=P_x(B)
=\wt P_x(B)
=\iprod{\wt P_Bx}{x},\quad x\in H.
$$
By polarisation, $P_B=\wt P_B$ for every $B\in\calB(K)$.
\end{proof}

We have seen that if $K\subseteq \C$ is compact and $P:\calB(K)\to\calL(H)$ is a projection-valued measure and $T = \int_{K} \la\ud P(\la)$, then
$\si(T) = \operatorname{supp}(P)$.
The next theorem computes the spectrum of an operator $f(T)$ in the bounded Borel calculus:

\begin{theorem}[Spectral mapping theorem]
\label{thm:SMT-pushforward-support}
Let $K\subseteq\C$ be compact and let $P:\mathcal B(K)\to\calL(H)$ be a projection-valued measure.
Set
$$
T:=\int_K \lambda\ud P(\lambda).
$$
For $f\in B_{\rm b}(K)$, put
$$
f(T):=\int_K f\ud P,
\qquad L:=\ov{f(K)},
$$
and regard $f_\#(P)$ as a projection-valued measure on $L$,
defined by
$$
(f_\#(P))_B:=P_{f^{-1}(B)},\quad B\in\mathcal B(L).
$$
Then
$$
\sigma(f(T))=\operatorname{supp}(f_\#(P)).
$$
Consequently,
$\sigma(f(T))\subseteq\ov{f(\sigma(T))}$.
If $f$ is continuous, then
$\sigma(f(T))=f(\sigma(T))$.
\end{theorem}
%

\begin{proof}

Since $f$ is bounded and Borel, Proposition~\ref{prop:projvalmeas-subst} (applied with $\Om':=L$)
implies that $Q:= f_\#(P)$ is a projection-valued measure on $L$ and that
$$
f(T) = \int_K f\ud P=\int_{L} z\ud Q(z).
$$
Applying Proposition~\ref{prop:TP} to the compact set $L$ and the projection-valued measure $Q$, we obtain
$$
\operatorname{supp}(f_\#(P)) = \operatorname{supp}(Q) = \sigma\Bigl(\int_{L} z\ud Q(z)\Bigr) = \sigma(f(T)).
$$

Let $F:=\ov{f(\si(T))}$. If $B\subseteq L$ is a Borel set with
$B\cap F=\emptyset$, then $f^{-1}(B)\cap \si(T)=\emptyset$, hence $P_{f^{-1}(B)\cap \si(T)}=0$.
Since $P$ is supported on $\si(T)$ by Proposition \ref{prop:TP}, this implies $P_{f^{-1}(B)}=0$. Therefore, $Q_B=P_{f^{-1}(B)}=0$ whenever $B\cap F=\emptyset$.
In particular, $Q_{L\setminus F}=0$. Since $\operatorname{supp}(Q)$ is the smallest closed set $E$ with $Q_{L\setminus E}=0$, this implies
$\operatorname{supp}(Q)\subseteq F=\ov{f(\si(T))}$.

If $f$ is continuous on $K$, then $f(\si(T))$ is compact.
For $\la_0\in\si(T)=\operatorname{supp}(P)$ and any relatively
open subset $V$ of $L$ containing $f(\la_0)$, the set
$f^{-1}(V)$ is a relatively open set in $K$ containing $\la_0$. Hence $Q_V=P_{f^{-1}(V)}\ne0$, so
$f(\la_0)\in\operatorname{supp}(Q)$.
Together with the preceding inclusion, this gives
$\si(f(T))=f(\si(T))$.
\end{proof}

\section{The Spectral Theorem for Bounded Normal Operators}\label{sec:spectral-thm-normal-bdd}

We are now ready to state and prove the spectral theorem for bounded normal operators.

\begin{theorem}[Spectral theorem for bounded normal operators]\label{thm:SMT-bdd}\index{theorem!spectral, for bounded normal operators}
Let $T\in \calL(H)$ be a normal operator. There exists a unique projection-valued measure $P$ on $\si(T)$
\index{projection-valued measure!of a bounded normal operator}
such that
$$ T= \int_{\si(T)} \la\ud P(\la).$$
\end{theorem}

For the proof of Theorem \ref{thm:SMT-bdd} we need the following elementary consequence of the Riesz representation theorem.

\begin{proposition}\label{prop:bilin-T}
Let $\aa:H\times H\to \C$ be a sesquilinear mapping with the property
that there exists a constant $C\ge 0$ such that
$$|\aa(x,y)| \le C \n x\n \n y\n, \quad x,y\in H.$$
Then there exists a unique operator $A\in\calL(H)$
such that $$\aa(x,y) = \iprod{Ax}{y}, \quad x,y\in H.$$
Moreover, $\n A\n \le C$, where $C$ is the boundedness constant of $\aa$.
\end{proposition}

\begin{proof}
For all $y \in H$, the mapping $x \mapsto \aa(x,y)$ is a
bounded functional on $H$, and the Riesz representation theorem gives a unique  $w = w(y) \in H$ satisfying
$$ \aa(x,y) = \iprod{x}{w(y)} ,\quad x \in  H.$$
Let  $B:y\mapsto w=w(y)$ be the resulting mapping. Then,
$$ \aa(x,y) = \iprod{x}{By},\quad x,y \in H.$$
We claim that $B: H \to H$ is a bounded operator. Indeed, if
$c_1,c_2\in\K$  and $y_1,y_2\in H$, for all $x \in H$ we have
\begin{align*} \iprod{x}{B(c_1 y_1 + c_2 y_2)}
 & = \aa(x,c_1 y_1 + c_2 y_2)
  = \overline c_1 \aa(x,y_1) + \overline c_2 \aa(x,y_2)
\\ & = \overline c_1 \iprod{x}{By_1} + \overline c_2 \iprod{x}{By_2}
= \iprod{x}{c_1 By_1 + c_2 By_2}.
\end{align*}
Since this equality holds for all $x \in H$, it follows that $B$ is linear. Furthermore,
$$\n By\n^2  = \iprod{By}{By} =|\iprod{By}{By}| = |\aa(By,y)| \le C\n By\n \n y\n.$$
Consequently, $\n By\n\le C\n  y\n$ for all $y\in H$,  so $B$ is bounded with $\n B\n\le C$.
The operator $A:=B^\star$ has the required properties.
\end{proof}

\begin{proof}[Proof of Theorem \ref{thm:SMT-bdd}]
We begin with existence. For $x,y\in H$ consider the linear mapping $\phi_{x,y}: C(\si(T))\to \C$, $$\phi_{x,y}(f):= \iprod{f(T)x}{y},$$
where $f(T)$ is given by the continuous functional calculus of $T$.
The bound $\n f(T)\n \le \n f\n_\infty$ implies that $\phi_{x,y}$ is bounded and $\n \phi_{x,y}\n  \le \n x\n \n y\n$.
By the Riesz representation theorem (Theorem \ref{thm:CK-dual}) there exists a unique complex Borel measure $P_{x,y}$ on $\si(T)$ such that
\begin{equation}\label{eq:mu-xy}
\iprod{f(T)x}{y} = \int_{\si (T)} f\ud P_{x,y}, \quad f\in C(\si(T)).
\end{equation}
The left-hand side is sesquilinear in $(x,y)$, hence so is the right-hand side.
By the uniqueness part of the Riesz representation theorem, $(x,y)\mapsto P_{x,y}$ is sesquilinear as a map into the complex measures. Moreover,
$$ \n P_{x,y}\n = \sup_{\n f\n_\infty\le 1}|\iprod{f(T)x}{y}| \le \n x\n\n y\n,$$
since $\n f(T)\n \le \n f\n_\infty$.

For $x=y$, $P_x:=P_{x,x}$ is a finite positive Borel measure.
Indeed, for nonnegative $f\in C(\si(T))$ we have $f(T)\ge 0$ by Corollary \ref{cor:normal-cont-fc-pos},
so
$$
\int_{\si(T)}f\ud P_x=\iprod{f(T)x}{x}\ge 0.
$$
Positivity of the measure $P_x$ now follows from the Riesz representation theorem.

By Proposition \ref{prop:bilin-T}, applied to $\aa(x,y):= \int_{\si(T)}f\ud P_{x,y}$, for all $f\in B_{\rm b}(\si(T))$
there exists a unique bounded operator on $H$, which we denote by $f(T)$, such that
\begin{align}\label{eq:spec-int-T}
\iprod{f(T)x}{y} = \int_{\si(T)}f\ud P_{x,y}
\end{align}
for all $x,y\in H$.
The bound on the norm of $P_{x,y}$ implies $$\n f(T)\n \le \n f\n_\infty, \quad f\in B_{\rm b}(\si(T)).$$
For $f\in C(\si(T))$,
\eqref{eq:spec-int-T} is consistent with \eqref{eq:mu-xy}.
Taking $f(\la) = \la$ and $x=y$,
and noting that $P_x:= P_{x,x}$ is a finite Borel measure on $\si(T)$,
we obtain the identity
$$  \iprod{Tx}{x} = \int_{\si(T)}\la\ud P_{x}(\la),$$
except that it remains to be proved that the measures $P_x$ come from a projection-valued measure.
The remainder of the proof is devoted to showing that this is indeed the case.

For all $f,g\in C(\si(T))$ and $x,y\in H$, by the multiplicativity of the continuous functional calculus of $T$ we have
$$\int_{\si(T)} f \ud P_{g(T)x,y} = \iprod{f(T)g(T)x}{y} = \iprod{(fg)(T)x}{y} = \int_{\si(T)} fg\ud P_{x,y}.$$
By the Riesz representation theorem (Theorem \ref{thm:CK-dual}),
this implies that
\begin{align}\label{eq:mainstep-mult-fc}
P_{g(T)x,y} = g P_{x,y}
\end{align}
 as finite Borel measures on $\si(T)$. This, in turn, implies that for all $f\in B_{\rm b}(\si(T))$, $g\in C(\si(T))$, and $x,y\in H$, we have
$$\iprod{f(T)g(T)x}{y}=\int_{\si(T)} f \ud P_{g(T)x,y} =\int_{\si(T)} fg\ud P_{x,y}= \iprod{(fg)(T)x}{y}.$$
Since also
$$ \iprod{f(T)g(T)x}{y}= \iprod{g(T)x}{f(T)^\star y} = \int_{\si(T)} g \ud P_{x,f(T)^\star y},$$ and since $g\in C(\si(T))$ was arbitrary, another application of the
Riesz representation theorem improves \eqref{eq:mainstep-mult-fc} to
$P_{x,f(T)^\star y} = f P_{x,y}$ as finite Borel measures on $\si(T)$.
It follows that for all $g\in B_{\rm b}(\si(T))$ we have
\begin{align}\label{eq:multipl-Bb}
\iprod{f(T)g(T)x}{y}=\int_{\si(T)} g \ud P_{x,f(T)^\star y} = \int_{\si(T)}fg \ud P_{x,y}
= \iprod{(fg)(T)x}{y}.
\end{align} Since $x,y\in H$ were arbitrary, this proves multiplicativity of the Borel calculus.

We can now complete the proof by showing that the measures $P_x$ come from a projection-valued measure.
For Borel sets $B\subseteq \sigma(T)$, define $P_B\in \calL(H)$ by
 $$ P_B  := \one_B(T).$$
This operator is contractive. It is also positive, since for all $x\in H$ we have
$$ \iprod{P_Bx}{x} = \int_{\si(T)} \one_B \ud P_x = P_x(B).$$
By \eqref{eq:multipl-Bb}, $P_B^2 = \one_B(T)\one_B(T) = \one_B(T) = P_B$.
It follows that $P_B$ is a projection, and this projection is orthogonal by selfadjointness.
Moreover, $P_{\si(T)}=\one_{\si(T)}(T)=I$,
and for every $x\in H$ the mapping
$ B\mapsto\iprod{P_Bx}{x}=P_x(B) $
is a finite measure on $\calB(\si(T))$.

Thus
$ P:\calB(\si(T))\to\calL(H)$
is a projection-valued measure. Let $\Phi$ be the bounded Borel calculus associated with $P$.
For all $f\in B_{\rm b}(\si(T))$ and $x\in H$ we have
$$
\iprod{\Phi(f)x}{x}
=\int_{\si(T)}f\ud P_x
=\iprod{f(T)x}{x}.
$$
By polarisation, $\Phi(f)=f(T)$.
By construction,
$$
T=\int_{\si(T)}\lambda\ud P(\lambda).
$$
Uniqueness follows from Proposition \ref{prop:spres-uniq}.
\end{proof}

\begin{example}\label{ex:position-operator}
On $L^2(0,1)$, the {\em position operator}\index{position operator} $X$ is defined by
$$
 X f(x):= xf(x),\quad x\in (0,1),\ f\in L^2(0,1).
$$
We have $\sigma(X) = [0,1]$, and the projection-valued measure of $X$ is given by
$$ P_B f = \one_B f$$ for all $f\in L^2(0,1)$ and Borel subsets $B$ of $[0,1]$.
In particular,
$\one_B(X)=0$ for any Borel null set $B$ of $ [0,1]$. As a consequence, for all $\phi \in L^{\infty}(0,1)$
the operator $\phi(X)$ is well defined as a bounded operator on $L^2(0,1)$. In fact
we have $\phi(X) = M_\phi$, where
$ M_\phi f(x):= \phi(x)f(x).$
\end{example}

The operators $\phi(X)$ of this example arise quite naturally:

\begin{proposition}\label{prop:commuteSX}
An operator $S\in \calL(L^2(0,1))$ commutes with the operator $X$ (that is, $SX = XS$)
if and only if there exists $\phi \in L^{\infty}(0,1)$ such that $S = \phi(X)$.
\end{proposition}
\begin{proof}
Every multiplication operator $\phi(X)$ with
$\phi\in L^\infty(0,1)$ commutes with $X$.

Conversely, suppose that $SX=XS$ and put
$\phi:=S\one\in L^2(0,1)$.
For every polynomial $p$ we have
$$
Sp=Sp(X)\one=p(X)S\one=p\phi.
$$

Let $f\in L^2(0,1)$.
By the Weierstrass approximation theorem and the density of
$C[0,1]$ in $L^2(0,1)$, there exist polynomials $p_n$ such that
$p_n\to f$ in $L^2(0,1)$.
Since $S$ is bounded, we also have $Sp_n\to Sf$ in $L^2(0,1)$.
Passing to a subsequence along which both convergences hold
almost everywhere, we deduce from the identities $Sp_n=p_n\phi$ that
$Sf=f\phi$ almost everywhere.

To prove that $\phi$ is essentially bounded, fix $a>\|S\|$
and put $F_a:=\{|\phi|>a\}$.
Applying the preceding identity to $\one_{F_a}$ gives
$$
a^2\|\one_{F_a}\|_2^2
\le
\|\phi\one_{F_a}\|_2^2
=
\|S\one_{F_a}\|_2^2
\le
\|S\|^2\|\one_{F_a}\|_2^2.
$$
Since $a>\|S\|$, the set $F_a$ has measure zero.
Thus $\phi\in L^\infty(0,1)$ with
$\|\phi\|_\infty\le\|S\|$, and the identity $Sf=f\phi$
for all $f\in L^2(0,1)$ proves that $S=\phi(X)$.
\end{proof}

We are now in a position to prove the assertions made in Section \ref{sec:sa-un-no}.

\begin{corollary}\label{cor:normal-spec-char} Let $T\in\calL(H)$ be a normal operator.
\begin{enumerate}[label={\rm(\arabic*)}, leftmargin=*]
 \item\label{it:normal-spec-char1} $T$ is selfadjoint if and only if $\si(T)\subseteq \R$;
 \item\label{it:normal-spec-char2} $T$ is positive if and only if $\si(T)\subseteq [0,\infty)$;
 \item\label{it:normal-spec-char3} $T$ is unitary if and only if $\si(T)\subseteq \mathbb{T}$;
 \item\label{it:normal-spec-char4} $T$ is an orthogonal projection if and only if $\si(T)\subseteq \{0,1\}$.
 \end{enumerate}
Furthermore, if $T$ is a projection, then it is an orthogonal projection.
\end{corollary}
\begin{proof}
 The `only if' statements have already been proved in Chapter \ref{chap:Hilbert-operators}.
 The `if' statements follow from Proposition \ref{prop:TP} and the properties of the bounded functional calculus of a projection-valued measure.
\end{proof}

As the example of the Volterra operator $V$ shows (see Example \ref{ex:Volterra2} and Problem \ref{prob:Volterra}), normality cannot be omitted in
parts \ref{it:normal-spec-char1}, \ref{it:normal-spec-char2}, and \ref{it:normal-spec-char4}. The operator $I+V$ shows that  normality cannot be omitted in
part \ref{it:normal-spec-char3}.

For normal operators $T\in\calL(H)$ and functions $f\in B_{\rm b}(\si(T))$, the operator $\Phi(f)\in\calL(H)$ defined in terms of the projection-valued measure of $T$ is denoted by $f(T)$,
$$ f(T):= \Phi(f) =  \int_{\si(T)} f\ud P\!.$$
The properties of the bounded calculus for $\Phi$ translate into corresponding properties for the mapping $f\mapsto f(T)$:

\begin{theorem}[Bounded functional calculus for normal operators]\label{thm:bff-FC-normal} Let $T\in \calL(H)$ be normal. Then:
\begin{enumerate}[label={\rm(\roman*)}, leftmargin=*]
  \item\label{it:Borel-T1}  for $f(z) = z^m \ov z^n$ we have $f(T) = T^m T^{\star n}$;
  \item\label{it:Borel-T4} for all $f,g\in B_{\rm b}(\si(T))$ we have $(fg)(T) = f(T)g(T)$;
  \item\label{it:Borel-T3} for all $f\in B_{\rm b}(\si(T))$ we have $\ov f(T) = (f(T))^\star$;
  \item\label{it:Borel-T5} for all $f\in B_{\rm b}(\si(T))$ we have $\n f(T)\n \le \n f\n_\infty$;
  \item\label{it:Borel-T6} for all $f_n,f\in B_{\rm b}(\si(T))$, if $\sup_{n\ge 1}\n f_n\n_\infty < \infty$ and $f_n\to f$ pointwise on $\si(T)$,
  then  for all $x\in H$ we have $f_n(T)x\to f(T)x$.
 \end{enumerate}
 Moreover, for all $x\in H$ and $f\in B_{\rm b}(\si(T))$ we have
 \begin{align*} \iprod{f(T)x}{x} = \int_{\si(T)} f\ud P_x, \qquad
\n f(T)x\n^2 = \int_{\si(T)} |f|^2\ud P_x.
\end{align*}
The operators $f(T)$ are normal, and if $f$ is real-valued (respectively, nonnegative) they are selfadjoint (respectively, positive).
\end{theorem}
\begin{proof}
 Everything but \ref{it:Borel-T1} follows by applying Theorem \ref{thm:Borel-FC} to the projection-valued measure of $T$ obtained in Theorem \ref{thm:SMT-bdd}; \ref{it:Borel-T1} follows from \ref{it:Borel-T4} and \ref{it:Borel-T3} and the fact that for ${\rm id}(z) = z$ we have ${\rm id}(T) = T$.
\end{proof}

The bounded Borel calculus agrees with the continuous functional
calculus for functions $f\in C(\si(T))$. Indeed, part \ref{it:Borel-T1} gives
agreement on polynomials in $z$ and $\ov z$. Such polynomials
are dense in $C(\si(T))$ by Theorem
\ref{thm:Stone-Weierstrass}, and both calculi are continuous
for the supremum norm.

Our next aim is to prove versions of the spectral mapping theorem and the composition rule for the bounded Borel calculus.

\begin{theorem}[Spectral mapping theorem revisited]\label{thm:SMT-normal-bdd}
Let $T\in\calL(H)$ be normal with projection-valued measure $P$ on $\si(T)$.
Then for all $f\in B_{\rm b}(\si(T))$ we have
$$
\si(f(T))=\operatorname{supp}(f_\#(P)).
$$
As a consequence,
$\si(f(T))\subseteq\overline{f(\si(T))}$. If $f$ is continuous, then $\si(f(T))=f(\si(T))$.
\end{theorem}

This is an immediate consequence of Theorem \ref{thm:SMT-pushforward-support}, which can be applied to the bounded Borel calculus of the projection-valued measure associated with $T$.

\begin{theorem}[Composition]\label{thm:comp-c-normal}
Let $T\in\calL(H)$ be normal and let $f\in B_{\rm b}(\si(T))$.
Then for every $g\in B_{\rm b}(\si(f(T)))$ we have
$$
g(f(T)) =(\wt g\circ f)(T)
$$
for any bounded Borel function $\wt g:\C\to\C$ such that $\wt g|_{\si(f(T))}=g$.
\end{theorem}
\begin{proof}
Let $P$ be the projection-valued measure of $T$, put $S:=f(T)$,
and set
$K:=\{z\in\C:\, |z|\le\|f\|_\infty\}$.
Viewing $f$ as a Borel mapping from $\sigma(T)$ into $K$, let
$Q:=f_\#P$. By Proposition \ref{prop:projvalmeas-subst},
$$
S=\int_K z\ud Q(z).
$$
Proposition \ref{prop:TP} gives $\sigma(S)\subseteq K$.
By Proposition \ref{prop:spres-uniq}, $Q$ is therefore the
projection-valued measure of $S$, viewed as a
projection-valued measure on $K$.

Let $g\in B_{\rm b}(\sigma(S))$ and let $\wt g:\C\to\C$ be
a bounded Borel extension of $g$. Another application of
Proposition \ref{prop:projvalmeas-subst} gives
$$
g(S)
=\int_K\wt g(z)\ud Q(z)
=\int_{\sigma(T)}\wt g(f(\lambda))\ud P(\lambda)
=(\wt g\circ f)(T).
$$
This proves the asserted identity and its independence of
the chosen extension.
\end{proof}

The following proposition expresses the null space of $f(T)$ in terms of the projection-valued measure of $T$.

\begin{proposition}\label{prop:ker-fT-zero-set}
Let $T\in \calL(H)$ be normal, and let $P$ be its projection-valued measure. For every bounded Borel function $f:\si(T)\to \C$ we have
$$
\ker(f(T))=\Ran(P_{f^{-1}(\{0\})}).
$$
\end{proposition}

\begin{proof}
Put $Z:=f^{-1}(\{0\})$.
If $x\in\Ran(P_Z)$, then $x=P_Zx$, and hence
$$
        f(T)x=f(T)P_Zx=(f\one_Z)(T)x=0.
$$
Thus $\Ran(P_Z)\subseteq\ker(f(T))$.
Conversely, let $x\in\ker(f(T))$. Then
$$
        0=\|f(T)x\|^2
        =
        \int_{\sigma(T)} |f(\lambda)|^2\ud P_x(\lambda).
$$
It follows that $P_x(\sigma(T)\setminus Z)=0$, or equivalently,
$$
        \|(I-P_Z)x\|^2
        =
        \iprod{P_{\sigma(T)\setminus Z}x}{x}
        =
        P_x(\sigma(T)\setminus Z)=0.
$$
Hence $x=P_Zx$, so $x\in\Ran(P_Z)$.
\end{proof}

In particular, this proposition implies that if $T\in \calL(H)$ is normal,
$f:\si(T)\to\C$ is bounded Borel, and $f(T)$ is injective, then
$P_{f^{-1}(\{0\})}=0$.

It is of interest to revisit the case of a compact normal operator $T$. The spectrum of $T$ is a finite or infinite sequence $(\la_n)_{n\ge 1}$ with $0$ as its only possible limit point. By Theorem \ref{thm:sa-isolated-pt}, for any nonzero $\la\in\sigma(T)$ the orthogonal projection $P_{\la}$
onto the corresponding eigenspace equals the spectral projection $P^{\{\la\}}$ of Theorem \ref{thm:spec-proj}.

\begin{proposition} Let $T\in \calL(H)$ be a compact normal operator and let $P$ be its projection-valued measure.
Then for all nonzero $\la\in\si(T)$, $$P_{\{\la\}} = P^{\{\la\}} = P_\la.$$
\end{proposition}

\begin{proof}
If we put $\wt P_{\{\la\}}:= P_{\la}$ and extend this definition by defining $\wt P_{\{0\}}$ to be the orthogonal projection onto $\ker (T)$ if $0\in \si(T)$, the spectral theorem for compact normal operators implies that $\wt P$ defines a projection-valued measure and
\begin{align*}\iprod{Tx}{x}
=  \sum_{\la\in \si(T)\setminus\{0\}} \la \iprod{P_{\la} x}{x}
=  \int_{\si(T)\setminus\{0\}}\la \ud \wt P_{x}(\la)
= \int_{\si(T)}\la \ud \wt P_{x}(\la).
   \end{align*}
Hence, by the uniqueness theorem for projection-valued measures, $\wt P = P$.
\end{proof}

Further results along this line are given in Problem \ref{prob:threeP} and Proposition \ref{prop:normal-eigen}.

We conclude this section with two famous results due to von Neumann.

\begin{theorem}[von Neumann]\label{thm:vonNeumann}\index{theorem!von Neumann, on contractions}\index{contraction!von Neumann theorem}
 If $T\in\calL(H)$ is a contraction, then for all polynomials $p$ in one variable we have
 $$ \n p(T)\n \le \sup_{|z|=1} |p(z)|.$$
\end{theorem}
\begin{proof}
First suppose that $U$ is a unitary operator. Then $\sigma(U)$ is contained in the unit circle.  Since unitaries are normal, by Theorem \ref{thm:SMT-bdd} we have $U =\int_{\sigma(U)}\la \ud P(\la)$, where
$P$ is the projection-valued measure of $U$.
By Theorem \ref{thm:normal-cont-fc} and the fact that $\sigma(U)\subseteq \mathbb{T}$,
 \begin{align}\label{eq:vN-unitary}
 \n p(U)\n
 = \sup_{z\in \sigma(U)} |p(z)| \le \sup_{|z|=1} |p(z)|.
 \end{align}
Next let $T$ be a contraction. By the Sz.-Nagy dilation theorem (Theorem \ref{thm:Nagy}), $T$ has a unitary dilation $U$, so that $T^n = J^\star U^n J$ for some isometric operator $J$ and all $n\in \N$.
Then $p(T) = J^\star p(U) J$, so $\n p(T)\n \le \n p(U)\n$ and the result follows from \eqref{eq:vN-unitary}.
\end{proof}

As an application of the spectral theorem we prove von Neumann's theorem on pairs of commuting normal operators. Two projection-valued measures
$P:\calF\to\calL(H)$ and $P':\calF'\to\calL(H)$ are said to {\em commute}
if $P_F$ and $P_{F'}'$ commute for all $F\in \calF$ and $F'\in \calF'$.

\begin{lemma}\label{lem:prod-spec-meas} Let $P^1,\dots,P^k$ be commuting projection-valued measures on compact metric spaces $K_1,\dots,K_k$, respectively. There exists a unique projection-valued measure $P$ on $K=K_1\times\cdots\times K_k$ such that $$ P_{B_1\times\cdots\times B_k} = P^1_{B_1}\cdots P^k_{B_k} $$ for all Borel sets $B_j\subseteq K_j$, $j=1,\dots,k$. \end{lemma}

\begin{proof}
Let $V$ be the vector space of all finite linear combinations of indicator functions of Borel rectangles in $K=K_1\times\cdots\times K_k$.
For a rectangle put
$$
\Phi(\one_{B_1\times\cdots\times B_k})
:=
P^1_{B_1}\cdots P^k_{B_k},
$$
and extend by linearity.
This is well defined and contractive. Indeed, for any $f\in V$ one may choose finite Borel partitions $K_j=\bigcup_{r=1}^{m_j}A_{j,r}$
such that $f$ is constant on every product
$A_{1,r_1}\times\cdots\times A_{k,r_k}$.
The corresponding projections $P^1_{A_{1,r_1}}\cdots P^k_{A_{k,r_k}}$
are pairwise orthogonal. Hence, if $c_{r_1,\dots,r_k}$ is the value of $f$ on the corresponding product,
$$
\Phi(f)
=
\sum_{r_1,\dots,r_k}
c_{r_1,\dots,r_k}
P^1_{A_{1,r_1}}\cdots P^k_{A_{k,r_k}},
$$
and therefore $\n\Phi(f)\n\le\n f\n_\infty$.
The same representation shows that
$$
\Phi(fg)=\Phi(f)\Phi(g),
\qquad
\Phi(\ov f)=\Phi(f)^\star,
\qquad
\Phi(\one)=I.
$$

By uniform continuity, the uniform closure of $V$ in $B_{\rm b}(K)$
contains $C(K)$. Extending $\Phi$ to this closure and restricting it
to $C(K)$ gives a contractive mapping
$\Phi:C(K)\to\calL(H)$ with the same properties.
Repeating the construction in the proof of Theorem \ref{thm:SMT-bdd}, with this $\Phi$ in place of the continuous functional calculus of a normal operator, gives a bounded Borel calculus
$\widehat\Phi:B_{\rm b}(K)\to\calL(H)$
extending $\Phi$. Setting
$$
P_E:=\widehat\Phi(\one_E),\quad E\in\calB(K),
$$
we obtain a projection-valued measure on $K$. To prove that it takes the stated values on Borel rectangles, let
$\pi_j:K_1\times\cdots\times K_k\to K_j$ be the coordinate projection.
We first claim that
$$
P_{\pi_j^{-1}(B)}=P^j_B, \quad B\in\calB(K_j),\ j=1,\ldots,k.
$$
For $\varphi\in C(K_j)$, the construction of $\Phi$ gives
$$
\widehat\Phi(\varphi\circ\pi_j) = \Phi(\varphi\circ\pi_j)
= \int_{K_j}\varphi\ud P^j .
$$
Hence the two projection-valued measures $B\mapsto P_{\pi_j^{-1}(B)}$ and $B\mapsto P^j_B$
on $K_j$ have the same integrals against all continuous functions. Arguing
as in Proposition \ref{prop:spres-uniq}, we conclude that they agree on all Borel
subsets of $K_j$.

Now let $B_j\in\calB(K_j)$, $j=1,\ldots,k$. Since $
\one_{B_1\times\cdots\times B_k} = \prod_{j=1}^k \bigl(\one_{B_j}\circ\pi_j\bigr)$, the multiplicativity of the bounded Borel calculus gives
$$
P_{B_1\times\cdots\times B_k}
= \widehat\Phi(\one_{B_1\times\cdots\times B_k})
= \prod_{j=1}^k \widehat\Phi(\one_{B_j}\circ\pi_j)
= P^1_{B_1}\circ\cdots\circ P^k_{B_k}.
$$
Thus $P$ has the required values on Borel rectangles.

Finally, suppose that $Q$ is another projection-valued measure on $K$ with
the same values on Borel rectangles. For fixed $x\in H$, the finite measures
$E\mapsto \iprod{P_Ex}{x}$ and $E\mapsto \iprod{Q_Ex}{x}$
agree on all Borel rectangles $B_1\times\cdots\times B_k$.
This collection contains $K$, is closed under finite intersections, and
generates $\calB(K)$. By Dynkin's lemma (Lemma \ref{lem:unique-app}), the two
finite measures agree on $\calB(K)$. Hence
$\iprod{P_Ex}{x}=\iprod{Q_Ex}{x}$ for all $x\in H,\ E\in\calB(K)$,
and by polarisation $P_E=Q_E$ for all $E\in\calB(K)$. This proves
uniqueness.
\end{proof}

\begin{lemma}[Commutation of projection-valued measures]\label{lem:PVMs-commute}
\index{projection-valued measure!commutation}
Let $N_1,N_2\in\calL(H)$ be commuting normal operators, and let $P^1$ and $P^2$ be their associated projection-valued measures. Then for all $A\in\calB(\si(N_1))$ and $B\in \calB(\si(N_2))$ we have
$$
P^1_A\,P^2_B=P^2_B\, P^1_A.
$$
\end{lemma}

\begin{proof}
We prove the lemma in two steps.

\smallskip{\em Step 1} --
First let $N\in\calL(H)$ be a normal operator, with projection-valued measure $P$, and let
$T\in\calL(H)$ commute with $N$. Then, by the Fuglede--Putnam--Rosenblum theorem (Theorem \ref{thm:FPR}), $T$ commutes with $N^\star$. It follows that $T$ commutes with each polynomial $p(N,N^\star)$. By the Stone--Weierstrass theorem, this implies that $T$ commutes with $f(N)$ for all $f\in C(\sigma(N))$.

For $x,y\in H$, let $P_{x,y}$ denote the finite complex Borel
measure given by $P_{x,y}(B):=\iprod{P_Bx}{y}$.
Since $T$ commutes with $f(N)$ for every $f\in C(\sigma(N))$,
we have
$$
\int_{\sigma(N)}f\ud P_{Tx,y}
=\iprod{f(N)Tx}{y}
=\iprod{f(N)x}{T^\star y}
=\int_{\sigma(N)}f\ud P_{x,T^\star y}.
$$
By the uniqueness part of the Riesz representation theorem
(Theorem \ref{thm:CK-dual}), $P_{Tx,y}=P_{x,T^\star y}$.
Consequently, for every $B\in\calB(\sigma(N))$,
$$
\iprod{P_BTx}{y}
=\iprod{P_Bx}{T^\star y}
=\iprod{TP_Bx}{y}.
$$
Since $x,y\in H$ were arbitrary, this proves that $P_BT=TP_B$.

\smallskip
{\em Step 2} --
Now let $N_1,N_2\in \calL(H)$ be commuting normal operators, and fix a set $A\in \calB(\si(N_1))$. By Step 1 we have
$ N_2P^1_A=P^1_AN_2$. Since $P^1_A$ commutes with $N_2$, we may apply Step 1 again and conclude that $P^1_AP^2_B =P^2_BP^1_A$ for all
$B\in\calB(\si(N_2))$.
\end{proof}

\begin{theorem}[von Neumann]\label{thm:Neumann-commute}
Let $N_1,N_2\in\calL(H)$ be commuting normal
operators on $H$. Then there exist a bounded selfadjoint operator $S\in\calL(H)$ and bounded continuous functions $\psi_1,\psi_2:\R\to\C$ such that
$$
N_j=\psi_j(S),\quad j=1,2.
$$
\end{theorem}

\begin{proof} The case $H=\{0\}$ is trivial, so we may assume that $H$ is nonzero.

Since $N_1$ and $N_2$ are commuting normal operators, their projection-valued measures commute by Lemma \ref{lem:PVMs-commute}. Hence, by Lemma \ref{lem:prod-spec-meas},
there exists a joint projection-valued measure $P$ on $\si(N_1)\times \si(N_2)$ such that
$$
N_1=\int_{\si(N_1)\times \si(N_2)} z_1\ud P(z_1,z_2),\qquad
N_2=\int_{\si(N_1)\times \si(N_2)} z_2\ud P(z_1,z_2).
$$
Put $K:=\sigma(N_1)\times\sigma(N_2)$.
By Theorem \ref{thm:Hausdorff-Alexandroff}, there exist a compact
set $L\subseteq[0,1]$ and a continuous surjection $\pi:L\to K$.
By Proposition \ref{prop:min-section}, there exists a Borel
mapping $\sigma:K\to L$ such that $\pi\circ\sigma=\mathrm{id}_K$.
Since $\sigma$ is bounded and real-valued, the operator
$$
S:=\int_K\sigma(z)\ud P(z)
$$
is bounded and selfadjoint.

For $j=1,2$, let $p_j:K\to\C$ be the coordinate projection.
The function $p_j\circ\pi$ is continuous on $L$.
By the Tietze extension theorem, applied to its real and
imaginary parts, it extends to a bounded continuous function
$\psi_j:\R\to\C$. For every $z\in K$ we have
$$
\psi_j(\sigma(z))
=p_j(\pi(\sigma(z)))
=p_j(z)=z_j.
$$

Let $Q:=\sigma_\#P$, viewed as a projection-valued measure
on $[0,1]$. By Proposition \ref{prop:projvalmeas-subst},
$$
S=\int_{[0,1]}t\ud Q(t).
$$
By Proposition \ref{prop:spres-uniq}, $Q$ is the
projection-valued measure of $S$, viewed on $[0,1]$.
Applying Proposition \ref{prop:projvalmeas-subst} once more gives
\[
\psi_j(S)
=\int_{[0,1]}\psi_j(t)\ud Q(t)
=\int_K\psi_j(\sigma(z))\ud P(z)
=\int_K z_j\ud P(z)
=N_j,\quad j=1,2.
\qedhere\]
\end{proof}

\section{The von Neumann Bicommutant Theorem}\label{sec:bicommutant}

In this section we prove a result of fundamental importance in the theory of operator algebras, von Neumann's celebrated bicommutant theorem.
We also show that the bicommutant of a single normal operator on a separable
Hilbert space is precisely its bounded functional calculus.

We begin by introducing the relevant terminology.

\begin{definition}[Commutant]
The {\em commutant}\index{commutant} of a subset $\mathscr{T}\subseteq \calL(H)$ is the set
$$ \mathscr{T}' := \{S\in \calL(H):\, ST = TS \ \ \textrm{for all} \ \ T\in \mathscr{T}\}.$$
The  {\em bicommutant}\index{bicommutant} of a subset $\mathscr{T}\subseteq \calL(H)$ is the set $\mathscr{T}'':= (\mathscr{T}')'$. Higher commutants are defined inductively.
\end{definition}

It is an immediate consequence of this definition that for any subset $\mathscr{T}\subseteq \calL(H)$ we have
$ \mathscr{T}' = \mathscr{T}'''.$

\begin{definition}[Strong and weak operator topologies]\label{def:weakoperatortop}
The {\em strong operator topology on $\calL(H)$}\index{strong!operator topology} is the smallest topology $\tau$ on $\calL(H)$ with the property that the linear mappings $T\mapsto Tx$ are continuous for all $x\in H$. The {\em weak operator topology on $\calL(H)$}\index{weak!operator topology} is the smallest topology $\tau$ on $\calL(H)$ with the property that the linear mappings $T\mapsto \iprod{Tx}{y}$ are continuous for all $x,y\in H$.
\end{definition}

Definition \ref{def:weakoperatortop} has natural counterparts for $\calL(X)$ with $X$ a Banach space, but these will not be needed.

As in Section \ref{sec:weak-star} for the weak and weak$^*$ topologies,
the strong operator topology is generated by the sets of the form
$$\{T\in \calL(H):\, \n (T-T_0)x\n < \eps\}$$
with $\eps>0$, $x\in H$, and $T_0\in\calL(H)$,
and likewise the weak operator topology is generated by the sets of the form
$$\{T\in \calL(H):\, |\iprod{(T-T_0)x}{y}| < \eps\}$$
with $\eps>0$, $x,y\in H$, and $T_0\in\calL(H)$.

For every set $\mathscr{T}\subseteq \calL(H)$, the commutant $\mathscr{T}'$ is closed in the weak operator topology. To see this,
suppose that $T_0\not\in \mathscr{T}'$\!. Then there exist an operator $S\in \mathscr{T}$, vectors $x,y\in H$, and a number $\delta>0$ such that $|\iprod{T_0Sx}{y} -  \iprod{ST_0x}{y}| = \delta$.
The set $$U:= \{T\in \calL(H): \, |\iprod{(T_0-T)Sx}{y}|<\delta/2, \, |\iprod{(T_0-T)x}{S^\star y}|<\delta/2\}$$
is open in the weak operator topology, contains $T_0$, and every $T\in U$ fails to commute with $S$. It follows that
$U\cap \mathscr{T}'=\emptyset$.

Recall that a {\em subalgebra} of $\calL(H)$ is a subspace of $\calL(H)$ which is closed under composition. A  {\em $\star$-subalgebra} of $\calL(H)$ is a subalgebra of $\calL(H)$ closed under taking Hilbert space adjoints. A subalgebra is said to be {\em unital} if it contains the identity operator.

\begin{theorem}[von Neumann bicommutant theorem]\label{thm:bicommutant}\index{theorem!von Neumann, bicommutant}\index{theorem!bicommutant} For a unital $\star$-subalgebra $\mathscr{A}$ of $\calL(H)$ the following assertions are equivalent:
 \begin{enumerate}[label={\rm(\arabic*)}, leftmargin=*]
  \item\label{it:bicommutant1} $\mathscr{A} = \mathscr{A}''$;
  \item\label{it:bicommutant3} $\mathscr{A}$ is closed in the weak operator topology;
  \item\label{it:bicommutant2} $\mathscr{A}$ is closed in the strong operator topology.
 \end{enumerate}
Moreover, $\mathscr A$ is strongly dense in $\mathscr A''$.
\end{theorem}
A $\star$-subalgebra $\mathscr A$ of $\calL(H)$ which is closed in
operator norm is called a {\em $C^\star$-algebra}.\index{C@$C^\star$-algebra} This is not the commonly used definition (the standard definition is mentioned in the Notes to Chapter \ref{ch:CompactOperators}), but one of the main theorems in the theory of $C^\star$-algebras asserts that this definition is equivalent to the standard one. A unital $\star$-subalgebra $\mathscr{A}$ of $\calL(H)$ satisfying the equivalent conditions
of Theorem \ref{thm:bicommutant} is called a {\em von Neumann algebra}.\index{von Neumann algebra}

\begin{proof}
The implication \ref{it:bicommutant1}$\Rightarrow$\ref{it:bicommutant3}
follows from the easy observation that commutants are closed in the weak operator topology, and the implication  \ref{it:bicommutant3}$\Rightarrow$\ref{it:bicommutant2} is trivial.

\smallskip
\ref{it:bicommutant2}$\Rightarrow$\ref{it:bicommutant1}: \ We proceed in two steps.

\smallskip{\em Step 1} -- Fix $x_0\in H$ and let $P$ denote the orthogonal projection of $H$ onto the closure $Y$ of the subspace $\{Tx_0:\, T\in \mathscr{A}\}$.
Since $I\in \mathscr{A}$, we have $Px_0 = x_0$.

We claim that $Y$ is invariant under all $T\in \mathscr{A}$: for if $T\in \mathscr{A}$ and $y\in Y$, say $y = \limn T_n x_0$ with $T_n\in \mathscr{A}$ for all $n\ge 1$, then $Ty = \limn  T T_n x_0$ with $TT_n\in \mathscr{A}$ for all $n\ge 1$, and therefore $Ty\in Y$.
Similarly, $Y^\perp$ is invariant under all $T\in \mathscr{A}$: for if $T\in \mathscr{A}$ and $x\in Y^\perp$\!,
then for all $y\in Y$ we have  $\iprod{Tx}{y} = \iprod{x}{T^\star y} =0$ since $T^\star \in \mathscr{A}$ and therefore $T^\star y\in Y$.

By the claim, for all $T\in \mathscr{A}$ and $x\in H$ we have $TPx\in Y$ and $T(I-P)x \in Y^\perp$\!, and therefore
$$ TPx = PTPx = PT (Px +(I-P)x) = PTx, \quad x\in H.$$
We conclude that $TP = PT$ for all $T\in \mathscr{A}$, that is, $P\in \mathscr{A}'$\!.

Let $T_0\in \mathscr{A}''$ be fixed.  Then $PT_0 = T_0P$ since $P\in \mathscr{A}'$\!, and this implies that $T_0x_0 = T_0Px_0 = PT_0x_0 \in Y$. By the definition of $Y$ this means that for all $\eps>0$ there exists an element $T\in \mathscr{A}$ such that $$\n T_0x_0 - Tx_0\n < \eps.$$

{\em Step 2} --
To show that every strongly open set
containing the operator $T_0\in  \mathscr{A}''$ intersects $\mathscr{A}$, it suffices to show that, for any choice of $x_1,\dots,x_k\in H$ and $\eps>0$, there exists $T\in \mathscr{A}$ such that
\begin{align}\label{eq:toprove-bicomm} \n (T_0-T)x_j\n < \eps, \quad j=1,\dots,k.
\end{align}
In what follows we set ${\bf x_0}:= (x_1,\dots,x_k)$.

For $S\in \calL(H)$ let $\rho(S)\in \calL(H^k)$ be given by $$\rho(S)(h_1,\dots,h_k):= (Sh_1,\dots,Sh_k).$$
We claim that $$\rho(T_0)\in (\rho(\mathscr{A}))''\!.$$
Indeed, suppose that ${\bf S} = (S_{ij})_{i,j=1}^k \in(\rho(\mathscr{A}))'$\!. This means that
$ {\bf S} \rho(T){\bf y} = \rho(T){\bf S} {\bf y}$ for all $T\in \mathscr{A}$ and ${\bf y} = (y_1,\dots,y_k)\in H^k$\!, that is,
$$ \sum_{j=1}^k S_{ij} Ty_j =  \sum_{j=1}^k T S_{ij} y_j, \quad i=1,\dots,k, \ y_1,\dots,y_k\in H,$$
which implies that for all $1\le i,j\le k$ we have $S_{ij}\in \{T\}'$ for all $T\in \mathscr{A}$, so $S_{ij}\in \mathscr{A}'$\!. But this clearly implies that
$\rho(T_0)$ commutes with ${\bf S}$.

We now apply Step 1, with $H$, $\mathscr{A}$, and $T_0$ replaced by $H^k$\!, $\rho(\mathscr{A})$, and $\rho(T_0)$, respectively. This gives an operator $T\in \mathscr{A}$ such that $\n (\rho(T_0) - \rho(T)){\bf x_0}\n<\eps$, that is,
$$\sum_{j=1}^k \n (T_0-T)x_j\n^2 < \eps^2\!.$$
In particular, \eqref{eq:toprove-bicomm} follows from this.

We have shown that every strongly open set containing an element
of $\mathscr A''$ intersects $\mathscr A$. The argument up to
this point did not require $\mathscr A$ to be strongly closed,
and therefore proves the additional density assertion.
If $\mathscr A$ is strongly closed, it follows that
$\mathscr A=\mathscr A''$, as required.
\end{proof}

The following result identifies the von Neumann algebra generated by a
projection-valued measure with the range of its bounded functional calculus.

\begin{theorem}\label{thm:vN-PVM}
Let $H$ be separable, let $(\Om,\calF)$ be a measurable space, and let
$P:\calF\to\calL(H)$ be a projection-valued measure. If $\mathscr P:=\{P_F:\, F\in\calF\}$,
then
$$
\mathscr P''
=
\Bigl\{\int_\Om f\ud P:\, f\in L^\infty(\Om,P)\Bigr\}.
$$
\end{theorem}

For the proof of this theorem we need the following lemma. We need the first part here; the second part is needed in Chapter \ref{chap:QM}.

\begin{lemma}[Commutative von Neumann algebras]
\label{lem:comm-vN}
Let $H$ be separable and let $\mathscr M\subseteq\calL(H)$ be a
commutative von Neumann algebra. Then:
\begin{enumerate}[label={\rm(\arabic*)}, leftmargin=*]
\item there exists a separating vector $h\in H$, that is,
$Ah=0$ implies $A=0$ for all $A\in\mathscr M$;
\item there exists $S\in\mathscr M$ with $0\le S\le I$ such that
$\mathscr M=\{S\}''$.
\end{enumerate}
\end{lemma}

\begin{proof}
Choose a dense sequence $(x_n)_{n\ge1}$ in $H$. Starting with
$G_0=\{0\}$, define inductively
$$
y_n:=(I-Q_{n-1})x_n,\qquad
H_n:=\ov{\mathscr M y_n},\qquad G_n:=G_{n-1}\oplus H_n,
$$
where $Q_{n-1}$ is the orthogonal projection onto $G_{n-1}$.
Inductively, $G_{n-1}$ is $\mathscr M$-invariant; since $\mathscr M$
is closed under adjoints, so is $G_{n-1}^\perp$.
Thus $H_n\subseteq G_{n-1}^\perp$ and
$H_n$ is $\mathscr M$-invariant.
Since $I\in\mathscr M$, we have $y_n\in H_n$, and therefore
$x_n=Q_{n-1}x_n+y_n\in G_n$. Hence $H=\bigoplus_{n\ge1}H_n$.
The series
$$
h:=\sum_{n\ge1}2^{-n}\frac{y_n}{1+\n y_n\n}
$$
converges in $H$. If $A\in\mathscr M$ and $Ah=0$, orthogonality
implies $Ay_n=0$ for every $n$. For $B\in\mathscr M$ we then have
$ABy_n=BAy_n=0$, so $A$ vanishes on every $H_n$ and hence on $H$.

The subspace $\ov{\mathscr M'h}$ is invariant under $\mathscr M'$
and its adjoints, so its orthogonal projection $Q$ belongs to
$\mathscr M''=\mathscr M$. Since $(I-Q)h=0$, we have $Q=I$ and $\ov{\mathscr M'h} = H$.
Choose projections $E_n\in\mathscr M$ so that
$\{E_nh:\,n\ge1\}$ is dense in
$\{Eh:\,E\in\mathscr M$ is an orthogonal projection $\}$;
this is possible since $H$ is separable.
For any orthogonal projection $E\in\mathscr M$, choose indices $n_k$
with $E_{n_k}h\to Eh$. For $B\in\mathscr M'$,
$$
(E_{n_k}-E)Bh=B(E_{n_k}-E)h\to 0.
$$
For $x\in H$ and $\eps>0$, choose $B\in\mathscr M'$ with
$\n x-Bh\n<\eps$. Since $\n E_{n_k}-E\n\le2$, we have
$$\n(E_{n_k}-E)x\n\le2\eps+\n B\n\,\n(E_{n_k}-E)h\n.$$
Letting $k\to\infty$ and then $\eps\downarrow0$ proves that
$E_{n_k}\to E$ strongly. Thus the strongly closed algebra
$\{E_n:\,n\ge1\}''$ contains every orthogonal projection in $\mathscr M$.

Put $\mathscr A:=\{E_n:\,n\ge1\}''$.
For selfadjoint $A\in\mathscr M$, its spectral projections belong to
$\mathscr M$: by the first step in the proof of Lemma
\ref{lem:PVMs-commute}, they commute with $\mathscr M'$ and hence
belong to $\mathscr M''=\mathscr M$. Hence, by the preceding paragraph, these spectral projections belong to $\mathscr A$.
Approximating the identity function uniformly on $\sigma(A)$ by
step functions shows that $A$ is a norm limit of finite linear
combinations of these projections. Since $\mathscr A$ is strongly
closed, it is also norm closed, and therefore $A\in\mathscr A$.
Every element of $\mathscr M$ is of the form $A+iB$ with
$A,B\in\mathscr M$ selfadjoint, so $\mathscr M\subseteq\mathscr A$.
The reverse inclusion follows because $\mathscr M$ is a von Neumann
algebra containing all $E_n$. We conclude that $\mathscr M = \mathscr A$.

Put $S:=\sum_{n\ge1}3^{-n}E_n$, with convergence in operator norm.
Then $S\in\mathscr M$ and $0\le S\le\tfrac12 I$.
Inductively, suppose that $E_1,\dots,E_{n-1}\in\{S\}''$, and put
$R_n:=S-\sum_{k<n}3^{-k}E_k=\sum_{k\ge n}3^{-k}E_k$.
By the induction hypothesis, $R_n\in\{S\}''$.
Since $E_n$ commutes with $R_n$, its range and kernel are invariant
under $R_n$. On the range $\Ran(E_n)$, the first term of this series equals
$3^{-n}I$ and all remaining terms are positive, so $R_n\ge3^{-n}I$.
On $\Ker(E_n)$, the first term vanishes, and therefore
$$0\le R_n\le\sum_{k>n}3^{-k}I=\tfrac12\,3^{-n}I.$$
The spectral theorem therefore shows that
$\one_{[3^{-n},\infty)}(R_n)$ acts as the identity on $\Ran(E_n)$
and as zero on $\Ker(E_n)$. Consequently,
$$
E_n=\one_{[3^{-n},\infty)}(R_n)\in\{S\}'',
$$
by the preceding observation applied to $\mathscr N=\{S\}''$.
Thus $\mathscr M\subseteq\{S\}''$; the reverse inclusion follows
from $S\in\mathscr M$.
\end{proof}

\begin{proof}[Proof of Theorem \ref{thm:vN-PVM}]
Let $\Phi:L^\infty(\Om,P)\to\calL(H)$ denote the bounded functional
calculus associated with $P$, and let $\mathscr M:=\mathscr P''$ be the von Neumann algebra generated by the range of $P$. Since the projections in $\mathscr P$ commute, we have
$\mathscr M=\mathscr P''\subseteq\mathscr P'$. If $S,T\in\mathscr M$,
then $S\in\mathscr P'$ and $T\in\mathscr P''=(\mathscr P')'$; hence
$ST=TS$. Thus $\mathscr M$ is commutative.
By Lemma \ref{lem:comm-vN}, there exists a vector $h\in H$ that is separating for $\mathscr M$.

We next prove that $\Ran(\Phi)=\mathscr M$. The inclusion
$\Ran(\Phi)\subseteq\mathscr M$ is immediate for simple functions.
For general $f\in L^\infty(\Om,P)$ it follows by uniform approximation
of a bounded representative of $f$ by simple functions.

For the converse, fix $S\in\mathscr M$.
The linear span $\mathscr A$ of $\mathscr P$ is a unital
$\star$-subalgebra, since $P_\Om=I$, $P_F^\star=P_F$, and
$P_FP_G=P_{F\cap G}$. Moreover,
$\mathscr A''=\mathscr P''=\mathscr M$.
The density assertion in Theorem \ref{thm:bicommutant}
therefore shows that $\mathscr A$ is strongly dense in
$\mathscr M$. For each $n\ge1$, strong density gives
$A_n\in\operatorname{span}\mathscr P$ with
$\n(A_n-S)h\n<1/n$. Each $A_n$ has the form $\Phi(f_n)$
for a simple function $f_n$, so $\Phi(f_n)h\to Sh$.
By \eqref{eq:intfP2},
$$
\n\Phi(f_n-f_m)h\n^2
=
\int_\Om |f_n-f_m|^2\ud P_h,
$$
so $(f_n)_{n\ge1}$ is a Cauchy sequence in $L^2(\Om,P_h)$. Let
$f\in L^2(\Om,P_h)$ be its limit.
We claim that $f$ is $P$-essentially bounded and
$\n f\n_{L^\infty(\Om,P)}\le\n S\n$.

Fix $F\in\calF$. Since $S$ and $\Phi(f_n)$ commute with $P_F$,
$$
\Phi(f_n)P_Fh=P_F\Phi(f_n)h\to P_FSh=SP_Fh.
$$
Consequently,
\begin{align*}
\int_F|f|^2\ud P_h =\limn\int_F|f_n|^2\ud P_h
& =\limn\n\Phi(f_n)P_Fh\n^2  \\
&=\n SP_Fh\n^2
 \le\n S\n^2\n P_Fh\n^2
 =\n S\n^2P_h(F).
\end{align*}
Since $F\in\calF$ was arbitrary, it follows that $|f|\le\n S\n$ $P_h$-almost everywhere.

Moreover, if $P_h(F)=0$, then $P_Fh=0$. Since $P_F\in\mathscr M$
and $h$ is separating for $\mathscr M$, this implies $P_F=0$.
Thus the $P_h$-null sets are precisely the $P$-null sets, and the
claim follows.

We may therefore regard $f$ as an element of $L^\infty(\Om,P)$.
Using \eqref{eq:intfP2} once more,
$$
\Phi(f_n)h\to\Phi(f)h.
$$
Hence $Sh=\Phi(f)h$. Both $S$ and $\Phi(f)$ belong to $\mathscr M$,
and $h$ is separating for $\mathscr M$. Therefore $S=\Phi(f)$.
This proves the reverse inclusion.
\end{proof}

As an immediate application we obtain:

\begin{theorem}[von Neumann, bicommutant of a normal operator]
\label{thm:bicomm-singleT}
\index{theorem!bicommutant of a normal operator}
Let $H$ be separable and let $T\in\calL(H)$ be normal. Then
$$
\{T\}''
=
\bigl\{f(T):\, f\in B_{\rm b}(\si(T))\bigr\}.
$$
\end{theorem}

\begin{proof}
Let $P$ be the projection-valued measure of $T$, and let
$\mathscr M_P$ be the von Neumann algebra generated by its range.

If $S\in\{T\}'$, then, by the first step in the proof of Lemma
\ref{lem:PVMs-commute}, $S$ commutes with every $P_F$. Hence
$P_F\in\{T\}''$ for every $F\in\calB(\si(T))$, and therefore
$\mathscr M_P\subseteq\{T\}''$.

Conversely, $T\in\mathscr M_P$. Since $T$ is normal, the
Fuglede--Putnam--Rosenblum theorem gives $\{T\}'=\{T,T^\star\}'$.
Thus $\{T\}''$ is the von Neumann algebra generated by $T$, and hence
$\{T\}''\subseteq\mathscr M_P$.

It follows that $\{T\}''=\mathscr M_P$. The result now follows from
Theorem \ref{thm:vN-PVM}.
\end{proof}

\section{Application to Orthogonal Polynomials}\label{sec:orthpol}\index{orthogonal!polynomials}

In this section we present an interesting application of the spectral theorem to orthogonal polynomials.
Let $\mu$ be a Borel measure on the real line satisfying the condition
\begin{align}\label{eq:mu-finite-moments} \int_{-\infty}^\infty |x|^n \ud \mu(x)< \infty, \quad n\in\N.
\end{align}
We further assume that the support of $\mu$ is not a finite set.
Suppose that $(p_n)_{n\in\N }$ is a sequence of polynomials with real coefficients satisfying the following two assumptions:
\begin{enumerate}[label={\rm(\roman*)}, leftmargin=*]
 \item\label{it:orthpol1} for all $n\in\N $ we have ${\rm deg}(p_n) = n$;
 \item\label{it:orthpol2} for all $m,n\in\N $ with $m\not=n$ we have the orthogonality relation
\begin{align}\label{eq:OP} \int_{-\infty}^\infty p_m(x)p_n(x) \ud \mu(x)= 0.
\end{align}
\end{enumerate}
For $n=0$, \ref{it:orthpol1} is understood to mean that $p_0\not\equiv 0$.
By linearity, \ref{it:orthpol1} and \ref{it:orthpol2} imply
$$ \int_{-\infty}^\infty x^mp_n(x) \ud \mu(x)= 0, \quad m<n.$$

\begin{proposition}\label{prop:orthpol}\index{three-term recurrence}
Let $\mu$ be a Borel measure on the real line with the properties stated above. For any sequence of polynomials $(p_n)_{n\in\N }$ with real coefficients satisfying the conditions \ref{it:orthpol1} and \ref{it:orthpol2},
there exist real numbers $A_n,B_n,C_n$ $(n\in\N )$ satisfying $C_0=0$ and $A_{n-1}C_{n}>0$  $(n\ge 1)$ such that,
 with $p_{-1}\equiv 0$,
 $$ xp_{n} =  A_np_{n+1}+B_np_{n}+ C_np_{n-1}, \quad n\in\N.
 $$
\end{proposition}
\begin{proof}
 Since $xp_{n}$ is a polynomial of degree $n+1$, it is of the form
$ xp_{n} = \sum_{j=0}^{n+1} c_{j,n} p_j$
with all coefficients $c_{j,n}$ real.
Setting $c_{j,n} := 0$ for $j\ge n+2$, \eqref{eq:OP} implies
$$\int_{-\infty}^\infty xp_n(x) p_{m}(x) \ud \mu(x)= c_{m,n} N_m, \quad \hbox{where} \ N_m =  \int_{-\infty}^\infty p_m^2(x) \ud \mu(x).$$
Since the support of $\mu$ is not finite, we have $N_m\not=0$ for all $m\in\N$.
For $n\in\N $ the polynomial $p_{n}(x)$ is orthogonal to $xp_m(x)$ for all $m= 0,1,\dots,n-2$.
This forces $c_{m,n} =0$ for all $m =0,1,\dots,n-2$. This, in turn,
implies $$  xp_{n} = c_{n+1,n}p_{n+1}+ c_{n,n}p_{n}+ c_{n-1,n} p_{n-1}, \quad n\in\N .$$
This gives the three-term recurrence relation with  $A_n = c_{n+1,n}$, $B_n = c_{n,n}$, and $C_n = c_{n-1,n}$, with the convention $C_0 = c_{-1,0}=0$.
Since the degree of $ xp_{n}$ is $n+1$, we have $A_n\not=0$. Also, for $n\ge 1$,
\begin{equation}\label{eq:orthpol-AC}
\begin{aligned} A_{n-1}N_{n} = c_{n,n-1}N_{n} & = \int_{-\infty}^\infty  xp_{n}(x)p_{n-1}(x) \ud \mu(x)
\\ & = \int_{-\infty}^\infty  xp_{n-1}(x) p_{n}(x)\ud \mu(x) = c_{n-1,n}N_{n-1} = C_n N_{n-1},
\end{aligned}
\end{equation}
and therefore $A_{n-1}C_n>0$ for $n\ge 1$.
\end{proof}

The polynomial $p_n$ has norm one in $L^2(\R,\mu)$ if and only if $N_n =1$.  Hence if the $p_n$ are orthonormal, \eqref{eq:orthpol-AC} gives $0\not = A_{n-1} = C_n$ for all $n\ge 1$. As an application of the spectral theorem we show that, conversely, for every sequence of polynomials satisfying the three-term recurrence relation subject to the conditions $0\not = A_{n-1} = C_n$ for all $n\ge 1$, and satisfying the additional boundedness assumption $$\sup_{n\in\N }\max\{|A_n|, |B_n|, |C_n|\} < \infty,$$ there exists a finite Borel measure $\mu$ on the real line with respect to which the polynomials are orthonormal.

\begin{theorem}[Three-term recurrence]\label{thm:orthpol}\index{theorem!three-term recurrence} For every sequence $(p_n)_{n\in\N }$ of nonzero polynomials,
with $p_0$ a nonzero constant, satisfying the three-term recurrence relation $$ xp_{n} =  A_np_{n+1}+B_np_{n}+ C_np_{n-1}, \quad n\in\N,$$ with $p_{-1}\equiv 0$, subject to the conditions
$A_n, B_n, C_n \in\R$ and
$0\not = A_{n} = C_{n+1}$ for all $n\in \N$, and $$\sup_{n\in\N }\max\{|A_n|, |B_n|, |C_n|\} < \infty,$$ there exists a finite Borel measure $\mu$ on the real line which satisfies \eqref{eq:mu-finite-moments} and such that the sequence $(p_n)_{n\in\N }$ is an orthonormal basis for $L^2(\R,\mu)$.
\end{theorem}
\begin{proof}
Write $p_0\equiv c$, where $c\neq0$, and let $(e_n)_{n\in\N}$
be the standard orthonormal basis of $\ell^2(\N)$.
Putting $e_{-1}:=0$, define
$$
Te_n:=A_ne_{n+1}+B_ne_n+C_ne_{n-1},\quad n\in\N.
$$
The boundedness of the coefficients implies that $T$ extends to a
bounded operator on $\ell^2(\N)$. Since the coefficients are real
and $A_n=C_{n+1}$, from
\begin{align*} \iprod{Te_n}{e_m} & = C_{n+1}\delta_{m,n+1} + B_{n}\delta_{m,n} + C_{n}\delta_{m,n-1}
\\ & = C_{m}\delta_{m-1,n} + B_{m}\delta_{m,n} + C_{m+1}\delta_{m+1,n}
=  \iprod{e_n}{T e_m}
\end{align*}
we see that $T$ is selfadjoint.
Let $P$ be its projection-valued measure and put
$$
\mu:=|c|^{-2}P_{e_0},
$$
viewed as a Borel measure on $\R$. This measure is finite and
supported on the compact set $\sigma(T)$, so it satisfies
\eqref{eq:mu-finite-moments}.

We claim that
$$
p_n(T)e_0=ce_n,\quad n\in\N.
$$
For $n=0$ this follows from $p_0\equiv c$.
If the identity holds at indices $n$ and $n-1$, the recurrence gives
\begin{align*}
A_np_{n+1}(T)e_0
&=(T-B_nI)p_n(T)e_0-C_np_{n-1}(T)e_0\\
&=c(Te_n-B_ne_n-C_ne_{n-1})
=cA_ne_{n+1}.
\end{align*}
Since $A_n\neq0$, the claim follows by induction, with the terms
at index $-1$ understood to be zero.

By the functional calculus of $T$, we now obtain
$$
\int_\R p_m(x)\overline{p_n(x)}\ud\mu(x)
=
|c|^{-2}\iprod{p_m(T)e_0}{p_n(T)e_0}
=
\delta_{mn}.
$$
Thus the polynomials are orthonormal in $L^2(\R,\mu)$.
The recurrence and the fact that $A_n\neq0$ also imply
$\deg(p_n)=n$, so their linear span is the space of all polynomials.
By the Stone--Weierstrass theorem, this space is dense in
$C(\sigma(T))$, which is dense in $L^2(\sigma(T),\mu)$ by
Remark \ref{rem:CKdenseLpK}.
Since $\mu$ is supported on $\sigma(T)$, the sequence
$(p_n)_{n\in\N}$ is therefore an orthonormal basis for
$L^2(\R,\mu)$.
\end{proof}

\begin{example}
 We have already encountered two examples of orthogonal polynomials.

\begin{enumerate}[label={\rm(\roman*)}, leftmargin=*]
\item The Hermite polynomials $H_n$, $n\in\N$, have been introduced in Section \ref{subsec:Hermite}.
They are orthogonal with respect to the Gaussian measure $\frac1{\sqrt{2\pi}}\exp(-\frac12 x^2)\ud x$ on the real line and satisfy the three-term recurrence relation $H_0(x) = 1$, $H_1(x) = x$, and
$$ H_{n+2}(x) = x H_{n+1}(x)-(n+1)H_n(x), \quad n\in\N.$$
\item
The monic Laguerre polynomials $L_n$, $n\in\N$, have been introduced in Problem \ref{prob:Laguerre} as a scaled version of the Laguerre polynomials.

They are orthogonal with respect to the measure $\one_{\R_+}(x)\exp(-x)\ud x$ on the real line and satisfy the three-term recurrence relation $L_0(x) = 1$, $L_1(x) = x-1$, and
$$L_{n+2}(x) = (x-(2n+3)) L_{n+1}(x) - (n+1)^2L_{n}(x), \quad n\in\N.$$
\end{enumerate}
\end{example}

\section{The Ideal of Compact Operators Revisited}\label{sec:compact-revisited}

In Chapter \ref{ch:CompactOperators} we proved that the compact operators on a Banach space form a closed two-sided ideal in the algebra of all bounded operators. In the Hilbert space setting one can say considerably more. The aim of this section is to show that, on a separable infinite-dimensional Hilbert space $H$, the ideal of compact operators is maximal (in the sense explained by the theorem below) among proper two-sided ideals of $\calL(H)$.

The proof is based on the spectral theorem for bounded selfadjoint operators. More precisely, we will show that if $\calI$ is a two-sided ideal in $\calL(H)$
containing a noncompact operator, then $\calI$ contains the identity operator and therefore $\calI = \calL(H)$.

\begin{theorem}\label{thm:proper-ideal-in-KH}
Let $H$ be a separable Hilbert space and let $\calI$ be a two-sided ideal in
$\calL(H)$. If $\calI\neq \calL(H)$, then
$$
\calI\subseteq \mathscr K(H).
$$
\end{theorem}

\begin{proof}
Arguing by contradiction, suppose that $\calI$ is not contained in
$\mathscr K(H)$. Then there exists $T\in \calI$ which is not compact. Since $\calI$ is a two-sided
ideal, we have $S:=T^\star T\in \calI$.
Moreover, $S$ is not compact. Indeed, if it were, then by Lemma \ref{lem:absTcompact} $S^{1/2} = |T|$ would be compact, and the polar decomposition $T=U|T|$ would imply that $T$ is compact.

Let $P$ be the projection-valued measure of the positive operator $S$. We claim that
there exists $\la>0$ such that the closed subspace
$Y:=\Ran(P_{[\la,\infty)})$
is infinite-dimensional.
Indeed, if $\Ran(P_{[\la,\infty)})$ were finite-dimensional for every
$\la>0$, then for each $n\ge 1$ the projection $P_n:=P_{[1/n,\infty)}$
would have finite rank. Since
$$
S-P_nS = (I-P_n)S = P_{[0,1/n)}S = f_n(S),
$$
where $f_n(\la) := \la\one_{[0,1/n)}(\la)$, the spectral theorem gives
$$
\n S-P_n S\n = \n f_n(S)\n \le \n f_n\n_\infty \le \frac1n.
$$
Hence $S$ would be the norm limit of finite-rank operators, and therefore
compact, a contradiction.

Fix $\la>0$ such that
$Y:=\Ran(P_{[\la,\infty)})$ is infinite-dimensional.
Since $H$ is separable, there exists an isometry
$U\in\calL(H)$ with $\Ran(U)=Y$.
Put $B:=U^\star SU$. For every $x\in H$, the spectral theorem gives
$$
\iprod{Bx}{x}
=\iprod{SUx}{Ux}
\ge\la\|Ux\|^2
=\la\|x\|^2.
$$
Thus $B\ge\la I$, so $B$ is invertible.
Since $S\in\calI$ and $\calI$ is a two-sided ideal,
we have $B\in\calI$ and hence $I=B^{-1}B\in\calI$.
It follows that $\calI=\calL(H)$, a contradiction.
\end{proof}

\begin{problems}

\item
Let $T\in \calL(H)$ be a compact normal operator with spectral decomposition $$T=\sum_{n\ge 1} \lambda_n P_n,$$
where $(\la_n)_{n\ge 1}$ is the (finite or infinite) sequence of distinct eigenvalues of $T$.
Prove that if $f:\sigma(T)\to \C$ is a bounded function, then for all $x\in H$ we have
$$ f(T)x = \sum_{n\ge 1} f(\lambda_n) P_n x$$
with convergence in $H$.
When does the sum $\sum_{n\ge 1} f(\lambda_n) P_n$ converge to $f(T)$ in $\calL(H)$?

\item\label{prob:spectrum-compact-normal}
Let $T\in\calL(H)$ be a compact normal operator.
Give a direct proof (that is, without invoking Theorem \ref{thm:comp-spec}) of the following two statements:
\begin{enumerate}[\rm(a), leftmargin=*]
\item If $\la$ is a nonzero element of $\sigma(T)$, then $\la$ is an eigenvalue.

\noindent{\em Hint:}\ Use the fact that $\ran(\la-T)$ is closed (Lemma \ref{lem:kerran}) to establish the equivalences $\ker(\la-T)=\{0\} \Leftrightarrow \ker(\ov\la-T^\star)=\{0\} \Leftrightarrow  \Ran(\la-T) =H$.
\item If $T$ has infinitely many distinct eigenvalues $\la_n$, then $\limn \la_n = 0$.

\noindent{\em Hint:}\ Choose eigenvectors $Th_n = \la_nh_n$, and define $H_0:=\{0\}$ and
$H_n:={\rm span}\{h_1,\dots,h_n\}$, $n\ge 1$.
Show that $H_{n-1}\subsetneq H_n$ and choose norm-one vectors $x_n \in H_n \cap H_{n-1}^\perp$.
Show that if $n>m$, then $ \n Tx_m - Tx_n\n \ge |\la_n|.$
\end{enumerate}

\item\label{prob:sa-StoneFormula}
Let $T\in\calL(H)$ be a selfadjoint operator with projection-valued measure $P$. Prove the following results.

\begin{enumerate}[\rm(a), leftmargin=*]
 \item If $t\in \R$, then for all $x\in H$ we have
 $$ \lim_{\eps\to 0} i\eps R(t+i\eps,T)x = P_{\{t\}}x.$$
 \item If $a,b\in\R$ with $a<b$, then for all $x\in H$ we have {\em Stone's formula}\index{formula!Stone}\index{Stone's formula}
 $$\lim_{\eps\downarrow 0} \frac{1}{2\pi i} \int_a^b R(t-i\eps,T)x - R(t+i\eps,T)  x \ud t = \frac12(P_{[a,b]}x + P_{(a,b)}x).$$
 {\em Hint:}\ Show first that
 $$ \lim_{\eps\downarrow 0} \frac{1}{2\pi i} \int_a^b \frac{1}{\la-i\eps -t} - \frac{1}{\la+i\eps-t} \ud t =
 \begin{cases}
  0, &  \la\not\in [a,b],\\
  \frac12, & \la\in \{a,b\},\\
  1, & \la\in (a,b).
 \end{cases}
$$
\end{enumerate}

\item
Show that for every positive operator $T\in\calL(H)$ one has $T \le \|T\|I$.

\noindent{\em Hint:}\ If $f$ is bounded and real-valued, then $f \leq \|f\|_\infty $ pointwise.

\item
Let $S,T\in\calL(H)$ satisfy $0\le S\le T$.
\begin{enumerate}[\rm(a), leftmargin=*]
\item Show that $\n S\n \le \n T\n$.

\noindent{\em Hint:} \ Use the result of the preceding problem.

\item
Show that if $S$ and $T$ are invertible, then $0\le T^{-1} \le S^{-1}$.

\noindent{\em Hint:} \ Observe that $T^{-1/2}ST^{-1/2} \le I$ by Problem \ref{prob:S1leqS2}. Deduce from this that $S^{1/2}T^{-1}S^{1/2}\le I$ and use  Problem \ref{prob:S1leqS2} once more.
\end{enumerate}

\item
Let $T\in\calL(H)$ be a normal operator with projection-valued measure $P$. Let $B$ be a Borel subset of $\sigma(T)$.
\begin{enumerate}[\rm(a), leftmargin=*]
  \item Show that $T$ leaves the range of $P_B$ invariant.
  \item Show that $\sigma(T|_{\Ran(P_B)})\subseteq \ov B$.
\end{enumerate}

\item
Prove that if $H \not=\{0\}$ and $T\in \calL(H)$ is normal, then
$$ \n T\n = \sup\{|\iprod{Tx}{x}|:\, \n x\n=1\}.$$

\noindent{\em Hint:}\ Fix $\la_0\in\si(T)$ with $|\la_0| = \n T\n$ (why does such $\la_0$ exist?) and $\eps>0$,
and consider the projection $P_{B(\la_0;\eps)}$, where $P$ is the projection-valued measure of $T$.
Show that if $x\in\ran(P_{B(\la_0;\eps)})$ is a norm-one vector, then
$|\iprod{Tx}{x}| \ge \n T\n - \eps$.

\item\label{prob:threeP}
Let $T\in \calL(H)$ be a normal operator with projection-valued measure $P$.
\begin{enumerate}[\rm(a), leftmargin=*]
\item Show that if $\la\in\si(T)$ is an eigenvalue, then $P_{\{\la\}}$ equals the orthogonal projection $P_\la$ onto the eigenspace $E_\la$.

\noindent{\em Hint:}\ Use Proposition \ref{prop:ker-fT-zero-set}.
\item Conclude that if $\la\in\si(T)$ is an isolated point, then $\la$ is an eigenvalue and $P_{\{\la\}} = P_\la$ equals the spectral projection associated with $\{\la\}$.
\end{enumerate}

\item
Show that the space
$L^\infty(\Om,P)$ introduced in Definition \ref{def:LinftyP} is a Banach space.

\item\label{prob:findPVM}
Find the projection-valued measures of the following unitary operators:
\begin{enumerate}[\rm(a), leftmargin=*]
  \item\label{it:findPVM1} the right shift on $\ell^2(\Z)$;
  \item\label{it:findPVM2} translation by a fixed $t\in\R$ on $L^2(\R)$;
  \item\label{it:findPVM3} rotation by a fixed $\theta\in\R$ on $L^2(\mathbb{T})$;
  \item\label{it:findPVM4} the Fourier transform on $L^2(\R^d)$.
\end{enumerate}

\item\label{prob:sa-diff-2pos}
Let $T\in \calL(H)$ be selfadjoint.
\begin{enumerate}[\rm(a), leftmargin=*]
  \item Show that there exist positive operators $T^+$ and $T^-$ such that $T = T^+-T^-.$

  \noindent{\em Hint:}\ Consider the functions $f^+(t):= t^+$ and $f^{-}(t):= t^-$.
  \item Show that these operators are unique if we also ask that $T^+T^- = T^-T^+=0$.
\end{enumerate}

\item\label{prob:UT}
Prove that if $U\in \calL(H)$ is unitary, there exists a selfadjoint operator $T$ such that $U = e^{iT}$ and $\sigma(T)\subseteq [-\pi,\pi]$. Is this operator $T$ unique?

\noindent{\em Hint:}\ Write $U = \int_{\sigma(U)} \la \ud P(\la)$
as in Theorem \ref{thm:SMT-bdd}.
Find a projection-valued measure $Q$ on $[-\pi,\pi]$ whose image under $t\mapsto e^{it}$ is $P$.

\item
This problem outlines another proof of the spectral theorem for normal operators.
\begin{enumerate}[\rm(a), leftmargin=*]
 \item Explain how the proof of the spectral theorem for normal operators simplifies for selfadjoint operators.
 \item Deduce the spectral theorem for normal operators $T$ from the selfadjoint case by considering the selfadjoint operators $\frac12(T+T^\star)$ and $\frac1{2i}(T-T^\star)$, and
 applying Lemma \ref{lem:prod-spec-meas} to their projection-valued measures.
\end{enumerate}

\item
Show that if $S,T\in \calL(H)$ are contractions and $\frac12(S+T) = I$, then $S = T = I$.
Deduce from this that $I$ is an extreme point of the closed unit ball of $\calL(H)$.

\noindent{\em Hint:}\ First consider the case where $S$ and $T$ are selfadjoint.
For the general case, observe that $\frac12(\frac12(S+S^\star)+ \frac12(T+T^\star))=I$.

\item\index{functional calculus!joint}
Let $H$ be separable and let $T_1,T_2\in\calL(H)$ be commuting
selfadjoint operators, with projection-valued measures $P^{(1)}$ and
$P^{(2)}$. Let $P$ be their product as in Lemma
\ref{lem:prod-spec-meas}. For $f\in B_{\rm b}(\R^2)$ put
$$
f(T_1,T_2):=\int_{\R^2}f(s,t)\ud P(s,t).
$$
Show that
$$
\{T_1,T_2\}''
=
\{f(T_1,T_2):\,f\in B_{\rm b}(\R^2)\}.
$$

\noindent{\em Hint:}\
Let $\mathscr M_P$ be the von Neumann algebra generated by the range of
$P$. Identify $\mathscr M_P$ with the bounded
functional calculus of $P$, and prove that
$\mathscr M_P=\{T_1,T_2\}''$.

\item\label{prob:cyclic-spectral}
Let $T\in\calL(H)$ be normal with projection-valued measure $P$, and let
$h\in H$ be a norm-one vector such that
$\{f(T)h:\, f\in C(\si(T))\}$
is dense in $H$. Such a vector is called {\em cyclic}\index{cyclic vector} for $T$.
Put $\mu:=P_h$.

\begin{enumerate}[\rm(a), leftmargin=*]
\item Show that $U: f(T)h \mapsto f$
defines an isometry $U$ from $\{f(T)h:\, f\in C(\si(T))\}$ into
$L^2(\si(T),\mu)$.

\item Using the density of $C(\si(T))$ in $L^2(\si(T),\mu)$, show that
$U$ extends uniquely to a unitary operator
$U:H\to L^2(\si(T),\mu).$

\item Show that $UTU^\star=M_z$, where $(M_zf)(z):=z f(z).$
\end{enumerate}

\noindent{\em Hint:}\ For part (a), use \eqref{eq:intfP2}.

\item
Find $\mathscr{K}(H)'$ and $\mathscr{K}(H)''$, where $\mathscr{K}(H)$ is the space of compact operators on $H$.

\item
Let $H$ be separable and let $(h_n)_{n\ge1}$ be an orthonormal basis.
If $P_n$ denotes the orthogonal projection onto
${\rm span}\{h_1,\ldots,h_n\}$, show that for every $T\in\calL(H)$ we have
$P_nTP_n\to T$ in the strong operator topology. Deduce that $\mathscr K(H)$ is strongly
dense in $\calL(H)$, although it is closed in the operator norm.

\item
Let $\mathscr A\subseteq\calL(H)$ be a unital $\star$-subalgebra and
let $h\in H$. We call $h$ {\em cyclic} for $\mathscr A$ if
$\ov{\{Ah:\,A\in\mathscr A\}}=H$,
and {\em separating} for $\mathscr A$ if
$Ah=0$ with $A\in\mathscr A$ implies $A=0$.
Show that an element $h\in H$ is cyclic for $\mathscr A$ if and only if $h$ is
separating for $\mathscr A'$.

\noindent{\em Hint:}\
In one direction, consider the orthogonal projection $P$ onto
$\ov{\{Ah:\,A\in\mathscr A\}}$ and show that $P\in\mathscr A'$.

\item\label{prob:comp-ONS-null-converse}
This problem establishes the following converse to Problem \ref{prob:comp-ONS-null}: if $H$ is a separable infinite-dimensional Hilbert space and $T\in \calL(H)$ satisfies $\limn \n Th_n\n = 0$ for every orthonormal basis
$(h_n)_{n\ge 0}$ of $H$, then $T$ is compact.

Consider the collection $\calI$ of all operators that have the stated property.
\begin{enumerate}[\rm(a), leftmargin=*]
\item Show that $\calI$ is a closed subspace of $\calL(H)$.
\item Show that if $T\in \calI$, then
$ST\in \calI$ for all $S\in \calL(H)$, that is, $\calI$ is a left ideal.
\item Show that if $T\in \calI$ and $U\in \calL(H)$ is unitary, then $TU\in \calI$\!.
\item Show that every contraction is a linear combination of four unitaries.

\noindent{\em Hint:}\ First show that every selfadjoint
contraction $A$ is the average of two unitaries. Indeed, use the continuous
functional calculus to define
$$
        U:=A+i(I-A^2)^{1/2}.
$$
Then $U$ is unitary and
$A=\frac12(U+U^\star).$
Apply this to the real and imaginary parts of an arbitrary contraction.

\item Conclude that $\calI$ is a two-sided closed ideal.
\item Observe that $I\notin\calI$, and hence $\calI\neq\calL(H)$. Now use Theorem
\ref{thm:proper-ideal-in-KH} to conclude that every $T\in\calI$ is compact.
\end{enumerate}

\end{problems}

%% file: ch10-UnboundedOperators.tex
\chapter{The Spectral Theorem for Unbounded Normal Operators}\label{ch:unbdd}

\blfootnote{\rm \small \begin{quote}This book has been published by Cambridge University Press in the series `Cambridge Studies in Advanced Mathematics'. The present version is free to view and download for personal use only. Not for re-distribution, re-sale or use in derivative works.\newline {\copyright} Jan van Neerven\end{quote}}

\noindent
Up to this point we have been dealing exclusively with bounded operators. In order to make the functional-analytic apparatus applicable to the study of partial differential equations, we need to accommodate differential operators within the theory. This leads to the notion of an unbounded operator as a linear operator defined only on a suitable subspace, the domain of the operator. Of special interest are unbounded selfadjoint and normal operators, and the main goal of this chapter is to extend the spectral theorems of the preceding chapter to these classes of operators.

\section{Unbounded Operators}\label{sec:unbounded}

Throughout this chapter, $X$ and $Y$ are Banach spaces.

\begin{definition}[Linear operators] A {\em linear operator}\index{linear operator} from $X$ to $Y$ is a pair $(A, \Dom(A))$, where $\Dom(A)$ is a subspace of $X$ and $A: \Dom(A)\to Y$ is a linear operator. The subspace $\Dom(A)$ is called the {\em domain}\index{domain!of a linear operator} of $A$. A linear operator is
{\em densely defined}\index{densely defined} when $\Dom(A)$ is a dense subspace of $X$.
\end{definition}

When no confusion is likely to arise, we omit the domain from the notation and write $A$ instead of $(A,\Dom(A))$. If we wish to emphasise that a linear operator is only defined on its domain, when $X=Y$ we will speak about operators defined {\em in} $X$ (as opposed to operators defined {\em on} $X$).

It is of course perfectly possible that $\Dom(A)=X$, so in particular every bounded operator $A:X\to Y$ is a linear operator in the sense of the above definition. More generally, it may happen that there exists a constant $C\ge 0$
such that $\n Ax\n \le C\n x\n $ for all $x\in \Dom(A)$. In this situation, $A$ admits a unique extension to a bounded operator (of norm at most $C$) defined on the closure of $\Dom(A)$. The interest in the above definition
arises from the fact that many interesting examples of {\em unbounded}\index{unbounded operator} linear operators exist, that is, linear operators for which such a constant $C$ does not exist. Typical examples, treated in more detail below, include differential operators and multiplication operators with unbounded multipliers.

The terms `linear operators' and `unbounded operators' are often used interchangeably. With the latter terminology, however, it becomes somewhat ambiguous whether bounded operators are to be considered as special cases of unbounded operators. To avoid such trivial issues, we generally prefer the terminology `linear operator', which is neutral in this respect.

\subsection{Closed Operators}\label{subsec:closed-operators}

A (globally defined) linear operator from $X$ to $Y$ is bounded if and only if its graph is closed in $X\times Y$, the `if' part being the content of the closed graph theorem. This motivates the following definition.

\begin{definition}[Closed operators]
 A linear operator $A$ from $X$ to $Y$ with domain $\Dom(A)$ is called {\em closed}\index{closed!linear operator} when its {\em graph}\index{graph!of a linear operator} $$\Gr(A):= \{(x,Ax):\, x\in \Dom(A)\}$$
  is closed in $X\times Y$.
\end{definition}

Every bounded operator from $X$ to $Y$ is closed, and
by the closed graph theorem a closed operator with domain $\Dom(A) = X$ is bounded.

If $A$ is a linear operator with domain $\Dom(A)$, then $A$ is bounded (in fact, contractive) as an operator defined on the normed space $\Dom(A)$ endowed with the {\em graph norm}\index{graph!norm}
$$\n x\n_{\Dom(A)} := \n x\n + \n Ax\n, \quad x\in \Dom(A).$$
This follows from the trivial inequality
$$ \n Ax \n \le \n x\n + \n Ax\n  =\n x\n_{\Dom(A)}.$$
The following proposition gives a necessary and sufficient
condition for closedness in terms of the graph norm.

\begin{proposition}\label{prop:closed-op}
 A linear operator is closed if and only if its domain
 is a Banach space with respect to its graph norm.
\end{proposition}

\begin{proof}
Endowing $X\times Y$ with the product norm
$\|(x,y)\|:=\|x\|+\|y\|$, the mapping $x\mapsto (x,Ax)$
is an isometric isomorphism from $\Dom(A)$ endowed with the graph norm onto
$\Gr(A)$.
Hence $\Dom(A)$ is complete if and only if $\Gr(A)$ is complete. Since
$X\times Y$ is a Banach space, this is equivalent to $\Gr(A)$ being closed.
\end{proof}

The definition of a closed linear operator can be restated sequentially as follows.

\begin{proposition} \label{prop:closed-sequential1}
 A linear operator $A$ with domain $\Dom(A)$ is closed if and only if the following holds:
 whenever $x_n\to x$ in $X$, with $x_n\in \Dom(A)$ for all $n$, and $Ax_n\to y$ in $Y$, then $x\in \Dom(A)$ and $Ax = y$.
\end{proposition}

This criterion is used to prove closedness in the next two examples.

\begin{example}\label{ex:closed1}
The derivative operator, as a linear operator in $C[0,1]$ with domain $C^1[0,1]$, is densely defined and closed.
The density of $C^1[0,1]$ in $C[0,1]$ is clear (by the Weierstrass approximation theorem we can even approximate continuous functions with polynomials). To prove closedness,
suppose that $f_n \to f$ in $C[0,1]$, with $f_n\in C^1[0,1]$ for all $n$, and $f'_n \to g$ in $C[0,1]$. We must prove that
$f\in C^1[0,1]$ and $f'=g$. For all $x\in [0,1]$ we have
$$ f(x)-f(0) = \limn f_n(x)-f_n(0) = \limn \int_0^x f'_n(y)\ud y = \int_0^x g(y)\ud y,$$
using the uniform convergence of $f_n'$ to $g$ in the last step.
The right-hand side is a continuously differentiable function, with derivative $g$. This proves that
$f\in C^1[0,1]$ and $f'=g$.
\end{example}

An analogous result holds for weak derivatives in $L^p(D)$, where $D$ is an open subset of $\R^d$; see
Section \ref{sec:weakder}.

\begin{example}\label{ex:closed2}
Let $(\Om,\calF\!,\mu)$ be a measure space and let $1\le p\le \infty$. Given a measurable function $m:\Om\to \K$, define
\begin{align*}\Dom(A_m) & := \{f\in L^p(\Om):\, mf \in L^p(\Om)\}, \\ A_m f &:=  mf, \quad f\in \Dom(A_m).
\end{align*}
 We claim that the linear operator $A_m$ is  closed. Moreover, if $1\le p<\infty$, then $A_m$ is densely defined.

To prove closedness, let $f_n\to f$ in $L^p(\Om)$ with each $f_n$ in $\Dom(A_m)$ and $A_mf_n \to g$ in $L^p(\Om)$. We must show that $f\in \Dom(A_m)$ and $A_m f = g$.
By passing to a subsequence, we may assume that both convergences also hold pointwise $\mu$-almost everywhere. Then, for $\mu$-almost all $\om\in\Om$,
$$ g(\om) = \limn A_m f_n(\om) = \limn m(\om)f_n(\om) = m(\om) f(\om).$$
This proves that $mf\in L^p(\Om)$ and $mf = g$ $\mu$-almost everywhere, hence as elements of $L^p(\Om)$.
This is equivalent to saying that $f\in \Dom(A_m)$ and $A_m f = g$.

Now let $1\le p<\infty$. By dominated convergence,  $\lim_{N\to\infty} \one_{\{|m|\le N\}} f = f$ for all $f\in L^p(\Om)$, with convergence in the norm of $L^p(\Om)$. Since $\one_{\{|m|\le N\}} f \in \Dom(A_m)$, this shows that $\Dom(A_m)$ is dense in $L^p(\Om)$.
\end{example}

\begin{example}\label{ex:bddpert-closed} If $A$ is a closed operator and $B$ is bounded, then the operator $A+B$ with domain $\Dom(A+B) := \Dom(A)$
defined by $(A+B)x:= Ax+Bx$ for $x\in\Dom(A)$ is closed. The easy proof is left as an exercise.
\end{example}

\begin{example}\label{ex:inverse-closed} If $A$ is an injective closed operator (in particular, if $A$ is an injective bounded operator), its inverse $A^{-1}$\!, with domain $\Dom(A^{-1}) = \Ran(A)$, is closed. This is immediate by noting that the graph of $A^{-1}$ equals
$$ \{(y,A^{-1}y): \, y\in \Dom(A^{-1})\} = \{(Ax,x): \, x\in \Dom(A)\}$$
and that the latter is closed in $Y\times X$ since $\{(x,Ax): \, x\in \Dom(A)\}$ is closed in $X\times Y$.
\end{example}

Further examples will be given later on. We highlight two of them:

\begin{example} The adjoint $A^*$ of a densely defined linear operator $A$ acting between Banach spaces is closed by Proposition \ref{prop:adj-dd-wsclosed}.
Likewise, by Proposition \ref{prop:graph-As}, the Hilbertian adjoint $A^\star$ of a densely defined linear operator $A$ acting between Hilbert spaces is closed.
\end{example}

\begin{example} Generators of $C_0$-semigroups are closed by Proposition \ref{prop:sg-prop}.
\end{example}

It frequently happens that linear operators are initially defined on a `too small' domain to be closed, but can be extended to a closed operator on a larger domain.
Typical examples of this situation arise in connection with differential operators, which initially can be defined on compactly supported smooth functions only.

When $A$ and $B$ are linear operators satisfying $\Dom(A)\subseteq \Dom(B)$ and $Ax = Bx$ for all $x\in \Dom(A)$, we call $B$ an {\em extension}\index{extension!of a linear operator} of $A$ and write\index{$A\subseteq B$}
$$ A\subseteq B.$$

\begin{definition}[Closability and closure]
 A linear operator is said to be {\em closable}\index{closable!operator} if it has a closed extension, or equivalently, if the closure of its graph is the graph of a linear operator.
The unique linear operator $\ov A$ whose graph is the closure of the graph of a closable operator $A$ is called the {\em closure}\index{closure!of a closable operator} of $A$; it is the smallest closed extension of $A$.
\end{definition}

We have the following analogue of Proposition \ref{prop:closed-sequential1}:

\begin{proposition} \label{prop:closed-sequential2}
A linear operator $A$ with domain $\Dom(A)$ is closable if and only if the following holds:
whenever $x_n\to 0$ in $X$ and $Ax_n\to y$ in $Y$, with all $x_n$ in $\Dom(A)$, then $y=0$.
\end{proposition}
\begin{proof}
We only need to prove the `if' part, the `only if' part being trivial. Denote the closure of $\Gr(A)$ by $\ov G$. We must prove that, under the stated condition, $\ov G$ is the graph of a linear operator $B$. This is the case if and only if
$(x,y_1)\in \ov G,\, (x,y_2)\in \ov G$ implies $y_1=y_2$, for in that case we may define $\Dom(B)$ to be the set of all $x\in X$ such that $(x,y)\in \ov G$; for $x\in \Dom(B)$ we may then define $Bx:= y$, where $y\in Y$ is the unique element such that $(x,y)\in \ov G$.
By a limiting argument, the linearity of $A$ implies that the operator $B$ thus defined is linear. Clearly it extends $A$ and its graph $\Gr(B) = \ov G$ is closed.

Suppose, therefore, that $(x,y_1)\in \ov G$ and $(x,y_2)\in \ov G$. Then $(0,y_1-y_2) \in \ov G$ since $\ov G$ is a linear subspace of $X\times Y$, and this means that there exists a sequence $(x_n,Ax_n)\to (0,y_1-y_2)$ in $X\times Y$. But then $x_n\to 0$ in $X$ and $Ax_n\to y_1-y_2$ in $Y$. By our assumption, this forces $y_1-y_2=0$.
\end{proof}

\begin{example}
 In the setting of Examples \ref{ex:closed1} and \ref{ex:closed2}, a closable operator is obtained by replacing the domain $\Dom(A)$ of the operator $A$ by any smaller subspace $Y$. The closure of the operator thus obtained equals $A$ if and only if $Y$ is dense in $\Dom(A)$ with respect to the graph norm.
\end{example}

\begin{example}
It is shown in Proposition \ref{prop:symmetry} in the next section that every densely defined symmetric operator acting in a Hilbert space is closable.
\end{example}

\begin{example}\label{ex:Lp-weak-der-closable} Let $D$ be a nonempty open subset of $\R^d$\!, let $1\le p\le \infty$, and let $\alpha\in \N^d$ be a multi-index. In $L^p(D)$ we consider the linear operator $A$ with domain $C_{\rm c}^\infty(D)$ defined by $$ Af := \partial^\al f, \quad f\in C_{\rm c}^\infty(D),$$
where $\partial^\al = \partial_1^{\al_1}\circ\cdots\circ \partial_d^{\al_d}$, with $\partial_j = \partial/\partial x_j$ denoting the $j$th directional derivative.
We claim that $A$ is closable. Indeed, suppose the functions $f_n\in C_{\rm c}^\infty(D)$ satisfy $f_n\to 0$ and $Af_n\to g$ in $L^p(D)$. Integrating by parts,
for all $\phi\in C_{\rm c}^\infty(D)$ we obtain
\begin{align}\label{eq:Lp-weak-der-closable} \int_D g \phi\ud x = \limn  \int_D (\partial^\al f_n) \phi\ud x = (-1)^{|\al|} \limn \int_D f_n \partial^\al\phi\ud x = 0,
\end{align}
where $|\al| = \al_1+\cdots+\al_d$; the last step follows by H\"older's inequality (cf. Corollary \ref{cor:Holder}). It is shown in Proposition \ref{prop:fundvarcalculus} in the next chapter that \eqref{eq:Lp-weak-der-closable} implies $g=0$ almost everywhere.
\end{example}

\begin{example}\label{ex:weak-Laplace-closable}
 Let $D$ be a nonempty open subset of $\R^d$ and let $1\le p\le \infty$. In $L^p(D)$ we consider the linear operator $A$ with domain $C_{\rm c}^\infty(D)$ defined by $$ Af := \Delta f, \quad f\in C_{\rm c}^\infty(D),$$
where $\Delta f = \partial_1^2 f+\cdots+\partial_d^2f$ is the Laplacian of $f$. As in the previous example one shows that this operator is closable.

Various closed realisations of the Laplacian extending this operator
will be discussed in Chapters
\ref{chap:boundaryvalueproblems}--\ref{chap:semigroups}; see in
particular Sections \ref{sec:Hs}, \ref{sec:examples-forms}, and
\ref{subsec:heat-sgr}.
\end{example}

\subsection{The Adjoint Operator}

When $A$ is a densely defined linear operator from $X$ to $Y$,
we may uniquely define a linear operator $A^*$ from $Y^*$ to $X^*$
by defining its domain $\Dom(A\s)$ to be
the set of all $y\s\in Y^*$ with the property that
there exists an element $x\s\in X^*$ such that
$$\lb x, x\s\rb=\lb Ax, y\s\rb,\quad x\in \Dom(A).$$ Since $\Dom(A)$
is dense in $X$, the element $x\s\in X^*$ (if it exists) is unique
and we can set $$A^*y\s:=x^*\!, \quad y\s\in \Dom(A\s).$$
Thus, by definition, we have the identity
$$ \lb Ax, y\s\rb = \lb x,A^*y\s\rb,\quad x\in \Dom(A), \ y\s\in \Dom(A\s).$$

\begin{definition}[Adjoint operator] The operator $A^*$ is called the
{\em adjoint}\index{adjoint!of a densely defined operator} of $A$.
\end{definition}

The adjoint of a closable densely defined operator $A$ equals the adjoint of the closure $\ov A$.
Indeed, if $x^*\in X^*$ and $y^*\in Y^*$ are such that
$\lb x,x^*\rb=\lb Ax,y^*\rb$ for all $x\in\Dom(A)$, then by continuity
this identity extends to all $x\in\Dom(\ov A)$.

\begin{proposition}\label{prop:adj-dd-wsclosed}
If $A$ is a densely defined linear operator from $X$ to $Y$, then $A^*$ is weak$^*$ closed in the sense that its graph is weak$^*$ closed in $Y^*\times X^*$ and we have $$\Gr(A\s) =(J(\Gr(A)))^\perp\!,$$
where $J: X\times Y\to Y\times X$ is defined by $J(x,y)= (-y,x)$.
If $A$ is densely defined and closed, then $A^*$ is weak$^*$ densely defined in the sense that its domain is weak$^*$ dense.
\end{proposition}
\begin{proof}
The pairing  $$\lb (y,x), (y^*\!, x\s)\rb:=\lb
y,y\s\rb + \lb x, x\s\rb$$
allows us to identify $Y\s\times X^*$ with the dual of $Y\times
X$. By the definition of the adjoint operator we have $(y^*\!, x\s)\in
\Gr(A\s)$ if and only if $$\lb (-Ax, x), (y^*\!, x\s)\rb =0,
\quad x\in \Dom(A).$$
This proves the identity $\Gr(A\s) =(J(\Gr(A)))^\perp\!$.
Since annihilators are weak$^*$ closed, this proves that $A^*$ is weak$^*$ closed.

Suppose now that $A$ is closed.
Since $J(\Gr(A))$ is closed, the graph identity just proved
and Proposition \ref{prop:perpperp} give
$$
{}^\perp(\Gr(A^*))
={}^\perp\bigl((J(\Gr(A)))^\perp\bigr)
=J(\Gr(A)).
$$
If $y\in Y$ annihilates $\Dom(A^*)$, then
$(y,0)\in{}^\perp(\Gr(A^*))$.
Thus $(y,0)=(-Ax,x)$ for some $x\in\Dom(A)$,
which forces $x=0$ and $y=0$.
Hence ${}^\perp\Dom(A^*)=\{0\}$, and Proposition
\ref{prop:perpperp} implies that $\Dom(A^*)$ is weak$^*$ dense
in $Y^*$.
\end{proof}

If the linear operators
$A$ and $B$ act from $X$ to $Y$, we define the operator $A+B$ acting from $X$ to $Y$ by
\begin{align*}
\Dom(A+B) & := \Dom(A)\cap \Dom(B), \\ (A+B)x & := Ax +Bx, \quad x\in \Dom(A+B).
\intertext{If $A$ acts from $X$ to $Y$ and $B$ acts from $Y$ to another Banach space $Z$, we define
the operator $BA$ acting from $X$ to $Z$ by}
\Dom(BA) & := \{x\in \Dom(A):\, Ax\in \Dom(B)\}, \\ BAx &:= B(Ax), \quad x\in \Dom(BA).
\end{align*}
There is of course {\em a priori} no guarantee that $\Dom(A+B)$ and $\Dom(BA)$ contain any nontrivial elements even when both $A$ and $B$ are densely defined.

\begin{proposition}\label{prop:AB-properties}
Let $A$ and $B$ be densely defined operators acting as indicated above. Then:
 \begin{enumerate}[label={\rm(\arabic*)}, leftmargin=*]
  \item\label{it:AB-properties1} if $A\subseteq B$, then $B^*\subseteq A^*$;
  \item\label{it:AB-properties2} if $A+B$ is densely defined, then $A^*+B^*\subseteq (A+B)^*$, with equality if $B$ is bounded;
  \item\label{it:AB-properties3} if $BA$ is densely defined, then $A^* B^*\subseteq (BA)^*$, with equality if $B$ is bounded.
 \end{enumerate}
\end{proposition}
\begin{proof}
 Part \ref{it:AB-properties1} is immediate from the definitions.

 If $y^* \in
 \Dom(A^*+B^*) = \Dom(A^*)\cap \Dom(B^*)$, then for all $x\in \Dom(A+B) = \Dom(A)\cap \Dom(B)$
 we have $\lb (A+B)x, y^* \rb  =\lb Ax, y^* \rb +\lb Bx, y^* \rb  =  \lb x, A^* y^* + B^* y^* \rb $,
 so $y^*\in \Dom((A+B)^*)$ and $(A+B)^* y^* = A^* y^* + B^* y^*$\!.
 If $B$ is bounded and $y^*\in \Dom((A+B)^*)$,
 then for all $x\in \Dom(A+ B) = \Dom(A)$ we have
 $\lb Ax, y^* \rb  = \lb (A+B)x, y^* \rb  - \lb Bx, y^* \rb
 = \lb x, (A+B)^* y^* \rb - \lb x, B^* y^* \rb ,$
 so that $y^*\in \Dom(A^*) = \Dom(A^*+B^*)$ and $A^* y^* = (A+B)^* y^* - B^* y^*$\!.
 This gives \ref{it:AB-properties2}.

 If $z^*\in \Dom(A^* B^*)$, then $z^*\in \Dom(B^*)$ and $B^* z^* \in \Dom(A^*)$, and
 for all $x\in \Dom(BA)$ we have
 $\lb (BA)x, z^* \rb  = \lb Ax,  B^* z^* \rb = \lb x, A^* B^* z^* \rb $, so that $z^*\in \Dom((BA)^*)$
 and $(BA)^* z^* = A^* B^* z^*$\!. If $B$ is bounded and $z^*\in \Dom((BA)^*)$,
 then for all $x\in \Dom(BA) = \Dom(A)$ we have
 $\lb Ax, B^* z^* \rb  = \lb BAx, z^* \rb  = \lb x,(BA)^* z^* \rb $,
 so $B^* z^*\in \Dom(A^*)$ and $z^*\in \Dom(A^* B^*)$ and $A^* B^* z^* = (BA)^* z^*.$
 This gives \ref{it:AB-properties3}.
\end{proof}

We have the following useful duality criterion for deciding whether an element belongs to the domain of an operator.

\begin{proposition}\label{prop:DomAweak}
Let $A$ be a densely defined closed operator from $X$ to $Y$.
If $x\in X$ and $y\in Y$ are such that
$\lb y,y\s\rb = \lb x,A^*y\s\rb$ for all $y\s\in
\Dom(A\s)$, then $x\in\Dom(A)$ and $Ax = y$.
\end{proposition}
\begin{proof}
Put $J(u,v):=(-v,u)$.
By the assumption, for every $y^*\in\Dom(A^*)$ we have
$$
\lb J(x,y),(y^*\!,A^*y^*)\rb
=-\lb y,y^*\rb+\lb x,A^*y^*\rb
=0.
$$
Thus $J(x,y)\in{}^\perp(\Gr(A^*))$.
By Propositions \ref{prop:adj-dd-wsclosed} and
\ref{prop:perpperp}, and the closedness of $A$,
$$
{}^\perp(\Gr(A^*))
={}^\perp\bigl((J(\Gr(A)))^\perp\bigr)
=J(\Gr(A)).
$$
Since $J$ is injective, it follows that $(x,y)\in\Gr(A)$.
Thus $x\in\Dom(A)$ and $Ax=y$.
\end{proof}

In what follows we let $H$ and $K$ be Hilbert spaces. When $A$ is a densely defined operator acting from $H$ to $K$,
the Riesz representation theorem may be used to identify the adjoint $A^*$, which acts from $K^*$ to $H^*$, with a linear operator $A^\star$ acting from $K$ to $H$.
Thus, by definition, an element $k\in K$ belongs to $\Dom(A^\star)$ if
there exists a (necessarily unique) element $h\in H$ such that
$$\iprod{x}{h} = \iprod{Ax}{k},\quad x\in \Dom(A),$$
and in that case $A^\star k = h$.
Thus we have the identity
$$ \iprod{x}{A^\star k} = \iprod{Ax}{k}, \quad x\in \Dom(A), \ k\in \Dom(A^\star).$$

\begin{definition}[Hilbert space adjoint]
The operator $A^\star$ is called the {\em Hilbert space adjoint} of $A$.\index{adjoint!of a densely defined Hilbert space operator}
\end{definition}

Proposition \ref{prop:adj-dd-wsclosed} admits the following Hilbertian version.
We denote by $H\oplus K$ the Hilbert space obtained by endowing the product $H\times K$ with the inner product
$$\iprod{(h, k)}{(h', k')} := \iprod{h}{h'}+ \iprod{k}{k'}.$$

\begin{proposition}\label{prop:graph-As}
If $A$ is a densely defined linear operator from $H$ to $K$, then $A^\star$ is closed and we have $$\Gr(A^\star) =(J(\Gr(A)))^\perp\!,$$
where $J: H\oplus K\to K\oplus H$ is defined by $J(x,y) := (-y,x)$.
If $A$ is densely defined and closed, then $A^\star$ is densely defined.
\end{proposition}
\begin{proof}
In Hilbert spaces the weak topology and the weak$^*$ topology agree, and a subspace is weakly dense (respectively, weakly closed) if and only if it is dense (respectively, closed).
Hence everything follows from Proposition \ref{prop:adj-dd-wsclosed}, except the statement that $\Gr(A^\star)$ is the orthogonal complement
of $J(\Gr(A))$ in $K\oplus H$. This follows from
\begin{align*}
(k,h) \in \Gr(A^\star)
  & \Leftrightarrow \iprod{Ax}{k} = \iprod{x}{h}\ \hbox{for all}  \ x\in \Dom(A)
\\ &  \Leftrightarrow (k,h)\perp (-Ax, x) \ \hbox{for all}  \ x\in \Dom(A)
  \Leftrightarrow
 (k,h)\perp J(\Gr(A)).
 \qedhere
\end{align*}
\end{proof}

{\em Mutatis mutandis},
Proposition \ref{prop:AB-properties} admits a Hilbertian version as well; we leave this as an exercise for the reader.

\begin{proposition}\label{prop:adj-dd}
If $A$ is a densely defined closed operator acting from $H$ to $K$, then
$A = A^{\star\star}$
with equality of domains.
\end{proposition}
\begin{proof}
If $Z$ is a subspace of $H\oplus K$, then
$$(k,h)\perp J(Z) \Leftrightarrow (h,-k)\in Z^\perp \Leftrightarrow
(k,h) \in J(Z^\perp),$$
which shows that $(J(Z))^\perp = J(Z^\perp).$
By Proposition \ref{prop:graph-As},
$ \Gr(A^\star) = (J(\Gr(A)))^\perp
= J((\Gr(A))^\perp)$
and
\[
\Gr(A^{\star\star}) = (J(\Gr(A^\star)))^\perp =(JJ((\Gr(A))^\perp))^\perp
= (\Gr(A))^{\perp\perp} = \Gr(A).
\qedhere\]
\end{proof}

For operators acting in Hilbert spaces we have the following extension of
Proposition \ref{prop:injective-denserange}, the proof of which is almost {\em verbatim} the same:

\begin{proposition}\label{prop:injective-denserange-unbdd}
If $A$ is a densely defined closed operator from $H$ into another Hilbert space $K$,
then $H$ and $K$ admit orthogonal decompositions
$$
H = \Ker(A) \oplus \overline{\Ran(A^\star)}, \qquad
K = \Ker(A^\star) \oplus \overline{\Ran(A)}.
$$
In particular,
\begin{enumerate}[label={\rm(\arabic*)}, leftmargin=*]
 \item $A$ is injective if and only if $A^\star$ has dense range;
 \item $A$ has dense range if and only if $A^\star$ is injective.
\end{enumerate}
\end{proposition}

\subsection{The Spectrum}\label{subsec:spectrum-unbdd}

The spectrum of a linear operator is defined in the same way as for bounded operators, except that explicit attention has to be paid to domains.
Throughout this section, all vector spaces are complex.

\begin{definition}[Spectrum and resolvent]
 The {\em resolvent set}\index{resolvent!set} of a linear operator $A$ acting in a Banach space $X$ is the set  $\varrho(A)$ consisting of all $\la\in\C$ for which
 the operator $\la I - A$ has a two-sided inverse, that is, there exists a bounded operator $U$ on $X$
 such that:
 \begin{enumerate}[label={\rm(\roman*)}, leftmargin=*]
  \item for all $x\in \Dom(A)$ we have $U(\la I-A)x = x$;
  \item for all $x\in X$ we have $Ux\in \Dom(A)$ and $(\la I-A)Ux = x$.
 \end{enumerate}
The {\em spectrum}\index{spectrum} of $A$ is defined as the complement of the resolvent set of $A$: $$\si(A):= \C\setminus\varrho(A).$$
\end{definition}

We emphasise that, although $A$ is allowed to be unbounded, the two-sided inverse $U = (\la I-A)^{-1}$ is
required to be bounded.
It is customary to write\index{$Raa$@$R(\la,A)$}
$$ R(\la,A) := (\la-A)^{-1}$$
for the {\em resolvent operator}\index{resolvent!operator} of $A$ at $\la$.
As in the bounded case, the {\em resolvent identity} holds: if $\la,\mu\in \varrho(A)$, then \index{resolvent!identity}
\begin{align}\label{eq:res-id-unbdd} R(\la,A) - R(\mu,A)= (\mu-\la)R(\la,A) R(\mu,A).
\end{align}

By the observations in Examples \ref{ex:bddpert-closed} and \ref{ex:inverse-closed}, a linear operator $A$ in $X$ with non\-empty resolvent set is closed.
The proofs of Lemma \ref{lem:res-open-bdd}, the holomorphy assertion in Lemma \ref{lem:spectr-holo-bdd}, and Propositions \ref{prop:res-blowup-bdd} and \ref{prop:approx-eigenvalue-bdry} carry over {\em verbatim}, and Proposition \ref{prop:closed-range} carries over with an obvious adaptation of the proof. For the reader's convenience we state the results here:

\begin{proposition}\label{prop:closed-range-unbdd} If $A$ is closed and satisfies $\n Ax \n \ge C \n x\n$ for some $C>0$ and all $x\in \Dom(A)$, then
 $A$ is injective and has closed range.
\end{proposition}

\begin{proposition}\label{prop:res-open}
The spectrum
$\sigma(A)$ is a closed subset of $\C$. More precisely, if
$X\neq\{0\}$ and $\lambda\in\varrho(A)$, then
 $B(\la;r)\subseteq\varrho(A)$ with $r = 1/\n R(\la,A)\n$. Moreover, if $|\la-\mu|\le \delta r$ with $0\le \delta<1$,
 then $$\n R(\mu,A)\n \le \frac{1}{1-\delta}\n R(\la,A)\n.$$
\end{proposition}

\begin{proposition}\label{prop:spectr-holo} The function $\la\mapsto R(\la,A)$ is holomorphic on $\varrho(A)$, and its complex derivative is given by  $-R(\la,A)^2$\!.
\end{proposition}

\begin{proposition}\label{prop:res-blowup}
 If $\la_n\to \la$ in $\C$, with each $\la_n\in \varrho(A)$ and with $\la\in \partial\varrho(A)$, then
 $$\limn \n R(\la_n,A)\n = \infty.$$
\end{proposition}

The following proposition gives a simple but powerful uniqueness result:

\begin{proposition}\label{prop:semigroupsAB}
 Let $A$ and $B$ be linear operators acting in a Banach space $X$. If $\varrho(A)\cap \varrho(B)\not=\emptyset$ and $B$ is an extension of $A$, then $A = B$ with equality of domains.
\end{proposition}
\begin{proof}
Fix $\la\in \varrho(A)\cap \varrho(B)$. Then for all $x\in X$ we have
$R(\la,A)x\in \Dom(A)\subseteq \Dom(B)$ and
$$ (\la-B)R(\la,A)x = (\la-A)R(\la,A)x = x.$$
Multiplying both sides from the left by $R(\la,B)$ gives
$ R(\la,A)x = R(\la,B)x$.
Since $x\in X$ was arbitrary, we conclude that $R(\la,A) = R(\la,B)$ and therefore $\Dom(A) = \Dom(B)$.
\end{proof}

The proof of the following result is the same as those of
Propositions \ref{prop:spectrum-dual} and \ref{prop:spectrum-HS-dual}.

\begin{proposition}\label{prop:spectrum-dual-unbdd} If $A$ is a densely defined closed operator in a Banach space $X$, then
$$\sigma(A^*) = \sigma(A).$$
If $A$ is a densely defined closed operator in a Hilbert space $H$, then
$$\sigma(A^\star) = \ov{\sigma(A)}.$$
\end{proposition}

For later use we compute the spectrum of a simple diagonal operator.

\begin{proposition}\label{prop:sigmaAdiag} Let $A$ be a densely defined closed operator in a separable Hilbert space $H$ with an orthonormal basis
$(h_n)_{n\ge 1}$ of eigenvectors. If $\varrho(A)\not=\emptyset$ and the corresponding eigenvalue sequence $(\la_n)_{n\ge 1}$ satisfies $\limn |\la_n| = \infty$, then
$$\sigma(A) = \{\la_n:\, n\ge 1\}.$$
\end{proposition}

The assumption $\varrho(A)\not=\emptyset$ cannot be omitted; see Problem \ref{prob:diag-empty-res}.

\begin{proof}
Since every $\la_n$ is an eigenvalue of $A$, we have
$\{\la_n:\,n\ge1\}\subseteq\sigma(A)$.
Choose $\mu_0\in\varrho(A)$ and put $R_0:=R(\mu_0,A)$.
Then
$$
R_0h_n=\frac1{\mu_0-\la_n}h_n,\quad n\ge1.
$$

Fix $\mu\notin\{\la_n:\,n\ge1\}$ and define
$$
C_\mu:=I+(\mu-\mu_0)R_0.
$$
This is a bounded operator satisfying
$$
C_\mu h_n
=\frac{\mu-\la_n}{\mu_0-\la_n}h_n,\quad n\ge1.
$$
Since $|\la_n|\to\infty$, the coefficients on the right-hand side converge
to $1$. None of them is zero, so their absolute values are
bounded away from zero. Consequently, $C_\mu$ is invertible
in $\calL(H)$.

For every $x\in\Dom(A)$ we have
$$
C_\mu(\mu_0-A)x
=(\mu_0-A)x+(\mu-\mu_0)x
=(\mu-A)x.
$$
Thus $\mu-A$ has the bounded inverse $R_0C_\mu^{-1}$,
and hence $\mu\in\varrho(A)$.
This proves the reverse inclusion
$\sigma(A)\subseteq\{\la_n:\,n\ge1\}$.
\end{proof}

We conclude this section with two useful elaborations on
Proposition \ref{prop:res-open}. In the first, we write $\Sigma_\varphi := \{z\in \C\setminus\{0\}:\, |\arg(z)|<\varphi\},$
the argument being taken in $(-\pi,\pi)$.

\begin{lemma}\label{lem:res-Taylor} Let $A$ be a linear operator acting in a nonzero Banach space $X$.
If the open half-line $(0,\infty)$ is contained in $\varrho(A)$ and
$$\sup_{\lambda>0} \n \lambda R(\lambda,A)\n =: M<\infty,$$
then $M\ge 1$,
and for all $\varphi\in (0,\frac12\pi)$ with $\sin \varphi<1/M$ we have $\Sigma_\varphi\subseteq \varrho(A)$
and $$\sup_{\lambda\in\Sigma_\varphi} \n \lambda R(\lambda,A)\n \le \frac{M}{1-M\sin\varphi}.$$
\end{lemma}

\begin{proof}
For nonzero $x\in \Dom(A)$ we have $\lambda R(\lambda,A)x=x+R(\lambda,A)Ax\to x$ as $\lambda\to\infty$, from which it follows that $M\geq 1$.

Let $\varphi\in(0,\frac12\pi)$ satisfy $\sin\varphi<1/M$, and fix $\la\in\Sigma_\varphi$. Write $\la=|\la|e^{i\theta}$, where $|\theta|<\varphi$, and put $\mu:=|\la|/\cos\theta$. Then $\mu>0$ and
$$
\frac{|\la-\mu|}{\mu}=|\sin\theta|\leq\sin\varphi<\frac1M,
$$
where the first identity follows by direct computation. Since $\n R(\mu,A)\n\leq M/\mu$, it follows that $|\la-\mu|<1/\n R(\mu,A)\n$. By Proposition \ref{prop:res-open}, $\la\in\varrho(A)$.

The Neumann series expansion around $\mu$ gives
\begin{align*}
\n R(\la,A)\n
&\leq \n R(\mu,A)\n
   \sum_{n=0}^\infty |\la-\mu|^n\n R(\mu,A)\n^n  \\
&\leq \frac{M}{\mu}
   \sum_{n=0}^\infty \Bigl(\frac{|\la-\mu|}{\mu}M\Bigr)^n
 \leq \frac{M}{\mu}\sum_{n=0}^\infty (M\sin\varphi)^n  \\
&= \frac{M}{1-M\sin\varphi}\cdot\frac1{\mu}
 \leq \frac{M}{1-M\sin\varphi}\cdot\frac1{|\la|}.\qedhere
\end{align*}
\end{proof}

A typical application of this lemma is the second part of the next lemma, which
 extends a uniform bound on $\la R(\la,A)$ on a half-plane to a larger sector. For reasons of completeness we also include its counterpart for uniform bounds on $R(\la,A)$.

\begin{lemma}\label{lem:res-Taylor-cor}
Suppose that the half-plane
$\C_+= \{\Re\la>0\} = \{|\arg(\la)|< \frac12\pi\}$ is contained in the resolvent set $\varrho(A)$ of the linear operator $A$ acting in a Banach space $X$. Then:
\begin{enumerate}[label={\rm(\arabic*)}, leftmargin=*]
 \item\label{it:res-Taylor1} if $\sup_{\la\in\C_+} \n R(\la,A)\n < \infty,$ then there exists a $\delta>0$ such that
 $\{\Re\la>-\delta\}\subseteq \varrho(A)$ and  $$\sup_{\Re\la>-\delta} \n R(\la,A)\n < \infty;$$
 \item\label{it:res-Taylor2} if $\sup_{\la\in\C_+} \n \la R(\la,A)\n < \infty,$ then there exists a $\delta>0$ such that
 $\{|\arg\la| < \frac12\pi +\delta\}\subseteq \varrho(A)$ and $$\sup_{|\arg\la| < \frac12\pi +\delta} \n \la R(\la,A)\n < \infty.$$
\end{enumerate}
\end{lemma}
\begin{proof} Proposition \ref{prop:res-blowup} implies that
in case \ref{it:res-Taylor1} we have $i\R\subseteq\varrho(A)$,
and that in case \ref{it:res-Taylor2} we have $i\R\setminus\{0\}\subseteq\varrho(A)$.
The result now follows from Proposition \ref{prop:res-open} applied to the points $\la \in i\R$ (in case \ref{it:res-Taylor1})
and Lemma \ref{lem:res-Taylor} applied to the operators $\pm iA$ (in case \ref{it:res-Taylor2}).
\end{proof}

\section{Unbounded Selfadjoint Operators}\label{sec:unbdd-sa}

We now specialise the general theory of unbounded operators to symmetric and selfadjoint operators.

In what follows we let $H$ be a complex Hilbert space.

\begin{definition}[Symmetric and positive operators]\label{def:A-pos-symm}  A linear operator $A$ acting in $H$ is called:
\begin{itemize}
 \item {\em symmetric},\index{symmetric!operator}  if for all $x,y\in \Dom(A)$ we have
$ \iprod{Ax}{y} = \iprod{x}{Ay}.$
 \item {\em positive},\index{positive!operator}  if for all $x\in \Dom(A)$ we have
$ \iprod{Ax}{x}\ge 0.$
\end{itemize}
\end{definition}

Over the complex scalars,
positive operators are symmetric.
Indeed, if $A$ is positive, then for all $x\in \Dom(A)$ we have
$ \iprod{Ax}{x} = \ov{\iprod{Ax}{x}} = \iprod{x}{Ax}.$
By polarisation (as in the proof of Proposition \ref{prop:polarisation}, this requires working over the complex scalars), this implies
$ \iprod{Ax}{y} = \iprod{x}{Ay}$ for all $x,y\in \Dom(A)$.

It is an immediate consequence of Definition \ref{def:A-pos-symm} and the definition of $A^\star$ that if $A$ is densely defined and symmetric, then $\Dom(A)\subseteq \Dom(A^\star)$ and
$Ax = A^\star x$ for all $x\in \Dom(A)$, that is, $A^\star$ is an extension of $A$. Since $A^\star$ is closed, we have shown:

\begin{proposition}\label{prop:symmetry} If $A$ is a densely defined symmetric operator, then $A$ is closable and $A^\star$ is a closed extension of $A$.
\end{proposition}

In general, $\Dom(A)$ may be strictly smaller than $\Dom(A^\star)$. A simple example is the Laplace operator\index{$D$@$\Delta$} $$\Delta =  \partial_1^2 +\cdots+\partial_d^2$$ on
$L^2(\R^d)$ with domain $C_{\rm c}^\infty(\R^d)$: this operator is densely defined and symmetric but not closed, and therefore $\Delta^\star$ is a proper extension of $\Delta$. This motivates the following definition.

\begin{definition}[Selfadjoint operators] A densely defined operator $A$ in $H$ is called
{\em selfadjoint}\index{selfadjoint operator}
 if $A= A^\star$\!, that is, if $\Dom(A) = \Dom(A^\star)$ and $Ah = A^\star h$ for all $h\in \Dom(A) = \Dom(A^\star)$.
The operator $A$ is called {\em essentially selfadjoint}\index{essentially!selfadjoint} if it is closable and its closure $\ov A$ is selfadjoint.
\end{definition}

By Propositions \ref{prop:adj-dd} and \ref{prop:symmetry}, a densely defined symmetric operator $A$ is essentially selfadjoint if and only if $A^\star$ is symmetric.

\begin{example}[Multipliers]\label{ex:spectr-meas-multipliers}
Let $(\Om,\calF\!,\mu)$ be a $\sigma$-finite measure space and let
$m:\Om\to\C$ be a measurable function.
It has been shown in Example \ref{ex:closed2} that the linear operator $M_m$
in $L^2(\Om)$ defined by
\begin{align*}
 \Dom(M_m) & := \{f\in L^2(\Om):\, mf\in L^2(\Om)\}, \\
 M_m f & := mf, \quad f\in \Dom(M_m),
\end{align*}
is densely defined and closed.
It is immediate from the definition of the Hilbert space adjoint that
$M_m^\star=M_{\ov m}$, with equality of domains
$\Dom(M_{\ov m})=\Dom(M_m)$. Hence, if $m$ is real-valued
$\mu$-almost everywhere, then $M_m$ is selfadjoint.

Conversely, if $M_m$ is selfadjoint, then
$(m-\ov m)f=0$ for all $f\in\Dom(M_m)$.
By $\sigma$-finite\-ness, choose measurable sets
$F_1\subseteq F_2\subseteq\cdots$ of finite measure whose union is
$\Om$. Then we have
$\one_{F_n\cap\{|m|\le n\}}\in\Dom(M_m)$,
and therefore
$(m-\ov m)\one_{F_n\cap\{|m|\le n\}}=0$
almost everywhere. Since the sets
$F_n\cap\{|m|\le n\}$ increase to $\Om$, it follows that
$m=\ov m$ almost everywhere.
\end{example}

\begin{example}[Fourier multipliers]\label{ex:sa-FM} Let $m:\R^d\to \R$ be real-valued and measurable,
and let $T_m$ denote the (possibly unbounded) Fourier
multiplier in $L^2(\R^d)$ defined by
\begin{equation}\label{eq:DomFM}
\begin{aligned}\Dom(T_m) & :=
\{f\in L^2(\R^d):\, m \wh f\in L^2(\R^d)\}, \\
T_m f &:= (m \wh f)\widecheck{\phantom{b}}, \quad f\in \Dom(T_m).
\end{aligned}
\end{equation}
Let $\calF$ denote the Fourier transform on $L^2(\R^d)$.
Since
$\Dom(T_m)=\calF^{-1}\Dom(M_m)$
and $M_m$ is densely defined, so is $T_m$.
For $g,h\in L^2(\R^d)$, the identity
$$
\iprod{T_mf}{g}=\iprod{f}{h},
\quad f\in\Dom(T_m),
$$
is equivalent, upon writing $u=\calF f$, to
$$
\iprod{M_mu}{\calF g}=\iprod{u}{\calF h},
\quad u\in\Dom(M_m).
$$
By the definition of the adjoint, it follows that
$\Dom(T_m^\star)=\calF^{-1}\Dom(M_m^\star)$ and
$$
T_m^\star=\calF^{-1}M_m^\star\calF.
$$
Since $M_m$ is selfadjoint by the preceding example,
we conclude that $T_m^\star=T_m$, with equality of domains.
\end{example}

Two cases are of special interest.

\begin{example}[The momentum operator]\label{ex:momentum}\index{momentum operator}
The multiplier $m(\xi) = \xi$ gives rise to the operator $\frac1i \frac{{\rm d}}{{\rm d}x}$ on $L^2(\R)$. With the domain given by \eqref{eq:DomFM},
this operator is selfadjoint.
With the notation and techniques developed in Section \ref{sec:Hs}, this domain is seen to be the Sobolev space $H^1(\R) = W^{1,2}(\R)$.
\end{example}

\begin{example}[The Laplacian]\label{ex:Laplacian-sa}
The multiplier $m(\xi) = - |\xi|^2$ gives rise to the Laplace operator $\Delta = \sum_{j=1}^d \frac{\partial^2}{\partial x_j^2}$ on $L^2(\R^d)$. With the domain given by \eqref{eq:DomFM}, this operator is selfadjoint.
With the notation and techniques developed in Section \ref{sec:Hs}, this domain is seen to be the Sobolev space
$H^2(\R^d) = W^{2,2}(\R^d)$.
\end{example}

The following version of Theorem \ref{thm:spect-sa} holds.

\begin{proposition}\label{prop:sa-pos-spectrum}
 If $A$ is selfadjoint, then $\si(A)\subseteq \R$. If, in addition, $A$ is positive, then $\si(A)\subseteq[0,\infty)$.
\end{proposition}
\begin{proof} This may be established by repeating parts of the proof of Theorem \ref{thm:spect-sa},
using Propositions \ref{prop:closed-range-unbdd} and \ref{prop:injective-denserange-unbdd} instead of Propositions \ref{prop:closed-range} and \ref{prop:injective-denserange}, respectively.
\end{proof}

The following proposition provides a sufficient condition for selfadjoint\-ness.

\begin{proposition}\label{prop:symm-sgr-sa}
 If $A$ is densely defined and symmetric and $\varrho(A)\cap \R\not=\emptyset$, then $A$ is selfadjoint.
\end{proposition}
\begin{proof} The nonemptiness of the resolvent set implies that $A$ is closed (cf. Example \ref{ex:inverse-closed}). The operator $A^\star$ is closed as well.  The symmetry of $A$ implies $\Dom(A)\subseteq \Dom(A^\star)$,
and if $\la\in \varrho(A)\cap \R$, then $\la \in \varrho(A^\star)$ in view of Proposition \ref{prop:spectrum-dual-unbdd}. The identity $A = A^\star$ with equality of domains therefore follows from Proposition \ref{prop:semigroupsAB}.
\end{proof}

We next show that a densely defined inverse of a selfadjoint
operator is selfadjoint.

\begin{proposition}\label{prop:sa-inv} If the linear operator $A$ in $H$ is selfadjoint, injective, and has dense range, then its inverse $A^{-1}$
with domain $\Dom(A^{-1}) = \Ran(A)$ is selfadjoint.
\end{proposition}
\begin{proof} Since $\Ran(A)$ is dense, the inverse $A^{-1}$, with domain
$\Dom(A^{-1})=\Ran(A)$, is densely defined. We claim that
$$
   (A^{-1})^\star=(A^\star)^{-1}.
$$
Indeed, let $y\in\Dom((A^{-1})^\star)$. Then there exists $z\in H$ such
that $\iprod{A^{-1}u}{y}=\iprod{u}{z}$ for all $u\in\Ran(A)$.
Writing $u=Ax$, with $x\in\Dom(A)$, gives
$\iprod{x}{y}=\iprod{Ax}{z}$ for all $x\in\Dom(A)$.
Thus $z\in\Dom(A^\star)$ and $A^\star z=y$, so
$y\in\Ran(A^\star)$ and
$(A^{-1})^\star y=z=(A^\star)^{-1}y.$
The converse follows by reversing the argument. Hence
$(A^{-1})^\star=(A^\star)^{-1}$ as claimed.
Since $A=A^\star$, we obtain $(A^{-1})^\star=A^{-1}$.
Therefore $A^{-1}$ is selfadjoint.
\end{proof}

As an application of Proposition \ref{prop:symm-sgr-sa},
we record the following result.

\begin{corollary}\label{cor:sa-c0sgr}
Let $A$ be a densely defined closed positive operator in $H$.
If $I+A$ has dense range, then $A$ is selfadjoint.
\end{corollary}
\begin{proof}
For all $x\in\Dom(A)$, positivity gives
$$
\|(I+A)x\|\,\|x\|
\ge \iprod{(I+A)x}{x}
=\|x\|^2+\iprod{Ax}{x}
\ge\|x\|^2.
$$
Thus $\|(I+A)x\|\ge\|x\|$.
Since $I+A$ is closed, Proposition \ref{prop:closed-range-unbdd}
implies that it is injective and has closed range.
Its range is dense by assumption, so it equals $H$.
The same estimate shows that $(I+A)^{-1}$ is bounded,
and hence $-1\in\varrho(A)$.
Since positivity implies symmetry, Proposition
\ref{prop:symm-sgr-sa} now gives that $A$ is selfadjoint.
\end{proof}

As an application of Corollary \ref{cor:sa-c0sgr}, we have the following sufficient condition for selfadjointness.

\begin{theorem}[Selfadjointness of $A^\star A$]\label{thm:AstarA}
 If $A$ is a densely defined closed operator from $H$ into another Hilbert space $K$, then:
 \begin{enumerate}[label={\rm(\arabic*)}, leftmargin=*]
 \item\label{it:AstarA1} the operator $A^\star A$ is densely defined, selfadjoint, and positive;
 \item\label{it:AstarA2} $\Dom(A^\star A)$ is dense in $\Dom(A)$ with respect to the graph norm.
\end{enumerate}
\end{theorem}

This operator will be revisited in Proposition \ref{prop:AstarA-via-forms} in connection with the theory of forms.

\begin{proof}
By Proposition \ref{prop:graph-As} we have the orthogonal decomposition
$$ K\oplus H = \Gr(A^\star)\oplus J(\Gr(A)),$$
where $J(x,y)= (-y,x)$. Hence, for any $u\in H$, we can find $x\in \Dom(A)$ and $y\in \Dom(A^\star)$ such that
$$ (0,u) = (y,A^\star y)+J(x,Ax) = (y-Ax, A^\star y+x).$$
It follows that $y=Ax$, which implies $x\in \Dom(A^\star A)$, and $u = A^\star y +x = (I+A^\star A)x.$ This proves that $I+A^\star A$ is surjective. Moreover, for $x\in\Dom(A^\star A)$,
$$\n x\n^2+\n Ax\n^2=\iprod{(I+A^\star A)x}{x}\le \n (I+A^\star A)x\n\,\n x\n,$$
so $I+A^\star A$ is injective and the inverse $(I+A^\star A)^{-1}$ is bounded. Hence $I+A^\star A$ is closed, and therefore so is $A^\star A$. The identity $\iprod{A^\star Ax}{x} = \n Ax\n^2$ shows that $A^\star A$ is positive.

To prove density of $\Dom(A^\star A)$ in $\Dom(A)$ with respect to the graph norm, suppose that $x\in \Dom(A)$ is such that
$\iprod{x}{y}_{\Dom(A)} =0$ for all $y\in \Dom(A^\star A)$, where
$\iprod{x}{y}_{\Dom(A)}:= \iprod{x}{y}+\iprod{Ax}{Ay}$ is the inner product of $\Dom(A)$, viewed as a Hilbert space
with respect to this inner product (completeness being a consequence of the closedness of $A$; see Proposition \ref{prop:closed-op}). Then, for all $y\in \Dom(A^\star A)$,
$$ 0= \iprod{x}{y}+\iprod{x}{A^\star Ay} =  \iprod{x}{(I + A^\star A)y}.$$
Since $I+A^\star A$ is surjective, it follows that
$\iprod{x}{z}=0$ for all $z\in H$, so $x=0$.

In particular, $A^\star A$ is densely defined. Selfadjointness now follows from Corollary \ref{cor:sa-c0sgr}.
\end{proof}

We conclude this section with another useful criterion for selfadjointness.

\begin{theorem}\label{thm:sa-criterion}
For a densely defined symmetric operator $A$ in $H$ the following assertions are equivalent:
\begin{enumerate}[label={\rm(\arabic*)}, leftmargin=*]
 \item\label{it:sa-criterion1} $A$ is selfadjoint;
 \item\label{it:sa-criterion2} $A$ is closed and $\Ker(A^\star+i) = \Ker(A^\star-i) =\{0\}$;
 \item\label{it:sa-criterion3} $\Ran(A+i) = \Ran(A-i) = H$.
\end{enumerate}
\end{theorem}
\begin{proof}
\ref{it:sa-criterion1}$\Rightarrow$\ref{it:sa-criterion2}:\
If $A$ is selfadjoint, then $A = A^\star$ is closed by Proposition \ref{prop:graph-As}. If $x\in \Dom(A^\star)$ satisfies $(A^\star + i)x=0$, then $x\in \Dom(A)$ and $Ax = A^\star x = -ix$, and
$$-i\iprod{x}{x} = \iprod{Ax}{x} = \iprod{x}{A^\star x} = i\iprod{x}{x}$$ implies $x=0$.
Similarly, $(A^\star - i)x=0$ implies $x=0$.

\smallskip
\ref{it:sa-criterion2}$\Rightarrow$\ref{it:sa-criterion3}:\ By the same argument as in the proof of Proposition \ref{prop:injective-denserange},
the injectivity of $(A\pm i)^\star = A^\star \mp i$ implies that $A\pm i$ has dense range (and conversely; this will be used in the proof of the next implication).
By the same argument as in the proof of Theorem \ref{thm:spect-sa},
the symmetry of $A$ implies $\n (A\pm i)x\n\ge \n x\n$ for all $x\in \Dom(A)$, and since $A$ is closed, Proposition \ref{prop:closed-range-unbdd}
implies that the ranges of $A\pm i$ are closed. We conclude that both ranges equal $H$.

\smallskip
\ref{it:sa-criterion3}$\Rightarrow$\ref{it:sa-criterion1}:\
Fix $h\in\Dom(A^\star)$.
The assumption $\Ran(A-i) = H$ implies that there exists an $h'\in \Dom(A)$ such that $(A-i)h' = (A^\star-i)h$. Since $A^\star$ extends $A$, we have $h'\in \Dom(A^\star)$ and $A^\star h' = Ah'$\!.
It follows that $(A^\star - i)h' = (A^\star-i)h$. As was noted in the proof of the previous implication,
the assumption $\Ran(A+i) = H$ implies that
$A^\star - i$ is injective. It follows that $h = h'$. Since $h'\in \Dom(A)$, this implies that $h\in \Dom(A)$
 and $Ah = A^\star h$, the latter since $A^\star$ extends $A$.

This shows that $A$ extends $A^\star$\!. Since
 $A^\star$ extends $A$, these operators are equal.
\end{proof}

The theory of selfadjoint operators is taken up again
in Section \ref{subsec:symm-pos-forms} in connection with the theory of forms.

\section{Unbounded Normal Operators}

Having dealt with unbounded selfadjoint operators, we now turn to unbounded normal operators. We fix a complex Hilbert space $H$.

\subsection{Definition and General Properties}

\begin{definition}[Normal operators]\index{normal!operator}
A linear operator $A$ in $H$ is said to be {\em normal} if
it is densely defined, closed, and satisfies $$ A^\star A = AA^\star\!.$$
\end{definition}

The equality $A^\star A = AA^\star$ is shorthand for equality of the domains
\begin{align*}
 \Dom(A^\star A) & :=  \{x\in \Dom(A):\, Ax\in \Dom(A^\star)\}, \\
 \Dom(AA^\star) & :=  \{x\in \Dom(A^\star):\, A^\star x\in \Dom(A)\},
\end{align*}
along with equality  $$A^\star A x= AA^\star x$$ for all $x$ in this common domain.

Since $A^{\star\star} = A$ by Proposition \ref{prop:adj-dd},
a densely defined closed operator $A$ is normal if and only if its adjoint $A^\star$ is normal.

  \begin{proposition}\label{prop:normal-ext}
   If $A$ is a normal operator, then:
   \begin{enumerate}[label={\rm(\arabic*)}, leftmargin=*]
    \item\label{it:normal-ext1} $\Dom(A) = \Dom(A^\star)$;
    \item\label{it:normal-ext2} $\n Ax \n = \n A^\star x\n $ for all $x\in \Dom(A) = \Dom(A^\star)$;
    \item\label{it:normal-ext3} if $A\subseteq B$ with $B$ normal, then $A=B$.
   \end{enumerate}
  \end{proposition}
  \begin{proof}
\ref{it:normal-ext1} and \ref{it:normal-ext2}: \   The normality of $A$ implies that if $x\in \Dom(A^\star A)= \Dom(AA^\star)$, then
$x\in \Dom(A)$, $x\in \Dom(A^\star)$, and
\begin{align*}\n Ax\n^2 =  \iprod{A^\star A x}{x} = \iprod{AA^\star x}{x} = \n A^\star x\n^2\!.
\end{align*}
By Theorem \ref{thm:AstarA},
$\Dom(A^\star A)$ is dense in $\Dom(A)$, so for any $x\in \Dom(A)$ we may choose a
sequence $x_n\to x$ with each $x_n\in \Dom(A^\star A)=\Dom(AA^\star)$ and with convergence in the graph norm of $\Dom(A)$.
Then $x_m,x_n\in \Dom(A^\star)$, and from $$\lim_{n,m\to\infty} \n A^\star x_n - A^\star x_m\n = \lim_{n,m\to\infty}\n A x_n - A x_m\n =0$$
we infer that $A^\star x_n\to y$ for some $y\in H$. From the closedness of $A^\star$ (Proposition \ref{prop:graph-As})
we infer that $x\in \Dom(A^\star)$ and $A^\star x = y$. This argument shows that
$\Dom(A) \subseteq \Dom(A^\star)$.

Since $A^\star$ is normal, what we just proved can be applied to $A^\star$\!. This, together with Proposition \ref{prop:adj-dd}, gives the reverse inclusion
$\Dom(A^\star) \subseteq \Dom(A^{\star\star}) = \Dom(A)$.

For $x\in\Dom(A) = \Dom(A^\star)$, the sequence $(x_n)_{n\ge1}$ chosen above
satisfies $Ax_n\to Ax$ and $A^\star x_n\to A^\star x$.
Hence
$$
\n A^\star x\n
= \limn\n A^\star x_n\n
= \limn\n Ax_n\n
= \n Ax\n.
$$
This proves \ref{it:normal-ext2}.

\smallskip
\ref{it:normal-ext3}: \   If $A\subseteq B$ with $A$ and $B$ normal, then by \ref{it:normal-ext1}, Proposition \ref{prop:AB-properties}\ref{it:AB-properties1}, and another application of \ref{it:normal-ext1},
$$ \Dom(B) = \Dom(B^\star)\subseteq \Dom(A^\star) = \Dom(A).$$
Together with the assumption $ \Dom(A)\subseteq  \Dom(B)$, this implies $ \Dom(A)= \Dom(B)$.
\end{proof}

In general, the spectrum of an unbounded closed operator may be empty; an example is given in Problem \ref{prob:empty-spectrum}.
We finish this subsection with the elementary observation that every normal operator has nonempty spectrum:

\begin{proposition}\label{prop:sp-normal-nonempty}
Let $A$ be a normal operator in a nonzero Hilbert space $H$. Then the spectrum $\sigma(A)$ is nonempty.
\end{proposition}

\begin{proof}
If $A$ is not invertible, then $0\in\sigma(A)$ and we are done. Therefore we may assume that $A$ is invertible.
Then $A^\star$ is invertible and
$(A^{-1})^\star=(A^\star)^{-1}$. Moreover,
\begin{align*}
(A^{-1})^\star A^{-1}
&=(A^\star)^{-1}A^{-1}
=(AA^\star)^{-1}
=(A^\star A)^{-1}
=A^{-1}(A^\star)^{-1}
=A^{-1}(A^{-1})^\star.
\end{align*}
Thus $A^{-1}$ is normal.

For $\lambda\neq 0$ we have
$A-\lambda I=A\bigl(I-\lambda A^{-1}\bigr)$,
so $A-\lambda I$ is invertible if and only if $I-\lambda A^{-1}$ is invertible. Stated differently, for $\la\not=0$ we have
\begin{equation*}
\lambda\in\sigma(A)\ \Longleftrightarrow\  \frac{1}{\lambda}\in\sigma(A^{-1}).
\end{equation*}

Suppose that $\sigma(A)=\varnothing$. Then the above equivalence implies $\sigma(A^{-1})\subseteq\{0\}$.
Since $A^{-1}$ is bounded, its spectrum is nonempty, and we conclude that $\sigma(A^{-1})=\{0\}$.
But then $\n A^{-1}\n = r(A^{-1}) = 0$ by Proposition \ref{prop:normal-spectralradius}. Since $A^{-1} = 0$ is a bijection of $H$ onto $\Dom(A)$, this forces $\Dom(A) = \{0\}$, and since $A$ is densely defined, it follows that $H=\{0\}$.
\end{proof}

\subsection{The Measurable Functional Calculus}\label{subsec:meas-calc-unbdd}

Projection-valued measures give rise to normal operators:

\begin{theorem}[Measurable functional calculus]\label{thm:Borel-FC-unbddnormal}\index{functional calculus!measurable}
Let $(\Om,\calF)$ be a measurable space, let $P:\calF\to \calL(H)$ be a projection-valued measure,  and let $f:{\Om}\to \C$ be a measurable function. There exists a unique normal operator $\Phi(f)$ in $H$ satisfying
 \begin{align*} \Dom(\Phi(f)) & = \Bigl\{x\in H:\, \int_{\Om} |f|^2\ud P_x <\infty \Bigr\}, \\
 \iprod{\Phi(f) x}{x}&  =  \int_{\Om} f\ud P_x, \quad x\in \Dom(\Phi(f)).
 \end{align*}
For all $x\in \Dom(\Phi(f))$ we have
\begin{align*}
\n \Phi(f) x\n^2 = \int_{\Om} |f|^2\ud P_x.
\end{align*}
Furthermore, if $f_n, f, g:\Om\to\C$ are measurable functions and $a,b \in \C$, then:
\begin{enumerate}[label={\rm(\arabic*)}, leftmargin=*]
\item\label{it:Borel-FC-unbddnormal-add}
$\Dom(\Phi(f))\cap\Dom(\Phi(g))\subseteq \Dom(\Phi(a f+b g))$
and
$$
   a\,\Phi(f)x+b\,\Phi(g)x = \Phi(a f+b g)x,
   \quad x\in \Dom(\Phi(f))\cap\Dom(\Phi(g)).
$$
\item\label{it:Borel-FC-unbddnormal1} $\Dom(\Phi(f)\Phi(g)) = \Dom(\Phi(fg))\cap \Dom(\Phi(g))$ and
$$\Phi(f)\Phi(g)x = \Phi(fg)x, \quad x\in \Dom(\Phi(f)\Phi(g));$$
\item\label{it:Borel-FC-unbddnormal2}
$\Dom(\Phi(f))=\Dom(\Phi(\ov f))$ and
$\Phi(f)^\star = \Phi(\ov f)$;
\item\label{it:Borel-FC-unbddnormal3} if
$0\le |f_n|\le |f|$ and $\limn f_n = f$ pointwise on $\Om$, then $\Dom(\Phi(f))\subseteq \Dom(\Phi(f_n))$ and
$$\limn \Phi(f_n)x = \Phi(f)x, \quad x\in \Dom(\Phi(f)).$$
\end{enumerate}
The operator $\Phi(f)$ is selfadjoint if and only if $f$ is real-valued $P_x$-almost everywhere for all $x\in H$.
\end{theorem}

It follows from \ref{it:Borel-FC-unbddnormal1} that
\begin{align}\label{eq:domPhifg}
\Phi(f)\Phi(g)= \Phi(fg) \ \Leftrightarrow \ \Dom(\Phi(fg))\subseteq \Dom(\Phi(g)).
\end{align}
This is trivially the case if $g$ is bounded, for then $\Dom(\Phi(g)) = H$. In that case
$\Phi(g)$ is bounded and equals the operator given by the bounded calculus of Theorem \ref{thm:Borel-FC}. This fact, and the properties of
the bounded calculus, will be frequently used below. For examples of strict operator inclusion, see
Problem \ref{prob:meas-calc-domains}.

\begin{example}\label{normal-powers} Under the above assumptions, it follows from \eqref{eq:domPhifg} that
 $$ \Phi(f^n) = (\Phi(f))^n, \quad n\ge 1.$$
To prove this by induction, it suffices to check that $ \Phi(f^{k+1}) = \Phi(f^k)\Phi(f)$ for all $k\ge 1$.
By  \eqref{eq:domPhifg}, this operator identity holds if and only if $$\Dom(\Phi(f^{k+1})) \subseteq \Dom(\Phi(f)).$$ If $x\in \Dom(\Phi(f^{k+1})) $, then $$\int_{\Om} |f|^{2k+2}\ud P_x <\infty.$$ Since $P_x$ is a finite measure, this implies $\int_{\Om} |f|^{2}\ud P_x <\infty$, that is,
$x\in \Dom(\Phi(f))$.
\end{example}

\begin{proof}[Proof of Theorem \ref{thm:Borel-FC-unbddnormal}]
For bounded measurable functions, let $\Phi$ denote the bounded Borel calculus
of Theorem~\ref{thm:Borel-FC}. For a measurable function
$f:\Om\to\C$ and all $n\ge 1$, put
$$
E_n:=\{|f|\le n\},
\qquad
f_n:=\one_{E_n}f,
$$
and set
$$
\Dom_f:=\Big\{x\in H:\,\int_\Om |f|^2\ud P_x<\infty\Big\}.
$$

For $x\in H$ and $m,n\ge1$, the bounded
functional calculus gives
$$
\n\Phi(f_n)x-\Phi(f_m)x\n^2
=
\int_\Om|f_n-f_m|^2\ud P_x.
$$
It follows that
$x\in\Dom_f$ if and only if the sequence $(\Phi(f_n)x)_{n\ge1}$ converges in $H$. Indeed, the forward implication follows by dominated convergence. Conversely, if the sequence $(\Phi(f_n)x)_{n\ge1}$ converges, then
 $$
 \sup_{n\ge 1}\int_\Om|f_n|^2\ud P_x= \sup_{n\ge 1}\n\Phi(f_n)x\n^2 <\infty
 $$
and monotone convergence gives $\int_\Om|f|^2\ud P_x<\infty$.

The domain characterisation implies that $\Dom_f$ is a linear subspace of $H$, and we may define a linear operator $\Phi(f):\Dom_f\to H$ by
$$
\Phi(f)x:=\lim_{n\to\infty}\Phi(f_n)x.
$$
We next show that $\Dom_f$ is dense. Indeed, for $x\in H$ we have
$x_n := P_{E_n}x\in\Dom_f$, since
$$
P_{x_n}(B)=\iprod{P_{B}x_n}{x_n}  = \iprod{P_B P_{E_n}x}{P_{E_n}x} = \iprod{P_{B\cap E_n}x}{x} =  P_x(B\cap E_n),\quad B\in\calF\!,
$$
and therefore
$$
\int_\Om|f|^2\ud P_{x_n}
=
\int_{E_n}|f|^2\ud P_x
\le n^2\n x\n^2\!.
$$
Moreover, $x_n = P_{E_n}x\to x$ by the strong countable additivity of $P$.

We record the following useful approximation observation. If $x\in\Dom_f$
and $u:\Om\to\C$ is bounded and measurable, then
\begin{equation}\label{eq:meas-FC-approx}
\n\Phi(u)x-\Phi(f)x\n^2
=
\int_\Om|u-f|^2\ud P_x.
\end{equation}
Indeed, for every $n$ the bounded functional calculus gives
$$
\n\Phi(u)x-\Phi(f_n)x\n^2
=
\int_\Om|u-f_n|^2\ud P_x,
$$
and letting $n\to\infty$ proves \eqref{eq:meas-FC-approx}. In particular,
if $(u_n)_{n\ge1}$ is a sequence of bounded measurable functions converging
to $f$ in $L^2(\Om,P_x)$, then
\begin{equation}\label{eq:meas-FC-L2-approx}
\Phi(u_n)x\to\Phi(f)x.
\end{equation}

Taking limits in the corresponding identities for the bounded functions
$f_n$, we obtain, for all $x\in\Dom(\Phi(f))$,
$$
\n\Phi(f)x\n^2=\int_\Om|f|^2\ud P_x.
$$
Since $P_x$ is finite and $f_n\to f$ in $L^2(\Om,P_x)$, we also have
$f_n\to f$ in $L^1(\Om,P_x)$, and hence
$$
\iprod{\Phi(f)x}{x}
=
\lim_{n\to\infty}\int_\Om f_n\ud P_x
=
\int_\Om f\ud P_x.
$$

To prove \ref{it:Borel-FC-unbddnormal-add}, let
$x\in\Dom(\Phi(f))\cap\Dom(\Phi(g))$, and let
$$
f_n:=\one_{\{|f|\le n\}}f,
\qquad
g_n:=\one_{\{|g|\le n\}}g.
$$
Then $af_n+bg_n\to af+bg$ in $L^2(\Om,P_x)$. Therefore
$x\in\Dom(\Phi(af+bg))$, and
\eqref{eq:meas-FC-L2-approx}, together with the linearity of the bounded
functional calculus, gives
$$
\Phi(af+bg)x
=
\lim_{n\to\infty}\Phi(af_n+bg_n)x
=
a\,\Phi(f)x+b\,\Phi(g)x.
$$

We next prove \ref{it:Borel-FC-unbddnormal3}. If
$0\le |f_n|\le |f|$ and $f_n\to f$ pointwise, then
$\Dom(\Phi(f))\subseteq\Dom(\Phi(f_n))$ by the definition of the domains. Moreover, for $x\in\Dom(\Phi(f))$, dominated convergence gives
$$
\n\Phi(f)x-\Phi(f_n)x\n^2
=
\int_\Om|f-f_n|^2\ud P_x\to 0.
$$

We now prove \ref{it:Borel-FC-unbddnormal1}. We first note that if $b$ is
bounded and $g$ is measurable, then
\begin{equation}\label{eq:meas-FC-bdd-mult}
\Phi(b)\Phi(g)x=\Phi(bg)x, \quad x\in\Dom(\Phi(g)).
\end{equation}
Indeed, for $g_n:=\one_{\{|g|\le n\}}g$, the bounded functional calculus
gives
$$
\Phi(b)\Phi(g_n)x=\Phi(bg_n)x.
$$
Letting $n\to\infty$ and using \eqref{eq:meas-FC-L2-approx}, we obtain
\eqref{eq:meas-FC-bdd-mult}.

Consequently, for bounded measurable $b$ and
$x\in\Dom(\Phi(g))$,
$$
\int_\Om|b|^2\ud P_{\Phi(g)x}
=
\n\Phi(b)\Phi(g)x\n^2
=
\n\Phi(bg)x\n^2
=
\int_\Om|b|^2|g|^2\ud P_x.
$$
Applying this identity to
$b_n:=\one_{\{|f|\le n\}}f$ and using monotone convergence, we obtain
\begin{equation*}
\int_\Om|f|^2\ud P_{\Phi(g)x}
= \int_\Om|fg|^2\ud P_x, \quad x\in\Dom(\Phi(g)),
\end{equation*}
where either side may be infinite. It follows that
$$
\Dom(\Phi(f)\Phi(g))
=
\Dom(\Phi(fg))\cap\Dom(\Phi(g)).
$$

Put
$f_n:=\one_{\{|f|\le n\}}f$.
For $x$ in this common domain, by
\ref{it:Borel-FC-unbddnormal3} we have
$\Phi(f_n)\Phi(g)x\to\Phi(f)\Phi(g)x$ and $\Phi(f_ng)x\to\Phi(fg)x$, and by \eqref{eq:meas-FC-bdd-mult} we have
$\Phi(f_n)\Phi(g)x=\Phi(f_ng)x$. Therefore
$$
\Phi(f)\Phi(g)x=\Phi(fg)x.
$$

We turn to \ref{it:Borel-FC-unbddnormal2}. Clearly,
$\Dom(\Phi(f))=\Dom(\Phi(\ov f))$.
For $x,y$ in this common domain, let
$f_n:=\one_{\{|f|\le n\}}f$. By the bounded functional calculus we have
$\iprod{\Phi(f_n)x}{y} = \iprod{x}{\Phi(\ov{f_n})y} $, and
letting $n\to\infty$ we obtain
$\iprod{\Phi(f)x}{y} = \iprod{x}{\Phi(\ov f)y}$.
Thus $\Phi(\ov f)\subseteq\Phi(f)^\star$.

For the reverse inclusion, let
$y\in\Dom(\Phi(f)^\star)$ and put $z:=\Phi(f)^\star y$.
For $n\ge1$, set
$f_n:=\one_{E_n}f$ with
$E_n:=\{|f|\le n\}$.
The product rule gives
$\Phi(f)P_{E_n}=\Phi(f_n)$ on $H$.
Hence, by the bounded functional calculus, for every $x\in H$
we have
\begin{align*}
\iprod{x}{\Phi(\ov{f_n})y}
=\iprod{\Phi(f_n)x}{y} =\iprod{\Phi(f)P_{E_n}x}{y}
=\iprod{P_{E_n}x}{z} =\iprod{x}{P_{E_n}z}.
\end{align*}
Consequently,
$$
\Phi(\ov{f_n})y=P_{E_n}z\to z,
$$
since $E_n\uparrow\Om$.
By the domain characterisation at the beginning of the proof,
this implies that $y\in\Dom(\Phi(\ov f))$ and
$\Phi(\ov f)y=z$.
This proves the reverse inclusion $\Phi(f)^\star\subseteq\Phi(\ov f)$. We conclude that
$$
\Phi(f)^\star=\Phi(\ov f).
$$
Applying this identity to $\ov f$, we find that $\Phi(f)$, being an adjoint operator, is closed.

Moreover,
$\Dom(\Phi(|f|^2)) \subseteq \Dom(\Phi(f)) =
\Dom(\Phi(\ov f))$, the inclusion being a consequence of the finiteness of $P_x$.
Part \ref{it:Borel-FC-unbddnormal1} therefore gives
$$
\Phi(f)^\star\Phi(f)
=
\Phi(\ov f)\Phi(f)
=
\Phi(|f|^2)
=
\Phi(f)\Phi(\ov f)
=
\Phi(f)\Phi(f)^\star.
$$
Hence $\Phi(f)$ is normal.

If $N$ is a normal operator with domain $\Dom_f$ and
$$
\iprod{Nx}{x}=\int_\Om f\ud P_x,
\quad x\in\Dom_f,
$$
then
$$
\iprod{Nx}{x}=\iprod{\Phi(f)x}{x},
\quad x\in\Dom_f.
$$
Polarisation gives $N=\Phi(f)$, proving uniqueness.

Finally, suppose first that $f$ is real-valued $P_x$-almost everywhere for
all $x\in H$. Then $\Phi(f)=\Phi(\ov f)$, and hence
$\Phi(f)=\Phi(f)^\star$.

Conversely, suppose that $\Phi(f)$ is selfadjoint. For
$x\in\Dom(\Phi(f))$, parts
\ref{it:Borel-FC-unbddnormal-add} and
\ref{it:Borel-FC-unbddnormal2} give
$$
0
=
\n\Phi(f-\ov f)x\n^2
=
\int_\Om|f-\ov f|^2\ud P_x.
$$
Now let $x\in H$ be given, and put
$E_n:=\{|f|\le n\}$. Since $P_{E_n}x\in\Dom(\Phi(f))$,
$$
0
=
\int_\Om|f-\ov f|^2\ud P_{P_{E_n}x}
=
\int_{E_n}|f-\ov f|^2\ud P_x.
$$
Letting $n\to\infty$, we conclude that
$f=\ov f$ $P_x$-almost everywhere. This completes the proof.
\end{proof}

The following substitution rule extends Proposition \ref{prop:projvalmeas-subst}.

\begin{proposition}[Substitution]\label{prop:normal-subst} Let $(\Om,\calF)$ and $(\Om'\!,\calF')$ be measurable spaces and let
$f:{\Om}\to \Om'$ be a measurable mapping. If $P:\calF\to \calL(H)$ is a projection-valued measure, and if $\Phi$ and $\Psi$ denote the measurable functional calculi of $P$ and $f_\#(P)$, for all measurable functions
$g:\Om'\to \C$ we have
$$\Phi(g\circ f) = \Psi(g)$$
with equality of domains.
\end{proposition}
\begin{proof}
The proof is the same as that of Proposition \ref{prop:projvalmeas-subst}, except that some domain issues have to be taken care of.
Following this proof, for all nonnegative measurable functions $g$ on $\Om'$ and $x\in H$ we obtain, with $Q:= f_\#(P)$,
\begin{align}\label{eq:normal-subst} \int_{\Om} g\circ f\ud P_x = \int_{\Om'} g\ud Q_x,
\end{align}
the finiteness of one of these integrals implying the finiteness of the other. Applying this with $g$ replaced by $|g|^2$\!, we see that
$x\in \Dom(\Phi(g\circ f))$ if and only if $x\in \Dom(\Psi(g))$.  On this common domain, the Cauchy--Schwarz inequality gives integrability of $g\circ f$ and $g$. Applying \eqref{eq:normal-subst} to the positive and negative parts of $\Re g$ and $\Im g$, we obtain
$$\iprod{\Phi(g\circ f)x}{x} = \iprod{\Psi(g)x}{x}.$$
By polarisation, for all $x,y$ in this common domain,
$$
\iprod{\Phi(g\circ f)x}{y}
=
\iprod{\Psi(g)x}{y}.
$$
This common domain is dense in $H$, since it is the domain of the normal
operator $\Psi(g)$, and therefore
$$
\Phi(g\circ f)x=\Psi(g)x
$$
for all $x$ in the common domain.
\end{proof}

Let $P$ be a projection-valued measure on $\calB(\C)$.
As in Section \ref{sec:Borel-FC}, its support is the smallest
closed set $F\subseteq\C$ such that $P_{\C\setminus F}=0$.
This set exists by a standard argument using that every open subset of $\C$ is the union of a countable family of open balls with rational centres and radii.

\begin{proposition}\label{prop:supp-P-sigmaA}
Let $P$ be a projection-valued measure on $\calB(\C)$ and consider the normal operator
$$ A := \int_\C z\ud P(z)$$
with domain $\Dom(A) = \{x\in H:\, \int_\C |z|^2\ud P_x(z) <\infty\}$.
The spectrum of this operator is the support of $P$,
$$\si(A) = \operatorname{supp}(P).$$
Moreover,
if $B\in\calB(\C)$ is such that $\sigma(A)\cap B= \emptyset$, then $P_B=0$.
\end{proposition}

\begin{proof}
For $H = \{0\}$ the result is trivial, so we may assume that $H$ is nonzero.

The inclusion $\sigma(A)\subseteq\operatorname{supp}(P)$ and the second assertion are proved as in Proposition \ref{prop:TP}.

The proof of the inclusion $\operatorname{supp}(P)\subseteq\sigma(A)$ also follows the pattern of Proposition \ref{prop:TP}, except that some domain issues have to be taken care of. Let $z\in\varrho(A)$. Choose $r_0>0$ such
that
$\overline{B(z;r_0)}\subseteq\varrho(A)$,
and put
$r:=\min\{r_0,\,1/(2\|R(z,A)\|)\}$.
We claim that $P_{B(z;r)}=0$. Suppose otherwise, and choose
$x\in\Ran(P_{B(z;r)})$ with $\|x\|=1$. Then $P_x$ is supported on
$B(z;r)$, and since $P_x$ is a probability measure,
$$
        \int_\C |\lambda|^2\ud P_x(\lambda)
        = \int_{B(z;r)}|\lambda|^2\ud P_x(\lambda)
        \le (|z|+r)^2<\infty,
$$
so $x\in\Dom(A)$. Therefore
\begin{align*}
        \|(z-A)x\|^2 = \int_\C |z-\lambda|^2\ud P_x(\lambda)
        = \int_{B(z;r)} |z-\lambda|^2\ud P_x(\lambda)
        \le r^2 .
\end{align*}
Since $x=R(z,A)(z-A)x$, it follows that
$$1=\|x\| \le \|R(z,A)\|\,\|(z-A)x\| \le \|R(z,A)\|\,r \le \frac12,$$
a contradiction. Hence $P_{B(z;r)}=0$. It follows that $\operatorname{supp}(P)$ is contained in the complement of $B(z;r)$. In particular, $z\not \in \operatorname{supp}(P)$.
\end{proof}

\section{The Spectral Theorem for Unbounded Normal Operators}

The proof of the spectral theorem for unbounded normal operators proceeds by a
reduction to the bounded case. The basic idea is to exploit the fact that
the mapping
\begin{align*}\zeta: z\mapsto \frac{z}{(1+|z|^2)^{1/2}}\end{align*}
maps the complex plane bijectively onto the open unit disc $\mathbb{D}$.
This suggests that if $A$ is a normal operator, then $$Z_A := A(I+A^\star A)^{-1/2}$$ is a normal contraction on $H$.
This is indeed the case, as will be proved in Proposition \ref{prop:normal-TA}. It follows that $\si(Z_A)\subseteq \ov{\mathbb{D}}$. By the spectral theorem for bounded normal operators,
there exists a projection-valued measure $Q$ on $\ov{\mathbb{D}}$ such that
$$ Z_A = \int_{\ov{\mathbb{D}}} \la \ud Q(\la).$$
We now define a projection-valued measure $P$ on $\C$ by setting $P_B := Q_{\zeta(B)}$ for Borel sets $B\subseteq\C$, and use Proposition \ref{prop:normal-subst}
to show that
$$ A = \int_{\C} \la \ud P(\la).$$
Similarly, the uniqueness of $P$
is reduced to the uniqueness of $Q$.

Some technical details need to be addressed to turn this simple idea into a rigorous proof: one has
to deal with subtle domain issues and with the fact that $\zeta$ maps $\C$ onto the open unit disc, whereas $Q$ is supported on the closed unit disc.

\smallskip
We start with the proof that $Z_A$ is well defined as a contractive normal operator on $H$. This is accomplished in
Proposition \ref{prop:normal-TA}, for which we need two lemmas.

\begin{lemma}\label{lem:Ta-approx}
 Let $A$ be a closed operator in $H$, let $T\in\calL(H)$, and let $f\in C(\si(T))$.

\begin{enumerate}[label={\rm(\arabic*)}, leftmargin=*]
 \item if $T$ is selfadjoint and $TA\subseteq AT$, then $f(T)A \subseteq Af(T)$;
 \item if $T$ is normal and $TA\subseteq AT$ and $T^\star A\subseteq AT^\star$, then $f(T)A \subseteq Af(T)$.
\end{enumerate}
\end{lemma}
\begin{proof}
We only prove the second assertion, the first being an immediate consequence.

We have
$T^2A = T(TA) \subseteq T(AT) = (TA)T \subseteq (AT)T = AT^2\!.$
Continuing by induction, we see that $T^kA\subseteq AT^k$ for all $k\in\N$.
Similarly, $(T^\star)^kA\subseteq A(T^\star)^k$ for all $k\in\N$.
These inclusions imply that $$p(T,T^\star)A\subseteq Ap(T,T^\star)$$ for all polynomials $p$
in the variables $z$ and $\ov z$; notation is as in Section \ref{subsec:contFCnorma}.
By the Stone--Weierstrass
theorem there exist polynomials $p_n$ in the variables $z$ and $\ov z$
such that $p_n(z,\ov z)\to f(z)$ uniformly for $z\in\sigma(T)$.
Then, by the properties of the continuous functional calculus for normal operators
(Theorem \ref{thm:normal-cont-fc}),
$$ \n p_n(T,T^\star)- f(T)\n = \sup_{z\in \sigma(T)} |p_n(z,\ov z)-f(z)| \to 0 \ \hbox{ as } \
n\to\infty.$$
The inclusion $p_n(T,T^\star)A\subseteq Ap_n(T,T^\star)$ implies that if $x\in \Dom(A)$, then $p_n(T,T^\star)x\in\Dom(A)$ and
$$\limn A p_n(T,T^\star)x = \limn p_n(T,T^\star) Ax = f(T)Ax.$$ Since also $\limn  p_n(T,T^\star) x = f(T)x$,
the closedness of $A$ implies that $f(T)x\in \Dom(A)$ and $Af(T)x = f(T)Ax$. This proves the result.
\end{proof}

If $A$ is a densely defined closed operator in $H$,
then $A^\star A$ is selfadjoint by Theorem \ref{thm:AstarA}, and by Proposition \ref{prop:sa-pos-spectrum} we have $\sigma(A^\star A)\subseteq [0,\infty)$.
This allows us to define
$$T_A := (I+A^\star A)^{-1}\!.$$
This operator is a positive contraction, and if $A$ is normal, we have $T_A =T_{A^\star}$.

\begin{lemma}\label{lem:normal-AAs}
 If $A$ is normal, then for all $x\in \Dom(A)$ we have
  $T_A x\in \Dom(A)$ and $$A T_A x = T_A  A x.$$
\end{lemma}

\begin{proof}
Let $x\in \Dom(A)$. Then $y :=  T_A x \in \Dom(A^\star A)\subseteq \Dom(A)$, $Ay = AT_A  x \in \Dom(A^\star)$, and
 $A^\star Ay = x- T_A x\in \Dom(A)$, so
 $Ay\in \Dom(AA^\star) = \Dom(A^\star A)$. Combining this with
$$(I+AA^\star )Ax = Ax+(AA^\star )Ax =  Ax+A(A^\star A)x =A(I+A^\star A)x, \quad x\in\Dom(AA^\star A),$$
we obtain
\[
AT_A x = [T_A(I+AA^\star)]AT_Ax
=T_A A[(I+A^\star A)T_A]x = T_A Ax.
\qedhere\]
\end{proof}

\begin{proposition}\label{prop:normal-TA}
If $A$ is a normal operator, then:
\begin{enumerate}[label={\rm(\arabic*)}, leftmargin=*]
 \item\label{it:normal-TA1} the range of $T_A^{1/2}$ is dense in $\Dom(A)$ with respect to the graph norm;
 \item\label{it:normal-TA2} the operator $Z_A:= A T_A^{1/2}$
 is contractive and we have $T_A = I -Z_A^\star Z_A$;
 \item\label{it:normal-TA4} $Z_A$ is normal, commutes with $T_A$, and satisfies $Z_A^\star = Z_{A^\star}$.
\end{enumerate}
\end{proposition}

\begin{proof}
\ref{it:normal-TA1} and \ref{it:normal-TA2}: \
Put $T:=T_A$. Since
$\Ran(T)=\Dom(A^\star A)\subseteq\Dom(A)$,
the operator $AT$ is defined on all of $H$.

Since $T$ is positive and injective, so is $T^{1/2}$. By selfadjointness and Proposition \ref{prop:injective-denserange}, the range of $T^{1/2}$ is dense in $H$. Define the linear operator $Z_0: \Ran(T^{1/2})\to H$ by $Z_0T^{1/2}x:=ATx$. This operator is well defined by the injectivity of $T^{1/2}$. For $y=T^{1/2}x$ we have
\begin{align*}
\n Z_0y\n^2
&= \n ATx\n^2
 = \iprod{A^\star ATx}{Tx}
 = \iprod{x-Tx}{Tx}  \\
&= \iprod{x}{Tx}-\iprod{Tx}{Tx}
 = \n y\n^2-\n T^{1/2}y\n^2
 \leq \n y\n^2\!.
\end{align*}
Thus $Z_0$ is contractive on the dense subspace $\Ran(T^{1/2})$, and therefore it has a unique extension to a contraction on $H$. We denote this extension by $Z_A$.

We claim that $Z_A=AT^{1/2}$. Let $y\in H$, and choose $x_n\in H$ such that $T^{1/2}x_n\to y$. Then $Tx_n=T^{1/2}(T^{1/2}x_n)\to T^{1/2}y$, while $ATx_n=Z_0T^{1/2}x_n\to Z_Ay$. Since $A$ is closed, it follows that $T^{1/2}y\in\Dom(A)$ and $AT^{1/2}y=Z_Ay$. Since $y\in H$ was arbitrary, this proves the claim.

For $y=T^{1/2}x\in\Ran(T^{1/2})$, the above computation gives
$$\iprod{Z_A^\star Z_Ay}{y} =
\n Z_Ay\n^2=\n y\n^2-\n T^{1/2}y\n^2=\iprod{(I-T)y}{y}.
$$
By polarisation, $\iprod{Z_A^\star Z_Ay}{z}=\iprod{(I-T)y}{z}$ for all $y,z\in\Ran(T^{1/2})$. Since $\Ran(T^{1/2})$ is dense in $H$ and both $Z_A^\star Z_A$ and $I-T$ are bounded operators on $H$, we obtain $Z_A^\star Z_A=I-T$, that is, $T_A=T = I-Z_A^\star Z_A$.

It remains to prove the density assertion in \ref{it:normal-TA1}. We have just shown that $\Ran(T_A^{1/2})\subseteq\Dom(A)$. Moreover, $\Ran(T_A)=\Dom(A^\star A)$ is dense in $\Dom(A)$ with respect to the graph norm by Theorem \ref{thm:AstarA}. Since $\Ran(T_A)\subseteq\Ran(T_A^{1/2})$, it follows that $\Ran(T_A^{1/2})$ is dense in $\Dom(A)$ with respect to the graph norm.

\smallskip
\ref{it:normal-TA4}: \
By Lemma \ref{lem:normal-AAs}, for $x\in\Dom(A)$ we have $T_Ax\in\Dom(A)$ and $AT_Ax=T_AAx$. Thus $T_AA\subseteq AT_A$, and Lemma \ref{lem:Ta-approx} gives
$$
T_A^{1/2}A\subseteq AT_A^{1/2}.
$$
Since $A$ is normal, $T_{A^\star}=T_A$. Applying \ref{it:normal-TA2} to $A^\star$, we see that the operator $Z_{A^\star}=A^\star T_A^{1/2}$ is bounded. Furthermore,
$$
Z_{A^\star}
= A^\star(T_A^{1/2})^\star
= (T_A^{1/2}A)^\star
\supseteq (AT_A^{1/2})^\star
= Z_A^\star,
$$
where the equality and the inclusion follow from Proposition \ref{prop:AB-properties}. Since both $Z_{A^\star}$ and $Z_A^\star$ are bounded operators on $H$, this inclusion is an equality. Hence $Z_{A^\star}=Z_A^\star$.

Applying \ref{it:normal-TA2} to $A$ and to $A^\star$, and using again that $T_{A^\star}=T_A$, we obtain
$$
I-Z_A^\star Z_A=T_A=T_{A^\star}
=I-Z_{A^\star}^\star Z_{A^\star}
=I-Z_AZ_A^\star.
$$
It follows that $Z_A^\star Z_A=Z_AZ_A^\star$, and therefore $Z_A$ is normal.
Commutation with $T_A$ is now immediate.
\end{proof}

We let $Q$ denote the projection-valued measure on $\sigma(Z_A)$ associated with $Z_A$. Since $Z_A$ is a contraction, we have $\sigma(Z_A)\subseteq\ov{\mathbb D}$. Thus $Q$ may be viewed as a projection-valued measure on $\ov{\mathbb D}$.

In what follows we let ${\rm id}: \C\to \C$\index{$I$@id} be the identity mapping,
${\rm id}(\la) := \la.$

\begin{lemma}\label{lem:Qnull-bdry} $Q_{\partial \mathbb D} = 0$.
\end{lemma}
\begin{proof}
The mapping
$\zeta: z\mapsto {z}/{(1+|z|^2)^{1/2}}$
maps the complex plane
bijectively onto the open unit disc $\mathbb{D}$, with inverse
$ \zeta^{-1}:w\mapsto {w}/{(1-|w|^2)^{1/2}}.$
Since $T_A=I-Z_A^\star Z_A=(1-|{\rm id}|^2)(Z_A)$,
the spectral theorem for bounded normal operators gives
$$
T_A=\int_{\ov{{\mathbb{D}}}} (1-|\mu|^2)\ud Q(\mu).
$$
Since $1-|\mu|^2=0$ if and only if $\mu\in \partial{\mathbb{D}}$, it follows from Proposition \ref{prop:ker-fT-zero-set} that $\ker(T_A)=\Ran(Q_{\partial{\mathbb{D}}})$.
Since $T_A=(I+A^\star A)^{-1}$ is injective, we have $\ker(T_A)=\{0\}$, and therefore $Q_{\partial{\mathbb{D}}}=0$.
\end{proof}

Now we are ready to state and prove the main result of this section.

\begin{theorem}[Spectral theorem for normal
operators]\index{theorem!spectral, for normal operators}\label{thm:ST-unboundedn-normal}
For every normal operator $A$ there exists a unique projection-valued measure $P$ on $\sigma(A)$ such that\index{projection-valued measure!of a normal operator}
$$ A = \int_{\sigma(A)} \la \ud P(\la).$$
\end{theorem}
\begin{proof} We proceed in two steps.

\smallskip
{\em Step 1} -- We first prove existence.
Let $Q$ be the projection-valued measure of $Z_A$, and let $\Psi$ denote
its measurable functional calculus.

Define, for $\mu\in\overline{\mathbb D}$,
$$
c(\mu):=(1-|\mu|^2)^{1/2},
\qquad
\eta(\mu):=
\begin{cases}
\mu(1-|\mu|^2)^{-1/2},& \mu\in\mathbb D,\\
0,& \mu\in\partial\mathbb D.
\end{cases}
$$
Since $c$ is bounded and nonnegative and
$c^2=1-|{\rm id}|^2$, the bounded functional calculus and
the uniqueness of positive square roots give
$\Psi(c)=T_A^{1/2}$.

By Lemma \ref{lem:Qnull-bdry}, we have
$\eta c={\rm id}$ $Q$-almost everywhere.
Since $c$ is bounded, the product rule of
Theorem \ref{thm:Borel-FC-unbddnormal} therefore gives,
with domain $H$,
$$
\Psi(\eta)T_A^{1/2}
=\Psi(\eta)\Psi(c)
=\Psi(\eta c)
=\Psi({\rm id})
=Z_A
=AT_A^{1/2}.
$$
Thus $\Psi(\eta)$ agrees with $A$ on $\Ran(T_A^{1/2})$.
By Proposition \ref{prop:normal-TA}\ref{it:normal-TA1},
this subspace is dense in $\Dom(A)$ with respect to the
graph norm. Since $\Psi(\eta)$ is closed, it follows that
$A\subseteq\Psi(\eta)$.
Both operators are normal, so Proposition
\ref{prop:normal-ext}\ref{it:normal-ext3} implies
$A=\Psi(\eta)$.

Now let
$$
P:=\eta_\#(Q).
$$
Since $\eta=\zeta^{-1}$ on $\mathbb D$ and
$Q_{\partial\mathbb D}=0$, for every $B\in\mathscr B(\C)$ we have
$$
P_B
=Q_{\eta^{-1}(B)}
=Q_{\eta^{-1}(B)\cap\mathbb D}
=Q_{\zeta(B)}.
$$
By Proposition~\ref{prop:normal-subst}, the measurable functional calculus
$\Phi$ associated with $P$ satisfies
$\Phi({\rm id})=\Psi(\eta)=A$,
and Proposition~\ref{prop:supp-P-sigmaA} then gives
$\supp(P)=\sigma(A)$. Thus, viewing the restriction of $P$ as a
projection-valued measure on $\sigma(A)$, we obtain
$$
A=\int_{\sigma(A)}\lambda\ud P(\lambda).
$$
This completes the existence proof.

\smallskip
{\em Step 2} -- Turning to uniqueness, suppose that $\wt P$ is a
projection-valued measure on $\si(A)$ such that
$$
A=\int_{\si(A)}\la\ud\wt P(\la),
$$
and let $\wt\Phi$ denote the associated measurable functional calculus.
Thus $\wt\Phi({\rm id})=A$ and
$\wt\Phi(\ov{\rm id})=A^\star$.

By the product rule,
$$
A^\star A
=\wt\Phi(\ov{\rm id})\wt\Phi({\rm id})
=\wt\Phi(|{\rm id}|^2)
$$
with common domain
$\Dom(\wt\Phi(|{\rm id}|^2))$.
In view of the inequalities
$
|{\rm id}|^4
\le
(1+|{\rm id}|^2)^2
\le
2(1+|{\rm id}|^4)
$, this common domain equals
$ \Dom(\wt\Phi(1+|{\rm id}|^2))$.
The addition rule therefore gives, with equality of domains,
$$
\wt\Phi(1+|{\rm id}|^2)=I+A^\star A.
$$

Put
$$
g:=(1+|{\rm id}|^2)^{-1}=(c\circ\zeta)^2\!,
$$
the equality being evident by direct computation. Since $(1+|{\rm id}|^2)g=\one$, another application of the product rule gives
$$
(I+A^\star A)\wt\Phi(g)
=\wt\Phi(1+|{\rm id}|^2)\wt\Phi(g)
=I.
$$
Consequently,
$$
\wt\Phi(g)=(I+A^\star A)^{-1}=T_A.
$$
The function $c\circ\zeta$ is bounded and nonnegative, and hence the
uniqueness of positive square roots gives
$\wt\Phi(c\circ\zeta)=T_A^{1/2}$.
Since $\zeta={\rm id}\cdot(c\circ\zeta)$, the product rule now gives
$$
\wt\Phi(\zeta)
=\wt\Phi({\rm id})\wt\Phi(c\circ\zeta)
=AT_A^{1/2}
=Z_A.
$$

Define
$$
\wt Q:=\zeta_\#(\wt P).
$$
By Proposition~\ref{prop:normal-subst},
$$
\int_{\ov{\mathbb D}}\mu\ud\wt Q(\mu)
=\wt\Phi(\zeta)
=Z_A.
$$
Thus $Q$ and $\wt Q$ are projection-valued measures on the
compact set $\overline{\mathbb D}$ with the same integral
of the coordinate function. By Proposition
\ref{prop:spres-uniq}, we have $\wt Q=Q$.

Finally, since $\eta\circ\zeta={\rm id}$ on $\C$, this implies
$$
\wt P
=\eta_\#(\zeta_\#(\wt P))
=\eta_\#(\wt Q)
=\eta_\#(Q)
=P.
$$
This proves uniqueness.
\end{proof}

For normal operators $A$ and measurable functions $f:\si(A)\to \C$, the operator $\Phi(f)$ defined in terms of the projection-valued measure $P$ of $A$ by the calculus of Theorem \ref{thm:Borel-FC-unbddnormal} will be denoted by $f(A)$:
$$ f(A):= \Phi(f) = \int_{\si(A)} f\ud P.$$
As for bounded normal operators in Theorem \ref{thm:bff-FC-normal}, the properties of the measurable calculus
$\Phi$ translate into corresponding properties for the mapping $f\mapsto f(A)$. The result of Example \ref{normal-powers} says that for any normal operator $A$ and measurable function $f:\si(A)\to \C$,
$$ (f^n)(A) = (f(A))^n, \quad n\ge 1.$$

Recall from Section \ref{sec:Borel-FC} that if $P$ is a projection-valued measure on a measurable space $(\Omega,\calF)$ and
$f:\Omega\to\C$ is measurable, then the image of $P$ under $f$ is the projec\-tion-valued measure $f_\#(P)$ on $\C$ defined by
$$
(f_\#(P))_B:=P_{f^{-1}(B)},\quad B\in\calB(\C).
$$

Theorem~\ref{thm:SMT-normal-bdd} extends as follows.

\begin{theorem}[Spectral mapping theorem]\label{thm:SMT-normal-unbdd}\index{theorem!spectral mapping, normal operators}
Let $A$ be normal with projection-valued measure $P$, and let $f:\sigma(A)\to \C$ be
measurable. Then $$\sigma(f(A)) = \supp(f_\#(P)).$$
As a consequence, $\sigma(f(A))\subseteq \ov{f(\sigma(A))}.$
If $f$ is continuous, then $$\sigma(f(A)) = \ov{f(\sigma(A))}.$$
\end{theorem}

\begin{proof}
Put $Q:=f_\#(P)$.
By Propositions \ref{prop:normal-subst} and
\ref{prop:supp-P-sigmaA},
$$
f(A)=\int_\C z\ud Q(z)
$$
and $\si(f(A))=\supp(Q)$.
Since $Q_{\C\setminus\ov{f(\si(A))}}=0$, we also have
$\supp(Q)\subseteq\ov{f(\si(A))}$.

If $f$ is continuous and $\mu\in \si(A)$ is such that $f(\mu)\not\in \si(f(A)) = \supp(f_\#(P))$, then for $\eps>0$ small enough the relatively open set
$$N:=\{\la\in\si(A):\, |f(\la)-f(\mu)|<\eps\}$$ satisfies $P_N = 0$.
The continuity of $f$ implies that $B(\mu;\delta)\cap \si(A)\subseteq N$ for some  $\delta>0$. But then $\mu\not\in \supp(P) = \sigma({\rm id}(A)) = \si(A)$. This proves the inclusion
$f(\si(A))\subseteq\si(f(A))$, and hence, since $\si(f(A))$ is closed,
$\ov{f(\si(A))}\subseteq\si(f(A))$.
\end{proof}

The following result extends Proposition \ref{prop:ker-fT-zero-set} to the unbounded case, and the proof extends {\em mutatis mutandis}.

\begin{proposition}\label{prop:normal-eigen}
Let $A$ be a normal operator, and let $P$ be its projection-valued measure. For every Borel function $f:\si(A)\to \C$ we have
$$
\ker(f(A))=\Ran(P_{f^{-1}(\{0\})}).
$$
\end{proposition}

\begin{corollary}\label{cor:normal-eigen-mapping} Let $A$ be a normal operator and let $f:\sigma(A)\to \C$ be a Borel function. If $x\in\Dom(A)$ is nonzero and satisfies $Ax = \la x$, then $x\in \Dom(f(A))$ and $f(A)x = f(\la)x$.
\end{corollary}

\begin{proof}
Since $Ax=\la x$, the number $\la$ is an eigenvalue of $A$. By Proposition \ref{prop:normal-eigen}, applied to the function $z\mapsto \la-z$, we have $P_{\{\la\}}x=x$. It follows that the measure $P_x$ is supported on $\{\la\}$. Hence
$$
\int_{\si(A)}|f|^2\ud P_x=|f(\la)|^2\n x\n^2<\infty,
$$
so $x\in\Dom(f(A))$. Moreover,
$$
\n f(A)x-f(\la)x\n^2
=\int_{\si(A)}|f-f(\la)|^2\ud P_x=0.
$$
Thus $f(A)x=f(\la)x$.
\end{proof}

As in the bounded case, we can use the functional calculus to define square roots:

\begin{proposition}\label{prop:pos-sqrt}
 If $A$ is a positive selfadjoint operator, there exists a unique positive selfadjoint operator $B$ such that $B^2 = A$ with equality of domains.
\end{proposition}

\begin{proof}
The operator $f(A)$ with $f(\la) = \la^{1/2}$ is selfadjoint and positive,
and squares to $A$ by the result of Example \ref{normal-powers}. This proves existence.

To prove uniqueness, suppose that $B$ is a positive selfadjoint operator satisfying $B^2=A$. Let $P$ and $Q$ be the projection-valued measures of $A$ and $B$, respectively; both are supported on $[0,\infty)$, since $A$ and $B$ are selfadjoint and positive. Let $\phi:[0,\infty)\to[0,\infty)$ be given by $\phi(\la)=\la^2$, and let $R:=\phi_\#(Q)$. By the substitution formula for projection-valued measures and the result of Example \ref{normal-powers}, we have
$$
\int_{[0,\infty)}\la\ud R(\la)
 = \int_{[0,\infty)}\phi(\la)\ud Q(\la)
 = B^2=A.
$$
Thus $R$, viewed as a projection-valued measure on $\C$ by zero extension, represents $A$. By Proposition \ref{prop:supp-P-sigmaA}, both $R$ and $P$ are supported on $\si(A)$. The uniqueness part of Theorem \ref{thm:ST-unboundedn-normal} therefore gives $R=P$. Hence
\begin{align*}
A^{1/2}
&  = \int_{[0,\infty)}\la^{1/2}\ud P(\la)
 = \int_{[0,\infty)}\la^{1/2}\ud R(\la)
\\ & = \int_{[0,\infty)}(\la^2)^{1/2}\ud Q(\la)
 = \int_{[0,\infty)}\la\ud Q(\la)=B.
\end{align*}
This proves uniqueness.
\end{proof}

We will denote the square root of a positive selfadjoint operator $A$ by $A^{1/2}$.

If $A$ is normal, we may use the measurable calculus to define
$|A| :=  f(A)$, where $f(\la) = |\la| $. Furthermore, $A^\star A$ is selfadjoint and positive, so it has a unique selfadjoint and positive square root $(A^\star A)^{1/2}$ by Proposition \ref{prop:pos-sqrt}. The next corollary extends Corollary \ref{cor:sqrt-normal} to unbounded normal operators:

\begin{corollary} For every normal operator $A$ we have $\Dom(|A|) = \Dom(A)$
and $$(A^\star A)^{1/2} = |A|.$$
\end{corollary}
\begin{proof}
 With ${\rm id}(\la) = \la$ we have $\Dom(|A|) = \Dom(\Phi(|{\rm id}|)) = \Dom(\Phi({\rm id})) = \Dom(A)$,
 the middle identity being immediate from the definition of these domains.
With $f(\la)=|\la|$, we obtain
$$|A|^2 =f(A)f(A) = f^2(A)
= (\ov{\rm id}\cdot{\rm id})(A) =  \ov{\rm id}(A) {\rm id}(A) =A^\star A;$$
the identities $f(A)f(A)=f^2(A)$ and
$(\overline{\rm id}\cdot{\rm id})(A)=\overline{\rm id}(A)\,{\rm id}(A)$
follow from Theorem \ref{thm:Borel-FC-unbddnormal}
\ref{it:Borel-FC-unbddnormal1}, while
$\overline{\rm id}(A)=A^\star$ follows from
Theorem \ref{thm:Borel-FC-unbddnormal}\ref{it:Borel-FC-unbddnormal2}.
Since $|A|$ is positive selfadjoint and $|A|^2=A^\star A$, the uniqueness
part of Proposition \ref{prop:pos-sqrt} gives
$(A^\star A)^{1/2}=|A|$.
 \end{proof}

We proceed with some examples.

\begin{example}[Multiplication operators]\label{ex:normal-ex1}
Let $(\Om,\calF\!,\mu)$ be a $\sigma$-finite measure space and let
$m:\Om\to\C$ be a measurable function.
The linear operator $A_m$ defined by
\begin{align*}
 \Dom(A_m)  := \{f\in L^2(\Om,\mu):\, mf\in L^2(\Om,\mu)\},\qquad A_m f  := mf,
\end{align*}
is normal, its spectrum $\sigma(A_m)$ equals the $\mu$-essential range of $m$,
and its projection-valued measure is given by
$P_B f= \one_{m^{-1}(B)}f$ for $B\in \mathscr{B}(\sigma(A_m))$ and $f\in L^2(\Om,\mu).$
\end{example}

\begin{example}[Fourier multiplier operators]\label{ex:normal-ex2}
 Let $m:\R^d\to \C$ be a measurable function.
The linear operator $T_m$ defined by
\begin{align*}
 \Dom(T_m) := \{f\in L^2(\R^d):\, \calF f \in \Dom(A_m)\}, \qquad T_m f := \calF^{-1} A_m \calF f,
\end{align*}
where $A_m$ is the operator of the preceding example,
is normal, has spectrum $\sigma(T_m)=\sigma(A_m)$, and has projection-valued measure given by
$P_B f = T_{\one_{m^{-1}(B)}}f$
for $B\in \mathscr{B}(\sigma(T_m))$ and $f\in L^2(\R^d).$
Thus, the projections in the range of the projection-valued measure of the Fourier multiplier operator $T_m$
are Fourier multiplier operators.
\end{example}

We conclude with a precursor to the Weyl calculus developed in
Chapter \ref{ch:SymmQuant}, which motivates the connection between the Borel
calculus for the position operator and the Weyl calculus for functions
depending on position and momentum.

\begin{theorem}\label{thm:Fourier-functional-calculus}
Let $A$ be a selfadjoint operator in a Hilbert space $H$, and let
$f:\R\to\C$ be a bounded continuous function for which there exists a function
$g\in L^1(\R)$ such that
$$
f(\la)=\frac{1}{(2\pi)^{1/2}}\int_\R g(\xi)e^{i\xi\la}\ud\xi,
\quad \la\in\R.
$$
Then $f(A)$ is bounded, and for all $h\in H$ we have
\begin{equation*}
f(A)h=\frac{1}{(2\pi)^{1/2}}\int_\R g(\xi)e^{i\xi A}h\ud\xi,
\end{equation*}
where the right-hand side is a Bochner integral in $H$.
\end{theorem}

\begin{proof}
Since $A$ is selfadjoint, $f(A)$ is bounded by the bounded functional
calculus. Let $P$ be the projection-valued measure of $A$.

For fixed $h\in H$, the mapping $\xi\mapsto e^{i\xi A}h$
is bounded and continuous; this follows from a direct estimate using the bounded Borel calculus. Hence $\xi\mapsto g(\xi)e^{i\xi A}h$
is strongly measurable.
Since $g\in L^1(\R)$, the Bochner integral
$$
Bh:=\frac1{(2\pi)^{1/2}}\int_\R g(\xi)e^{i\xi A}h\ud\xi
$$
exists for all $h\in H$, and $B\in\calL(H)$.
Since
$$
\int_\R\int_\R |g(\xi)|\ud P_{h}(\lambda)\ud\xi<\infty,
$$
Fubini's theorem gives
\begin{align*}
 \iprod{Bh}{h}
 &=
 \frac1{(2\pi)^{1/2}}
 \int_\R g(\xi)\iprod{e^{i\xi A}h}{h}\ud\xi
=
 \frac1{(2\pi)^{1/2}}
 \int_\R g(\xi)
      \int_\R e^{i\xi\lambda}\ud P_{h}(\lambda)\ud\xi        \\
 &=
 \int_\R
 \Bigl(
 \frac1{(2\pi)^{1/2}}
 \int_\R g(\xi)e^{i\xi\lambda}\ud\xi
 \Bigr)
 \ud P_{h}(\lambda)
 =
 \int_\R f(\lambda)\ud P_{h}(\lambda)
 =
 \iprod{f(A)h}{h}.
\end{align*}
Since this holds for all $h\in H$, by polarisation we have $B=f(A)$.
\end{proof}

\begin{problems}

\item\label{prob:dd-bdd}
Let $A$ be a densely defined linear operator in a Banach space $X$ such that $\n Ax\n \le C\n x\n$ for some constant $C\ge 0$ and all $x\in \Dom(A)$. Prove that $A$ is closable, $\Dom(\ov A) = X$, and $\ov A$ is bounded with $\n \ov A x\n \le C\n x\n$ for all $x\in X$.

\item
Define the linear operator $A$ in $L^2(0,1)$ by $\Dom(A)  := C[0,1]$ and
\begin{align*}
Af := f(0) \one, \quad f\in \Dom(A).
\end{align*}
Show that $A$ is densely defined but nonclosable.

\item
Let $A$ be a densely defined closed operator in a Banach space $X$, and suppose there is a subspace $Y$, contained in $\Dom(A)$ and dense in $X$, such that $Ay = 0$ for all $y\in Y$. Does it follow that $Ax = 0$ for all $x\in \Dom(A)$?
What happens if the closedness assumption is dropped?

\item
In $C[0,1]$ consider the linear operators $A_1f := f'$ and $A_2 f:=f'$ with domains $\Dom(A_1):=C^{\infty}[0,1]$ and $\Dom(A_2):= C_{\rm c}^\infty(0,1)$.
\begin{enumerate}[\rm(a), leftmargin=*]
  \item Show that $A_1$ is closable and find the domain of its closure.
  \item Show that $A_2$ is closable and find the domain of its closure.
\end{enumerate}

\item
Provide the details for Example \ref{ex:bddpert-closed}.

\item
Let $A$ be a densely defined closed linear operator from a Banach space $X$ to a Banach space $Y$.
Prove or disprove:
\begin{enumerate}[\rm(a), leftmargin=*]
  \item for all $T\in \calL(X)$, the operator $AT$ with domain $\Dom(AT) = \{x\in X:\, Tx\in \Dom(A)\}$ is closed;
  \item for all $T\in \calL(Y)$ the operator $TA$ with domain $\Dom(TA) = \Dom(A)$ is closed.
\end{enumerate}

\item\label{prob:Hille}
Let $(\Om,\calF\!,\mu)$ be a measure space, $X$ be a Banach space, and suppose that
$f:\Om\to X$ is Bochner integrable with respect to $\mu$.

\begin{enumerate}[\rm(a), leftmargin=*]
  \item\label{it:Hille2} Prove {\em Hille's theorem}:\index{theorem!Hille}
  if $A$ is a closed linear operator in a Banach space $X$,
  $f$ takes values in $\Dom(A)$ $\mu$-almost everywhere, and
  the $\mu$-almost everywhere defined function
  $Af:\Om\to X$ is Bochner integrable, then $\int_\Om f\ud\mu\in \Dom(A)$ and
  $$A\int_\Om f\ud\mu = \int_\Om Af\ud\mu.$$

  \noindent{\em Hint:}\ Show that $\om\mapsto (f(\om), Af(\om))$ is Bochner integrable as a function with values
  in $X\times X$ and hence, by the result of Problem \ref{prob:Bochner-int-subspace}, as a function with values in the graph $\Gr(A)$.

  \item\label{it:Hille3} Justify the identity
  $$ \frac{\partial }{\partial t} \int_0^1 f(t,s)\ud s = \int_0^1 \frac{\partial}{\partial t}f(t,s)\ud s$$
  by providing conditions on $f$ so that the result of part \ref{it:Hille2}
  can be applied.
\end{enumerate}

\item
Give a proof of Proposition \ref{prop:injective-denserange-unbdd}.

\item
Let $A$ be a densely defined closed linear operator from a Banach space $X$ to a Banach space $Y$. Prove or disprove:
\begin{enumerate}[\rm(a), leftmargin=*]
  \item if $T\in \calL(X)$ and $AT$ is densely defined, then $\Dom((AT)\s) = \Dom(T^*A\s)$;
  \item if $T\in \calL(Y)$ and $A^*T\s$ is densely defined, then $\Dom((TA)\s) = \Dom(A^*T\s)$.
\end{enumerate}

\item
Let $A$ be a linear operator in a Banach space $X$.
Show that $\la\mapsto R(\la,A)$ is holomorphic on $\varrho(A)$ as an $\calL(X,\Dom(A))$-valued mapping, where $\Dom(A)$ carries the graph norm.

\item\label{prob:Arhancet}
Let $A$ and $B$ be linear operators acting in a Banach space $X$ such that $A \subseteq B$. Suppose that there exists $\lambda\in \C$ such that $\la+A$ is surjective and $\lambda+B$ is injective. Show that $A = B$.

\item\label{prob:full-spectrum}
Let $M$ be a nonempty closed subset of $\C$. In the Banach space $C_{\rm b}(M)$ of bounded continuous functions on $M$ consider the linear operator $A$ given by
\begin{align*}
\Dom(A) := \{f\in C_{\rm b}(M):\, z\mapsto zf(z) \in C_{\rm b}(M)\}, \qquad
Af(z) := zf(z).
\end{align*}
\begin{enumerate}[\rm(a), leftmargin=*]
  \item Show that $A$ is a closed operator.
  \item Show that $\sigma(A) = M$.
\end{enumerate}

\item \label{prob:diag-empty-res}
We show that the assumption $\varrho(A)\not=\emptyset$ in Proposition \ref{prop:sigmaAdiag} cannot be omitted.

Let $H:=\ell^2$, with standard orthonormal basis
$(e_n)_{n\ge1}$, and define
$$
\Dom(B):=\{a=(a_n)_{n\ge 1}:\,(na_n)_{n\ge 1}\in\ell^2\},
\qquad
Ba:=(na_n)_{n\ge 1}.
$$
Put $x:=(1/n)_{n\ge1}$, and define
$$
\Dom(A):=\Dom(B)\dotplus\operatorname{span}\{x\},
\qquad
A(u+c x):=Bu.
$$

\begin{enumerate}[\rm(a), leftmargin=*]
 \item Show that $A$ is densely defined and satisfies $Ae_n=ne_n$ for all $n\ge1$.
 \item Show that the graph of $A$ is given by $\Gr(A)=\Gr(B)+\operatorname{span}\{(x,0)\}$, and conclude that $A$ is closed.
\end{enumerate}
Thus $A$ has an
orthonormal basis of eigenvectors whose eigenvalues tend to infinity.
We next show that every $\mu\in\C$ is an eigenvalue of $A$; this implies $\sigma(A)=\C$.

\begin{enumerate}[\rm(a), leftmargin=*]\setcounter{enumii}{2}
 \item Show that every $\mu\in \N$ is an eigenvalue.
 \item Show that every $\mu\notin\N$ is an eigenvalue by completing the following argument. Define $u=(u_n)_{n\ge1}$ by
$u_n:=\frac{\mu}{n(n-\mu)}$. Show that
$
(nu_n)_{n\ge1} \in\ell^2(\N)$, and conclude that $u\in\Dom(B)$.
\item Show that
$(u+x)_n = \frac1{n-\mu}$,
so $u+x\ne 0$.
\item Show that $u+x\in\Dom(A)$ and
$A(u+x) = Bu = \mu(u+x)$.
\end{enumerate}

\item
Show that if $A$ and $B$ are linear operators in a complex Hilbert space
such that $\Dom(A)=\Dom(B)$ is dense and $ \iprod{Ax}{x} = \iprod{Bx}{x}$ for all $x\in \Dom(A) = \Dom(B)$, then $A = B$.

\item\label{prob:symm-dd}
Let $A$ be any nonclosable operator in a Hilbert space $H$.
Show that the operator $B$ on the Hilbert space direct sum $H\oplus H$ defined by $\Dom(B)  := \Dom(A) \oplus \{0\}$ and
\begin{align*}
B (x,0) :=  (0, Ax), \quad  (x , 0) \in \Dom(B),
\end{align*}
is symmetric and nonclosable. This example shows that
the assumption that the operator be densely defined cannot be omitted from Proposition \ref{prop:symmetry}.

\item\label{prob:meas-calc-domains}
Let $H=\ell^2$ and let $P$ be the projection-valued measure
on $\N$ given by
$$
P_F(x_n)_{n\ge1}:=(\one_F(n)x_n)_{n\ge1},
\quad F\subseteq\N.
$$
Let $\Phi$ be its measurable functional calculus, and put $f(0) = g(0):=0$, and
$f(n):=1/n$ and $g(n):=n$ for $n\ge 1$.

\begin{enumerate}[\rm(a), leftmargin=*]
\item Determine the domains and actions of $\Phi(f)$ and $\Phi(g)$.

\item Show that
$\Phi(f)\Phi(g)\subsetneq\Phi(fg)=I$ and $\Phi(g)\Phi(f)=I$.
Explain the difference between the two products using
\eqref{eq:domPhifg}.

\item Show that
$\Phi(g)+\Phi(-g)\subsetneq\Phi(g-g)=0$.
\end{enumerate}
In parts {\rm(b)} and {\rm(c)}, exhibit a vector on which
the right-hand side is defined but the proper restriction
on the left is not.

\item
Prove the claims in Examples \ref{ex:normal-ex1} and \ref{ex:normal-ex2}.

\item\label{prob:PVM-Laplace}
Combining Examples \ref{ex:Laplacian-sa} and  \ref{ex:normal-ex2},
find the projection-valued measure of the Laplace operator $\Delta$ on $L^2(\R^d)$, viewed as a selfadjoint operator in this space with domain $\Dom(\Delta) = H^2(\R^d)$.

\item
Let $A$ be a normal operator with projection-valued measure $P$, and let $B\subseteq\sigma(A)$ be a bounded Borel subset. Show that
$P_B x\in \Dom(A)$ for all $x\in H$.

\item
Let $A$ be a normal operator with projection-valued measure $P$, and let $f,g\in B_{\rm b}(\sigma(A))$.
Show that if $f = g$ $P$-almost everywhere (in the sense that there is a Borel set $N$ such that $P_N = 0$ and
$f=g$ on $\complement N$), then $f(A) = g(A)$.

\item
Extend the result of Problem \ref{prob:sa-StoneFormula} to unbounded
selfadjoint operators.

\item
Show that if $A$ is a positive selfadjoint operator, then $e^{-A}$ is bounded and injective. Is $e^{-A}$ always invertible?

\item
Show that normal $A$ is bounded if and only if $\sigma(A)$ is bounded.

\item
Let $A$ be a normal operator with projection-valued measure $P$ and let $f:\sigma(A)\to \C$ be Borel. Show that $f(A)$ is injective with dense range if and only if $ f\not=0$ $P$-almost everywhere, where $P$ is the projection-valued measure of $A$, and that in this case the densely defined inverse $(f(A))^{-1}$ is equal to $(1/f)(A)$.

\noindent{\em Hint:}\ Define $(1/f)(A)$ via the measurable functional calculus and use Theorem \ref{thm:Borel-FC-unbddnormal} to check that $\Dom((1/f)(A)f(A)) = \Dom(f(A))$.

\end{problems}

%% file: ch11-BoundaryValueProblems.tex
\chapter{Boundary Value Problems}\label{chap:boundaryvalueproblems}

\blfootnote{\rm \small \begin{quote}This book has been published by Cambridge University Press in the series `Cambridge Studies in Advanced Mathematics'. The present version is free to view and download for personal use only. Not for re-distribution, re-sale or use in derivative works.\newline {\copyright} Jan van Neerven\end{quote}}

\noindent
Having developed some of the core results of functional analysis, we now turn to applications to partial differential equations. Boundary value problems require us to solve a partial differential equation while prescribing a boundary condition. Solutions may lack the regularity needed to interpret the equation pointwise. Sobolev spaces supply a functional analytic setting in which derivatives and boundary conditions can be understood in a weaker sense. Variational methods connect the resulting problems with the geometry of Hilbert spaces and identify weak solutions as minimisers of a suitable energy functional.

\section{Sobolev Spaces}\label{sec:Sobolev}

We begin by developing some elements of the theory of Sobolev spaces.
Throughout this chapter we assume that $d\ge 1$ is an integer and $D$ is a nonempty open subset of $\R^d\!$.

\paragraph{Multi-Index Notation}

A $d$-tuple $\alpha = (\alpha_1,\hdots,\alpha_d)\in \N^d$ is called a {\em multi-index}\index{multi-index} of dimension $d$. Its {\em order}\index{order of a multi-index} is the nonnegative integer $$|\alpha| := \alpha_1+\cdots+\alpha_d.$$ We also define
$$ \al!:= \al_1!\cdots \al_d!.$$
We write $\al\le \be$ if $\al_j\le \be_j$ for all $j=1,\dots,d$, and in such cases we define
$\beta-\al := (\beta_1-\alpha_1,\hdots,\beta_d-\alpha_d)$ and
$$ \binom{\beta}{\al} := \frac{\beta!}{\al!(\beta-\al)!}.$$
For $x\in\R^d$ and $\al\in\N^d$ we set\index{$X$@$x^\al$} $$ x^\al := x_1^{\al_1}\cdots x_d^{\al_d}\!.$$
Similarly, we define\index{$D$@$\partial^\alpha$} $$\partial^{\alpha} := \partial_1^{\alpha_1}\circ\cdots\circ\partial_d^{\alpha_d}\!,$$
where $\partial_j $\index{$D$@$\partial_j$} is the partial derivative in the $j$th direction. By a standard result from calculus, for $C^{|\alpha|}$-functions the order in which the derivatives are taken is unimportant.

\paragraph{Test Functions}

By $C^\infty(D)$\index{$Caa$@$C^\infty(D)$} we denote the space of functions $f:D\to \K$ having continuous
derivatives $\partial^\al f$ on $D$ of all orders $\al\in\N^d$\!, and by $C_{\rm c}^\infty(D)$\index{$Cac$@$C_{\rm c}^\infty(D)$} its subspace consisting of all functions compactly supported in $D$; recall that this means that the closure of the set $\{x\in D:\, f(x)\not=0\}$ is a compact set contained in $D$. Elements of $C_{\rm c}^\infty(D)$ are referred to as {\em test functions} on $D$.\index{test function}

Test functions with suitable additional properties always exist:

\begin{lemma}[Test functions]\label{lem:smooth-cutoff}
The following assertions hold.
\begin{enumerate}[\rm(a), leftmargin=*]
\item\label{it:basic-test1}
For every $r>0$ there exists a function
$\eta\in C_{\rm c}^\infty(B(0;r))$ such that
$\eta\ge0$ and
$$
\int_{\R^d}\eta(x)\ud x=1.
$$

\item\label{it:basic-test2}
There exist functions $\psi_n\in C_{\rm c}^\infty(\R^d)$, $n\ge1$,
such that $0\le\psi_n\le1$, $\supp(\psi_n)\subseteq\ov B(0;2n)$, $\psi_n\equiv1$ on $B(0;n)$,
and
$$
\|\nabla\psi_n\|_{L^\infty(\R^d)}\le C_0/n
$$
for some constant $C_0\ge0$ independent of $n$.
\end{enumerate}
\end{lemma}

\begin{proof}
The function
$$
\rho(t):=
\begin{cases}
e^{-1/t},&t>0,\\
0,&t\le0,
\end{cases}
$$
belongs to $C^\infty(\R)$: by induction, for $t>0$ and $n\ge1$ we have $\rho^{(n)}(t)=p_n(1/t)e^{-1/t}$ for a polynomial $p_n$, and since $\lim_{t\downarrow0}\rho^{(n)}(t)=0$, each derivative extends continuously to $t=0$.

For \ref{it:basic-test1}, after scaling it suffices to take $r=1$.
The function $x\mapsto \rho(1-4|x|^2)$
is nonnegative, smooth, not identically zero, and has compact support
in $B(0;1)$. Multiplying it by a suitable positive constant gives
the required function.

For \ref{it:basic-test2}, put
$$
\theta(t):=
\frac{\rho(2-t)}
{\rho(2-t)+\rho(t-1)},\quad t\in\R.
$$
Then $\theta\in C^\infty(\R)$, $0\le\theta\le1$,
$\theta\equiv1$ on $(-\infty,1]$, and $\theta\equiv0$ on $[2,\infty)$.
Set
$$
\psi_n(x):=\theta(|x|/n).
$$
Since $\theta$ is constant near $0$, the function $\psi_n$ is smooth
also at the origin. The asserted support and pointwise properties are
immediate, and
$$
|\nabla\psi_n(x)|
\le \frac{\|\theta'\|_\infty}{n},
$$
which proves the result with $C_0:=\|\theta'\|_\infty$.
\end{proof}

Similarly, one defines the space $C^k(D)$\index{$Cad$@$C^k(D)$}, $k\in\N$, as the space of functions $f:D\to \K$ having continuous derivatives $\partial^\al f$ on $D$ of all orders $\al\in\N^d$ satisfying $|\al|\le k$ (with the convention that $C^0(D) = C(D)$), and $C_{\rm c}^k(D)$\index{$Cae$@$C_{\rm c}^k(D)$}
as its subspace of compactly supported functions.

By $C^\infty(\ov D)$\index{$Cab$@$C^\infty(\ov D)$} we denote the space of all functions
$f\in C^\infty(D)$ such that $\partial^\al f$ has a continuous extension to $\ov D$ for all $\al\in\N^d$\!. The spaces $C^k(\ov D)$\index{$Caf$@$C^k(\ov D)$} are defined similarly, by considering only the multi-indices satisfying $|\al|\le k$.

A measurable function $f:D\to \K$ is called
{\em locally integrable}\index{locally!integrable}\index{$L^1_{\rm loc}(D)$}
if its restriction to every open set $U$ with compact closure contained in $D$ is integrable. The space of all locally integrable
functions $f:D\to \K$ is denoted by $L^1_{\rm loc}(D)$; as always we identify functions that are equal almost everywhere.
In our study of weak derivatives we need the following result on convolutions.

\begin{proposition}\label{prop:YoungCk} Let $k$ be a nonnegative integer.
 If $f\in L^1_{\rm loc}(\R^d)$ and $g\in C_{\rm c}^k(\R^d)$, then the convolution $f*g$ is pointwise well defined and belongs to $C^k(\R^d)$, and
 we have
 $$\partial ^\alpha(f*g) = f*(\partial^\alpha g)$$
 for all multi-indices $\alpha\in\N^d$ satisfying $|\alpha|\le k$.
\end{proposition}
\begin{proof}
Fix $x_0\in\R^d$, and choose a compact set $K$ containing
$B(x_0;1)-\supp g$.
For every $|\alpha|\le k$ and $x\in B(x_0;1)$,
$$
|f(y)\partial^\alpha g(x-y)|
\le
\|\partial^\alpha g\|_\infty\,|f(y)|\one_K(y).
$$
The right-hand side is integrable. Hence the convolutions
$f*(\partial^\alpha g)$ are pointwise well defined, and dominated
convergence shows that they are continuous.

For $|\alpha|<k$, let $e_j$ denote the $j$th coordinate vector.
The fundamental theorem of calculus gives, for sufficiently small
$h\ne 0$,
$$
\Bigl|\frac1h(\partial^\alpha g(x_0+he_j-y)
                 -\partial^\alpha g(x_0-y))f(y)
\Bigr|
\le
\|\partial^{\alpha+e_j}g\|_\infty\,|f(y)|\one_K(y).
$$
Dominated convergence therefore permits differentiation under
the integral, giving
$$
\partial_j\bigl(f*(\partial^\alpha g)\bigr)(x_0)
=
\bigl(f*(\partial^{\alpha+e_j}g)\bigr)(x_0).
$$
Since $x_0$ was arbitrary, induction on $|\alpha|$ proves that
$f*g\in C^k(\R^d)$ and that its derivatives have the asserted form.
\end{proof}

The following version of Theorem \ref{thm:partition-unity} will be useful.

\begin{proposition}[Smooth partition of unity]\label{prop:partition-unity-smooth}\index{partition of unity!smooth}
 Let $F\subseteq U_1\cup\cdots\cup U_k,$
 where $F\subseteq \R^d$ is compact and the sets $U_j\subseteq\R^d$ are open for all $j=1,\dots,k$.
 Then there exist nonnegative functions $f_j\in C_{\rm c}^\infty(U_j)$, $j=1,\dots,k$, such that
 $$ f_1 + \cdots +f_k \equiv 1 \ \hbox{on $F$}.$$
\end{proposition}

Here, we think of the functions $f_j$ as elements of $C_{\rm c}^\infty(\R^d)$ with support in $U_j$.

\begin{proof}
After intersecting each $U_j$ with an open ball containing $F$, we may assume that the sets $U_j$ are bounded.
Since $F$ is compact and the sets $U_j$ are open, there exists a $\delta>0$ such that $F_\delta\subseteq U_1^\delta\cup\cdots\cup U_k^\delta$, where
$F_\delta := \{x\in\R^d :\, d(x,F) \le \delta\}$ and $U_j^\delta = \{x\in U_j :\, d(x,\complement U_j) > \delta\}$.
Theorem \ref{thm:partition-unity} provides us with nonnegative continuous functions $g_j:\R^d\to [0,1]$ supported in $U_j^\delta$ such that
 $$ g_1 + \cdots +g_k \equiv 1 \ \hbox{on $F_\delta$}.$$
Choose a nonnegative test function $\phi\in C_{\rm c}^\infty(\R^d)$ with compact support in the open ball $B(0;\delta)$ and satisfying $\int_{\R^d} \phi\ud x = 1$.
Put $f_j := g_j * \phi$. Then $f_j\in C^\infty(\R^d)$ by Proposition \ref{prop:YoungCk}, and $f_j\ge0$.
Moreover,
$$
  \supp(f_j)\subseteq \supp(g_j)+\supp(\phi)\subseteq U_j,
$$
since $\supp(g_j)\subseteq U_j^\delta$ and $\supp(\phi)\subseteq B(0;\delta)$.
Thus $f_j\in C_{\rm c}^\infty(U_j)$. Finally, if $x\in F$, then
$x-y\in F_\delta$ for all $y\in\supp(\phi)$, and therefore
\[
 \sum_{j=1}^k f_j(x) = \int_{\R^d}\sum_{j=1}^k g_j(x-y)\phi(y)\ud y
 = \int_{\R^d}\phi(y)\ud y = 1.
\qedhere\]
\end{proof}

\subsection{Weak Derivatives}\label{sec:weakder}

In order to make the theorems of functional analysis applicable to the theory of partial differential
equations, it is desirable to be able to discuss derivatives of functions in $L^p(D)$. The difficulty is that for such functions, the classical pointwise definition of differentiability through limits of difference quotients does not make sense since their values are well defined only almost everywhere. This necessitates an approach that is insensitive to redefining functions on sets of measure zero. Such an approach is provided by the notion of a {\em weak derivative}.
With the help of weak derivatives we then introduce the class of Sobolev spaces, which provides the $L^p$-analogues of
the classical spaces of continuously differentiable functions.

If $f\in C^k(D)$, then for all test functions $\phi\in C^{\infty}_{\rm c}(D)$ and multi-indices $\alpha\in \N^d$ with $|\alpha|\le k$
we have the integration by parts formula
\begin{align}\label{eq:def-weak-der-class}
\int_{D} f(x)\partial^{\alpha}\phi(x)\ud x = (-1)^{|\alpha|} \int_{D} g(x)\phi(x)\ud x,
\end{align}
where $g = \partial^\alpha f$. Using a smooth partition of unity (Proposition \ref{prop:partition-unity-smooth}), we can reduce the proof of this identity to the situation where the support of $\phi$ is contained in an open rectangle contained within $D$;
for such $\phi$, the formula follows by Fubini's theorem and integration by parts on intervals in dimension one.

This motivates the following definition.

\begin{definition}[Weak derivatives]\label{def:weak-der} Let $f\in L^1_{\rm loc}(D)$. A function $g\in L^1_{\rm loc}(D)$ is said to be a {\em weak derivative of order $\alpha\in\N^d$}\index{weak!derivative}
of $f$ if for all $\phi\in C^{\infty}_{\rm c}(D)$ we have
\begin{align}\label{eq:def-weak-der}
\int_{D} f(x)\partial^{\alpha}\phi(x)\ud x = (-1)^{|\alpha|} \int_{D} g(x)
\phi(x)\ud x.
\end{align}
A function $ f\in L^1_{\rm loc}(D)$ is said to be {\em weakly differentiable of order $k$} if it has weak derivatives $\partial^\alpha f\in L^1_{\rm loc}(D)$ for all multi-indices satisfying $|\alpha|\le k$.
\end{definition}

Here we adopt the convention that test functions enter bilinearly, without
complex conjugation.

\begin{remark}\label{rem:equiv-Ck-weakder}
The definition of a weak derivative of order $\al$ can equivalently be stated by using functions $\phi\in C_{\rm c}^{k}(D)$ for any integer $k\ge |\al|$. To see this, suppose that $g\in L^1_{\rm loc}(D)$ is a weak derivative of order $\al$ for the function $f\in L^1_{\rm loc}(D)$. We wish to prove that the integration by parts formula \eqref{eq:def-weak-der-class} holds for functions $\phi\in C_{\rm c}^{k}(D)$.
To this end we claim that there exist functions $\phi_n\in C_{\rm c}^\infty(D)$ such that $\phi_n\to \phi$ and $\partial^\alpha\phi_n\to \partial^\alpha\phi$ uniformly. Once this has been shown,  \eqref{eq:def-weak-der} follows from
\begin{align*}
\int_{D} f(x)\partial^{\alpha}\phi(x)\ud x
& = \limn \int_{D} f(x)\partial^{\alpha}\phi_n(x)\ud x
\\* & = (-1)^{|\alpha|} \limn \int_{D} g(x) \phi_n(x)\ud x = (-1)^{|\alpha|} \int_{D} g(x)\phi(x)\ud x.
\end{align*}
To prove the claim, let $\eta\in C_{\rm c}^\infty(\R^d)$ be supported in the unit ball $B(0;1)$ of $\R^d$ and satisfy $\int_{\R^d} \eta\ud x = 1$.
For $n\ge 1$ let $\eta^{(n)}(x) = n^{d}\eta(nx)$.
We extend $\phi$ by zero outside $D$ and define, for $y\in D$,
$$ \phi_n(y) := \eta^{(n)} * \phi (y) = \int_{\R^d} \eta^{(n)}(y-x)\phi(x) \ud x.
$$
Since $\eta^{(n)}$ is supported in $B(0;\frac1n)$, for sufficiently large $n$ the functions $ \phi_n$ are compactly supported in $D$. They are also smooth and hence belong to $C_{\rm c}^\infty(D)$ by Proposition \ref{prop:YoungCk},
and the desired convergence properties follow by elementary calculus arguments.
\end{remark}

The following proposition implies that weak derivatives, if they exist, are necessarily uniquely defined
up to a null set. This allows us to speak of {\em the} weak derivative of order $\alpha$ of a function $f$,
and denote it by $ \partial^\alpha f.$

In what follows we write\index{$AB$@$A\Subset B$}
$$U\Subset D$$
to express that the closure of $U$ is compact and contained in $D$.

\begin{proposition}\label{prop:fundvarcalculus}
If a function $g\in L^1_{\rm loc}(D)$ satisfies
$$\int_{D} g(x) \phi(x) \ud x = 0$$
for all $\phi\in C_{\rm c}^\infty(D)$, then $g = 0$ almost everywhere on $D$.
\end{proposition}

\begin{proof}
Let $B$ be any open ball such that $B\Subset D$ and let $\psi\in C_{\rm c}^\infty(D)$ satisfy $\psi \equiv 1$ on $B$.
Then $\psi g\in L^1(\R^d)$ after extending it by zero outside $D$.
Pick a mollifier $\eta\in C_{\rm c}^\infty(\R^d)$ satisfying $\int_{\R^d}\eta\ud x = 1$.
For $y\in\R^d$ set $\eta_{\eps,y} (x):= \eta_\eps(y-x)$, where $\eta_\eps(x) = \eps^{-d}\eta(\eps^{-1} x)$ for $x\in \R^d\!$.
Extending $\psi$ and $g$ by zero outside $D$ and noting that $\eta_{\eps,y}\psi \in C_{\rm c}^\infty(D)$, we infer from the assumption that for all $y\in \R^d$ we have
$$
 \eta_{\eps} * (\psi g) (y)
 = \int_{\R^d} \eta_{\eps}(y-x)\psi(x)g(x) \ud x
 = \int_{D} \eta_{\eps,y}(x)\psi(x) g(x) \ud x
 = 0.
$$
By Theorem \ref{thm:approx-identity} we have $\eta_{\eps} * (\psi g) \to \psi g$ in $L^1(\R^d)$ as $\eps\downarrow 0$.
It follows that $\psi g = 0$ almost everywhere on $\R^d$ and therefore $g = 0$ almost everywhere on $B$.
Since this is true for every open ball $B\Subset D$, the result follows.
\end{proof}

The following simple observation will be used repeatedly without further comment.

\begin{lemma}\label{lem:Sob-elementary} If $f\in L_{\rm loc}^1(D)$ and $\al\in\N^d$ is a multi-index, then:
\begin{enumerate}[label={\rm(\arabic*)}, leftmargin=*]
\item\label{it:Sob-elementary1} if $f$ has a weak derivative of order $\alpha$ and $D'$ is a nonempty open subset of $D$, then
$f|_{D'}$ has a weak derivative of order $\alpha$ given by $$\partial^\alpha (f|_{D'}) = (\partial^\alpha f)|_{D'};$$

\item\label{it:Sob-elementary2}
if $f$ has a weak derivative $g$ of order $\alpha$ and $g$ has a weak derivative $h$ of order $\beta$, then $f$ has a weak derivative of order $\alpha+\beta$ given by $h$, that is,
$$ \partial^\beta(\partial^\alpha f) = \partial^{\alpha+\beta} f.$$
\end{enumerate}
\end{lemma}
\begin{proof}
\ref{it:Sob-elementary1}: \ We consider only test functions $\phi\in C_{\rm c}^\infty(D')$ in \eqref{eq:def-weak-der}.
Extending them by zero to test functions defined on all of $D$, we obtain
\begin{align*}
 \int_{D'} f(x)\partial^{\alpha}\phi(x)\ud x & = \int_{D} f(x)\partial^{\alpha}\phi(x)\ud x
\\ & = (-1)^{|\alpha|} \int_{D} g(x)\phi(x)\ud x
= (-1)^{|\alpha|} \int_{D'} g(x)\phi(x)\ud x.
\end{align*}

\ref{it:Sob-elementary2}: \
For all $\phi,\psi\in C_{\rm c}^\infty(D)$ we have
$$\int_{D} f(x)\partial^{\alpha}\phi(x)\ud x = (-1)^{|\alpha|} \int_{D} g(x)
\phi(x)\ud x $$
and
$$ \int_{D} g(x)\partial^{\beta}\psi(x)\ud x = (-1)^{|\beta|} \int_{D} h(x)
\psi(x)\ud x,
$$
and the result follows by applying the first identity with $\phi = \partial^{\beta}\psi$.
\end{proof}

\begin{example}
The classical integration by parts formula \eqref{eq:def-weak-der-class} says that functions in $C^k(D)$ are weakly differentiable of order $k$, with weak derivatives given by their classical derivatives.
\end{example}

\begin{example}\label{ex:modulus-x-weakder}
 The function $f(x) = |x|$ has a weak derivative on $(-1,1)$, given by
 $${\rm sign}(x) =
 \begin{cases}
 \phantom{-}1, & x\in (0,1), \\ -1, & x\in (-1,0).
 \end{cases}
$$ This follows from
\begin{align*}
 \int_{-1}^1 |x|\phi'(x)\ud x & = \int_0^1 x\phi'(x)\ud x - \int_{-1}^0 x\phi'(x)\ud x
\\ & = - \int_0^1 \phi(x)\ud x + \int_{-1}^0 \phi(x)\ud x = -\int_{-1}^1 {\rm sign}(x)\phi(x)\ud x.
\end{align*}
A far-reaching generalisation of this example is given in Theorem \ref{thm:H1-lattice}.

In contrast, the function $f(x) = \sign(x)$ does not have a weak derivative on $(-1,1)$.
For if it did, by restriction the weak derivative $g$ would coincide with the classical derivative on both $(-1,0)$ and $(0,1)$, forcing $g = 0$ almost everywhere on $(-1,1)$. For $\phi\in C_{\rm c}^\infty(-1,1)$ with $\phi(0)\not=0$, this leads to the contradiction
$$
0 = \int_{-1}^1 g(x)\phi(x)\ud x = -\int_{-1}^1 \sign(x)\,\phi'(x)\ud x = 2\phi(0).
$$
\end{example}

We have the following local version of Proposition \ref{prop:YoungCk}.

\begin{proposition}\label{prop:YoungCk2}
Let $f\in L_{\rm loc}^1(D)$ have a weak derivative of order $\alpha$ on
$D$.
Suppose that $\eta\in C_{\rm c}^\infty(\R^d)$ has support in
$B(0;r)$ for some $r>0$, and let the open set $U\Subset D$ satisfy
$U+B(0;r)\subseteq D$.
For $x\in U$ define
$$
  (\eta*f)(x):=\int_{B(0;r)}\eta(y)f(x-y)\ud y.
$$
Then $\eta*f$ has weak and classical derivatives of order $\alpha$ on
$U$, and both are given by
\begin{align*}
\partial^\alpha (\eta* f)
= (\partial^\alpha \eta)*f =
\eta *(\partial^\alpha f) \ \hbox{ on } U.
\end{align*}
\end{proposition}

\begin{proof}
We first prove that $\eta*f$ is $C^\infty$ on $U$. Let $V\Subset U$. Then
$V-\supp\eta\Subset D$, and therefore we may choose
$\chi\in C_{\rm c}^\infty(D)$ such that $\chi\equiv1$ on an open set containing
$V-\supp\eta$. Extending $\chi f$ by zero outside $D$, on $V$ we have
$\eta*f=(\chi f)*\eta$.
By Proposition \ref{prop:YoungCk}, the right-hand side is $C^\infty$ on
$V$, and on this set we have
$$
  \partial^\alpha(\eta*f) = (\chi f)*(\partial^\alpha\eta) = (\partial^\alpha\eta)*f.
$$
Since $V\Subset U$ was arbitrary, this proves that $\eta*f$ is $C^\infty$
on $U$ and that on this set
$$
  \partial^\alpha(\eta*f)=(\partial^\alpha\eta)*f.
$$

It remains to identify this derivative with $\eta*(\partial^\alpha f)$. Fix
$x\in U$.
Since $\supp\eta\subseteq B(0;r)$ and $U+B(0;r)\subseteq D$, the function
$z\mapsto \eta(x-z)$ belongs to $C_{\rm c}^\infty(D)$.
Applying the weak derivative identity to
this test function gives
\begin{align*}
 \int_D f(z)\partial^\alpha\eta(x-z)\ud z
 &=
 \int_D (\partial^\alpha f)(z)\eta(x-z)\ud z .
\end{align*}
Indeed, differentiating $z\mapsto \eta(x-z)$ contributes the factor
$(-1)^{|\alpha|}$, which cancels the factor in the definition of the weak
derivative. Hence on $U$ we have
$$
  (\partial^\alpha\eta)*f = \eta*(\partial^\alpha f).
$$
This identity shows that the classical derivative
$\partial^\alpha(\eta*f)$ is equal to $\eta*(\partial^\alpha f)$. Since
classical derivatives are weak derivatives, the same identity holds in the
weak sense.
\end{proof}

In the proof of the next proposition we will use the fact that for all $1\le p\le \infty$, the operator
$f\mapsto \partial^\alpha f$ is closed as a linear operator in $L^p(D)$ with domain
\begin{align*}\Dom(\partial^\alpha):= \{f\in L^p(D):\, f \hbox{ has a weak derivative of order $\alpha$ in }L^p(D)\}.
\end{align*}
This domain of course depends on $p$, but we suppress this from the notation.
To prove that $\partial^\alpha$ is a closed operator, suppose that $f_n \to f$ in $L^p(D)$, with $f_n\in \Dom(\partial^\alpha)$ for all $n$, and
$\partial^\alpha f_n \to g$ in $L^p(D)$. We must prove that $f\in\Dom(\partial^\alpha)$ and $\partial^\alpha f=g$.
For all $\phi\in C_{\rm c}^\infty(D)$ we have
$$ \int_D f_n \partial^\al\phi\ud x = (-1)^{|\al|} \int_D \partial^\al f_n \phi\ud x .$$
Passing to the limit $n\to\infty$ in this formula (which is possible by H\"older's inequality, thanks to the fact that
test functions belong to $L^q(D)$ with $\frac1p+\frac1q=1$), we obtain
$$ \int_D f \partial^\al\phi\ud x = (-1)^{|\al|} \int_D g \phi\ud x .$$
This means that the function $g\in L^p(D)$ is a weak derivative of $f$ of order $\al$.

By $L_{\rm loc}^p(D)$ we denote the space of measurable functions whose restrictions to all sets $U\Subset D$ belong to $L^p(U)$, identifying functions that are equal almost everywhere on $D$.

\begin{proposition}\label{prop:L1loc-approx} Let $1\le p<\infty$.
 A function $f\in L_{\rm loc}^p(D)$ admits a weak derivative of order $\al$ in $L_{\rm loc}^p(D)$ if and only if there exist a sequence of functions $f_n\in C_{\rm c}^\infty(D)$ and a function $g\in L_{\rm loc}^p(D)$ such that for all open sets $U\Subset D$ we have:
 \begin{enumerate}[label={\rm(\roman*)}, leftmargin=*]
 \item\label{it:L1loc-approx1} $f_n \to f$ in $L^p(U)$;
 \item\label{it:L1loc-approx2} $\partial^\al f_n \to g$ in $L^p(U)$.
 \end{enumerate}
In this situation we have $\partial^\al f = g.$
\end{proposition}

\begin{proof}
`If':\ Let $\phi\in C_{\rm c}^\infty(D)$, and choose an open set
$U\Subset D$ with $\supp\phi\subseteq U$. Then
$$
 \int_D g\phi\ud x
 =
 \lim_{n\to\infty}\int_U \partial^\al f_n\,\phi\ud x
 =
 (-1)^{|\al|}\lim_{n\to\infty}\int_U f_n\,\partial^\al\phi\ud x
 =
 (-1)^{|\al|}\int_D f\,\partial^\al\phi\ud x .
$$
Hence $g=\partial^\al f$ weakly on $D$.

\smallskip
`Only if':\ Put $g:=\partial^\al f$. Let $(K_n)_{n\ge1}$ be an increasing
sequence of compact subsets of $D$ such that every compact subset of $D$ is
contained in $K_n$ for all sufficiently large $n$. Choose
$\chi_n\in C_{\rm c}^\infty(D)$ such that $0\le\chi_n\le 1$ and $\chi_n\equiv 1$ on an open set $W_n$
containing $K_n$.
Let $\eta\in C_{\rm c}^\infty(\R^d)$ satisfy
$\supp\eta\subseteq B(0;1)$ and $\int_{\R^d}\eta\ud x=1$,
and put $\eta_\eps(x):=\eps^{-d}\eta(\eps^{-1}x)$. Choose a sequence
$\eps_n\downarrow0$ so that for each $n$ we have
$$
  \supp(\chi_n)+\overline{B(0;\eps_n)}\Subset D,
  \qquad
  K_n+\overline{B(0;\eps_n)}\subseteq W_n .
$$
Extending $\chi_n f$ by $0$ outside $D$, define
$$
  f_n:=\eta_{\eps_n}*(\chi_n f).
$$
Then $f_n\in C_{\rm c}^\infty(D)$.

Let $U\Subset D$ be open. Choose an open set $\Omega$ such that
$U\Subset\Omega\Subset D$. For all sufficiently large $n$ we have
$\overline U\subseteq K_n$ and $U+B(0;\eps_n)\subseteq\Omega$.
Moreover, by the choice of $\eps_n$, we have $\chi_n\equiv1$ on
$U+B(0;\eps_n)$. Hence, for $x\in U$,
$$
  f_n(x)
  =
  \int_{\R^d}\eta_{\eps_n}(y)f(x-y)\ud y
  =
  \bigl(\eta_{\eps_n}*(\one_\Omega f)\bigr)(x).
$$
Since $\one_\Omega f\in L^p(\R^d)$ and $\one_\Omega f = f$ on $U$,
Theorem \ref{thm:approx-identity} gives $f_n\to f$ in $L^p(U)$.

It remains to prove the convergence of the weak derivatives. By Lemma
\ref{lem:Sob-elementary}\ref{it:Sob-elementary1}, the restriction
$f|_\Omega$ has weak derivative $g|_\Omega$ of order $\al$. Since
$U+B(0;\eps_n)\subseteq\Omega$, Proposition \ref{prop:YoungCk2}, applied
on the open set $\Omega$, gives
$$
  \partial^\al f_n
  =
  \eta_{\eps_n}*(\one_\Omega g)
  \ \hbox{ on }\, U
$$
for all sufficiently large $n$. Since $\one_\Omega g\in L^p(\R^d)$ and $\one_\Omega g = g$ on $U$, another
application of Theorem \ref{thm:approx-identity} gives
$\partial^\al f_n\to g$ in $L^p(U)$.
This proves \ref{it:L1loc-approx1} and \ref{it:L1loc-approx2}.
\end{proof}

As an application of this proposition we prove the following result on the existence of lower-order weak derivatives.

\begin{theorem}[Lower-order weak derivatives]\label{thm:lower-order}
Let $f\in L^1_{\rm loc}(D)$ be given and fix an integer $k\ge 1$.
If the weak derivative $\partial^\alpha f \in L^1_{\rm loc}(D)$ exists for every multi-index
$\alpha\in\N^d$ with $|\alpha|=k$, then $f$ admits weak derivatives $\partial^\beta f \in L^1_{\rm loc}(D)$ for all
$\beta\in\N^d$ with $|\beta|\le k$.
\end{theorem}

The proof uses the following one-dimensional estimate, which we will
apply along coordinate lines with $p=1$.

\begin{lemma}[Localised estimate]\label{lem:lower-order-local}
Let $I,J\subseteq\R$ be open intervals with $I\Subset J\Subset\R$,
and let $k\ge1$. Then there exists a constant $C\ge0$,
depending only on $k$ and $\delta:=d(I,\R\setminus J)$, such that for all
$f\in C^k(J)$ and $1\le p\le\infty$ we have
$$
\|f^{(k-1)}\|_{L^p(I)}
\le C\bigl(\|f\|_{L^p(J)}+\|f^{(k)}\|_{L^p(J)}\bigr).
$$
\end{lemma}

\begin{proof}
Choose a fixed function $\eta\in C_{\rm c}^\infty(-1,1)$ equal to
$1$ in an open set containing $0$, and put $\zeta(t):=\eta(t/\delta)$.

For $x\in I$, a first integration by parts gives
$$
\int_0^\delta\zeta(t)f^{(k)}(x+t)\ud t
=
-f^{(k-1)}(x)
-\int_0^\delta\zeta'(t)f^{(k-1)}(x+t)\ud t,
$$
since $\zeta(0)=1$ and $\zeta$ vanishes near $\delta$.
In each subsequent integration by parts, both boundary terms
vanish: all derivatives of $\zeta$ vanish near $\delta$, and
$\zeta^{(r)}(0)=0$ for $r\ge1$. Thus, for $r=1,\dots,k-1$,
$$
\int_0^\delta\zeta^{(r)}(t)f^{(k-r)}(x+t)\ud t
=
-\int_0^\delta\zeta^{(r+1)}(t)f^{(k-r-1)}(x+t)\ud t.
$$
Substituting these identities successively into the first one
and rearranging gives
$$
f^{(k-1)}(x)
=
-\int_0^\delta\zeta(t)f^{(k)}(x+t)\ud t
+(-1)^k\int_0^\delta\zeta^{(k)}(t)f(x+t)\ud t,
\quad x\in I.
$$
Since $I+t\subseteq J$ for $0\le t<\delta$, taking $L^p(I)$-norms
gives
$$
\|f^{(k-1)}\|_{L^p(I)}
\le
\|\zeta\|_{L^1(\R)}\|f^{(k)}\|_{L^p(J)}
+
\|\zeta^{(k)}\|_{L^1(\R)}\|f\|_{L^p(J)}.
$$
The identities
$\|\zeta\|_1=\delta\|\eta\|_1$ and
$\|\zeta^{(k)}\|_1=\delta^{1-k}\|\eta^{(k)}\|_1$
give the asserted dependence of the constant.
\end{proof}

\begin{proof}[Proof of Theorem \ref{thm:lower-order}]
{\em Step 1} --
We first show that if $u\in L_{\rm loc}^1(D)$ admits a weak
derivative $\partial_j^m u\in L_{\rm loc}^1(D)$ for some
$1\le j\le d$ and $m\ge1$, then it admits weak derivatives
$\partial_j^r u\in L_{\rm loc}^1(D)$ for all $0\le r\le m$.

For $m=1$ there is nothing to prove. Suppose that $m\ge2$, and
choose the sequence $(u_n)_{n\ge1}$ supplied by
Proposition \ref{prop:L1loc-approx} with $\alpha=me_j$.
Let
$$
Q=\prod_{i=1}^d I_i
\Subset
Q'=\prod_{i=1}^d J_i
\Subset D
$$
be open rectangles. Applying Lemma \ref{lem:lower-order-local}
with $p=1$ in the variable $x_j$, and integrating over the remaining
variables, gives
\begin{align*}
\|\partial_j^{m-1}(u_n-u_\ell)\|_{L^1(Q)}
&\le C\bigl(
\|u_n-u_\ell\|_{L^1(Q')}
+
\|\partial_j^m(u_n-u_\ell)\|_{L^1(Q')}
\bigr)
\to 0
\end{align*}
as $n,\ell\to\infty$.
Thus $\partial_j^{m-1}u_n$ converges in $L^1(Q)$ to some
$v_Q\in L^1(Q)$. Integrating by parts and passing to the limit, we obtain
$$
\int_Q v_Q\phi\ud x
=
(-1)^{m-1}\int_Q u\,\partial_j^{m-1}\phi\ud x,
\quad \phi\in C_{\rm c}^\infty(Q).
$$
These local weak derivatives agree on overlapping rectangles.
They therefore define a function in $L_{\rm loc}^1(D)$, which,
by a partition of unity argument, is the weak derivative
$\partial_j^{m-1}u$ on $D$. Repeating the argument proves the
assertion.

\smallskip
{\em Step 2} --
We prove the theorem by induction on $k$. The case $k=1$ is
immediate. Suppose that $k\ge2$ and that the result holds for
$k-1$.

Since $\partial_j^k f$ exists for every $j$, Step 1 shows that
$g_j:=\partial_j f$ belongs to $L_{\rm loc}^1(D)$.
For every multi-index $\alpha$ with $|\alpha|=k-1$ and every
$\phi\in C_{\rm c}^\infty(D)$, we have
\begin{align*}
\int_D g_j\,\partial^\alpha\phi\ud x
&=
-\int_D f\,\partial^{\alpha+e_j}\phi\ud x
=
(-1)^{k-1}\int_D(\partial^{\alpha+e_j}f)\phi\ud x.
\end{align*}
Thus $g_j$ has all weak derivatives of order $k-1$.
By the induction hypothesis, it has all weak derivatives of
orders at most $k-1$.

For any $\beta$ with $1\le|\beta|\le k$, choose $j$ such that
$\beta_j\ge1$. Lemma
\ref{lem:Sob-elementary}\ref{it:Sob-elementary2} now gives
$\partial^\beta f = \partial^{\beta-e_j}g_j\in L_{\rm loc}^1(D)$.
This completes the induction step.
\end{proof}

One may ask whether the following stronger version of
Theorem \ref{thm:lower-order} is true: if
$f\in L_{\rm loc}^1(D)$ admits a weak derivative $\partial^\al f$ for some $\al\in\N^d$\!, does it admit weak derivatives of all orders $0\le \beta\le\al$?
The following example gives a negative answer for $d\ge 2$.

\begin{example}\label{ex:lower-order-fails}
Let $D:=(-1,1)\times(-1,1)$ and define
$$
f(x):=\sign(x_1),\quad x = (x_1,x_2)\in D.
$$
Since $f$ is independent of $x_2$, it has weak derivative $\partial_2 f=0$, and hence also $\partial_1\partial_2 f=0$ in the sense of weak derivatives.

We next show that $f$ has no weak derivative $\partial_1 f$.
If there were $g\in L^1_{\rm loc}(D)$ such that
$$
\int_D f\partial_1 \phi \ud x = -\int_D g\phi\ud x,\quad \phi\in C_{\rm c}^\infty(D),
$$
then choosing $\phi(x_1,x_2)=\psi(x_1)\eta(x_2)$ with $\psi\in C_{\rm c}^\infty(-1,1)$ and $\eta\in C_{\rm c}^\infty(-1,1)$
satisfying $\int_{-1}^1\eta\ud x =1$ would give
\begin{align*}
 \int_{-1}^1 \sign(x_1)\psi'(x_1)\ud x_1
 & = \int_{-1}^1 \eta(x_2)\Bigl(\int_{-1}^1 \sign(x_1)\psi'(x_1)\ud x_1\Bigr) \ud x_2
\\ & = \int_{D} \sign(x_1)\partial_1 \phi(x)\ud x = \int_{D} f(x)\partial_1 \phi(x)\ud x
\\ & = -\!\int_D g(x)\phi(x)\ud x
 = -\!\int_{-1}^1 \Big(\int_{-1}^1 g(x_1,x_2)\eta(x_2)\ud x_2\Big)\psi(x_1)\ud x_1.
\end{align*}
This would mean that $x\mapsto \sign x$ has a weak derivative in $L^1_{\rm loc}(-1,1)$, contradicting Example \ref{ex:modulus-x-weakder}.
\end{example}

We have the following product rule.

\begin{proposition}[Product rule]\label{prop:productrule}
Let $G\subseteq \R^d$ be open, let $k\in\N$, and let $a\in C^k(G)$.
Suppose that $u\in L_{\rm loc}^1(G)$ has weak derivatives
$\partial^\beta u\in L_{\rm loc}^1(G)$ for all $|\beta|\le k$. Then $au$
has weak derivatives up to order $k$, and for every $|\alpha|\le k$,
$$
  \partial^\alpha(au)
  =
  \sum_{\beta\le\alpha}
  \binom{\alpha}{\beta}
  (\partial^{\alpha-\beta}a)(\partial^\beta u)
$$
in the weak sense.
\end{proposition}

\begin{proof}
The case $k=0$ is immediate. We first prove the first-order product rule.
Let $1\le j\le d$, and let $\phi\in C_{\rm c}^\infty(G)$. Since
$a\phi\in C_{\rm c}^1(G)$, Remark \ref{rem:equiv-Ck-weakder} permits us
to use $a\phi$ as a test function in the weak derivative identity. Hence
\begin{align*}
 \int_G au\,\partial_j\phi\ud x
 &= \int_G u\,\partial_j(a\phi)\ud x
  - \int_G u\,(\partial_j a)\phi\ud x
= -\int_G(\partial_j u)a\phi\ud x
  - \int_G u(\partial_j a)\phi\ud x .
\end{align*}
Thus
$\partial_j(au)=a\,\partial_j u+(\partial_j a)u$
weakly on $G$.

The general Leibniz formula follows by induction on $|\alpha|$, using this first-order formula at each step and the binomial identity for multi-indices.
\end{proof}

We next record a lemma on zero extensions.

\begin{lemma}[Zero extensions]\label{lem:der-zeroextension}
Let $D\subseteq\wt D$ be open sets, let $F\subseteq D$ be closed in
$\wt D$, and let $f\in L^1_{\rm loc}(D)$ vanish almost everywhere on
$D\setminus F$. If $f$ admits a weak derivative $\partial^\al f$ on
$D$, then the zero extension $\wt f$ of $f$ to $\wt D$ admits a weak
derivative $\partial^\al\wt f$ on $\wt D$, given by the zero extension
of $\partial^\al f$.

Consequently, if $1\le p<\infty$ and $f\in W^{k,p}(D)$, then
$\wt f\in W^{k,p}(\wt D)$ and
$$
\n\wt f\n_{W^{k,p}(\wt D)}
=
\n f\n_{W^{k,p}(D)}.
$$
\end{lemma}

\begin{proof}
Let $\phi\in C_{\rm c}^\infty(\wt D)$. Since
$F\cap\supp(\phi)$ is a compact subset of $D$, we may choose
$\eta\in C_{\rm c}^\infty(D)$ such that $\eta\equiv1$ on an open set containing $F\cap\supp(\phi)$. Then
$\partial^\al\phi=\partial^\al(\eta\phi)$ on an open set containing $F$,
and therefore
$$
\int_{\wt D}\wt f\,\partial^\al\phi\ud x
=
\int_D f\,\partial^\al(\eta\phi)\ud x
=
(-1)^{|\al|}\int_D(\partial^\al f)\eta\phi\ud x.
$$
By Lemma~\ref{lem:Sob-elementary}\ref{it:Sob-elementary1},
the restriction of $\partial^\al f$ to $D\setminus F$ is the weak
derivative of the zero function there, and hence vanishes almost
everywhere. Thus the last integral equals
$$
(-1)^{|\al|}
\int_{\wt D}\wt{\partial^\al f}\,\phi\ud x.
$$
This proves the first assertion. The final assertions follow by applying
this result to all $|\al|\le k$.
\end{proof}

The {\em gradient}\index{gradient} of a weakly differentiable function $f$ is the function $\nabla f\in L_{\rm loc}^1(D;\K^d)$ defined by \index{$D$@$\nabla$}
 $$\nabla f := (\partial_1 f, \dots,\partial_d f).$$

An open subset $U$ of $\R^d$ is said to be {\em (pathwise) connected}\index{connected} if any two points $x,y\in U$
can be joined by a continuous curve in $U$, that is, there exists a continuous mapping $\varphi:[0,1]\to U$
such that $\varphi(0)=x$ and $\varphi(1) = y$.

\begin{proposition}\label{prop:der-const}
Let $f\in L^1_{\rm loc}(D)$ be weakly differentiable. If $\nabla f = 0$ almost everywhere on an open connected subset $U$ of $D$, then $f$ is almost everywhere equal to a constant on $U$.
\end{proposition}
\begin{proof}
Let $B\Subset B'\Subset U$ be open balls, and let
$g\in L^1(\R^d)$ be the zero extension of $f|_{B'}$.
Choose a nonnegative function $\eta\in C_{\rm c}^\infty(B(0;1))$
with $\int_{\R^d}\eta\ud x=1$, and put
$\eta_\eps(x):=\eps^{-d}\eta(x/\eps)$ and $f_\eps:=\eta_\eps*g$.
For sufficiently small $\eps>0$, Proposition
\ref{prop:YoungCk2}, applied on $B'$, shows that $f_\eps$ is smooth
and has zero gradient on $B$. Thus $f_\eps$ is constant on $B$.

By Theorem \ref{thm:approx-identity},
$f_\eps\to f$ in $L^1(B)$. Consequently, the constants
$f_\eps|_B$ converge, and $f$ is almost everywhere constant on $B$.
Since $B\Subset U$ was arbitrary, and the constants on overlapping
balls agree, connectedness of $U$ gives the result.
\end{proof}

\subsection{The Sobolev Spaces \texorpdfstring{$W^{k,p}(D)$}{W\textasciicircum\{k,p\}(D)}}\label{sec:Wkp}

By H\"older's inequality, every $f\in L^p(D)$ with $1\le p\le \infty$ belongs to $L^1_{\rm loc}(D)$. This suggests the following definition.

\begin{definition}[Sobolev spaces]\label{def:Sobolev}
For $k\in\N$ and $1\le p\le \infty$ the {\em Sobolev space}\index{Sobolev space}
$W^{k,p}(D)$\index{$WKP$@$W^{k,p}(D)$}
is the space of all $f\in L^p(D)$ whose weak derivatives of all orders $|\alpha|\leq k$
exist and belong to $L^p(D)$.
\end{definition}

Endowed with the norm $$\|f\|_{W^{k,p}(D)} :=
\sum_{|\alpha|\le k} \|\partial^\alpha f\|_{L^p(D)},$$
$W^{k,p}(D)$ is a Banach space. Indeed, if $(f_n)_{n\ge 1}$
is a Cauchy sequence in $W^{k,p}(D)$, then for all $|\alpha|\leq k$ the sequence
of weak derivatives $(\partial^\alpha f_n)_{n\ge 1}$
is a Cauchy sequence in $L^{p}(D)$ and
hence convergent to some $f^{(\alpha)}\in L^p(D)$. Set $f :=f^{(0)}$.
Using H\"older's inequality as in Corollary \ref{cor:Holder}, we may pass to the limit $n\to \infty$
in the identity
$$
\int_{D} f_n(x)\partial^{\alpha}\phi(x)\ud x = (-1)^{|\alpha|} \int_{D} \partial^{\al} f_n(x)
\phi(x)\ud x, \quad \phi\in C_{\rm c}^\infty(D).
$$
It follows that $f$ has weak derivatives of all orders $|\alpha|\leq k$ given by $\partial^\alpha f = f^{(\alpha)}$\!.
Having observed this, we see that $f\in W^{k,p}(D)$ and $f_n\to f$ in $W^{k,p}(D)$.

We record the following simple consequence of Proposition \ref{prop:productrule}.

\begin{proposition}\label{prop:Cm-mult-Wmp}
Let $k\in\N$, let $1\le p\le\infty$, and let $a\in C^k(D)$.
If $\partial^\alpha a\in L^\infty(D)$ for all $|\alpha|\le k$,
then multiplication by $a$ defines a bounded operator on
$W^{k,p}(D)$.
\end{proposition}

\begin{proof}
By Proposition \ref{prop:productrule}, for every $u\in W^{k,p}(D)$
and every $|\alpha|\le k$,
$$
\partial^\alpha(au)
=
\sum_{\beta\le\alpha}
\binom{\alpha}{\beta}
(\partial^{\alpha-\beta}a)(\partial^\beta u).
$$
Consequently,
$$ \|au\|_{W^{k,p}(D)}
\le
\sum_{|\alpha|\le k}
\sum_{\beta\le\alpha}
\binom{\alpha}{\beta}
\|\partial^{\alpha-\beta}a\|_\infty
\|\partial^\beta u\|_{p}\le
C\max_{|\gamma|\le k}\|\partial^\gamma a\|_\infty
\|u\|_{W^{k,p}(D)},
$$
where $C$ is a constant depending only on $k$ and the dimension $d$.
\end{proof}

For $1\le p<\infty$, the linear operators $\partial^\alpha$ are closable as operators in $L^p(D)$ with initial domain $C_{\rm c}^\infty(D)$; this follows by an argument similar to the one just given. By Proposition \ref{prop:Cc-dense}, these operators are densely defined.
This prompts the question of which functions in $W^{k,p}(D)$ can be approximated, in the norm of this space, by test functions. A first result in this direction is the following proposition.

\begin{proposition}[Local approximation]\label{prop:local-approx-W1p}
 For all $f\in W^{1,p}(D)$ with $1\le p<\infty$ there exists a sequence of functions $f_n\in C_{\rm c}^\infty(D)$ such that
 for every open set $U\Subset D$ we have
 $$\limn f_n|_U = f|_U \ \ \hbox{in} \ \ W^{1,p}(U).$$
\end{proposition}

\begin{proof}
The functions $f_n$ constructed in the proof of Proposition \ref{prop:L1loc-approx}
have the required properties.
Indeed, their construction depends only on $f$ and
not on the multi-index $\alpha$. Applying the conclusion of that proof with
$\alpha=e_j$, $1\le j\le d$, gives $\partial_j f_n\to \partial_j f$
in $L^p(U)$
for every open set $U\Subset D$. The same construction also gives
$f_n\to f$ in $L^p(U)$. Hence $f_n|_U\to f|_U$ in $W^{1,p}(U)$.
\end{proof}

\begin{proposition}[Chain rule]\label{prop:chainrule}\index{chain rule}
Let $\rho:\R\to\R$ be a $C^1$-function such that $\rho'$ is bounded and $\rho(0)=0$, and let $1\le p\le \infty$.
For every real-valued $f\in W^{1,p}(D)$ we have
$\rho\circ f\in W^{1,p}(D)$ and
$$
\partial_j(\rho\circ f) = (\rho'\circ f)\partial_j f, \quad j=1,\dots,d.
$$
\end{proposition}

\begin{proof}
Put $M:=\|\rho'\|_\infty$.
Since $\rho(0)=0$, we have $|\rho(f)|\le M|f|$.
Also,
$|(\rho'\circ f)\partial_jf|\le M|\partial_jf|$.
Thus these functions belong to $L^p(D)$, and it remains to verify the identity
for the weak derivatives.

Fix $\phi\in C_{\rm c}^\infty(D)$, and choose open sets
$U,V$ such that
$\supp\phi\subseteq U\Subset V\Subset D$.
Since $V$ has finite measure, $f|_V\in W^{1,1}(V)$.
By Proposition \ref{prop:local-approx-W1p}, applied on $V$, there
exist real-valued functions $f_n\in C_{\rm c}^\infty(V)$ such that
$f_n\to f$ in $W^{1,1}(U)$.
Passing to a subsequence, we may also assume that $f_n\to f$
almost everywhere on $U$.

The Lipschitz estimate for $\rho$ gives
$\rho\circ f_n\to\rho\circ f$ in $L^1(U)$.
Moreover,
\begin{align*}
&\|(\rho'\circ f_n)\partial_jf_n
      -(\rho'\circ f)\partial_jf\|_{L^1(U)}
\\
&\qquad\le
M\|\partial_jf_n-\partial_jf\|_{L^1(U)}
+
\int_U|\rho'(f_n)-\rho'(f)|\,|\partial_jf|\ud x
\to 0.
\end{align*}
Here the first term tends to zero by the choice of $f_n$, and the
second by dominated convergence.
Passing to the limit in the classical chain rule
$$
\int_U(\rho\circ f_n)\partial_j\phi\ud x
=
-\int_U(\rho'\circ f_n)(\partial_jf_n)\phi\ud x
$$
therefore gives
$$
\int_D(\rho\circ f)\partial_j\phi\ud x
=
-\int_D(\rho'\circ f)(\partial_jf)\phi\ud x.
$$
Since $\phi$ was arbitrary, this proves the assertion.
\end{proof}

We write $W_{\rm loc}^{1,1}(D)$ for the space of all
$f\in L_{\rm loc}^1(D)$ whose weak first-order derivatives belong to
$L_{\rm loc}^1(D)$.

\begin{proposition}[Substitution rule]\label{prop:substitutionrule}
Let $D,D'\subseteq\R^d$ be nonempty open sets, and let
$\varrho:D\to D'$ be a $C^1$-diffeomorphism. If
$f\in W_{\rm loc}^{1,1}(D)$, then
$f\circ\varrho^{-1}\in W_{\rm loc}^{1,1}(D')$.
Writing $g:=f\circ\varrho^{-1}$, on $D'$ we have
$$
  \partial_j g = \sum_{i=1}^d \frac{\partial(\varrho^{-1})_i}{\partial y_j}
 \,(\partial_i f)\circ\varrho^{-1}.
$$
\end{proposition}

\begin{proof}
The idea of the proof is that, for smooth functions $f$, this is the substitution rule from calculus; the general case follows by local approximation with smooth functions.

Let $G'\Subset D'$. Choose an open set $H'$ such that
$G'\Subset H'\Subset D'$,
and put $G:=\varrho^{-1}(G')$, $H:=\varrho^{-1}(H')$. Then
$G\Subset H\Subset D$. Since $\varrho$ and $\varrho^{-1}$ are $C^1$,
their first derivatives are bounded on $\overline H$ and
$\overline{H'}$, and the Jacobian determinant of $\varrho$ is bounded
above and below away from $0$ there.

Choose $\chi\in C_{\rm c}^\infty(D)$ such that $\chi\equiv1$ on an open set
containing $\overline H$. Then $\chi f\in W^{1,1}(D)$. By
Proposition \ref{prop:local-approx-W1p}, there exist functions
$f_n\in C_{\rm c}^\infty(D)$ such that
$f_n|_H\to (\chi f)|_H=f|_H$ in $W^{1,1}(H)$.
Put $g_n:=f_n\circ\varrho^{-1}$. Then $g_n\in C^1(H')$, and by the
classical chain rule,
$$
  \partial_j g_n
  =
  \sum_{i=1}^d
  \frac{\partial(\varrho^{-1})_i}{\partial y_j}
  \,(\partial_i f_n)\circ\varrho^{-1}
  \quad\hbox{on }H'.
$$
By the change of variables formula and the boundedness of the Jacobian
determinant on $\overline H$, we have $g_n\to f\circ\varrho^{-1}$ in $L^1(G')$.
Similarly, using also the boundedness of the coefficients
$\partial(\varrho^{-1})_i/\partial y_j$ on $\overline{G'}$, the right-hand
side converges in $L^1(G')$ to
$$
  \sum_{i=1}^d
  \frac{\partial(\varrho^{-1})_i}{\partial y_j}
  \,(\partial_i f)\circ\varrho^{-1}.
$$
Passing to the limit in the weak derivative identity for $g_n$ on $G'$
therefore gives
$$
  \partial_j(f\circ\varrho^{-1})
  =
  \sum_{i=1}^d
  \frac{\partial(\varrho^{-1})_i}{\partial y_j}
  \,(\partial_i f)\circ\varrho^{-1}
  \quad\hbox{weakly on }G'.
$$
Thus $f\circ\varrho^{-1}\in W^{1,1}(G')$. Since $G'\Subset D'$ was
arbitrary, the result follows.
\end{proof}

\subsection{The Sobolev Spaces \texorpdfstring{$W_0^{1,p}(D)$}{W\textunderscore 0\textasciicircum\{1,p\}(D)}}\label{subsec:W01p}

We now introduce a class of spaces which play an important role in the theory of partial differential equations, where they provide the $L^p$-setting for studying boundary value problems subject to Dirichlet boundary conditions.

\begin{definition}[The spaces $W_0^{1,p}(D)$]\label{def:W01p} For $1\le p< \infty$ we define $W_0^{1,p}(D)$\index{$WZ$@$W_0^{1,p}(D)$} to be the closure of $C_{\rm c}^\infty(D)$ in $W^{1,p}(D)$.
\end{definition}

The following result gives a simple sufficient condition for membership of $W_0^{1,p}(D)$:

\begin{proposition}\label{prop:local-approx} Let $U\Subset D$ be open and let $1\le p<\infty$.
If $f\in W^{1,p}(D)$ vanishes outside $U$, then $f\in W_0^{1,p}(D)$.
\end{proposition}

\begin{proof}
Choose $\chi\in C_{\rm c}^\infty(D)$ equal to $1$ on an open set containing $\ov U$, and choose an open set $V\Subset D$ containing
$\supp\chi$. By Proposition \ref{prop:local-approx-W1p}, there
exist functions $f_n\in C_{\rm c}^\infty(D)$ such that
$f_n|_V\to f|_V$ in $W^{1,p}(V)$.
Since $f$ vanishes outside $U$, we have $\chi f=f$.
The support of $\chi$ is contained in $V$, so Proposition
\ref{prop:Cm-mult-Wmp} gives
\begin{align*}
\|\chi f_n-f\|_{W^{1,p}(D)} = \|\chi(f_n-f)\|_{W^{1,p}(V)} \le C\|f_n-f\|_{W^{1,p}(V)} \to 0.
\end{align*}
As $\chi f_n\in C_{\rm c}^\infty(D)$, this proves that
$f\in W_0^{1,p}(D)$.
\end{proof}

It is immediate that if $f\in W^{1,p}(D)$, then
$\Re f,\, \Im f \in W^{1,p}(D)$, and similarly if $f\in W_0^{1,p}(D)$, then
$\Re f,\, \Im f \in W_0^{1,p}(D)$. The next theorem asserts that the positive part $f^+$ of a
real-valued function $f$ in $W^{1,p}(D)$ belongs to $W^{1,p}(D)$ and provides an explicit expression for its weak derivatives; moreover, if $f\in W_0^{1,p}(D)$, then $f^+\in W_0^{1,p}(D)$. The analogous results hold for the negative part $f^- = f^+-f$ and the absolute value $|f| = f^+ + f^-$.

\begin{theorem}[Positive parts]\label{thm:H1-lattice} If $1\le p<\infty$ and $f\in W^{1,p}(D)$ is real-valued, then:
\begin{enumerate}[label={\rm(\arabic*)}, leftmargin=*]
 \item\label{it:H1-lattice1}
we have $f^+\in W^{1,p}(D)$ and
$$\partial_j f^{+} = \one_{\{f>0\}}\partial_j f, \quad j=1,\dots,d;$$
\item\label{it:H1-lattice2} the mapping $f\mapsto f^+$ is continuous in the norm of $W^{1,p}(D)$;
\item\label{it:H1-lattice3} if $f\in W_0^{1,p}(D)$, then $f^+\in W_0^{1,p}(D)$.
\end{enumerate}
\end{theorem}

\begin{proof} The idea of the proof is to approximate the function $\rho: t\mapsto t^+ = \max\{t,0\}$ with $C^1$-functions $\rho_n$ in such a way that for all $t\in \R$ we have
\begin{enumerate}[label={\rm(\roman*)}, leftmargin=*]
  \item\label{it:H1-lattice-i} $0\le \rho_n(t)\uparrow t^+$;
  \item\label{it:H1-lattice-ii} $0\le \rho_n'(t)\uparrow \one_{(0,\infty)}(t)$.
\end{enumerate}
For instance, the choice $\rho_n(t) = (t^2+\frac1{n^2})^{1/2} - \frac1n$ for $t>0$ and $\rho_n(t) = 0$ for $t\le 0$ will do; see Figure \ref{fig:rho}.

\begin{figure}
\begin{center}
 \begin{tikzpicture}[domain=0:2.6, scale = 1.5]
 \draw[very thin,color=gray] (-0.9,-0.9) grid (2.9,2.9);
 \draw[->] (-0.9,0.001) -- (2.9,0.001) node[right] {$t$};
 \draw[->] (0.004,-0.9) -- (0.004,2.9) node[right] {};
\draw[thick] plot (\x,\x)      node[right] {$\rho(t) =t^+$};
\draw[thick] (-0.9,0.005) -- (0,0.005);
\draw[thin] plot (\x, {sqrt(\x^2+1) - 1}) node[right] {$\rho_1(t) = (t^2+1)^{1/2}-1$};
\draw[thin] plot (\x, {sqrt(\x^2+1/4) - 1/2}) node[right] {$\rho_2(t)= (t^2+\frac14)^{1/2}-\frac12$};
\draw[thin] plot (\x, {sqrt(\x^2+1/9) - 1/3}) node[right] {};
\draw[thin] plot (\x, {sqrt(\x^2+1/16) - 1/4}) node[right] {};
 \end{tikzpicture}
 \caption{The functions $\rho_n(t) = (t^2+\frac1{n^2})^{1/2} - \frac1n$ on the positive half-line \label{fig:rho}}
\end{center}
\end{figure}

\smallskip \ref{it:H1-lattice1}: \
Let $f$ be a real-valued function in $W^{1,p}(D)$. By Proposition \ref{prop:chainrule} we have
$\rho_n \circ f \in W^{1,p}(D)$ and,
for test functions $\phi\in C_{\rm c}^\infty(D)$,
\begin{align*}
 - \int_D (\rho_n \circ f)(x) \partial_j \phi(x)\ud x = \int_D (\rho_n' \circ f)(x) \partial_j f(x)\phi(x)\ud x.
\end{align*}
By \ref{it:H1-lattice-i}, \ref{it:H1-lattice-ii}, and dominated convergence we obtain
$$ - \int_D f^+(x) \partial_j \phi(x)\ud x = \int_D\one_{\{f>0\}}(x) \partial_j f(x)\phi(x)\ud x.$$
This proves \ref{it:H1-lattice1}.

\smallskip\ref{it:H1-lattice2}:\
Applying \ref{it:H1-lattice1} to $-f$ gives $f^-:=(-f)^+\in W^{1,p}(D)$ and
$$
\partial_j f^-=-\one_{\{f<0\}}\partial_j f,\quad j=1,\dots,d.
$$
Since $f=f^+-f^-$, it follows that
$$
\partial_j f=\partial_j f^+-\partial_j f^-
=\one_{\{f>0\}}\partial_j f+\one_{\{f<0\}}\partial_j f
=\one_{\{f\neq0\}}\partial_j f.
$$
Hence $\one_{\{f=0\}}\partial_j f=0$ almost everywhere on $D$.

Now let $f_n\to f$ in $W^{1,p}(D)$, with all functions real-valued. Clearly
$f_n^+\to f^+$ in $L^p(D)$. It remains to prove that
$\partial_j f_n^+\to \partial_j f^+$ in $L^p(D)$ for all $1\le j\le d$. By \ref{it:H1-lattice1} this amounts to proving that
$\one_{\{f_n>0\}}\partial_j f_n \to \one_{\{f>0\}}\partial_j f$ in $L^p(D)$.

We have
\begin{align*}
\|\one_{\{f_n>0\}}\partial_j f_n
   -\one_{\{f>0\}}\partial_j f\|_p
\le
 \|\partial_j f_n-\partial_j f\|_p + \|(\one_{\{f_n>0\}}-\one_{\{f>0\}})\partial_j f\|_p .
\end{align*}
The first term tends to $0$. To handle the second term, by a standard subsequence argument it suffices to consider a subsequence along which
$f_n\to f$ almost everywhere. Along such a subsequence we have
$(\one_{\{f_n>0\}}-\one_{\{f>0\}})\partial_j f\to 0$ almost everywhere.
Indeed, at points where $f\neq0$ the characteristic functions converge, while
at points where $f=0$ we have $\partial_j f=0$ almost everywhere, and dominated convergence gives
the convergence to $0$ in $L^p(D)$. Thus $f_n^+\to f^+$ in $W^{1,p}(D)$.

\smallskip\ref{it:H1-lattice3}: \
Suppose now that $f\in W_0^{1,p}(D)$. Choose functions
$f_n\in C_{\rm c}^\infty(D)$ such that
$f_n\to f$ in $W^{1,p}(D)$. Replacing every $f_n$ by $\Re f_n$, we
may assume that every $f_n$ is real-valued. Then
$f_n^+\to f^+$ in $W^{1,p}(D)$ by
\ref{it:H1-lattice2}.
For each fixed $n$, the function $f_n^+$ belongs to $W^{1,p}(D)$ by \ref{it:H1-lattice1}
and vanishes outside $\supp(f_n)\Subset D$. Hence
$f_n^+\in W_0^{1,p}(D)$ by Proposition \ref{prop:local-approx}. Since
$W_0^{1,p}(D)$ is closed in $W^{1,p}(D)$, it follows that
$f^+\in W_0^{1,p}(D)$.
\end{proof}

The following theorem connects the space $W_0^{1,p}(D)$ with Dirichlet
boundary conditions. The notion of a $C^1$-boundary is defined at the
beginning of Section \ref{subsec:extension}, and the proof uses the
local invariance result proved there.

\begin{theorem}\label{thm:H01-cont} Let $D$ be bounded and let $1\le p<\infty$. For all $f\in W^{1,p}(D) \cap C(\ov D)$
the following assertions hold:
\begin{enumerate}[label={\rm(\arabic*)}, leftmargin=*]
 \item\label{it:H01-cont1} if $f|_{\partial D} \equiv 0$, then $f\in W_0^{1,p}(D)$;
 \item\label{it:H01-cont2} if $\partial D$ is $C^1$ and $f\in W_0^{1,p}(D)$, then $f|_{\partial D} \equiv 0$.
\end{enumerate}
\end{theorem}
\begin{proof}
\ref{it:H01-cont1}: \ By considering real and imaginary parts separately, we may assume that $f$ is real-valued, and by Theorem \ref{thm:H1-lattice} we may even assume that $f$ is nonnegative.
Since $f-\frac1k \to f$ in $W^{1,p}(D)$ as $k\to\infty$ and taking positive parts is continuous in $W^{1,p}(D)$ by Theorem \ref{thm:H1-lattice}, we see that
$f_k := (f-\frac1k)^+ \to f^+ = f$ in $W^{1,p}(D)$ as $k\to\infty$. Hence, to prove that $f\in W_0^{1,p}(D)$ it suffices to prove that $f_k\in W_0^{1,p}(D)$ for all $k\ge 1$.

By continuity, each $f_k$ vanishes in an open set containing $\partial D$. Then Proposition \ref{prop:local-approx} shows that $f_k$ can be approximated in $W^{1,p}(D)$ by functions $f_{k,n}\in C_{\rm c}^\infty(D)$.

\smallskip
\ref{it:H01-cont2}: \ We first reduce to the case of the half-space $\R_+^d=\{x\in\R^d:\,x_d>0\}$.
Let $x_0\in\partial D$. By the definition of a $C^1$-boundary, there are open sets $U$ containing $x_0$ and $V$ containing $0$, and a $C^1$-diffeomorphism
$\varrho:U\to V$, such that
$$
  \varrho(U\cap D)=V\cap\R_+^d,
  \qquad
  \varrho(U\cap\partial D)=V\cap\partial\R_+^d .
$$
Choose $\psi\in C_{\rm c}^\infty(U)$ with $\psi\equiv1$ on an open set containing $x_0$.
Since $f\in W_0^{1,p}(D)$, we may choose $f_n\in C_{\rm c}^\infty(D)$
with $f_n\to f$ in $W^{1,p}(D)$. By Proposition \ref{prop:Cm-mult-Wmp},
$\psi f_n\to\psi f$ in $W^{1,p}(D)$, and each $\psi f_n$ belongs to
$C_{\rm c}^\infty(D)$. Hence $\psi f\in W_0^{1,p}(D)$.

Choose an open set $G$ such that
$$
\supp(\psi)\Subset G\Subset U,
$$
and put $G_+:=G\cap D$. Then
$$
\varrho(G_+)=\varrho(G)\cap\R_+^d.
$$
By Proposition \ref{prop:diffeo-Wkp}\ref{it:diffeo-Wkp2}, each
$(\psi f_n)\circ\varrho^{-1}$ belongs to
$W_0^{1,p}(\varrho(G_+))$. Hence its zero extension belongs to
$W_0^{1,p}(\R_+^d)$. These zero extensions converge in
$W^{1,p}(\R_+^d)$ to the zero extension of
$(\psi f)\circ\varrho^{-1}$. Consequently, this zero extension belongs
to $W_0^{1,p}(\R_+^d)$.
It is also continuous on $\ov{\R_+^d}$. Thus it remains to prove the
assertion for functions
$f\in W_0^{1,p}(\R_+^d)\cap C(\ov{\R_+^d})$.

Let $\phi_n\to f$ in $W_0^{1,p}(\R_+^d) $ with
$\phi_n\in C_{\rm c}^\infty(\R_+^d)$
for all $n\ge 1$. For all
$x = (x',x_d)\in \R_+^d = \R^{d-1}\times (0,\infty)$,
\begin{align*}
 |\phi_n(x',x_d)| \le \int_0^{x_d} |\partial_d\phi_n (x',y)|\ud y.
\end{align*}
Integrating over $|x'|<1$ and averaging over $0<x_d<\eps$ with $0<\eps<1$, we obtain
\begin{align*}
\frac1\eps \int_0^\eps\int_{\{|x'|<1\}} |\phi_n(x',x_d)|\ud x'\ud x_d
& \le  \frac1\eps \int_0^\eps\int_{\{|x'|<1\}}\int_0^{x_d} |\partial_d \phi_n (x',y)|\ud y \ud x' \ud x_d
\\ & =  \int_0^{\eps}\int_{\{|x'|<1\}} \frac1\eps\int_y^\eps |\partial_d \phi_n (x',y)|\ud x_d \ud x' \ud y
\\ & \le \int_0^{\eps}\int_{\{|x'|<1\}} |\partial_d \phi_n (x',y)| \ud x' \ud y.
\end{align*}
Letting $n\to\infty$, we obtain
\begin{align*}
\frac1\eps \int_0^\eps\int_{\{|x'|<1\}} |f(x',x_d)|\ud x'\ud x_d \le
 \int_0^{\eps}\int_{\{|x'|<1\}} |\partial_d f(x',y)| \ud x' \ud y.
\end{align*}
By uniform continuity of $f$ on the compact set
$\{(x',x_d):\, |x'|\le1,\ 0\le x_d\le1\}$,
the left-hand side converges to $\int_{\{|x'|<1\}}|f(x',0)|\ud x'$
as $\eps\downarrow0$.
The right-hand side tends to $0$ by the absolute continuity of the
integral, using H\"older's inequality when $p>1$. Hence
\begin{align*}
\int_{\{|x'|<1\}} |f(x',0)|\ud x' = 0.
\end{align*}

This implies that $f(x',0)=0$ for almost all $x'$ with $|x'|<1$, and hence,
by continuity, for all such $x'$.

Now let $a\in \R^{d-1}$ be arbitrary, and define
$$
g_a(x',x_d):=f(x'+a,x_d), \quad (x',x_d)\in \ov{\R_+^d}.
$$
Then $g_a\in W_0^{1,p}(\R_+^d)\cap C(\ov{\R_+^d})$. Applying the preceding
argument to $g_a$, we obtain
$$
g_a(x',0)=0, \quad |x'|<1,
$$
or equivalently,
$$
f(y',0)=0, \quad |y'-a|<1.
$$
Since $a\in \R^{d-1}$ was arbitrary, it follows that
$f(x',0)=0$ for all $x'\in \R^{d-1}$.
\end{proof}

\begin{theorem}\label{thm:W01pRd} For all $1\le p<\infty$ we have
$$ W_0^{1,p}(\R^d) = W^{1,p}(\R^d).$$
\end{theorem}

\begin{proof} Clearly, $W_0^{1,p}(\R^d)\subseteq W^{1,p}(\R^d)$. To prove the reverse inclusion, we must show that every $f\in W^{1,p}(\R^d)$ can be approximated in $W^{1,p}(\R^d)$ by functions in $C_{\rm c}^\infty(\R^d)$.

Let $(\psi_n)_{n\ge1}$ be the functions of Lemma
\ref{lem:smooth-cutoff}. By Proposition \ref{prop:Cm-mult-Wmp},
$\psi_nf\in W^{1,p}(\R^d)$ and
$
\partial_j(\psi_nf) = \psi_n\partial_jf+f\partial_j\psi_n$.
By dominated convergence,
$\psi_nf\to f$ and $\psi_n\partial_jf\to\partial_jf$
in $L^p(\R^d)$. Moreover,
$$
\|f\partial_j\psi_n\|_p
\le
\frac{C_0}{n}\|f\|_p
\to0.
$$
Hence $\psi_nf\to f$ in $W^{1,p}(\R^d)$.

Accordingly, it suffices to prove that every compactly supported function in $W^{1,p}(\R^d)$ can be approximated, in the norm of $W^{1,p}(\R^d)$, by test functions. This was accomplished in Proposition \ref{prop:local-approx}.
\end{proof}

With the same method of proof one obtains the following density result.

\begin{theorem}\label{thm:test-dense-Wkp}
 For all $k\in\N$ and $1\le p<\infty$ the space $C_{\rm c}^\infty(\R^d)$ is dense in $W^{k,p}(\R^d)$.
\end{theorem}
\begin{proof}
The cutoffs constructed in Lemma \ref{lem:smooth-cutoff}
satisfy $\psi_n(x)=\psi_1(x/n)$, and hence
$\n\partial^\gamma\psi_n\n_\infty
=n^{-|\gamma|}\n\partial^\gamma\psi_1\n_\infty$.
The Leibniz rule and dominated convergence therefore give
$\psi_nf\to f$ in $W^{k,p}(\R^d)$ for every
$f\in W^{k,p}(\R^d)$.

For compactly supported $f\in W^{k,p}(\R^d)$, its standard
mollifications $\eta_\eps*f$ belong to $C_{\rm c}^\infty(\R^d)$.
By Proposition \ref{prop:YoungCk2} and Theorem \ref{thm:approx-identity},
$$
\partial^\alpha(\eta_\eps*f)
=\eta_\eps*\partial^\alpha f
\to\partial^\alpha f
\ \ \hbox{in }L^p(\R^d),\quad |\alpha|\le k.
$$
Thus $\eta_\eps*f\to f$ in $W^{k,p}(\R^d)$, proving density.
\end{proof}

\subsection{Extension Operators}\label{subsec:extension}

Let $k\ge 1$ be an integer.
We say that $D$ has a {\em $C^k$-boundary}\index{boundary!$C^k$} if for every $x_0\in\partial D$
there exist open sets $U,V\subseteq\R^d\!$, with $x_0\in U$,
and a $C^k$-diffeomorphism $\varrho:U\to V$ with the following properties:
 \begin{enumerate}[label={\rm(\roman*)}, leftmargin=*]
 \item\label{it:Ckdomain1} $\varrho(D\cap U) = \{x\in V:\, x_d>0\}$;
 \item\label{it:Ckdomain2} $\varrho(\partial D\cap U) = \{x\in V:\, x_d=0\}$;
 \item\label{it:Ckdomain3} there exists a constant $C>0$ such that
 $$ C^{-1} \le |\det(D\varrho(x))|\le C, \quad x\in U,$$
 where $D\varrho$ is the total derivative of $\varrho$.
 \end{enumerate}
See Figure \ref{fig:DU}.

\begin{proposition}[Local invariance under $C^k$-diffeomorphisms]
\label{prop:diffeo-Wkp}
Let $U,V\subseteq\R^d$ be open, let $k\ge 1$ be an integer, and let
$\varrho:U\to V$ be a $C^k$-diffeomorphism. Let $G\Subset U$ be open, and set
$G':=\varrho(G)\subseteq V$. For $u:G\to\K$ define
$v:G'\to\K$ by $v:=u\circ\varrho^{-1}$. Fix $1\le p<\infty$.
\begin{enumerate}[label={\rm(\arabic*)}, leftmargin=*]
\item\label{it:diffeo-Wkp1}
We have $u\in W^{k,p}(G)$ if and only if $v\in W^{k,p}(G')$, and there exists a constant $C\ge 1$ such that
\begin{equation}\label{eq:diffeo-p-local}
C^{-1}\|v\|_{W^{k,p}(G')}\le \|u\|_{W^{k,p}(G)}\le C\|v\|_{W^{k,p}(G')}.
\end{equation}
The constant $C$ depends only on $p,k,G$ and bounds for the derivatives of $\varrho$ and $\varrho^{-1}$ up to order $k$
on $\overline G$ and $\overline{G'}$, respectively, and on positive lower and upper bounds for $|\det D\varrho|$ on $\overline G$.

\item\label{it:diffeo-Wkp2}
The composition operator $T:u\mapsto u\circ\varrho^{-1}$ restricts to an isomorphism
$$
T:W^{1,p}_0(G)\,\cong\,W^{1,p}_0(G').
$$
\end{enumerate}
\end{proposition}

\begin{proof}
Since $G\Subset U$, the set $\overline G$ is compact. As $\varrho$ is a
homeomorphism, $\overline{G'}=\varrho(\overline G)$ is compact and contained in
$V$. Therefore all partial derivatives of $\varrho$ up to order $k$ are bounded
on $\overline G$, and all partial derivatives of $\varrho^{-1}$ up to order $k$
are bounded on $\overline{G'}$. Moreover, by continuity there exist constants
$0<m\le M$ such that
$$
m\le |\det D\varrho(x)|\le M, \quad x\in\overline G,
$$
where $D\varrho$ is the total derivative of $\varrho$.

\smallskip
\ref{it:diffeo-Wkp1}:\ We argue by induction on $k$.

\smallskip
{\em Step 1} --
By the change of variables $y=\varrho(x)$ and the bounds on $|\det D\varrho|$ we obtain
$$
\|u\|_{L^p(G)}^p=\int_G |v(\varrho(x))|^p\ud x
=\int_{G'} |v(y)|^p\,|\det D\varrho^{-1}(y)|\ud y
\le m^{-1}\|v\|_{L^p(G')}^p,
$$
and similarly $\|v\|_{L^p(G')}\lesssim \|u\|_{L^p(G)}$.

\smallskip
{\em Step 2} --
By the substitution rule for weak derivatives (Proposition~\ref{prop:substitutionrule}), applied respectively to $\varrho$ on $G$ and to $\varrho^{-1}$ on $G'$, we have
$$
\partial_i u=\sum_{j=1}^d a_{ij}\,(\partial_j v)\circ\varrho,
\qquad a_{ij}:=\frac{\partial\varrho_j}{\partial x_i},
$$
and
$$
\partial_j v=\sum_{i=1}^d b_{ji}\,(\partial_i u)\circ\varrho^{-1},
\qquad b_{ji}:=\frac{\partial(\varrho^{-1})_i}{\partial y_j}.
$$
Since the coefficients $a_{ij}$ are bounded on $G$, the displayed formula gives
$$
\|\partial_i u\|_{L^p(G)} \lesssim \sum_{j=1}^d \|(\partial_j v)\circ\varrho\|_{L^p(G)}
\lesssim \sum_{j=1}^d \|\partial_j v\|_{L^p(G')},
$$
where the last estimate is again obtained by a change of variables. Together with Step~1 this gives
$\|u\|_{W^{1,p}(G)}\lesssim \|v\|_{W^{1,p}(G')}$.
The reverse inequality follows in the same way from the formula for $\partial_j v$.

\smallskip
{\em Step 3} --
Assume now that \eqref{eq:diffeo-p-local} holds with $k$ replaced by $k-1$.
By the characterisation that $u\in W^{k,p}(G)$ if and only if $u\in W^{1,p}(G)$ and $\partial_i u\in W^{k-1,p}(G)$ for $i=1,\dots,d$,
it suffices to estimate $\|\partial_i u\|_{W^{k-1,p}(G)}$ in terms of $v$.
From Step~2 we have $\partial_i u=\sum_j a_{ij}(\partial_j v)\circ\varrho$.
Each $a_{ij}$ belongs to $C^{k-1}(\overline G)$, and therefore has bounded derivatives up to order $k-1$ on $G$.
By Proposition \ref{prop:Cm-mult-Wmp}, multiplication by each coefficient
$a_{ij}\in C^{k-1}(\overline G)$ is bounded on $W^{k-1,p}(G)$. Combining
this with the induction hypothesis applied to $\partial_j v$ gives
$$
\|\partial_i u\|_{W^{k-1,p}(G)} \lesssim \sum_{j=1}^d \|(\partial_j v)\circ\varrho\|_{W^{k-1,p}(G)}
\lesssim \sum_{j=1}^d \|\partial_j v\|_{W^{k-1,p}(G')}.
$$
Together with the already proved $W^{1,p}$-estimate this gives $\|u\|_{W^{k,p}(G)}\lesssim \|v\|_{W^{k,p}(G')}$.
The converse estimate is obtained analogously, using the representation of $\partial_j v$ in terms of $(\partial_i u)\circ\varrho^{-1}$.
This completes the proof of \ref{it:diffeo-Wkp1}.

\smallskip
\ref{it:diffeo-Wkp2}:\ Let $u\in W^{1,p}_0(G)$, and choose functions $u_n\in C_{\rm c}^\infty(G)$ with
$u_n\to u$ in $W^{1,p}(G)$. The functions
$ v_n:=u_n\circ\varrho^{-1}$ belong to $W^{1,p}(G')$ by \ref{it:diffeo-Wkp1}, and we have
$\supp(v_n)=\varrho(\supp u_n)\Subset G'$.
Choose an open set $G_n'\Subset G'$ such that $\supp(v_n)\subseteq G_n'$.
Then $v_n$ vanishes outside $G_n'$, and therefore
$v_n\in W^{1,p}_0(G')$ by Proposition \ref{prop:local-approx}.

Since $u_n\to u$ in $W^{1,p}(G)$, another application of
\ref{it:diffeo-Wkp1} gives
$u_n\circ\varrho^{-1}\to u\circ\varrho^{-1}$ in $W^{1,p}(G')$.
As $W^{1,p}_0(G')$ is closed in $W^{1,p}(G')$, it follows that
$u\circ\varrho^{-1}\in W^{1,p}_0(G')$. The converse inclusion follows by
applying the same argument to $\varrho^{-1}$.
\end{proof}

\begin{figure}
\begin{center}
\begin{tikzpicture}[scale=0.7]

\filldraw[thick, nearly transparent] (-6,0) .. controls (-1.5,5) and (0,3) .. (-1.5,0);
\filldraw[thick, nearly transparent] (-1.5,0) .. controls (-3,-3) and (-8.5,-3) .. (-6,0);
\draw (-4.5,2.3) node {$D$};
\filldraw[thick, nearly transparent] (-1.2,2.8) circle (40pt);
\draw (0.3,3.8) node {$U$};
\draw (-1.7,2.4) node {$D\cap U$};
\draw (3,3) node {$\varrho$};
\draw [->] (0.8,3) arc (86:42:150pt);
\draw (4,-1.6) -- (4,2.6);
\draw (2.4,0) -- (7.6,0);
\filldraw[thick, nearly transparent] (6,0) circle (30pt);
\filldraw[thick, nearly transparent] (7.05,0) arc (0:180:30pt) -- cycle;
\draw (7.2,-0.8) node {$V$};
\draw (6,0.3) node {$\varrho(D\cap U)$};
\end{tikzpicture}
\caption{The definition of a $C^k$-boundary \label{fig:DU}}
\end{center}
 \end{figure}

Our next aim is to construct an extension operator for bounded domains with a
$C^k$-boundary. We start with two results for the half-space
$$
\R_+^d:=\{x\in\R^d:\, x_d>0\}.
$$

\begin{lemma}[Approximation on the half-space]\label{lem:halfspace-approx}
Let $k\in\N$ and let $1\le p<\infty$. For every $f\in W^{k,p}(\R_+^d)$ there
exists a sequence $f_n\in C_{\rm c}^\infty(\R^d)$ such that
$$
f_n|_{\R_+^d}\to f\quad\hbox{in }W^{k,p}(\R_+^d).
$$
\end{lemma}

\begin{proof}
Let $\psi\in C_{\rm c}^\infty(\R^d)$ satisfy
$\one_{B(0;1)}\le\psi\le\one_{B(0;2)}$, and put
$\psi_n(x):=\psi(x/n)$.
By the product rule and dominated convergence,
$\psi_nf\to f$ in $W^{k,p}(\R_+^d)$.
We may therefore assume that $f$ has bounded support.

For $|\alpha|\le k$, let $g_\alpha\in L^p(\R^d)$ denote the zero
extension of $\partial^\alpha f$.
Choose a nonnegative function $\eta\in C_{\rm c}^\infty(B(0;1))$
with $\int_{\R^d}\eta\ud x=1$, put $\eta_\eps(x):=\eps^{-d}\eta(x/\eps)$, and define
$$
 f_\eps(x):=(\eta_\eps*g_0)(x+2\eps e_d),
$$
where $e_d$ is the $d$th coordinate vector.
Then $f_\eps\in C_{\rm c}^\infty(\R^d)$.
By Proposition \ref{prop:YoungCk2},
applied locally,
$$
\partial^\alpha f_\eps(x)
=
(\eta_\eps*g_\alpha)(x+2\eps e_d),
\quad x\in\R_+^d,\ |\alpha|\le k.
$$
It follows that
\begin{align*}
\|\partial^\alpha f_\eps-\partial^\alpha f\|_{L^p(\R_+^d)}
&\le
\|(\eta_\eps*g_\alpha)(\cdot+2\eps e_d)-g_\alpha\|_{L^p(\R^d)}
\\
&\le
\|\eta_\eps*g_\alpha-g_\alpha\|_{L^p(\R^d)}
+
\|g_\alpha(\cdot+2\eps e_d)-g_\alpha\|_{L^p(\R^d)}
\to 0
\end{align*}
as $\eps\downarrow0$, by Theorem \ref{thm:approx-identity}
and the $L^p$-continuity of translations.
Thus $f_\eps|_{\R_+^d}\to f$ in $W^{k,p}(\R_+^d)$,
which proves the desired approximation.
\end{proof}

With this result at hand, we can prove the following proposition.
Here and below, the support of a locally integrable function means its
essential support.

\begin{proposition}[Reflection extension on the half-space]
\label{prop:halfspace-extension}
Let $k\ge 1$ be an integer.
Denote by
$L_{\rm loc}^1(\overline{\R_+^d})$ the space of functions on
$\R_+^d$ which are integrable on $K\cap\R_+^d$ for every compact
set $K\subseteq\R^d$. There exists a linear mapping
$$
E_0:L_{\rm loc}^1(\overline{\R_+^d})
  \to L_{\rm loc}^1(\R^d)
$$
such that:
\begin{enumerate}[label={\rm(\roman*)}, leftmargin=*]
\item $(E_0f)|_{\R_+^d}=f$ for all
$f\in L_{\rm loc}^1(\overline{\R_+^d})$;

\item for every $1\le p<\infty$ and every
$\ell=0,1,\dots,k$ there exists a constant
$C_{\ell,p}\ge0$ such that
$$
\|E_0f\|_{W^{\ell,p}(\R^d)}
\le C_{\ell,p}\|f\|_{W^{\ell,p}(\R_+^d)}
$$
for all $f\in W^{\ell,p}(\R_+^d)$;

\item if $f\in L_{\rm loc}^1(\overline{\R_+^d})$ and $\supp(f)\subseteq W\cap \ov{\R_+^d}$
for some open set $W\subseteq\R^d$, then
$$
\supp(E_0f)\subseteq \ov{W^\sharp},
$$
where, with $s_j(y',y_d):=(y',-y_d/j)$,
$$
W^\sharp:=(W\cap \ov{\R_+^d})\cup \bigcup_{j=1}^{k+1} s_j(W\cap \ov{\R_+^d}).
$$
If, moreover, $f\in C^k(\ov{\R_+^d})$ and $W^\sharp\Subset\R^d$, then
$E_0f\in C_{\rm c}^k(\R^d)$.
\end{enumerate}
\end{proposition}

The reflection formula in the proof below may enlarge the support in the normal direction:
a point in the support above the boundary may give rise to finitely many reflected
points below it. The sets $W^\sharp$ keep track of this;
later we will choose the boundary charts so that the reflected support remains inside
the coordinate neighbourhood.

\begin{proof}
Choose scalars $c_j$ in such a way that
$$ \sum_{j=1}^{k+1} (-j)^m c_j = 1, \quad m = 0,1,\dots,k.$$
By the theory of Vandermonde determinants, this system of $k+1$
equations has a unique solution $(c_1,\dots,c_{k+1})$.

Define the linear mapping $E_0: L_{\rm loc}^1 (\overline{\R_+^d})\to L_{\rm loc}^1(\R^d)$ by
\begin{align*}
 E_0 f(x):= \begin{cases}
       f(x), & x_d>0, \\
       \displaystyle\sum_{j=1}^{k+1}c_j f(x_1,\dots,x_{d-1},-jx_d), & x_d<0.
      \end{cases}
\end{align*}
The value on the set $\{x_d=0\}$ is immaterial; in the $C^k$ case below
we choose the continuous representative.

If $f\in L_{\rm loc}^1(\ov{\R_+^d})$ and
$\supp(f)\subseteq W\cap \ov{\R_+^d}$, the defining formula for $E_0$ immediately gives
$\supp(E_0f)\subseteq\ov{W^\sharp}$.
We next check the asserted regularity under the additional assumption that
$f\in C^k(\overline{\R_+^d})$.
Indeed, for $x_d<0$ and
$|\alpha|\le k$,
$$
  \partial^\alpha E_0f(x)
  =
  \sum_{j=1}^{k+1}
   c_j(-j)^{\alpha_d}
   (\partial^\alpha f)(x_1,\ldots,x_{d-1},-jx_d).
$$
Hence, as $x_d\uparrow0$,
$$
  \partial^\alpha E_0f(x',x_d)
  \to
  \Bigl(\sum_{j=1}^{k+1}c_j(-j)^{\alpha_d}\Bigr)
  \partial^\alpha f(x',0)
  =
  \partial^\alpha f(x',0),
$$
because $\alpha_d\le|\alpha|\le k$. Thus the derivatives up to order $k$
match across the hyperplane $\{x_d=0\}$, and $E_0f\in C^k(\R^d)$. If, moreover, $W^\sharp\Subset\R^d$, then
the support inclusion above will imply $E_0f\in C_{\rm c}^k(\R^d)$.

Fix $0\le\ell\le k$ and let
$f\in C^k(\ov{\R_+^d})\cap W^{\ell,p}(\R_+^d)$.
For every multi-index $\alpha$ with $|\alpha|\le\ell$,
\begin{align*}
\|\partial^\alpha E_0f\|_{L^p(\R^d)}^p
& = \int_{\R_+^d} |\partial^\alpha f|^p\ud x
 + \int_{\complement\R_+^d}
 \Bigl|\sum_{j=1}^{k+1} (-j)^{\alpha_d}c_j
 \partial^\alpha f(x_1,\dots,x_{d-1},-jx_d)\Bigr|^p\ud x
\\* & \le C_{\ell,p}^p \|\partial^\alpha f\|_{L^p(\R_+^d)}^p.
\end{align*}
Summing over all multi-indices $|\alpha|\le \ell$, we obtain the bound
\begin{align}\label{eq:boundE0f}
\|E_0f\|_{W^{\ell,p}(\R^d)} \le C_{\ell,p}\|f\|_{W^{\ell,p}(\R_+^d)}.
\end{align}
By Lemma \ref{lem:halfspace-approx}, this estimate extends to all
$f\in W^{\ell,p}(\R_+^d)$. Indeed, let $f_n\in C_{\rm c}^\infty(\R^d)$ be such that
$f_n|_{\R_+^d}\to f$ in $W^{\ell,p}(\R_+^d)$. Then the functions $E_0(f_n|_{\R_+^d})$
form a Cauchy sequence in $W^{\ell,p}(\R^d)$ by the estimate just proved. Passing to a subsequence if necessary, we may assume that $f_n|_{\R_+^d}\to f$ almost
everywhere on $\R_+^d$, and then by the explicit formula for $E_0$ also
$E_0(f_n|_{\R_+^d})\to E_0f$ almost everywhere on $\R^d$. Hence the limit in $W^{\ell,p}(\R^d)$ is $E_0f$.
It follows that the bound \eqref{eq:boundE0f} extends to $f$.
\end{proof}

We are now ready to state and prove the main theorem of this section.

\begin{theorem}[Extension operator]\label{thm:Sobolev-extension}\index{extension!operator}Let $k\ge 1$ be an integer, and let $D$ be bounded and have a $C^k$-boundary. Then there exists a linear
mapping
$$
 E:L^1(D)\to L_{\rm loc}^1(\R^d)
$$
with the following properties:
\begin{enumerate}[label={\rm(\roman*)}, leftmargin=*]
 \item for all $f\in L^1(D)$ we have $(Ef)|_D=f$;
 \item for all $1\le p<\infty$ and $\ell=0,1,\dots,k$ there exists a
 constant $C\ge0$ such that for all $f\in W^{\ell,p}(D)$ we have
 $Ef\in W^{\ell,p}(\R^d)$ and
 $$
  \n Ef\n_{W^{\ell,p}(\R^d)}
  \le C\n f\n_{W^{\ell,p}(D)}.
 $$
\end{enumerate}
\end{theorem}

\begin{proof}
For each $x\in\partial D$, choose open sets $U_x,V_x$ and a
$C^k$-diffeomorphism $\varrho_x:U_x\to V_x$ as in the definition of a
$C^k$-boundary. Thus, in particular,
$\varrho_x(D\cap U_x)=V_x\cap\R_+^d$ and $\varrho_x(\partial D\cap U_x)=V_x\cap\partial\R_+^d$.

Writing
$$
s_j(y',y_d):=(y',-y_d/j),\quad j=1,\dots,k+1,
$$
and using the fact that the maps $s_j$ fix the hyperplane
$\{y_d=0\}$, we may choose $W_x$ to be an open ball centred at $\varrho_x(x)$
with $\ov W_x\subseteq V_x$.
Each $s_j$ fixes this centre and does not increase distances,
so $W_x^\sharp\subseteq W_x\Subset V_x$.

Put $U'_x:=\varrho_x^{-1}(W_x)$. By compactness of $\partial D$, finitely
many of the sets $U'_x$ cover $\partial D$;
relabel the corresponding objects as $U_m,V_m,W_m,U'_m,\varrho_m$, $m=1,\dots,M$. Choose open sets
$V'_m$ such that
$W_m^\sharp\Subset V'_m\Subset V_m.$
Adding an open set $U'_0\Subset D$, we may arrange that
$\ov D\subseteq\bigcup_{m=0}^M U'_m.$
Choose nonnegative functions $\eta_m\in C_{\rm c}^\infty(U'_m)$ such that
$\sum_{m=0}^M\eta_m\equiv1$ on $D$.

Let $f\in L^1(D)$ and put $f_m:=\eta_mf$. Let $h_0$ be the zero
extension of $f_0$ to $\R^d$. For $m=1,\dots,M$, let $g_m$ be the zero
extension to $\R_+^d$ of
$f_m\circ\varrho_m^{-1}$ on $W_m\cap\R_+^d$,
and put $G_m:=\varrho_m^{-1}(V'_m)$. If $E_0$ is the operator of
Proposition~\ref{prop:halfspace-extension}, let $h_m$ be the zero
extension to $\R^d$ of $(E_0g_m)\circ\varrho_m$ on $G_m$.
Since $\supp(g_m)\subseteq W_m\cap\ov{\R_+^d}$, Proposition
\ref{prop:halfspace-extension} gives
$$
\supp(E_0g_m)\subseteq \ov{W_m^\sharp}\Subset V'_m.
$$
Thus $(E_0g_m)\circ\varrho_m$ is supported compactly inside $G_m$;
this will justify the zero extension in the Sobolev estimates below.
 Moreover, on $D$ we have $h_m=f_m$ almost everywhere
for $m=0,\dots,M$.
Indeed, for $m\ge1$, on $D\cap G_m$ this follows from
$(E_0g_m)|_{\R_+^d}=g_m$ and the definition of $g_m$.
On $D\setminus G_m$ both functions vanish, since
$D\cap U'_m\subseteq G_m$ and $f_m=0$ outside $U'_m$.
Hence
$$
Ef:=\sum_{m=0}^M h_m
$$
defines a linear mapping from $L^1(D)$ into $L_{\rm loc}^1(\R^d)$ and
$(Ef)|_D=f$.

It remains to prove the Sobolev estimates. Fix $1\le p<\infty$,
$0\le\ell\le k$, and $f\in W^{\ell,p}(D)$. By the product rule,
$f_m\in W^{\ell,p}(D)$. Since $\supp(f_0)\Subset D$, Lemma
\ref{lem:der-zeroextension} gives
$$
\n h_0\n_{W^{\ell,p}(\R^d)}
\lesssim \n f\n_{W^{\ell,p}(D)}.
$$
For $m\ge1$, Proposition~\ref{prop:diffeo-Wkp}, with the
change-of-variables formula when $\ell=0$, gives
$$
\n f_m\circ\varrho_m^{-1}\n_{W^{\ell,p}(W_m\cap\R_+^d)}
\lesssim \n f\n_{W^{\ell,p}(D)}.
$$
The support of this function is contained in
$F_m:=\varrho_m(\supp\eta_m)\cap\R_+^d$.
This set is closed in $\R_+^d$ and contained in
$W_m\cap\R_+^d$. Hence Lemma~\ref{lem:der-zeroextension} shows that
$g_m\in W^{\ell,p}(\R_+^d)$ and
$$
\n g_m\n_{W^{\ell,p}(\R_+^d)}
=
\n f_m\circ\varrho_m^{-1}
  \n_{W^{\ell,p}(W_m\cap\R_+^d)}
\lesssim
\n f\n_{W^{\ell,p}(D)}.
$$
Proposition~\ref{prop:halfspace-extension} now gives
$$
\n E_0g_m\n_{W^{\ell,p}(\R^d)}
\lesssim \n f\n_{W^{\ell,p}(D)}.
$$
Since $\supp(E_0g_m)\subseteq\ov{W_m^\sharp}\Subset V'_m$,
Proposition~\ref{prop:diffeo-Wkp} and Lemma
\ref{lem:der-zeroextension} imply, again using change of variables when
$\ell=0$,
$$
\n h_m\n_{W^{\ell,p}(\R^d)}
\lesssim \n f\n_{W^{\ell,p}(D)}.
$$
Summing over $m=0,\dots,M$ gives
\[
\n Ef\n_{W^{\ell,p}(\R^d)}
\le C\n f\n_{W^{\ell,p}(D)}.
\qedhere\]
\end{proof}

The boundedness of $D$ is used only to obtain finitely many boundary
charts and hence uniform global estimates. The same construction applies
locally to arbitrary $C^k$-domains; for global estimates on unbounded
domains one needs suitable uniformity assumptions on the boundary charts.

\begin{theorem}[Density of smooth functions up to the boundary]\label{thm:trace-W1p}
Let $1\le p<\infty$ and let $k\ge 1$. If $D$ is bounded with $C^k$-boundary,
then for every $f\in W^{k,p}(D)$ there exists a sequence
$f_n\in C^\infty(\R^d)$ such that
$$
f_n|_D\to f\quad\hbox{in }W^{k,p}(D).
$$
In particular, $C^\infty(\ov D)$ is dense in $W^{k,p}(D)$.
\end{theorem}

\begin{proof}
Let $E$ be the extension operator of
Theorem~\ref{thm:Sobolev-extension}. For $f\in W^{k,p}(D)$ we have
$Ef\in W^{k,p}(\R^d)$. By Theorem~\ref{thm:test-dense-Wkp}, there exist
$f_n\in C_{\rm c}^\infty(\R^d)$ such that
$f_n\to Ef$ in $W^{k,p}(\R^d)$. Restriction to $D$ gives
$f_n|_D\to f$ in $W^{k,p}(D)$.
\end{proof}

\subsection{Bessel Potential Spaces}\label{sec:Hs}

A function $f\in L^p(D)$, with $1\le p\le \infty$, is said to admit a {\em weak $L^p$-Laplacian}\index{weak!Laplacian}\index{Laplace operator!weak} if there exists a function
 $g\in L^p(D)$ such that for all test functions $\phi\in C_{\rm c}^\infty(D)$ we have
 $$ \int_{D} f \Delta \phi\ud x = \int_{D} g\phi\ud x.$$
This defines a linear operator $\Delta_{p,D}$ in $L^p(D)$, the {\em weak $L^p$-Laplacian on $D$}, with domain
$$\Dom(\Delta_{p,D}):= \{f\in L^p(D) :\, \hbox{$f$ admits a weak $L^p$-Laplacian on $D$}\}.$$
By the argument of Example \ref{ex:weak-Laplace-closable}, this operator is closed.
In what follows we will denote weak Laplacians simply by $\Delta$.

As an application of the Plancherel theorem, the next theorem shows,
among other things, that the domain of the weak Laplacian in $L^2(\R^d)$ equals $W^{2,2}(\R^d)$.

\begin{theorem}[Fourier-analytic characterisation, weak Laplacian]\label{thm:H1-Sob} Let $k\ge 0$ be a nonnegative integer. For a function $f \in L^2(\R^d)$ the following assertions are equivalent:
\begin{enumerate}[label={\rm(\arabic*)}, leftmargin=*]
\item\label{it:H1-Sob1} $f$ belongs to $W^{k,2}(\R^d)$;
\item\label{it:H1-Sob2} $\xi \mapsto (1 + |\xi|^2)^{k/2} \wh{f}(\xi)$ belongs to $ L^2(\R^d).$
\end{enumerate}
Moreover, $$f\mapsto \big\n \xi \mapsto (1 + |\xi|^2)^{k/2} \wh{f}(\xi)\big\n_2$$ defines an equivalent norm on $W^{k,2}(\R^d)$.
For $k=2$, \ref{it:H1-Sob1} and \ref{it:H1-Sob2} are equivalent to:
\begin{enumerate}[label={\rm(\arabic*)}, leftmargin=*]\setcounter{enumi}{2}
\item\label{it:H1-Sob3} $f$ admits a weak Laplacian in $L^2(\R^d)$.
\end{enumerate}
Moreover,
$\Dom(\Delta) = W^{2,2}(\R^d)$ and
$f\mapsto \n f\n_2+\n \Delta f\n_2$ defines an equivalent norm on $W^{2,2}(\R^d)$.
\end{theorem}

\begin{definition}[The Bessel potential spaces $H^s(\R^d)$]
For real numbers $s\ge 0$
the subspace of all $f \in L^2(\R^d)$ such that
$$\xi \mapsto (1 + |\xi|^2)^{s/2} \wh{f}(\xi)$$ belongs to $ L^2(\R^d)$
is denoted by $H^s(\R^d)$ and is called the {\em Bessel potential space}\index{Bessel potential space} with smoothness exponent $s$.\index{$HS$@$H^s(\R^d)$}
With respect to the norm
\begin{align}\label{eq:norm-Hs} \n f\n_{H^s(\R^d)} := \big\n \xi \mapsto (1 + |\xi|^2)^{s/2} \wh{f}(\xi)\big\n_{L^2(\R^d)}
\end{align}
this space is a Hilbert space. The easy verification is left as an exercise.
\end{definition}

For noninteger values of $s$, the spaces $H^s(\R^d)$ play an important role in the regularity theory for partial differential equations.

The equivalence of \ref{it:H1-Sob1} and \ref{it:H1-Sob2} of Theorem \ref{thm:H1-Sob} asserts:

\begin{theorem}[$W^{k,2}(\R^d) = H^{k}(\R^d)$]\label{thm:W=H}
If $k\ge 0$ is a nonnegative integer, then
$$ W^{k,2}(\R^d) = H^{k}(\R^d)$$
with equivalence of norms.
\end{theorem}

The proof of Theorem \ref{thm:H1-Sob} relies on the following lemma, where for vectors $\xi\in \R^d$ and multi-indices $\alpha\in \N^d$ we use the shorthand notation
$\xi^\al := \xi_1^{\al_1}\cdots \xi_d^{\al_d}.$

 \begin{lemma}\label{lem:FT-partial} A function $f\in L^2(\R^d)$ admits a weak
 derivative $\partial^\alpha f$ belonging to $L^2(\R^d)$ if and only if $\xi\mapsto \xi^\alpha \wh f(\xi)$ belongs to
 $L^2(\R^d)$, and in that case we have
 $$\wh{\partial^\alpha f}(\xi) = i^{|\alpha|}\xi^\alpha \wh f(\xi)$$
 for almost all $\xi\in \R^d\!$.
 \end{lemma}
\begin{proof}
For test functions $\phi \in C_{\rm c}^\infty(\R^d)$ and $\xi\in \R^d$, integration by parts gives
\begin{equation}\label{eq:FT-test-alpha}
\begin{aligned}
\widehat{\partial^\alpha \phi}(\xi) &= \frac1{(2\pi)^{d/2}}
\int_{\R^d} \exp(- i x\cdot\xi) \partial^\alpha \phi(x) \ud x
\\ & = \frac{i^{|\alpha|}\xi^{\alpha}}{(2\pi)^{d/2}}\int_{\R^d} \exp(- i x\cdot\xi) \phi(x)\ud x  =
 i^{|\alpha|}\xi^{\alpha} \widehat{\phi}(\xi).
\end{aligned}
\end{equation}
Also, by differentiating under the integral, we see that $\wh\phi$ is smooth and
\begin{equation}\label{eq:FT-test-alpha1}
\partial^\alpha\wh\phi(\xi)
=
(-i)^{|\alpha|}\widehat{x^\alpha\phi}(\xi).
\end{equation}

`If': \ Suppose that $\xi\mapsto \xi^{\alpha} \widehat{f}(\xi)$ belongs to $L^2(\R^d)$
and denote its inverse Fourier--Plancherel transform by $g_\alpha$.

Using the Plancherel isometry and \eqref{eq:FT-test-alpha}, for real-valued
$\phi\in C_{\rm c}^\infty(\R^d)$ we obtain
\begin{align*}
(-1)^{|\alpha|}\int_{\R^d} f(x)\partial^\alpha\phi(x)\ud x
&= (-1)^{|\alpha|}\iprod{f}{\partial^\alpha\phi}
 = (-1)^{|\alpha|}\iprod{\wh f}{\wh{\partial^\alpha\phi}}
\\ &=
i^{|\alpha|}\int_{\R^d}\xi^\alpha\wh f(\xi)\,
\ov{\wh\phi(\xi)}\ud\xi
= i^{|\alpha|}\int_{\R^d} g_\alpha(x)\phi(x)\ud x.
\end{align*}
Considering real and imaginary parts separately, we extend the identity to complex-valued $\phi \in C_{\rm c}^\infty(\R^d)$.
This shows that $f$ has weak derivative $\partial^\alpha f= i^{|\alpha|}g_\alpha$ in $L^2(\R^d)$.

`Only if': \
Fix a test function $\phi\in C_{\rm c}^\infty(\R^d)$. By \eqref{eq:FT-test-alpha1}, the function $\wh\phi$ and all its
derivatives up to order $|\alpha|$ belong to $L^2(\R^d)$. Hence Theorem
\ref{thm:test-dense-Wkp} permits us to extend the weak derivative
identity from test functions to $\wh\phi$.
 Using Proposition \ref{prop:FT-fFg} and \eqref{eq:FT-test-alpha1}, we obtain
 \begin{align*}
 \int_{\R^d} \wh{\partial^\alpha f}(\xi) \phi(\xi) \ud \xi
 & = \int_{\R^d}  (\partial^\alpha f)(\xi) \wh\phi(\xi) \ud \xi
 = (-1)^{|\al|} \int_{\R^d} f (\xi)\partial^\alpha \wh\phi(\xi)\ud \xi
 \\ & = i^{|\al|}\int_{\R^d}f(\xi)(x^\alpha \phi(x))\wh{\phantom{a}}(\xi)\ud \xi
 = i^{|\alpha|}\int_{\R^d} \xi^\alpha\wh f(\xi) \phi(\xi)\ud \xi.
 \end{align*}
 By Proposition \ref{prop:fundvarcalculus} this implies
 $\wh{\partial^\alpha f}(\xi) = i^{|\alpha|} \xi^\alpha \wh f(\xi)$ for almost all $\xi\in\R^d\!$.
 Since by assumption $\partial^\alpha f\in L^2(\R^d)$, the Plancherel theorem implies that $\wh{\partial^\alpha f}\in L^2(\R^d)$.
This shows that $\xi\mapsto \xi^\alpha \wh f(\xi)$ belongs to $L^2(\R^d)$.
\end{proof}

\begin{proof}[Proof of Theorem \ref{thm:H1-Sob}]
We begin with the proof of the equivalence \ref{it:H1-Sob1}$\Leftrightarrow$\ref{it:H1-Sob2} for general integers $k\ge 0$.
The heart of the matter is contained in the two-sided estimate
\begin{align}\label{eq:est-Wk2}
(1 + |\xi|^2)^{k/2} \eqsim_{d,k} \Bigl(\sum_{|\alpha| \leq k} |\xi^{\alpha}|^2\Bigr)^{1/2}\!,
\end{align}
which follows from the binomial identity. Here, the notation $A\eqsim_{d,k} B$ is shorthand for the existence
of constants $0<c_{d,k}, c_{d,k}'<\infty$,
both depending only on $d$ and $k$, such that $A \le c_{d,k}B$ and $ B \le c_{d,k}'A$.

\smallskip
\ref{it:H1-Sob1}$\Rightarrow$\ref{it:H1-Sob2}:\
First let $f\in W^{k,2}(\R^d)$.
Then $\partial^\alpha f\in L^2(\R^d)$ for all $|\alpha|\le k$.
By Lemma \ref{lem:FT-partial} we have, for almost all $\xi\in\R^d$,
$\widehat{\partial^\alpha f}(\xi)=i^{|\al|}\xi^\alpha \widehat f(\xi)$.
Since $\widehat{\partial^\alpha f}\in L^2(\R^d)$ by Plancherel's theorem, it follows that
$\xi^\alpha \widehat f(\xi)\in L^2(\R^d)$ for all $|\alpha|\le k$.
By \eqref{eq:est-Wk2},
this shows that $\xi\mapsto (1 + |\xi|^2)^{k/2} \wh{f}(\xi)$ belongs to $L^2(\R^d)$, and that its $L^2$-norm is
equivalent to the norm of $f$ in $W^{k,2}(\R^d)$.

\smallskip
\ref{it:H1-Sob2}$\Rightarrow$\ref{it:H1-Sob1}: \ Suppose that $f\in L^2(\R^d)$ is a function with the property that $\xi\mapsto
(1 + |\xi|^2)^{k/2} \wh{f}(\xi)$ belongs to $L^2(\R^d)$.
Then $\xi\mapsto \xi^\alpha \wh{f}(\xi)$ belongs to $L^2(\R^d)$ for all multi-indices
$\alpha\in\N^d$ satisfying $|\alpha|\le k$,
and therefore $f$ has weak derivatives $\partial^\alpha f$ in $L^2(\R^d)$ for all $|\alpha|\le k$
by Lemma \ref{lem:FT-partial}. This shows that $f\in W^{k,2}(\R^d)$.

\smallskip
For $k=2$, \ref{it:H1-Sob1} obviously implies \ref{it:H1-Sob3}.
Conversely, if \ref{it:H1-Sob3} holds (with $\Delta f = g$), then by the same argument as in the `only if' part of Lemma \ref{lem:FT-partial} we find that $\xi\mapsto |\xi|^2 \wh{f}(\xi)$ belongs to $L^2(\R^d)$ (and equals $-\wh g(\xi)$ for almost all $\xi\in\R^d$), and therefore \ref{it:H1-Sob2} holds. The equivalence of norms again follows from the Plancherel theorem and \eqref{eq:est-Wk2}.
\end{proof}

\subsection{Sobolev Inequalities}\label{subsec:Sobolev-inequalities}

We conclude this section with a treatment of the basic Sobolev inequalities. They express the principle that integrability of derivatives improves the integrability, and eventually the regularity, of a function.

\begin{lemma}[Gagliardo]\label{lem:Gagliardo}
Let $d\ge2$ and let $g_1,\ldots,g_d\in L^1(\R^{d-1})$
be nonnegative. Denote by $\widehat x_j$ the vector obtained
from $x$ by deleting its $j$th coordinate. Then
$$
\int_{\R^d}\prod_{j=1}^d g_j(\widehat x_j)^{1/(d-1)}\ud x
\le
\prod_{j=1}^d
\Bigl(\int_{\R^{d-1}}g_j(y)\ud y\Bigr)^{1/(d-1)}.
$$
\end{lemma}

\begin{proof}
For $d=2$, Fubini's theorem gives equality.
Assume the result in dimension $d$ and consider dimension
$d+1$. Write $x=(x',t)\in\R^d\times\R$ and set
$$
G_j(y):=\int_\R g_j(y,t)\ud t, \quad y\in\R^{d-1},\ j=1,\ldots,d.
$$
By Fubini's theorem, $G_j\in L^1(\R^{d-1})$ and
$\n G_j\n_1=\n g_j\n_1$.
Integrating first in $t$ and applying H\"older's inequality
with $d$ factors gives
\begin{align*}
I:=\int_{\R^{d+1}}
       \prod_{j=1}^{d+1}g_j(\widehat x_j)^{1/d}\ud x
&\le\int_{\R^d}g_{d+1}(x')^{1/d}
       \prod_{j=1}^dG_j(\widehat{x'}_j)^{1/d}\ud x'\\
&\le\n g_{d+1}\n_1^{1/d}
 \Bigl(\int_{\R^d}\prod_{j=1}^d
 G_j(\widehat{x'}_j)^{1/(d-1)}\ud x'\Bigr)^{(d-1)/d}\\
&\le\n g_{d+1}\n_1^{1/d}
       \prod_{j=1}^d\n G_j\n_1^{1/d}
 =\prod_{j=1}^{d+1}\n g_j\n_1^{1/d}.
\end{align*}
The second inequality uses H\"older's inequality with
exponents $d$ and $d/(d-1)$, and the third uses the
induction hypothesis.
\end{proof}

\begin{theorem}[Sobolev inequalities on $\R^d$]\label{thm:Sobolev-inequalities-Rd}\index{theorem!Sobolev}
Let $1\le p<\infty$.
\begin{enumerate}[label={\rm(\arabic*)}, leftmargin=*]
\item\label{it:Sobolev-inequalities-Rd1}
If $1\le p<d$ and $p^*:=dp/(d-p)$, then every $f\in W^{1,p}(\R^d)$
belongs to $L^{p^*}\!(\R^d)$ and
$$
\n f\n_{p^*}\le C\n\nabla f\n_p,
\quad f\in W^{1,p}(\R^d),
$$
where $C$ is a constant depending only on $d$ and $p$.
\item\label{it:Sobolev-inequalities-Rd2}
If $p>d$, or if $d=p=1$, then every $f\in W^{1,p}(\R^d)$ is equal almost
everywhere to a function in $C_0(\R^d)$ and
$$\n f\n_\infty\le C\n f\n_{W^{1,p}(\R^d)},$$
where $C$ is a constant depending only on $d$ and $p$.
\end{enumerate}
\end{theorem}

\begin{proof}
In the proof, $C$ is a constant depending on $p$ and $d$ whose numerical value changes from line to line.

\smallskip
\ref{it:Sobolev-inequalities-Rd1}: \
First consider $p=1$.
For $f\in C_{\rm c}^\infty(\R^d)$, the fundamental theorem of
calculus gives
$$
f(x)
=
\int_{-\infty}^{x_j}
\partial_jf(x_1,\dots,x_{j-1},t,x_{j+1},\dots,x_d)\ud t,
$$
and hence
$$
|f(x)|
\le
\int_\R
|\partial_jf(x_1,\dots,x_{j-1},t,x_{j+1},\dots,x_d)|
\ud t.
$$
Multiplying these inequalities and taking the power $1/(d-1)$, we obtain
$$
|f(x)|^{d/(d-1)}
\le
\prod_{j=1}^d
\Bigl(
\int_\R|\partial_jf(x_1,\ldots,x_{j-1},t,x_{j+1},\ldots,x_d)|\ud t
\Bigr)^{1/(d-1)}.
$$
By Lemma \ref{lem:Gagliardo},
$$
\n f\n_{d/(d-1)}^{d/(d-1)}
\le
\prod_{j=1}^d\n\partial_jf\n_1^{1/(d-1)}.
$$
Taking the power $(d-1)/d$ and using the arithmetic--geometric mean
inequality, we conclude that
$$
\n f\n_{d/(d-1)}
\le
\prod_{j=1}^d\n\partial_jf\n_1^{1/d}
\le
\frac1d\sum_{j=1}^d\n\partial_jf\n_1
\le C\n\nabla f\n_1.
$$
By density, this estimate holds for all $f\in W^{1,1}(\R^d)$.

Now let $1<p<d$ and put $r:=p(d-1)/(d-p)$. Considering real and
imaginary parts separately, it suffices to prove the estimate for
real-valued functions $f\in C_{\rm c}^\infty(\R^d)$.
For such functions, since $r>1$, the function $\rho(t):=|t|^r$ is continuously differentiable, with
$\rho'(t)=r|t|^{r-2}t$ and $\rho'(0)=0$. Hence $\rho\circ f=|f|^r$ belongs to
$C_{\rm c}^1(\R^d)\subseteq W^{1,1}(\R^d)$ and
$$
|\nabla(|f|^r)|\le r|f|^{r-1}|\nabla f|
$$
almost everywhere.
Applying the case $p=1$ to $|f|^r$, and using that
$rd/(d-1)=p^*$, we obtain
$$
\n f\n_{p^*}^r
=\n |f|^r\n_{d/(d-1)}
\le Cr\int_{\R^d}|f|^{r-1}|\nabla f|\ud x.
$$
H\"older's inequality with exponents $p'$ and $p$, together with
$(r-1)p'=p^*$, gives
$$
\int_{\R^d}|f|^{r-1}|\nabla f|\ud x
\le \n f\n_{p^*}^{r-1}\n\nabla f\n_p.
$$
If $f\ne0$, dividing by $\n f\n_{p^*}^{r-1}$ proves the assertion;
the case $f=0$ is trivial.
The general case follows by density.

\smallskip
\ref{it:Sobolev-inequalities-Rd2}: \
Let $f\in C_{\rm c}^\infty(\R^d)$, and choose a nonnegative function
$\phi\in C_{\rm c}^\infty(B(0;1))$ with $\int_{\R^d}\phi(u)\ud u=1$. Then
\begin{align}\label{eq:Sobolec-2terms}
f(x)
=
\int_{B(0;1)}\phi(u)f(x-u)\ud u
+
\int_{B(0;1)}\phi(u)(f(x)-f(x-u))\ud u.
\end{align}
Since $\phi$ is bounded, the first term in \eqref{eq:Sobolec-2terms} is bounded by
$C\int_{B(x;1)}|f(y)|\ud y$.
By the fundamental theorem of calculus,
$$
|f(x)-f(x-u)|
\le |u|\int_0^1|\nabla f(x-tu)|\ud t.
$$
The second term in \eqref{eq:Sobolec-2terms} is therefore bounded by
$
C\int_0^1\int_{B(0;1)}
  |u|\,|\nabla f(x-tu)|\ud u\ud t.
$
For fixed $t\in(0,1)$, the change of variables $y=x-tu$ gives
$u=(x-y)/t$, $\ud u=t^{-d}\ud y$, and the condition $|u|<1$
becomes $|x-y|<t$. Hence, by changing the order of integration,
\begin{align*}
C\int_0^1\int_{B(0;1)}
  |u|\,|\nabla f(x-tu)|\ud u\ud t
&=
C\int_0^1\int_{B(x;t)}
  |x-y|t^{-d-1}|\nabla f(y)|\ud y\ud t\\
&=
C\int_{B(x;1)}|x-y||\nabla f(y)|
  \int_{|x-y|}^1t^{-d-1}\ud t\ud y.
\end{align*}
Finally,
$$
\int_{|x-y|}^1t^{-d-1}\ud t
=\frac1d\bigl(|x-y|^{-d}-1\bigr)
\le C|x-y|^{-d},
$$
and therefore the last expression is bounded by
$
C\int_{B(x;1)}|x-y|^{1-d}|\nabla f(y)|\ud y.
$

Combining the two estimates gives the pointwise bound
$$
|f(x)|
\le C\int_{B(x;1)}
\big(|f(y)|+|x-y|^{1-d}|\nabla f(y)|\big)\ud y.
$$

If $p>d$, the function
$y\mapsto |y|^{1-d}\one_{B(0;1)}(y)$ belongs to $L^{p'}(\R^d)$.
The same is true for $d=p=1$, with $p'=\infty$. H\"older's inequality then gives
the asserted estimate for $\n f\n_\infty$. By density, every $f\in W^{1,p}(\R^d)$ is the uniform limit of a sequence in $C_{\rm c}^\infty(\R^d)$, and hence is equal almost
everywhere to a function in $C_0(\R^d)$.
\end{proof}

The extension theorem gives the corresponding results on bounded domains.

\begin{corollary}[Sobolev embeddings on bounded domains]
\label{cor:Sobolev-inequalities-domain}
Let $D\subseteq\R^d$ be bounded with $C^1$-boundary, and let $1\le p\le \infty$.
\begin{enumerate}[label={\rm(\roman*)}, leftmargin=*]
\item\label{it:Sobolev-ineq1} If $1\le p<d$, then $W^{1,p}(D)$ embeds continuously into
$L^q(D)$ for all $1\le q\le p^*$.
\item\label{it:Sobolev-ineq2} If $p=d\ge2$, then $W^{1,p}(D)$ embeds continuously into
$L^q(D)$ for all $1\le q<\infty$.
\item\label{it:Sobolev-ineq3} If $p>d$, or if $d=p=1$, then every $f\in W^{1,p}(D)$ is equal
almost everywhere to a function in $C(\ov D)$, and
$\n f\n_\infty\le C\n f\n_{W^{1,p}(D)}$.
\end{enumerate}
In \ref{it:Sobolev-ineq1} and \ref{it:Sobolev-ineq2}, we have
$\n f\n_q\le C\n f\n_{W^{1,p}(D)}$, where $C$ depends
only on $D,d,p,q$. The constant in \ref{it:Sobolev-ineq1} depends
only on $D,d,p$.
\end{corollary}

\begin{proof}
First suppose that $1\le p<\infty$.
Let $E$ be the extension operator of
Theorem \ref{thm:Sobolev-extension}.
For $f\in W^{1,p}(D)$, we have
$$
\|Ef\|_{W^{1,p}(\R^d)}
\le C\|f\|_{W^{1,p}(D)}.
$$
Applying Theorem \ref{thm:Sobolev-inequalities-Rd} to $Ef$
and restricting to $D$ gives the first assertion for
$q=p^*$ and the third assertion for $d<p<\infty$
and for $d=p=1$.
The remaining values of $q$ in the first
assertion follow from H\"older's inequality.

Suppose that $p=d\ge2$. For $q>d$, put $r:=dq/(d+q)$.
Then $1<r<d$ and $r^*=q$.
Since $D$ has finite measure,  the
continuous inclusion $W^{1,d}(D)\subseteq W^{1,r}(D)$ follows from H\"older's inequality.
Applying the first assertion with exponent $r$, we obtain
$$
\|f\|_{L^q(D)}
\le C\|f\|_{W^{1,r}(D)}
\le C'\|f\|_{W^{1,d}(D)}.
$$
The cases $1\le q\le d$ follow directly from H\"older's inequality.

Finally, let $p=\infty$ and fix a finite exponent $q>d$.
Since $D$ has finite measure, H\"older's inequality gives
the continuous inclusion
$W^{1,\infty}(D)\subseteq W^{1,q}(D)$.
Applying the third assertion with exponent $q$, we obtain
a representative of $f$ in $C(\ov D)$ and the estimate
$$
\n f\n_\infty
\le C\n f\n_{W^{1,q}(D)}
\le C'\n f\n_{W^{1,\infty}(D)}.
$$
This proves the third assertion also for $p=\infty$.
\end{proof}

The Fourier transform gives a particularly short higher-order result.

\begin{theorem}[Sobolev embedding for Bessel potential spaces]
\label{thm:Hs-embedding-C0}
If $s>d/2$, then $H^s(\R^d)$ is continuously embedded in $C_0(\R^d)$.
Consequently, if $k>d/2$ is an integer and $D\subseteq\R^d$ is bounded with
$C^k$-boundary, then $W^{k,2}(D)$ is continuously embedded in $C(\ov D)$.
\end{theorem}

\begin{proof}
For $f\in H^s(\R^d)$, the Cauchy--Schwarz inequality gives
$$
\n\wh f\n_1
\le
\Bigl(\int_{\R^d}(1+|\xi|^2)^{-s}\ud\xi\Bigr)^{1/2}
\n f\n_{H^s(\R^d)}.
$$
The integral is finite if and only if $s>d/2$. The inverse Fourier transform
of $\wh f$ therefore belongs to $C_0(\R^d)$ by the Riemann--Lebesgue lemma and equals $f$ in $L^2(\R^d)$
by the Fourier--Plancherel inversion theorem.

For the final assertion, let $E$ be the extension operator of
Theorem~\ref{thm:Sobolev-extension}. If $f\in W^{k,2}(D)$ with $k>d/2$, then
$Ef\in W^{k,2}(\R^d)$ and
$$
\n Ef\n_{W^{k,2}(\R^d)}
\le C\n f\n_{W^{k,2}(D)}.
$$
By Theorem~\ref{thm:W=H}, we have $Ef\in H^k(\R^d)$, with
$$
\n Ef\n_{H^k(\R^d)}
\le C\n f\n_{W^{k,2}(D)}.
$$
Since $k>d/2$, the first assertion shows that $Ef$ is equal almost
everywhere to a function $F\in C_0(\R^d)$ satisfying
$$
\n F\n_\infty
\le C\n Ef\n_{H^k(\R^d)}
\le C\n f\n_{W^{k,2}(D)}.
$$
As $(Ef)|_D=f$, the restriction $F|_{\ov D}$ belongs to $C(\ov D)$ and is
equal almost everywhere to $f$ on $D$. This proves the asserted continuous
inclusion of $W^{k,2}(D)$ into $C(\ov D)$.
\end{proof}

More generally, if $m\in\N$ and $s>m+d/2$, then every
$f\in H^s(\R^d)$ has a $C^m$-represen\-tative. Indeed, for
$|\al|\le m$ the function $(i\xi)^\al\wh f(\xi)$ belongs to
$L^1(\R^d)$, and Fourier inversion and differentiation under the
integral give the assertion.

\section{The Poisson Problem}\label{sec:Poisson-problem}\index{Poisson!problem}

The results developed above are now applied to the study of the Poisson problem.

\subsection{Dirichlet Boundary Conditions}\label{sec:Poisson-D}\index{Dirichlet!boundary conditions}

We recall from Theorem \ref{thm:W=H} that $H^k(\R^d)= W^{k,2}(\R^d)$ with equivalent norms.
For nonempty open subsets $D$ of $\R^d$ it is customary to {\em define}\index{$H_0^1(D)$}\index{$H^k(D)$}
\begin{align*} H^k(D):= W^{k,2}(D), \qquad H_0^1(D):= W_0^{1,2}(D),
\end{align*}
where $W_0^{1,2}(D)$
is the closure of $C_{\rm c}^\infty(D)$ in $W^{1,2}(D)$.
We endow $H^k(D)$ with the norm
$$ \n f\n_{H^k(D)}^2 := \sum_{|\alpha|\le k} \n \partial^\alpha f\n_2^2,$$
where the summation extends over all multi-indices $\alpha\in\N^d$ of order at most $k$.
This norm is associated with the inner product
$$ \iprod{f}{g}_{H^k(D)} =\sum_{|\alpha|\le k} \iprod{\partial^\alpha f}{\partial^\alpha g}_2$$
and it turns $H^k(D)$ into a Hilbert space.
As a closed subspace of $H^1(D)$, the space $H_0^1(D)$ is a Hilbert space as well.

Let us now consider the {\em Poisson problem with Dirichlet boundary
conditions:}\index{Poisson!problem, Dirichlet boundary conditions}
\begin{equation}\label{eq:Poisson}
\begin{cases} -\Delta u & = f \ \ \hbox{on $D$}, \\ u|_{\partial D} & = 0,
\end{cases}
\end{equation}
where $f\in L^2(D)$ is a given function and $\Delta = \sum_{i=1}^d \partial_i^2$ is the Laplace operator.
Multiplying both sides of the equation \eqref{eq:Poisson} by a test function $\phi\in C_{\rm c}^\infty(D)$
and integrating, we obtain the following integrated version of the problem:
$$ \int_D (\Delta u)\phi\ud x = -\int_D f\phi\ud x,$$
which, after a formal integration by parts (which can be rigorously justified if $u\in C^2(D)$), can be rewritten as
\begin{align}\label{eq:Poisson-int}
\int_D \nabla u \cdot\nabla\phi\ud x = \int_D f\phi\ud x.
\end{align}

This formal derivation justifies the following definition.

\begin{definition}[Weak solutions]\label{def:Poisson-weaksol}
A function $u\in H_0^{1}(D)$ is called a {\em weak solution}\index{weak!solution} of the Poisson problem
with Dirichlet boundary conditions \eqref{eq:Poisson} if
$$
\int_D \nabla u \cdot\nabla \phi\ud x = \int_D f\phi\ud x, \quad \phi\in C_{\rm c}^\infty(D).
$$
\end{definition}

The notion of weak solution makes sense for inhomogeneities $f\in L^2(D)$. In the special case $f\in C(\ov D)$, a {\em classical solution} may be defined as a function $u\in C^2(D)\cap C(\ov D)$ satisfying the equations of the Poisson problem,
$-\Delta u = f$ on $D$ and $u|_{\partial D} = 0$ pointwise. Classical solutions may not exist, however, even when $f\in C_{\rm c}(D)$; a counterexample is presented at the end of this section.

The advantage of working with weak solutions is that the integrated
equation in \eqref{eq:Poisson-int} involves only first derivatives and both integrals in \eqref{eq:Poisson-int} can be interpreted as inner products. This makes the problem of finding weak solutions amenable to Hilbert space methods.
The requirement that $u$ be in $H_0^1(D)$ implements the boundary condition $u|_{\partial D} = 0$, as evidenced by
Theorem \ref{thm:H01-cont}.

Having adopted membership in $H_0^1(D)$ as a way of imposing
Dirichlet boundary conditions, one may still try to look for
solutions in $H_0^1(D)\cap H^2(D)$. It can be shown that this works (that is, a unique weak solution in the sense of Definition \ref{def:Poisson-weaksol} belonging to this space exists) if $D$ is assumed to be bounded with $C^2$-boundary (see Remark \ref{rem:MR-d1}). Without additional assumptions on $D$, however, such solutions need not exist.

Before attacking the problem \eqref{eq:Poisson}
using Hilbert space methods, we pause to emphasise that in some special cases this problem is simple enough to admit a direct
elementary solution. For instance, for $d=1$ and $D = (a,b)$ we may integrate the equation twice and determine the two integration constants
by substituting the boundary conditions, which take the form $u(a)=u(b)=0$. After some computations one arrives at
$$ u(x) = \int_a^b k(x,y)f(y)\ud y,\quad x\in (a,b),$$
where
\begin{align*} k(x,y) = \frac1{b-a}\begin{cases}
       (b-y)(x-a), & \ x\le y, \\
       (b-x)(y-a), & \ y\le x,
      \end{cases}
\end{align*}
is the so-called {\em Green function}\index{Green function} for the Poisson problem on $(a,b)$ with Dirichlet boundary conditions.
The reader may check (see Problem \ref{prob:Green}) that the function $u$ thus defined belongs to $H_0^1(a,b)$ and
is indeed a weak solution of \eqref{eq:Poisson}.
It is unclear, however, how to extend this elementary method to higher dimensions.
In contrast, the Hilbert space
method adopted here works in arbitrary dimensions.
Our main tool is the following inequality,
which is of interest in its own right.

\begin{theorem}[Poincar\'e inequality]\label{thm:Poincare}\index{inequality!Poincar\'e} Let $D$ be a bounded open set
contained in $R:=(0,r)\times \R^{d-1}$ and let $1\le p< \infty$.
Then the
following estimate holds:
$$ \n f\n_{p} \le rp^{-1/p} \Big\n \frac{\partial f}{\partial x_1}\Big\n_p, \quad f\in W_0^{1,p}(D).$$
As a consequence,
$\nn f\nn_{W_0^{1,p}(D)} := \n \nabla f\n_{L^p(D;\K^d)}$ defines an equivalent norm on $W_0^{1,p}(D)$.
\end{theorem}

By translation, the second assertion holds for any bounded open set $D\subseteq \R^d$.

\begin{proof}
First assume that $f\in C_{\rm c}^\infty(D)$.
Extending $f$ by zero to $R$, we may view $f$ as a function in
$C_{\rm c}^\infty(R)$.

For $x_1\in [0,r]$ we have, by H\"older's inequality and the fact that $f(0,x_2,\dots,x_d) = 0$,
\begin{align*}|f(x_1,x_2,\dots,x_d)|^p
 & = \Big|\int_{0}^{x_1} \frac{\partial f}{\partial x_1}(t,x_2,\dots,x_d)\ud t\Big|^p
 \le x_1^{p/q}\int_{0}^{x_1} \Big|\frac{\partial f}{\partial x_1}(t,x_2,\dots,x_d)\Big|^p\ud t,
\end{align*}
where $\frac1p+\frac1q=1$ and $x_1^{p/q}:=1$ if $p=1$.
Integrating both sides over $R$, we obtain
\begin{align*}
 \n f\n_p^p &  \le \int_R \Bigl(x_1^{p/q} \int_0^{r} \Big|\frac{\partial f}{\partial x_1}(t,x_2,\dots,x_d)\Big|^p\ud t\Bigr)\ud x
\\ & = \int_0^r x_1^{p/q}\ud x_1 \int_{\R^{d-1}}\int_{0}^{r} \Big|\frac{\partial f}{\partial x_1}(t,x_2,\dots,x_d)\Big|^p\ud t\ud x_2\cdots\ud x_d
 = \frac{r^{p}}{p} \Big\n \frac{\partial f}{\partial x_1}\Big\n_p^p.
\end{align*}
Since $C_{\rm c}^\infty(D)$ is dense in $W_0^{1,p}(D)$, the estimate for a general $f\in W_0^{1,p}(D)$ follows by approximation.
The equivalence of the norms follows from
$$\n \nabla f\n_p \le \n f\n_p+ \n \nabla f\n_p \le (rp^{-1/p}+1)\n \nabla f\n_p,$$
where we used the trivial estimate $\n \frac{\partial f}{\partial x_1}\n_p \le\n \nabla f\n_p$.
\end{proof}

We now take $p=2$ and recall the notation $H^1(D) = W^{1,2}(D)$ and $H_0^1(D) = W_0^{1,2}(D)$.
Poincar\'e's inequality then states that, if $D$ is bounded,
\begin{align*}\nn f\nn_{H_0^1(D)} := \n \nabla f\n_2, \quad f\in H_0^1(D),
\end{align*}
defines an equivalent
norm on $H_0^1(D)$. With respect to this norm,
$H_0^1(D)$ is again a Hilbert space, this time with respect to the inner product
$$ (\!(f_1|f_2)\!)_{H_0^1(D)} := \iprod{\nabla f_1}{\nabla f_2}_{2}.$$

\begin{theorem}[Poisson problem, Dirichlet boundary conditions]\label{thm:Poisson-D}\index{existence and uniqueness!for the Poisson problem}
If $D$ is bounded, then for every $f\in L^2(D)$ the Poisson problem \eqref{eq:Poisson}
admits a unique weak solution $u\in H_0^1(D)$. Moreover, there exists a constant $C\ge 0$, independent of $f$,
such that $$ \n u\n_{H_0^1(D)} \le C \n f\n_2.$$
\end{theorem}

\begin{proof}
By the Cauchy--Schwarz inequality and Poincar\'e's inequality, the linear mapping
$L : g \mapsto \int_D g\ov f\ud x $ is bounded from $H_0^1(D)$
to $\K$:
$$ |L(g)| \le \n g\n_2 \n f\n_2 \le C \n \nabla g\n_2 \n f\n_2 = C \nn g\nn_{H_0^1(D)}\n f\n_2.$$
Therefore $L$ defines a bounded functional on $H_0^1(D)$.
Hence, by the Riesz representation theorem, there exists a unique $u\in H_0^1(D)$ such that
\begin{align}\label{eq:Poiss-ip} L(g) = (\!\iprod{g}{u}\!)_{H_0^1(D)}, \quad g\in H_0^1(D),
\end{align}
and it satisfies
$\nn u\nn_{H_0^1(D)} = \n L\n \le C\n f\n_2$.
Writing out the identity \eqref{eq:Poiss-ip}, we obtain
\begin{align}\label{eq:H2-u} \int_D \nabla g\cdot \ov{\nabla u}\ud x =  \int_D g\ov{f}\ud x, \quad g\in H_0^1(D).
\end{align}
In particular, \eqref{eq:H2-u} holds for all $g\in C_{\rm c}^\infty(D)$, since $C_{\rm c}^\infty(D)$ is contained in $H_0^1(D) $.
Replacing $g$ by $\ov g$ and taking conjugates on both sides, we see that $u$ is a weak solution of \eqref{eq:Poisson}.

If $v\in H_0^1(D)$ is another weak solution, then
$$ \int_D (\nabla u-\nabla v)\cdot{\nabla \phi}\ud x = 0, \quad \phi\in C_{\rm c}^\infty(D).$$
Since $C_{\rm c}^\infty(D)$ is dense in $H_0^1(D) $, it follows that
$$ \int_D (\nabla u-\nabla v)\cdot{\nabla g}\ud x = 0, \quad g\in H_0^1(D),$$
and applying this to $\ov g$ gives $(\!\iprod{u-v}{g}\!)_{H_0^1(D)} = 0$ for all $g\in H_0^1(D)$. This implies $u-v = 0$
in $H_0^1(D)$.
\end{proof}

We next prove that the weak solution actually belongs to $H_0^1(D)\cap H_{\rm loc}^2(D)$, where $H_{\rm loc}^2(D)$
is the space of all $f\in L_{\rm loc}^1(D)$ with the property that $f|_U \in H^2(U)$ for all open sets $U\Subset D$.
With $L_{\rm loc}^2(D)$ defined similarly, this will follow from the following lemma.

\begin{lemma}\label{lem:MR-d1}
If $u\in H^1(D)$ admits a weak Laplacian in $L_{\rm loc}^2(D)$, then $u\in H_{\rm loc}^2(D)$.
\end{lemma}

\begin{proof}
Let $U,U'$ be bounded open sets such that
$U\Subset U'\Subset D$, and choose
$\psi\in C_{\rm c}^\infty(U')$ such that $\psi\equiv1$ on $U$.
Let $w$ be the zero extension of $\psi u$ to $\R^d$.

We claim that $w$ has weak Laplacian
$$
G
:=
\psi\Delta u
+
2\sum_{j=1}^d(\partial_j\psi)(\partial_j u)
+
(\Delta\psi)u,
$$
where the right-hand side is extended by zero outside $D$.
Indeed, for $\phi\in C_{\rm c}^\infty(\R^d)$ we have
$\psi\phi\in C_{\rm c}^\infty(D)$, and therefore
\begin{align*}
\int_{\R^d}w\Delta\phi\ud x
&=
\int_Du\psi\Delta\phi\ud x
\\
&=
\int_Du\Delta(\psi\phi)\ud x
-
2\sum_{j=1}^d
\int_Du(\partial_j\psi)(\partial_j\phi)\ud x
-
\int_Du(\Delta\psi)\phi\ud x
\\
&=
\int_D(\Delta u)\psi\phi\ud x
+
2\sum_{j=1}^d
\int_D(\partial_j u)(\partial_j\psi)\phi\ud x
+
\int_Du(\Delta\psi)\phi\ud x
=
\int_{\R^d}G\phi\ud x.
\end{align*}
All terms defining $G$ belong to $L^2(\R^d)$. Hence Theorem
\ref{thm:H1-Sob} gives $w\in H^2(\R^d)$. Since $w=u$ on $U$, it follows
that $u|_U\in H^2(U)$.
\end{proof}

\begin{theorem}\label{thm:MR-d1} Let $D$ be bounded. The weak solution $u$ of the Poisson problem $-\Delta u =f$ with $f\in L^2(D)$, subject to Dirichlet boundary conditions, belongs to $H_0^1(D)\cap H_{\rm loc}^2(D)$.
\end{theorem}
\begin{proof}
The very definition of a weak solution implies that
$u$ admits a weak Laplacian belonging to $L^2(D)$, given by $\Delta u = -f$. The result now follows from Lemma \ref{lem:MR-d1}.
\end{proof}

\begin{remark}\label{rem:MR-d1}
If $D$ is bounded and has a $C^2$-boundary, then in fact the weak solution
belongs to $H^2(D)$. This is a standard elliptic boundary regularity result for
the Dirichlet Laplacian, the proof of which is beyond the scope of this book.
\end{remark}

The existence of a weak solution does not imply the existence of a classical
solution, even when the right-hand side is continuous and compactly supported.
We illustrate this for the Dirichlet problem on the unit disc
$\mathbb D\subseteq\R^2$.

\begin{example}[Nonexistence of classical solutions]\label{ex:noclass-Poisson}
We now show that there exists a function $f\in C_{\rm c}(\mathbb D)$ such that the problem
$$
\begin{cases} -\Delta u & = f \ \ \hbox{on $\mathbb D$}, \\ u|_{\partial \mathbb D} & = 0,
\end{cases}
$$
has no classical solution $u\in C^2(\mathbb D)\cap C(\ov{\mathbb D})$.

For $0<r<\frac12$ put $\ell(r):=\log|\log r|$ and define
$$
v(x,y):=(x^2-y^2)\ell((x^2+y^2)^{1/2}),
\quad 0<x^2+y^2<\frac14.
$$
Writing $r=(x^2+y^2)^{1/2}$ and differentiating, we obtain
$$
v_x=2x\ell(r)+\frac{x(x^2-y^2)}{r^2\log r},
\qquad
v_y=-2y\ell(r)+\frac{y(x^2-y^2)}{r^2\log r}.
$$
Since $r^2\ell(r)\to0$, $r\ell(r)\to0$, and $r/|\log r|\to0$ as
$r\downarrow0$, the function $v$ extends to an element of
$C^1(B(0;\frac12))$ by setting $v(0)=0$ and $\nabla v(0)=0$.

A further calculation gives
$$
\Delta v(x,y)
=
\frac{x^2-y^2}{r^2}
\Bigl(\frac4{\log r}-\frac1{(\log r)^2}\Bigr).
$$
The right-hand side tends to $0$ as $r\downarrow0$, and therefore extends
to a function $g\in C(B(0;\frac12))$ by putting $g(0)=0$. Moreover,
$\Delta v=g$ in the sense of weak derivatives. Indeed, applying Green's
formula on $B(0;\frac12)\setminus\ov{B(0;\eps)}$ and letting
$\eps\downarrow0$ gives
$$
\int_{B(0;\frac12)}v\Delta\phi\ud x
=
\int_{B(0;\frac12)}g\phi\ud x,
\quad \phi\in C_{\rm c}^\infty(B(0;\tfrac12)),
$$
because on $\partial B(0;\eps)$,
$$
|v|\le C\eps^2\log|\log\eps|,
\qquad
|\nabla v|
\le C\eps\Bigl(\log|\log\eps|+\frac1{|\log\eps|}\Bigr),
$$
and hence the boundary terms tend to zero.

Choose $\eta\in C_{\rm c}^\infty(B(0;\frac12))$ such that
$\eta=1$ on $B(0;\frac14)$, and put $u:=-\eta v$. Then
$u\in C_{\rm c}^1(\mathbb D)\subseteq H_0^1(\mathbb D)$ and, in the
sense of weak derivatives,
$-\Delta u
=
\eta g+2\nabla\eta\cdot\nabla v+(\Delta\eta)v=:f$.
Thus $f\in C_{\rm c}(\mathbb D)$ and $u$ is a weak solution.

For $x\ne0$ sufficiently small, differentiation gives
$$
v_{xx}(x,x)
=
2\log|\log{(\sqrt2|x|)}|
+2/{\log{(\sqrt2|x|)}},
$$
which tends to infinity as $x\to0$. Since $u=-v$ near the origin, we have
$u\notin C^2(\mathbb D)$.

Suppose now, for a contradiction, that
$\wt u\in C^2(\mathbb D)\cap C(\ov{\mathbb D})$ is a classical solution,
and put $w:=\wt u-u$. Then $w\in H^1_{\rm loc}(\mathbb D)$ and
$\Delta w=0$ weakly. Iterating Lemma~\ref{lem:MR-d1}, using that weak
derivatives commute with the Laplacian, gives
$w\in H^4_{\rm loc}(\mathbb D)$.

Choose a disc $B$ centred at the origin with
$\ov B\subseteq B(0;\frac14)$ and a function
$\chi\in C_{\rm c}^\infty(\mathbb D)$ equal to $1$ on an open set containing $\ov B$. Then $\chi w\in H^4(\mathbb D)$, and
Lemma~\ref{lem:der-zeroextension} gives $\chi w\in H^4(\R^2)$.
Since $4>2+\frac d2=3$, the higher-order Sobolev embedding noted after
Theorem~\ref{thm:Hs-embedding-C0} shows that $\chi w$ has a
$C^2$-representative.
On $B$, this representative agrees almost everywhere with the continuous
function $w$, and therefore everywhere. Thus $w\in C^2(B)$, so
$u=\wt u-w\in C^2(B)$,
contradicting the fact that $u=-v$ near the origin and
$v_{xx}(x,x)\to\infty$.
We conclude that no classical solution exists.
\end{example}

\subsection{Neumann Boundary Conditions}\label{subsec:Neumann}\index{Neumann!boundary conditions}

Having dealt with Dirichlet boundary conditions, we now consider the
{\em Poisson problem with Neumann boundary conditions:}\index{Poisson!problem, Neumann boundary conditions}
\begin{equation}\label{eq:Poisson2}
\begin{cases} -\Delta u & = f \ \ \hbox{on $D$}, \\ \frac{\partial u}{\partial \nu} & = 0 \ \ \hbox{on $\partial D$},
\end{cases}
\end{equation}
where $\frac{\partial u}{\partial\nu}:=\nabla u\cdot\nu$
is the derivative in the outward normal direction along $\partial D$ and $f\in L^2(D)$ is a given function. The treatment of this example in higher dimensions requires some familiarity with standard techniques from partial differential equations, as the notion of an
outward normal is meaningful only under some regularity assumptions on the boundary $\partial D$. For the present treatment
$C^1$-regularity of the boundary suffices.

To motivate the notion of weak solutions we need Green's theorem:\index{theorem!Green}
{\em If $D$ is bounded with $C^1$-boundary, then for all $u\in C^2(\ov D)$ and $v\in C^1(\ov D)$ we have
 $$ \int_D \nabla u\cdot \nabla v \ud x = -\int_D v\Delta u\ud x + \int_{\partial D} v \frac{\partial u}{\partial \nu}\ud S,
 $$
where $\ud S$ is the surface measure on $\partial D$.}
We temporarily disregard the boundary condition and ask what information
is conveyed by the integrated equation
\begin{align}\label{eq:Poisson-N} \int_D \nabla u\cdot \nabla\phi\ud x = \int_D f\phi\ud x
\end{align}
if it is to hold for all $\phi\in C^\infty(\ov D)$ (and not just for all $\phi\in C_{\rm c}^\infty(D)$, since that would ignore what happens at the boundary).
By Green's theorem,
\begin{align}\label{eq:Poisson-int2a}
 \int_D \phi\Delta u\ud x - \int_{\partial D} \phi \frac{\partial u}{\partial \nu}\ud S= - \int_D f\phi\ud x.
\end{align}
Since we wish to solve $-\Delta u= f$, we substitute this relation into \eqref{eq:Poisson-int2a} and find that
$$ \int_{\partial D} \phi \frac{\partial u}{\partial \nu}\ud S = 0.$$
This can only hold for all $\phi\in C^\infty(\ov D)$ if
$$ \frac{\partial u}{\partial \nu} = 0 \ \ \hbox{on $\partial D$},$$
that is, {\em if Neumann boundary conditions hold.}

By considering $\phi\equiv 1$ we see that the identity \eqref{eq:Poisson-N} can only hold for all $\phi\in C^\infty(\ov D)$ if $f$ satisfies the compatibility condition
$$ \int_D f\ud x = 0.$$
More generally, the integral of $f$ against every locally constant function should vanish.
Likewise, solutions of \eqref{eq:Poisson-N} cannot be unique: if
$u$ is a solution (in whatever weak sense), then $u+C$ is also a solution, for any locally constant function $C$.
In order to simplify matters, we henceforth assume that $D$ is connected, so that the only locally constant functions are the constant functions. Under this assumption we eliminate integration constants by imposing the constraint
$$ \int_D u\ud x = 0,$$
that is, the average of $u$ over $D$ should vanish. Accordingly, we define\index{$H_{\rm av}^1(D)$}
$$H_{\rm av}^1(D):= \Bigl\{u\in H^1(D) :\, \int_D u\ud x = 0\Bigr\}.$$
Since $D$ is bounded, the mapping
$u\mapsto\int_Du\ud x$
is bounded on $H^1(D)$. Thus $H_{\rm av}^1(D)$ is a closed subspace of
$H^1(D)$ and hence is a Hilbert space.

The above discussion leads to the following weak formulation of the problem \eqref{eq:Poisson2}.

\begin{definition}[Weak solutions] Let $D$ be bounded and connected with $C^1$-boundary
and let $f\in L^2(D)$ satisfy $\int_D f\ud x = 0$.
A function $u\in H_{\rm av}^1(D)$ is called a {\em weak solution}\index{weak!solution} of the Poisson problem \eqref{eq:Poisson2} if
$$ \int_D \nabla u \cdot\nabla\phi\ud x = \int_D f\phi\ud x , \quad \phi\in {C^\infty}(\ov D).$$
\end{definition}

Although the weak formulation is defined using test functions
$\phi\in C^\infty(\overline D)$, it is equivalent here to the same identity
for all $\phi\in H^1(D)$: by Theorem \ref{thm:trace-W1p}
$C^\infty(\overline D)$ is dense in $H^1(D)$, and
both sides of \eqref{eq:Poisson-N} are bounded linear functionals of $\phi$
on $H^1(D)$.

Our treatment of the Poisson problem with Dirichlet boundary conditions crucially depended on
the Poincar\'e inequality for $H_0^1(D)$.
The treatment of Neumann boundary conditions proceeds analogously,
with $H_{\rm av}^1(D)$ replacing $H_0^1(D)$. We will prove a version of the Poincar\'e inequality for this space in Theorem \ref{thm:Poincare-mean}. Its proof depends on the following compactness result.

\begin{theorem}[Rellich--Kondrachov compactness theorem]\label{thm:Rellich}\index{theorem!Rellich--Kondrachov}
 If $D$ is bounded and $1\le p<\infty$, then:
\begin{enumerate}[label={\rm(\arabic*)}, leftmargin=*]
 \item\label{it:Rellich1} the inclusion mapping from $W_0^{1,p}(D)$ into $L^p(D)$ is compact;
 \item\label{it:Rellich2} if $D$ has $C^1$-boundary, the inclusion mapping from $W^{1,p}(D)$ into $L^p(D)$ is compact.
\end{enumerate}
\end{theorem}

\begin{proof} We first prove \ref{it:Rellich1} and deduce \ref{it:Rellich2} from
it
by means of the extension theorem.

\smallskip
\ref{it:Rellich1}:\
For $f\in L^p(D)$, let $\wt f$ denote its zero extension to $\R^d$.
If $f\in C_{\rm c}^\infty(D)$, then $\wt f$ is smooth and the
fundamental theorem of calculus gives
$$
\wt f(x+h)-\wt f(x)
=
\int_0^1 h\cdot\nabla\wt f(x+th)\ud t.
$$
Writing $\tau_hg(x):=g(x+h)$, we obtain
$$
\|\tau_h\wt f-\wt f\|_{L^p(\R^d)}
\le
|h|\int_0^1
\|\nabla\wt f(\cdot+th)\|_{L^p(\R^d)}\ud t
=
|h|\,\|\nabla f\|_{L^p(D)}.
$$
By density of $C_{\rm c}^\infty(D)$ in $W_0^{1,p}(D)$, and
continuity of zero extension and translation in $L^p$, this estimate
extends to every $f\in W_0^{1,p}(D)$.

Taking the supremum over the unit ball $B$ of $W_0^{1,p}(D)$ gives
\begin{align*}
\lim_{|h|\to0}\sup_{f\in B}
\|\tau_h\wt f-\wt f\|_{L^p(\R^d)}=0.
\end{align*}
The zero extensions form a bounded subset of $L^p(\R^d)$ supported
in the fixed compact set $\ov D$.
Corollary \ref{cor:RF} therefore shows that $B$ is relatively
compact in $L^p(D)$.

\smallskip
\ref{it:Rellich2}: \
Choose a bounded open set $D'$ with $D\Subset D'$ and
$\psi\in C_{\rm c}^\infty(D')$ such that $\psi\equiv1$ on $D$.
If $E_D:W^{1,p}(D)\to W^{1,p}(\R^d)$ is the extension operator of
Theorem~\ref{thm:Sobolev-extension}, then
$$
f\mapsto\psi E_Df
$$
is bounded from $W^{1,p}(D)$ into $W_0^{1,p}(D')$. The inclusion of the
latter space into $L^p(D')$ is compact by \ref{it:Rellich1}, and restriction
from $D'$ to $D$ is bounded. Since the composition of these three mappings
is the inclusion of $W^{1,p}(D)$ into $L^p(D)$, this inclusion is compact.
\end{proof}

As an application of the Rellich--Kondrachov theorem, we have the following variant of Poincar\'e's inequality.
Define
$$W_{\rm av}^{1,p}(D):= \Bigl\{u\in W^{1,p}(D) :\, \int_D u\ud x = 0\Bigr\}.$$

\begin{theorem}[Poincar\'e--Wirtinger inequality]\label{thm:Poincare-mean}\index{inequality!Poincar\'e--Wirtinger} Let $D$ be bounded and connected with $C^1$-boundary. Let $1\le p<\infty$.
Then there exists a constant $C = C_{p,D}$ such that for all $f\in W_{\rm av}^{1,p}(D)$ the following estimate holds:
$$ \n f\n_p \le C \n \nabla f\n_p.$$
In particular,
$\nn f\nn_{W_{\rm av}^{1,p}(D)} := \n \nabla f\n_{p}$ defines an equivalent norm on $W_{\rm av}^{1,p}(D)$.
\end{theorem}
\begin{proof}
We argue by contradiction. If the theorem were false, we could find a sequence $(f_n)_{n\ge 1}$ in $W_{\rm av}^{1,p}(D)$
such that $\n f_n\n_p \ge n \n \nabla f_n\n_p$ for $n\ge 1$. By scaling, we may assume that
$\n f_n\n_p=1$, so that $\n \nabla f_n\n_p \le \frac1n$.

Since $(f_n)_{n\ge 1}$ is bounded in $W_{\rm av}^{1,p}(D)$,
we may use the Rellich--Kondrachov theorem to extract a subsequence $(f_{n_k})$ that converges in $L^p(D)$ to some $f\in L^p(D)$.
Also $\n \nabla f_{n_k}\n_p \le \frac1{n_k} \to 0$ as
$k\to \infty$,
and therefore
the closedness of $\nabla$ as an operator from $L^p(D)$ to $L^p(D;\K^d)$ implies that
$f\in \Dom(\nabla) = W^{1,p}(D)$ and $\nabla f = 0$.
In view of
\begin{align}\label{eq:Poinc-ave}\int_D f\ud x = \limk \int_D f_{n_k}\ud x = 0,
\end{align}
we even have $f\in W_{\rm av}^{1,p}(D)$. But $\nabla f = 0$
implies, via Proposition \ref{prop:der-const}, that $f$ is a constant almost everywhere. In view of \eqref{eq:Poinc-ave}
this is only possible if $f = 0$ almost everywhere. We thus arrive at the contradiction
$ 0 = \n f\n_p =\limk \n f_{n_k}\n_p = 1.$
\end{proof}

We are now in a position to solve the Poisson problem with Neumann boundary
conditions.

\begin{theorem}[Poisson problem, Neumann boundary conditions]\label{thm:Poisson-N}\index{existence and uniqueness!for the Poisson problem} Let $D$ be bounded and connected with $C^1$-boundary. For every
$f\in L^2(D)$ satisfying $$ \int_D f\ud x = 0,$$
 the Poisson problem \eqref{eq:Poisson2} admits a unique weak solution $u\in H_{\rm av}^1(D)$.
 Moreover, there exists a constant $C\ge 0$, independent of $f$,
such that $$ \n u\n_{H_{\rm av}^1(D)} \le C \n f\n_2.$$
\end{theorem}
\begin{proof}
The argument is similar to the proof of Theorem \ref{thm:Poisson-D},
with minor adjustments.
By the Poincar\'e--Wirtinger inequality,
$$
\nn g\nn_{H_{\rm av}^1(D)}:=\n\nabla g\n_2
$$
is an equivalent Hilbert space norm on $H_{\rm av}^1(D)$. The functional
$$
L(g):=\int_Dg\ov f\ud x,\quad g\in H_{\rm av}^1(D),
$$
is bounded with respect to this norm. Hence the Riesz representation theorem
gives a unique $u\in H_{\rm av}^1(D)$ such that
$$
\int_D\nabla g\cdot\ov{\nabla u}\ud x
=
\int_Dg\ov f\ud x,\quad g\in H_{\rm av}^1(D),
$$
and
$$
\n\nabla u\n_2\le C\n f\n_2.
$$

Both sides of the preceding integral identity vanish when $g$ is constant. Since every
$g\in H^1(D)$ is the sum of a constant and an element of
$H_{\rm av}^1(D)$, the identity holds for all $g\in H^1(D)$. Taking
$g=\ov\phi$ and conjugating shows that $u$ is a weak solution.

If $u$ and $v$ are weak solutions, their difference satisfies the homogeneous
identity for all test functions in $C^\infty(\ov D)$, and hence, by density,
for all test functions in $H^1(D)$. Taking $\phi=\ov u-\ov v$ gives
$\n\nabla(u-v)\n_2=0$. Since $u-v$ has mean zero, the
Poincar\'e--Wirtinger inequality implies that $u=v$.
The asserted estimate follows from the equivalence of norms.
\end{proof}

The analogue of Theorem \ref{thm:MR-d1} holds, with the same proof:

\begin{theorem}\label{thm:MR-d1-N} Let $D$ be bounded and connected with $C^1$-boundary
and let $f\in L^2(D)$ satisfy $\int_D f\ud x = 0$.
The weak solution of the Poisson problem $-\Delta u =f$, subject to Neumann boundary conditions, belongs to $H^1(D)\cap H_{\rm loc}^2(D)$.
\end{theorem}

If $D$ is bounded and has a $C^2$-boundary, it can be shown that the weak solution $u$ belongs to $H^2(D)$.

\section{Variational Methods}\label{sec:VariantionalMethods}

In this section we develop a variational framework for elliptic boundary value problems. The guiding principle is that such equations arise as the Euler--Lagrange equations for certain minimisation problems.

The terminology used in this section comes from the calculus of variations, particularly its applications to classical mechanics. In that setting one studies functionals whose argument is not a point in a vector space but a \emph{curve} (a trajectory). The basic example is the \emph{action functional} in Lagrangian mechanics,
$$
\mathcal A(q):=\int_{t_0}^{t_1} L(q(t),\dot q(t),t)\ud t,
$$
where $q:[t_0,t_1]\to\R^n$ is a sufficiently regular curve and $L$ is a given function (the \emph{Lagrangian}).
Hamilton's principle asserts that the physical trajectories are those curves $q$ for which the action is stationary under perturbations that keep the endpoints fixed.
More precisely, one considers variations of the form $q_\varepsilon:=q+\varepsilon\eta$, where $\eta$ is a smooth curve with $\eta(t_0)=\eta(t_1)=0$. Stationarity of $\mathcal A$ at $q$ means that for such $\eta$,
$$
\frac{\ud}{\ud\varepsilon}\Big|_{\varepsilon=0}\mathcal A(q_\varepsilon)=0.
$$
A standard integration by parts argument shows that this condition is equivalent to the \emph{Euler--Lagrange equation}\index{Euler--Lagrange equation}
$$
\frac{\ud}{\ud t} \frac{\partial L}{\partial \dot q}(q(t),\dot q(t),t)
-\frac{\partial L}{\partial q} (q(t),\dot q(t),t)=0,
\quad t\in(t_0,t_1).
$$
In the present section, the functionals are defined on a Hilbert space
and we are interested in their minimisers. We first consider energy
functionals associated with bounded symmetric coercive forms, and then
apply the resulting characterisation to the Poisson problem with
Dirichlet and Neumann boundary conditions.
The underlying idea of taking variations reappears in the proof of
Theorem~\ref{thm:quadratic-energy-Euler-Lagrange}: given a minimiser $v_0$, we perturb it along lines $v_0+\eps v$ and obtain the variational identity from the vanishing of the first derivative at $\eps=0$.

\subsection{The Direct Method}

We start by presenting the direct method for proving existence and uniqueness of
minimisers for energy functionals associated with certain sesquilinear forms,
and for characterising these minimisers by an abstract Euler--Lagrange
equation.

\begin{definition}[Forms]\label{def:form-bounded}
A {\em form}\index{form} on a Hilbert space $V$ is a sesquilinear mapping
$\aa:V\times V\to\K$. A form $\aa:V\times V\to\K$ is called:
\begin{itemize}
\item[$\bullet$] \emph{accretive},\index{accretive!form} if
$\Re\aa(v,v)\ge 0$ for all $v\in V$;
\item[$\bullet$] \emph{coercive},\index{coercive form} if there exists $\alpha>0$ such that
$\Re\aa(v,v)\ge \alpha\n{v}\n_V^2$ for all $v\in V$;
\item[$\bullet$] \emph{symmetric},\index{symmetric!form} if
$\aa(u,v)=\ov{\aa(v,u)}$ for all $u,v\in V$;
\item[$\bullet$] \emph{bounded},\index{bounded!form} if there exists $C\ge 0$ such that
 $|\aa(u,v)|\le C\n{u}\n_V\n{v}\n_V$ for all $u,v\in V$.
\end{itemize}
\end{definition}

In the complex case, a form $\aa:V\times V\to\C$ is accretive
and symmetric if and only if $\aa$ is positive in the sense that
$\aa(v,v)\ge 0$ for all $v\in V$; see Proposition \ref{prop:form-continuity} for the easy proof.

\begin{remark}
The reason for denoting the Hilbert spaces by $V$ in this section can be explained as follows. Taking the Poisson problem
$\Delta u = -f$ with Dirichlet boundary conditions as our guiding example, let us revisit the identity \eqref{eq:Poisson-int},
$$ \int_D \nabla u\cdot \nabla \phi \ud x = \int_D f \phi \ud x,$$
where $\phi\in C_{\rm c}^\infty(D)$ is a test function.
Replacing $\phi$ by its complex conjugate and taking complex conjugates on both sides of the equation, we may write this abstractly as $$ \aa(\phi,u) = L\phi,$$
where
$$ \aa(\phi,u) = \int_D \nabla \phi\cdot \ov{\nabla u} \ud x, \qquad L\phi = \int_D \phi \ov f\ud x.$$
To give this a rigorous meaning, one interprets $\aa$ as a form on the Hilbert space $V = H_0^1(D)$. In this setting, the symbol $H$ is reserved for the ambient space $L^2(D)$ where the forcing function $f$ lives.
\end{remark}

Let $\aa:V\times V\to\K$ be a symmetric form, let $L:V\to\K$ be a bounded linear functional, and consider the {\em energy functional}\index{energy functional}
$$
E(v):=\frac12\aa(v,v)-\Re L(v), \quad v\in V.
$$
For accretive symmetric forms we have the following variational characterisation of minimisers.

\begin{theorem}[Minimisers and the Euler--Lagrange equation]
\label{thm:quadratic-energy-Euler-Lagrange}
Let $V$ be a Hilbert space, let $\aa:V\times V\to\K$ be an accretive symmetric form, let $L:V\to\K$ be a bounded linear functional, and let
$$
E(v):=\frac12\aa(v,v)-\Re L(v),\quad v\in V.
$$
For $v_0\in V$, the following are equivalent:
\begin{enumerate}[label={\rm(\roman*)}, leftmargin=*]
\item\label{it:quadratic-energy-Euler-Lagrange1} $v_0$ is a minimiser of $E$;
\item\label{it:quadratic-energy-Euler-Lagrange2} the {\em Euler--Lagrange equation}
$$\aa(v,v_0)=L(v), \quad v\in V,$$
is satisfied.
\end{enumerate}
If these equivalent conditions hold, then
\begin{align}\label{eq:E0-difference}
E(v)-E(v_0)=\frac12\aa(v-v_0,v-v_0), \quad v\in V.
\end{align}
If, in addition, $\aa$ is coercive, then the minimiser is unique.
\end{theorem}

\begin{proof}
Assume first that \ref{it:quadratic-energy-Euler-Lagrange1} holds, that is, $v_0$ is a minimiser. Fix $v\in V$ and define $g:\R\to\R$ by $g(\eps):=E(v_0+\eps v)$.
Using the symmetry of $\aa$, we obtain
$$
g(\eps)=E(v_0)+\eps\Re(\aa(v,v_0)-L(v))+\frac{\eps^2}{2}\aa(v,v).
$$
Since $\eps=0$ minimises $g$, we have $\frac{{\rm d}}{{\rm d}\eps}g(0)=0$, and therefore
$$
\Re(\aa(v,v_0)-L(v))=0,\quad v\in V.
$$
If $\K=\R$, this is exactly \ref{it:quadratic-energy-Euler-Lagrange2}. If $\K=\C$, applying the same identity to $iv$ gives
$$
0=\Re(\aa(iv,v_0)-L(iv))=-\Im(\aa(v,v_0)-L(v)),
$$
so both the real and imaginary parts of $\aa(v,v_0)-L(v)$ vanish, and \ref{it:quadratic-energy-Euler-Lagrange2} follows.

Assume next that \ref{it:quadratic-energy-Euler-Lagrange2} holds. For any $v\in V$, writing $w:=v-v_0$ and using symmetry, we obtain
\begin{align*}
E(v)-E(v_0)
&=\frac12\aa(v_0+w,v_0+w)-\Re L(v_0+w)-\frac12\aa(v_0,v_0)+\Re L(v_0)
\\
&=\Re(\aa(w,v_0)-L(w))+\frac12\aa(w,w) =\frac12\aa(w,w).
\end{align*}
This proves the identity \eqref{eq:E0-difference}, and accretivity shows that the right-hand side is nonnegative.
Hence $v_0$ is a minimiser. This shows that \ref{it:quadratic-energy-Euler-Lagrange1} holds.

If $v_1$ is another minimiser, then applying \eqref{eq:E0-difference} with $v=v_1$ gives
$$
0=E(v_1)-E(v_0)=\frac12\aa(v_1-v_0,v_1-v_0).
$$
If $\aa$ is coercive, this implies $v_1=v_0$.
\end{proof}

For bounded, coercive, symmetric forms, existence of minimisers can be shown. For the proof we need the following lemma.

\begin{lemma}\label{lem:quadratic-energy-properties}
Let $\aa: V\times V \to \K$ be a bounded symmetric form, let $L:V\to \K$ be bounded and linear, and let
$$ E(v):= \frac12\aa(v,v) - \Re L(v).$$
\begin{enumerate}[label={\rm(\roman*)}, leftmargin=*]
\item\label{it:quadratic-energy-properties1} if $\aa$ is accretive, then $E$ is sequentially weakly lower semicontinuous on $V$,\index{lower semicontinuous!energy functional} in the sense that whenever $v_n\to v$ weakly in $V$, then
$$
E(v)\le \liminf_{n\to\infty}E(v_n);
$$
\item\label{it:quadratic-energy-properties2} if $\aa$ is coercive, then $E$ is bounded below and coercive in the sense that
$$\inf_{v\in V}E(v)>-\infty, \qquad E(v) \to\infty \ \ \hbox{as}\ \ \n{v}\n_V\to\infty.$$
\end{enumerate}
\end{lemma}

\begin{proof}
\ref{it:quadratic-energy-properties1}: \
Let $v_n\to v$ weakly in $V$. Since $L$ is linear and bounded, it is weakly continuous and
therefore $L(v_n)\to L(v)$. Moreover, since $\aa(v,v)\in \R$ by symmetry,
$$
\aa(v_n,v_n)-\aa(v,v)=\aa(v_n-v,v_n-v)+2\Re\aa(v_n-v,v).
$$
By accretivity, the first term on the right-hand side is nonnegative, while the second term tends to $0$
because $v_n-v\to 0$ weakly and $u\mapsto \aa(u,v)$ is a bounded linear functional on $V$, by the boundedness of $\aa$. Hence
$$
\aa(v,v)\le \liminf_{n\to\infty}\aa(v_n,v_n),
$$
and therefore
$$
E(v)=\frac12\aa(v,v)-\Re L(v)
\le \liminf_{n\to\infty}\Bigl(\frac12\aa(v_n,v_n)-\Re L(v_n)\Bigr)
=\liminf_{n\to\infty}E(v_n).
$$
So $E$ is sequentially weakly lower semicontinuous.

\smallskip
\ref{it:quadratic-energy-properties2}:\ If $\alpha>0$ is the coercivity constant of $\aa$, then for all $v\in V$ we have
$$
E(v)\ge \frac{\alpha}{2}\n{v}\n_V^2-\n{L}\n_{V^*}\n{v}\n_V,
$$
and both assertions follow.
\end{proof}

\begin{theorem}[Existence and uniqueness]
\label{thm:quadratic-energy-minimiser}
Let $V$ be a Hilbert space, let $\aa:V\times V\to\K$ be a bounded, coercive, symmetric form,
and let $L:V\to\K$ be a bounded linear functional. Define
$$
E(v):=\frac12\aa(v,v)-\Re L(v),\quad v\in V.
$$
Then $E$ admits a unique minimiser $v_0\in V$.
\end{theorem}

\begin{proof}
By the first assertion in part \ref{it:quadratic-energy-properties2} of Lemma~\ref{lem:quadratic-energy-properties}, the functional $E$ is bounded from below. Let $(v_n)_{n\ge 1}$ be a minimising
sequence. By the second assertion in part \ref{it:quadratic-energy-properties2} of Lemma \ref{lem:quadratic-energy-properties}, $(v_n)_{n\ge 1}$ is bounded in $V$. By Proposition \ref{prop:Hilbert-weakconvergence}, upon passing to
a subsequence, we may assume that $v_n\to v_0$ weakly for some $v_0\in V$. By part \ref{it:quadratic-energy-properties1} of Lemma \ref{lem:quadratic-energy-properties},
$$
E(v_0)\le \liminf_{n\to\infty}E(v_n)=\inf_{v\in V}E(v),
$$
so $v_0$ is a minimiser.

Uniqueness follows from Theorem \ref{thm:quadratic-energy-Euler-Lagrange}.
\end{proof}

\subsection{The Dirichlet Poisson Problem Revisited}\label{subsec:Poisson-Dir-revisited}

The \emph{Dirichlet principle} asserts that a
harmonic function (or, more generally, a solution of the Poisson equation) may be characterised as the function which
minimises a natural energy among functions with the same boundary data. In its most classical form, one prescribes a boundary value and seeks a function $u$ with the correct boundary value for which the \emph{Dirichlet energy}
$$
\frac12\int_D |\nabla u|^2\ud x
$$
is minimal. In the treatment of the Poisson equation $-\Delta u = f$ with Dirichlet boundary conditions, the boundary condition is encoded by choosing the `energy space' $V$ to be $H_0^1(D)$. The forcing term $f$ leads to the energy functional
$$
\frac12\int_D |\nabla u|^2\ud x-\Re\int_D u\,\overline f\ud x
$$
on this space. If we write this functional abstractly as
$$\frac12\aa(u,u)- \Re L(u)$$ for $u\in V := H_0^1(D)$,
the Euler--Lagrange equation $\aa(u,u_0) = L(u)$ for the unique minimiser $u_0$ of this functional is exactly the equation for the weak solution of the Poisson equation with Dirichlet boundary conditions. This will be apparent from the proof of the following theorem. Throughout, $D$ is a nonempty open subset of $\R^d$.

\begin{theorem}[Variational characterisation]\label{thm:Dirprinc}
Let $D$ be bounded and let $f\in L^2(D)$. For a function $u_0\in H_0^1(D)$ the following assertions are equivalent:
\begin{enumerate}[label={\rm(\arabic*)}, leftmargin=*]
\item\label{it:VM-Dirprinc1} $u_0$ is the weak solution of the Poisson problem $-\Delta u = f$ on $D$ subject to Dirichlet boundary conditions;
\item\label{it:VM-Dirprinc2} $u_0$ minimises the energy functional $E:H_0^1(D)\to \R$ defined by
$$
E(u):= \frac12\int_D |\nabla u|^2\ud x -\Re \int_D u \, \ov f \ud x.
$$
\end{enumerate}
\end{theorem}

\begin{proof}
Set $V:=H_0^1(D)$, equipped with the norm inherited from $H^1(D)$,
and define
\begin{align}\label{eq:form-Dir-revisited}
\aa(u,v):=\int_D\nabla u\cdot\ov{\nabla v}\ud x,
\qquad
L(v):=\int_Dv\ov f\ud x.
\end{align}
Then $\aa$ is bounded, accretive, and symmetric, $L$ is bounded,
and $E(u)=\frac12\aa(u,u)-\Re L(u)$.
Since $D$ is bounded, Theorem \ref{thm:Poincare} also
shows that $\aa$ is coercive on $V$.
By Theorem \ref{thm:quadratic-energy-Euler-Lagrange}, $u_0$
minimises $E$ if and only if
$$
\aa(v,u_0)=L(v),\quad v\in V.
$$
By density of $C_{\rm c}^\infty(D)$ in $V$, this identity is
equivalent to its restriction to $v=\ov\phi$ with
$\phi\in C_{\rm c}^\infty(D)$. Taking complex conjugates gives
$$
\int_D\nabla u_0\cdot\nabla\phi\ud x
=
\int_Df\phi\ud x,
\quad \phi\in C_{\rm c}^\infty(D).
$$
This is precisely the weak formulation of the Dirichlet problem.
Existence and uniqueness of the minimiser follow from
Theorem \ref{thm:quadratic-energy-minimiser}.
\end{proof}

\subsection{The Neumann Poisson Problem Revisited}\label{subsec:Poisson-Neum-revisited}

A variational characterisation of weak solutions in the spirit of Theorem~\ref{thm:Dirprinc} can also be given in the case of Neumann boundary conditions. Before stating the result, it is worth noting that, under the assumption $\int_D f\ud x=0$, the energy functional
$$
E(u):=\frac12\int_D |\nabla u|^2\ud x-\Re\int_D u\overline f\ud x
$$
is invariant under addition of constants: $E(u+c)=E(u)$.
Consequently, minimisers of $E$ over $H^1(D)$ are not unique: if $u$ is a minimiser, then so is $u+c$. For this reason, we can only expect a \emph{unique} minimiser in
$H_{\rm av}^1(D)=\{u\in H^1(D):\,\int_D u\ud x=0\}$.

\begin{theorem}[Variational characterisation]\label{thm:Dirprinc-N}
Let $D$ be bounded and connected with $C^1$-boundary, and let $f\in L^2(D)$
satisfy $\int_D f\ud x=0$. For a function $u_0\in H_{\rm av}^1(D)$ the following assertions are equivalent:
\begin{enumerate}[label={\rm(\arabic*)}, leftmargin=*]
\item\label{it:VM-DirprincN-1} $u_0$ is the weak solution of the Poisson problem $-\Delta u = f$ on $D$ subject to Neumann boundary conditions;
\item\label{it:VM-DirprincN-2} $u_0$ minimises the energy functional $E:H^1(D)\to \R$ defined by
$$
E(u):= \frac12\int_D |\nabla u|^2\ud x -\Re \int_D u \, \ov f \ud x.
$$
\end{enumerate}
\end{theorem}

\begin{proof}
Set $V:=H_{\rm av}^1(D)$, equipped with the equivalent norm
$\n v\n_V:=\n\nabla v\n_2$, and define $\aa:V\times V \to \C$ and $L:V\to \C$ by
\begin{align}\label{eq:form-Neum-revisited}
\aa(u,v):=\int_D\nabla u\cdot\ov{\nabla v}\ud x,
\qquad
L(v):=\int_Dv\ov f\ud x.
\end{align}
By the Poincar\'e--Wirtinger inequality, $\aa$ is a bounded, coercive,
symmetric form on $V$, and $L$ is bounded on $V$.

For $u\in H^1(D)$ put
$$
u_{\rm av}:=
u-\frac1{|D|}\int_Du\ud x.
$$
Since $\int_Df\ud x=0$, we have
$E(u)=E(u_{\rm av})$.
Thus a function $u_0\in V$ minimises $E$ on $H^1(D)$ if and only if it
minimises the restriction of $E$ to $V$.
By Theorem~\ref{thm:quadratic-energy-Euler-Lagrange}, the latter is equivalent
to $\aa(v,u_0)=L(v),\quad v\in V$.
Since every element of $H^1(D)$
is the sum of a constant and an element of $V$, the identity is
equivalent to
$$
\aa(v,u_0)=L(v),\quad v\in H^1(D).
$$
By density of $C^\infty(\ov D)$ in $H^1(D)$, and by taking complex
conjugates as in the Dirichlet case, this is precisely the weak
formulation of the Neumann problem.
Existence and uniqueness of the mean-zero minimiser now follow from
Theorem~\ref{thm:quadratic-energy-minimiser}.
\end{proof}

\section{The Lax--Milgram Theorem}\label{sec:LM}

The considerations of the previous two subsections depended crucially on the symmetry of the forms $\aa$ defined by \eqref{eq:form-Dir-revisited} and \eqref{eq:form-Neum-revisited}, while the methods of Section \ref{sec:Poisson-problem} depended on the Riesz representation theorem.
The link is that an inner product is a symmetric sesquilinear form. The methods of Section \ref{sec:Poisson-problem}
can be extended to broader classes of boundary value problems by using a more flexible version of the Riesz representation theorem, the so-called Lax--Milgram theorem.

\subsection{The Theorem}\label{subsec:LM}

In what follows, $V$ is a Hilbert space.
The reason for using the letter $V$ is that in applications, typical choices are $V = H_0^1(D)$ and $V = H^1(D)$,
where $D$ is an open subset of $\R^d\!$. In such settings the letter $H$ will be reserved for the space $L^2(D)$. To avoid possible confusion, the inner product and norm of $V$ will be denoted by $\iprod{\cdot}{\cdot}_V$ and $\n \cdot\n_V$, respectively.

In the language of forms, Proposition \ref{prop:bilin-T} asserts that if $\aa:V\times V\to \K$ is a bounded form, then there exists a unique bounded operator $A$ on $V$
such that $$\aa(v,v') = \iprod{Av}{v'}_V, \quad v,v'\in V.$$
Moreover, $\n A\n_{\calL(V)} \le C$, where $C$ is the boundedness constant of $\aa$. For bounded coercive forms we have the following version of Proposition \ref{prop:bilin-T}.

\begin{theorem}[Lax--Milgram]\label{thm:LaxMilgram}\index{theorem!Lax--Milgram}
If $\aa$ is bounded and coercive on $V$, then:
\begin{enumerate}[label={\rm(\arabic*)}, leftmargin=*]
 \item\label{it:LaxMilgram1} the bounded operator $A$ associated with $\aa$ is boundedly invertible and $\n A^{-1}\n_{\calL(V)} \le \al^{-1}$\!, where $\alpha$ is the coercivity constant of $\aa$;
 \item\label{it:LaxMilgram2} for every bounded functional $L: V \to \K$ there exists a unique $v' \in V$ such that
$$ L(v) = \aa(v,v'), \quad v \in V.$$
Moreover, $\n v'\n_V \le \n A^{-1}\n_{\calL(V)} \n L\n$.
\end{enumerate}
\end{theorem}
\begin{proof} We proceed in two steps.

\smallskip
{\em Step 1} --
Let $A$ be the bounded operator provided by Proposition \ref{prop:bilin-T}.
The estimate
\begin{align*}
\alpha\n v\n_V^2 \le \Re\aa(v,v) \le |\aa(v,v)| = |\iprod{Av}{v}_V |\le \n Av\n_V \n v\n_V
\end{align*}
implies that $\alpha\n v\n_V \le \n Av\n_V$ for all $v\in V$. From this we infer that
$A$ is one-to-one and has closed range $\ran(A)$ in $V$ (see Proposition \ref{prop:closed-range}).
The operator $A$ is also surjective, for other\-wise there exists a nonzero element
 $v^\perp \in (\Ran(A))^\perp$ and we arrive at
$$0<\alpha\n v^\perp\n_V^2 \le\Re \aa(v^\perp\!, v^\perp) =\Re \iprod{Av^\perp}{v^\perp}_V = 0,$$
which is impossible.
The estimate $\alpha\n v\n_V \le \n Av\n_V$ now implies that the linear inverse $A^{-1}$ is bounded, with norm
$\n A^{-1}\n_{\calL(V)}\le \al^{-1}$\!.

\smallskip
{\em Step 2} --
Let $L:V\to\K$ be a bounded functional.
By the Riesz representation theorem, there exists a unique
$v'\in V$ such that $(L\circ A^{-1})(w)=\iprod{w}{v'}_V$ for all $w\in V$,
with norm
$\n v'\n_V = \n L\circ A^{-1}\n \le \n A^{-1}\n_{\calL(V)}\n L\n$.
Taking $w=Av$ gives
$$
L(v)=\iprod{Av}{v'}_V=\aa(v,v'),\quad v\in V.
$$
Since $A$ is surjective, uniqueness follows from the uniqueness part
in the Riesz representation theorem.
\end{proof}

Part \ref{it:LaxMilgram2} of the theorem provides a generalisation of the Riesz representation theorem with the inner product replaced by a bounded coercive form $\aa$.
If $\aa$ is {\em symmetric},\index{symmetric!sesquilinear mapping} that is, $\aa(v,v') = \overline{\aa(v'\!,v)}$ for all $v,v'\in V$
(some authors refer to this as $\aa$ being {\em Hermitian}), then
$\aa(v,v')$ defines an inner product on $V$ generating an equivalent norm. In this situation the Lax--Milgram theorem
is an immediate consequence of the Riesz representation theorem. The principal interest in the theorem lies in the nonsymmetric case.

\subsection{The Elliptic Problem \texorpdfstring{$\la u-\Delta u = f$}{lambda u - Delta u = f}}\index{elliptic!problem}

With the Lax--Milgram theorem at our disposal, the results of Subsections \ref{subsec:Poisson-Dir-revisited} and \ref{subsec:Poisson-Neum-revisited} admit straightforward modifications for the
{\em elliptic problem}\index{elliptic!problem}
$$\la u-\Delta u = f$$
for $\Re\la>0$ and $f\in L^2(D)$, subject to Dirichlet or Neumann boundary conditions.
As always, we assume that $D$ is a nonempty bounded open subset of $\R^d\!$, and in the case of Neumann boundary conditions we furthermore assume that $D$ has $C^1$-boundary.

To treat the case of Dirichlet boundary conditions we define a {\em weak solution}\index{weak!solution} to be a function $u\in H_0^1(D)$ such that
$$ \la \int_D u\phi \ud x + \int_D \nabla u \cdot\nabla \phi\ud x = \int_D f\phi\ud x, \quad \phi\in C_{\rm c}^\infty(D).$$
Repeating the steps of the proofs of Theorems \ref{thm:Poisson-D} and
\ref{thm:Dirprinc}, but using the Lax--Milgram theorem in place of the
Riesz representation theorem, one obtains:

\begin{theorem}[Elliptic problem, Dirichlet boundary conditions]\label{thm:elliptic-D}
If $D$ is bounded, then for all $\Re\la > 0$ and $f\in L^2(D)$ the elliptic problem $\la u - \Delta u = f$ subject to Dirichlet boundary conditions admits a unique weak solution. For $\la>0$ this weak solution is the unique minimiser of the energy functional\index{energy functional!Dirichlet conditions} $E:H_0^1(D)\to \R$ defined by $$E(u):= \frac12\int_D |\nabla u|^2 +\la |u|^2 \ud x -\Re \int_D u \ov f \ud x.$$
\end{theorem}

In the case of Neumann boundary conditions, one now works directly on $H^1(D)$: the form
$$
\aa_\la(u,v):=\int_D\nabla u\cdot\ov{\nabla v}\ud x
+\la\int_Du\ov v\ud x, \quad u,v\in H^1(D),
$$
satisfies
$$
\Re\aa_\la(u,u)
\ge \min\{1,\Re\la\}\n u\n_{H^1(D)}^2
$$
with respect to the usual norm $\n u\n_{H^1(D)}^2 = \n u\n_2^2 + \n \nabla u\n_2^2$;
the Poincar\'e--Wirtinger inequality, the mean-zero condition, and the connectedness assumption are not needed any longer. Repeating the heuristic arguments for the Poisson problem, we now define a {\em weak
solution}\index{weak!solution} to be an element $u\in H^1(D)$ such that
$$
  \lambda\int_D u\phi\ud x+\int_D\nabla u\cdot\nabla\phi\ud x
  =
  \int_D f\phi\ud x,
  \quad \phi\in C^\infty(\overline D).
$$
Repeating the steps of the proofs of Theorems \ref{thm:Poisson-N} (again
replacing the Riesz representation theorem with the Lax--Milgram theorem) and \ref{thm:Dirprinc-N}, one obtains:

\begin{theorem}[Elliptic problem, Neumann boundary conditions]\label{thm:elliptic-N}
 If $D$ is bounded with $C^1$-boundary, then for all $\Re\la > 0$ and $f\in L^2(D)$ the elliptic problem $\la u - \Delta u = f$ subject to Neumann boundary conditions admits a unique weak solution. For $\la>0$ this weak solution is the unique minimiser of the energy functional\index{energy functional!Neumann conditions} $E:H^1(D)\to \R$ defined by $$E(u):= \frac12\int_D |\nabla u|^2 +\la |u|^2 \ud x -\Re \int_D u \ov f \ud x.$$
\end{theorem}

We limit the present treatment of the elliptic problem to the above two theorems. In the next two chapters
we will develop powerful techniques that allow us, in a more general setting, to give precise $L^2$-estimates for the solutions $u$ in terms of the data $f$ and to extend these estimates to $L^p$ for $1\le p<\infty$.

\subsection{The Sturm--Liouville Problem}\label{subsec:SL}

Generalising the results for the elliptic problem of the preceding subsection, we now consider the Sturm--Liouville problem\index{Sturm--Liouville problem}
 \begin{equation}\label{eq:Sturm-Liouville-11}
 \begin{cases} -\operatorname{div}(a\nabla u) +qu & = f \ \ \hbox{on $D$}, \\ u|_{\partial D} & = 0.
 \end{cases}
 \end{equation}
Here, as always, $D$ is a nonempty bounded open subset of $\R^d$. We make the following assumptions:
 \begin{itemize}
 \item the function $f$ belongs to $L^2(D)$;
 \item the matrix-valued function $a:D\to M_d(\K)$ has bounded measurable coefficients and is \emph{uniformly elliptic}, in the sense that there exists $\alpha>0$ such that for almost all $x\in D$ and all $\xi\in\K^d$,
$$
\Re(a(x)\xi\cdot\overline{\xi})\ge \alpha|\xi|^2;
$$
 \item the function $q:D\to\K$ is bounded and measurable and satisfies $\Re q(x)\ge 0$ for almost all $x\in D$.
\end{itemize}
 A function $u\in H_0^1(D)$ is called a {\em weak solution} of \eqref{eq:Sturm-Liouville-11}
 if for all $\phi\in C_{\rm c}^\infty(D)$ we have
 \begin{align*}
 \int_D a\nabla u\cdot\nabla\phi\ud x + \int_D qu\phi\ud x
 = \int_D f\phi\ud x.
 \end{align*}

 \begin{theorem}[Sturm--Liouville problem]\label{thm:SturmLiouville}
Let $D$ be bounded. Under the above assumptions on $a$, $q$, and $f$, \eqref{eq:Sturm-Liouville-11} admits a unique weak solution $u$ in $H_0^1(D)$. Moreover, there exists a constant $C\ge 0$ independent of $f$ such that
$$\n u\n_{H_0^1(D)} \le C\n f\n_2.$$
\end{theorem}
\begin{proof}
The proof is a straightforward adaptation of the proof of existence and uniqueness for the Poisson problem with Dirichlet boundary conditions. This time we apply
the Lax--Milgram theorem to the form $\aa: H_0^1(D)\times H_0^1(D)\to \K$,
 \begin{align*}\aa(u,v):= \int_D a^\star \nabla u\cdot\ov{\nabla v}\ud x + \int_D \ov q u\ov{v}\ud x,
\end{align*}
where $a_{ij}^\star = \overline{a_{ji}}$. This form is bounded and coercive. Indeed,
boundedness follows from
\begin{align*}
|\aa(u,v)|
&\le
\|a\|_\infty\|\nabla u\|_2\|\nabla v\|_2
+
\|q\|_\infty\|u\|_2\|v\|_2
\le
\bigl(\|a\|_\infty+C_D^2\|q\|_\infty\bigr)
\|u\|_{H_0^1(D)}\|v\|_{H_0^1(D)},
\end{align*}
where $C_D$ is the Poincar\'e constant.
To prove coercivity we take $\xi=\nabla v(x)$ in the uniform ellipticity assumption on $a$. Integrating and using the assumption $\Re q\ge 0$, we obtain
\begin{align*}
\Re \aa(v,v) & = \Re\int_D a^\star \nabla v \cdot\ov{\nabla v}\ud x + \Re \int_D \ov q |v|^2\ud x
 \\* & \ge \int_D \Re (a^\star \nabla v \cdot\ov{\nabla v}) \ud x
 = \int_D \Re (a \nabla v \cdot\ov{\nabla v}) \ud x\ge \al \|\nabla v\|_2^2 =
 \al \nn v\nn_{H_0^1(D)}^2,
\end{align*}
where $\nn v\nn_{H_0^1(D)} = \n \nabla v\n_2$ is the equivalent norm on $H_0^1(D)$ introduced earlier.

With these ingredients, the proof may now be completed along the lines of the proof for the Poisson problem with Dirichlet boundary conditions.
\end{proof}

If $a=a^\star$ almost everywhere and $q: D\to [0,\infty) $ is nonnegative,
then the weak solution has again a variational interpretation: it
is the unique minimiser over $H_0^1(D)$ of the energy functional
$$
E(v)
=
\frac12\int_D a\nabla v\cdot\overline{\nabla v}\ud x
+
\frac12\int_D q|v|^2\ud x
-
\Re\int_D f\overline v\ud x.
$$
Indeed, the Euler--Lagrange equation for this functional is precisely the
weak formulation of the Sturm--Liouville problem.

The case of Neumann boundary conditions can be handled similarly and is left to the reader (see Problem \ref{prob:SL}).

\begin{problems}

\item
Determine whether the function $f\in L^1((-1,1)\times(-1,1))$ given by
$$ f(x,y):= |xy|, \quad (x,y)\in (-1,1)\times(-1,1)$$
has weak derivatives of order one. If `no', provide a proof; if `yes', compute the
weak derivatives $\partial_1 f$ and $\partial_2 f$.

\item\label{prob:mainthm-Sob}
Show that for all $f\in L^p(0,1)$ with $1\le p\le \infty$ the function
$$ I_f(x) :=\int_0^x f (y) \ud y,\quad x \in (0, 1), $$
belongs to $W^{1,p}(0,1)$ and its weak derivative is given by $\partial I_f= f$. Moreover,
the mapping $f \mapsto I_f$ from $L^p(0, 1)$ to $W^{ 1,p}(0, 1)$ is bounded.

\item\label{prob:cont-version-W1p}
Let $f\in W^{1,p}(0,1)$ with $1\le p\le \infty$.
\begin{enumerate}[\rm(a), leftmargin=*]
 \item Show that $f$ is equal almost everywhere to a (unique) continuous function $\wt f \in C[0, 1]$.

 \noindent {\em Hint:}\ the function
 $f- \int_0^\cdot f'(y) \ud y$ has weak derivative $0$.
 \item Show that the resulting mapping $f \mapsto \wt f $ from $W^{1,p}(0,1)$ to $C[0,1]$ is bounded.
 \item Show that a function $f \in W^{1,p}(0, 1)$ with $1 \le p <\infty$ belongs to
 $W_0^{1,p}(0,1)$
 if and only if its continuous version $\wt f$ satisfies
 $\wt f (0) = \wt f (1) = 0.$
\end{enumerate}

\item
Give a direct proof that $C^\infty[0, 1]$ is dense in $W^{1,p}(0, 1)$ for all $1 \le p< \infty$.

\item
Fix $1<p<\infty$ and $f \in W^{1,p}(0,1)$, and let $\widetilde{f} \in C[0,1]$ be its continuous version (see Problem \ref{prob:cont-version-W1p}).
\begin{enumerate}[\rm(a), leftmargin=*]
 \item Suppose that $\widetilde{f}(0)=0$. Show that $x\mapsto\frac{f(x)}{x}$ belongs to $L^p(0,1)$ with
 \begin{equation*}
  \Big\n x\mapsto \frac{f(x)}{x}\Big\n_{p} \leq \frac{p}{p-1} \n f'\n_{p}.
 \end{equation*}
 {\em Hint:}\ Use Young's inequality for $(\R_+,\frac{{\rm d} x}{x})$ from Problem \ref{prob:Young-mult}
 with
 $$g(x) = \begin{cases}
  x^{1/p} f'(x), & \ x \in [0,1],\\
  0,       & \ x \in (1,\infty).
 \end{cases}$$
 \item Suppose that $x\mapsto\frac{f(x)}{x}$ belongs to $L^p(0,1)$. Show that $\widetilde{f}(0)=0$.

 \noindent{\em Hint:}\ Argue by contradiction.
 \item Define
 \begin{equation*}
  f(x) = \frac{1}{1-\log x}, \quad x\in (0,1).
 \end{equation*}
 Show that $f \in W^{1,1}(0,1)$ and $\widetilde{f}(0)=0$, but $\frac{f(x)}{x} \notin L^1(0,1)$.
\end{enumerate}

\item
Show that if a real-valued function $f\in L_{\rm loc}^1(D)$ admits weak derivatives $\partial_j f$ and $\rho:\R\to \K$ is a $C^1$-function with bounded derivative, then $\rho\circ f$ admits weak derivatives given by
$$ \partial_j(\rho\circ f) = (\rho'\circ f)\partial_j f. $$
How does this compare with Proposition \ref{prop:chainrule}?

\item
When the scalar field is real, is $H^1(D)$ a Banach lattice?

\item\label{prob:productrule}
Let $k\in\N$, and let $1\le p,q,r\le \infty$ satisfy $\frac1p+\frac1q=\frac1r$.
Prove that if $f\in W^{k,p}(D)$ and $g\in W^{k,q}(D)$, then $fg\in W^{k,r}(D)$, and for all multi-indices $\al$ with $|\al|\le k$ we have the Leibniz formula
\begin{align*}
 \partial^\al (fg) =
\sum_{0\le \beta\le \al} \binom{\al}{\beta} (\partial^\beta f) (\partial ^{\al-\beta}g).
\end{align*}

\item
For $1\le p< \infty$ we consider the weak derivative $\partial$
as a linear operator in $L^p(0,1)$
with domain $\Dom(\partial):= C_{\rm c}^\infty(0,1)$.
\begin{enumerate}[\rm(a), leftmargin=*]
 \item Show that $\partial$ is closable.
 \item Show that the domain of the closure of $\partial$ equals $W_0^{1,p}(0,1)$.
 \item Show that this closure has a proper closed extension, given by the weak derivative with domain $W^{1,p}(0,1)$.
 \item Why doesn't this contradict the result of Proposition \ref{prop:semigroupsAB}?
\end{enumerate}

\item
Show that if $f\in W_0^{1,p}(D)$ with $1\le p< \infty$ and $D'$ is an open set containing $D$,
then the function $\wt f:D'\to \K$ defined by
$$\wt f(x) := \left\{\begin{aligned}
         & f(x), && x\in D, \\
         & 0,  && x\in D'\setminus D,
         \end{aligned}
\right.$$
belongs to $W_0^{1,p}(D')$.

\item
Show that for all $s\ge 0$ the space $H^s(\R^d)$ is a Hilbert space with respect to the norm given by \eqref{eq:norm-Hs}.

\item
Let $D$ be bounded and fix $1\le p<\infty$. Let $r\in \R$ satisfy $r>1-\frac{d}{p}$.
Show that $f(x):= |x|^r$ belongs to $W^{1,p}(D)$, and compute its weak partial derivatives.

\item
Show that for all $f\in W^{2,2}(\R^d)$ we have
$$ \n \nabla f\n_2 \le \n f\n_2^{1/2} \n \Delta f\n_2^{1/2}\!.$$

\noindent{\em Hint:}\
For $f\in C_{\rm c}^\infty(\R^d)$, apply integration by parts to
$\n\nabla f\n_2^2$. Then use the density of $C_{\rm c}^\infty(\R^d)$ in $W^{2,2}(\R^d)$.

\item
Let $f\in W^{1,\infty}(D)$.
\begin{enumerate}[\rm(a), leftmargin=*]
 \item Let $\eta\in C_{\rm c}^\infty(\R^d)$ be nonnegative, with support in the unit ball $B(0;1)$ of $\R^d$ and satisfying $\int_{\R^d}\eta(x)\ud x=1$. For $\eps>0$ denote $\eta_\eps(x):=\eps^{-d}\eta(\eps^{-1}x)$ and, for $x\in D_\eps$, put
 $$ f_\eps(x):=\int_{\R^d}\eta_\eps(y)f(x-y)\ud y,$$
 where $D_\eps := \{x\in D: \, d(x,\partial D)>\eps\}$. Show that $f_\eps\in C^\infty(D_\eps)$ and satisfies the pointwise bound
 $$ | \nabla f_\eps(x)| \le \n \nabla f \n_\infty, \quad x\in D_\eps,$$
 where $\nabla f$ is the weak gradient of $f$.
 \item Show that if $D$ is convex, then
 $$ | f_\eps(x) - f_\eps(y)| \le \n \nabla f \n_\infty |x-y|, \quad x,y\in D_\eps.$$
 \item \label{prob:convex-Lip-a} Deduce that if $D$ is convex, then for every $f\in W^{1,\infty}(D)$ there exists a Lipschitz continuous function $g:D\to \K$ such that $f=g$ almost everywhere, with Lipschitz constant $L_g \le \n \nabla f\n_\infty$.
 \item Show that the result of part \ref{prob:convex-Lip-a} fails for the nonconvex open set in $\R^2$ obtained by removing the nonnegative part of the $x$-axis from $B(0;1)$.
\end{enumerate}

\item
Does Theorem \ref{thm:Rellich}\ref{it:Rellich2} extend to the case $p=\infty$?

\item
For $t\in \R\setminus\{0\}$ let
$$ D^t_j f(x) := \frac1t({f(x+ t e_j) - f(x)}), \quad 1\leq j\leq d,$$
where $e_j$ is the $j$th standard unit vector of $\R^d\!$.
\begin{enumerate}[\rm(a), leftmargin=*]
 \item Prove that if $f\in W^{1,p}(\R^d)$ with $1\le p\le \infty$, then
 $$ \n D^t_j f \n_p \le \n \partial_j f\n_p, \quad 1\leq j\leq d, \ t\in \R\setminus\{0\},$$
 where $\partial_j f$ denotes the $j$th partial derivative of $f$.
 \item Prove that if $1< p<\infty$ and $f\in L^p(\R^d)$ satisfies
 $$ \n D^t_j f \n_p \le C, \quad 1\le j\le d, \ t\in \R\setminus\{0\},$$
 then $f\in W^{1,p}(\R^d)$ and $\n \partial_j f\n_p \le C$ for all $1\le j\le d.$
\end{enumerate}

\item \label{prob:comp-Sob-d1}
Let $1\le p\le \infty$ and consider the inclusion mapping $f\mapsto \wt f$ from $W^{1,p}(0,1)$ to $C[0,1]$ of Problem \ref{prob:cont-version-W1p}.
\begin{enumerate}[\rm(a), leftmargin=*]
 \item\label{it:comp-Sob-d1-3} Show that for
 $1<p\le \infty$
 the inclusion $W^{1,p}(0,1)\subseteq C[0,1]$ is compact.

 \noindent{\em Hint:}\ Use the Arzel\`a--Ascoli theorem.

 \item\label{it:comp-Sob-d1-6} Show that the inclusion $W^{1,1}(0,1)\subseteq C[0,1]$ fails to be compact.

 \noindent{\em Hint:}\ Approximate $\one_{(\frac12,1)}$ pointwise
 by a sequence of piecewise-linear functions that is bounded in $W^{1,1}(0,1)$.
\end{enumerate}

\item
Show that the inclusion mapping $W^{1,2}(\R)\subseteq L^2(\R)$ fails to be compact.

\item
\label{prob:Green}
Consider the Green function on the unit interval $[0,1]$ (see Section \ref{sec:Poisson-D}):
$$ k(x,y) := \begin{cases}
          (1-x)y, & 0\le y\le x, \\
          (1-y)x, & x\le y\le 1.
          \end{cases}
$$
\begin{enumerate}[\rm(a), leftmargin=*]
 \item Show that the associated integral operator
 $$T_k f(x):= \int_0^1 k(x,y)f(y)\ud y $$
 on $L^2(0,1)$ is compact and has eigenvalues $1/(n\pi)^2$\!, $n\ge 1$, with corresponding eigenfunctions $x\mapsto \sin(n\pi x)$.
\end{enumerate}
Now let $f\in L^2(0,1)$ and define $u\in L^2(0,1)$ by $$ u(x):= T_k f(x), \quad x\in (0,1).$$
\begin{enumerate}[\rm(a), leftmargin=*]\setcounter{enumii}{1}
 \item Show that $u\in H_0^1(0,1)$.
 \item Show that $u$ is a weak solution of the Poisson problem with Dirichlet boundary conditions
 \begin{equation*}
 \begin{cases} -u'' = f \ \hbox{ on $(0,1)$}, \\  u(0) =u(1) = 0.
 \end{cases}
 \end{equation*}
\end{enumerate}

\item
Show that if a function $f\in L_{\rm loc}^2(D)$ admits a weak Laplacian in $L_{\rm loc}^2(D)$, then $f$ belongs to $H_{\rm loc}^2(D)$.

 \noindent{\em Hint:} \ First prove that $\psi f\in H^1(D)$ for every test function $\psi\in C_{\rm c}^\infty(D)$. Then use Lemma \ref{lem:MR-d1}.

\item\label{prob:weak-comparison}
Taking the scalar field to be real, let $f,g\in L^2(D)$ and let
$u_f,u_g\in H_0^1(D)$ be the corresponding weak solutions of the
Poisson problem with Dirichlet boundary conditions. Prove the
{\em weak comparison principle}
$$
f\le g\ \hbox{ almost everywhere}
\quad\Longrightarrow\quad
u_f\le u_g\ \hbox{ almost everywhere}.
$$
Deduce that $f\ge0$ almost everywhere implies $u_f\ge0$ almost everywhere.

\noindent{\em Hint:}\ The weak formulation extends by density to all
test functions in $H_0^1(D)$. Apply it to
$(u_f-u_g)^+\in H_0^1(D)$.

\item
Let $V$ be a Hilbert space, let $C\subseteq V$ be a nonempty closed
convex set, and let $E:V\to\R$ be sequentially weakly lower
semicontinuous and coercive, in the sense that
$E(v)\to\infty$ as $\n v\n_V\to\infty$.
Show that $E$ attains its minimum on $C$, and show that the minimiser is unique if $E$ is strictly convex on
$C$.

\noindent{\em Hint:}\ A closed convex subset of a Hilbert space is
weakly closed.

\item
Take the scalar field to be real. Let $D$ be bounded, let
$f\in L^2(D)$ and $\psi\in H_0^1(D)$, and put
$$
K:=\{v\in H_0^1(D):\,v\ge\psi\ \hbox{almost everywhere}\}.
$$
Show that $K$ is closed and convex and that the energy functional
$$
E(v):=\frac12\int_D|\nabla v|^2\ud x-\int_Dfv\ud x
$$
has a unique minimiser $u\in K$. Show that $u\in K$ is this
minimiser if and only if
$$
\int_D\nabla u\cdot\nabla(v-u)\ud x
\ge \int_D f(v-u)\ud x,\quad v\in K.
$$

\noindent{\em Hint:}\ Use the direct method and Poincar\'e's inequality.
For the last assertion, consider $u+t(v-u)$ with $0\le t\le1$.

\item\label{prob:elliptic-LM}
Prove Theorems \ref{thm:elliptic-D} and \ref{thm:elliptic-N}.

\item
Show that if in Theorem \ref{thm:SturmLiouville} we make the stronger assumption
$\Re q\ge q_0>0$ almost everywhere on $D$,
then coercivity holds on $H_0^1(D)$ (with the original norm of $H^1(D)$) without invoking Poincar\'e's inequality, and boundedness of $D$ is no longer needed.

\item\label{prob:SL}
Let $D$ be bounded with $C^1$-boundary, and let $a$, $q$, and $f$
satisfy the assumptions preceding Theorem \ref{thm:SturmLiouville}.
Formulate the weak Neumann problem and prove:
\begin{enumerate}[\rm(a), leftmargin=*]
\item\label{it:prob:SL1} If
$\Re q\ge q_0>0$
almost everywhere on $D$, then the Neumann problem admits a unique weak
solution in $H^1(D)$.
\item\label{it:prob:SL2} Suppose that $q=0$ and that $D$ is connected. Show that a weak
solution can exist only if
$\int_Df\ud x=0$.
Under this compatibility condition, prove that there is a unique weak
solution in the space
$$
H_{\rm av}^1(D)
=
\Bigl\{
u\in H^1(D):\,\int_Du\ud x=0
\Bigr\}.
$$
\item Suppose in addition that $a=a^\star$ almost everywhere. In case
\ref{it:prob:SL1}, assume that $q$ is real-valued and show that the weak solution
is the unique minimiser over $H^1(D)$ of the energy functional
$$
E(v)
=
\frac12\int_D a\nabla v\cdot\overline{\nabla v}\ud x
+
\frac12\int_D q|v|^2\ud x
-
\Re\int_D f\overline v\ud x.
$$
In case \ref{it:prob:SL2}, show that the unique weak solution in
$H_{\rm av}^1(D)$ is the unique minimiser over $H_{\rm av}^1(D)$ of
$$
E(v)
=
\frac12\int_D a\nabla v\cdot\overline{\nabla v}\ud x
-
\Re\int_D f\overline v\ud x.
$$
\end{enumerate}

\item\label{prob:Cea}
Let $V$ be a Hilbert space and consider a bounded coercive form
$\aa: V\times V \to\K$. Let $L: V\to\K$ be a bounded functional.
By the Lax--Milgram theorem there is a unique $u_V\in V$ satisfying
$\aa(v,u_V) = L(v)$ for all $v\in V$.

Suppose now that $W$ is a closed subspace of $V$.
By the Lax--Milgram theorem, applied to the restriction
of $\aa$ to $W\times W$, there is a unique $u_W\in W$ satisfying
$\aa(w,u_W) = L(w)$ for all $w\in W$.
\begin{enumerate}[\rm(a), leftmargin=*]
 \item Show that $\aa(w,u_V - u_W) = 0$ for all $w\in W$.
 \item Show that
 \begin{align}\label{eq:Cea}
 \n u_V-u_W\n \le C\alpha^{-1} \inf_{w\in W}\n u_V-w\n,
 \end{align}
 with $C\ge 0$ and $\alpha>0$ the boundedness and coercivity constants of $\aa$.
 \item Show that if $\aa$ is symmetric, that is, $\aa(v_1,v_2) = \ov{\aa(v_2,v_1)}$ for all $v_1,v_2\in V$, then
 \begin{align}\label{eq:Cea-symm}
 \n u_V -u_W\n \le \sqrt{C\alpha^{-1}} \inf_{w\in W}\n u_V-w\n.
 \end{align}
\end{enumerate}
The {\em quasi-optimality estimates}\index{quasi-optimality estimate} in \eqref{eq:Cea} and
\eqref{eq:Cea-symm} are known as {\em C\'ea's lemma}.\index{lemma!C\'ea}

\item\label{prob:FEM}
In this problem we outline an application to the so-called {\em finite element method}\index{finite!element method} for
the Poisson problem \eqref{eq:Poisson} on the unit interval $(0,1)$ with datum $f\in L^2(0,1)$.

In what follows we endow $V:= H_0^1(0,1)$ with the norm $\n v\n_{H_0^1(0,1)}:= \n v'\n_2$. By the Poincar\'e inequality,
this norm is equivalent to the Sobolev norm
$\n v\n_2+ \n v'\n_2$.

Consider a partition $\pi = \{x_0,\dots,x_N\}$ of the interval $[0,1]$, that is, we assume that $0= x_0<x_1< \dots < x_{N-1} < x_N = 1$. Let $V_\pi$ denote the closed subspace consisting of all $v\in V$ that are linear on each of
the intervals $[x_{n-1},x_n]$.

Define
$$
    \aa(v,u):=\int_0^1 v'(t)\overline{u'(t)}\ud t,\quad u,v\in V.
$$
\begin{enumerate}[\rm(a), leftmargin=*]
 \item\label{it:FEM1} Show that for every $f\in L^2(0,1)$ there exist unique elements $u\in H_0^1(0,1)$ and $u_\pi\in V_\pi $ such that for all $v\in V$ and $ v_\pi\in V_\pi$,
 \begin{equation}\label{eq:Vpi}
 \begin{aligned}
 \aa(v,u)  = \int_0^1 v(x)\ov{f(x)}\ud x, \qquad
 \aa(v_\pi,u_\pi)  = \int_0^1 v_\pi(x)\ov{f(x)}\ud x.
 \end{aligned}
 \end{equation}
 Using the results of Problem \ref{prob:Cea}, prove the quasi-optimality estimate
 \begin{equation*}
 \n u-u_\pi \n_{H_0^1(0,1)} \le \inf_{v_\pi\in V_\pi} \n u-v_\pi\n_{H_0^1(0,1)}.
 \end{equation*}
 \item\label{it:FEM2} Show that for all $v\in H_0^1(0,1)\cap H^2(0,1)$ we have $\pi v \in  V_\pi$ and
 $$ \n v - \pi v\n_{H_0^1(0,1)} \le h_\pi \n v''\n_{L^2(0,1)},$$
 where $\pi v\in V_\pi$ is obtained by piecewise-linear interpolation of the values $v(x_n)$, $0\le n\le N$, and $h_\pi:= \max_{1\le n\le N}|x_{n}-x_{n-1}|$ is the mesh of $\pi$.

 \noindent{\em Hint:}\ Fix $x\in [0,1]\setminus\pi$ and choose $1\le n\le N$ such that $x_{n-1}< x<x_n$.
 Then $$ (\pi v)'(x) = \frac{v(x_n) - v(x_{n-1})}{x_n - x_{n-1}}$$ and, since $v\in H^2(0,1)$,
 \begin{align*} v(x_n)
 = v(x_{n-1}) + (x_n - x_{n-1})v'(x_{n-1}) +
 \int_{x_{n-1}}^{x_n} (x_n-y) v''(y) \ud y.
 \end{align*}
 Rewriting the latter as
 $$ v'(x_{n-1}) - \frac{v(x_n) - v(x_{n-1})}{x_n - x_{n-1}} = -
 \int_{x_{n-1}}^{x_n} \frac{x_n-y}{x_n - x_{n-1}} v''(y) \ud y
 $$
 and using the inequalities
$$
\Big|\frac{y-x_{n-1}}{x_n-x_{n-1}}\Big|\le 1,\qquad
\Big|\frac{x_n-y}{x_n-x_{n-1}}\Big|\le 1,
$$
show that
 \begin{align*} | v'(x) - (\pi v)'(x)|
 \le \int_{x_{n-1}}^{x_n} |v''(y)| \ud y.
 \end{align*}
 \item\label{it:FEM3}

Let the assumptions of Theorem \ref{thm:Poisson-D}
be satisfied with $d=1$ and $D=(0,1)$, and let
$u\in H_0^1(0,1)$ be the weak solution of
\eqref{eq:Poisson}. By Remark \ref{rem:MR-d1},
$u\in H^2(0,1)$. Prove that
 $$ \n u - u_\pi\n_{H_0^1(0,1)} \le h_\pi \n u''\n_{L^2(0,1)}.$$
\end{enumerate}

Since the norm $\n v\n_{H_0^1(0,1)}$ is equivalent to the Sobolev norm
$\n v\n_2+ \n v'\n_2$, the result of part \ref{it:FEM3}
shows that $u_\pi$ and its weak derivative $u_\pi'$ provide
good approximations of $u$ and its weak derivative $u'$ in the $L^2(0,1)$-norm
if $h_\pi$ is small.

The approximate solution $u_\pi$ can be constructed explicitly as follows.
In the remaining explicit matrix description we take the scalar field to be real.
Every $u\in V_\pi$
can be written uniquely as a finite linear combination
$ u = \sum_{n=1}^{N-1} c_n \psi_n,$
where $\psi_n\in V_\pi$ is the piecewise-linear function given by the conditions
$$ \psi_n(x_m) = \begin{cases}
                1, &  \ n=m, \\
                0,&  \ n\not=m,
                \end{cases}
$$
since these functions form a basis for $V_\pi$.
By definition, $u_\pi$ is the unique element of $V_\pi$ solving
\eqref{eq:Vpi}, which, for our boundary value problem, takes the form
$$ \int_0^1 v' u_\pi'\ud x = \int_0^1 vf\ud x, \quad v\in V_\pi.$$
Since the functions $\psi_n$ form a basis for $V_\pi$, this holds if and only if
$$ \int_0^1 \psi_n' u_\pi'\ud x = \int_0^1 \psi_nf\ud x, \quad n = 1,\dots,N-1.$$
Writing $u_\pi=\sum_{m=1}^{N-1}c_m\psi_m$, we reduce the problem to
determining the coefficients $c_1,\dots,c_{N-1}$ from the system of $N-1$ linear equations
$$ \sum_{m=1}^{N-1} c_m \int_0^1 \psi_m' \psi_n'\ud x = \int_0^1 \psi_n f\ud x, \quad n = 1,\dots, N-1.$$
The functions $\psi_n'$ take nonzero constant values on the intervals $(x_{n-1},x_n)$ and $(x_{n},x_{n+1})$
and vanish on the remaining subintervals.
It follows from this that $\int_0^1 \psi_m' \psi_n'\ud x=0 $ unless $m-n \in \{-1,0,1\}$.
Therefore the computation of the coefficients $c_m$ reduces to a matrix problem of the form
$Sc = d$, where the so-called {\em stiffness matrix}\index{stiffness matrix} $S$ is the $(N-1)\times(N-1)$ matrix whose coefficients $$s_{nm} =\int_0^1 \psi_n' \psi_m'\ud x$$
vanish off the diagonal and the two neighbouring off-diagonals, and $$d_n = \int_0^1 \psi_n f\ud x.$$
This problem is easily solved using numerical linear algebra.

\end{problems}

%% file: ch12-LaplaceOperator.tex
\chapter{Forms}\label{chap:forms}

\blfootnote{\rm \small \begin{quote}This book has been published by Cambridge University Press in the series `Cambridge Studies in Advanced Mathematics'. The present version is free to view and download for personal use only. Not for re-distribution, re-sale or use in derivative works.\newline {\copyright} Jan van Neerven\end{quote}}

\noindent
This chapter develops elements of the theory of sesquilinear forms and uses it
to define and study certain bounded and unbounded operators, including second-order differential operators such as the Laplace operator subject to Dirichlet and Neumann boundary conditions. It culminates in a proof of Weyl's theorem on counting the eigenvalues of the Dirichlet Laplacian and a discussion of the landscape
function, which gives spectral information from solving a single Poisson problem.

\section{Forms}\label{sec:forms}

In the previous chapter we proved existence and uniqueness of weak solutions of the Poisson problem $-\Delta u = f$ on nonempty bounded open subsets
$D\subseteq \R^d$ for functions $f\in H = L^2(D)$ by exploiting the properties
of the symmetric sesquilinear mapping $\aa: V\times V\to \K$,
\begin{align*}
\aa(u,v) = \int_{D} \nabla u \cdot \ov{\nabla v}\ud x,
\end{align*}
where $V$ is a suitable closed subspace of $H^1(D)$, the choice of which depends on the imposed boundary conditions.
For bounded measurable uniformly elliptic matrix-valued functions
$a:D\to M_d(\K)$, the sesquilinear mapping
\begin{align}\label{eq:forms-ex2}\aa(u,v) = \int_{D} a\nabla u \cdot \ov{\nabla v}\ud x
\end{align}
played the same role in solving the Sturm--Liouville problem.
In each of these cases, the key ingredient was the Poincar\'e inequality, which can be phrased in terms of $\aa$ as
$$ \Re\aa(v,v) \ge \al\n v\n^2\!, \quad v\in V,$$
where $\al >0$ is a positive constant and the norm is taken in $H=L^2(D)$.

In order to study these matters from an abstract point of view, it will be useful to
interpret a form $\aa$ defined on a subspace $V$ of a Hilbert space $H$ as one in $H$ with domain $\Dom(\aa)=V$, in the same way as the notion of a bounded operator
has been generalised to that of a linear operator $A$ defined on a domain $\Dom(A)$.

\begin{definition}[Forms, accretivity, and coercivity]\label{def:form-Dom}
A {\em form} in a Hilbert space $H$\index{form} is a pair $(\aa, \Dom(\aa))$, where $\Dom(\aa)$ is a subspace of $H$,
the {\em domain} of $\aa$,\index{domain!of a form}
and $\aa:\Dom(\aa)\times \Dom(\aa)\to \K$ is a sesquilinear mapping.
A form $(\aa, \Dom(\aa))$ is called {\em accretive}\index{accretive!form}
if $$ \Re \aa(x,x) \ge 0, \quad x\in \Dom(\aa),$$
and
{\em coercive}\index{coercive form}
if there exists a constant $\al> 0$ such that
$$ \Re \aa(x,x) \ge \alpha \n x\n^2\!, \quad x\in \Dom(\aa).$$
\end{definition}

In what follows, $H$ always denotes a Hilbert space.
The notions of a form, accretivity, and coercivity in
Definition \ref{def:form-bounded} are recovered in the special case
$\Dom(\aa)=H$. A more general compatibility result is proved in Proposition \ref{prop:closed-corecive-2defs}.

When no confusion is likely to arise, we omit $\Dom(\aa)$ from the notation and denote the form by $\aa$.

\begin{example}\label{ex:forms}
Let $D$ be a bounded open subset of $\R^d$. If $a:D\to M_d(\K)$ is measurable, bounded, and uniformly elliptic, then the form in $H = L^2(D)$ defined by \eqref{eq:forms-ex2} is
accretive and bounded on the domain $\Dom(\aa)=H^1(D)$, and coercive on the
domain $\Dom(\aa)=H_0^1(D)$. If, in addition, $D$ is a
bounded connected $C^1$ domain, then this form is also coercive on the
domain $\Dom(\aa)=H^1_{\rm av}(D)$.
\end{example}

If $(\aa, \Dom(\aa))$ is an accretive form in $H$, then
\begin{equation*}
  \iprod{x}{y}_\aa:=\Re\aa(x,y)+\iprod{x}{y}, \quad x,y\in \Dom(\aa),
\end{equation*}
defines an inner product on $\Dom(\aa)$;
here,
$\iprod{x}{y}$ is the inner product of $x$ and $y$ in $H$ and $$ \Re\aa := \frac12(\aa + \aa^\star)$$  is the {\em symmetric part} of $\aa$,\index{symmetric!part} where $\aa^\star(x,y):= \overline{\aa(y,x)}$.
The norm associated with this inner product equals
$$\n x\n_\aa = \iprod{x}{x}_\aa^{1/2}, \quad x\in \Dom(\aa).$$
{\em Warning:} One should not confuse $\Re\aa(x,y) = \frac12(\aa(x,y) + \ov{\aa(y,x)})$ with $\Re(\aa(x,y)) =\frac12(\aa(x,y) + \ov{\aa(x,y)}).$
The former defines a sesquilinear form, the {\em real part} of $\aa$,\index{real part!of a form} but the latter generally does not.
It is true, however, that
$\Re \aa(x,x) = \Re (\aa(x,x))$ for all $x\in \Dom(\aa)$.

\medskip
From Definition \ref{def:form-bounded} we recall that a sesquilinear form $\aa:V\times V\to \K$ is {\em bounded} if there exists a constant $C\ge 0$ such that $$|\aa(u,v)| \le C \n u\n_V \n v\n_V, \quad u,v\in V.$$
This definition can be extended to forms in $H$ as follows.

\begin{definition}[Bounded forms]\label{def:form-continuous}
An accretive
form $(\aa,\Dom(\aa))$ in $H$ is called {\em bounded}\index{bounded!form} if there exists a constant $C\ge 0$ such that
$$|\aa(x,y)| \le C\n x\n_\aa \n y\n_\aa, \quad x,y\in \Dom(\aa).$$
\end{definition}

A sufficient condition for boundedness will be given in Proposition \ref{prop:sect-cont-forms}.

\subsection{Closed Forms}\label{subsec:assoc-op}

The following definition is motivated by the simple fact, observed in
Proposition \ref{prop:closed-op}, that
a linear operator $A$ is closed if and only if its domain $\Dom(A)$ is a Banach space with respect to the graph norm.

\begin{definition}[Closed forms]\label{def:form-closed}
An accretive form $(\aa,\Dom(\aa))$ in $H$ is called {\em closed}\index{closed!form} if $\Dom(\aa)$ is a
 Hilbert space with respect to the norm $\n \cdot\n_\aa$.
\end{definition}

The following proposition clarifies the relation between
Definition \ref{def:form-Dom}, where coercivity of forms in $H$
was defined through the condition
\begin{align*}
 \Re \aa (x,x) &\ge \al \n x\n^2\!,\quad x\in\Dom(\aa),
\intertext{and Definition \ref{def:form-bounded},
where coercivity of a form on a Hilbert space $V$ was defined through the condition
          }
 \Re \aa (x,x) & \ge \al \n x\n_V^2,\quad x\in V.
\end{align*}
In the former case, one could view $\aa$ as a form on $V = \Dom(\aa)$ and ask why norms are taken in $H$ rather than in $V$. As it turns out,
except for the numerical value of the constant, this leads to the same definition.

\begin{proposition}\label{prop:closed-corecive-2defs}
A closed form $(\aa,\Dom(\aa))$ in $H$ is accretive (respectively coercive, bounded) if and only if $\aa$, as a form on the Hilbert space $V := (\Dom(\aa), \n\cdot\n_\aa)$, is accretive (respectively coercive, bounded).
\end{proposition}
\begin{proof} Only the assertion concerning coercivity needs proof. We must prove that
there exists a constant $\al>0$ such that
\begin{align}\label{eq:form-accr1}
\Re \aa(x,x) \ge \al\n x\n^2\!, \quad x\in\Dom(\aa),
\end{align}
if and only if there exists a constant $\beta>0$ such that
\begin{align}\label{eq:form-accr2}
\Re \aa(x,x) \ge \beta\n x\n_\aa^2, \quad x\in\Dom(\aa).
\end{align}
If \eqref{eq:form-accr1} holds, then for all $x\in\Dom(\aa)$ we have
$$(1+\al)\Re \aa(x,x) \ge \al \Re\aa(x,x) + \alpha \n x\n^2 = \al \n x\n_\aa^2$$ and therefore
\eqref{eq:form-accr2} holds with $\beta = \frac{\al}{1+\al}$.

Conversely, if \eqref{eq:form-accr2} holds, then
$$
\Re\aa(x,x)\ge\beta\n x\n_\aa^2\ge\beta\n x\n^2,
\quad x\in\Dom(\aa),
$$
so \eqref{eq:form-accr1} holds with $\alpha=\beta$.
\end{proof}

From now on we will omit the domain $\Dom(\aa)$ from the terminology if no confusion is likely to arise. Thus, the phrase `$\aa$ is a form in $H$' will always mean `$(\aa,\Dom(\aa))$ is a form in $H$'.

\begin{remark}\label{rem:Gelfand-triple}
Let $V$ and $H$ be Hilbert spaces, with $V$ \emph{densely embedded}
in $H$, that is, there is a bounded injective operator
$\iota:V\to H$ with dense range.
After identifying $V$ with its image $\iota(V)$, a bounded
accretive form
$$
\aa:V\times V\to\K
$$
may be viewed as an accretive form in $H$ with domain $\iota(V)$.
It is closed as a form on $V$, but as a form in $H$ it need not be
closed in general; see Problem \ref{prob:closedVnonclosedH}.
This motivates the study of {\em Gelfand triples},\index{Gelfand triple}
that is, triples $(\iota,V,H)$ with the properties just mentioned. Some properties of such triples will be explored in the problems at the end of this chapter.

We mention two standard applications of Gelfand triples to partial
differential equations. In both applications the space $V$ describes
the spatial regularity of the solution.
Let $V'$ denote the
\emph{conjugate dual}\index{conjugate dual} of $V$, that is, the vector space of all mappings $\phi:V\to \K$ that are conjugate-linear in the sense that
$$ \phi(u+v) = \phi(u)+\phi(v), \qquad \phi(cv) = \ov c \phi(v), \quad u,v\in V, \ c\in\K,$$
and bounded in the sense that
$$|\phi(v)|\le C\n v\n_V, \quad v\in V.$$ The space $V'$ is a Banach space in a natural way with norm $$\n \phi\n_{V'} := \sup_{\n v\n_V\le 1} |\phi(v)|.$$
Every element $h\in H$ naturally defines an element $\phi_h\in V'$ by
$$
\phi_h(v):=\iprod{h}{\iota(v)},\quad v\in V.
$$
The mapping $\phi:H\to V'$, $h\mapsto\phi_h$, is linear.

First, weak formulations of boundary value problems often make sense for forcing terms $f$ which are not in $H$, but belong to $V'$.
For example, for the Dirichlet problem one naturally takes $V=H^1_0(D)$ and $H=L^2(D)$, and in this setting $V'=H^{-1}(D)$ accommodates data $f$ with negative Sobolev regularity. One may then consider equations of the form
$$
\aa(u,v)=\langle f,v\rangle,\quad v\in V,
$$
where $\langle\cdot,\cdot\rangle$ denotes the duality pairing between $V'$ and $V$.

Second, for inhomogeneous time-dependent problems of the form
$$
u'(t)+Au(t)=f(t)
$$
one often looks for solutions $t\mapsto u(t)$ taking values in the smaller space $V$, allowing the time derivative $t\mapsto u'(t)$ to take values in the larger space $V'$.
Such problems can be treated in a weak sense, requiring $u\in L^2(0,T;V)$ and $u'\in L^2(0,T;V')$.
\end{remark}

In many applications, one starts out by defining a form on an initial space and then proceeds by extending the form to a larger domain in which it enjoys certain closedness properties. This idea is captured by the following definition.

\begin{definition}[Closable forms]\label{def:form-closable}
An accretive form $\aa$ in $H$ is called {\em closable}\index{closable!form} if there exists a closed accretive form $\wt\aa$ in $H$
 extending $\aa$, that is, $\wt\aa$ is closed and accretive, $\Dom(\aa)\subseteq\Dom(\wt\aa)$, and $\wt\aa(u,v) =\aa(u,v)$ for all $u,v\in\Dom(\aa).$
\end{definition}

The following proposition, in which we view $\Dom(\aa)$ as a (not necessarily complete) normed space with norm $\n\cdot\n_\aa$, gives a useful necessary and sufficient condition for a form $\aa$ in $H$ to be closable. It should be compared with Proposition \ref{prop:closed-sequential2}.

\begin{proposition}\label{prop:form-closable} For a bounded accretive form $\aa$ in $H$ the following assertions are equivalent:
\begin{enumerate}[label={\rm(\arabic*)}, leftmargin=*]
 \item\label{it:form-closable1}  $\aa$ is closable;
 \item\label{it:form-closable2}  every Cauchy sequence in $\Dom(\aa)$ converging to $0$ in $H$ converges to $0$ in $\Dom(\aa)$.
\end{enumerate}
\end{proposition}

\begin{proof}
Put $V:=\Dom(\aa)$, equipped with the norm
$\n\cdot\n_V:=\n\cdot\n_\aa$.

\smallskip
\ref{it:form-closable1}$\Rightarrow$\ref{it:form-closable2}:\
Let $\wt\aa$ be a closed accretive extension of $\aa$. If
$(v_n)_{n\ge1}$ is Cauchy in $V$ and $v_n\to0$ in $H$, then it is Cauchy
in $\Dom(\wt\aa)$ with respect to $\n\cdot\n_{\wt\aa}$. Hence it converges
in this norm to some $v\in\Dom(\wt\aa)$. Since $\wt\aa$ is a form in $H$,
we have $\Dom(\wt\aa)\subseteq H$, and accretivity gives
$$
\n w\n^2\le \n w\n_{\wt\aa}^2,
\quad w\in\Dom(\wt\aa).
$$
Thus $v_n\to v$ in $H$. Since also $v_n\to0$ in $H$, we have $v=0$.
As $\wt\aa$ extends $\aa$, it follows that
$\n v_n\n_V=\n v_n\n_{\wt\aa}\to0$, hence $v_n\to 0$ in $\Dom(\aa)$.

\smallskip
\ref{it:form-closable2}$\Rightarrow$\ref{it:form-closable1}:\
Let $\widehat V$ be the completion of $V$. Since
$\n v\n\le\n v\n_V$, the inclusion mapping
$j:V\to H$ extends to a bounded operator
$\widehat j:\widehat V\to H$. By boundedness, $\aa$ extends uniquely to
a bounded accretive sesquilinear form $\widehat\aa$ on $\widehat V$.

We claim that $\widehat j$ is injective. Indeed, suppose that
$\widehat j\,\widehat v=0$, and choose $(v_n)_{n\ge1}$ in $V$ such that
$v_n\to\widehat v$ in $\widehat V$. Then $(v_n)_{n\ge1}$ is Cauchy in
$V$, and the continuity of $\widehat j$ implies that $v_n=jv_n = \widehat jv_n \to\widehat j\,\widehat v=0$ in $H$. By \ref{it:form-closable2}, $v_n\to0$ in $V$, and therefore
$\widehat v=0$.

Put $\ov V:=\widehat j(\widehat V)\subseteq H$ and define
$$
\ov\aa(\widehat j\,\widehat u,\widehat j\,\widehat v)
:=\widehat\aa(\widehat u,\widehat v),
\quad \widehat u,\widehat v\in\widehat V.
$$
By injectivity of $\widehat j$, this is well defined. It is accretive and
extends $\aa$. Moreover,
$$
\n\widehat j\,\widehat v\n_{\ov\aa}^2
=
\n\widehat j\,\widehat v\n^2
+\Re\widehat\aa(\widehat v,\widehat v)
=
\n\widehat v\n_{\widehat V}^2,
$$
where the last identity follows by continuity from the corresponding
identity on the dense subspace $V$ of $\widehat V$. Thus $\widehat j$ is
an isometric isomorphism from $\widehat V$ onto
$(\ov V,\n\cdot\n_{\ov\aa})$. The latter space is therefore complete,
so $\ov\aa$ is closed.
\end{proof}

The form $\ov\aa$ constructed in the above proof is called the {\em closure}\index{closure!of a closable form} of $\aa$. This extension is minimal among closed bounded accretive extensions. Indeed, let $\wt\aa$ be a closed
bounded accretive extension of $\aa$ with domain $\wt V$.
The inclusion $V\to\wt V$ is isometric for the form norms
and therefore extends to an isometry
$\iota:\widehat V\to\wt V$. Its composition with the inclusion
$\wt V\to H$ equals $\widehat j$, since these maps agree on
the dense subspace $V$. Hence $\ov V\subseteq\wt V$, and
continuity gives
$$
\wt\aa(\widehat j\,\widehat u,\widehat j\,\widehat v)
=\widehat\aa(\widehat u,\widehat v)
=\ov\aa(\widehat j\,\widehat u,\widehat j\,\widehat v),
\quad \widehat u,\widehat v\in\widehat V.
$$
Thus $\wt\aa$ extends $\ov\aa$.

\subsection{Coercive Forms}

One of the main purposes of the theory of forms is to provide a unified method to rigorously treat domain issues for certain classes of operators.
\begin{definition}[The linear operator associated with a densely defined form]
\label{def:op-form}
The {\em operator associated with a densely defined form $\aa$} in $H$
is the linear operator $(A,\Dom(A))$ defined by
\begin{align*}
\Dom(A) & :=
\big\{u\in\Dom(\aa):\,\text{ there exists $h_u\in H$ such that }
\iprod{h_u}{v}=\aa(u,v) \text{ for all } v\in\Dom(\aa)\big\}, \\
Au & :=h_u,\quad u\in\Dom(A).
\end{align*}
\end{definition}

Since $\Dom(\aa)$ is dense in $H$, the element $h_u \in H$, if it exists, is uniquely defined, and thus
$A$ is well defined as a linear operator in $H$ (linearity being a consequence of uniqueness).

Without imposing further properties on $\aa$ this definition is not very useful.
Under appropriate additional assumptions on $\aa$, the next theorem provides some interesting properties of the associated operator.

\Needspace{6\baselineskip}
\begin{theorem}[Resolvent estimates for accretive bounded forms]\label{thm:form-accretive-resolvent}
The operator $A$ associated with a densely
defined, closed, accretive, bounded form $\aa$ in $H$ is densely defined and closed, and for all $\la\in \C$ with $\Re\lambda>0$ we have $\lambda\in\varrho(-A)$ and
\begin{align*}
\|(\lambda+A)^{-1}\| & \le \frac{1}{\Re\lambda}.
\end{align*}
\end{theorem}

In combination with arguments from semigroup theory, the theorem says that $-A$ is the generator of a $C_0$-semigroup of contractions.

A resolvent estimate for coercive bounded forms will be given in Corollary \ref{cor:form-coercive-resolvent}.

\begin{proof}
Put $V:=(\Dom(\aa),\n\cdot\n_\aa)$ and fix $\lambda\in\C$ with
$\Re\lambda>0$. The form
$$
\bb_\lambda(u,v):=\aa(u,v)+\lambda\iprod{u}{v},
\quad u,v\in V,
$$
is bounded on $V$, and
$$
\Re\bb_\lambda(v,v)
=\Re\bb_\lambda^\star(v,v)
=\Re\aa(v,v)+(\Re\lambda)\n v\n^2
\ge\min\{1,\Re\lambda\}\n v\n_V^2.
$$
Thus both $\bb_\lambda$ and $\bb_\lambda^\star$ are coercive
on $V$. For $f\in H$, the functional $v\mapsto\iprod{v}{f}$
is bounded on $V$.
Applying the Lax--Milgram theorem to
$\bb_\lambda^\star$ and taking complex conjugates, we obtain
a unique $u\in V$ such that
\begin{equation}\label{eq:resolvent-variational}
\aa(u,v)+\lambda\iprod{u}{v}=\iprod{f}{v},
\quad v\in V.
\end{equation}
By the definition of the associated operator, this is equivalent
to saying that $u\in\Dom(A)$ and $(\lambda+A)u=f$. In particular, $\lambda+A$ is
bijective.

Taking $v=u$ in \eqref{eq:resolvent-variational} gives
$$
(\Re\lambda)\n u\n^2+\Re\aa(u,u)
=\Re\iprod{f}{u}
\le\n f\n\,\n u\n.
$$
Accretivity therefore implies that $(\lambda+A)^{-1}$ is bounded and
$$
\n(\lambda+A)^{-1}\n\le\frac1{\Re\lambda}.
$$
It also follows that $\lambda+A$, and hence $A$, is closed.
To prove that $\Dom(A)$ is dense, let $x\in\Dom(A)^\perp$
and put $u:=(I+A)^{-1}x$. The preceding identity with
$\lambda=1$ and $f=x$ gives
$$
\n u\n^2+\Re\aa(u,u)=\Re\iprod{x}{u}=0.
$$
Thus $u=0$ by accretivity, and therefore $x=(I+A)u=0$. This proves density.
\end{proof}

An elegant application is the following duality result. Recall that if $\aa$ is a form in $H$, we defined
$\aa^\star(x,y):= \overline{\aa(y,x)}$ for $x,y\in\Dom(\aa)$.

\begin{corollary}[$A^\star$ is associated with $\aa^\star$]\label{cor:form-dual}
Let $A$ be the densely defined closed operator associated with a densely defined closed bounded accretive form $\aa$ in $H$.
Then $A^\star$ is the densely defined closed operator associated with the form $\aa^\star$.
\end{corollary}

\begin{proof}
Since $\Dom(\aa^\star)=\Dom(\aa)$ by definition, $\aa^\star$ is densely defined.
Moreover, for $v\in\Dom(\aa)$,
$$
\|v\|_{\aa^\star}^2=\|v\|^2+\Re\aa^\star(v,v)
=\|v\|^2+\Re\ov{\aa(v,v)}
=\|v\|^2+\Re\aa(v,v)=\|v\|_{\aa}^2,
$$
so $\aa^\star$ is closed, accretive, and bounded.

Let $B$ be the operator associated with $\aa^\star$.
If $x\in\Dom(B)$, then for all $y\in\Dom(A)$ we have
$$
\iprod{y}{Bx}
=\ov{\iprod{Bx}{y}}
=\ov{\aa^\star(x,y)}
=\aa(y,x)
=\iprod{Ay}{x}.
$$
Thus $x\in\Dom(A^\star)$ and $A^\star x=Bx$. This proves the inclusion $B\subseteq A^\star$.

By Theorem \ref{thm:form-accretive-resolvent}, applied to
$\aa^\star$ and $\aa$, both $I+B$ and $I+A$ are boundedly
invertible. Hence so is $I+A^\star$, with inverse
$((I+A)^{-1})^\star$. Thus
$-1\in\varrho(B)\cap\varrho(A^\star)$, and
Proposition \ref{prop:semigroupsAB}, together with
$B\subseteq A^\star$, gives $B=A^\star$.
\end{proof}

\subsection{Sectorial Forms}

For $\omega\in (0,\pi)$ we let
$$
\Sigma_\omega:=\{z\in\C\setminus\{0\}:\, |\arg z|<\omega\}
$$
denote the open sector of angle $\omega$ (arguments are taken in $(-\pi,\pi)$).

\begin{definition}[Sectorial forms]
Let $0<\om\le\frac12\pi$. A form $\aa$ in $H$ is called
\emph{$\omega$-sectorial}\index{sectorial form}  if
$ \aa(v):= \aa(v,v) \in \overline{\Sigma_\omega}$ for all $v\in \Dom(\aa).$
\end{definition}

Sectorial forms of angle less than $\frac12\pi$ are bounded:

\begin{proposition}\label{prop:sect-cont-forms}
Let $\aa$ be an $\omega$-sectorial form in $H$ with $0<\omega<\frac{\pi}{2}$. Then $\aa$ is accretive and bounded and, for all $u,v\in\Dom(\aa)$,
$$
  |\aa(u,v)| \le (1+\tan\omega)\,
  \bigl(\Re\aa(u,u)\bigr)^{1/2}\,\bigl(\Re\aa(v,v)\bigr)^{1/2}.
$$
\end{proposition}

\begin{proof}
Accretivity is clear. By the Cauchy--Schwarz inequality applied to the form $\Re\aa=\frac12(\aa+\aa^\star)$ with $\aa^\star(u,v):=\overline{\aa(v,u)}$,
$$
  \bigl|\Re\aa(u,v)\bigr| \le \bigl(\Re\aa(u,u)\bigr)^{1/2}\,\bigl(\Re\aa(v,v)\bigr)^{1/2}.
$$

Consider the form $\Im \aa = \frac1{2i}(\aa-\aa^\star)$.
Fix $u,v\in \Dom(\aa)$.
Replacing $v$ by $e^{i\theta}v$ if necessary, we may assume that
$\Im \aa(u,v)\in\R$. Then, by direct computation,
\begin{align*}
  \Im\aa(u,v) &
= \frac14\Bigl(\Im\aa(u+v,u+v) - \Im\aa(u-v,u-v)\Bigr)
\end{align*}
and therefore, by $\om$-sectoriality,
\begin{align*}
|\Im\aa(u,v)| & \le \frac14\tan\om \bigl(\Re\aa(u+v,u+v) + \Re\aa(u-v,u-v)\bigr)
 \\ & = \frac12\tan\om \bigl(\Re\aa(u) + \Re\aa(v)\bigr).
 \end{align*}
 Replacing $u$ and $v$ with $\sqrt \eps u$ and $v/\sqrt{\eps}$ and taking the infimum over $\eps>0$, we obtain
$$
|\Im\aa(u,v)|
\le
\tan\omega
\big(\Re\aa(u,u)\big)^{1/2}
\big(\Re\aa(v,v)\big)^{1/2}.
$$
This remains valid when either factor on the right-hand side vanishes,
as can be seen by letting $\eps\to\infty$ or $\eps\downarrow0$.

Combining the estimates for the real and imaginary parts, we obtain the result.
\end{proof}

\begin{theorem}[Resolvent estimates for sectorial bounded forms]\label{thm:form-sectorial-resolvent}
Let $\aa$ be a densely defined closed $\om$-sectorial form in a Hilbert space $H$, with $\om\in (0,\frac12\pi)$. Let $A$ be the operator associated with $\aa$.
Then $A$ is densely defined and closed, and the following assertions hold:

\begin{enumerate}[label={\rm(\arabic*)}, leftmargin=*]
\item\label{it:form-sectorial-resolvent-master} $\sigma(A)\subseteq \ov{\Sigma_\omega}$,  and for all $\lambda\in\Sigma_{\pi-\omega}$ we have
$\displaystyle \|(\lambda+A)^{-1}\|
\le \frac{1}{d(-\lambda,\ov{\Sigma_\omega})}$;
\item\label{it:sectorial-resolvent3}
if $\lambda\in \Sigma_{\frac12\pi-\omega}$, then
$\displaystyle\|(\lambda+A)^{-1}\|\le \frac1{|\lambda|}$;

\item\label{it:form-sectorial-resolvent2}
if $\lambda\in \Sigma_{\pi-\phi}$ with $\phi\in(\omega,\pi)$, then
$\displaystyle \|(\lambda+A)^{-1}\|
\le \frac{1}{\sin(\phi-\omega)}\frac1{|\lambda|}$.
\end{enumerate}
\end{theorem}

Note that \ref{it:form-sectorial-resolvent2} only provides new information for
angles $\phi\in (\om,\om+\frac12\pi]$; for angles $\phi \in (\om+\frac12\pi,\pi)$, part \ref{it:sectorial-resolvent3} gives already a better estimate.

In combination with arguments from semigroup theory, part
\ref{it:sectorial-resolvent3}, together with the Hille--Yosida theorem (Theorem \ref{thm:HilleYosida}) applied to the rotated operators, implies that $-A$ generates an
analytic $C_0$-semigroup of contractions on
$\Sigma_{\frac12\pi-\omega}$; see Theorem
\ref{thm:lumerphilips} and its proof.

\begin{proof}
By Proposition \ref{prop:sect-cont-forms} and Theorem
\ref{thm:form-accretive-resolvent}, $A$ is densely defined
and closed. By Corollary \ref{cor:form-dual}, $A^\star$
is the operator associated with $\aa^\star$. Since
$\aa^\star(u,u)=\ov{\aa(u,u)}$, this form is also
$\omega$-sectorial.

Let $\lambda\in\Sigma_{\pi-\omega}$. Then $-\lambda\notin \ov{\Sigma_\omega}$,
so $d_\lambda:=d(-\lambda,\ov{\Sigma_\omega})>0$.
For $u\in\Dom(A)$ with $\n u\n=1$, by $\om$-sectoriality we have
$$
\n(\lambda+A)u\n
\ge|\iprod{(\lambda+A)u}{u}|
=|\lambda+\aa(u,u)|
\ge d_\lambda.
$$
By homogeneity,
$$
\n(\lambda+A)u\n\ge d_\lambda\n u\n,
\quad u\in\Dom(A).
$$
Proposition \ref{prop:closed-range-unbdd} therefore shows
that $\lambda+A$ is injective and has closed range.
Applying the same argument to $A^\star$ and $\ov\lambda$
gives that $\ov\lambda+A^\star$ is injective. Hence
$$
\Ran(\lambda+A)^\perp
=\Ker(\ov\lambda+A^\star)
=\{0\}.
$$
Thus $\lambda+A$ has dense range, and therefore is surjective, and
$$
\n(\lambda+A)^{-1}\n
\le \frac1{d_\la} = \frac1{d(-\lambda,\ov{\Sigma_\omega})},
\quad \lambda\in\Sigma_{\pi-\omega}.
$$
In particular, $\Sigma_{\pi-\omega}\subseteq\varrho(-A)$,
which proves \ref{it:form-sectorial-resolvent-master}.

\smallskip
\ref{it:sectorial-resolvent3}:\
If $\lambda\in\Sigma_{\frac12\pi-\omega}$, then
$\Re(\lambda\ov w)\ge0$ for every $w\in \ov{\Sigma_\omega}$.
Consequently,
$$
|\lambda+w|^2
=|\lambda|^2+2\Re(\lambda\ov w)+|w|^2
\ge|\lambda|^2.
$$
Since equality holds for $w=0$, we have
$d(-\lambda,\ov{\Sigma_\omega})=|\lambda|$. The preceding resolvent
estimate gives the assertion.

\smallskip
\ref{it:form-sectorial-resolvent2}:\
Let $\phi\in(\omega,\pi)$ and
$\lambda\in\Sigma_{\pi-\phi}$. Elementary geometry gives
$$
d(-\lambda,\ov{\Sigma_\omega})
\ge|\lambda|\sin(\phi-\omega).
$$
Combining this with the same resolvent estimate gives
\[
\n(\lambda+A)^{-1}\n
\le\frac1{\sin(\phi-\omega)}\frac1{|\lambda|}.
\qedhere\]
\end{proof}

\begin{corollary}[Resolvent estimate for coercive bounded forms]
\label{cor:form-coercive-resolvent}
Let $\aa$ be a densely defined closed bounded coercive form
in $H$, and let $A$ be its associated operator. Then
$$
\|(\lambda+A)^{-1}\|
\le \frac{C}{\alpha}\frac1{|\lambda|},
\quad \Re\lambda>0,
$$
where $C$ and
$\alpha$ are the boundedness and coercivity constants of $\aa$ as a form on
$V:=(\Dom(\aa),\|\cdot\|_\aa)$.
\end{corollary}

\begin{proof}
For $u\in V$,
$$
|\aa(u,u)|\le C\|u\|_V^2
\le \frac{C}{\alpha}\Re\aa(u,u).
$$
If $C>\alpha$, the form is $\omega$-sectorial with
$\cos\omega=\alpha/C$. Taking $\phi=\frac12\pi$ in
part \ref{it:form-sectorial-resolvent2} of Theorem
\ref{thm:form-sectorial-resolvent} gives the estimate.
If $C=\alpha$, apply that theorem with any
$\omega\in(0,\frac12\pi)$ and let $\omega\downarrow0$.
\end{proof}

In the language of semigroup theory, this corollary says that $-A$ is the generator of a bounded analytic $C_0$-semigroup.

\subsection{Symmetric Forms}\label{subsec:symm-pos-forms}

Recall that a sesquilinear form $\aa$ in $H$ is accretive if $\Re\aa(u,u)\ge 0$ for all $u\in\Dom(\aa)$. In this section we study forms satisfying the stronger condition $\aa(u,u)\ge 0$ for all $u\in\Dom(\aa)$.

\begin{definition}[Symmetry and positivity]
Let $\aa$ be a sesquilinear form in $H$.
\begin{enumerate}[label={\rm(\roman*)}, leftmargin=*]
\item $\aa$ is \emph{positive}\index{positive!form} if $\aa(u,u)\ge 0$ for all $u\in\Dom(\aa)$.
\item $\aa$ is \emph{symmetric}\index{symmetric!form} if $\aa(u,v)=\ov{\aa(v,u)}$ for all $u,v\in\Dom(\aa)$.
\end{enumerate}
\end{definition}

\begin{proposition}[Positive forms]\label{prop:form-continuity}
For every sesquilinear form
$\aa:\Dom(\aa)\times\Dom(\aa)\to\C$
the following assertions are equivalent:
\begin{enumerate}[label={\rm(\arabic*)}, leftmargin=*]
 \item\label{it:form-continuity1} $\aa$ is positive;
 \item\label{it:form-continuity2} $\aa$ is accretive and symmetric.
\end{enumerate}
If these equivalent conditions are satisfied, then $|\aa(u,v)| \le \n u\n_\aa \n v\n_\aa$; in particular, $\aa$ is bounded.
\end{proposition}
\begin{proof}
\ref{it:form-continuity1}$\Rightarrow$\ref{it:form-continuity2}: \
It is clear that positivity implies accretivity. To prove symmetry, put
$$
        q(w):=\aa(w,w),\quad w\in\Dom(\aa).
$$
By positivity, $q(w)\in\R$ for all $w\in\Dom(\aa)$. Since the form is
linear in the first variable, the polarisation identity gives
$$
        \aa(u,v)
        =
        \frac14\bigl(
        q(u+v)-q(u-v)
        +i q(u+iv)-i q(u-iv)
        \bigr),
        \quad u,v\in\Dom(\aa).
$$
Since $q$ is real-valued and $q(c w)=|c|^2q(w)$, this
identity implies
$\aa(u,v)=\overline{\aa(v,u)}$.

\smallskip
\ref{it:form-continuity2}$\Rightarrow$\ref{it:form-continuity1}: \
Applying symmetry with $u=v\in \Dom(\aa)$, we obtain $\aa(u,u)=\overline{\aa(u,u)}$, so
$\aa(u,u)\in \R$. Since $\Re \aa(u,u)\ge0 $ by accretivity, it follows that $\aa(u,u)\ge0 $ for all $u\in\Dom(\aa)$, that is, $\aa$ is positive.

\smallskip
It remains to prove the final assertion. Since $\aa$ is positive and symmetric, for all $u,v\in\Dom(\aa)$ we have the Cauchy--Schwarz inequality
$$
|\aa(u,v)|^2\le \aa(u,u)\,\aa(v,v).
$$
Since trivially $\aa(w,w) \le \|w\|^2+\aa(w,w)=\|w\|_{\aa}^2$ for $w\in\Dom(\aa)$,
it follows that $|\aa(u,v)|\le \|u\|_{\aa}\,\|v\|_{\aa}$. This says that $\aa$ is bounded.
\end{proof}

We have the following result for the operator associated with a densely defined, closed, positive form.

\begin{proposition}\label{prop:repr-closed-symm-form}
Let $\aa$ be a densely defined, closed, positive form in $H$. Then:

\begin{enumerate}[label={\rm(\arabic*)}, leftmargin=*]
\item\label{it:form-repr1} the associated operator $A$ is positive and selfadjoint;

\item\label{it:form-repr2} $\Dom(A)$ is dense in $\Dom(\aa)$ with respect to the form norm.
\end{enumerate}
\end{proposition}

In Theorem \ref{thm:Kato} we will prove that, under the same assumptions,
$\Dom(\aa)=\Dom(A^{1/2}).$

\begin{proof}
\ref{it:form-repr1}: \ By Proposition~\ref{prop:form-continuity}, positivity of $\aa$ implies that
$\aa$ is symmetric, accretive, and bounded. Hence
$\aa^\star=\aa$, and Corollary~\ref{cor:form-dual} gives $A^\star=A$.
Moreover, for $u\in\Dom(A)$ we have
$\iprod{Au}{u}=\aa(u,u)\ge 0$,
so $A$ is positive.

\smallskip
\ref{it:form-repr2}:\
Let $u\in\Dom(\aa)$ be orthogonal to $\Dom(A)$ with respect
to the form inner product. By symmetry of $\aa$, for every
$v\in\Dom(A)$ we have
$$
0=\iprod{u}{v}_\aa
=\iprod{u}{v}+\aa(u,v)
=\iprod{u}{(I+A)v}.
$$
Since $I+A$ is surjective by
Theorem \ref{thm:form-accretive-resolvent}, it follows that
$u=0$. Thus $\Dom(A)$ is dense in $\Dom(\aa)$ with respect
to the form norm.
\end{proof}

As an application, we obtain the following uniqueness result.

\begin{proposition}\label{prop:uniq-closed-symm-form}
If $\aa$ and $\bb$ are densely defined, closed, positive forms with the same associated operator $A$, then $\aa = \bb$ with equal domains.
\end{proposition}

\begin{proof}
By Proposition~\ref{prop:repr-closed-symm-form}, $\Dom(A)$ is dense in
both form domains. Moreover, for $u\in\Dom(A)$,
$$
\n u\n_\aa^2
=
\n u\n^2+\iprod{Au}{u}
=
\n u\n_\bb^2.
$$
Thus the identity mapping on $\Dom(A)$ extends uniquely to a unitary operator
$$
U:(\Dom(\aa),\n\cdot\n_\aa)
\to
(\Dom(\bb),\n\cdot\n_\bb).
$$

Let $j_\aa:\Dom(\aa)\to H$ and $j_\bb:\Dom(\bb)\to H$ be the inclusion
mappings. These mappings are bounded and
$j_\bb U=j_\aa$ on the dense subspace $\Dom(A)$, hence on all of
$\Dom(\aa)$. Consequently, $Uu=u$ as an element of $H$. It follows
that $\Dom(\aa)=\Dom(\bb)$
and $\n u\n_\aa=\n u\n_\bb$ for all $u\in\Dom(\aa)$.
In particular, it follows that $\aa(u,u)=\bb(u,u)$ on the common domain. Since both forms are symmetric, polarisation gives $\aa=\bb$.
\end{proof}

It has been shown in Corollary \ref{cor:sa-c0sgr} that if
 $A$ is a densely defined closed operator which is positive in the sense that
 $\iprod{Ax}{x} \ge 0$ for all $x\in \Dom(A)$ and has the property that
$I+A$ has dense range, then $A$ is selfadjoint.
The next theorem states that if we drop the dense range condition,
selfadjoint extensions still exist.

\begin{theorem}[Friedrichs extension]\label{thm:Friedrichs}\index{theorem!Friedrichs}\index{Friedrichs extension}
 Let $A$ be a densely defined positive operator acting in a complex Hilbert space $H$.
Then:

\begin{enumerate}[label={\rm(\arabic*)}, leftmargin=*]
 \item\label{it:Friedrichs1} the form $\aa$ in $H$ given by $\Dom(\aa):=\Dom(A)$ and
$$ \aa(x,y):=\iprod{Ax}{y}, \quad x,y\in \Dom(A),$$
is densely defined, positive, bounded, and closable;
 \item\label{it:Friedrichs2} the operator associated with the
 closure of $\aa$ is a positive selfadjoint extension of $A$.
\end{enumerate}
\end{theorem}

\begin{proof}
\ref{it:Friedrichs1}: \
It is clear that $\aa$ is densely defined and positive, and boundedness of $\aa$ follows from Proposition \ref{prop:form-continuity}.
To prove closability, let $(v_n)_{n\ge1}$ be Cauchy in the
form norm and suppose that $v_n\to0$ in $H$. Put
$M:=\sup_n\n v_n\n_\aa<\infty$. Since $\aa$ is symmetric,
$$
\iprod{v_n}{v_m}_\aa
=\iprod{v_n}{v_m}+\aa(v_n,v_m)
=\iprod{v_n}{(I+A)v_m}.
$$
Writing $v_n=(v_n-v_m)+v_m$ in the second argument of the
form inner product and applying the Cauchy--Schwarz inequality,
we obtain
$$
\n v_n\n_\aa^2
\le M\n v_n-v_m\n_\aa
   +|\iprod{v_n}{(I+A)v_m}|.
$$
For fixed $m$, the last term tends to zero as $n\to\infty$.
Taking the limit superior as $n\to\infty$ and then letting
$m\to\infty$, we obtain $\n v_n\n_\aa\to0$ by the Cauchy property.
Proposition \ref{prop:form-closable} therefore shows that
$\aa$ is closable.

\smallskip\ref{it:Friedrichs2}: \
By \ref{it:Friedrichs1} the form $\aa$ is densely defined, bounded, closable, and satisfies $\aa(v,v)\ge 0$ for all $v\in\Dom(\aa)$. Its closure $\ov \aa$ is densely defined, closed, and positive, and therefore the operator $B$ associated with $\ov\aa$ is positive and selfadjoint by Proposition \ref{prop:repr-closed-symm-form}.

It remains to check that $B$ extends $A$.
Let $x\in\Dom(A)=\Dom(\aa)$. We claim that $x\in\Dom(B)$ and $Bx=Ax$.
Fix $v\in\Dom(\ov\aa)$. By definition of the closure, there exists a sequence $(v_n)_{n\ge1}\subseteq\Dom(\aa)$ such that
$v_n\to v$ in the form norm. In particular, $v_n\to v$ in $H$. Using the definition of $\ov\aa$, we obtain
$$
\ov\aa(x,v)=\lim_{n\to\infty}\aa(x,v_n)
=\lim_{n\to\infty}\iprod{Ax}{v_n}
=\iprod{Ax}{v}.
$$
Since this holds for all $v\in\Dom(\ov\aa)$, Definition~\ref{def:op-form} shows that $x\in\Dom(B)$ and $Bx=Ax$.
\end{proof}

If $A$ is a densely defined closed operator from $H$ to another Hilbert space $K$, then by
Theorem \ref{thm:AstarA} the operator $A^\star A$
with domain $$\Dom(A^\star A) = \{x\in \Dom(A):\, Ax \in \Dom(A^\star)\}$$ is positive and selfadjoint.
The next result relates this operator to the theory of forms.

\begin{proposition}\label{prop:AstarA-via-forms}
Let $A$ be a densely defined closed operator from $H$ to another Hilbert space $K$.
The form $\aa$ given by $\Dom(\aa):=\Dom(A)$ and
$$ \aa(u, v) := \iprod{Au}{Av}, \quad u,v \in \Dom(A),$$ is densely defined, closed, bounded, and accretive,
and the operator associated with $\aa$ equals the selfadjoint operator $A^\star A$.
\end{proposition}
\begin{proof} That $\aa$ is densely defined, bounded, and accretive is clear.  For $v\in \Dom(\aa)$ we have
$$ \n v\n_\aa^2 = \n v\n^2 + \aa(v,v) = \n v\n^2 + \n Av\n^2\!,$$
from which we deduce that $\n\cdot\n_\aa$ is equivalent to the graph norm of $A$. Since $A$ is closed, $\Dom(\aa)=\Dom(A)$ is  complete with respect to $\n\cdot\n_\aa$
and the closedness of $\aa$ follows.

Let $B$ be the operator associated with $\aa$. Then
$u\in\Dom(B)$ if and only if $u\in\Dom(A)$ and there exists
$f\in H$ such that
$\iprod{f}{v}=\iprod{Au}{Av}$ for all $v\in\Dom(A)$.
Taking complex conjugates, we obtain the equivalent condition
$$\iprod{Av}{Au}=\iprod{v}{f}, \quad v\in\Dom(A).$$
This, in turn, is equivalent to saying that
$Au\in\Dom(A^\star)$ and $A^\star Au=f$.

Consequently,
$\Dom(B)=\Dom(A^\star A)$
and $Bu=A^\star Au$ for all $u\in\Dom(B).$
The selfadjointness of $A^\star A$ has already been proved in Theorem \ref{thm:AstarA}.
\end{proof}

\begin{theorem}[Kato]\label{thm:Kato}
Let $\aa$ be a densely defined closed positive form in a complex Hilbert space $H$, and let $A$ be the associated positive selfadjoint operator. Then
$$
\Dom(\aa)=\Dom(A^{1/2})
$$
with equality of norms $\n u\n_\aa^2=\n u\n^2+\n A^{1/2}u\n^2$, and $\Dom(A)$ is dense in $\Dom(\aa)$ with respect to the form norm.
\end{theorem}

\begin{proof}
Since $A$ is selfadjoint and positive, the square root $A^{1/2}$ is well defined as a selfadjoint and positive operator.

The form
$$
\bb(u,v):=\iprod{A^{1/2}u}{A^{1/2}v},\qquad \Dom(\bb):=\Dom(A^{1/2}),
$$
is positive and symmetric.
Proposition \ref{prop:AstarA-via-forms}, applied to $A^{1/2}$, shows that
$\bb$ is densely defined, closed, and bounded, and that its associated
operator equals $(A^{1/2})^\star A^{1/2}$.
Since $A^{1/2}$ is selfadjoint and $(A^{1/2})^2=A$ with equality of
domains by Proposition~\ref{prop:pos-sqrt}, we have
$$
   (A^{1/2})^\star A^{1/2} = A^{1/2}A^{1/2} = A .
$$
Thus both $\aa$ and $\bb$ are closed, symmetric, nonnegative forms with the same associated operator $A$.
Proposition \ref{prop:uniq-closed-symm-form} therefore gives $\aa=\bb$.
In particular, $\Dom(\aa)=\Dom(\bb)=\Dom(A^{1/2})$. Moreover,
$$\|u\|_\aa^2=\|u\|^2+\aa(u,u)=\|u\|^2+\bb(u,u)=\|u\|^2+\|A^{1/2}u\|^2\!, \quad u\in \Dom(\aa)=\Dom(A^{1/2}).$$

Finally, density of $\Dom(A)$ in $\Dom(\aa)$ follows from Proposition \ref{prop:repr-closed-symm-form}.
\end{proof}

\section{The Laplace Operator}\label{sec:examples-forms}

We now turn to some examples that connect the theory developed in the preceding section to the boundary value problems studied in the previous chapter.

\subsection{The Laplace Operator on \texorpdfstring{$\R^d$}{R\textasciicircum d}}\label{subsec:heat-forms}\index{Laplace operator}
Consider the sesquilinear form $\aa$ with domain $\Dom(\aa) = H^1(\R^d)$ defined by
\begin{equation}\label{eq:form-aa-Delta}
 \aa(u,v):=  \int_{\R^d} \nabla u\cdot \overline{\nabla v}\ud x, \quad u,v\in H^1(\R^d).
\end{equation}
This form is bounded, positive, and closed. Indeed,
$$
\n u\n_\aa^2
=
\n u\n_2^2+\n\nabla u\n_2^2
=
\n u\n_{H^1(\R^d)}^2,
\quad u\in H^1(\R^d).
$$
We claim that the densely defined closed operator $A$ in $L^2(\R^d)$ associated with $\aa$ equals $-\Delta$,
where $\Delta$ is the weak Laplacian in $L^2(\R^d)$ with domain $\Dom(-\Delta) = H^2(\R^d)$ (cf. Theorem \ref{thm:H1-Sob}).

To prove the claim we begin by noting that if
$u\in  H^2(\R^d)$, then  $\partial_j u\in H^1(\R^d)= W^{1,2}(\R^d)$ by Theorems \ref{thm:H1-Sob} and  \ref{thm:W=H}, and therefore
\begin{equation*} \aa(u,v) = -\sum_{j=1}^d \int_{\R^d} \partial_j^2 u(x)\overline{v(x)}\ud x
= - \iprod{\Delta u}{v}
\end{equation*}
for all $v\in C_{\rm c}^\infty(\R^d)$.
By approximation this identity extends to all $v\in H^1(\R^d)$.
This means that $u\in \Dom(A)$ and $Au = -\Delta u$.

Conversely, if $u\in \Dom(A)$, then $u\in H^1(\R^d)$ and for all $v\in C_{\rm c}^\infty(\R^d)$ we have
\begin{align*}
\int_{\R^d} Au(x)v(x)\ud x
= \iprod{Au}{\ov{v}}  = \aa(u,\ov{v})  = \int_{\R^d} \nabla u\cdot \nabla v\ud x = - \int_{\R^d} u(x){\Delta v(x)}\ud x
\end{align*}
by the definition of weak derivatives.
This shows that $u$ admits a weak Laplacian given by $\Delta u = -Au$ in the sense of Theorem \ref{thm:H1-Sob},
and therefore $u\in H^2(\R^d)$ by this theorem.

Another description of the operator $\Delta$ can be given on the basis of Theorem \ref{thm:AstarA} and Proposition \ref{prop:AstarA-via-forms}.
These results identify the operator associated with the form $\aa$ defined by \eqref{eq:form-aa-Delta} as $\nabla^\star\nabla$, with domain $ \Dom(\nabla^\star\nabla) = \{f\in\Dom(\nabla):\, \nabla f\in\Dom(\nabla^\star)\}$, where $\Dom(\nabla) = H^1(\R^d)$.

Summarising this discussion, we have proved:

\begin{theorem}\label{thm:LaplacianRd} The following operators in $L^2(\R^d)$ are equal, with equal domains:
 \begin{enumerate}[label={\rm(\arabic*)}, leftmargin=*]
  \item the weak Laplacian $\Delta$ with domain
  $$ \Dom(\Delta) =  \{f\in L^2(\R^d):\, \hbox{$f$ admits a weak Laplacian in $L^2(\R^d)$}\};$$
  \item the operator $-A$, where $A$ is the operator in $L^2(\R^d)$ associated with the form $\aa$ with domain $\Dom(\aa) = H^1(\R^d)$ given by $$\aa(u,v):=  \int_{\R^d} \nabla u\cdot \overline{\nabla v}\ud x;$$
  \item the operator $-\nabla^\star\nabla$ with domain $$\Dom(\nabla^\star\nabla) = \{f\in \Dom(\nabla):\, \nabla f\in \Dom(\nabla^\star)\},$$
  where $\nabla$ is the weak gradient, viewed as a densely defined closed operator from $L^2(\R^d)$ to $L^2(\R^d\!,\K^d)$ with domain $\Dom(\nabla) = H^1(\R^d).$
 \end{enumerate}
\end{theorem}

A fourth description of $\Delta$ will be added to this list in Section \ref{subsec:heat-sgr}, namely, as the generator of the heat semigroup on $L^2(\R^d)$.

For the sake of completeness we also mention a maximal regularity result for the Poisson problem $-\Delta u = f$ on the full space $\R^d\!$.

\begin{proposition}[Maximal regularity on $\R^d$]
\label{prop:maxregH2}
Let $f,u\in L^2(\R^d)$ and suppose that
$-\Delta u=f$
in the sense of weak derivatives. Then $u\in H^2(\R^d)$.
\end{proposition}

\begin{proof}
The function $u$ has weak Laplacian $-f\in L^2(\R^d)$, and the result
follows from Theorem~\ref{thm:H1-Sob}.
\end{proof}

\subsection{The Dirichlet Laplace Operator}\label{subsec:heat-sg-domains-Dirichlet}

Let $D$ be a nonempty bounded open subset of $\R^d\!$. As before we write $$H_0^{1}(D):= W_0^{1,2}(D).$$
Consider the form $\aa_{\rm Dir}$ with domain $\Dom(\aa_{\rm Dir}) = H_0^{1}(D)$ given by
\begin{equation}\label{eq:app:aaa-Dir}
 \aa_{\rm Dir}(u,v):= \int_D \nabla u\cdot \overline{\nabla v}\ud x, \quad u,v\in H_0^{1}(D).
\end{equation}
This form is bounded, positive, coercive, and closed. Indeed,
$H_0^1(D)$ is a closed subspace of $H^1(D)$, and
$$
\n u\n_{\aa_{\rm Dir}}^2
=
\n u\n_2^2+\n\nabla u\n_2^2
=
\n u\n_{H^1(D)}^2.
$$
Coercivity follows from the Poincar\'e inequality.

The densely defined closed operator in $L^2(D)$  associated with $\aa_{\rm Dir}$  is denoted by $-\Delta_{\rm Dir}$.\index{$D$@$\Delta_{\rm Dir}$}
The operator $\Delta_{\rm Dir}$ is called the {\em Dirichlet Laplacian}\index{Dirichlet!Laplacian}
on $L^2(D)$.

By the definition of the associated operator,
$u\in\Dom(\Delta_{\rm Dir})$ if and only if $u\in H_0^1(D)$ and there
exists $f\in L^2(D)$ such that
$$
\int_D\nabla u\cdot\ov{\nabla v}\ud x
=
\int_Df\ov v\ud x,
\quad v\in H_0^1(D).
$$
Testing this identity with $v=\ov\phi$, where
$\phi\in C_{\rm c}^\infty(D)$, shows that $u$ has weak Laplacian
$\Delta u=-f$. Conversely, if $u\in H_0^1(D)$ has weak Laplacian in
$L^2(D)$, the definition of weak derivatives and the density of
$C_{\rm c}^\infty(D)$ in $H_0^1(D)$ give the preceding identity with
$f=-\Delta u$. Consequently,
\begin{align*}
\Dom(\Delta_{\rm Dir})
& =
\{u\in H_0^1(D):\,\text{ $u$ has weak Laplacian $\Delta u$ in $L^2(D)$}\},
\\
\Delta_{\rm Dir}u & =\Delta u, \quad u\in \Dom(\Delta_{\rm Dir}).
\end{align*}
Theorem \ref{thm:MR-d1} now gives the equivalent description
\begin{equation}\label{eq:dom-Dir-Laplacian}
\Dom(\Delta_{\rm Dir})
=
\{u\in H_0^1(D)\cap H_{\rm loc}^2(D):\,\Delta u\in L^2(D)\}.
\end{equation}
In particular, if $u\in C^2(\ov D)$ and $u|_{\partial D}=0$, then
Theorem~\ref{thm:H01-cont} gives $u\in H_0^1(D)$, and the preceding
characterisation shows that $u\in\Dom(\Delta_{\rm Dir})$ and
$\Delta_{\rm Dir}u=\Delta u$.

The same reasoning, together with Remark \ref{rem:MR-d1}, shows that if
$D$ has a $C^2$-boundary, this domain characterisation improves to
\begin{align*}
\Dom(\Delta_{\rm Dir})= H_0^1(D)\cap H^2(D).
\end{align*}

This identifies the $L^2$-gradient with an operator associated with a form. Indeed,
for the Dirichlet energy
$$
E(u):=\frac12\n\nabla u\n_2^2,\quad u\in H_0^1(D),
$$
we have
$$
\operatorname{grad}_{L^2(D)}E(u)=-\Delta_{\rm Dir}u,
\quad u\in\Dom(\Delta_{\rm Dir}).
$$
Thus the form domain $H_0^1(D)$ is the domain of the energy,
whereas its $L^2$-gradient exists precisely on the operator domain
$\Dom(\Delta_{\rm Dir})$. The corresponding gradient flow equation is
$$
u'(t)=\Delta_{\rm Dir}u(t),\quad t>0.
$$
In Chapter \ref{chap:semigroups} we will construct this evolution
by proving that $\Delta_{\rm Dir}$ generates an analytic
$C_0$-semigroup.

\subsection{The Neumann Laplace Operator}\label{subsec:heat-sg-domains-Neumann}

Let $D$ be a nonempty open subset of $\R^d$. We now consider the form
$\aa_{\rm Neum}$ with domain $\Dom(\aa_{\rm Neum})=H^1(D)$ given by
$$
 \aa_{\rm Neum}(u,v):= \int_D \nabla u \cdot\overline{\nabla v}\ud x,
 \quad u,v\in H^1(D).
$$
The only difference from \eqref{eq:app:aaa-Dir} is the choice of the
domain of the form.
This form is bounded, positive, and closed, since
$$
\n u\n_{\aa_{\rm Neum}}^2 = \n u\n_2^2+\n\nabla u\n_2^2 = \n u\n_{H^1(D)}^2.
$$
The densely defined
closed operator in $L^2(D)$ associated with it is denoted by
$-\Delta_{\rm Neum}$.\index{$D$@$\Delta_{\rm Neum}$}
The operator $\Delta_{\rm Neum}$ is called the {\em Neumann Laplacian}
\index{Neumann!Laplacian} on $L^2(D)$.

By the definition of the associated operator,
$u\in\Dom(\Delta_{\rm Neum})$ if and only if $u\in H^1(D)$, has weak
Laplacian $\Delta u\in L^2(D)$, and
$$
\int_D(\Delta u)\phi\ud x
=
-\int_D\nabla u\cdot\nabla\phi\ud x,
\quad \phi\in H^1(D).
$$
If $|D|<\infty$, then $\one\in H^1(D)$, and taking
$\phi=\one$ gives
$$
\int_D\Delta u\ud x=0.
$$
By Lemma \ref{lem:MR-d1}, every
$u\in\Dom(\Delta_{\rm Neum})$ belongs to $H_{\rm loc}^2(D)$.
Consequently, for every nonempty open set $D\subseteq\R^d$,
\begin{equation}\label{eq:dom-Neum-Laplacian}
\begin{aligned}
\Dom(\Delta_{\rm Neum})
 & =
\Big\{u\in  H^1(D)\cap H_{\rm loc}^2(D):\,
\Delta u\in L^2(D),\\
& \qquad\qquad \int_D(\Delta u)\phi\ud x
=
-\int_D\nabla u\cdot\nabla\phi\ud x,
\quad \phi\in H^1(D)
\Big\}.
\end{aligned}
\end{equation}
If $D$ is bounded with a $C^1$-boundary, then
$\{ u\in C^2(\ov D):\, \frac{\partial u}{\partial\nu}=0 \text{ on }\partial D\}
$ is contained in $\Dom(\Delta_{\rm Neum})$.
Indeed, if $u\in C^2(\ov D)$ satisfies
$\frac{\partial u}{\partial\nu}=0$ on $\partial D$, then Green's identity gives
the identity in \eqref{eq:dom-Neum-Laplacian} for
$\phi\in C^2(\ov D)$. By Theorem~\ref{thm:trace-W1p}, this space is dense
in $H^1(D)$. Hence $u\in\Dom(\Delta_{\rm Neum})$ and
$\Delta_{\rm Neum}u=\Delta u$.

\subsection{Selfadjointness}

The following result is an immediate consequence of Proposition
\ref{prop:repr-closed-symm-form}, noting that the forms involved are
densely defined, closed, and positive:

\begin{theorem}[Selfadjointness of the Laplacian]\label{thm:Laplace-sa} Let $\Delta$ denote the Laplacian on $L^2(\R^d)$ or the Dirichlet or Neumann Laplacian on $L^2(D)$ with $D\subseteq\R^d$ nonempty, bounded, and open. Then
$-\Delta$ is positive and selfadjoint.
\end{theorem}

These operators are understood to have the domains specified in
Sections \ref{subsec:heat-forms}--\ref{subsec:heat-sg-domains-Neumann}.
They also fall within the setting of Theorem \ref{thm:AstarA}.
Indeed, by Proposition \ref{prop:AstarA-via-forms}, each of these three
positive operators is of the form $\nabla^\star\nabla$, where $\nabla$ is the gradient
viewed as a densely defined closed operator from $L^2(U)$ to
$L^2(U;\K^d)$, with domain
$H^1(\R^d)$ in the case $U=\R^d$,
$H_0^1(D)$ in the Dirichlet case,
and $H^1(D)$ in the Neumann case.
This gives an alternative proof of their selfadjointness.

By Theorem \ref{thm:Kato}, $$\Dom((-\Delta)^{1/2}) = V,$$ where $V$ equals $H^1(\R^d)$, $H_0^1(D)$, and $H^1(D)$ respectively.

\subsection{Riesz Transforms}
\label{subsec:Riesz-transforms}

The square root of the Laplacian connects the theory of forms with
the Fourier multiplier operators studied in Chapter \ref{ch:operators}.
In dimension one this leads back to the Hilbert transform; in higher
dimensions it gives a family of operators associated with the
coordinate derivatives.

Let $d\ge1$ and let $\Delta$ be the Laplace operator with domain $\Dom(\Delta) = H^2(\R^d)$ in $L^2(\R^d)$. By Theorem \ref{thm:Kato},
$$
\Dom((-\Delta)^{1/2})=H^1(\R^d),\qquad
\n (-\Delta)^{1/2}u\n_2=\n\nabla u\n_2,\quad u\in H^1(\R^d).
$$
Under the Fourier transform, $-\Delta$ becomes multiplication by
$|\xi|^2$. The spectral theorem therefore gives
$$
\widehat{(-\Delta)^{1/2}u}(\xi)=|\xi|\wh u(\xi),\quad u\in H^1(\R^d).
$$
In particular, $(-\Delta)^{1/2}$ is injective. Since $(-\Delta)^{1/2}$ is selfadjoint, we have
$\Ran((-\Delta)^{1/2})^\perp=\Ker((-\Delta)^{1/2})=\{0\}$, so its range is dense in
$L^2(\R^d)$.

For $j=1,\ldots,d$, define
$$
R_j^{(0)}f:=-\partial_j (-\Delta)^{-1/2}f,\quad f\in\Ran((-\Delta)^{1/2}),
$$
where $(-\Delta)^{-1/2}$ is the inverse of $(-\Delta)^{1/2}$ on its range.
Since $(-\Delta)^{-1/2}f\in H^1(\R^d)$, this is well defined, and the weak derivative in
this definition belongs to $L^2(\R^d)$.
Note that $R_j^{(0)}$ maps $(-\Delta)^{1/2}u$ to $-\partial_j u$ for $u\in H^1(\R^d)$.

\begin{proposition}[Riesz transforms]
\label{prop:Riesz-transforms}
Each $R_j^{(0)}$ extends uniquely to a bounded operator $R_j$ on
$L^2(\R^d)$ with $\n R_j\n\le1$. Its Fourier multiplier is
$$
m_j(\xi):=
\begin{cases}
-i\xi_j/|\xi|, & \xi\ne0,\\
0, & \xi=0.
\end{cases}
$$
Moreover,
\begin{equation}\label{eq:Riesz-square-root}
\begin{aligned}
R_j(-\Delta)^{1/2}u&=-\partial_j u,\quad u\in H^1(\R^d),
\\
\sum_{j=1}^d\n R_jf\n_2^2&=\n f\n_2^2,
\quad f\in L^2(\R^d).
\end{aligned}
\end{equation}
For $d=1$, the operator $R_1$ is the Hilbert transform $H$.
\end{proposition}

The operators $R_1,\ldots,R_d$ are called the
{\em Riesz transforms}.\index{Riesz transform}

\begin{proof}
For $f=(-\Delta)^{1/2}u$ with $u\in H^1(\R^d)$, we have
$$
\sum_{j=1}^d\n R_j^{(0)}f\n_2^2
=\n\nabla u\n_2^2
=\n (-\Delta)^{1/2}u\n_2^2
=\n f\n_2^2.
$$
Since $\Ran((-\Delta)^{1/2})$ is dense, each $R_j^{(0)}$ has a unique bounded
extension, and the displayed identity passes to these extensions.
The identity \eqref{eq:Riesz-square-root} follows from the
definition. Furthermore, the Fourier formula for weak derivatives
gives, for $\xi\ne0$,
$$
\widehat{R_j^{(0)}f}(\xi)
=-i\xi_j\wh u(\xi)
=-i\frac{\xi_j}{|\xi|}\wh f(\xi).
$$
By continuity, this multiplier formula holds for all
$f\in L^2(\R^d)$. If $d=1$, the multiplier is
$-i\,{\rm sign}(\xi)$, so Theorem \ref{thm:HT} gives $R_1=H$.
\end{proof}

\subsection{Operators in Divergence Form}\label{subsec:divergence-form-operator}

Let $D$ be a nonempty bounded open subset of $\R^d$ and consider a matrix-valued function
$ a: D \to M_d(\K)$ satisfying the following conditions:
\begin{enumerate}[label={\rm(\roman*)}, leftmargin=*]
 \item\label{it:divergence-form1}
 the coefficients $a_{ij}:D\to \K$ are measurable and bounded;
 \item\label{it:divergence-form2}
 there exists a constant $\alpha>0$ such that, for almost all $x\in D$
and all $\xi\in\K^d$,
$$\Re\sum_{i,j=1}^d a_{ij}(x)\xi_j\ov{\xi_i}
\ge\alpha|\xi|^2\!.
$$
\end{enumerate}
Condition \ref{it:divergence-form2} is the uniform ellipticity
assumption on $a$.

Under the assumptions \ref{it:divergence-form1} and \ref{it:divergence-form2}, the form
$$
\aa_a(u,v):=\int_D a\nabla u\cdot\ov{\nabla v}\ud x
=\int_D\sum_{i,j=1}^d a_{ij}(x)\partial_j u\,\ov{\partial_i v}\ud x
$$
with domain $\Dom(\aa_a)= H^1(D)$ is a densely defined,
closed, bounded, accretive form in $L^2(D)$.
Accretivity follows from \ref{it:divergence-form2}, which implies
$\Re \aa_a(u,u) \ge \alpha
\n\nabla u\n_2^2$.
By boundedness of the coefficients there exists a constant $C\ge0$ such that
$|\aa_a(u,v)|\le C\n\nabla u\n_2\n\nabla v\n_2$ for all $u,v\in H^1(D)$.
Taking $v=u$, we obtain
$\alpha\n\nabla u\n_2^2 \le \Re\aa_a(u,u)
\le |\aa_a(u,u)| \le C\n\nabla u\n_2^2$.
Consequently,
$$
\n u\n_2^2+\alpha\n\nabla u\n_2^2
\le \n u\n_{\aa_a}^2
\le \n u\n_2^2+C\n\nabla u\n_2^2,
$$
so the form norm $\n\cdot\n_{\aa_a}$ is equivalent to the $H^1(D)$-norm.
Hence $\aa_a$ is closed. It also follows that
$|\aa_a(u,v)| \le \frac C\alpha\n u\n_{\aa_a}\n v\n_{\aa_a}$,
so $\aa_a$ is bounded.

The operator in $L^2(D)$ associated with $\aa_{a}$ is usually denoted by $$- {\rm div}(a \nabla)$$
in recognition of the fact that (at least formally) the Hilbert space adjoint of $\nabla$ equals $-{\rm div}$.
This operator is often referred to as a second-order differential operator in {\em divergence form}\index{divergence form operator}
subject to Neumann boundary conditions.
If the coefficients satisfy the symmetry condition
$$a_{ij}(x)=\ov{a_{ji}}(x) \ \hbox{for almost all} \  x\in D,$$
then $\aa_a$ is symmetric and the associated operator $-{\rm div}(a\nabla)$ is positive and selfadjoint.

The analogue for Dirichlet boundary conditions is obtained by replacing the form domain by $H_0^1(D)$. This results in a densely defined,
closed, bounded, coercive form in $L^2(D)$; coercivity follows from Poincar\'e's inequality.

If, in addition, $D$ is connected and has $C^1$ boundary, there is a
coercive analogue for Neumann boundary conditions. Put
$$
L^2_{\rm av}(D)
:=\Bigl\{
f\in L^2(D):\,
\int_D f\ud x=0
\Bigr\}.
$$
Then, viewed as a form in $L^2_{\rm av}(D)$ with domain
$H^1_{\rm av}(D)$, the form $\aa_a$ is densely defined, closed, bounded,
and coercive. Indeed, coercivity follows from the Poincar\'e inequality
on $H^1_{\rm av}(D)$.

\section{The Poisson Problem Revisited}\label{sec:Poisson-revisited}

We now revisit the
Poisson problem $-\Delta u = f$ by viewing it as a special instance of the abstract problem
$$ Au = x,$$
where $A$ is assumed to be a closed operator acting in a Banach space $X$, $x\in X$ is a given element, and $u\in X$
is the unknown. One could define a {\em strong solution} as an element $u\in \Dom(A)$ such that $Au = x$, but this is not what we did in Section \ref{sec:Poisson-problem}. Instead, we considered {\em weak solutions} defined in terms of the
sesquilinear form with which $A$ is associated. If $A$ is also densely defined, a third option is to use duality to define a {\em scalar solution} to be an element $u\in X$ with the property that
$$ \lb u, A^* x^*\rb = \lb x,x^*\rb, \quad x^*\in \Dom(A^*).$$
In line with standard functional-analytic terminology it would be more appropriate to call {\em this} a weak solution, but the usage of the term `weak solution' in connection with integration by parts using test functions is well established.

\begin{proposition}\label{prop:weakvsstrongsol}
 Let $A$ be a densely defined closed linear operator acting in a Banach space $X$ and let $x\in X$ be given. For any $u\in X$ the following assertions are equivalent:
 \begin{enumerate}[label={\rm(\arabic*)}, leftmargin=*]
  \item\label{it:weakvsstrongsol1} $u$ is a strong solution of $Au = x$, that is, $u\in \Dom(A)$ and $Au = x$;
  \item\label{it:weakvsstrongsol2} $u$ is a scalar solution of $Au = x$, that is, $\lb u, A^* x^*\rb = \lb x,x^*\rb$ for all $x^*\in \Dom(A^*)$.
 \end{enumerate}
 If $X = H$ is a Hilbert space and $A$ is the operator associated with a densely defined form
 $\aa$ in $H$, then \ref{it:weakvsstrongsol1} and \ref{it:weakvsstrongsol2} are equivalent to:
 \begin{enumerate}[label={\rm(\arabic*)}, leftmargin=*]\setcounter{enumi}{2}
  \item\label{it:weakvsstrongsol3} $u$ is a weak solution of $Au = x$, that is,
  $u\in \Dom(\aa)$ and $\aa(u,v) = \iprod{x}{v}$ for all $v\in\Dom(\aa)$.
  \end{enumerate}
\end{proposition}

\begin{proof}
 The implications \ref{it:weakvsstrongsol1}$\Rightarrow$\ref{it:weakvsstrongsol2} and \ref{it:weakvsstrongsol1}$\Rightarrow$\ref{it:weakvsstrongsol3}
 are trivial. The implication \ref{it:weakvsstrongsol2}$\Rightarrow$\ref{it:weakvsstrongsol1} is an immediate consequence of Proposition \ref{prop:DomAweak}.
 Finally, if \ref{it:weakvsstrongsol3} holds, then by the definition of the associated operator we have
 $u\in \Dom(A)$ and $Au = x$, so \ref{it:weakvsstrongsol1} holds.
\end{proof}

In the special case where $A=-\Delta$, with $\Delta$ the Dirichlet or
Neumann Laplacian, Proposition~\ref{prop:weakvsstrongsol} says that every
weak solution of the Poisson problem with right-hand side in $L^2(D)$ is
a strong solution for the corresponding Laplace operator. The domain
descriptions \eqref{eq:dom-Dir-Laplacian} and
\eqref{eq:dom-Neum-Laplacian} express the regularity results of
Theorems~\ref{thm:MR-d1} and \ref{thm:MR-d1-N} in operator-theoretic language.

\section{Weyl's Theorem}\label{sec:Weyl}

This section is a digression from the main line of development and is
devoted to a proof of Weyl's celebrated asymptotic formula for the number of eigenvalues of the Dirichlet Laplacian.

\subsection{Spectrum of the Dirichlet and Neumann Laplacians}\label{sec:spec-Dir-Neuum}

By Theorem \ref{thm:Laplace-sa} and Proposition
\ref{prop:sa-pos-spectrum}, $-\Delta_{\rm Dir}$ and $-\Delta_{\rm Neum}$ are positive
selfadjoint operators with spectra contained in $[0,\infty)$. As a warm-up, let us compute the spectra
for $d=1$ and $D = (0,1)$:

\begin{example}\label{ex:ev-Dir}
We will use the fact that if $u\in C^2[0,1]$ satisfies
$$
u(0)=u(1)=0,
$$
then $u\in \Dom(\Delta_{\rm Dir})$ and $\Delta_{\rm Dir}u=u''$, whereas if
$u\in C^2[0,1]$ satisfies
$$
u'(0)=u'(1)=0,
$$
then $u\in \Dom(\Delta_{\rm Neum})$ and $\Delta_{\rm Neum}u=u''$.
We will also use the fact that a function $u\in C^2[0,1]$ belongs to $H_0^1(0,1)$ if
and only if $u(0)=u(1)=0$; we leave the elementary proof to the reader (see
Problem \ref{prob:cont-version-W1p} for a more precise result).

The functions $u_n(\theta) = \sin(\pi n \theta)$, $n\ge 1$, satisfy $-\Delta_{\rm Dir}u_n = -u_n'' = \pi^2 n^2 u_n$ and obey Dirichlet
boundary conditions.
Moreover, by Theorem \ref{thm:L2bases-sin-cos}, after normalisation these functions form an
orthonormal basis for $L^2(0,1)$.
By Proposition \ref{prop:sigmaAdiag},
this implies
$$\si(-\Delta_{\rm Dir}) = \bigl\{\pi^2 n^2:\, n\ge 1\bigr\},$$
with eigenfunctions $u_n(\theta) = \sin(\pi n \theta)$.
Likewise, the functions $v_n(\theta) = \cos(\pi n \theta)$, $n\ge 0$, satisfy
$-\Delta_{\rm Neum}v_n = -v_n'' = \pi^2 n^2 v_n$
and obey Neumann boundary conditions. Again by Theorem \ref{thm:L2bases-sin-cos}, after normalisation
they form an orthonormal basis for $L^2(0,1)$. Hence
$$\si(-\Delta_{\rm Neum}) = \bigl\{\pi^2 n^2:\, n =0,1,2, \dots\bigr\}.$$
\end{example}

Turning to higher dimensions, we begin with a simple observation.

\begin{proposition}\label{prop:DeltaFredholm} Let $D$ be a nonempty bounded open subset of $\R^d\!$.
\begin{enumerate}[label={\rm(\arabic*)}, leftmargin=*]
 \item\label{it:DeltaFredholm1} $\Delta_{\rm Dir}$ is both injective and surjective, and hence invertible;
 \item\label{it:DeltaFredholm2} if, in addition, $D$ is connected and has $C^1$-boundary, the null space of
 $\Delta_{\rm Neum}$ consists of the constant functions and its range is the orthogonal complement of the constant functions. In particular, its range is closed and
 $$\dim \Ker(\Delta_{\rm Neum}) = \codim \Ran(\Delta_{\rm Neum}) = 1.$$
\end{enumerate}
\end{proposition}
\begin{proof}
\ref{it:DeltaFredholm1}: \
The existence and uniqueness results for the Dirichlet Poisson problem, together with
Proposition~\ref{prop:weakvsstrongsol}, show that
$\Delta_{\rm Dir}$ is bijective. Since it is closed, the closed graph
theorem shows that its inverse is bounded. Thus
$0\in\varrho(\Delta_{\rm Dir})$.

\smallskip
\ref{it:DeltaFredholm2}: \
The existence and uniqueness results for the Neumann Poisson
problem show that
$\Ker(\Delta_{\rm Neum})=\operatorname{span}\{\one\}$ and
$\Ran(\Delta_{\rm Neum})=\{\one\}^\perp$, from which the assertions follow.
\end{proof}

Extending the corresponding definition for bounded operators, we define a {\em  Fredholm operator}\index{Fredholm!operator}
as a closed operator whose null space is finite-dimensional and whose range has finite codimension.
The same proof as in the bounded case shows that the second condition
implies that the range is closed. The {\em index}\index{index} of such an operator $A$ is defined as
$$ \ind(A):= \dim \Ker(A) - \codim \Ran(A).$$
Proposition \ref{prop:DeltaFredholm} implies that both $\Delta_{\rm Dir}$ and $\Delta_{\rm Neum}$ (the latter under the stated more restrictive assumptions on $D$) are Fredholm operators with index $0$.

\begin{theorem}\label{thm:spectra-Laplacian}
Let $D$ be a nonempty bounded open subset of $\R^d$.

\begin{enumerate}[label={\rm(\arabic*)}, leftmargin=*]
\item There exist an orthonormal basis $(h_n)_{n\ge1}$ of $L^2(D)$ and
numbers
$$
0<\lambda_1\le\lambda_2\le\cdots,\qquad \lambda_n\to\infty,
$$
such that
$-\Delta_{\rm Dir}h_n=\lambda_nh_n$ for all $n\ge 1$.
 In particular, $\sigma(-\Delta_{\rm Dir}) = \{\la_n:\, n\ge 1\}$.

\item If, in addition, $D$ is connected and has $C^1$-boundary, there
exist an orthonormal basis $(k_n)_{n\ge1}$ of $L^2(D)$ and numbers
$$
0=\mu_1\le\mu_2\le\cdots,\qquad \mu_n\to\infty,
$$
such that
$-\Delta_{\rm Neum}k_n=\mu_nk_n$ for all $n\ge1$, and
the eigenvalue $0$ has multiplicity one. In particular,
$\sigma(-\Delta_{\rm Neum}) = \{\mu_n:\, n\ge 1\}$.
\end{enumerate}

In both cases the eigenvalues are repeated according to multiplicity.
\end{theorem}

\begin{proof} By Proposition \ref{prop:sigmaAdiag}, the formulas for $\sigma(-\Delta_{\rm Dir})$ and $\sigma(-\Delta_{\rm Neum})$ follow from the assertions preceding them.

Let $A=-\Delta_{\rm Dir}$ and $V=H_0^1(D)$, or let
$A=-\Delta_{\rm Neum}$ and $V=H^1(D)$. In the latter case we impose the assumptions stated in the theorem. In either case, since $A$ is positive and selfadjoint, we have
$-1\in\varrho(A)$.

Put $B:=(I+A)^{-1}$. If $u=Bf$, then
$$
\n u\n_V^2
=
\n u\n_2^2+\aa(u,u)
=
\iprod{f}{u}
\le\n f\n_2\n u\n_2
\le\n f\n_2\n u\n_V.
$$
Thus $B$ is bounded from $L^2(D)$ into $V$. The inclusion of $V$ into
$L^2(D)$ is compact by the Rellich--Kondrachov theorem, and hence $B$ is
compact on $L^2(D)$.

The operator $B$ is positive, selfadjoint, compact, and injective.
The spectral theorem for compact selfadjoint operators therefore gives
an orthonormal basis $(h_n)_{n\ge1}$ and eigenvalues, repeated according
to multiplicity and arranged so that
$$
1\ge\beta_1\ge\beta_2\ge\cdots>0,
\qquad
\beta_n\to0,
$$
such that
$Bh_n=\beta_nh_n$ for every $n\ge 1$.
Since $\Ran(B)=\Dom(A)$, we have $h_n\in\Dom(A)$ and
$$
Ah_n=(\beta_n^{-1}-1)h_n.
$$
Putting
$\nu_n:=\beta_n^{-1}-1$, we have
$$
0\le\nu_1\le\nu_2\le\cdots,
\qquad
\nu_n\to\infty.
$$
In the Dirichlet case Proposition \ref{prop:DeltaFredholm} gives
$\nu_1>0$. In the Neumann case it gives $\nu_1=0$, with multiplicity
one. Renaming the eigenvectors and eigenvalues accordingly proves the
two assertions.
\end{proof}

As a variant on the min-max theorem for compact positive Hilbert space operators (Theorem \ref{thm:minmax}),
we prove an explicit formula for the Dirichlet and Neumann eigenvalues of the Laplace operator on a nonempty bounded open set $D\subseteq \R^d$; in the case of Neumann boundary conditions we make the additional assumption that
$D$ is connected and has $C^1$-boundary. We denote by
$0<\la_1\le \la_2 \le \dots$ and $0 = \mu_1 \le \mu_2\le  \dots$ the sequences of eigenvalues of
$-\Delta_{\rm Dir}$ and $-\Delta_{\rm Neum}$, respectively, taking multiplicities into account.

\begin{theorem}[Courant--Fischer]\label{thm:CF}\index{theorem!min-max}\index{theorem!Courant--Fischer}
With the notation just introduced,
\begin{enumerate}[label={\rm(\arabic*)}, leftmargin=*]
 \item
for all $n\ge 1$ we have
\begin{align*}
\la_n & = \inf_{\substack{Y\subseteq H_0^1(D) \\ \dim(Y) = n}} \ \sup_{\substack{y \in Y\\y\not=0}} \ \frac{\n \nabla y\n_{L^2(D)}^2}{\n y\n_{L^2(D)}^2},
\end{align*}
where the infima are taken over all subspaces $Y$ of dimension $n$;
\item if, in addition, $D$ is connected and has a $C^1$-boundary, then for all $n\ge 1$ we have
\begin{align*}
\mu_n & = \inf_{\substack{Y\subseteq H^1(D) \\ \dim(Y) = n}} \ \sup_{\substack{y \in Y\\ y\not=0}} \ \frac{\n \nabla y\n_{L^2(D)}^2}{\n y\n_{L^2(D)}^2},
\end{align*}
where the infima are taken over all subspaces $Y$ of dimension $n$.
\end{enumerate}
\end{theorem}

\begin{proof}
We present the case of Dirichlet eigenvalues, the proof for Neumann eigenvalues being identical. Put
$A:=-\Delta_{\rm Dir}$, and let $(h_j)_{j\ge1}$ be the orthonormal basis
of eigenvectors from Theorem~\ref{thm:spectra-Laplacian}, with
$Ah_j=\lambda_jh_j$.

By Kato's theorem (Theorem~\ref{thm:Kato}),
$$
H_0^1(D)=\Dom(A^{1/2}), \qquad
\n\nabla f\n_2^2=\n A^{1/2}f\n_2^2.
$$
Thus, if $f=\sum_{j\ge1}c_jh_j$, then
$f\in H_0^1(D)$ if and only if $\sum_{j\ge1}\lambda_j|c_j|^2<\infty$,
and in this case
$$
\n\nabla f\n_2^2=\sum_{j\ge1}\lambda_j|c_j|^2\!.
$$

Let $Y\subseteq H_0^1(D)$ have dimension $n$. If
$f\in Y$ is orthogonal to $h_1,\ldots,h_{n-1}$, then
$$
\n\nabla f\n_2^2
=
\sum_{j\ge n}\lambda_j|c_j|^2
\ge
\lambda_n\n f\n_2^2.
$$
It follows that
$$
\sup_{\substack{f\in Y\\f\ne0}}
\frac{\n\nabla f\n_2^2}{\n f\n_2^2}
\ge\lambda_n.
$$

Conversely, for $Y_n:=\operatorname{span}\{h_1,\ldots,h_n\}$ we have
$$
\n\nabla f\n_2^2
=
\sum_{j=1}^n\lambda_j|c_j|^2
\le
\lambda_n\n f\n_2^2,
\quad f\in Y_n.
$$
Combining the two inequalities proves the formula.
\end{proof}

\begin{corollary}\label{cor:DirNeumEV} If $D$ is connected
with $C^1$-boundary, then $$\mu_n \le \la_n, \quad n\ge 1.$$
\end{corollary}
\begin{proof}
Since $H_0^1(D)\subseteq H^1(D)$, the infimum in the Neumann
Courant--Fischer formula is taken over a larger collection of
$n$-dimensional subspaces.
\end{proof}

\begin{corollary}[Optimal Poincar\'e constants]
\label{cor:optimal-Poincare}
Let $D$ be a nonempty bounded open subset of $\R^d$.

\begin{enumerate}[label={\rm(\arabic*)}, leftmargin=*]
\item If $\lambda_1$ is the first Dirichlet eigenvalue, then
$$
\lambda_1
=
\inf_{\substack{u\in H_0^1(D)\\u\ne0}}
\frac{\n\nabla u\n_2^2}{\n u\n_2^2}.
$$
Consequently, the optimal constant in the Poincar\'e inequality
$$
\n u\n_2\le C\n\nabla u\n_2,
\quad u\in H_0^1(D),
$$
is $C=\lambda_1^{-1/2}$.

\item Suppose, in addition, that $D$ is connected and has
$C^1$-boundary. If $\mu_2$ is the first nonzero Neumann eigenvalue, then
$$
\mu_2
=
\inf_{\substack{u\in H^1(D),\ u\ne0\\
                 \int_Du\ud x=0}}
\frac{\n\nabla u\n_2^2}{\n u\n_2^2}.
$$
Consequently, the optimal constant in the Poincar\'e--Wirtinger
inequality
$$
\n u-u_D\n_2\le C\n\nabla u\n_2,
\quad u\in H^1(D),
$$
where $u_D:=\frac1{|D|}\int_Du\ud x$,
is $C=\mu_2^{-1/2}$.
\end{enumerate}
\end{corollary}

\begin{proof}
The first identity is the case $n=1$ of the
Courant--Fischer theorem for the Dirichlet Laplacian. Equality is attained by an eigenfunction
corresponding to $\lambda_1$, and the assertion about the optimal
constant follows.

For the Neumann problem, let $(h_n)_{n\ge1}$ be an orthonormal basis of
eigenfunctions, with $h_1=|D|^{-1/2}\one$ and ordered by their eigenvalues $0=\mu_1<\mu_2\le\cdots$, repeated according to multiplicity. If $u\in H^1(D)$ has integral zero, then
$\iprod{u}{h_1}=0$. Writing
$u=\sum_{n\ge2}c_nh_n$,
we obtain
$$
\n\nabla u\n_2^2
=
\sum_{n\ge2}\mu_n|c_n|^2
\ge
\mu_2\sum_{n\ge2}|c_n|^2
=
\mu_2\n u\n_2^2.
$$
Equality is attained by an eigenfunction corresponding to $\mu_2$.
Applying this result to $u-u_D$ proves the final assertion.
\end{proof}

\subsection{Weyl's Theorem}

Weyl's theorem gives an asymptotic formula for the number
of Dirichlet eigenvalues in $[0,r]$ as $r\to\infty$.

\begin{theorem}[Weyl]\label{thm:Weyl}\index{theorem!Weyl}
Let $D$ be a nonempty bounded open subset of $\R^d$ satisfying $|\partial D| = 0$, let $0<\la_1\le \la_2\le \dots$ be the sequence of eigenvalues of $-\Delta_{\rm Dir}$ on $L^2(D)$, taking multiplicities into account, and for $r>0$ let
$$N_{\rm Dir}(r):= \max\bigl\{n\ge 1:\, \la_n\le r\bigr\},$$
with the convention $\max\emptyset= 0$. Then
$$  \lim_{r\to\infty} \frac{N_{\rm Dir}(r)}{r^{d/2}} =\frac{\om_d}{(2\pi)^d}|D| ,$$
where $\omega_d = \pi^{d/2}/\Gamma(1+ \frac12d) $ is the volume of the unit ball in $\R^d\!$.
\end{theorem}

For $d=2$, the condition $|\partial D|=0$ is satisfied if the boundary
is a rectifiable curve. More generally, the condition holds for domains with
Lipschitz boundary.

Before turning to the proof, it is instructive to revisit Example \ref{ex:ev-Dir}.
For the Dirichlet Laplacian in $L^2(0,1)$ we obtain
$$N_{\rm Dir}(r) = \max\bigl\{n\ge 1:\, \pi^2 n^2 \le r\bigr\}.$$
On the other hand, $\omega_1 = |(-1,1)| = 2$ and $|D| = |(0,1)| =1$. It follows that
$$
\lim_{r\to\infty} \frac{N_{\rm Dir}(r)}{r^{1/2}} = \frac1\pi, \qquad \frac{\om_1}{2\pi}|D| = \frac{2}{2\pi}\cdot 1 = \frac1\pi.
$$

The main lemma needed for the proof of Weyl's theorem is a monotonicity result.

\begin{lemma}\label{lem:Weyl-monotone2}
 Let $D_1$ and $D_2$ be nonempty bounded open subsets of $\R^d$ with $D_1\subseteq D_2$.
 Then the corresponding Dirichlet eigenvalues, taking multiplicities into account, satisfy
$$ \la_{n,D_1} \ge \la_{n,D_2}, \quad n\ge 1.$$
As a consequence,
 $N_{D_1}(r) \le N_{D_2}(r)$ for all $r>0$.
\end{lemma}
\begin{proof}
This follows from the Courant--Fischer theorem, because the zero extensions of functions in $H_0^1(D_1)$ belong to $H_0^1(D_2)$.
\end{proof}

The analogue of this lemma fails for Neumann boundary conditions. It is for this reason that we only present Weyl's theorem for Dirichlet eigenvalues. The case of Neumann boundary conditions is discussed in the Notes to this chapter.

We also need a lemma on counting lattice points in dilations of a given set:

\begin{lemma}[Counting lattice points]\label{lem:lattice-dilate}
Let $E\subseteq\R^d$ be bounded and measurable, and assume that $|\partial E|=0$.
Then, as $t\to\infty$,
$$
\#(tE\cap\Z^d)=t^d|E|+o(t^d).
$$
\end{lemma}
\begin{proof}
For $t>0$ and $m\in\Z^d$ put
$Q_{m,t}:=t^{-1}(m+[0,1)^d)$
and define
$$
F_t:=\sum_{m\in\Z^d}\one_E(m/t)\one_{Q_{m,t}}.
$$
Then
$$
\int_{\R^d}F_t(x)\ud x
=
t^{-d}\#(tE\cap\Z^d).
$$

If $x\notin\partial E$, then $x$ belongs either to the interior of $E$
or to the interior of its complement. The distance between $x$ and the
point $m/t$ belonging to the same cube $Q_{m,t}$ is at most
$\sqrt d/t$. Consequently,
$F_t(x)\to\one_E(x)$
for every $x\notin\partial E$.

Since $E$ is bounded, the functions $F_t$, for $t\ge1$, are supported in
one fixed bounded set and satisfy $0\le F_t\le1$. As
$|\partial E|=0$, dominated convergence gives
$$
t^{-d}\#(tE\cap\Z^d)
=
\int_{\R^d}F_t(x)\ud x
\to
\int_{\R^d}\one_E\ud x
=
|E|
$$
as $t\to\infty$.
\end{proof}

\begin{proof}[Proof of Theorem \ref{thm:Weyl}]
For an open subset $U$ of $\R^d$ we denote the Dirichlet Laplacian in $L^2(U)$ by $\Delta_U$.

\smallskip
{\em Step 1} --
We first prove the theorem for open cubes with sides parallel
to the coordinate axes. By translation, we may assume that
$D=(0,\ell)^d$ for some $\ell>0$. For
$n=(n_1,\ldots,n_d)\in\N_1^d$, where $\N_1:=\{1,2,\ldots\}$,
put
$$
u_n(x):=\Bigl(\frac2\ell\Bigr)^{d/2}
\prod_{j=1}^d\sin\Bigl(\frac{\pi n_jx_j}{\ell}\Bigr),
\quad x\in D.
$$
Each $u_n$ belongs to $C^\infty(\ov D)$ and vanishes on
$\partial D$, so $u_n\in\Dom(\Delta_D)$. Direct differentiation
gives
$$
-\Delta_Du_n=\frac{\pi^2}{\ell^2}|n|^2u_n.
$$
By rescaling Example \ref{ex:ev-Dir} and using
Section \ref{subsec:tensorbases}, the functions $u_n$ form an
orthonormal basis of $L^2(D)$.
Put $\lambda_n:=\pi^2\sum_{j=1}^d n_j^2/\ell^2$.
For every $R>0$, only finitely many $n\in\N_1^d$ satisfy
$\lambda_n\le R$. If $-\Delta_Dh=\lambda h$,
selfadjointness gives
$(\lambda-\lambda_n)\iprod{h}{u_n}=0$ for all $n\in\N_1^d$.
Completeness of the orthogonal family $(u_n)$ therefore implies
$\Ker(-\Delta_D-\lambda)=\operatorname{span}\{u_n:\,\lambda_n=\lambda\}$.
Thus the values $\lambda_n$ exhaust the eigenvalues,
each repeated according to its multiplicity.
Consequently,
$$
N_{\rm Dir}(r)
=\#\Bigl\{n\in\N_1^d:\,|n|\le\frac{\ell\sqrt r}{\pi}\Bigr\}.
$$

Let
$$
E:=\{x\in(0,\infty)^d:\,|x|\le1\}.
$$
Its boundary is contained in the union of the unit sphere
and the coordinate hyperplanes, so $|\partial E|=0$, and
$|E|=2^{-d}\omega_d$. Applying Lemma \ref{lem:lattice-dilate}
with $t=\ell\sqrt r/\pi$, we obtain
$$
N_{\rm Dir}(r)
=\frac{\omega_d}{2^d}
 \Bigl(\frac{\ell\sqrt r}{\pi}\Bigr)^d+o(r^{d/2})
=c_d|D|r^{d/2}+o(r^{d/2}),
$$
where $c_d:=\omega_d/(2\pi)^d$. This proves the theorem
for cubes.

\smallskip
{\em Step 2} -- In this step we prove that if $U_1$ and $U_2$ are bounded disjoint open subsets of $\R^d$ and
$U=U_1\cup U_2$, then for all $r>0$ we have
$N_U(r)=N_{U_1}(r)+N_{U_2}(r).$

For $v\in H_0^1(U_j)$, let $\widetilde v$ denote its extension by zero to $U$.
Since $U_1$ and $U_2$ are disjoint, the mapping from $H_0^1(U_1)\oplus H_0^1(U_2)$ to $H_0^1(U)$ defined by
$(v_1,v_2)\mapsto \widetilde v_1+\widetilde v_2$
is an isometric isomorphism, with inverse given by restriction to $U_1$ and
$U_2$. Likewise we have an isometric isomorphism $
L^2(U)=L^2(U_1)\oplus L^2(U_2)$.
Under these identifications, the Dirichlet form on $U$
is the direct sum of the Dirichlet forms on $U_1$ and $U_2$.
By Definition \ref{def:op-form}, testing with $(v_1,0)$
and $(0,v_2)$, where $v_j\in H_0^1(U_j)$, gives
$$
-\Delta_U=(-\Delta_{U_1})\oplus(-\Delta_{U_2}).
$$
Each eigenspace is therefore the direct sum of the
corresponding eigenspaces of the two summands. Thus their
eigenvalue multiplicities add, and
$$
N_U(r)=N_{U_1}(r)+N_{U_2}(r),
\quad r>0.
$$

{\em Step 3} --
By Steps 1 and 2, the theorem holds for finite disjoint unions
of open cubes with sides parallel to the coordinate axes.

We first prove the lower bound. Let $R_1,\ldots,R_k$ be
disjoint such cubes contained in $D$, and put
$R:=\bigcup_{j=1}^kR_j$. By Lemma \ref{lem:Weyl-monotone2}
and Steps 1 and 2,
$$
c_d |R| = \lim_{r\to\infty}\frac{N_R(r)}{r^{d/2}}
\le \liminf_{r\to\infty}\frac{N_{\rm Dir}(r)}{r^{d/2}}.
$$
Since $D$ is open, it can be approximated from within,
up to a set of arbitrarily small measure, by finite disjoint
unions of such cubes. Hence,
$$
c_d |D| \le \liminf_{r\to\infty}\frac{N_{\rm Dir}(r)}{r^{d/2}} .
$$

For the upper bound, choose an open cube $Q$ with sides
parallel to the coordinate axes and $\ov D\subseteq Q$,
and put $V:=Q\setminus\ov D$. Since $V$ is open, the lower
bound just proved gives
$$
c_d |V| \le \liminf_{r\to\infty}\frac{N_V(r)}{r^{d/2}}.
$$
Moreover, $D$ and $V$ are disjoint subsets of $Q$, and therefore,
by Step 2 and Lemma \ref{lem:Weyl-monotone2},
$$
N_{\rm Dir}(r)+N_V(r) = N_{D\cup V}(r) \le N_Q(r),\quad r>0 .
$$
It follows that
\begin{align*}
 \limsup_{r\to\infty}\frac{N_{\rm Dir}(r)}{r^{d/2}}
 \le  \lim_{r\to\infty}\frac{N_Q(r)}{r^{d/2}}
 - \liminf_{r\to\infty}\frac{N_V(r)}{r^{d/2}}
 \le c_d |Q| - c_d |V| = c_d |D|,
\end{align*}
where in the last equality we used $|\partial D|=0$. Together with the
lower bound this gives
\[
\lim_{r\to\infty}\frac{N_{\rm Dir}(r)}{r^{d/2}} = c_d |D|.
\qedhere\]
\end{proof}

If one imagines a bounded open set $D$ in $\R^d$ as a $d$-dimensional `drum', the eigenvalues of the negative Dirichlet Laplacian
on $L^2(D)$ can be interpreted as the `frequencies' of the drum. This prompted the famous question of Mark Kac:
``Can one hear the shape of a drum?'' In its mathematical formulation, the question is whether the shape of $D$, up to an isometry of $\R^d\!$, is determined by
its sequence of frequencies. Without further assumptions on $D$, the answer is negative in general. Nevertheless, Weyl's theorem implies that one can hear the $d$-dimensional volume of the drum.

\section{The Landscape Function}\label{sec:landscape}

Weyl's theorem describes the distribution of eigenvalues at high
energies. We now consider a different way of obtaining spectral
information, based on the solution of a single Poisson problem. Besides
giving bounds for the first eigenvalue, this solution controls the
distribution of the mass of eigenfunctions and provides an upper bound
for the number of eigenvalues below a prescribed energy.

Throughout this section, $D\subseteq\R^d$ is a nonempty bounded open set. By Proposition
\ref{prop:DeltaFredholm}, the function
$u:=-\Delta_{\rm Dir}^{-1}\one $
is well defined and belongs to $\Dom(\Delta_{\rm Dir})\subseteq H_0^1(D)$. It is the
weak solution of the Poisson problem $\Delta u = -\one$ with Dirichlet boundary conditions, so
\begin{align}\label{eq:landscape-weak}
\int_D\nabla u\cdot\ov{\nabla \phi}\ud x
=\int_D\ov \phi\ud x,
\quad \phi\in H_0^1(D).
\end{align}
In connection with eigenvalue problems it is called the
\emph{landscape function}.\index{landscape function}

\begin{lemma}\label{lem:landscape-positive}
The operator $-\Delta_{\rm Dir}^{-1}$ maps nonnegative functions to nonnegative
functions. Moreover, $u\in L^\infty(D)$ and $u\ge0$.
\end{lemma}

Here, and in the remainder of this section,
all pointwise inequalities are understood in the almost everywhere sense.

\begin{figure}[h]
 \centering
 \includegraphics[scale=0.45]{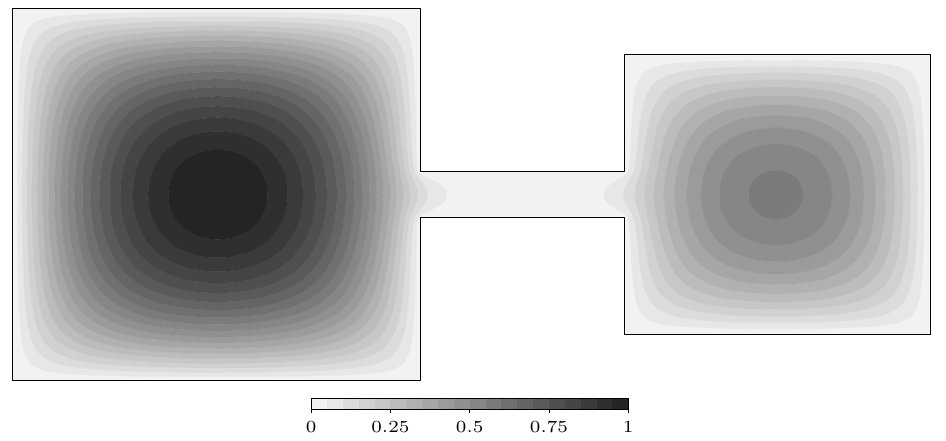}
 \caption{The $L^\infty$-normalised landscape function of an asymmetric dumbbell domain}
 \label{fig:landscape-dumbbell}
 \end{figure}

\begin{proof} Let $g\in L^2(D)$ be nonnegative and put
$w:=-\Delta_{\rm Dir}^{-1}g$. Since $w$ is a weak solution of $-\Delta w = g$ with Dirichlet boundary conditions, for all $\phi\in H_0^1(D)$, we have
\begin{align*}
\int_D\nabla\ov w\cdot{\nabla\phi}\ud x
=\ov{\int_D\nabla w\cdot\nabla\ov\phi\ud x}
=\ov{\int_Dg\ov\phi\ud x}
=\int_Dg\phi\ud x,
\end{align*}
where we used that $g$ is real-valued. Thus $\ov w$ is another
weak solution. Uniqueness of weak solutions gives
$\ov w=w$, so $w$ is real-valued.

By Theorem
\ref{thm:H1-lattice}, the negative part $w^-$ belongs to $H_0^1(D)$.
Testing the Poisson equation for $w$ with $\phi = w^-$ gives
$$
-\int_D|\nabla w^-|^2\ud x
=\int_Dg w^-\ud x\ge0.
$$
Hence $\nabla w^- = 0$ almost everywhere. Poincar\'e's inequality implies that $w^-=0$, and hence $w\ge0$.
Applying this with $g = \one$, in particular we obtain $u\ge0$.

Choose $R>0$ such that $\ov D\subseteq B(0;R)$ and put
$$b(x):=(R^2-|x|^2)/(2d).$$ Then $b>0$ on $\ov D$ and $-\Delta b=1$.
The function $(u-b)^+$ belongs to $H_0^1(D)$. Indeed, if the real-valued functions
$u_n\in C_{\rm c}^\infty(D)$ converge to $u$ in $H^1(D)$, then
$(u_n-b)^+\to(u-b)^+$ in $H^1(D)$ by Theorem \ref{thm:H1-lattice}.
Each approximating function has compact support in $D$ and therefore
belongs to $H_0^1(D)$ by Proposition \ref{prop:local-approx}.
Since $b$ is smooth and $-\Delta b=1$, integration by parts gives
$$
\int_D\nabla b\cdot\ov{\nabla\phi}\ud x
=-\int_D\Delta b\,\ov\phi\ud x
=\int_D\ov\phi\ud x,
\quad \phi\in C_{\rm c}^\infty(D).
$$
By density, this identity holds for every $\phi\in H_0^1(D)$.
Thus $b$ satisfies the same weak equation as $u$, even though
$b$ does not satisfy the zero boundary condition.
Subtracting this identity from \eqref{eq:landscape-weak} and
taking $\phi=(u-b)^+$ gives
$$
\int_D\nabla(u-b)\cdot\nabla ((u-b)^+)\ud x=0.
$$
On the other hand, the differentiation rule for positive parts (Theorem \ref{thm:H1-lattice})
gives $\nabla ((u-b)^+)=\one_{\{u>b\}}\nabla(u-b)$
almost everywhere. Therefore,
$$
\int_D|\nabla ((u-b)^+)|^2\ud x
=\int_D\one_{\{u>b\}}|\nabla(u-b)|^2\ud x
=\int_D\nabla(u-b)\cdot\nabla ((u-b)^+)\ud x
=0,
$$
and Poincar\'e's inequality gives $(u-b)^+=0$. Thus
$0\le u\le b\le R^2/(2d)$ on $D$.
\end{proof}

The following inequality explains how the landscape function controls
spectral quantities.

\begin{theorem}[Landscape inequality]\label{thm:landscape-inequality}
The landscape function satisfies $u>0$ almost everywhere on $D$, and
\begin{align}\label{eq:landscape-inequality}
\int_D\frac{|v|^2}{u}\ud x
\le\int_D|\nabla v|^2\ud x,
\quad v\in H_0^1(D).
\end{align}
\end{theorem}
\begin{proof}
First let $v\in C_{\rm c}^\infty(D)$ and $\eps>0$. The product rule and chain rule
(the latter applied to $\rho-\rho(0)$, where
$\rho(t)=1/(t+\eps)$ for $t\ge0$, extended to a
$C^1$-function on $\R$ with bounded derivative)
show that $\psi_\eps:=|v|^2/(u+\eps)$ belongs to $H_0^1(D)$,
so it is an admissible test function in \eqref{eq:landscape-weak}. The product and chain rules give
$$
\nabla\psi_\eps
=
\frac{\ov v\nabla v+v\,\ov{\nabla v}}{u+\eps}
-\frac{|v|^2\nabla u}{(u+\eps)^2}.
$$
Since $\psi_\eps$ and $u$ are real-valued, testing
\eqref{eq:landscape-weak} with $\psi_\eps$ gives
\begin{align*}
\int_D\frac{|v|^2}{u+\eps}\ud x
=\int_D\nabla u\cdot\nabla\psi_\eps\ud x
=2\Re\int_D\frac{\ov v\nabla v\cdot\nabla u}{u+\eps}\ud x
-\int_D\frac{|v|^2|\nabla u|^2}{(u+\eps)^2}\ud x.
\end{align*}
By completing squares, we obtain
\begin{align*}
\int_D\frac{|v|^2}{u+\eps}\ud x
=\int_D|\nabla v|^2\ud x
-\int_D\Bigl|\nabla v-\frac{v\nabla u}{u+\eps}\Bigr|^2\ud x
\le\int_D|\nabla v|^2\ud x.
\end{align*}
In particular,
$$\eps^{-1}\int_{\{u=0\}}|v|^2\ud x\le\n\nabla v\n_2^2.$$
Letting $\eps\downarrow0$ shows that $v=0$ almost everywhere on
$\{u=0\}$. Choosing countably many test functions whose nonzero sets
cover $D$, we conclude that $u>0$ almost everywhere.

Monotone convergence now gives \eqref{eq:landscape-inequality} for
$v\in C_{\rm c}^\infty(D)$. For general $v\in H_0^1(D)$, choose such
functions $v_n$ converging to $v$ in $H^1(D)$ and, after passing to a
subsequence, almost everywhere. Fatou's lemma gives
\[
\int_D\frac{|v|^2}{u}\ud x
\le\liminf_{n\to\infty}\int_D\frac{|v_n|^2}{u}\ud x
\le\lim_{n\to\infty}\n\nabla v_n\n_2^2
=\n\nabla v\n_2^2.
\qedhere\]
\end{proof}

We call $W:=1/u$ the \emph{effective potential}.\index{effective potential}
The landscape inequality says that the Dirichlet energy controls the
weighted integral $\int_DW|v|^2\ud x$. In particular, if
$-\Delta_{\rm Dir}\phi=\lambda\phi$, then
$$
\int_DW|\phi|^2\ud x
\le\n\nabla\phi\n_2^2
=\lambda\n\phi\n_2^2.
$$
Consequently, for every $a>1$,
\begin{align*}
\int_{\{x\in D:\,W(x)\ge a\lambda\}}|\phi(x)|^2\ud x
\le\frac1a\n\phi\n_2^2.
\end{align*}
Thus an eigenfunction has little $L^2$-mass where the effective
potential is large compared with its eigenvalue. This estimate does
not require any boundedness or continuity of the eigenfunction.

There are also simple bounds for the first eigenvalue:
\begin{align}\label{eq:landscape-first-eigenvalue}
\frac1{\n u\n_\infty}
\le\lambda_1
\le\frac{\int_Du\ud x}{\int_Du^2\ud x}.
\end{align}
The lower bound follows from \eqref{eq:landscape-inequality} and the
Courant--Fischer theorem (Theorem \ref{thm:CF}). For the upper bound,
the same theorem gives
$$
\lambda_1
=\inf_{\substack{v\in H_0^1(D)\\v\ne0}}
  \frac{\n\nabla v\n_2^2}{\n v\n_2^2}
\le \frac{\n\nabla u\n_2^2}{\n u\n_2^2},
$$
since $u\not=0$ in $H_0^1(D)$.
Taking $u$ as the test function in \eqref{eq:landscape-weak} gives
$$
\n\nabla u\n_2^2=\int_Du\ud x.
$$
Substituting this into the preceding inequality, we obtain the upper bound.

\begin{remark}[A pointwise bound]\label{rem:landscape-pointwise}
For a bounded eigenfunction $\phi\in\Dom(-\Delta_{\rm Dir})\cap L^\infty(D)$ with
$-\Delta_{\rm Dir}\phi=\lambda\phi$, we also have $|\phi|\le\lambda u\n\phi\n_\infty$.
%
For real-valued $\phi$, this follows by applying the positivity-preserving operator
$-\Delta_{\rm Dir}^{-1}$ to
$-\n\phi\n_\infty\one\le\phi\le\n\phi\n_\infty\one$
and using $-\Delta_{\rm Dir}^{-1}\phi=\lambda^{-1}\phi$.
For complex-valued $\phi$, apply the real-valued estimate to
$\Re(\zeta\phi)$ and take the supremum over a countable dense subset
of $\{\zeta\in\C:\,|\zeta|=1\}$.
\end{remark}

We next use the landscape to count eigenvalues. The argument will
combine \eqref{eq:landscape-inequality} with the following elementary
Poincar\'e inequality on cubes. In particular, it applies to restrictions of functions on $\R^d$, as needed below.

\begin{lemma}[Poincar\'e inequality for cubes]\label{lem:landscape-cube-Poincare}\index{inequality!Poincar\'e}
Let $Q = (0,\ell)^d$. For $v\in H^1(Q)$, put
$v_Q:=|Q|^{-1}\int_Qv\ud x$. Then
\begin{align*}
\int_Q|v-v_Q|^2\ud x
\le \frac{\ell^2}{\pi^2}\int_Q|\nabla v|^2\ud x.
\end{align*}
Equality holds for $v(x)=\cos(\pi x_1/\ell)$.
\end{lemma}
\begin{proof}
Expanding
$v=\sum_{n\in\{0,1,\ldots\}^d}c_ne_n$ in
the normalised tensor product cosine basis,
integrating by parts, and using Parseval's identity, we obtain
$$
\int_Q|\nabla v|^2\ud x
=\frac{\pi^2}{\ell^2}\sum_n|n|^2|c_n|^2
\ge\frac{\pi^2}{\ell^2}\int_Q|v-v_Q|^2\ud x,
$$
with equality if $v(x)=\cos(\pi x_1/\ell)$.
\end{proof}

For $E>0$, put
$$ \ell_E:=\frac{\pi}{2\sqrt E},\qquad
\mathscr Q_E:=\{\ell_E(m+(0,1)^d):\,m\in\Z^d\}.
$$
The cubes in $\mathscr Q_E$ are disjoint and cover $\R^d$ up to a set
of measure zero. Let $\widetilde u$ be the zero extension of the landscape function $u$ and
define
\begin{align*}
\mathscr G_E :=\Bigl\{Q\in\mathscr Q_E:\,
\esssup_{x\in Q}\widetilde u(x)\ge\frac1{4E}\Bigr\},
\qquad
M_u(E):=\#\mathscr G_E.
\end{align*}
This number is finite because $D$ is bounded. Thus $M_u(E)$ counts
the cubes, at a scale proportional to $E^{-1/2}$, on which the
landscape function has sufficiently large essential supremum.

\begin{theorem}[A landscape bound for the eigenvalue count]
\label{thm:landscape-count}
Let $N_{\rm Dir}(E)$ be the number of Dirichlet eigenvalues at most $E$,
counted with multiplicity. Then
$$
N_{\rm Dir}(E)\le M_u(E),\quad E>0.
$$
\end{theorem}
\begin{proof}
Let $F_E$ be the span of the eigenfunctions with eigenvalues at most
$E$. By Theorem \ref{thm:spectra-Laplacian}, this space has dimension
$N_{\rm Dir}(E)$, and
\begin{align}\label{eq:landscape-spectral-energy}
\n\nabla v\n_2^2\le E\n v\n_2^2,
\quad v\in F_E.
\end{align}
For $v\in H_0^1(D)$, its zero extension $\widetilde v$ belongs to
$H^1(\R^d)$, with
$\n\nabla\widetilde v\n_{L^2(\R^d)}=\n\nabla v\n_{L^2(D)}$;
this follows by extending the approximating functions from
$C_{\rm c}^\infty(D)$ by zero.

Consider the linear map which assigns to $v\in F_E$ the averages
of $\widetilde v$ over the cubes in $\mathscr G_E$. We will show
that this map is injective. Suppose that all these averages vanish.
By Lemma \ref{lem:landscape-cube-Poincare},
$$
\sum_{Q\in\mathscr G_E}\int_Q|\widetilde v|^2\ud x
\le \frac{\ell_E^2}{\pi^2}\int_{\R^d}|\nabla\widetilde v|^2\ud x
=\frac1{4E}\n\nabla v\n_2^2.
$$
On each remaining cube we have $\widetilde u<1/(4E)$ almost
everywhere, so $W\ge4E$ on its intersection with $D$. The landscape
inequality therefore gives
$$
\sum_{Q\in\mathscr Q_E\setminus\mathscr G_E}
\int_Q|\widetilde v|^2\ud x
\le\frac1{4E}\int_DW|v|^2\ud x
\le\frac1{4E}\n\nabla v\n_2^2.
$$
Adding the two estimates and using
\eqref{eq:landscape-spectral-energy}, we obtain
$$
\n v\n_2^2
\le\frac1{2E}\n\nabla v\n_2^2
\le\frac12\n v\n_2^2.
$$
Thus $v=0$. The averaging map is injective and takes values in a
space of dimension $M_u(E)$, so $N_{\rm Dir}(E)=\dim F_E\le M_u(E)$.
\end{proof}

The theorem gives an upper bound at every energy, with the
correct order but not Weyl's sharp constant. To compare
leading terms, consider the expression
\begin{align}\label{eq:landscape-effective-Weyl}
\mathcal N_u(E)
:=\frac{\omega_d}{(2\pi)^d}
\int_D(E-W(x))_+^{d/2}\ud x,
\quad E>0.
\end{align}
At a point $x$, the quantity $\omega_d(E-W(x))_+^{d/2}$ is the
volume of the set of $\xi\in\R^d$ satisfying
$|\xi|^2+W(x)\le E$.

\begin{proposition}[Effective Weyl approximation]\label{prop:landscape-Weyl}
For every nonempty bounded open set $D\subseteq\R^d$, we have
$$
\lim_{E\to\infty}\frac{\mathcal N_u(E)}{E^{d/2}}
=\frac{\omega_d}{(2\pi)^d}|D|.
$$
If, in addition, $|\partial D|=0$, then
$\mathcal N_u(E)\sim N_{\rm Dir}(E)$ as $E\to\infty$.
\end{proposition}

This comparison is asymptotic and gives no bound at fixed $E$.

\begin{proof}
Since $u>0$ almost everywhere, $W$ is finite almost everywhere and we have
$(1-W/E)_+^{d/2}\to 1$ almost everywhere as $E\to\infty$.
These functions take values in $[0,1]$. Hence dominated convergence
gives
$$
\frac{\mathcal N_u(E)}{E^{d/2}}
=\frac{\omega_d}{(2\pi)^d}
\int_D\Bigl(1-\frac{W(x)}E\Bigr)_+^{d/2}\ud x
\to\frac{\omega_d}{(2\pi)^d}|D|.
$$
The final assertion follows from Weyl's theorem
(Theorem \ref{thm:Weyl}).
\end{proof}

\begin{example}\label{ex:landscape-interval}
For $D=(0,L)$, the landscape and its reciprocal are
$$u(x)=\frac12x(L-x), \qquad W(x)=2/(x(L-x)).$$
Since
$\n u\n_\infty=\frac{L^2}8$, $\int_0^Lu\ud x=\frac{L^3}{12}$, and $\int_0^Lu^2\ud x=\frac{L^5}{120}$,
the bounds \eqref{eq:landscape-first-eigenvalue} become
$$
\frac8{L^2}\le\lambda_1\le\frac{10}{L^2}.
$$

Rescaling Example \ref{ex:ev-Dir}, we see that the exact values of the Dirichlet eigenvalues
on $(0,L)$ are
$
\lambda_n={\pi^2n^2}/{L^2}$,
$n\ge1,
$
each with multiplicity one, with corresponding eigenfunctions
$x\mapsto\sin(\pi nx/L)$. In particular,
$\lambda_1=\pi^2/L^2$. Counting the eigenvalues at most $E$ gives
$$
N_{\rm Dir}(E)
=\#\bigl\{n\ge1:\, n\le L\sqrt E/\pi\bigr\}
=\lfloor L\sqrt E/\pi\rfloor.
$$

For the landscape expression, the volume of the unit ball in
$\R$ is $\omega_1=2$, so the prefactor in
\eqref{eq:landscape-effective-Weyl} equals $1/\pi$.
Substituting $W(x)=2/(x(L-x))$, we obtain
$$
\mathcal N_u(E)
=\frac1\pi\int_0^L
\Bigl(E-\frac2{x(L-x)}\Bigr)_+^{1/2}\ud x.
$$
Both expressions have leading term $L\sqrt E/\pi$ as
$E\to\infty$, the first by the formula involving the integer
part and the second by Proposition \ref{prop:landscape-Weyl}.
\end{example}

\begin{problems}

\item
Let $\aa$ be a closed accretive form in $H$. Show that if $V:=\Dom(\aa)$ admits an inner product $\iprod{\cdot}{\cdot}_V$ turning $V$ into a Hilbert space such that the inclusion mapping from $V$ into $H$ is bounded, then the associated norm $\n \cdot\n_V$ is equivalent to $\n \cdot\n_\aa$.

\noindent{\em Hint:}\ Both norms make $V$ a Banach space,
and convergence in either norm implies convergence in $H$.
Use this to show that the identity
$I:(V,\n\cdot\n_V)\to(V,\n\cdot\n_\aa)$ has closed graph.
Apply the closed graph theorem to $I$, and then the open
mapping theorem to obtain boundedness of $I^{-1}$.

\item \label{prob:closedVnonclosedH}
Give an example of an accretive form $\aa$ in $H$, which is bounded (and hence closed) as a form on $V:= \Dom(\aa)$ but not closed as a form in $H$.

\item
Let $H:=L^2(0,1)$, and define
$$
\aa(u,v):=u(0)\ov{v(0)},\quad u,v\in \Dom(\aa):= C([0,1]).
$$
Show that $\aa$ is densely defined and positive, but not closable.

\noindent{\em Hint:} \ For nonclosability, consider the functions
$u_n(t):=\max(1-nt,0)$.

\item\label{prob:derivative-form}
Let $H:=L^2(\R)$, and define the sesquilinear form
$$
\aa(u,v):=\iprod{u'}{v},\quad u,v\in\Dom(\aa):=H^1(\R).
$$
\begin{enumerate}[\rm(a), leftmargin=*]
\item Show that $\aa$ is densely defined, accretive, closed, and bounded as a form on $H^1(\R)$ in the norm of $H^1(\R)$,
$$
|\aa(u,v)|\le \|u\|_{H^1(\R)}\|v\|_{H^1(\R)},\quad u,v\in H^1(\R).
$$

\item Using integration by parts, show that for all $u,v\in H^1(\R)$ we have $\aa^\star(u,v)=-\aa(u,v)$,
and deduce that $\Re\aa(u,u)=0$ for all $u\in H^1(\R)$.
In particular, $\aa$ is not coercive.

\item Show that for all $u\in H^1(\R)$ we have
$$
\|u\|_\aa^2=\Re\aa(u,u)+\|u\|^2=\|u\|_{L^2(\R)}^2.
$$
Deduce that $\aa$ is not closed as a form in $H$.

\item Show that $\aa$ is not bounded as a form in $H$.

\noindent {\em Hint:}\ Take $u_n(x):=e^{inx}\phi(x)$ with $\phi\in C_{\rm c}^\infty(\R)$ nonzero.
\end{enumerate}

\item\label{prob:shifted-form}
Let $\aa$ be a densely defined closed accretive form in a Hilbert space $H$, and let $A$ be the operator associated with $\aa$. Define a form $\aa_1$ with domain $\Dom(\aa_1):=\Dom(\aa)$ by
$$
\aa_1(u,v):=\aa(u,v)+\iprod{u}{v},\quad u,v\in\Dom(\aa).
$$

\begin{enumerate}[label={\rm(\alph*)}, leftmargin=*]
\item Show that $\aa_1$ is densely defined, closed, and coercive, more precisely,
$$
\Re\aa_1(u,u)\ge \|u\|^2\!,\quad u\in\Dom(\aa).
$$

\item\label{it:prob:shifted-form-d} Show that if $\aa$ is bounded, then $\aa_1$ is bounded as well.

\item Show that the operator associated with $\aa_1$ is $A+I$.
\end{enumerate}

\item
Let $(\iota,V,H)$ be a Gelfand triple. Show that if $\iota: V \to H$ has dense range, then the mapping $\phi: H \to V'$ described in Remark \ref{rem:Gelfand-triple} is injective and has dense range.

\item\label{prob:form-elliptic1}
Let $(\iota,V,H)$ be a Gelfand triple. A form $\aa$ on $V$ is said to be
{\em elliptic} if there exist $\lambda>0$ and $\alpha>0$ such that
$$
  \Re\aa(v,v)+\lambda\|\iota(v)\|_H^2\ge \alpha\|v\|_V^2,
  \quad v\in V.
$$

\begin{enumerate}[\rm(a), leftmargin=*]
\item Show that the form $\aa$ on $V$ is elliptic if and only if
there exists $\lambda>0$ such that the form
$$
\aa_\lambda(u,v):=\aa(u,v)
+\lambda\iprod{\iota(u)}{\iota(v)}_H
$$
is coercive on $V$.

\item Let $\aa$ be a bounded coercive form on
 $V$, and let $B \in\calL(V,H)$. Show that the form $\aa_B$ on $V$ defined by
 $$
         \aa_B(u,v):=\aa(u,v)+\iprod{Bu}{\iota(v)},\quad u,v\in V,
 $$
 is bounded and elliptic.

\item
Let $\aa$ be a bounded elliptic form on $V$. State and prove
versions of Theorem \ref{thm:form-accretive-resolvent} and Corollary \ref{cor:form-coercive-resolvent} for the operator
$A$ in $H$ associated with $\aa$.
\end{enumerate}

\item\label{prob:form-elliptic3}
Let $D\subseteq \R^d$ be a nonempty bounded open set, and let
$a:D\to M_d(\K)$ have bounded measurable coefficients satisfying the uniform ellipticity condition
$$ \Re\sum_{i,j=1}^d a_{ij}(x)\xi_j\ov{\xi_i}
\ge\alpha|\xi|^2\!, \quad \xi\in \K^d\!,$$
for almost all $x\in D$, with $\alpha>0$.
Let $b:D\to \K^d$ have bounded measurable coordinate functions and let $c:D\to \K$ be bounded and measurable.
Show that
$$ \aa(u,v) := \int_D a \nabla u \cdot\ov{\nabla v}\ud x + \int_D b\cdot \nabla u \, \ov v\ud x + \int_D cu\ov v\ud x, \quad u,v\in H_0^1(D),$$
defines an elliptic form $\aa$.

\item
Let $H=L^2(0,1)$ and define
$$
\aa(u,v):=\int_0^1x\,u(x)\ov{v(x)}\ud x,
\quad u,v\in H.
$$
Show that $\aa$ is densely defined, closed, positive, and therefore
sectorial of every angle $\omega\in(0,\frac12\pi)$. Show also that its
associated operator is injective, but that $\aa$ is not coercive.
Conclude that sectoriality, even together with injectivity of the
associated operator, does not imply coercivity.

\item\label{prob:sectorial-form-closable}
Let $A$ be a densely defined linear operator in a Hilbert space $H$. Suppose that, for some $\omega\in(0,\frac12\pi)$, we have $\iprod{Ax}{x}\in\ov{\Sigma_\omega}$ for all $x\in\Dom(A)$.
Show that $$\aa(x,y):=\iprod{Ax}{y}, \quad x,y\in \Dom(\aa):=\Dom(A),$$
defines a closable form $\aa$ in $H$, and that its closure is $\omega$-sectorial.

\item
\label{prob:positive-form-inner-product}
Let $\aa:H\times H\to\C$ be a {\em positive definite} form, that is, $\aa$ is positive and nondegenerate in the sense that $\aa(u,u)=0$ implies $u=0$.
\begin{enumerate}[\rm(a), leftmargin=*]
\item Show that $\aa$ is an inner product on $H$.
\item Show that if $\aa$ is also closed and coercive, then $u\mapsto \aa(u,u)^{1/2}$ is an equivalent norm on $H$.
\end{enumerate}

\item\label{Prob:dom-Dir-Neum-Lapl}
Prove the following facts for $d=1$ and $D = (a,b)$:
\begin{enumerate}[\rm(a), leftmargin=*]
  \item $\Dom(\Delta_{\rm Dir}) = \{f\in H^2(D):\, f(a) = f(b) = 0\}$;
  \item $\Dom(\Delta_{\rm Neum}) = \{f\in H^2(D):\, f'(a) = f'(b) = 0\}$.
\end{enumerate}

\item\label{prob:compact-resolvent}
Let $A$ be a closed operator acting in a Banach space $X$.

\begin{enumerate}[\rm(a), leftmargin=*]
\item
Show that for all $\lambda\in\varrho(A)$,
$$
\sigma(R(\lambda,A))\setminus\{0\}
=
\Bigl\{\frac1{\lambda-\mu}:\, \mu\in\sigma(A)\Bigr\}.
$$
Show also that for $\mu\in\sigma(A)$, the eigenspaces of $A$
corresponding to $\mu$ and of $R(\lambda,A)$ corresponding to
$(\lambda-\mu)^{-1}$ coincide.

\item
Suppose that $\varrho(A)\ne\emptyset$ and that
$R(\lambda_0,A)$ is compact for some
$\lambda_0\in\varrho(A)$. Show that $R(\lambda,A)$ is compact for
every $\lambda\in\varrho(A)$.

\item
Under the assumptions of part {\rm(b)}, show that every point of
$\sigma(A)$ is an isolated eigenvalue with finite-dimensional
eigenspace, and that $\sigma(A)$ is finite or countable. If it is
infinite, show that its elements can be enumerated as a sequence
$(\mu_n)_{n\ge1}$ satisfying $|\mu_n|\to\infty$.
\end{enumerate}

\item
Let $B_R:=\{x\in\R^d:\, |x|<R\}$, where $R>0$.
Compute its landscape function and use
\eqref{eq:landscape-first-eigenvalue} to prove that its first Dirichlet eigenvalue satisfies
$$ \frac{2d}{R^2}\le\lambda_1(B_R)\le\frac{d(d+4)}{2R^2}. $$

\end{problems}

%% file: ch13-Semigroups.tex
\chapter{Semigroups of Linear Operators}\label{chap:semigroups}

\blfootnote{\rm \small \begin{quote}This book has been published by Cambridge University Press in the series `Cambridge Studies in Advanced Mathematics'. The present version is free to view and download for personal use only. Not for re-distribution, re-sale or use in derivative works.\newline {\copyright} Jan van Neerven\end{quote}}

\noindent
In this chapter we set up a functional-analytic framework for the study of linear and nonlinear initial value problems. This includes the treatment of parabolic problems such as the heat equation and hyperbolic problems such as the wave equation.  From the operator-theoretic perspective the main challenge is to arrive at a thorough understanding of linear equations. This is achieved through the theory of $C_0$-semigroups developed in the present chapter. Once this is done, nonlinear equations are handled by perturbation techniques.

\section{\texorpdfstring{$C_0$}{C\textunderscore 0}-Semigroups}\label{sec:semigroups}

Equations of mathematical physics describing systems involving time evolution can often be cast in the abstract form
\begin{equation*}
\left\{
\begin{aligned}
u'(t) &=Au(t) + f(t,u(t)),\quad t\in [0,T],\\
u(0) &=u_0,
\end{aligned}
\right.
\end{equation*}
where the unknown is a function $u$ from the time interval $[0,T]$ into a Banach space $X$, the operator $A$ is a linear, usually unbounded, operator acting in $X$, $f:[0,T]\times X\to X$ is a given function, and the initial value $u_0$ is assumed to be an element of $X$.
This initial value problem is referred to as
the {\em abstract Cauchy problem} associated with $A$ and $f$.
In applications, typically $X$ is a Banach space of
functions suited for the particular problem and $A$ is a partial differential
operator. For instance, for the heat equation on a bounded open subset $D$ of $\R^d$ subject to Dirichlet boundary conditions one could choose $X=L^2(D)$ and take $A$ to be the Dirichlet Laplacian studied in the previous chapter.

If $A$ is a {\em bounded} operator, the unique solution $u$ of the linear abstract Cauchy problem\index{Cauchy problem!linear}
\begin{equation*}\label{ACP}\tag{ACP}
\left\{
\begin{aligned}
u'(t) &=Au(t),\quad t\in [0,T],\\
u(0) &=u_0,
\end{aligned}
\right.
\end{equation*}
is given by
$$ u(t) = e^{tA}u_0 =  \sum_{n=0}^\infty \frac{t^n}{n!}A^n u_0,\quad t\in [0,T].$$
The operators $e^{tA}$ may be thought of as `solution
operators' mapping the initial value $u_0$ to the solution
$e^{tA}u_0$ at time $t$. For unbounded operators $A$, this simple strategy does not work, since we run into convergence and domain issues. In the case of selfadjoint operators $A$, and, more generally, normal operators $A$ acting in a Hilbert space, one could instead use the functional calculus of Chapter \ref{ch:unbdd} to define the exponentials $e^{tA}$. This would still limit the scope and applicability of the theory considerably. In order to set up a more general and flexible framework, we take a more abstract approach which is motivated by the properties of the exponentials $e^{tA}$ for bounded operators $A$: they satisfy
$e^{0A} = I$ and $e^{tA}e^{sA} = e^{(t+s)A}$, and the mapping $t\mapsto e^{tA}$
is continuous with respect to the operator norm.

\subsection{Definition and General Properties}

Throughout this chapter, $X$ is a Banach space and $H$ is a Hilbert space.

The preceding discussion suggests the following definition.

\begin{definition}[$C_0$-Semigroups]\label{def:semigroup}
A family $S=(S(t))_{t\ge 0}$ of bounded operators acting
on $X$ is called
a $C_0$-{\em semigroup}\index{C0semigroup@$C_0$-semigroup}\index{semigroup, $C_0$-} if the following three properties are satisfied:
\begin{enumerate}[leftmargin=*,label=(S\arabic*)]
\item\label{S1} $S(0)=I$;

\item\label{S2} (semigroup property) $S(t)S(s)=S(t+s)$ for all $t,s\ge 0$;

\item\label{S3} (strong continuity) $\lim_{t\downarrow 0} \n S(t)x-x\n =0$ for all $x\in X$.
\end{enumerate}
Its {\em infinitesimal generator}\index{infinitesimal generator}, or briefly the {\em generator}\index{generator},
is the linear operator $A$ defined by
\begin{align*}
\Dom(A) & := \Bigl\{x\in X:\, \lim_{t\downarrow 0} \frac{1}{t} (S(t)x-x)\hbox{ exists in $X$}\Bigr\},\\
     Ax & := \lim_{t\downarrow 0} \frac{1}{t} (S(t)x-x),\quad x\in \Dom(A).
\end{align*}
\end{definition}

The idea is to interpret the orbit
$ u(t):= S(t)u_0$
as the `solution' of the linear problem \eqref{ACP}.
To find a precise way of making this idea rigorous, and subsequently to
cover nonlinear initial value problems as well, is among the main objectives
of this chapter.

\begin{remark}[Strong convergence versus uniform convergence]
The properties of $e^{tA}$ suggest replacing \ref{S3} by the stronger condition
$\lim_{t\downarrow 0} \n S(t)-I\n =0$. As it turns out, however,
this condition forces the generator $A$ to be bounded (see Problem \ref{prob:bddgen}). This renders the theory useless, as it would fail to cover equations in which $A$ is a differential operator acting in a Banach space $X$ of functions. In a sense the strong convergence imposed in \ref{S3} is also more natural, as it ensures that the `solution' $u(t)=S(t)u_0$ is continuous in the norm of $X$ (see Proposition \ref{prop:sg-prop}).
\end{remark}

The next two propositions collect
some elementary properties of $C_0$-semi\-groups and their generators.

\begin{proposition}\label{prop:boundS}
Let $S$ be a $C_0$-semigroup on $X$.
There exist $M\ge 1$ and $\omega\in \R$ such that $\n S(t)\n
\le Me^{\omega t}$ for all $t\ge 0$.
\end{proposition}
\begin{proof}
There exists a number $\delta>0$ such that
$\sup_{t\in [0,\delta]}\n S(t) \n =:\sigma <\infty$. Indeed, otherwise we could find a
sequence $t_n\downarrow 0$ such that $\limn \n S(t_n)\n = \infty$. By the uniform
boundedness theorem, this implies the existence of an $x\in X$ such that
$\sup_{n\ge 1} \n S(t_n)x\n = \infty$, contradicting the strong continuity assumption
\ref{S3}.

By the semigroup property \ref{S2}, if
$t\in[(k-1)\delta,k\delta]$, then $\n S(t)\n\le \sigma^k\le \sigma^{1+t/\delta}$, where
the second inequality uses that $\sigma\ge 1$ by \ref{S1}. This proves the
proposition, with $M = \sigma$ and $\omega = \frac1\delta \log\sigma$.
\end{proof}

We will frequently use the trivial observation that
if $A$ generates the $C_0$-semigroup $(S(t))_{t\ge 0}$, then for all scalars $\mu$
the linear operator $A - \mu$ generates the $C_0$-semigroup $(e^{-\mu t}S(t))_{t\ge 0}$.
For $\mu>\om$, with $\om$ as in Proposition \ref{prop:boundS}, this rescaled semigroup has exponential decay in operator norm.

\begin{proposition}\label{prop:sg-prop}
Let $S$ be a $C_0$-semigroup on $X$ with generator $A$. The following assertions hold:
\begin{enumerate}[label={\rm(\arabic*)}, leftmargin=*]
\item\label{it:sg-prop1} for all $x\in X$ the orbit $t\mapsto S(t)x$ is continuous for
$t\ge 0$;
\item\label{it:sg-prop2} for all $x\in \Dom(A)$ the orbit $t\mapsto S(t)x$ is continuously
differentiable for $t\ge 0$, we have $S(t)x\in \Dom(A)$, and
$$\frac{{\rm d}}{{\rm d}t}S(t)x = AS(t)x = S(t)Ax,\quad t\ge 0;$$
\item\label{it:sg-prop3} for all $x\in X$ and $t\ge 0$ we have $\int_0^t S(s)x\ud s\in \Dom(A)$ and
$$A\int_0^t S(s)x\ud s = S(t)x-x,$$ and if $x\in \Dom(A)$, then both sides are equal
to $\int_0^t S(s)Ax\ud s $;
\item\label{it:sg-prop4} the generator $A$ is a densely defined closed operator.
\end{enumerate}
\end{proposition}

\begin{proof}
The proof uses the calculus rules for Banach space-valued Riemann integrals (Proposition \ref{prop:der-f-zero}).

\smallskip
\ref{it:sg-prop1}: \ Right continuity of $t\mapsto S(t)x$ follows from the right continuity
at $t=0$ in \ref{S3} and the semigroup property \ref{S2}. For left continuity,
observe that by the semigroup property, for $0<h<t$ we have
$$\n S(t)x - S(t-h)x\n \le \n S(t-h)\n \n S(h)x-x\n
\le \sup_{s\in [0,t]}\n S(s)\n\n S(h)x-x\n,$$ where the supremum is finite by Proposition \ref{prop:boundS}.

\smallskip
\ref{it:sg-prop2}: \ Fix $x\in \Dom(A)$ and $t\ge 0$. By the semigroup property we have
$$ \lim_{h\downarrow 0} \frac1h (S(t+h)x - S(t)x) =  S(t) \lim_{h\downarrow 0}
\frac1h (S(h)x - x) = S(t)Ax.$$
This proves all assertions except left differentiability.
For $t>0$ we note that
$$\lim_{h\downarrow 0} \frac1{h}(S(t)x-S(t-h)x) = \lim_{h\downarrow 0} S(t-h)\Big(\frac1h(S(h)x-x)\Big) = S(t)Ax, $$
where we used the fact that $x\in \Dom(A)$ and the fact that the convergence $\lim_{h\downarrow 0} S(t-h)y = S(t)y$ for all $y\in X$ implies convergence uniformly on compact sets by Proposition \ref{prop:uniform-limits-on-K}.
\smallskip

\ref{it:sg-prop3}: \ The first identity follows from
\begin{align*} \lim_{h\downarrow 0} \frac1h(S(h)-I)\int_0^t S(s)x\ud s
&  = \lim_{h\downarrow 0} \frac1h \Big(\int_0^t S(s+h)x\ud s - \int_0^t S(s)x\ud s\Big)
\\ & =  \lim_{h\downarrow 0}\frac1h\Big(\int_t^{t+h} S(s)x\ud s -\int_0^h
S(s)x\ud s \Big)
 = S(t)x-x,
\end{align*}
where we first made a substitution and then used the continuity of $t\mapsto S(t)x$.
 The identity for $x\in \Dom(A)$ follows by integrating the identity of part \ref{it:sg-prop2}, or by
noting that
\begin{align*}
 \lim_{h\downarrow 0} \frac1h(S(h)-I)\int_0^t S(s)x\ud s
& = \lim_{h\downarrow 0}  \int_0^t S(s) \Bigl(\frac1h(S(h)x-x)\Bigr)\ud s
 = \int_0^t S(s)Ax\ud s,
\end{align*}
where the convergence under the integral is justified by the fact that the convergence of the difference quotients $\frac1h(S(h)x-x)$ to $Ax$
implies uniform convergence of the integrands on $[0,t]$.
\smallskip

\ref{it:sg-prop4}:  \ Denseness of $\Dom(A)$ follows from \ref{it:sg-prop1} and the first part of \ref{it:sg-prop3}: for any $x\in X$, the latter implies that
$\int_0^t S(s)x\ud s \in \Dom(A)$ for all $t>0$, while the former implies that
$\lim_{t\downarrow 0} \frac1t \int_0^t S(s)x\ud s = x$.

To prove that $A$ is closed we must check that the graph $\Gr(A)= \{(x,Ax):\, x\in\Dom(A)\}$ is closed in $X\times X$. Suppose that $(x_n)_{n\ge 1}$ is a sequence in
$\Dom(A)$ such that $\limn x_n = x$   and $\limn Ax_n = y$ in $X$.
Then, by the second part of \ref{it:sg-prop3},
\begin{align*} \frac1h (S(h)x-x) &= \limn \frac1h (S(h)x_n-x_n)
 = \limn  \frac1h \int_0^h S(s)Ax_n\ud s
 = \frac1h \int_0^h S(s)y\ud s.
\end{align*}
Passing to the limit as $h\downarrow 0$, we obtain $x\in \Dom(A)$ and $Ax = y.$
\end{proof}

We have just seen that the generator of a $C_0$-semigroup
is always densely defined and closed. As a consequence of the latter, $\Dom(A)$ is a Banach space with respect to its graph norm.
In various applications it is of interest to know when a
subspace $Y$\!, which is dense in $X$ and contained in $\Dom(A)$, is dense in $\Dom(A)$.
If this is the case, $Y$ is called a {\em core}\index{core} for $A$.
The next result gives a simple sufficient condition.

\begin{proposition}\label{prop:core}
Let $S$ be a $C_0$-semigroup with generator $A$ on $X$. If $Y$ is a subspace of $\Dom(A)$
 which is dense in $X$ and invariant under each operator $S(t)$, then $Y$ is dense in $\Dom(A)$.
\end{proposition}
\begin{proof}
Let $Z$ be the closure of $Y$ in $\Dom(A)$ with respect to the
graph norm. For $t>0$ define
$$
J_t x:=\frac1t\int_0^t S(s)x\ud s,\quad x\in X.
$$
By Proposition \ref{prop:sg-prop}, $J_t x\in\Dom(A)$ and
$$
AJ_t x=\frac1t(S(t)x-x).
$$
These identities show that $J_t$ is bounded from $X$ into
$\Dom(A)$ with its graph norm.

For $y\in Y$, the function $s\mapsto S(s)y$ is continuous in
the graph norm, since $AS(s)y=S(s)Ay$. The Riemann sums of the integral $\int_0^t S(s)y\ud s$
belong to $Y$ by the assumed invariance of $Y$ under the semigroup, and therefore $J_t y\in Z$.
Since $Y$ is dense in $X$ and $J_t:X\to\Dom(A)$ is bounded,
it follows that $J_tX\subseteq Z$.

Finally, for $x\in\Dom(A)$, strong continuity gives
$$
J_t x\to x,\qquad AJ_t x=J_tAx\to Ax
$$
as $t\downarrow0$. Thus $J_t x\to x$ in the graph norm.
Since $J_t x\in Z$ and $Z$ is closed in that norm, we obtain
$x\in Z$. This proves that $Z=\Dom(A)$.
\end{proof}

This proposition is often helpful in determining the domain of the generator explicitly when
the semigroup is given; see for instance Section \ref{subsec:translation-sg}.

The proof of the next proposition uses the following version of the product rule. It is proved in the same way as the product rule in calculus; uniform convergence on compact sets follows from Proposition \ref{prop:uniform-limits-on-K}.

\begin{lemma}\label{lem:diff}\index{product!rule}
Let $I\subseteq\R$ be an interval of positive length, let
$T:I\to\calL(X)$ be strongly continuous, and let $u:I\to X$.
Let $t_0\in I$. Suppose that
\begin{enumerate}[label={\rm(\roman*)}, leftmargin=*]
\item $u$ is differentiable at $t_0$;
\item the function $t\mapsto T(t)u(t_0)$ is differentiable at $t_0$.
\end{enumerate}
Then the function $t\mapsto T(t)u(t)$ is differentiable at $t_0$, and
$$
\frac{\ud}{\ud t}\Big|_{t=t_0}T(t)u(t)
=
\frac{\ud}{\ud t}\Big|_{t=t_0}T(t)u(t_0)
+
T(t_0)u'(t_0).
$$
At an endpoint of $I$, the assertion is understood in the one-sided
sense.
\end{lemma}

\begin{proof}
For $t\neq t_0$ we have
\begin{align*}
\frac{T(t)u(t)-T(t_0)u(t_0)}{t-t_0}
&=
T(t)\frac{u(t)-u(t_0)}{t-t_0}
+
\frac{T(t)u(t_0)-T(t_0)u(t_0)}{t-t_0}.
\end{align*}
The second term has the asserted limit by assumption. Strong
continuity implies that the operators $T(t)$ are uniformly bounded for
$t$ in a neighbourhood of $t_0$. Hence
\begin{align*}
&
\Bigl\|
T(t)\frac{u(t)-u(t_0)}{t-t_0}
-
T(t_0)u'(t_0)
\Bigr\|
\\
&\qquad\le
\|T(t)\|
\Bigl\|
\frac{u(t)-u(t_0)}{t-t_0}-u'(t_0)
\Bigr\|
+
\|(T(t)-T(t_0))u'(t_0)\|,
\end{align*}
which tends to zero.
\end{proof}

A $C_0$-semigroup is uniquely determined by its generator:

\begin{proposition}
 If $A$ is the generator of the $C_0$-semigroups $S$ and $T$, then $S(t) = T(t)$ for all $t\ge 0$.
\end{proposition}
\begin{proof}
By Lemma \ref{lem:diff}, for all $t>0$ and $x\in \Dom(A)$ the function $\phi_t(s):= S(t-s)T(s)x$ is continuously differentiable
 on $[0,t]$ with derivative $\phi_t'(s) = -A S(t-s)T(s)x +  S(t-s)AT(s)x = 0$, and therefore $\phi_t$ is constant by
 Proposition \ref{prop:der-f-zero}. Hence, $S(t)x = \phi_t(0) = \phi_t(t) = T(t)x$.
 Since this holds for all $x$ in the dense subspace $\Dom(A)$ of $X$,
it follows that $S(t)=T(t)$.
\end{proof}

The next proposition identifies the resolvent of the generator as the Laplace transform\index{Laplace transform} of the semigroup.

\begin{proposition}\label{prop:resolv}
Let $A$ be the generator of a $C_0$-semigroup $S$ on $X$, and fix constants $M\ge 1$ and $\omega\in\R$
such that $\n S(t)\n\le Me^{\omega t}$ for all $t\ge 0$.
Then $\{\la\in \C:\, \Re\la >\omega\}\subseteq \varrho(A)$, and on this set
the resolvent of $A$ is given by
$$ R(\la, A) x = \int_0^\infty e^{-\la t}S(t)x\ud t, \quad x\in X.$$
As a consequence, for $\Re\la>\omega$ we have
$$ \n R(\la,A)\n \le \frac{M}{\Re\la-\omega}.$$
\end{proposition}

\begin{proof}
Fix $x\in X$ and $\la\in\C$ with $\Re\la>\omega$, and define $R_\la  x:=  \int_0^\infty e^{-\la t}S(t)x\ud t$.
Using the semigroup property \ref{S2} and a substitution,
we obtain the identity
\begin{align*}
 \lim_{h\downarrow 0}\frac1h (S(h)-I)R_\la  x
 &   = \lim_{h\downarrow 0} \frac1h \Bigl(e^{\la h}\int_h^\infty\! e^{-\la t}S(t)x\ud t - \int_0^\infty e^{-\la t}S(t)x\ud t\Bigr)
\\ & = \lim_{h\downarrow 0} \frac{e^{\la h}-1}{h} \int_h^\infty\! e^{-\la t}S(t)x\ud t - \lim_{h\downarrow 0}
\frac1h \int_0^h e^{-\la t}S(t)x\ud t
  = \la R_\la  x - x,
\end{align*}
from which it follows that $R_\la  x \in\Dom(A)$
and $A R_\la  x = \la R_\la  x-x$. This shows that the bounded operator
$R_\la$ is a right inverse for $\la-A$.

Integrating by parts and using the identity
$\frac{{\rm d}}{{\rm d}t} S(t)x = S(t)Ax $ for $x\in\Dom(A)$, we obtain
$$ \la \int_0^T e^{-\la t}S(t)x\ud t = -e^{-\la T}S(T)x +x + \int_0^T e^{-\la
t}S(t)Ax\ud t.$$
Since $\Re\la>\omega$, sending $T\to\infty$ gives $\la R_\la  x = x + R_\la  Ax$. This
shows that $R_\la$ is also a left inverse.

The estimate for the resolvent follows from
\[
\Bigl\n \int_0^\infty e^{-\la t}S(t)x\ud t \Bigr\n \le
 \int_0^\infty e^{-\Re\la t}\n S(t) x\n\ud t
 \le M \n x\n \int_0^\infty e^{(\omega-\Re\la) t} \ud t= \frac{M}{\Re\la - \omega}\n x\n.
\qedhere\]
\end{proof}

The resolvent bound can be strengthened to
$$
\|R(\lambda,A)^k\|\le \frac{M}{(\Re\lambda-\omega)^k}, \quad k\ge 1,
$$
for all $\Re\lambda>\omega$.
This can be proved by expanding $R(\la,A)^k$ as a repeated integral and
using the semigroup property, or by the following elegant argument. With
respect to the equivalent norm
$\|x\|_\omega:=\sup_{t\ge0} e^{-\omega t}\|S(t)x\|$
we have
$\|x\|\le \|x\|_\omega\le M\|x\|$.
If $\|\cdot\|_\omega$ also denotes the induced operator norm on
$\calL(X)$, then
$\|S(t)\|_\omega\le e^{\omega t}$ and hence
$$
  \|R(\la,A)\|_\omega\le \frac1{\Re\la-\omega},\quad \Re\la>\omega.
$$
Taking $k$-th powers and switching back to the original norm gives the
result.

Combining Proposition \ref{prop:resolv} with Proposition \ref{prop:semigroupsAB},
we obtain the following maximality property of generators.

\begin{proposition}\label{prop:AB-sgrs} If $A$ and $B$ generate $C_0$-semigroups on  $X$, and if $B$ is an extension of $A$, then $A=B$.
\end{proposition}
\begin{proof}
Proposition \ref{prop:resolv} implies that the resolvent sets of $A$ and $B$ share a common half-plane.
The equality $A=B$ then follows from Proposition \ref{prop:semigroupsAB}.
\end{proof}

For operators satisfying the resolvent estimate of Proposition \ref{prop:resolv} for real $\la$, we have the following convergence result.

\begin{proposition}\label{prop:resolv-conv}
Let $A$ be a densely defined closed operator acting in $X$, and suppose that for some $M\ge1$ and $\omega\in\R$ we have
$\{\la>\omega\}\subseteq\varrho(A) $ and  $$ \n R(\la,A)\n \le \frac{M}{\la-\omega}, \quad \lambda>\omega.$$
Then for all $x\in X$ we have
$$ \lim_{\lambda\to\infty} \lambda R(\la,A) x = x.$$
\end{proposition}
\begin{proof}
For $x\in\Dom(A)$ we have
$\la R(\la,A)x-x=R(\la,A)Ax$,
and therefore
$$
\n\la R(\la,A)x-x\n
\le\frac{M}{\la-\omega}\n Ax\n\to0
$$
as $\la\to\infty$.
Since $\Dom(A)$ is dense in $X$ and the operators
$\la R(\la,A)$ are uniformly bounded for sufficiently
large $\la$, the convergence extends to every $x\in X$.
\end{proof}

The final result of this subsection gives a useful sufficient condition for a semigroup of operators to be strongly continuous.
We need the following terminology.
A family of bounded operators  $S= (S(t))_{t\ge 0}$ on $X$ is said to be a {\em weakly continuous semigroup} if conditions
\ref{S1} and \ref{S2} in Definition \ref{def:semigroup} hold and \ref{S3} is replaced by the condition that for all $x\in X$ and $x^*\in X^*$ one has $\lim_{t\downarrow 0} \lb S(t) x,x^*\rb  = \lb x,x^*\rb$.

\begin{theorem}[Phillips]\label{thm:Phillips}\index{theorem!Phillips}
Every weakly continuous semigroup is strongly continuous.
\end{theorem}
\begin{proof}
Let $$X_0:= \Bigl\{x\in X:\, \lim_{t\downarrow 0}\n S(t)x-x\n = 0\Bigr\}.$$ It is evident that $X_0$ is a linear subspace of $X$.
Arguing as in the proof of Proposition \ref{prop:boundS}, we see that
the family
$\{S(t):\, 0\le t\le1\}$
is uniformly bounded. It follows that $X_0$ is a closed subspace of
$X$. We wish to show that $X_0=X$.

Fix $x\in X$. For every $x^*\in X^*$, the scalar function
$t\mapsto\lb S(t)x,x^*\rb$
is right-continuous at every $t\ge0$, since
$S(t+h)x=S(t)S(h)x$
and $S(h)x$ converges weakly to $x$ as $h\downarrow0$. In particular,
every orbit is weakly measurable.

Let
$$
Y_x:=
\overline{\operatorname{span}\{S(q)x:\, q\in\Q_+\}}.
$$
For every $t\ge0$, choose rational numbers $q_n\downarrow t$. Then
$S(q_n)x$ converges weakly to $S(t)x$. Since the norm-closed subspace
$Y_x$ is weakly closed, it follows that $S(t)x\in Y_x$. Thus the orbit
of $x$ is contained in the separable closed subspace $Y_x$.
It follows that we can apply the
Pettis measurability theorem (Theorem \ref{thm:Pettis-secondversion}) and conclude that every orbit $t\mapsto S(t)x$ is strongly measurable.
It follows from these considerations that the Bochner integrals
$x_t:= \frac1t\int_0^t S(s)x\ud s$ are well defined.

Fix $x\in X$
and $0<t<\frac12$. For $0< s< t$,
\begin{align*}
\n S(s)x_t-x_t\n
& = \frac1t\Bigl\n \int^t_0 S(s+r)x\ud r -\int^t_0 S(r)x \ud r \Bigr\n
\\ & = \frac1t\Bigl\n \int^{t+s}_t S(r)x\ud r -\int^s_0 S(r)x\ud r \Bigr\n
\le 2s \cdot \frac1t\Bigl(\sup_{0\le r\le 1}\n S(r)\n\Bigr)\,\n x\n.
\end{align*}
This shows that $x_t\in X_0$.

Suppose now, for a contradiction, that $X_0\not=X$. Then there exists an $x\in X\setminus X_0$ and
by the Hahn--Banach theorem
we can find an $x^{*}\in X^{*}$ which vanishes on $X_0$ but not on $x$. Then, with $x_t=\frac1t\int_0^t S(s)x\ud s$ as before,
\begin{align*} 0  = \lim_{t\downarrow 0}\, \lb  x_t, x^{*}\rb
 = \lim_{t\downarrow 0} \frac1t\int_0^t \lb S(s)x, x\s\rb \ud s
 = \lb x,x\s\rb \not=0,
\end{align*}
 a contradiction.
\end{proof}

\subsection{\texorpdfstring{$C_0$}{C\textunderscore 0}-Groups}\label{subsec:C0groups}

Instead of considering only forward time, we could also include backward time. This leads to the notion of a {\em $C_0$-group}.

\begin{definition}[$C_0$-groups]
 A {\em $C_0$-group}\index{C0group@$C_0$-group}\index{group, $C_0$-} is a family $S = (S(t))_{t\in\R}$ of bounded operators acting on $X$
 with the following properties:
\begin{enumerate}[leftmargin=*,label=(G\arabic*)]
\item\label{G1} $S(0)=I$;

\item\label{G2} $S(t)S(s)=S(t+s)$ for all $t,s\in \R$;

\item\label{G3} $\lim_{t\to 0} \n S(t)x-x\n =0$ for all $x\in X$.
\end{enumerate}
Its {\em infinitesimal generator}, or briefly its {\em generator}, is the linear operator $A$ defined by
\begin{align*}
\Dom(A) &:= \Bigl\{x\in X:\, \lim_{t\to 0} \frac{1}{t} (S(t)x-x)\hbox{ exists in $X$}\Bigr\},\\
Ax &:= \lim_{t\to 0} \frac{1}{t} (S(t)x-x),\quad x\in \Dom(A).
\end{align*}
\end{definition}

It is evident from the definition that if $A$ generates a $C_0$-group $(S(t))_{t\in\R}$, then
both $(S(t))_{t\ge 0}$ and $(S(-t))_{t\ge 0}$ are $C_0$-semigroups. Denoting their generators by
$A_+$ and $A_-$, it is evident that $\Dom(A) \subseteq \Dom(A_+)\cap \Dom(A_-)$ and that for all $x\in \Dom(A)$ we have
$Ax = A_+x = -A_-x$. In fact, more is true:

\begin{proposition}\label{prop:C0group}
A linear operator $A$ in $X$ generates a $C_0$-group $(S(t))_{t\in\R}$ if and only if both $A$ and $-A$ generate
$C_0$-semigroups. These semigroups are $(S(t))_{t\ge 0}$ and $(S(-t))_{t\ge 0}$, respectively.
\end{proposition}

\begin{proof}
Let $A$ generate a $C_0$-group $(S(t))_{t\in\R}$, and let $A_+$ be the
generator of the restricted semigroup $(S(t))_{t\ge0}$. It is clear that
$\Dom(A)\subseteq\Dom(A_+)$ and that $Ax=A_+x$ for $x\in\Dom(A)$.
Conversely, let $x\in\Dom(A_+)$. For $h>0$ we have
$$
        \frac{S(-h)x-x}{-h}
        =
        S(-h)\frac{S(h)x-x}{h}.
$$
As $h\downarrow0$, the quotient on the right-hand side tends to $A_+x$,
while $S(-h)\to I$ strongly, and by Proposition \ref{prop:uniform-limits-on-K} the convergence is uniform on compact sets. It follows that the left derivative at $0$ exists and
is equal to $A_+x$. Thus $x\in\Dom(A)$ and $Ax=A_+x$. Therefore
$A=A_+$. Applying the same argument to the restricted semigroup
$(S(-t))_{t\ge0}$, we find that its generator is $-A$.

For the converse, suppose that $A$ and $-A$ generate $C_0$-semigroups
$(S_+(t))_{t\ge 0}$ and $(S_-(t))_{t\ge 0}$, respectively. By Lemma \ref{lem:diff}, for $x\in \Dom(A) = \Dom(-A)$ the function $t\mapsto S_-(t)S_+(t)x $ is continuously differentiable and
$$ \frac{{\rm d}}{{\rm d}t} S_-(t)S_+(t)x = -AS_-(t)S_+(t)x + S_-(t)AS_+(t)x
= 0,$$
where we used the fact that $S_-(t)$ commutes with $A$.
It follows from Proposition \ref{prop:der-f-zero} that $t\mapsto S_-(t)S_+(t)x$
is constant, and evaluation at $t=0$ shows that $S_-(t)S_+(t)x = x$ for all $t\ge 0$.
Since $\Dom(A)$ is dense, this identity extends to arbitrary $x\in X$. This proves that
$S_-(t)$ is a left inverse for $S_+(t)$. Interchanging the roles of $S_-(t)$ and $S_+(t)$, we find that
$S_-(t)$ is also a right inverse for $S_+(t)$. As a result, $S_+(t)$ is invertible and $(S_+(t))^{-1} = S_-(t)$
for all $t\ge 0$.
For $t\in \R$ define
$$ S(t) :=
\begin{cases}
 S_+(t), & \ t \ge 0,\\
 S_-(-t), & \ t \le 0.
\end{cases}
$$
By what we have proved, it is immediate that $(S(t))_{t\in\R}$ is a
$C_0$-group and that $A$ is its generator.
\end{proof}

Proposition \ref{prop:resolv}, applied to the semigroups generated by $\pm A$, implies:

\begin{corollary}
If $A$ generates a uniformly bounded $C_0$-group on $X$, then $\sigma(A)\subseteq i\R$.
\end{corollary}

The spectrum of the generator of a $C_0$-semigroup may be empty (an example is given in Problem \ref{prob:empty-spectrum}).
This contrasts with the second part of the following result. For a uniformly bounded $C_0$-group $S$
on $X$ and $f\in L^1(\R)$ we define $S(f)\in \calL(X)$ by
\begin{align}\label{eq:Sf}S(f)x := \int^\infty_{-\infty} f(t) S(t) x\ud t,\quad x\in X.
\end{align}

\begin{theorem}\label{thm:sp-bdd-group}
If $A$ generates a uniformly bounded $C_0$-group $S$ on $X$, then:
 \begin{enumerate}[label={\rm(\arabic*)}, leftmargin=*]
  \item\label{it:sp-bdd-group1} if the Fourier transform of a function $f\in L^1(\R)$
is compactly supported and vanishes in a neighbourhood of $i\sigma(A)$, then $S(f) =0$;
  \item\label{it:sp-bdd-group2} if $X\not=\{0\}$, then $\sigma(A)\not=\emptyset$.
\end{enumerate}
\end{theorem}
\begin{proof} \ref{it:sp-bdd-group1}:\ For all $\delta>0$ and $s\in\R$ we have $\pm\delta -is\in\varrho(A)$, and for all $x\in X$ we have
$$R(\delta-is,A)x =\int^\infty_0 e^{-(\delta-is)t} S(t)x\ud t$$ and $$
R(-\delta-is,A)x= -R(\delta+is,-A)x=-\int^\infty_0 e^{-(\delta+is)t} S(-t)x\ud t.$$
Hence by dominated convergence, Fourier inversion (Theorem \ref{thm:FT-inversion}),
Fubini's theorem, and Propositions \ref{prop:resolv} and \ref{prop:C0group},
\begin{align*}
S(f)x & = \lim_{\delta\downarrow 0} \int^\infty_{-\infty} e^{-\delta |t|} f(t)S(t) x\ud t
\\ & = \frac1{\sqrt{2\pi}} \lim_{\delta\downarrow 0} \int^\infty_{-\infty}
e^{-\delta |t|} \Bigl(\int^\infty_{-\infty} e^{ist}\wh f(s)\ud s\Bigr)\,S(t) x\ud t
\\ & =
\frac1{\sqrt{2\pi}} \lim_{\delta\downarrow 0} \int^\infty_{-\infty}
\wh f(s)\Bigl(\int^\infty_{-\infty}
e^{-\delta|t|} e^{is t}S(t) x\ud t\Bigr)\ud s
\\ & =
\frac1{\sqrt{2\pi}}\lim_{\delta\downarrow 0}  \int^\infty_{-\infty}
\wh f(s) (R(\delta-is,A)-R(-\delta-is,A))x\ud s.
\end{align*}
Put $K:=\supp(\wh f)$. Since $\wh f$ vanishes in a neighbourhood of
$i\si(A)$, we have
$-is\in\varrho(A)$ for all $s\in K$.
By continuity of the resolvent,
$R(\delta-is,A)-R(-\delta-is,A)\to 0$
uniformly for $s\in K$ as $\delta\downarrow0$. The resolvents occurring
here are uniformly bounded for $s\in K$ and sufficiently small
$\delta>0$. Since $K$ is compact and $\wh f$ is integrable, dominated
convergence gives $S(f)=0$.

\smallskip
\ref{it:sp-bdd-group2}:\ Suppose that $\sigma(A)=\emptyset.$
Part \ref{it:sp-bdd-group1} implies that
$S(f) = 0$ for all $f\in L^1(\R)$ whose
Fourier transform has compact support.
We claim that such functions are dense in $L^1(\R)$.
Indeed, choose $\phi\in C_{\rm c}^\infty(\R)$ with
$\phi(0)\ne0$, and put $\psi:=\widecheck{\phi}$.
Then $\psi\in L^1(\R)$ and $\int_{\R}\psi\ud x\ne0$. Multiplying
$\psi$ by a scalar, we may suppose that $\int_{\R}\psi\ud x=1$.
By Theorem \ref{thm:approx-identity} we then have $\lim_{\eps\downarrow 0}\psi_\eps*f = f$ in $L^1(\R)$, where
$\psi_\eps(x):=\eps^{-1}\psi(\eps^{-1}x)$, and the Fourier transforms $\wh{\psi_\eps*f} = \sqrt{2\pi} \wh{\psi_\eps} \wh{f}$ are compactly supported. This proves the claim.

By approximation we find that $S(f) = 0$ for all $f\in L^1(\R)$.
In particular, taking $f_0(t):=e^{-t}$ for $t\geq 0$ and $f_0(t):=0$ for $t<0$, we obtain
$R(1,A)=S(f_0)=0$ from Proposition \ref{prop:resolv}. Since $\Ran(R(1, A)) = \Dom(A)$ is dense in $X$, this implies $X=\{0\}.$
\end{proof}

We will use this theorem to give a proof of Wiener's Tauberian theorem (Theorem \ref{thm:Wiener}).
Recall that this theorem asserts that if the Fourier transform
of a function $f\in L^1(\R)$ is zero-free, then
the span of the set of all translates of $f$
is dense in $L^1(\R)$.

We start with some preparations. In what follows, $X$ is a complex Banach space.
If $S$ is a uniformly bounded $C_0$-group on $X$, we define
$$I_S:=\{f\in L^1(\R):\, S(f)=0\},$$
where $S(f)$ is given by \eqref{eq:Sf}.
The {\it Arveson spectrum}\index{Arveson spectrum} of $S$ is the set
$$  {\rm Sp}(S) := \{\om\in \R:\, \wh f(\om) =0 \ \hbox{ for all } \ f\in I_S\}.$$
The key to proving Wiener's Tauberian theorem is the following result, which is of independent interest.

\begin{theorem}\label{thm:Jorgensen} Let $S$ be a uniformly bounded $C_0$-group on $X$ with generator $A$. Then $${\rm Sp}(S) = i\sigma(A).$$
\end{theorem}
\begin{proof}
First let $\om\in\R$ satisfy $\om\not\in i\sigma(A)$. Since $i\sigma(A)\subseteq\R$ is closed, we may choose an
open set $U$ containing $i\sigma(A)$ such that
$\omega\notin\ov U$. Pick $\psi\in C^\infty_{\rm c}(\R)$
with $\supp(\psi)\subseteq \R\setminus U$ and $\psi(\omega)\ne 0$. Then
$ f:=\widecheck{\psi} \in L^1(\R)$ and $\wh f=\psi$, so
$\wh f$ vanishes on a neighbourhood of $i\sigma(A)$ but $\wh f(\omega)\ne 0$.
By Theorem \ref{thm:sp-bdd-group}, $S(f) = 0$, so $f\in I_S$. But then $\wh f(\om)
\not = 0$ implies that $\om\not \in {\rm Sp}(S)$.

Conversely, let $\om\in i\si(A)$. Since
$\si(A)\subseteq i\R$, every point of
$\si(A)$ is a boundary point. Since the boundary of the spectrum consists of approximate eigenvalues (see Section \ref{subsec:spectrum-unbdd}, where it was observed that the corresponding result for bounded operators, Proposition  \ref{prop:approx-eigenvalue-bdry}, extends to unbounded operators), $-i\om$ is an approximate eigenvalue of $A$. Hence we may
choose a sequence $(x_n)_{n\ge 1}$ of norm-one vectors in $X$, with
$x_n\in \Dom(A)$ for all $n\ge 1$, such that
$\n A x_n + i \om x_n\n\to 0$.
In view of
$$S(t)x_n- e^{-i\om t}x_n =e^{-i\omega t}\int^t_0 e^{i\om s}S(s) (A+i\om)x_n \ud s \to 0 \ \ \hbox{as} \ n\to\infty,$$
$(x_n)_{n\ge 1}$ is an approximate eigensequence of $S(t)$ with
approximate eigenvalue $e^{-i\om t}$.

Let $f\in L^1(\R)$. By dominated convergence,
$$\limn \int^\infty_{-\infty}
f(t) (S(t) x_n - e^{-i\om t}x_n )\ud t =0.$$
Thus, since $\n x_n\n=1$,
$$
\n S(f)\n \ge \limn \n S(f)x_n\n
=\limn \Bigl\n\int^\infty_{-\infty} f(t) S(t) x_n\ud t\Bigr\n
=\Bigl|\int^\infty_{-\infty} e^{-i\om t}f(t)\ud t\Bigr| = \sqrt{2\pi} |\wh f(\om)|.$$
This inequality implies that $\wh f(\om) = 0$ for all $f\in I_S$.
Therefore, $\om\in {\rm Sp}(S)$.
\end{proof}

The {\it right translation group} is the $C_0$-group $U = (U(t))_{t\in\R}$ on $L^1(\R)$ defined by
$$ U(t)f(s) :=f(s-t),\quad s,t\in \R.$$
Note that $U(f)g = f * g$ for all $f, g\in L^1(\R)$, where $*$
denotes convolution.

We are now ready for the proof of Wiener's Tauberian theorem.

\begin{proof}[Proof of Theorem \ref{thm:Wiener}]
Let $f\in L^1(\R)$ be a function whose Fourier transform is zero-free and let
$X:=\overline{\hbox{span}\{U(t) f:\, t\in \R\}}$. We wish to prove that
$X=L^1(\R)$. Consider the quotient space $Y:= L^1(\R)/X$ and let $U_Y = (U_Y(t))_{t\in\R}$ denote the
associated quotient translation group on $Y$. This group is strongly
continuous and bounded. For all $g\in L^1(\R)$ we have
$U(f)g = f*g = g * f =  U(g)f.$
By the translation invariance
of $X$, $U(g)f\in X$. Hence
$U(f)g\in X$, so $U_Y(f)(g+X) = 0$ for all
$g\in L^1(\R)$. It follows that $U_Y(f)=0$. On the other hand, by assumption we have
$\wh f(\om)\not=0$ for all $\om\in\R$.
Therefore, ${\rm Sp}(U_Y)=\emptyset$.
We conclude that $Y=\{0\}$ and $X=L^1(\R)$.
\end{proof}

The Fourier argument used in the above proof also gives the following
stability result for semigroups.

\begin{theorem}[Strong stability]\label{thm:strong-stability}
Let $A$ generate a uniformly bounded $C_0$-semigroup $S$ on $X$. If $\sigma(A)\cap i\R=\emptyset$, then
$S$ is {\em strongly stable}\index{strongly stable}, that is,
$$
\lim_{t\to\infty}S(t)x=0,\quad x\in X.
$$
In fact, $\lim_{t\to\infty}\n S(t)A^{-1}\n=0$.
\end{theorem}

\begin{proof}
Put $M:=\sup_{t\ge0}\n S(t)\n$. For $f\in L^1(\R)$ and $t\ge0$,
define $T_f(t)\in\calL(X)$ by
$$
T_f(t)x:=\int_0^\infty f(r-t)S(r)x\ud r,\quad x\in X.
$$
Then $\n T_f(t)\n\le M\n f\n_1$ for all $t\ge0$.

First suppose that $\wh f$ has compact support. As in the proof
of Theorem \ref{thm:sp-bdd-group}, dominated convergence, Fourier
inversion, Fubini's theorem, and Proposition \ref{prop:resolv} give
\begin{align*}
T_f(t)x
&=\lim_{\delta\downarrow0}
  \int_0^\infty e^{-\delta r}f(r-t)S(r)x\ud r\\
  &=\frac1{\sqrt{2\pi}}\lim_{\delta\downarrow0}
  \int_0^\infty e^{-\delta r}
  \left(\int_\R e^{is(r-t)}\wh f(s)\ud s\right)S(r)x\ud r\\
&=\frac1{\sqrt{2\pi}}\lim_{\delta\downarrow0}
  \int_{\R}e^{-ist}\wh f(s)R(\delta-is,A)x\ud s
=\frac1{\sqrt{2\pi}}
  \int_{\R}e^{-ist}\wh f(s)R(-is,A)x\ud s.
\end{align*}
Here the last passage to the limit follows from the uniform
convergence of the resolvents on the compact support of $\wh f$.

The Riemann--Lebesgue lemma now implies that
$\n T_f(t)\n\to0$ as $t\to\infty$. To justify convergence in
operator norm, we approximate the continuous function
$s\mapsto R(-is,A)$ uniformly on $\supp(\wh f)$ by step functions
with values in $\calL(X)$, and apply the scalar
Riemann--Lebesgue lemma to each term.

As shown in the proof of Theorem \ref{thm:sp-bdd-group}, the
functions in $L^1(\R)$ whose Fourier transforms have compact
support are dense in $L^1(\R)$. Since
$$
\sup_{t\ge0}\n T_f(t)-T_g(t)\n\le M\n f-g\n_1,
$$
we conclude that $\n T_f(t)\n\to0$ for every $f\in L^1(\R)$.

Take $f_0(r):=e^{-r}$ for $r\ge0$ and $f_0(r):=0$ for $r<0$.
By the semigroup property and Proposition \ref{prop:resolv},
$$
T_{f_0}(t)x=\int_t^\infty e^{-(r-t)}S(r)x\ud r
=S(t)\int_0^\infty e^{-u}S(u)x\ud u = S(t)R(1,A)x.
$$
Thus $\n S(t)R(1,A)\n\to0$. Since
$\Ran(R(1,A))=\Dom(A)$ is dense in $X$ and $S$ is uniformly
bounded, this proves strong stability.

Finally, $0\in\varrho(A)$, and the identity
$A^{-1}=R(1,A)(A^{-1}-I)$
gives $\n S(t)A^{-1}\n\to0$.
\end{proof}

\section{The Hille--Yosida Theorem}\label{sec:HY}

The main theorem on generation of $C_0$-semigroups is the {\em Hille--Yosida theorem}, which gives necessary and sufficient
conditions in terms of resolvent growth. We only need the version for contraction semigroups, which is somewhat easier to state and prove. Its extension to general semigroups can be done via reductions similar to those used in the proof of Theorem \ref{thm:sgr-perturbation}, with the resolvent replacing the semigroup.

\begin{theorem}[Hille--Yosida]\label{thm:HilleYosida}\index{theorem!Hille--Yosida}
 For a densely defined closed linear operator $A$ in $X$ the following assertions are equivalent:
\begin{enumerate}[label={\rm(\arabic*)}, leftmargin=*]
 \item\label{it:HilleYosida1} $A$ generates a $C_0$-semigroup of contractions on $X$;
 \item\label{it:HilleYosida2} $\{\lambda\in\C:\, \Re\lambda>0\}\subseteq\varrho(A)$ and
$$ \n R(\la,A)\n \le \frac1{\Re\lambda}, \quad \Re\lambda>0;$$
 \item\label{it:HilleYosida3} $\{\lambda\in\R:\, \lambda>0\}\subseteq\varrho(A)$ and
$$ \n R(\la,A)\n \le \frac1{\lambda}, \quad \lambda>0.$$
\end{enumerate}
\end{theorem}
\begin{proof}
The implication \ref{it:HilleYosida1}$\Rightarrow$\ref{it:HilleYosida2} follows from Propositions \ref{prop:sg-prop} and
\ref{prop:resolv}, and the implication \ref{it:HilleYosida2}$\Rightarrow$\ref{it:HilleYosida3} is trivial.

Assume now that \ref{it:HilleYosida3} holds.
For the bounded operators $A_n := nA R(n,A) = n^2 R(n,A)-nI$, $n\ge 1$,  Proposition \ref{prop:resolv-conv} implies
$\limn A_n x = Ax$ for all $x\in \Dom(A)$. Also,
\begin{align}\label{eq:Yos-approx-contr} \n e^{tA_n}\n \le e^{n^2\n R(n,A)\n t }e^{-nt} \le e^{nt}e^{-nt} = 1.
\end{align}
Fix $x\in \Dom(A)$ and $t\ge 0$. The identity
\begin{align*} e^{tA_n}x - e^{tA_m} x
& = \int_0^t \frac{{\rm d}}{{\rm d}s} [e^{(t-s)A_m}e^{sA_n}x] \ud s
 = \int_0^t e^{(t-s)A_m}e^{sA_n}(A_n x-A_m x)\ud  s
\end{align*}
and the contractivity estimate \eqref{eq:Yos-approx-contr}
imply
$$ \n e^{tA_n}x - e^{tA_m} x \n \le t \n  A_n x-A_m x\n.$$
Thus, for every $T>0$ and $x\in\Dom(A)$,
the sequence $(e^{tA_n}x)_{n\ge 1}$ is Cauchy in
$C([0,T];X)$. By density of $\Dom(A)$ and contractivity,
the same holds for every $x\in X$.
Thus
$$
S(t)x:=\limn e^{tA_n}x
$$
exists, with convergence uniform on bounded time intervals.
Each $S(t)$ is linear and contractive, and every orbit
$t\mapsto S(t)x$ is continuous.

Clearly $S(0)=I$. Contractivity also allows us to pass
to the limit in the semigroup identity:
$$
S(t)S(s)x
=\limn e^{tA_n}e^{sA_n}x
=\limn e^{(t+s)A_n}x
=S(t+s)x.
$$
Hence $S$ is a $C_0$-semigroup of contractions.

Let $B$ denote its generator. For $x\in\Dom(A)$,
the convergence $A_nx\to Ax$ and the locally uniform
convergence just proved give
$$
S(t)x-x
=\limn\int_0^t e^{sA_n}A_nx\ud s
=\int_0^t S(s)Ax\ud s.
$$
Dividing by $t$ and letting $t\downarrow0$ shows that
$x\in\Dom(B)$ and $Bx=Ax$. Thus $A\subseteq B$.
Both resolvent sets contain $(0,\infty)$, by assumption
for $A$ and by Proposition \ref{prop:resolv} for $B$.
Proposition \ref{prop:semigroupsAB} therefore gives $A=B$.
\end{proof}

As an application, we obtain the following perturbation result.

\begin{theorem}[Perturbation]\label{thm:sgr-perturbation}\index{perturbation}
Let $A$ be the generator of a $C_0$-semigroup $S$ on $X$, and let $B$ be a bounded operator on $X$.
Then $A+B$ generates a $C_0$-semigroup on $X$.
\end{theorem}
Here it is understood that $\Dom(A+B)= \Dom(A)$ and $(A+B)x = Ax + Bx$ for $x\in \Dom(A)$.
The proof of the theorem shows that if  $\n S(t)\n \le Me^{\omega t}$,
then $\n S_B(t)\n \le Me^{(\omega+M\n B\n) t}$.
\begin{proof}
We prove the theorem in two steps.

\smallskip
{\em Step 1} -- We first prove the theorem for contraction semigroups.

Fix $\lambda\in\C$ with $\Re\lambda>0$. Then $\lambda\in\varrho(A)$ and $\n R(\la,A)\n \le (\Re\lambda)^{-1}$\! by the Hille--Yosida theorem. Because
$$(\lambda-(A+B))=(I-BR(\la,A))(\lambda-A)$$ and $\n BR(\la,A)\n \le (\Re\lambda)^{-1}\n B\n $,
for $\Re\lambda>\n B\n $ the operator $I-BR(\la,A)$ is
invertible, and the Neumann series for its inverse gives
\begin{align*}\n R(\la,A)\n \n (I-BR(\la,A))^{-1}\n & \le
(\Re\lambda)^{-1}(1-(\Re\lambda)^{-1}\n B\n )^{-1}  = (\Re\lambda-\n B\n)^{-1}\!.
\end{align*} Hence,
for $\Re\lambda> \n B\n$, the operator $\lambda-(A+B)$ is invertible and
its inverse satisfies $$\n R(\lambda,A+B)\n \le (\Re\lambda-\n
B\n)^{-1}\!.$$
The operator $\la - (A+B-\n B\n)$ is then invertible for $\Re\la>0$ and its inverse satisfies
$$\n R(\lambda, A+B-\n B\n)\n \le (\Re\lambda)^{-1}\!.$$
Thus, by the Hille--Yosida theorem, $A+B-\n B\n$ generates a $C_0$-semigroup
 $T$ of contractions. Then $A+B$  generates the $C_0$-semigroup
given by $S_B(t) := e^{t\n B\n}T(t)$.

\smallskip
{\em Step 2} -- In the general case, choose $M\ge1$ and
$\omega\in\R$ such that $\n S(t)\n\le Me^{\omega t}$.
The equivalent norm
$$
\n x\n_\omega:=\sup_{t\ge0}e^{-\omega t}\n S(t)x\n
$$
satisfies
$\n x\n\le\n x\n_\omega\le M\n x\n$.
With respect to this norm, the semigroup
$(e^{-\omega t}S(t))_{t\ge0}$ is contractive.
Writing $\n B\n_\omega$ for the corresponding operator
norm, we have
$$
\n Bx\n_\omega
\le M\n Bx\n
\le M\n B\n\n x\n_\omega,
$$
so $\n B\n_\omega\le M\n B\n$.

By Step 1, $A-\omega+B$ generates a $C_0$-semigroup $T$
satisfying $\n T(t)\n_\omega\le e^{t\n B\n_\omega}$.
Hence $A+B$ generates $S_B(t):=e^{\omega t}T(t)$, and
\[
\n S_B(t)x\n
\le e^{\omega t}\n T(t)x\n_\omega
\le e^{(\omega+\n B\n_\omega)t}\n x\n_\omega
\le M e^{(\omega+M\n B\n)t}\n x\n.
\qedhere\]
\end{proof}

In general there is no closed-form expression for $S_B(t)$, but we do have the so-called {\em variation of constants identity}
$$ S_B(t) x = S(t)x + \int_0^t S(t-s)BS_B(s)x\ud s.$$
To prove this identity, fix $x\in\Dom(A)=\Dom(A+B)$ and set $\phi(s):= S(t-s)S_B(s)x$. By Lemma \ref{lem:diff}, $\phi$ is differentiable and $$\phi'(s) =  -AS(t-s)S_B(s)x + S(t-s)(A+B)S_B(s)x = S(t-s)BS_B(s)x.$$
Integrating this identity over the interval $[0,t]$ gives the identity for $x\in\Dom(A)$. For fixed $t\ge0$, both
sides define bounded linear functions of $x\in X$. Since $\Dom(A)$ is
dense in $X$, the identity extends to every $x\in X$.

As a consequence of the variation of constants identity,
$$ \n S(t) - S_B(t)\n = O(t) \ \hbox{ as } \ t\downarrow 0.$$

We next prove an Euler formula for general $C_0$-semigroups.
It allows us to deduce properties of a semigroup from properties
of its resolvent along the positive real line. Later we will
use it to prove positivity of the heat semigroups with Dirichlet
and Neumann boundary conditions.

We recall Euler's formula for the exponential, which entails that for all $a\in \R$ and $t>0$,
$$ e^{ta} = \limn \Big(1-\frac{t}{n}a\Big)^{-n} =  \limn \Bigl(\frac{n}{t} (\frac{n}{t}- a)^{-1}\Bigr)^n\!.$$

\begin{theorem}[Euler formula, semigroup version]\label{thm:Euler}\index{theorem!Euler formula}
Let $A$ be the generator of a $C_0$-semigroup $S$ on  $X$. Then for all $x\in X$ and $t> 0$ we have
$$ S(t)x = \limn \Bigl(\frac{n}{t} R(\frac{n}{t},A)\Bigr)^n x.$$
\end{theorem}

\begin{proof}
Choose $M\geq1$ and $\om\in\R$ such that $\n S(s)\n\leq Me^{\om s}$ for all $s\geq0$.

\smallskip
{\em Step 1} -- We first derive an integral representation. By Proposition \ref{prop:spectr-holo}, the resolvent of $A$ is holomorphic, with complex derivative given by $\frac{{\rm d}}{{\rm d}\la}R(\la,A)=-R(\la,A)^2$. By induction, for $n\geq1$,
$$
\frac{{\rm d}^{\,n-1}}{{\rm d}\la^{\,n-1}}R(\la,A)
=(-1)^{n-1}(n-1)!R(\la,A)^n .
$$
On the other hand, for $\Re\la>\om$, repeated differentiation under the integral in the Laplace transform representation $R(\la,A)x=\int_0^\infty e^{-\la s}S(s)x\ud s$ gives
$$
\frac{{\rm d}^{\,n-1}}{{\rm d}\la^{\,n-1}}R(\la,A)x
=(-1)^{n-1}\int_0^\infty s^{n-1}e^{-\la s}S(s)x\ud s .
$$
Substituting $s=rt$ and putting $\la=n/t$, for $n$ so large that $n/t>\om$, we obtain
$$
\Bigl(\frac ntR(\frac nt,A)\Bigr)^nx
=\frac{n^n}{(n-1)!}\int_0^\infty r^{n-1}e^{-nr}S(rt)x\ud r .
$$

\smallskip
{\em Step 2} -- We prove that the right-hand side tends to $S(t)x$. Let $\nu_n$ be the probability measure on $[0,\infty)$ defined by
$$
\ud\nu_n(r):=\frac{n^n}{(n-1)!}r^{n-1}e^{-nr}\ud r .
$$
We shall use the elementary identities
$$
\int_0^\infty (r-1)^2\ud\nu_n(r)=\frac1n,
\qquad
\int_0^\infty e^{cr}\ud\nu_n(r)=\Bigl(1-\frac cn\Bigr)^{-n},\quad 0\leq c<n.
$$
Both follow from the formula $\int_0^\infty r^me^{-nr}\ud r=m!/n^{m+1}$.

Fix $x\in X$ and $\eps>0$. By strong continuity, choose
$0<\delta<1$ such that
$$
\n S(rt)x-S(t)x\n<\eps,\quad |r-1|<\delta.
$$
Put $E_\delta:=\{r\ge0:\, |r-1|\ge\delta\}$ and
$\beta:=\max\{\om t,0\}$. Then
$$
\nu_n(E_\delta)
\le \frac1{\delta^2}\int_{E_\delta}(r-1)^2\ud\nu_n(r)
\le \frac1{\delta^2}\int_0^\infty(r-1)^2\ud\nu_n(r)
=\frac1{n\delta^2}.
$$
By the Cauchy--Schwarz inequality, for $n>2\beta$ we have
\begin{align*}
\int_{E_\delta}e^{\beta r}\ud\nu_n(r)
&\le
\nu_n(E_\delta)^{1/2}
\Bigl(\int_0^\infty e^{2\beta r}\ud\nu_n(r)\Bigr)^{1/2}
\le
\frac1{\delta\sqrt n}
\Bigl(1-\frac{2\beta}{n}\Bigr)^{-n/2}
\to0.
\end{align*}
Since $\n S(rt)x\n\le Me^{\beta r}\n x\n$, splitting the
integral over $E_\delta$ and its complement gives
\begin{align*}
\Bigl\n\int_0^\infty (S(rt)x-S(t)x)\ud\nu_n(r)\Bigr\n
\le \eps+M\n x\n\int_{E_\delta}e^{\beta r}\ud\nu_n(r)+\n S(t)x\n\,\nu_n(E_\delta).
\end{align*}
The last two terms tend to zero. Since $\eps>0$ was arbitrary,
the integral representation from Step 1 therefore gives
\[
\limn\Bigl(\frac ntR(\frac nt,A)\Bigr)^nx=S(t)x.
\qedhere\]
\end{proof}

For sufficiently large $n$,
$\frac ntR(\frac nt,A)=(I-\frac tnA)^{-1}$.
Thus the theorem says that $n$ steps of the implicit Euler
scheme for $u'=Au$, with step size $t/n$ and initial value
$x$, converge to $S(t)x$.

We continue with a simple result about compact semigroups.

\begin{proposition}[Compact semigroups]\label{prop:compact-sgr}\index{compact!semigroup}
Let $A$ be the generator of a $C_0$-semi\-group $S$ on  $X$. If $S(t)$ is a compact operator for every $t>0$,
 then:
\begin{enumerate}[label={\rm(\arabic*)}, leftmargin=*]
 \item\label{it:compact-sgr1} the semigroup is norm continuous for $t>0$;
 \item\label{it:compact-sgr2} the resolvent operator $R(\la,A)$ is compact for every $\la\in\varrho(A)$;
 \item\label{it:compact-sgr3} the spectrum of $A$ is finite or countable and consists of isolated eigenvalues,
 and the corresponding eigenspaces are finite-dimensional;
 \item\label{it:compact-sgr4} for all $t>0$ we have the {\em spectral mapping formula}\index{spectral!mapping, compact semigroups}
 $$ \sigma(S(t))\setminus\{0\} = \exp(t\sigma(A)).$$
Moreover, the eigenspace corresponding to $\la\in \sigma(A)$ is contained in the eigenspace corresponding to $e^{\la t}\in\sigma(S(t))$.
\end{enumerate}
\end{proposition}
\begin{proof} \ref{it:compact-sgr1}: \
Fix $s>0$. For $t>s/2$ we have
\begin{align*}
\n S(t)-S(s)\n =  \sup_{\n x\n\le 1} \n(S(t-s/2)-S(s/2))S(s/2)x\n.
\end{align*}
Since $S(s/2)\ov B_X$ is relatively compact,
this implies
$\lim_{t\to s}\n S(t)-S(s)\n = 0$ by Proposition \ref{prop:uniform-limits-on-K}.

\smallskip
\ref{it:compact-sgr2}: \ Choose $M\ge 1$ and $\om\in\R$ such that $\n S(t)\n\le Me^{\om t}$ for all $t\ge 0$.
We first claim that the compactness of the semigroup operators $S(t)$ for $t>0$ implies that $R(\mu,A)$ is compact
for all $\Re\mu>\om$. For all $x\in X$ and $t>0$ we have
$$ S(t)R(\mu,A)x - R(\mu,A)x = \int_0^t S(s)AR(\mu,A)x\ud s,$$
and therefore $$\n S(t)R(\mu,A)-R(\mu,A)\n \le t \sup_{s\in [0,t]}\n S(s)\n \n AR(\mu,A)\n,$$
where $AR(\mu,A) = \mu R(\mu,A)-I$ is a bounded operator.
Since $S(t)R(\mu,A)$ is compact for every $t>0$,
Proposition \ref{prop:limcomp} implies that $R(\mu,A)$ is compact.
This proves the claim. The compactness of $R(\la,A)$ for arbitrary $\la\in\varrho(A)$ now follows from the resolvent identity \eqref{eq:res-id-unbdd}.

\smallskip
\ref{it:compact-sgr3}: \
Fix $\mu\in\varrho(A)$ and put $K:=R(\mu,A)$, which is compact by
part \ref{it:compact-sgr2}. For $\lambda\ne\mu$ we have
$ (\mu-\lambda)^{-1}I-K =(\mu-\lambda)^{-1}(\lambda-A)K$.
Since $K$ maps $X$ bijectively onto $\Dom(A)$, it follows that
$\lambda-A$ is injective, respectively surjective, if and only if
$(\mu-\lambda)^{-1}I-K$ is injective, respectively surjective.
Consequently,  $\Ker(\lambda-A) = \Ker((\mu-\lambda)^{-1}I-K)$ and
$$
\lambda\in\sigma(A)
\quad\Longleftrightarrow\quad
(\mu-\lambda)^{-1}\in\sigma(K).
$$

By the Riesz--Schauder theorem, the nonzero spectrum of $K$ is finite
or countable, consists of isolated eigenvalues, and has no possible
accumulation point other than $0$; the corresponding eigenspaces are
finite-dimensional. The preceding correspondence therefore implies
that $\sigma(A)$ is finite or countable, consists of isolated
eigenvalues with finite-dimensional eigenspaces, and, if infinite, can
accumulate only at infinity.

\smallskip
\ref{it:compact-sgr4}: \
Since $S(t)$ is compact, every nonzero element of $\sigma(S(t))$ is an eigenvalue
of $S(t)$ by the Riesz--Schauder theorem. By part \ref{it:compact-sgr3}, every element of $\sigma(A)$ is an eigenvalue of $A$. The spectral mapping formula now
follows from the next proposition.
\end{proof}

In the next proposition we denote by $\sigma_{\rm p}(B)$ the {\em point spectrum}\index{point!spectrum} of a bounded or unbounded operator $B$, that is, the set of its eigenvalues.

\begin{proposition}[Spectral mapping theorem for the point spectrum]\index{theorem!spectral mapping, point spectrum}
Let $A$ be the generator of a $C_0$-semi\-group $S$ on $X$. Then for all $t>0$ we have
$$\sigma_{\rm p}(S(t))\setminus \{0\} = \exp( {t\sigma_{\rm p}(A)}).$$
Moreover, for every $\la\in\sigma_{\rm p}(A)$,
$$
   \Ker(\la-A)\subseteq \Ker(e^{\la t}-S(t)).
$$
\end{proposition}

\begin{proof}
If $x\in \Dom(A)$ is an eigenvector of $A$ corresponding to the eigenvalue $\la$, the identity
$$
\int^t_0 e^{-\la s}S(s)(\la-A)x\ud s=x-e^{-\la t}S(t)x
$$
shows that $S(t)x=e^{\la t}x$. This proves the inclusion $\sigma_{\rm p}(S(t))\backslash\{0\}\supseteq\exp(t\sigma_{\rm p}(A))$ and the inclusion of the corresponding eigenspaces.

Conversely, fix $t>0$, and suppose that $x\in X$ is an eigenvector of $S(t)$ corresponding to a nonzero eigenvalue $\mu$. Choose $\la\in\C$ such that $\mu=e^{\la t}$. Then the map $s\mapsto e^{-\la s}S(s)x$ is periodic with period $t$. Since this map is not identically zero, the uniqueness theorem for the Fourier transform and the Hahn--Banach theorem imply that, after scaling the interval $[0,t]$ to $[0,2\pi]$, at least one of its Fourier coefficients is nonzero. Thus there exists an integer $k\in\Z$ such that, with $\la_k:=\la+2\pi i k/t$,
$$
x_k:=\frac1t\int^t_0 e^{-\la_k s}S(s)x\ud s\not=0.
$$
We will show that $x_k$ is an eigenvector of $A$ corresponding to $\la_k$.

By the observation preceding Proposition \ref{prop:sg-prop}, the operator $A-\la_k$ generates the $C_0$-semigroup $S_k$ given by $S_k(s):=e^{-\la_k s}S(s)$. Hence Proposition \ref{prop:sg-prop}, applied to $S_k$, gives $x_k\in\Dom(A)$ and
$$
(A-\la_k)x_k
 =\frac1t\bigl(e^{-\la_k t}S(t)x-x\bigr)
 =\frac1t\bigl(e^{-\la_k t}e^{\la t}x-x\bigr)=0,
$$
since $e^{-\la_k t}e^{\la t}=e^{-2\pi i k}=1$. Thus $\la_k\in\sigma_{\rm p}(A)$, and $\mu=e^{\la t}=e^{\la_k t}\in\exp(t\sigma_{\rm p}(A))$.
\end{proof}

\section{The Abstract Cauchy Problem}\label{sec:ACP}

Now that we have set up the general theory of $C_0$-semigroups, it is time to put them to use
in solving abstract Cauchy problems.

\subsection{The Inhomogeneous Cauchy Problem}\label{subsec:inhomACP}

If $A$ is the generator of a $C_0$-semigroup $S$ on $X$, then by Proposition \ref{prop:sg-prop} for initial values $u_0\in \Dom(A)$
the function
\begin{align}\label{eq:orbit} u(t):= S(t)u_0, \quad t\ge 0,
\end{align}
solves the initial value problem \eqref{ACP},
\begin{equation*}
\left\{\begin{aligned}
u'(t) &=Au(t),\quad t\in [0,T],\\
u(0) &=u_0,
\end{aligned}
\right.
\end{equation*}
in the sense that $u$ is continuously differentiable, takes values in $\Dom(A)$, and satisfies the equation
pointwise in time. A function $u$ with these properties is called a {\em classical solution}.\index{classical solution}
However, the definition \eqref{eq:orbit} makes sense for arbitrary $u_0\in X$, not just for
$u_0\in \Dom(A)$, and for all $u_0\in X$ the function
$u(t) = S(t)u_0$ solves the following integrated version of \eqref{ACP}:
\begin{align}\label{eq:integrated-sol} u(t) = u_0 + A\int_0^t u(s)\ud s, \quad t\in [0,T].
\end{align}
Indeed, by Proposition \ref{prop:sg-prop}\ref{it:sg-prop3}, for arbitrary $u_0\in X$ we have $\int_0^t S(s)u_0\ud s \in \Dom(A)$
and $A \int_0^t S(s)u_0\ud s = S(t)u_0-u_0$, confirming that \eqref{eq:integrated-sol} holds for $u(t) = S(t)u_0$.
This observation leads to the notion of
{\em strong solution} which we develop next in the more general context
of the inhomogeneous Cauchy problem\index{Cauchy problem!inhomogeneous}
\begin{equation*}\label{iACP}\tag{IACP}
\left\{\begin{aligned}
 u'(t) & = Au(t)+f(t), \quad t\in [0,T], \\
  u(0) & = u_0,
\end{aligned}
\right.
\end{equation*}
with initial value $u_0\in X$.
We assume that $f$ belongs to $L^1(0,T;X)$,\index{$L^1(0,T;X)$} the
space of all strongly measurable functions $f:(0,T)\to X$ such that $$\n f\n_1:= \int_0^T \n f(t)\n\ud t < \infty,$$
where functions that are equal almost everywhere are identified. As in the scalar-valued case, one shows that $L^1(0,T;X)$ is a Banach space.

\begin{definition}[Strong solutions]
A {\em strong solution}\index{strong!solution} of \eqref{iACP} is a continuous function $u:[0,T]\to X$ such that
for all $t\in [0,T]$ we have $\int_0^t u(s)\ud s \in \Dom(A)$ and
$$
u(t) = u_0 + A\int_0^t u(s)\ud s
+ \int_0^t f(s)\ud s.$$
\end{definition}

We proceed with an existence and uniqueness result for strong solutions
of \eqref{iACP}. It is based on the following lemma.

\begin{lemma}\label{lem:cont-mild}
Let $f\in L^1(0,T;X)$. Then:
\begin{enumerate}[label={\rm(\arabic*)}, leftmargin=*]
 \item\label{it:cont-mild1} for all $t\in [0,T]$ the function $s\mapsto S(t-s)f(s)$ has a strongly measurable representative and is integrable on $[0,t]$;
 \item\label{it:cont-mild2} the function $t\mapsto \int_0^t S(t-s)f(s)\ud s$ is continuous on $[0,T]$.
\end{enumerate}
\end{lemma}
\begin{proof}
\ref{it:cont-mild1}: \ Choose a strongly measurable representative for $f$, which we denote again by $f$, as well as a sequence of simple functions $f_n$ converging to $f$ pointwise.
Each function $s\mapsto S(t-s)f_n(s)$
is strongly measurable, since it is a linear combination of functions of the form $s\mapsto \one_B(s)S(t-s)x$ with $B\subseteq [0,T]$ a Borel subset,
and such functions are strongly measurable because continuous functions on an interval are strongly measurable and
if $g$ is strongly measurable and $B$ is a Borel set, then $\one_B g$ is strongly measurable.
By Corollary \ref{cor:ptw-lim-strmeas}, the pointwise limit $s\mapsto S(t-s)f(s)$ is strongly measurable.
Integrability follows from the estimate $\n S(t-s)f(s)\n \le M \n f(s)\n$, where $M = \sup_{t\in [0,T]}\n S(t)\n$, and the integrability of $f$.

\smallskip
\ref{it:cont-mild2}: \ Let $0\le t\le t'\le T$. Then
\begin{align*} \ &  \Big\n \int_0^{t'} S(t'-s)f(s)\ud s - \int_0^t S(t-s)f(s)\ud s\Big\n
 \\* & \qquad \le  \Big\n \int_t^{t'} S(t'-s)f(s)\ud s \Big\n
+  \Big\n \int_0^{t} S(t'-s)f(s)\ud s - \int_0^t S(t-s)f(s)\ud s\Big\n.
\end{align*}
The first term on the right-hand side can be bounded above by
$M \int_t^{t'}\n f(s)\n \ud s$, which tends to $0$ by dominated convergence as  $t'-t\to 0$.
The second term tends to $0$ by dominated convergence as well: for simple functions $f$ this follows from the strong continuity and local boundedness of
the semigroup, and for general $f\in L^1(0,T;X)$ this follows by approximation by simple functions in the $L^1(0,T;X)$-norm.
\end{proof}

\begin{theorem}[Existence and uniqueness]\label{thm:iACP}\index{existence and uniqueness!for inhomogeneous problems}
For all $u_0\in X$ and $f\in L^1(0,T;X)$ the problem \eqref{iACP} admits a unique
strong solution $u$. It is given by the convolution formula
\begin{equation}\label{eq:ACPconvol} u(t) = S(t)u_0 + \int_0^t S(t-s)f(s)\ud s.
\end{equation}
\end{theorem}
The function $ t\mapsto S(t)u_0 + \int_0^t S(t-s)f(s)\ud s$
is called the {\em mild solution} of \eqref{iACP}.\index{mild solution}
\begin{proof}
Let $u$ be given by \eqref{eq:ACPconvol}. By Lemma
\ref{lem:cont-mild}, this function is continuous.
Fix $0<t\le T$ and put
$$
F_t(r):=\int_0^{t-r}S(s)f(r)\ud s,\quad 0<r<t.
$$
By Proposition \ref{prop:sg-prop}, for almost all $r\in(0,t)$
we have $F_t(r)\in\Dom(A)$ and
$$
AF_t(r)=S(t-r)f(r)-f(r).
$$
Writing $M:=\sup_{0\le s\le T}\n S(s)\n$, we obtain
$$
\n F_t(r)\n_{\Dom(A)}
\le (tM+M+1)\n f(r)\n.
$$
Moreover, $F_t$ is strongly measurable with values in $\Dom(A)$ with its graph norm. Indeed, for simple functions $f$
this follows from strong continuity, and the general case follows by pointwise
approximation by simple functions.
The estimate therefore shows that $F_t$ is Bochner integrable as a
$\Dom(A)$-valued function.

By Fubini's theorem,
$$
\int_0^t u(s)\ud s
=
\int_0^t S(s)u_0\ud s+\int_0^t F_t(r)\ud r
\in\Dom(A).
$$
Since $A:\Dom(A)\to X$ is bounded for the graph norm, we may
apply it under the second integral. Proposition
\ref{prop:sg-prop} then gives
\begin{align*}
A\int_0^t u(s)\ud s
&=
S(t)u_0-u_0
+\int_0^t\bigl(S(t-r)f(r)-f(r)\bigr)\ud r =
u(t)-u_0-\int_0^t f(r)\ud r.
\end{align*}
The same identity holds trivially at $t=0$.
This shows that $u$ is a strong solution.

To prove uniqueness, suppose that $u$ and $\wt u$ are strong solutions of
\eqref{iACP}. It follows from the definition that $u$ and $\wt u$ are continuous.
Set $v:= u-\wt u$. Then $v$ is continuous, $\int_0^t v(s)\ud s \in \Dom(A)$, and
$ v(t) = A\int_0^t v(s)\ud s$
for all $t\in [0,T]$.

Fix $0\le t\le T$ and define $w:[0,t]\to X$ by $w(s) = S(t-s)\int_0^s v(r)\ud r$.
By Lemma \ref{lem:diff}, this function is differentiable with derivative
$$ w'(s) = S(t-s)v(s)- S(t-s)A\int_0^s v(r)\ud r
= S(t-s)v(s)- S(t-s)v(s) = 0.$$
It follows that $w$ is constant.
Hence
$\int_0^t v(r)\ud r  = w(t) = w(0) = 0.$
Since this is true for all $0\le t\le T$ and $v$ is continuous, it follows that $v=0$ on $[0,T]$.
\end{proof}

The solution $u$ depends continuously on $u_0$ in the norm of $C([0,T];X)$, the Banach space of all continuous functions from $[0,T]$ to $X$.
Indeed, if $\wt u_0$ is another initial value and the corresponding unique strong solution is denoted by $\wt u$, then
\begin{align*}
 \n u(t) - \wt u(t)\n
 \le \n S(t)\n \n u_0-\wt u_0\n
  \le M_T \n u_0-\wt u_0\n,
\end{align*}
where $M_T := \sup_{t\in [0,T]}\n S(t)\n$,
and therefore
$  \n u - \wt u\n _\infty \le M_T \n u_0-\wt u_0\n.$

Unique solvability together with continuous dependence on the initial
value is usually summarised as {\em well-posedness}.\index{well posed}
Thus, the inhomogeneous problem \eqref{iACP} is well posed for strong solutions.

\subsection{The Semilinear Cauchy Problem}\label{subsec:semilinACP}

In this section we study a class of nonlinear evolution equations of the form
\begin{equation*}\label{eq:SDE} \tag{SCP}
\left\{
\begin{aligned} u'(t) & = Au(t) + f(t,u(t)), \quad t\in [0,T], \\
      u(0) & = u_0.
\end{aligned}
\right.
\end{equation*}
Equations of this form are referred to as {\em semilinear}\index{Cauchy problem!semilinear} equations.
We assume that $A$ generates a $C_0$-semigroup $S$ on  $X$
and that the initial value $u_0$ lies in $X$.
We make the following assumptions on the function $f:[0,T]\times X\to X$:
\begin{enumerate}[label={\rm(\roman*)}, leftmargin=*]
 \item\label{it:eeq-standing1} (Strong measurability) for all $x\in X$ the function $t\mapsto f(t,x)$ is strongly measurable on $[0,T]$;
 \item\label{it:eeq-standing2} (Linear growth) there exists a constant $C\ge 0$ such that
 $$ \n f(t,x)\n
 \le C(1+\n x\n), \quad t\in [0,T], \ x\in X;$$
 \item\label{it:eeq-standing3} (Lipschitz continuity)
 there exists a constant $L\ge 0$ such that
 $$ \n f(t,x)-f(t,x')\n \le L\n x-x'\n, \quad t\in [0,T], \ x,x'\in X.$$
\end{enumerate}
Under these assumptions, in force throughout this section, we will prove existence, uniqueness, and continuous dependence on the initial conditions of mild solutions. Thus \eqref{eq:SDE} is well posed for mild solutions.

\begin{definition}[Mild solutions]\label{def:mild-semilinear}
A function $u:[0,T]\to X$ is called a {\em mild solution}\index{mild solution} of \eqref{eq:SDE}
if it is continuous and satisfies
\begin{equation*}
u(t) = S(t)u_0 + \int_0^t S(t-s)f(s,u(s))\ud s, \quad t\in [0,T].
\end{equation*}
\end{definition}

To see that this is well defined we must check that the integral converges as a Bochner integral in $X$.
Lemma \ref{lem:mild-welldef} takes care of this. Taking the lemma for granted for the moment, let us first motivate the definition.

First of all, it generalises the formula of Theorem \ref{thm:iACP} for the strong solution
of the inhomogeneous problem.
Perhaps more importantly, every classical solution is a mild solution.
Here, by a {\em classical solution}\index{classical solution} we mean a function
$$
u\in C^1([0,T];X)\cap C([0,T];\Dom(A)),
$$
which satisfies \eqref{eq:SDE} pointwise. Here, $\Dom(A)$ is treated as a Banach space endowed with the graph norm.

Suppose now that $u$ is a classical solution, and fix $t\in[0,T]$. By Lemma \ref{lem:diff}, the function
$s\mapsto S(t-s)u(s)$
is continuously differentiable on $[0,t]$, with
\begin{align*}
\frac{\ud}{\ud s}\bigl(S(t-s)u(s)\bigr)
&=
-S(t-s)Au(s)+S(t-s)u'(s)
=
S(t-s)f(s,u(s)).
\end{align*}
Integrating from $0$ to $t$ gives
$$
u(t)-S(t)u_0
=
\int_0^tS(t-s)f(s,u(s))\ud s.
$$
Thus every classical solution is a mild solution.

In general there is no reason to expect the existence of classical solutions, but, under the standing assumptions \ref{it:eeq-standing1}--\ref{it:eeq-standing3} formulated above, a unique mild solution always exists. Unlike a classical solution, a mild solution need not be differentiable or take values in $\Dom(A)$.

As promised, we now check that the integral in Definition \ref{def:mild-semilinear}
is well defined as a Bochner integral in $X$. The next result extends Lemma \ref{lem:cont-mild} to the present situation.

\begin{lemma}\label{lem:mild-welldef} Let $f:[0,T]\times X\to X$ satisfy the conditions \ref{it:eeq-standing1}--\ref{it:eeq-standing3} and suppose that $u:[0,T]\to X$ is continuous. Then:
\begin{enumerate}[label={\rm(\arabic*)}, leftmargin=*]
\item\label{it:mild-welldef1} for every $t\in[0,T]$, the functions $s\mapsto f(s,u(s))$ on $[0,T]$ and $s\mapsto S(t-s)f(s,u(s))$ on $[0,t]$ have strongly measurable representatives and are integrable;
\item\label{it:mild-welldef2} the function $t\mapsto \int_0^t S(t-s)f(s,u(s))\ud s$ is continuous on $[0,T]$.
\end{enumerate}
\end{lemma}
\begin{proof}
Choose step functions $u_n:[0,T]\to X$ converging uniformly
to $u$. By assumption \ref{it:eeq-standing1}, each function
$s\mapsto f(s,u_n(s))$ is strongly measurable. Moreover,
Lipschitz continuity gives
$$
\n f(s,u_n(s))-f(s,u(s))\n
\le L\n u_n-u\n_\infty\to0,\quad s\in[0,T].
$$
Thus $g(s):=f(s,u(s))$ is strongly measurable by Corollary
\ref{cor:ptw-lim-strmeas}. The linear growth assumption gives
$$
\n g(s)\n\le C(1+\n u\n_\infty),\quad s\in[0,T],
$$
so $g\in L^1(0,T;X)$. Both assertions now follow from
Lemma \ref{lem:cont-mild}, applied to $g$.
\end{proof}

We are now ready to state and prove our main result:

\begin{theorem}[Well-posedness of the semilinear problem]\label{thm:mainexistenceL}\index{well posed}\index{existence and uniqueness!for semilinear problems}
Under the assumptions \ref{it:eeq-standing1}--\ref{it:eeq-standing3} formulated at the beginning of the section, the semilinear problem \eqref{eq:SDE} admits a unique mild solution $u\in C([0,T];X)$.
This solution depends continuously, in the norm of $C([0,T];X)$, on the initial condition $u_0\in X$.
\end{theorem}

\begin{proof}
To obtain existence and uniqueness we define a nonlinear mapping $\Phi$ from $C([0,T];X)$ to itself
by
$$ (\Phi(v))(t) := S(t)u_0 + \int_0^t S(t-s)f(s,v(s))\ud s, \quad t\in [0,T], \ v\in C([0,T];X).$$
We have already observed in Lemma \ref{lem:mild-welldef} that the integrand is integrable, and the continuity of $\Phi(v)$ follows
from the strong continuity of the semigroup and Lemma \ref{lem:cont-mild}. It follows that
$\Phi$ is well defined as a mapping of $C([0,T];X)$ into itself. We now reuse the idea from the proof of Theorem \ref{thm:DV-Lip}
and set, for $\la>0$ to be chosen in a moment,
$$\n g\n_\lambda := \sup_{t\in [0,T]} e^{-\lambda t} \n g(t)\n.$$
This defines an equivalent norm on $C([0,T];X)$.
By the Lipschitz continuity assumption, for all $v,w\in C([0,T];X)$ and $t\in [0,T]$ we have
\begin{align*}
\n (\Phi(v))(t) - (\Phi(w))(t)\n
& \le  \int_0^{t} \n S(t-s)(f(s,v(s))-f(s,w(s)))\n \ud s
\\ & \le LM \int_0^{t} e^{\la s}e^{-\la s}\n v(s)-w(s)\n \ud s
\\ & \le LM   \int_0^{t} e^{\la s}\n v-w\n_\lambda\ud s
 = \frac{LM}{\la} (e^{\la t}-1) \n v-w\n_\lambda ,
\end{align*}
where $M = \sup_{t\in [0,T]}\n S(t)\n$.
It follows that $$
\n \Phi(v) - \Phi(w)\n_\la  \le \frac{LM}{\la}\n v-w\n_\la.$$
If we choose $\la>LM$, the mapping $\Phi$ is a uniform contraction on $C([0,T];X)$ and therefore
has a unique fixed point $u\in C([0,T];X)$ by the Banach fixed point theorem (Theorem \ref{thm:fixed-point-banach}).\index{fixed point!argument} Then,
$$ u(t) = (\Phi(u))(t) = S(t)u_0 + \int_0^t S(t-s)f(s,u(s))\ud s, \quad t\in [0,T],$$
so $u$ is a mild solution. Conversely, any mild solution is a fixed point of $\Phi$, and since $\Phi$ has
a unique fixed point, the mild solution $u$ is unique.

To complete the proof we check the continuous dependence of the mild solution on the initial value $u_0$. If $\wt u_0$ is another
initial value and the corresponding unique mild solution is denoted by $\wt u$, estimating as before, we obtain
\begin{align*}
 \n u(t) - \wt u(t)\n
 & \le \n S(t)\n \n u_0-\wt u_0\n + \int_0^t \n S(t-s)(f(s,u(s))- f(s,\wt u(s)))\n \ud s
 \\ & \le M \n u_0-\wt u_0\n  + \frac{LM}{\la} (e^{\la t}-1) \n u-\wt u\n_\la
\end{align*}
and therefore
$$  \n u - \wt u\n _\la \le  M \n u_0-\wt u_0\n + \frac{LM}{\la} \n u-\wt u\n_\la.$$
Choosing $\la = \max\{1,2LM\}$ gives
$ \frac12 \n u - \wt u\n _{\max\{1,2LM\}} \le  M \n u_0-\wt u_0\n$,
and the desired continuity follows, as $\n \cdot\n_{\max\{1,2LM\}}$ is an equivalent norm on $C([0,T];X)$.
\end{proof}

For this result to be useful, one needs a way of translating
nonlinearities arising in concrete partial differential equations into
our abstract framework. We illustrate this with an example.
Consider the following semilinear heat equation on a bounded open subset $D$ of $\R^d$:
\index{heat!equation, semilinear}
\begin{equation*}
\left\{\begin{aligned}
\frac{\partial u}{\partial t}(t,\xi) & = \Delta u(t,\xi)
+  b(u(t,\xi)),
&& \xi\in D, && t\in
[0,T],
\\  u(t,\xi) &= 0, && \xi\in \partial D, && t\in [0,T],
\\ u(0,\xi) &= u_0(\xi), && \xi\in D.
\end{aligned}
\right.
\end{equation*}
We work over the real scalars and assume that  $b:\R\to \R$ is Lipschitz continuous, with Lipschitz constant $L$:
$$ |b(\xi) - b(\xi')| \le L|\xi-\xi'|, \quad \xi,\xi'\in \R.$$ The assumption that $b$ depends only on the solution
is made for simplicity; the more general case where $b$ also depends on time can be treated in the same way.

To cast this problem into a semilinear abstract Cauchy problem we assume that the initial value $u_0$ belongs to $L^2(D)$.
The above problem may then be written in the form
\begin{equation*}
\left\{
\begin{aligned}
u'(t) & = A u(t) + B(u(t)), \\
 u(0) & = u_0,
\end{aligned}
\right.
\end{equation*}
where $A$ is the Dirichlet Laplacian on $L^2(D)$. Since $-A$ is
positive and selfadjoint, we have
$$
(0,\infty)\subseteq\varrho(A),
\qquad
\n R(\la,A)\n\le\frac1\la,\quad \la>0.
$$
The Hille--Yosida theorem therefore shows that $A$ generates a
contractive $C_0$-semigroup on $L^2(D)$.
The mapping $B:L^2(D)\to L^2(D)$ is the {\em Nemytskii mapping}\index{Nemytskii mapping} associated
with $b$: $$(B(x))(\xi):= b(x(\xi)), \quad x\in L^2(D).$$
The next proposition checks that this mapping is well defined, of
linear growth, and Lipschitz continuous on $L^2(D)$ (and, with the same proof, on
$L^p(D)$ with $1\le p<\infty$).

\begin{proposition}\label{prop:Nemytskii}
Under the above assumptions on $b$, the
Nemytskii mapping
 $B: L^2(D)\to L^2(D)$
is well defined, of linear growth, and Lipschitz continuous (in the sense that $f(t,x):= B(x)$ satisfies conditions \ref{it:eeq-standing2} and \ref{it:eeq-standing3}
at the beginning of this section).
\end{proposition}
\begin{proof}
Let us first check that $B(x)\in L^2(D)$ for all $x\in L^2(D)$.
Using the triangle inequality in $L^2(D)$, for all $x\in L^2(D)$ we have
\begin{align*} \| B(x)\|_{L^2(D)} & =\Bigl(\int_D |b(x(\xi))|^2\ud \xi\Bigr)^{1/2}
\\ & \le \Bigl(\int_D |b(x(\xi))-b(0)|^2\ud \xi\Bigr)^{1/2}+ \Bigl(\int_D |b(0)|^2\ud \xi\Bigr)^{1/2}
\\ & \le L \Bigl(\int_D |x(\xi)-0|^2\ud \xi\Bigr)^{1/2} + |b(0)|\Bigl(\int_D 1\ud \xi\Bigr)^{1/2}
 = L\n x\n_2 + |D|^{1/2}|b(0)|,
\end{align*}
where $|D|$ stands for the Lebesgue measure of $D$.
This proves that $B$ is well defined and of linear growth.
Lipschitz continuity follows by a similar estimate. For all $x,y\in L^2(D)$,
\begin{align*}
 \|B(x) - B(y)\|_{L^2(D)}
 & = \Bigl(\int_D |b(x(\xi)) - b(y(\xi))|^2\ud \xi\Bigr)^{1/2}
\\ & \le L \Bigl(\int_D |x(\xi) - y(\xi)|^2\ud \xi\Bigr)^{1/2}
= L \n x-y\n_{L^2(D)}.\qedhere
\end{align*}
\end{proof}

We have thus shown that all assumptions of Theorem \ref{thm:mainexistenceL} are fulfilled. It follows that the corresponding semilinear abstract Cauchy problem
admits a unique mild solution.

\section{Analytic Semigroups}\label{sec:analytic-semigroups}

Analytic semigroups provide an abstract framework for discussing a class of initial value problems, referred to in the partial differential equations literature as {\em parabolic}. An important characteristic of this class of problems is that solutions are smooth.

\subsection{The Main Result}\label{subsec:analytic-main}

For $\om\in (0,\pi)$ consider the open sector \index{sector}\index{$S$@$\Sigma_\omega$}
$$\Sigma_\om := \{z\in \C\setminus\{0\}: \, |\arg(z)|<\om\},$$
where the argument is taken in $(-\pi,\pi)$.

\begin{definition}[Analytic $C_0$-semigroups] \label{def:analytic-semigroup}
A $C_0$-semigroup $S$ on  $X$ is called an
{\em analytic $C_0$-semigroup}\index{analytic semigroup}
on $\Sigma_\om$ if it extends to a holomorphic mapping $S:\Sigma_\om\to\calL(X)$ such that,
for every $x\in X$,
$$ \lim_{\substack{z\in \Sigma_\om,\\ z\to 0}} S(z)x = x.$$
We call $S$ an {\em analytic $C_0$-semigroup} if it is analytic on $\Sigma_\om$
for some $\om\in (0,\pi)$.
\end{definition}

If $S$ is an analytic $C_0$-semigroup on $\Sigma_\om$, then for all $z_1,z_2\in \Sigma_\om$ with $z_1+z_2 \in \Sigma_\omega$ we
have $$S(z_1)S(z_2) = S(z_1+z_2).$$ Indeed, for each $x\in X$ the functions
 $z_1\mapsto S(z_1)S(t)x$ and $S(z_1+t)x$ are
holomorphic extensions of $s\mapsto S(s+t)x$
and are therefore equal. Repeating this argument, we find that the functions
$z_2\mapsto S(z_1)S(z_2)x$ and $S(z_1+z_2)x$
are equal.

As in the proof of Proposition \ref{prop:boundS}, the uniform boundedness theorem
implies that if $S$ is an analytic $C_0$-semigroup on $\Sigma_\om$,
then the operators $S(z)$ are uniformly bounded on
$\Sigma_{\om'}\cap B(0;r)$ for every $0<\om'<\om$ and
$r\ge 0$. The same argument as in Proposition \ref{prop:boundS} gives exponential boundedness
on $\Sigma_{\om'}$ for all $0<\om'<\om$, in the sense that there are constants $M'\ge 1$ and $c' = c_{\om'}\in\R$
such that $$\n S(z)\n\le M'e^{c'|z|}, \quad z\in \Sigma_{\om'}.$$
We say that $S$ is a
{\em bounded analytic $C_0$-semigroup on $\Sigma_\om$}\index{bounded!analytic semigroup}
if $S$ is an analytic $C_0$-semigroup on $\Sigma_\om$ and the operators $S(z)$ are uniformly bounded on $\Sigma_\om$.
{\em Contractive analytic $C_0$-semigroups on $\Sigma_\om$}\index{contractive analytic semigroup}
are defined similarly.
There is a rather subtle point here: boundedness and contractivity are imposed on a sector,
not just on the positive real line.
That this makes a difference is shown by the simple example of the rotation group on $\C^2$\!,
 $$ S(t) \ = \
\begin{pmatrix}
\cos t & -\sin t  \\
\sin t & \cos t
\end{pmatrix}.
$$
 For each $t\in \R$ we have $\n S(t)\n =1$. Upon replacing $t$ by a complex parameter
$z$ the group extends holomorphically to the entire complex plane, but it is unbounded on every sector $\Sigma_\om$
with $0<\om<\pi$.

\begin{theorem}[Bounded analytic semigroups, complex characterisation]\label{thm:analytic}
For a dense\-ly defined closed operator $A$ in $X$ the following assertions are
equivalent:
\begin{enumerate}[label={\rm(\arabic*)}, leftmargin=*]
\item\label{it:analytic-semigroup1}
$A$ generates a bounded analytic $C_0$-semigroup on $\Sigma_\eta$ for some $\eta\in (0,\frac12\pi)$;
\item\label{it:analytic-semigroup2}
there exists $\theta\in (\frac12\pi,\pi)$ such that
$\dps\Sigma_\theta\subseteq\varrho(A)$ and
$$\sup_{\la\in \Sigma_\theta} \n \la R(\la,A)\n <\infty;$$
\item\label{it:analytic-semigroup3}
we have $\C_+ \subseteq \varrho(A)$ and
$$\sup_{\la\in \C_+} \n \la R(\la,A)\n <\infty.$$
\end{enumerate}
Denoting the suprema of all admissible $\eta$ and $\theta$ in \ref{it:analytic-semigroup1}
and \ref{it:analytic-semigroup2} by $\om_{\rm\, holo}(A)$ and $\om_{\rm\, res}(A)$, respectively, we have
$$ \om_{\rm\, res}(A) = \frac12\pi + \om_{\rm\, holo}(A).$$
Under the equivalent conditions \ref{it:analytic-semigroup1} and \ref{it:analytic-semigroup2} we have the inverse Laplace transform\index{Laplace transform!inverse} representation
\begin{equation*}
S(t)x = \frac1{2\pi i}\int_\Gamma e^{\la t}R(\la,A)x\ud \la, \quad t>0, \ x\in X,
\end{equation*}
where $\Gamma = \Gamma_{\!\theta'\!,B}$ is the upwards oriented boundary of $\Sigma_{\theta'}\setminus B$, for any
$\theta'\in (\frac12\pi,\theta)$ and any closed ball $B$ of positive radius centred at the origin.
\end{theorem}

Note that \ref{it:analytic-semigroup2} implies $\si(A)\cap i\R \subseteq\{0\}$.

\begin{proof}
By Cauchy's theorem, if the integral representation holds for some
$\theta'\in (\frac12\pi,\theta)$ and some closed ball $B$ centred at the origin,
then it holds for any such $\theta'$ and $B$.

\smallskip
\ref{it:analytic-semigroup1}$\Rightarrow$\ref{it:analytic-semigroup2}: \ We start with the preliminary observation that if a linear operator $\wt A$ generates a uniformly bounded $C_0$-semigroup $\wt S$
on $X$, then, by
Proposition \ref{prop:resolv}, the open right half-plane $\C_+=\{\Re\la>0\}$
is contained in the resolvent set of $\wt A$ and we have the bound  $\n R(\la,\wt A)\n\le {M}/{\Re\la}$
for all $\la\in \C_+.$ Moreover, for all $\theta\in (0,\frac12\pi)$ and $\la\in\Sigma_\theta$ we have $\Re \la \ge |\la|\cos(\theta)$ and therefore
$$\sup_{\lambda\in \Sigma_{\theta}} \n \la R(\la,\wt A)\n  \le \frac{M}{\cos\theta}.$$

Now if $A$ generates a $C_0$-semigroup which is bounded on a sector
$\Sigma_\eta$ with $\eta\in(0,\frac12\pi)$, say by a constant $M$, we
can apply the above reasoning to the bounded semigroups
$S(e^{i\eta'}t)$, with $|\eta'|<\eta$.
Note that the generator of $t\mapsto S(e^{i\eta'}t)$ is
$e^{i\eta'}A$.
Letting $\eta'$ vary in this
symmetric interval gives the resolvent estimate on the corresponding
symmetric sector. Optimising the various choices of angles, we obtain
$$
  \omega_{\rm\, res}(A)\ge \frac12\pi+\omega_{\rm\, holo}(A).
$$

\smallskip
\ref{it:analytic-semigroup2}$\Rightarrow$\ref{it:analytic-semigroup1}: \ The idea is to
define the semigroup operators by the integral representation given in the statement of the theorem,
and prove that they define a bounded $C_0$-semigroup which has the properties stated in part \ref{it:analytic-semigroup1}.
The construction gives \ref{it:analytic-semigroup1} for every $\eta<\theta-\frac12\pi$, and hence the inequality $$\om_{\rm\, holo}(A)\ge\om_{\rm\, res}(A)-\frac12\pi.$$

Let $\eta:= \theta - \frac12\pi$. For any $\zeta\in \Sigma_\eta$ let
$$ S(\zeta)x := \frac1{2\pi i}\int_\Gamma e^{\la \zeta}R(\la,A)x\ud \la, \quad  x\in X, $$
where $\Gamma$ is the boundary of $\Sigma_{\theta'}\setminus B$ with $\theta'\in
(\frac12\pi, \theta)$ any number such that $\frac12\pi +
|\arg(\zeta)| < \theta'< \theta$
and $B$ is any closed ball centred at the origin;
see Figure \ref{fig:Gamma}.
This integral converges absolutely, defines a bounded operator $S(\zeta)$ on $X$,
and the function $\zeta\mapsto S(\zeta)x$ is holomorphic on $\Sigma_\eta$.

\begin{figure}
\begin{center}
\begin{tikzpicture}[scale=2]
\draw (-1.73,0) -- (1.73,0);
\draw (0,-1.73) -- (0,1.73);

\draw[thick] (1,0) arc (0:120:1cm);
\draw[thick] (1,0) arc (0:-120:1cm);

\draw[->,thick] (-0.5,0.866) -- (-0.75,1.298);
\draw[thick] (-0.5,-0.866) -- (-0.75,-1.298);

\draw[thick] (-0.75,1.298) -- (-1,1.73);
\draw[->,thick] (-1,-1.73) -- (-0.75,-1.298);

\draw[dashed] (0,0) -- (-0.5,0.866);
\draw[dashed] (0,0) -- (-0.5,-0.866);

\draw (0,0) -- (-1.285,1.53);
\draw (0,0) -- (-1.285,-1.53);

\draw[->,thin] (0.3,0) arc (0:120:0.3cm);
\draw[->,thin] (0.4,0) arc (0:130:0.4cm);

\draw (0.11,0.11) node [fill=white]{$\theta'$};
\draw (0.37,0.37) node [fill=white]{$\theta$};

\draw (1,-0.55) node [fill=white]{$\Gamma$};
\end{tikzpicture}
\caption{The contour $\Gamma$\label{fig:Gamma}}
\end{center}
\end{figure}

The proof of the semigroup property proceeds much as the proof of the
multiplicativity of the holomorphic calculus. Fix $\zeta,\zeta'\in \Sigma_{\eta}$ and choose contours $\Gamma$ and $\Gamma'$ as above, with
$\Gamma$ to the left of $\Gamma'$\!.
Then, by the resolvent identity \eqref{eq:res-id-unbdd}, Cauchy's theorem, Fubini's theorem, and the Cauchy integral formula,
\begin{align*}
 S(\zeta')S(\zeta)x
   & =  \frac1{(2\pi i)^2}\int_{\Gamma}\int_{\Gamma'} e^{\la\zeta+\mu\zeta'}R(\la,A)R(\mu,A)x\ud \mu\ud \la
\\ & =  \frac1{(2\pi i)^2}\int_{\Gamma}\int_{\Gamma'} e^{\la\zeta+\mu\zeta'}\frac{R(\la,A)x-R(\mu,A)x}{\mu-\la}\ud \mu\ud \la
\\ & =  \frac1{(2\pi i)^2}\int_{\Gamma}\int_{\Gamma'} e^{\la\zeta+\mu\zeta'}\frac{R(\la,A)x}{\mu-\la}\ud \mu\ud \la
\\ & =  \frac1{2\pi i}\int_{\Gamma} e^{\la\zeta+\la\zeta'} R(\la,A)x\ud\la = S(\zeta+\zeta')x.
\end{align*}

Put $M:= \sup_{\la\in \Sigma_\theta}\|\la R(\la,A)\|$. Fix
$0<\eta_0<\eta$. We will prove that $S$ is uniformly bounded on $\Sigma_{\eta_0}$.
Choose $\theta'$ such that $\frac12\pi+\eta_0<\theta'<\theta$. For $\zeta\in\Sigma_{\eta_0}$, consider the contour
$\Gamma=\Gamma_{\!\theta'\!,B_r}$, where $B_r:=\ov B(0;r)$. Along the rays
$\la=\rho e^{\pm i\theta'}$, $\rho\ge r$, we have
$$
  \Re(\la\zeta)\le -\delta \rho|\zeta|,
  \qquad
  \delta:=-\cos(\theta'-\eta_0)>0 .
$$
Thus the arc contributes at most
$\frac{\theta'M}{\pi}e^{r|\zeta|}$,
and the two rays together contribute at most
$\frac{M}{\pi r\delta|\zeta|}$.
Taking $r=|\zeta|^{-1}$, we obtain
$$
  \|S(\zeta)\|
  \le
  \frac{M}{\pi}\Bigl(\theta'e+\frac1\delta\Bigr),
  \quad \zeta\in\Sigma_{\eta_0}.
$$
Hence $S$ is bounded on the sector $\Sigma_{\eta_0}$ as claimed.

Next we prove strong continuity on each sector
$\Sigma_{\eta'}$ with $0<\eta'<\eta$. Fix $\eta'\in(0,\eta)$ and choose
$\theta'$ such that $\frac12\pi+\eta'<\theta'<\theta$.
Let $\Gamma=\Gamma_{\theta',B}$ be the corresponding contour, and choose
$\mu>0$ larger than the radius of $B$.
Let $x\in\Dom(A)$ and write $x=R(\mu,A)y$. Then, by the resolvent
identity and Cauchy's theorem,
$$
S(\zeta)x
 = \frac1{2\pi i}\int_\Gamma e^{\la\zeta}(\mu-\la)^{-1}R(\la,A)y\ud\la,
 \quad \zeta\in\Sigma_{\eta'}.
$$
Indeed, the additional scalar integral arising from the resolvent
identity vanishes, since $\mu$ lies to the right of $\Gamma$.

Let $r>0$ denote the radius of $B$. For $\zeta\in\Sigma_{\eta'}$ with
$r|\zeta|\le1$, along the rays of $\Gamma$ we have
$$
\|(\mu-\la)^{-1}R(\la,A)y\|\lesssim |\la|^{-2}\|y\|,
$$
while on the circular arc $|e^{\la\zeta}|\le e$. Hence the integrand is
dominated by an integrable function on $\Gamma$, uniformly for
$\zeta\in\Sigma_{\eta'}$ with $r|\zeta|\le1$. Therefore,
$$
S(\zeta)x \to \frac1{2\pi i}\int_\Gamma (\mu-\la)^{-1}R(\la,A)y\ud\la
=R(\mu,A)y=x
$$
as $\zeta\to 0$ in $\Sigma_{\eta'}$
by dominated convergence and Cauchy's theorem, closing
$\Gamma$ to the right in $\Sigma_{\theta'}\setminus B$
for the contour identity.
This proves that $S(\zeta)x\to x$ for $x\in \Dom(A)$ as $\zeta\to 0$ inside $\Sigma_{\eta'}$.
In view of the uniform boundedness of $S(\zeta)$ on $\Sigma_{\eta'}$,
convergence for general $x\in X$ follows by approximation.

It remains to show that $A$ is the generator of the semigroup
$(S(t))_{t\ge0}$. Let us denote this generator by $G$; thus we wish to prove that $G = A$.
Let $x\in\Dom(A)$. For $t>0$, differentiation under the contour
integral is justified by absolute convergence, and gives
\begin{align*}
\frac{\ud}{\ud t}S(t)x
&=
\frac1{2\pi i}
\int_\Gamma
e^{\la t}\la R(\la,A)x\ud\la
\\
&=
\frac1{2\pi i}
\int_\Gamma e^{\la t}x\ud\la
+
\frac1{2\pi i}
\int_\Gamma e^{\la t}R(\la,A)Ax\ud\la
= S(t)Ax.
\end{align*}
Here we used the identity
$\la R(\la,A)x=x+R(\la,A)Ax$
and the fact that
$\int_\Gamma e^{\la t}\ud\la=0.$
The latter follows by closing $\Gamma \cap \{|z|\le R\}$ with a circular arc of radius $R$ in the left half-plane, and then letting $R\to\infty$.

It follows that
$$
S(t)x-x=\int_0^tS(s)Ax\ud s,\quad t\ge0.
$$
Hence
$
\lim_{t\downarrow0}\frac1t(S(t)x-x)=Ax,
$
so $A\subseteq G$.

The semigroup $(S(t))_{t\ge0}$ is uniformly bounded, and therefore
$(0,\infty)\subseteq\varrho(G)$. By assumption,
$(0,\infty)\subseteq\varrho(A)$. Proposition
\ref{prop:semigroupsAB} now gives
$G=A.$

For every $0<\eta_0<\theta-\frac12\pi$, the preceding construction
therefore gives a bounded analytic $C_0$-semigroup on
$\Sigma_{\eta_0}$ with generator $A$. This proves
\ref{it:analytic-semigroup1} as well as the inequality
$\om_{\rm\, holo}(A)\ge \om_{\rm\, res}(A)-\frac12\pi$.

\smallskip
\ref{it:analytic-semigroup2}$\Leftrightarrow$\ref{it:analytic-semigroup3}:\ This is an immediate consequence of the second part of Lemma \ref{lem:res-Taylor-cor}.
\end{proof}

The following result characterises analytic $C_0$-semigroups directly in terms of the semigroup and
its generator, without reference to the resolvent. Its importance lies in the smoothing property
revealed by condition \ref{it:analytic-real2} below: the semigroup operators $S(t)$ map every $x\in X$ into the smaller subspace $\Dom(A)$ for all $t>0$.

\begin{theorem}[Bounded analytic semigroups, real characterisation]\label{thm:analytic-real} Let $A$ be the generator of a uniformly bounded $C_0$-semigroup $S$ on $X$.
The following assertions are equivalent:
\begin{enumerate}[label={\rm(\arabic*)}, leftmargin=*]
 \item\label{it:analytic-real1} $S$ is bounded analytic;
 \item\label{it:analytic-real2} $S(t)x \in \Dom(A)$ for all $x\in X$ and $t>0$, and
$\sup_{t>0} \, t \n AS(t)\n <\infty.$
\end{enumerate}
\end{theorem}

\begin{proof}
First note that assertion \ref{it:analytic-real2}, together
with the semigroup property and the commutation of $A$ with
$S(t)$ on $\Dom(A)$, implies $S(t)X\subseteq\Dom(A^n)$ and
$$
A^nS(t)=(AS(t/n))^n,\quad t>0,\ n\ge1.
$$
Consequently, $\sup_{t>0}t^n\n A^nS(t)\n<\infty$ for every
$n\ge1$.

\smallskip
\ref{it:analytic-real1}$\Rightarrow$\ref{it:analytic-real2}: \
Choose $0<\eta<\frac12\pi$ such that $S$ is bounded analytic
on $\Sigma_\eta$, and put
$M_\eta:=\sup_{z\in\Sigma_\eta}\n S(z)\n$.
For $t>0$ and $x\in X$, analyticity and the semigroup property
give
$$
\lim_{h\downarrow0}\frac1h(S(h)S(t)x-S(t)x)=S'(t)x.
$$
Thus $S(t)x\in\Dom(A)$ and $AS(t)x=S'(t)x$.
Fix $0<c<\sin\eta$. The closed disc with centre $t$ and
radius $ct$ lies in $\Sigma_\eta$,
so Cauchy's integral formula gives
\begin{align*}
\n AS(t)x\n
&=\Bigl\|\frac1{2\pi i}\int_{|z-t|=ct}
       \frac{S(z)x}{(z-t)^2}\ud z\Bigr\|
\le\frac{2\pi ct}{2\pi(ct)^2}M_\eta\n x\n
=\frac{M_\eta}{ct}\n x\n
\end{align*}
for $t>0$ and $x\in X$.
This proves \ref{it:analytic-real2}.

\smallskip
\ref{it:analytic-real2}$\Rightarrow$\ref{it:analytic-real1}: \
Put
$M:=\sup_{t\ge0}\n S(t)\n$ and $C:=\sup_{t>0}t\n AS(t)\n$.
By the observation at the beginning of the proof, every orbit $t\mapsto S(t)x$ is twice
continuously differentiable on $(0,\infty)$, with derivatives
$AS(t)x$ and $A^2S(t)x$, and we have
$$
\n A^2S(t)\n
=\n(AS(t/2))^2\n
\le\frac{4C^2}{t^2},\quad t>0.
$$

Let $\Re\la>0$ and put $\tau:=|\la|^{-1}$.
By Proposition \ref{prop:resolv}, $\la\in\varrho(A)$.
Splitting the Laplace integral at $\tau$ and integrating by parts
the second part twice, we obtain
\begin{align*}
R(\la,A)x
&=\int_0^\tau e^{-\la s}S(s)x\ud s
  +\frac{e^{-\la\tau}}{\la}S(\tau)x
  +\frac{e^{-\la\tau}}{\la^2}AS(\tau)x
  +\frac1{\la^2}\int_\tau^\infty
       e^{-\la s}A^2S(s)x\ud s.
\end{align*}
The boundary terms at infinity vanish by the preceding
bounds and the assumption $\Re\la>0$. Taking norms gives
\begin{align*}
\n R(\la,A)\n
&\le M\tau+\frac M{|\la|}
     +\frac C{|\la|^2\tau}
     +\frac{4C^2}{|\la|^2}\int_\tau^\infty s^{-2}\ud s
=\frac{2M+C+4C^2}{|\la|}.
\end{align*}
Theorem \ref{thm:analytic} and uniqueness of the semigroup
generated by $A$ now imply that $S$ is bounded analytic.
\end{proof}

\subsection{The Lumer--Phillips Theorem}\label{subsec:analytic-LP}

In Hilbert spaces, we have the following characterisation of contractive analytic $C_0$-semigroups (for an extension to Banach spaces see Problem \ref{prob:LuPhi}).

\begin{theorem}[Lumer--Phillips, contractive analytic semigroups] \label{thm:lumerphilips}\index{theorem!Lumer--Phillips}
Let $A$ be a densely defined closed operator in a Hilbert space $H$ and
let $0< \eta< \frac12\pi$. The following assertions are
equivalent:
\begin{enumerate}[label={\rm(\arabic*)}, leftmargin=*]
  \item\label{it:lumerphilips1} $A$ generates a contractive analytic $C_0$-semigroup on $H$ on the sector $\Sigma_\eta$;
  \item\label{it:lumerphilips2} $\mu-A$ has dense range for some $\mu>0$ and
  $-\iprod{Ax}{x}\in \ov{\Sigma_{\frac12\pi-\eta}}$  for all $x\in \Dom(A)$.
\end{enumerate}
\end{theorem}

\begin{proof}
\ref{it:lumerphilips1}$\Rightarrow$\ref{it:lumerphilips2}:
Since the $C_0$-semigroup $S$ generated by $A$ is contractive, Proposition \ref{prop:resolv} implies
$\{\Re\lambda>0\}\subseteq\varrho(A)$; in particular $1\in\varrho(A)$, so $I-A$ is surjective and has dense range.

Let $|\eta'| < \eta$, $x\in \Dom(A)$,
and consider the function $f(t) := \Re\iprod{S(te^{i\eta'})x}{x}$.
Observe that $f(0)=\|x\|^2$. Furthermore, for all $t\ge0$ we have
$$
|f(t)|
\le
|\iprod{S(te^{i\eta'})x}{x}|
\le
\|S(te^{i\eta'})x\|\|x\|
\le
\|x\|^2\!,
$$
where we used the contractivity of $S$ on $\Sigma_\eta$.
From these two observations we infer
that $f'(0) \leq 0$. On the other hand, differentiation gives
$f'(t) = \Re\iprod{e^{i\eta'}S(te^{i\eta'})Ax}{x}.$
Evaluating at $t=0$ gives $$\Re(e^{i\eta'}\iprod{Ax}{x}) \leq 0.$$
This can only be true for all $|\eta'|<\eta$ if
$\iprod{Ax}{x} \in -\ov{\Sigma_{\frac12\pi-\eta}}$.

\ref{it:lumerphilips2}$\Rightarrow$\ref{it:lumerphilips1}: \
Put
$W:=-\ov{\Sigma_{\frac12\pi-\eta}}$ and $\beta:=\frac12\pi+\eta$,
so that $U:=\C\setminus W=\Sigma_\beta$.
For $\la\in U$ write $d(\la):=\operatorname{dist}(\la,W)>0$.
Since $\iprod{Ax}{x}/\n x\n^2\in W$ for every nonzero
$x\in\Dom(A)$, the Cauchy--Schwarz inequality gives
\begin{align*}
\n(\la-A)x\n
\ge\frac{|\iprod{(\la-A)x}{x}|}{\n x\n}
\ge d(\la)\n x\n.
\end{align*}
Thus, for every $\la\in U$, the operator $\la-A$ is injective
and has closed range by Proposition \ref{prop:closed-range-unbdd}.
Since $\mu>0$, we have $\mu\in U$, so this applies in particular
to $\mu-A$. Its range is dense by assumption and therefore
equals $H$. The preceding estimate also shows that its inverse
is bounded, and hence $\mu\in\varrho(A)$. Moreover,
$$
\n R(\la,A)\n\le\frac1{d(\la)},
\quad \la\in\varrho(A)\cap U.
$$

The set $\varrho(A)\cap U$ is relatively open in $U$.
It is also relatively closed: if resolvent points converge
to a point of $U$, the preceding estimate keeps their
resolvent norms bounded, since $d$ is continuous and positive
on $U$. Proposition \ref{prop:res-blowup} therefore excludes
the possibility that such limits belong to $\sigma(A)$.
Since $U$ is connected and
$\mu\in\varrho(A)\cap U$, we conclude that
$U\subseteq\varrho(A)$.

For $\frac12\pi<\theta<\beta$, elementary geometry gives
$d(\la)\ge|\la|\sin(\beta-\theta)$ on $\Sigma_\theta$,
and hence
$$
\n R(\la,A)\n
\le\frac1{\sin(\beta-\theta)}\frac1{|\la|},
\quad \la\in\Sigma_\theta.
$$
For every $0<\eta_0<\eta$, choose $\theta$ with
$\frac12\pi+\eta_0<\theta<\beta$.
By Theorem \ref{thm:analytic}, $A$ generates a bounded
analytic $C_0$-semigroup on $\Sigma_{\eta_0}$.
These semigroups agree on the positive real axis by
uniqueness, and therefore agree on overlapping sectors by
uniqueness of analytic continuation. They define a
holomorphic semigroup $S$ on $\Sigma_\eta$, bounded and
strongly continuous at zero on every smaller sector.

For $x\in\Dom(A)$ and $z\in\Sigma_\eta$, the semigroup identity
$S(z+h)x=S(z)S(h)x=S(h)S(z)x$, $h>0$, and analyticity give
$S(z)x\in\Dom(A)$ and $S'(z)x=S(z)Ax=AS(z)x$.
Fix $|\alpha|<\eta$ and $x\in\Dom(A)$. The numerical range
assumption \ref{it:lumerphilips2} implies
$\Re\iprod{e^{i\alpha}Ay}{y}\le0$ for all $y\in\Dom(A)$.
Consequently,
$$
\frac{{\rm d}}{{\rm d}t}\|S(te^{i\alpha})x\|^2
=
2\Re\iprod{e^{i\alpha}AS(te^{i\alpha})x}
                 {S(te^{i\alpha})x}
\le0,\quad t>0.
$$
Since $S(te^{i\alpha})x\to x$ as $t\downarrow0$, we obtain
$\|S(te^{i\alpha})x\|\le\|x\|$. By density of $\Dom(A)$,
$S$ is contractive on $\Sigma_\eta$.

Finally, for $x\in\Dom(A)$, integration along the segment from $0$
to $z\in\Sigma_\eta$ gives
$$
S(z)x-x=\int_0^z S(w)Ax\ud w,
$$ which gives the bound $\|S(z)x-x\|\le |z|\,\|Ax\|$.
Contractivity and density now imply $S(z)x\to x$ for every $x\in H$
as $z\to0$ through $\Sigma_\eta$. Thus $S$ is a contractive analytic
$C_0$-semigroup on that sector.
\end{proof}

The conditions of the theorem are satisfied if $-A$ is a positive selfadjoint operator in $H$. In that case, the
semigroup $S$ generated by $A$ is given by
$$ S(t) = \int_{\sigma(-A)} e^{-t\la}\ud P(\la),$$
where $P$ is the projection-valued measure associated with $-A$. More generally, this formula can be used to
associate a $C_0$-semigroup of contractions with every normal operator $-A$ whose
spectrum is contained in the
closed right half-plane; see Theorem \ref{thm:normal-sgr}.

With the same proof, Theorem \ref{thm:lumerphilips} extends to $\eta = 0$,
provided we interpret $\Sigma_0$ as the positive real line and replace `analytic
$C_0$-semigroup of contractions' by `$C_0$-semigroup of contractions':

\begin{theorem}\label{thm:contrsgr-accretive}
Let $A$ be a densely defined closed operator acting in a Hilbert space $H$. The following assertions are
equivalent:
\begin{enumerate}[label={\rm(\arabic*)}, leftmargin=*]
 \item $A$ generates a $C_0$-semigroup of contractions on $H$;
 \item $\mu-A$ has dense range for some $\mu>0$ and $-\Re\iprod{Ax}{x}\ge 0$ for all $x\in \Dom(A)$.
\end{enumerate}
\end{theorem}

The condition `$-\Re\iprod{Ax}{x}\ge 0$ for all $x\in \Dom(A)$' says that $-A$ is {\em accretive}.\index{accretive!operator}
Since the open half-line $(0,\infty)$ is contained in the resolvent set of any operator generating a $C_0$-semigroup of contractions,
in Theorems \ref{thm:lumerphilips} and \ref{thm:contrsgr-accretive}, the condition
`$\mu-A$ has dense range for some $\mu>0$' may be replaced by `$1\in \varrho(A)$'. An accretive operator $-A$ satisfying $1\in\varrho(A)$
is also called a {\em maximal accretive} operator, or briefly, an {\em $m$-accretive operator}.\index{maccretive@$m$-accretive}\index{accretive!maximal}

\subsection{Semigroups Associated with Forms}\label{sec:semigroups-forms}

Theorem \ref{thm:form-accretive-resolvent} shows that the assumptions of the Hille--Yosida theorem are fulfilled. For coercive bounded forms, Corollary \ref{cor:form-coercive-resolvent} also supplies the estimate needed in Theorem \ref{thm:analytic}. Thus we obtain the following result.

\begin{theorem}[Bounded analytic semigroups via forms]\label{thm:form-analytic-sgr}
Let $A$ be a linear operator in a Hilbert space $H$. Then:
\begin{enumerate}[label={\rm(\arabic*)}, leftmargin=*]
\item\label{it:form-analytic-sgr1}
if $-A$ is the operator associated with a densely defined closed
accretive bounded form in $H$, then $A$ generates a $C_0$-contraction
semigroup on $H$, and for all $\delta>0$ the operator $A-\delta$
generates a bounded analytic $C_0$-semigroup on $H$;
\item\label{it:form-analytic-sgr2}
if $-A$ is the operator associated with a densely defined closed
coercive bounded form in $H$, then $A$ generates a $C_0$-contraction
semigroup on $H$ that extends to a bounded analytic $C_0$-semigroup
on $H$.
\end{enumerate}
\end{theorem}

\begin{proof}
Only the second assertion in \ref{it:form-analytic-sgr1} needs proof.
Let $\aa$ denote the densely defined closed accretive bounded form
with which $-A$ is associated.

Fix $\delta>0$ and let $\aa_\delta$ be the form with domain
$\Dom(\aa_\delta):=\Dom(\aa)$ defined by
$$
\aa_\delta(u,v):=\aa(u,v)+\delta\iprod{u}{v}.
$$
Its form norm satisfies
$$
\n u\n_\aa^2
\le
\n u\n_{\aa_\delta}^2
=
\n u\n_\aa^2+\delta\n u\n^2
\le
(1+\delta)\n u\n_\aa^2.
$$
Thus the form norms are equivalent, and hence $\aa_\delta$ is closed.
If $C$ is a boundedness constant for $\aa$, then
$$
|\aa_\delta(u,v)|
\le
(C+\delta)\n u\n_{\aa_\delta}\n v\n_{\aa_\delta},
$$
so $\aa_\delta$ is bounded. Moreover,
$\Re\aa_\delta(u,u)\ge\delta\n u\n^2$, and therefore
$\aa_\delta$ is coercive. The operator associated with
$\aa_\delta$ is $-A+\delta$. Hence Corollary \ref{cor:form-coercive-resolvent}, applied to the form
$\aa_\delta$, gives a constant $M_\delta$ such that
$$
\n R(\lambda,A-\delta)\n
\le
\frac{M_\delta}{|\lambda|},
\quad \Re\lambda>0.
$$
Theorem \ref{thm:analytic} now implies that $A-\delta$ generates
a bounded analytic $C_0$-semigroup.
\end{proof}

\begin{example}[Operators in Divergence Form I]\label{ex:exmaples-revisited-II}
Let $D$ be a nonempty bounded open subset of $\R^d$, and assume
that the matrix-valued function $a:D\to M_d(\C)$ satisfies
\begin{enumerate}[label={\rm(\roman*)}, leftmargin=*]
\item the functions $a_{ij}:D\to\C$ are measurable and bounded;
\item the function $a$ is uniformly elliptic, that is,
there exists $\alpha>0$ such that for almost all $x\in D$,
$$
\Re\sum_{i,j=1}^d a_{ij}(x)\xi_j\ov{\xi_i}
\ge\alpha|\xi|^2\!,
\quad \xi\in\C^d.
$$
\end{enumerate}
As shown in Section \ref{subsec:divergence-form-operator}, the form
$$
\aa_a(u,v):=\int_D a\nabla u\cdot\ov{\nabla v}\ud x,
\quad u,v\in\Dom(\aa_a):=H^1(D),
$$
is densely defined, closed, bounded, and accretive in $L^2(D)$.
The operator $-A_a$ associated with this form is the Neumann
realisation of $-{\rm div}(a\nabla)$. By part
\ref{it:form-analytic-sgr1} of Theorem \ref{thm:form-analytic-sgr},
$A_a$ generates a $C_0$-contraction semigroup on $L^2(D)$, and
$A_a-\delta$ generates a bounded analytic $C_0$-semigroup for
every $\delta>0$.

If the form domain is replaced by $H_0^1(D)$, we obtain the
Dirichlet realisation in $L^2(D)$. The resulting form is densely
defined, closed, bounded, and coercive, so part
\ref{it:form-analytic-sgr2} of Theorem \ref{thm:form-analytic-sgr}
applies.
If, in addition, $D$ is connected and has a $C^1$-boundary, we may
instead view $\aa_a$ as a form in $L^2_{\rm av}(D)$ with domain
$H^1_{\rm av}(D)$. This form is densely defined, closed, bounded,
and coercive, so part \ref{it:form-analytic-sgr2} of Theorem
\ref{thm:form-analytic-sgr} also applies to the Neumann realisation
on $L^2_{\rm av}(D)$.
\end{example}

A sectoriality condition on the form gives contractivity on a
sector. We first state the general result and then return to the
divergence form operators of Example \ref{ex:exmaples-revisited-II}.

\begin{theorem}[Contractive analytic semigroups via forms]\label{thm:om-sect-form-anal-sgr}
Let $H$ be a Hilbert space and let $A$ be the densely defined closed
operator in $H$ associated with a densely defined closed form $\aa$
that is $\omega$-sectorial for some $0<\omega<\frac12\pi$.
Then $-A$ generates an analytic $C_0$-semigroup of contractions on
the sector $\Sigma_{\frac12\pi-\omega}$.
\end{theorem}

\begin{proof}
By Theorem \ref{thm:form-sectorial-resolvent}, $I+A$ is surjective.
Moreover,
$$
\iprod{Au}{u}=\aa(u,u)\in\ov{\Sigma_\omega},
\quad u\in\Dom(A).
$$
Theorem \ref{thm:lumerphilips}, applied to $-A$, gives the result.
\end{proof}

\begin{example}[Operators in Divergence Form II]\label{ex:exmaples-revisited-I}
Let $D$ and $a$ satisfy the assumptions of Example
\ref{ex:exmaples-revisited-II}, and let $A_a$ be the
Dirichlet or the Neumann realisation of ${\rm div}(a\nabla)$
in $L^2(D)$. The corresponding form $\aa_a$ has domain
$H_0^1(D)$ or $H^1(D)$, respectively.

Choose $M>0$ such that $|a(x)\xi|\le M|\xi|$ for almost all
$x\in D$ and all $\xi\in\C^d$. Then
$$
|\Im\aa_a(u,u)|
\le M\n\nabla u\n_2^2
\le \frac M\alpha\,\Re\aa_a(u,u),
\quad u\in\Dom(\aa_a).
$$
Thus $\aa_a$ is $\omega$-sectorial for
$\omega:=\arctan(M/\alpha)\in(0,\frac12\pi)$.
Since $-A_a$ is associated with $\aa_a$, Theorem
\ref{thm:om-sect-form-anal-sgr} shows that $A_a$ generates an
analytic $C_0$-semigroup of contractions on
$\Sigma_{\frac12\pi-\omega}$. This conclusion holds for both
realisations on $L^2(D)$.
If $D$ is connected and has a $C^1$-boundary, the same argument
applies to the Neumann realisation on $L^2_{\rm av}(D)$.

Suppose now that the coefficient assumptions are strengthened to:
\begin{enumerate}[label={\rm(\roman*)}, leftmargin=*]
\item\label{it:divergence-form1a}
the functions $a_{ij}:D\to\C$ are measurable and bounded;
\item\label{it:divergence-form2a}
for almost all $x\in D$, the matrix $a(x)$ is Hermitian, and
there exists a constant $\alpha>0$, independent of $x$, such that
$$
\sum_{i,j=1}^d a_{ij}(x)\xi_j\ov{\xi_i}
\ge\alpha|\xi|^2\!,
\quad \xi\in\C^d.
$$
\end{enumerate}
Condition \ref{it:divergence-form2a} adds Hermitian symmetry
to the uniform ellipticity assumption in Example
\ref{ex:exmaples-revisited-II}. The form $\aa_a$ is therefore
symmetric and positive. In particular, $\aa_a(u,u)\ge0$ for
all $u\in\Dom(\aa_a)$, so $\aa_a$ is $\omega$-sectorial for
every $\omega\in(0,\frac12\pi)$.
\end{example}

\subsection{Maximal Regularity}\label{subsec:maxreg}

In Section \ref{sec:ACP} we have seen that the mild solution $u$ of the inhomogeneous problem $u' = Au+f$ with initial condition $u(0)=u_0$, which is given by
$$ u(t) = S(t)u_0+ \int_0^t S(t-s)f(s)\ud s, \quad t\ge 0,$$
is also a strong solution, that is, for all $t\ge 0$ we have $\int_0^t u(s)\ud s \in \Dom(A)$ and
$$
u(t) = u_0 + A\int_0^t u(s)\ud s
+ \int_0^t f(s)\ud s.$$
In general it cannot be asserted, however, that
$$
u(t) = u_0 + \int_0^t Au(s)\ud s
+ \int_0^t f(s)\ud s,$$
the problem being that $u$ need not take values in $\Dom(A)$ almost everywhere, and even if it does, we cannot be certain that $s\mapsto Au(s)$ is integrable on intervals $(0,t)$. The aim of the present section is to prove that these things do hold for $f\in L^2(\R_+;H)$ with zero initial condition $u_0=0$, provided that $A$ generates a bounded analytic $C_0$-semigroup on a Hilbert space.

\begin{theorem}[Maximal regularity]\label{thm:MR}\index{maximal regularity}
 Let $A$ be the generator of a bounded analytic $C_0$-semigroup $S$ on a Hilbert space $H$.
 Then for all $f\in L^2(\R_+;H)$ the mild solution $u=u_f$ of the inhomogeneous problem $u' = Au+f$ with initial condition $u(0)=0$
 enjoys the following properties:
 \begin{enumerate}[label={\rm(\arabic*)}, leftmargin=*]
  \item\label{it:MR1}  $u(t)\in \Dom(A)$ for almost all $t\ge 0$, $Au$ belongs to $L^2(\R_+;H)$,
  and for
  all $t\ge 0$ we have $$u(t) = \int_0^t Au(s)\ud s + \int_0^t f(s)\ud s;$$
  \item\label{it:MR2} we have $$\n Au\n_2 \le C\n f\n_2,$$
 where $C = \sup_{\xi\in \R\setminus\{0\}}\n AR(i \xi,A)\n$.
 \end{enumerate}
\end{theorem}

By the Lebesgue differentiation theorem (Theorem \ref{thm:Leb-diff}, or rather, its vector-valued version which is proved in exactly the same way), the identity in \ref{it:MR1}
implies that $u$ is differentiable almost everywhere and that the pointwise identity
$$ u'(t) = Au(t) + f(t)$$
holds for almost all $t\ge 0$.
This, in combination with \ref{it:MR2}, implies that
$u'$ also belongs to $L^2(\R_+;H)$
(with estimate $\n u'\n_2\le (C+1)\n f\n_2$). This explains the name `maximal regularity' attached to the theorem.

We begin with a reduction to a class of `nice' functions $f$. To this end, for subspaces $F$ and $Y$ of $L^2(\R_+)$
and $H$, respectively, we introduce the notation
$F\otimes Y$ for the vector space of all linear combinations of functions $f:\R_+\to H$ of the form $f= \phi\otimes y$ with $\phi\in F$ and $y\in Y$, where
$$ (\phi\otimes y)(t):= \phi(t)y, \quad t\ge 0.$$
If $F$ and $Y$ are dense in $L^2(\R_+)$
and $H$, respectively, then $F\otimes Y$ is a dense subspace of $L^2(\R_+;H)$. This is because the simple functions are dense in $L^2(\R_+;H)$ and every such function is a linear combination of functions of the form $\one_B\otimes h$ with $B\subseteq \R_+$ a Borel set of finite measure and $h\in H$; we now approximate $\one_B$ by functions in $F$ (in the norm of $L^2(\R_+)$) and $h$ by elements of $Y$ (in the norm of $H$).

In what follows we consider the dense subspaces $F = C_{\rm c}^1(\R_+)$,
the space of continuously differentiable functions with compact support in $(0,\infty)$, and $Y = \Dom(A)$.
For functions $f\in C_{\rm c}^1(\R_+)\otimes \Dom(A)$
the mild solution of the problem $u' = Au+f$ with initial condition $u(0)=0$, given by
$$ u(t) = \int_0^t S(t-s)f(s)\ud s, \quad t\ge 0,$$
is continuously differentiable in $H$, takes values in $\Dom(A)$, and satisfies $u'(t) = Au(t) + f(t)$ for every $t\ge 0$. Indeed, $s\mapsto S(t-s)f(s)$ is continuous in the graph norm of $\Dom(A)$, so we may pass $A$ under the integral and differentiate using Proposition \ref{prop:sg-prop}\ref{it:sg-prop2}.

\begin{lemma}\label{lem:MR1} Let $A$ be the generator of a bounded analytic $C_0$-semigroup $S$ on a Hilbert space $H$.
 If there exists a constant $C\ge 0$ such that for all $f\in C_{\rm c}^1(\R_+)\otimes \Dom(A)$ the mild solution $u$ associated with $f$ satisfies $Au \in L^2(\R_+;H)$ and
 $$ \n Au\n_2 \le C\n f\n_2,$$
where $C\ge 0$ is a constant independent of $f$,
then the assertions \ref{it:MR1} and \ref{it:MR2} of Theorem \ref{thm:MR} hold
 for all $f\in L^2(\R_+;H)$, with the same constant $C$.
\end{lemma}
\begin{proof}
Since $C_{\rm c}^1(\R_+)\otimes \Dom(A)$ is dense in $L^2(\R_+;H)$, for any
$f\in L^2(\R_+;H)$ we may choose functions $f_n\in C_{\rm c}^1(\R_+)\otimes \Dom(A)$ converging to $f$ in $L^2(\R_+;H)$.
Writing $u$ and $u_n$ for the mild solutions corresponding to $f$ and $f_n$, respectively, from the assumed {\em a priori} bound we infer that the functions $Au_n$ form a Cauchy sequence in $L^2(\R_+;H)$ and therefore converge to a limit $v$ in $L^2(\R_+;H)$.

Let
$M_T:=\sup_{0\le t\le T}\n S(t)\n.$
The formula for the mild solution gives
$$
\n u_n-u\n_{C([0,T];H)}
\le M_T T^{1/2} \n f_n-f\n_{L^2(0,T;H)},
$$
and hence $u_n\to u$ uniformly on every interval $[0,T]$.
Passing to a subsequence, we may assume that
$Au_n(t)\to v(t)$
for almost all $t\ge0$. Since $u_n(t)\to u(t)$ for every $t\ge0$ and
$A$ is closed, it follows that
$u(t)\in\Dom(A)$ and $Au(t)=v(t)$
for almost all $t\ge0$.

For all $n\ge 1$ and $t\ge0$,
$$
u_n(t) = \int_0^tAu_n(s)\ud s + \int_0^tf_n(s)\ud s.
$$
Passing to the limit as $n\to\infty$ gives
$$
u(t) = \int_0^tAu(s)\ud s + \int_0^tf(s)\ud s.
$$
Finally,
$ \n Au\n_2 = \n v\n_2 \le C\n f\n_2$.
\end{proof}

The constant is sharp even for $H=\C$ and $A=-I$: here $C=1$, and taking $f_T=\one_{[0,T]}$ gives $\n Au_{f_T}\n_2^2/\n f_T\n_2^2=1-(1-e^{-T})/(T)\to1$ as $T\to\infty$.

The Fourier--Plancherel transform on $L^2(\R)$ extends to an isometry
of $L^2(\R;H)$ onto itself by first defining
$$
\calF(\one_B\otimes h):=(\calF\one_B)\otimes h
$$
for Borel sets $B\subseteq\R$ of finite measure and $h\in H$,
and then extending linearly and continuously. That this extension enjoys the stated properties can be proved in exactly the same way as in the scalar-valued case, by repeating the proof given for that case word by word with the obvious adjustments.

\begin{proof}[Proof of Theorem \ref{thm:MR}]
By Lemma \ref{lem:MR1}, it suffices to prove the asserted
estimate for $f\in C_{\rm c}^1(\R_+)\otimes\Dom(A)$.
For such $f$, the associated mild solution takes values
in $\Dom(A)$ and satisfies
$$
Au(t)=Vf(t):=\int_0^t S(t-s)Af(s)\ud s,\quad t\ge0.
$$
Extending $f$ and $Af$ by zero to $\R$, we may write
$$
Vf(t)=\int_0^\infty S(r)Af(t-r)\ud r,\quad t\in\R.
$$
This also extends $Vf$ by zero to the negative half-line.

Put
$$
m(\xi):=AR(i\xi,A)=i\xi R(i\xi,A)-I,\quad \xi\ne0,
$$
and set $m(0):=0$. By Theorem \ref{thm:analytic}, this
function is bounded, with
$$
C:=\sup_{\xi\ne0}\n m(\xi)\n<\infty.
$$
Indeed, Theorem \ref{thm:analytic} gives $\theta>\frac12\pi$ and $M':=\sup_{\la\in\Sigma_\theta}\n\la R(\la,A)\n<\infty$, so $i\R\setminus\{0\}\subseteq\Sigma_\theta$ and $C\le M'+1$. Moreover, $m$ is continuous in operator norm on $\R\setminus\{0\}$, so $m\calF f$ is strongly measurable; the value assigned to $m(0)$ is immaterial.
Thus the Fourier multiplier
$T_mf:=\calF^{-1}(m\calF f)$ is bounded on $L^2(\R;H)$
with norm at most $C$. We will show that $Vf=T_mf$.

For $\eps>0$ define
$$
V_\eps f(t):=\int_0^\infty
e^{-\eps r}S(r)Af(t-r)\ud r,\quad t\in\R.
$$
Since $Af\in L^1(\R;H)\cap L^2(\R;H)$ and $S$ is bounded,
Young's inequality gives
$V_\eps f\in L^1(\R;H)\cap L^2(\R;H)$.
Moreover, since $f$ is a finite sum of functions with
values in $\Dom(A)$, we have
$\wh f(\xi)\in\Dom(A)$ and
$\wh{Af}(\xi)=A\wh f(\xi)$.
Fubini's theorem and Proposition \ref{prop:resolv} therefore give
\begin{align*}
\wh{V_\eps f}(\xi)
&=\int_0^\infty e^{-(\eps+i\xi)r}
S(r)\wh{Af}(\xi)\ud r
=R(\eps+i\xi,A)\wh{Af}(\xi)
=AR(\eps+i\xi,A)\wh f(\xi).
\end{align*}

By Theorem \ref{thm:analytic} and the identity
$AR(\la,A)=\la R(\la,A)-I$, the operators
$AR(\eps+i\xi,A)$ are uniformly bounded for
$\eps>0$ and $\xi\in\R$. For every $\xi\ne0$,
continuity of the resolvent gives
$$
AR(\eps+i\xi,A)\to AR(i\xi,A)=m(\xi)
$$
in operator norm as $\eps\downarrow0.
$
Dominated convergence and Plancherel's theorem imply that
$$
V_\eps f\to T_mf\quad\hbox{in }L^2(\R;H).
$$
On the other hand, dominated convergence in the defining
integral gives $V_\eps f(t)\to Vf(t)$ for every $t\in\R$.
Passing to a sequence $\eps_n\downarrow0$ along which
the $L^2$-convergence also holds almost everywhere, we obtain
\begin{align*}
Vf=T_mf\quad\hbox{almost everywhere on }\R.
\end{align*}
Consequently,
\begin{align*}
\n Au\n_{L^2(\R_+;H)}
=\n Vf\n_{L^2(\R;H)}
=\n m\wh f\n_{L^2(\R;H)}
\le C\n f\n_{L^2(\R_+;H)}.
\end{align*}
Lemma \ref{lem:MR1} now completes the proof.
\end{proof}

We demonstrate the usefulness of maximal regularity by proving local existence for the time-dependent inhomogeneous Cauchy problem
\begin{equation}\label{eq:nonaut}
\left\{
\begin{aligned}
 u'(t)& = A(t)u(t) + f(t), \quad t\in [0,T],\\
 u(0) &= 0,
\end{aligned}
\right.
\end{equation}
where $(A(t))_{t\in [0,T]}$ is a family of densely defined closed operators in a Hilbert space $H$.

We make the following assumptions:
\begin{itemize}
\item there exists a Banach space $D$, continuously and densely embedded in $H$, such that $$\Dom(A(t))=D, \quad t\in[0,T];$$
\item for each $t\in[0,T]$, the operator $A(t)$ belongs to $\calL(D,H)$;
\item the mapping $t\mapsto A(t)\in \calL(D,H)$ is continuous on $[0,T]$;
\item the operator $A(0)$ is invertible and generates a bounded analytic
$C_0$-semigroup on $H$;
\item $f \in L^2(0,T;H)$.
\end{itemize}
The idea is to rewrite the problem in the form
$$
 u'(t)=A(0)u(t)+g_u(t) \ \hbox{ with }
 g_u(t):=(A(t)-A(0))u(t)+f(t).
$$
Now let $0<a\le T$ and consider a fixed function $u\in L^2(0,a;D)$.
Denote by $K_a(u)$ the unique mild solution of the inhomogeneous problem
\begin{equation}\label{eq:A0}
\left\{
\begin{aligned}
v'(t)&=A(0)v(t)+g_u(t),\quad t\in[0,a], \\ v(0)&=0.
\end{aligned}
\right.
\end{equation}
Then, at least formally, the solutions of \eqref{eq:nonaut} are the
fixed points of $K_a$.
The maximal regularity of $A(0)$ will now be used to show that
$K_a$ is a uniform contraction  (that is, its Lipschitz constant is strictly smaller than one) on $L^2(0,a;D)$, provided $0<a\le
T$ is small enough. The Banach fixed point theorem then gives the existence of a unique fixed point of $K_a$ in $L^2(0,a;D)$.
This fixed point will be called the {\em solution} on $(0,a)$.

Let us begin by checking that $K_a$ maps $L^2(0,a;D)$ into itself.
For $u\in L^2(0,a;D)$, continuity of
$t\mapsto A(t)$ in $\calL(D,H)$ implies that
$g_u\in L^2(0,a;H)$.
Maximal regularity gives
$$A(0)K_a(u)\in L^2(0,a;H).$$
Since $0\in\varrho(A(0))$, the operator $A(0):D\to H$ is bijective.
By the open mapping theorem there exists a constant $c_D\ge0$ such that
$$
\n x\n_D\le c_D\n A(0)x\n_H,\quad x\in D.
$$
It follows that
$K_a(u)\in L^2(0,a;D)$.

Next, if $ u_1,u_2\in L^2(0,a;D),$ then
$K_a{(u_1)}-K_a{(u_2)}$ equals the solution $u$ of
$$
u'(t)=A(0)u(t)+g_{u_1}(t)-g_{u_2}(t), \quad u(0)=0.
$$
Using the maximal regularity inequality for the problem \eqref{eq:A0} on the
interval $(0,a)$, we obtain
\begin{align*}
\|K_a(u_1)- K_a(u_2)\|_{L^2(0,a;D)}
 &\le c_D\|A(0)u\|_{L^2(0,a;H)}\\
 &\le c_D C\|g_{u_1}-g_{u_2}\|_{L^2(0,a;H)}\\
 &=c_D C\|[A(\cdot)-A(0)](u_1-u_2)\|_{L^2(0,a;H)}\\
 &\le c_D C \sup_{t\in [0,a]}\|A(t)-A(0)\|_{\calL(D,H)}
        \|u_1-u_2\|_{L^2(0,a;D)}.
\end{align*}
To justify the second inequality, we extend the inhomogeneity $g_{u_1}-g_{u_2}$ by zero on $[a,\infty)$,
and observe that the mild solution to the resulting inhomogeneous problem on $\R_+$ restricts to $u$ on the interval $(0,a)$.

If the constant $a>0$ is small enough, then
$$
  c_D C\sup_{t\in[0,a]}\|A(t)-A(0)\|_{\calL(D,H)}<1,
$$
and the Banach fixed point
theorem provides a unique solution for \eqref{eq:nonaut} on $(0,a)$.

\section{Stone's Theorem}\label{sec:Stone}

If $A$ is a positive selfadjoint operator in a Hilbert space $H$, then $-A$ satisfies the conditions of Theorem \ref{thm:lumerphilips}. Let $S$ be the contractive analytic $C_0$-semigroup generated by $-A$. It can be shown (see Problem \ref{prob:bdry-sgr}) that for all $t\in \R$ and $x\in H$ the limit
$$U(t)x:= \lim_{s\downarrow 0} S(s+it)x$$
exists
and that the family $(U(t))_{t\in \R}$ is a $C_0$-group of unitary operators with
generator $-iA$.
The main result of this section is Stone's theorem, which asserts that, for any selfadjoint operator $A$, the operator $-iA$ generates a $C_0$-group of unitary operators.
For the proof of this theorem we need the following auxiliary result.

\begin{proposition}\label{prop:sa-real-spectrum}
If $A$ is a selfadjoint operator in $H$, then $\sigma(A)\subseteq \R$ and
$$ \n R(\la,A)\n \le \frac1{|\Im \la|}, \quad  \la\in\C\setminus\R.$$
\end{proposition}
\begin{proof} The inclusion $\sigma(A)\subseteq \R$  has been shown in Proposition \ref{prop:sa-pos-spectrum}.
Let $\la =  \alpha+i\beta$ with $\beta \not=0$.
For all $x\in \Dom(A)$ we have
$\iprod{Ax}{x}\in \R$ and
 \begin{align*} \n (\la-A)x\n\n x\n \ge |\iprod{(\la-A)x}{x}| & = \big|\al \iprod{x}{x} - \iprod{Ax}{x} + i\beta\iprod{x}{x}\big| \ge |\beta|\n x\n^2
 \end{align*}
and therefore $\n (\la-A)x\n \ge |\beta|\n x\n$.
Since $\la\in \varrho(A)$, this implies $\n R(\la,A)\n \le 1/|\beta|$.
\end{proof}

\begin{theorem}[Stone]\label{thm:Stone}\index{theorem!Stone}
 For a densely defined operator $A$ in $H$, the following assertions are equivalent:
 \begin{enumerate}[label={\rm(\arabic*)}, leftmargin=*]
  \item\label{it:Stone1} $A$ is selfadjoint;
  \item\label{it:Stone2} $-iA$ is the generator of a $C_0$-group of unitary operators.\index{unitary!group}
 \end{enumerate}
\end{theorem}
\begin{proof}

\smallskip
\ref{it:Stone1}$\Rightarrow$\ref{it:Stone2}: \
By Proposition \ref{prop:sa-real-spectrum}, $\si(A)$ is contained in the real line
and for all $\la\in \C\setminus\R$ we have $\n R(\la, A)\n \le 1/|\Im \la|$.
Hence, by the Hille--Yosida theorem (Theorem \ref{thm:HilleYosida}), the operators $\pm iA$ generate $C_0$-contraction semigroups $S_\pm$.
Hence, by Proposition \ref{prop:C0group}, $-iA$ generates a $C_0$-group of contractions given by
$$ U(t) :=
\begin{cases}
 S_-(t), & \ t \ge 0,\\
 S_+(-t), & \ t \le 0.
\end{cases}
$$
Since $U(t)$ and $U(-t)=U(t)^{-1}$ are contractions, for every $x\in H$ we have
$$
\n U(t)x\n\le\n x\n=\n U(-t)U(t)x\n\le\n U(t)x\n.
$$
Thus $U(t)$ is a surjective isometry, and hence unitary.

\smallskip
\ref{it:Stone2}$\Rightarrow$\ref{it:Stone1}: \
Suppose that $-iA$ generates the unitary group $U$.
For $x,y\in\Dom(A)$, differentiation of
$\iprod{U(t)x}{U(t)y}=\iprod{x}{y}$ at $t=0$ gives
$$
0=-i\iprod{Ax}{y}+i\iprod{x}{Ay}.
$$
Thus $A$ is symmetric, that is, $A\subseteq A^\star$.
Moreover, $A$ is closed because $-iA$ is a generator.

Both $-iA$ and $iA$ generate contraction semigroups,
so Proposition \ref{prop:resolv} gives
$\pm i\in\varrho(A)$. By Proposition
\ref{prop:spectrum-dual-unbdd}, we also have
$i\in\varrho(A^\star)$.
The inclusion $A\subseteq A^\star$ and Proposition
\ref{prop:semigroupsAB} now imply $A=A^\star$.
\end{proof}

\section{Examples}\label{sec:examples-sg}

In this section we collect some important examples of $C_0$-semigroups and $C_0$-groups.

\subsection{Multiplication Semigroups}\label{subsec:multsgr}\index{multiplication!semigroup}
Let $(\Omega,\calF\!,\mu)$ be a measure space and let $m:\Omega\to \K$ be
measurable with real part essentially bounded from below:
$$ M:= \operatorname*{ess\,inf}_{\om\in \Om} \Re m(\om)> -\infty.$$
The operators
$$ S(t)f:= e^{-tm}f, \quad t\ge 0,$$
are bounded on $L^p(\Om)$, $1\le p\le \infty$, with norm $\n S(t)\n \le e^{-tM}$. We will prove that $S$ is a $C_0$-semigroup on $L^p(\Omega)$ for $1\le p<\infty$, with
generator $A$ given by
\begin{equation}\label{eq:multsgr-DA}
\begin{aligned}
\Dom(A) &:= \big\{f\in L^p(\Omega):\, mf\in L^p(\Omega)\big\},\\   Af &:= -mf, \quad f\in \Dom(A).
\end{aligned}
\end{equation}

Fix $1\le p<\infty$.
The semigroup properties \ref{S1} and \ref{S2} are clear and \ref{S3} follows by dominated convergence.
To prove \eqref{eq:multsgr-DA} let $f\in L^p(\Omega)$ be such that $mf\in L^p(\Omega)$.
For $\mu$-almost all $\omega\in \Omega$ we have
$$ S(t)f(\omega) - f(\omega) =  e^{-tm(\omega)} f(\omega) - f(\omega) = -m(\omega)f(\omega)\int_0^t e^{-sm(\omega)}\ud s.$$
After division by $t$, in $L^p(\Om)$ this identity reads
$$
\frac1t(S(t)f-f)
=-\frac1t\int_0^t S(s)(mf)\ud s.
$$
Since $mf\in L^p(\Om)$, strong continuity implies that
the right-hand side converges to $-mf$ as $t\downarrow0$.
Thus $f\in\Dom(A)$ and $Af=-mf$.

Conversely, let $f\in\Dom(A)$. By the definition of the generator,
$\frac1t(S(t)f-f)\to Af$ in $L^p(\Omega)$ as $t\downarrow0$.
Since $p<\infty$, there exists a sequence $t_n\downarrow0$
along which the difference quotients converge to $Af$
pointwise almost everywhere. On the other hand,
$\frac1{t_n}(e^{-t_nm(\omega)}-1)f(\omega)
\to -m(\omega)f(\omega)$
for almost all $\omega\in\Omega$ as $n\to\infty$. Hence $Af=-mf$ almost
everywhere, and therefore $mf\in L^p(\Omega)$. This proves
the reverse inclusion.

\subsection{The Translation Group}\label{subsec:translation-sg}\index{translation!group}
On the space $L^p(\R)$, $1\le p<\infty$,  the formula
$$ (S(t)f)(x) := f(x+t), \quad x\in \R, \ t\in\R,$$
defines a $C_0$-group $S$.
Its generator $A$ is given by
\begin{equation*}
 \begin{aligned}
\Dom(A) & = W^{1,p}(\R), \\ Af & = f'\!, \quad f\in \Dom(A).
 \end{aligned}
\end{equation*}
The group properties \ref{G1} and \ref{G2} are clear and \ref{G3}
follows from Proposition \ref{prop:Lp-transl}.

To prove that $\Dom(A) = W^{1,p}(\R)$ and $Af = f'$\!, we first note that for $f\in C_{\rm c}^1(\R)$
we have
$$S(t)f(x) - f(x) =  f(x+t) - f(x) = \int_0^t f'(x+s)\ud s = \int_0^t S(s)f'(x)\ud s.$$
It follows that
$S(t) f-f = \int_0^t S(s)f'\ud s$.
Also,
$ S(t)f-f = A\int_0^t S(s)f\ud s.$
It follows that $$A\int_0^t S(s)f\ud s = \int_0^t S(s)f'\ud s.$$
Next we note that $\lim_{t\to 0} \frac1t\int_0^t S(s)f\ud s = f$ and
$$\lim_{t\to 0}A\frac1t\int_0^t S(s)f\ud s = \lim_{t\to 0}\frac1t\int_0^t S(s)f'\ud s = f'$$ in $L^p(\R)$.
Since $A$ is closed, this implies
$ f\in \Dom(A)$ and $Af = f'\!.$

Since $C^1_{\rm c}(\R)$ is dense in $L^p(\R)$ and invariant under translations, from Proposition \ref{prop:core}
we infer that $C^1_{\rm c}(\R)$ is dense in $\Dom(A)$. Since $C^1_{\rm c}(\R)$ is also
dense in $W^{1,p}(\R)$ and $\n f\n_{\Dom(A)} = \n f\n + \n Af\n = \n f \n + \n f'\n = \n f\n_{W^{1,p}(\R)}$
for all $f\in C_{\rm c}^1(\R)$, it follows that $\Dom(A) = W^{1,p}(\R)$ and $Af = f'$ for all $f\in \Dom(A) = W^{1,p}(\R)$.

\subsection{The Heat Semigroup}\label{subsec:heat-sgr}\index{heat!semigroup, on $\R^d$}

\paragraph{The Heat Semigroup on \texorpdfstring{$\R^d$}{R\textasciicircum d}}\label{par:heat-sgr-Rd}\index{heat!semigroup, on $\R^d$}

For $1\le p<\infty$ and $t>0$ we define a bounded operator $H(t)$ on $L^p(\R^d)$ by
\begin{equation*}
 H(t)f = K_t *f, \quad f\in L^p(\R^d),
\end{equation*}
where $$K_t(x) := (4\pi t)^{-d/2}
e^{-|x|^2/4t}$$ is the {\em heat kernel}.\index{heat!kernel}
Since $K_t\in L^1(\R^d)$ with $\n K_t\n_1 =1$, it follows from
Young's inequality (Proposition \ref{prop:Young}) that for all $1\le p<\infty$ the operators $H(t)$ are
well defined and bounded on $L^p(\R^d)$ and satisfy
$$\|H(t) f\|_{L^p(\R^d)} \leq \|f\|_{L^p(\R^d)}.$$ We furthermore set $H(0):= I$, the identity operator
on $L^p(\R^d)$.
We will prove that the family $H=(H(t))_{t\ge 0}$ is a $C_0$-semigroup of contractions, the so-called {\em heat semigroup}, on
$L^p(\R^d)$ and that its generator $A$ is the weak $L^p$-Laplacian $\Delta$.
Thus the heat semigroup solves the linear heat equation\index{heat!equation, linear}
$$\left\{
\begin{aligned} \frac{\partial u}{\partial t}(t,x) & = \Delta u(t,x), && t\ge 0, \ x\in \R^d\!,\\
 u(0,x) & = f(x), && x\in \R^d\!,
\end{aligned}
\right.
$$
in the sense that its orbits are strongly continuous and satisfy $H(0)f=f$ and $\frac{{\rm d}}{{\rm d}t}H(t)f = \Delta H(t)f$ for $t>0$.

\smallskip
{\em Step 1} -- Fix $1\le p<\infty$.  First we prove that $H$ is a $C_0$-semigroup on $L^p(\R^d)$.
For all $t>0$,  by Lemma \ref{lem:Gauss} and a change of variables the Fourier transform of $K_t$
is given by
$$
\widehat{K_t}(\xi):=\frac1{(2\pi)^{d/2}}\int_{\R^d} K_t(x) e^{- i x \cdot\xi}\ud x = \frac1{(2\pi)^{d/2}}e^{-{t |\xi|^2}}, \quad \xi\in \R^d\!.
$$
It follows that $(2\pi)^{d/2}\widehat{K_t} \widehat{K_s} = \widehat{K_{t+s}}$ for all $t,s>0$, and by Proposition \ref{prop:FT-convol} this implies  $H(t+s)f = H(t)H(s)f$ for all $f\in L^2(\R^d)$.
In particular, this identity holds for functions $f\in L^p(\R^d)\cap L^2(\R^d)$, and since we have already seen that the operators $H(t)$ are contractive on $L^p(\R^d)$ the identity $H(t+s)f = H(t)H(s)f$ extends to general functions $f\in L^p(\R^d)$, by the density
of $L^p(\R^d)\cap L^2(\R^d)$ in $L^p(\R^d)$.

Strong continuity of the semigroup is an immediate consequence of Theorem \ref{thm:approx-identity}.

\smallskip
{\em Step 2} -- We now prove that $A=\Delta$, the weak $L^p$-Laplacian, with equal domains.

We begin by proving the inclusion $\Dom(\Delta)\subseteq \Dom(A)$ along with the fact that
$$
  Af=\Delta f, \quad f\in \Dom(\Delta).
$$
First let $ f\in C_{\rm c}^\infty(\R^d)$. For all $t>0$ we have the pointwise identities
$$
 \frac{\partial}{\partial t}K_t = \Delta K_t,\qquad \frac{\partial}{\partial t}(K_t *f) = \Delta K_t * f,
$$
and therefore
\begin{equation}\label{eq:heat-diff}
  H(t)f-f = K_t*f -f  = \int_0^t \Delta K_s * f\ud s = \int_0^t \Delta H(s)f\ud s.
\end{equation}
Since we are assuming that $f\in C_{\rm c}^\infty(\R^d)$, all identities can be rigorously justified by
elementary calculus arguments.

Now let $f\in\Dom(\Delta)$ and set $g:=\Delta f\in L^p(\R^d)$. Since $K_s$ is smooth and decays sufficiently rapidly, and since convolution
commutes with weak derivatives, we have
$$
  \Delta H(s)f = \Delta (K_s*f) = K_s * \Delta f = H(s)g \ \hbox{ in } L^p(\R^d)\ \hbox{for each } s>0.
$$
For $0<\eps<t$, the mapping
$
s\mapsto H(s)f
$
is differentiable on $[\eps,t]$, and
$$
\frac{\ud}{\ud s}H(s)f=\Delta H(s)f=H(s)g.
$$
Consequently,
$$
H(t)f-H(\eps)f=\int_\eps^tH(s)g\ud s.
$$
Letting $\eps\downarrow0$ and using strong continuity of $H$ gives
$$
H(t)f-f=
\int_0^tH(s)g\ud s.
$$
Dividing by $t$ and letting $t\downarrow 0$, from the strong continuity of $s\mapsto H(s)$ on $L^p(\R^d)$ we obtain
$$
  \lim_{t\downarrow 0}\frac1t\int_0^t H(s)g\ud s = g \ \hbox{ in } L^p(\R^d).
$$
By the definition of the generator, this shows $f\in\Dom(A)$ and $Af=g=\Delta f$.

To prove the converse inclusion $\Dom(A)\subseteq \Dom(\Delta)$, we must
show that every $f\in\Dom(A)$ admits a weak $L^p$-Laplacian. Applying \eqref{eq:heat-diff} to a test function
$\phi\in C_{\rm c}^\infty(\R^d)$, and using the symmetry
$K_t(x)=K_t(-x)$ of the heat kernel together with Fubini's theorem, we obtain
\begin{align*}
 \int_{\R^d} (H(t)f(x)-f(x))\phi(x)\ud x
 &= \int_{\R^d} f(x)(H(t)\phi(x)-\phi(x))\ud x \\
 &= \int_0^t \int_{\R^d} f(x)H(s)\Delta\phi(x)\ud x\ud s \\
 &= \int_0^t \int_{\R^d} H(s)f(x)\Delta\phi(x)\ud x\ud s .
\end{align*}
Dividing by $t$ and passing to the limit as $t\downarrow 0$, and using the assumption $f\in \Dom(A)$ and the strong continuity of $H$ on $L^p(\R^d)$, we obtain the identity
$$ \int_{\R^d} Af(x)\phi(x)\ud x= \int_{\R^d} f(x) \Delta\phi(x)\ud x,$$ which shows that $f$ admits a weak Laplacian, given by the function $Af\in L^p(\R^d)$.

By Theorem \ref{thm:H1-Sob}, for $p=2$ we have
\begin{align}\label{eq:heat-domL2}
\Dom(A)=\Dom(\Delta) = W^{2,2}(\R^d).
\end{align}

\begin{remark}
For $1<p<\infty$ one has the analogous equality
$$ \Dom(A) =\Dom(\Delta) = W^{2,p}(\R^d),$$
but this is highly nontrivial and depends on the $L^p$-boundedness of the Riesz transforms (see  the Notes to Chapter \ref{ch:operators}).

For $d=1$ there is a more elementary argument which also works
for $p=1$. Let $1\le p<\infty$. Since
$W^{2,p}(\R)\subseteq\Dom(\Delta)$, it suffices to prove the
reverse inclusion. Suppose that $f\in L^p(\R)$ admits a weak
second derivative $f''\in L^p(\R)$. By
Theorem \ref{thm:lower-order}, $f$ has a weak derivative
$f'\in L^1_{\rm loc}(\R)$, whose weak derivative is $f''$.
Thus $f$ and $f'$ have locally absolutely continuous
representatives. We claim that $f'\in L^p(\R)$. To prove this, we apply
the fundamental theorem of calculus twice to obtain
$$
f(x+1)-f(x)
=f'(x)+\int_0^1(1-s)f''(x+s)\ud s.
$$
The integral on the right-hand side
belongs to $L^p(\R)$. Indeed, this follows from
$$
\Bigl\n \int_0^1(1-s)f''(\cdot+s)\ud s\Bigr\n_p
\le \int_0^1(1-s)\n f''(\cdot+s)\n_p\ud s
=\tfrac12\n f''\n_p.
$$
Since $x\mapsto f(x+1)-f(x)$ belongs to $L^p(\R)$ as well, this proves the claim.
Consequently, $f\in W^{2,p}(\R)$ as required.
\end{remark}

\begin{remark}
There is a
slightly different route to the identification $A = \Delta$ for $p=2$ which depends on the fact that each of the operators
$H(t)$ is a Fourier multiplier operator associated with the multiplier $m_t(\xi) = \exp(-t |\xi|^2)$.
Defining $\wt H(t)g:= m_t g$ we obtain a multiplication semigroup $\wt H$ on $L^2(\R^d)$ which is strongly continuous
and whose generator $\wt A$ is given by
\begin{align*} \Dom(\wt A) &= \big \{g\in L^2(\R^d):\, \xi\mapsto |\xi|^2 g(\xi) \in L^2(\R^d)\big\} = \calF(H^2(\R^d)), \\
 \wt A g (\xi) & = - |\xi|^2 g(\xi), \quad g\in\Dom(\wt A), \ \xi\in\R^d;
\end{align*}
this follows from the results proved in Section \ref{subsec:multsgr}.
This semigroup is related to the heat semigroup through the identity
  $$ H(t) = \calF^{-1}\circ \wt H(t)\circ \calF\!, \quad t\ge 0,$$
from which it follows that a function $f\in L^2(\R^d)$
belongs to the domain of the generator $A$ of $H$
  if and only if $\calF f = \wh f$ belongs to the domain of the generator $\wt A$ of $\wt H$, and that
  $$ Af =  \calF^{-1}\circ \wt A\circ \calF f$$
in this case. As we have seen, this is the case if and only if $f\in H^2(\R^d)$. Since $H^2(\R^d) = W^{2,2}(\R^d)$ up to equivalence of norms, this implies \eqref{eq:heat-domL2}.
\end{remark}

\paragraph{The Heat Semigroup on Bounded Domains}\label{par:heat-sgr-D}\index{heat!semigroup, on domains}

Let $D\subseteq \R^d$ be a nonempty bounded open set.

\begin{proposition}\label{prop:LaplacianDsgr}
The Dirichlet and Neumann Laplacians on $L^2(D)$ generate analytic
$C_0$-semigroups $S_{\rm Dir}$ and $S_{\rm Neum}$ which are contractive on every sector of angle less
than $\frac12\pi$. For every $t\ge0$, both $S_{\rm Dir}(t)$ and $S_{\rm Neum}(t)$ are selfadjoint contractions.
\end{proposition}
\begin{proof}
Everything except selfadjointness of the real-time operators follows
from the Lumer--Phillips theorem.

For $\lambda>0$, the resolvent $R(\lambda,\Delta)$ is selfadjoint.
Euler's formula gives, for $t>0$,
$$
S(t)x
=
\limn
\Bigl(
\frac ntR(\frac nt,\Delta)
\Bigr)^nx.
$$
Every operator in this sequence is selfadjoint. Therefore, for all
$x,y\in L^2(D)$,
\begin{align*}
\iprod{S(t)x}{y}
&=
\limn
\iprod{
(\frac ntR(\frac nt,\Delta))^nx
}{y}
=
\limn
\iprod{x}{
(\frac ntR(\frac nt,\Delta))^ny
}
=
\iprod{x}{S(t)y}.
\end{align*}
Thus $S(t)$ is selfadjoint.

To check the conditions of the Lumer--Phillips theorem, let us denote the Dirichlet and Neumann Laplacians by $\Delta$. The operator $-\Delta$ is positive and selfadjoint by Theorem \ref{thm:Laplace-sa}, and therefore $I-\Delta$ is injective.
Dualising and using selfadjointness, we conclude in turn that $I-\Delta$ has dense range. This verifies the first condition of the Lumer--Phillips theorem; the second follows immediately from the positivity of $-\Delta$.
\end{proof}

Alternatively, Proposition \ref{prop:LaplacianDsgr} can be deduced from the
spectral theorem for selfadjoint operators (as such the result is a special
case of Theorem \ref{thm:normal-sgr}, where further details are
provided). If $P$ denotes the projection-valued measure associated with the
positive selfadjoint operator $-\Delta$, this gives the representation
(see Example \ref{ex:HSgr-ST})
$$
  S(z)=\int_{[0,\infty)} e^{-zt}\ud P(t),\quad \Re z>0.
$$

\paragraph{The Heat Equation as a Gradient Flow}
Let $D\subseteq\R^d$ be a nonempty bounded open set, and put
$H:=L^2(D)$ and $V:=H_0^1(D)$. Consider the Dirichlet energy
$$
E(v):=\frac12\int_D|\nabla v|^2\ud x,\quad v\in V.
$$
Its first variation at $v\in V$ in the direction $h\in V$ is
$$
\left.\frac{\ud}{\ud\eps}\right|_{\eps=0}E(v+\eps h)
=\Re\iprod{\nabla v}{\nabla h},
$$
where $\eps$ is real and the inner product on the right is taken
in $L^2(D;\K^d)$. We say that $E$ has an {\em $H$-gradient at $v$}
if there exists $g\in H$ such that
$$
\Re\iprod{\nabla v}{\nabla h}
=\Re\iprod{g}{h}_H,\quad h\in V.
$$
Since $V$ is dense in $H$, such a vector is unique; we denote it
by $\operatorname{grad}_H E(v)$. This definition uses the inner
product of the ambient space $H$, although the energy is defined
only on $V$.

The weak definition of the Laplacian and the density of
$C_{\rm c}^\infty(D)$ in $V$ show that this representing vector
exists precisely when $\Delta v\in L^2(D)$, and is then
$g=-\Delta v$. Thus, putting $A:=-\Delta_{\rm Dir}$, we have
$$
\Dom(\operatorname{grad}_H E)
=\Dom(A)
=\{v\in H_0^1(D):\,\Delta v\in L^2(D)\},
\qquad
\operatorname{grad}_H E(v)=Av.
$$
In particular, the $H$-gradient is an unbounded operator.
A {\em gradient flow} prescribes the velocity to be the negative
gradient. For the Dirichlet energy this gives the heat equation
$$
u'(t)=-\operatorname{grad}_H E(u(t))
=-Au(t)=\Delta_{\rm Dir}u(t),
$$
with the Dirichlet boundary condition expressed by $u(t)\in V$.

Put $S:=S_{\rm Dir}$. For every $u_0\in H$, the function
$u(t):=S(t)u_0$ belongs to
$C([0,\infty);H)\cap C^1((0,\infty);H)$, satisfies $u(0)=u_0$,
and, by analyticity, takes values in $\Dom(A)$ for $t>0$.
Thus the heat semigroup produces a gradient flow for every
initial value in $L^2(D)$.
This solution is unique in the stated class. Indeed, the difference
$w$ of two solutions with the same initial value satisfies
$$
\frac12\frac{\ud}{\ud t}\n w(t)\n_2^2
=-\Re\iprod{Aw(t)}{w(t)}_H
=-\n\nabla w(t)\n_2^2\le0,\quad t>0.
$$
Hence $\n w(t)\n_2\le\n w(s)\n_2$ for $0<s<t$; letting
$s\downarrow0$ gives $w(t)=0$.

For positive times, analyticity also gives continuous
differentiability with values in $\Dom(A)$ equipped with its
graph norm. Since
$E(v)=\frac12\iprod{Av}{v}_H$ for $v\in\Dom(A)$,
selfadjointness of $A$ permits differentiation and gives
$$
\frac{\ud}{\ud t}E(u(t))
=\Re\iprod{Au(t)}{u'(t)}_H
=-\n u'(t)\n_2^2,\quad t>0.
$$
Consequently,
$$
E(u(t))+\int_s^t\n u'(r)\n_2^2\ud r=E(u(s)),
\quad 0<s\le t.
$$
To include time zero, suppose that $u_0\in V$. By Kato's theorem
(Theorem \ref{thm:Kato}),
$$
V=\Dom(A^{1/2}),\qquad
E(v)=\frac12\n A^{1/2}v\n_2^2,\quad v\in V.
$$
The spectral representation of $S$ gives
$A^{1/2}S(s)u_0=S(s)A^{1/2}u_0$, so strong continuity implies
$E(u(s))\to E(u_0)$ as $s\downarrow0$. We therefore obtain
\begin{align*}
E(u(t))+\int_0^t\n u'(r)\n_2^2\ud r=E(u_0),
\quad t\ge0.
\end{align*}
Thus the decrease in Dirichlet energy is exactly the integral
of the squared speed of the flow in $L^2(D)$.

Finally, for a fixed $f\in H$, consider
$$
E_f(v):=E(v)-\Re\iprod{v}{f}_H,\quad v\in V.
$$
The same first-variation calculation gives
$\operatorname{grad}_H E_f(v)=Av-f$ on $\Dom(A)$.
Its gradient flow equation is the inhomogeneous heat equation
$$
u'(t)=-Au(t)+f=\Delta_{\rm Dir}u(t)+f,
$$
whose solution with initial value $u_0$ is
$$
u(t)=S(t)u_0+\int_0^tS(t-r)f\ud r.
$$
The stationary equation $Au=f$ is the Poisson problem with
Dirichlet boundary conditions. By Theorem \ref{thm:Dirprinc},
its weak solution is the unique minimiser of $E_f$ over $V$.
Thus the elliptic problem describes the equilibrium of the
evolution associated with the same energy.

\paragraph{Positivity}

The Dirichlet and Neumann heat semigroups on $L^2(D)$ are positivity preserving,\index{heat!equation, positivity}
that is, they map nonnegative functions to nonnegative functions.  From the physics point of view it is natural to expect that heat semigroups should have this property, as they are meant to describe the time evolution of heat distributions. Positivity of the heat semigroup on $L^2(\R^d)$ is evident from the explicit representation through convolution with the heat kernel.

\begin{theorem}[Positivity]\label{thm:heat-pos} Let $D$ be bounded. Then the $C_0$-semigroups on $L^2(D)$ generated by $\Delta_{\rm Dir}$ and $\Delta_{\rm Neum}$ are positivity preserving.
\end{theorem}
\begin{proof}
Let $A$ denote the Dirichlet or Neumann Laplacian on $L^2(D)$ and let $S$ be the
$C_0$-contraction semigroup generated by $A$ on $L^2(D)$.
We must prove that $S(t)f\ge 0$ for all $t\ge 0$ whenever $f\in L^2(D)$ satisfies $f\ge 0$. In what follows we fix such a function $f$.

\smallskip{\em Step 1} -- We first prove that for all real-valued
$f\in L^2(D)$ satisfying $f\ge 0$ and all $\la>0$ we have
$u:= R(\la,A) f \ge 0$. We do this for the Dirichlet Laplacian; the proof
for the Neumann Laplacian is the same provided $H_0^1(D)$ is replaced by
$H^1(D)$.

Let $\aa$ be the form associated with $-A$. We first note that $u$ is
real-valued. Indeed, $u\in\Dom(A)$ is characterised by
$$
        \la\iprod{u}{v}_{L^2(D)}+\aa(u,v)=\iprod{f}{v}_{L^2(D)},
        \quad v\in\Dom(\aa).
$$
Since the form $\aa$ is real, in the sense that
$\aa(\overline u,\overline v)=\overline{\aa(u,v)}$, and since $f$ is
real-valued, taking complex conjugates shows that $\overline u$ satisfies
the same identity. By uniqueness of the resolvent equation,
$\overline u=u$. Thus $u$ is real-valued.

Since $u\in \Dom(A)\subseteq \Dom(\aa) = H_0^1(D)$ (see Section \ref{sec:examples-forms}), also $u^-:=(-u)\vee 0$ belongs to $H_0^1(D)$. Moreover,
$$
\la \iprod{u}{u^-}_{L^2(D)}+\aa(u,u^-)=\iprod{f}{u^-}_{L^2(D)}.
$$
Now $\iprod{u}{u^-}_{L^2(D)}=-\|u^-\|_{L^2(D)}^2$ and $\aa(u,u^-)= -\aa(u^-,u^-)$,
while $\iprod{f}{u^-}_{L^2(D)}\ge 0$.
Hence
$$
-\la \|u^-\|_{L^2(D)}^2-\aa(u^-,u^-)\ge 0.
$$
Since $\aa(v,v)\ge 0$ for all $v\in H_0^1(D)$, it follows that
$$
\la \|u^-\|_{L^2(D)}^2+\aa(u^-,u^-)=0,
$$
and therefore $u^-=0$. Thus $u\ge 0$ almost everywhere on $D$.

\smallskip
{\em Step 2} -- The positivity of the operators $S(t)$ follows from the result of Step 1 via the Euler formula
(Theorem \ref{thm:Euler})
\[S(t)f = \limn \Bigl(\frac{n}{t} R(\frac{n}{t},A)\Bigr)^n f, \quad f\in L^2(D).\qedhere\]
\end{proof}

\paragraph{\texorpdfstring{$L^p$}{L\textasciicircum p}-Theory}

We have seen that the Laplace operator $\Delta$ generates a $C_0$-semi\-group of contractions on $L^p(\R^d)$ for all $1\le p<\infty$.
For bounded open subsets $D$ of $\R^d\!$, up to this point we have only considered the analogues of this semigroup on the space $L^2(D)$. We next prove that the Dirichlet and Neumann Laplacians also generate $C_0$-semigroups of contractions on the space $L^p(D)$
for $1\le p<\infty$.
This will be derived from an abstract result on $L^p$-boundedness of doubly submarkovian operators which we discuss first.

Let $(\Om,\calF\!,\mu)$ be a finite measure space.
A bounded operator $T$ on $L^2(\Om)$ is called {\em doubly submarkovian}\index{doubly submarkovian}
if it has the following properties:
\begin{enumerate}[label={\rm(\roman*)}, leftmargin=*]
 \item $T f\ge 0$ for all $f\ge 0$;
 \item $T\one\le \one$ and $T^\star\one  \le \one $.
\end{enumerate}
Such operators enjoy the following extension property.

\begin{theorem}[$L^p$-Boundedness of doubly submarkovian operators]\label{thm:doublysubMarkovian} Let $(\Om,\calF\!,\mu)$ be a finite measure space and let
$T$ be a doubly submarkovian operator on $L^2(\Om)$. Then
for all $1\le p\le \infty$,
the restriction of $T$ to $L^2(\Om)\cap L^p(\Om)$ has a unique
extension to a contraction on $L^p(\Om)$.
\end{theorem}
\begin{proof}
We first prove that for all $f\in L^2(\Om)$ we have
$|Tf|\le T|f|$ almost everywhere.
For real-valued $f$ this follows from
$-|f|\le f\le|f|$ and positivity.
For complex-valued $f$, let $\alpha$ range over a countable
dense subset of the unit circle. Since $T$ maps real-valued
functions to real-valued functions, positivity gives
$$
\Re(\alpha Tf)=T(\Re(\alpha f))\le T|f|
$$
outside one common null set. Taking the supremum over
$\alpha$ proves the inequality.

Since $\mu$ is finite, we have $\one\in L^2(\Om)$, and for
$f\in L^2(\Om)$ we obtain
\begin{align*}
\n Tf\n_1 = \n \,|Tf|\,\n_1
\le\iprod{T|f|}{\one}
=\iprod{|f|}{T^\star\one}
\le\iprod{|f|}{\one}
=\n f\n_1.
\end{align*}
By density, $T$ therefore extends uniquely to a contraction
on $L^1(\Om)$.

For $f\in L^\infty(\Om)\subseteq L^2(\Om)$, positivity gives
$$
|Tf|\le T|f|
\le\n f\n_\infty T\one
\le\n f\n_\infty\one.
$$
Thus $T$ restricts to a contraction on $L^\infty(\Om)$.
The Riesz--Thorin interpolation theorem
(Theorem \ref{thm:RieszThorin}) now gives the assertion
for $1<p<\infty$.
\end{proof}

Turning to the $L^p$-boundedness of the heat semigroup,
we begin with the case of Neumann boundary conditions.

\begin{theorem}[$L^p$-Boundedness, Neumann boundary conditions]\label{thm:heat-Neumann-Lp}
Let $D$ be bounded and let $S_{\rm Neum}$ denote the $C_0$-semigroup generated by $\Delta_{\rm Neum}$ on $L^2(D)$.
For all $1\le p<\infty$, the restriction of $S_{\rm Neum}$ to $L^2(D)\cap L^p(D)$ uniquely extends to a $C_0$-semigroup of positivity-preserving contractions on $L^p(D)$.
\end{theorem}
\begin{proof}
By Proposition \ref{prop:LaplacianDsgr} and Theorem \ref{thm:heat-pos}, for all $t\ge 0$ the operator $S_{\rm Neum}(t)$ is selfadjoint and positivity preserving.
From $\Delta_{\rm Neum}\one = 0$ it follows that $S_{\rm Neum}^\star(t)\one = S_{\rm Neum}(t)\one=\one$ and therefore the operators $S_{\rm Neum}(t)$ are doubly submarkovian. Applying Theorem \ref{thm:doublysubMarkovian}, we find that for all $1\le p<\infty$ and $t\ge 0$
the restriction of $S_{\rm Neum}(t)$ to $L^2(D)\cap L^p(D)$ has a unique
extension, also denoted by $S_{\rm Neum}(t)$, to a contraction on $L^p(D)$.

By H\"older's inequality, the strong continuity of $S_{\rm Neum}$ on $L^2(D)$ implies that for all $1\le p< 2$ and $f\in L^2(D)$ we have
$$\n S_{\rm Neum}(t)f-f\n_p \le |D|^{1/r}\n S_{\rm Neum}(t)f-f\n_2\to 0$$ as $t\downarrow 0$, where $\frac12+\frac1r = \frac1p$.
Since $\n S_{\rm Neum}(t)\n_p\le 1$ for all $t\ge 0$, the density of $L^2(D)$ in $L^p(D)$ implies that the semigroup property and strong continuity extend to all $f\in L^p(D)$.

For $2<p<\infty$, first let $f\in L^2(D)\cap L^\infty(D)$. By Lemma \ref{lem:Lp-logconvex},
$$
        \|S_{\rm Neum}(t)f-f\|_p
        \le
        \|S_{\rm Neum}(t)f-f\|_2^{2/p}
        \|S_{\rm Neum}(t)f-f\|_\infty^{1-2/p}.
$$
The first factor tends to $0$, while the second is bounded by
$2\|f\|_\infty^{1-2/p}$, because the operators $S_{\rm Neum}(t)$ are contractions
on $L^\infty(D)$. Hence $S_{\rm Neum}(t)f\to f$ in $L^p(D)$. Since
$L^2(D)\cap L^\infty(D)$ is dense in $L^p(D)$ and the operators
$S_{\rm Neum}(t)$ are contractions on $L^p(D)$, strong continuity follows
for all $f\in L^p(D)$.

Positivity on $L^p(D)$ follows from the positivity on $L^2(D)$ by a density argument.
\end{proof}

Our next aim is to prove the following analogue of Theorem \ref{thm:heat-Neumann-Lp} for the Dirichlet heat semigroup.

\begin{theorem}[$L^p$-Boundedness, Dirichlet boundary conditions]\label{thm:heat-Dirichlet-Lp}
Let $D$ be a bounded open subset of $\R^d$ and let $S_{\rm Dir}$ denote the $C_0$-semigroup generated by $\Delta_{\rm Dir}$ on $L^2(D)$.
For all $1\le p<\infty$, the restriction of $S_{\rm Dir}$ to $L^2(D)\cap L^p(D)$ uniquely extends to a $C_0$-semigroup of positivity-preserving contractions on $L^p(D)$.
\end{theorem}

The heart of the matter is to prove the following resolvent inequality.

\begin{lemma}\label{lem:dominDirNeum} For all $0\le f\in L^2(D)$ and $\la>0$ we have
 $$ 0\le R(\la,\Delta_{\rm Dir})f \le  R(\la,\Delta_{\rm Neum})f.$$
\end{lemma}

\begin{proof}
Fix $\la>0$ and $0\le f\in L^2(D)$.
Then $0\le u:= R(\la, \Delta_{\rm Dir})f \in \Dom(\Delta_{\rm Dir})$ and $0\le v:= R(\la, \Delta_{\rm Neum})f \in
\Dom(\Delta_{\rm Neum})$ and
\begin{align}\label{eq:Dir-Neum} \la u- \Delta_{\rm Dir} u = f = \la v- \Delta_{\rm Neum} v.
\end{align}
The lemma will be proved by showing that this implies $u\le v$.

Fix a nonnegative test function $\phi\in C_{\rm c}^\infty(D)$.
Multiplying \eqref{eq:Dir-Neum} on both sides by $\phi$ and integrating by parts,  we arrive at
\begin{align}\label{eq:dominDirNeum}
 \la \int_D u\phi \ud x +\int_D \nabla u \cdot\nabla\phi \ud x
 = \la \int_D v\phi \ud x +\int_D \nabla v \cdot \nabla\phi \ud x.
\end{align}
We claim that this equality extends to all nonnegative functions $\phi\in H_0^1(D)$. Indeed, if $\phi_n\to \phi$ in $H_0^1(D)$ with $\phi_n\in C_{\rm c}^\infty(D)$ for all $n\ge 1$, then $\phi_n^+\to \phi$ in $H_0^1(D)$ by Theorem \ref{thm:H1-lattice}.
Since each $\phi_n^+$ has compact support, mollification with a nonnegative compactly supported smooth mollifier allows us to approximate $\phi_n^+$ by nonnegative test functions in $H_0^1(D)$ as in Proposition \ref{prop:local-approx}.
Applying \eqref{eq:dominDirNeum} to these test functions and taking limits, we obtain the claim.

Next we claim that $(u-v)^+$ belongs to $H_0^1(D)$, the point here being that $u\in H_0^1(D)$ but $v\in H^1(D)$.
To prove the claim, let $u_k\to u$ in $H_0^1(D)$ with $u_k\in C_{\rm c}^\infty(D)$. Since $v\ge 0$, for each $k\ge 1$ the function
$(u_k-v)^+$ is supported in the compact support of $u_k$ and therefore it belongs to $H_0^1(D)$ by Proposition \ref{prop:local-approx}.
By another application of Theorem \ref{thm:H1-lattice} it then follows that
$(u-v)^+ = \limk (u_k-v)^+$  belongs to $H_0^1(D)$ as claimed.

By \eqref{eq:dominDirNeum}, which we may now apply to $\phi = (u-v)^+$,
\begin{align*}
 \la \int_D u(u-v)^+\ud x +\int_D \nabla u \cdot \nabla (u-v)^+ \ud x
 = \la \int_D v(u-v)^+\ud x +\int_D \nabla v \cdot \nabla (u-v)^+ \ud x.
\end{align*}
As a consequence,
\begin{align*}
 \la \int_D (u-v)^{+2}\ud x & = \la \int_D (u-v)(u-v)^+\ud x
 \\ & = \int_D \nabla(v-u)\cdot\nabla(u-v)^+\ud x
  = - \int_D |\nabla(u-v)^+|^2\ud x  \le 0,
\end{align*} where in the last step we argue as in the proof of Theorem \ref{thm:heat-pos}.
This implies that $(u-v)^+=0$, that is, $u\le v$.
\end{proof}

\begin{proof}[Proof of Theorem \ref{thm:heat-Dirichlet-Lp}]
By Lemma \ref{lem:dominDirNeum} and Euler's formula
(Theorem \ref{thm:Euler}), for all $0\le f\in L^2(D)$ and $t>0$ we have
$$
 S_{\rm Dir}(t)f
 = \limn \Bigl(\frac{n}{t} R(\frac{n}{t},\Delta_{\rm Dir})\Bigr)^n f
 \le
 \limn \Bigl(\frac{n}{t} R(\frac{n}{t},\Delta_{\rm Neum})\Bigr)^n f
 = S_{\rm Neum}(t)f .
$$
In particular, this implies
$S_{\rm Dir}(t)\one \le S_{\rm Neum}(t)\one=\one$.
Together with positivity and selfadjointness, this implies that the operators
$S_{\rm Dir}(t)$ are doubly submarkovian.
The proof can now be finished along the lines of that of Theorem
\ref{thm:heat-Neumann-Lp}.
\end{proof}

\subsection{The Poisson Semigroup}\label{subsec:Poisson-sgr}\index{Poisson!semigroup}

Let $1\le p<\infty$.
For $t>0$ we define the operator $P(t)$ on $L^p(\R^d)$ by convolution with the {\em Poisson kernel}\index{Poisson!kernel}
\begin{align}\label{eq:PoissonRd}
p_t(x) = \frac{c_d t}{(t^2 + |x|^2)^{\frac12(d+1)}},\quad t>0, \ x \in\R^d\!,
\end{align}
where $c_d = {\Gamma(\frac12(d+1))}/{\pi^{\frac12(d+1)}}$
with $\Gamma(t) = \int_0^\infty x^{t-1}e^{-x}\ud x$ the Euler Gamma function.
The change of variables $y = x/t$ gives the norm estimate
\begin{align*}
\n p_t \n_{L^1(\R^d)}
& = c_d \int_{\R^d} \frac{1}{(1 + |y|^2)^{\frac12(d+1)}}\ud y
\\* & = c_d \sigma_{d-1} \int_0^\infty \frac{1}{(1 + r^2)^{\frac12(d+1)}}r^{d-1}\ud r
  = c_d \sigma_{d-1} \cdot \frac12\sqrt{\pi} \frac{\Gamma(\frac12d)}{\Gamma(\frac12(d+1))} = 1,
\end{align*}
where $\sigma_{d-1} = {2\pi^{d/2}}/{\Gamma(\frac12d)}$ is the surface area of the unit sphere in $\R^d\!$.
Young's inequality guarantees that the operators $P(t)$ defined by $P(0):=I$ and
$$ P(t)f := p_t * f, \quad t>0, \ f\in L^p(\R^d), $$
are well defined and contractive on $L^p(\R^d)$.
For $d=1$, the formula \eqref{eq:PoissonRd} takes the simpler form
$$ p_t(x) = \frac1\pi\frac{t}{t^2 + x^2}, \quad t>0, \ x\in \R.$$

To see that the operators $P(t)$ satisfy the semigroup property we will show that the Fourier transform of $p_t$ is given by
\begin{align}\label{eq:Poisson-FT}
 \wh{p_t}(\xi)  = \frac1{(2\pi)^{d/2}}\exp(- t|\xi|).
\end{align}
To prove this identity we compute the inverse Fourier transform of the exponential on the right-hand side.
By a standard contour integration argument, for $\gamma>0$ we have
$$ e^{-\gamma} = \frac1\pi\int_{\R} \frac{e^{i\gamma y}}{1+y^2}\ud y.$$
Writing $\frac{1}{1+y^2} = \int_0^\infty e^{-(1+y^2)u}\ud u$, using Fubini's theorem to interchange the order of integration, and using Lemma \ref{lem:Gauss} and a change of variables to evaluate the inner integral, we obtain
\begin{align*}
e^{-\gamma}
& = \frac1\pi\int_{\R} \int_0^\infty e^{i\gamma y} e^{-(1+y^2)u}\ud u\ud y
\\ & = \frac1\pi\int_0^\infty e^{-u} \Bigl(\int_{\R}  e^{i\gamma y} e^{-y^2 u}\ud y\Bigr)\ud u
 = \frac1{\sqrt{\pi}} \int_0^\infty \frac{e^{-u}}{\sqrt u} e^{-\gamma^2/4u}\ud u.
\end{align*}
We apply this with $\gamma = t|\xi|$. Using Fubini's theorem, Lemma \ref{lem:Gauss}, and another substitution, we obtain
\begin{align*}
\ &  \frac1{(2\pi)^{d/2}}  \int_{\R^d} \frac1{(2\pi)^{d/2}}\exp(- t|\xi|) \exp( i x\cdot \xi)\ud \xi
\\ & \qquad = \frac1{(2\pi)^d}\int_{\R^d} \frac1{\sqrt{\pi}}
 \int_0^\infty \frac{e^{-u}}{\sqrt u} e^{- t^2|\xi|^2/4u} \exp(i x\cdot \xi)\ud u \ud \xi
\\ & \qquad = \frac1{\sqrt{\pi}}\int_0^\infty  \frac{e^{-u}}{\sqrt u}\frac1{(2\pi)^d}\int_{\R^d}
 e^{- t^2|\xi|^2/4u} \exp( i x\cdot \xi)\ud \xi\ud u
 \\ & \qquad =  \frac1{\sqrt{\pi}}  \int_0^\infty  \frac{e^{-u}}{\sqrt u}
 \big(\frac{u}{\pi t^2}\big)^{d/2}e^{-|x|^2u/t^2}\ud u
 \\ & \qquad = \frac{1}{\pi^{\frac12(d+1)}}\frac1{t^d} \int_0^\infty e^{-(t^2+|x|^2)u/t^2}u^{\frac12(d-1)}\ud u
 \\ & \qquad = \frac{1}{\pi^{\frac12(d+1)}} \frac{t}{(t^2+|x|^2)^{\frac12(d+1)}}\int_0^\infty  e^{-v} v^{\frac12(d-1)}\ud v
 \\ & \qquad = \frac{1}{\pi^{\frac12(d+1)}} \frac{t}{(t^2+|x|^2)^{\frac12(d+1)}}\Gamma(\frac12(d+1))
  = p_t(x).
\end{align*}
This completes the proof of \eqref{eq:Poisson-FT}. Thanks to this identity, for $t,s>0$ we obtain
\begin{align*}\wh{p_t * p_s}  = (2\pi)^{d/2}\wh{p_t} \wh{p_s} & = (2\pi)^{-d/2}\exp(- t|\xi|)\exp(- s|\xi|)
\\ & = (2\pi)^{-d/2}\exp(- (t+s)|\xi|) = \wh{p_{t+s}}
\end{align*} and therefore $p_t*p_s = p_{t+s}$.
It follows that for, say, $f\in C_{\rm c}(\R^d)$,
$$ P(t)P(s)f = p_t * (p_s* f) = (p_t * p_s)* f = p_{t+s}*f = P(t+s)f.$$
Since $C_{\rm c}(\R^d)$ is dense in $L^p(\R^d)$ for $1\le p<\infty$, this proves that $P(t)P(s) = P(t+s)$ for all $t,s>0.$
This identity trivially extends to $t,s\ge 0$.
Strong continuity of the family $P$ on $L^p(\R^d)$ for $1\le p<\infty$ is an immediate consequence of
Theorem \ref{thm:approx-identity}.

We will determine the generator $A$ of this semigroup for $p=2$. It will turn out that
$$ A = -(-\Delta)^{1/2}\!, \qquad \Dom(A) = H^1(\R^d).$$

Here, the square root can be defined either via the functional calculus of the selfadjoint operator
$-\Delta$ (see Proposition \ref{prop:pos-sqrt}), or via the following direct Fourier-analytic argument. Recall the definition of
$H^1(\R^d)$ as the space of all $f\in L^2(\R^d)$ for which
$$\xi \mapsto (1+|\xi|^2)^{1/2} \wh f(\xi)$$ belongs to $L^2(\R^d)$.
In view of the trivial inequality $|\xi| \le  (1+|\xi|^2)^{1/2}$, for all $f\in H^1(\R^d)$ the function
$ \xi \mapsto |\xi| \wh f(\xi)$ belongs to $L^2(\R^d)$.
Thus we can define an operator $(B, \Dom(B))$ by
$$ Bf := (\xi \mapsto |\xi| \wh f(\xi))\widecheck{\phantom{b}}, \qquad \Dom(B) = H^1(\R^d).$$
Note the formal analogy with the definition of Fourier multiplier operators; the only difference here is that the multiplier function $m(\xi) = |\xi|$ does not belong to $L^\infty(\R^d)$,
so a domain must be specified, as in Example \ref{ex:sa-FM}.

For $f\in H^2(\R^d)$ we similarly have
$$ -\Delta f = (\xi \mapsto |\xi|^2 \wh f(\xi))\widecheck{\phantom{b}},$$
so that $ Bf \in \Dom(B)$ and $B^2 f := B(Bf) = -\Delta f.$ This justifies the notation
\begin{align}\label{eq:domLaplhalf} (-\Delta)^{1/2} := B, \qquad \Dom((-\Delta)^{1/2}) = H^1(\R^d).
\end{align}

Since multiplication by the real-valued function $\xi\mapsto |\xi|$ is selfadjoint on its natural domain in $L^2(\R^d)$ and the Fourier--Plancherel transform is unitary, the operator $B$ is selfadjoint. Moreover,
$$ \iprod{Bf}{f} = \bigl(|\cdot|\wh f(\cdot)\big|\wh f(\cdot)\bigr) = \int_{\R^d} |\xi||\wh f(\xi)|^2\ud \xi \ge 0,$$
so $B$ is positive. Also,
$$\Dom(B^2)=\big\{f\in L^2(\R^d):\, \xi\mapsto |\xi|^2\wh f(\xi)\in L^2(\R^d)\big\}=H^2(\R^d),$$
and $B^2f=-\Delta f$ for all $f\in H^2(\R^d)$. Thus $B$ is a positive selfadjoint square root of $-\Delta$, and the uniqueness part of Proposition \ref{prop:pos-sqrt} implies that $B$ coincides with the
positive square root of $-\Delta$ obtained by means of the functional calculus of $-\Delta$.

In dimension $d=1$ the identity $|\xi| = -i\,{\rm sign}(\xi) \cdot  i\xi$ implies that
$$ (-{\rm d}^2/{\rm d}x^2)^{1/2} = H \circ \frac{\rm d}{{\rm d}x},
$$
where $H$ is the Hilbert transform\index{Hilbert transform!and the Poisson semigroup} (which, as we recall from Section \ref{sec:HT}, is the Fourier multiplier operator
corresponding to the multiplier $\xi \mapsto -i\,{\rm sign}(\xi) $).

\begin{theorem}[Poisson semigroup] The Poisson semigroup on $L^2(\R^d)$ is generated by the selfadjoint operator $-(-\Delta)^{1/2}$.
\end{theorem}

\begin{proof}
 We start by noting that a function $f\in L^2(\R^d)$ belongs to $H^1(\R^d)$ if and only if
  $\xi \mapsto |\xi| \wh f(\xi)$ belongs to $L^2(\R^d)$. The `only if' part has already been noted, and the
  `if' part follows similarly from the inequality $(1+|\xi|^2)^{1/2} \le 1+ |\xi|$.

  On $L^2(\R^d)$ we now define the multiplication semigroup $Q$ by
  $$ Q(t)g(\xi) := e^{- t|\xi|}g(\xi)$$ for $t\ge 0$ and $g\in L^2(\R^d)$.
  As we have shown in Section \ref{subsec:multsgr}, this is a $C_0$-semigroup whose generator $(C,\Dom(C))$ is given by
  $$ Cg(\xi) = - |\xi|g(\xi)$$ for $ g\in \Dom(C) = \{g\in L^2(\R^d):\, \xi\mapsto  |\xi|g(\xi) \in L^2(\R^d)\}.$
  Evidently,
  $$ P(t) = \calF^{-1}\circ Q(t)\circ \calF\!, \quad t\ge 0,$$
  from which it follows that a function $f\in L^2(\R^d)$ belongs to the domain of the generator $A$ of $P$
  if and only if $\calF f = \wh f$ belongs to the domain of the generator $C$ of $Q$, and that
  $$ Af =  \calF^{-1}\circ C\circ \calF f$$
  in this case. By the observation at the beginning of the proof and \eqref{eq:domLaplhalf}, $\calF f$ belongs to $\Dom(C)$ if and only if
  $f\in H^1(\R^d) = \Dom((-\Delta)^{1/2})$, and in that case we have
  $$-(-\Delta)^{1/2} f = (\xi \mapsto -|\xi| \wh f)\widecheck{\phantom{b}} =  \calF^{-1}\circ C\circ \calF f.$$
  These considerations prove that $A = -(-\Delta)^{1/2}$ with equality of their domains
  $\Dom(A)=\Dom((-\Delta)^{1/2})=H^1(\R^d)$.
\end{proof}

The Poisson
semigroup can be recovered from the heat semigroup $S$ on $L^2(\R^d)$ by
$$
P(t)f=\frac{t}{2\sqrt\pi}\int_0^\infty
s^{-3/2}e^{-t^2/(4s)}S(s)f\ud s,
\quad t>0,\ f\in L^2(\R^d).
$$
To verify this identity, one applies the
Fourier--Plancherel transform and uses the scalar integral formula
in the proof of \eqref{eq:Poisson-FT}, with $\gamma=t|\xi|$ and
the substitution $s=t^2/(4u)$. This construction is called
{\em subordination}.\index{subordination}
More generally, the same formula relates $e^{-tN^{1/2}}$ and
$e^{-sN}$ for every positive selfadjoint operator $N$, by the
spectral theorem.
There is also an interesting connection with the wave equation, which will be elaborated in Section \ref{subsec:wave}.

The Poisson semigroup solves a boundary value problem on the
upper half-space $\R^d\times(0,\infty)$. For $f\in L^2(\R^d)$ put
$$
u(x,t):=(P(t)f)(x),\quad x\in\R^d,\ t>0.
$$
For every integer $k\ge0$ and $t>0$, the function
$\xi\mapsto |\xi|^k e^{-t|\xi|}\wh f(\xi)$ belongs to
$L^1(\R^d)\cap L^2(\R^d)$, since the multiplier is bounded
and square integrable. By
Fourier inversion and differentiation it therefore follows that $u$ is
smooth on the upper half-space and satisfies
$$
\partial_t^2u+\Delta_xu=0,\qquad
\lim_{t\downarrow0}\n u(\cdot,t)-f\n_2=0.
$$
Thus $u$ is a harmonic extension of $f$, with boundary values
understood in $L^2(\R^d)$. If $f\in H^1(\R^d)$, the generator
identity gives
$$
-\partial_tu(\cdot,t)=(-\Delta)^{1/2}P(t)f
\to (-\Delta)^{1/2}f \ \hbox{ as } t\downarrow0,
$$
with convergence in $L^2(\R^d)$.
Since the outward normal points in the negative $t$-direction,
this identifies $(-\Delta)^{1/2}$ as the
{\em Dirichlet-to-Neumann operator}:\index{Dirichlet-to-Neumann operator}
it maps the boundary value to the outward normal derivative of
its harmonic extension. In dimension $d=1$, for real-valued $f$,
the function
$$
(x,t)\mapsto (P(t)f)(x)+i(P(t)Hf)(x)
$$
is holomorphic as a function of $x+it$, where $H$ is the Hilbert
transform. This is the connection with harmonic conjugates
discussed in Problem \ref{prob:harmonic}.

\subsection{The Ornstein--Uhlenbeck Semigroup}\label{subsec:OU}

In this section we assume some elementary knowledge about Gaussian random variables.
Let\index{$Gamma1$@$\gamma$} $$\!\ud\gamma(x) = \frac1{{(2\pi)^{d/2}}}\exp\Bigl(-\frac12 |x|^2\Bigr)\ud x$$ denote the standard Gaussian measure\index{Gaussian!measure, standard} on $\R^d\!$.
For $t\ge 0$ and  $f\in C_{\rm c}(\R^d)$, define the operator ${OU}(t)$ on $L^2(\R^d\!,\gamma)$ by
$$ {OU}(t)f(x):= \int_{\R^d} f(e^{-t}x + \sqrt{1-e^{-2t}}y)\ud \gamma(y), \quad x\in \R^d\!.$$
We will prove that $OU$ is a $C_0$-semigroup on $L^p(\R^d\!,\gamma)$ for all $1\le p<\infty$.
This semigroup is known as the {\em Ornstein--Uhlenbeck semigroup}.\index{Ornstein--Uhlenbeck!semigroup}
It plays a central role in the so-called Malliavin calculus,\index{Malliavin calculus}
an infinite-dimensional Gaussian version of calculus which finds applications in, for example, the theory of
stochastic (partial) differential equations and mathematical finance. Interestingly, this semigroup also makes
its appearance in quantum field theory,
where its negative generator plays the role of the so-called bosonic
number operator. This point of view will be taken up in Section \ref{sec:SQ}.

Let us first show that each operator ${OU}(t)$ extends to a bounded operator on $L^p(\R^d\!,\gamma)$ of norm at most $1$.
By H\"older's inequality, for all $f\in C_{\rm c}(\R^d)$ we have
\begin{align*}
\| {OU}(t) f \|_{L^p(\R^d\!,\gamma)}^p &= \int_{\R^d} \Bigl| \int_{\R^d} f(e^{-t} x + \sqrt{1 - e^{-2t}} y)\ud \gamma(y) \Bigr|^p\ud\gamma(x) \\
&\leq \int_{\R^d} \int_{\R^d} \Bigl| f(e^{-t} x + \sqrt{1-e^{-2t}} y) \Bigr|^p\ud\gamma(x) \ud\gamma(y)
\\ & =  \mathbb{E}\big| f(e^{-t} X + \sqrt{1-e^{-2t}} Y)\big|^p\!,
\end{align*}
where $X$ and $Y$ are independent standard Gaussian random variables defined on some probability space and $\mathbb{E}$ is the expectation with respect to the probability measure. Since $(e^{-t})^2 +( \sqrt{1-e^{-2t}})^2 = 1$,
the random variable $e^{-t} X + \sqrt{1-e^{-2t}} Y$ is standard Gaussian again, so it is equal in distribution to $X$. Hence
$$  \mathbb{E}\big| f(e^{-t} X + \sqrt{1-e^{-2t}} Y)\big|^p =  \mathbb{E}|f(X)|^p = \int_{\R^d} |f(x)|^p \ud \gamma(x) = \|f\|_{L^p(\R^d\!,\gamma)}^p\!. $$
This proves that $\n {OU}(t)f\n_p \le \n f\n_p$.

We claim that $C_{\rm c}^\infty(\R^d)$ is dense in
$L^p(\R^d\!,\gamma)$. To this end we first approximate $f$ by
$\one_{B(0;r)}f$ for sufficiently large $r>0$.
On $B(0;r)$ the Gaussian density is bounded above and below
by positive constants, so convergence in the
$L^p(B(0;r),\gamma)$-norm is equivalent to convergence in the
$L^p(B(0;r))$-norm. Since $C_{\rm c}^\infty(B(0;r))$ is dense
in $L^p(B(0;r))$, the desired result follows by extending
these approximations by zero.

Combining the above two steps, we find that the operators ${OU}(t)$ extend uniquely to contractions on $L^p(\R^d\!,\gamma)$.
The same estimate, using Fubini's theorem, applies to every
$f\in L^p(\R^d\!,\gamma)$. It shows that
the defining integral exists for almost all $x$ and defines
an element of $L^p(\R^d\!,\gamma)$ independently of the chosen
representative. Hence the defining formula for $OU(t)$
extends to arbitrary $f\in L^p(\R^d\!,\gamma)$.

Next we prove that $OU$ is a $C_0$-semigroup on $L^p(\R^d\!,\gamma)$. It is clear that \ref{S1} holds. To prove the semigroup property \ref{S2}, let us first fix a function $f\in C_{\rm c}(\R^d)$. Then ${OU}(t) f\in C_{\rm b}(\R^d)$ and
\begin{align*}
{OU}(t) {OU}(s) f(x)
& = \int_{\R^d}  {OU}(s)f(e^{-t} x+ \sqrt{1-e^{-2t}} y)  \ud\gamma(y)
\\ & = \int_{\R^d} \int_{\R^d}f( e^{-s} ( e^{-t} x + \sqrt{1-e^{-2t}} y) + \sqrt{1-e^{-2s}} z) \ud\gamma(z) \ud\gamma(y)
\\ & = \mathbb{E} f(e^{-(t+s)} x + e^{-s} \sqrt{1-e^{-2t}} Y + \sqrt{1-e^{-2s}} Z),
\end{align*}
where $Y$ and $Z$ are independent standard Gaussians.
In view of the identity
$$(e^{-s} \sqrt{1-e^{-2t}})^2 + (\sqrt{1-e^{-2s}})^2 = 1 - e^{-2(t+s)},$$
the random variable $e^{-s} \sqrt{1-e^{-2t}} Y + \sqrt{1-e^{-2s}} Z$ is equal in distribution to a Gaussian random variable with variance $1-e^{-2(t+s)}$.
Therefore
\begin{align*}
 \mathbb{E} f(e^{-(t+s)} x & + e^{-s} \sqrt{1-e^{-2t}} Y + \sqrt{1-e^{-2s}} Z)
 \\ & = \mathbb{E}f(e^{-(t+s)} x + \sqrt{1-e^{-2(t+s)}} Y)
 \\ & = \int_{\R^d}  f(e^{-(t+s)} x + \sqrt{1 - e^{-2(t+s)}} y) \ud\gamma(y)
 = {OU}(t+s) f(x).
 \end{align*}
 This proves the identity ${OU}(t){OU}(s)f = {OU}(t+s)f$ for $f\in C_{\rm c}(\R^d)$.
 By the density of these functions in $L^p(\R^d\!,\gamma)$, the identity extends to every
$f\in L^p(\R^d\!,\gamma)$.

To prove the strong continuity property \ref{S3} we again first consider a function $f\in C_{\rm c}(\R^d)$. By H\"older's inequality,
 \begin{align*}
 \| {OU}(t) f - f\|_{L^p(\R^d\!,\gamma)}^p
 &\leq \int_{\R^d} \int_{\R^d} \bigl| f(e^{-t} x + \sqrt{1-e^{-2t}} y) - f(x)\bigr|^p \ud\gamma(x)\ud\gamma(y).
 \end{align*}
 The right-hand side tends to $0$ as $t\downarrow 0$ by dominated convergence. This gives strong
 continuity for functions $f\in C_{\rm c}(\R^d)$. The general case follows again by approximation, keeping in mind that
 the operators ${OU}(t)$ are all contractive.

For $t>0$, substituting $\sqrt{1-e^{-2t}}y=u$ and $e^{-t}x+u=v$
into the definition of $OU(t)$ and writing out the Gaussian density,
we obtain
\begin{equation}\label{eq:OU-Mehler}
\begin{aligned}
{OU}(t)f(x)
& = \frac1{(2\pi)^{d/2}}\Bigl(\frac1{1-e^{-2t}}\Bigr)^{d/2}
\int_{\R^d}f(e^{-t}x+u)
\exp\Bigl(-\frac12\frac{|u|^2}{1-e^{-2t}}\Bigr)\ud u
\\ & = \frac1{(2\pi)^{d/2}}\Bigl(\frac1{1-e^{-2t}}\Bigr)^{d/2}
\int_{\R^d}f(v)
\exp\Bigl(-\frac12\frac{|e^{-t}x-v|^2}{1-e^{-2t}}\Bigr)\ud v
\\ & = \int_{\R^d}M_t(x,v)f(v)\ud v.
\end{aligned}
\end{equation}
This represents ${OU}(t)$ as an integral operator with kernel
\begin{align*}
M_t(x,y)
= \frac1{(2\pi)^{d/2}}\Bigl(\frac1{1-e^{-2t}}\Bigr)^{d/2}
\exp\Bigl(-\frac12\frac{|e^{-t}x-y|^2}{1-e^{-2t}}\Bigr).
\end{align*}
The kernel $M_t$ is called the {\em Mehler kernel}\index{Mehler kernel}
at time $t$.
We can express the Mehler kernel in terms of the heat kernel as
$M_t(x,y)=K_{\frac12(1-e^{-2t})}(e^{-t}x-y)$, and therefore, for
$f\in C_{\rm c}(\R^d)$ and $t>0$, we have the pointwise identity
\begin{align}\label{eq:OU-heat}
{OU}(t)f(x)=H(\tfrac12(1-e^{-2t}))f(e^{-t}x),
\end{align}
where $H$ is the heat semigroup on $L^p(\R^d)$.

Writing the integral with respect to the Gaussian measure $\gamma(y)$,
we obtain
$$
OU(t)f(x)=\int_{\R^d}K_t(x,y)f(y)\ud\gamma(y),
$$
where, with $r=e^{-t}$,
$$
K_t(x,y)
=(1-r^2)^{-d/2}
\exp\Bigl(
-\frac{r^2|x|^2-2r x\cdot y+r^2|y|^2}{2(1-r^2)}
\Bigr).
$$
This kernel is symmetric in $x$ and $y$, hence $OU(t)$ is selfadjoint
on $L^2(\R^d\!,\gamma)$.
Denoting by $L_p$ the generator of $({OU}(t))_{t\ge0}$ on
$L^p(\R^d\!,\gamma)$, the observation just made implies that $L_2$
is symmetric. Since $1\in\varrho(L_2)$ by contractivity, Proposition
\ref{prop:symm-sgr-sa} implies that $L_2$ is selfadjoint.

Returning to the $L^p$-setting, we show next that
$C_{\rm c}^\infty(\R^d)$ is contained in $\Dom(L_p)$ and
\begin{align}\label{eq:O-gen}
L_pf(x)=\Delta f(x)-x\cdot\nabla f(x),
\quad f\in C_{\rm c}^\infty(\R^d).
\end{align}

Writing $a(t):=\frac12(1-e^{-2t})$ and defining the dilation operator
$(D_tg)(x):=g(e^{-t}x)$, by \eqref{eq:OU-heat} we have, for
$f\in C_{\rm c}^\infty(\R^d)$ and $t>0$,
$$
OU(t)f=D_tH(a(t))f.
$$
On $L^p(\R^d)$ we have $\n D_t\n=e^{dt/p}$, and $D_tg\to g$
in $L^p(\R^d)$ as $t\downarrow0$.

For $f\in C_{\rm c}^\infty(\R^d)$,
\begin{align*}
\frac1t(OU(t)f-f)
&=D_t\frac1t(H(a(t))f-f)+\frac1t(D_tf-f).
\end{align*}
Since $a(t)\to0$, $a(t)/t\to1$, and $f\in\Dom(\Delta)$,
$$
\frac1t(H(a(t))f-f)
=\frac{a(t)}t\frac1{a(t)}(H(a(t))f-f)\to\Delta f
$$
in $L^p(\R^d)$. Local uniform boundedness and strong convergence
of $D_t$ therefore imply
$$
D_t\frac{H(a(t))f-f}{t}\to\Delta f
$$
in $L^p(\R^d)$. Also, by dominated convergence,
$$
\frac1t(D_tf-f)
=-\int_0^1e^{-st}x\cdot\nabla f(e^{-st}x)\ud s
\to-x\cdot\nabla f
$$
in $L^p(\R^d)$.
It follows that
$$
\frac1t(OU(t)f-f)\to\Delta f-x\cdot\nabla f
$$
in $L^p(\R^d)$. Since
$\n g\n_{L^p(\R^d\!,\gamma)}
\le(2\pi)^{-d/(2p)}\n g\n_{L^p(\R^d)}$,
the same convergence holds in $L^p(\R^d\!,\gamma)$.
Hence $f\in\Dom(L_p)$ and \eqref{eq:O-gen} holds.

Let us now prove that $C_{\rm c}^\infty(\R^d)$ is a core for $L_p$.

\begin{proposition}\label{prop:OU-core}
For all $1\le p<\infty$, $C_{\rm c}^\infty(\R^d)$ is a core for
the generator $L_p$ of the Ornstein--Uhlenbeck semigroup on
$L^p(\R^d\!,\gamma)$.
\end{proposition}

\begin{proof}
We introduce the subspace
$$
Y:=C_{\rm b}^\infty(\R^d):=
\bigl\{f\in C^\infty(\R^d):\, \hbox{$f$ and all its derivatives are bounded}\bigr\}.
$$
We first prove that $Y\subseteq\Dom(L_p)$ and that every element
of $Y$ can be approximated by test functions in the graph norm
of $L_p$.

Fix $f\in Y$ and put $g:=\Delta f-x\cdot\nabla f$.
Since $f$ has bounded derivatives and the Gaussian measure has
finite moments of every order, $g\in L^p(\R^d\!,\gamma)$.
Choose $\chi\in C_{\rm c}^\infty(\R^d)$ with $0\le\chi\le1$,
$\chi\equiv1$ on $\{|x|\le1\}$, and
$\supp(\chi)\subseteq\{|x|\le2\}$. For $n\ge1$ set
$\chi_n(x):=\chi(x/n)$.
Then $\chi_n f\in C_{\rm c}^\infty(\R^d)$ and $\chi_n f\to f$
in $L^p(\R^d\!,\gamma)$ by dominated convergence.

We claim that also $L_p(\chi_n f)\to g$ in
$L^p(\R^d\!,\gamma)$. By \eqref{eq:O-gen} and the product rule
for differentiation, after some rewriting we obtain
$$
L_p(\chi_n f)-g
=(\chi_n-1)g+f\,L_p\chi_n+2\nabla\chi_n\cdot\nabla f.
$$
Since $\chi_n\to1$ pointwise and $|\chi_n-1|\le1$, we have
$|(\chi_n-1)g|^p\le|g|^p$, and the right-hand side belongs to
$L^1(\R^d\!,\gamma)$. By dominated convergence,
$(\chi_n-1)g\to0$ in $L^p(\R^d\!,\gamma)$.

Next, the supports of both $\nabla\chi_n$ and $\Delta\chi_n$ are
contained in the set $\{n\le|x|\le2n\}$ and
$$
\|\nabla\chi_n\|_\infty\lesssim n^{-1},\qquad
\|\Delta\chi_n\|_\infty\lesssim n^{-2},\qquad
\|x\cdot\nabla\chi_n\|_\infty\lesssim1,
$$
so $\|L_p\chi_n\|_\infty\lesssim1$ with
$\supp(L_p\chi_n)\subseteq\{n\le|x|\le2n\}$.
Since $f$ and $\nabla f$ are bounded,
$$
\begin{aligned}
&\|f\,L_p\chi_n\|_{L^p(\R^d\!,\gamma)}
+\|\nabla\chi_n\cdot\nabla f\|_{L^p(\R^d\!,\gamma)}
\lesssim
\gamma(\{n\le|x|\le2n\})^{1/p}\to0
\quad\hbox{as }n\to\infty,
\end{aligned}
$$
since $\gamma(\{|x|\ge n\})\to0$.
Thus $L_p(\chi_n f)\to g$ in $L^p(\R^d\!,\gamma)$ as claimed.
By closedness of $L_p$, we conclude that $f\in\Dom(L_p)$ and
$L_pf=g$. Moreover, $\chi_n f\to f$ in $\Dom(L_p)$.

Differentiating under the integral in the defining formula for
$OU(t)$, we obtain
$$
\partial^\alpha OU(t)f
=e^{-|\alpha|t}OU(t)(\partial^\alpha f),
\quad f\in Y.
$$
It follows that $OU(t)f\in Y$ whenever $f\in Y$.
Furthermore, the density of $C_{\rm c}^\infty(\R^d)$ in
$L^p(\R^d\!,\gamma)$ implies that $Y$ is dense in this space.
Hence by Proposition \ref{prop:core}, $Y$ is a core for $L_p$.
Since we have also shown that $C_{\rm c}^\infty(\R^d)$ is dense
in $Y$ with respect to the graph norm, we conclude that
$C_{\rm c}^\infty(\R^d)$ is a core for $L_p$.
\end{proof}

Returning to $p=2$, we wish to point out an interesting feature of
the formula \eqref{eq:O-gen}. For $f\in C_{\rm c}^\infty(\R^d)$,
multiplying both sides by the conjugate of a test function $\phi$,
integrating with respect to $\gamma$, and integrating by parts
after having written out the Gaussian density, we obtain
\begin{align*}
\int_{\R^d}L_2f(x)\ov{\phi(x)}\ud\gamma(x)
&=\frac1{(2\pi)^{d/2}}\int_{\R^d}
(\Delta f(x)-x\cdot\nabla f(x))\ov{\phi(x)}
\exp\Bigl(-\frac12|x|^2\Bigr)\ud x
\\ &=-\frac1{(2\pi)^{d/2}}\int_{\R^d}\!
\nabla f\cdot\ov{\nabla\phi}
\exp\Bigl(-\frac12|x|^2\Bigr)\ud x
=-\int_{\R^d}\!\nabla f\cdot\ov{\nabla\phi}\ud\gamma(x).
\end{align*}

The form $\aa_{OU}$ given by
$$
\aa_{OU}(f,g)=\int_{\R^d}\nabla f\cdot\ov{\nabla g}\ud\gamma,
\quad f,g\in C_{\rm c}^\infty(\R^d),
$$
is densely defined, bounded, positive, and closable.
Indeed, closability follows from Theorem \ref{thm:Friedrichs},
applied to the positive operator
$-L_2|_{C_{\rm c}^\infty(\R^d)}$.
Its closure $\ov\aa_{OU}$ is densely defined, closed, bounded,
and positive, and the associated operator $A_{OU}$ is positive
and selfadjoint. We will prove that $-L_2=A_{OU}$.

The above argument gives
$-\iprod{L_2f}{\phi}=\aa_{OU}(f,\phi)$ for all
$f,\phi\in C_{\rm c}^\infty(\R^d)$.
The density of $C_{\rm c}^\infty(\R^d)$ in
$\Dom(\ov\aa_{OU})$ with respect to the form norm and the
continuity of $\ov\aa_{OU}$ imply that this identity extends to
$\phi\in\Dom(\ov\aa_{OU})$.
Then the definition of $A_{OU}$ as the associated operator implies
$f\in\Dom(A_{OU})$ and $-L_2f=A_{OU}f$.
Since $C_{\rm c}^\infty(\R^d)$ is a core for $L_2$ and $A_{OU}$
is closed, we obtain the operator inclusion $-L_2\subseteq A_{OU}$.
We have $1\in\varrho(L_2)$ by contractivity and
$-1\in\varrho(A_{OU})$ by positivity and selfadjointness.
Proposition \ref{prop:semigroupsAB} therefore gives $-L_2=A_{OU}$
as operators in $L^2(\R^d\!,\gamma)$, with equality of domains.

We next show that the analogue of Theorem \ref{thm:LaplacianRd}
holds for $L_2$. Let $G_0$ be the gradient with domain
$C_{\rm c}^\infty(\R^d)$, viewed as an operator from
$L^2(\R^d\!,\gamma)$ to $L^2(\R^d\!,\gamma;\C^d)$.
For $f\in C_{\rm c}^\infty(\R^d)$ we have
$$
\n f\n_{\aa_{OU}}^2=\n f\n_2^2+\n G_0f\n_2^2.
$$
Since $\aa_{OU}$ is closable, this identity implies that $G_0$
is closable. Indeed, if $f_n\to0$ and $G_0f_n\to h$ in the
respective $L^2$-spaces, then $(f_n)$ is Cauchy in the form norm.
By Proposition \ref{prop:form-closable},
$\n G_0f_n\n_2^2=\aa_{OU}(f_n,f_n)\to0$, so $h=0$.
Writing $\nabla:=\ov G_0$, the same norm identity shows that
$\Dom(\ov\aa_{OU})=\Dom(\nabla)$ and
$$
\ov\aa_{OU}(f,g)=\iprod{\nabla f}{\nabla g},
\quad f,g\in\Dom(\nabla).
$$
We have already identified $-L_2$ as the operator associated with
the closed form $\ov\aa_{OU}$. This, in combination with
Proposition \ref{prop:AstarA-via-forms}, implies that
$$
-L_2=\nabla^\star\nabla,
$$
with equality of domains, the adjoint being taken with respect to the
Gaussian $L^2$-spaces.

Finally, since $C_{\rm c}^\infty(\R^d)$ is a core for the
selfadjoint operator $L_2$, for $f,g\in L^2(\R^d\!,\gamma)$
we have $f\in\Dom(L_2)$ and $L_2f=g$ if and only if
$$
\iprod{f}{L_2\phi}=\iprod{g}{\phi},
\quad \phi\in C_{\rm c}^\infty(\R^d).
$$
Indeed, the core property extends this identity to all
$\phi\in\Dom(L_2)$, and selfadjointness then gives
$f\in\Dom(L_2)$ and $L_2f=g$.
This is a special instance of Proposition \ref{prop:weakvsstrongsol}.
It identifies $L_2$ with the weak $L^2$-Ornstein--Uhlenbeck
operator, with the weak equation expressed by the preceding
test identity.

For $1<p<\infty$ it can be shown that
$$
\Dom(L_p)=W^{2,p}(\R^d\!,\gamma),
$$
the Gaussian Sobolev space of all $f\in L^p(\R^d\!,\gamma)$
admitting weak derivatives up to order $2$, all of which belong
to $L^p(\R^d\!,\gamma)$.
The proof of this fact is beyond the scope of this work.
Some properties of the first-order Gaussian Sobolev space
$W^{1,2}(\R^d\!,\gamma)$, including its identification with
$\Dom(\nabla)$, are explored in Problem \ref{prob:gaussianW12}.

For $d=1$ one has the following representation for the
Ornstein--Uhlenbeck semigroup in terms of the Hermite polynomials
$H_n$ discussed in Section \ref{subsec:Hermite}:
\begin{align}\label{eq:O-Hermite}
{OU}(t)H_n=e^{-nt}H_n,\quad n\in\N.
\end{align}
The recurrence relations for the Hermite polynomials give
$H_n''-xH_n'=-nH_n$.
Since polynomials belong to $L^2(\R,\gamma)$, integration by
parts against test functions and the preceding weak
characterisation show that $H_n\in\Dom(L_2)$ and
$L_2H_n=-nH_n$. This proves \eqref{eq:O-Hermite}.
Hence by Proposition \ref{prop:sigmaAdiag},
$$
\sigma(-L_2)=\N.
$$
The formula \eqref{eq:O-Hermite} generalises to arbitrary dimension
$d$ once one has found an analogue of the Hermite basis for
$L^2(\R^d\!,\gamma)$. This is accomplished in Section
\ref{subsec:WI} (see Theorems \ref{thm:WI} and \ref{thm:OU-diag}).

The analyticity of the Ornstein--Uhlenbeck semigroup already
follows from the positivity and selfadjointness of $-L_2$
by the Lumer--Phillips theorem. We summarise this conclusion in the following
theorem.

\begin{theorem}\label{thm:OU-analytic}
The operator $L_2$ generates an analytic $C_0$-semigroup of
contractions on $L^2(\R^d\!,\gamma)$ on every sector
$\Sigma_\omega$ with $0<\omega<\frac12\pi$.
\end{theorem}
\begin{proof}
The operator $-L_2$ is positive and selfadjoint.
Theorem \ref{thm:lumerphilips}, applied to $L_2$, therefore gives
the result for every $0<\omega<\frac12\pi$.
\end{proof}

The holomorphic extension to the right half-plane may also be
expressed in terms of the functional calculus as
$$
OU(z)=\exp(zL_2),\quad\Re z>0.
$$
This point of view will be taken up in Theorem \ref{thm:normal-sgr}.

\subsection{The Schr\"odinger Group}\label{subsec:Schr}\index{Schr\"odinger!group}

 Let $0\le V \in L_{\rm loc}^1(\R^d)$ be given; in applications we think of
$V$ as a {\em potential}.\index{potential}
Consider the forms
  $\aa_\Delta,\aa_V$ with domains $\Dom(\aa_\Delta) = \Dom(\aa_V) = C_{\rm c}^\infty(\R^d)$ defined by
 \begin{align*}
  \aa_\Delta (u,v)  := \frac12\int_{\R^d} \nabla u(x)\cdot \ov{\nabla v(x)}\ud x, \qquad
 \aa_V (u,v)  := \int_{\R^d} V(x) u(x)\ov{v(x)}\ud x.
 \end{align*}
 The form $\aa_\Delta+\aa_V$ with domain $\Dom(\aa_\Delta+\aa_V) = C_{\rm c}^\infty(\R^d)$ is densely defined, symmetric, positive, and
closable. To prove closability, suppose that $u_n\to0$ in $L^2(\R^d)$ and that $(u_n)_{n\ge 1}$ is
Cauchy for the form norm. Then $\nabla u_n\to g$ in $L^2(\R^d;\K^d)$ and $V^{1/2}u_n\to h$ in $L^2(\R^d)$ for some $g\in L^2(\R^d;\K^d)$ and $h\in L^2(\R^d)$. Since $u_n\to0$ in $L^2(\R^d)$,
closedness of weak differentiation gives $g=0$. Passing to a
subsequence, we may assume that $u_n\to 0$ almost everywhere.
After passing to a further subsequence in the second convergence, we may also assume that $ V^{1/2}u_n\to h$ almost everywhere, and therefore $h=0$. Hence
$(\aa_\Delta+\aa_V)(u_n,u_n)\to 0$.
By Proposition \ref{prop:form-closable}, this proves the asserted closability.

By Proposition \ref{prop:repr-closed-symm-form}, the operator associated with the closure $\aa$ of $\aa_\Delta+\aa_V$ is a positive selfadjoint operator. We denote this operator somewhat suggestively by
$A = -\frac12\Delta+V$.

An interesting special case arises if we take $ V(x) = \frac12|x|^2$\!. This choice results in the selfadjoint operator
$H =  -\frac12\Delta + \frac12|x|^2$ on $L^2(\R^d)$, the so-called {\em Hermite operator}.\index{Hermite!operator}
In quantum mechanics, this operator is called the {\em quantum harmonic oscillator}\index{quantum harmonic oscillator}. This special case will be addressed in more detail in Section \ref{subsec:QHA}.

For the remainder of this subsection, let
$$
A=-\frac12\Delta+V
$$
be the positive selfadjoint operator introduced above. By Stone's
theorem, the operator $-iA$ generates a unitary $C_0$-group on
$L^2(\R^d)$, the so-called {\em Schr\"odinger group} with potential $V$.
It solves the {\em Schr\"odinger equation}
$$
  i\frac{\partial}{\partial t}u(t,x)
  =
  -\frac12\Delta u(t,x)+ V(x)u(t,x),
  \quad t\in\R,\ x\in\R^d\!.
$$

The special case $V \equiv 0$ is of special interest:

\begin{example}[Free Schr\"odinger group]
The $C_0$-group $(S(t))_{t\in \R}$ on $L^2(\R^d)$ generated by the operator
$\frac{i}{2}\Delta$ with domain $H^2(\R^d)$ is called the {\em free Schr\"odinger group}.

On a formal level, $S(t)=H(it)$, where $(H(z))_{\Re z>0}$ is the
holomorphic extension of the scaled heat semigroup generated by
$\frac12\Delta$. Indeed, by the Fourier representation of $H$,
$$
\calF(H(s+it)f)(\xi)
=
e^{-\frac12(s+it)|\xi|^2}\wh f(\xi),
\quad s>0.
$$
The multipliers converge pointwise to $e^{-it|\xi|^2/2}$ and have modulus
at most one. Hence dominated convergence and Plancherel's theorem give
$$
\lim_{s\downarrow0}H(s+it)f
=
\calF^{-1}\bigl(e^{-it|\cdot|^2/2}\wh f\bigr)
\quad\hbox{in }L^2(\R^d).
$$
The operators on the right form a unitary $C_0$-group with generator
$\frac i2\Delta$ on $H^2(\R^d)$, and therefore coincide with $S$. An alternative derivation could be based on the result of Problem \ref{prob:bdry-sgr}.

We next derive the kernel representation of this group. For
$f\in L^1(\R^d)\cap L^2(\R^d)$ and $t\ne0$, the heat kernel formula
suggests
\begin{align}\label{eq:free-schr}
S(t)f
=
\frac1{(2\pi(it))^{d/2}}
\int_{\R^d}
\exp\Bigl(i\frac{|\cdot-y|^2}{2t}\Bigr)f(y)\ud y.
\end{align}
To make this rigorous, for $\Re z>0$ we set
$z^{d/2} := \exp(\frac d2\operatorname{Log} z)$,
where $\operatorname{Log}$ is the holomorphic logarithm on the right half-plane which is real-valued on $(0,\infty)$,
and
$$
(it)^{d/2}
:=
\lim_{s\downarrow0}(s+it)^{d/2}
=
|t|^{d/2}
\exp\Bigl(\frac{i\pi d}{4}\operatorname{sgn}(t)\Bigr).
$$
By analytic continuation, the Gaussian Fourier transform
formula of Lemma \ref{lem:Gauss} remains valid for complex
parameters with positive real part. Consequently, the heat
kernel formula holds for $H(z)$ whenever $\Re z>0$.

Now fix $t\ne0$. For $s>0$ we have the bound
$$
\Bigl|
\frac1{(2\pi(s+it))^{d/2}}
\exp\Bigl(-\frac{|x-y|^2}{2(s+it)}\Bigr)f(y)\Bigr|
\le
\frac1{(2\pi|t|)^{d/2}}|f(y)|,
$$
since
$$
|s+it| \ge |t|, \qquad
\Re\frac1{s+it} = \frac{s}{s^2+t^2} \ge 0.
$$
The function on the right is integrable. Dominated convergence
therefore gives, for every $x\in\R^d$,
$$
\lim_{s\downarrow0} H(s+it)f(x)
=
\frac1{(2\pi(it))^{d/2}}
\int_{\R^d}
\exp\Bigl(i\frac{|x-y|^2}{2t}\Bigr)f(y)\ud y.
$$
On the other hand, from the $L^2$-convergence we may choose a
sequence $s_n\downarrow0$ such that
$$\limn H(s_n+it)f = S(t)f$$
almost everywhere.
Comparing the two limits shows that  \eqref{eq:free-schr} holds.

The representation \eqref{eq:free-schr} implies that $S(t)f\in L^\infty(\R^d)$ for all $f\in L^1(\R^d)\cap L^2(\R^d)$ and $t\in\R\setminus\{0\}$, with bounds
\begin{equation*}
\begin{aligned}
 \n S(t)f\n_\infty &\le \frac1{(2\pi|t|)^{d/2}}\n f\n_1, \qquad \n S(t)f\n_2 \le \n f\n_2,
\end{aligned}
\end{equation*}
the former by a direct estimate and the latter by Plancherel's theorem.
By the Riesz--Thorin interpolation theorem, for all $t\in\R\setminus\{0\}$ the operators $S(t)$, when restricted to $L^1(\R^d)\cap L^2(\R^d)$, extend to bounded operators from
$L^p(\R^d)$ to $L^q(\R^d)$ for all $1\le p\le 2$ and $\frac1p+\frac1q =1$, with bound
\begin{equation}\label{eq:Schrod-Lp-est}
\| S(t)\|_{\calL(L^p(\R^d), L^q(\R^d))} \le \frac{1}{(2\pi|t|)^{d(\frac1p-\frac12)}}, \quad t\not=0.
\end{equation}
\end{example}

The bound \eqref{eq:Schrod-Lp-est} is an example of a
{\em dispersive estimate}.
For $f\in L^1(\R^d)\cap L^2(\R^d)$, it shows that
$\n S(t)f\n_\infty\to0$ as $|t|\to\infty$, while unitarity gives
$\n S(t)f\n_2=\n f\n_2$ for all $t\in\R$.
This decay in amplitude, accompanied by conservation of the
$L^2$-norm, expresses the spreading of the solution in space,
a phenomenon known as {\em dispersion}.\index{dispersion}

Dispersion also gives rise to estimates in which integrability
in space and time are considered together. The following
{\em Strichartz estimate} bounds an $L^p$-norm over
$\R\times\R^d$ in terms of the $L^2$-norm of the initial value.
It applies to arbitrary initial data in $L^2(\R^d)$ and
expresses a gain in spatial integrability after integration
in time.

\begin{theorem}[Strichartz estimate]\label{thm:Strichartz}
\index{Strichartz estimates}\index{theorem!Strichartz}
Let $S$ be the free Schr\"odinger group on $L^2(\R^d)$ and put
$p_d:=2(d+2)/d$. There exists a constant $C_d$ depending only on $d$
such that
\begin{align}\label{eq:Strichartz}
\Bigl(\int_\R\int_{\R^d}|S(t)f(x)|^{p_d}\ud x\ud t\Bigr)^{1/p_d}
\le C_d\n f\n_2,
\quad f\in L^2(\R^d).
\end{align}
\end{theorem}

In dimensions $1$, $2$, and $3$, the exponent $p_d$ equals $6$, $4$, and $10/3$,  respectively. The exponent is forced by scaling considerations as follows. For $\la>0$ put
$f_\la(x):=f(\la x)$. The identity
$S(t)f_\la(x)=(S(\la^2t)f)(\la x)$ shows that the
$L^p(\R^{1+d})$-norm of the solution scales by
$\la^{-(d+2)/p}$, whereas
$\n f_\la\n_2=\la^{-d/2}\n f\n_2$.
Thus an estimate of the form \eqref{eq:Strichartz} with exponent
$p$ can hold for all $f$ only if $(d+2)/p=d/2$, that is, $p=p_d$.

To prove the theorem, we first establish a
fractional integration estimate on the time axis. Its proof uses
the Hardy--Littlewood maximal theorem.

\begin{lemma}[Fractional integration]\label{lem:Strichartz-fractional}
Let $2<p<\infty$ and $\frac1p+\frac1{q}=1$. For $0\le g\in L^{q}(\R)$,
put $$
Jg(t):=\int_\R\frac{g(s)}{|t-s|^{2/p}}\ud s.
$$
Then $Jg$ is finite almost everywhere, belongs to $L^p(\R)$, and
$$
\n Jg\n_p\le C_p\n g\n_{q},
$$
where $C_p$ depends only on $p$.
\end{lemma}

\begin{proof} In the proof, $C_p$ is a constant whose value may change from line to line.

Put $\alpha:=2/p$, so that $0<\alpha<1$, and let $Mg$ be the
Hardy--Littlewood maximal function of $g$. For $R>0$ we split the
integral defining $Jg(t)$ into the regions $|t-s|<R$ and
$|t-s|\ge R$.

For the first region we use the identity
$$
r^{-\alpha}=R^{-\alpha}
+\alpha\int_r^R\rho^{-\alpha-1}\ud\rho,
\quad 0<r<R.
$$
Fubini's theorem and the definition of $Mg$ give
\begin{align*}
\int_{|t-s|<R}\frac{g(s)}{|t-s|^\alpha}\ud s
&=R^{-\alpha}\int_{|t-s|<R}g(s)\ud s + \alpha\int_0^R\rho^{-\alpha-1}
       \Bigl(\int_{|t-s|<\rho}g(s)\ud s\Bigr)\ud\rho\\
&\le 2Mg(t)\Bigl(R^{1-\alpha}
       +\alpha\int_0^R\rho^{-\alpha}\ud\rho\Bigr)
 =\frac{2}{1-\alpha}R^{1-\alpha}Mg(t).
\end{align*}
For the second region, H\"older's inequality and $\alpha p=2$
give
\begin{align*}
\int_{|t-s|\ge R}\frac{g(s)}{|t-s|^\alpha}\ud s
&\le\Bigl(\int_{|r|\ge R}|r|^{-2}\ud r\Bigr)^{1/p}\n g\n_{q}
=\Bigl(\frac2R\Bigr)^{1/p}\n g\n_{q}.
\end{align*}
Consequently,
$$
Jg(t)\le C_p\bigl(R^{1-2/p}Mg(t)+R^{-1/p}\n g\n_{q}\bigr).
$$

If $g=0$, there is nothing to prove. Otherwise $Mg(t)>0$ for every
$t$, and $Mg(t)<\infty$ for almost every $t$ by Theorem
\ref{thm:HL-maximal}. For such $t$, choose
$$
R:=\Bigl(\frac{\n g\n_{q}}{Mg(t)}\Bigr)^{q}.
$$
The two terms in the preceding bound then have the same value,
and we obtain
$$
Jg(t)\le C_p\n g\n_{q}^{1-q/p}\bigl(Mg(t)\bigr)^{q/p}.
$$
Taking the $L^p$-norm and applying Theorem \ref{thm:HL-maximal}
once more gives
$$
\n Jg\n_p
\le C_p\n g\n_{q}^{1-q/p}\n Mg\n_{q}^{q/p}
\le C_p\n g\n_{q},
$$
which proves the lemma.
\end{proof}

\begin{proof}[Proof of Theorem \ref{thm:Strichartz}]
Let $p:= p_d =2(d+2)/d$ and $\frac1p+\frac1{q}=1$. Then $1<q<2<p<\infty$ and
$d(\frac1{q}-\frac12) =\frac{d}{d+2}=\frac2p$.
Thus, for $v\in L^{q}(\R^d)\cap L^2(\R^d)$,
\eqref{eq:Schrod-Lp-est} gives
$$
\n S(t-s)v\n_p
\le (2\pi|t-s|)^{-2/p}\n v\n_{q},
\quad t\ne s.
$$

{\em Step 1} --
Let $F\in C_{\rm c}^\infty(\R^{1+d})$, write
$F(t):=F(t,\cdot)$, and put
$$
BF:=\int_\R S(-s)F(s)\ud s.
$$
This is a well-defined $L^2(\R^d)$-valued Riemann integral, since
the integrand is continuous and vanishes outside a bounded time
interval. Expanding its squared norm and using that $S(-t)^\star = S(t)$ by unitarity, we find
\begin{align*}
\n BF\n_2^2
=\int_\R\int_\R
  \iprod{S(-s)F(s)}{S(-t)F(t)}\ud s\ud t
  =\int_\R\int_\R
  \iprod{S(t-s)F(s)}{F(t)}\ud s\ud t.
\end{align*}
These scalar integrals are absolutely convergent, because the
integrands vanish outside a bounded set and their absolute values
are at most $\n F(s)\n_2\n F(t)\n_2$.

Put $g(t):=\n F(t)\n_{q}$. For $t\ne s$, H\"older's inequality
and the dispersive estimate give
\begin{align*}
\bigl|\iprod{S(t-s)F(s)}{F(t)}\bigr|
\le\n S(t-s)F(s)\n_p\n F(t)\n_{q}
\le (2\pi|t-s|)^{-2/p}g(s)g(t).
\end{align*}
We may ignore the diagonal $s=t$, which has measure zero.
Applying Lemma \ref{lem:Strichartz-fractional} and H\"older's
inequality on $\R$, we obtain
\begin{align*}
\n BF\n_2^2
&\le (2\pi)^{-2/p}
  \int_\R g(t)\Bigl(\int_\R
       \frac{g(s)}{|t-s|^{2/p}}\ud s\Bigr)\ud t\\
&\le (2\pi)^{-2/p}\n g\n_{q}\n Jg\n_p
 \le (2\pi)^{-2/p}C_p\n g\n_{q}^2.
\end{align*}
Since
$$
\n g\n_{q}^{q}
=\int_\R\int_{\R^d}|F(t,x)|^{q}\ud x\ud t,
$$
this proves, with $C_d:=((2\pi)^{-2/p}C_p)^{1/2}$, that
\begin{align}\label{eq:Strichartz-dual}
\Bigl\n\int_\R S(-s)F(s)\ud s\Bigr\n_2
\le C_d\n F\n_{L^{q}(\R^{1+d})}.
\end{align}

{\em Step 2} --
We now obtain the desired estimate by duality. First let
$f\in L^2(\R^d)$ with $\wh f\in L^1(\R^d)$. The Fourier
representation
$$
u(t,x):=\frac1{(2\pi)^{d/2}}
\int_{\R^d}e^{ix\cdot\xi-it|\xi|^2/2}\wh f(\xi)\ud\xi
$$
gives a bounded continuous function $u$ on $\R^{1+d}$ satisfying
$u(t,\cdot)=S(t)f$ in $L^2(\R^d)$. For
$F\in C_{\rm c}^\infty(\R^{1+d})$, unitarity and
\eqref{eq:Strichartz-dual} give
\begin{align*}
\Bigl|\int_\R\int_{\R^d}
u(t,x)\ov{F(t,x)}\ud x\ud t\Bigr|
=\Bigl|\int_\R\iprod{S(t)f}{F(t)}\ud t\Bigr|=|\iprod{f}{BF}|
 \le C_d\n f\n_2\n F\n_{L^{q}(\R^{1+d})}.
\end{align*}

For $R>0$, put $Q_R:=(-R,R)^{1+d}$. Since $u$ is bounded, its
restriction to $Q_R$ belongs to $L^p(Q_R)$. Taking the supremum
over $F\in C_{\rm c}^\infty(Q_R)$ with
$\n F\n_{L^{q}(Q_R)}\le1$, and using $L^p$-duality and the
density of these test functions in $L^{q}(Q_R)$, we conclude that
$$
\n u\n_{L^p(Q_R)}\le C_d\n f\n_2.
$$
Letting $R\to\infty$ and applying the monotone convergence theorem,
we obtain \eqref{eq:Strichartz}.

For general $f\in L^2(\R^d)$, put
$f_n:=\calF^{-1}(\one_{B(0;n)}\wh f)$.
By the Cauchy--Schwarz inequality, $\wh f_n\in L^1(\R^d)$, and
Plancherel's theorem gives $f_n\to f$ in $L^2(\R^d)$. The estimate
just proved, applied to $f_n-f_m$, shows that the functions
$u_n(t,x):=S(t)f_n(x)$ converge in $L^p(\R^{1+d})$ to a function
$u$. By Fubini's theorem,
$$
\int_\R\n u_n(t,\cdot)-u(t,\cdot)\n_p^p\ud t
=\n u_n-u\n_{L^p(\R^{1+d})}^p\to 0.
$$
After passing to a subsequence, we therefore have
$u_n(t,\cdot)\to u(t,\cdot)$ in $L^p(\R^d)$ for almost every $t$.
On the other hand, unitarity gives
$$
\n u_n(t,\cdot)-S(t)f\n_2
=\n f_n-f\n_2\to 0,
\quad t\in\R.
$$
The limits in $L^p(\R^d)$ and $L^2(\R^d)$ must agree almost
everywhere, so $u(t,\cdot)=S(t)f$ for almost every $t$. Passing
to the limit in the estimate for $f_n$ proves
\eqref{eq:Strichartz}.
\end{proof}

\subsection{The Wave Group}\label{subsec:wave}

The wave equation is second order in time, so its evolution is
determined by an initial displacement and an initial velocity.
By keeping track of both displacement and velocity, we can
rewrite the equation as a first-order system and apply the
theory of $C_0$-groups. We first consider bounded domains with
Dirichlet boundary conditions, where conservation of energy
corresponds to unitarity of the wave group, and then turn to
$\R^d$, where the Fourier transform gives an explicit formula
for the evolution.

\paragraph{The Wave Group on Bounded Domains}

Let $D$ be a nonempty bounded open set in $\R^d\!$.
The space $H:=H_0^1(D)\times L^2(D)$ is a Hilbert space with the norm
\begin{align*}
\n(u,v)\n^2 = \nn u\nn_{H_0^1(D)}^2 + \n v\n_2^2,
\end{align*}
where we consider the norm on $H_0^1(D)$ given by
\begin{align}\label{eq:wave-normH01} \nn u\nn_{H_0^1(D)}^2 := \int_{D} |\nabla u|^2 \ud x, \quad u\in H_0^1(D).
\end{align}
This norm is equivalent to the usual Sobolev norm on $H_0^1(D)$ by Poincar\'e's inequality.
In $H$ we define the operator $A$ by
\begin{align*}
 A & := \begin{pmatrix} 0 & I \\ \Delta_{\rm Dir} & 0
        \end{pmatrix},
\qquad \Dom(A)  := \Dom(\Delta_{\rm Dir}) \times H_0^1(D),
\end{align*}
where $\Delta_{\rm Dir}$ is the Dirichlet Laplacian on $L^2(D)$.
 We will prove that $A$ is the generator of a unitary $C_0$-group\index{wave!group} $W$ on $H$.
 This group solves the linear wave equation\index{wave!equation}
 $$
 \left\{
\begin{aligned}
\frac{\partial^2 u}{\partial t^2}(t,x) & = \Delta u(t,x), && t\in\R, \ x\in D,
\\ u(0,x) & = u_0(x), && x\in D,
\\ \frac{\partial u}{\partial t}(0,x) & = v_0(x),  && x\in D,
\end{aligned}
\right.
$$ written as a system of abstract first-order ODEs
 $u'=v$, $v'=\Delta u$, with initial value $u(0)=u_0$, $v(0)=v_0$,
 subject to Dirichlet boundary conditions.

The operator $A$ is densely defined and an integration by parts gives,
for $u=\begin{pmatrix} u_1 \\ u_2 \end{pmatrix}\in \Dom(A)$ and $v= \begin{pmatrix} v_1 \\ v_2 \end{pmatrix} \in \Dom(A)$,
\begin{align*}
\iprod{Au}{v}  = \left( \begin{pmatrix} u_2 \\ \Delta u_1 \end{pmatrix} \Bigg| \begin{pmatrix} v_1 \\ v_2 \end{pmatrix}\right)
& = \int_D \nabla u_2 \cdot \ov{\nabla v_1} - \nabla u_1\cdot\ov{\nabla v_2}\ud x
\\ & = - \left(\begin{pmatrix} u_1 \\ u_2 \end{pmatrix} \Bigg| \begin{pmatrix} v_2 \\ \Delta v_1 \end{pmatrix}\right)
= -\iprod{u}{Av}.
\end{align*}
This implies that $iA$ is symmetric.

Since $0\in \varrho(\Delta_{\rm Dir})$ by Theorem \ref{thm:spectra-Laplacian}, we may consider the operator
$$ R:= \begin{pmatrix} 0 & \Delta_{\rm Dir}^{-1} \\ I & 0 \end{pmatrix}
$$ on $H$. For $(u,v)\in H$ we have
$$R \begin{pmatrix} u \\ v \end{pmatrix} = \begin{pmatrix} \Delta_{\rm Dir}^{-1} v \\ u \end{pmatrix}
\in \Dom(\Delta_{\rm Dir})\times H_0^1(D) = \Dom(A).$$
This operator is bounded; the boundedness of $\Delta_{\rm Dir}^{-1}: L^2(D)\to H_0^1(D)$ follows from the closed graph theorem or by noting that the Cauchy--Schwarz and Poincar\'e inequalities imply the bound
$\n\nabla (\Delta_{\rm Dir}^{-1})v\n_2
\le C_D\n v\n_2$, with $C_D$ the Poincar\'e constant.

It is immediate to check that $AR = I$ and $RAh=h$ for $h\in \Dom(A)$. This proves that $A$ is boundedly invertible, that is, $0\in\varrho(A)$. An application of
Proposition \ref{prop:symm-sgr-sa} now gives that $iA$ is selfadjoint.
Therefore, by Stone's theorem (Theorem \ref{thm:Stone}) we obtain:

\begin{theorem}[Wave group on bounded domains]\label{thm:wave-groupD}
The operator $A$ generates a unitary $C_0$-group $(W(t))_{t\in \R}$ on $H_0^1(D)\times L^2(D)$, provided $H_0^1(D)$
is endowed with the equivalent norm given by \eqref{eq:wave-normH01}.
\end{theorem}

\paragraph{The Wave Group on \texorpdfstring{$\R^d$}{R\textasciicircum d}}\index{wave!group}

Let us next consider the case $D = \R^d\!$, which
is not covered by the above considerations since it was assumed that
$D$ be bounded.
We have $H_0^1(\R^d) = H^1(\R^d)$ by Theorem \ref{thm:W01pRd}
and Theorem \ref{thm:W=H} and $\Delta_{\rm Dir} = \Delta$ with $\Dom(\Delta) = H^2(\R^d)$.
This suggests considering in $H^1(\R^d)\times L^2(\R^d)$ the operator
\begin{align*} A & := \begin{pmatrix}
                            0 & I \\ \Delta & 0
                      \end{pmatrix}, \qquad
\Dom(A)  := H^2(\R^d)\times H^1(\R^d).
\end{align*}
We will use the theory of Fourier multipliers to prove that $A$ generates a $C_0$-group on  $H^1(\R^d)\times L^2(\R^d)$  and give an explicit expression for it.

To motivate the expressions below we first consider the matrix $$ A_a := \left(\begin{matrix} 0 & 1 \\ -a & 0\end{matrix}\right),$$ where $a\ge 0$ is a nonnegative scalar, and compute its exponentials $e^{tA_a}$. The powers of $A_a$ are given by
$$A_a^{2k} =
\left(\begin{matrix} (-a)^{k} & 0 \\ 0 & (-a)^{k} \end{matrix}\right), \qquad A_a^{2k+1} =
\left(\begin{matrix} 0 & (-a)^k \\ (-a)^{k+1} & 0\end{matrix}\right),
$$
so that with $b:= a^{1/2}$\!,
\begin{align*}
e^{tA_a}
& = \sum_{k=0}^\infty \left(\begin{matrix} \displaystyle(-1)^k\frac{t^{2k}}{(2k)!}b^{2k} & \displaystyle(-1)^k\frac{t^{2k+1}}{(2k+1)!}b^{2k}\\ \\ \displaystyle -(-1)^k\frac{t^{2k+1}}{(2k+1)!}b^{2k+2} & \displaystyle (-1)^k\frac{t^{2k}}{(2k)!}b^{2k} \end{matrix}\right)
= \left(\begin{matrix} \cos tb & b^{-1}\sin tb\\ -b\sin tb &  \cos tb \end{matrix}\right)\!,
\end{align*}
where $b^{-1}\sin(tb)$ is understood to be $t$ when $b=0$.
Substituting $(-\Delta)^{1/2}$ for $b$, we arrive at the following guess for the expression for the wave group:
$$W(t) = \left(\begin{matrix} \cos (t(-\Delta)^{1/2}) & (-\Delta)^{-1/2}\sin (t(-\Delta)^{1/2})\\ \\
-(-\Delta)^{1/2} \sin (t(-\Delta)^{1/2}) & \cos (t(-\Delta)^{1/2})\end{matrix}\right), \quad t\in\R.$$
We need to give a meaning to the operators occurring in this matrix, which can be accomplished by the functional calculus of $-\Delta$, or by interpreting them as
Fourier multiplier operators as follows.
The operators on the diagonal can be interpreted as
Fourier multiplier operators on $L^2(\R^d)$ associated with the multipliers
$$m_{1,1;t}(\xi)=  m_{2,2;t}(\xi)=  \cos( t|\xi|);
$$
this function belongs to $L^\infty(\R^d)$ with norm $1$ for every $t\in\R$. Recalling
the characterisation of $H^1(\R^d)$ as those functions $f$ in $L^2(\R^d)$ for which
$\xi\mapsto (1+|\xi|^2)^{1/2}\wh f(\xi)$ is in $L^2(\R^d)$, we see moreover that $\cos (t(-\Delta)^{1/2})$ maps
$W^{1,2}(\R^d)$ into itself. Recalling the norm of $H^1(\R^d)$ given by \eqref{eq:norm-Hs}, this argument also gives the estimates
$$ \Big\n \cos (t(-\Delta)^{1/2})\Big\n_{\calL(L^2(\R^d))}\le 1, \qquad \Big \n \cos (t(-\Delta)^{1/2})\Big\n_{\calL(H^1(\R^d))}\le 1.$$
Similarly, we can associate a bounded operator
from $H^1(\R^d)$ to $L^2(\R^d)$
with the multiplier $$m_{2,1;t}(\xi)=  -|\xi|\sin( t|\xi|).$$ Indeed, if $f\in H^1(\R^d)$, then
 \begin{align*}
  \big\n  m_{2,1;t}(\xi) \wh f(\xi)\big\n_{L^2(\R^d)}
 & \le \Big(\sup_{\xi\in \R^d} \frac{|\xi|}{(1+| \xi|^2)^{1/2}}\bigl|\sin( t|\xi|)\bigr| \Big)\n f\n_{H^1(\R^d)}
  \le \n f\n_{H^1(\R^d)}
 \end{align*}
 and therefore
$$ \Big\n -(-\Delta)^{1/2} \sin (t(-\Delta)^{1/2}) \Big\n_{\calL(H^1(\R^d),L^2(\R^d))} \le 1, \quad t\in \R.$$
Similarly, the operators $(-\Delta)^{-1/2} \sin (t(-\Delta)^{1/2})$ are interpreted as bounded
operators from $L^2(\R^d)$ into $H^1(\R^d)$ given by the multipliers $m_{1,2;t}(\xi)=\sin(t|\xi|)/|\xi|$, where the value at $\xi=0$ is understood to be $t$. These multipliers satisfy (distinguish the cases $|\xi|\le 1$ and $|\xi|>1$)
\begin{align*}
  \big\n \calF^{-1}(m_{1,2;t}\wh f)\big\n_{H^1(\R^d)}
 &  \le \Big(\sup_{\xi\in \R^d} \frac{|\sin( t|\xi|)|}{|\xi|} (1+|\xi|^2)^{1/2}\Big)\n \wh f\n_{L^2(\R^d)}
\\ & \le C(1+|t|)\n \wh f\n_{L^2(\R^d)}= C(1+|t|)\n f\n_{L^2(\R^d)},
 \end{align*}
 where $C$ is a universal constant.
 Therefore
$$  \Bigl\n (-\Delta)^{-1/2} \sin(t(-\Delta)^{1/2}) \Bigr\n_{\calL(L^2(\R^d),H^1(\R^d))} \le C(1+|t|), \quad t\in \R.$$

\begin{theorem}[Wave group on $\R^d$]\label{thm:wave-group} The operator $A$ generates a $C_0$-group $(W(t))_{t\in \R}$ on $H^1(\R^d)\times L^2(\R^d)$, given by
$$W(t) = \left(\begin{matrix} \cos (t(-\Delta)^{1/2}) & (-\Delta)^{-1/2}\sin (t(-\Delta)^{1/2})\\
-(-\Delta)^{1/2} \sin (t(-\Delta)^{1/2}) & \cos (t(-\Delta)^{1/2})\end{matrix}\right), \quad t\in \R.$$
Moreover, there is a constant $C\ge 0$ such that $$ \n W(t)\n \le C(1+|t|), \quad t\in\R.$$
\end{theorem}

\begin{proof} We proceed in two steps.

\smallskip
{\em Step 1} --
The group property follows from formal matrix multiplication, which can be made rigorous by noting that in the Fourier domain we are just multiplying matrices of scalar-valued multipliers, much like what we did in the treatment of the heat and Poisson semigroups.

To prove strong continuity we note that $\n m_{1,1;t}\n_\infty\le 1$
and $\lim_{t\to 0} m_{1,1;t}(\xi)= 1$ pointwise,
implying that $\lim_{t\to 0}  m_{1,1;t}\wh f  =  \wh f$ in $L^2(\R^d)$ by dominated convergence. Hence
 $\lim_{t\to 0}\cos (t(-\Delta)^{1/2})f = f$ in $L^2(\R^d)$ for all $f\in L^2(\R^d)$ by Plancherel's theorem. The strong convergence of the other three terms is proved similarly.

 The $1+|t|$ bound follows from the computations preceding the theorem.

 \smallskip
{\em Step 2} -- We now identify the generator domain.
Let $G$ be the generator of the $C_0$-group $(W(t))_{t\in\R}$ on
$X:=H^1(\R^d)\times L^2(\R^d)$,
endowed with its usual Hilbertian norm.

Set $B:=(-\Delta)^{1/2}$ on $L^2(\R^d)$, so that
$\Dom(B)=H^1(\R^d)$, $\Dom(B^2)=H^2(\R^d)$, and $B^2 u=-\Delta u$ for $u\in \Dom(B^2)$.
Let
$$
A_0:=\begin{pmatrix}0&I\\ \Delta&0\end{pmatrix}
$$
with domain $\Dom(A_0):=H^2(\R^d)\times H^1(\R^d)$.

We first show that $A_0\subseteq G$.
Let $(u,v)\in H^2(\R^d)\times H^1(\R^d)$.
The inequalities $|\cos z-1|\le |z|$, $|\sin z|\le |z|$,
and $|\sin z/z|\le1$ for $z\in\R$ (with the quotient interpreted
as $1$ at $z=0$) justify dominated convergence in the Fourier
representation and we obtain
$$
\lim_{t\to 0}\frac{\cos(tB)u-u}{t} = 0 \quad \hbox{in }H^1(\R^d),
\qquad
\lim_{t\to 0}\frac{B^{-1}\sin(tB)v}{t} = v \quad \hbox{in }H^1(\R^d),
$$
and
$$
\lim_{t\to 0}\frac{-B\sin(tB)u}{t} = -B^2u=\Delta u \quad \hbox{in }L^2(\R^d),
\qquad
\lim_{t\to 0}\frac{\cos(tB)v-v}{t} = 0 \quad \hbox{in }L^2(\R^d).
$$
Hence
$$
\lim_{t\to 0}
\frac{W(t)(u,v)-(u,v)}{t}= (v,\Delta u)
\quad \hbox{in } X,
$$
so $(u,v)\in\Dom(G)$ and $G(u,v)=A_0(u,v)$.
The first and third limits used
$u\in H^2(\R^d)$, and the second and fourth used
$v\in H^1(\R^d)$.
Thus $A_0\subseteq G$.

To prove the converse inclusion, fix $\la>1$. Since $(W(t))_{t\in\R}$ is a
$C_0$-group with bound $ \n W(t)\n \le C(1+|t|)$ for $t\in \R$, we have $\la\in\varrho(G)$. We claim that $\la-A_0$ is surjective on $X$.
Indeed, given $(f,g)\in H^1(\R^d)\times L^2(\R^d)$, define
$$
\wh u(\xi):=\frac{\la \wh f(\xi)+\wh g(\xi)}{\la^2+|\xi|^2},
\qquad
v:=\la u-f.
$$
Then $u\in H^2(\R^d)$ and $v\in H^1(\R^d)$, and a direct computation shows that
$(\la-A_0)(u,v)=(f,g)$.
Thus $\la-A_0$ is surjective.

Now let $(u,v)\in\Dom(G)$ and put $(f,g):=(\la-G)(u,v)\in X$.
By the surjectivity of $\la-A_0$, there exists $(u_0,v_0)\in\Dom(A_0)$ such that
$(\la-A_0)(u_0,v_0)=(f,g)$.
Since $A_0\subseteq G$, we also have
$$
(\la-G)(u_0,v_0)=(f,g)=(\la-G)(u,v).
$$
Since $\la\in\varrho(G)$, the operator $\la-G$ is injective and we have
$(u,v)=(u_0,v_0)\in\Dom(A_0)$.
Therefore $G=A_0$, that is,
$\Dom(G)=H^2(\R^d)\times H^1(\R^d)$ and
$G(u,v)=(v,\Delta u)$ for $(u,v)\in \Dom(G)$.
This proves that the generator of $(W(t))_{t\in\R}$ equals the operator
$A = \begin{pmatrix}0&I\\ \Delta&0\end{pmatrix}$ with the claimed domain.
\end{proof}

\begin{remark}
Notwithstanding the linear growth bound for $W(t)$, the {\em energy functional}
$$
E(t):=\frac12(\|v(t)\|_2^2+\|\nabla u(t)\|_2^2),\qquad (u(t),v(t)):=W(t)(u_0,v_0),
$$
is constant in time. For $(u_0,v_0)\in\Dom(A)$ one has $E'(t)=0$ by direct
calculation.
For a differentiable Hilbert space-valued function $w$, the function $t\mapsto\frac12\|w(t)\|^2$ has derivative $\Re\iprod{\dot w(t)}{w(t)}$. Thus, since $u'(t)=v(t)$ and $v'(t)=\Delta u(t)$,
\begin{align*}
E'(t)
& = \Re\iprod{\dot v(t)}{v(t)} + \Re\iprod{\nabla \dot u(t)}{\nabla u(t)}
 = \Re\iprod{\Delta u(t)}{v(t)} + \Re\iprod{\nabla v(t)}{\nabla u(t)}.
\end{align*}
Integrating the first term on the right-hand side by parts and using $\Re\iprod{\nabla u(t)}{\nabla v(t)}=\Re\iprod{\nabla v(t)}{\nabla u(t)}$, we obtain
$$
E'(t) = -\Re\iprod{\nabla u(t)}{\nabla v(t)} + \Re\iprod{\nabla v(t)}{\nabla u(t)}=0.
$$

For general $(u_0,v_0)\in H^1(\R^d)\times L^2(\R^d)$ one obtains $E(t)=E(0)$ by approximating
$(u_0,v_0)$ in $H^1\times L^2$ with elements of $\Dom(A)$ and using strong continuity of $W(t)$.
\end{remark}

We conclude with some informal remarks establishing a connection with the Poisson semigroup of Section \ref{subsec:Poisson-sgr}.
Since $(-\Delta)^{1/2}$ with domain $H^1(\R^d)$ is selfadjoint, $-i(-\Delta)^{1/2}$ is the generator of a $C_0$-group $(U(t))_{t\in \R}$ of unitary operators on $L^2(\R^d)$ by Stone's theorem.
For $f\in H^2(\R^d)$ we have $(-\Delta)^{1/2}f\in H^1(\R^d)$  and
$$ \frac{{\rm d}^2}{{\rm d}t^2}U(t) f  = [-i(-\Delta)^{1/2}]^2 U(t)f = \Delta U(t)f.$$
In this sense, $t\mapsto u(t,x) = U(t)f(x)$ satisfies the wave equation with initial condition $u(0,x) = f(x)$. Its initial velocity is $u'(0)=-i(-\Delta)^{1/2}f$. In fact, we could run the same argument for $i(-\Delta)^{1/2}$\!, which generates $(U(-t))_{t\in\R}$, and conclude that its orbits also solve the wave equation with initial condition $u(0,x) = f(x)$. Interpreting $U(t)$ and $U(-t)$ as Fourier multipliers one obtains
$$ \frac12(U(t) + U(-t)) = \cos\bigl(t (-\Delta)^{1/2}\bigr),$$
which is the first entry in the matrix representation for the wave group.
For $f\in H^2(\R^d)$ and $g\in H^1(\R^d)$, set
$$
u(t):=\cos(tB)f+B^{-1}\sin(tB)g,\qquad B:=(-\Delta)^{1/2}.
$$
Then
$$
u'(t)=-B\sin(tB)f+\cos(tB)g,\qquad
u''(t)=-B^2u(t)=\Delta u(t),
$$
with derivatives in $H^1(\R^d)$ and $L^2(\R^d)$, respectively.
Thus $u(0)=f$ and $u'(0)=g$; the averaged group above gives
the case $g=0$.

\section{Semigroups Generated by Normal Operators}\label{sec:spectaltheoremrevisited}

We begin with a general observation about semigroup generation by normal operators. Recall from
Section \ref{sec:analytic-semigroups} the notation
$$\Sigma_\om := \{z\in \C\setminus\{0\}: \, |\arg(z)|<\om\}$$
for the open sector of angle $\om\in (0,\pi)$, where arguments are taken in $(-\pi,\pi)$.

\begin{theorem}[Semigroups generated by normal operators]\label{thm:normal-sgr} Let $N$ be a normal operator in a Hilbert space $H$ with associated projec\-tion-valued measure $P$. Then:
 \begin{enumerate}[label={\rm(\arabic*)}, leftmargin=*]
  \item\label{it:normal-sgr1}  if $\sigma(N)$ is contained in the closed right half-plane, then
  $-N$ is the generator of a $C_0$-semigroup of contractions $S$, given
  by
  $$ S(t) = \int_{\sigma(N)} e^{-\la t} \ud P(\la),\quad t\ge 0;$$
  \item\label{it:normal-sgr2}  if $\sigma(N)$ is contained in a closed  sector of angle $\theta\in (0,\frac12\pi)$, then
  $-N$ is the generator of an analytic $C_0$-semigroup of contractions $(S(z))_{z\in \Sigma_{\frac12\pi -\theta}}$, given
  by
  $$ S(z) = \int_{\sigma(N)} e^{-\la z} \ud P(\la), \quad z\in \Sigma_{\frac12\pi-\theta}.$$
 \end{enumerate}
\end{theorem}

\begin{proof}
We give a detailed proof of \ref{it:normal-sgr1}; the proof of \ref{it:normal-sgr2} is entirely similar.

The operators $S(t)$ are well defined and contractive by Theorem \ref{thm:Borel-FC}(iv), since $|e^{-t\la}|\leq1$ for $t\geq0$ and $\la\in\si(N)$. The identity $S(0)=I$ follows from $P_{\si(N)}=I$, and the semigroup property follows from the multiplicativity of the bounded Borel calculus.

We next prove strong continuity at $0$. For $x\in H$, dominated convergence gives
$$
\n S(t)x-x\n^2
 = \int_{\si(N)} |e^{-t\la}-1|^2\ud P_x(\la)\to0
$$
as $t\downarrow0$, since $|e^{-t\la}-1|\leq2$ on the closed right half-plane and $P_x$ is a finite measure. It follows from the semigroup property that $(S(t))_{t\geq0}$ is a $C_0$-semigroup of contractions.

It remains to identify its generator. Let $A$ denote the generator of $(S(t))_{t\geq0}$. We first prove that $-N\subseteq A$. Let $x\in\Dom(N)$. For $\la\in\si(N)$ and $t>0$ we have
$\frac1t (e^{-t\la}-1) = -\la\int_0^1 e^{-st\la}\ud s$,
and therefore
$$
\Bigl|\frac{e^{-t\la}-1}{t}+\la\Bigr|
 \leq |\la|\int_0^1 |1-e^{-st\la}|\ud s
 \leq 2|\la|.
$$
By dominated convergence, using $x\in\Dom(N)$, we obtain
$$
\Big\n\frac{S(t)x-x}{t}+Nx\Big\n^2
 =\int_{\si(N)}
   \Bigl|\frac{e^{-t\la}-1}{t}+\la\Bigr|^2\ud P_x(\la)
 \to0
$$
as $t\downarrow0$. Hence $x\in\Dom(A)$ and $Ax=-Nx$.

Thus $-N\subseteq A$. To prove equality, fix $\mu>0$.
Since $\sigma(N)$ is contained in the closed right
half-plane we have $\mu\in\varrho(-N)$. Also, $\mu\in\varrho(A)$ by
Proposition \ref{prop:resolv}, since $A$ generates
a contraction semigroup.
Proposition \ref{prop:semigroupsAB} therefore gives $A=-N$.

For \ref{it:normal-sgr2} the same argument applies for $z\in\Sigma_{\frac12\pi-\theta}$. Indeed, if $\la\in\si(N)$, then $\Re(z\la)\geq0$, so $|e^{-z\la}|\leq1$. We leave some details to the reader. The semigroup law, contractivity, and strong continuity follow as above, and holomorphy follows by applying dominated convergence to the difference quotients of $z\mapsto e^{-z\la}$ on compact sub-sectors of $\Sigma_{\frac12\pi-\theta}$. The generator on the positive real axis is again $-N$, by the result just proved.
\end{proof}

Some of the semigroup examples of the previous section can be constructed rather easily using the spectral theorem.

\begin{example}[Heat semigroup, Poisson semigroup, free Schr\"odinger group, wave group revisited]\label{ex:HSgr-ST}
 Let $Q$ be the projection-valued measure on $\R$ associated with the negative Laplace operator
 $-\Delta$, viewed as a selfadjoint operator in $L^2(\R^d)$ (see  Problem \ref{prob:PVM-Laplace}).
 The heat semigroup $H$ is then given by
 $$ H(t) = \int_{\R} e^{-\la t}\ud Q(\la),\quad t\ge 0,$$
 and the free Schr\"odinger group by
$$
        S(t)=\int_{\R} e^{-i\la t/2}\ud Q(\la),\quad t\in\R.
$$

The positive square root $(-\Delta)^{1/2}$ defined through Proposition \ref{prop:pos-sqrt} coincides with the unbounded Fourier multiplier operator
corresponding to the multiplier $m(\xi) = |\xi|$. Using Proposition \ref{prop:normal-subst} to switch between the  projection-valued measure $Q$ of $-\Delta$ and $R$ of $(-\Delta)^{1/2}$\!,
we see that the Poisson semigroup generated by $-(-\Delta)^{1/2}$ is given by
$$ P(t) = \int_{[0,\infty)} e^{-\la t}\ud R(\la)= \int_{[0,\infty)} e^{-\la^{1/2}t }\ud Q(\la),\quad t\ge 0.$$

Similarly, the operators $\cos(t(-\Delta)^{1/2})$
and $\sin(t(-\Delta)^{1/2})$ featuring in the wave group are given by
\begin{align*}\cos\bigl(t(-\Delta)^{1/2}\bigr)
& = \int_{[0,\infty)} \cos(t\la)\ud R(\la)
  = \int_{[0,\infty)} \cos\bigl(t\la^{1/2}\bigr)\ud Q(\la), \\
  \sin\bigl(t(-\Delta)^{1/2}\bigr)
  & =\int_{[0,\infty)} \sin(t \la)\ud R(\la)
  = \int_{[0,\infty)} \sin\bigl(t\la^{1/2}\bigr)\ud Q(\la).
 \end{align*}
\end{example}

\begin{example}[Stone's theorem  revisited]\index{theorem!Stone}
 Let $P$ be the projection-valued measure on $\R$ associated with a selfadjoint operator $A$
 acting in a Hilbert space $H$.
 Then the unitary $C_0$-group $(U(t))_{t\in \R}$ generated by $-iA$ is given by
 $$ U(t) = \int_{\R} e^{-i\la t}\ud P(\la),\quad t\in\R.$$
\end{example}

We conclude with an application of the spectral theorem to the
behaviour of selfadjoint semigroups as $t\to\infty$.

\begin{proposition}[Convergence to equilibrium, spectral gap]
\label{prop:spectral-gap-sgr}\index{convergence to equilibrium}\index{spectral gap}
Let $N$ be a positive selfadjoint operator in a Hilbert space $H$,
and let $S$ be the $C_0$-semigroup generated by $-N$.
Denote by $\Pi$ the orthogonal projection onto $\Ker(N)$.
Then
$$
\lim_{t\to\infty}S(t)x=\Pi x,\quad x\in H.
$$
Moreover, $\lim_{t\to\infty}\n S(t)-\Pi\n=0$ if and only if
there exists $\delta>0$ such that
$$
\sigma(N)\subseteq\{0\}\cup[\delta,\infty).
$$
In this case,
$$
\n S(t)-\Pi\n\le e^{-\delta t},\quad t\ge0.
$$
\end{proposition}

\begin{proof}
Let $P$ denote the projection-valued measure of $N$ on
$[0,\infty)$. The spectral theorem gives $\Pi=P_{\{0\}}$.
For every $x\in H$ we therefore have
$$
\n S(t)x-\Pi x\n^2
=
\int_{(0,\infty)}e^{-2t\lambda}\ud P_x(\lambda)
\to0
$$
as $t\to\infty$, by dominated convergence.

If $\sigma(N)\subseteq\{0\}\cup[\delta,\infty)$, the integral
in the preceding display is bounded by
$e^{-2\delta t}\n x\n^2$. This proves the asserted estimate
and, in particular, convergence in the operator norm.

Conversely, suppose that no such $\delta$ exists.
Since the support of $P$ equals $\sigma(N)$, we have
$P_{(0,\eps)}\ne0$ for every $\eps>0$.
Choose a unit vector $x_\eps$ in the range of $P_{(0,\eps)}$.
For every $t>0$ we then have
$$
\n S(t)-\Pi\n^2
\ge
\n (S(t)-\Pi)x_\eps\n^2
=
\int_{(0,\eps)}e^{-2t\lambda}\ud P_{x_\eps}(\lambda)
\ge e^{-2t\eps}.
$$
Letting $\eps\downarrow0$ gives $\n S(t)-\Pi\n\ge1$.
The reverse inequality follows from the first display,
so $\n S(t)-\Pi\n=1$ for every $t>0$.
In particular, convergence in operator norm fails.
\end{proof}

\begin{example}[Dirichlet and Neumann heat semigroups]
\label{ex:heat-spectral-gap}
Let $D$ be a nonempty bound\-ed open subset of $\R^d$, and
let $S_{\rm Dir}$ and $S_{\rm Neum}$ denote the heat
semigroups on $L^2(D)$ generated by $\Delta_{\rm Dir}$
and $\Delta_{\rm Neum}$, respectively.

By Theorem \ref{thm:spectra-Laplacian}, the smallest
eigenvalue $\lambda_1$ of $-\Delta_{\rm Dir}$ is strictly
positive. In particular, $\Ker(\Delta_{\rm Dir})=\{0\}$.
Proposition \ref{prop:spectral-gap-sgr} gives the upper
bound in
$$
\n S_{\rm Dir}(t)\n=e^{-\lambda_1t},\quad t\ge0.
$$

Suppose now, in addition, that $D$ is connected and has a
$C^1$-boundary. By Theorem \ref{thm:spectra-Laplacian},
the eigenvalue $0$ of $-\Delta_{\rm Neum}$ has multiplicity
one and the smallest positive eigenvalue $\mu_2$ is
strictly positive. The null space consists of the constant
functions, so the orthogonal projection onto it is
$$
\Pi f
=
\Bigl(\frac1{|D|}\int_D f\ud x\Bigr)\one,
\quad f\in L^2(D).
$$
Proposition \ref{prop:spectral-gap-sgr} now gives the upper
bound in
$$
\n S_{\rm Neum}(t)-\Pi\n=e^{-\mu_2t},\quad t\ge0.
$$
Equality in both norm estimates follows by applying the
operators to norm-one eigenfunctions corresponding to
$\lambda_1$ and $\mu_2$, respectively.

Thus the Dirichlet heat flow tends exponentially to zero,
and the Neumann heat flow tends exponentially to the
spatial average of its initial value.
\end{example}

\begin{problems}

\item\label{prob:bddgen}
The aim of this problem is to prove that if $A$ generates a $C_0$-semigroup $S$
on  $X$,
then $A$ is bounded if and only if $$\lim_{t\downarrow 0} \n S(t)-I\n = 0.$$
\begin{enumerate} [\rm(a), leftmargin=*]
  \item Show that if $A$ is bounded, then $S(t) = e^{tA}$ and  $\lim_{t\downarrow 0} \n S(t)-I\n = 0$.
  \item Use a Neumann series argument to prove that if $\lim_{t\downarrow 0} \n S(t)-I\n = 0$, then
  for small enough $t>0$ the operators $T_t:= \int_0^t S(s)\ud s$
  are invertible, and show that for such $t>0$ we have $A =   T_t^{-1}(S(t)-I).$
\end{enumerate}

\item\label{prob:etA}
Let $A$ be the generator of a $C_0$-semigroup $S$
on  $X$. This problem gives a rigorous interpretation to the `formula' `$S(t) = e^{tA}$'.

\smallskip
For each $h>0$ consider the bounded operator $$A(h)x := \frac1h(S(h)x-x).$$
\begin{enumerate} [\rm(a), leftmargin=*]
  \item\label{it:etA1} Choosing $M\ge 1$ and $\omega\ge 0$ so that $\n S(t)\n \le Me^{\omega t}$ for all $t\ge 0$, show that
  $$ \n e^{tA(h)}\n \le M\exp\Bigl(\frac{t}{h}(e^{\omega h} -1)\Bigr).$$
  Deduce that for all $0<h\le 1$ we have $  \n e^{tA(h)}\n\le Me^{t (e^{\omega} -1)}\!.$
  \item\label{it:etA2} Using the identity $$S(t)x - e^{tA(h)}x = \int_0^t \frac{{\rm d}}{{\rm d}s}[e^{(t-s)A(h)}S(s)x]\ud s$$
  deduce from part \ref{it:etA1}
  that for all $x\in \Dom(A)$ and $0<h<1$ we have
  $$ \n S(t)x - e^{tA(h)} x \n \le tM^2 e^{t (\omega+e^{\omega} -1)}\n Ax - A(h)x\n.$$
  \item\label{it:etA3} Prove that for all $x\in X$ and $t\ge 0$  we have
  $ \lim_{h\downarrow 0} e^{tA(h)}x = S(t)x.$
  \item For $n\in \N$ with $n>\om$, let $A_n:= nAR(n,A)$ as in the proof of the Hille--Yosida theorem. Prove that
  for all $x\in X$ and $t\ge 0$  we have
  $ \limn e^{tA_n}x = S(t)x.$
\end{enumerate}

\item\label{prob:adjoint-sg}
Let $A$ be the generator of a $C_0$-semigroup $S$ on  $X$.
The {\em adjoint semigroup}\index{adjoint!semigroup} on $X\s$ is the family
$S\s =(S\s(t))_{t\ge 0}$, where $S\s(t) = (S(t))\s$ for $t\ge 0$.
\begin{enumerate}[\rm(a), leftmargin=*]
  \item Show that the adjoint semigroup has the semigroup properties \ref{S1} and \ref{S2} but may fail \ref{S3}.
\end{enumerate}
Set $$ X^\odot :=\{x\s\in X\s:\, \lim_{t\downarrow 0} \n S\s(t)x\s-x\s\n =0\}.$$
\begin{enumerate}[\rm(a), leftmargin=*]\setcounter{enumii}{1}
  \item Show that $X^\odot$ is a closed subspace of $X^*$ and that $S^*$ restricts to a $C_0$-semigroup $S^\odot$ on $X^\odot$.
  \item Show that $X^\odot = \overline{\Dom(A\s)}$ and deduce that  $X^\odot$ is weak$^*$ dense in $X^*$\!.
  \item Show that $$\Dom(A\s) = \Bigl\{x\s\in X\s:\, \limsup_{t\downarrow 0} \frac1t \n S\s(t)x\s-x\s\n < \infty\Bigr\}.$$
  \item Show that the generator $A^\odot$ of $S^\odot$ is the part of $A^*$ in $X^\odot$, that is,
  $$ \Dom(A^\odot) = \bigl\{x^* \in \Dom(A^*):\, A^* x^* \in X^\odot\bigr\}$$
  and $A^\odot x^* = A^* x^*$ for all $x^*\in \Dom(A^\odot)$.
\end{enumerate}

\item
By using Proposition \ref{prop:sg-prop}\ref{it:sg-prop3}, show that if a $C_0$-semigroup satisfies $\lim_{t\to\infty} S(t)x=0$ for all $x\in X$, then its generator $A$ has dense range. Do we have $0\in\varrho(A)$?

\item\label{prob:empty-spectrum}
 Let $A$ be the generator
of the $C_0$-semigroup of left translations on $L^2(0,1)$, inserting zeroes from the right. Show that $\sigma(A) = \emptyset$.

\noindent{\em Hint:}\ Apply Proposition \ref{prop:resolv}.

\item
Let $A$ be the generator of a $C_0$-semigroup $S$ on $X$.
Prove the following {\em spectral inclusion formula}\index{spectral inclusion formula}: for all $t\ge 0$ we have
$$  \exp({t\sigma(A)})\subseteq\sigma(S(t)) .$$
{\em Hint:}\
First show that for all $\la\in\C$, $t\ge 0$, and $x\in \Dom(A)$ we have
$$ e^{\la t}x-S(t)x = \int^t_0 e^{\la(t-s)} S(s) (\la-A) x\ud s
=(\la-A)\int^t_0 e^{\la(t-s)} S(s)  x\ud s.$$

\item\label{prob:TrotterKato}
This problem proves the {\em Trotter--Kato approximation theorem}.
\index{theorem!Trotter--Kato}
Let $A$ and $A_n$ generate $C_0$-semigroups $S$ and $S_n$, $n\ge1$,
on $X$ satisfying
$$
\n S(t)\n,\ \n S_n(t)\n\le Me^{\omega t},\quad t\ge0,\ n\ge1,
$$
with $M\ge1$ and $\omega\in\R$.
Fix $\lambda>\omega$ and put $R:=R(\lambda,A)$ and
$R_n:=R(\lambda,A_n)$.
\begin{enumerate}[\rm(a), leftmargin=*]
\item For $x\in X$, show that $S_n(t)x\to S(t)x$ for all $t\ge0$
implies $R_nx\to Rx$.

\item Prove that for all $x\in X$ and $t\ge0$,
$$
S_n(t)R_nRx-R_nRS(t)x
=\int_0^t S_n(t-s)(R_n-R)S(s)x\ud s.
$$
\noindent{\em Hint:}\ Differentiate $S_n(t-s)R_nRS(s)x$ with respect to $s$.

\item Assuming that $R_nx\to Rx$ for all $x\in X$, deduce that
$$
\limn\sup_{0\le t\le T}\n S_n(t)x-S(t)x\n=0,
\quad x\in X,\ T>0.
$$
\noindent{\em Hint:}\ Use Proposition \ref{prop:uniform-limits-on-K}
and the density of $\Ran(R^2)$ in $X$.

\item Show that these conclusions hold if some core $D$ for $A$
is contained in every $\Dom(A_n)$ and $A_nx\to Ax$ for all $x\in D$.

\item On $L^2(\R)$ define the bounded operators
$$
B_nf:=n^2\bigl(f(\cdot+n^{-1})-2f+f(\cdot-n^{-1})\bigr).
$$
Show that $e^{tB_n}f\to H(t)f$ uniformly for $t$ in bounded intervals
of $[0,\infty)$, where $H$ is the heat semigroup and $f\in L^2(\R)$.

\noindent{\em Hint:}\ Show that $\n e^{tB_n}\n\le1$ and use the core
$C_{\rm c}^\infty(\R)$.
\end{enumerate}

\item\label{prob:semigroup-square-root}
Let $S$ be a uniformly bounded $C_0$-semigroup on a Banach
space $X$, with generator $A$, and put
$M:=\sup_{t\ge0}\|S(t)\|$. For $x\in\Dom(A)$, define
$$
B_0x:=-\frac1{2\sqrt\pi}\int_0^\infty
\frac1{t^{3/2}} (S(t)x-x)\ud t.
$$
\begin{enumerate}[\rm(a), leftmargin=*]
\item Show that this Bochner integral converges absolutely.
\item Prove that $B_0$ is closable, and that its closure $B$ satisfies
$$
\Dom(B^2)=\Dom(A),\qquad B^2x=-Ax,\quad x\in\Dom(A).
$$

\noindent{\em Hint:}\ Fix $\eps>0$ and set
$$
Q_\eps x:=\frac1{\sqrt\pi}\int_0^\infty
t^{-1/2}e^{-\eps t}S(t)x\ud t.
$$
First show that $Q_\eps^2=R(\eps,A)$.
Then show that $Q_\eps^{-1}-B_0$, initially defined on
$\Dom(A)$, extends to a bounded operator $K_\eps$ satisfying
$\|K_\eps\|\le M\sqrt\eps$ and commuting with $Q_\eps$.
Use that $\Ran(Q_\eps^2)=\Dom(A)$ is a core for
$Q_\eps^{-1}$, and let $\eps\downarrow0$.

\item Show that if $X$ is a Hilbert space and $-A$ is positive
and selfadjoint, then $B = (-A)^{1/2}$, the square root of $-A$.
\end{enumerate}

\item
Let $A$ be the generator of a $C_0$-semigroup on $X$,
let $f\in L^1(0,T;X)$ and $u_0\in X$ be fixed.
A function $u\in L^1(0,T; X)$ is said to be a {\em weak solution}\index{weak!solution}
of the inhomogeneous abstract Cauchy problem
\begin{equation*}
\left\{\begin{aligned}
 u'(t) & = Au(t)+f(t), \quad t\in [0,T], \\
  u(0) & = u_0,
\end{aligned}
\right.
\end{equation*}
if for all $x^*\in \Dom(A^*)$ we have
$$
 \lb u(t),x^*\rb = \lb u_0,x^*\rb
 + \int_0^t \lb u(s), A^* x^*\rb\ud s
 + \int_0^t \lb f(s),x^*\rb\ud s
$$
for almost all $t\in[0,T]$.
Prove that the inhomogeneous abstract Cauchy problem has a unique weak
solution and that this solution is equal almost everywhere to the
unique strong solution.

\item\label{prob:semilinear-Schauder}
Let $A$ generate a $C_0$-semigroup $S$ on a Banach space $X$, and
suppose that $S(t)$ is compact for every $t>0$.
Fix $u_0\in X$, let
$f:[0,T]\times X\to X$ be continuous, and suppose that, for some
$c\ge0$,
$$
\n f(t,x)\n\le c(1+\n x\n),
\quad t\in[0,T],\ x\in X.
$$

\begin{enumerate}[\rm(a), leftmargin=*]

\item
With $M:=\sup_{0\le t\le T}\n S(t)\n$ and $r(t):=(M\n u_0\n+1)e^{Mct}-1$, show that
$$
r(t)
=
M\n u_0\n
+
Mc\int_0^t(1+r(s))\ud s,
\quad 0\le t\le T.
$$

\item
Show that
$$
K:=
\bigl\{
u\in C([0,T];X):\,
\n u(t)\n\le r(t),\ 0\le t\le T
\bigr\}
$$
is a nonempty closed convex subset of $C([0,T];X)$.
\item Show that $\Phi(K)\subseteq K$, where for $u\in K$ we define
$$
(\Phi u)(t)
:=
S(t)u_0+
\int_0^tS(t-s)f(s,u(s))\ud s.
$$

\item
Show that $\Phi:K\to K$ is continuous with respect to the
norm of $C([0,T];X)$.

\item
Show that for every $t\in[0,T]$ the set
$\{(\Phi u)(t):\,u\in K\}$
is relatively compact.

\noindent{\em Hint:}\
For $t>0$ and $0<\eps<t$, split the integral at $t-\eps$.

\item
Show that the family $\Phi(K)$ is equicontinuous. Deduce from the
vector-valued Arzel\`a--Ascoli theorem that $\Phi(K)$ is relatively
compact in $C([0,T];X)$.

\noindent{\em Hint:}\
Use the norm continuity of a compact
$C_0$-semigroup on $(0,\infty)$ together with the bound on $f$.

\item
Apply the Schauder fixed point theorem of Problem
\ref{prob:Schauder} to conclude that the semilinear Cauchy problem \eqref{eq:SDE}
admits a mild solution on $[0,T]$.

\item
Explain why this argument gives no uniqueness assertion in general.
\end{enumerate}

\item
Let $B$ be a bounded operator on $X$. Show that if $A$ generates a bounded analytic $C_0$-semigroup on  $X$, then so does $A+B-c I $ for large enough $c$. Conclude that if $A$ generates an analytic $C_0$-semigroup, then so does $A+B$.

\item
Let $A$ generate a bounded $C_0$-group $(U(t))_{t\in \R}$ on $X$.
\begin{enumerate}[\rm(a), leftmargin=*]
\item Show that $\Dom(A^2)$ is dense in $X$.
\item Show that, for $\Re\lambda>0$, the operator $\lambda^2-A^2$
is bijective, with bounded inverse
$(\lambda^2-A^2)^{-1}
=\frac1{2\lambda}(R(\lambda,A)-R(-\lambda,A))$.
Deduce that $A^2$ is closed.
\item Deduce that $A^2$ generates a $C_0$-semigroup $(S(t))_{t\ge 0}$ which is analytic and bounded
on every sector $\Sigma_\eta$ with $0<\eta<\frac12\pi$.
\item Prove the formula
$$
S(t)x=\frac1{\sqrt{4\pi t}}
\int_\R e^{-s^2/(4t)}U(s)x\ud s,
\quad t>0,\ x\in X.
$$
{\em Hint:} \ First take $x\in\Dom(A^2)$.
Differentiate the integral in $t$, integrate by parts twice in $s$,
and use uniqueness of classical solutions.
\end{enumerate}

\item
Let $A$ be the generator of a contractive analytic $C_0$-semigroup $S$ on a Hilbert space $H$ of angle $\omega\in (0,\frac12\pi)$. Show that the form $\aa$ in $H$ defined by
$$\aa(x,y):= -\iprod{Ax}{y}, \quad x,y\in \Dom(\aa):= \Dom(A),$$ is accretive and closable, that its closure is $(\frac12\pi-\omega)$-sectorial, and that the
Hille--Yosida estimate, when applied to the semigroups $S(e^{i\theta} t)$ with $|\theta|<\omega$, reproduces the resolvent estimates of Theorem \ref{thm:form-sectorial-resolvent}.

\noindent{\em Hint:}\ First use the contractivity of the rotated semigroups
to show that $-\iprod{Ax}{x}\in\ov{\Sigma_{\frac12\pi-\omega}}$
for $x\in\Dom(A)$, as in the proof of Theorem \ref{thm:lumerphilips}.
Then apply Problem \ref{prob:sectorial-form-closable} to $-A$.

\item\label{prob:bdry-sgr}
Let $A$ be the generator of a $C_0$-semigroup $S$ on $X$ which is bounded analytic on the open right half-plane.
Show that for all $t\in \R$ and $x\in X$ the limit
$$T(t)x:= \lim_{s\downarrow 0} S(s+it)x$$
exists and that the family $(T(t))_{t\in \R}$ is a uniformly bounded $C_0$-group of operators with
generator $iA$.

\noindent{\em Hint:}\ Begin by observing that if $0<s'<s<\infty$ and $t\in\R$, then $$\n S(s+it)x - S(s'+it)x\n \le M \n S(s-s')x-x\n,$$
where $M = \sup_{\Re z>0} \n S(z)\n$. Deduce from this the existence of the limits. Deduce the semigroup properties and strong continuity in a similar manner. Finally, show that if $x\in\Dom(A)$, then $\int_0^tT(s)iAx\ud s=T(t)x-x$.
Deduce that $iA\subseteq B$,
where $B$ is the generator of $(T(t))_{t\in\R}$.
Show also that $T(t)$ leaves $\Dom(A)$ invariant. Proposition
\ref{prop:core} then shows that $\Dom(A)$ is a core for $B$.
Since $iA$ is closed, conclude that $B=iA$.

\item\label{prob:ues}
Let $A$ be the generator of an analytic $C_0$-semigroup $S$ on  $X$. Prove that if $$\sigma (A)\subseteq \{z\in\C:\, \Re z<0\},$$ then $S$ is {\em uniformly exponentially stable},
\index{exponentially stable, uniformly} that is, there exist $M\ge 1$ and $\om>0$ such that $ \n S(t)\n \le Me^{-\om t}$ for all $t\ge 0$.

\noindent {\em Hint:}\ Verify the assumptions of Theorem \ref{thm:analytic} for $\delta+A$ for small enough $\delta>0$.

\item\label{prob:DP}
Prove the {\em Datko--Pazy theorem}:\index{theorem!Datko--Pazy} for any given $1\le p<\infty$, a $C_0$-semigroup $S$ on $X$ is uniformly exponentially stable (see  Problem \ref{prob:ues}) if and only if the orbit $t\mapsto S(t)x$ belongs to $L^p(\R_+;X)$ for all $x\in X$.

\noindent{\em Hint:}\ Apply the uniform boundedness theorem and argue by contradiction.

\item\label{prob:GP}
This problem proves the {\em Gearhart--Pr\"uss theorem}\index{theorem!Gearhart--Pr\"uss}: a $C_0$-semigroup $(S(t))_{t\geq 0}$ with generator $A$ on a Hilbert space $H$ is uniformly exponentially stable (see  Problem \ref{prob:ues}) if and only if
$\{\la\in \C:\, \Re \la>0\} \subseteq \varrho(A)$ and $$\sup_{\Re \la >0} \n R(\la,A)\n <\infty.$$
Suppose that $\n S(t)\n \le Me^{\om t}$.  Complete the following steps:
\begin{enumerate}[\rm(a), leftmargin=*]
  \item Extend the Fourier--Plancherel theorem to $L^2(\R;H)$.
  \item Prove that $t\mapsto R(s+it,A)x$ belongs to $L^2(\R;H)$ for all $x\in H$ and $s>\om$.
\end{enumerate}
By Lemma \ref{lem:res-Taylor-cor} there exists a $\delta>0$ such that $\{\la\in \C:\, \Re \la>-\delta\} \subseteq \varrho(A)$ and $$\sup_{\Re \la >-\delta} \n R(\la,A)\n <\infty.$$
\begin{enumerate}[\rm(a), leftmargin=*]\setcounter{enumii}{2}
  \item Use the resolvent identity to prove that the function $t\mapsto R(it,A)x$ belongs to $L^2(\R;H)$ for all $x\in H$.

\noindent{\em Hint:}\
Fix $s_0>\max\{\omega,0\}$. Use the resolvent identity to show that
$$
\sup_{0<s\le s_0}\int_{\R}\|R(s+i\xi,A)x\|^2\ud\xi<\infty.
$$
Now let $s\downarrow0$ and apply Fatou's lemma.

  \item Conclude that $t\mapsto S(t)x$ belongs to $L^2(\R_+;H)$ for all $x\in H$.

\noindent{\em Hint:}\ For $s>0$ let $g_s$ be the inverse Fourier--Plancherel
transform of $\xi\mapsto R(s+i\xi,A)x$. Show that for almost all $t>0$,
$$
g_s(t)=\sqrt{2\pi}\,\one_{(0,\infty)}(t)e^{-st}S(t)x.
$$
Apply Plancherel's theorem and pass to the limit $s\downarrow0$. Finally, apply Problem \ref{prob:DP}.
\end{enumerate}

 \item\label{prob:GoldysvanNeerven}
 This problem gives a two-dimensional example of a bounded analytic
 $C_0$-semi\-group which is uniformly exponentially stable, contractive on $\R_+$,
 and fails to be contractive on any open sector containing $\R_+$.

 On $\C^2$ define
 $\iprod{x}{y}_Q:= \iprod{Qx}{y},$ where
 $ Q = \begin{pmatrix} 1 & 1 \\ 1 & 2 \end{pmatrix}.$
 \begin{enumerate}[\rm(a), leftmargin=*]
   \item Show that $\iprod{\cdot}{\cdot}_Q$ defines an inner product on $\C^2$\!.
 \end{enumerate}
 Let $\n \cdot\n_Q$ be the associated norm. On $(\C^2\!,\n\cdot\n_Q)$ we consider the $C_0$-semigroup $S$,
 $$ S(t) = e^{-t/2}\begin{pmatrix} 1 & \,t \\ 0 & \,1\end{pmatrix}, \quad t\ge 0.$$
 \begin{enumerate}[\rm(a), leftmargin=*]\setcounter{enumii}{1}
   \item Show that $$\n S(t)\n_Q^2 = \frac12 e^{-t }(t^2 + 2+t \sqrt{t^2+4})$$
   and conclude that $S(t)$ is contractive with respect to $\n \cdot\n_Q$ for all $t\ge 0$.

 \noindent
 {\em Hint:}\ Use the fact that
 $\n S(t)\n_Q^2$ equals the largest eigenvalue of $S(t)S^\star(t)$ (see  Problem \ref{prob:hermit}), where the adjoint refers to the inner product $\iprod{\cdot}{\cdot}_Q$.
   \item Show that $S$ extends to an entire $C_0$-semigroup
   which is uniformly bounded on the open sector $\Sigma_\eta$ for all $0<\eta<\frac12\pi$.
   \item Show that $S$ fails to be contractive on any open sector $\Sigma_\eta$.
 \end{enumerate}

\item\label{prob:LuPhi}
Let $X$ be a nonzero Banach space.
For $x\in X$ define the {\em subdifferential} of $x$\index{subdifferential}\index{$D$@$\partial(x)$} by
$$
\partial(x):=\{x^*\in X^*:\, \|x^*\|\le 1,\ \langle x,x^*\rangle=\|x\|\}.
$$
\begin{enumerate}[\rm(a), leftmargin=*]

\item Show that $\partial(x)\not=\emptyset$.

\item Show that if $X$ is a Hilbert space, then under the usual identification of $X$ with its dual, for all nonzero $x\in X$ we have $$\partial(x) = \{x/\n x\n\}.$$
\item Let $1<p<\infty$.
Show that for $X=L^p(\Om)$ and nonzero $f\in L^p(\Om)$ we have
$$
\partial(f)=\{g_{f,p}/\|f\|_p^{p-1}\},
$$ where
$f(\om) = e^{i\theta(\om)}|f(\om)|$ and
$ g_{f,p}(\om) =  e^{-i\theta(\om)}|f(\om)|^{p-1}\!.$
\item Is the subdifferential always a single\-ton?
\item Extend the Lumer--Phillips theorem (Theorem \ref{thm:lumerphilips}) to Banach spaces.\index{theorem!Lumer--Phillips}
\end{enumerate}

\item\label{prob:gaussianW12}

Define the Gaussian Sobolev space $W^{1,2}(\R^d\!,\gamma)$ to consist of all $f\in L^2(\R^d\!,\gamma)$ whose weak first-order derivatives exist and belong to $L^2(\R^d\!,\gamma)$.

\begin{enumerate}[\rm(a), leftmargin=*]
 \item Show that $W^{1,2}(\R^d\!,\gamma)$ is a Hilbert space in a natural way.
 \item Via truncation and mollification, show that $C_{\rm c}^\infty(\R^d)$ is dense in $W^{1,2}(\R^d\!,\gamma)$.

 \item Show that the gradient $\nabla$ is closed as a densely defined operator
 acting from $L^2(\R^d\!,\gamma)$ to $L^2(\R^d\!,\gamma;\C^d)$ with domain $\Dom(\nabla) := W^{1,2}(\R^d\!,\gamma)$.
 \item Show that the form domain of the closed form $\ov \aa_{OU}$ is precisely $W^{1,2}(\R^d\!,\gamma)$.
\end{enumerate}

\end{problems}

%% file: ch14-HS-TC.tex
\chapter{Trace Class Operators}\label{chap-HilbertSchmidt-TraceClass}

\blfootnote{\rm \small \begin{quote}This book has been published by Cambridge University Press in the series `Cambridge Studies in Advanced Mathematics'. The present version is free to view and download for personal use only. Not for re-distribution, re-sale or use in derivative works.\newline {\copyright} Jan van Neerven\end{quote}}

\noindent
This chapter is devoted to the study of trace class operators and the related class of Hilbert--Schmidt operators. In a sense that will be explained in the next chapter, we can think of positive trace class operators as analogues of finite measures, with the trace corresponding to total mass. After proving some general properties of
trace class operators, we compute traces in a number of interesting examples.

\section{Hilbert--Schmidt Operators}\label{sec:Hilbert-Schmidt}

Throughout this chapter we assume that all Hilbert spaces are separable and complex.

Let $H$ and $K$ be separable Hilbert spaces.

\begin{definition}[Hilbert--Schmidt operators]
A bounded operator $T\in\calL(H,K)$ is called a
{\em Hilbert--Schmidt operator}\index{operator!Hilbert--Schmidt}
if the quantity
$$
\n T\n_{\calL_2(H,K)}
:=
\Bigl(\sum_{n\ge1}\n Th_n\n_K^2\Bigr)^{1/2}
$$
is finite for some (equivalently, for every) orthonormal basis
$(h_n)_{n\ge1}$ of $H$.
The space of these operators is denoted by $\calL_2(H,K)$.
When $K=H$, we write $\calL_2(H)$.
\end{definition}

This definition is independent of the choice of orthonormal basis.
Indeed, let $(h_n)_{n\ge1}$ and $(k_j)_{j\ge1}$ be orthonormal
bases of $H$ and $K$, respectively.
Parseval's identity and Tonelli's theorem give
\begin{align*}
\sum_{n\ge1}\n Th_n\n_K^2
&=\sum_{n\ge1}\sum_{j\ge1}
  |\iprod{Th_n}{k_j}_K|^2
=\sum_{j\ge1}\sum_{n\ge1}
  |\iprod{h_n}{T^\star k_j}_H|^2
=\sum_{j\ge1}\n T^\star k_j\n_H^2.
\end{align*}
For a fixed basis of $K$, the right-hand side is independent
of the choice of basis of $H$, proving the claim.

It is clear that
$\n T\n_{\calL(H,K)}\le\n T\n_{\calL_2(H,K)}$.
We have also shown:

\begin{proposition}\label{prop:HilbSchm-adjoint}
A bounded operator $T\in\calL(H,K)$ is a Hilbert--Schmidt
operator if and only if its adjoint
$T^\star\in\calL(K,H)$ is a Hilbert--Schmidt operator.
In this case,
$$
\n T\n_{\calL_2(H,K)}
=
\n T^\star\n_{\calL_2(K,H)}.
$$
\end{proposition}

For $g\in K$ and $h\in H$ we recall the notation $g\,\bar\otimes\, h$ for the operator from $H$ to $K$ defined by
\begin{align*}
(g\,\bar\otimes\, h)x:= \iprod{x}{h}g, \quad x\in H.
\end{align*}

\begin{example}[Finite-rank operators]\label{ex:FinRank-HS}
Every finite-rank operator from $H$ to $K$ is a Hilbert--Schmidt operator. Indeed,
by a Gram--Schmidt argument we may represent any such operator $T$ as $\sum_{j=1}^k g_j\,\bar\otimes\,h_j$ with $g_1,\dots,g_k\in K$ orthonormal and $h_1,\dots,h_k\in H$.
Taking an orthonormal basis $(e_n)_{n\ge 1}$ of $H$, we have
\begin{align*}\sum_{n\ge 1} \n Te_n\n^2 = \sum_{n\ge 1} \sum_{j=1}^k |\iprod{e_n}{h_j}|^2
= \sum_{j=1}^k \sum_{n\ge 1} |\iprod{e_n}{h_j}|^2
= \sum_{j=1}^k \|h_j\|^2\!.
\end{align*}
\end{example}

\begin{example}[Integral\, operators\, with\, square-integrable\, kernels]\label{ex:kernel-HS}\
Let\, $(\Om,\mu)$\, and $(\Om',\mu')$ be $\sigma$-finite measure spaces such that
$L^2(\Om,\mu)$ and $L^2(\Om',\mu')$ are separable and let
$k\in L^2(\Om'\times\Om,\mu'\times\mu)$. Then, for every
$f\in L^2(\Om,\mu)$, the function $t\mapsto k(s,t)f(t)$
belongs to $L^1(\Om,\mu)$ for $\mu'$-almost all $s\in\Om'$, and for such $s$ we may define the bounded operator $T$ by
$$
Tf(s):=\int_\Om k(s,t)f(t)\ud\mu(t).
$$
This operator is Hilbert--Schmidt from $L^2(\Om,\mu)$ to $L^2(\Om',\mu')$, with norm
$$
\|T\|_{\calL_2(L^2(\Om,\mu),L^2(\Om',\mu'))}
=\|k\|_{L^2(\Om'\times\Om)}.
$$
Indeed, if $(h_n)_{n\ge 1}$ is an orthonormal basis of $L^2(\Om,\mu)$, then
\begin{align*}
\sum_{n\ge 1} \int_{\Om'}\Big|\int_{\Om} k(s,t)h_n(t)\ud \mu(t)\Big|^2\ud\mu'(s)
 & = \int_{\Om'}\sum_{n\ge 1}\Big|\int_{\Om} k(s,t)h_n(t)\ud \mu(t)\Big|^2\ud\mu'(s)
\\ & = \int_{\Om'}\n k(s,\cdot)\n_{L^2(\Om,\mu)}^2\ud\mu'(s)
 = \n k\n_{L^2(\Om'\times\Om,\mu'\times\mu)}^2.
\end{align*}
As a special case, any $m\times d$ matrix
$A = (a_{jk})_{1\le j\le m,\,1\le k\le d}$ is a Hilbert--Schmidt operator
as a linear operator from $\C^d$ to $\C^m$\!, and
$$
\n A\n_{\calL_2(\C^d,\C^m)}^2
= \sum_{j=1}^m\sum_{k=1}^d |a_{jk}|^2\!.
$$
\end{example}

A converse to this example will be stated at the end of this section.

\begin{proposition}\label{prop:HS-inner}
The space $\calL_2(H,K)$\index{$L$@$\calL_2(H,K)$} of all Hilbert--Schmidt operators from $H$ to $K$ is a Hilbert space with respect to the inner product
 $$\iprod{T_1}{T_2} := \sum_{n\ge 1} \iprod{T_1h_n}{T_2h_n},$$
where $(h_n)_{n\ge 1}$ is any orthonormal basis of $H$.
\end{proposition}

\begin{proof}
For $T_1,T_2\in\calL_2(H,K)$, the series
$\sum_{n\ge1}\iprod{T_1h_n}{T_2h_n}$
converges absolutely by the Cauchy--Schwarz inequality. It therefore
defines a sesquilinear form whose associated norm is the
Hilbert--Schmidt norm. Since this norm is independent of the chosen
orthonormal basis, the polarisation identity shows that the
sesquilinear form is independent of the basis as well. It is plainly
an inner product.

The triangle inequality in $\ell^2$
implies that $\calL_2(H,K)$ is a normed space. To prove completeness, suppose that $(T_n)_{n\ge 1}$ is a Cauchy sequence in $\calL_2(H,K)$.
Then $(T_n)_{n\ge 1}$ is a Cauchy sequence in $\calL(H,K)$. Let $T\in \calL(H,K)$ be its limit.
If $(h_j)_{j\ge 1}$ is an orthonormal basis for $H$, then for all $n\ge 1$ we have
$$ \sum_{j=1}^n \n Th_j\n^2 = \limk \sum_{j=1}^n \n T_k h_j\n^2 \le \limk \n T_k\n_{\calL_2(H,K)}^2 <\infty.$$
Letting $n\to\infty$, we see that $T$ is a Hilbert--Schmidt operator and $$ \n T\n_{\calL_2(H,K)}\le \limk \n T_k\n_{\calL_2(H,K)}.$$
Also,
$$\sum_{j=1}^n \n (T_k-T)h_j\n^2 = \lim_{m\to\infty} \sum_{j=1}^n \n (T_k -T_m)h_j\n^2 \le \limsup_{m\to\infty} \n T_k-T_m\n_{\calL_2(H,K)}^2.$$
It follows that $\n T_k-T\n_{\calL_2(H,K)} \le \limsup_{m\to\infty} \n T_k-T_m\n_{\calL_2(H,K)}$.
Since the latter tends to $0$ as $k\to\infty$, it follows that $\limk T_k = T$ in $\calL_2(H,K)$.
This proves completeness.
\end{proof}

\begin{proposition}\label{prop:HS-compact} Every Hilbert--Schmidt operator from $H$ to $K$ is compact and can be approximated, in the Hilbert--Schmidt norm, by finite-rank operators.
\end{proposition}
 \begin{proof} Let $T$ be a Hilbert--Schmidt operator from $H$ to $K$,
  let $(h_n)_{n\ge 1}$ be an orthonormal basis for $H$, and denote by $P_N$ the orthogonal projection onto the span of $\{h_1,\dots,h_N\}$.
  Then $TP_N$ is a finite-rank operator and hence Hilbert--Schmidt, and we have
  $$\limsup_{N\to\infty} \n TP_N -T\n^2 \le \limsup_{N\to\infty} \n TP_N -T\n_{\calL_2(H,K)}^2 = \limsup_{N\to\infty}\sum_{n\ge N+1} \n Th_n\n^2=0.$$
  Each $T P_N$ is a finite-rank operator, hence compact. Since uniform limits of compact operators are compact, it follows that $T$ is compact.
\end{proof}

Hilbert--Schmidt operators have the following ideal property:

\begin{proposition}\label{prop:ideal-HS}\index{ideal property!of Hilbert--Schmidt operators} If $T$ is a Hilbert--Schmidt operator from $H$ to $K$ and $S\in\calL(K)$ and $U\in\calL(H)$ are bounded, then $STU$ is a Hilbert--Schmidt operator
and $$\n STU\n_{\calL_2(H,K)} \le \n S\n \n T\n_{\calL_2(H,K)} \n U\n.$$
\end{proposition}
\begin{proof}
It is clear that $ST$ is a Hilbert--Schmidt operator
and $$\n ST\n_{\calL_2(H,K)} \le \n S\n \n T\n_{\calL_2(H,K)}.$$
Applying the preceding estimate to $U^\star$ and $T^\star$ and using
Proposition \ref{prop:HilbSchm-adjoint}, we find that
$U^\star T^\star$ is Hilbert--Schmidt and
$$\n U^\star T^\star\n_{\calL_2(K,H)} \le \n U^\star\n \n T^\star\n_{\calL_2(K,H)}=\n U\n \n T\n_{\calL_2(H,K)}.$$
Hence $TU=(U^\star T^\star)^\star$ is Hilbert--Schmidt
and $$\n TU\n_{\calL_2(H,K)}= \n U^\star T^\star\n_{\calL_2(K,H)} \le \n U\n \n T\n_{\calL_2(H,K)}.$$ Using the first step once more, we conclude that $STU$ is a Hilbert--Schmidt operator and
satisfies the estimate in the statement of the proposition.
\end{proof}

\begin{theorem}\label{thm:Mercer-HS}
Let $(\Om,\mu)$ and $(\Om',\mu')$ be $\sigma$-finite measure spaces such that
$L^2(\Om,\mu)$ and $L^2(\Om',\mu')$ are separable.
If $T\in\calL_2(L^2(\Om,\mu),L^2(\Om',\mu'))$, there exists a unique
$k\in L^2(\Om'\times\Om,\mu'\times\mu)$ such that for all $f\in L^2(\Om,\mu)$ we have
$$ Tf(\om') = \int_\Om k(\om',\om) f(\om)\ud\mu(\om)$$
for $\mu'$-almost all $\om'\in\Om'$.
\end{theorem}

\begin{proof}
The mapping
$\Phi:L^2(\Om'\times\Om,\mu'\times\mu)\to\calL_2(L^2(\Om,\mu),L^2(\Om',\mu'))$, $\Phi(k):=T_k$,
where $T_k$ is the integral operator with kernel $k$ as in the example, is a linear
isometry by Example \ref{ex:kernel-HS}. Since $L^2(\Om'\times\Om,\mu'\times\mu)$ is complete, $\Phi$ has closed range.

For $g\in L^2(\Om',\mu')$ and $h\in L^2(\Om,\mu)$, the function
$k(\om',\om):=g(\om')\ov{h(\om)}$ belongs to
$L^2(\Om'\times\Om,\mu'\times\mu)$ and satisfies
$T_k=g\,\bar\otimes\,h$. Hence the range of $\Phi$ contains every
finite-rank operator. These
operators are dense in $\calL_2(L^2(\Om,\mu),L^2(\Om',\mu'))$ by Proposition \ref{prop:HS-compact}, so $\Phi$ is surjective.
Uniqueness of the kernel follows from the fact that $\Phi$ is an isometry.
\end{proof}

\section{Trace Class Operators}\label{sec:trace}

Recall that a bounded operator $T\in \calL(H)$ is called {\em positive}
if $\iprod{Th}{h}\ge 0$ for all $h\in H$. Since the scalar field is complex, every bounded positive operator is selfadjoint.

\begin{definition}[Trace, of a positive operator] The {\em trace}\index{trace} of a positive operator $T\in \calL(H)$ is
the nonnegative extended-real number defined by
 $$ \tr(T) := \sum_{n\ge 1} \iprod{Th_n}{h_n},$$
 where $(h_n)_{n\ge 1}$ is any orthonormal basis of $H$.
\end{definition}

To see that $\tr(T)$ is well defined, suppose that $(h_n)_{n\ge 1}$ and $(h_n')_{n\ge 1}$
are orthonormal bases of $H$ and that $\sum_{n\ge 1}\iprod{Th_n}{h_n} < \infty$. Then $T^{1/2}$ is Hilbert--Schmidt and
$$ \sum_{n \ge 1} \iprod{Th_n'}{h_n'} =
   \sum_{n \ge 1} \n T^{1/2}h_n'\n^2  =
   \sum_{n\ge 1} \n T^{1/2}h_n\n^2  =
   \sum_{n\ge 1}\iprod{Th_n}{h_n}.
$$
By interchanging the roles of the two bases, the same conclusion is reached when $\sum_{n\ge 1}\iprod{Th_n'}{h_n'} < \infty$, and the result follows.

\begin{definition}[Trace class operators] A bounded operator $T\in\calL(H)$ is called a {\em trace class operator}\index{operator!trace class}
if its modulus $|T|:= (T^\star T)^{1/2}$ has finite trace.
\end{definition}

\begin{proposition} If $T\in \calL(H)$ is a trace class operator, then
$ \n T\n \le \tr(| T|)$.
\end{proposition}

\begin{proof}
By Proposition \ref{prop:normal-spectralradius} we have $\n T\n^2 = \n T^\star T\n = \n (T^\star T)^{1/2}\n^2$, and therefore $\|T\|=\||T|\|$. In view of this, it is enough to prove that $\|\,|T|\,\|\le \tr(|T|)$.

Let $h\in H$ with $\|h\|=1$, and choose an orthonormal basis
$(h_n)_{n\ge1}$ of $H$ with $h_1=h$. Since $|T|$ is positive,
$$
        \iprod{|T|h}{h}
        \le
        \sum_{n\ge1}\iprod{|T|h_n}{h_n}
        =
        \tr(|T|),
$$
and therefore
$\|\,|T|\,\| = \sup_{\|h\|=1}\iprod{|T|h}{h} \le \tr(|T|)$.
\end{proof}

\begin{example}[Finite-rank operators]\label{ex:FinRank-TC}
Every finite-rank operator $T$ is a trace class operator.
Indeed, the proof of Theorem \ref{thm:sing-value-comp} gives a representation
$$|T|= \sum_{n\ge 1} \mu_n h_n\,\bar\otimes\,h_n,$$
where $(\mu_n)_{n\ge 1}$ is the sequence of nonzero eigenvalues of $|T|$
repeated according to multiplicities and the orthonormal sequence $(h_n)_{n\ge 1}$ consists of eigenvectors of $|T|$.
Then $$T^\star T = \sum_{n\ge 1} \mu_n^2 h_n\,\bar\otimes\,h_n,$$
and since $T^\star T$ has finite rank, this sum must be
a finite sum. Therefore the same is true for the sum representing $|T|$.
This implies that $\tr(|T|) =\sum_{n=1}^N \mu_n$ is finite.
\end{example}

\begin{proposition} Every trace class operator is compact.
\end{proposition}
\begin{proof}
 First assume that $T$ is a positive  trace class operator.  Let $(h_n)_{n\ge 1}$
be an orthonormal basis of $H$.
 Then from
 $$\sum_{n\ge 1} \n T^{1/2} h_n\n^2 = \sum_{n\ge 1} \iprod{Th_n}{h_n} = \tr(T) <\infty$$
 we see that $T^{1/2}$ is a Hilbert--Schmidt operator and therefore compact.  Hence also $T = (T^{1/2})^2$ is compact.

 In the general case let $T = U|T|$ be the polar decomposition of $T$, with $U$ an isometry from $\ov{\Ran(|T|)}$ onto $\ov{\Ran(T)}$
 (see Theorem \ref{thm:polar}).
 Since the positive operator $|T|$ is a trace class operator, $|T|$ is compact, hence so is $T$.
\end{proof}

Let $T$ be a compact operator, with polar decomposition $T = U|T|$.
Viewing $U$ as an isometry from $\ov{\Ran(|T|)}$ onto $\ov{\Ran(T)}$,
its adjoint $U^\star$ is an isometry from $\ov{\Ran(T)}$ onto $\ov{\Ran(|T|)}$
satisfying $U^\star U = P_{\,\ov{\Ran(|T|)}}$, and consequently $|T| = U^\star T$.
It follows that $|T|$ is compact.

\begin{definition}[Singular values]
 The {\em singular values}\index{singular value} of a compact operator $T$
on $H$ are the nonzero
eigenvalues of the compact operator $|T|$.
\end{definition}

Since $|T|$ is positive, every singular value is a strictly positive real number, and since $|T|$ is compact,
the set of singular values is finite or countable with $0$ as its only possible accumulation point. We may therefore think of the set of singular values as a
nonincreasing (finite or infinite) sequence $(\mu_n)_{n\ge 1}$. This sequence, where each $\mu_n$ is repeated according to its multiplicity,
is called the {\em singular value sequence}.\index{singular value!sequence}
The singular value sequence $(\mu_n)_{n\ge 1}$ of a compact normal operator $T\in \calL(H)$
satisfies $\mu_n=|\la_n|$, provided eigenvalues are repeated according to multiplicity and
the sequences are ordered by decreasing absolute value;
this is immediate from the spectral theory of these operators.

To understand the role of the singular values, it is useful to recall how the singular value decomposition is constructed. Since $|T|$ is compact and selfadjoint, its nonzero eigenspaces
are finite-dimensional and mutually orthogonal. Choosing an orthonormal basis in each of these eigenspaces, and listing the resulting vectors in the same order as the singular values,
gives an orthonormal family $(h_n)_{n\ge 1}$ satisfying
$$
|T|h_n=\mu_nh_n.
$$
Here the repetition of each singular value according to its multiplicity accounts for all the basis vectors in the corresponding eigenspace.

By the spectral theorem for compact selfadjoint operators (Theorem \ref{thm:spect-thm-comp}), the closed linear span of these eigenvectors is the orthogonal complement of $\ker(|T|)$.
Thus
$$
H=\ker(|T|)\oplus
\overline{\operatorname{span}}\{h_n: \, n\ge 1\}.
$$
Choosing an orthonormal basis of $\ker(|T|)$ and adjoining it to $(h_n)_{n\ge 1}$ therefore gives an orthonormal basis $(e_j)_{j\ge 1}$ of $H$. The additional basis vectors correspond to the eigenvalue zero, which is not included in the singular value sequence.
Moreover, the identity
$$
\n Tx\n^2
=\iprod{T^\star Tx}{x}
=\iprod{|T|^2x}{x}
=\n |T|x\n^2,
\quad x\in H,
$$
shows that $\ker(T)=\ker(|T|)$. Hence both $T$ and $|T|$ vanish on these additional basis vectors.
On the remaining basis vectors, $|T|$ acts by multiplication by the corresponding singular values. Consequently,
$$
|T|=\sum_n\mu_n h_n\,\bar\otimes\,h_n,
$$
with convergence in operator norm. In the infinite case, this convergence follows from $\mu_n\to0$: the norm of the remainder after the first $N$ terms is $\mu_{N+1}$.

To recover $T$, write $T=U|T|$ and put $g_n:=Uh_n$. Since $\mu_n>0$, the identity
$h_n=\mu_n^{-1}|T|h_n$ shows that every $h_n$ belongs to $\Ran(|T|)$. The operator $U$ is isometric on $\ov{\Ran(|T|)}$, so the family $(g_n)_{n\ge 1}$ is orthonormal.
Furthermore,
$$
Th_n=\mu_ng_n
$$
and $\n Th_n\n=\mu_n$. Multiplying the expansion of $|T|$ by $U$ gives
$$
T=\sum_n\mu_n g_n\,\bar\otimes\,h_n,
$$
again with convergence in operator norm. This is the singular value decomposition of Theorem \ref{thm:sing-value-comp}. If $T$ has finite rank, the families and sums in this
discussion are finite.

In the next theorem we use the notation
$$\n T\n_{\calL_1(H)} := \tr(|T|).$$
Later in this section we prove that the set $\calL_1(H)$ of all trace class operators on $H$ is a vector space, that $\tr(|T|)$ defines a norm on it, and that the resulting normed space is a Banach space.

\begin{theorem}[Singular value decomposition]\label{thm:tc-ell1}\index{theorem!singular value decomposition}\index{singular value!decomposition}
Let $T\in\calL(H)$ be compact, and let $(\mu_n)_{n\ge 1}$ be its singular value sequence. Then:
\begin{enumerate}[label={\rm(\arabic*)}, leftmargin=*]
\item\label{it:tc-ell1-1} $T$ is a trace class operator if and only if
$\sum_{n\ge 1} \mu_n <\infty$. In this case we have $$\n T\n_{\calL_1(H)} = \sum_{n\ge 1} \mu_n.$$
\item\label{it:tc-ell1-2} $T$ is a Hilbert--Schmidt operator if and only if
$\sum_{n\ge 1} \mu_n^2 <\infty$. In this case we have $$\n T\n_{\calL_2(H)}^2  = \sum_{n\ge 1} \mu_n^2.$$
\end{enumerate}
Putting $$ T_N:=\sum_{n=1}^N \mu_n g_n\,\bar\otimes\,h_n, $$
where $(g_n)_{n\ge 1}$ and $(h_n)_{n\ge 1}$ are the orthonormal
sequences in a singular value decomposition of $T$,
in the situation of \ref{it:tc-ell1-1} we have
$\n T-T_N\n_{\calL_1(H)}\to 0$ and in the situation of \ref{it:tc-ell1-2} we have $\n T-T_N\n_{\calL_2(H)}\to 0$.
If $T$ is positive we may take $(g_n)_{n\ge 1} = (h_n)_{n\ge 1}$.
\end{theorem}
\begin{proof}
Choose a singular value decomposition
$$
T=\sum_{n\ge 1}\mu_n g_n\,\bar\otimes\,h_n, \qquad |T|=\sum_{n\ge 1}\mu_n h_n\,\bar\otimes\,h_n,
$$
with $|T|h_n=\mu_nh_n$ for $n\ge 1$, and complete $(h_n)_{n\ge 1}$
by an orthonormal basis of $\ker(|T|)$ to an orthonormal basis
$(e_j)_{j\ge 1}$ of $H$. Since $\n Tx\n=\n |T|x\n$ for all $x\in H$,
the operator $T$ vanishes on $\ker(|T|)$ and
$\n Th_n\n=\mu_n$. Consequently,
$$
\tr(|T|)=\sum_{n\ge 1}\mu_n,
\qquad
\sum_j\n Te_j\n^2=\sum_{n\ge 1}\mu_n^2,
$$
both in the sense that if one side is finite, then so is the other side and equality holds.
This proves the two characterisations and their norm identities.

For the truncations $T_N$ in the statement, the nonzero singular
values of $T-T_N$ are precisely the remaining terms
$(\mu_n)_{n>N}$. Applying the identities just proved to this
tail operator gives, in the respective cases,
\begin{align*}
\n T-T_N\n_{\calL_1(H)}
=\sum_{n>N}\mu_n\to0,\qquad
\n T-T_N\n_{\calL_2(H)}^2
=\sum_{n>N}\mu_n^2\to0.
\end{align*}
If $T$ is positive, we may take $g_n=h_n$ in its singular value
decomposition.
\end{proof}

We now extend the trace from positive operators to arbitrary trace
class operators.

\begin{theorem}[Trace]\label{thm:tc}
If $T\in\calL(H)$ is a trace class operator, for any orthonormal basis $(h_n)_{n\ge 1}$ of $H$ the sum
$$\tr(T) := \sum_{n\ge 1} \iprod{Th_n}{h_n}$$
converges absolutely, the sum is independent of the choice of the basis,
and $$|\tr(T)|\le \tr(|T|).$$
\end{theorem}

\begin{proof}
Let $T=U|T|$ be the polar decomposition, and put $R:=|T|^{1/2}$.
Then $R$ is Hilbert--Schmidt, with
$\n R\n_{\calL_2(H)}^2=\tr(|T|)$.
Since $\n U\n\le1$, Proposition \ref{prop:ideal-HS} implies that
$RU^\star$ is Hilbert--Schmidt and
$\n RU^\star\n_{\calL_2(H)}\le\n R\n_{\calL_2(H)}$.

For any orthonormal basis $(h_n)_{n\ge1}$ of $H$, we have
$\iprod{Th_n}{h_n}=\iprod{Rh_n}{RU^\star h_n}$.
Hence the Cauchy--Schwarz inequality gives
\begin{align*}
\sum_{n\ge1}|\iprod{Th_n}{h_n}|
&=\sum_{n\ge1}|\iprod{Rh_n}{RU^\star h_n}|\le
\n R\n_{\calL_2(H)}\n RU^\star\n_{\calL_2(H)}\le \tr(|T|).
\end{align*}
Moreover,
$
\sum_{n\ge1}\iprod{Th_n}{h_n}
= \iprod{R}{RU^\star}_{\calL_2(H)}
$ is independent of the chosen basis by
Proposition \ref{prop:HS-inner}. This proves the theorem.
\end{proof}

\begin{definition}[Trace, of a trace class operator]
The {\em trace}\index{trace} of a trace class operator $T\in \calL(H)$ is defined by
$$\tr(T) := \sum_{n\ge 1} \iprod{Th_n}{h_n},$$
where $(h_n)_{n\ge 1}$ is any orthonormal basis of $H$.
\end{definition}

By Theorem \ref{thm:tc}, the trace is well defined.

\begin{example}[Finite-rank operators, continued]\label{ex:trace-finiterank}
The trace of a rank-one operator $g\,\bar\otimes\, h$ is equal to $\iprod{g}{h}$. To see this, without loss of generality we may assume that $h\ne 0$. Choose an orthonormal basis $(h_n')_{n\ge 1}$ such that $h_1' = h/\n h\n$; then
$$ \tr(g\,\bar\otimes\,h) = \sum_{n\ge 1} \iprod{g}{h_n'}\iprod{h_n'}{ h} = \iprod{g}{h}. $$
Hence, by linearity, the trace of a finite-rank operator $T = \sum_{n=1}^N g_n\,\bar\otimes\, h_n$
is given by
$$ \tr (T) = \sum_{n=1}^N \tr(g_n\,\bar\otimes\, h_n) =  \sum_{n=1}^N \iprod{g_n}{h_n}.$$

If $P$ is a (not necessarily orthogonal) projection onto an $N$-dimensional subspace, then $\tr (P) = N$.
Indeed, choose an orthonormal basis $g_1,\dots,g_N$ of the range of $P$, and set $h_n := P^\star g_n$. Then $P=\sum_{n=1}^N g_n\bar\otimes h_n$ and
$$
\tr(P) = \sum_{n=1}^N \iprod{g_n}{h_n}
= \sum_{n=1}^N \iprod{g_n}{P^\star g_n}
= \sum_{n=1}^N \iprod{P g_n}{g_n}
= \sum_{n=1}^N \iprod{g_n}{g_n}
= N.
$$
\end{example}

More interesting examples will be given in Section \ref{sec:trace-formulas}.

We next prove that the set $\calL_1(H)$\index{$L$@$\calL_1(H)$} of all trace class operators on $H$ is a vector space and, endowed with norm
$$\n T\n_{\calL_1(H)} := \tr(|T|),$$ a Banach space.
We first prove that $\calL_1(H)$ is a vector space. It is evident that if $T$ is a trace class operator, then so is $cT$ for all $c\in\C$ and
$\tr(|cT|) = |c|\tr(|T|)$. Additivity is less trivial and is based on a characterisation
of trace class operators which we prove first. The crucial ingredient is the following lemma.

\begin{lemma}\label{lem:tc}
Let $T\in \calL(H)$ be compact and let $(\mu_n)_{n\ge 1}$ be its singular value sequence. Then
$$ \sum_{j\ge 1} \mu_j = \sup_{g,h} \Bigl|\sum_{j} \iprod{Tg_j}{h_j}\Bigr| = \sup_{g,h} \sum_{j} |\iprod{Tg_j}{h_j}|,$$
where the suprema are taken over all integers $k\ge 1$ and all finite orthonormal sequences $g=(g_j)_{j=1}^k$ and $h=(h_j)_{j=1}^k$ of length $k$ in $H$.
\end{lemma}

Here we allow the possibility that all three expressions are infinite.

\begin{proof} Without loss of generality, we assume that $T\neq0$.
Choose an orthonormal sequence $(e_j)_{j\ge1}$ of eigenvectors of $|T|$
corresponding to the strictly positive eigenvalues $(\mu_j)_{j\ge1}$, repeated
according to multiplicity, so that $|T|e_j=\mu_j e_j$. Consider a polar
decomposition
$T=U|T|$, with $U$ a partial isometry which is isometric from
$\overline{\ran(|T|)}$ onto $\overline{\ran(T)}$.
Then,
\begin{align*}\sum_{j=1}^n \mu_j = \sum_{j=1}^n \iprod{|T|e_j}{e_j} = \sum_{j=1}^n \iprod{Te_j}{U e_j} = \Bigl|\sum_{j=1}^n \iprod{Te_j}{Ue_j}\Bigr|,
\end{align*}
using the fact that $\iprod{x}{y} = \iprod{Ux}{Uy}$ for all $x,y\in  \ov{\ran(|T|)}$
and that all $\mu_j$ are positive. For the same reason, $(Ue_j)_{j=1}^n$ is an orthonormal sequence. This proves the inequalities
$$ \sum_{j\ge 1} \mu_j \le \sup_{g,h} \Bigl|\sum_{j} \iprod{Tg_j}{h_j}\Bigr| \le \sup_{g,h} \sum_{j} |\iprod{Tg_j}{h_j}|.
$$
If $ \sum_{j\ge 1} \mu_j =\infty$, this immediately gives the equalities of the lemma.

Suppose next that $\sum_{k\ge1}\mu_k<\infty$.
Put $R:=|T|^{1/2}$. As in the proof of Theorem \ref{thm:tc}, $R$ and $RU^\star$ are Hilbert--Schmidt, with $\n RU^\star\n_{\calL_2(H)} \le \n R\n_{\calL_2(H)}$ and
$$
\n R\n_{\calL_2(H)}^2 =\tr(|T|)=\sum_{k\ge1}\mu_k.
$$
For any two orthonormal sequences $(g_j)_{j=1}^n$ and
$(h_j)_{j=1}^n$ in $H$, the Cauchy--Schwarz inequality gives
\begin{align*}
\sum_{j=1}^n|\iprod{Tg_j}{h_j}|
= \sum_{j=1}^n|\iprod{Rg_j}{RU^\star h_j}|
& \le
\Bigl(\sum_{j=1}^n\n Rg_j\n^2\Bigr)^{1/2}
\Bigl(\sum_{j=1}^n\n RU^\star h_j\n^2\Bigr)^{1/2}
\\* & \le
\n R\n_{\calL_2(H)}\n RU^\star\n_{\calL_2(H)}
\le \sum_{k\ge1}\mu_k.
\end{align*}
It follows that
\[
\sup_{g,h} \sum_{j} |\iprod{Tg_j}{h_j}| \le \sum_{k\ge1}\mu_k.
\qedhere\]
\end{proof}

\begin{theorem}[Trace class operators]\label{thm:tc-char}
For a bounded operator $T\in \calL(H)$ the following assertions are equivalent:
\begin{enumerate}[label={\rm(\arabic*)}, leftmargin=*]
\item\label{it:tc-char1} $T$ is a trace class operator;
\item\label{it:tc-char2} we have $$\sup_{g,h;k} \,\Big|\sum_{j= 1}^k\iprod{Tg_j}{h_j}\Big| <\infty,$$
the supremum being taken over all integers $k\ge1$ and all orthonormal families $g=(g_j)_{j=1}^k$ and $h=(h_j)_{j=1}^k$ in $H$;
\item\label{it:tc-char3} we have
$$\sup_{g,h}\, \sum_{j\ge 1}|\iprod{Tg_j}{h_j}| <\infty,$$
the supremum being taken over all finite or infinite orthonormal sequences
$g=(g_j)_{j\ge 1}$ and $h=(h_j)_{j\ge 1}$ in $H$ having the same number of terms.
\end{enumerate}
In this situation, the supremum in \ref{it:tc-char3} is in fact a maximum, and we have
$$\n T\n_{\calL_1(H)} =  \sup_{g,h} \,\Big|\sum_{j\ge 1}\iprod{Tg_j}{h_j}\Big|
= \sup_{g,h} \,\sum_{j\ge 1}|\iprod{Tg_j}{h_j}|,
$$
the supremum being taken over all (finite or infinite) orthonormal sequences $g=(g_j)_{j\ge 1}$ and $h=(h_j)_{j\ge 1}$ in $H$.
\end{theorem}
\begin{proof} The implication \ref{it:tc-char1}$\Rightarrow$\ref{it:tc-char3} follows from Lemma \ref{lem:tc}, which also gives the equalities in the final assertion of the theorem once the other equivalences have been proved. The implication \ref{it:tc-char3}$\Rightarrow$\ref{it:tc-char2} is trivial. To prove the implication \ref{it:tc-char2}$\Rightarrow$\ref{it:tc-char1} we argue as follows.
Denote the supremum in \ref{it:tc-char2} by $M$.
Let $T=U|T|$ be the polar decomposition of $T$, where $U$ is a partial isometry from
$\ov{\Ran(|T|)}$ onto $\ov{\Ran(T)}$. Choose an orthonormal basis $(e_n)_{n\ge1}$ of $\ov{\Ran(|T|)}$.
Then $(Ue_n)_{n\ge1}$ is an orthonormal sequence in $H$, and for every $N\ge1$ we have
$$
\sum_{n=1}^N \bigiprod{|T|e_n}{e_n}
=\sum_{n=1}^N \bigiprod{U|T|e_n}{Ue_n}
=\sum_{n=1}^N \iprod{Te_n}{Ue_n}.
$$
Since each summand on the left is real and nonnegative, the right-hand side equals its absolute value, so by the
definition of $M$,
$$
\sum_{n=1}^N \bigiprod{|T|e_n}{e_n}
=\Bigl|\sum_{n=1}^N \iprod{Te_n}{Ue_n}\Bigr|
\le M.
$$
Letting $N\to\infty$ and using monotone convergence, we obtain $$\sum_{n\ge1}\iprod{|T|e_n}{e_n}\le M<\infty.$$
Extending $(e_n)_{n\ge1}$ to an orthonormal basis of $H$ by adding an orthonormal basis of $\ker(|T|)$, we see that $\tr(|T|)<\infty$, that is, $T$ is of trace class.

To see that the supremum in \ref{it:tc-char3} is in fact a maximum, consider the sequences given by the singular value decomposition of Theorem \ref{thm:tc-ell1}.
\end{proof}

We record the following immediate consequence:

\begin{proposition}
 An operator $T\in\calL(H)$ is a trace class operator if and only if $T^\star$ is a trace class operator, in which case
$\tr(T^\star) = \ov{\tr(T)}$ and $\n T^\star\n_{\calL_1(H)} = \n T\n_{\calL_1(H)}$.
\end{proposition}

The trace class condition is expressed in terms of the summability of the
singular value sequence. As a first application of Theorem \ref{thm:tc-char} we show that if an operator is of trace class, then its eigenvalue sequence is absolutely summable:

 \begin{proposition}\label{prop:Schur}
 For any trace class operator $T\in\calL(H)$, let
$(\la_n)_{n\ge1}$ list its nonzero eigenvalues, repeated according
to algebraic multiplicity, with zeros appended if the list is finite.
Then
$$\sum_{n\ge1}|\la_n|\le\n T\n_{\calL_1(H)}.$$
 \end{proposition}
 \begin{proof}
We prove the proposition in two steps.

\smallskip
{\em Step 1} -- In this step we let $T$ be any linear operator acting on a $d$-dimensional Hilbert space $H$, with eigenvalue sequence $(\lambda_j)_{j=1}^d$ repeated according to algebraic multiplicities.
Our aim is to prove that there exists an orthonormal basis
$(h_j)_{j=1}^d$ of $H$ such that
$\iprod{Th_j}{h_j}=\lambda_j$ for all $j=1,\ldots,d$.

As a first step we prove that there exists an orthonormal basis
$(h_j)_{j=1}^d$ of $H$ such that the matrix representation
$(t_{ij})_{i,j=1}^d$ of $T$ with respect to this basis is upper triangular,
that is, it satisfies $t_{ij}=0$ whenever $i>j$.
We prove this by induction on the dimension $d$. For $d=1$ there is nothing to be proved.

Assume now that the claim has been proved for all dimensions less than $d$,
and suppose that $\dim H=d$. Let $\lambda$ be an eigenvalue of $T$.
Then $\overline\lambda$ is an eigenvalue of $T^\star$, so we may choose
a norm-one vector $h_d\in\ker(T^\star-\overline\lambda)$.
Put $G:=\{h_d\}^\perp$.
Then $G$ is invariant under $T$. Indeed, if $x\in G$, then
$$\iprod{Tx}{h_d}
   =\iprod{x}{T^\star h_d}
   = \iprod{x}{\overline\lambda h_d}
   = 0.
$$
By the induction hypothesis applied to $T|_G$, there is an orthonormal
basis $h_1,\ldots,h_{d-1}$ of $G$ such that the matrix of $T|_G$ in this
basis is upper triangular. Since $G$ is invariant under $T$, we have
$$\iprod{Th_i}{h_d}=0,\quad i=1,\ldots,d-1.$$
Moreover,
$$\iprod{Th_d}{h_d} = \iprod{h_d}{T^\star h_d} =\lambda.$$
Thus, with respect to the orthonormal basis
$h_1,\ldots,h_{d-1},h_d$, the matrix of $T$ is upper triangular. This
proves the claim.

With the orthonormal basis $(h_j)_{j=1}^d$ as in the claim, we have
$$ \det(\la-T) = (\la - \iprod{Th_1}{h_1})\dots (\la - \iprod{Th_d}{h_d}).$$
But this implies that $(\iprod{Th_j}{h_j})_{j=1}^d$ is the eigenvalue sequence of $T$ repeated according to algebraic multiplicities.

\smallskip
{\em Step 2} -- Now let $H$ and $T$ be as in the statement of the proposition. Since $T$ is compact, its nonzero spectrum consists of at most
countably many eigenvalues of finite algebraic multiplicity.
Let $(\nu_j)_{j\ge 1}$ be the distinct nonzero eigenvalues of $T$, and let
$m_j$ denote the algebraic multiplicity of $\nu_j$.
Then $m_j = \dim(X_j)$ is the dimension of the range $X_j = \Ran(P_{\nu_j})$, where $P_{\nu_j}$ is the spectral projection associated with $\nu_j$ (see Definition \ref{def:alg-multipl}).

For each admissible $N\ge 1$, let
$H_N\subseteq H$ be the finite-dimensional $T$-invariant subspace spanned by $X_1,\dots,X_N$. The identities $P_{\nu_j}P_{\nu_k}=0$ for $j\ne k$ show that
$H_N=X_1\oplus\cdots\oplus X_N$ is an algebraic direct sum.
Moreover, $(T-\nu_j I)|_{X_j}$ is nilpotent, so $T|_{H_N}$
has precisely the eigenvalues $\nu_1,\ldots,\nu_N$, with
algebraic multiplicities $m_1,\ldots,m_N$.
Applying Step 1 to $T|_{H_N}$, we obtain an orthonormal basis
$(h_n)_{n=1}^{d_N}$ of $H_N$
such that
$$
\sum_{j=1}^N m_j|\nu_j|
 = \sum_{n=1}^{d_N} |\iprod{Th_n}{h_n}|,
$$
where $d_N:=\dim H_N=\sum_{j=1}^N m_j$.
By Theorem \ref{thm:tc-char},
$$
\sum_{n=1}^{d_N} |\iprod{Th_n}{h_n}|
 \le \n T\n_{\calL_1(H)}.
$$
The proof is completed by letting $N\to\infty$ if there are infinitely many distinct nonzero eigenvalues; if there are only finitely many, we take them all, and if there are none, the assertion is immediate.
\end{proof}

We continue with some results about the structure of the set of trace class operators.

\begin{theorem}[Sums]\label{thm:sum-of-TC}
 If $T_1$ and $T_2$ are trace class operators on $H$, then so is their sum $T_1+T_2$, and we have
$\tr(T_1+T_2) = \tr(T_1)+\tr(T_2)$ and
 $$\n T_1+T_2\n_{\calL_1(H)} \le\n T_1\n_{\calL_1(H)}+\n T_2\n_{\calL_1(H)}.$$
\end{theorem}
\begin{proof}
By Theorem \ref{thm:tc-char},
\begin{align*}
\ & \sup_{g,h;N} \Bigl|\sum_{n=1}^N \iprod{(T_1+T_2)g_n}{h_n}\Bigr|
\\ & \qquad \le \sup_{g,h;N} \Bigl|\sum_{n=1}^N \iprod{T_1g_n}{h_n}\Bigr|
+ \sup_{g,h;N} \Bigl|\sum_{n=1}^N \iprod{T_2g_n}{h_n}\Bigr|
 =
\|T_1\|_{\calL_1(H)}+\|T_2\|_{\calL_1(H)}.
\end{align*}
Here the suprema are taken over all $N\ge1$ and all orthonormal families
$g_1,\ldots,g_N$ and $h_1,\ldots,h_N$ in $H$. Applying
Theorem \ref{thm:tc-char} once more shows that $T_1+T_2$ is of trace class and
that $\|T_1+T_2\|_{\calL_1(H)} \le\|T_1\|_{\calL_1(H)}+\|T_2\|_{\calL_1(H)}$.
Additivity of the traces is clear.
\end{proof}

\begin{theorem}[Completeness]\label{thm:tc-Banach} The normed space $\calL_1(H)$ is a Banach space, and the finite-rank operators
are dense in this space.
\end{theorem}
\begin{proof}
Density of the finite-rank operators follows from Theorem \ref{thm:tc-ell1}.

Suppose $(T_n)_{n\ge 1}$ is a Cauchy sequence in $\calL_1(H)$. Then $(T_n)_{n\ge 1}$ is a Cauchy sequence in $\calL(H)$.
Let $T\in \calL(H)$ be its limit. Since each $T_n$ is compact, so is $T$. Moreover, for all
orthonormal sequences $(g_j)_{j\ge 1}$ and $(h_j)_{j\ge 1}$ in $H$ and all $k\ge 1$,
$$ \sum_{j=1}^k |\iprod{Tg_j}{h_j}| = \limn \sum_{j=1}^k |\iprod{T_n g_j}{h_j}| \le \sup_{n\ge 1} \n T_n  \n_{\calL_1(H)} < \infty.$$
Therefore Theorem \ref{thm:tc-char} implies that $T$ is a trace class operator. Also,
$$ \sum_{j=1}^k |\iprod{(T_n-T)g_j}{h_j}| = \lim_{m\to\infty}  \sum_{j=1}^k |\iprod{(T_n-T_m )g_j}{h_j}| \le \limsup_{m\to\infty}\n T_n - T_m\n_{\calL_1(H)}.$$
Again by Theorem \ref{thm:tc-char}, this implies that
$$ \n T_n - T\n_{\calL_1(H)}\le \limsup_{m\to\infty}\n T_n - T_m\n_{\calL_1(H)}.$$
Since the latter tends to $0$ as $n\to\infty$, it follows that $\limn T_n = T$ in $\calL_1(H)$.
This proves that $\calL_1(H)$ is complete.
\end{proof}

 Trace class operators have the following ideal property:

 \begin{theorem}[Ideal property]\label{thm:ideal-tc}\index{ideal property!of trace class operators}
If $T$ is a trace class operator and if $S$ and $U$ are bounded, then $STU$ is a trace class operator and
$$
\|STU\|_{\calL_1(H)} \le \|S\|\,\|T\|_{\calL_1(H)}\,\|U\|.
$$
\end{theorem}

\begin{proof}
By Theorem \ref{thm:tc-ell1}, we may write
$$T=\sum_{n\ge 1}\mu_n\, g_n\,\bar\otimes\, h_n,$$
where $(g_n)_{n\ge 1}$ and $(h_n)_{n\ge 1}$ are orthonormal sequences in $H$, with
$\mu_n\ge 0$ for all $n\ge 1$ and $\sum_{n\ge 1}\mu_n=\|T\|_{\calL_1(H)}$.

For each $n\ge 1$ we have
$S(g_n\,\bar\otimes\, h_n)U = Sg_n\,\bar\otimes\, U^\star h_n$.
Since rank-one operators satisfy
$\|x\,\bar\otimes\,y\|_{\calL_1(H)}=\|x\|\,\|y\|$, we have
$$
\sum_{n\ge 1}\mu_n
\|Sg_n\,\bar\otimes\,U^\star h_n\|_{\calL_1(H)}
\le
\|S\|\,\|U\|\sum_{n\ge1}\mu_n<\infty.
$$
Hence $\sum_{n\ge1}\mu_n\,Sg_n\,\bar\otimes\,U^\star h_n $
converges in $\calL_1(H)$, and therefore also in $\calL(H)$.
Its partial sums are $ST_NU$, which converge in operator norm to $STU$.
Thus the trace-norm sum is $STU$, so $STU$ is of trace class, and
\[
\|STU\|_{\calL_1(H)}
\le \sum_{n\ge 1}\mu_n \|Sg_n\|\,\|U^\star h_n\|
\le \|S\|\,\|U\| \sum_{n\ge 1}\mu_n
= \|S\|\,\|T\|_{\calL_1(H)}\,\|U\|.
\]
\end{proof}

We continue with a proposition describing the relationship between trace class operators and Hilbert--Schmidt operators.

\begin{proposition}\label{prop:tr-HSHS}
 A bounded operator on $H$ is a trace class operator if and only if it is the product of two Hilbert--Schmidt operators. If $T = S_2 S_1$
 is such a decomposition, then
 $$\n T\n_{\calL_1(H)} \le \n S_1\n_{\calL_2(H)}\n S_2\n_{\calL_2(H)}.$$
\end{proposition}
\begin{proof}
`If': \ If $S_1$ and $S_2$ are Hilbert--Schmidt and $T:= S_2S_1$,
then for all orthonormal sequences $(g_j)_{j=1}^n$ and $(h_j)_{j=1}^n$ in $H$ and all $n\ge 1$,
\begin{align*}
\sum_{j=1}^n |\iprod{Tg_j}{h_j}|
& \le  \sum_{j=1}^n \n S_1 g_j\n \n S_2^\star h_j\n \le \Bigl( \sum_{j=1}^n \n S_1 g_j\n^2\Bigr)^{1/2}  \Bigl(  \sum_{j=1}^n \n S_2^\star h_j\n ^2\Bigr)^{1/2}
\\ & \le \n S_1\n_{\calL_2(H)} \n S_2^\star \n_{\calL_2(H)}
= \n S_1\n_{\calL_2(H)} \n S_2 \n_{\calL_2(H)}.
\end{align*}
Letting $n\to\infty$ and applying Theorem \ref{thm:tc-char}, we conclude that $T$ is a trace class operator and satisfies the inequality $\n T\n_{\calL_1(H)} \le \n S_1\n_{\calL_2(H)}\n S_2\n_{\calL_2(H)}.$

\smallskip
`Only if': \ Using a polar decomposition $T = U|T|$ with $U$ a partial isometry, take $S_1 = |T|^{1/2}$ and $S_2 = U|T|^{1/2}$\!.
Since $|T|$ is a trace class operator, $|T|^{1/2}$ is a Hilbert--Schmidt operator, and hence so are $S_1$ and $S_2$.
\end{proof}

 \begin{proposition}\label{prop:trace-commute}
  Let $S,T\in \calL(H)$, and suppose that at least one of the following three conditions holds:
   \begin{enumerate}[label={\rm(\arabic*)}, leftmargin=*]
   \item\label{it:trace-commute1a} $T\in \calL_1(H)$;
   \item\label{it:trace-commute1b} $S, T\in \calL_2(H)$;
   \item\label{it:trace-commute2} $ST$ and $TS$ have finite rank.
  \end{enumerate}
 Then $ST$ and $TS$ are trace class operators and $\tr(ST) = \tr(TS)$.
 \end{proposition}

\begin{proof}
That $ST$ and $TS$ are trace class operators follows from Theorem \ref{thm:ideal-tc} (for \ref{it:trace-commute1a}), the preceding proposition
(for \ref{it:trace-commute1b}),
and the fact that finite-rank operators are of trace class (for \ref{it:trace-commute2}).

We first prove the trace identity in case
\ref{it:trace-commute2}. The spaces
$E:=\Ran(ST)$ and $F:=\Ran(TS)$ are finite-dimensional, and
$S(F)\subseteq E$ and $T(E)\subseteq F$.
Let $S_0:F\to E$ and $T_0:E\to F$ be the corresponding
restrictions. Then
$$
S_0T_0=ST|_E,
\qquad
T_0S_0=TS|_F.
$$
Extending orthonormal bases of $E$ and $F$ to orthonormal bases
of $H$, and using the range inclusions, we obtain
$$
\tr(ST)=\tr_E(S_0T_0)=\tr_F(T_0S_0)=\tr(TS),
$$
where the middle equality is the elementary trace identity
for rectangular matrices.

We next prove the trace identity in case \ref{it:trace-commute1a}.
Choose finite-rank operators $T_n$ such that $T_n\to T$ in $\calL_1(H)$. By the case just proved, $\tr(ST_n)=\tr(T_nS)$ for all $n\geq1$. Moreover, by the ideal property,
$\n S(T_n-T)\n_{\calL_1(H)}
\leq \n S\n\,\n T_n-T\n_{\calL_1(H)}$ and $\n (T_n-T)S\n_{\calL_1(H)}
\leq \n S\n\,\n T_n-T\n_{\calL_1(H)}$.
Thus we see that $ST_n\to ST$ and $T_nS\to TS$ in $\calL_1(H)$. Since the trace is continuous with respect to the trace norm, we obtain
$$
\tr(ST)=\limn\tr(ST_n)=\limn\tr(T_nS)=\tr(TS).
$$

The proof of \ref{it:trace-commute1b} proceeds in the same way, choosing
finite-rank operators $T_n\to T$ in $\calL_2(H)$ and using the bound of
Proposition \ref{prop:tr-HSHS}.
\end{proof}

As a consequence of these two propositions, the inner product of $\calL_2(H)$ can be written in
terms of the trace:
$$
\iprod{T_1}{T_2}_{\calL_2(H)}
=
\tr(T_2^\star T_1)
=
\tr(T_1T_2^\star).
$$
Indeed, both products are trace class by Proposition
\ref{prop:tr-HSHS}, and the identities follow by computing the trace
with respect to an orthonormal basis and using Proposition
\ref{prop:trace-commute}.

\section{Trace Duality}\label{sec:trace-duality}

We have already noted that the space $\calL_2(H)$ of Hilbert--Schmidt operators on $H$ is a Hilbert space with respect to the inner product
given by
$$ \iprod{T_1}{T_2} = \tr(T_1T_2^\star),\quad T_1,T_2\in\calL_2(H).$$
The next theorem establishes that
$\calL_1(H)$ can be identified isometrically with the dual of
$\mathscr K(H)$, the closed subspace of $\calL(H)$ consisting of all compact operators on $H$,
and $\calL(H)$ with the dual of $\calL_1(H)$.

\begin{theorem}[Trace duality]\label{thm:traceduality}\index{theorem!trace duality}
By trace duality we have isometric isomorphisms $$ (\mathscr{K}(H))^* \simeq \calL_1(H) \ \ \hbox{ and } \ \ (\calL_1(H))^* \simeq \calL(H).$$
More precisely, the following results hold:
\begin{enumerate}[label={\rm(\arabic*)}, leftmargin=*]
 \item\label{it:traceduality1}  for every $T\in \calL_1(H)$ the mapping $\phi_T:\mathscr{K}(H)\to \C$ given by
 $\phi_T(S):= \tr(ST)$ is linear and bounded
and satisfies $$\n \phi_T\n_{(\mathscr{K}(H))^*} = \n T\n_{\calL_1(H)};$$ conversely, for every $\phi\in (\mathscr{K}(H))^*$ there is a unique
$T\in \calL_1(H)$ such that $\phi=\phi_T$;
 \item\label{it:traceduality2} for every $T\in \calL(H)$ the mapping $\psi_T: \calL_1(H)\to \C$ given by $\psi_T(S):= \tr(ST)$ is linear and bounded
and satisfies $$\n \psi_T\n_{(\calL_1(H))^*} = \n T\n_{\calL(H)};$$   conversely, for every $\psi\in (\calL_1(H))^*$ there is a unique
$T\in \calL(H)$ such that $\psi=\psi_T$.
\end{enumerate}
\end{theorem}

The point of working with $\tr(ST)$ rather than $\tr(ST^\star)$ is that
the resulting identifications are linear rather than conjugate-linear.

\begin{proof}
Linearity of $\phi_T$ and $\psi_T$ is clear. The ideal property
and the trace estimate give
$$
\n\phi_T\n_{(\mathscr K(H))^*}\le\n T\n_{\calL_1(H)},
\qquad
\n\psi_T\n_{(\calL_1(H))^*}\le\n T\n.
$$

Let $\Phi$ be a bounded functional on either $\mathscr K(H)$ or
$\calL_1(H)$. Since
$$
\n g\,\bar\otimes\,h\n
=
\n g\,\bar\otimes\,h\n_{\calL_1(H)}
=
\n g\n\,\n h\n,
$$
the sesquilinear form
$b(g,h):=\Phi(g\,\bar\otimes\,h)$ satisfies
$|b(g,h)|\le\n\Phi\n\,\n g\n\,\n h\n$.
By Proposition \ref{prop:bilin-T}, there exists a unique
$T\in\calL(H)$ with $\n T\n\le\n\Phi\n$ such that
$$
\Phi(g\,\bar\otimes\,h)
=
\iprod{Tg}{h}
=
\tr((g\,\bar\otimes\,h)T),
\quad g,h\in H.
$$

\smallskip
\ref{it:traceduality1}:\
Let $\phi\in(\mathscr K(H))^*$, and let $T$ be the operator
constructed above. For finite orthonormal families
$(g_j)_{j=1}^k$ and $(h_j)_{j=1}^k$, the operator
$F:=\sum_{j=1}^k g_j\,\bar\otimes\,h_j$
has finite rank and satisfies $\n F\n\le1$. Hence
$$
\Bigl|\sum_{j=1}^k\iprod{Tg_j}{h_j}\Bigr|
=
|\phi(F)|
\le\n\phi\n.
$$
By Theorem \ref{thm:tc-char}, $T$ is of trace class and
$\n T\n_{\calL_1(H)}\le\n\phi\n$.

The identity on rank-one operators and linearity give
$\phi(F)=\tr(FT)$ for every finite-rank operator $F$.
Since these operators are dense in $\mathscr K(H)$ by
Proposition \ref{prop:fr-dense-compH}, we obtain $\phi=\phi_T$.
The identity on rank-one operators also gives uniqueness of $T$.
Combining the two norm bounds proves
$\n\phi_T\n_{(\mathscr K(H))^*}=\n T\n_{\calL_1(H)}$.

\smallskip
\ref{it:traceduality2}:\
Let $\psi\in(\calL_1(H))^*$, and let $T$ be the operator
constructed above. Then $\n T\n\le\n\psi\n$, and
$\psi(F)=\tr(FT)$ for every finite-rank operator $F$.
By Theorem \ref{thm:tc-Banach}, these operators are dense in
$\calL_1(H)$, so $\psi=\psi_T$.
Again uniqueness follows from the identity on rank-one operators,
and the two norm bounds give
$\n\psi_T\n_{(\calL_1(H))^*}=\n T\n$.
\end{proof}

\section{The Partial Trace}\label{subsec:partialtrace}

The partial trace is the noncommutative analogue of integrating out one
variable. For a trace class operator $T$ on $H\otimes K$, its partial trace
over $K$ is the trace class operator $\tr_K(T)$ on $H$ characterised by
$$
   \tr(\tr_K(T)S)=\tr(T(S\otimes I)),\quad S\in\calL(H).
$$
Thus $\tr_K(T)$ is obtained by tracing out the $K$-component. In
quantum terminology, if $T$ is a state, then $\tr_K(T)$ is the reduced
state of the subsystem described by $H$.

\newcommand{\ott}{\ot}

Using the notation of Appendix \ref{sec:tensor}, we introduce the following definition.

\begin{definition}[Hilbert space tensor product]\label{def:Hilbert-tensor}
The {\em Hilbert space tensor product}\index{tensor!product, of Hilbert spaces} of the Hilbert spaces $H_1,\dots,H_N$ is the completion of the algebraic tensor product $H_1\ot\cdots\ot H_N$ with respect to the norm obtained from the inner product
$$ \Bigl(\sum_{i=1}^k g_1^{(i)}\ot\cdots\ot g_N^{(i)}\Big| \sum_{j=1}^\ell h_1^{(j)}\ot\cdots\ot h_N^{(j)}\Bigr) := \sum_{i=1}^k  \sum_{j=1}^\ell \prod_{n=1}^N \iprod{g_n^{(i)}}{h_n^{(j)}}.$$
\end{definition}

With a slight abuse of notation the Hilbert space tensor product of $H_1,\dots,H_N$ is denoted again by $H_1\ot\cdots\ot H_N$.
We leave it to the reader to check that if $H_1,\dots,H_N$ are separable, with an orthonormal basis $(h_j^{(n)})_{j\ge 1}$ for each $H_n$, then the tensors
$h_{j_1}^{(1)}\ot\cdots\ot h_{j_N}^{(N)}$ form an orthonormal basis for $H_1\ot\cdots \ot H_N$.

If $(\Om_1,\mu_1), \dots, (\Om_N,\mu_N)$ are $\sigma$-finite measure spaces, then the linear mapping
from $L^2(\Om_1,\mu_1)\ot\cdots\ot L^2(\Om_N,\mu_N)$ into $L^2(\Om_1\times\cdots\times\Om_N,\mu_1\times\cdots\times\mu_N)$ defined by
$$f_1\ot\cdots\ot f_N \mapsto \Bigl[(\om_1,\dots,\om_N) \mapsto \prod_{n=1}^N f_n(\om_n)\Bigr]$$
extends uniquely to an isometric isomorphism
\begin{align}\label{eq:L2-tensor} L^2(\Om_1,\mu_1)\ot\cdots\ot L^2(\Om_N,\mu_N) \simeq L^2(\Om_1\times\cdots\times\Om_N,\mu_1\times\cdots\times\mu_N).
\end{align}

Now let $H$ and $K$ be Hilbert spaces.
If $S\in\calL(H)$, the operator $S\ot I$, defined on the algebraic tensor
product by
$$
(S\ot I)(h\ot k):=Sh\ot k
$$
and extended linearly, has a unique bounded extension to $H\ot K$, and
$$
\|S\ot I\|=\|S\|
$$
when $K\not=\{0\}$. Indeed, every element $u$ of the algebraic tensor product can
be written as
$u=\sum_{k=1}^N h_k\otimes e_k$,
where $e_1,\dots,e_N$ are orthonormal in $K$. To obtain such
a representation, choose an orthonormal basis of the
finite-dimensional span of the second components in any
representation of $u$ as a finite sum of elementary tensors.
The definition of the tensor product inner product then gives
\begin{align*}
\n(S\ot I)u\n^2
&=\sum_{k=1}^N\n Sh_k\n^2
\le\n S\n^2\sum_{k=1}^N\n h_k\n^2
=\n S\n^2\n u\n^2.
\end{align*}
Since the algebraic tensor product is dense in $H\ot K$,
$S\ot I$ extends uniquely to a bounded operator with
$\n S\ot I\n\le\n S\n$.
If $K\ne\{0\}$, fix a unit vector $k\in K$. For every
$h\in H$ with $\n h\n\le1$ we have
$$
\n S\ot I\n\ge\n(S\ot I)(h\otimes k)\n=\n Sh\n.
$$
Taking the supremum over these $h$ proves the reverse inequality.

\begin{theorem}[Partial trace]\label{thm:partial-trace}\index{partial!trace}
 Let $H$ and $K$ be separable Hilbert spaces and let $T\in \calL_1(H\ott K)$. There exists
 a unique operator $\tr_K(T)\in \calL_1(H)$ such that for all $S\in \calL(H)$ we have
\begin{align*}
\tr(\tr_K(T)S) = \tr(T(S\ot I)).
\end{align*}
\end{theorem}

The mapping $T\mapsto \tr_K(T)$ is called the {\em partial trace} with respect to $K$ and is obtained by {\em tracing out} $K$.

\begin{proof}
For $k\in K$ let $U_k:H\to H\ot K$ be given by $$U_k h:= h\otimes k.$$
Its Hilbert space adjoint is given by $U_k^\star (h\ot k') = \iprod{k'}{k}h$.

 We claim that if  $(k_n)_{n\ge 1}$ is an orthonormal basis of $K$,
the sum
\begin{align}\label{eq:partial-trace2}  \tr_K(T) := \sum_{n\ge 1} U_{k_n}^\star T U_{k_n}
\end{align}
converges in $\calL_1(H)$ and its sum has the required properties.

By Theorem \ref{thm:ideal-tc},
each operator $U_{k_n}^\star T U_{k_n}$ is a trace class operator.
Hence by Theorem \ref{thm:tc-char}, for each $n\ge 1$ there exist orthonormal sequences $(g_j^{(n)})_{j\ge 1}$ and $(h_j^{(n)})_{j\ge 1}$ in $H$ such that
$$ \n U_{k_n}^\star T U_{k_n}\n_{\calL_1(H)} =
\sum_{j\ge 1}|\iprod{ U_{k_n}^\star T U_{k_n} g_j^{(n)}}{h_j^{(n)}}|.$$
It follows that
\begin{align*}
\sum_{n\ge 1} \n U_{k_n}^\star T U_{k_n}\n_{\calL_1(H)}
& =
\sum_{n\ge 1} \sum_{j\ge 1}|\iprod{ U_{k_n}^\star T U_{k_n} g_j^{(n)}}{h_j^{(n)}}|
 =
\sum_{n\ge 1} \sum_{j\ge 1}|\iprod{T ( g_j^{(n)}\ot k_n)}{h_j^{(n)}\ot k_n}|,
\end{align*}
and the right-hand side is finite since $T$ is of trace class
and the sequences $(g_j^{(n)}\ot k_n)_{j,n\ge 1}$ and $(h_j^{(n)}\ot k_n)_{j,n\ge 1}$ are orthonormal in $H\ott K$.

Next we check the required identity.
If $(h_m)_{m\ge 1}$ is an orthonormal basis for $H$, then
\begin{align*}
 \tr(\tr_K(T)S)
 &= \sum_{n\ge 1} \tr(U_{k_n}^\star T U_{k_n}S)\\
 &= \sum_{n\ge 1}\sum_{m\ge 1}
    \iprod{T(U_{k_n}S h_m)}{U_{k_n}h_m} \\
 &= \sum_{n\ge 1}\sum_{m\ge 1}
    \iprod{T((S\ot I)(h_m\otimes k_n))}{h_m\otimes k_n} = \tr(T(S\ot I)).
\end{align*}
In the first identity we used the fact that the series in \eqref{eq:partial-trace2} converges absolutely in
$\calL_1(H)$, and therefore multiplication on the right by $S$ and application of the trace may be performed term by term.

It remains to prove uniqueness. If $A$ is a trace class operator on $H$
such that $\tr(AS)= \tr(\tr_K(T)S)$
for all $S\in \calL(H)$, then Theorem \ref{thm:traceduality} implies that $A = \tr_K(T)$.
\end{proof}

\begin{example}[Partial trace of a rank-one projection]\label{ex:partialtrace} If $h\in H$ and $k\in K$ have norm one and $T = (h\otimes k)\,\bar\otimes\, (h\otimes k)$
is the corresponding rank-one projection on $H\ot K$, then
$\tr_K(T)$ is the rank-one projection $h\,\bar\otimes\,h$ on $H$:
$$ \tr_K(T) = h\,\bar\otimes\, h.$$
Indeed, for all $S\in \calL(H)$ we have
\begin{align*}
  \tr(\tr_K(T)S) & = \tr(T(S\ot I))  =
   \iprod{(S\ot I)(h\ot k)}{h\ot k} \\ & = \iprod{Sh\ot k}{h\ot k}
  = \iprod{Sh}{h}\iprod{k}{k} =  \iprod{Sh}{h} =  \tr((h\,\bar\otimes\,h)S).
\end{align*}
The result now follows from the uniqueness part of Theorem \ref{thm:partial-trace}.
\end{example}

In the terminology of the next chapter, the following proposition states that the partial trace of a state is again a state.

\begin{proposition}\label{prop:part-trace-props} Let $H$ and $K$ be separable Hilbert spaces and let $T\in \calL_1(H\ott K)$. Then:
\begin{enumerate}[label={\rm(\arabic*)}, leftmargin=*]
 \item\label{it:part-trace-props1} if $T$ has unit trace, then so has $\tr_K(T)$;
 \item\label{it:part-trace-props2} if $T$ is positive, then so is $\tr_K(T)$.
\end{enumerate}
\end{proposition}

\begin{proof}
For \ref{it:part-trace-props1}, we apply the defining identity of the partial trace with $S=I$. This gives
$$
\tr(\tr_K(T))=\tr(\tr_K(T)I)=\tr(T(I\ot I))=\tr(T).
$$
Thus, if $T$ has unit trace, then so has $\tr_K(T)$.

For \ref{it:part-trace-props2}, let $(k_n)_{n\geq1}$ be an orthonormal basis of $K$. By \eqref{eq:partial-trace2}, for all $h\in H$ we have
$$
\iprod{\tr_K(T)h}{h}
 =\sum_{n\geq1}\iprod{U_{k_n}^\star TU_{k_n}h}{h}
 =\sum_{n\geq1}\iprod{T(h\ot k_n)}{h\ot k_n}\geq0.
$$
Here the last inequality follows from the positivity of $T$. Hence $\tr_K(T)$ is positive.
\end{proof}

\section{Trace Formulas}\label{sec:trace-formulas}

In this final section we illustrate the preceding theory by computing traces in a number of interesting situations.

\subsection{Lidskii's Theorem}

If $T$ is a linear operator acting on $\C^d$ whose matrix representation with respect to the orthonormal basis $(e_n)_{n=1}^d$ is in Jordan normal form, then
$$ \tr(T) = \sum_{n= 1}^d \iprod{Te_n}{e_n} = \sum_{n= 1}^d \la_n ,$$
where $\la_1,\dots,\la_d$ are the eigenvalues of $T$ repeated according to their algebraic multiplicities; see Example \ref{ex:Jordan}.
By the result of Step 1 of the proof of Proposition \ref{prop:Schur}, this identity extends to arbitrary linear operators $T$ acting on $\C^d$.

If $T$ is a normal trace class operator on a Hilbert space $H$, the spectral theorem for compact normal operators allows us to select an orthonormal basis for $H$ consisting of eigenvectors in the following way. For each of the eigenspaces corresponding to the eigenvalues of
$T$ we select an orthonormal basis. These eigenspaces are mutually orthogonal and the union of these bases, after a relabelling, is an orthonormal basis $(h_n)_{n\ge 1}$ for $H$. For this basis we have
$$ \tr(T) = \sum_{n\ge 1} \iprod{Th_n}{h_n} =\sum_{n\ge 1} \la_n ,$$
where in the final sum the zero eigenvalues have been omitted and the remaining eigenvalues, repeated according to their multiplicities, have been relabelled as $(\la_n)_{n\ge 1}$.

The following deep result asserts that these formulas for the trace extend to general trace class operators:

\begin{theorem}[Lidskii]\label{thm:Lidskii}\index{theorem!Lidskii}
 For every trace class operator $T$ we have
 $$ \tr(T) =\sum_{n\ge 1} \la_n ,$$
 where  $(\la_n)_{n\ge 1}$ is the sequence of nonzero eigenvalues of $T$ repeated according to their algebraic multiplicities.
\end{theorem}

Here the sum is understood to be zero if $T$ has no nonzero
eigenvalues.
The absolute summability of the eigenvalue sequence has already been established in Proposition \ref{prop:Schur}.

We present the beautiful proof of this theorem due to Simon, which is based on the theory of Fredholm determinants.
In order to introduce these, we need some notation from multilinear algebra. We refer to Appendix \ref{sec:tensor} for the definitions.
The $n$-fold exterior product\index{exterior product} of a vector space $V$ is denoted by $\Lambda^n V$.
If $T$ is a linear operator on $V$, then the restriction
$$\Lambda^n(T):=T^{\ot n}\big|_{\Lambda^n(V)}$$
is a linear operator on $\Lambda^n(V)$.
If $S$ is another linear operator on $V$, then
$\Lambda^n (ST) = \Lambda^n (S)\Lambda^n (T)$.

If $H$ is a Hilbert space and $T$ is bounded on $H$, then $\Lambda^n (T)$ is bounded on $\Lambda^n (H)$, which is a Hilbert space in a natural way. More precisely, we view $\Lambda^n(H)$ as the range of the orthogonal antisymmetrisation projection $P_\Lambda$ from Appendix
\ref{sec:tensor}, with the inherited inner product, and use the normalised
notation
$$
v_1\wedge\cdots\wedge v_n
:=
\sqrt{n!}\,P_\Lambda(v_1\ot\cdots\ot v_n).
$$
With this convention, wedges of distinct vectors from an orthonormal basis
form an orthonormal basis of $\Lambda^n(H)$.
The adjoint of $\Lambda^n(T)$ satisfies $(\Lambda^n (T))^\star = \Lambda^n (T^\star)$.
Moreover, $\Lambda^n(|T|)=
\Lambda^n(|T|^{1/2})^\star\Lambda^n(|T|^{1/2})$ is positive.
From this we infer that $|(\Lambda^n (T))| = \Lambda^n (|T|)$. If $T$ is of trace class with singular value sequence $(\mu_j)_{j\ge 1}$, the singular values of $\Lambda^n(T)$ are $\mu_{j_1}\cdots\mu_{j_n}$ with $j_1<\dots<j_n$.
Indeed, choose an orthonormal basis $(e_j)_{j\ge1}$ of
eigenvectors of $|T|$, and write $|T|e_j=\nu_j e_j$.
The nonzero terms of $(\nu_j)$, counted with multiplicity,
are precisely the singular values $(\mu_j)$.
 Then the vectors
$e_{j_1}\wedge\cdots\wedge e_{j_n}$ with $j_1<\dots<j_n$
form an orthonormal basis of $\Lambda^nH$, and
$$
   \Lambda^n(|T|) (e_{j_1}\wedge\cdots\wedge e_{j_n})
   = \nu_{j_1}\cdots\nu_{j_n} (e_{j_1}\wedge\cdots\wedge e_{j_n}).
$$
Moreover,
$$
   |\Lambda^n(T)|=\Lambda^n(|T|),
$$
since $(\Lambda^n(T))^\star\Lambda^n(T) = \Lambda^n(T^\star T)
   = \Lambda^n(|T|^2) = \Lambda^n(|T|)^2 $.
Thus the nonzero singular values of $\Lambda^n(T)$ are precisely the products
$\mu_{j_1}\cdots\mu_{j_n}$.

It follows that $\Lambda^n(T)$ is a trace class operator and
\begin{equation}\label{eq:LaT}
\|\Lambda^n(T)\|_{\calL_1(\Lambda^n H)}
= \sum_{j_1<\dots< j_n}\mu_{j_1}\cdots\mu_{j_n}
\le \frac{1}{n!}\Big(\sum_{j\ge 1}\mu_j\Big)^n
= \frac{1}{n!}\,\|T\|_{\calL_1(H)}^{\,n}.
\end{equation}

For $n\times n$ matrices $A$ we have the following identity relating the determinant to traces and exterior products, known as
{\em MacMahon's formula}:\index{MacMahon formula}
$$ \det (I+A) =\sum_{k=0}^n \tr (\Lambda^k A),$$
with the convention $\Lambda^0(A)=I$.
A proof is sketched in Problem \ref{prob:MacMahon}. Since $\Lambda^k (V) = \{0\}$ when $k> \dim(V)$, MacMahon's formula suggests the following definition.

\begin{definition}[Fredholm determinant]\label{def:Fredholm-det}
Let $T\in \calL(H)$ be a trace class operator. The {\em Fredholm determinant}\index{Fredholm!determinant} of $I+T$ is defined as
$$ \det(I+T):= \sum_{n\in\N } \tr (\Lambda^n (T)),$$
where $\Lambda^0(H):=\C$ and $\Lambda^0(T):=I_\C$, so the
term with $n=0$ equals $1$.
\end{definition}

The sum on the right-hand side is absolutely convergent since
\begin{equation}\label{eq:abssumLa}
\begin{aligned}\sum_{n\in\N } |\tr (\Lambda^n (T))| &\le  \sum_{n\in\N} \n(\Lambda^n (T))\n_{\calL_1(\Lambda^n H)}
 \le \sum_{n\in\N } \frac1{n!}\n T\n_{\calL_1(H)}^n = \exp(\n T\n_{\calL_1(H)}).
\end{aligned}
\end{equation}
For fixed $T\in\calL_1(H)$, the function
$z\mapsto \det(I+z T)$
is entire. Indeed,
$$
   \det(I+z T)
   =
   \sum_{n=0}^{\infty}z^n\tr(\Lambda^n(T)),
$$
and by \eqref{eq:LaT} this series converges uniformly for $z$ in compact
subsets of $\mathbb C$. In particular, we note that
$$
   \det(I+z T)
   =
   1+z\tr(T)+O(z^2),
   \quad z\to0.
$$

The crucial step in the proof of Lidskii's theorem is to establish the validity of the following identity for all trace class operators $T$ and all $z\in\C$:
$$\det(I+z T) = \prod_{n\ge 1} (1+z\la_n).$$
Here $(\la_n)_{n\ge 1}$ is the sequence of nonzero eigenvalues of $T$ repeated according to algebraic multiplicities.
Notice that Proposition \ref{prop:Schur} guarantees the convergence of the infinite product. Once this formula has been obtained, Lidskii's theorem
is immediate by comparing the linear term of this product with the linear term in the definition $\det(I+z T)= \sum_{n\in\N } z^n\,\tr (\Lambda^n (T))$.

The remainder of this subsection is devoted to proving Lidskii's theorem. We fix a separable Hilbert space $H$ and start with some preliminary results.

\begin{lemma}\label{lem:detI+T}
 Let $T\in \calL(H)$ be a trace class operator and let $(\mu_n)_{n\ge 1}$ be its singular value sequence, repeated according to multiplicities. Then
 $$ |\det(I+T)| \le \prod_{n\ge 1} (1+ \mu_n).$$
\end{lemma}

\begin{proof}
It follows from \eqref{eq:LaT} and \eqref{eq:abssumLa} that
\begin{align*}
|\det(I+T)| & \le \sum_{n\in\N } |\tr(\Lambda^n(T))| \le \sum_{n\in\N } \n \Lambda^n(T)\n_{\calL_1(\Lambda^n H)}
\\* & = \sum_{n\in\N } \sum_{j_1< \dots< j_n}\mu_{j_1}\cdots\mu_{j_n} = \prod_{n\ge 1} (1+ \mu_n).\qedhere
\end{align*}
\end{proof}

\begin{lemma}\label{lem:det-exp-bd} Let $T\in \calL(H)$ be a trace class operator.
For all $\eps>0$ there exists a constant $C_\eps\ge 0$ such that for all $z\in \C$ we have $$|\det(I+z T)| \le C_\eps \exp(\eps |z|).$$
\end{lemma}
\begin{proof}
Using the inequality $|1+t|\le \exp(|t|)$,  Lemma \ref{lem:detI+T} implies, for any $N\ge
1$,
\begin{align*}
 |\det(I+z T)| & \le\prod_{n\ge 1} (1+ |z|\mu_n)
 \le \prod_{n= 1}^N (1+ |z|\mu_n) \exp\Bigl(\sum_{n\ge N+1}  |z|\mu_n\Bigr).
\end{align*}
Fix $\eps>0$. If we choose $N\ge 1$ so large that $\sum_{n\ge N+1} \mu_n <\frac12\eps$, the desired
estimate is obtained, with
\[
C_\eps := \sup_{z\in\C}\prod_{n= 1}^N (1+ |z|\mu_n) \exp\Bigl(-\frac12\eps |z|\Bigr).
\qedhere\]
\end{proof}

\begin{lemma}\label{lem:det-cont}
 The map $T \mapsto \det(I+T)$ is continuous from $\calL_1(H)$ to $\C$.
\end{lemma}
\begin{proof}
By tensoring singular value decompositions, we see that
$$
\n A_1\otimes\cdots\otimes A_n\n_{\calL_1(H^{\otimes n})}
=
\prod_{j=1}^n\n A_j\n_{\calL_1(H)}
$$
for trace class operators $A_1,\dots,A_n$.
Hence, if $\n S\n_{\calL_1(H)},\n T\n_{\calL_1(H)}\le C$,
a telescoping expansion gives
$$
\n S^{\otimes n}-T^{\otimes n}\n_{\calL_1(H^{\otimes n})}
\le
n\,C^{n-1}\n S-T\n_{\calL_1(H)}.
$$
Compressing to $\Lambda^n(H)$ and using the ideal property, we obtain
$$
\n\Lambda^n(S)-\Lambda^n(T)\n_{\calL_1(\Lambda^n H)}
\le n\,C^{n-1}\n S-T\n_{\calL_1(H)}.
$$
Thus $T\mapsto\tr(\Lambda^n(T))$ is continuous for every $n\ge1$.

By \eqref{eq:LaT}, on the ball $\n T\n_{\calL_1(H)}\le C$ we have
$$
|\tr(\Lambda^n(T))|\le\frac{C^n}{n!}.
$$
The series defining $\det(I+T)$ therefore converges uniformly
on every bounded set in $\calL_1(H)$. Since its summands are continuous,
its sum is continuous.
\end{proof}

\begin{lemma}\label{lem:det-tr}
Let $T\in\calL(H)$, and let $Y\subseteq H$ be a finite-dimensional
subspace containing $\Ran(T)$. Then
$$
\det(I+T)=\det(I_Y+T|_Y),
$$
where the determinant on the right is the ordinary determinant
on $Y$.
\end{lemma}
\begin{proof}
Put $m:=\dim(Y)$ and extend an orthonormal basis of $Y$ to one
of $H$. The induced orthonormal basis of $\Lambda^n(H)$ contains
an orthonormal basis of $\Lambda^n(Y)$.
Since $\Lambda^n(T)$ has range in $\Lambda^n(Y)$, all diagonal
terms corresponding to basis vectors outside $\Lambda^n(Y)$
vanish. Hence
$$
\tr(\Lambda^n(T))
=
\tr(\Lambda^n(T|_Y)),
\quad 0\le n\le m.
$$
For $n>m$ we have $\Lambda^n(T)=0$. Therefore
$$
\det(I+T)
=
\sum_{n=0}^m\tr(\Lambda^n(T|_Y))
=
\det(I_Y+T|_Y)
$$
by the finite-dimensional determinant formula.
\end{proof}

\begin{lemma}\label{lem:det-mult}
 If $S,T\in\calL(H)$ are trace class operators, then
 $$ \det(I+T)\det(I+S) = \det(I+(T+S+TS)).$$
\end{lemma}

\begin{proof}
First suppose that $T$ and $S$ have finite rank. Choose a
finite-dimensional subspace $Y$ containing their ranges, and
put $T_0:=T|_Y$ and $S_0:=S|_Y$.
Applying Lemma \ref{lem:det-tr} to $T$, $S$, and $T+S+TS$,
and using multiplicativity of the ordinary determinant,
\begin{align*}
\det(I+T)\det(I+S)
&=\det(I_Y+T_0)\det(I_Y+S_0)\\
&=\det((I_Y+T_0)(I_Y+S_0)) =\det(I+(T+S+TS)).
\end{align*}
This proves the finite-rank case.

For general trace class $T$ and $S$, choose finite-rank operators $T_n$ and $S_n$ such that $T_n\to T$ and $S_n\to S$ in $\calL_1(H)$. Then $T_nS_n\to TS$ in $\calL_1(H)$, since
$$
\n T_nS_n-TS\n_{\calL_1(H)}
\leq \n S_n\n\,\n T_n-T\n_{\calL_1(H)}
     +\n T\n\,\n S_n-S\n_{\calL_1(H)}.
$$
Thus $T_n+S_n+T_nS_n\to T+S+TS$ in $\calL_1(H)$. By the finite-rank case,
$$
\det(I+T_n)\det(I+S_n)=\det(I+(T_n+S_n+T_nS_n)).
$$
Letting $n\to\infty$ and using Lemma \ref{lem:det-cont}, we obtain
the desired identity.
\end{proof}

\begin{proposition}\label{prop:Fredholm-det}
If $T\in\calL_1(H)$, then $I+T$ is invertible if and only if
$\det(I+T)\neq0$.
\end{proposition}

\begin{proof}
Suppose first that $I+T$ is invertible and let
$S:=(I+T)^{-1}-I=-(I+T)^{-1}T$.
Then $S\in\calL_1(H)$ and
$(I+T)(I+S)=I$.
By Lemma~\ref{lem:det-mult},
$\det(I+T)\det(I+S) = \det((I+T)(I+S)) = \det(I) = 1$,
so $\det(I+T)\neq0$.

Suppose next that $I+T$ is not invertible.
By the Riesz--Schauder theorem, $-1$ is an eigenvalue of the
compact operator $T$, so there exists a unit vector $e\in H$
such that $Te=-e$.
Let $P:=e\,\bar\otimes\,e$ be the orthogonal projection onto
the span of $e$. Then $TP=-P$, and therefore
$(I+T)(I-P)=I+T$.
By Lemma \ref{lem:det-tr}, $\det(I-P)=0$.
Multiplicativity of the Fredholm determinant
(Lemma \ref{lem:det-mult}) now gives
\[
\det(I+T) = \det((I+T)(I-P)) = \det(I+T)\det(I-P) = 0.
\qedhere\]
\end{proof}

\begin{proposition}\label{prop:det-zero-multiplicity}
Let $T\in\calL_1(H)$, and let $-1/z_0$ be a nonzero eigenvalue
of $T$ of algebraic multiplicity $\nu$. Then
$F(z):=\det(I+zT)$ has a zero at $z_0$ of multiplicity $\nu$.
\end{proposition}

\begin{proof}
Put $\lambda_0:=-1/z_0$,
and let $P$ be the spectral projection of $T$ associated with
$\lambda_0$. Then $P$ has finite rank, $\dim\ran(P)=\nu$,
and $P$ commutes with $T$. Hence
$H=\ran(P)\oplus\ran(I-P)$
is an algebraic direct sum reducing $T$.

Since $P$ commutes with $T$ and $P(I-P)=0$,
Lemma \ref{lem:det-mult} gives, for every $z\in\C$,
\begin{align*}
        \det(I+z T)
        =
        \det\bigl((I+z TP)(I+z T(I-P))\bigr)
        =
        \det(I+z TP)\det(I+z T(I-P)).
\end{align*}

The restriction of $T$ to $\ran(I-P)$ has spectrum
$$
        \sigma(T|_{\ran(I-P)})
        \subseteq
        \sigma(T)\setminus\{\lambda_0\}.
$$
Since $\lambda_0\neq0$, it follows that $I+z_0T(I-P)$ is invertible; indeed
$I+z_0T(I-P)$ is the identity on $\Ran(P)$ and equals $I_Z + z_0 T|_Z$ on $Z:=\Ran(I-P)$.
Hence, by Proposition~\ref{prop:Fredholm-det},
$\det(I+z_0T(I-P))\neq0$.
Since $z\mapsto\det(I+z T(I-P))$
is entire, it is nonzero in a neighbourhood of $z_0$. Therefore the
multiplicity of the zero of $\det(I+z T)$ at $z_0$ is the same as the
multiplicity of the zero of $\det(I+z TP)$ at $z_0$.

By Lemma \ref{lem:det-tr}, applied to $z TP$ and
$Y=\ran(P)$, the determinant $\det(I+z TP)$ is the ordinary
determinant of $I_{\ran(P)}+z T|_{\ran(P)}$.
The operator $T|_{\ran(P)}$ has
only the eigenvalue $\lambda_0=-1/z_0$, counted with algebraic multiplicity
$\nu=\dim\ran(P)$. Hence, by the finite-dimensional determinant formula,
$$
  \det(I+z TP)
    =(1+z\lambda_0)^\nu
    =\Bigl(1-\frac{z}{z_0}\Bigr)^\nu.
$$
Thus $\det(I+z TP)$ has a zero of multiplicity $\nu$ at $z_0$, and
so does $\det(I+z T)$.
\end{proof}

We need the following standard lemma from the theory of entire functions; for
convenience we recall the short final argument needed in addition to the standard Hadamard factorisation.

\begin{lemma}\label{lem:Hadamard}
 Let $F$ be an entire function whose zeroes $z_1,z_2,\dots$ (counting multiplicities) satisfy $\sum_{n\ge 1} 1/|z_n| <\infty.$
 Assume furthermore that $F(0)=1$ and that for all $\eps>0$ there exists a constant $C_\eps\ge 0$ such that $|F(z)| \le C_\eps \exp(\eps |z|)$. Then
 $$ F(z) = \prod_{n\ge 1} \Bigl(1- \frac{z}{z_n}\Bigr), \quad z\in\C.$$
\end{lemma}
\begin{proof}
We use without proof the following standard special case of Hadamard's
factorisation theorem. If $F$ is an entire function satisfying $F(0)\not=0$ and $|F(z)| \le Ce^{|z|}$, and whose nonzero zeroes $(z_n)$ satisfy
$\sum_{n\ge 1}\frac1{|z_n|}<\infty$, then
$$
F(z)=e^{az+b}\prod_{n\ge 1}\Bigl(1-\frac z{z_n}\Bigr)
$$
for suitable $a,b\in\C$.

Since $\sum_{n\ge1}1/{|z_n|}<\infty$, the (finite or infinite) product
$$
P(z):=\prod_{n\ge1}\Bigl(1-\frac{z}{z_n}\Bigr)
$$
converges uniformly on compact subsets of $\C$ and therefore defines an entire
function. Moreover, for every $\eps>0$
there exists a constant $C_\eps$ such that
$$
|P(z)|\le C_\eps e^{\eps |z|},\quad z\in\C.
$$
Indeed, choose $N\ge1$ so large that
$\sum_{n>N}|z_n|^{-1}<\eps/2$. Using $|1+w| \le e^{|w|}$, we see that
$$\prod_{n> N}\Bigl|1-\frac{z}{z_n}\Bigr| \le e^{\eps |z|/2}.$$ The product of the first $N$ factors has polynomial growth, and hence is bounded by $C_\eps e^{\eps |z|/2}$ for a suitable
constant $C_\eps$. This gives the asserted estimate.

By Hadamard's factorisation theorem and the fact that $F(0)=1$, there exists a constant
$a\in\C$ such that
$$
F(z)=e^{az}P(z),\quad z\in\C.
$$
It remains to prove that $a=0$. Fix $\eps>0$. By the growth
estimates for $F$ and $P$, there exists $C_\eps\ge1$ such that
$$
|F(z)|,\ |P(z)|\le C_\eps e^{\eps|z|},
\quad z\in\C.
$$
Since $P(z)=e^{-az}F(z)$, for every $r>0$ such that $P$ has no
zero on the circle $|z|=r$ we obtain
$$
\log|P(re^{i\theta})|
\le
\log C_\eps+\eps r-\max\{0,\Re(ar e^{i\theta})\}.
$$
We now recall Jensen's formula: if $G$ is holomorphic in a neighbourhood of the closed
disc $\overline{B}(0;r)$, satisfies $G(0)\ne0$, and has zeros
$a_1,\dots,a_m$ in $B(0;r)$, counted with multiplicity, then
$$
\log |G(0)|
+\sum_{j=1}^m \log \frac{r}{|a_j|}
 =
 \frac1{2\pi}\int_0^{2\pi}\log |G(re^{i\theta})|\ud\theta .
$$
Jensen's formula and $P(0)=1$ therefore give
\begin{align*}
0
&\le \sum_{|z_n|<r}\log\frac{r}{|z_n|}
 =\frac1{2\pi}\int_0^{2\pi}
       \log|P(re^{i\theta})|\ud\theta\\
&\le
\log C_\eps+\eps r
-\frac1{2\pi}\int_0^{2\pi}
       \max\{0,\Re(ar e^{i\theta})\}\ud\theta
       = \log C_\eps+\eps r-\frac{|a|r}{\pi}.
\end{align*}
Dividing by $r$ and letting $r\to\infty$ through such radii gives
$|a|\le\pi\eps$. Since $\eps>0$ was arbitrary, we conclude that
$a=0$.
\end{proof}

\begin{theorem}\label{thm:Fredholm-det}
If $T\in\calL(H)$ is a trace class operator, with sequence
$(\la_n)_{n\ge 1}$ of nonzero eigenvalues repeated according to algebraic
multiplicities, then for all $z\in\C$ we have
$$ \det(I+z T) = \prod_{n\ge 1} (1+z\la_n).$$
\end{theorem}
\begin{proof}
By Lemma~\ref{lem:det-exp-bd}, the function
$F(z):=\det(I+zT)$ satisfies the growth assumption of
Lemma~\ref{lem:Hadamard}. Moreover, $F(0)=1$.

For $z\ne0$, Proposition~\ref{prop:Fredholm-det} gives
$F(z)=0$ if and only if $-1/z\in\sigma(T)\setminus\{0\}$.
Since $T$ is compact, every such spectral value is an eigenvalue.

If $T$ has no nonzero eigenvalues, then Proposition~\ref{prop:Fredholm-det}
implies that $F$ has no zeroes. Lemma~\ref{lem:Hadamard}, with an empty zero
set, gives $F\equiv1$, which is the desired empty-product formula.
Otherwise, by Propositions~\ref{prop:Fredholm-det} and
\ref{prop:det-zero-multiplicity}, the zeroes of $F$, counted with
multiplicity, are precisely the points
$ -1/{\lambda_n}$.
By Proposition~\ref{prop:Schur}, $\sum_{n\ge1}|\lambda_n|<\infty$,
so the zeroes satisfy the summability assumption in Lemma~\ref{lem:Hadamard}.
Therefore
\[
   F(z)
   =\prod_{n\ge1}\Bigl(1-\frac{z}{-1/\lambda_n}\Bigr)
   = \prod_{n\ge1}(1+z\lambda_n).
\qedhere\]
\end{proof}

\begin{proof}[Proof of Theorem \ref{thm:Lidskii}]
By Theorem \ref{thm:Fredholm-det}, for all $z\in\C$ we have
$$
\det(I+z T)=\prod_{n\ge1}(1+z\la_n).
$$
The product converges absolutely and uniformly on compact subsets of $\C$, because
$\sum_{n\ge1}|\lambda_n| \le \|T\|_{\calL_1(H)}<\infty$ by Proposition~\ref{prop:Schur}.

On the other hand, by the definition of the Fredholm determinant,
$$
\det(I+z T) = 1+z\,\tr(T)+O(z^2)
$$
as $z\to 0$.
Comparing the derivatives at $z=0$ of the two sides gives
$$
\tr(T)=\sum_{n\ge1}\la_n .
$$
The comparison of derivatives is justified by the absolute convergence
$\sum_{n\ge1}|\la_n|<\infty$, which allows termwise differentiation of the
product at $z=0$.
 \end{proof}

\subsection{Trace Formula for Integral Operators}\label{subsec:ktt}

The trace of an integral operator with continuous kernel can be computed as follows.

\begin{theorem}[Mercer]\label{thm:Mercer}\index{theorem!Mercer}
Let $\mu$ be a finite Borel measure on a compact metric space $K$.
Let $T$ be an integral operator on $L^2(K,\mu)$ of the form
$$ Tf(s) = \int_K k(s,t)f(t)\ud \mu(t)$$
with continuous kernel $k\in C(K\times K)$. Then:
\begin{enumerate}[label={\rm(\arabic*)}, leftmargin=*]
 \item\label{it:Mercer1} if $T$ is a trace class operator, then its trace is given by
\begin{align*}
\tr(T) = \int_K k(t,t)\ud \mu(t);
\end{align*}
 \item\label{it:Mercer2} if $T$ is positive, that is, if $\iprod{Tf}{f}\ge 0$ for all $f\in L^2(K,\mu)$, then $T$ is of trace class.
\end{enumerate}
\end{theorem}

If $\supp(\mu)=K$, then $T$ is positive if and only if the kernel
$k$ is {\em positive definite}\index{positive definite}, that is,
for all integers $N\ge1$, points $t_1,\dots,t_N\in K$, and
scalars $z_1,\dots,z_N\in\C$ we have
$$
\sum_{n,m=1}^N k(t_m,t_n)z_m\ov z_n\ge0.
$$
For the `only if' part, the assumption $\supp(\mu)=K$ ensures that every
ball $B(t_j;r)$ has positive measure, so we may test positivity
of $T$ on linear combinations of the normalised indicators
$\mu(B(t_j;r))^{-1}\one_{B(t_j;r)}$.
Letting $r\downarrow0$ then gives the required inequality by
continuity of $k$.

For the `if' part, choose finite Borel partitions whose meshes tend
to zero, select one point in each of the sets of the partitions, and approximate $k$
on products of these sets by its values at the corresponding pairs
of chosen points. The resulting kernels converge uniformly
to $k$ and, by positive definiteness, define positive integral
operators converging to $T$ in operator norm, so $T$ is positive.

\begin{proof}
If $\mu(K)=0$, both assertions are immediate. Assume that
$\mu(K)>0$, and put $H:=L^2(K,\mu)$, which is separable by
Remark \ref{rem:CKdenseLpK}.

For each $m\ge1$, choose a finite Borel partition of $K$ into
sets of diameter at most $1/m$, and let $\pi_m$ denote the
collection of the sets in this partition with positive measure. These sets cover
$K$ up to a $\mu$-null set. For $B\in\pi_m$ put
$$
e_B:=\mu(B)^{-1/2}\one_B,
$$
and let $P_m$ be the orthogonal projection onto the span of
$\{e_B:\,B\in\pi_m\}$.

For $f\in C(K)$, averaging over the sets in $\pi_m$ gives
$$
\n P_mf-f\n_2
\le
\mu(K)^{1/2}
\sup_{\substack{s,t\in K\\d(s,t)\le1/m}}|f(s)-f(t)|
\to 0.
$$
Since $C(K)$ is dense in $H$ and the operators $P_m$ are
contractions, it follows that $P_m\to I$ strongly on $H$.

The operator $P_mTP_m$ has finite rank, and
\begin{align*}
a_m:=\tr(P_mTP_m) =\sum_{B\in\pi_m}\iprod{Te_B}{e_B}
= \sum_{B\in\pi_m}\frac1{\mu(B)}
  \int_B\int_B k(s,t)\ud\mu(t)\ud\mu(s).
\end{align*}
Put $D:=\int_K k(s,s)\ud\mu(s)$. Comparing the integrand
$k(s,t)$ with $k(s,s)$ gives
$$
|a_m-D|
\le
\mu(K)
\sup_{\substack{s,t\in K\\d(s,t)\le1/m}}
|k(s,t)-k(s,s)|
\to0
$$
by the uniform continuity of $k$.

\smallskip
\ref{it:Mercer1}:\
Suppose that $T$ is of trace class. For a rank-one operator
$F=g\,\bar\otimes\,h$, strong convergence of $P_m$ gives
$$
P_mFP_m=(P_mg)\,\bar\otimes\,(P_mh)\to F
\quad\hbox{in }\calL_1(H),
$$
where we use
$\n u\,\bar\otimes\,v\n_{\calL_1(H)}=\n u\n\,\n v\n$.
By linearity the same holds for every finite-rank operator.
The density of these operators in $\calL_1(H)$
and the ideal property therefore imply
$P_mTP_m\to T$ in $\calL_1(H)$.
By continuity of the trace,
$$
\tr(T)=\lim_{m\to\infty}a_m
=\int_K k(s,s)\ud\mu(s).
$$

\ref{it:Mercer2}:\
Suppose that $T$ is positive. Then every $P_mTP_m$ is positive.
For any finite orthonormal family $u_1,\dots,u_N$ in $H$,
strong convergence of $P_m$ gives
\begin{align*}
\sum_{j=1}^N\iprod{Tu_j}{u_j}
=
\lim_{m\to\infty}
\sum_{j=1}^N\iprod{P_mTP_mu_j}{u_j}\le
\lim_{m\to\infty}\tr(P_mTP_m)
=D.
\end{align*}
By the definition of the trace of a positive operator, it follows
that $\tr(T)\le D<\infty$. Hence $T$ is of trace class, and
part \ref{it:Mercer1} gives the trace formula.
\end{proof}

\subsection{Trace Formula for Fredholm Operators}\label{subsec:TraceFormulas}

The following theorem gives a trace formula for the index of a Fredholm operator.

\begin{theorem}[Fedosov]\label{thm:Fedosov}\index{theorem!Fedosov}\index{trace formula!and Fredholm operators}
Let $T\in\calL(H)$ be Fredholm, and let
$S\in\calL(H)$ be such that $I-ST$ and $I-TS$ have finite rank. Then
$$\ind(T) = \tr([T,S]),$$
where $[T,S] = TS - ST$ is the commutator\index{commutator} of $T$ and $S$.
\end{theorem}
\begin{proof}
Since $T$ is Fredholm, $\Ran(T)$ is closed.
The restricted operator $T_0:=T|_{\ker(T)^\perp}:\ker(T)^\perp\to\Ran(T)$
is bijective, so it has a bounded inverse.
With respect to the orthogonal decomposition
$H=\Ran(T)\oplus\Ran(T)^\perp$,
define $S_0\in \calL(H)$ by $S_0(y+y^\perp):=T_0^{-1}y$.
Then $I-S_0T=P_{\ker(T)}$ and $I-TS_0=P_{\Ran(T)^\perp}$,
the orthogonal projections onto $\ker(T)$ and $\Ran(T)^\perp$.
Therefore
$$
\tr(I-S_0T)-\tr(I-TS_0)
=\dim\ker(T)-\dim\Ran(T)^\perp
=\ind(T).
$$

Now let $S$ be any operator satisfying the assumptions of the theorem, and put $R:=S-S_0$.
Then $RT=ST-S_0T$ and $TR=TS-TS_0$ have finite rank. By Proposition \ref{prop:trace-commute}, applied to $R$ and $T$, we have
$\tr(RT)=\tr(TR)$.
Consequently,
\begin{align*}
\tr([T,S]) & = \tr(I-ST)-\tr(I-TS)
\\ & =\tr(I-S_0T)-\tr(I-TS_0) -\tr(RT)+\tr(TR) =\ind(T).
\qedhere
\end{align*}
\end{proof}

\subsection{Trace Formula for Commutators of Toeplitz Operators}\label{subsec:Toeplitz-trace}

From Section \ref{subsec:Toeplitz} we recall
that $H^2(\mathbb{D})$ is the vector space of all holomorphic functions on $\mathbb{D}$ of the form $\sum_{n\in\N} c_n z^n$ with $\sum_{n\in\N} |c_n|^2 < \infty$. Under its identification with the closed subspace of $L^2(\mathbb{T})$ consisting of all functions whose negative Fourier coefficients vanish,
$H^2(\mathbb{D})$ is the range of the Riesz projection
$$ P: \sum_{n\in \Z}\wh f(n)e_n \mapsto \sum_{n\in\N}\wh f(n)e_n$$
on $L^2(\mathbb{T})$, where $e_n(\theta) = e^{ in\theta}$. This projection discards the terms in the sequence of Fourier coefficients $(\wh f(n))_{n\in \Z}$ of $f$ with negative indices $n = -1,-2,\dots$.

Given a function $\phi\in L^\infty(\mathbb{T})$, the {\em Toeplitz operator} with symbol $\phi$
has been defined as the bounded operator $T_\phi$ on $H^2(\mathbb{D})$ given by
$$T_\phi f := P(\phi f), \quad f\in H^2(\mathbb{D}).$$
It follows from Lemma \ref{lem:Toep-comp} that for all
$\phi,\psi\in C(\mathbb{T})$ the commutator $[T_\phi, T_\psi] = T_\phi T_\psi - T_\psi T_\phi$ is compact.
For functions $\phi,\psi\in C^2(\mathbb{T})$ we have the following stronger result.

\begin{theorem}[Helton--Howe]\label{thm:Toeplitz-trace}\index{theorem!Helton--Howe}\index{trace formula!and Toeplitz operators} For all
$\phi,\psi\in C^2(\mathbb{T})$ the commutator
$[T_\phi, T_\psi]$ is a trace class operator and
\begin{align}\label{eq:Toeplitz-trace} \tr([T_\phi, T_\psi]) = \frac1{2\pi i} \int_{-\pi}^\pi \phi(\theta)\psi'(\theta)\ud\theta.
\end{align}
\end{theorem}
\begin{proof}
For $n,m\in\Z$, the commutator
$[T_{e_n},T_{e_m}]$ has rank at most
$\min\{|m|,|n|\}$.
Indeed, for $n,p>0$, a direct calculation gives
$$
[T_{e_n},T_{e_{-p}}]e_j=
\begin{cases}
-e_{j+n-p},&\max\{0,p-n\}\le j<p,\\
0,&\text{otherwise}.
\end{cases}
$$
The $\min\{n,p\}$ nonzero images are linearly independent.
Interchanging the indices handles the other opposite-sign case.
When $nm\ge0$, the commutator $[T_{e_n},T_{e_m}]$ vanishes.

We have
$$
\|[T_{e_n},T_{e_m}]\|
\le
\|T_{e_n}T_{e_m}\|
+
\|T_{e_m}T_{e_n}\|
\le2.
$$
Since a rank-$r$ operator $R$ satisfies
$\|R\|_{\calL_1(H)}\le r\|R\|$,
we obtain
\begin{equation}\label{eq:HH}
\|[T_{e_n},T_{e_m}]\|_{\calL_1(H^2(\mathbb D))}
\le
2\min\{|m|,|n|\}.
\end{equation}
Second, for all $n,m\in \Z$ and $j\in \N$ we have
$T_{e_n}T_{e_m}e_j = \la_j^{nm} e_{n+m+j}$ with $\la_j^{nm}\in\{0,1\}$, so that
\begin{align}\label{eq:comm-Toep} \tr([T_{e_n}, T_{e_m}]) = \sum_{j\ge 0} (\la_j^{nm}-\la_j^{mn})\iprod{e_{n+m+j}}{e_{j}}.
\end{align}

\smallskip
{\em Case 1:} \  $n+m\not=0$. In that case \eqref{eq:comm-Toep} gives
\begin{align*} \tr([T_{e_n}, T_{e_m}]) = 0 = \frac1{2\pi i} \int_{-\pi}^\pi e_n(\theta)e_m'(\theta)\ud\theta.
\end{align*}

\smallskip
{\em Case 2:} \ $n+m=0$ with $n\ge 0$. In that case $\la_j^{n,-n} =0$ if $j<n$ and  $ \la_j^{n,-n} =1$ if $j\ge n$,
whereas $\la_j^{-n,n}=1$ for all $j$, and  \eqref{eq:comm-Toep} gives
\begin{align*} \tr([T_{e_n}, T_{e_m}]) = -n = -\frac{n}{2\pi} \int_{-\pi}^\pi e_n(\theta)e_{-n}(\theta)\ud\theta= \frac1{2\pi i} \int_{-\pi}^\pi e_n(\theta)e_{m}'(\theta)\ud\theta.
\end{align*} The case $n+m=0$ with $n<0$ is entirely similar.
This completes the proof of \eqref{eq:Toeplitz-trace} for $\phi = e_n$ and $\psi = e_m$.

If $\phi,\psi\in C^2(\mathbb T)$, their Fourier coefficients are
absolutely summable. Hence the corresponding Toeplitz series converge
in operator norm. Estimate \eqref{eq:HH} and the calculation below show
that the double commutator series converges absolutely in trace norm.
Since both the left- and right-hand sides of \eqref{eq:Toeplitz-trace} are linear in both $\phi$ and $\psi$,
for $\phi = \sum_{n\in \Z} a_n e_n$ and $\psi = \sum_{n\in \Z} b_n e_n$ we have
\begin{align*} [T_\phi, T_\psi] = \sum_{m,n\in\Z} a_n b_m [T_{e_n}, T_{e_m}]
\end{align*}
provided the sum converges in $\calL_1(H^2(\mathbb{D}))$. In view of \eqref{eq:HH}, this can be guaranteed if we assume that $\phi$ and $\psi$ are $C^2$\!, for then $|a_n|$ and $|b_n|$ are of order $n^{-2}$
as $n\to\infty$
and
\begin{align*}
\sum_{m,n\in\Z}
\frac{\min\{|m|,|n|\}}{(1+m^2)(1+n^2)}
\le
\sum_{m,n\in\Z}
\frac{\sqrt{|m||n|}}{(1+m^2)(1+n^2)}
=
\Bigl(\sum_{n\in\Z}\frac{\sqrt{|n|}}{1+n^2}\Bigr)^2
<\infty.
\end{align*}

Taking traces,
\begin{align*}
\tr( [T_\phi, T_\psi])
= \sum_{m,n\in\Z} a_n b_m \frac1{2\pi i} \int_{-\pi}^\pi e_n(\theta)e_m'(\theta)\ud\theta
= \frac1{2\pi i} \int_{-\pi}^\pi \phi(\theta)\psi'(\theta)\ud\theta.
\end{align*}
The Fourier series of $\psi'$ converges in $L^2(\mathbb T)$,
with coefficients $im b_m$; together with the uniform convergence
of the Fourier series of $\phi$, this justifies termwise integration.
\end{proof}

\subsection{Trace Formula for the Dirichlet Heat Semigroup}\label{subsec:Dirichlet-heat-trace}\index{trace formula!for the Dirichlet heat semigroup}

Let $D$ be a bounded open subset of $\R^d$ and let
$S_{\rm Dir}$ be the $C_0$-semigroup on $L^2(D)$ generated by the Dirichlet Laplacian $\Delta_{\rm Dir}$ associated with $D$.

If  $|\partial D| = 0$, the eigenvalue asymptotics given by Weyl's theorem give a quick proof that the operators $S_{\rm Dir}(t)$ are of trace class for $t>0$. The next theorem proves that this conclusion holds for general bounded open domains $D$, with a sharp upper bound for their traces. The proof combines the spectral theorem with the Fourier transform and scalar convexity.

\begin{theorem}[Kac]\label{thm:Kac-bound}\index{theorem!Kac}
Let $D\subseteq\R^d$ be nonempty, bounded, and open, and let
$S_{\rm Dir}$ be the $C_0$-semigroup generated by
$\Delta_{\rm Dir}$ on $L^2(D)$. Then $S_{\rm Dir}(t)$ is
of trace class for every $t>0$, and
\begin{align*}
\tr(S_{\rm Dir}(t))
\le\frac{|D|}{(4\pi t)^{d/2}},\quad t>0.
\end{align*}
\end{theorem}

\begin{proof}
By Theorem \ref{thm:spectra-Laplacian}, there is an orthonormal
basis $(h_n)_{n\ge1}$ of $L^2(D)$ with
$$-\Delta_{\rm Dir}h_n=\lambda_nh_n,$$ where
$0<\lambda_1\le\lambda_2\le\cdots$ and eigenvalues are counted
according to multiplicity. Since $h_n\in H_0^1(D)$, its zero
extension $\wt h_n$ belongs to $W^{1,2}(\R^d)$ and its weak
first derivatives are the zero extensions of those of $h_n$.
Indeed, this follows by approximating $h_n$ in $H_0^1(D)$
by functions in $C_{\rm c}^\infty(D)$ and extending these by zero.

Plancherel's theorem and the form definition of the Dirichlet
Laplacian therefore give
$$
\int_{\R^d}|\wh{\wt h_n}(\xi)|^2\ud\xi=1,
\qquad
\int_{\R^d}|\xi|^2|\wh{\wt h_n}(\xi)|^2\ud\xi
=\n\nabla h_n\n_2^2=\lambda_n.
$$
Thus $|\wh{\wt h_n}(\xi)|^2\ud\xi$ is a probability measure
on $\R^d$. Applying Jensen's inequality to the convex function
$s\mapsto e^{-ts}$ gives
$$
e^{-t\lambda_n}
\le\int_{\R^d}e^{-t|\xi|^2}
       |\wh{\wt h_n}(\xi)|^2\ud\xi.
$$

For fixed $\xi\in\R^d$, put
$e_\xi(x):=(2\pi)^{-d/2}e^{ix\cdot\xi}$ for $x\in D$.
The zero extensions $\wt h_n$ are integrable, since $D$ has finite measure,
so their Fourier transforms have the representatives
$$
\wh{\wt h_n}(\xi)
=(2\pi)^{-d/2}\int_D h_n(x)e^{-ix\cdot\xi}\ud x
=\iprod{h_n}{e_\xi},
$$
and Parseval's identity in $L^2(D)$ gives
$$
\sum_{n\ge1}|\wh{\wt h_n}(\xi)|^2
=\n e_\xi\n_2^2=\frac{|D|}{(2\pi)^d}.
$$
Summing the Jensen inequalities and applying Fubini's theorem,
we obtain
\begin{align*}
\sum_{n\ge1}e^{-t\lambda_n}
&\le\int_{\R^d}e^{-t|\xi|^2}
       \sum_{n\ge1}|\wh{\wt h_n}(\xi)|^2\ud\xi
=\frac{|D|}{(2\pi)^d}\int_{\R^d}e^{-t|\xi|^2}\ud\xi
=\frac{|D|}{(4\pi t)^{d/2}}.
\end{align*}
Since $S_{\rm Dir}(t)h_n=e^{-t\lambda_n}h_n$, the semigroup
operator is positive, and finiteness of the sum above shows that it is
of trace class. Its trace is $\sum_{n\ge1}e^{-t\lambda_n}$,
which proves the assertion.
\end{proof}

No regularity assumption on $\partial D$ is needed for this bound.
When $|\partial D|=0$, Theorem \ref{thm:heat-trace} below shows that
the upper bound is asymptotically exact as $t\downarrow0$.
In particular, the constant $(4\pi)^{-d/2}$ is optimal.

\begin{theorem}[Trace formula for the Dirichlet heat semigroup]\label{thm:heat-trace}
Let $D\subseteq\R^d$ be nonempty, bounded, and open, and assume that $|\partial D| = 0$. For all $t>0$ the operator $S_{\rm Dir}(t)$ is a trace class operator on $L^2(D)$ with
\begin{align*} \lim_{t\downarrow 0} t^{d/2}\tr(S_{\rm Dir}(t)) = \frac{|D|}{(4\pi)^{d/2}}.
\end{align*}
\end{theorem}

For the proof of this formula we need the following lemma.

\begin{lemma}\label{lem:Abel} Let $\mu$ be a Borel measure on $[0,\infty)$ whose Laplace transform\index{Laplace transform!of a measure} satisfies
 $$ \calL \mu(t):= \int_0^\infty e^{-t x}\ud \mu(x)<\infty$$
 for all $t>0$. If for some $r\ge 0$ and $a\in\R$ we have $ \lim_{x\to\infty} x^{-r} \mu([0,x]) = a,$
 then $$ \lim_{t\downarrow 0} t^r \calL\mu(t) = a \Gamma(1+r),$$
where $\Gamma(s) = \int_0^\infty x^{s-1}e^{-x}\ud x$, $s>0$, is the Euler Gamma function.
\end{lemma}
\begin{proof}
Put $\nu(x):=\mu([0,x])$. Finiteness of $\calL\mu(1)$ implies
that $\mu$ is finite on bounded intervals.
By Fubini's theorem and a change of variable,
\begin{align*}
\calL\mu(t)
=\int_0^\infty\Bigl(t\int_s^\infty e^{-tx}\ud x\Bigr)\ud\mu(s)=
t\int_0^\infty e^{-tx}\nu(x)\ud x =
t^{-r}
\int_0^\infty e^{-y}t^r\nu(y/t)\ud y.
\end{align*}
For every $y>0$, the assumed asymptotics imply
$$
t^r\nu(y/t)
=
y^r(y/t)^{-r}\nu(y/t)
\to ay^r
\quad\hbox{as }t\downarrow0.
$$
Moreover, local boundedness of $\nu$ and the same asymptotics
give $\nu(x)\le C(1+x)^r$ for all $x\ge0$.
Thus, for $0<t\le1$ and $y\ge0$ we have $0\le t^r\nu(y/t)\le C(t+y)^r\le C(1+y)^r$.
Dominated convergence now gives
\[
\lim_{t\downarrow0}t^r\calL\mu(t)
=
a\int_0^\infty e^{-y}y^r\ud y
=
a\Gamma(1+r).
\qedhere\]
\end{proof}

\begin{proof}[Proof of Theorem \ref{thm:heat-trace}]
The case $D=\emptyset$ is immediate, so we may assume that $D$ is nonempty.
By Theorem \ref{thm:spectra-Laplacian}, there is an orthonormal
basis $(h_n)_{n\ge1}$ of $L^2(D)$ such that
$-\Delta_{\rm Dir}h_n=\lambda_nh_n$ for all $n\ge 1$,
where $0<\lambda_1\le\lambda_2\le\cdots$ are the eigenvalues
repeated according to multiplicity. This implies
$$
S_{\rm Dir}(t)h_n=e^{-t\lambda_n}h_n,
\quad t\ge0,\ n\ge1.
$$
By Weyl's theorem (Theorem \ref{thm:Weyl}), there exists a constant $c>0$
such that $\lambda_n\ge c n^{2/d}$ for all sufficiently large $n$.
Hence $\sum_{n\ge1}e^{-t\lambda_n}<\infty$ for every $t>0$.
The diagonal representation therefore shows that $S_{\rm Dir}(t)$
is positive and of trace class, with
$$
\tr(S_{\rm Dir}(t))=\sum_{n\ge1}e^{-t\lambda_n}.
$$

We now apply Lemma \ref{lem:Abel} to the Borel measure $\mu = \sum_{n\ge 1} \delta_{\la_n}.$
Setting $N(x):= \max\{n\ge 1:\, \la_n\le x\}$ (with $\max\emptyset = 0$), by Weyl's theorem we have
\begin{align*}
\lim_{x\to\infty} x^{-d/2} \mu([0,x])
= \lim_{x\to\infty} x^{-d/2} N(x) = \frac{\om_d}{(2\pi)^d}|D|,
\end{align*}
where $\omega_d = \pi^{d/2}/\Gamma(1+ \frac12d) $ is the volume of the unit ball in $\R^d$\!.
Lemma  \ref{lem:Abel} allows us to conclude that
\begin{align*}
 t^{d/2} \tr(S_{\rm Dir}(t)) & =  t^{d/2} \sum_{n\ge 1} e^{-\la_n t}  =
 t^{d/2} \calL \mu(t) \to   \frac{\om_d}{(2\pi)^d}|D|\Gamma\bigl(1+\frac{d}{2}\bigr)
= \frac{|D|}{(4\pi)^{d/2}}
\end{align*}
as $t\downarrow 0$.
\end{proof}

\subsection{Euler's Identity Revisited}\label{subsec:Euler}

Consider the Dirichlet Laplacian $\Delta_{\rm Dir}$ on $L^2(0,1)$.
As shown in Example \ref{ex:ev-Dir},
$$
\sigma(-\Delta_{\rm Dir})
=
\{\pi^2n^2:\, n\ge 1\},
$$
and the functions
$f_n(t):=\sqrt{2}\sin(n\pi t)$, $n\ge 1$,
form an orthonormal basis of eigenfunctions. Consequently, the positive
compact operator
$T:=(-\Delta_{\rm Dir})^{-1}$
satisfies
$Tf_n=(\pi^2n^2)^{-1}f_n$, $n\ge 1$.
Since $(f_n)_{n\ge1}$ is an orthonormal basis, from
$Tf_n=(\pi^2n^2)^{-1}f_n$ we infer that $T$ is injective and has
no eigenvalues other than $(\pi^2n^2)^{-1}$, $n\ge1$.
By the Riesz--Schauder theorem (Theorem \ref{thm:comp-spec}),
every nonzero spectral value of $T$ is an eigenvalue.
Moreover, $0\in\sigma(T)$ since $T$ is compact and $L^2(0,1)$
is infinite-dimensional. Thus,
$$
\sigma(T)
=\{0\}\cup
\Bigl\{\frac1{\pi^2n^2}:\, n\ge1\Bigr\},
$$
where $0$ belongs to the spectrum but is not an eigenvalue.

Since $T$ is the positive
integral operator with continuous kernel $k$, Mercer's theorem gives
$$
\sum_{n=1}^\infty\frac1{\pi^2n^2}
=
\tr(T)
=
\int_0^1k(t,t)\ud t,
$$
where $k$ is Green's function for the Poisson problem on the unit
interval with Dirichlet boundary conditions.
From Section \ref{sec:Poisson-D} we recall that it is given by
\begin{align*} k(s,t) = \begin{cases}
             (1-t)s, &  s\le t, \\
             (1-s)t, &  t\le s.
            \end{cases}
\end{align*}
In view of
$$ \int_0^1 k(t,t)\ud t = \int_0^1 (1-t)t\ud t = \frac16,$$
we recover Euler's identity\index{Euler's identity $\sum_{n=1}^\infty \frac1{n^2} = \frac{\pi^2}{6}$}\index{trace formula!and Euler's identity}
$$ \sum_{n=1}^\infty \frac1{n^2} = \frac{\pi^2}{6}.$$

\begin{problems}

\item
\begin{enumerate}[\rm (a), leftmargin=*]
\item
Find a Hilbert--Schmidt operator on $\ell^2$ that is not of trace class.
\item
Find a compact operator on $\ell^2$ that is not Hilbert--Schmidt.
\end{enumerate}

\item
Show that every trace class operator $T$ is Hilbert--Schmidt and that
$$
\n T\n_{\calL_2(H)}^2
\le \n T\n\,\n T\n_{\calL_1(H)}.
$$
\noindent{\em Hint:}\ Use the singular values of $T$.

\item\label{prob:HaaseL2Linfty}
Let $T:L^2(0,1)\to L^\infty(0,1)$ be a bounded operator, and let $(h_n)_{n\ge 1}$ be an orthonormal basis for $L^2(0,1)$.
For each $j\ge1$, fix a measurable representative of the equivalence class $Th_j\in L^\infty(0,1)$; below, $Th_j(t)$ refers to this representative.
\begin{enumerate}[\rm(a), leftmargin=*]
\item\label{it:HaaseL2Linfty}
Show that for every $k\ge 1$ there exists a null set $N_k\subseteq (0,1)$ such that
for all $c=(c_j)_{j=1}^k\in \C^k$,
$$
\Big| \sum_{j=1}^k c_j\, Th_j(t)\Big| \le \|T\|\,\Big(\sum_{j=1}^k |c_j|^2\Big)^{1/2},\quad t\in (0,1)\setminus N_k.
$$
\item Deduce from part \ref{it:HaaseL2Linfty} that $ \sum_{j=1}^k |Th_j(t)|^2 \le \n T\n^2$ for all $t\in (0,1)\setminus N_k$.
\end{enumerate}
Let $i:L^\infty(0,1)\to L^2(0,1)$
be the inclusion mapping.
\begin{enumerate}[\rm(a), leftmargin=*, resume]
\item Show that $i\circ T$ is Hilbert--Schmidt on $L^2(0,1)$ and
 $\n i\circ T\n_{\calL_2(L^2(0,1))}\le \n T\n$.
\end{enumerate}

\item
Prove the following analogue of Proposition \ref{prop:Schur}:
If $T\in \calL(H)$ is a Hilbert--Schmidt operator, with sequence $(\la_n)_{n\ge 1}$ of nonzero eigenvalues repeated according to algebraic multiplicity, then
$$\sum_{n\ge 1}|\la_n|^2 \le \n T\n_{\calL_2(H)}^2.$$

\item
Let $T\in\calL(H)$ be compact and let
$\mu_1\ge\mu_2\ge\cdots\ge0$ be its singular values, extended by
zeros when necessary. If $H$ is finite-dimensional, let
$1\le n\le\dim H$; otherwise let $n\ge1$.
Show that for all such $n$ we have
$$ \mu_n  = \inf_{\substack{Y\subseteq H \\ \dim(Y) = n-1}} \sup_{\substack{\n y\n=1 \\ y \perp Y}} \n Ty\n,$$
where the infima are taken over all subspaces $Y$ of $H$ of dimension $n-1$.

\noindent
{\em Hint:}\ Use Theorem \ref{thm:minmax}.

\item
As was observed in the main text, the nonzero eigenvalues $\la_n$ and singular values $\mu_n$ of a compact normal operator, repeated according to multiplicities and ordered by decreasing absolute value, are related by $|\la_n| = \mu_n$. Show that this relation breaks down in the absence of normality, by computing the eigenvalues and singular values of the matrix $\begin{pmatrix} 1 & 1 \\ 0 & 1\end{pmatrix}$.

\item
Show that if $T\in \calL(H)$ is a trace class operator and $U\in \calL(H)$ is unitary, then $UTU^\star$ is a trace class operator and $\tr(T) = \tr(UTU^\star)$.

\item\label{prob:postraceclass}
Show that if $A\in \calL(H)$ is a positive trace class operator, then
$ A \le \tr(A)I$.

\item
Show that if $S = x\,\bar\otimes\, x$ and $T = y\,\bar\otimes\, y$ with $\n x\n = \n y\n = 1$ are two rank-one orthogonal projections, then
$$\n S-T\n^2 = 1- |\iprod{x}{y}|^2 = 1 - \tr(ST).$$

\item
Prove that if $T\in \calL(H)$ is self\-adjoint and $S\in \calL(H)$ is compact, and if the commutator $[T,S]$ is a trace class operator, then $\tr [T,S]=0$.

\noindent{\em Hint:}\ Compute the traces of $[T,S\pm S^\star]$ relative to orthonormal bases which diagonalise
$S\pm S^\star$\!.

\item
Let $A$ be a complex $d\times d$ matrix and let $(\la_n)_{n= 1}^d$ and $(\mu_n)_{n= 1}^d$ be the sequences of eigenvalues
of $A$ and $|A|$, respectively, repeated according to
algebraic multiplicities. Show that
$$ \prod_{n=1}^d |\la_n| = \prod_{n=1}^d \mu_n .$$

\noindent{\em Hint:}\ Use the result of Step 1 of the proof of Proposition \ref{prop:Schur}
to see that $\det(A)= \prod_{n=1}^d \la_n.$
Apply this to $|A|$.

\item
Prove the following properties of the partial trace.
\begin{enumerate}[\rm(a), leftmargin=*]
  \item If
  $T\in \calL_1(H\ott K)$ and
  $A_1,A_2\in \calL(H)$, then $(A_1\ot I)T(A_2\ot I) \in \calL_1(H\ott K)$ and
  $$ \tr_K ((A_1\ot I)T(A_2\ot I)) = A_1\tr_K (T)A_2.$$
  \item If $A\in \calL_1(H)$ and $B\in \calL_1(K)$, then $A\ot B\in \calL_1(H\ott K)$ and
  $$\tr_K (A\otimes B) = \tr(B)A.$$
  \item The adjoint of the mapping $T \mapsto \tr_K T$ from $\calL_1(H\ot K)$ to $\calL_1(H)$ is
  the mapping $S \mapsto S \otimes I$ from $\calL(H)$ to $\calL(H\ot K)$.
\end{enumerate}

\item\label{prob:partial-trace-contractive}
Let $H$ and $K$ be separable Hilbert spaces. Show that for every
$T\in\calL_1(H\otimes K)$,
$$
\n\tr_K(T)\n_{\calL_1(H)}
\le \n T\n_{\calL_1(H\otimes K)}.
$$

\noindent{\em Hint:}\ Use trace duality and $\tr(\tr_K(T)A)=\tr(T(A\otimes I))$.

\item\label{prob:MacMahon}
Complete the following outline of a proof of MacMahon's formula
$$ \det (I+A) =\sum_{k=0}^d \tr (\Lambda^k(A))$$
for complex $d\times d$ matrices $A$.

\begin{enumerate}[\rm(a), leftmargin=*]
  \item Prove the formula for the special case when $A$ is diagonalisable, by showing that in this case the formula reduces to the identity
  $$ \prod_{n=1}^d (1+\la_n) = \sum_{k=0}^d \sum_{1\le i_1 < \dots < i_k \le d} \la_{i_1}\cdots \la_{i_k},$$
  where $(\la_n)_{n=1}^d$ is the sequence of eigenvalues of $A$ repeated according to algebraic multiplicities.
  \item Show that the diagonalisable matrices are dense in $M_d(\C)$.
\end{enumerate}

\item\label{prob:tensor-HS-identification}
Let $H$ and $K$ be separable Hilbert spaces, and let $H\otimes K$ be their Hilbert space tensor product. Denote by
$\overline K$ the conjugate Hilbert space of $K$, that is, the Hilbert
space obtained by replacing scalar multiplication $(c,k)\mapsto ck$ by
$(c,k)\mapsto \overline c\,k$ and defining
$\iprod{k}{\ell}_{\overline K}:=\iprod{\ell}{k}_K$.
\begin{enumerate}[\rm(a), leftmargin=*]
\item\label{it:Schmidt0}
Viewing $h\,\bar\otimes\,k$ as the rank-one operator from $\ov K$ to $H$ given by $\overline v\mapsto \iprod{k}{v}_K h$, show that the linear mapping
$$\sum_{j=1}^N h_j\ot k_j
  \mapsto
  \sum_{j=1}^N h_j\,\bar\otimes\,k_j$$
extends uniquely to an isometric isomorphism from $H\ott K$ onto
$\calL_2(\overline K,H)$.

\item\label{prob:Schmidt-decomposition}
Using the singular value decomposition of compact operators
(Theorem \ref{thm:sing-value-comp}) and the operator constructed in part \ref{it:Schmidt0}, prove the {\em Schmidt decomposition}\index{Schmidt decomposition}: for all $u\in H\ot K$ there exist a (possibly empty) finite or
infinite sequence $(s_n)_{n\ge1}$ of strictly positive numbers, an
orthonormal sequence $(h_n)_{n\ge1}$ in $H$, and an orthonormal
sequence $(k_n)_{n\ge1}$ in $K$, such that
$$
  u=\sum_{n\ge1} s_n h_n\ot k_n,
$$
with convergence in $H\otimes K$, and show that $\|u\|_{H\ott K}^2=\sum_{n\ge1}s_n^2$.
Conclude that the coefficients $s_n$, counted with multiplicity and
arranged in nonincreasing order, are uniquely determined.

\item\label{it:Schmidt3}
Assume that $\|u\|=1$, and let $P_u:=u\,\bar\otimes\,u$ be the rank-one projection onto the span of $u$. Prove the partial trace identities
$$
  \tr_K(P_u)=\sum_{n\ge1}s_n^2\,h_n\,\bar\otimes\,h_n, \qquad
  \tr_H(P_u)=\sum_{n\ge1}s_n^2\,k_n\,\bar\otimes\,k_n,
$$
with convergence in trace norm.
\end{enumerate}

\item\label{prob:MacMahon-symm}
Prove the following symmetric analogue of MacMahon's formula: for complex $d\times d$ matrices $A$ with $\n A\n <1$ one has
$$ \frac1{\det (I-A)} =\sum_{n=0}^\infty \tr (\Gamma^n(A)),$$
where $\Gamma^n(A)$ is the natural extension of $A$ to the
$n$-fold symmetric tensor product $\Gamma^n(\C^d)$ (cf. Appendix \ref{sec:tensor}; we use the conventions $\Gamma^0(\C^d)=\C$ and $\Gamma^0(A)=I_{\C}$).

\item\label{prob:ConnesConsani}
Let $A, B \in \calL(H)$ be bounded and positive. Show that if $AB$ is of trace class, then its trace satisfies $\tr(AB) \ge 0.$

\noindent{\em Hint:}\
Use Problem \ref{prob:sigmaSTvsTS} to show that every nonzero
eigenvalue of $AB$ is a positive real number, since
$ \sigma(AB)\setminus\{0\} = \sigma(A^{1/2}BA^{1/2})\setminus\{0\}$.
Then apply Lidskii's theorem to the trace class operator $AB$.

\item
Let $T\in\calL_1(H)$ be positive. Using Theorem \ref{thm:Fredholm-det}, Fubini's theorem, and the scalar analogue of the formula,
prove the trace formula
$$
\tr(T^{1/2})
=
\frac1{2\pi}\int_0^\infty
\la^{-1/2}\log\det(I+T/\la)\ud\la,
$$
where both sides are allowed to take the value $\infty$.
Deduce that $T^{1/2}$ is of trace class if and only if the integral on the
right-hand side is finite.

\item\label{prob:Murphy1}
Using Fedosov's theorem, prove that if $T$ is a Fredholm
operator on $H$ and $T = U|T|$ is its polar decomposition, then $$\ind(T) = \tr(UU^\star - U^\star U).$$

\noindent{\em Hint:}\ Show that $I - U^\star U$ and $I - UU^\star$ are the
projections onto the null spaces of $T$ and $T^\star$, respectively, and that $\codim\Ran(T) = \dim\ker(T^\star)$.

\item\label{prob:Murphy2}
Use Fedosov's theorem to give an alternative proof of the identity
$$\ind(T_1 T_2) = \ind(T_1) + \ind(T_2)$$
for Fredholm operators $T_1$ and $T_2$ acting on $H$.

\item\label{prob:HeltonHowe}
Let $\phi,\psi$ be smooth real-valued functions on the unit circle and let $f,g: \overline{\mathbb{D}}\to \R$ denote their harmonic extensions. Applying Green's theorem to $(f \frac{\partial g}{\partial x}, f \frac{\partial g}{\partial y})$,
show that Theorem \ref{thm:Toeplitz-trace} implies the identity
$$ \tr([T_\phi, T_\psi]) = \frac1{2\pi i} \int_{\mathbb{D}} \Bigl(\frac{\partial f}{\partial x}\frac{\partial g}{\partial y}- \frac{\partial g}{\partial x}\frac{\partial f}{\partial y}\Bigr) \ud x\ud y.$$

\item\label{prob:spectral-shift}
Let $0\le A\le B$ be trace class operators on $H$, and let
$M>\n B\n$. For $t>0$, let $N_T(t)$ count the eigenvalues  $\la$
of a positive operator $T$ satisfying $\la \ge t$, repeated according
to multiplicities, and put $\xi(t):=N_B(t)-N_A(t)$.
\begin{enumerate}[\rm(a), leftmargin=*]
\item\label{it:spectral-shift1} Prove that $\xi\ge0$ and
$\tr(B-A)=\int_0^M\xi(t)\ud t$.
\item\label{it:spectral-shift2} Prove that, for every $f\in C^1([0,M])$, the operator
$f(B)-f(A)$ is of trace class and
$$
\tr\bigl(f(B)-f(A)\bigr)=\int_0^M f'(t)\xi(t)\ud t.
$$
\end{enumerate}
\noindent{\em Hint:}\ For \ref{it:spectral-shift1}, compare the ordered eigenvalues using Theorem \ref{thm:minmax}. For \ref{it:spectral-shift2}, first replace $f$ by $f-f(0)$.

\item
Let $(\lambda_n)_{n\ge1}$ be the Dirichlet eigenvalues of a bounded open set $D\subseteq\R^d$, counted with multiplicity. Prove {\em Berezin's inequality}:\index{inequality!Berezin}
$$
\sum_{n\ge1}(E-\lambda_n)_+\le\frac{2\omega_d|D|}{(d+2)(2\pi)^d}E^{1+d/2},\quad E>0,
$$
where $a_+:=\max\{a,0\}$ and $\omega_d$ is the volume of the unit ball.

\noindent{\em Hint:}\ Adapt the proof of Theorem \ref{thm:Kac-bound}, applying Jensen's inequality to $s\mapsto(E-s)_+$.

\item
Let $(OU(t))_{t\ge0}$ be the Ornstein--Uhlenbeck semigroup on
$L^2(\R^d\!,\gamma)$.
Put $$H_\alpha(x):=\prod_{j=1}^d H_{\alpha_j}(x_j)$$ for
$\alpha\in\N^d$ and $x\in\R^d$, where the $H_n$ are the
one-dimensional Hermite polynomials of Section \ref{subsec:Hermite}.
Using the normalised tensor products $(\alpha!)^{-1/2}H_\alpha$,
with $\alpha!:=\prod_{j=1}^d\alpha_j!$, show that
$$
OU(t)H_\alpha=e^{-|\alpha|t}H_\alpha,
\quad
\alpha\in\N^d.
$$
Deduce that $OU(t)$ is of trace class for every $t>0$ and that
$$
\|OU(t)\|_{\calL_1(L^2(\R^d\!,\gamma))}
=
\tr(OU(t))
=
(1-e^{-t})^{-d}.
$$

\end{problems}

%% file: ch15-QM.tex
\chapter{States and Observables}\label{chap:QM}

\blfootnote{\rm \small \begin{quote}This book has been published by Cambridge University Press in the series `Cambridge Studies in Advanced Mathematics'. The present version is free to view and download for personal use only. Not for re-distribution, re-sale or use in derivative works.\newline {\copyright} Jan van Neerven\end{quote}}

\noindent
In this chapter we apply some of the ideas developed in the preceding chapters to set up a functional-analytic framework for quantum mechanics. More specifically, we will show how the replacement of Borel sets in classical mechanics by orthogonal projections on a Hilbert space leads, in a natural way, to the quantum mechanical formalism for states and observables.

\section{States and Observables in Classical Mechanics}\label{sec:states-observables-classical}

We start by taking a brief look at the notions of state and observable in classical mechanics from a rather abstract measure-theoretic point of view.

\subsection{States}\label{subsec:states-classical}

In classical mechanics, the {\em state space} of a physical system is a measurable space
$(X,\X)$, typically a manifold with its Borel $\sigma$-algebra. For example, the state space of a system of $N$ labelled point particles moving freely in $\R^3$ is $\R^{3N}\times \R^{3N}$ (three position coordinates $x_j$ and three momentum coordinates $p_j$ for each particle) and the state space of the one-dimensional harmonic oscillator is
$\R\times\R$, with position coordinate $x$ and momentum coordinate $p$;
after normalising the physical constants, the circles $\frac12(x^2+p^2)=E$ are the
fixed-energy surfaces.

\begin{definition}[States, pure states]\label{def:clas-state} Let $(X,\X)$ be a measurable space.
\begin{enumerate}[label={\rm(\roman*)}, leftmargin=*]
 \item
A {\em state}\index{state} is a probability measure $\nu$ on $(X,\X)$.
 \item
A {\em pure state}\index{pure state} is an extreme point of the set of probability measures on $(X,\X)$.
\end{enumerate}
For a measurable set $B\in\X$, the number $\nu(B)$ is thought of as `the probability that the system occupies a point of $B$'.
\end{definition}

Thus we identify the `state' of a system with the collection of probabilities
assigned to certain questions about the system.
For example, although in general one has incomplete information about
the exact positions and momenta of all particles in a gas container at
a given time, one might ask about the probability of finding a certain portion of the gas in a certain subset of the container.

Recall that a probability measure $\nu$ on $(X,\X)$ is called
\emph{two-valued}\index{two-valued measure} if
$\nu(B)\in\{0,1\}$ for all $B\in\X$.

\begin{proposition}
The pure states are precisely the two-valued probability measures.
\end{proposition}

This was shown in Example \ref{ex:extreme-points-K}.
In particular, every point $x\in X$ defines a pure state $\delta_x$. Such a state describes a sharply specified phase-space point, whereas a general probability measure describes a statistical ensemble.

\subsection{Observables}\label{subsec:observables:classical}

\begin{definition}[Observables] Let $(X,\mathscr{X})$ and $(\Om,\calF)$
be measurable spaces.
An {\em $\Om$-valued observable}\index{observable} is a measurable function $f: X\to \Om$.
An {\em elementary observable}\index{elementary observable} is a $\{0,1\}$-valued observable.
\end{definition}

For example, the three position coordinates $x_j$ and momentum coordinates $p_j$
of a freely moving particle in $\R^3$ are real-valued observables on the state space $X = \R^3\times \R^3$, and so are the kinetic energy $|p|^2/2m$ (where the mass $m$ is treated as a constant) and potentials $V(x)$.

If the system is in state $\nu$ and $f:X\to\Om$ is an observable, the
distribution of its measured values is the image measure $f_\#\nu$,
defined by
$$
(f_\#\nu)(F):=\nu(f^{-1}(F)),
\quad F\in\calF\!.
$$
Thus $(f_\#\nu)(F)$ is the probability that `measuring $f$ gives a value
in $F$'.

\subsection{From Classical to Quantum}\label{subsec:class-quantum}

Every elementary classical observable is an indicator function $\one_B$
and therefore satisfies
$\one_B^2=\one_B$.
The selfadjoint operators in $\calL(H)$ satisfying $P^2 = P$ are precisely the orthogonal
projections.
The range of an elementary observable $\one_B$ equals $\{0,1\}$ unless $B = \emptyset$ or $B = X$, in which case one has $\one_\emptyset \equiv 0$ and $\one_X \equiv 1$. Orthogonal projections on a complex Hilbert space enjoy similar
properties {\em spectrally}: if $P$ is an orthogonal projection on
a Hilbert space $H$, its spectrum equals
$\sigma(P) = \{0,1\}$ unless $P=0$ or $P=I$; in these cases one has
$\sigma(0) = \{0\}$ and $\sigma(I) = \{ 1\}$.

The basic idea that underlies quantum mechanics is to {\em replace elementary
observables by orthogonal projections}.
The set of all orthogonal projections on a complex Hilbert space $H$ is denoted by $\PP(H)$.\index{$P$@$\PP(H)$} This set is partially ordered in
a natural way by declaring $P_1\le P_2$ to mean that the range of $P_1$ is contained in the range of $P_2$; this is equivalent to the statement that the operator $P_2-P_1$ is positive.
With respect to this partial ordering, $\PP(H)$ is a lattice in the sense of Definition \ref{def:lattice};
for $P_1$ and $P_2$ in $\PP(H)$ the greatest lower bound $$P_1\wedge P_2$$ in $\PP(H)$ is the orthogonal projection
onto $\Ran(P_1)\cap \Ran(P_2)$,
and the least upper bound $$P_1\vee P_2$$ in $\PP(H)$ is the orthogonal projection
onto the closed subspace spanned by $\Ran(P_1)$ and $\Ran(P_2)$. In addition to these operators, the {\em negation}\index{negation} of an orthogonal projection
$P\in \calP(H)$ is the orthogonal projection $$\neg P := I-P$$ onto the orthogonal complement of $\Ran(P)$.
One has the associative laws $$(P_1\wedge P_2)\wedge P_3 = P_1\wedge (P_2\wedge P_3), \qquad
(P_1\vee P_2)\vee P_3 = P_1\vee (P_2\vee P_3)$$ and the identities
$$ \neg(P_1\wedge P_2) = \neg P_1\vee \neg P_2, \qquad  \neg(P_1\vee P_2) = \neg P_1\wedge \neg P_2.$$
The important difference from the classical setting is that the distributive laws
\begin{align*}P_1\wedge (P_2\vee P_3) & = (P_1\wedge P_2)\vee (P_1\wedge P_3)\\
 P_1\vee (P_2\wedge P_3) & = (P_1\vee P_2)\wedge (P_1\vee P_3)
\end{align*}
generally fail. For {\em commuting} orthogonal projections $P,Q$, however, one has
$P\wedge Q=PQ$ and $P\vee Q=P+Q-PQ$. Hence if $P,Q,R$ are pairwise commuting
orthogonal projections, then
\begin{align*} P \wedge (Q\vee R) & = P(Q+R-QR) = PQ+PR-PQR
\\ & = PQ+PR - (PQ)(PR) = (P \wedge Q)\vee (P\wedge R)
\end{align*}
and similarly $P \vee (Q\wedge R) = (P \vee Q)\wedge (P\vee R)$.

\begin{example}
In $\C^2$ consider the orthogonal projections $P_1$, $P_2$, and $P_3$ onto the first and second coordinate axes and the diagonal, respectively. Then $P_1\vee P_2 = I$, $P_1\wedge P_3 = P_2\wedge P_3 =0$, and
$$ (P_1\vee P_2)\wedge P_3 = P_3, \qquad (P_1\wedge P_3)\vee (P_2\wedge P_3) =0.$$
\end{example}

\section{States and Observables in Quantum Mechanics}\label{sec:states-observables}

Motivated by the preceding discussion, we now formulate states and
observables in the Hilbert space setting. States will be represented by
positive trace class operators of trace one, pure states by rank-one orthogonal
projections, and observables by projection-valued measures.

Throughout the remainder of this chapter, unless stated otherwise, $H$ denotes a {\em separable complex Hilbert space.}

\subsection{States}\label{subsec:states}

Upon replacing indicator functions of measurable sets by orthogonal projections on $H$, one is led to the idea of defining a {\em state} as a mapping $\nu:\calP(H)\to [0,1]$ that satisfies $\nu(0)=0$, $\nu(I)=1$, and is {\em countably additive} in the sense that $$\sum_{n\ge 1} \nu(P_n) = \nu(P)$$ whenever $(P_n)_{n\ge 1}$ is a (finite or infinite) sequence of pairwise disjoint orthogonal projections and $P$ is their least upper bound, that is, $P$ is the orthogonal projection onto the closure of the span of the ranges of $P_n$, $n\ge 1$. Here, two orthogonal projections are called {\em disjoint}\index{disjoint projections} if their ranges are mutually orthogonal, and a family of projections is said to be {\em disjoint} if distinct members of this family are disjoint.

Although this definition is quite satisfactory in many ways, countable
additivity only controls sums of pairwise disjoint orthogonal projections. In particular,
if $S=\sum_{n=1}^N c_nP_n$,
where $c_n\ge0$ and the $P_n$ are pairwise disjoint, then the spectral
decomposition of $S$ shows that
$\nu(S):=\sum_{n=1}^Nc_n\nu(P_n)$
is unambiguously defined. What countable additivity does not ensure is
compatibility with arbitrary linear relations between orthogonal projections. Thus, if
$$
\sum_{n=1}^N c_nP_n=\sum_{m=1}^M d_mQ_m
$$
for nonnegative coefficients and orthogonal projections which need not be pairwise
disjoint, there is no reason {\em a priori} why
$$
\sum_{n=1}^Nc_n\nu(P_n)=\sum_{m=1}^Md_m\nu(Q_m)
$$
should hold. Example \ref{ex:cex-Gleason} shows that it may indeed
fail.
The next definition proposes a way of avoiding this difficulty.

Recall that the {\em convex hull} of a subset $S$ of a vector space $V$ is the smallest convex set in $V$ containing $S$ and is denoted by ${\rm co}(S)$.

\begin{definition}[Affine mappings] Let $S$ be a subset of a real vector space $V$, and let $W$ be a real vector space. A mapping $\nu:S\to W$ is called
\emph{affine} if whenever
$$
\sum_{j=1}^m\lambda_j v_j=\sum_{k=1}^n\mu_k w_k,\qquad
\sum_{j=1}^m\lambda_j=\sum_{k=1}^n\mu_k=1,
$$
with all $\lambda_j,\mu_k\ge0$ and $v_j,w_k\in S$, one has
$$
\sum_{j=1}^m\lambda_j\nu(v_j)
=
\sum_{k=1}^n\mu_k\nu(w_k).
$$
\end{definition}
The case $\nu:S\to[0,1]\subseteq\R$ will be used in our
description of states.

If $S$ is convex, a mapping
$\nu:S\to W$ is affine if and only if it preserves convex combinations, that is, if $s = \sum_{j=1}^n \la_j s_j$ with $\la_j\ge 0$ and $\sum_{j=1}^n \la_j  =1$, then $$\nu(s) = \sum_{j=1}^n \la_j \nu(s_j).$$

The following lemma shows that a real-valued affine mapping
extends linearly to the real-linear span of its domain,
provided that the domain contains zero and the mapping
vanishes there.

\begin{lemma}[Linear extension of affine mappings]
\label{lem:affine-linear-extension}
Let $S$ be a subset of a real vector space satisfying $0\in S$, and let
$\nu:S\to\R$ be affine and satisfy $\nu(0)=0$. Then $\nu$ extends
uniquely to a real-linear mapping on the real-linear span of $S$.
\end{lemma}

\begin{proof}
Suppose that $\sum_{j=1}^na_js_j=0$ for certain
$a_j\in\R$ and $s_j\in S$, $j=1,\dots, n$.
Write $a_j=a_j^+-a_j^-$, put
$$
a^+:=\sum_{j=1}^na_j^+,
\qquad
a^-:=\sum_{j=1}^na_j^-,
$$ and let
$ M:=\max\{a^+, a^-\}$.
If $M=0$, there is nothing to prove. Otherwise,
$$\sum_{j=1}^n\frac{a_j^+}{M}s_j
+\Bigl(1-\frac{a^+}{M}\Bigr)0
=
\sum_{j=1}^n\frac{a_j^-}{M}s_j
+\Bigl(1-\frac{a^-}{M}\Bigr)0.
$$
By affinity and the identity $\nu(0)=0$,
$$
\sum_{j=1}^na_j^+\nu(s_j)
=
\sum_{j=1}^na_j^-\nu(s_j).
$$
Thus
$$
\wt\nu\Bigl(\sum_{j=1}^na_js_j\Bigr)
:=
\sum_{j=1}^na_j\nu(s_j)
$$
is well defined. Clearly it is the unique real-linear extension of
$\nu$.
\end{proof}

Let us denote by $\calP_{\rm fin}(H)$\index{$P$@$\calP_{\rm fin}(H)$} the set of all {\em finite-rank projections} in $\calP(H)$, that is, the set of all orthogonal projections with finite-dimensional ranges. To prepare for the definition of a state, we prove the following result.

\begin{proposition}\label{prop:Gleason-light} Let $\nu:\calP_{\rm fin}(H)\to [0,1]$ be affine and satisfy $\nu(0)=0$.
Then there exists a unique positive trace class operator $T$ on $H$ such that
 \begin{align*} \nu(P) = \tr(PT), \quad P\in \calP_{\rm fin}(H).
 \end{align*}
It satisfies $$\tr(T) = \sup_{P\in \calP_{\rm fin}(H)} \nu(P).$$
Conversely, if $T$ is a positive trace class operator on $H$ with $\tr(T) \le 1$, then
 \begin{align*} \nu(P) := \tr(PT), \quad P\in \calP(H),
 \end{align*}
defines an affine mapping $\nu:\calP(H)\to [0,1]$ satisfying $\nu(0)=0$ and
$$ \sup_{P\in \calP_{\rm fin}(H)} \nu(P) = \sup_{P\in \calP(H)} \nu(P) =  \nu(I) = \tr(T).$$
Moreover, $\nu$ is countably additive.
\end{proposition}

\begin{proof}
By Lemma~\ref{lem:affine-linear-extension}, $\nu$ extends uniquely to a
real-linear functional on the real-linear span of $\calP_{\rm fin}(H)$. By the
spectral theorem, this span is the space of selfadjoint finite-rank
operators. Complexification therefore gives a linear functional, still
denoted by $\nu$, on the space $\mathscr F(H)$ of all finite-rank
operators on $H$.
This functional is positive. Indeed, if $S\ge0$ has finite rank, its
spectral decomposition has the form
$S=\sum_{j=1}^n\lambda_jP_j$ with
$\lambda_j\ge0$ and $P_j\in\calP_{\rm fin}(H)$,
and hence
$$
\nu(S)=\sum_{j=1}^n\lambda_j\nu(P_j)\ge0.
$$
If $S=S^\star$ has finite rank and $P$ is the orthogonal projection onto its range,
then
$$
-\n S\n P\le S\le\n S\n P.
$$
Positivity then gives $|\nu(S)|\le  M\n S\n$, where $M:=\sup_{P\in\calP_{\rm fin}(H)}\nu(P)
$.

Writing a finite-rank operator as the sum of its real and
imaginary parts, we obtain
$$
|\nu(S)|\le 2M\n S\n,
\quad S\in\mathscr F(H).
$$
Since the finite-rank operators on $H$ are dense in the space $\mathscr K(H)$ of compact operators on $H$ by Proposition \ref{prop:fr-dense-compH},
$\nu$ extends uniquely to a bounded functional on
$\mathscr K(H)$ of norm at most $2M$.
By trace duality, there exists a unique $T\in\calL_1(H)$ such that
$$
\nu(S)=\tr(ST),
\quad S\in\mathscr K(H).
$$
For every unit vector $h\in H$,
$$
\iprod{Th}{h} = \tr ((h\,\bar\otimes\,h)\circ T) = \nu(h\,\bar\otimes\, h)\ge0,
$$
so $T$ is positive. If $(P_n)_{n\ge1}$ is an increasing sequence in $\mathscr P_{\rm fin}(H)$ converging strongly to $I$, then
$$
\tr(T)
=
\lim_{n\to\infty}\tr(P_nT)
=
\lim_{n\to\infty}\nu(P_n)
\le M.
$$
Conversely, for every $P\in\calP_{\rm fin}(H)$,
$$
\nu(P)=\tr(PT)\le\tr(T),
$$
and therefore $M\le \tr(T)$.

We conclude that  $\tr(T)=M$.
Uniqueness of $T$ follows by testing against the rank-one projections
$h\,\bar\otimes\,h$ with $\n h\n=1$ and using polarisation.

Conversely, let $T\ge0$ be trace class with $\tr(T)\le 1$, and put
$$
\nu(P):=\tr(PT),
\quad P\in\calP(H).
$$
Since $\tr(PT) = \tr((PT^{1/2})T^{1/2}) = \tr(T^{1/2}PT^{1/2})$ and $0\le P\le I$, we have $$0\le\nu(P) = \sumn \iprod{PT^{1/2}h_n}{PT^{1/2}h_n}
\le \sumn \iprod{T^{1/2}h_n}{T^{1/2}h_n} = \tr(T)\le 1$$
for any orthonormal basis $(h_n)_{n\ge 1}$ of $H$.
Thus $\nu$ takes its values in $[0,1]$.
Linearity of the trace shows that $\nu$ is affine. If $(P_n)_{n\ge1}$
are pairwise disjoint and $P$ is their least upper bound, then, writing
$P^{(N)}:=\sum_{n=1}^NP_n$, monotone convergence in an eigenbasis of
$T$ gives
$$
\nu(P)
=
\lim_{N\to\infty}\tr(P^{(N)}T)
=
\sum_{n\ge1}\nu(P_n).
$$
This proves that $\nu$ is countably additive. Finally,
$$
\sup_{P\in\calP(H)}\nu(P)=\nu(I)=\tr(T),
$$
and, since we may approximate $I$ by finite-rank orthogonal projections, the same
supremum is obtained if we replace $\calP(H)$ by $\calP_{\rm fin}(H)$.
\end{proof}

In what follows we denote by $\mathscr{S}(H)$\index{$S$@$\mathscr{S}(H)$} the convex set of all positive trace class operators on $H$ with unit trace. We will see below (Proposition \ref{prop:states-SH}) that
this set is the closed convex hull of its set of extreme points and that these extreme points are precisely the rank-one orthogonal projections on $H$.

A functional $\phi:\calL(H)\to \C$ is called {\em positive}\index{positive!functional} if $\phi(T)\ge 0$ for every positive $T\in \calL(H)$, and {\em normal}\index{normal!functional} if $$\sum_{n\ge 1}\phi(P_n) = \phi(P)$$ for every sequence $(P_n)_{n\ge 1}$ of pairwise disjoint orthogonal projections on $H$, where $P$ denotes their least upper bound.
The same definition of positivity applies to functionals
$\phi:\mathscr K(H)\to\C$.

\begin{theorem}\label{thm:state-posfc} The following six sets are in one-to-one correspondence:
\begin{enumerate}[label={\rm(\arabic*)}, leftmargin=*]
 \item\label{it:state-posfc1} affine mappings $\nu:\calP_{\rm fin}(H)\to [0,1]$ satisfying $\nu(0)=0$ and $ \sup_{P\in \calP_{\rm fin}(H)}\nu(P)=1$;
 \item\label{it:state-posfc1a} countably additive affine mappings
$\nu:\calP(H)\to [0,1]$ satisfying $\nu(0)=0$ and $\nu(I)=1$;
 \item\label{it:state-posfc2} positive trace class operators $T$ on $H$ satisfying $\tr(T)=1$, via
 $$ \nu(P) = \tr(PT), \quad P\in \calP_{\rm fin}(H);$$
 \item\label{it:state-posfc2a} positive trace class operators $T$ on $H$ satisfying $\tr(T)=1$, via
 $$ \nu(P) = \tr(PT), \quad P\in \calP(H);$$
 \item\label{it:state-posfc3} positive  functionals $\phi:\mathscr{K}(H)\to \C$ satisfying $ \sup_{P\in \calP_{\rm fin}(H)}\phi(P)=1$, via
 $$ \phi(S) = \tr(ST), \quad S\in \mathscr{K}(H);$$
 \item\label{it:state-posfc4} positive normal functionals $\phi:\calL(H)\to \C$ satisfying $\phi(I)=1$, via
 $$ \phi(S) = \tr(ST), \quad S\in \calL(H).$$
\end{enumerate}
\end{theorem}

\begin{proof}
We show that each description amounts to specifying a unique
positive trace class operator $T$ with $\tr(T)=1$.

First, such an operator defines a functional
$$
\phi_T(S):=\tr(ST),\quad S\in\calL(H).
$$
This functional is positive, since for $S\ge0$ we have
$\phi_T(S)=\tr(S^{1/2}TS^{1/2})\ge0$.
Proposition \ref{prop:Gleason-light} shows that the restriction of $\phi_T$ to
$\calP(H)$ is affine and countably additive, with value $1$ at $I$
and supremum $1$ on $\calP_{\rm fin}(H)$. Thus $\phi_T$ is normal,
and restriction gives all the objects listed in the theorem.

Conversely, an affine mapping as in \ref{it:state-posfc1} is
represented by a unique such $T$ by Proposition
\ref{prop:Gleason-light}. If $\nu$ is as in
\ref{it:state-posfc1a}, countable additivity and $\nu(I)=1$ imply
that its restriction to $\calP_{\rm fin}(H)$ has supremum $1$.
Let $T$ represent this restriction. For any $P\in\calP(H)$,
choose finite-rank orthogonal projections $P_n\uparrow P$.
Countable additivity of both projection mappings gives
$$
\nu(P)=\limn\nu(P_n)
=\limn\tr(P_nT)=\tr(PT).
$$
This establishes the correspondences involving
\ref{it:state-posfc1}, \ref{it:state-posfc1a},
\ref{it:state-posfc2}, and \ref{it:state-posfc2a}.

\smallskip
Now let $\phi$ be as in \ref{it:state-posfc3}.
We first check that $\phi$ is bounded. Otherwise, decomposing
compact operators into real and imaginary parts and then into
positive and negative parts, we could choose positive compact
operators $S_n$ with $\n S_n\n\le2^{-n}$ and $\phi(S_n)\ge1$.
The norm-convergent sum $S:=\sum_{n\ge1}S_n$ would then be compact,
and positivity would give
$$
\phi(S)\ge\sum_{n=1}^N\phi(S_n)\ge N,\quad N\ge1,
$$
a contradiction.
By trace duality (Theorem \ref{thm:traceduality}), there is therefore
a unique $T\in\calL_1(H)$ such that
$\phi(S)=\tr(ST)$ for all $S\in\mathscr K(H)$.
Testing on rank-one orthogonal projections shows that $T\ge0$,
and
$$
\tr(T)=\sup_{P\in\calP_{\rm fin}(H)}\tr(PT)
=\sup_{P\in\calP_{\rm fin}(H)}\phi(P)=1.
$$

Finally, let $\phi$ be as in \ref{it:state-posfc4}.
Its restriction to $\calP(H)$ satisfies
\ref{it:state-posfc1a}, so the preceding argument gives a unique
$T\in\mathscr S(H)$ such that $\phi(P)=\tr(PT)$ for every
orthogonal projection $P$.
For selfadjoint $S$, positivity and
$-\n S\n I\le S\le\n S\n I$ give $|\phi(S)|\le\n S\n$.
Approximating $S$ in operator norm by finite real linear
combinations of its spectral projections therefore gives
$\phi(S)=\tr(ST)$ for every selfadjoint $S\in\calL(H)$.
By complex linearity the identity holds for every
$S\in\calL(H)$. This completes all the correspondences.
\end{proof}

The countable additivity assumption in \ref{it:state-posfc1a} cannot be
omitted if $H$ is infinite-\-dimen\-sional:

\begin{example}[Failure of countable additivity] Let $(e_n)_{n\ge1}$ be an orthonormal basis of $H$. Using the real Hahn--Banach theorem, let
$L_\R:\ell_\R^\infty\to\R$ be a linear functional of norm one extending the ordinary limit on convergent real sequences. Its complexification $L:\ell^\infty\to\C$ satisfies $L(\one)=1$ and is positive: for if $0\le x\le\n x\n_\infty\one$ with $x\not=0$, then $x\in \ell_\R^\infty$ and $$ \Bigl\n\one-\frac{x}{\n x\n_\infty}\Bigr\n_\infty\le1, $$
and therefore
$$ \Bigl|L\Bigl(\one - \frac{x}{\n x\n_\infty}\Bigr)\Bigr| = \Bigl|L_\R\Bigl(\one - \frac{x}{\n x\n_\infty}\Bigr)\Bigr| = \Bigl|1 - \frac{L_\R(x)}{\n x\n_\infty}\Bigr|\le 1,$$
which forces $L(x) = L_\R(x) \ge 0$.

Define
$$
        \phi(T):=L\bigl((\iprod{Te_n}{e_n})_{n\ge1}\bigr),
        \quad T\in\calL(H).
$$
Then $\phi(T)\ge 0$ for all $T\ge 0$, and $\phi(0)= 0$ and $\phi(I)=1$.
If
$K\in\mathscr K(H)$, then $\iprod{Ke_n}{e_n} \to 0$ as $n\to \infty$,
and therefore $\phi(K)=0$.

Now set
$$
        \nu(P):=\phi(P),\quad P\in\calP(H).
$$
Then $\nu(0)=0$, $\nu(I)=1$, and the linearity of $\phi$ implies directly that $\nu$ is affine; it also implies that $\nu$ is finitely additive on mutually orthogonal projections. However, if $P$ has finite rank, then $P\in\mathscr K(H)$, and
therefore $\nu(P)=0$.

If $P_n = e_n\,\bar\otimes\,e_n$ denotes the rank-one orthogonal projection onto the one-dimensional subspace spanned by $e_n$, then
$\sum_{n=1}^\infty P_n = I$ in the strong operator topology. On the other hand,
$$
\sum_{n=1}^\infty \nu(P_n) = 0 \neq 1 = \nu(I),
$$
so $\nu$ is not countably additive.
Proposition \ref{prop:Gleason-light} can be applied to the restriction
of $\nu$ to $\calP_{\rm fin}(H)$; the representing operator is $T=0$.
\end{example}

We may now define the set of {\em states} as any one of these six sets. For definiteness we take the sixth:

\begin{definition}[States]\label{def:state}
 A {\em state}\index{state} is a positive normal functional $\phi:\calL(H)\to \C$ satisfying $\phi(I)=1$.
\end{definition}

This definition captures what is generally called a {\em normal state} in the
mathematical literature; the term {\em state} is often
reserved for arbitrary positive functionals $\phi:\calL(H)\to\C$ satisfying
$\phi(I)=1$. A positive functional on $\calL(H)$ is normal if and only if
its restriction to the closed operator unit ball is continuous
for the weak operator topology; see Theorem \ref{thm:normal-continuity}.
All states considered in this chapter are normal, and for this reason {\em by `state' we will always understand a normal state.}

\begin{theorem}[Normality and operator topologies]
\label{thm:normal-continuity}
Let $\phi:\calL(H)\to\C$ be a positive linear functional.
The following assertions are equivalent:
\begin{enumerate}[label={\rm(\arabic*)}, leftmargin=*]
\item\label{it:normal-continuity1}
$\phi$ is normal;
\item\label{it:normal-continuity2}
the restriction of $\phi$ to the closed unit ball of $\calL(H)$
is continuous for the weak operator topology;
\item\label{it:normal-continuity3}
the restriction of $\phi$ to the closed unit ball of $\calL(H)$
is continuous for the strong operator topology.
\end{enumerate}
\end{theorem}

\begin{proof}
\ref{it:normal-continuity1}
$\Rightarrow$\ref{it:normal-continuity2}: \
If $\phi(I)=0$, positivity implies $\phi=0$, and the assertion
is immediate. Otherwise, applying Theorem~\ref{thm:state-posfc} to $\phi/\phi(I)$ gives a positive
trace class operator $T_0$ with unit trace. Putting $T:=\phi(I)T_0$,
we obtain
$$
\phi(S)=\tr(ST),\quad S\in\calL(H).
$$
Choose finite-rank operators $T_n$ converging to $T$ in trace
norm. Each functional $S\mapsto\tr(ST_n)$ is a finite sum of
functionals of the form $S\mapsto\iprod{Sx}{y}$, with $x,y\in H$,
and is therefore continuous for the weak operator topology. Moreover,
$$
\sup_{\|S\|\le1}|\tr(ST)-\tr(ST_n)|
\le\|T-T_n\|_1\to0.
$$
Thus the restriction of $\phi$ to the closed unit ball is a
uniform limit of weak operator continuous functions, and is
therefore continuous in this topology.

\smallskip
\ref{it:normal-continuity2}
$\Rightarrow$\ref{it:normal-continuity3}: \
The strong operator topology is finer than the weak operator
topology.

\smallskip
\ref{it:normal-continuity3}
$\Rightarrow$\ref{it:normal-continuity1}: \
Let $(P_n)_{n\ge1}$ be pairwise disjoint orthogonal projections,
and let $P$ be their least upper bound. The partial sums
$\sum_{n=1}^N P_n$ belong to the closed unit ball and converge
strongly to $P$. Hence
$$
\phi(P)
=\lim_{N\to\infty}\phi\Bigl(\sum_{n=1}^N P_n\Bigr)
=\sum_{n\ge1}\phi(P_n),
$$
which proves normality.
\end{proof}

\begin{remark}[Density operators]
In the physics literature, the positive trace class operator $T$ with unit trace associated with a state $\phi$ is called the {\em density operator}\index{density!operator} associated with $\phi$.
\end{remark}

As the following example shows, a countably additive mapping $\nu:\calP(\C^2)\to [0,1]$ satisfying $\nu(0)=0$ and $\nu(I)=1$ need not be affine (and therefore need not define a state).

\begin{example}[Failure of affinity in two dimensions]
\label{ex:cex-Gleason}
Let $H=\C^2$, fix an orthonormal basis $(e_1,e_2)$, and define
$\nu:\calP(H)\to[0,1]$ by $\nu(0)=0$, $\nu(I)=1$, and, for every
rank-one orthogonal projection $P$,
$$
\nu(P):=
\begin{cases}
      1,&\text{if }\iprod{Pe_1}{e_1}>\frac12,\\
\frac12,&\text{if }\iprod{Pe_1}{e_1}=\frac12,\\
      0,&\text{if }\iprod{Pe_1}{e_1}<\frac12.
\end{cases}
$$
If $P$ has rank one, then $I-P$ has rank one and
$\iprod{(I-P)e_1}{e_1}=1-\iprod{Pe_1}{e_1}$.
It follows that
$$
\nu(P)+\nu(I-P)=1.
$$
As the only nontrivial decomposition of an orthogonal projection on $\C^2$ into
disjoint nonzero orthogonal projections is $I=P+(I-P)$, the mapping
$\nu$ is countably additive.

We show that $\nu$ is not affine. Fix $0<a<1$ and put
$$
h_\pm
:=
\Bigl(\frac{1+a}{2}\Bigr)^{1/2}e_1
\pm
\Bigl(\frac{1-a}{2}\Bigr)^{1/2}e_2.
$$
Let $P_\pm:=h_\pm\bar\otimes h_\pm$ and
$P_j:=e_j\bar\otimes e_j$, $j=1,2$. Direct calculation gives
\begin{align}\label{ex:2-dim-Gleason}
\frac12P_++\frac12P_-
=
\frac{1+a}{2}P_1+\frac{1-a}{2}P_2.
\end{align}
On the other hand,
$$
\nu(P_+)=\nu(P_-)=\nu(P_1)=1,
\qquad
\nu(P_2)=0.
$$
For example,
$ \iprod{P_+e_1}{e_1} = |\iprod{e_1}{h_+}|^2 = \frac12(1+a) > \frac12$,
hence $\nu(P_+)=1$; the other identities are proved similarly.
Thus, if $\nu$ were affine, the two sides of \eqref{ex:2-dim-Gleason} would be assigned
the respective values $1$ and $\frac{1+a}{2}$.
Since these are unequal, $\nu$ is not affine.
\end{example}

It is not a coincidence that this counterexample lives in two dimensions: a celebrated theorem due to Gleason asserts that if $\dim(H)\ge 3$, then every countably additive mapping $\nu:\calP(H)\to [0,1]$ satisfying $\nu(I)=1$ is affine and hence defines a state.

The following operator version of the monotone convergence theorem will be needed in Section \ref{sec:POVM}.

\begin{proposition}[Monotone convergence for normal states]
\label{prop:normal-MCT}
Let $\phi:\calL(H)\to\C$ be a normal state.
Suppose that $(Q_n)_{n\in\N}$ is an increasing sequence of positive operators such that $0\le Q_1\le Q_2\le \dots  \uparrow Q \le I$ and $Q_n\to Q$ weakly as $n\to \infty$.
Then
$$
  \phi(Q_n)\uparrow \phi(Q)\quad\text{as }n\to\infty.
$$
\end{proposition}

\begin{proof}
By the normality of $\phi$, there exists a unique positive trace class operator
$T \in\calL_1(H)$ with $\tr(T)=1$ such that $\phi(S) = \tr(ST)$ for all $S\in\calL(H)$. Choose an orthonormal basis $(h_k)_{k\ge 1}$ of $H$ consisting of eigenvectors of $T$, that is, $T  h_k = \lambda_k h_k$ with $\lambda_k\ge 0$ for all $k\ge 1$, and $\sum_{k=1}^\infty \lambda_k = 1$.
Since $0\le Q_n\uparrow Q$ and $Q_n\to Q$ weakly,
for all $k\ge 1$ we have
$$\iprod{Q_n h_k}{h_k} \uparrow \iprod{Q h_k}{h_k} \ \ \hbox{as $n\to\infty$}.$$
Hence, by the monotone convergence theorem (viewing sums as integrals with respect to the counting measure),
\[
\limn \phi(Q_n) = \limn \tr(Q_n T)
   = \limn  \sum_{k=1}^\infty \lambda_k \iprod{Q_n h_k}{h_k}
   = \sum_{k=1}^\infty \lambda_k \iprod{Q h_k}{h_k}
   = \tr(Q T)
   = \phi(Q).
\qedhere\]
\end{proof}

\subsection{Pure States}

Theorem \ref{thm:state-posfc} establishes six equivalent ways of looking at the convex set of all states.
Since the correspondences between them preserve convex combinations and hence extreme points, the following definition makes sense from each of these points of view:

\begin{definition}[Pure states]
A {\em pure state}\index{pure state} is an extreme point of the convex set of states.
\end{definition}

\begin{proposition}\label{prop:pure} A state $\phi:\calL(H)\to \C$ is pure if and only if it is a {\em vector state},\index{vector!state} that is,
there exists a unit vector $h\in H$ such that
$$\phi(S) = \iprod{Sh}{h}, \quad S\in \calL(H). $$
This unit vector is unique up to a scalar multiple of modulus one.
\end{proposition}

The first assertion can be equivalently stated as saying that the extreme points of the set of all positive trace class operators with unit trace are precisely the orthogonal projections of rank one.

\begin{proof}
`Only if': \ Let $\phi$ be a pure state and let $T$ be the associated positive trace class operator on $H$ with unit trace.
By the singular value decomposition (Theorem \ref{thm:tc-ell1}) we have
$ T =\sum_{n\ge 1} \la_n h_n\,\bar\otimes\, h_n$ for some orthonormal basis $(h_n)_{n\ge 1}$ of $H$ and a nonnegative scalar sequence $(\la_n)_{n\ge 1}$ such that $\sum_{n\ge 1} \la_n = \tr(T) = 1$.
If at least two of the numbers $\lambda_n$ are nonzero, choose an
index $j$ such that
$0<\lambda_j<1$, put
$$
P:=h_j\,\bar\otimes\,h_j,
\qquad
R:=
\frac{T-\lambda_jP}{1-\lambda_j},
$$
and note that $R$ is a positive trace class operator of trace one.
Moreover, $R\neq P$, and
$$
T=\lambda_jP+(1-\lambda_j)R
$$
is a nontrivial convex decomposition into two distinct states.
Consequently, a pure state can have only one nonzero eigenvalue.
Since $\tr(T)=1$, this eigenvalue must equal one, and hence
$T=h\,\bar\otimes\,h$
for some unit vector $h\in H$. This shows that $\phi$ is a vector state.

\smallskip
`If': \ If $\phi$ is a vector state, then the associated operator in $\mathscr S(H)$ is of the form $T = h\,\bar\otimes\,h$ with $\n h\n =1$. Suppose that $T = (1-\la) T_0 + \la T_1$ with $0<\la<1$ and $T_0,T_1 \in \mathscr S(H)$.
If $g\perp h$, then
$$
0 = \iprod{Tg}{g} = (1-\lambda)\iprod{T_0g}{g} + \lambda\iprod{T_1g}{g}.
$$
Since $\iprod{T_0g}{g}\ge 0$ and $\iprod{T_1g}{g}\ge 0$, it follows that
 $\iprod{T_0g}{g} = \iprod{T_1g}{g} = 0$
for all $g\perp h$, so $T_0$ and $T_1$ vanish on $\{h\}^\perp$. By selfadjointness their ranges are contained in
$\operatorname{span}\{h\}$, and therefore there exist $c_0,c_1\ge 0$ with $T_0=c_0 T$ and $T_1=c_1 T$. Taking traces gives $c_0=c_1=1$, so $T_0=T_1=T$. It follows that $T$ is an extreme point of the convex set $\mathscr S(H)$. This means that $\phi$ is a pure state.

It remains to prove the converse uniqueness assertion. Suppose that $h,h'\in H$
are unit vectors and that
$$
\iprod{Sh}{h}=\iprod{Sh'}{h'},\quad S\in\calL(H).
$$
Then the associated rank-one operators are equal, that is,
$h\,\bar\otimes\,h=h'\,\bar\otimes\, h'$.
Applying both sides to $h$ gives
$$
h=\iprod{h}{h}h=\iprod{h}{h'}h'.
$$
Thus $h$ is a scalar multiple of $h'$. Since both vectors have norm one, the
scalar has modulus one. Conversely, for all $\theta\in\R$ one has
$$
(e^{i\theta}h)\bar\otimes(e^{i\theta}h)=h\,\bar\otimes\,h,
$$
and hence $e^{i\theta}h$ and $h$ define the same vector state.
\end{proof}

\begin{remark}[Bras, kets, superpositions, mixed states]
In physics notation, a vector $h\in H$ is written $\ket{h}$; the symbol $\ket{\cdot}$ is referred to as a {\em ket}\index{ket}.
For a unit vector $h$, the associated pure state is the {\em ray}\index{ray} $\C h$, or
equivalently, the rank-one projection $\ket{h}\bra{h}$ onto this ray (denoted here by $h\,\bar\otimes\,h$); in this context the symbol $\bra{\cdot}$ is referred to as a {\em bra}\index{bra}.
In this notation, our bracket $\iprod{Th}{h}$ becomes the `bra-ket' $\bra{h}\!T\!\ket{h}$. We retain our convention that the inner product is linear in
its first variable. Thus the physics bra--ket $\lb h|k\rb$
corresponds to $\iprod{k}{h}$ in our notation.

The ket notation $\ket{h}$ is also used for the associated pure state when no ambiguity can arise. If $H$ is realised as an $L^2$-space, the representing
function $h$ is called a {\em wave function}\index{wave!function}.

The addition in $H$ can be used to define, for orthogonal unit vectors $h_1,h_2\in H$ and scalars $\al_1,\al_2\in \C$ satisfying $|\al_1|^2+|\al_2|^2=1$, the pure state $$| \al_1 h_1+\al_2 h_2\rb$$
or, in physics notation, $\al_1\ket{h_1} + \al_2 \ket{h_2}$.
Such states are referred to as (coherent) {\em superpositions}\index{superposition} of the states  $\ket{h_1}$ and $\ket{h_2}$.
Such states should be carefully distinguished from states that can be built by using the addition of $\calL_1(H)$.
Indeed, if $h_1$ and $h_2$ are orthogonal unit vectors in $H$
and if $\la\in [0,1]$, then the convex combination
$$ (1-\la) h_1\,\bar\otimes\,h_1 + \la h_2\,\bar\otimes\,h_2$$ or, in physics notation, $(1-\la) \ket{h_1}\bra{h_1} + \la \ket{h_2}\bra{h_2}$,
is a positive trace class operator with unit trace on $H$, and hence defines a state.
Such states, which are not pure unless $\la=0$ or $\la=1$, are called {\em mixed states}\index{mixed state} or, more precisely, {\em mixtures} of $\ket{h_1}$ and $\ket{h_2}$.
\end{remark}

\begin{proposition}\label{prop:states-SH} The set $\mathscr{S}(H)$ is the closed convex hull of its extreme points, the closure being taken with respect to the trace norm.
The extreme points are precisely the rank-one orthogonal projections of the form $h\,\bar\otimes\,h$ with $h\in H$ of norm one.
\end{proposition}
\begin{proof}
By the singular value decomposition of Theorem \ref{thm:tc-ell1}, every
$T\in\mathscr S(H)$ is of the form
$T=\sum_{n\ge1}\lambda_nh_n\,\bar\otimes\,h_n$,
with convergence in trace norm, where $(h_n)_{n\ge1}$ is an orthonormal
basis of $H$ and $(\lambda_n)_{n\ge1}$ is a nonnegative sequence satisfying
$\sum_{n\ge1}\lambda_n=1$.
Put $s_N:=\sum_{n=1}^N\lambda_n$
and, for $N$ such that $s_N>0$, define
$$
T_N:=\frac1{s_N}\sum_{n=1}^N
\lambda_nh_n\,\bar\otimes\,h_n.
$$
Then $T_N$ is a finite convex combination of rank-one orthogonal projections.
Moreover, since all the operators involved are diagonal with respect to
$(h_n)_{n\ge1}$,
$$
\n T_N-T\n_1
=
\sum_{n=1}^N\Bigl|\frac{\lambda_n}{s_N}-\lambda_n\Bigr|
+\sum_{n>N}\lambda_n
=
2(1-s_N)\to 0.
$$
Thus every element of $\mathscr S(H)$ belongs to the closed
convex hull in $\calL_1(H)$ of the rank-one orthogonal projections.
This proves the first assertion.

The second assertion follows from Theorem \ref{thm:state-posfc}, which
identifies operators of the form $h\,\bar\otimes\,h$, with $\|h\|=1$,
with vector states. By Proposition \ref{prop:pure}, these are precisely
the extreme points of the convex set of all states.
\end{proof}

\subsection{Observables}

Let $(\Om,\calF)$ be a measurable space. Classically, an $\Om$-valued observable on the state space $(X,\mathscr{X})$ is a measurable function $f:X\to \Om$. By definition of measurability,
$f$ induces a mapping from $\calF$ to $\X$ given by $$F \mapsto f^{-1}(F), \quad F\in \calF\!,$$ and this mapping is countably additive, in the sense that if the sets $F_n\in\calF$ are pairwise disjoint, then
$f^{-1}(\bigcup_{n\ge 1} F_n) =\bigcup_{n\ge 1}  f^{-1}(F_n)$.
Identifying sets in $\X$ by their indicator functions and replacing them by orthogonal projections on a Hilbert space $H$, and interpreting `countable additivity' as strong countable additivity in the sense of Section \ref{sec:PVM}, we arrive at the following definition of an observable in quantum mechanics.

Now let $H$ be a separable complex Hilbert space.

\begin{definition}[Observables]\label{def:observable}\index{observable} Let $(\Om,\calF)$ be a measurable space.
An {\em $\Om$-valued observable} is a projection-valued measure (PVM) $P:\calF\to \mathscr{P}(H)$.
An {\em elementary observable}\index{elementary observable} is a $\{0,1\}$-valued observable.
\end{definition}

Identifying a projection-valued measure $P$ on $\{0,1\}$ with the orthogonal projection $P(\{1\})$, the elementary observables are in one-to-one correspondence with the orthogonal projections.
This should be compared with the classical situation, in which elementary
observables are the indicator functions of measurable sets.

Observables defined in this way are sometimes called {\em sharp observables}\index{sharp observable}, as opposed to {\em unsharp observables}, which will be introduced in Section \ref{subsec:POVM}.

Following the notation introduced in Chapter \ref{chap:spectral-theorem}, we write $P_F := P(F)$ for $F\in\calF$\!. For vectors $h\in H$, we denote by $P_h$ the finite measure on $(\Om,\calF)$ given by $$P_h(F) := \iprod{P_F h}{h}, \quad F\in\calF\!.$$
In the language of Chapter \ref{chap:spectral-theorem}, a real-valued observable is nothing but
a projection-valued measure $P$ on $\R$, and by the spectral theorem (Theorem \ref{thm:ST-unboundedn-normal})
we can associate a unique selfadjoint operator $A$ with $P$,
determined by
$$ \Dom(A) = \Big\{h\in H:\, \int_\R |\la|^2 \ud P_h(\la)< \infty\Big\}$$
and
$$ \iprod{Ah}{h} = \int_\R \lambda\ud P_h(\la), \quad h\in\Dom(A).$$
In the converse direction, the spectral theorem asserts that every selfadjoint operator
$A$ arises from a projection-valued measure on $\R$ in this way and hence defines an observable.

Thus we arrive at the conclusion that {\em real-valued} observables are in one-to-one correspondence with self\-adjoint operators. In most treatments of quantum mechanics this is simply taken as a postulate. Defining real-valued observables as projection-valued measures
makes their correspondence with selfadjoint operators a consequence
of the spectral theorem. This connects the operator description
with the passage from classical events to orthogonal projections.
A further advantage of the present approach is that, in the same vein,
the spectral theorem for normal operators can be reinterpreted as establishing a one-to-one correspondence between complex-valued observables and normal operators, and likewise between observables with values in the unit circle and unitary operators.

We return to the abstract setting of observables $P:\calF\to \calP(H)$ with values in $\Om$. If $\phi:\calL(H)\to\C$ is a pure state represented by the unit vector $h\in H$, then
$$\phi(P_F) = \iprod{P_F h}{h} = P_h (F), \quad F\in \calF\!,$$
so the assignment $F\mapsto \phi(P_F)$ defines a probability measure.
The following proposition is an immediate consequence of the fact that states are normal
and that the least upper bound of a sequence of disjoint orthogonal projections is given by their sum in the strong operator topology.
It is the mathematical counterpart of the so-called Born rule in quantum mechanics
and  allows us to interpret the number $\phi(P_F)$
as `the probability that measuring $P$ produces a value in $F\in\F$
when the system is in state $\phi$'.

\begin{proposition}[Born rule]\label{prop:observ-mapping}
 If $\phi:\calL(H)\to\C$ is a state and $P:\calF\to \calP(H)$ an $\Om$-valued observable, the mapping $$ F \mapsto \phi(P_F), \quad F\in \calF\!,$$
 defines a probability measure on $(\Om,\calF)$.
\end{proposition}

If $P$ is a real-valued (or complex-valued) observable represented by a selfadjoint (or, more generally, a normal) operator $A$,
then, as a projection-valued measure, $P$ is supported on the spectrum $\sigma(A)$ (see Proposition \ref{prop:supp-P-sigmaA}) and therefore $P$ can be thought of
as a
$\sigma(A)$-valued observable. The physical interpretation is that `with probability one,
a measurement of $A$ produces a value belonging to $\sigma(A)$'.

\subsection{The Uncertainty Principle}

If $P$ is a real-valued observable represented by a bounded selfadjoint operator $A$, the {\em expected value of $P$ in state $\phi$}\index{expected value!of an observable} is defined as the number
$$ \lb A\rb_\phi := \phi(A).$$
If $\phi = \ket{h}$ is a pure state, this definition reduces to $ \lb A\rb_{\ket{h}} = \iprod{Ah}{h}$. The latter may be generalised to unbounded selfadjoint operators $A$ by defining, for a pure state
$\ket h$ with $h\in\Dom(A)$,
$$
\lb A\rb_{\ket h}:=\iprod{Ah}{h}.
$$
In this situation, for $h\in \Dom(A)$ we can define the {\em variance}\index{variance} by
$$\var_{\ket{h}}(A) := \n (A - \lb A\rb_{\ket h})h\n^2\!.$$
For $h\in \Dom(A^2)$ we have
$$\var_{\ket{h}}(A) = \lb A^2 \rb_{\ket h} - \lb A\rb_{\ket{h}}^2 = \big\lb (A - \lb A\rb_{\ket{h}})^2\big\rb_{\ket{h}}.$$
The {\em uncertainty} of $A$ in state ${\ket{h}}$ is defined by
$$\Delta_{\ket{h}}(A) :=  (\var_{\ket{h}}(A))^{1/2}\!.$$

\begin{theorem}[Uncertainty principle]\label{thm:uncertainty}\index{uncertainty principle} Let $\ket{h}$ be a pure state associated with the unit vector $h\in H$,
and consider two real-valued observables with associated selfadjoint operators
$A$ and $B$.
If $h\in \Dom([A,B]):= \{h\in \Dom(A)\cap \Dom(B): \, Ah\in \Dom(B), \,  Bh\in \Dom(A)\}$ and $[A,B]h:= ABh  - BAh$, then
$$ \Delta_{\ket{h}}(A)\Delta_{\ket{h}}(B) \ge \frac12 |\iprod{[A,B]h}{h}|.$$
\end{theorem}
\begin{proof} The operators $\wt A:= A -  \iprod{Ah}{h}$ and $\wt B:= B -  \iprod{Bh}{h}$
with domains $\Dom(\wt A) = \Dom(A)$ and $\Dom(\wt B) =  \Dom(B)$ are selfadjoint. In particular,
we note that $\wt A h\in \Dom(\wt B)$, $\wt B h\in \Dom(\wt A)$, and we have $[\wt A,\wt B]h = [A,B]h$.
The Cauchy--Schwarz inequality implies
\begin{align*}
\Delta_{\ket{h}}(A)\Delta_{\ket{h}}(B)
  & = \n \wt A h\n\n \wt B h\n
 \ge |\iprod{\wt A h}{\wt B h}|
  \ge |\Im \iprod{\wt A h}{\wt B h}|
 \\ & = \frac12 |\iprod{\wt A h}{\wt B h} - \iprod{\wt B h}{\wt A h}|
  = \frac12 |\iprod{[\wt A, \wt B]h}{h}| = \frac12 |\iprod{[A, B]h}{h}|.
  \qedhere
\end{align*}
\end{proof}

The physical interpretation of the next result is that a measurement of $A$ in a pure state $\ket{h}$ gives the expected value $\iprod{Ah}{h}$ with probability one if and only if the representing unit vector $h$ is an eigenvector of $A$, and in this case the eigenvalue equals $\iprod{Ah}{h}$.

\begin{proposition}\label{prop:pure-state} Let $P$ be a real-valued observable, represented by the selfadjoint operator $A$, and let $h\in \Dom(A)$ satisfy $\n h\n=1$. The following assertions are equivalent:
\begin{enumerate}[label={\rm(\arabic*)}, leftmargin=*]
\item\label{it:pure-state1} $A$ has zero uncertainty in the state $\ket{h}$;
\item\label{it:pure-state2} $h$ is an eigenvector for $A$.
\end{enumerate}
If these equivalent conditions hold, then for the corresponding eigenvalue $\la$ we have $$
 \la = \iprod{Ah}{h} \ \ \hbox{and} \ \  \iprod{P_{\{\la\}}h}{h} = 1.$$
\end{proposition}

\begin{proof}
By definition,
$\var_{\ket h}(A)=0$ if and only if $(A-\iprod{Ah}{h})h=0$.
The latter is equivalent to saying that $h$ is an eigenvector of $A$, with
eigenvalue $\lambda=\iprod{Ah}{h}$. In this case
Corollary~\ref{cor:normal-eigen-mapping}, applied to
$\one_{\{\lambda\}}$, gives
$P_{\{\lambda\}}h=h$,
and hence $\iprod{P_{\{\lambda\}}h}{h}=1$.
\end{proof}

\subsection{The Qubit}\label{subsec:qubit}

It is instructive to take a closer look at the simplest genuinely quantum mechanical system, the qubit. It is the quantum version of the {\em bit}\index{bit} $\{0,1\}$, which we think of as equipped with the counting measure $\mu$ giving mass $1$ to each of the two elements of $\{0,1\}$. A physical realisation is the spin degree of freedom of a spin-$1/2$ particle.

We write
$\one_{\{0\}}$ and $\one_{\{1\}}$ for the standard basis vectors of the Hilbert space $L^2(\{0,1\})$ and denote the pure states associated with them by $\ket{0}$ and $\ket{1}$. Every pure state is then of the form  $\alpha\ket{0} +\beta\ket{1}$ with $\alpha,\beta\in\C$ satisfying $|\alpha|^2 + |\beta|^2 = 1$.
Since multiplying a representing unit vector by a scalar of modulus one
does not change the pure state, except at the poles $\theta=0$ and $\theta=\pi$, every
pure state can be uniquely written in the form
\begin{align*}
\cos(\theta/2)\ket{0} +e^{i\varphi }\sin(\theta/2)\ket{1}
\end{align*}
for suitable $0\leq \theta \leq \pi$ and $0\leq \varphi <2\pi$. In spherical coordinates, the variables $\theta $ and $\varphi$
uniquely determine a point
\begin{align*}
(\sin \theta \cos \varphi ,\,\sin \theta \sin \varphi ,\,\cos \theta)
\end{align*}
on the unit sphere $S^2$ of $\R^3$. Conversely, the point determines
$\theta$ uniquely, and also $\varphi$ unless $\theta=0,\pi$. This representation of pure states is frequently referred to as the {\em Bloch sphere}.\index{Bloch!sphere}

In what follows, we identify $L^2(\{0,1\})$ isometrically with $\C^2$\!.
Under this identification,
linear operators on $L^2(\{0,1\})$ correspond to $2\times 2$ matrices with complex coefficients.
States can be identified with points in the closed unit ball of $\R^3$ as follows. Any selfadjoint operator $T = (t_{ij})_{i,j=1}^2$ on $\C^2$ with unit trace $\tr(T) = t_{11}+t_{22} = 1$ is of the form
\begin{align}\label{eq:T-qubit} T = \frac12\begin{pmatrix}
        1+c_3 & \, c_1-ic_2 \\ \, c_1+ic_2 & \, 1-c_3
       \end{pmatrix}
\end{align}
with $c_1,c_2,c_3\in \R$. The vector $c = (c_1,c_2,c_3) \in\R^3$ is called the {\em Bloch vector}\index{Bloch!vector} of $T$. It is easily checked that the eigenvalues of $T$ are $\frac12(1\pm |c|)$. From this we see that $T\ge 0$ if and only if $|c|\le 1$.

A routine computation shows that the pure state $\ket{h} = \cos(\theta/2)\ket{0} +e^{i\varphi}\sin(\theta/2)\ket{1}$
corresponds to the operator
\begin{align*}
T = h\,\bar\otimes\, h
  = \frac12\begin{pmatrix}
    1+ \cos \theta & \, e^{-i\varphi }\sin\theta \\ \, e^{i\varphi }\sin \theta & \, 1-\cos\theta
    \end{pmatrix}
\end{align*}
with Bloch vector $$(\sin\theta\cos\varphi, \,\sin\theta\sin\varphi, \,\cos\theta).$$
Thus the Bloch sphere representation of the pure state $\ket{h}$ coincides with the Bloch vector of the associated operator
$h\,\bar\otimes\, h$.

Equation \eqref{eq:T-qubit} can be written as
\begin{align*} T & = \frac12 \begin{pmatrix}\, 1 & \, 0\, \\ \,0 & 1  \end{pmatrix}
+ \frac{c_1}2  \begin{pmatrix}\, 0 & \, 1\, \\ 1 & 0  \end{pmatrix}
+ \frac{c_2}2  \begin{pmatrix}\, 0 & \,-i\, \\ i & 0  \end{pmatrix}
+ \frac{c_3}2  \begin{pmatrix}\, 1 & \, 0\, \\ 0 & -1 \end{pmatrix}
\\ & = \frac12(I + c_1\sigma_1 + c_2 \sigma_2+ c_3 \sigma_3),
\end{align*}
where
$$\sigma_1 = \begin{pmatrix}\, 0 & \, 1\, \\ 1 & 0  \end{pmatrix}, \qquad
  \sigma_2 = \begin{pmatrix}\, 0 & \,-i\, \\ i & 0  \end{pmatrix}, \qquad
  \sigma_3 = \begin{pmatrix}\, 1 & \, 0\, \\ 0 & -1 \end{pmatrix}$$
are the three {\em Pauli matrices}.\index{Pauli matrices}
These matrices are selfadjoint and their spectra equal $\{\pm 1\}$. Therefore they are associated with $\{\pm 1\}$-valued observables, also denoted by $\sigma_1$, $\sigma_2$, and $\sigma_3$. The corresponding eigenstates of $\sigma_j$ are called the {\em `spin up/spin down' states along the $j$th axis}.
Every selfadjoint operator $A$ on $\C^2$ is of the form
$$ A =   \begin{pmatrix}\, a & c-id\, \\ \,c+id & b  \end{pmatrix}
=  \begin{pmatrix}\, c_0+c_3 & c_1-ic_2\, \\ \,c_1+ic_2 & c_0-c_3  \end{pmatrix}
= c_0 I + c_1\sigma_1 + c_2 \sigma_2+ c_3 \sigma_3
$$
for certain $a,b,c,d\in \R$ and $c_0,c_1,c_2,c_3\in \R$ with $a = c_0+c_3$, $b = c_0-c_3$, $c = c_1$, and $d = c_2$. It follows that
the quadruple $\{I,\sigma_1,\sigma_2,\sigma_3\}$ is a basis of the real-linear vector space of selfadjoint operators on $\C^2$\!.

\subsection{Entanglement}\label{subsec:entaglement}

The natural choice for the state space of a system of $N$ classical point particles in $\R^3$ is $\R^{3N}\times \R^{3N}$, the idea being that six coordinates are needed (three for position, three for momentum) to describe the state of each particle.
In quantum mechanics, the natural choice of Hilbert space is $L^2(\R^{3N})$. Labelling the points of $\R^{3N}$ as $x = (x_j^{(n)})_{1\le j\le 3,\,1\le n\le N}$, this choice suggests the following natural definition of observables $\wh x_j^{(n)}$ describing the $j$th coordinate of the $n$th particle:
$$ \wh x_j^{(n)} f(x) := x_j^{(n)}f(x),\quad f\in \Dom(\wh x_j^{(n)}),\ x\in\R^{3N},$$
where $\Dom(\wh x_j^{(n)})= \{f\in L^2(\R^{3N}):\, x_j^{(n)} f\in L^2(\R^{3N})\}$.
Later we will see that the corresponding momentum operators are given by
$$  \wh p_j^{(n)} f := \frac1i \frac{\partial}{\partial x_j^{(n)}}f,\quad f\in \Dom(\wh p_j^{(n)})$$
with their natural domains.

The space $L^2(\R^{3N})$ is isometric in a natural way to the $N$-fold Hilbert space tensor product (see Definition \ref{def:Hilbert-tensor} and the discussion following it):
$$ L^2(\R^{3N}) \simeq \underbrace{L^2(\R^3)\ot\cdots\ot L^2(\R^3)}_{N \ {\rm times}}.$$
This suggests that if the Hilbert spaces $H_1,\hdots,H_N$ describe the states of $N$ quantum mechanical systems, then
their Hilbert space tensor product $$H_1\ot\cdots\ot H_N$$ serves to describe the system composed of these $N$ subsystems. In what follows we focus on the case $N=2$, but everything we say extends to general $N$ without difficulty.

Let $H$ and $K$ be separable Hilbert spaces.
A pure state in $H\otimes K$ is called a \emph{product state} if it
can be represented by a vector
$h\otimes k$
with $h\in H$ and $k\in K$ unit vectors.
In the physics literature, the associated product state $\ket{h\otimes k}$ is usually written as $|h\rb |k\rb$.

A pure state which is not a product state is called an {\em entangled state}.\index{entangled state}
For example, let $h_1,h_2\in H$ and $k_1,k_2\in K$ be orthonormal and let
$|c_1|^2+|c_2|^2=1$.
If both $c_1\not=0$ and $c_2\not=0$, then the pure state
$$\psi := \ket{c_1h_1\otimes k_1+c_2h_2\otimes k_2}$$ (in physics notation: $c_1\ket{h_1}\ket{k_1}+ c_2 \ket{h_2}\ket{k_2}$) is not a product state. Indeed, under the identification of
$\operatorname{span}\{h_1,h_2\}\otimes\operatorname{span}\{k_1,k_2\}$ with the space of
$2\times2$ matrices, product vectors correspond to matrices of rank
one, whereas the unit vector representing $\psi$ corresponds to the rank-2 matrix
$\operatorname{diag}(c_1,c_2)$.

The partial trace (see Section \ref{subsec:partialtrace}) can be used to define states of subsystems
starting from the state of a composite system. More concretely, suppose that $T\in \calL_1(H\ot K)$ is a positive trace class operator with unit trace.
By Proposition \ref{prop:part-trace-props}, $\tr_K(T)$ and $\tr_H(T)$ are positive trace class operators in
$\calL_1(H)$ and $\calL_1(K)$, respectively, each with unit trace. If we think of $T$ as describing the state of a system with Hilbert space $H\ot K$, $\tr_K(T)$ and $\tr_H(T)$ can be thought of as describing the states of the two constituent subsystems with Hilbert spaces $H$ and $K$, respectively. If the operator $T$ corresponds to a unit vector $h\ot k$ in $H\ot K$, the states corresponding to $\tr_K(T)$ and $\tr_H(T)$ are the pure states $\ket{h}$ and $\ket{k}$, or more precisely,
$$ \tr_K(T) = \ket{h}\bra{h}, \qquad  \tr_H(T) = \ket{k}\bra{k}.$$
More generally, if $T$ corresponds to
$c_1h_1\otimes k_1+c_2h_2\otimes k_2$, where $h_1,h_2$ and $k_1,k_2$ are orthonormal pairs and
$|c_1|^2+|c_2|^2=1$, then
\begin{align*}
\tr_K(T)
& = |c_1|^2\ket{h_1}\bra{h_1} + |c_2|^2\ket{h_2}\bra{h_2}, \\
\tr_H(T)
& = |c_1|^2\ket{k_1}\bra{k_1} + |c_2|^2\ket{k_2} \bra{k_2}.
\end{align*}
Thus an entangled pure state of the
composite system may induce mixed states on its subsystems.

\section{Positive Operator-Valued Measures}\label{sec:POVM}

We next enlarge the class of observables considered so far. A
projection-valued measure assigns to each measurable set $F$ of outcomes
an orthogonal projection $P_F$, and the probability of obtaining an outcome
in the set $F$ when the system is in the state $T$ equals
$\tr(P_FT)$.
This framework describes ideal measurements in which the events are
represented by orthogonal projections. More general measurement procedures,
including measurements affected by noise, coarse-graining, or an
indirect coupling to an auxiliary system, lead instead to operators
$Q_F$ satisfying
$$0\le Q_F\le I.$$
Such operators are called {\em effects} and are studied in the next subsection, and
their introduction leads to the notion of a {\em positive operator-valued measure}; projection-valued measures are recovered when every operator in their range is an orthogonal projection.

\subsection{Effects}

Let $(X,\X)$ be a measurable space, and let $B_{\rm b}(X)$ be the Banach space of bounded measurable functions on $X$.
As a warm-up we show:

\begin{proposition}\label{prop:convexhullBbM}
The closed convex hull in $B_{\rm b}(X)$ of the set of elementary observables
$\{\one_B:\,B\in \X\}$ equals
$$ \calE(X):= \{f\in B_{\rm b}(X):\, 0\le f \le \one \ \hbox{pointwise}\}.$$
The extreme points of $\calE(X)$ are precisely the elementary observables $\one_B$, $B\in \calX$.
\end{proposition}
\begin{proof}
Denote by $E(X)$ the closed convex hull of the set of elementary
observables. The inclusion $E(X)\subseteq\calE(X)$ is immediate.

Conversely, let $f\in\calE(X)$. For $n\ge1$ and
$k=1,\ldots,2^n$, put
$$
A_{k,n}:=\{x\in X:\, f(x) \ge k2^{-n}\}
$$
and
$$
f_n:=2^{-n}\sum_{k=1}^{2^n}\one_{A_{k,n}}.
$$
Each $f_n$ is a convex combination of elementary observables. Moreover,
$0\le f-f_n<2^{-n}$ pointwise, and hence
$\n f-f_n\n_\infty\le 2^{-n}$.
Thus $f_n\to f$ in $B_{\rm b}(X)$, and therefore $f\in E(X)$.

If $g\in \calE(X)$ is an elementary observable and
$g = \lambda f_0+(1-\lambda)f_1$ with $0< \la< 1$ and
$0\le f_j\le \one$ for $j=0,1$, then for every $\xi\in X$ with
$g(\xi)=0$ we have $f_0(\xi)=f_1(\xi)=0$, while for every
$\xi\in X$ with $g(\xi)=1$ we have $f_0(\xi)=f_1(\xi)=1$.
Thus $f_0=f_1=g$ pointwise.
It follows that every elementary observable is an extreme point
of $\calE(X)$.

If $g\in \calE(X)$ is not an elementary observable, then $g$ does not
take only the values $0$ and $1$, so there exists $\eps>0$ such that
the set $A:=\{\eps\le g\le 1-\eps\}$
is nonempty. The functions $f_0:=g+\eps\one_A$ and $f_1:=g-\eps\one_A$ are distinct elements of $\calE(X)$ and we have $g = \frac12f_0+\frac12 f_1$. It follows that $g$ is not an extreme point of $\calE(X)$.
\end{proof}

The quantum mechanical counterparts of the elementary observables are
the orthogonal projections. We now give an analogous characterisation of the closed convex hull of $\calP(H)$ in $\calL(H)$.
We write $S\le T$ to express that $T-S$ is a positive operator.

\begin{proposition}\label{prop:effect-convexhull} The closed convex hull in $\calL(H)$ of $\calP(H)$ equals\index{$E$@$\calE(H)$}
$$ \calE(H):= \{ E\in \calL(H): \, 0\le E\le I\}.$$
The extreme points of $\calE(H)$ are precisely the orthogonal projections.
\end{proposition}

\begin{proof}
Every element of the convex hull of $\calP(H)$ belongs to $\calE(H)$; this remains true for the closed convex hull. In the converse direction, let $E\in\calE(H)$, and let $P$ be its projection-valued
measure. For $n\ge1$ define
$$
f_n(t):=
2^{-n}\sum_{k=1}^{2^n}\one_{[k2^{-n},1]}(t),
\quad t\in[0,1],
$$
and put
$$
E_n:=f_n(E)
=
2^{-n}\sum_{k=1}^{2^n}P_{[k2^{-n},1]}.
$$
Each $E_n$ is a convex combination of orthogonal projections. Since
$$
0\le t-f_n(t)<2^{-n},
\quad t\in[0,1],
$$
the functional calculus gives, with ${\rm id}(t)=t$,
$$
\n E-E_n\n
\le
\n{\rm id}-f_n\n_\infty
\le2^{-n}.
$$
Thus $E_n\to E$ in operator norm, and $E$ belongs to the closed convex
hull of $\calP(H)$.

We next identify the extreme points of $\calE(H)$ as orthogonal projections.
If $E\in \calE(H)$ is an orthogonal projection and $E =   \lambda E_0+(1-\lambda)E_1$ with $0< \la< 1$ and $0\le E_j\le I$ for $j=0,1$, then for all $x\in \ker(E)$ we have $ \lambda \iprod{E_0 x}{x} +(1-\lambda)\iprod{E_1 x}{x}=0$
with $\iprod{E_i x}{x} \ge 0$ for $i=0,1$, and this is possible only if $\iprod{E_0 x}{x} = \iprod{E_1 x}{x} = 0$.
For all norm-one vectors $x\in \ran(E)$ we have $\iprod{Ex}{x} = \iprod{x}{x} = 1$ and consequently $\lambda \iprod{E_0 x}{x} +(1-\lambda)\iprod{E_1 x}{x} = 1$.
Since $\iprod{E_i x}{x} \le 1$ for $i=0,1$, this is possible only if $\iprod{E_0 x}{x} = \iprod{E_1 x}{x} = 1$.
It follows that $E_0=E_1 = 0$ on $\ker(E)$ and $E_0=E_1=I$ on $\ran(E)$,
and therefore $E_0=E_1=E$.
It follows that $E$ is an extreme point of $\calE(H)$.

Conversely, suppose that $E$ is not an orthogonal projection. Since $0\le E\le I$, we have
$\sigma(E)\subseteq[0,1]$. If $\sigma(E)\subseteq\{0,1\}$, then by the spectral
mapping theorem we have $E^2=E$, contradicting our assumption. Hence there exists
$\lambda\in\sigma(E)\cap(0,1)$.

Choose $\delta>0$ such that $[\lambda-\delta,\lambda+\delta]\subseteq(0,1)$.
Let $\psi\in C([0,1])$ be nonzero, with
$\operatorname{supp}\psi\subseteq[\lambda-\delta,\lambda+\delta]$ and
$$
0\le \psi(t)\le \min\{t,1-t\},\quad t\in[0,1],
$$
and such that $\psi(\lambda)\ne0$. Then $\psi(E)\ne0$ by the continuous functional calculus. Moreover, $0\le E-\psi(E)\le I$ and $0\le E+\psi(E)\le I$.
Thus $E-\psi(E)$ and $E+\psi(E)$ are effects, and
$$
E=\frac12\bigl(E-\psi(E)\bigr)+\frac12\bigl(E+\psi(E)\bigr)
$$
is a nontrivial convex decomposition. Hence $E$ is not an extreme point.
\end{proof}

\begin{definition}[Effects] An {\em effect}\index{effect}
is an element of the set $\calE(H)$.
\end{definition}

Effects are selfadjoint, and it follows from Theorem \ref{thm:spect-sa} that
a selfadjoint operator on $H$ is an effect if and only if
its spectrum is contained in the unit interval $[0,1]$.
If $T$ is an arbitrary nonzero positive operator, then for all $0\le c \le \n T\n^{-1}$ the operator $cT$ is an effect. Indeed, this is clear for $c=0$, and if $c>0$ the operator
$c^{-1} I - T$ is positive since it is selfadjoint and has nonnegative spectrum.

A mapping $\nu: \calE(H)\to [0,1]$ is said to be {\em finitely additive}\index{finite!additivity} if
$$\sum_{n=1}^N \nu(E_n) = \nu(E)$$ whenever $E_1,\dots,E_N, E\in \calE(H)$ satisfy $E_1+\cdots+E_N = E$.

\begin{theorem}[Busch]\index{theorem!Busch}\label{thm:Busch}
Every finitely additive mapping $\nu:\calE(H)\to [0,1]$ satisfying $\nu(I)=1$
restricts to an affine mapping $\nu:\calP(H)\to [0,1]$ satisfying
$\nu(0)=0$ and $\nu(I)=1$. If this restriction is countably additive, it
defines a state in the sense of Definition \ref{def:state}.
\end{theorem}

\begin{proof}
By assumption we have $\nu(I)=1$ and
from $1 = \nu(I) = \nu(I+0) = \nu(I)+\nu(0) = 1 + \nu(0)$ it follows that $\nu(0)=0$.
The function $f_P(c):=\nu(cP)$ is additive on $[0,1]$, bounded, and satisfies $f_P(0)=0$. Hence, by the lemma below, $\nu(cP) =c\nu(P)$ for all $c\in[0,1]$.
If $P_1,\dots,P_N\in\calP(H)$ and $c_1,\dots,c_N\in[0,1]$ with $\sum_{i=1}^N c_i=1$, then $E:=\sum_{i=1}^N c_i P_i$ belongs to $\calE(H)$ and, by finite additivity,
$$\nu(E)=\sum_{i=1}^N \nu(c_i P_i)=\sum_{i=1}^N c_i\nu(P_i).$$
Thus the restriction of $\nu$ to $\calP(H)$ is affine and satisfies
$\nu(0)=0$ and $\nu(I)=1$. If this restriction is countably additive, then by
Theorem \ref{thm:state-posfc} it defines a state.
\end{proof}

The proof of Busch's theorem is completed by the following real-variable lemma.

\begin{lemma}\label{lem:additive-bounded-linear}
Let $f:[0,1]\to[0,\infty)$ satisfy
$f(s+t)=f(s)+f(t)$
whenever $s,t\ge0$ and $s+t\le1$. Then
$$
f(t)=tf(1),\quad t\in[0,1].
$$
\end{lemma}

\begin{proof}
The function $f$ is increasing, since
$f(t)=f(s)+f(t-s)\ge f(s)$ whenever $0\le s\le t\le1$.
Additivity gives
$$
f\Bigl(\frac mn\Bigr)=\frac mnf(1),
\quad 0\le m\le n.
$$
Approximating an arbitrary $t\in[0,1]$ from below and above by rational
numbers and using monotonicity proves the result.
\end{proof}

\subsection{Positive Operator-Valued Measures}\label{subsec:POVM}

The next definition generalises the notion of a projection-valued measure by replacing orthogonal projections with effects.

\begin{definition}[Positive operator-valued measures]\label{def:POVM} A {\em positive operator-valued measure} (POVM)\index{POVM}
on a measurable space $(\Om,\calF)$
is a mapping $Q: \calF\to \calE(H)$ that assigns to every set $F\in \calF$ an effect $Q_F:=Q(F) \in \calE(H)$ with
the following properties:
\begin{enumerate}[label={\rm(\roman*)}, leftmargin=*]
 \item\label{it:POVM1} $Q_\Om = I$;
 \item\label{it:POVM2}  for all $x\in H$ the mapping $F\mapsto \iprod{Q_F x}{x}$
defines a measure on $(\Om,\calF)$.
\end{enumerate}
\end{definition}

The measure defined by \ref{it:POVM2} is denoted by $Q_{x}$. Thus, for all $F\in \calF$ and $x\in H$, by definition we have
$$ \iprod{Q_F x}{x}  = Q_{x}(F) =\int_{\Om} \one_F \ud Q_{x}.$$
It follows from
$ Q_x(\Om) = \iprod{Q_\Om x}{x} = \iprod{x}{x} =\n x\n^2$ that the measures $Q_x$ are finite.

Every projection-valued measure is a POVM. The following proposition gives a necessary and sufficient condition for a POVM to be projection-valued.

\begin{proposition}\label{prop:POVM-PVM} A POVM $Q: \calF\to \calE(H)$ is a projection-valued measure
 if and only if $Q_FQ_{F'} = Q_{F\cap F'}$ for all $F,F'\in \calF$\!.
\end{proposition}
\begin{proof}
The `only if' part has already been established in Section \ref{sec:PVM}. The `if' part is evident from $Q_F^2 = Q_{F\cap F} = Q_F$, which shows that each $Q_F$ is a projection. Since $Q_F$ is also positive, it is an orthogonal projection.
\end{proof}

POVMs retain the measure-theoretic features of projection-valued measures, but in general lose their `projection-theoretic' features. More precisely, one still has $Q_\emptyset=0$, monotonicity ($F\subseteq G$ implies $Q_F \le Q_G$), and strong countable additivity. For the latter, the proof for projection-valued measures
carries over, provided one replaces the identity $\n Px\n^2 = \iprod{Px}{x}$
for orthogonal projections $P$ by the inequality $\n Qx\n^2 \le \iprod{Qx}{x}$ for operators $0\le Q\le I$; this inequality follows from
$$
\|Qx\|^2 = \iprod{Qx}{Qx} = \iprod{Q Q^{1/2}x}{Q^{1/2} x} \le \iprod{Q^{1/2}x}{Q^{1/2} x} = \iprod{Qx}{x}.
$$
In particular, strong continuity from below and above
remain valid. The above proposition shows that the identity $Q_FQ_G=Q_{F\cap G}$ generally fails; in particular, for
disjoint sets $F,G$ one need not have $Q_FQ_G=0$. The support property carries over, provided one replaces the condition $x\in \Ran(P_F)$ by $Q_F x=x$.

A POVM which is not projection-valued is sometimes called an {\em unsharp observable}.\index{unsharp observable} An interesting example will be discussed in
Section \ref{subsec:number-phase}.

We have seen in Proposition \ref{prop:observ-mapping} that if $P:\calF\to \calP(H)$ is a projection-valued measure, then
for every $T\in \mathscr{S}(H)$ the mapping $$ F \mapsto \tr(P_F T), \quad F\in \calF\!,$$
is a probability measure on $(\Om,\calF)$. This sets up an affine mapping from the convex set $\mathscr{S}(H)$ to the convex set
$M_1^+(\Om)$\index{$M$@$M_1^+(\Om)$} of probability measures on $(\Om,\calF)$.
Inspection of this argument shows that it extends to POVMs $Q$; countable additivity
of $ F \mapsto \tr(Q_F T)$ follows from Proposition \ref{prop:normal-MCT}. The following theorem shows that, conversely, every affine mapping from $\mathscr{S}(H)$ to $M_1^+(\Om)$ arises from a unique POVM.

\begin{theorem}[POVMs as unsharp observables]\label{thm:POVM-unsharp} Let $(\Om,\calF)$ be a measurable space.
If $\Phi: \mathscr{S}(H) \to M_1^+(\Om)$ is an affine mapping, then there exists a unique POVM $Q:\calF\to\calE(H)$ such that
for all $T\in \mathscr{S}(H)$ we have
$$ (\Phi(T))(F) = \tr(Q_F T), \quad F\in \calF\!.$$
\end{theorem}

No continuity or normality assumption on $\Phi$ is needed: affinity alone suffices, because $\Phi$ takes values in probability measures, so countable additivity is built in.

\begin{proof}
Fix $F\in\calF$. For a nonzero positive operator
$T\in\calL_1(H)$ define
$$
\ell_F(T):=\tr(T)\,
\Bigl(\Phi\Bigl(\frac{T}{\tr(T)}\Bigr)\Bigr)(F),
$$
and put $\ell_F(0):=0$. This mapping is positively homogeneous
and satisfies
$$
0\le\ell_F(T)\le\tr(T),\quad 0\le T\in\calL_1(H).
$$
It is also additive. Indeed, if $S,T$ are nonzero and positive,
put $a:=\tr(S)$ and $b:=\tr(T)$. The identity
$$
\frac{S+T}{a+b}
=\frac{a}{a+b}\frac Sa+\frac{b}{a+b}\frac Tb
$$
and the affinity of $\Phi$ give
$\ell_F(S+T)=\ell_F(S)+\ell_F(T)$.
If either operator is zero, the same identity is immediate.

By additivity and positive homogeneity, $\ell_F$ is affine
on the cone of positive trace class operators and satisfies
$\ell_F(0)=0$.
By Lemma \ref{lem:affine-linear-extension}, $\ell_F$ extends
uniquely to a real-linear functional on the selfadjoint trace
class operators, and then by complexification to a linear
functional on $\calL_1(H)$.
For selfadjoint $T$, its positive and negative parts give
$$
|\ell_F(T)|
\le\ell_F(T_+)+\ell_F(T_-)
\le\tr(T_++T_-)=\n T\n_1.
$$
For arbitrary $T$, its real and imaginary parts each have trace
norm at most $\n T\n_1$, so $|\ell_F(T)|\le2\n T\n_1$.
Thus $\ell_F$ is bounded.

By trace duality (Theorem \ref{thm:traceduality}), there exists
a unique $Q_F\in\calL(H)$ such that
$$
\ell_F(T)=\tr(Q_FT),\quad T\in\calL_1(H).
$$
In particular, for $T\in\mathscr S(H)$ this gives the identity
asserted in the theorem. For every unit vector $h\in H$,
$$
\iprod{Q_Fh}{h}
=\ell_F(h\,\bar\otimes\,h)
=(\Phi(h\,\bar\otimes\,h))(F)\in[0,1].
$$
Hence $0\le Q_F\le I$.
These identities also show that $Q_\Om=I$ and that, for every
unit vector $h$, the mapping $F\mapsto\iprod{Q_Fh}{h}$ is a
probability measure. By scaling, it is a finite positive
measure for every $h\in H$. Thus $F\mapsto Q_F$ is a POVM.

Finally, the required representation determines
$\iprod{Q_Fh}{h}$ for every unit vector $h$ and every $F\in\calF$.
By scaling and polarisation, it therefore determines each $Q_F$
uniquely.
\end{proof}

\begin{remark}\label{rem:convex-reasonable}
 The assumption that $\Phi$ should be affine is a reasonable one in the light of the following argument. Suppose we have two quantum mechanical systems at our disposal, represented by the operators $T_1$ and $T_2$ in $\mathscr{S}(H)$ describing their states.
 We use a classical coin to decide which state is observed: with
probability $p$ we observe the system corresponding to $T_1$, and with
probability $1-p$ the system corresponding to $T_2$. This experiment can be described as observing the state corresponding to the convex combination $pT_1 + (1-p)T_2$. If $\Phi$ is the observable to be measured, we expect the probability distribution of the outcomes
 $\Phi(pT_1 + (1-p)T_2)$
 to be given by
 $p\Phi(T_1)+(1-p)\Phi(T_2)$.
\end{remark}

POVMs admit a bounded functional calculus, but an important difference from the bounded functional calculus for projection-valued measures of Theorem \ref{thm:Borel-FC} is that this calculus fails to be multiplicative
(see, however, \eqref{eq:mult-discalgebra}).

\begin{proposition}[Bounded functional calculus for POVMs]\label{prop:POVM-calculus} Let $(\Omega,\calF)$ be a measurable space and let $Q: \calF\to \calE(H)$ be a POVM.
There exists a unique linear mapping $\Psi:B_{\rm b}(\Om)\to \calL(H)$
satisfying
 $$ \Psi(\one_F) = Q_F, \quad F\in \calF\!,$$
  and
 $$ \n \Psi(f)\n \le \n f\n_\infty, \quad f\in B_{\rm b}(\Om).$$
 It satisfies $$\Psi(f)^\star = \Psi(\ov f), \quad f\in B_{\rm b}(\Om).$$
 Moreover,
$\Psi(\one)=I$,
and $\Psi(f)$ is positive whenever $f$ is nonnegative.
\end{proposition}

\begin{proof}
For $x,y\in H$ consider the complex measure $Q_{x,y}$ defined by
$$ Q_{x,y} (F):= \iprod{Q_F x}{y}, \quad F\in \calF\!.$$
A polarisation argument, using the countable additivity of the measures
$Q_x = Q_{x,x}$, shows that this indeed defines a measure.
For any measurable partition $\Om = F_1\cup\cdots\cup F_k$ we have, by the Cauchy--Schwarz inequality applied twice,
\begin{align*} \sum_{j=1}^k
|Q_{x,y}(F_j)|  = \sum_{j=1}^k |\iprod{Q_{F_j} x}{y}|
& \le \sum_{j=1}^k \iprod{Q_{F_j} x}{x}^{1/2}\iprod{Q_{F_j} y}{y}^{1/2}
\\ & =  \sum_{j=1}^k Q_x(F_j)^{1/2}Q_y(F_j)^{1/2}
\\ &  \le \Bigl(\sum_{j=1}^k Q_x(F_j)\Bigr)^{1/2} \Bigl(\sum_{j=1}^k Q_y(F_j)\Bigr)^{1/2}
\\ & = Q_x(\Om)^{1/2}Q_y(\Om)^{1/2}
  = \n x\n\n y\n,
\end{align*}
from which it follows that $Q_{x,y}$ has finite variation $|Q_{x,y}|(\Om) \le \n x\n\n y\n.$

For $f\in B_{\rm b}(\Om)$ define
 $$\aa_f(x,y):= \int_\Om f\ud Q_{x,y}, \quad x,y\in H.$$
The form $\aa_f$ is sesquilinear
and bounded and defines a bounded operator $\Psi(f)$ on $H$ by Proposition \ref{prop:bilin-T}.
It is clear that $\Psi(\one_F) = Q_F$ for all $F\in \calF$ and
$$ |\iprod{\Psi(f)x}{y}| = \Big|\int_\Om f\ud Q_{x,y}\Big| \le \int_\Om |f|\ud |Q_{x,y}| \le \n f\n_\infty\n x\n\n y\n.$$
In particular, $\Psi(\one)= Q_\Om = I$.
If $f$ is nonnegative,
then $$\iprod{\Psi(f)x}{x} = \aa_f(x,x) = \int_\Om f\ud Q_x \ge 0,$$ so $\Psi(f)$ is positive.
The identity $(\Psi(f))^\star = \Psi(\ov f)$ follows from
\begin{align*}
 \iprod{(\Psi(f))^\star x}{y} & =\iprod{x}{\Psi(f) y} = \ov{\iprod{\Psi(f) y}{x}}
 = \ov{\aa_f (y,x)} \\ & = \ov{\int_\Om f\ud Q_{y,x}} = \int_\Om \ov f \ud Q_{x,y} = \aa_{\ov f}(x,y) =\iprod{\Psi(\ov f)x}{y}.
\end{align*}
Since the simple functions are dense in $B_{\rm b}(\Om)$, uniqueness
follows by approximation from the identities $\Psi(\one_F)=Q_F$,
$F\in\calF$.
\end{proof}

For a POVM $Q$ on $(\Om,\calF)$ we write
$\int_\Om f\ud Q:=\Psi(f)$ for $f\in B_{\rm b}(\Om)$,
where $\Psi$ is the mapping of Proposition \ref{prop:POVM-calculus}.

\subsection{Naimark's Theorem}\label{subsec:Naimark}

If $J$ is an isometry from $H$ into another Hilbert space $\wt H$ and $\wt P$ is an orthogonal projection on $\wt H$,
then
$J^\star \wt PJ$ is an effect on $H$: for all $x\in H$ we have
$$0\le \iprod{\wt PJx}{Jx} = \n \wt P Jx\n^2 \le \n x\n^2 = \iprod{x}{x}$$ and therefore $0\le J^\star \wt PJ \le I$.
This gives a method of producing POVMs from projection-valued measures:

\begin{proposition}[Compression]\label{prop:POVM-compression}
Let $J$ be an isometry from $H$ into another Hilbert space $\wt H$, and let $(\Omega,\calF)$ be a measurable space. If
$\wt P: \calF\to \calP(\wt H)$ is a projection-valued measure, then
$Q:= J^\star \wt P J: \calF\to \calE(H)$ is a POVM.
\end{proposition}
\begin{proof}
By what we just observed, $Q$ maps each set $F\in\calF$ to an element
of $\calE(H)$.
It is clear that $Q_{\Om} = J^\star J = I$. To see that $Q$ is a POVM, it remains to observe that
for all $x\in H$ and $F\in \calF$ we have $$Q_x(F) = \iprod{Q_F x}{x} = \iprod{\wt P_FJx}{Jx} = \wt P_{Jx}(F),$$ from which it follows that $Q_x$ is a finite measure on $(\Om,\calF)$.
\end{proof}

The main result of this section is {\em Naimark's theorem}, which asserts that, conversely, every POVM arises in this way.

\begin{theorem}[Naimark]\label{thm:Naimark}\index{theorem!Naimark}\index{dilation!of a POVM}
Let $(\Om,\calF)$ be a measurable space and let
$Q:\calF\to \calE(H)$ be a POVM.
Then there exist a Hilbert space $\wt H$, a projection-valued measure
$\wt P:\calF\to\calP(\wt H)$, and an isometry $J:H\to\wt H$ such that
$$
Q_F=J^\star \wt P_FJ,\quad F\in\calF\!.
$$
Moreover, one may choose this dilation to be {\em minimal} in the sense that
$$
\wt H=\overline{\mathrm{span}}\{\wt P_FJh:\, F\in\calF\!, \, h\in H\},
$$
and any two such minimal dilations are unitarily equivalent.
\end{theorem}

\begin{proof}
Let $V$ be the vector space of all $H$-valued simple functions on
$(\Om,\calF)$. Given $f,g\in V$, choose a finite measurable partition
$(F_j)_{j=1}^N$ of $\Om$ such that
$f=\sum_{j=1}^N\one_{F_j}x_j$ and $g=\sum_{j=1}^N\one_{F_j}y_j$,
and put
$$
\iprod{f}{g}_Q
:=
\sum_{j=1}^N\iprod{Q_{F_j}x_j}{y_j}.
$$
Finite additivity of $Q$ shows that this definition is independent of
the chosen common partition.
From the definition
and the corresponding properties of the inner product on $H$, we then
find that $\iprod{\cdot}{\cdot}_Q$ is sesquilinear and satisfies
$\iprod{g}{f}_Q=\ov{\iprod{f}{g}_Q}$ for all $f,g\in V$.
Moreover, if
$f=\sum_{j=1}^N\one_{F_j}x_j$
with $(F_j)_{j=1}^N$ a finite measurable partition of $\Om$, then
$$
\iprod{f}{f}_Q
=
\sum_{j=1}^N\iprod{Q_{F_j}x_j}{x_j}
\ge0,
$$
since every $Q_{F_j}$ is positive.
The Cauchy--Schwarz inequality for such forms implies that
$$
N:=\{f\in V:\,\iprod{f}{f}_Q=0\}
$$
is a linear subspace and that $\iprod{f}{g}_Q=0$ whenever $f\in N$ and
$g\in V$. Hence the form induces an inner product on $V/N$. Let $\wt H$
be the completion of this inner product space.

Define
$$
Jx:=\one_\Om x,\quad x\in H.
$$
Then $\iprod{Jx}{Jy}_Q=\iprod{Q_\Om x}{y}=\iprod{x}{y}$,
so $J:H\to\wt H$ is an isometry.

Define
$$
\wt P_Af:=\one_Af, \quad A\in\calF\!.
$$
If $f=\sum_j\one_{F_j}x_j$, then
$$
\n\wt P_Af\n_Q^2
=
\sum_j\iprod{Q_{A\cap F_j}x_j}{x_j}
\le
\sum_j\iprod{Q_{F_j}x_j}{x_j}
=
\n f\n_Q^2.
$$
Thus $\wt P_A$ maps $N$ into $N$ and induces a contraction on $V/N$, which extends uniquely to a contraction on $\wt H$. It is routine to check that this extension, again denoted by $\wt P_A$, is
selfadjoint and satisfies $\wt P_A^2=\wt P_A$, so it is an orthogonal projection, and it satisfies $\wt P_A\wt P_B=\wt P_{A\cap B}$.

If $A=\bigcup_{n\ge1}A_n$ is a disjoint union and
$f=\sum_{j=1}^N\one_{F_j}x_j$, then
$$
\Bigl\|
\wt P_Af-\sum_{n=1}^m\wt P_{A_n}f
\Bigr\|_Q^2
=
\sum_{j=1}^N
\iprod{
Q_{(A\setminus\bigcup_{n=1}^mA_n)\cap F_j}x_j
}{x_j}
\to 0.
$$
This proves strong countable additivity on the dense subspace $V/N$.
Moreover,
$$ \n \wt P_A-\sum_{n=1}^m\wt P_{A_n}\n
=\n\wt P_{A\setminus\bigcup_{n=1}^mA_n}\n \le 1$$
on $V/N$ for all $m\ge 1$. The
convergence therefore extends from $V/N$ to all of $\wt H$. Hence
$\wt P:\calF\to\calP(\wt H)$
is a projection-valued measure.

Finally, for $A\in\calF$ and $x,y\in H$,
$$
\iprod{J^\star\wt P_AJx}{y}
=
\iprod{\one_Ax}{\one_\Om y}_Q
=
\iprod{Q_Ax}{y}.
$$
Thus $Q_A=J^\star\wt P_AJ$.

Every simple function is a finite sum of functions of the form
$\one_Ax=\wt P_AJx$. Consequently,
$$
\wt H
=
\ov{\operatorname{span}}\{\wt P_AJx:\, A\in\calF\!,\, x\in H\},
$$
so the dilation is minimal.

To prove the final statement,
let $(\wt H^{(1)},P^{(1)},J^{(1)})$ and $(\wt H^{(2)},P^{(2)},J^{(2)})$ be two
minimal Naimark dilations of the same POVM $Q:\calF\to\calE(H)$.
On the linear span of the vectors $P_A^{(1)}J^{(1)}h$ define
$$
U\Big(\sum_{k=1}^nP_{A_k}^{(1)}J^{(1)}h_k\Big)
:=
\sum_{k=1}^nP_{A_k}^{(2)}J^{(2)}h_k.
$$
This mapping is well defined and isometric. Indeed, for
$A,B\in\calF$ and $h,g\in H$ we have
\begin{align*}
\iprod{P_A^{(i)}J^{(i)}h}{P_B^{(i)}J^{(i)}g}
= \iprod{P_{A\cap B}^{(i)}J^{(i)}h}{J^{(i)}g}
= \iprod{Q_{A\cap B}h}{g}, \quad i=1,2.
\end{align*}
It follows by expanding inner products of finite sums that their
norms are preserved. By minimality, the domain of $U$ is dense in
$\wt H^{(1)}$ and its range is dense in $\wt H^{(2)}$, so $U$ extends
uniquely to a unitary operator from $\wt H^{(1)}$ onto $\wt H^{(2)}$.

Since $P_\Om^{(i)}=I$, we have $UJ^{(1)}=J^{(2)}$. Moreover, on the dense
linear span of the vectors $P_B^{(1)}J^{(1)}h$,
$$
UP_A^{(1)}P_B^{(1)}J^{(1)}h
=
P_{A\cap B}^{(2)}J^{(2)}h
=
P_A^{(2)}UP_B^{(1)}J^{(1)}h,
$$
and hence $UP_A^{(1)}=P_A^{(2)}U$.
Finally, these identities determine $U$ on the same dense linear
span. This proves the asserted uniqueness up to unitary equivalence.
\end{proof}

\begin{theorem}\label{thm:contr-POVM}\index{contraction!represented by a POVM} For every contraction $T\in\calL(H)$, there exists a unique POVM $Q:\calB(\T)\to\calE(H)$
 such that
 $$ T^n = \int_{\mathbb{T}} z^n \ud Q(z), \quad n\in\N.$$
If $Q$ is a POVM with the above property, then $T$ is unitary if and only if $Q$ is a projection-valued measure.
\end{theorem}
\begin{proof}
By the Sz.-Nagy dilation theorem (Theorem \ref{thm:Nagy}) there
exist a Hilbert space $\wt H$, a unitary operator $U\in\calL(\wt H)$, and an
isometry $J:H\to\wt H$ such that
$$
T^n=J^\star U^nJ,\quad n\in\N.
$$
Let $\wt P:\mathscr B(\mathbb T)\to\calP(\wt H)$ be the projection-valued
measure associated with $U$ by the spectral theorem. Define
$$
Q_B:=J^\star\wt P_BJ,\quad B\in\mathscr B(\mathbb T).
$$
Then $Q:\mathscr B(\mathbb T)\to\calE(H)$ is a POVM by Proposition
\ref{prop:POVM-compression}, and for all $n\in\N$ we have
$$
\int_{\mathbb T} z^n\ud Q(z)
= J^\star\Bigl(\int_{\mathbb T}z^n\ud\wt P(z)\Bigr)J
= J^\star U^nJ = T^n.
$$

To prove uniqueness, suppose that for all $x\in H$ we have
$$
  \iprod{T^n x}{x} = \int_{\mathbb{T}} z^n \ud Q_x(z) = \int_{\mathbb{T}} z^n \ud \wt Q_x(z), \quad n\in\N,
$$
where $\wt Q$ is another POVM on $\T$.
This means that the nonnegative Fourier coefficients of the finite measures
$Q_x$ and $\wt Q_x$ agree. By the injectivity of the Fourier--Stieltjes transform (Theorem \ref{thm:uniq-FT-T} and the observation following it) we conclude that $Q_x = \wt Q_x$ for all $x\in H$.
Thus $\iprod{Q_B x}{x} = \iprod{\wt Q_B x}{x}$ for all $B \in\calB(\mathbb T)$ and $x\in H$,
and this implies $Q_B=\wt Q_B$ for all $B \in\calB(\mathbb T)$.

For the final statement, if $Q$ is projection-valued, then its functional calculus is multiplicative and
$$
T^\star T
=
\Bigl(\int_\mathbb T\ov z\ud Q(z)\Bigr)
\Bigl(\int_\mathbb Tz\ud Q(z)\Bigr)
=
\int_\mathbb T|z|^2\ud Q(z)
=
I.
$$
Similarly, $TT^\star=I$, so $T$ is unitary. This proves the `if' part.

The `only if' part follows from uniqueness. Indeed, suppose that $\wt Q$ is a POVM with the above moment property. If $T$ is unitary, we may also represent $T$ in terms of its associated projection-valued measure $P$. Then
$$
T^n=\int_{\mathbb T}z^n\ud P(z),\quad n\in\N,
$$
and by uniqueness, $\wt Q=P$.
\end{proof}

If $Q$ is as in Theorem \ref{thm:contr-POVM}, for all polynomials $p\in C(\mathbb{T})$ of the form
$p(z) =\sum_{n=0}^N c_n z^n$ we have
$$\Psi(p) =  \int_{\mathbb{T}} p \ud Q = p(T),$$ where $ T = \int_{\mathbb{T}} z \ud Q(z)$.
By Proposition \ref{prop:POVM-calculus} and the maximum modulus
principle,
$$
\n p(T)\n=\n\Psi(p)\n
\le\sup_{z\in\mathbb T}|p(z)|
=\sup_{z\in\ov{\mathbb D}}|p(z)|.
$$
The contractive mapping $p\mapsto p(T)$ for polynomials $z\mapsto p(z)$ thus extends uniquely to a contractive functional calculus for the {\em disc algebra}\index{disc algebra}
 $A(\mathbb{D})$,\index{$A(\mathbb{D})$}
the Banach space of functions continuous on
$\ov{\mathbb D}$ and holomorphic on $\mathbb D$.
Under the identification of functions in $A(\mathbb D)$ with their restrictions to $\mathbb T$, these are precisely the functions belonging to the closure in $C(\mathbb{T})$ of the polynomials of the form just considered. For $f,g\in A(\mathbb D)$, choose polynomials $p_n\to f$
and $q_n\to g$ uniformly on $\ov{\mathbb D}$.
Since $\Psi(p_n)\Psi(q_n)=\Psi(p_nq_n)$, norm continuity
gives that the bounded functional calculus associated with this POVM is multiplicative on the disc algebra $A(\mathbb{D})$:
\begin{align}\label{eq:mult-discalgebra}\Psi(f) \Psi(g)=  \Psi(fg), \quad f,g\in A(\mathbb{D}).
\end{align}
This multiplicativity need not extend to $C(\mathbb T)$ in general.
Indeed, we have $\Psi(z)=T$, $\Psi(\ov z)=T^\star$, and
$\Psi(z\ov z)=I$, whereas
$$
\Psi(z)\Psi(\ov z)=TT^\star,\qquad
\Psi(\ov z)\Psi(z)=T^\star T.
$$
If $T$ is not unitary, at least one of these products differs
from $I$.

\subsection{The Phase--Number Pair}\label{subsec:number-phase}

A convenient (single-mode) model for the {\em number operator}, the selfadjoint
operator used in quantum optics to represent the photon-number observable,
can be given
on the Hardy space $H^2(\mathbb{D})$ considered in Section \ref{subsec:Toeplitz}. Recall that $H^2(\mathbb{D})$
is the Hilbert space of all holomorphic functions on $\mathbb{D}$ of the form $f(z) = \sum_{n\in\N} c_n z^n$ with
$$\n f\n^2 := \sum_{n\in\N} |c_n|^2 < \infty.$$
As we have seen in that section, the mapping $$\sum_{n\in\N} c_n z_n\mapsto \sum_{n\in\N} c_n e_n,$$
where $z_n(z):= z^n$ and $e_n(\theta):= e^{in\theta}$, defines an isometry from $H^2(\mathbb{D})$ onto the closed subspace of $L^2(\mathbb{T})$ consisting of all
functions whose negative Fourier coefficients vanish. This allows us to identify $H^2(\mathbb{D})$ with the range of the
Riesz projection
$ \sum_{n\in\Z} c_n e_n \mapsto \sum_{n\in\N} c_n e_n$ on $L^2(\mathbb{T})$.

In $H^2(\mathbb{D})$ we consider the unbounded selfadjoint operator $N$
with domain
$$\Dom(N) = \Bigl\{f = \sum_{n\in\N} c_n z_n\in H^2(\mathbb{D}):\, \sum_{n\in\N} n^2|c_n|^2 <\infty\Bigr\},$$
given by
$$ N z_n =
nz_n, \quad n\in\N.$$
The sequence $(z_n)_{n\in \N}$ is an orthonormal basis of eigenvectors of $N$ and accordingly we have $\N\subseteq\si(N)$.
On the other hand, if $\la\in\C\setminus\N$, then for every $f = \sum_{n\in\N}c_n z_n$ in $H^2(\mathbb{D})$ the equation $(\la - N )u = f$
is uniquely solved by $u = \sum_{n\in \N}\frac{c_n}{\la-n}z_n \in \Dom(N)$. This implies that $\la\in \varrho(N)$. We conclude that
\begin{align}\label{eq:specN} \sigma(N) = \N.\end{align}
This is a special case of Proposition \ref{prop:sigmaAdiag}, but the proof could be simplified here because we
have precise information about the domain of the operator.
We think of the eigenfunctions $z_n$ of $N$ as the pure states describing the $n$-photon states of an electromagnetic field. In this interpretation, \eqref{eq:specN} tells us that the number of photons observed is a nonnegative integer.

The projection-valued measure associated with $N$ (denoted by the same letter) is given by $N_{\{n\}} = \pi_n$, the orthogonal projection on $H^2(\mathbb{D})$
onto the one-dimensional subspace spanned by $z_n$, so that
$$ \iprod{N f}{f} = \int_{\N} n \ud N_f(n) = \sum_{n\in \N} n \iprod{N_{\{n\}}f}{f},
\quad f\in \Dom(N). $$

To define {\em phase}\index{phase} as a $\T$-valued unsharp observable in the sense of POVMs, we proceed as follows.
Let $S$ be the unilateral shift on $H^2(\mathbb{D})$, that is,
$$Sf(z)=zf(z), \quad f\in H^2(\mathbb D).$$
Under the identification of $H^2(\mathbb D)$ with the range of the Riesz projection on $L^2(\T)$, a unitary dilation of $S$ is given by
the operator of multiplication by $z$ on $L^2(\T)$, equivalently described as the bilateral shift $\wt S$ on $L^2(\T)$,
$$
 \wt S \sum_{n\in\Z} c_n e_n := \sum_{n\in \Z} c_n e_{n+1},
$$
with $e_n(\theta) = e^{in\theta}$ as before.
The projection-valued measure $P: \mathscr{B}(\T)\to \calP(L^2(\T))$ associated with $\wt S$ is easily checked to be given by
\begin{align}\label{eq:P-phase}
 P_B f = \one_B f, \quad B\in \mathscr{B}(\T), \ f\in L^2(\T).
\end{align}
Its compression to $H^2(\mathbb{D})$ is a POVM $\Phi :\mathscr{B}(\T)\to \calE(H^2(\mathbb{D}))$, which is called the {\em phase observable}. It satisfies
$$
S^n = \int_{\mathbb T} z^n\ud \Phi(z), \quad n\in\N.
$$
The covariance property expressed in the following theorem identifies the POVM $\Phi$
as the `complementary unsharp observable' to the number observable $N$. The unilateral shift $S$ is not unitary, and Theorem
\ref{thm:contr-POVM} therefore shows that $\Phi$ is not
projection-valued.

The notions of covariance and complementarity will be developed in more detail in the next chapter.

\begin{theorem}[Covariance of phase]\label{thm:quantum-phase}
The phase observable $\Phi$ is {\em covariant} with respect to the unitary $C_0$-group $U$ generated by $-iN$, that is,
\begin{align*}
U(t)\Phi_B U^\star(t) = \Phi_{e^{it}B}, \quad t\in \R, \ B\in\mathscr B(\mathbb T),
\end{align*}
where $e^{it}B = \{e^{it}z: \, z\in B\}$ is the rotation of $B$ through the angle $t$.
\end{theorem}
\begin{proof}
Let $J:H^2(\mathbb D)\to L^2(\T)$ be the canonical inclusion, so that
$\Phi_B=J^\star P_BJ$. On $L^2(\T)$ consider the unitary operators
$$
\wt U(t)f(z):=f(e^{-it}z),\quad t\in\R.
$$
Since $\wt U(t)e_n=e^{-int}e_n$ and $U(t)e_n=e^{-int}e_n$
for $n\ge0$, we have $\wt U(t)J=JU(t)$.

By \eqref{eq:P-phase}, conjugating $P_B$ by $\wt U(t)$ gives
multiplication by
$\one_B(e^{-it}z)=\one_{e^{it}B}(z)$. Thus
$$\wt U(t)P_B\wt U^\star(t)=P_{e^{it}B},$$
and compressing this identity to $H^2(\mathbb D)$ we obtain
\[
U(t)\Phi_BU^\star(t)
=J^\star\wt U(t)P_B\wt U^\star(t)J
=J^\star P_{e^{it}B}J
=\Phi_{e^{it}B}.
\qedhere\]
\end{proof}

\section{Markov Kernels and Effective Models}\label{sec:hidden}

In the preceding sections, observables were described by
projection-valued measures (`sharp observables') or, more generally,
by positive operator-valued measures (`unsharp observables'). The
classical analogue of a sharp observable is a measurable function from
the state space to the outcome space. We now argue that its unsharp
analogue is a Markov kernel, which assigns to each point of the classical state space a probability
measure on the outcome space. We then show that a family of POVMs can always be represented by a family of Markov kernels on a common classical state space.
This provides a natural framework for {\em effective models}\index{effective model} of quantum statistics, which we illustrate for the spin observables of a qubit and whose strengths and limitations we discuss.

\subsection{Markov Kernels as Classical Unsharp Observables}\label{subsec:Markov-kernels-POVMs}

We begin with the definition of a Markov kernel.

\begin{definition}[Markov kernels]\label{def:Markov-kernel}
Let $(X,\calX)$ and $(\Omega,\calF)$ be measurable spaces. A mapping
$K:X\times\calF\to[0,1]$
is called a \emph{Markov kernel}\index{Markov kernel} if:
\begin{enumerate}[label={\rm(\roman*)}, leftmargin=*]
\item for every $x\in X$ the set function $F\mapsto K(x,F)$ is a probability measure on $(\Omega,\calF)$;
\item for every $F\in\calF$ the function $x\mapsto K(x,F)$ is $\calX$-measurable.
\end{enumerate}
\end{definition}

Given a Markov kernel $K$ and a probability measure $\mu\in M_1^+(X)$, one defines the induced measure $K_\#\mu$ by
\begin{equation}\label{eq:kernel-pushforward}
(K_\#\mu)(F):=\int_X K(x,F)\ud\mu(x),\quad F\in\calF\!.
\end{equation}
Then $K_\#\mu$ is a probability measure. Indeed, $K(x,\Omega)=1$ for all $x\in X$, hence $(K_\#\mu)(\Omega)=\int_X 1\ud\mu=1$; countable additivity follows from the monotone convergence theorem. Moreover, the map $\mu\mapsto K_\#\mu$ is affine.

\begin{example}[Deterministic observables as Markov kernels]
If $f:X\to\Omega$ is measurable, then
$$
K_f(x,F):=\one_F(f(x)),\quad x\in X,
$$
defines a Markov kernel, and \eqref{eq:kernel-pushforward} reduces to
$(K_f)_\#\mu=f_\#\mu$, the image measure of $\mu$ under $f$. Thus Markov kernels extend the notion of classical
observables.
\end{example}

Recall from Proposition \ref{prop:convexhullBbM}
that the closed convex hull in $B_{\rm b}(X)$ of the set of elementary observables
$\{\one_B:\,B\in \X\}$ equals
the set $ \calE(X)$ of all real-valued bounded measurable functions $f:X\to\R$ such that $0\le f \le \one$ pointwise on $X$. We will prove three theorems providing justification for viewing Markov kernels as the classical analogue of POVMs.

We define a {\em classical positive operator-valued measure} as a mapping $q:\calF\to \calE(X)$ satisfying:
\begin{itemize}
 \item  $q_\Omega=\one$;
 \item for all $x\in X$ the mapping $F\mapsto q_F(x)$ defines a measure on $(\Omega,\calF)$.
\end{itemize}
Here we write $q_F:=q(F)$.
When $q_F$ is an elementary observable for every $F\in \calF$, we say that $q$ is a {\em classical projection-valued measure}.

\begin{theorem}[Markov kernels as classical POVMs]\label{thm:classical-POVM-kernels}
Let $(X,\calX)$ and $(\Omega,\calF)$ be measurable spaces. The set of Markov kernels
$K:X\times\calF\to[0,1]$ is in one-to-one correspondence with the set of classical
positive operator-valued measures $q:\calF\to\calE(X)$, the correspondence
being given by
$$
        q_F(x)=K(x,F),\quad x\in X,\ F\in\calF\!.
$$
\end{theorem}

\begin{proof}
If $K$ is a Markov kernel, then $q_F(x):=K(x,F)$ defines a function
$q_F\in\calE(X)$ for every $F\in\calF$, since $0\le K(x,F)\le 1$ for all
$x\in X$. Moreover, $q_\Omega=\one$, and for fixed $x\in X$ the mapping
$F\mapsto q_F(x)$ is the probability measure $F\mapsto K(x,F)$. Conversely,
if $q:\calF\to\calE(X)$ has these properties, then $K(x,F):=q_F(x)$ is a
Markov kernel: the measure property is built into the assumption, and the
measurability of $x\mapsto K(x,F)$ follows from $q_F\in B_{\rm b}(X)$.
\end{proof}

The correspondence also extends to bounded functional calculi.
A Markov kernel $K$ defines a unital positive linear mapping
$$
\Psi_K:B_{\rm b}(\Omega)\to B_{\rm b}(X),\qquad
(\Psi_K g)(x):=\int_\Omega g(\omega)\,K(x,\ud\omega).
$$
It satisfies
$$
\Psi_K(\one_F)=q_F,\qquad
\|\Psi_Kg\|_\infty\le\|g\|_\infty.
$$
Moreover, if $0\le g_n\uparrow g$ pointwise and
$g\in B_{\rm b}(\Omega)$, then
$\Psi_Kg_n\uparrow\Psi_Kg$ pointwise.

Conversely, every unital positive linear mapping
$\Psi:B_{\rm b}(\Omega)\to B_{\rm b}(X)$ with this monotone
convergence property is of the form $\Psi_K$ for a unique Markov
kernel $K$. Indeed, set
$$
K(x,F):=(\Psi\one_F)(x),\quad x\in X,\quad F\in\calF.
$$
Positivity and unitality give nonnegative finitely additive
probability measures, and the monotone convergence property gives
countable additivity. Measurability in $x$ follows from
$\Psi\one_F\in B_{\rm b}(X)$. The identity $\Psi=\Psi_K$ follows
first for simple functions and then by uniform approximation.

Thus both classical and quantum unsharp observables admit
unital positive integration maps. The classical map takes values
in the commutative algebra $B_{\rm b}(X)$, whereas the quantum
map takes values in $\calL(H)$. The relation with normal maps
between von Neumann algebras is discussed in the Notes.

\begin{theorem}[Classical POVMs versus classical PVMs]\label{thm:classical-PVM-POVM-kernels}
Let $(X,\calX)$ and $(\Omega,\calF)$ be measurable spaces, and assume that
 $\calF$ is countably generated. Let $K: X\times \calF \to [0,1]$ be a Markov kernel, and set $$q_F(x):= K(x,F), \quad x\in X, \ F\in \calF\!.$$
The following assertions are equivalent:
\begin{enumerate}[label={\rm(\arabic*)}, leftmargin=*]
\item\label{it:classical-PVM-POVM-kernels1} $K$ is represented by a measurable function $f:X\to\Omega$, in the sense that
$$K(x,F)=\one_F(f(x)), \quad x\in X, \ F\in \calF;$$

\item\label{it:classical-PVM-POVM-kernels4} for all $f,g\in B_{\rm b}(\Omega)$ we have $$\Psi_K(fg) = \Psi_K(f)\Psi_K(g);$$

\item\label{it:classical-PVM-POVM-kernels3} for all $F,G\in\calF$ we have $$q_Fq_G=q_{F\cap G};$$

\item\label{it:classical-PVM-POVM-kernels2} for all $F\in\calF$ the function $q_F$ takes values in $\{0,1\}$.
\end{enumerate}
In this case we have $q_F=\one_{f^{-1}(F)}$ for all $F\in\calF$.
\end{theorem}

Assertion~\ref{it:classical-PVM-POVM-kernels1} says that $K$
is deterministic. Since $\Psi_K$ is unital and preserves
complex conjugation,
assertion~\ref{it:classical-PVM-POVM-kernels4} says precisely
that $\Psi_K$ is a $\star$-homomorphism.
Assertions~\ref{it:classical-PVM-POVM-kernels3}
and~\ref{it:classical-PVM-POVM-kernels2} are the classical
counterparts of the equivalent conditions
$Q_FQ_G=Q_{F\cap G}$ and $Q_F\in\calP(H)$ for a POVM
to be projection-valued; see Proposition~\ref{prop:POVM-PVM}.

For the proof of the theorem we isolate a lemma.

\begin{lemma}\label{lem:01-measure-Dirac} Every $\{0,1\}$-valued probability measure on a countably generated measurable space $(\Omega,\calF)$ is a Dirac measure.
\end{lemma}

\begin{proof} Let $(A_n)_{n\ge1}$ generate $\calF$, and let $\nu:\calF\to \{0,1\}$ be a probability measure. Put
$$
B_n:=
\begin{cases} A_n,&\nu(A_n)=1,\\ \complement A_n,&\nu(A_n)=0. \end{cases}
$$
Then
$C:=\bigcap_{n\ge1}B_n$ satisfies $\nu(C)=1$ and is therefore nonempty. The collection of all $F\in\calF$ such that either $C\subseteq F$ or $C\cap F=\emptyset$ is a $\sigma$-algebra containing every $A_n$, and hence equals $\calF$. Choose $\omega\in C$. It follows that $\nu(F)=\one_F(\omega)$ for all $F\in\calF$.
\end{proof}

\begin{proof}[Proof of Theorem \ref{thm:classical-PVM-POVM-kernels}]
\ref{it:classical-PVM-POVM-kernels1}
$\Rightarrow$\ref{it:classical-PVM-POVM-kernels4}: \
If $K=K_f$, then $\Psi_Kg=g\circ f$ for every
$g\in B_{\rm b}(\Omega)$, so $\Psi_K$ is multiplicative.

\smallskip
\ref{it:classical-PVM-POVM-kernels4}
$\Rightarrow$\ref{it:classical-PVM-POVM-kernels3}: \
Since $\one_F\one_G=\one_{F\cap G}$, multiplicativity gives
$$
q_Fq_G
=\Psi_K(\one_F)\Psi_K(\one_G)
=\Psi_K(\one_{F\cap G})
=q_{F\cap G}.
$$

\smallskip
\ref{it:classical-PVM-POVM-kernels3}
$\Rightarrow$\ref{it:classical-PVM-POVM-kernels2}: \
Taking $F=G$ gives $q_F^2=q_F$. Hence $q_F$ takes only the values
$0$ and $1$.

 \smallskip
 \ref{it:classical-PVM-POVM-kernels2}
 $\Rightarrow$\ref{it:classical-PVM-POVM-kernels1}: \
For $x\in X$, define the $\{0,1\}$-valued probability measure
$$
\nu_x(F):=q_F(x),
\quad F\in\calF\!.
$$
By Lemma~\ref{lem:01-measure-Dirac}, choose $f(x)\in\Omega$ such that
$$
q_F(x)=\one_F(f(x)),
\quad F\in\calF\!.
$$
The function $f$ thus defined is measurable, since
$f^{-1}(F)=\{x\in X:\,q_F(x)=1\}$ belongs to $\calX$.
Thus $K=K_f$, and $q_F=\one_{f^{-1}(F)}$ for all $F\in\calF$.
\end{proof}

Every Markov kernel therefore induces an affine mapping
$K_\#:M_1^+(X)\to M_1^+(\Omega)$.
Conversely, Problem~\ref{prob:Markov} shows that, for complete
separable metric spaces, every weakly continuous affine mapping
between these spaces of probability measures is induced by a
unique Markov kernel. This continuity assumption is sufficient
for the converse; an arbitrary Markov kernel need not induce
a weakly continuous mapping.

By Proposition \ref{prop:POVM-PVM}, a POVM $Q: \calF\to \calE(H)$ is a projection-valued measure
if and only if $Q_FQ_{F'} =  Q_{F'}Q_F = Q_{F\cap F'}$ for all $F,F'\in \calF$.
The final theorem of this subsection represents POVMs on the real line which only satisfy the commutation condition $Q_FQ_{F'} = Q_{F'}Q_F$ for all $F,F'\in \mathscr{B}(\R)$ as
`smearings', by a Markov kernel, of a projection-valued measure on an auxiliary compact space.

\begin{theorem}[Commutative POVMs as smearings of PVMs]
\label{thm:commutative-POVM-kernel}
Let $Q:\calB(\R)\to\calE(H)$ be a POVM on a separable Hilbert space $H$.
The following assertions are equivalent:
\begin{enumerate}[label={\rm(\arabic*)}, leftmargin=*]
\item\label{it:commutative-POVM-kernel1}
$Q$ is commutative, that is, $Q_FQ_G=Q_GQ_F$ for all $F,G\in \calB(\R)$;
\item\label{it:commutative-POVM-kernel2}
There exist a compact set $X\subseteq[0,1]$, a projection-valued
measure $P:\calB(X)\to\mathscr P(H)$,
and a Markov kernel $K:X\times\calB(\R)\to[0,1]$
such that
\begin{equation}\label{eq:commutative-POVM-smearing}
Q_F=\int_X K(x,F)\ud P(x),\quad F\in\calB(\R).
\end{equation}
\end{enumerate}
Moreover, $P$ may be chosen so that its range generates the same
von Neumann algebra as the range of $Q$.
\end{theorem}

\begin{proof}
Suppose first that \ref{it:commutative-POVM-kernel2} holds. For all
$F\in\calB(\R)$, the operator $Q_F$ belongs to the bounded functional
calculus of $P$. Hence for all $F,G\in\calB(\R)$, by multiplicativity of the bounded functional calculus,
$$
Q_FQ_G
=
\int_X K(x,F)K(x,G)\ud P(x)
=
Q_GQ_F.
$$
Thus \ref{it:commutative-POVM-kernel1} holds.

Conversely, suppose that $Q$ is commutative, and let $\mathscr M$
be the von Neumann algebra generated by its range. Since the
effects $Q_B$, $B\in\calB(\R)$, commute, $\mathscr M$ is commutative.
By Lemma \ref{lem:comm-vN}, there exists $S\in\mathscr M$ with
$0\le S\le I$ such that $\mathscr M=\{S\}''$.

Let $P$ be the projection-valued measure of $S$ and put
$X:=\sigma(S)\subseteq[0,1]$. For $r\in\Q$ put
$A_r:=Q_{(-\infty,r]}$. Then $A_r\in \calM$, and by Theorem \ref{thm:bicomm-singleT},
there exists a bounded Borel function $f_r$ on $X$ such that
$$
A_r=f_r(S)=\int_X f_r\ud P.
$$
Since $0\le A_r\le I$, we may choose the functions $f_r$ with
$0\le f_r\le1$ outside a $P$-null set. Moreover, if $r<s$, then
$A_r\le A_s$, and hence $f_r\le f_s$ outside a $P$-null set.
Since there are only countably many rational pairs $r<s$, the
functions may be modified on one common Borel $P$-null set so that
these inequalities hold everywhere outside that set.

Strong continuity of the POVM from above and below
allows us, after enlarging the same null set if necessary, to assume
also that outside this null set we have
$$
f_r(x)=\inf_{\substack{s\in\Q\\s>r}}f_s(x),
\qquad
\lim_{r\to-\infty}f_r(x)=0,
\qquad
\lim_{r\to+\infty}f_r(x)=1.
$$
Fix a probability measure $\nu_0$ on $\R$. On the exceptional Borel
$P$-null set $N$, redefine, simultaneously for all $r\in\Q$,
$$
f_r(x):=\nu_0((-\infty,r]).
$$
The functions $f_r$ remain Borel representatives of the same elements
of $L^\infty(X,P)$, and the preceding monotonicity, right-continuity,
and limit properties now hold for every $x\in X$.

For $x\in X$ define
$$
F_x(t):=\inf_{\substack{r\in\Q\\r>t}}f_r(x),\quad t\in\R.
$$
Then $F_x$ is a distribution function of a unique
probability measure $K(x,\cdot)$ on $\R$ with
$K(x,(-\infty,t])=F_x(t)$. For every $t\in\R$, the map
$$
x\mapsto K(x,(-\infty,t])
=\inf_{\substack{r\in\Q\\r>t}}f_r(x)
$$
is Borel measurable. The sets $B\in\calB(\R)$ for which
$x\mapsto K(x,B)$ is Borel measurable form a Dynkin system containing
the half-lines and $\R$. Thus $K$ is a Markov kernel.

For $r\in\Q$ we have, for all $x\in X$,
$$
K(x,(-\infty,r])=f_r(x),
$$
and therefore
$$
 Q_{(-\infty,r]}
=
\int_X K(x,(-\infty,r])\ud P(x).
$$
The collection of all $B\in\calB(\R)$ for which
$$
 Q_B=\int_XK(x,B)\ud P(x)
$$
is again a Dynkin system containing the rational half-lines. Thus
the identity holds for every $B\in\calB(\R)$.
This proves \eqref{eq:commutative-POVM-smearing}.

By construction, the range of $P$ generates
$\{S\}''=\mathscr M$, the von Neumann algebra generated by the
range of $Q$.
\end{proof}

\subsection{Effective Models for Families of POVMs}
\label{subsec:effective-models}

The preceding discussion explained why Markov kernels may be viewed as the classical counterparts of POVMs.
We now explain how, starting from a family of quantum unsharp
observables, one obtains a common classical state space on which the
corresponding outcome statistics are described by Markov kernels. In the language of the
foundations of quantum mechanics, its points are the `effective states' discernible by the given set of observables.

Let $Q^{(i)}:\calF_i\to\calE(H)$, $i\in I$, be a family of POVMs on
measurable outcome spaces $(\Omega_i,\calF_i)$. Thus each $Q^{(i)}$
defines an affine observable
$$
        \Phi_{Q^{(i)}}:\mathscr S(H)\to M_1^+(\Omega_i),
        \qquad
        \Phi_{Q^{(i)}}(T)(F):=\tr(Q_F^{(i)}T),
$$
where $T\in\mathscr S(H)$ and $F\in\calF_i$.

Let
$$
        \mathscr P_1(H):=\{ h\,\bar\otimes\, h :\, h\in H,\, \|h\|=1\}
$$
be the set of pure normal states, identified with the rank-one orthogonal
projections on $H$. In physics notation, $h\,\bar\otimes\, h = \ket{h}\bra{h}$. Define an equivalence relation $\sim$ on $\mathscr P_1(H)$ by
$$
  P\sim R \ \Longleftrightarrow \ \tr(Q_F^{(i)}P)=\tr(Q_F^{(i)}R)
    \ \hbox{ for all } i\in I \hbox{ and } F\in\calF_i.
$$
Thus two pure states are equivalent if none of the given POVMs
distinguishes them statistically. We now set
$$
        X:=\mathscr P_1(H)/\!\sim
$$
and denote the quotient map by $q:\mathscr P_1(H)\to X$.
On $X$ we consider the
$\sigma$-algebra $\calX$ generated by the functions $k_F^{(i)}: X\to[0,1]$, $i\in I$, $F\in\calF_i$, defined by
$$
        k_F^{(i)}(q(P)):=\tr(Q_F^{(i)}P).
$$
These functions are well defined by the definition of $\sim$.

The following theorem, which is a variant of a classical representation
theorem for statistical models due to Holevo, associates a measurable
space $X$ and a family of Markov kernels $K^{(i)}:X\times\calF_i\to[0,1]$ with a family of POVMs.
One may think of the points $x\in X$ as {\em effective states}: they are not directly observed, but they
determine the probability laws of the observed measurement outcomes through the {\em response kernels}
$K^{(i)}$. A probability measure $\mu_T$ on $X$ representing a state $T$ describes the distribution of the effective states in a preparation of $T$. A triple consisting of such a classical state space, response kernels, and representing measures is called an {\em effective model}\index{effective model}.

\begin{theorem}[Effective model] \label{thm:effective-model}
Let $Q^{(i)}:\calF_i\to\calE(H)$, $i\in I$, be a family of POVMs, and let $(X,\calX)$ be the measurable quotient space constructed above. For each $i\in I$, the formula
$$ K^{(i)}(x,F):=k_F^{(i)}(x), \quad x\in X,\ F\in\calF_i, $$
defines a Markov kernel.
Moreover, let $T\in\mathscr S(H)$ and choose a spectral decomposition
$T=\sum_{n\ge1}\lambda_nP_n$ with each $\lambda_n\ge0$, $\sum_{n\ge1}\lambda_n=1$, and each $P_n$ a rank-one projection.
Then the probability measure
$$ \mu_T:=\sum_{n\ge1}\lambda_n\delta_{q(P_n)} $$
on $X$ satisfies
\begin{align}\label{eq:muTrepresenting} \tr(Q_F^{(i)}T) = \int_XK^{(i)}(x,F)\ud\mu_T(x), \quad i\in I,\ F\in\calF_i.
\end{align}
\end{theorem}

\begin{proof} For fixed $i\in I$ and $F\in\calF_i$, the function $x\mapsto K^{(i)}(x,F)$ is measurable by the definition of $\calX$. If $x=q(P)$ with $P\in\mathscr P_1(H)$, then $ K^{(i)}(x,F)=\tr(Q_F^{(i)}P)$.
Thus, for fixed $x$, the mapping $F\mapsto K^{(i)}(x,F)$ is a probability measure. This is independent of the chosen representative $P$ by the definition of the equivalence relation. Hence $K^{(i)}$ is a Markov kernel.

Finally, let $T = \sum_{n\ge 1}\lambda_nP_n$ as in the theorem, and put $T_N: = \sum_{n=1}^N\lambda_nP_n$.
Then $T_N\to T$ in the trace norm, and $
\tr(Q_F^{(i)}T_N)\to\tr(Q_F^{(i)}T)
$ since $Q_F^{(i)}$ is bounded. Hence,
\[
\int_XK^{(i)}(x,F)\ud\mu_T(x) = \sum_{n\ge1}\lambda_nK^{(i)}(q(P_n),F)
= \sum_{n\ge1}\lambda_n\tr(Q_F^{(i)}P_n)= \tr(Q_F^{(i)}T).
\]
\end{proof}

The representing measure $\mu_T$ constructed above may depend on the
chosen spectral decomposition of $T$; it will be shown below that this is the case for the spin observables of the qubit.

\begin{proposition}
An affine assignment $T\mapsto\mu_T$ of measures such that \eqref{eq:muTrepresenting} holds
exists if and only if there is a POVM $G$ on $X$ generating all the
$Q^{(i)}$ in the sense that
\begin{align}\label{eq:POVMs-oneG}
Q_F^{(i)}
=\int_X K^{(i)}(x,F)\ud G(x),
\quad i\in I,\ F\in\calF_i.
\end{align}
\end{proposition}

If no such POVM $G$ exists, then this
effective model does not admit an affine assignment of representing measures.

\begin{proof}
Suppose first that $T\mapsto\mu_T$ is affine and satisfies
\eqref{eq:muTrepresenting}. By Theorem \ref{thm:POVM-unsharp}, there
exists a POVM $G$ on $X$ such that
$\mu_T(B)=\tr(G_BT)$ for all $T\in\mathscr S(H)$ and $B\in\calX$.
This identity and the bounded POVM functional calculus give
$$
\int_X f(x)\ud\mu_T(x)
=
\tr\Bigl(\Bigl(\int_X f(x)\ud G(x)\Bigr)T\Bigr)
$$
for every bounded measurable function $f$ on $X$. Taking
$f(x)=K^{(i)}(x,F)$ and using \eqref{eq:muTrepresenting}, we obtain
$$
\tr(Q_F^{(i)}T)
=
\tr\Bigl(\Bigl(\int_X K^{(i)}(x,F)\ud G(x)\Bigr)T\Bigr)
$$
for all $T\in\mathscr S(H)$. Since states separate selfadjoint
operators, this gives \eqref{eq:POVMs-oneG}.

Conversely, suppose that a POVM $G$ satisfying \eqref{eq:POVMs-oneG} exists and define $$ \mu_T(B):=\tr(G_BT), \quad B\in\calX. $$ Then $\mu_T$ is a probability measure and the mapping $T\mapsto\mu_T$ is affine. Again by the bounded POVM functional calculus,
$$ \int_X K^{(i)}(x,F)\ud\mu_T(x) = \tr\Bigl(\Bigl(\int_XK^{(i)}(x,F)\ud G(x)\Bigr)T\Bigr) = \tr(Q_F^{(i)}T), $$
so $\mu_T$ is a representing measure for every $T\in\mathscr S(H)$.
\end{proof}

\subsection{The Qubit Revisited}\label{subsec:qubit-effective}

We now apply the construction of Theorem \ref{thm:effective-model} to the simplest genuinely quantum case, namely, the Pauli
observables $\sigma_1, \sigma_2, \sigma_3$ of the qubit.
The example makes explicit that the Pauli observables, which are projection-valued on the quantum side, are represented on the space of effective states by classical unsharp observables, that is, by Markov kernels which are not deterministic.

It is readily checked that for $p\in\{\pm1\}$, the POVM $P_j$ of $\sigma_j$ is given by
\begin{align}\label{eq:sp-proj-Pauli}
   (P_j)_{\{p\}}=\frac12(I+p\sigma_j), \quad p = \pm 1, \ j=1,2,3.
\end{align}
On the other hand, parametrising the pure states of $\mathscr S(\C^2)$ by the Bloch sphere $S^2\subseteq\R^3$, the point
$\xi=(\xi_1,\xi_2,\xi_3)\in S^2$ corresponds to the rank-one orthogonal projection
$$
   T_\xi=\frac12(I+\xi_1\sigma_1+\xi_2\sigma_2+\xi_3\sigma_3).
$$
Using the identities $\tr(\sigma_i)=0$ and
$\tr(\sigma_i\sigma_k)=2\delta_{ik}$, for $j=1,2,3$ we obtain
\begin{align*}
   \tr((P_j)_{\{p\}}T_\xi)
   &=
   \frac14
   \tr\Bigl((I+p\sigma_j)
   \Bigl(I+\sum_{k=1}^3\xi_k\sigma_k\Bigr)\Bigr)       \\
   &=
   \frac14
   \Bigl(
      \tr(I)
      +p\sum_{k=1}^3\xi_k\tr(\sigma_j\sigma_k)
   \Bigr)
   =
   \frac14(2+2p\xi_j)
    =
   \frac12(1+p\xi_j).
\end{align*}
Thus the probability of obtaining the value $p$ when measuring
$\sigma_j$ in the pure state $T_\xi$ is $\frac12(1+p\xi_j)$.
In the language of Markov kernels, this shows that the kernel
corresponding to the observable $\sigma_j$ is given by
$$
   K_j(\xi,A) =
   \sum_{p\in A}\frac12(1+p\xi_j),
   \quad
   \xi\in S^2\!,\ A\subseteq\{\pm1\}.
$$

Since the three coordinate functions $\xi\mapsto \xi_j$, $j=1,2,3$,
separate the points of $S^2$, the quotient space $X$ of Theorem
\ref{thm:effective-model} may in this case be identified with $S^2$. Under this identification, the $\sigma$-algebra $\calX$ is the Borel
$\sigma$-algebra of $S^2$, since it is generated by the three
coordinate functions.

For a probability measure $\mu$ on $S^2$, let
$$
m:=\int_{S^2}\xi\ud\mu(\xi)\in\R^3
$$
be its barycentre, and put
$$
T_\mu:=\frac12\Bigl(I+\sum_{k=1}^3m_k\sigma_k\Bigr).
$$
Since $|m|\le1$, the Bloch vector characterisation shows that
$T_\mu$ defines a state. For $j=1,2,3$ and $p\in\{-1,1\}$ we have
\begin{align*}
((K_j)_\#\mu)(\{p\})
= \int_{S^2}K_j(\xi,\{p\})\ud\mu(\xi)
& = \frac12\Bigl(1+p\int_{S^2}\xi_j\ud\mu(\xi)\Bigr)
\\ & = \frac12(1+pm_j)
= \tr((P_j)_{\{p\}}T_\mu).
\end{align*}
Thus the classical mixture $\mu$ reproduces the outcome distributions
of the three Pauli observables in the state $T_\mu$.

For this family of three Pauli observables, there is no affine
assignment $T\mapsto\mu_T$ of representing measures on $S^2$. For $c\in\R^3$ with $|c|\le 1$,
a probability measure $\mu$ on $S^2$ represents the state
$T_c=\frac12(I+c\cdot\sigma)$
if and only if
$$
\int_{S^2}\xi\ud\mu(\xi)=c.
$$
If $|c|=1$, then
$$
\int_{S^2}|\xi-c|^2\ud\mu(\xi)
=
2-2c\cdot\int_{S^2}\xi\ud\mu(\xi)
=
2-2|c|^2
=
0.
$$
Since the integrand is nonnegative, $\xi=c$ for $\mu$-almost every $\xi$, and hence $\mu=\delta_c$.
Suppose that the representing measures could be chosen affinely.
Since $T_0=\frac12T_\xi+\frac12T_{-\xi}$ for all $\xi\in S^2$,
we would have
$$
\mu_0=\frac12\delta_\xi+\frac12\delta_{-\xi}
$$
for every $\xi\in S^2$, which is impossible. Thus this effective model admits no affine assignment of representing measures.

\section{Classical and Quantum Correlations}\label{sec:Grothendieck}

The effective model of the preceding section reproduces the outcome distribution of each given unsharp observable by a Markov kernel. We now investigate whether the pair distributions obtained from
measurements on two subsystems can arise as marginals of a single probability distribution for all the individual outcomes. This is a special case of the problem of prescribed marginals discussed in Subsection \ref{subsec:marginals}.

We restrict attention to measurements with outcomes $\pm1$. Their correlations satisfy certain inequalities whenever a common joint distribution exists. Quantum measurements can violate these
inequalities, as we will see in a two-qubit example attaining the Tsirel'son bound. We then prove Grothendieck's inequality, which gives a universal factor comparing quantum and classical correlation bounds, independently of the number of observables and the dimensions of the Hilbert spaces.

\subsection{Classical Correlations}\label{subsec:epi-classical-correlations}

Fix integers $m,n\ge1$. A joint distribution of $m$ binary outcomes for the first subsystem and $n$ binary outcomes for the second is a probability measure $\mu$ on $\{-1,1\}^{m+n}$. Writing its coordinates as $(\varepsilon,\delta)$, we call the matrix $E^\mu=(E^\mu_{ij})_{ij}$ defined by
\begin{align}\label{eq:Eij}
E^\mu_{ij}:=\int_{\{-1,1\}^{m+n}}\varepsilon_i\delta_j
       \ud\mu(\varepsilon,\delta),
\quad i=1,\ldots,m,\ j=1,\ldots,n,
\end{align}
a {\em classical correlation matrix}.\index{correlation matrix} Here and below, integrals over this finite set are simply finite sums.

Let us write $\mathscr C_{m,n}$ for the set of classical correlation matrices,
$$
\mathscr C_{m,n}:=
\operatorname{conv}
\{(\varepsilon_i\delta_j)_{ij}:\,
  (\varepsilon,\delta)\in\{-1,1\}^{m+n}\}.
$$
The set $\mathscr C_{m,n}$ is compact and convex, and symmetric about $0$.

For a real matrix $M=(m_{ij})_{ij}\in\R^{m\times n}$, define its {\em classical correlation bound} by
\begin{align*}
\beta(M):=
\max_{(\varepsilon,\delta)\in\{-1,1\}^{m+n}}
\Bigl|\sum_{i=1}^m\sum_{j=1}^n m_{ij}\varepsilon_i\delta_j\Bigr|.
\end{align*}
If $E\in\mathscr C_{m,n}$, say $E= E^\mu$, then by \eqref{eq:Eij} we have
$$
\Bigl|\sum_{i=1}^m\sum_{j=1}^n m_{ij}E_{ij}\Bigr|\le\beta(M).$$
Conversely, the validity of these inequalities for every $M\in \R^{m\times n}$ characterises the classical correlation matrices.

\begin{proposition}\label{prop:epi-correlation-duality}
A real matrix $E\in\R^{m\times n}$ belongs to $\mathscr C_{m,n}$
if and only if
$$
\Bigl|\sum_{i,j}m_{ij}E_{ij}\Bigr|\le\beta(M),
\quad M\in\R^{m\times n}.
$$
\end{proposition}

\begin{proof}
Only sufficiency remains to be proved. If $E\notin\mathscr C_{m,n}$, the Hahn--Banach separation theorem gives a real linear functional $\ell$ on $\R^{m\times n}$ such that
$$
\ell(E)>\max_{C\in\mathscr C_{m,n}}\ell(C).
$$
Write $\ell(A)=\sum_{i,j}m_{ij}A_{ij}$ for $A \in \R^{m\times n}$. A linear functional takes
its maximum over the convex hull of a finite set at one of the
points of that set. Since $\mathscr C_{m,n}$ is symmetric, every $C\in\mathscr C_{m,n}$
is accompanied by $-C\in\mathscr C_{m,n}$, and
$\ell(-C)=-\ell(C)$. Thus the maximum of $\ell$ equals the
maximum of its absolute value:
$$\sum_{i,j}m_{ij}E_{ij} = \ell(E)>
\max_{C\in\mathscr C_{m,n}}\ell(C)
=\max_{C\in\mathscr C_{m,n}}|\ell(C)|
=\beta(M).
$$
This contradiction completes the proof.
\end{proof}

For the matrix
$M_{\rm CHSH}:=\begin{pmatrix}1&1\\1&-1\end{pmatrix}$,
the classical bound is particularly simple. For any four signs,
$$
\bigl|\varepsilon_1(\delta_1+\delta_2)
      +\varepsilon_2(\delta_1-\delta_2)\bigr|=2,
$$
because exactly one of $\delta_1+\delta_2$ and $\delta_1-\delta_2$
is zero, while the other has absolute value $2$. Thus
$\beta(M_{\rm CHSH})=2$, and every classical correlation matrix
$E\in\mathscr C_{2,2}$ satisfies the CHSH inequality
of Subsection \ref{subsec:marginals},
\begin{align}\label{eq:epi-classical-CHSH}
|E_{11}+E_{12}+E_{21}-E_{22}|\le 2.
\end{align}

\subsection{Quantum Correlations}\label{subsec:epi-quantum-correlations}

Let $H_A,H_B$ be nonzero separable complex Hilbert spaces and put
$H:=H_A\otimes H_B$. We consider a state $T\in\mathscr S(H)$ and
two finite families of selfadjoint contractions
$$
A_1,\ldots,A_m\in\calL(H_A),\qquad
B_1,\ldots,B_n\in\calL(H_B).
$$
Each $A_i$ determines a binary POVM with effects
$$
(E_{i})_{\{p\}}:=\frac12(I+pA_i),\quad p\in\{-1,1\}.
$$
Indeed, both effects are positive and $(E_i)_{\{1\}}+(E_i)_{\{-1\}}=I$.
Conversely, every POVM on $\{-1,1\}$ is of this form, with
$A_i=(E_i)_{\{1\}}-(E_i)_{\{-1\}}$. Define the effects
$(F_j)_{\{q\}}:=\frac12(I+qB_j)$ in the same way.
If $A_i^2=I$, its two effects are orthogonal projections, so that
the corresponding observable is sharp; the same holds for $B_j$.

For any given state $T$,
the pair of choices $(i,j)$ gives a probability measure $\mu_{ij}^T$
on $\{-1,1\}^2$ by
\begin{align}\label{eq:epi-quantum-pair-law}
\mu_{ij}^T(\{(p,q)\})
:=\tr\bigl(((E_i)_{\{p\}}\otimes (F_j)_{\{q\}})T\bigr),
\quad p,q\in\{-1,1\}.
\end{align}
These numbers are nonnegative and sum to one, since the tensor
products of the effects form a POVM. Their one-coordinate mar\-ginals
satisfy
\begin{align*}
\sum_{q\in\{-1,1\}}\mu_{ij}^T(\{(p,q)\})
&=\tr\bigl(((E_i)_{\{p\}}\otimes I)T\bigr),\\
\sum_{p\in\{-1,1\}}\mu_{ij}^T(\{(p,q)\})
&=\tr\bigl((I\otimes (F_j)_{\{q\}})T\bigr).
\end{align*}
The first expression is independent of $j$, and the second is
independent of $i$. In quantum information theory, this is the {\em no-signalling condition}:\index{no-signalling}
the outcome distribution on either subsystem is independent of
the choice of measurement on the other subsystem. It can be
equivalently stated by saying that the pair distributions $\mu_{ij}^T$ agree on their
common coordinates; this is the necessary compatibility condition
for the marginal problem discussed in Subsection \ref{subsec:marginals}.

The correlations $E_{ij}^T$ are defined by
\begin{align}\label{eq:epi-quantum-correlation}
E_{ij}^T:=\sum_{p,q\in\{-1,1\}}pq\,\mu_{ij}^T(\{(p,q)\})
=\tr\bigl((A_i\otimes B_j)T\bigr).
\end{align}
We call a matrix obtained in this way a {\em quantum correlation
matrix}.\index{correlation matrix!quantum}

\begin{lemma}\label{lem:epi-quantum-vectors}
There exists a real Hilbert space $\HH(H)$ such that for every quantum correlation matrix of the form $E^T \in \R^{m\times n}$, with $T\in \calS(H)$, there are vectors
$u_1,\ldots,u_m$ and $v_1,\ldots,v_n$ in $\HH(H)$ of norm at most one such that
$$
E_{ij}^T=\iprod{u_i}{v_j}_{\HH(H)},\quad i=1,\ldots,m,\ j=1,\ldots,n.
$$
\end{lemma}

\begin{proof} Let $E^T$ be of the form \eqref{eq:epi-quantum-correlation}.
Since $T\ge0$ and $\tr(T)=1$, its positive square root $T^{1/2}$
is a Hilbert--Schmidt operator with
$\n T^{1/2}\n_{\calL_2(H)}=1$. We regard $\calL_2(H)$ as a real
Hilbert space, with inner product
$\iprod{U}{V}_{\calL_2(H)}:=\Re\tr(V^\star U)$.

Put
$$
u_i:=(A_i\otimes I)T^{1/2},\qquad
v_j:=(I\otimes B_j)T^{1/2}.
$$
By the ideal property of Hilbert--Schmidt operators, all these
operators have Hilbert--Schmidt norm at most one. The trace identity
of Proposition \ref{prop:trace-commute} gives
\begin{align*}
\iprod{u_i}{v_j}_{\calL_2(H)}
&=\Re\tr\bigl(T^{1/2}(I\otimes B_j)(A_i\otimes I)T^{1/2}\bigr)\\*
&=\Re\tr\bigl(T^{1/2}(A_i\otimes B_j)T^{1/2}\bigr)
=\Re\tr\bigl((A_i\otimes B_j)T\bigr) = \tr\bigl((A_i\otimes B_j)T\bigr)
=E_{ij}^T,
\end{align*}
noting that $\tr((A_i\otimes B_j)T)$ is real by \eqref{eq:epi-quantum-correlation}.
\end{proof}

Write $\mathscr Q_{m,n}$ for the set of all quantum correlation
matrices, allowing the local Hilbert spaces, state, and observables
to vary. For a real matrix $M=(m_{ij})_{ij}\in\R^{m\times n}$,
define its {\em quantum correlation bound} by
\begin{align*}
\beta_{\rm q}(M):=
\max_{\substack{u_i,v_j\in\R^{m+n}\\
                \n u_i\n,\n v_j\n\le1}}
\Bigl|\sum_{i=1}^m\sum_{j=1}^n
m_{ij}\iprod{u_i}{v_j}\Bigr|.
\end{align*}
If $E\in\mathscr Q_{m,n}$, then
Lemma \ref{lem:epi-quantum-vectors} gives
$$
\Bigl|\sum_{i=1}^m\sum_{j=1}^n m_{ij}E_{ij}\Bigr|
\le\beta_{\rm q}(M).
$$
Here the span of the $m+n$ representing vectors can be identified
isometrically with a subspace of $\R^{m+n}$.
Conversely, the validity of these inequalities for every
$M\in\R^{m\times n}$ characterises the quantum correlation
matrices. This gives a quantum analogue of
Proposition \ref{prop:epi-correlation-duality}.

\begin{proposition}\label{prop:epi-quantum-correlation-duality}
A real matrix $E\in\R^{m\times n}$ belongs to $\mathscr Q_{m,n}$
if and only if
$$
\Bigl|\sum_{i,j}m_{ij}E_{ij}\Bigr|
\le\beta_{\rm q}(M),
\quad M\in\R^{m\times n}.
$$
\end{proposition}

\begin{proof}
Only sufficiency remains to be proved. Put $d:=m+n$ and let
$$
\mathscr V_{m,n}:=
\left\{
(\iprod{u_i}{v_j})_{i,j}:\,
u_i,v_j\in\R^d,\
\n u_i\n,\n v_j\n\le1
\right\}.
$$
This is the image, under $G\mapsto(G_{i,m+j})_{i,j}$,
of the compact convex set of real positive semidefinite
$d\times d$ matrices with diagonal entries at most one.
Indeed, these are precisely the Gram matrices of $d$ vectors
of norm at most one in $\R^d$.
Moreover, $\mathscr V_{m,n}$ is symmetric about $0$.
The separation argument of
Proposition \ref{prop:epi-correlation-duality} therefore shows
that the assumed inequalities imply $E\in\mathscr V_{m,n}$.

To realise such a matrix quantum mechanically, put $D:=2^d$
and, on $K:=(\C^2)^{\otimes d}$, define
$$
\Gamma_k:=
\sigma_3^{\otimes(k-1)}\otimes\sigma_1
\otimes I^{\otimes(d-k)},
\quad k=1,\ldots,d.
$$
These real symmetric matrices satisfy
$$
\Gamma_k^2=I,\qquad
\Gamma_k\Gamma_l=-\Gamma_l\Gamma_k\quad(k\ne l),
\qquad
\tr(\Gamma_k\Gamma_l)=D\delta_{kl}.
$$
Thus $C(x):=\sum_{k=1}^d x_k\Gamma_k$ satisfies
$C(x)^2=\n x\n^2I$ and $\n C(x)\n=\n x\n$.
Choose vectors representing $E$ and take
$$
A_i:=C(u_i),\qquad B_j:=C(v_j),\qquad
\psi:=D^{-1/2}\sum_{a=1}^D e_a\otimes e_a,
$$
where $(e_a)_{a=1}^D$ is the standard orthonormal basis of $K$.
The operators $A_i,B_j$ are selfadjoint contractions.
Since they are real symmetric, the pure state
$T:=\psi\,\bar\otimes\,\psi$ satisfies
$$
\tr((A_i\otimes B_j)T)
=D^{-1}\tr(A_iB_j)
=\iprod{u_i}{v_j}
=E_{ij}.
$$
Hence $E\in\mathscr Q_{m,n}$.
\end{proof}

The proof establishes Tsirel'son's characterisation
$$
\mathscr Q_{m,n}
=
\left\{
(\iprod{u_i}{v_j})_{i,j}:\,
u_i,v_j\in\R^{m+n},\
\n u_i\n,\n v_j\n\le1
\right\}.
$$
In particular, $\mathscr Q_{m,n}$ is compact, convex, and
symmetric about $0$, and every quantum correlation matrix
admits a realisation on finite-dimensional local Hilbert
spaces with a pure state.

\begin{theorem}[Tsirel'son bound]\label{thm:epi-tsirelson}
\index{theorem!Tsirel'son}
If $E\in\R^{m\times n}$ is a quantum correlation matrix, with $m,n\ge 2$, then
\begin{align}\label{eq:epi-tsirelson}
|E_{11}+E_{12}+E_{21}-E_{22}|\le2\sqrt2.
\end{align}
Equality can occur for a pure state and
sharp binary observables.
\end{theorem}

\begin{proof}
Choose the vectors given by Lemma \ref{lem:epi-quantum-vectors}.
The Cauchy--Schwarz inequality and the parallelogram identity give
\begin{align*}
|E_{11}+E_{12}+E_{21}-E_{22}|
&=|\iprod{u_1}{v_1+v_2}+\iprod{u_2}{v_1-v_2}|\\
&\le\n v_1+v_2\n+\n v_1-v_2\n\\
&\le\sqrt{2\bigl(\n v_1+v_2\n^2+\n v_1-v_2\n^2\bigr)}
=2\sqrt{\n v_1\n^2+\n v_2\n^2}
\le2\sqrt2.
\end{align*}

To prove sharpness, let $e_0,e_1$ be the standard orthonormal basis
of $\C^2$ and put
$$
\psi:=\frac1{\sqrt2}(e_0\otimes e_0+e_1\otimes e_1),
\qquad T:=\psi\,\bar\otimes\,\psi.
$$
Using the Pauli matrices of Subsection \ref{subsec:qubit}, choose
$$
A_1:=\sigma_3,\qquad A_2:=\sigma_1,\qquad
B_1:=\frac{\sigma_3+\sigma_1}{\sqrt2},\qquad
B_2:=\frac{\sigma_3-\sigma_1}{\sqrt2}.
$$
All four matrices are selfadjoint and have square $I$. For $B_1,B_2$
this follows from
$\sigma_1\sigma_3+\sigma_3\sigma_1=0$ and
$\sigma_1^2=\sigma_3^2=I$. Thus their binary observables are sharp.

Now $\sigma_3e_0=e_0$, $\sigma_3e_1=-e_1$, and $\sigma_1$ interchanges
$e_0,e_1$. Hence
$$
(\sigma_3\otimes\sigma_3)\psi=\psi,\qquad
(\sigma_1\otimes\sigma_1)\psi=\psi.
$$
The operator associated with the CHSH expression therefore satisfies
\begin{align*}
\bigl(A_1\otimes(B_1+B_2)+A_2\otimes(B_1-B_2)\bigr)\psi
=\sqrt2(\sigma_3\otimes\sigma_3+\sigma_1\otimes\sigma_1)\psi
=2\sqrt2\,\psi.
\end{align*}
Taking the inner product with $\psi$, we obtain equality in
\eqref{eq:epi-tsirelson}.
\end{proof}

In this example, both $(\sigma_3\otimes\sigma_1)\psi$ and
$(\sigma_1\otimes\sigma_3)\psi$ are orthogonal to $\psi$.
Consequently, the correlation matrix is
$$
E=\frac1{\sqrt2}\begin{pmatrix}1&1\\1&-1\end{pmatrix}.
$$
The expectations of $A_i\otimes I$ and $I\otimes B_j$ in the
state $T$ are zero. Expanding the effects in
\eqref{eq:epi-quantum-pair-law} gives the four pair laws explicitly:
$$
\mu_{ij}^T(\{(p,q)\})=\frac14(1+pqE_{ij}),
\quad p,q\in\{-1,1\}.
$$
Their one-coordinate mar\-ginals are all uniform. Nevertheless, they
cannot be the mar\-ginals of one probability measure on
$\{-1,1\}^4$: such a measure would make $E$ a classical correlation
matrix, contradicting \eqref{eq:epi-classical-CHSH}, since
$E_{11}+E_{12}+E_{21}-E_{22}=2\sqrt2>2$.
This is the example announced in Subsection \ref{subsec:marginals}:
the pair distributions agree on their common coordinates, yet admit
no joint distribution.

It is interesting to compare the CHSH example with the effective model of Theorem
\ref{thm:effective-model}. Applied to the pair POVMs with effects $(E_i)_{\{p\}}\otimes (F_j)_{\{q\}}$, the latter
gives Markov kernels $K^{(i,j)}$ on a common measurable space
and, for any given state $T$, a representing measure $\mu_T$ such that
$$
\mu_{ij}^T(\{(p,q)\})
=\int_X K^{(i,j)}(x,\{(p,q)\})\ud\mu_T(x).
$$
These kernels depend on both measurement choices; the theorem does not assert that they factor into kernels for the individual subsystems.

In fact, for this example no such factorisation exists. Suppose that on some probability space $(Y,\mathscr Y,\nu^T)$ there were Markov kernels $k_i,l_j$ with outcome space $\{-1,1\}$ such that
\begin{align}\label{eq:epi-kernel-factorisation}
\mu_{ij}^T(\{(p,q)\})
=\int_Y k_i(y,\{p\})l_j(y,\{q\})\ud\nu^T(y)
\end{align}
for all $i,j,p,q$, with $k_i$ independent of $j$ and $l_j$ independent of $i$. Then
$$
\mu^T(\{(\varepsilon,\delta)\})
:=\int_Y\prod_{i=1}^m k_i(y,\{\varepsilon_i\})
             \prod_{j=1}^n l_j(y,\{\delta_j\})\ud\nu^T(y)
$$
would define a joint probability measure with pair marginals $\mu_{ij}^T$: summing over each omitted coordinate replaces its factor by $1$. Conversely, any joint distribution gives such a factorisation by taking the coordinate outcomes as deterministic kernels. Thus \eqref{eq:epi-kernel-factorisation} is equivalent to the existence of a joint distribution, and the CHSH example rules it out.

In the language of effective models, \eqref{eq:epi-kernel-factorisation} is Bell's factorisability
condition (local causality), combined with the requirement that the distribution $\nu^T$ of the effective states does not depend on the measurements chosen; the CHSH example therefore shows that no effective model of quantum theory can have these properties (Bell's theorem), although the effective model of Theorem~\ref{thm:effective-model} exists.

This complements Subsection \ref{subsec:qubit-effective}. There, the representing measures could not be chosen affinely in the state, or equivalently, no single POVM realised the prescribed kernels for all states. Here, even for one fixed state, a representation by products of kernels for the individual
choices may be impossible. A common effective model therefore need not have either of these additional properties, namely preparation noncontextuality and locality.

\subsection{Grothendieck's Inequality}\label{subsec:epi-grothendieck}

For the CHSH expression, the quantum bound exceeds the classical bound by the factor $\sqrt2$. We now ask whether the increase can be arbitrarily large for other linear combinations of correlations. By Lemma \ref{lem:epi-quantum-vectors}, it is enough to compare expressions involving inner products of vectors of norm at most one with the corresponding expressions involving scalar signs. Grothendieck's inequality gives a bound independent of the number of vectors and of the Hilbert space.

\begin{theorem}[Grothendieck's inequality]\label{thm:epi-grothendieck}
\index{theorem!Grothendieck}\index{inequality!Grothendieck}
Let $H$ be a real Hilbert space and let $u_1,\ldots,u_m$ and $v_1,\ldots,v_n\in H$ have norm at most one. Then, for every $M=(m_{ij})\in\R^{m\times n}$,
\begin{align}\label{eq:epi-grothendieck}
\Bigl|\sum_{i=1}^m\sum_{j=1}^n m_{ij}\iprod{u_i}{v_j}\Bigr|
\le \frac{\pi}{2\log{(}1+\sqrt2{)}} \beta(M).
\end{align}
\end{theorem}

The least constant which can replace the constant in \eqref{eq:epi-grothendieck}, for every real Hilbert space and every finite matrix, is called the {\em real Grothendieck constant},
\index{Grothendieck constant} and is denoted by $K_G$. Thus the theorem asserts that
 $$K_G\le \frac{\pi}{2\log{(}1+\sqrt2{)}} \approx 1.7822.$$
Below we will also obtain the lower bound $K_G \ge \sqrt 2$.

We begin with a Gaussian lemma.

\begin{lemma}[Gaussian signs]\label{lem:epi-gaussian-signs}
Let $u,v$ be unit vectors in $\R^d$, and let $\gamma_d$ be the standard Gaussian measure on $\R^d$. Then
\begin{align*}
\int_{\R^d}{\rm sign}(x\cdot u)\,{\rm sign}(x\cdot v)\ud\gamma_d(x)
=\frac2\pi\arcsin(u\cdot v).
\end{align*}
\end{lemma}

 \begin{proof}
Since proper subspaces of $\R^d$ have Gaussian measure zero, the value assigned to
${\rm sign}(0)$ does not affect the integral. If $v=u$ or $v=-u$, the integrand equals $1$ or $-1$ almost everywhere, respectively, and the assertion follows.

Suppose now that $v\ne\pm u$, and put $\theta:=\arccos(u\cdot v)\in(0,\pi)$. Choose an orthonormal basis
whose first two vectors are
$$
e_1:=u,\qquad
e_2:=\frac{v-(u\cdot v)u}{\sin\theta}.
$$
Then $v=\cos\theta\,e_1+\sin\theta\,e_2$. Making the corresponding orthogonal change of coordinates and integrating out the remaining coordinates, we obtain
\begin{align*}
&\int_{\R^d}{\rm sign}(x\cdot u)\,{\rm sign}(x\cdot v)\ud\gamma_d(x)\\
&\qquad =
\frac1{2\pi}\int_{\R^2} {\rm sign}(y_1)\, {\rm sign}(y_1\cos\theta+y_2\sin\theta) e^{-(y_1^2+y_2^2)/2}\ud y_1\ud y_2.
\end{align*}
Here we used orthogonal invariance of the Gaussian density and Fubini's theorem.

We now turn to the polar coordinates $y_1=r\cos\varphi$, $y_2=r\sin\varphi$, with
$r>0$ and $-\frac12\pi<\varphi<\frac32\pi$.
Since
$$
y_1\cos\theta+y_2\sin\theta=r\cos(\varphi-\theta),
\qquad
\int_0^\infty e^{-r^2/2}r\ud r=1,
$$
the preceding integral equals
$$
\frac1{2\pi}
\int_{-\pi/2}^{3\pi/2}
{\rm sign}(\cos\varphi)\,
{\rm sign}(\cos(\varphi-\theta))\ud\varphi.
$$

On the chosen angular interval, $\cos\varphi$ is positive precisely when $\varphi\in(-\frac12\pi,\frac12\pi)$, whereas $\cos(\varphi-\theta)$ is positive precisely when
$\varphi\in(\theta-\frac12\pi,\theta+\frac12\pi)$.
Apart from the endpoints, the two signs therefore differ precisely on
$$
D_\theta:=
\Bigl(-\frac12\pi,\theta-\frac12\pi\Bigr)
\cup
\Bigl(\frac12\pi,\theta+\frac12\pi\Bigr).
$$
These intervals are disjoint and each has length $\theta$.
The product of the signs is consequently $-1$ on a set of length $2\theta$ and $1$ on a set of length $2\pi-2\theta$. Thus
$$
\int_{\R^d}{\rm sign}(x\cdot u)\,
              {\rm sign}(x\cdot v)\ud\gamma_d(x)
=\frac{(2\pi-2\theta)-2\theta}{2\pi}
=1-\frac{2\theta}{\pi}.
$$
Finally, $\arcsin(\cos\theta)=\frac12\pi-\theta$ for
$0\le\theta\le\pi$, and $\cos\theta=u\cdot v$.
This proves the lemma.
\end{proof}

The following proof of Grothendieck's inequality is due to Krivine.

\begin{proof}[Proof of Theorem \ref{thm:epi-grothendieck}]
We may replace $H$ by the finite-dimensional span of the vectors
$u_i$ and $v_j$. We may also assume that all these vectors have norm
one. Indeed, in $H\oplus\R\oplus\R$ replace them by
$$
\widetilde u_i:=(u_i,\sqrt{1-\n u_i\n^2},0),\qquad
\widetilde v_j:=(v_j,0,\sqrt{1-\n v_j\n^2}).
$$
These vectors have norm one and
$\iprod{\widetilde u_i}{\widetilde v_j}=\iprod{u_i}{v_j}$.
We therefore retain the notation $H,u_i,v_j$ and assume from now
on that $H$ is finite-dimensional and $u_i,v_j$ are unit vectors.

Put
$c:=\log{(}1+\sqrt2{)}$ and
$a_k:={c^{2k+1}}/{(2k+1)!}$ for $k\ge 0$.
Then $0<c<\frac12\pi$ and
$$\sum_{k\ge0}a_k=\sinh c=1.$$ Consider the real Hilbert space
$\mathscr H:=\bigoplus_{k\ge0}\Gamma^{2k+1}(H)$,
the direct sum of the odd symmetric tensor powers of $H$
(see Appendix \ref{sec:tensor}). For a unit vector $u\in H$, define
$$
U(u):=\bigl(\sqrt{a_k}\,u^{\otimes(2k+1)}\bigr)_{k\ge0},
\qquad
V(u):=\bigl((-1)^k\sqrt{a_k}\,u^{\otimes(2k+1)}\bigr)_{k\ge0}.
$$
The Hilbert space tensor product satisfies
$\iprod{u^{\otimes r}}{v^{\otimes r}}=\iprod{u}{v}^r$ and $\n u^{\otimes r}\n=\n u\n^r$.
Consequently, $U(u),V(u)\in\mathscr H$ and
$\n U(u)\n^2=\n V(u)\n^2=\sum_{k\ge0}a_k=1$.
For unit vectors $u,v\in H$, the power series of the sine gives
\begin{align}\label{eq:epi-krivine-sine}
\iprod{U(u)}{V(v)}
=\sum_{k\ge0}(-1)^k
  \frac{c^{2k+1}\iprod{u}{v}^{2k+1}}{(2k+1)!}
=\sin\bigl(c\iprod{u}{v}\bigr).
\end{align}

Although $\mathscr H$ is infinite-dimensional, only the finitely
many vectors $U(u_i)$ and $V(v_j)$ are needed. Their span is isometric to $\R^d$ for some finite $d$.
On this space put
$$
\varepsilon_i(x):={\rm sign}\iprod{x}{U(u_i)},\qquad
\delta_j(x):={\rm sign}\iprod{x}{V(v_j)}.
$$
All these functions take values in $\{-1,1\}$ outside one
$\gamma_d$-null set. Assigning the value $1$ on that set, if necessary,
we may regard their joint law as a probability measure on
$\{-1,1\}^{m+n}$.
By Lemma \ref{lem:epi-gaussian-signs} and
\eqref{eq:epi-krivine-sine},
\begin{align*}
\int_{\R^d}\varepsilon_i\delta_j\ud\gamma_d
&=\frac2\pi\arcsin\bigl(\sin(c\iprod{u_i}{v_j})\bigr)
=\frac{2c}\pi\iprod{u_i}{v_j}.
\end{align*}
In the last equality we used
$|c\iprod{u_i}{v_j}|\le c<\frac12\pi$.
It follows that
\[
\Bigl|\sum_{i,j}m_{ij}\iprod{u_i}{v_j}\Bigr|
=\frac\pi{2c}
  \Bigl|\int_{\R^d}\sum_{i,j}m_{ij}\varepsilon_i\delta_j\ud\gamma_d\Bigr|
\le\frac\pi{2c}
  \int_{\R^d}\Bigl|\sum_{i,j}m_{ij}\varepsilon_i\delta_j\Bigr|\ud\gamma_d
\le\frac\pi{2c}\beta(M).
\qedhere\]
\end{proof}

We can now apply the inequality to quantum correlations.

\begin{theorem}[Grothendieck bound for quantum correlations]
\label{thm:epi-quantum-grothendieck}
For every quantum correlation matrix $E\in\R^{m\times n}$ and
every $M\in\R^{m\times n}$ we have
$$
\Bigl|\sum_{i,j}m_{ij}E_{ij}\Bigr|
\le\frac{\pi}{2\log{(}1+\sqrt2{)}}\beta(M).
$$
Moreover, $K_G^{-1}E$ is a classical correlation matrix.
\end{theorem}

\begin{proof} The inequality follows by
applying Grothendieck's inequality
to the vectors of Lemma \ref{lem:epi-quantum-vectors}.
For the second assertion, the definition of $K_G$ gives $|\sum_{i,j}m_{ij}E_{ij}|\le K_G\beta(M)$. Dividing by $K_G$ and applying Proposition \ref{prop:epi-correlation-duality} gives the result.
\end{proof}

Finally, let $E$ be the correlation matrix of the two-qubit example
attaining the Tsirel'\-son bound. We have $\beta(M_{\rm CHSH})=2$ and
$\sum_{i,j}(M_{\rm CHSH})_{ij}E_{ij}=2\sqrt2$.
Its Hilbert space representation from Lemma
\ref{lem:epi-quantum-vectors} therefore shows that $K_G\ge\sqrt2$.
Together with Krivine's bound, this gives
$$
\sqrt2\le K_G\le\frac{\pi}{2\log{(}1+\sqrt2{)}}.
$$

\begin{problems}

\item
Prove the assertions about orthogonal projections in Section \ref{subsec:class-quantum}.

\item\label{prob:orthomodular}
Prove that if $P$ and $Q$ are orthogonal projections on a Hilbert space such that $\ran(P)$ is contained in $\Ran(Q)$, then $  Q = P \vee (Q \wedge \neg P).$

\item
Let $\phi:\calL(H)\to \C$ be a state (recall our convention that states are normal), let $(h_n)_{n\ge 1}$ be an orthonormal basis of $H$, and let $P_n$ denote the orthogonal projection of $H$ onto the span of the set $\{h_1,\dots,h_n\}$. Prove that for all $T\in \calL(H)$
we have
$$ \phi(T) = \limn \phi(P_nT).$$

\noindent{\em Hint:}\ Use Theorem \ref{thm:state-posfc},
or apply the Cauchy--Schwarz inequality to the positive
sesquilinear form $(S,U)\mapsto\phi(SU^\star)$,
with $S:=I-P_n$ and $U:=T^\star$.
Use normality to show that $\phi(I-P_n)\to0$.

\item\label{prob:phi-bicomm}
The aim of this problem is to prove that for a linear mapping $\phi: \calL(H)\to \C$ the following assertions are equivalent:
\begin{enumerate}[label={\rm(\arabic*)}, leftmargin=*]
  \item\label{it:lem:weakoperatortop1}  $\phi(T) = \sum_{j=1}^k \iprod{Tx_j}{y_j}$ for suitable $k\ge 1$ and $x_1,\dots,x_k,y_1,\dots,y_k\in H$;
  \item\label{it:lem:weakoperatortop3}  $\phi$ is continuous with respect to the weak operator topology;
  \item\label{it:lem:weakoperatortop2}  $\phi$ is continuous with respect to the strong operator topology.
\end{enumerate}
\begin{enumerate}[\rm(a), leftmargin=*]
  \item Prove the implications \ref{it:lem:weakoperatortop1}$\Rightarrow$\ref{it:lem:weakoperatortop3}$\Rightarrow$\ref{it:lem:weakoperatortop2}.
\end{enumerate}
The remainder of the problem is devoted to a proof of the implication \ref{it:lem:weakoperatortop2}$\Rightarrow$\ref{it:lem:weakoperatortop1}.
\begin{enumerate}[\rm(a), leftmargin=*]\setcounter{enumii}{1}
  \item Show that \ref{it:lem:weakoperatortop2}
  implies that there exist $x_1,\dots,x_k\in H$ such that
  $$ |\phi(T)| \le \max_{1\le j\le k} \n Tx_j\n, \quad T\in \calL(H).$$
  \item Let $K$ be the closure of the subspace
  $\{(Tx_1,\dots,Tx_k)\in H^k:\, T\in \calL(H)\}$ in $H^k$\!. Show that the linear mapping $\psi: K\to \C$ defined by
  $$ \psi(Tx_1,\dots,Tx_k):= \phi(T)$$
  is well defined and bounded.
  \item Using the Riesz representation theorem, show that $\phi$ is of the form as in \ref{it:lem:weakoperatortop1}.
\end{enumerate}
\item\label{prob:phi-state}
Let $H$ be a separable complex Hilbert space and let
$\phi:\calL(H)\to\C$ be a bounded linear functional.
Prove that the following assertions are equivalent:
\begin{enumerate}[label={\rm(\arabic*)}, leftmargin=*]
\item\label{it:phi-state-trace}
There exists $T\in\calL_1(H)$ such that
$$
\phi(S)=\tr(ST),\quad S\in\calL(H).
$$
\item\label{it:phi-state-series}
There exist sequences $(x_n)_{n\ge1}$ and $(y_n)_{n\ge1}$ in $H$
satisfying
$$
\sum_{n\ge1}\|x_n\|^2<\infty,\qquad
\sum_{n\ge1}\|y_n\|^2<\infty,
$$
such that
$$
\phi(S)=\sum_{n\ge1}\iprod{Sx_n}{y_n},
\quad S\in\calL(H).
$$
\item\label{it:phi-state-weak}
The restriction of $\phi$ to the closed unit ball of $\calL(H)$
is continuous for the weak operator topology.
\item\label{it:phi-state-strong}
The restriction of $\phi$ to the closed unit ball of $\calL(H)$
is continuous for the strong operator topology.
\end{enumerate}
Show that, in this case, the operator $T$ is unique and
$\|\phi\|=\|T\|_1$.

If $\phi$ is also positive, show that $T$ is positive and that
the representation in \ref{it:phi-state-series} can be chosen
with $y_n=x_n$ and the vectors $x_n$ mutually orthogonal.
If, in addition, $\phi(I)=1$, show that
$$
\sum_{n\ge1}\|x_n\|^2=1.
$$

\noindent
{\em Hint:}\ For
\ref{it:phi-state-trace}$\Rightarrow$\ref{it:phi-state-series},
use the singular value decomposition.
For
\ref{it:phi-state-strong}$\Rightarrow$\ref{it:phi-state-trace},
choose finite-rank orthogonal projections $P_n$ converging
strongly to $I$. By Problem~\ref{prob:phi-bicomm}, the
functionals $\phi_n(S):=\phi(SP_n)$ are represented by
finite-rank operators. Use continuity of $\phi$ at $0$ on
the closed unit ball and the estimate
$$
\sup_{\|S\|\le1}\|S(I-P_n)x\|
\le\|(I-P_n)x\|,\quad x\in H,
$$
to show that $\|\phi-\phi_n\|\to0$, and then use trace duality.


\item
Consider a qubit in state $\al\ket{0} + \beta\ket{1}$, where $\al,\beta\in\C $ satisfy $|\al|^2+|\beta|^2 =1$. Compute the probabilities that a measurement of the spin in direction
$j\in\{1,2,3\}$ yields `up' or `down', respectively.

\item
We take a closer look at the Pauli matrices $\sigma_1$, $\sigma_2$, $\sigma_3$.

\begin{enumerate}[\rm(a), leftmargin=*]
\item Prove the identity \eqref{eq:sp-proj-Pauli}.
\item Show that the complex exponentials of the Pauli matrices are given by
$$ \exp(i\theta\sigma_j) = (\cos\theta)I + i(\sin\theta)\sigma_j, \quad j=1,2,3.$$
\item Show that if $v$ is a unit vector in $\R^3$, then for all $n\in\N$ we have
$$ (v\cdot \sigma)^n = \begin{cases}
                        I, & \hbox{$n$ even},\\
                        v\cdot \sigma, & \hbox{$n$ odd}.
                       \end{cases}
$$
Use this to prove the identity
\begin{align*}
\exp(i\theta v\cdot\sigma)  =  (\cos\theta)I + i(\sin\theta)v\cdot\sigma.
\end{align*}
Furthermore, show that $\det(\exp(i\theta v\cdot\sigma))=1$.
\item Let $\mathrm{SU}(2)$ denote the group of unitary matrices acting on $\C^2$ with determinant $1$.
Show that for every $U\in \mathrm{SU}(2)$ there exist $\theta\in [0,2\pi]$
and $v\in \R^3$ with $|v|=1$ such that
$U=\exp(i\theta\, v\cdot \sigma)$.
\end{enumerate}

\item\label{prob:state-discrimination}
A quantum system is prepared in one of two states
$\rho,\sigma\in\mathscr S(H)$, each with probability $\frac12$.
We wish to determine which state was prepared.
A measurement with two outcomes is described by the effects
$E$ and $I-E$.
If the first outcome occurs, we guess $\rho$; otherwise we
guess $\sigma$. The probability of a correct guess is
$$
p(E):=\frac12\tr(E\rho)+\frac12\tr((I-E)\sigma).
$$

\begin{enumerate}[\rm(a), leftmargin=*]
\item Put $D:=\rho-\sigma$. Show that
$$
\sup_{0\le E\le I}\tr(ED)=\frac12\n D\n_1.
$$

\noindent{\em Hint:}\ Use the spectral decomposition of $D$
and the identity $\tr(D)=0$.

\item Deduce the {\em Holevo--Helstrom formula}
\index{theorem!Holevo--Helstrom}
$$
\sup_{0\le E\le I}p(E)
=\frac12+\frac14\n\rho-\sigma\n_1.
$$
Show that the supremum is attained by taking $E$ to be
the spectral projection of $D$ corresponding to $(0,\infty)$.
Thus an optimal measurement can be chosen to be
projection-valued.

\item Suppose that $\rho=h\,\bar\otimes\,h$ and
$\sigma=k\,\bar\otimes\,k$, where $\n h\n=\n k\n=1$.
Show that
$$
\n\rho-\sigma\n_1
=2\sqrt{1-|\iprod{h}{k}|^2}.
$$
Deduce that these pure states can be distinguished with
probability $1$ if and only if $h$ and $k$ are orthogonal.
\end{enumerate}

\item\label{prob:local-discrimination}
Let $H,K$ be separable complex Hilbert spaces and let
$\rho,\sigma\in\mathscr S(H\ot K)$.
\begin{enumerate}[\rm(a), leftmargin=*]
\item
Interpret the inequality $\n\tr_K(\rho)-\tr_K(\sigma)\n_1
\le\n\rho-\sigma\n_1$ (see Problem \ref{prob:partial-trace-contractive})
in terms of distinguishing states
by measurements on $H$ alone.

\item Now take $H=K=\C^2$, with standard basis $e_1,e_2$, and put
$$
u_\pm:=\frac1{\sqrt2}(e_1\ot e_1\pm e_2\ot e_2),
\qquad
\rho_\pm:=u_\pm\,\bar\otimes\,u_\pm.
$$
Compute both partial traces of these states.
Show that a measurement of the composite system can
distinguish them with probability $1$, whereas access to
either subsystem alone gives no information about which
state was prepared.
\end{enumerate}

\item
Let $Q:\calF\to\calE(H)$ be a POVM and let
$\Psi:B_{\rm b}(\Om)\to\calL(H)$ be its bounded functional calculus.
Prove that the following assertions are equivalent:
\begin{enumerate}[\rm(a), leftmargin=*]
\item $Q$ is projection-valued;
\item $\Psi(\one_F\one_G)=\Psi(\one_F)\Psi(\one_G)$ for all
$F,G\in\calF$;
\item $\Psi(fg)=\Psi(f)\Psi(g)$ for all bounded simple functions
$f,g$;
\item $\Psi(fg)=\Psi(f)\Psi(g)$ for all
$f,g\in B_{\rm b}(\Om)$.
\end{enumerate}

\item\label{prob:Naimark-finite}
Let $(e_j)_{j=1}^n$ be the standard orthonormal basis of $\C^n$. Define $J:H\to H\otimes\C^n$ by
$$
Jh:=\sum_{j=1}^n Q_j^{1/2}h\otimes e_j,
$$
where $Q_1,\ldots,Q_n\in\calE(H)$ satisfy
$\sum_{j=1}^n Q_j=I$.
\begin{enumerate}[\rm(a), leftmargin=*]
\item Show that $J$ is an isometry.
\item If
$P_j:=I\otimes(e_j\,\bar\otimes\,e_j)$, show that the $P_j$ are pairwise
orthogonal projections with $\sum_{j=1}^nP_j=I$ and
that $Q_j=J^\star P_jJ$ for $1\le j\le n$.
\item Conclude that this gives an explicit representation of a finite-outcome POVM as the compression of a projection-valued measure on a larger Hilbert space.
\end{enumerate}

\item \label{prob:Markov}
Let $X$ and $\Omega$ be complete separable metric spaces, and let
$\Phi:M_1^+(X)\to M_1^+(\Omega)$
be continuous for the weak topologies.
\begin{enumerate}[\rm(a), leftmargin=*]
 \item  Prove that $\Phi$ is affine if and only if there exists a Markov kernel
$K:X\times\calB(\Omega)\to[0,1]$ such that
$\Phi(\mu)=K_\#\mu$ for all $\mu\in M_1^+(X)$.
 \item Show that in this case $K(x,\cdot)=\Phi(\delta_x)$ for all $x\in X,$
and the mapping $x\mapsto K(x,\cdot)$ is continuous for the weak
topology of $M_1^+(\Omega)$.
\end{enumerate}

\noindent
{\em Hint:}\ First show that the finitely supported measures are weakly dense in $M_1^+(X)$.

\item Explain how the result of Problem \ref{prob:mult-fc-CK} is related to Theorem \ref{thm:classical-PVM-POVM-kernels}.

\end{problems}

%% file: ch16-SymmQuant.tex
\chapter{Symmetry and Quantisation}\label{ch:SymmQuant}

\blfootnote{\rm \small \begin{quote}This book has been published by Cambridge University Press in the series `Cambridge Studies in Advanced Mathematics'. The present version is free to view and download for personal use only. Not for re-distribution, re-sale or use in derivative works.\newline {\copyright} Jan van Neerven\end{quote}}

\noindent
Having described quantum states and observables, we now examine their interplay with symmetries implemented by unitary operators. Covariance describes how observables transform under their action. For position and momentum, these ideas lead to the Weyl commutation relations and the Stone--von Neumann theorem. They also provide a starting point for quantisation, which associates operators with functions on classical phase space and relates quantum commutators to classical Poisson brackets. Second quantisation extends the Hilbert space description to systems with a variable number of particles. Its construction using symmetric tensor products reveals close connections between Gaussian analysis and quantum mechanics, in which the Fourier transform and semigroups encountered earlier in the book acquire new interpretations.

Unless otherwise stated, $H$ is a separable complex Hilbert space.

\section{Transformations and Symmetries}\label{sec:symmetries}

The description of states and observables also allows us to describe
transformations of a quantum system. We begin with transformations
implemented by unitary operators and examine their effect on measurement
probabilities. This will lead to the notion of a symmetry of the pure
state space and to its application to time evolution.

\subsection{Unitary Transformations}\label{subsec:unitary-transformations}

We denote the set of all unitary operators on $H$ by
$\mathscr U(H)$.\index{$U$@$\mathscr U(H)$}
A unitary operator $U$ transforms a unit vector $h$ into the unit vector
$Uh$, and hence transforms the corresponding pure state by
$$
h\,\bar\otimes\,h\mapsto
Uh\,\bar\otimes\,Uh=U(h\,\bar\otimes\,h)U^\star.
$$
More generally, it acts on states by
$$
\Theta_U(T):=UTU^\star,
\quad T\in\mathscr S(H).
$$
Conjugation by $U$ preserves positivity and trace, so $\Theta_U$ maps
$\mathscr S(H)$ into itself. It is affine and bijective, with inverse
$\Theta_{U^\star}$.

Let $Q:\calF\to\calE(H)$ be a POVM on a measurable space
$(\Om,\calF)$. If the state $T$ is transformed into $\Theta_U(T)$,
while the observable is kept fixed, then the probability of an outcome
in $F\in\calF$ becomes
$$
\tr(Q_F\Theta_U(T))
=\tr(Q_FUTU^\star)
=\tr((U^\star Q_FU)T).
$$
Thus the same probabilities are obtained by keeping the state $T$
fixed and replacing $Q$ by the POVM $F\mapsto U^\star Q_FU$.

One may also transform both the state and the observable by the same
unitary operator. If
$$
T':=UTU^\star,\qquad Q'_F:=UQ_FU^\star,
$$
then
$$
\tr(Q'_FT')
=\tr(UQ_FTU^\star)
=\tr(Q_FT).
$$
All outcome probabilities are therefore unchanged.

More generally, the same construction applies to a unitary operator
$U:H\to K$ between different Hilbert spaces, transporting states and
observables on $H$ to states and observables on $K$. Such changes of
representation will occur in the Stone--von Neumann theorem and in our
discussion of second quantisation.

\subsection{Symmetries and Wigner's Theorem}\label{subsec:Wigner}

A transformation is called a symmetry when it preserves a specified
structure. On the pure state space $\mathscr P_1(H)$, the structure
of interest is given by the transition probabilities
$(P,R)\mapsto\tr(PR)$. We call a bijection of $\mathscr P_1(H)$ a
{\em symmetry}\index{symmetry} if it preserves
these probabilities. Every unitary operator induces such a symmetry.
To describe all of them, we also need conjugate-linear operators.

The {\em adjoint} of a bounded conjugate-linear mapping $U:H\to H$ (that is, a mapping satisfying $U(x+y) = U(x)+U(y)$ and $U(cx) = \ov c Ux$)
is the unique conjugate-linear mapping $U^\star:H\to H$ defined by
$$ \iprod{x}{U^\star y} = \ov{\iprod{Ux}{y}}, \quad x,y\in H.$$
A bounded conjugate-linear mapping $U:H\to H$ is called
{\em antiunitary}\index{antiunitary} if $UU^\star = U^\star U = I.$
It is straightforward to check that if $U:H\to H$ is unitary or
antiunitary, the mapping  $\Theta(T) = UTU^\star$
maps $\mathscr S(H)$ onto itself, sends pure states to pure states, and
preserves transition probabilities,
$$
   \tr(\Theta(P)\Theta(T)) = \tr(PT),
        \quad P\in\mathscr P_1(H),\ T \in \mathscr S(H).
$$ Indeed, denoting by $\mathscr P_1(H)$ the set of rank-one orthogonal
projections, for
$P=h\,\bar\otimes\,h\in\mathscr P_1(H)$ with $\|h\|=1$ we have $UPU^\star=Uh\,\bar\otimes\, Uh$ and hence, for all $T\in\mathscr S(H)$,
\begin{align*}
\tr(\Theta(P)\Theta(T)) = \iprod{(UTU^\star)Uh}{Uh} =
\iprod{Th}{h} = \tr(PT).
\end{align*}
In the antiunitary case, the middle equality uses that
$\iprod{Th}{h}$ is real.

Wigner's theorem gives the following converse: every bijection of the pure state space preserving transition probabilities is induced by a unitary or antiunitary operator.

\begin{theorem}[Wigner's theorem]\label{thm:Wigner}\index{theorem!Wigner}
Let $H$ be a separable complex Hilbert space and let
$\Theta:\mathscr P_1(H)\to\mathscr P_1(H)$
be a bijection such that
$$
        \tr(\Theta(P)\Theta(Q))=\tr(PQ),
        \quad P,Q\in\mathscr P_1(H).
$$
Then there exists either a unitary or an antiunitary operator $U:H\to H$
such that
$$
        \Theta(h\,\bar\otimes\,h)= Uh\,\bar\otimes\,Uh,
        \quad \|h\|=1.
$$
If $\dim H\ge2$, this operator is unique up to multiplication
by a scalar of modulus one.
\end{theorem}

A complete proof would take us too far afield.
In the one-dimensional case, say $H=\C e$ with $\|e\|=1$, the set
$\mathscr P_1(H)$ consists of the single orthogonal projection $e\,\bar\otimes\, e$. The
unitary and antiunitary operators on $\C e$ are precisely the operators of the form
$ze\mapsto c ze$, respectively $ze\mapsto c\ov z e$, with $|c|=1$.
For every given $c\in \C$ with $|c|=1$, these operators are different as mappings on $\C e$, but they induce the same transformation on states.

Wigner's theorem tells us that symmetries of the pure state space, understood
as bijections preserving transition probabilities, are implemented on the
underlying Hilbert space by operators that are either unitary or antiunitary.

Preservation of transition probabilities does not require individual
states or observables to remain unchanged. Invariance of a particular
observable is an additional condition, to which we return in the next
section.

\subsection{One-Parameter Groups and Time Evolution}\label{subsec:unitary-time-evolution}

Time evolution provides a principal example of a family of transformations.
In dynamics one is primarily interested in one-parameter groups
of symmetries of the pure state space. Suppose that
$(\Theta(t))_{t\in\R}$ is such a group, and assume that
$$
t\mapsto\tr(\Theta(t)(P)Q)
$$
is continuous for every $P,Q\in\mathscr P_1(H)$.
This expresses continuity in time of all transition probabilities.

If $\dim H\le1$, its action on pure states is trivial and can be
implemented by the constant unitary group $I$. We therefore assume
that $\dim H\ge2$. By the uniqueness part of Wigner's theorem, the type of
an implementation, unitary or antiunitary, is determined by
$\Theta(t)$. The type therefore defines a group homomorphism from
$\R$ into the two-element group $\{\pm1\}$.
Such a homomorphism $\chi$ satisfies
$\chi(t)=\chi(t/2)^2=1$ for every $t\in\R$.
Thus the implementing operators $U(t)$ may be chosen
unitary. The identity
$\Theta(t)\Theta(s)=\Theta(t+s)$
then implies that there exists a scalar $c(t,s)$ of modulus one such that
$$
        U(t)U(s)=c(t,s)U(t+s).
$$
From the associative law $U(t)(U(s)U(r))=(U(t)U(s))U(r)$ we obtain the
{\em cocycle identity}\index{cocycle identity}
$$
        c(s,r)c(t,s+r)=c(t,s)c(t+s,r).
$$
The {\em unitary operator ray} of a unitary operator $U$ is
the set $\{cU:\,c\in\C,\ |c|=1\}$.
Bargmann's lifting theorem, specialised to the additive
group $\R$, states that a one-parameter group of unitary operator
rays admits a strongly continuous unitary group of representatives
provided that
$$
t\mapsto|\iprod{U(t)h}{k}|
$$
is continuous for all $h,k\in H$.
This condition is independent of the phases chosen for $U(t)$.
In our setting it follows from the assumed continuity of
transition probabilities, since for unit vectors $h,k$ we have
$$
|\iprod{U(t)h}{k}|^2
=
\tr\bigl(\Theta(t)(h\,\bar\otimes\,h)
              (k\,\bar\otimes\,k)\bigr).
$$

Consequently, there exists a function $d:\R\to\C$ with
$|d(t)|=1$ such that
$$
\wt U(t):=d(t)^{-1}U(t),\quad t\in\R,
$$
forms a strongly continuous unitary group.
Its group law implies
$$
c(t,s)=\frac{d(t)d(s)}{d(t+s)}.
$$
Changing these phases does not change the action on pure states,
so
$$
\Theta(t)(h\,\bar\otimes\,h)
=
\wt U(t)h\,\bar\otimes\,\wt U(t)h,
\quad \|h\|=1.
$$
Hence, by Stone's theorem,
there exists a selfadjoint operator $E$, the {\em Hamiltonian}
associated with the group $(\wt U(t))_{t\in\R}$, such that
$$
        \wt U(t)=e^{-itE},
        \quad t\in\R.
$$
For $h\in\Dom(E)$, the function $u(t)= \wt U(t)h$ satisfies $u(0) = h$ and
$$
        \frac{{\rm d}}{{\rm d}t}u(t) =
        -iE u(t), \quad t\in\R.
$$
This is an abstract version of the Schr\"odinger
equation,\index{Schr\"odinger!equation}
which, in concrete examples, corresponds to the case $H=L^2(\R^d)$ and $E=-\Delta+V$, where $V$ is a potential.

For an initial state $T_0\in\mathscr S(H)$, the corresponding evolution
of states is
$$
T(t):=\wt U(t)T_0\wt U(t)^\star,
\quad t\in\R.
$$
If $Q$ is a POVM, the trace identity from the first subsection gives
$$
\tr(Q_FT(t))
=\tr((\wt U(t)^\star Q_F\wt U(t))T_0).
$$
In the {\em Schr\"odinger picture},\index{Schr\"odinger!picture}
the state evolves and the observable is kept fixed. In the
{\em Heisenberg picture},\index{Heisenberg!picture}
the state is kept fixed and the observable evolves according to
$Q_F(t):=\wt U(t)^\star Q_F\wt U(t)$. The identity shows that these descriptions
give the same outcome probabilities.

Once the Hamiltonian $E$ has been specified, a unitary transformation
$W$ preserves it if $WEW^\star=E$, where the equality includes equality
of domains. The spectral calculus then gives
$$
W\wt U(t)W^\star
=We^{-itE}W^\star
=e^{-itWEW^\star}
=e^{-itE}
=\wt U(t),
\quad t\in\R.
$$
Thus such a transformation commutes with the time evolution.
Preserving the Hamiltonian is an additional requirement on $W$;
it does not follow merely from unitarity.

\section{Covariance}\label{sec:covariance}

We now examine how unitary transformations of states are reflected
in the outcome distributions of observables. For position, a spatial
translation shifts the outcome distribution by the same amount.
More generally, covariance expresses compatibility between a unitary
transformation and a transformation of the outcome space.
Conservation is the special case in which the outcome distribution
remains unchanged for every initial state.

It will be useful to introduce the following definition.

\begin{definition}[Symmetry of a measure space] A {\em symmetry} of the measure space $(\Om,\calF\!,\mu)$ is a measurable bijection $g:\Om\to\Om$, with measurable inverse, which leaves $\mu$ invariant, that is,
$$ (g(\mu))(F):=  \mu(g^{-1}(F)) = \mu(F), \quad  F\in \calF\!.$$
\end{definition}

If $g$ is a  symmetry of $(\Om,\calF\!,\mu)$, then for all $F\in \calF$ the set $g(F)$ is measurable and $\mu(F) = \mu(g^{-1}(g(F))) = \mu(g(F))$, that is, $\mu$ is also invariant under $g^{-1}$\!.

\begin{definition}[Conservation and covariance]\index{covariant}\index{conserved} Let $(\Om,\calF\!,\mu)$ be a measure space and let $U\in \calL(H)$ be unitary.
\begin{enumerate}[label={\rm(\roman*)}, leftmargin=*]
\item An observable $P:\calF\to \PP(H)$ is said to be {\em conserved} under $U$ if
$U P_F U^\star  = P_F$ for all $F\in \calF\!.$
\item An observable $P:\calF\to \PP(H)$ is said to be {\em covariant} under the pair $(g,U)$, where $g$ is a symmetry of $(\Om,\calF\!,\mu)$, if
$U P_F U^\star   = P_{g(F)}$ for all $F\in \calF\!,$
that is, if the following diagram commutes:
\begin{figure}[ht]
 \begin{center}
 \begin{tikzcd}[scale cd=1.1, sep=huge]
 \calF \arrow[r, "F\mapsto g(F)"] \arrow[d, "P"]
& \calF \arrow[d, "P"] \\
   \PP(H) \arrow[r, " A\mapsto UAU^\star"]
 & \PP(H)
\end{tikzcd}
\end{center}
\end{figure}
\end{enumerate}
\end{definition}
For example, position is covariant with respect to translations (an object at position $x$ appears at position $x-x'$ if the origin is translated by $x'$), and similarly momentum is covariant with respect to boosts (translations `in momentum space').

\subsection{Locally Compact Abelian Groups}

An interesting special case arises when observables take values in a locally compact abelian (LCA) group $G$.
Our treatment borrows some results from the theory of LCA groups that will not be proved here. For the special cases $\R^d$ and $\T$ the presentation is self-contained, as all missing details can be supplied using the results of
Chapter \ref{ch:operators}.
A fuller treatment should also cover the case of Lie groups such as $SO(3)$ and $SU(2)$, but this would take us too far afield.

Every LCA group $G$ admits a {\em Haar measure}\index{Haar measure},
that is, a nonzero Borel measure $\mu$ that is finite on compact
sets, satisfies $\mu(gB)=\mu(B)$ for all $g\in G$ and
$B\in\calB(G)$, and has the regularity properties as in the definition of a Radon measure, namely,
\begin{align*}
\mu(B)&=\inf\{\mu(O):\, B\subseteq O,\, O\hbox{ open}\}, \\*
\mu(O)&=\sup\{\mu(K):\, K\subseteq O,\, K\hbox{ compact}\},
\end{align*}
for Borel sets $B$ and open sets $O$.
This measure is unique up to multiplication by a positive scalar.
With respect to a Haar measure $\mu$, every $g\in G$ induces a symmetry on $G$ by
$$ g: g'\mapsto gg'\!, \quad g'\in G.$$
This induces a unitary operator $U_g$ on $L^2(G) := L^2(G,\mu)$ given by
$$U_g f = f \circ g^{-1}\!.$$
We refer to $U_g$ as the {\em translation by $g$}.\index{translation}\index{$U_g$}
We have
$$U_{g_1}U_{g_2}f = (f\circ g_2^{-1})\circ g_1^{-1} = f \circ (g_1 g_2)^{-1} = U_{g_1 g_2}f,$$
so $U_{g_1}U_{g_2} = U_{g_1 g_2}$. This means that the mapping $U:G\to\calL(L^2(G))$, $g\mapsto U_g$, is multiplicative and hence defines a unitary representation. Proceeding as in the proof of Proposition \ref{prop:Lp-transl}, one checks that this representation is strongly continuous.

A {\em character}\index{character} of $G$ is a continuous group homomorphism $\gamma:G\to \T$.
The set $\Gamma$ of all characters of $G$ is called the {\em Pontryagin
dual of $G$}. Endowed with the compact-open topology, it is a locally
compact abelian group.
The Pontryagin duality theorem asserts that every character of $\Gamma$
is of the form $\gamma\mapsto\gamma(g)$ for a unique $g\in G$, thereby
identifying the Pontryagin dual of $\Gamma$ with $G$ as an LCA group.

Every character $\gamma\in\Gamma$ induces a unitary operator on $L^2(G)$ via
$$V_\gamma f(g') = \gamma(g')f(g'), \quad g'\in G, \ f\in L^2(G).$$
We refer to $V_\gamma$ as the {\em boost by $\gamma$}.\index{boost}\index{$V_\gamma$}
In the language of Chapter \ref{ch:operators}, $V_\gamma$ is the pointwise multiplier by $\gamma$.

\smallskip
In the remainder of this section $G$ will always denote an LCA group with Pontryagin dual $\Gamma$.
The above definitions imply the so-called {\em Weyl commutation relations}\index{Weyl!commutation relation}:

\begin{proposition}[Weyl commutation relation]\label{prop:WeylCR}
For all $ g\in G$ and $\gamma\in \Gamma$
we have $$ V_\gamma U_g  =  \gamma(g) U_g V_\gamma .$$
\end{proposition}
\begin{proof} For $f\in L^2(G)$ and $g'\in G$ we have
\begin{align*}
 V_\gamma U_g f(g') & = \gamma(g')U_g f(g') = \gamma(g') f(g^{-1}g')
 \\ & = \gamma(g) \gamma(g^{-1}g') f(g^{-1}g')
 = \gamma(g) V_\gamma f(g^{-1}g') = \gamma(g) U_g V_\gamma f(g').
 \qedhere
\end{align*}
\end{proof}

We now turn to a pair of canonical observables associated with $G$ and
its dual $\Gamma$.
For the statement of the second part of the theorem we need the {\em Plancherel theorem},\index{theorem!Plancherel} which asserts that the Fourier transform $\calF:L^1(G)\to L^\infty(\Gamma)$ defined by
$$ \F f(\gamma):= \int_G f(g)\ov{\gamma(g)}\ud \mu(g), \quad \gamma\in\Gamma\!,$$
where $\mu$ is a Haar measure on $G$,
maps $L^1(G)\cap L^2(G)$ into $L^2(\Gamma)$ and there is a unique normalisation of the Haar measure of $\Gamma$ such that $\F$ extends to a unitary operator, called the Fourier--Plancherel transform, from $L^2(G)$ onto $L^2(\Gamma)$.
For later reference we
note that $\ov{\gamma(g)} = \gamma(g^{-1})$. For all $f\in L^1(G)$ and $g\in G$ we have $U_g f \in L^1(G)$ (by invariance of $\mu$) and
\begin{equation}\label{eq:transl-Fourier}
\begin{aligned}
\F U_g f(\gamma) &= \int_G f(g^{-1}g')\ov{\gamma(g')}\ud \mu(g')
 = \int_G f(g')\ov{\gamma(g g')}\ud \mu(g') =\gamma(g^{-1})\calF f(\gamma),
\end{aligned}
\end{equation}
where the second identity follows by substitution and invariance of $\mu$.
For $g\in G$ we write
$$
\chi_g(\gamma):=\gamma(g^{-1})=\ov{\gamma(g)},\quad \gamma\in\Gamma.
$$
Thus \eqref{eq:transl-Fourier} says that
\begin{align}\label{eq:Ug-chig}
\calF U_g f=\chi_g\,\calF f.
\end{align}

Let us call a $G$-valued observable
$P:\calB(G)\to\calP(H)$ {\em Radon}\index{Radon!observable}
if each of the measures $P_h$, $h\in H$, is Radon.
This is the case, for instance, if $P_h$ is regular for every
$h\in H$ and the topology of $G$ is countably generated
(inner regularity by compacts then follows from
Lemma \ref{lem:sigma-compact}).
Injectivity of the Fourier--Stieltjes transform (Theorem \ref{thm:uniq-FT-Rd}) extends to this setting; more precisely, if the finite Radon measures $\nu$ and $\nu'$ on $G$ have the same Fourier--Stieltjes transforms, that is, if for all $\gamma\in \Gamma$ we have
$ \int_G \gamma \ud \nu = \int_G \gamma \ud \nu',$
then $\nu = \nu'$.

\begin{theorem}[Position and momentum]\label{thm:Ugamma}
Let $G$ be an LCA group, let $\Gamma$ be its Pontryagin dual, and let $\calB(G)$ and $\calB(\Gamma)$ be their Borel $\sigma$-algebras. Then:
\begin{enumerate}[label={\rm(\arabic*)}, leftmargin=*]
\item\label{it:Ugamma1} there exists a unique $G$-valued Radon observable $X:\calB(G)\to \PP(L^2(G))$ such that for all $\gamma\in \Gamma$ we have
$$ V_\gamma =\int_G \gamma \ud X,$$
and it is given by $X_Bf = \one_B f$ for $f\in L^2(G)$ and $B\in \calB(G)$;

\item\label{it:Ugamma2} there exists a unique $\Gamma$-valued Radon observable
 $\Xi:\calB(\Gamma)\to \PP(L^2(G))$ such that for all $g\in G$ we have
$$
 U_g=\int_\Gamma \chi_g\ud\Xi,
$$
where $\chi_g(\gamma):=\gamma(g^{-1})$, and it is given by
$$
\Xi_B f=\calF^{-1}\one_B\calF f,\quad f\in L^2(G),\ B\in\calB(\Gamma).
$$
\end{enumerate}
\end{theorem}
\begin{proof}
\ref{it:Ugamma1}: \
Consider the $G$-valued observable $X:\calB(G)\to \calL(L^2(G))$ defined by
$$ X_B f:= \one_B f, \quad B\in \calB(G), \ f\in L^2(G).$$
Since $X_f(B)=\int_B|f|^2\ud\mu$ and $|f|^2\in L^1(G,\mu)$,
the finite measure $X_f$ is Radon by the regularity properties
of Haar measure.

Since $\int_G \one_B \ud X = X_B$ is the pointwise multiplier by $\one_B$,
the bounded functional calculus of $X$ gives that $\int_G\gamma\ud X$ is the pointwise multiplier by $\gamma$, that is, \ref{it:Ugamma1} holds.

If the $G$-valued Radon observable $\wt X$ has the same property, then for every $f\in L^2(G)$ the
finite Radon measures $X_f$ and $\wt X_f$ have the same
Fourier--Stieltjes transforms. Hence they are equal, and polarisation
gives $X_B=\wt X_B$ for every $B$.

\smallskip
\ref{it:Ugamma2}: \
Let $X^\Gamma$ be the multiplication projection-valued measure on
$L^2(\Gamma)$ and put
$$
\Xi_B:=\calF^{-1}X_B^\Gamma\calF\!.
$$
The same argument as above proves that $\Xi_f$ is Radon and $\int_\Gamma\chi_g\ud\Xi$
equals $\calF^{-1}M_{\chi_g}\calF$,
where $M_{\chi_g}$ is multiplication by $\chi_g$. By \eqref{eq:Ug-chig},  $\calF^{-1}M_{\chi_g}\calF = U_g$.
Uniqueness follows as in the first part.
\end{proof}

\begin{definition}[Position and momentum]\label{def:posmom} The $G$-valued observable $X$ and
the $\Gamma$-valued observable $\Xi$ of the  theorem are called the {\em position}\index{position operator} and {\em momentum}\index{momentum operator} observables of $G$, respectively.
\end{definition}

The special cases where $G = \R^d$ and $G = \mathbb{T}$ will be discussed in Sections \ref{subsec:pos-mom} and \ref{subsec:angular-momentum}, respectively.

\begin{proposition}\label{prop:LCA-pos-mom} Let $X$ and $\Xi$ be the position and momentum observables of $G$, respectively.
\begin{enumerate}[label={\rm(\arabic*)}, leftmargin=*]
 \item\label{it:LCA-covar1} $X$ is covariant with respect to every pair $(g,U_g)$ and conserved under every $V_\gamma$,
 $$ U_g X_B U_g^\star = X_{g B}, \qquad V_\gamma X_B V_\gamma^\star = X_B .$$
 \item\label{it:LCA-covar2} $\Xi$ is conserved under every $U_g$ and covariant with respect to every pair $(\gamma,V_\gamma)$,
 $$ U_g \Xi_B U_g^\star  =  \Xi_B , \qquad V_\gamma \Xi_B V_\gamma^\star = \Xi_{\gamma B}.$$
\end{enumerate}
\end{proposition}

\begin{proof}
\ref{it:LCA-covar1}: \  For all $g,g'\in G$, $f\in L^2(G)$, and $B\in \calB(G)$ we have
\begin{align*}
U_g X_B U_g^\star f(g')
= [\one_BU_{g^{-1}}f](g^{-1}g')
 = \one_B(g^{-1}g')f(g')  = X_{g B}f(g'),
\end{align*}
proving covariance with respect to $(g,U_g)$. Conservation under $V_\gamma$ is even simpler:
\begin{align*}
V_\gamma X_B V_\gamma^\star f =  \gamma(\cdot) \one_B(\cdot) \gamma^{-1}(\cdot)f(\cdot) = \one_B f = X_B f.
\end{align*}
The proof of \ref{it:LCA-covar2} is entirely similar.
\end{proof}

\begin{remark}[Complementarity]
In the physics literature, the `duality' between the position and momentum observables is referred to as {\em complementarity}.\index{complementarity}
\end{remark}

It is not an accident that position and momentum and the strongly continuous groups $U$ and $V$ associated with them satisfy both the Weyl commutation relations of Proposition \ref{prop:WeylCR} and the covariance properties of Proposition \ref{prop:LCA-pos-mom}.
The next proposition shows that this is an instance of a general
phenomenon.

\begin{proposition}[Covariance versus Weyl relations]\label{prop:covar-Weyl}
Let $P:\calB(G)\to\PP(H)$ and $\Pi:\calB(\Gamma)\to\PP(H)$ be $G$-valued and $\Gamma$-valued Radon observables, respectively, and let $U:G\to\mathscr U(H)$ and $V: \Gamma \to \mathscr U(H)$ be the associated strongly continuous unitary representations, that is,

$$
 U_g:=\int_\Gamma \chi_g(\gamma)\ud \Pi(\gamma),
 \qquad
 V_\gamma:=\int_G \gamma(g)\ud P(g),
$$
where $\chi_g(\gamma):=\gamma(g^{-1})$.
Then the following are equivalent:
\begin{enumerate}[label={\rm(\arabic*)}, leftmargin=*]
\item\label{it:covar-Weyl1} {\rm(Covariance on $G$)}\ $U_g P_B U_g^\star = P_{gB}$ for all $g\in G$ and $B\in\calB(G)$;
\item\label{it:covar-Weyl2} {\rm(Covariance on $\Gamma$)}\ $V_\gamma \Pi_B V_\gamma^\star = \Pi_{\gamma B}$ for all $\gamma\in \Gamma$ and $B\in\calB(\Gamma)$;
\item\label{it:covar-Weyl3} {\rm(Weyl relations)}\ $V_\gamma U_g = \gamma(g)\, U_g V_\gamma$ for all $g\in G$ and $\gamma\in\Gamma$.
\end{enumerate}
\end{proposition}
\begin{proof}
We prove the equivalence
\ref{it:covar-Weyl1}$\Leftrightarrow$\ref{it:covar-Weyl3};
the equivalence \ref{it:covar-Weyl2}$\Leftrightarrow$\ref{it:covar-Weyl3} is proved in the same way by interchanging the roles of $G$ and
$\Gamma$ and using that $\chi_g(\gamma^{-1}\gamma')=\gamma(g)\chi_g(\gamma')$.

Suppose first that
$U_gP_BU_g^\star=P_{gB}$ for all $B\in\calB(G)$.
By the bounded functional calculus, for every
bounded Borel function $f$ on $G$ we obtain
$$
U_g\Bigl(\int_Gf(g')\ud P(g')\Bigr)U_g^\star
=
\int_Gf(g^{-1}g')\ud P(g').
$$
Taking $f=\gamma$ gives
\begin{align*}
U_gV_\gamma U_g^\star =
\int_G\gamma(g^{-1}g')\ud P(g') = \gamma(g^{-1})V_\gamma,
\end{align*}
that is, the Weyl relations hold.

Conversely, suppose that the Weyl relations hold. Fix $g\in G$ and
define
$$
Q_B:=U_gP_BU_g^\star,
\qquad
R_B:=P_{gB}.
$$
Then
$$
\int_G\gamma(g')\ud Q(g')
=
U_gV_\gamma U_g^\star
=
\gamma(g^{-1})V_\gamma,
$$
while a change of variables gives
$$
\int_G\gamma(g')\ud R(g')
=
\gamma(g^{-1})V_\gamma.
$$
For every $h\in H$, the finite Radon measures $Q_h$ and $R_h$ therefore
have the same Fourier--Stieltjes transforms. By injectivity they are
equal. Polarisation gives $Q_B=R_B$ for every $B$, which is the
covariance relation.
\end{proof}

\subsection{The Case \texorpdfstring{$G = \R^d$}{G = R\textasciicircum d}: Position and Momentum}\label{subsec:pos-mom}

We now specialise to $G = \R^d$ with normalised Lebesgue measure ${\rm d}m(x) = (2\pi)^{-d/2}\ud x$ as the Haar measure.
Every character $\gamma:\R^d\to \T$ is of the form
$$\gamma(x) = e_\xi(x) := e^{ ix\cdot\xi}$$ for some $\xi\in \R^d\!$.
Under the identification $\gamma\leftrightarrow\xi$ we have $\Gamma = \R^d\!$, the Haar measure being the normalised Lebesgue measure $m$.
To distinguish $G=\R^d$ from its dual $\Gamma = \R^d$ we use Roman letters for elements of $G$ and Greek letters for elements of its dual $\Gamma$\!.

For every $y\in \R^d\!$, the translation $y: x\mapsto x+y$ is a symmetry of $\R^d\!$. The induced unitary operator $U_y$ on $L^2(\R^d\!,m)$ is right translation by $y$:
 $$ U_y f(x) = f(x-y), \quad x\in \R^d\!, \ f\in L^2(\R^d\!,m).$$
For every $\xi\in \R^d$ the boost $V_\xi$ on $L^2(\R^d\!,m)$ is given by
$$ V_\xi f(x) =  e^{i x\cdot\xi}f(x), \quad x\in \R^d\!,\ f\in L^2(\R^d\!,m).$$
The Weyl commutation relation takes the form
$$ V_\xi U_x  = e^{i x\cdot \xi}  U_x V_\xi, \quad x,\xi\in \R^d\!.   $$

The position observable $X:\mathscr{B}(\R^d)\to \PP(L^2(\R^d\!,m))$ and momentum observable $\Xi:\mathscr{B}(\R^d)\to \PP(L^2(\R^d\!,m))$ of Theorem \ref{thm:Ugamma} can be described as selfadjoint operators as follows.

For $1\le j\le d$, define the real-valued coordinate observables
$X_j$ and $\Xi_j$ by
$$
(X_j)_B
:=
X_{\{x\in\R^d:\, x_j\in B\}},
\qquad
(\Xi_j)_B
:=
\Xi_{\{\xi\in\R^d:\, \xi_j\in B\}},
\quad B\in\calB(\R).
$$
The projection-valued measure $X_j$ is the multiplication
projection-valued measure associated with the coordinate function
$x\mapsto x_j$. Its associated selfadjoint operator is therefore
$$
\wh x_jf(x):=x_jf(x),
$$
with domain
$
\Dom(\wh x_j)
=
\{f\in L^2(\R^d\!,m):\, x_jf\in L^2(\R^d\!,m)\}$.

Similarly, we have
$\Xi_j=\calF^{-1}X_j\calF$
after identifying the dual copy of $\R^d$ with $\R^d$, and its
associated selfadjoint operator is
$$
\wh\xi_j
=
\calF^{-1}M_{\xi_j}\calF
=
\frac1i\partial_j,
$$
with domain
$\Dom(\wh\xi_j) =\{f\in L^2(\R^d\!,m):\,\partial_jf\in L^2(\R^d\!,m)\}$.

Direct calculation gives the
{\em Heisenberg commutation relation}\index{Heisenberg!commutation relation}
\begin{align*}
 \wh x_j\wh \xi_k  - \wh\xi_k\wh x_j   =  i \delta_{jk} I,
\end{align*}
on $C_{\rm c}^\infty(\R^d)$. By the weak derivative identity
$\partial_k(x_jf)=\delta_{jk}f+x_j\partial_kf$,
the same relation holds for every
$f\in\Dom([\wh x_j,\wh\xi_k])$.

For pure states $\phi$ represented by a norm-one function $f\in L^2(\R^d\!,m)$ such that $f\in \Dom(\wh x_j)\cap \Dom(\wh\xi_j)$ and $\wh x_jf, \wh\xi_jf\in \Dom(\wh x_j)\cap \Dom(\wh\xi_j)$, the uncertainty principle of Theorem \ref{thm:uncertainty} takes the form
$$ \Delta_\phi(\wh x_j)\Delta_\phi(\wh\xi_j) \ge \frac1{2}.$$

It follows from Proposition \ref{prop:LCA-pos-mom} that the position observable
$X$ is covariant with respect to translations and conserved under boosts, and that the momentum observable $\Xi$ is conserved under translations and covariant with respect to boosts. The next theorem shows that this essentially characterises these observables.
For the proof we will use the following lemma.

\begin{lemma}\label{lem:transl-boost-commutant}
If $T\in\calL(L^2(\R^d\!,m))$ commutes with all translations $U_x$ and
all boosts $V_\xi$, then $T=cI$ for some $c\in\C$.
\end{lemma}

\begin{proof}
Since $T$ commutes with every boost $V_\xi=M_{e_\xi}$,
Lemma \ref{lem:char-mult} gives $T=M_b$ for some
$b\in L^\infty(\R^d)$. Commutation with translations then implies
$$
b(\cdot-x)=b
\quad\hbox{in }L^\infty(\R^d),\quad x\in\R^d.
$$
Choose $\rho\in C_{\rm c}^\infty(\R^d)$ with $\rho\ge0$ and
$\int_{\R^d}\rho\ud x=1$, and put $\rho_n(x):=n^d\rho(nx)$.
Each function $b*\rho_n$ is continuous and invariant under every
translation, and is therefore constant. Moreover,
$b*\rho_n\to b$ in $L^1_{\rm loc}(\R^d)$ by mollification, applied after localisation.
Convergence on any fixed ball of positive measure shows that
these constants converge to some $c\in\C$. Hence $b=c$ almost
everywhere, and therefore $T=cI$.
\end{proof}

\begin{theorem}[Covariance characterisation of position and momentum]\label{thm:pos-mom}
Up to conjugation by a translation and a boost, respectively, position and momentum are characterised by their covariance and conservation properties. More precisely, denoting by $X$ and $\Xi$ the position and momentum observables, the following assertions hold.
\begin{enumerate}[label={\rm(\arabic*)}, leftmargin=*]
\item\label{it:pos-mom1} if the observable $P: \calB(\R^d)\to \PP(L^2(\R^d\!,m))$ is covariant with respect to all pairs $(x,U_x)$, $x\in \R^d\!$, and conserved under all boosts $V_\xi$, $\xi\in \R^d\!$, then there exists a unique $y\in \R^d$ such that
 $$P_B  = U_{y} X_B U_y^\star\!, \quad B\in \calB(\R^d);$$
\item\label{it:pos-mom2} if the observable $P: \calB(\R^d)\to \PP(L^2(\R^d\!,m))$ is conserved under all translations $U_x$, $x\in \R^d\!$, and covariant with respect to all pairs $(\xi,V_\xi)$, $\xi\in \R^d\!$, then there exists a unique $\eta\in \R^d$ such that
 $$P_B  = V_{\eta} \Xi_B V_{\eta}^\star,  \quad B\in \calB(\R^d).$$
\end{enumerate}
\end{theorem}

\begin{proof}
\ref{it:pos-mom1}:\ For $\xi\in\R^d$ put
$$
W_\xi:=\int_{\R^d}e^{ix\cdot\xi}\ud P(x).
$$
The properties of the functional calculus readily imply that the resulting mapping $W: \R^d \to \calL(L^2(\R^d\!,m))$ is a strongly
continuous unitary representation.
Since $P$ is conserved under every boost, $W_\xi$ commutes with every
$V_\eta$. Moreover, covariance under translations gives
$$
U_yW_\xi U_y^\star=e^{-iy\cdot\xi}W_\xi,
$$
or equivalently,
$$
W_\xi U_y=e^{iy\cdot\xi}U_yW_\xi.
$$
The same Weyl relation holds for $V_\xi$. Consequently,
$C_\xi:=V_\xi^\star W_\xi$
commutes with every $U_y$ and with every $V_\eta$. Hence, by Lemma~\ref{lem:transl-boost-commutant}, there exists
$c(\xi)\in\C$ such that $C_\xi=c(\xi)I$, or equivalently
$$
W_\xi=c(\xi)V_\xi.
$$
Since $C_\xi=V_\xi^\star W_\xi$ is unitary, we have $
I=C_\xi^\star C_\xi=|c(\xi)|^2I$.
It follows that $|c(\xi)|=1$.

Since both $W$ and $V$ are strongly continuous representations,
$c:\R^d\to\T$ is a continuous character. Hence there exists a unique
$y\in\R^d$ such that
$$
c(\xi)=e^{-iy\cdot\xi},
\quad \xi\in\R^d\!.
$$
It follows that
$$
W_\xi
=
e^{-iy\cdot\xi}V_\xi
=
U_yV_\xi U_y^\star.
$$

Put
$$
Q_B:=U_yX_BU_y^\star,
\quad B\in\calB(\R^d).
$$
Then
$$
\int_{\R^d}e^{ix\cdot\xi}\ud Q(x)
=
U_yV_\xi U_y^\star
=
W_\xi
=
\int_{\R^d}e^{ix\cdot\xi}\ud P(x).
$$
Thus, for every $h\in L^2(\R^d\!,m)$, the finite measures $P_h$ and
$Q_h$ have the same Fourier--Stieltjes transform. By injectivity of the
Fourier--Stieltjes transform, $P_h=Q_h$. Polarisation gives $P_B=Q_B$
for every $B\in\calB(\R^d)$. Hence
$$
P_B=U_yX_BU_y^\star, \quad B \in \mathscr B(\R^d).
$$
The uniqueness of $y$ follows from the uniqueness of the character
$\xi\mapsto e^{-iy\cdot\xi}$. This completes the proof of the first assertion.

\smallskip
\ref{it:pos-mom2}:\ The observable $\widetilde P$ defined by
$$
\widetilde P_B:=\calF P_B\calF^{-1}, \quad B \in \mathscr B(\R^d),
$$
is conserved under boosts and covariant
under translations. By \ref{it:pos-mom1}, there exists a unique
$\eta\in\R^d$ such that $\widetilde P_B=U_\eta X_BU_\eta^\star$.
Conjugating with $\calF^{-1}$ gives
$P_B = V_\eta\Xi_BV_\eta^\star$.
\end{proof}

\subsection{The Case \texorpdfstring{$G = \T$}{G = T}: Angle and Angular Momentum}\label{subsec:angular-momentum}

The results of the preceding subsection have natural analogues for the unit circle $\T$.
We identify $\T$ with the unit circle of $\C$ and take the normalised Lebesgue measure on $\T$ as the Haar measure. Every character $\gamma:\T\to \T$ is of the form
$$ \gamma(z) = z^k\!, \quad z\in \T,$$ for some $k\in \Z$.
Under this identification we have $\Gamma = \Z$, its normalised Haar measure being the counting measure.

For every $w\in \mathbb{T}$ the rotation
$z\mapsto wz$ is a symmetry of $\T$.
The induced unitary operator $U_w$ on $L^2(\T,m)$ is given by
\begin{align}\label{def-Uw} U_w f(z) = f(w^{-1}z) , \quad z\in \T, \ f\in L^2(\T).
\end{align}
For every $k\in \Z$ the boost $V_k$ on $L^2(\T)$ is given by
\begin{align}\label{def-Vk} V_k f(z) = z^{k}f(z), \quad z\in \T, \  f\in L^2(\T).
\end{align}

The Weyl commutation relation takes the form
$$ V_k U_w  = w^k U_w V_k, \quad w\in \T, \ k\in\Z. $$

The position and momentum observables for $\T$ associated with the unitary operators $U_z$ and $V_k$ are denoted by $\theta$ and $L$ and are called the
{\em angle}\index{angle}
and (orbital) {\em angular momentum}\index{angular momentum}
observables. They take values in $\T$ and $\Z$, respectively; in particular, angular momentum can only assume discrete values. In the physics literature one speaks about `quanta' of angular momentum.

By Proposition \ref{prop:LCA-pos-mom}, $\theta$ is covariant under every pair $(z,U_z)$ and conserved under every $V_k$,
and $L$ is conserved under every $U_z$ and covariant under every pair $(k,V_k)$.
Repeating the proof of Theorem \ref{thm:pos-mom} almost {\em verbatim} gives the following result.

\begin{theorem}[Covariance characterisation of angle and angular momentum]\label{thm:angle} Up to conjugation by a rotation and a boost, respectively, angle and angular momentum are characterised by their covariance and conservation properties. More precisely, denoting by
$\theta$ and $L$ the position and momentum observables associated with the rotations $U_z$ and boosts $V_k$ given by \eqref{def-Uw} and \eqref{def-Vk}, the following assertions hold.
\begin{enumerate}[label={\rm(\arabic*)}, leftmargin=*]
 \item if the observable $P: \calB(\T)\to \PP(L^2(\T))$ is covariant with respect to all pairs $(z,U_z)$, $z\in\T$, and conserved under all $V_k$, $k\in \Z$, then there exists a unique $w\in \T$ such that
 $$P_B  = U_w \theta_B U_w^\star,  \quad B\in \calB(\T); $$
 \item if the observable $P: \calB(\Z)\to \PP(L^2(\T))$ is conserved under all $U_z$, $z\in \T$, and covariant with respect to all pairs $(k,V_k)$, $k\in\Z$, then there exists a unique $j\in\Z $ such that
 $$P_B  = V_j L_B V_j^\star\!,  \quad B\in \calB(\Z).$$
\end{enumerate}
\end{theorem}

\subsection{The Time--Energy Pair}

In the physics literature, `time' and `energy' are treated as complementary observables much in the same way as position and momentum. For a system with a semibounded Hamiltonian, however, there is
no selfadjoint time operator whose spectral measure is covariant
under time evolution in the sense described below.
Nevertheless, both the Heisenberg commutation relation and the uncertainty principle are often applied to this pair in an informal operational way.

In this section we formulate this obstruction in terms of covariance. As in Proposition \ref{prop:covar-Weyl}, covariance leads to a Weyl commutation relation, which in general cannot
be deduced from a Heisenberg commutation relation on a dense domain. Given a (possibly unbounded) selfadjoint `energy' operator $E$ acting in a Hilbert space $H$, let $U(t):= e^{-itE}$, $t\in \R$, be the unitary $C_0$-group on $H$ generated by $-iE$. In physical terms, we think of this group as describing the time evolution of a quantum system with state space $H$ governed by the operator $E$, which we interpret as the Hamiltonian of the system. A projection-valued measure $T: \calB(\R)\to \PP(H)$ may then be called a `time observable' if it is covariant with the evolution $U$ in the sense that
\begin{align}\label{eq:cov}
U(t)T_B U^\star(t)=T_{B+t}, \quad t\in\R, \ B\in\calB(\mathbb R).
\end{align}
Mathematically, this is just saying that we treat energy and time in the same way as position and momentum. Physically, \eqref{eq:cov} can be motivated as follows.
For a unit vector $h\in H$ (or rather, a pure state $\ket{h}$), we interpret the numbers
$$
  \nu_h(B):=\iprod{T_Bh}{h}
$$
as the probability that a clock reading produces an outcome belonging to the Borel set $B\in\mathscr{B}(\R)$.
Waiting $t$ time units should shift the distribution of outcomes by $t$. Indeed,
if one prepares $h$ at time $0$ and performs the same measurement after a waiting time $t$, then the probability at time $0$ of obtaining an outcome in the Borel set $B$ should equal the probability at time $t$
of obtaining an outcome in the shifted set $B+t$, that is,
\begin{equation}\label{eq:dist-shift}
  \nu_h(B) = \nu_{U(t)h}(B+t),\quad t\in\mathbb R,\ B\in\calB(\mathbb R).
\end{equation}
Using the unitarity of $U(t)$, we can rewrite the right-hand side of \eqref{eq:dist-shift} as
$$
  \nu_{U(t)h}(B+t)
  =\iprod{T_{B+t}U(t)h}{U(t)h}
  =\iprod{U^\star(t)T_{B+t}U(t)h}{h}.
$$
Thus \eqref{eq:dist-shift} is equivalent to
$$
  \iprod{U^\star(t)T_{B+t}U(t)h}{h}
  =\iprod{T_Bh}{h}.
$$
Since this holds for all unit vectors $h\in H$, we obtain the covariance relation
\begin{equation*}
  U^\star (t)T_{B+t}U(t)=T_B,\quad t\in\mathbb R,\ B\in\calB(\mathbb R),
\end{equation*}
which is equivalent to \eqref{eq:cov}.

A celebrated no-go theorem due to Pauli, presented here in the language of covariance, gives the mathematical
formulation of the nonexistence of a time operator in the above sense. It implies that if the selfadjoint operator $E$ is {\em semibounded} (by which we mean that for some $\om\in\R$ we have $E \ge \om I$, or equivalently, $\sigma(E)\subseteq [\om,\infty)$), there is no projection-valued measure $T$ covariant with $E$.

\begin{theorem}[Pauli's no-go theorem]\label{thm:Pauli-PVM}\index{theorem!Pauli, no-go} Let $E$ be a selfadjoint operator acting in $H$.
If there exists a \emph{projection-valued} measure $T:\calB(\R)\to\PP(H)$ such that
$$
U(t)\,T_{B}\,U^\star(t) =  T_{B+t},\quad t\in\R,\ B\in\calB(\R),
$$
then the spectrum $\sigma(E)$ is translation-invariant, that is, $\sigma(E)=\sigma(E)+s$ for every $s\in\R$.
As a consequence, if $H\not=\{0\}$, then $\sigma(E) = \R$.
\end{theorem}

\begin{proof}
By the spectral theorem, there is a unique selfadjoint operator $A$
associated with $T$. Put $V(t):=e^{itA}$. The covariance of $T$
and the bounded functional calculus give
$$
U(s)V(t)U^\star(s)=e^{-ist}V(t),
$$
or equivalently,
$$
V^\star(t)U(s)V(t)
=e^{-its}U(s)
=e^{-is(E+tI)},\quad s\in\R.
$$
The generators of the groups on the left and right are
$-iV^\star(t)EV(t)$ and $-i(E+tI)$, respectively.
By uniqueness of the generator,
$$
V^\star(t)EV(t)=E+tI
$$
with equality of domains. Since $V(t)$ is unitary, it follows that
$$
\sigma(E)
=\sigma(V^\star(t)EV(t))
=\sigma(E+tI)
=\sigma(E)+t.
$$
The only subsets of $\R$ invariant under all translations are
the empty set and $\R$ itself. If $H\ne\{0\}$, then
$\sigma(E)\ne\emptyset$ by Proposition
\ref{prop:sp-normal-nonempty}, and therefore $\sigma(E)=\R$.
\end{proof}

Pauli's theorem rules out covariant {\em sharp} time observables for semibounded Hamiltonians. If projection-valued measures are replaced by POVMs, covariant time observables can in fact be constructed in many situations, often by compressing a covariant projection-valued measure acting on a larger Hilbert space.

A periodic version of this phenomenon has already appeared in the phase--number pair. Up to unitary equivalence, the number operator $N$ may be viewed as the Hamiltonian of the one-dimensional quantum harmonic oscillator after subtraction of its ground-state energy. Adding a scalar multiple of the identity to a Hamiltonian changes its unitary evolution only by a scalar phase and therefore does not affect covariance relations for observables. Thus the phase POVM $\Phi$ of Section \ref{subsec:number-phase}, which was seen to be covariant with respect to the number operator,
may equally well be interpreted as a periodic time observable covariant with respect to the harmonic oscillator: it records the position within an oscillation cycle, with time understood modulo the period.

The following example is a canonical model of a covariant real-valued time observable; depending on the physical system, the measured time may represent a clock reading, an event
time, or another operationally defined temporal quantity.

\begin{example}[A covariant time POVM]
Let $\wt H=L^2(\R)$, let $X$ and $\Xi$ be the position and momentum
observables, and let $\wh\xi$ be the selfadjoint operator associated with $\Xi$.
Put $U(t):=e^{-it\wh\xi}$ and $P:=\Xi_{[0,\infty)}$.
Let $H:=P\wt H$, and let $J:H\to\wt H$ be the inclusion.
Then $H$ is invariant under every $U(t)$, and the restriction
$E:=\wh\xi|_H$, with domain $\Dom(E)=\Dom(\wh\xi)\cap H$,
is selfadjoint with $\sigma(E)=[0,\infty)$. Define
$$
T_B:=J^\star X_BJ,\quad B\in\calB(\R).
$$
By Proposition \ref{prop:POVM-compression}, $T$ is a POVM. Since
$U(t)H=H$ and $U(t)X_BU^\star(t)=X_{B+t}$, we have
$$
e^{-itE}T_Be^{itE}=T_{B+t}, \quad t\in \R, \ B \in \calB(\R).
$$
Thus $T$ is a covariant time observable for the semibounded
Hamiltonian $E$. By Theorem \ref{thm:Pauli-PVM}, it cannot be
projection-valued.
\end{example}

\section{The Stone--von Neumann Theorem}\label{subsec:StonevonNeumann}

We have seen in Theorem \ref{thm:pos-mom} that the $\R^d$-valued position and momentum observables are uniquely determined, up to conjugation by a translation and a boost, respectively, by their transformation properties under translations and boosts. We have already observed in  Proposition \ref{prop:covar-Weyl} that both the covariance relation for position
$$
U_x X_B U_{x}^\star = X_{x+B}, \quad x\in \R^d\!,\ B\in\calB(\R^d),
$$
and the covariance relation for momentum
$$
V_\xi \Xi_B V_{\xi}^\star = \Xi_{\xi+B}, \quad \xi\in \R^d\!,\ B\in\calB(\R^d)
$$ imply the Weyl commutation relation and vice versa.
In view of this, it is reasonable to ask the following: if $H$ is a Hilbert space and $\wt X: \calB(\R^d)\to \PP(H)$ and $\wt \Xi:\calB(\R^d)\to \PP(H)$ are projection-valued measures such that the associated strongly continuous unitary groups
$\wt U,\wt V: \R^d \to \mathscr{U}(H)$ on $H$ satisfy the Weyl commutation relations, does it follow that $\wt U$ and $\wt V$ are unitarily equivalent to the translation representation $U$ and the boost representation $V$? Modulo an irreducibility assumption the answer to this question is affirmative; this is a celebrated theorem due to Stone and von Neumann (Theorem \ref{thm:Stone-vonNeumann} and its corollary). Proving this theorem is our next objective.

We start with some preparations. Suppose that $\wt U, \wt V:\R^d\to \calL(H)$ are strong\-ly continuous unitary representations of $\R^d$ on a Hilbert space $H$ such that the Weyl commutation relation holds, that is,
\begin{align}\label{eq:WCR}
\wt V_\xi \wt U_x  = e^{ ix\cdot \xi} \wt U_x \wt V_\xi, \quad x,\xi\in \R^d\!.
\end{align}
The relation \eqref{eq:WCR} states that $\wt U$ and $\wt V$ `commute up to a multiplicative scalar of modulus one'. This suggests interpreting \eqref{eq:WCR} as a `projective' unitary representation\index{projective representation} of $\R^d\times\R^d$ on $H$.
The following construction turns this projective representation into a
genuine unitary representation.
Consider the unitary operators
\begin{align}\label{eq:WWeyl} \wt W(x,\xi):= e^{\frac12 ix\cdot \xi}\wt U_x \wt V_\xi = e^{-\frac12 ix\cdot \xi}\wt V_\xi \wt U_x, \quad x,\xi\in \R^d\!.
\end{align}
The operators $\wt W(x,\xi)$ defined by \eqref{eq:WWeyl} satisfy
\begin{equation}\label{eq:WW}
\begin{aligned}
\wt W(x,\xi)\wt W(x'\!,\xi')
 & =  e^{-\frac12 i(x\cdot\xi+x'\cdot\xi')} \wt V_\xi \wt U_x\wt V_{\xi'} \wt U_{x'}
 \\ & =  e^{-\frac12 i(x\cdot\xi+x'\cdot\xi') - ix\cdot\xi'}\wt V_\xi \wt V_{\xi'}\wt U_x \wt U_{x'}
 \\ & = e^{\frac12 i(x'\cdot\xi - x\cdot\xi')} e^{-\frac12 i(x+x')\cdot(\xi+\xi')}\wt V_{\xi+\xi'} \wt U_{x+x'}
\\ & = e^{\frac12 i(x'\cdot\xi - x\cdot\xi')} \wt W(x+x'\!,\xi+ \xi').
\end{aligned}
\end{equation}
From this it follows that the unitary operators defined by
\begin{align*}
\wt W(x,\xi,t) := e^{it} \wt W(x,\xi)
\end{align*}
satisfy
\begin{equation}\label{eq:Heisenberg-repr}
\begin{aligned} \wt W(x,\xi,t) \wt W(x'\!,\xi'\!,t')
 & =  e^{i(t+t')} \wt W(x,\xi)\wt W(x'\!,\xi')
\\ & = e^{i(t+t' + \frac12(x'\cdot\xi-x\cdot\xi'))} \wt W(x+x'\!,\xi+\xi')
\\ & = \wt W((x,\xi,t)\circ (x'\!,\xi'\!,t')),
\end{aligned}
\end{equation}
where
\begin{align}\label{eq:Heisenberg-mult} (x,\xi,t)\circ (x'\!,\xi'\!,t'):= \Bigl(x+x'\!,\ \xi+\xi'\!,\ t+t' + \frac12(x'\cdot\xi-x\cdot\xi')\Bigr).
\end{align}
 One easily checks that the operation $\circ$ turns $\H^d:=\R^d\times\R^d\times \R$\index{$Hh$@$\H^d$} into a group.

\begin{definition}[Heisenberg group]
The {\em Heisenberg group}\index{Heisenberg!group} in dimension $d$ is the group $\H^d:=\R^d\times\R^d\times \R$ with composition law \eqref{eq:Heisenberg-mult}.
\end{definition}

The identity \eqref{eq:Heisenberg-repr} shows that $\wt W$ defines a
unitary representation of the Heisenberg group $\H^d$ on $H$. It is strongly continuous and it satisfies
\begin{align*}
\wt W(0,0,t) = e^{it}I, \quad t\in \R.
\end{align*}

\begin{definition}[Schr\"odinger representation]\label{def:Schroedinger-representation} The {\em Schr\"odinger representation}\index{Schr\"odinger!representation} is the unitary representation $W: \H^d\to \calL(L^2(\R^d\!,m))$ arising
in the special case where the unitary representations $U,V:\R^d\to \calL(L^2(\R^d\!,m))$ are given by translations and boosts, respectively.
\end{definition}

Explicitly, the Schr\"odinger representation is given by
\begin{align}\label{eq:Schrod-repr} W(x,\xi,t)f(x') =
 e^{ it} e^{-\frac12 ix\cdot \xi} e^{i x'\cdot \xi}f(x'-x).
\end{align}

\begin{proposition}\label{prop:Schrod-irred}
The\, Schr\"odinger\, representation\, is\, {\em irreducible},\index{irreducible} that\! is, the\, only\, closed subspaces of $L^2(\R^d\!,m)$ invariant under the action of $W$ are the trivial subspaces $\{0\}$ and $L^2(\R^d\!,m)$.
\end{proposition}
 \begin{proof} Let $Y$ be a closed $W$-invariant
subspace. Then the orthogonal projection $P_Y$ onto $Y$ commutes with every $W(x,\xi,t)$, hence with all $U_x$ and $V_\xi$.
By Lemma \ref{lem:transl-boost-commutant}, $P_Y=cI$ for some $c\in\{0,1\}$, so $Y=\{0\}$ or $Y=L^2(\R^d\!,m)$.
\end{proof}

We can now state and prove the main result of this section.

\begin{theorem}[Stone--von Neumann]\label{thm:Stone-vonNeumann}\index{theorem!Stone--von Neumann}
Let $\wt W: \H^d\to \calL(H)$ be a strongly continuous unitary representation of $\H^d$ on a nonzero separable Hilbert space $H$. If $\wt W$ is irreducible and satisfies $\wt W(0,0,t) = e^{it}I$ for all $t\in \R$, then
$\wt W$ is unitarily equivalent to the Schr\"odinger representation $W$.
More precisely, there exists a unitary operator $S: L^2(\R^d\!,m)\to H$ such that
$$\wt W(x,\xi,t) = SW(x,\xi,t)S^\star\!, \quad (x,\xi,t)\in \H^d\!.$$
The operator $S$ is unique up to a multiplicative scalar of modulus one.
\end{theorem}

Here, we use the term {\em unitary operator}\index{unitary!operator} for an operator $S: H\to K$, where $H$ and $K$ are Hilbert spaces, such that $S^\star S = I$ (the identity operator on $H$) and $SS^\star = I$ (the identity operator on $K$).

We have the following immediate corollary for representations arising from  pairs of unitary representations satisfying the Weyl commutation relation.

\begin{corollary} Let $\wt U, \wt V:\R^d\to\calL(H)$ be strongly continuous unitary representations on a nonzero separable Hilbert space $H$ satisfying the Weyl commutation relation
$$ \wt V_\xi \wt U_x  = e^{ix\cdot \xi} \wt U_x \wt V_\xi, \quad x,\xi\in \R^d\!.$$
Suppose furthermore that the family $\{\wt U_x, \wt V_\xi:\, x\in \R^d\!,\, \xi\in\R^d\} $ acts irreducibly on $H$
in the sense that the only closed subspaces invariant under all operators $\wt U_x$ and $\wt V_\xi$, $x\in \R^d\!,\, \xi\in\R^d\!$, are the trivial subspaces $\{0\}$ and $H$.
Then there exists a unitary operator $S:L^2(\R^d\!,m)\to H$ such that
\begin{align*}
   \wt U_x   & = S U_x S^\star\!,\quad x\in\R^d\!,
\\ \wt V_\xi & = S V_\xi S^\star\!, \quad \xi\in\R^d\!,
\end{align*}
where $U$ and $V$ are the translation and boost representations on $L^2(\R^d\!,m)$, respectively.
The operator $S$ is unique up to a multiplicative constant of modulus one.
\end{corollary}

\begin{proof}
Define $\wt W(x,\xi):=e^{\frac12 i x\cdot\xi}\wt U_x\wt V_\xi$ and $\wt W(x,\xi,t):=e^{it}\wt W(x,\xi)$.
The Weyl relation gives the product formula
$$
   \wt W(x,\xi,t)\wt W(x',\xi',t')
   =
   \wt W(x+x',\xi+\xi',
      t+t'+\frac12(x'\cdot\xi-x\cdot\xi')),
$$
so $\wt W$ is a unitary representation of $\H^d$. Its strong continuity
follows from the strong continuity of $\wt U$ and $\wt V$, and we have
$\wt W(0,0,t)=e^{it}I$ for $t\in\R$.
Moreover, a closed subspace is invariant under all $\wt W(x,\xi,t)$ if and
only if it is invariant under all $\wt U_x$ and $\wt V_\xi$. Thus $\wt W$
is irreducible. Hence Theorem \ref{thm:Stone-vonNeumann} gives a unitary $S:L^2(\R^d\!,m)\to H$ such that
$$
   \wt W(x,\xi,t)=S W(x,\xi,t)S^\star .
$$
Taking $(x,\xi,t)=(x,0,0)$ and $(0,\xi,0)$ gives the stated identities for
$\wt U_x$ and $\wt V_\xi$. The uniqueness follows from the uniqueness part
of the theorem.
\end{proof}

We now fix a strongly continuous unitary representation
$\wt W:\H^d\to\calL(H)$ satisfying
$$
        \wt W(0,0,t)=e^{it}I,\quad t\in\R .
$$
Put
$$
   \wt U_x:=\wt W(x,0,0),\qquad
   \wt V_\xi:=\wt W(0,\xi,0),\qquad
   \wt W(x,\xi):=\wt W(x,\xi,0).
$$
Then the group law of $\H^d$ gives
$$
   \wt W(x,\xi)=e^{\frac12 i x\cdot\xi}\wt U_x\wt V_\xi
              =e^{-\frac12 i x\cdot\xi}\wt V_\xi\wt U_x,
$$
and hence \eqref{eq:WCR}--\eqref{eq:Heisenberg-repr} hold again.

We write $m$ for both the normalised Lebesgue measures on $\R^d$ and $\R^{2d}$.

\begin{definition}[Weyl transform]\index{Weyl!transform} For $a\in L^1(\R^{2d}\!,m)$
we define the operator $\wt W(a)\in \calL(H)$ by
\begin{align*}
\wt W(a)h
 & := \int_{\R^{2d}} a(x,\xi) \wt W(x,\xi)h \ud m(x)\ud m(\xi), \quad h\in H,
\end{align*}
where the integral is a Bochner integral in $H$.
\end{definition}

The next two lemmas state some properties of the Weyl transform associated with the Schr\"odinger representation $W:\H^d\to \calL(L^2(\R^d\!,m))$.

\begin{lemma}\label{lem:Moyal-injective}
For all $a\in L^1(\R^{2d}\!,m)\cap L^2(\R^{2d}\!,m)$ the operator $W(a)$ is Hilbert--Schmidt on $L^2(\R^{d}\!,m)$ and
$$ \n W(a)\n_{\calL_2(L^2(\R^{d}\!,m))} = \n a\n_{L^2(\R^{2d}\!,m)}.$$
\end{lemma}
\begin{proof}
Define
$$
k(x,x')
:=
\int_{\R^d}a(x'-x,\xi)
e^{\frac12 i(x+x')\cdot\xi}\ud m(\xi).
$$
Under the change of variables
$$
u:=\frac{x+x'}2,\qquad v:=x'-x,
$$
whose Jacobian has absolute value one, we have
$$
k\Bigl(u-\frac v2,u+\frac v2\Bigr)
=
\int_{\R^d}a(v,\xi)e^{iu\cdot\xi}\ud m(\xi).
$$
For almost every $v$, the right-hand side is the inverse Fourier
transform of $a(v,\cdot)$. Hence Plancherel's theorem and Fubini's
theorem give
\begin{align*}
\int_{\R^{2d}}|k(x,x')|^2\ud m(x)\ud m(x')
&=
\int_{\R^{2d}}
\Bigl|\int_{\R^d}a(v,\xi)e^{iu\cdot\xi}\ud m(\xi)\Bigr|^2
\ud m(u)\ud m(v)\\
&=
\int_{\R^{2d}}|a(v,\xi)|^2\ud m(v)\ud m(\xi)
=\n a\n_2^2.
\end{align*}
Thus $k$ defines a Hilbert--Schmidt integral operator, by
Example \ref{ex:kernel-HS} with the kernel variables interchanged.
For $f\in C_{\rm c}(\R^d)$, the formula
\eqref{eq:Schrod-repr}, Fubini's theorem, and the substitution
$x=x'-u$ give
\begin{equation}\label{eq:W(a)f}
\begin{aligned}
(W(a)f)(x')
&=
\int_{\R^{2d}}
a(u,\xi)e^{i(x'-u/2)\cdot\xi}f(x'-u)
\ud m(u)\ud m(\xi) \\ & =
\int_{\R^d}k(x,x')f(x)\ud m(x).
\end{aligned}
\end{equation}
Since $C_{\rm c}(\R^d)$ is dense in $L^2(\R^d\!,m)$ and both
operators are bounded, this identifies $W(a)$ with the integral
operator defined by $k$. Consequently,
\[
\n W(a)\n_{\calL_2(L^2(\R^d\!,m))}
=\n k\n_{L^2(\R^{2d}\!,m)}
=\n a\n_{L^2(\R^{2d}\!,m)}.
\qedhere\]
\end{proof}

Since $L^1(\R^{2d}\!,m)\cap L^2(\R^{2d}\!,m)$  is dense in $L^2(\R^{2d}\!,m)$, the lemma implies that the mapping $W: a\mapsto W(a)$ has a unique extension to an isometry from $L^2(\R^{2d}\!,m)$ into the space of Hilbert--Schmidt operators  $\calL_2(L^2(\R^{d}\!,m))$. This extension is again denoted by $W$.

A special role is played by the functions
\begin{align*}a_0 (x,\xi) & := \exp\Bigl(-\frac14(|x|^2+|\xi|^2)\Bigr), \quad x,\xi\in \R^d\!,\\
\phi_0(x) &:= 2^{d/4}  \exp\Bigl(-\frac12|x|^2\Bigr), \quad x\in \R^d\!.
\end{align*}
Note that $\n a_0\n_{L^2(\R^{2d}\!,m)} = \n \phi_0\n_{L^2(\R^d\!,m)} = 1$.

\begin{lemma}\label{lem:Schrod-rankone}
The operator $W(a_0)$ equals the rank-one orthogonal projection $\phi_0\,\bar\otimes\, \phi_0$.
\end{lemma}
\begin{proof}
By \eqref{eq:W(a)f}, the operator $W(a_0)$ is the integral operator with
kernel
$$
k_0(x,x')
=
\int_{\R^d}a_0(x'-x,\xi)
e^{\frac12i(x+x')\cdot\xi}\ud m(\xi).
$$
By Lemma~\ref{lem:Gauss},
\begin{align*}
k_0(x,x')
&=
\exp\Bigl(-\frac14|x'-x|^2\Bigr)
\int_{\R^d}
\exp\Bigl(-\frac14|\xi|^2\Bigr)
e^{\frac12i(x+x')\cdot\xi}\ud m(\xi)\\
&=
2^{d/2}
\exp\Bigl(-\frac14|x'-x|^2-\frac14|x+x'|^2\Bigr)\\
&=
2^{d/2}\exp\Bigl(-\frac12|x|^2\Bigr)
\exp\Bigl(-\frac12|x'|^2\Bigr)
=
\phi_0(x)\phi_0(x').
\end{align*}
This is the kernel of the rank-one orthogonal projection
$\phi_0\,\bar\otimes\,\phi_0$, and hence
$W(a_0)=\phi_0\,\bar\otimes\,\phi_0$.
\end{proof}

\begin{lemma}\label{lem:Moyal}
For all $a,b\in L^1(\R^{2d}\!,m)$ we have
$$ \wt W(a)\wt W(b) = \wt W(a\,\#\,b),  $$
where the so-called {\em twisted convolution}\index{twisted convolution}\index{convolution!twisted} $a\,\#\, b \in L^1(\R^{2d}\!,m)$ is defined by
$$ {a \,\#\, b}(x,\xi) := \int_{\R^{2d}} e^{\frac12 i(x'\cdot\xi - x\cdot \xi' )} a(x-x'\!, \xi-\xi') b(x'\!,\xi')\ud m(x')\ud m(\xi').  $$
\end{lemma}

Young's inequality implies that $a \,\#\, b$ does indeed belong to $L^1(\R^{2d}\!,m)$.

\begin{proof} Fix $h\in H$.
By \eqref{eq:WW}, a change of variables, and Fubini's theorem,
\begin{align*}
 \wt W(a)\wt W(b)h&  = \int_{\R^{4d}} a(x,\xi)b(x'\!,\xi') \wt W(x,\xi)\wt W(x'\!,\xi')h \ud m(x)\ud m(\xi)\ud m(x')\ud m(\xi')
\\ & = \int_{\R^{4d}} e^{\frac12 i(x'\cdot\xi - x\cdot \xi')} a(x,\xi)b(x'\!,\xi')
\\ & \hskip3cm \times \wt W(x+x'\!,\xi+ \xi')h \ud m(x)\ud m(\xi)\ud m(x')\ud m(\xi')
 \\ &  = \int_{\R^{4d}} e^{\frac12 i(x'\cdot\xi - x\cdot \xi')} a(x-x'\!,\xi-\xi')b(x'\!,\xi')
 \\ & \hskip3cm \times \wt W(x,\xi) h \ud m(x)\ud m(\xi)\ud m(x')\ud m(\xi')
 \\ &  = \int_{\R^{2d}} a\,\#\,b(x,\xi) \wt W(x,\xi)h \ud m(x)\ud m(\xi)
= \wt W(a\,\#\,b)h.
\qedhere
\end{align*}
\end{proof}

This lemma is used to establish the following technical fact.

\begin{lemma}\label{lem:WWW} We have
$$\wt W(a_0)\wt W(x,\xi)\wt W(a_0) =  a_0(x,\xi)\wt W(a_0), \quad x,\xi\in \R^d\!.$$
\end{lemma}
\begin{proof}
Repeating the steps in the proof of Lemma \ref{lem:Moyal}, for all $h\in H$ we obtain
\begin{align*}
 \wt W(x,\xi)\wt W(a_0)h
&= \int_{\R^{2d}} a_0(x'\!,\xi')\wt W(x,\xi)\wt W(x'\!,\xi')h
   \ud m(x')\ud m(\xi')                                      \\
&= \int_{\R^{2d}} e^{\frac12 i(x'\cdot\xi - x\cdot \xi')}
   a_0(x'\!,\xi')\wt W(x+x'\!,\xi+\xi')h
   \ud m(x')\ud m(\xi')                                      \\
&= \int_{\R^{2d}} e^{\frac12 i(x'\cdot\xi - x\cdot \xi')}
   a_0(x'-x,\xi'-\xi)\wt W(x'\!,\xi')h
   \ud m(x')\ud m(\xi')                                      \\
&= \wt W(a_{x,\xi})h,
\end{align*}
where $a_{x,\xi}(x'\!,\xi')
:= e^{\frac12 i(x'\cdot\xi - x\cdot \xi')}
   a_0(x'-x,\xi'-\xi)$.
Hence, by Lemma~\ref{lem:Moyal},
$$
\wt W(a_0)\wt W(x,\xi)\wt W(a_0)
=
\wt W(a_0)\wt W(a_{x,\xi})
=
\wt W(a_0\,\#\,a_{x,\xi}).
$$
We claim that
$$
a_0\,\#\,a_{x,\xi}=a_0(x,\xi)a_0.
$$
Indeed, the preceding calculation, applied to the Schr\"odinger
representation, gives
$$
W(a_{x,\xi})=W(x,\xi)W(a_0).
$$
Consequently,
$$
W(a_0\,\#\,a_{x,\xi})
=
W(a_0)W(x,\xi)W(a_0).
$$
By Lemma~\ref{lem:Schrod-rankone},
$$
W(a_0)W(x,\xi)W(a_0)
=
\iprod{W(x,\xi)\phi_0}{\phi_0}W(a_0).
$$
Moreover, by \eqref{eq:Schrod-repr} and an elementary Gaussian
calculation,
\begin{align*}
\iprod{W(x,\xi)\phi_0}{\phi_0}
&=
2^{d/2}e^{-\frac12ix\cdot\xi}
\int_{\R^d}e^{iy\cdot\xi}
\exp\Bigl(-\frac12|y-x|^2-\frac12|y|^2\Bigr)\ud m(y)\\
&=
\exp\Bigl(-\frac14|x|^2-\frac14|\xi|^2\Bigr)
=
a_0(x,\xi).
\end{align*}
Thus
$$
W(a_0\,\#\,a_{x,\xi})
=
a_0(x,\xi)W(a_0).
$$
Since $a_0,a_{x,\xi}\in L^1(\R^{2d}\!,m)\cap
L^2(\R^{2d}\!,m)$, Young's inequality gives
$a_0\,\#\,a_{x,\xi}\in L^2(\R^{2d}\!,m)$. The injectivity of the
$L^2$-extension of $W$ therefore gives
\[
a_0\,\#\,a_{x,\xi}=a_0(x,\xi)a_0.
\qedhere\]
\end{proof}

We are now ready for the proof of the Stone--von Neumann theorem.

\begin{proof}[Proof of Theorem \ref{thm:Stone-vonNeumann}]
We proceed in three steps.

\smallskip
{\em Step 1} --
We begin by showing that $P:= \wt W(a_0)$ is a nonzero orthogonal projection.
Taking $x=\xi=0$ in Lemma~\ref{lem:WWW} and using the identities
$\wt W(0,0)=I$ and $a_0(0,0)=1$, we obtain
$P^2=P$.
Moreover, \eqref{eq:WW} implies
$\wt W(x,\xi)^\star=\wt W(-x,-\xi)$. Since $a_0$ is real-valued
and even, a change of variables in the defining integral gives
$$
P^\star
=
\int_{\R^{2d}}a_0(x,\xi)\wt W(-x,-\xi)
\ud m(x)\ud m(\xi)
=
P.
$$

Suppose that $P=0$, and fix a unit vector $h\in H$.
Applying \eqref{eq:WW} twice, we obtain, for all $x,\xi\in\R^d$,
\begin{align*}
0 =\iprod{\wt W(x,\xi)P\wt W(-x,-\xi)h}{h}
=\int_{\R^{2d}} e^{i(u\cdot\xi-x\cdot v)} a_0(u,v)\iprod{\wt W(u,v)h}{h} \ud m(u)\ud m(v).
\end{align*}
Thus the Fourier transform of the continuous integrable function
$$
F(u,v):=a_0(u,v)\iprod{\wt W(u,v)h}{h}
$$
vanishes at every point $(-\xi,x)\in\R^{2d}$.
Injectivity of the Fourier transform implies that $F=0$ almost
everywhere, and continuity then gives $F\equiv0$. This contradicts
$F(0,0)=\n h\n^2=1$. Hence $P\ne0$.

\smallskip
{\em Step 2} --
Choose a unit vector $h_0\in\Ran(P)$. By Lemma \ref{lem:WWW},
$$
\iprod{\wt W(x,\xi)h_0}{h_0}
=\iprod{P\wt W(x,\xi)Ph_0}{h_0}
=a_0(x,\xi).
$$
The same identity holds for $W$ and $\phi_0$, since
$W(a_0)=\phi_0\,\bar\otimes\,\phi_0$ by
Lemma \ref{lem:Schrod-rankone}. Combining these identities with
\eqref{eq:WW}, we obtain
\begin{equation}\label{eq:iprodWW}
\begin{aligned}
\iprod{\wt W(x,\xi)h_0}{\wt W(x',\xi')h_0}
&=e^{\frac12 i(x'\cdot\xi-x\cdot\xi')}
  a_0(x-x',\xi-\xi')\\
&=\iprod{W(x,\xi)\phi_0}{W(x',\xi')\phi_0}.
\end{aligned}
\end{equation}
Consequently, the prescription
$$
S\Bigl(\sum_{j=1}^N c_jW(x_j,\xi_j)\phi_0\Bigr)
:=
\sum_{j=1}^N c_j\wt W(x_j,\xi_j)h_0
$$
defines a linear isometry on the span of the vectors
$W(x,\xi)\phi_0$: expanding the squared norms of the two sums
and using \eqref{eq:iprodWW} shows that they are equal, which also
proves that $S$ is well defined.

The closed spans of the orbits of $\phi_0$ and $h_0$ are nonzero
and invariant under $W$ and $\wt W$, respectively, by the product
formula \eqref{eq:Heisenberg-repr}. By the irreducibility of these representations (Proposition \ref{prop:Schrod-irred} for
$W$), these spans are $L^2(\R^d\!,m)$ and $H$.
Thus $S$ extends to a unitary operator from $L^2(\R^d\!,m)$ onto $H$.

To verify the intertwining property, the product formula gives
\begin{align*}
SW(y,\eta,t)W(x,\xi)\phi_0
&=e^{it+\frac12 i(x\cdot\eta-y\cdot\xi)}
  \wt W(x+y,\xi+\eta)h_0\\
&=\wt W(y,\eta,t)SW(x,\xi)\phi_0.
\end{align*}
By linearity and density, it follows that
$$
SW(y,\eta,t)=\wt W(y,\eta,t)S,
\quad (y,\eta,t)\in\H^d.
$$
Since $S$ is unitary, this is the required identity.

\smallskip
{\em Step 3} --
Let $T:L^2(\R^d\!,m)\to H$ be another unitary operator with the
same intertwining property. Then $S^\star T$ commutes with every
$W(x,\xi,t)$, and hence with all translations and boosts.
By Lemma \ref{lem:transl-boost-commutant},
$S^\star T=cI$ for some $c\in\C$. Since $S^\star T$ is unitary,
$|c|=1$, and therefore $T=cS$.
\end{proof}

\begin{remark}
Separability need not actually be assumed in Theorem~\ref{thm:Stone-vonNeumann}. Once a nonzero vector $h\in H$ is chosen, irreducibility says that the closed span of its orbit is $H$; strong continuity and separability of $\R^{2d}$ then imply that this orbit span is separable.
\end{remark}

 The final result of this section describes the Ornstein--Uhlenbeck semigroup in terms of the Weyl calculus. We introduce the unitary operator
\begin{align}\label{eq:U-unit}
U:L^2(\R^d\!,\gamma)\to L^2(\R^d\!,m),
\qquad
U:=D\circ E,
\end{align}
where
$$
(Df)(x):=2^{d/4}f(2^{1/2}x),
\qquad
(Ef)(x):=\exp\bigl(-|x|^2/4\bigr)f(x).
$$
Indeed,
$$
\ud\gamma(x)=\exp\bigl(-|x|^2/2\bigr)\ud m(x),
$$
so $E$ is unitary, and a change of variables shows that $D$ is
unitary.

\begin{theorem}\label{thm:exptL} For all $t > 0$ we have,
with $s:= \frac{1-e^{-t}}{1+e^{-t}}$,
\begin{equation*}
 OU(t) = (1+s)^d\,  U^\star  W(\wh{a_{s}})U,
\end{equation*}
where $W$ denotes the Weyl transform associated with the Schr\"odinger representation and $$ a_s(x,\xi):= \exp(-s(|x|^2+|\xi|^2)), \quad x,\xi\in\R^d\!.$$
\end{theorem}

\begin{proof}
By density, it suffices to prove the identity for $f\in C_{\rm c}(\R^d)$.

Put $b_s:=\wh{a_s}$.
By Lemma \ref{lem:Gauss} (in dimension $2d$),
$$
b_s(u,v)=\frac1{(2s)^d}\exp\Bigl(-\frac{|u|^2+|v|^2}{4s}\Bigr),\quad u,v\in\R^d\!.
$$
In particular, $b_s\in L^1(\R^{2d}\!,m)\cap L^2(\R^{2d}\!,m)$.
Therefore, by \eqref{eq:W(a)f}, the operator $W(b_s)=W(\wh a_s)$ is the
integral operator
$$
W(\wh a_s)g(y)=\int_{\R^d} k_s(x,y)g(x)\ud m(x),
$$
where
$$
k_s(x,y)=\int_{\R^d} b_s(y-x,\xi)\exp\Bigl(\frac{i}{2}(x+y)\cdot\xi\Bigr)\ud m(\xi).
$$
Substituting the above formula for $b_s$ and applying Lemma \ref{lem:Gauss}
once more (now in dimension $d$), we obtain
\begin{align*}
k_s(x,y)
&= \frac1{(2s)^d}\exp\Bigl(-\frac{|x-y|^2}{4s}\Bigr)
   \int_{\R^d}\exp\Bigl(-\frac1{4s}|\xi|^2\Bigr)
          \exp\Bigl(\frac{i}{2}(x+y)\cdot\xi\Bigr)\ud m(\xi)
\\
&= \frac1{(2s)^{d/2}}
   \exp\Bigl(-\frac{|x-y|^2}{4s}\Bigr)
   \exp\Bigl(-\frac{s}{4}|x+y|^2\Bigr).
\end{align*}

Since
$$
Uf(x)=2^{d/4}\exp\Bigl(-\frac12|x|^2\Bigr)f(\sqrt2\,x), \qquad
U^\star g(y)=2^{-d/4}\exp\Bigl(\frac14|y|^2\Bigr)g\Bigl(\frac{y}{\sqrt2}\Bigr),
$$
it follows by the change of variables $x\mapsto x/\sqrt2$ that
\begin{align*}
U^\star W(\wh a_s)U f(y)
   & = 2^{-d/2}\exp\Bigl(\frac14|y|^2\Bigr)
       \int_{\R^d} k_s\Bigl(\frac{x}{\sqrt2},\frac{y}{\sqrt2}\Bigr)
          \exp\Bigl(-\frac14|x|^2\Bigr) f(x)\ud m(x)
\\ & = \int_{\R^d} K_{a_s}(y,x)f(x)\ud m(x),
\end{align*}
where
\begin{align*}
K_{a_s}(y,x)
& = \frac1{2^d s^{d/2}}
   \exp\Bigl(-\frac1{8s}|x-y|^2-\frac{s}{8}|x+y|^2
             -\frac14|x|^2+\frac14|y|^2\Bigr)
\\ & = \frac1{2^d s^{d/2}} \exp\Bigl(-\frac{|(1+s)x-(1-s)y|^2}{8s}\Bigr),
\end{align*}
completing squares in the second identity.
With
$$
r:=\frac{1-s}{1+s}=e^{-t},
\qquad
1-r^2=\frac{4s}{(1+s)^2},
$$
it follows that
\begin{align*}
(1+s)^d U^\star W(\wh a_s)Uf(y)
&=
\frac1{(1-r^2)^{d/2}}
\int_{\R^d}
\exp\Bigl(-\frac{|x-ry|^2}{2(1-r^2)}\Bigr)
f(x)\ud m(x)\\
&=
\int_{\R^d}f\bigl(ry+\sqrt{1-r^2}\,z\bigr)\ud\gamma(z)
=OU(t)f(y).
\end{align*}
Here the second equality follows by substituting
$x=ry+\sqrt{1-r^2}\,z$, and the last equality is the Mehler
formula \eqref{eq:OU-Mehler}.
\end{proof}

\section{Quantisation}\label{sec:first-quantisation}

In this section we interpret the Weyl calculus as a {\em quantisation method},\index{quantisation!first} that is, a prescription for associating operators on a Hilbert space with functions on classical phase space.

\subsection{Commutators and Poisson Brackets}

In classical mechanics, observables are functions on phase space $\R^{2d}$
and dynamics is encoded by Poisson brackets; in quantum mechanics, observables
are selfadjoint operators and the corresponding operation is the commutator.
The main result of this section shows that Weyl calculus recovers the
classical Poisson bracket from commutators in the semiclassical limit.

\smallskip
Throughout this section we work with the Schr\"odinger representation $W$ on
$L^2(\R^d\!,m)$ introduced in Definition \ref{def:Schroedinger-representation}.

Fix a small parameter $h>0$ (in physics, this is the reduced Planck constant)\index{Planck constant} and
define the rescaled Weyl operators by
$$
W_h(x,\xi):=W(h^{1/2}x,h^{1/2}\xi),
\quad x,\xi\in\R^d\!.
$$
For $a\in L^1(\R^{2d}\!,m)$ define the {\em semiclassical Weyl
transform}\index{Weyl!transform}\index{$O$@$\Op_h$}
by
$$
\Op_h(a)
:=
\int_{\R^{2d}}a(x,\xi)W_h(x,\xi)\ud m(x)\ud m(\xi).
$$
For $a\in L^1(\R^{2d}\!,m)\cap L^2(\R^{2d}\!,m)$, put
$$
a_h(z):=h^{-d}a(h^{-1/2}z),
\quad z\in\R^{2d}.
$$
The change of variables $z=h^{1/2}w$ gives
$$
\Op_h(a)=W(a_h).
$$
Consequently, Lemma~\ref{lem:Moyal-injective} gives
$$
h^{d/2}
\|\Op_h(a)\|_{\calL_2(L^2(\R^d\!,m))}
=
h^{d/2}\|a_h\|_{L^2(\R^{2d}\!,m)}
=
\|a\|_{L^2(\R^{2d}\!,m)}.
$$
By density, $\Op_h$ extends uniquely to a linear mapping from
$L^2(\R^{2d}\!,m)$ into the Hilbert--Schmidt operators, and the extension satisfies the same identity.

For $a,b\in L^1(\R^{2d}\!,m)$, the
proof of Lemma~\ref{lem:Moyal} gives
$$
\Op_h(a)\Op_h(b)=\Op_h(a\,\#_h\, b),
$$
where
$$
(a\,\#_h\, b)(x,\xi)
:=
\int_{\R^{2d}}
e^{\frac12 ih(x'\cdot\xi-x\cdot\xi')}
a(x-x',\xi-\xi')b(x',\xi')\ud m(x')\ud m(\xi')
$$
is the {\em Moyal product}\index{Moyal!product} of $a$ and $b$ with parameter $h$.

It will be convenient to regard the symbols $a$ and $b$ as functions on phase space itself rather than on the Fourier side. To this end,
for $f\in L^1(\R^{2d}\!,m)$ define the {\em symplectic Fourier transform}
by
$$
(F_\sigma f)(x,\xi)
:=
\int_{\R^{2d}} f(y,\eta)e^{-i(y\cdot\xi-\eta\cdot x)}\ud m(y)\ud m(\eta).
$$
If, in addition, $F_\sigma f\in L^1(\R^{2d}\!,m)$, we set
$$
\Op_h[f]:=\Op_h(F_\sigma f)
=
\int_{\R^{2d}} (F_\sigma f)(x,\xi)W_h(x,\xi)
\ud m(x)\ud m(\xi).
$$
In particular, this applies to functions
$f\in C_{\rm c}^\infty(\R^{2d})$.
The operator $\Op_h[f]$ will
be called the {\em Weyl quantisation} of $f$.
\index{Weyl!quantisation}\index{quantisation!Weyl}

\begin{theorem}[Poisson brackets as semiclassical limits of commutators]\label{thm:commutator-Poisson}
For all $f,g\in C_{\rm c}^\infty(\R^{2d})$ we have
$$
h^{d/2}
\Big\|\frac{1}{ih}\big[\Op_h[f],\Op_h[g]\big]-\Op_h[\{f,g\}]\Big\|_{\calL_2(L^2(\R^d\!,m))}
\le C h^2\!,\quad h\in(0,1],
$$
where
$$
\{f,g\}:=
\sum_{j=1}^d
\bigl(
   \partial_{y_j}f\,\partial_{\eta_j}g
   -
   \partial_{\eta_j}f\,\partial_{y_j}g
\bigr)
$$
is the Poisson bracket of $f$ and $g$ on $\mathbb R^{2d}$,\index{Poisson!bracket}  and the constant $C$ depends only on $d$ and on the $L^1$-norms of
$\partial^\alpha f$, $|\alpha|\le2d+4$, and the $L^2$-norms of  $\partial^\beta g$, $|\beta|\le3$.
\end{theorem}

\begin{proof}
By the product formula for $\Op_h$ we have
\begin{equation}\label{eq:commutator-reduction}
\frac{1}{ih}[\Op_h[f],\Op_h[g]]
=
\Op_h\Bigl(
\frac{1}{ih}\bigl(
F_\sigma  f\,\#_h\,F_\sigma  g
-
F_\sigma  g\,\#_h\,F_\sigma  f
\bigr)
\Bigr).
\end{equation}
Write $a:=F_\sigma  f$ and $b:=F_\sigma  g$, and abbreviate
$$
\omega(x,\xi;x',\xi'):=x'\cdot\xi-x\cdot\xi'.
$$
Changing variables
$$
(x',\xi')\mapsto(x-x',\xi-\xi')
$$
in the integral defining $b\,\#_h\,a$, and observing that
$$
\omega(x,\xi;x-x',\xi-\xi')
=
-\omega(x,\xi;x',\xi'),
$$
we see that the {\em Moyal commutator}\index{Moyal!commutator} of
$a$ and $b$ equals
\begin{align*}
\frac{1}{ih}(a\,\#_h\, b-b\,\#_h\, a)
&=\frac{1}{ih}\Bigl(\int_{\R^{2d}}\bigl(e^{\frac12 ih\omega}-e^{-\frac12 ih\omega}\bigr)
\,a(x-x',\xi-\xi')\,b(x',\xi')\ud m(x')\ud m(\xi')\Bigr)
\\
&=\frac{2}{h}\int_{\R^{2d}}\sin\bigl(\frac{h}{2}\omega\bigr)\,a(x-x',\xi-\xi')\,b(x',\xi')\ud m(x')\ud m(\xi').
\end{align*}

Define
$\omega_h := \frac{2}{h^3}(\sin(\frac{h}{2}\omega)-\frac{h}{2}\omega)$,
so that
$$
\frac{2}{h}\sin\Bigl(\frac{h}{2}\omega\Bigr)
=
\omega+h^2\omega_h.
$$
By Taylor's theorem,
$$
|\omega_h|\le\frac1{24}|\omega|^3.
$$
This bound implies
\begin{equation}\label{eq:leading-term-plus-remainder}
\frac{1}{ih}(a\,\#_h\, b-b\,\#_h\, a)=c+h^2 r_h,
\end{equation}
where
\begin{align*}
c(x,\xi)
&:=\int_{\R^{2d}} \omega(x,\xi;x',\xi')\,a(x-x',\xi-\xi')\,b(x',\xi')\ud m(x')\ud m(\xi'),
\\
r_h(x,\xi)
&:=\int_{\R^{2d}} \omega_h(x,\xi;x',\xi')\,a(x-x',\xi-\xi')\,b(x',\xi')\ud m(x')\ud m(\xi').
\end{align*}

We claim that $c=F_\sigma \{f,g\}$. To prove this claim,
we start by noting that $F_\sigma$ satisfies the same product--convolution
identity as the ordinary Fourier transform:
$$
   F_\sigma(uv)=(F_\sigma u)*(F_\sigma v),
   \quad u,v\in C_{\rm c}^\infty(\R^{2d}),
$$
where convolution is taken with respect to the normalised Lebesgue measure on
$\R^{2d}$; this follows by writing $F_\sigma$ as the ordinary Fourier
transform composed with the linear change of variables
$(x,\xi)\mapsto(\xi,-x)$. Also, we note that multiplication by the coordinate functions corresponds,
under $F_\sigma $, to differentiation:
\begin{equation}\label{eq:symplFT-derivatives}
F_\sigma (\partial_{\eta_j}f)=-i x_j\,F_\sigma  f,\qquad
F_\sigma (\partial_{y_j}f)=i \xi_j\,F_\sigma  f,
\end{equation}
and similarly for $g$.

Using \eqref{eq:symplFT-derivatives} and this convolution identity, we obtain
\begin{align*}
F_\sigma (\partial_{y_j}f\,\partial_{\eta_j}g)
& = (F_\sigma (\partial_{y_j}f))*(F_\sigma (\partial_{\eta_j}g))
= (i\xi_j a)*(-ix_j b) = (\xi_j a)*(x_j b),
\\
F_\sigma (\partial_{\eta_j}f\,\partial_{y_j}g)
& =
(F_\sigma (\partial_{\eta_j}f))*(F_\sigma (\partial_{y_j}g))
= (-ix_j a)*(i\xi_j b) = (x_j a)*(\xi_j b).
\end{align*}
Therefore
\begin{align*}
F_\sigma \{f,g\}
&=
\sum_{j=1}^d
\bigl((\xi_j a)*(x_j b)-(x_j a)*(\xi_j b)\bigr)
\\
&=
\int_{\mathbb R^{2d}}
\Bigl(
   x'\cdot(\xi-\xi')-(x-x')\cdot\xi'
\Bigr)
a(x-x',\xi-\xi')b(x',\xi')\ud m(x')\ud m(\xi')
\\
&=
\int_{\mathbb R^{2d}}
(x'\cdot\xi-x\cdot\xi')
a(x-x',\xi-\xi')b(x',\xi')\ud m(x')\ud m(\xi')
= c.
\end{align*}
This proves the claim.

We now turn to the $O(h^2)$ bound in the normalised
Hilbert--Schmidt norm.
By \eqref{eq:commutator-reduction} and \eqref{eq:leading-term-plus-remainder},
$$
\frac{1}{ih}[\Op_h[f],\Op_h[g]]-\Op_h[\{f,g\}]
=h^2\Op_h(r_h).
$$
Using
$$
h^{d/2}\|\Op_h(a)\|_{\calL_2(L^2(\R^d\!,m))}=\|a\|_{L^2(\R^{2d}\!,m)},
$$
we obtain
$$
h^{d/2}
\Big\|\frac{1}{ih}[\Op_h[f],\Op_h[g]]-\Op_h[\{f,g\}]\Big\|_{\calL_2(L^2(\R^d\!,m))}
=
h^2\|r_h\|_{L^2(\R^{2d}\!,m)}.
$$

Finally, writing $z=(x,\xi)$ and $z'=(x',\xi')$, letting
$J(x',\xi'):=(-\xi',x')$, and using that $\omega(x,\xi;x',\xi') = \iprod{Jz'}{z-z'}_{\R^{2d}}$,
we have
$$
|\omega(x,\xi;x',\xi')|
\le |z-z'||z'|.
$$
It follows that
$$
|r_h(z)|
\le
\frac1{24}
\int_{\R^{2d}}
|z-z'|^3|a(z-z')||z'|^3|b(z')|\ud m(z'),
$$
and hence, by Young's inequality,
$$
\|r_h\|_2
\le
\frac1{24}
\bigl\||\cdot|^3a\bigr\|_1\,
\bigl\||\cdot|^3b\bigr\|_2.
$$

By the Plancherel theorem and
\eqref{eq:symplFT-derivatives},
$$
\||\cdot|^3F_\sigma g\|_2
\le
C\sum_{|\alpha|\le3}\|\partial^\alpha g\|_2.
$$
To estimate $\|\,|\cdot|^3F_\sigma f\,\|_1$, let $\Delta$ denote
the Laplacian on $\R^{2d}$. By
\eqref{eq:symplFT-derivatives},
$$
(1+|z|^2)^{d+2}F_\sigma f(z)
=
F_\sigma\bigl((1-\Delta)^{d+2}f\bigr)(z).
$$
The $L^1$-$L^\infty$ bound for the symplectic Fourier transform
therefore gives
$$
|F_\sigma f(z)|
\le
\frac{\|(1-\Delta)^{d+2}f\|_1}{(1+|z|^2)^{d+2}}.
$$
Consequently,
\begin{align*}
\bigl\|\,|\cdot|^3F_\sigma f\,\bigr\|_1
&\le \|(1-\Delta)^{d+2}f\|_1
\int_{\R^{2d}}\frac{|z|^3}{(1+|z|^2)^{d+2}}\ud m(z)
\le C_d\sum_{|\alpha|\le2d+4}\|\partial^\alpha f\|_1.
\end{align*}
The integral is finite, since its integrand is bounded near zero
and decays like $|z|^{-2d-1}$ at infinity.

Collecting the estimates, we see that $\|r_h\|_2$ is bounded independently of $h\in(0,1]$, with the
dependence stated in the theorem. This completes the proof.
\end{proof}

\subsection{Position and Momentum as Quantised Operators}

The estimate in Theorem~\ref{thm:commutator-Poisson} comes from the
Taylor expansion
$$
\frac{2}{h}\sin\Bigl(\frac{h}{2}\omega\Bigr)
=
\omega+O(h^2\omega^3),
$$
where
$\omega(x,\xi;x',\xi'):=x'\cdot\xi-x\cdot\xi'$.
Heuristically, the higher-order terms correspond to derivatives of
order at least three. This suggests that, after extending the Weyl
calculus suitably to polynomial symbols, the identity
$$
\frac1{ih}[\Op_h[f],\Op_h[g]]
=
\Op_h[\{f,g\}]
$$
is exact whenever one of the symbols is a polynomial of degree at most
$2$. Polynomial symbols do not belong to
$L^2(\R^{2d}\!,m)$, so a rigorous formulation requires approximation by
sufficiently regular square-integrable symbols.  On a suitable core, the
resulting operators converge strongly to the corresponding symmetrically
ordered polynomial in the position and momentum operators.

For the remainder of this subsection we take $h=1$. With this normalisation in place,
the position and momentum operators, as well as the quantum
harmonic oscillator, arise in this generalised sense from the coordinate
functions and from the classical harmonic oscillator
$\frac12(|x|^2+|\xi|^2)$.

We briefly sketch the argument for the position operator in dimension
$d=1$; the momentum case is analogous, and the harmonic oscillator can be
treated similarly.
In what follows, $W$ is the Schr\"odinger representation on $L^2(\R,m)$.

Recall that the boost operator for $\xi\in\R$ is
$$
(V_\xi f)(x)=e^{ix\xi}f(x), \quad x\in\R,
$$
and that $V_\xi=W(0,\xi)$.
Since $\wh x $ is the multiplication operator $(\wh x f)(x)=xf(x)$, it follows that
$$
e^{i\xi \wh x }=V_\xi=W(0,\xi), \quad \xi\in\R.
$$
For $f\in C_{\rm c}^\infty(\R)$, Theorem
\ref{thm:Fourier-functional-calculus} gives the Fourier integral
representation
$$
f(\wh x )=\frac{1}{(2\pi)^{1/2}}\int_\R \widehat f(\xi)\,e^{i\xi \wh x }\ud\xi
=\frac{1}{(2\pi)^{1/2}}\int_\R \widehat f(\xi)\,W(0,\xi)\ud\xi.
$$
On the other hand, the Weyl transform is given by
$$
W(a)=\int_{\R^2} a(x,\xi)\,W(x,\xi)\ud m(x)\ud m(\xi).
$$
Formally, if one takes the `symbol'
$a(x,\xi):=\widehat f(\xi)\,\delta_0(x)$,
where $\delta_0$ denotes the Dirac delta distribution relative to $m$,
characterised by
$$
\int_\R \varphi(x)\delta_0(x)\ud m(x)=\varphi(0),
$$
then the above identities suggest that $W(a)=f(\wh x )$.
Since $\delta_0$ is not an $L^1$-function, this identity has to be interpreted through approximation.
Replacing $\delta_0$ by an approximate identity in the $x$-variable, one then proves strong convergence of the corresponding Weyl operators to $f(\wh x )$.

To recover the position operator itself, take
$\chi\in C_{\rm c}^\infty(\R)$ such that $0\le\chi\le1$ and
$\chi=1$ near $0$, and put
$f_R(\lambda):=\lambda\chi(\lambda/R)$.
For every $u\in\Dom(\wh x )$, dominated convergence gives
$$
\n f_R(\wh x )u-\wh x u\n_2^2
=
\int_\R |x|^2|\chi(x/R)-1|^2|u(x)|^2\ud m(x)
\to0.
$$
Together with the preceding approximate-identity argument, this
recovers $\wh x $ by successive strong limits on its domain.

\section{Second Quantisation}\label{sec:SQ}

Up to this point we have been concerned with the problem of (first) quantisation, namely, how to define the quantum analogues of classical observables. Many quantum systems involve processes in which particles are
created or annihilated. To describe such systems, we need a Hilbert
space that accommodates different particle numbers. If one uses a Hilbert space $H$ to describe the pure states of a single particle,
one postulates that the $n$-fold Hilbert space tensor product $$H^{\hot n}:= \underbrace{H\otimes\cdots\ot H}_{n\ {\rm times}}$$
describes the pure states of a system of $n$ such particles. Since permutation operators extend to unitaries on $H^{\hot n}$,
the symmetrisation and antisymmetrisation operators of Appendix
\ref{sec:tensor} extend to orthogonal projections on this space.
Their ranges, denoted by $\Gamma^n(H)$ and $\Lambda^n(H)$, are the
symmetric and antisymmetric Hilbert space tensor products and are closed. A {\em boson}\index{boson} is a particle whose $n$-particle states are given by elements of $\Gamma^n(H)$, and a {\em fermion}\index{fermion} is a particle whose $n$-particle states are given by elements of $\Lambda^n(H)$. We restrict attention to the bosonic theory because its connections
with Gaussian analysis bring together several topics developed
earlier in the book.

The elements of the Hilbert space direct sum
$$\Gamma(H):= \bigoplus_{n\in\N} \Gamma^n(H) $$
correspond to superpositions of states carrying different numbers of bosons. The process of passing from $H$ to $\Gamma(H)$ is called {\em (bosonic, or symmetric) second quantisation}.\index{second quantisation}\index{quantisation!second}

Every contraction $T$ on $H$ extends to a contraction $\Gamma(T)$ on $\Gamma(H)$ (see Section \ref{subsec:secq}). For $H=\K^d$, Theorem \ref{thm:OU-2Q} identifies
$\Gamma(e^{-t}I)$ with the Ornstein--Uhlenbeck semigroup:
$$
OU(t)=\Gamma(e^{-t}I).
$$
Over the complex scalars, under this identification the negative of the generator of the Ornstein--Uhlenbeck semigroup is unitarily equivalent to the
number operator studied in Section \ref{subsec:number-phase} for $d=1$. Corollary \ref{cor:FP-Gamma} shows that the Fourier--Plancherel transform is unitarily equivalent to $\Gamma(-iI)$.

We restrict attention to the case in which the Hilbert
space describing the pure states of a single particle is finite-dimensional. {\em Mutatis mutandis}, the theory generalises to arbitrary Hilbert spaces $H$ if one replaces the Gaussian measure by a so-called $H$-isonormal process, a central object in Malliavin calculus.\index{Malliavin calculus} Although this generalisation does not pose any mathematical difficulties, we will not pursue it, as it adds a layer of abstraction that would only obscure the various connections just described.

Unless otherwise stated, the scalar field $\K$ is allowed to be either real or complex.

\subsection{The Wiener--It\^o Chaos Decomposition}\label{subsec:WI}

For $h\in \R^d$
we define $\phi_h\in L^2(\R^d\!,\gamma) = L^2(\R^d\!,\gamma;\K)$, where $\K=\R$ or $\K=\C$,\index{$F$@$\phi_h$} by
$$ \phi_h (x):= \iprod{x}{h} = x\cdot h = \sum_{j=1}^d x_j h_j , \quad x\in \R^d\!.$$
Let $(H_n)_{n\in\N}$ be the sequence of Hermite polynomials introduced in Section \ref{subsec:Hermite}.
For $n\in\N$
we define\index{$Hz$@$\HH_n$} $$\HH_n:=\overline{\hbox{span}}\,\{{H_n}(\phi_h):\, h\in \R^d\!, \, |h|=1\},$$
the closure being taken in $L^2(\R^d\!,\gamma)$. Here, $({H_n}(\phi_h))(x):= {H_n}(\phi_h(x)) = H_n(\iprod{x}{h})$ for $x\in\R^d$\!.
The space $\HH_n$ is sometimes referred to as the
{\em Gaussian chaos} of order $n$.\index{Gaussian!chaos} Note that $\HH_0 = \K\one$ is the one-dimensional space of constant functions.

We are going to prove that the
subspaces $\HH_n$ are pairwise orthogonal and induce a direct sum decomposition
$$\bigoplus_{n\in\N}  \HH_n = L^2(\R^d\!,\gamma).$$
In a second step we will identify orthonormal bases of the summands $\HH_n$.

To start our analysis, for $h\in \R^d$ we define the functions\index{$K$@$K_h$}
\begin{align*}K_h := \exp\Bigl(\phi_h-\frac12 |h|^2\Bigr).
\end{align*}
From $\exp(\phi_h)\in L^2(\R^d\!,\gamma)$ we see that
$K_h$ is well defined as an element of $L^2(\R^d\!,\gamma)$.

From \eqref{eq:Hermite-gen-fc} we see that $K_0=1$ and
\begin{align}\label{eq:Kh-powerseries}
K_h = \exp(|h|\phi_{h/|h|} - \frac12|h|^2) = \sum_{n\in\N}  \frac{|h|^n}{n!}{H_n}( \phi_{h/|h|}), \quad h\in \R^d\!\setminus\{0\}.
\end{align}
In particular,
\begin{align}\label{eq:Kh-powerseries-1}
K_h = \sum_{n\in\N}  \frac1{n!}{H_n}(\phi_h), \quad |h|=1.
\end{align}
The series converges in $L^2(\R^d\!,\gamma)$. Indeed, since
$\phi_{h/|h|}$ is standard Gaussian, Theorem
\ref{thm:Hermite-ONB} gives
$$
\sum_{n\in\N}
\frac{|h|^n}{n!}
\n H_n(\phi_{h/|h|})\n_2
=
\sum_{n\in\N}\frac{|h|^n}{\sqrt{n!}}
<\infty.
$$

\begin{lemma}\label{lem:Kh} The functions $K_h$, $h\in \R^d\!$, span a dense subspace in $L^2(\R^d\!,\gamma)$.
\end{lemma}
\begin{proof}
Suppose that $f\in L^2(\R^d\!,\gamma)$ satisfies
$\iprod{f}{K_h}=0$ for every $h\in\R^d$. Fix $y\in\R^d$ and define
$$
F_y(z):=\int_{\R^d}f(x)e^{z x\cdot y}\ud\gamma(x),
\quad z\in\C.
$$
The Cauchy--Schwarz inequality and Gaussian integrability justify
differentiation under the integral and show that $F_y$ is entire.
For every $t\in\R$,
$$
F_y(t)=e^{t^2|y|^2/2}\iprod{f}{K_{ty}}=0.
$$
The identity theorem therefore gives $F_y\equiv0$. In particular,
$$
\int_{\R^d}f(x)e^{-ix\cdot y}\ud\gamma(x)=0,
\quad y\in\R^d.
$$
This says that the Fourier transform of the integrable function
$x\mapsto f(x)e^{-|x|^2/2}$ vanishes. By injectivity of the Fourier
transform, this function is zero almost everywhere, and hence
$f=0$ in $L^2(\R^d\!,\gamma)$. Thus the functions $K_h$ have dense
linear span.
\end{proof}

\begin{theorem}[Wiener--It\^o decomposition]\label{thm:WI}
We have the orthogonal decomposition
$$L^2(\R^d\!,\gamma)=\bigoplus_{n\in\N}  \HH_n.$$
\end{theorem}
\begin{proof}
Fix $h,h'\in\R^d$ with $|h| = |h'| = 1$ and $s,t\in\R$. Repeating the steps of
\eqref{eq:Hermite-c},
for all $s,t\in\R$ we have
\begin{align*}
\ & \iprod{H(s,\phi_h)}{H(t,\phi_{h'})}_{L^2(\R^d\!,\gamma)}
= \exp(st\,\iprod{h}{h'}).
\end{align*}
Substituting
$ H(r,\cdot) = \sum_{k=0}^\infty \frac{r^k}{k!}H_k$ and $\exp(st\,\iprod{h}{h'}) =  \sum_{k=0}^\infty \frac{(st)^k}{k!}(\iprod{h}{h'})^k$, and taking the partial derivative
$\frac{\partial^{n+m}}{\partial s^m \partial t^n}$
at $s=t=0$
on both sides of the resulting identity,
we obtain
\begin{align*}
 \iprod{H_m(\phi_h)}{H_n(\phi_{h'})}_{L^2(\R^d\!,\gamma)} = \delta_{mn} n! \iprod{h}{h'}^n\!.
\end{align*}
Noting that $\delta_{mn} n! = \delta_{mn} \sqrt{m!}\sqrt{n!}$, we can rewrite this as
\begin{align}\label{eq:ortho-Hermite}
 \Bigl(\frac{H_m(\phi_h)}{\sqrt{m!}} \Big| \frac{H_n(\phi_{h'})} {\sqrt{n!}}\Bigr) = \delta_{mn} \iprod{h}{h'}^n\!.
\end{align}
For $m\not=n$ this implies $\HH_m \perp\HH_n$.

If $f\perp \HH_n$ for all $n\in\N$, then $\iprod{f}{H_n(\phi_h)}_{L^2(\R^d\!,\gamma)} = 0$ for all $h\in \R^d$ with $|h|=1$,
and therefore \eqref{eq:Kh-powerseries}
implies that
$ \iprod{f}{K_h}_{L^2(\R^d\!,\gamma)} =0$ for all $h\in \R^d\!$, and therefore $f=0$ by Lemma \ref{lem:Kh}.
\end{proof}

The next result shows that the Wiener--It\^o decomposition diagonalises the Ornstein--Uhlenbeck semigroup $OU$ on $L^2(\R^d\!,\gamma)$ introduced in Section \ref{subsec:OU}:

\begin{theorem}\label{thm:OU-diag}
The following identities hold:
\begin{enumerate}[label={\rm(\arabic*)}, leftmargin=*]
 \item\label{it:OU-diag1} for all $h\in \R^d$ and $t\ge 0$,
 $$ OU(t)K_h = K_{e^{-t}h};$$
 \item\label{it:OU-diag2} for all $n\in \N$, $F\in \HH_n$, and $t\ge 0,$ $$OU(t)F = e^{-n t}F\!.$$
 \end{enumerate}
\end{theorem}
\begin{proof}
Completing squares in the exponential, for all $h\in\R^d$ and $t\ge 0$ we have
\begin{align*} \! \int_{\R^d}\exp(\sqrt{1-e^{-2t}}\iprod{y}{h})\ud \gamma(y)
& = \prod_{j=1}^d \frac1{\sqrt{2\pi}}\int_\R \exp\Bigl(\sqrt{1-e^{-2t}}y_jh_j - \frac12|y_j|^2\Bigr)\ud y_j
\\ &  = \prod_{j=1}^d \exp\Bigl(\frac12(1-e^{-2t})|h_j|^2\Bigr)
= \exp\Bigl(\frac12(1-e^{-2t})|h|^2\Bigr)\!.
\end{align*}
Hence, from the definitions of $OU(t)$, $\phi_h$, and $K_h$,
\begin{align*}
OU(t)K_h(x) & = \int_{\R^d}  \exp\Bigl(\phi_h(e^{-t}x+ \sqrt{1-e^{-2t}}y)-\frac12 |h|^2\Bigr)\ud \gamma(y)
\\ & = \exp\Bigl(e^{-t}\iprod{x}{h} -\frac12 |h|^2\Bigr)\int_{\R^d}\exp(\sqrt{1-e^{-2t}}\iprod{y}{h})\ud \gamma(y)
\\ & = \exp\Bigl(e^{-t}\iprod{x}{h} -\frac12 |h|^2\Bigr) \exp\Bigl(\frac12(1-e^{-2t})|h|^2\Bigr)
\\ & =  \exp\Bigl(\iprod{x}{e^{-t}h} -\frac12 |e^{-t}h|^2\Bigr) \
= K_{e^{-t}h}(x).
\end{align*}

If $|h|=1$ and $s\ge 0$, it follows from \eqref{eq:Kh-powerseries} and the preceding calculation that
\begin{align*}
OU(t)\sum_{n\in\N}  \frac{s^n}{n!}{H_n}( \phi_{h})= OU(t)K_{sh} = K_{e^{-t}sh} = \sum_{n\in\N}  \frac{s^n e^{-nt}}{n!}{H_n}( \phi_{h}).
\end{align*}
Taking $n$th derivatives in $s$ and evaluating at $s=0$, we obtain the identity
$$OU(t) {H_n}( \phi_{h}) = e^{-nt}H_n( \phi_{h}).$$ By linearity and continuity, we obtain \ref{it:OU-diag2}.
\end{proof}

Part \ref{it:OU-diag2} of the theorem implies that each summand $\HH_n$ is contained in $\Dom(L)$,
where $L$ is the generator of $(OU(t))_{t\ge 0}$, and
$$ LF = -nF, \quad F\in\HH_n, \ n\in\N.$$
Over the complex scalars, each $\HH_n$ is finite-dimensional, since it is contained in the finite-dimensional space of polynomials on $\R^d$ of degree at most $n$. Choose an orthonormal basis of each $\HH_n$. By the Wiener--It\^o decomposition, the union of these bases is an orthonormal basis of $L^2(\R^d\!,\gamma)$ consisting of eigenvectors of $-L$, and the corresponding eigenvalues tend to infinity. Moreover, $-L$ is selfadjoint, so its resolvent set is nonempty. Proposition \ref{prop:sigmaAdiag} therefore applies, and we obtain:

\begin{corollary}\label{cor:siL}
$\sigma(-L) = \N$.
\end{corollary}

\subsection{The Wiener--It\^o Isometry}

Our next aim is to find an orthonormal basis of
each summand $\HH_n$. This will be achieved in Theorem \ref{thm:basisHHn} by means of multivariate Hermite polynomials.

For $\bfn = (n_1,\dots,n_k)\in \N^k$ we write $$ |{\bf n}| :=\sum_{j=1}^k  n_j, \qquad \bfn!:= \prod_{j=1}^k   n_j!.$$
For orthonormal systems ${\bf h}=(h_j)_{j=1}^k$ in $\R^d$ and $\bfn\in\N^k$ we define
$$H_{\bfn}(\phi_{\,\bf h}) :=\prod_{j=1}^k  H_{n_j}(\phi_{h_j}).$$

\begin{theorem}\label{thm:basisHHn}
Let ${\bf h}=(h_j)_{j=1}^d$ be an orthonormal basis of $\R^d$\!.
For each $n\in\N$ the family $$\Bigl\{\frac1{\sqrt{{\bfn}!}} H_{\bfn}(\phi_{\,\bf h}): \, {\bfn}\in {\N^d},\,|{\bfn}|=n
\Bigr\}$$
defines an orthonormal basis of $\HH_n$. As a consequence,
the family
$$\Bigl\{\frac1{\sqrt{{\bfn}!}} H_{\bfn}(\phi_{\,\bf h}): \, {\bfn}\in {\N^d}\Bigr\}$$ defines an orthonormal basis of
$L^2(\R^d\!,\gamma)$.
\end{theorem}
\begin{proof} The proof is divided into two steps.

\smallskip
{\em Step 1} --
We first prove that the family $\{\frac1{\sqrt{{\bfn}!}} H_{\bfn}(\phi_{\,\bf h}): \, {\bfn}\in {\N^d}\}$ is an orthonormal
system in $L^2(\R^d\!,\gamma)$.
Let ${\bfm},{\bfn}\in {\N^d}$. Since $(h_1,\dots,h_d)$ is an orthonormal basis,
the random variables $\phi_{h_1},\dots,\phi_{h_d}$ are independent standard
Gaussian random variables. Therefore, by separation of variables and
\eqref{eq:ortho-Hermite},
\begin{align*}
\Bigl(
\frac1{\sqrt{{\bfm}!}}H_{\bfm}(\phi_{\,\bf h})
\Big|
\frac1{\sqrt{{\bfn}!}}H_{\bfn}(\phi_{\,\bf h})
\Bigr)
& =\Bigl(\prod_{i=1}^d\frac1{\sqrt{m_i!}}H_{ m_i}(\phi_{h_i})\Big| \prod_{j=1}^d\frac1{\sqrt{n_j!}}H_{n_j}(\phi_{h_j})\Bigr)
\\ &=
\prod_{j=1}^d
\Bigl(
\frac1{\sqrt{m_j!}}H_{m_j}(\phi_{h_j})
\Big|
\frac1{\sqrt{n_j!}}H_{n_j}(\phi_{h_j})
\Bigr)
\\
&=
\prod_{j=1}^d \delta_{m_j n_j}
=
\delta_{{\bfm}{\bfn}}.
\end{align*}

{\em Step 2} --
The final step is to prove that
$$
\Bigl\{\frac1{\sqrt{{\bfn}!}}H_{\bfn}(\phi_{\,\bf h}): \, \bfn\in{\N^d},\,|{\bfn}|=n\Bigr\}
$$
is an orthonormal basis of $\HH_n$. Denote by $\GG_n$ the closed linear span of the set $\{H_{\bfn}(\phi_{\,\bf h}):\, \bfn\in{\N^d},\,|{\bfn}|=n\}$.
By Step~1, the spaces $\GG_n$ are mutually orthogonal.

We first show that $H_m(\phi_h)\in\GG_m$ for all $m\in\N$ and $h\in\R^d$ with $\n h\n=1$. Write $h=\sum_{j=1}^d c_jh_j$. Using the generating function for the Hermite polynomials, we have
$$
\sum_{m=0}^\infty \frac{t^m}{m!}H_m(\phi_h)
 =
\exp\Bigl(t\sum_{j=1}^d c_j\phi_{h_j}-\frac12t^2\sum_{j=1}^d c_j^2\Bigr)
 =
\prod_{j=1}^d
\sum_{k_j=0}^\infty \frac{(tc_j)^{k_j}}{k_j!}H_{k_j}(\phi_{h_j}).
$$
Comparing coefficients of $t^m$, we obtain
$$
H_m(\phi_h)
 =
\sum_{|{\bf k}|=m}
\frac{m!}{{\bf k}!}c^{\bf k}H_{\bf k}(\phi_{\,\bf h}).
$$
Thus $H_m(\phi_h)\in\GG_m$, and therefore $\HH_m\subseteq\GG_m$ for all $m\in\N$.

Since also $L^2(\R^d\!,\gamma)=\bigoplus_{m\in\N}\HH_m$, these inclusions are equalities. Indeed, if $g\in\GG_m$ is orthogonal to $\HH_m$, then $g$ is orthogonal to $\HH_j$ for all $j\not=m$, because $\HH_j\subseteq\GG_j$ and $\GG_m\perp\GG_j$. Hence $g$ is orthogonal to all of $L^2(\R^d\!,\gamma)$, and so $g=0$. Thus $\GG_m=\HH_m$ for all $m\in\N$. Together with Step~1 this proves the assertion.
\end{proof}

The orthogonal projection of $L^2(\R^d\!,\gamma)$ onto $\HH_n$ will be denoted by
$J_n$.\index{$J$@$J_n$}

\begin{corollary}\label{cor:secq-Jn}
For all $n\in\N$ and $h\in \R^d$ with $|h|=1$ we have
$$J_n(\phi_h^n) = H_n(\phi_h).$$
More generally, if $(h_j)_{j=1}^k$ is orthonormal in $\R^d$ and ${\bf n} = (n_1,\dots,n_k)\in \N^k$\!,
then
$$J_{|{\bf n}|}(\phi_{h_1}^{n_1}\cdots \phi_{h_k}^{n_k}) =
H_{n_1}(\phi_{h_1})\cdots H_{n_k}(\phi_{h_k}).$$
\end{corollary}
\begin{proof}
We have
\begin{align}\label{eq:Jproj} H_{n_1}(\phi_{h_1})\cdots H_{n_k}(\phi_{h_k})
=\phi_{h_1}^{n_1}\cdots \phi_{h_k}^{n_k} + q(\phi_{h_1},\dots,\phi_{h_k}),
\end{align}
where $q$ is a polynomial of
degree strictly less than $|{\bf n}|$. Hence Theorem \ref{thm:basisHHn} implies that $q(\phi_{h_1},\dots ,\phi_{h_k})\in \bigoplus_{j=0}^{|{\bf n}|-1} \HH_j$,
and the corollary follows
by projecting \eqref{eq:Jproj} onto $\HH_{|{\bf n}|}$. Since $H_{n_1}(\phi_{h_1})\cdots H_{n_k}(\phi_{h_k})\in \HH_{|\bf n|}$, the left-hand side remains unchanged; the right-hand side is mapped to
$J_{|{\bf n}|}(\phi_{h_1}^{n_1}\cdots \phi_{h_k}^{n_k}) $.
\end{proof}

\begin{corollary}\label{cor:KhJn}
For all $h\in \R^d$  with $|h|=1$,
$$K_h= \sum_{n\in\N}
\frac{1}{n!} J_n(\phi_h^n).$$
\end{corollary}
\begin{proof}
This follows from the previous corollary and \eqref{eq:Kh-powerseries-1}.
\end{proof}

\begin{proposition}\label{prop:secq-unitary}
For all $g_1,\dots ,g_n\in \R^d$ and $h_1,\dots ,h_n\in \R^d$
we have
$$\iprod{J_n(\phi_{g_1}\cdots \phi_{g_n})}{J_n(\phi_{h_1}\cdots
\phi_{h_n})}
=\!\sum_{\sigma\in S_n} \iprod{g_1}{h_{\sigma(1)}}\cdots \iprod{g_n}{h_{\sigma(n)}},$$
where $S_n$ is the group of permutations of $\{1,\dots ,n\}$.
\end{proposition}
\begin{proof}
Both sides are symmetric $n$-linear functions of
$g_1,\dots,g_n$ and of $h_1,\dots,h_n$. By polarisation, it is
therefore enough to prove the identity when
$$
g_1=\cdots=g_n=g
\quad\hbox{and}\quad
h_1=\cdots=h_n=h.
$$
By Corollary~\ref{cor:secq-Jn} and
\eqref{eq:ortho-Hermite}, with scaling if $g$ or $h$ is not a unit
vector,
$$
\iprod{J_n(\phi_g^n)}{J_n(\phi_h^n)}
=
n!\iprod{g}{h}^n.
$$
On the other hand,
\[
\sum_{\sigma\in S_n}
\iprod{g}{h}\cdots\iprod{g}{h}
=
n!\iprod{g}{h}^n.
\qedhere\]
\end{proof}

The main result of this section relates the spaces $\HH_n$ to the $n$-fold symmetric tensor product of $\K^d\!$.
The $n$-fold tensor product\index{$Kkk$@$(\K^d)^{\ot n}$}
$$ (\K^d)^{\otimes n} := \underbrace{\K^d \ot\cdots \otimes \K^d }_{n \ \mathrm{ times}}$$
is a Hilbert space with respect to the inner product
\begin{align*}
 \Bigl(\sum_{j=1}^\ell g_1^{(j)}\ot\cdots\ot g_n^{(j)}\Big| \sum_{k=1}^m h_1^{(k)}\ot\cdots\ot h_n^{(k)}
\Bigr)
:= \sum_{j=1}^\ell \sum_{k=1}^m (g_1^{(j)}| h_1^{(k)})\cdots (g_n^{(j)}| h_n^{(k)}).
\end{align*}
We identify $(\K^d)^{\ot 0}$ with the scalar field $\K$.
From Appendix \ref{sec:tensor} we recall that the {\em $n$-fold symmetric tensor product}\index{symmetric!tensor product} of $\K^d $, denoted by $\Gamma^n(\K^d)$,\index{$G$@$\Gamma^n(\K^d)$} is defined
as the range of the orthogonal projection $P_{\Gamma^n} \in \calL((\K^d)^{\ot
n})$ given by
 \begin{align*}
 P_{\Gamma^n}(h_1 \ot \cdots \ot h_n) :=
   \frac{1}{n!} \sum_{\sigma \in S_n} h_{\sigma(1)} \ot \cdots \ot
   h_{\sigma(n)}, \quad h_1,\dots,h_n\in \K^d\!,
 \end{align*}
and extended by linearity, where $S_n$ is the group of permutations of $\{1,\ldots,n\}$.
Stated differently, $\Gamma^n(\K^d)$ is the subspace of those elements of $(\K^d)^{\ot n}$ that are invariant
under the action of $S_n$.

For the formulation of the next theorem we introduce the following notation.
Let ${\bf h} = (h_j)_{j=1}^k$ be a finite sequence in $\K^d$.
For $\bfn\in \N^k$ with $|\bfn|=n$ let
$$
{\bf h}^{\ot \bfn}:=
h_1^{\ot n_1} \ot \cdots \ot h_k^{\ot n_k},
$$
where terms of the form $h_j^{\otimes 0}$ are omitted. If, in addition,
$h_1,\ldots,h_k\in\R^d$, we write
$$
\phi_{\bf h}^{\ot \bfn}:=
\phi_{h_1}^{n_1}\cdots \phi_{h_k}^{n_k}.
$$

\begin{theorem}[Wiener--It\^o isometry]\label{thm:Wiener-Ito-isometry}
There exists a unique isometric isomorphism $$W: \Gamma(\K^d) \to L^2(\R^d\!,\gamma;\K)$$ with the following property:
For every $1\le k\le d$, every orthonormal system ${\bf h} = (h_j)_{j=1}^k$ in $\R^d$\!, every $n\in\N$, and every multi-index $\bfn\in \N^k$ with $|\bfn|=n$,
\begin{align*}
W(P_{\Gamma^n}({\bf h}^{\ot \bfn})) =
\frac{1}{\sqrt{n!}}H_{\bfn}(\phi_{\,\bf h}), \quad \bfn\in\N^k\!.
\end{align*}
\end{theorem}

This mapping $W$ will be referred to as the {\em Wiener--It\^o isometry}.\index{theorem!Wiener--It\^o isometry}

\begin{proof}
For complex scalars, we observe that
$\Gamma^n(\C^d)$ is naturally identified with the complexification of
$\Gamma^n(\R^d)$, and $L^2(\R^d\!,\gamma;\C)$ with the complexification of
$L^2(\R^d\!,\gamma;\R)$. It is therefore enough to prove the theorem in
the real case.

For $h_1,\dots,h_n\in\R^d$, define
$$
        W_n\bigl(P_{\Gamma^n}(h_1\otimes\cdots\otimes h_n)\bigr)
        :=
        \frac1{\sqrt{n!}}\,J_n(\phi_{h_1}\cdots\phi_{h_n}).
$$
By Proposition \ref{prop:secq-unitary}, for all
$g_1,\dots,g_n,h_1,\dots,h_n\in\R^d$ we have
$$
\iprod{
        J_n(\phi_{g_1}\cdots\phi_{g_n})
      }{
        J_n(\phi_{h_1}\cdots\phi_{h_n})
      }
 =
\sum_{\sigma\in S_n}
        \iprod{g_1}{h_{\sigma(1)}}\cdots
        \iprod{g_n}{h_{\sigma(n)}} .
$$
On the other hand, by the definition of the orthogonal projection
$P_{\Gamma^n}$,
$$
\iprod{
        P_{\Gamma^n}(g_1\otimes\cdots\otimes g_n)
      }{
        P_{\Gamma^n}(h_1\otimes\cdots\otimes h_n)
      } =
\frac1{n!}
\sum_{\sigma\in S_n}
        \iprod{g_1}{h_{\sigma(1)}}\cdots
        \iprod{g_n}{h_{\sigma(n)}} .
$$
Thus $W_n$ preserves inner products on the linear span of the symmetrised
elementary tensors. Consequently, $W_n$ is well defined and extends to an
isometry from $\Gamma^n(\R^d)$ into $\HH_n$.

We next check that its range is all of $\HH_n$. Indeed, by Corollary \ref{cor:secq-Jn}, applied
to the standard basis of $\R^d$, the range of $W_n$ contains the functions
$$
        \frac1{\sqrt{n!}}H_{\bfn}(\phi_{\bf e}),
        \quad |\bfn|=n,
$$
and these span $\HH_n$ by Theorem \ref{thm:basisHHn}. This proves the claim.

We conclude that for each $n\in \N$, $$W_n: \Gamma^n(\R^d)\to \HH_n$$ is an
isometric isomorphism.
Taking the orthogonal direct sum over $n\in\N$ and using the
Wiener--It\^o decomposition gives an isometric isomorphism
$$
        W:\Gamma(\R^d)\to L^2(\R^d\!,\gamma;\R).
$$
Finally, if ${\bf h}=(h_j)_{j=1}^k$ is an orthonormal system in $\R^d$
and $\bfn\in\N^k$ with $|\bfn|=n$, then Corollary
\ref{cor:secq-Jn} gives
$$
        W(P_{\Gamma^n}({\bf h}^{\ot\bfn}))
        =
        \frac1{\sqrt{n!}}J_n(\phi_{\bf h}^{\ot\bfn})
        =
        \frac1{\sqrt{n!}}H_{\bfn}(\phi_{\bf h}).
$$
This proves the desired formula. Uniqueness follows from the fact that for each $n\in \N$ the
symmetrised elementary tensors span $\Gamma^n(\R^d)$.
\end{proof}

As a special case, note that for all $h\in \R^d$ with $|h|=1$,
\begin{align}\label{eq:Wiener-Ito-isometry} W (h^{\ot n}) = \frac1{\sqrt{n!}}H_n(\phi_h).\end{align}

\subsection{Second Quantised Operators}\label{subsec:secq}

For $T\in\calL(H)$, the formula
$$
T^{\ot n}(h_1\otimes\cdots\otimes h_n)
:=
Th_1\otimes\cdots\otimes Th_n
$$
and linearity define $T^{\ot n}$ on the algebraic tensor product.
For $n=0$ we understand that $H^{\ot n}=\K$ and $T^{\ot 0}=I_\K$, where $\K$ is the scalar field.

\begin{proposition}\label{prop:tensorTn}
The operator $T^{\ot n}$ extends uniquely to a bounded operator on the
Hilbert space tensor product $H^{\ot n}$, and
$$
\n T^{\ot n}\n=\n T\n^n.
$$
\end{proposition}

\begin{proof} For $n=0$ the assertion is immediate, so we may assume that $n\ge1$.
For $j=1,\ldots,n$, on the algebraic tensor product we define
$$
T_j:=I\otimes\cdots\otimes I\otimes T
          \otimes I\otimes\cdots\otimes I,
$$
where $T$ acts on the $j$-th tensor factor.
As in Section \ref{subsec:partialtrace} one shows that each $T_j$ extends uniquely to a bounded operator on $H^{\ot n}$ with norm $\n T_j\n=\n T\n$.

Since on the algebraic tensor product we have
$T^{\ot n}=T_1T_2\cdots T_n$, the bounded operator on the right-hand side provides the unique bounded extension of $T^{\ot n}$, and
$$
\n T^{\ot n}\n
\le\prod_{j=1}^n\n T_j\n
\le\n T\n^n.
$$

For the reverse inequality, let $h\in H$ satisfy $\n h\n\le1$.
Since $\n h^{\ot n}\n=\n h\n^n\le1$, we have
$$
\n T^{\ot n}\n
\ge
\n T^{\ot n}(h^{\ot n})\n
=
\n Th\n^n.
$$
Taking the supremum over $\n h\n\le1$ gives
$\n T^{\ot n}\n\ge\n T\n^n$.
\end{proof}

If $T\in\calL(H)$ is a contraction, the restriction of $T^{\ot n}$ to
the symmetric part $\Gamma^n(H)$ defines a contraction $\Gamma^n(T)$
on $\Gamma^n(H)$. By taking direct sums,\index{$G$@$\Gamma(T)$}
$$ \Gamma(T) := \bigoplus_{n\in\N} \Gamma^n(T)$$
defines a contraction on $\bigoplus_{n\in\N} \Gamma^n(H)$.

\begin{definition}[Symmetric second quantisation]
The Hilbert space $\bigoplus_{n\in\N}\Gamma^n(H) =: \Gamma(H)$ is called the {\em symmetric Fock space}
\index{symmetric!Fock space}\index{Fock space} over $H$. When $T$ is a contraction on $H$,
the contraction $\Gamma(T)$ on $\Gamma(H)$ is called the {\em symmetric second quantisation}\index{symmetric!second quantisation} of $T$.
\end{definition}

{\em Antisymmetric second quantisation} can be defined similarly but will not be studied here. For this reason, we will omit the adjective `symmetric' from now on and simply talk about {\em second quantisation}.

If $S$ and $T$ are contractions on $H$, their second quantisations satisfy
\begin{align}\label{eq:secq-1} \Gamma(I) = I, \qquad \Gamma(ST) = \Gamma(S)\Gamma(T), \qquad \Gamma(T^\star) = (\Gamma(T))^\star\!.
\end{align}
In what follows we again take $H=\K^d$\!. If $T$ is a contraction on $\K^d\!$, the operator $\Gamma(T)$ induces, via the Wiener--It\^o isometry (Theorem \ref{thm:Wiener-Ito-isometry}), a contraction on $L^2(\R^d\!,\gamma) = L^2(\R^d\!,\gamma;\K)$. By a slight abuse of notation, it will be denoted by $\Gamma(T)$ as well.
It is easily checked that \eqref{eq:secq-1} holds again.

For $h\in \K^d$ define the {\em exponential vector}\index{exponential vector}
$$
\epsilon(h):=\sum_{n\in\N}\frac{1}{\sqrt{n!}}\,h^{\ot n}\in \Gamma(\K^d).
$$
This series converges in $\Gamma(\K^d)$, and
$$
\|\epsilon(h)\|^2=\sum_{n\in\N}\frac{|h|^{2n}}{n!}=e^{|h|^2}.
$$
We now extend the notation $K_h$ to all $h\in \K^d$ by setting
$$
K_h:=W(\epsilon(h)),\quad h\in \K^d.
$$
If $h\in \R^d$, this agrees with the earlier definition
$$
K_h(x)=\exp\Bigl(\phi_h(x)-\frac12|h|^2\Bigr),
$$
by \eqref{eq:Kh-powerseries} and \eqref{eq:Wiener-Ito-isometry}.

\begin{lemma}\label{lem:Mall-Kh}
If $T$ is a contraction on $\K^d$, then for all $h\in \K^d$,
$$
\Gamma(T)K_h=K_{Th}.
$$
\end{lemma}

\begin{proof}
For $h\in \K^d$ we have
$$
\Gamma(T)\epsilon(h)
 = \sum_{n\in\N}\frac{1}{\sqrt{n!}}\,\Gamma^n(T)h^{\ot n}
 = \sum_{n\in\N}\frac{1}{\sqrt{n!}}\,(Th)^{\ot n}
 = \epsilon(Th).
$$
Applying the Wiener--It\^o isometry $W$, we obtain
\[
\Gamma(T)K_h
 = W(\Gamma(T)\epsilon(h))
 = W(\epsilon(Th))
 = K_{Th}.
\qedhere\]
\end{proof}

As a special case, for the Ornstein--Uhlenbeck semigroup we obtain:

\begin{theorem}[Ornstein--Uhlenbeck semigroup and second quantisation]\label{thm:OU-2Q} Under the Wiener--It\^o isometry, for all $t\ge 0$ we have $$ OU(t) = \Gamma(e^{-t} I).$$
\end{theorem}

\begin{proof} This follows from Theorem \ref{thm:OU-diag} and Lemma \ref{lem:Mall-Kh}, which give
$$OU(t)K_h = K_{e^{-t}h}
= \Gamma(e^{-t}I)K_h,$$
and the density of the span of the functions $K_h$ in $L^2(\R^d\!,\gamma)$ shown in Lemma  \ref{lem:Kh}.
\end{proof}

For a real contraction $T$ on $\R^d$, its second quantisation extends
complex-linearly to $L^2(\R^d\!,\gamma;\C)$, and we have the following
positivity result on the real part.

\begin{theorem}[Positivity]\label{thm:sec-quant-pos}
If $T$ is a contraction on $\R^d\!$, then $\Gamma(T)$
is a positivity-preserving contraction on $L^2(\R^d\!,\gamma)$.
\end{theorem}

For a genuinely complex-linear contraction $T$ on $\C^d$, positivity need
not be preserved: even if $h\in\R^d$, the vector $Th$ need not be real, so
the corresponding first-chaos function is no longer real-valued.

\begin{proof}
Put $C:=(I-T^\star T)^{1/2}$. For $f\in L^2(\R^d\!,\gamma)$ define
$$
P_Tf(x):=\int_{\R^d} f(T^\star x+Cy)\ud\gamma(y).
$$
This formula first makes sense for bounded Borel functions $f$, and then extends to all $f\in L^2(\R^d\!,\gamma)$. Indeed, if $X$ and $Y$ are independent standard Gaussian random variables in $\R^d$, then $T^\star X+CY$ is again standard Gaussian, since its covariance equals $T^\star T+C^2=I$. Hence, by Jensen's inequality,
$$
\int_{\R^d}|P_Tf|^2\ud\gamma
 \leq \int_{\R^d}\int_{\R^d}|f(T^\star x+Cy)|^2\ud\gamma(y)\ud\gamma(x)
 = \int_{\R^d}|f|^2\ud\gamma.
$$
Thus $P_T$ is a contraction on $L^2(\R^d\!,\gamma)$. It is clear from its definition that $P_T$ is positivity preserving.

It remains to identify $P_T$ with $\Gamma(T)$. Let $h\in\R^d$. Using the Gaussian identity $$\int_{\R^d}\exp(\phi_g(y))\ud\gamma(y)=\exp(\frac12|g|^2),$$ we obtain
\begin{align*}
P_TK_h(x)
&= \int_{\R^d}
   \exp\Bigl(\iprod{T^\star x+Cy}{h}-\frac12|h|^2\Bigr)
   \ud\gamma(y)                                      \\
&= \exp\Bigl(\iprod{x}{Th}-\frac12|h|^2\Bigr)
   \exp\Bigl(\frac12|Ch|^2\Bigr)                     \\
&= \exp\Bigl(\iprod{x}{Th}-\frac12|Th|^2\Bigr)
 = K_{Th}(x).
\end{align*}
By Lemma \ref{lem:Mall-Kh} we also have $\Gamma(T)K_h=K_{Th}$. Hence $P_T$ and $\Gamma(T)$ agree on the linear span of the functions $K_h$, $h\in\R^d$. This span is dense in $L^2(\R^d\!,\gamma)$ by Lemma \ref{lem:Kh}, and both $P_T$ and $\Gamma(T)$ are bounded. Therefore $P_T=\Gamma(T)$. Since $P_T$ is positivity preserving, so is $\Gamma(T)$.
\end{proof}

\subsection{The Segal--Plancherel Transform}\label{subsec:Segal-Plancherel}

In this section we discuss the Segal--Plancherel transform
$\mathscr W$ on $L^2(\R^d\!,\gamma)$, a Gaussian analogue of the
Fourier--Plancherel transform. With $\ud m(x):=(2\pi)^{-d/2}\ud x$,
the Fourier transform
$$
 \calF f(\xi)=\int_{\R^d}e^{-ix\cdot\xi}f(x)\ud m(x)
$$
extends from
$L^1(\R^d\!,m)\cap L^2(\R^d\!,m)$ to an isometry on $L^2(\R^d\!,m)$.

As in Section \ref{subsec:StonevonNeumann} we let $U:= D\circ E$,
where  $D: L^2(\R^d\!,m)\to L^2(\R^d\!,m)$ and $E:L^2(\R^d\!,\gamma)\to L^2(\R^d\!,m)$
are the unitary operators
\begin{align*} Df(x) := 2^{d/4} f\bigl({\sqrt 2}x\bigr), \qquad Ef(x) := e(x)f(x),
\end{align*}
with $e(x) :=\exp(-\frac14|x|^2).$

 \begin{theorem}\label{thm:Segal}
 The mapping $\mathscr{W}: f\mapsto \mathscr{W}f$, defined for multivariate polynomials
$f:\R^d\to\C$ by analytic continuation as
$$
\mathscr Wf(x):= \int_{\R^d} f(-ix+\sqrt 2 y)\ud \gamma(y), \quad x \in \R^{d}\!,
$$
extends to a unitary operator on $L^2(\R^d\!,\gamma)$, and we have
\begin{align}\label{eq:FWT}
\mathscr{W}  = U^\star  \circ \calF \circ U.
\end{align}
\end{theorem}
\begin{proof}
Fix $h\in\R^d$ with $|h|=1$. Since $\phi_h$ is standard Gaussian,
for $s\in\R$ and $x\in\R^d$ we have
\begin{align*}
\int_{\R^d}
\exp\Bigl(s\bigl(-i\phi_h(x)+\sqrt2\,\phi_h(y)\bigr)
-\frac12s^2\Bigr)\ud\gamma(y)
=\exp\Bigl(-is\phi_h(x)+\frac12s^2\Bigr).
\end{align*}
Differentiating $n$ times at $s=0$ and using the generating
function \eqref{eq:Hermite-gen-fc} gives
$$
\mathscr W H_n(\phi_h)=(-i)^nH_n(\phi_h).
$$
Differentiation under the integral is justified by Gaussian
integrability. The Wiener--It\^o decomposition and the Hermite
bases of the spaces $\HH_n$ therefore show that $\mathscr W$
extends to a unitary operator on $L^2(\R^d\!,\gamma)$, acting
as multiplication by $(-i)^n$ on $\HH_n$.

To identify this extension, fix $h\in\R^d$. Completing the square
in the definition of $UK_h$ gives
$$
(UK_h)(x)
=
2^{d/4}\exp\Bigl(\frac12|h|^2\Bigr)
\exp\Bigl(-\frac12|x-\sqrt2\,h|^2\Bigr).
$$
Therefore, by Lemma \ref{lem:Gauss},
$$
(\calF UK_h)(\xi)
=
2^{d/4}\exp\Bigl(\frac12|h|^2\Bigr)
e^{-i\sqrt2\,\xi\cdot h}
\exp\Bigl(-\frac12|\xi|^2\Bigr).
$$
Applying $U^\star$, we obtain
$$
(U^\star\calF UK_h)(x)
=
\exp\Bigl(-i\phi_h(x)+\frac12|h|^2\Bigr).
$$
On the other hand, the Hermite expansion
\eqref{eq:Kh-powerseries} and the action of $\mathscr W$ on each
chaos give precisely the same expression:
$$
\mathscr WK_h
=
\exp\Bigl(-i\phi_h+\frac12|h|^2\Bigr).
$$
Since the functions $K_h$, $h\in\R^d$, have dense linear span by
Lemma \ref{lem:Kh}, the two bounded operators agree on all of
$L^2(\R^d\!,\gamma)$. This proves \eqref{eq:FWT}.
\end{proof}

We have the following representation of $\mathscr{W}$ in terms of second quantisation:

\begin{theorem}[Segal]\label{thm:WT}\index{theorem!Segal} $\mathscr{W} = \Gamma(-i I).$
\end{theorem}
\begin{proof}
By definition, $\Gamma^n(-iI)=(-i)^nI$ on $\Gamma^n(\C^d)$.
Under the Wiener--It\^o isometry, this means that $\Gamma(-iI)$
acts as multiplication by $(-i)^n$ on $\HH_n$. The preceding
proof shows that $\mathscr W$ has the same action on every
$\HH_n$. The Wiener--It\^o decomposition therefore gives
$\mathscr W=\Gamma(-iI)$.
\end{proof}

Theorems \ref{thm:Segal} and \ref{thm:WT} combine to give the following result. Recall the definition of the unitary operator $U$ of \eqref{eq:U-unit}.

\begin{corollary}\label{cor:FP-Gamma}
The Fourier--Plancherel transform is unitarily equivalent to $\Gamma(-i I)$. More precisely, we have
$$ U^\star \circ \calF\circ U = \Gamma(-iI).$$
\end{corollary}

This gives a neat `explanation' of the identity $\calF^4=I$: by
multiplicativity of second quantisation, it follows from $(-i)^4=1$.

\subsection{Creation and Annihilation}\label{subsec:annih-creat}

For $h\in\R^d$ and $n\in\N$, the {\em creation operator}
\index{creation operator}
$a_n^\dagger(h):\Gamma^n(\R^d)\to\Gamma^{n+1}(\R^d)$
is defined by
$$
a_n^\dagger(h)x
:=
\sqrt{n+1}\,P_{\Gamma^{n+1}}(h\otimes x),
\quad x\in\Gamma^n(\R^d).
$$
The {\em annihilation operator}\index{annihilation operator}
$a_{n+1}(h):\Gamma^{n+1}(\R^d)\to\Gamma^n(\R^d)$
is defined by
$$
a_{n+1}(h):=(a_n^\dagger(h))^\star,
$$
and we put $a_0(h):=0$.
Here, as before,
 \begin{align*}
 P_{\Gamma^{n+1}}(h_1 \ot \cdots \ot h_{n+1}) :=
   \frac{1}{(n+1)!} \sum_{\sigma \in S_{n+1}} h_{\sigma(1)} \ot \cdots \ot
   h_{\sigma(n+1)}, \quad h_1,\dots,h_{n+1}\in \R^d\!.
 \end{align*}
On symmetrised elementary tensors these definitions give
\begin{align*}
a_n^\dagger(h)
P_{\Gamma^n}(h_1\otimes\cdots\otimes h_n)
= \sqrt{n+1}\, P_{\Gamma^{n+1}} (h\otimes h_1\otimes\cdots\otimes h_n)
\end{align*}
and
\begin{align*}
 a_{n+1}(h)
P_{\Gamma^{n+1}}(h_1\otimes\cdots\otimes h_{n+1})
= \frac1{\sqrt{n+1}} \sum_{j=1}^{n+1} \iprod{h_j}{h}\,
P_{\Gamma^n} (h_1\otimes\cdots\otimes\widehat{h_j}\otimes\cdots \otimes h_{n+1}),
\end{align*}
where a hat over a term indicates that the term is omitted.

These operators satisfy
$$
\n a_n^\dagger(h)\n
=
\n a_{n+1}(h)\n
=
\sqrt{n+1}\,|h|,
$$
and the canonical commutation relations
$$
a_{n+1}(h)a_n^\dagger(g)
-
a_{n-1}^\dagger(g)a_n(h)
=
\iprod{g}{h}I_{\Gamma^n(\R^d)},
\quad n\ge1.
$$
For $n=0$ the corresponding identity is
$ a_1(h)a_0^\dagger(g) = \iprod{g}{h}I_{\Gamma^0(\R^d)}$.

Let $(e_j)_{j=1}^d$ denote the standard basis of $\R^d$.
Using the Wiener--It\^o isometry, we define $A_0$ on vectors of the form $Wx$, with $x\in\Gamma^n(\R^d)$ and $n\in\N$, by
$$
A_0(Wx):=\bigl(W(a_n(e_1)x),\dots,W(a_n(e_d)x)\bigr)\in L^2(\R^d\!,\gamma;\K^d).
$$
Likewise, on vectors $(Wx_1,\dots,Wx_d)$ with $x_1,\dots,x_d\in\Gamma^n(\R^d)$ and $n\in\N$, we define
$$
A_0^\dagger\bigl(Wx_1,\dots,Wx_d\bigr)
 := W\bigl(a_n^\dagger(e_1)x_1+\cdots+a_n^\dagger(e_d)x_d\bigr)\in  L^2(\R^d\!,\gamma).
$$
Thus
$$
A_0(\HH_{n+1})\subseteq \HH_n^d,
\qquad
A_0^\dagger(\HH_n^d)\subseteq \HH_{n+1}.
$$

Put
$\mathscr P := {\rm Pol}(\R^d)$ for brevity. Since this space is the span of the spaces $\HH_n$, $n\in\N$, the preceding definitions
extend linearly to densely defined operators
$$
A_0:\mathscr P\to L^2(\R^d\!,\gamma;\K^d), \qquad
A_0^\dagger:\mathscr P^d\to L^2(\R^d\!,\gamma).
$$
The orthogonality of the Wiener chaoses implies
$$
\iprod{A_0f}{g}
=
\iprod{f}{A_0^\dagger g},
\quad
f\in\mathscr P,\ g\in\mathscr P^d.
$$
Thus
$$
A_0^\dagger\subseteq A_0^\star
\quad\hbox{and}\quad
A_0\subseteq(A_0^\dagger)^\star.
$$
Since $\mathscr P^d$ is dense in $L^2(\R^d\!,\gamma; \K^d)$, the domain
$\Dom(A_0^\star)$ is dense, and therefore $A_0$ is closable.
Similarly, $\Dom((A_0^\dagger)^\star)$ is dense, and hence
$A_0^\dagger$ is closable.
We denote their closures by $A$ and $A^\dagger$, respectively, and reserve the notation $\Dom(A)$ and $\Dom(A^\dagger)$ for the domains of their closures.

Let $\nabla = (\partial_1,\dots,\partial_d)$ be the gradient, viewed as a densely defined closed operator from
$L^2(\R^d\!,\gamma)$ to $L^2(\R^d\!,\gamma;\K^d)$ with its natural domain $\Dom(\nabla) = H^{1}(\R^d\!,\gamma)$,
the Hilbert space of all functions in $L^2(\R^d\!,\gamma)$ whose weak
first-order derivatives belong to $L^2(\R^d\!,\gamma)$.

\begin{lemma}\label{lem:H1gamma-dense}
The space ${\rm Pol}(\R^d)$ of polynomials in the real variables $x_1,\dots,x_d$ is dense in $H^{1}(\R^d\!,\gamma)$.
\end{lemma}
\begin{proof}
For $\bfn=(n_1,\dots,n_d)\in\N^d$ put $p_{\bfn}(x):=H_{n_1}(x_1)\cdots H_{n_d}(x_d)$.
The Hermite recurrence relations of
Proposition~\ref{prop:Hermite-pols} give
$$
-\Delta p_{\bfn}+x\cdot\nabla p_{\bfn}
=|\bfn|p_{\bfn}.
$$

Suppose that $f\in H^1(\R^d\!,\gamma)$ is orthogonal to all
polynomials with respect to the usual inner product of this
space. Integration by parts with respect to the Gaussian
measure gives
\begin{align*}
0
=\int_{\R^d}f\ov{p_{\bfn}}\ud\gamma
  +\int_{\R^d}\!\nabla f\cdot\ov{\nabla p_{\bfn}}\ud\gamma
=\int_{\R^d}
  f\,\ov{p_{\bfn}-\Delta p_{\bfn}
                 +x\cdot\nabla p_{\bfn}}\ud\gamma
=(1+|\bfn|)\int_{\R^d}f\ov{p_{\bfn}}\ud\gamma.
\end{align*}
To justify the integration by parts for weak derivatives,
insert a smooth cutoff $\chi_R(x):=\chi(x/R)$, where
$\chi\in C_{\rm c}^\infty(\R^d)$ equals $1$ near the origin,
and let $R\to\infty$. The additional term satisfies
$$
\left|
\int_{\R^d}f\,\nabla\chi_R\cdot
                  \ov{\nabla p_{\bfn}}\ud\gamma
\right|
\lesssim R^{-1}\n f\n_2\n\nabla p_{\bfn}\n_2\to0,
$$
and the remaining terms converge by dominated convergence,
since every polynomial belongs to $L^2(\R^d\!,\gamma)$.

Thus $f$ is orthogonal to every $p_{\bfn}$ in
$L^2(\R^d\!,\gamma)$. By Theorem~\ref{thm:basisHHn},
these polynomials span a dense subspace of
$L^2(\R^d\!,\gamma)$, so $f=0$.
Hence the orthogonal complement of ${\rm Pol}(\R^d)$ in
$H^1(\R^d\!,\gamma)$ is trivial, which proves the assertion.
\end{proof}

\begin{proposition}\label{prop:grad} We have $A = \nabla $ with equality of domains.
\end{proposition}
\begin{proof}
First let $f = \frac1{\sqrt{n!}}H_n(\phi_h)$ with $n\in \N$ and $h\in \R^d$ with $|h|=1$.
For $n=0$ the identity $Af = \nabla f$ is trivial, so in what follows we take $n\ge 1$.
The function $f$ is the image under the Wiener--It\^o isometry of the element $h^{\ot n}\in \Gamma^n(\R^d)$ and we have
\begin{align*}
(Af )_j = (AW(h^{\ot n}))_j & = W\bigl(a_n(e_j)\underbrace{h\ot\cdots\ot h}_{n \ {\rm times}}\bigr)
 =  \frac1{\sqrt{n}} \sum_{\ell = 1}^n \iprod{h}{e_j}W\bigl(\underbrace{h\ot\cdots\ot h}_{n-1 \ {\rm times}}\bigr)
\\ & =  \sqrt{n}\iprod{h}{e_j}W\bigl(\underbrace{h\ot\cdots\ot h}_{n-1 \ {\rm times}}\bigr)
= \frac{\sqrt{n}}{\sqrt{(n-1)!}}\iprod{h}{e_j}H_{n-1}(\phi_h).
\end{align*}
On the other hand, since $H_n' = n H_{n-1}$ and $\partial_j \phi_h = \iprod{e_j}{h}$,
\begin{align*}
(\nabla f)_j  =  \frac1{\sqrt{n!}}\partial_j H_n(\phi_h) =\frac{\sqrt{n}}{\sqrt{(n-1)!}}\iprod{e_j}{h} H_{n-1}(\phi_h).
\end{align*}
By Theorem \ref{thm:basisHHn}, the functions $H_n(\phi_h)$,
$n\in\N$, $h\in\R^d$ with $|h|=1$, span ${\rm Pol}(\R^d)$.
The preceding computation therefore shows that
$$
A_0=\nabla|_{{\rm Pol}(\R^d)}.
$$
By Lemma \ref{lem:H1gamma-dense}, polynomials are dense in
$\Dom(\nabla)=H^1(\R^d\!,\gamma)$ with respect to the graph norm
of $\nabla$. Since $\nabla$ is closed and $A$ is the closure of
$A_0$, it follows that
$$
A=\ov{A_0}
=\ov{\nabla|_{{\rm Pol}(\R^d)}}
=\nabla
$$
with equality of domains.
\end{proof}

Working over the complex scalars,
for $j=1,\dots,d$ and $f\in {\rm Pol}(\R^d)$ we let
$$a_j f:= (Af)_j, \qquad a_j^\dagger f:= A^\dagger(0,\dots,0,f,0,\dots)$$ with $f$ in the $j$th place.
These operators are dual to each other, in the sense that with respect to the inner product of $L^2(\R^d\!,\gamma)$ we have
$$ \iprod{a_jf}{g} = \iprod{f}{a_j^\dagger g}, \quad f,g \in {\rm Pol}(\R^d).$$

\subsection{Position and Momentum}
Recall the selfadjoint position and momentum operators $\wh x_j$ and
$\wh\xi_j=\frac1i\partial_j$ on $L^2(\R^d\!,m)$ introduced in
Section~\ref{subsec:pos-mom}. Using the unitary operator
$$
U:L^2(\R^d\!,\gamma)\to L^2(\R^d\!,m)
$$
of Sections~\ref{subsec:StonevonNeumann} and
\ref{subsec:Segal-Plancherel}, we define the {\em position
operator}\index{position operator} $Q=(q_1,\dots,q_d)$ and the
{\em momentum operator}\index{momentum operator} $P=(p_1,\dots,p_d)$
on $L^2(\R^d\!,\gamma)$ by
$$
q_j:=U^\star\wh x_jU,
\qquad
p_j:=U^\star\wh\xi_jU.
$$
Thus $q_j$ and $p_j$ are selfadjoint, with domains
$\Dom(q_j)=U^\star\Dom(\wh x_j)$ and $\Dom(p_j)=U^\star\Dom(\wh\xi_j)$.

Let us identify their restrictions to ${\rm Pol}(\R^d)$. By the
definition $U=D\circ E$ we have
$$
Uf(x)
=
2^{d/4}\exp\bigl({-\frac12|x|^2}\bigr)f(\sqrt2x).
$$
A direct computation therefore gives, for $f\in{\rm Pol}(\R^d)$,
$$
U^\star\wh x_jUf
=
\frac1{\sqrt2}x_jf, \qquad
U^\star\wh\xi_jUf
=
\frac1{i\sqrt2}(2\partial_jf-x_jf).
$$
By Proposition~\ref{prop:grad}, $a_j=\partial_j$ on
${\rm Pol}(\R^d)$. Moreover, integration by parts with respect to the
Gaussian measure gives $a_j^\dagger f=x_jf-\partial_jf$.
Consequently, for $f\in{\rm Pol}(\R^d)$,
$$
q_jf
=
\frac1{\sqrt2}(a_j+a_j^\dagger)f, \qquad
p_jf
=
\frac1{i\sqrt2}(a_j-a_j^\dagger)f.
$$
The normalising constant $1/\sqrt2$ in these identities is chosen so
that, under the unitary operator $U$, the operators $q_j$ and $p_j$
correspond to the standard position and momentum operators
$\wh x_j$ and $\wh\xi_j$, respectively.

The identities below are understood in the sense that they hold when
applied to functions in ${\rm Pol}(\R^d)$. Some additional details are
addressed in Problems~\ref{prob:crea-anni2} and \ref{prob:crea-anni3}.
From the commutation relation $[a_j,a_j^\dagger] = I$ we have
\begin{equation*}
\begin{aligned}
 [p_j,q_j] & =
 \frac1{i}(a_ja^\dagger_j - a^\dagger_ja_j ) = \frac1{i}I
\end{aligned}
\end{equation*}
as well as the identity
\begin{equation}\label{reason:2}
\begin{aligned}
\frac12(p_j^2 + q_j^2)
 & = \frac12(a_j^\dagger a_j + a_ja_j^\dagger)
= a_j^\dagger a_j + \frac12[a_j,a_j^\dagger]
= a_j^\dagger a_j + \frac12 I.
\end{aligned}
\end{equation}

\subsection{The Quantum Harmonic Oscillator}\label{subsec:QHA}

The Hamiltonian of the classical harmonic oscillator with unit mass and unit angular frequency equals
$$
\frac12\bigl(|\xi|^2+|x|^2\bigr), \quad (x,\xi)\in \R^d\times \R^d = \R^{2d}.
$$
In quantum mechanics, the position and momentum variables are replaced
by the corresponding operators $Q=(q_1,\dots,q_d)$ and
$P=(p_1,\dots,p_d)$. This leads to the operator
$$
\frac12(P^2+Q^2)
=
\sum_{j=1}^d\frac12(p_j^2+q_j^2),
$$
which is called the {\em quantum harmonic oscillator}.
\index{quantum harmonic oscillator}

An easy computation in terms of the annihilation and creation operators
gives
\begin{align*}
\nonumber
-L
&=
\nabla^\star\nabla
=
A^\dagger A
=
\sum_{j=1}^d a_j^\dagger a_j,
\end{align*} and hence, by \eqref{reason:2},
\begin{align}
\label{eq:QHO2}
-L & =
-\frac d2+\frac12(P^2+Q^2)=
-\frac d2+\sum_{j=1}^d\frac12(p_j^2+q_j^2).
\end{align}
These identities are understood on ${\rm Pol}(\R^d)$.
Thus the negative Ornstein--Uhlenbeck generator $-L$ is
obtained from the quantum harmonic oscillator by subtracting
its ground-state energy $d/2$.

\begin{theorem}[The Ornstein--Uhlenbeck operator and the quantum harmonic oscillator]
\label{thm:OU-harmonic-oscillator}
The operator
$$
-\frac d2+\frac12(P^2+Q^2),
$$
considered as a densely defined operator in $L^2(\R^d\!,\gamma)$ with
domain ${\rm Pol}(\R^d)$, is essentially selfadjoint, and its closure
equals $-L$.
\end{theorem}

\begin{proof}
By \eqref{eq:QHO2}, we have
$$
\Bigl(-\frac d2+\frac12(P^2+Q^2)\Bigr)f=-Lf,
\quad f\in{\rm Pol}(\R^d).
$$
Thus the operator in the statement is the restriction of the
selfadjoint operator $-L$ to ${\rm Pol}(\R^d)$. In particular, it is
closable. It remains to prove that ${\rm Pol}(\R^d)$ is a core for
$L$.

Let $f\in\Dom(L)$. By Theorem~\ref{thm:OU-diag} and the spectral
theorem,
$$
\sum_{n\in\N}n^2\n J_nf\n_2^2<\infty.
$$
Put $f_N:=\sum_{n=0}^N J_nf$.
By Theorem~\ref{thm:basisHHn}, each $f_N$ is a polynomial. Moreover,
$$
\n f-f_N\n_2^2 = \sum_{n>N}\n J_nf\n_2^2 \to 0
$$
and, since $LJ_nf=-nJ_nf$,
$$
\n Lf-Lf_N\n_2^2
= \sum_{n>N}n^2\n J_nf\n_2^2 \to 0.
$$
Hence $f_N\to f$ in the graph norm of $L$. Thus
${\rm Pol}(\R^d)$ is a core for $L$, and the closure of the operator in
the statement is $-L$.
\end{proof}

To recover the usual realisation of the quantum harmonic oscillator on
Lebesgue space, define
$$
H:=U\Bigl(-L+\frac d2I\Bigr)U^\star
$$
on $L^2(\R^d\!,m)$. By the definitions of $q_j$ and $p_j$ and
Theorem~\ref{thm:OU-harmonic-oscillator}, $H$ is the closure of
$$
\frac12\sum_{j=1}^d
\bigl(\wh\xi_j^{\,2}+\wh x_j^{\,2}\bigr)
=
-\frac12\Delta+\frac12|x|^2
$$
on the dense subspace $U({\rm Pol}(\R^d))$. This selfadjoint operator is
also called the {\em Hermite operator}\index{Hermite!operator}.

The identity $U(-L)U^\star=H-\frac d2I$ implies
$$
e^{-tH}=e^{-td/2}U\,OU(t)\,U^\star,
\quad t\ge0.
$$
The semigroup $(e^{-tH})_{t\ge0}$ generated by $-H$ is called
the {\em Hermite semigroup}.\index{Hermite!semigroup}
Thus, up to the scalar factor $e^{-td/2}$, it is unitarily
equivalent to the Ornstein--Uhlenbeck semigroup. In particular,
it extends to an analytic $C_0$-semigroup of contractions on
every sector $\Sigma_\omega$ with $\omega<\frac12\pi$.

Theorem~\ref{thm:OU-diag} and Corollary~\ref{cor:siL} give
$$
\sigma(H)=\Bigl\{n+\frac d2:\, n\in\N\Bigr\},
$$
with eigenspace $U\HH_n$ associated with $n+\frac d2$.
The lowest eigenvalue $d/2$ is the
{\em ground-state energy}\index{ground-state energy},
and the corresponding eigenspace is spanned by the Gaussian
$$
U\one(x)
=
2^{d/4}\exp\bigl(-\frac12|x|^2\bigr).
$$

We can also use the Hermite semigroup to construct states in
the sense of Chapter~\ref{chap:QM}.
Transporting the orthonormal Hermite basis of Theorem~\ref{thm:basisHHn} by $U$ gives an eigenbasis for $H$, with eigenvalues $|\bfn|+d/2$, indexed by $\bfn\in\N^d$.
Summing in this basis shows that $e^{-\beta H}$ is trace class
for every $\beta>0$, with
\begin{align*}
Z(\beta):=\tr(e^{-\beta H})
=e^{-\beta d/2}\prod_{j=1}^d\sum_{n=0}^\infty e^{-\beta n}
=\bigl(2\sinh(\beta/2)\bigr)^{-d}.
\end{align*}
The number $Z(\beta)$ is called the {\em partition function}.\index{partition function}
Normalising by this number gives the positive trace class operator
$$
\rho_\beta:=Z(\beta)^{-1}e^{-\beta H},
$$
which has trace one and therefore defines a state. This is the {\em Gibbs state}\index{Gibbs state}
of the harmonic oscillator, with $\beta$ interpreted as inverse temperature in units in which Boltzmann's constant is one. The operator $H\rho_\beta$ is trace class, and termwise
differentiation of the exponentially convergent spectral series gives the mean energy of this state:
$$
\tr(H\rho_\beta)
=-\frac{\ud}{\ud\beta}\log Z(\beta)
=\frac d2\coth(\beta/2).
$$
As $\beta\to\infty$, the mean energy tends to the ground-state energy $d/2$, and the Gibbs state tends to the pure ground state in trace norm. Indeed, if $P_0$ is the orthogonal projection onto
$\operatorname{span}\{U\one\}$, its weight in $\rho_\beta$ is $(1-e^{-\beta})^d$, and diagonalisation gives
$$
\n\rho_\beta-P_0\n_{\calL_1}
=2\bigl(1-(1-e^{-\beta})^d\bigr)\to0.
$$

\begin{problems}

\item
Determine the position and momentum observables on $L^2(\Z/2\Z)\simeq \C^2$\!.

\item
Find the projection-valued measures associated with
the selfadjoint operators $\wh x_j$ and $\wh \xi_j$ discussed in Section \ref{subsec:pos-mom}.

\item
Show that position and momentum are covariant with respect to rotations $R_\rho$ on $L^2(\R^d\!,m)$ given by $ R_\rho f(x) = f(\rho^{-1} x)$, where $\rho\in SO(d)$, the group of orthogonal transformations on $\R^d$ with determinant $1$.

\item
In this problem we prove {\em Wintner's theorem}:\index{theorem!Wintner}
There is no pair of {\em bounded} operators $S,T$ on a
nonzero Hilbert space $H$ satisfying the Heisenberg commutation relation
$ST-TS=iI$.
 We may absorb the imaginary constant $i$ into one of the two operators and consider the identity $ST-TS = I$. Assuming that $S,T\in \calL(H)$ satisfy $ST-TS = I$, obtain a contradiction as follows:
\begin{enumerate}[\rm(a), leftmargin=*]
  \item Show that for all $n\ge 1$ we have $S^n T -TS^n = n S^{n-1}$\!.
  \item Deduce that $S^{n-1}\not=0$ and $n\n S^{n-1}\n \le2\n S^{n-1}\n \n S\n\n T\n.$
\end{enumerate}

\item
Let $f\in C_{\rm c}^\infty(\R^{2d})$ and $h>0$, and suppose that
$\Op_h[f]$ is of trace class. Show that
$$
\tr(\Op_h[f])
=
h^{-d}\int_{\R^{2d}}f\ud m
=
\frac1{(2\pi h)^d}\int_{\R^{2d}}f(z)\ud z.
$$
\noindent{\em Hint:}\
Use \eqref{eq:W(a)f} and Fourier inversion to compute the diagonal of the integral kernel. Compress the operator to  $L^2(\ov B(0;R),m)$, apply Theorem \ref{thm:Mercer}, and let $R\to\infty$, using convergence in trace norm.

\item\label{prob:secq-compact}
Let $T$ be a contraction on $\R^d\!$.
\begin{enumerate}[\rm(a), leftmargin=*]
\item
Show that if $T$ is an isometry, then for every
$f\in L^2(\R^d\!,\gamma)$,
$$
\Gamma(T)f(x)=f(T^\star x)
$$
for almost all $x\in\R^d\!$.

\item
Show that for every $1\le p<\infty$, the operator
$\Gamma(T)$ extends to a contraction on
$L^p(\R^d\!,\gamma)$.

\item
Show that $\n\Gamma^n(T)\n=\n T\n^n$ for $n\ge1$.
Deduce that $\Gamma(T)$ is compact if and only if
$\n T\n<1$.

\item
Show that $\Gamma(rI)f\to f$ in $L^2(\R^d\!,\gamma)$
as $r\uparrow1$ for every $f\in L^2(\R^d\!,\gamma)$, but
$$
\n\Gamma(rI)-I\n=1,\quad 0\le r<1.
$$
\end{enumerate}

\item\label{prob:gaussian-poincare}
Let $\gamma$ be the standard Gaussian measure on $\R^d$,
and put $\Pi f:=(\int_{\R^d}f\ud\gamma)\one$.
\begin{enumerate}[\rm(a), leftmargin=*]
\item Show that
$$
\n\nabla f\n_2^2=\sum_{n\in\N}n\n J_nf\n_2^2,
\quad f\in H^1(\R^d\!,\gamma).
$$

\noindent{\em Hint:}\ First use the Hermite
basis for polynomials, then apply Lemma \ref{lem:H1gamma-dense}.

\item Deduce the sharp {\em Gaussian Poincar\'e inequality}
\index{inequality!Gaussian Poincar\'e}
$$
\n f-\Pi f\n_2^2\le\n\nabla f\n_2^2,
\quad f\in H^1(\R^d\!,\gamma),
$$
and show that equality holds if and only if
$f\in\HH_0\oplus\HH_1$.
\item Show that $\n OU(t)-\Pi\n=e^{-t}$ for every $t\ge0$.
\end{enumerate}

\item\label{prob:crea-anni2}
Show that, on ${\rm Pol}(\R^d)$, the following commutation relations hold:
\begin{enumerate}[label={\rm(\arabic*)}, leftmargin=*]
\item $[a_j,a_j^\dagger]=I$ for all $j=1,\dots,d$,
and $[a_j,a_k^\dagger]=0$ if $j\ne k$.
\item $[a_j,a_k]=[a_j^\dagger,a_k^\dagger]=0$ for all
$j,k = 1,\dots,d$.
\item $[p_j,q_k]=0$ if $j\ne k$.
\end{enumerate}

\item\label{prob:crea-anni3}
Let $\mathscr{W}$ be the Segal--Plancherel transform.
Show that, on ${\rm Pol}(\R^d)$, the position and momentum operators
$q_j$ and $p_j$ satisfy the relations
\begin{align*}
  q_j\circ \mathscr{W} = \mathscr{W} \circ p_j,\qquad
  p_j\circ \mathscr{W}  = -\mathscr{W}\circ  q_j.
\end{align*}

\end{problems}

%% file: appendix.tex
\chapter{Zorn's Lemma}

\blfootnote{\rm \small \begin{quote}This book has been published by Cambridge University Press in the series `Cambridge Studies in Advanced Mathematics'. The present version is free to view and download for personal use only. Not for re-distribution, re-sale or use in derivative works.\newline {\copyright} Jan van Neerven\end{quote}}

\noindent
Zorn's lemma provides a sufficient condition for the existence of maximal elements in
partially ordered sets. Its formulation uses some terminology which we introduce first.
A {\em relation}\index{relation} on a set $S$ is a subset $R$ of $S\times S$. Instead of $(x,y)\in R$ we often write $x\,R\,y$.

\begin{definition}[Partially ordered sets]
 A {\em partially ordered}\index{partial!order} set is a pair $(S,\le)$, where $S$ is a set and
 $\le$ is a relation on $S$ such that for all $x,y,z\in S$ we have:
 \begin{enumerate}[label={\rm(\roman*)}, leftmargin=*]
  \item (reflexivity): $x\le x$;
  \item (antisymmetry): if $x\le y$ and $y\le x$, then $x=y$;
  \item (transitivity): if $x\le y$ and $y\le z$, then $x\le z$.
 \end{enumerate}
 A {\em totally ordered}\index{totally!ordered} set is a partially ordered set $(S,\le)$ with the property that for all $x,y\in S$ we have $x\le y$ or $y\le x$ (or both, in which case $x=y$).
 \end{definition}

 \begin{definition}[Maximal elements, upper bounds] Let $(S,\le)$ be a partially ordered set. An element $x\in S$ is said to be
 {\em maximal}\index{maximal element}
 if $x\le y$ implies $y=x$. An element $x\in S$ is said to be an {\em upper bound}\index{upper!bound}
 for the subset $S'\subseteq S$
 if $x'\le x$ holds for all $x'\in S'$\!.
\end{definition}

Assuming the Axiom of Choice,\index{Axiom!of Choice} one has the following result.

\begin{theorem}[Zorn's lemma]\label{lem:Zorn}\index{lemma!Zorn}
 If $(S,\le)$ is a nonempty partially ordered set with the property that each of its nonempty totally ordered subsets has an upper bound in $S$, then $S$ has a maximal element.
\end{theorem}

The proof uses the following terminology.
A subset $C$ of $S$ is said to be {\em well ordered} if every nonempty subset of $C$ has a least element. If $C$ is well-ordered and $x\in C$, we write
$C_x:=\{y\in C:\,  y<x\}$.
Such a set is called a {\em proper initial segment} of $C$. More generally,
$I\subseteq C$ is called an {\em initial segment} of $C$ if it follows from
$x\in I$, $y\in C$, and $y<x$ that $y\in I$.

\begin{proof}
Write $x<y$ if $x\le y$ and $x\ne y$. We argue by contradiction and
assume that $S$ has no maximal element.

Let $C\subseteq S$ be nonempty and totally ordered. By assumption, $C$ has an upper
bound $u\in S$. Since $u$ is not maximal, there exists $v\in S$ with
$u<v$. Hence $v$ is a strict upper bound for $C$, that is, $c<v$ for all
$c\in C$. Applying the Axiom of Choice to the family of all nonempty sets
$$
   \{s\in S:\,  c<s\hbox{ for all }c\in C\},
$$
where $C$ varies over all totally ordered subsets of $S$, we obtain a map
$\Phi$ which assigns to every totally ordered $C\subseteq S$ a strict
upper bound $\Phi(C)$ of $C$.

Call a subset $C\subseteq S$ a {\em $\Phi$-chain}
if it is well-ordered by the order inherited from $S$ and for every $x\in C$ one has $x=\Phi(C_x)$, where $C_x:=\{y\in C:\,  y<x\}$.

We claim that any two $\Phi$-chains are comparable by initial segments.
Indeed, let $C$ and $D$ be $\Phi$-chains, and let $I$ be the union of
all common initial segments of $C$ and $D$. If $I\ne C$ and $I\ne D$,
let $c$ and $d$ be the least elements of $C\setminus I$ and
$D\setminus I$, respectively. Then $C_c=I=D_d$, and therefore
$$
   c=\Phi(I)=d,
$$
contradicting the definition of $I$. Thus either $C$ is an initial
segment of $D$, or $D$ is an initial segment of $C$.

Let $M$ be the union of all $\Phi$-chains. By the claim, $M$ is again a
well-ordered totally ordered subset of $S$. Moreover, if $x\in M$, then
$x$ belongs to some $\Phi$-chain $C$, and the initial segment of $M$
below $x$ agrees with the initial segment of $C$ below $x$. Hence
$$
   x=\Phi(M_x),\qquad M_x:=\{y\in M:\,  y<x\},
$$
so $M$ is itself a $\Phi$-chain.

But then $\Phi(M)$ is a strict upper bound for $M$, and
$$
   M\cup\{\Phi(M)\}
$$
is again a $\Phi$-chain. This contradicts the definition of $M$ as the
union of all $\Phi$-chains. Therefore our assumption was false, and $S$
has a maximal element.
\end{proof}

\cleardoublepage  

\cleardoublepage  

\chapter{Tensor Products}\label{sec:tensor}

\blfootnote{\rm \small \begin{quote}This book has been published by Cambridge University Press in the series `Cambridge Studies in Advanced Mathematics'. The present version is free to view and download for personal use only. Not for re-distribution, re-sale or use in derivative works.\newline {\copyright} Jan van Neerven\end{quote}}

\noindent
Let $V$ and $W$ be vector spaces. We write $\mathscr F(V\times W)$ for the vector space of all finitely supported functions from $V\times W$ into $\K$. For $(v,w)\in V\times W$ let $\delta_{(v,w)}$ denote the function which takes the value $1$ at $(v,w)$ and the value $0$ at all other points of $V\times W$. Thus every element of $\mathscr F(V\times W)$ has a unique representation as a finite sum $\sum_{j=1}^n c_j\delta_{(v_j,w_j)}$.

Let $\mathscr N$ be the linear subspace of $\mathscr F(V\times W)$ generated by all elements of the following four types:
\begin{align*}
\delta_{(v+v',w)}-\delta_{(v,w)}-\delta_{(v',w)},
\qquad
\delta_{(v,w+w')}-\delta_{(v,w)}-\delta_{(v,w')},
\end{align*}
\begin{align*}
\delta_{(cv,w)}-c\delta_{(v,w)},
\qquad
\delta_{(v,cw)}-c\delta_{(v,w)},
\end{align*}
where $c\in\K$, $v,v'\in V$, and $w,w'\in W$.

\begin{definition}[Algebraic tensor product] The (algebraic) {\em tensor product}\index{tensor!product, algebraic}\index{$V\ot W$} of $V$ and $W$ is the quotient vector space
$$ V\otimes W:=\mathscr F(V\times W)/\mathscr N.$$
For $v\in V$ and $w\in W$ we write
$$v\otimes w := \delta_{(v,w)}+\mathscr N.$$
\end{definition}

The elements $v\otimes w$ are called {\em pure tensors}. From the definition it is clear that
\begin{align*}c(v\otimes w) & = (cv)\otimes w = v\otimes(cw)\end{align*}
for all $c\in \K$, $v\in V$, and $w\in W$, and
\begin{align*}
 (v+v')\otimes w & = v\ot w+v'\ot w, \qquad
 v\otimes (w+w') = v\ot w +v\ot w'
\end{align*}
for all $v,v'\in V$ and $w,w'\in W$. Also, every element of $V\ot W$ admits a representation as a finite sum $\sum_{j=1}^k v_j\otimes w_j$.

The tensor product is characterised by the following universal property.

\begin{proposition}[Universal property]\label{prop:tensor-universal}
Let $U$ be a vector space. If $b:V\times W\to U$ is bilinear, then there exists a unique linear mapping $\widetilde b:V\otimes W\to U$ satisfying
$$\widetilde b(v\otimes w)=b(v,w),\quad v\in V,\ w\in W.$$
\end{proposition}
\begin{proof}
Define the linear mapping $\widehat b:\mathscr F(V\times W)\to U$ by
$$\widehat b\Big(\sum_{j=1}^n c_j\delta_{(v_j,w_j)}\Big):=\sum_{j=1}^n c_j b(v_j,w_j).$$
The bilinearity of $b$ implies that $\mathscr N\subseteq \ker (\widehat b)$. Hence $\widehat b$ passes to the quotient and defines a linear mapping $\widetilde b:V\otimes W\to U$ with the required property. Since the pure tensors span $V\otimes W$, this property determines $\widetilde b$ uniquely.
\end{proof}

We have natural isomorphisms of vector spaces
$$\K\otimes V \simeq V\otimes \K \simeq V,$$
given on pure tensors by $c\otimes v\mapsto cv$ and $v\otimes c\mapsto cv$.

If $V$ and $W$ have bases $(v_i)_{i\in I}$ and $(w_j)_{j\in J}$, respectively, then $(v_i\otimes w_j)_{(i,j)\in I\times J}$ is a basis for $V\otimes W$. Indeed, these tensors span $V\otimes W$ by expanding each vector appearing in a pure tensor with respect to the given bases. To prove linear independence, suppose that
$$ \sum_{(i,j)\in F} c_{ij}\,v_i\otimes w_j=0 $$
for some finite set $F\subseteq I\times J$. Let $\phi_i$ and $\psi_j$ be the coordinate functionals associated with the two bases. Applying the linear mapping induced by the bilinear form
$$ \zeta_{ij}(v,w):=\phi_i(v)\psi_j(w) $$
to the above identity gives $c_{ij}=0$ for each $(i,j)\in F$. This proves the assertion.

As a corollary, if $V$ and $W$ are finite-dimensional, with bases $(v_i)_{i=1}^{d_V}$
and $(w_j)_{j=1}^{d_W}$, then $(v_i\otimes w_j)_{i,j=1}^{d_V, d_W}$ is a basis for $V\otimes W$. In particular,
$$ \dim(V\otimes W) = \dim(V) \dim(W).$$

For vector spaces $U,V,W$, the mapping $(u,v,w)\mapsto u\ot(v\ot w)$ is trilinear from $U\times V\times W$ into $U\ot(V\ot W)$. Applying the universal property twice, we obtain a linear mapping
$$A:(U\ot V)\ot W\to U\ot(V\ot W), \qquad A((u\ot v)\ot w)=u\ot(v\ot w).$$
Similarly, the mapping $(u,v,w)\mapsto (u\ot v)\ot w$ induces a linear mapping
$$B:U\ot(V\ot W)\to (U\ot V)\ot W, \qquad B(u\ot(v\ot w))=(u\ot v)\ot w.$$
Since the pure tensors span, $A$ and $B$ are inverse to each other. Stated differently, taking tensor products is associative.
This allows us to define the tensor product
$$ U\ot V\ot W$$
as either one of the spaces in this isomorphism; for the sake of concreteness we will use the space on the left-hand side. With this in mind we can define the tensor product $V_1\otimes\cdots\ot V_N$ of vector spaces $V_n$, $n=1,\dots,N$, inductively by
$$V_1\otimes\cdots\ot V_N:= (V_1\otimes\cdots\ot V_{N-1})\ot V_N.$$
Thus, by repeated use of the universal property, linear mappings from $V_1\otimes\cdots\ot V_N$ into a vector space $Z$ may be identified with $N$-linear mappings from $V_1\times\cdots\times V_N$ into $Z$.
In what follows we write $V^{\ot n}:= V\otimes\cdots\otimes V$ for the $n$-fold tensor product of $V$.

Each permutation $\sigma\in S_n$ induces a linear operator on $V^{\ot n}$ by
$$v_1\ot\cdots\ot v_n\mapsto v_{\sigma(1)}\ot\cdots\ot v_{\sigma(n)}.$$
This is well defined by the preceding multilinear version of the universal property.
The {\em $n$-fold symmetric tensor product}\index{symmetric!tensor product} $\Gamma^n(V)$\index{$G$@$\Gamma^n(V)$} is defined
as the range of the projection $P_\Gamma$\index{$P$@$P_\Gamma$} on $V^{\ot n}$ given by
 \begin{align*}
P_\Gamma: v_1 \ot \cdots \ot v_n\mapsto
   \frac{1}{n!} \sum_{\sigma \in S_n} v_{\sigma(1)} \ot \cdots \ot
   v_{\sigma(n)}, \quad v_1, \dots, v_n \in V .
 \end{align*}
Likewise, one defines the {\em $n$-fold antisymmetric tensor product}\index{antisymmetric tensor product} $\Lambda^n(V)$,\index{$L$@$\Lambda^n(V)$}
also known as the {\em $n$-fold exterior product}\index{exterior product}, of a vector space $V$
as the range of the projection on $V^{\ot n}$ given by
 \begin{align*}
P_\Lambda: v_1 \ot \cdots \ot v_n \mapsto
   \frac{1}{n!} \sum_{\sigma \in S_n} \sign(\sigma) v_{\sigma(1)} \ot \cdots \ot
   v_{\sigma(n)}, \quad v_1, \dots, v_n \in V.
 \end{align*}
In connection with second quantisation, these spaces are sometimes denoted by $V^{\hbox{\tiny\textcircled{s}}n}$\index{$V^{\hbox{\tiny\textcircled{s}}n}$} and $V^{\hbox{\tiny\textcircled{a}}n}$,\index{$V^{\hbox{\tiny\textcircled{a}}n}$}  respectively.

We conclude with the observation that if $S$ and $T$ are linear operators on the vector spaces $V$ and $W$, respectively, then\index{$Sz$@$S\otimes T$}
$$ (S\otimes T) : v\otimes w \mapsto Sv\otimes Tw$$
uniquely defines a linear operator $S\otimes T$ on the tensor product $V\ot W$.
Indeed, the mapping $(v,w)\mapsto Sv\otimes Tw$ from $V\times W$ to $V\otimes W$ is bilinear, and the assertion follows from the universal property.

\cleardoublepage 

\cleardoublepage 

\chapter{Topological Spaces}\label{sec:topology}

\blfootnote{\rm \small \begin{quote}This book has been published by Cambridge University Press in the series `Cambridge Studies in Advanced Mathematics'. The present version is free to view and download for personal use only. Not for re-distribution, re-sale or use in derivative works.\newline {\copyright} Jan van Neerven\end{quote}}

\noindent
This appendix offers a brief treatment of topological spaces. Only notions used in the main text are covered. Several other notions are needed only in the more concrete setting of
metric spaces and will be discussed there.

\section*{Definition and General Properties}

\begin{definition}[Topological spaces]\label{def:top}
 A {\em topological space}\index{topological space} is a pair $(X,\tau)$, where $X$ is a set and $\tau$ is a {\em topology}\index{topology} on $X$, that is, $\tau$ is a collection of subsets of $X$ with the following properties:
 \begin{enumerate}[label={\rm(\roman*)}, leftmargin=*]
  \item $\emptyset\in \tau$ and $X \in \tau$;
  \item $\tau$ is closed under taking arbitrary unions;
  \item $\tau$ is closed under taking finite intersections.
 \end{enumerate}
A subset $U$ of $X$ is said to be {\em open}\index{open!set} if $U\in \tau$, and {\em closed}\index{closed!set} if its complement is open.
The {\em interior}\index{interior} $S^\circ$ of a subset $S$ is the union of all open subsets $U$ in $X$ contained in $S$.
The {\em closure}\index{closure} $\ov S$ of a subset $S$ is the intersection of all closed subsets $F$ of $X$ containing $S$.
Note that $S^\circ$ is the largest open subset of $X$ contained in $S$ and $\ov S$ is the smallest closed subset of $X$ containing $S$. A set $S$ is {\em dense}\index{dense} in a closed set $S'$ if $\ov S =S'$.

If $\mathscr{C}$ is a collection of subsets of $X$, the {\em topology generated by $\mathscr{C}$} is the intersection of all topologies on $X$
containing $\mathscr{C}$.
The topology generated by a collection $\mathscr{C}$
is the smallest topology containing every element of $\mathscr{C}$.
\end{definition}

We always assume that the scalar field $\K$ (that is, $\R$ or $\C$) carries the
topology generated by the Euclidean open balls.

\begin{remark}
The {\em relative topology}\index{relative topology} on a
subset $Y\subseteq X$ consists of the sets $Y\cap U$ with
$U$ open in $X$. For $K\subseteq Y$, compactness in $Y$
and in $X$ are equivalent.
\end{remark}

We often omit the topology $\tau$ from our notation and write $X$ instead of the more cumbersome $(X,\tau)$ to denote a
topological space, except in situations in which confusion could arise. In such situations we may speak of {\em $\tau$-open} and {\em $\tau$-closed} sets
instead of open and closed sets in order to emphasise the role of $\tau$.

 A topological space $X$ is said to be {\em Hausdorff}\index{topological space!Hausdorff} if for any two distinct points $x_1,x_2\in X$
 there exist disjoint open sets $U_1,U_2\in \tau$ such that $x_1\in U_1$ and $x_2\in U_2$.

\begin{proposition}
Finite subsets of a Hausdorff topological space are closed.
\end{proposition}
\begin{proof}
Since finite unions of closed sets are closed, it suffices to prove that every singleton
 $\{x\}$ in a Hausdorff space $(X,\tau)$ is closed. For any $y\in X\setminus\{x\}$ choose an open set $U_y$ such that $y\in U_y$ and $x\not\in U_{y}$. This is possible by the Hausdorff assumption (and actually uses less than that). We then have $\complement \{x\} = \bigcup_{y\in X\setminus\{x\}} U_y$,
 and this set is open since $\tau$ is closed under taking arbitrary unions. It follows that $\{x\}$ is closed.
\end{proof}

An important class of Hausdorff topological spaces is the class of metric spaces; they are discussed in more detail in Appendix \ref{sec:metric_spaces}. Further examples relevant to functional analysis are Banach spaces with their weak topology, dual Banach spaces with their weak$^*$ topology, and spaces of bounded operators acting between Banach spaces with their strong and weak operator topologies. For their definitions we refer to the main text.

\section*{Continuity}

Let $(X,\tau_X)$ and $(Y,\tau_Y)$ be topological spaces and consider a mapping $f:X\to Y$.

\begin{definition}
We call $f$ \emph{continuous at the point $x_0\in X$}\index{continuous!at a point}
 if for every open set $V\in \tau_Y$ containing $f(x_0)$ there exists an open set $U\in \tau_X$ containing $x_0$ such that $f(U)\subseteq V$.
We call $f$ \emph{continuous}\index{continuous}
if $f$ is continuous at every point of $X$.

A bijection $f:X\to Y$ between topological spaces is called
a {\em homeomorphism}\index{homeomorphism} if both $f$ and
its inverse $f^{-1}:Y\to X$ are continuous. The spaces $X$
and $Y$ are called {\em homeomorphic} if there exists a
homeomorphism between them.
\end{definition}

As an immediate consequence of the definition, we note that if $(X,\tau_X)$, $(Y,\tau_Y)$, $(Z,\tau_Z)$ are topological spaces and
$f:X\to Y$ is continuous at the point $x_0\in X$ and $g:Y\to Z$ is continuous at the point $f(x_0)\in Y$,
then the composition $g\circ f:X\to Z$ is continuous at the point $x_0\in X$. In particular,
the composition of two continuous mappings is continuous.

\begin{proposition}\label{prop:open-cont-top}
Let $(X,\tau_X)$ and $(Y,\tau_Y)$ be topological spaces.
For a mapping $f:X\to Y$ the following assertions are equivalent:
\begin{enumerate}[label={\rm(\arabic*)}, leftmargin=*]
\item\label{it:open-cont1-top} $f$ is continuous;
\item\label{it:open-cont2-top} $f^{-1}(V)$ is open for every
open subset~$V$ of~$Y$;
\item\label{it:open-cont3-top} $f^{-1}(F)$
is closed for every
closed subset~$F$ of~$Y$.
\end{enumerate}
\end{proposition}

\begin{proof}
\ref{it:open-cont1-top}$\Leftrightarrow$\ref{it:open-cont2-top}:\
If $f$ is continuous and $V\subseteq Y$ is open, every
$x\in f^{-1}(V)$ is contained in an open set $U_x$ such that
$f(U_x)\subseteq V$. Hence
$f^{-1}(V)=\bigcup_{x\in f^{-1}(V)}U_x$ is open.
Conversely, if \ref{it:open-cont2-top} holds and $V$ is an
open set containing $f(x_0)$, then $f^{-1}(V)$ is an
open set containing $x_0$ which is mapped into $V$.

\smallskip
\ref{it:open-cont2-top}$\Leftrightarrow$\ref{it:open-cont3-top}:\
This follows by taking complements, since
$f^{-1}(Y\setminus A)=X\setminus f^{-1}(A)$ for every
$A\subseteq Y$.
\end{proof}

Let $X$ be a topological space. A function
$f:X\to(-\infty,\infty]$ is called \emph{lower semicontinuous
at the point $x_0\in X$}\index{lower semicontinuous!function} if for every real number $a<f(x_0)$ there exists an open set $U$ containing $x_0$ such that $a<f(x)$ for all $x\in U$. The function is called {\em lower semicontinuous} if it is lower semicontinuous at every point of $X$, or equivalently, if for every $a\in\R$ the set
$$
\{x\in X:\,  f(x)>a\}
$$
is open in $X$. The function $f:X\to[-\infty,\infty)$ is called \emph{upper semicontinuous at the point $x_0\in X$}\index{upper!semicontinuous} if for every real number $a>f(x_0)$ there exists an open set $U$ containing $x_0$ such that $a>f(x)$ for all $x\in U$. The function is called {\em upper semicontinuous} if it is upper semicontinuous at every point of $X$, or equivalently,
if for every $a\in\R$ the set
$$
\{x\in X:\,  f(x)<a\}
$$
is open in $X$. If $U\subseteq X$ is open and $c\ge 0$, then the function $x\mapsto c\one_U(x)$ is lower semicontinuous; if $F\subseteq X$ is closed and $c\ge 0$, then the function $x\mapsto c\one_F(x)$
is upper semicontinuous. It is easy to prove that a function $f:X\to\R$ is continuous if and only if it is both lower semicontinuous and upper semicontinuous. Also, if $(f_i)_{i\in I}$ is a family of lower semicontinuous functions, then
$f:=\sup_{i\in I} f_i$ is lower semicontinuous, and if $(g_i)_{i\in I}$ is a family of upper semicontinuous functions, then $g:=\inf_{i\in I} g_i$ is upper semicontinuous.

\section*{Compactness}

Let $X$ be a topological space and let $S$ be a subset of $X$.
A collection $\mathscr{U} $ of open subsets
of $X$ is called an \emph{open cover}\index{open!cover}\index{cover} of $S$
if  $S\subseteq\bigcup_{U\in \mathscr{U}}U$.
A {\em subcover}\index{subcover} is a cover $\mathscr{U}'$ of $S$ contained in $\mathscr{U}$.
The set $S$ is called {\em compact}\index{compact}
if every open cover of $S$ has a finite subcover. A set is called {\em relatively compact}\index{relatively compact} if its closure is compact.

\begin{proposition}\label{prop:Hausdorff-closed} Let $X$ be a topological space. Then:
\begin{enumerate}[label={\rm(\arabic*)}, leftmargin=*]
 \item\label{it:Hausdorff-closed1}  every closed subset of $X$ contained in a compact subset of $X$ is compact;
 \item\label{it:Hausdorff-closed2}  if $X$ is Hausdorff, then every compact subset of $X$ is closed.
 \end{enumerate}
\end{proposition}

\begin{proof} \ref{it:Hausdorff-closed1}: \
Let the closed set $F$ be contained in the compact subset $S$ of $X$.
Let $\mathscr{U}_F$ be an open cover of $F$, and extend it to an open cover $\mathscr{U}$ of $S$ by adjoining the open set $\complement F$. The resulting cover of $S$ has a finite subcover, and this subcover also covers $F$. Removing the set $\complement F$ from this subcover, we are left with a finite subcover of $\mathscr{U}_F$ for $F$.
It follows that $F$ is compact.

\smallskip
\ref{it:Hausdorff-closed2}: \ Let $S$ be a compact subset of the Hausdorff space $X$. If $S$ is empty, there is nothing to prove, so we assume that $S$ is nonempty.

We first claim that for every $x\in\complement S$ there is an open set $U_x$ containing $x$ and disjoint from $S$. Indeed, for every $y\in S$, the Hausdorff property provides us with two disjoint open sets $U_{x,y}$ and
$V_{x,y}$ such that $x\in U_{x,y}$ and $y\in V_{x,y}$. The open cover $\mathscr{V}_x = \{V_{x,y}:\, y\in S\}$ of $S$
has a finite subcover, say $\mathscr{V}_x' = \{V_{x,y_1}, \dots, V_{x,y_{k_x}}\}$, where $k_x\ge 1$ is an integer depending on $x$. The set $$U_x:= \bigcap_{j=1}^{k_x} U_{x,y_j}$$ is open, contains $x$, and is disjoint from $S$. This proves the claim. But now we see that $\complement S = \bigcup_{x\in \complement S} U_x$, so $\complement S$ is open and $S$ is closed.
\end{proof}

A collection of subsets of a topological space has the {\em finite intersection property}\index{finite!intersection property} if every finite subcollection has nonempty intersection.

\begin{proposition}\label{prop:FIP}
 A  topological space $X$ is compact if and only if every collection of closed subsets of
 $X$ with the finite intersection property has nonempty intersection.
\end{proposition}
\begin{proof}
Taking complements gives a one-to-one correspondence between open
covers of $X$ and families of closed subsets with empty
intersection. An open cover has a finite subcover if and
only if the corresponding family has a finite subfamily
with empty intersection. The assertion follows.
\end{proof}

Compactness is preserved under continuous mappings:

\begin{proposition}\label{prop:compact-cont-image}
Let $X$ and $Y$ be topological spaces.
Let $f: X\to Y$ be a continuous mapping.
If $S$ is a compact subset of $X$, then $f(S)$ is compact in $Y$.
\end{proposition}

\begin{proof}
Let $\mathscr{U}$ be an open cover of $f(S)$.
Then $\{f^{-1}(U): \, U\in\mathscr{U}\}$ is an open cover
of $S$ by Proposition \ref{prop:open-cont-top}.
Since $S$ is compact, it has
a  finite subcover
$\{f^{-1}(U_1),\dots,f^{-1}(U_n)\}$.
The collection $\{U_1,\dots,U_n\}$ is then a  finite subcover of
$f(S)$.
\end{proof}

We continue with a simple, but often useful observation.

\begin{proposition}\label{prop:compact-homeomorphism}
Let $X$ be a compact topological space and let $Y$ be a
Hausdorff space. Every continuous bijection $f:X\to Y$
is a homeomorphism.
\end{proposition}

\begin{proof}
Every closed subset $F$ of $X$ is compact. Its image $f(F)$
is compact by Proposition \ref{prop:compact-cont-image},
and hence closed in $Y$ by Proposition
\ref{prop:Hausdorff-closed}. Thus the inverse mapping
$f^{-1}:Y\to X$ is continuous by Proposition
\ref{prop:open-cont-top}.
\end{proof}
Every continuous function $f:[a,b]\to\R$ has a
global maximum and a  global minimum on~$[a,b]$.
More generally, we have:

\begin{theorem}[Global maxima and minima]\label{thm:compact-max}
Let $X$ be a nonempty compact topological space and let $f: X\to\R$
be continuous.
Then $f$ attains a global maximum and a global minimum.
\end{theorem}
\begin{proof}
By Proposition \ref{prop:compact-cont-image}, the image
$f(X)$ is a nonempty compact subset of $\R$.
The open cover $\{(-n,n):\, n\ge1\}$ has a finite subcover,
so $f(X)$ is bounded. It is also closed by Proposition
\ref{prop:Hausdorff-closed}, and therefore contains its
infimum and supremum.
\end{proof}

\section*{Urysohn's lemma}

A Hausdorff space is called {\em normal}\index{normal}\index{topological space!normal} if for any two disjoint closed subsets $F$ and $G$ there exist disjoint open subsets $U$ and $V$ such that $F\subseteq U$ and $G\subseteq V$.
In Appendix \ref{sec:metric_spaces} it is shown that every metric space is normal.

\begin{proposition}\label{prop:compact-normal}
Every compact Hausdorff space $X$ is normal.
\end{proposition}
\begin{proof}
We may assume that $F$ and $G$ are nonempty disjoint closed subsets of
the compact Hausdorff space $X$, since the case where one of them is empty
is trivial. Then $F$ and $G$ are compact by Proposition \ref{prop:Hausdorff-closed}.
Fix a point $x\in F$. Since $X$ is Hausdorff, for all $y\in G$
 there exist disjoint open subsets $U_{x,y}$ and $V_{x,y}$ such that $x\in U_{x,y}$ and $y\in V_{x,y}$.
 By compactness of $G$, we find an open cover $V_{x,y_1},\dots,V_{x,y_{k_x}}$ of $G$. Set
 $$U_x:= \bigcap_{j=1}^{k_x} U_{x,y_j}, \qquad V_x:=
 \bigcup_{j=1}^{k_x} V_{x,y_j}.$$
 Then $x\in U_x$, $G\subseteq V_x$, and $U_x\cap V_x = \emptyset$. By compactness of $F$, we find an open cover
 $U_{x_1},\dots,U_{x_\ell}$ of $F$. The sets $U:= \bigcup_{j=1}^\ell U_{x_j}$ and $V:= \bigcap_{j=1}^\ell V_{x_j}$
 are open and satisfy $F\subseteq U$, $G\subseteq V$, and $U\cap V = \emptyset$.
\end{proof}

\Needspace{6\baselineskip}
\begin{proposition}\label{prop:Ury}
Let $X$ be a normal space. If $F\subseteq U\subseteq X$,
with $F$ closed and $U$ open, then there exists an open
set $V$ such that
$$F\subseteq V\subseteq\ov V\subseteq U.$$
\end{proposition}

\begin{proof}
By normality, choose disjoint open sets $V$ and $W$ such that
$F\subseteq V$ and $X\setminus U\subseteq W$.
Since $X\setminus W$ is closed and contains $V$, we have
$F\subseteq V\subseteq\ov V\subseteq X\setminus W\subseteq U$.
\end{proof}

Let $X$ be a topological space. The {\em support}\index{support} of a continuous function $f:X\to\K$ is defined as the complement of the largest open set $U\subseteq X$ such that $f\equiv 0$ on $U$ or, equivalently, as the closure of the set
$\{x\in X: \,f(x)\not=0\}$. The support of $f$ is denoted by supp$(f)$.

\begin{proposition}[Urysohn's lemma]\label{prop:Urysohn}\index{lemma!Urysohn}
Let $X$ be a normal space. If $F\subseteq U\subseteq X $ with $F$ closed and $U$ open, there exists
a continuous function $f:X\to [0,1]$ such that $f \equiv 1$ on $F$ and ${\rm supp}(f)\subseteq U$.
\end{proposition}

\begin{proof}
A rational number $q\in [0,1]$ is called {\em dyadic}\index{dyadic} if it is of the form $\frac{k}{2^n}$, where $k,n\in \N$
and $0\le k\le 2^n$\!. We will construct, for every dyadic $q\in [0,1]$, an open set $U_q$ such that
$$ F\subseteq U_q \subseteq \ov {U_q}\subseteq U$$ and, for all dyadic $r\in [0,1]$,
$$\hbox{$q>r$ implies $\ov{U_q}\subseteq U_r$}.$$

By Proposition \ref{prop:Ury} (applied twice) there exist open sets $U_0$ and $U_1$ such that
$$F\subseteq U_1\subseteq \ov{U_1} \subseteq U_0\subseteq \ov{U_0} \subseteq U.$$
Reasoning by induction, suppose that for some $n\in\N$, the open sets $U_{k/2^n}$, $k=0,\dots,2^n$, have been chosen such that $q>r$ implies $\ov{U_q}\subseteq U_r$. Using Proposition \ref{prop:Ury}, for all $k=0,\dots,2^n-1$ we find an open set
$U_{\frac{2k+1}{2^{n+1}}}$ such that
$$ \ov{U_{\frac{k+1}{2^n}}}\subseteq  U_{\frac{2k+1}{2^{n+1}}}\subseteq \ov{U_{\frac{2k+1}{2^{n+1}}}}\subseteq  U_{\frac{k}{2^n}}.$$
Then $\ov{U_q}\subseteq U_r$ holds for all dyadic $q>r$ of the form $\frac{k}{2^{n+1}}$ with $0\le k\le 2^{n+1}$\!.

Now define
$$ f_q(x):=
\begin{cases}
    q, & \hbox{if $x\in U_q$},
 \\ 0, & \hbox{otherwise},
\end{cases} \qquad
g_r(x):=
\begin{cases}
    1, & \hbox{if $x\in \ov{U_r}$},
 \\ r, & \hbox{otherwise},
\end{cases}
$$
and put $f(x):= \sup_q f_q(x)$ and $g(x):= \inf_r g_r(x).$
Then $f$ is lower semicontinuous, $g$ is upper semicontinuous, $0\le f\le 1$, $f\equiv 1$ on $F$, and ${\rm supp}(f)\subseteq U$. To conclude the proof, we show that $f=g$.

If $f_q(x)>g_r(x)$, then we must have $q>r$, $x\in U_q$, and $x\not\in \ov{U_r}$. But $q>r$ implies $U_q\subseteq U_r$.
This contradiction shows that $f_q(x)\le g_r(x)$ for all dyadic $q,r\in [0,1]$ and $x\in X$. This implies $f\le g$.

If $f(x)<g(x)$, there are dyadic numbers $q,r\in [0,1]$ such that $f(x)<r<q<g(x)$. But $f(x)<r$ implies that $x\not\in U_r$ and $g(x)>q$ implies $x\in \ov{U_q}$. This contradicts the fact that $q>r$ implies $\ov{U_q}\subseteq U_r$. It follows that $f=g$.
\end{proof}

As a first application, we prove an extension theorem for continuous functions.

\begin{theorem}[Tietze extension theorem]\label{thm:Tietze}\index{theorem!Tietze extension}
Let $F$ be a closed subset of a normal space $X$ and let $f:F\to [0,1]$ be continuous.
Then there exists a continuous function $g:X\to [0,1]$ such that $g|_F = f$.
\end{theorem}
\begin{proof}
As a preliminary step, let $M>0$ be arbitrary and let
$u:F\to[0,M]$ be continuous. The sets
$$
A:=\{x\in F:\, u(x)\le M/3\},
\qquad
B:=\{x\in F:\, u(x)\ge2M/3\}
$$
are disjoint and closed in $X$. By Urysohn's lemma there exists
a continuous function $v:X\to[0,M/3]$ such that $v=0$ on $A$
and $v=M/3$ on $B$. These properties imply that
$0\le u-v\le2M/3$ on $F$.

Starting with $r_0:=f$, apply this approximation step recursively
to $u=r_{n-1}$ and $M=(2/3)^{n-1}$. Denote the resulting
continuous function on $X$ by $g_n$ and put
$r_n:=r_{n-1}-g_n|_F$. By induction,
\begin{align*}
0\le g_n
\le \frac13\Bigl(\frac23\Bigr)^{n-1}
\quad\hbox{on }X,\qquad
0\le r_n=f-\sum_{j=1}^n g_j
\le \Bigl(\frac23\Bigr)^n
\quad\hbox{on }F.
\end{align*}
The geometric bound shows that the series
$g:=\sum_{n\ge1}g_n$ converges uniformly on $X$, so $g$ is
continuous. Moreover,
$$
0\le g\le
\sum_{n\ge1}\frac13\Bigl(\frac23\Bigr)^{n-1}=1.
$$
Since $r_n\to0$ on $F$, we have $g|_F=f$.
\end{proof}

We now turn to locally compact spaces.

\begin{definition}[Locally compact spaces]\label{def:locally-compact}
A topological space $X$ is called {\em locally compact}\index{locally!compact}\index{topological space!locally compact} if every point $x\in X$ is contained in an open set with compact closure.
\end{definition}

For example, the spaces $\K^d$ are locally compact.

For a locally compact Hausdorff space $X$, we write $C_{\rm c}(X)$ for the vector space of compactly supported continuous functions $f:X\to\K$. We let $C_0(X)$\index{$C_0(X)$} denote the space of continuous functions $f:X\to\K$ {\em vanishing at infinity}\index{vanishing at infinity}, that is, for every $\eps>0$ there exists a compact set $K\subseteq X$ such that $|f(x)|<\eps$ for all $x\in X\setminus K$. Both spaces are endowed with the supremum norm.

\begin{proposition}[Urysohn's lemma for locally compact spaces]\label{prop:Urysohn2}
Let $X$ be a locally compact Hausdorff space. If $K\subseteq  U\subseteq  X $ with $K$ compact and $U$ open, then there exists
a function $f\in C_{\rm c}(X)$ with compact support contained in $U$ such that
$0\le f\le \one$ pointwise on $X$ and $f \equiv 1$ on $K$.
\end{proposition}
\begin{proof}
If $K=\emptyset$, take $f=0$. Otherwise, for each $x\in K$ choose an open $U_x$ with compact closure contained in $U$; such $U_x$ exists by Proposition \ref{prop:Ury} applied to the compact closure of a suitable open set containing $x$.
Choose finitely many $U_{x_1},\dots,U_{x_k}$ covering $K$ and put
$W:=U_{x_1}\cup\cdots\cup U_{x_k}$. Then $K\subseteq W$ and $\overline W$ is compact and contained in $U$.
Urysohn's lemma, applied in $\overline W$ to the closed set
$K$ and the relatively open set $W$, gives a continuous
function $f:\overline W\to[0,1]$ with $f\equiv1$ on $K$,
and $S:={\rm supp}_{\overline W}(f)$ is contained in $W$.
We extend $f$ by zero outside $\overline W$.
Since $S$ is compact, it is closed in $X$. The zero extension
is continuous on the open set $W$ and identically zero
on the open set $X\setminus S$; these sets cover $X$.
Thus the extension belongs to $C_{\rm c}(X)$ and is
supported in $S\subseteq U$.
\end{proof}

As an application, we obtain the following approximation result.

\begin{proposition}\label{prop:Cc-dense-in-C0}
If $X$ is locally compact Hausdorff, then $C_{\rm c}(X)$ is dense in $C_0(X)$.
\end{proposition}
\begin{proof}
For given $f\in C_0(X)$ and $\eps>0$, let $K$ be a compact set such that $|f|<\eps$ outside $K$, and apply Proposition \ref{prop:Urysohn2} to obtain a function $g\in C_{\rm c}(X)$ such that $0\le g\le \one$ pointwise on $X$ and $g\equiv 1$ on $K$. Then $fg\in C_{\rm c}(X)$ and $\n f - fg\n_\infty \le \eps$.
\end{proof}

We next use the two versions of Urysohn's lemma to construct partitions of unity.

\begin{theorem}[Partition of unity]\label{thm:partition-unity}\index{partition of unity}
Let $X$ be a normal space and let
$$F\subseteq U_1\cup\cdots\cup U_k,$$
where $F$ is compact and the sets $U_j$ are open in $X$ for all $j=1,\dots,k$. Then there exist continuous functions $f_j:X\to[0,1]$ with support in $U_j$, $j=1,\dots,k$, such that
$$f_1+\cdots+f_k\equiv1\quad\hbox{on }F.$$
The same result holds if $X$ is a locally compact Hausdorff space; in this case the functions $f_j$ may be chosen with compact support.
\end{theorem}

\begin{proof}
If $F=\emptyset$, take $f_j=0$ for all $j$. Otherwise, every $x\in F$ belongs to some $U_j$, and we may choose an open neighbourhood $V$ of $x$ whose closure is contained in $U_j$. In the normal case this follows from Proposition \ref{prop:Ury}. In the locally compact Hausdorff case, apply Proposition \ref{prop:Urysohn2} to $\{x\}\subseteq U_j$ and take $V:=\{h>\frac12\}$, where $h$ is the resulting function. Then
$\ov V\subseteq\{h\ge\frac12\}\subseteq{\rm supp}(h)\subseteq U_j,$
and $\ov V$ is compact.

By compactness of $F$, finitely many such open sets $V_1,\dots,V_n$ cover $F$, with $\ov{V_m}\subseteq U_{j_m}$ for some $1\le j_m\le k$. Put
$$
F_j:=\bigcup_{m:\,j_m=j}\ov{V_m},\quad j=1,\dots,k.
$$
Each $F_j$ is closed and contained in $U_j$; in the locally compact Hausdorff case it is compact. By Proposition \ref{prop:Urysohn} in the normal case, and Proposition \ref{prop:Urysohn2} in the locally compact Hausdorff case, there are continuous functions $g_j:X\to[0,1]$ supported in $U_j$ such that $g_j\equiv1$ on $F_j$. In the latter case these functions have compact support.

Set $f_1:=g_1$ and
$$
f_j:=(1-g_1)\cdots(1-g_{j-1})g_j,\quad j=2,\dots,k.
$$
Then $0\le f_j\le1$, and the support of $f_j$ is contained in that of $g_j$. By induction,
\begin{align}\label{eq:pou}
f_1+\cdots+f_k=1-(1-g_1)\cdots(1-g_k).
\end{align}
Every $x\in F$ belongs to some $F_j$, so at least one $g_j(x)$ equals $1$. Thus \eqref{eq:pou} gives the required identity on $F$.
\end{proof}

Finally, local compactness and countable generation of the topology give the following useful observation.

\begin{lemma}\label{lem:sigma-compact}
Let $\Omega$ be a locally compact Hausdorff space whose topology is countably generated.
Then there exist compact subsets $K_1\subseteq K_2\subseteq\hdots$ such that
$\bigcup_{n\ge1}K_n=\Omega$.
\end{lemma}

\begin{proof}
We may assume that $\Om$ is nonempty. Let $\calG = (G_m)_{m\ge1}$ be a
countable collection of nonempty open sets generating the topology $\tau$ of $\Omega$.
Let $\calB = (B_n)_{n\ge 1}$ be an enumeration of all nonempty finite intersections of
sets in $\calG$, including the `empty intersection' $\Omega$. We claim that
$\calB$ is a {\em basis} for the topology, that is, for every open set $U$ and every
$\om\in U$ there is an $n\ge 1$ such that $\om\in B_n\subseteq U$.

By construction, every topology containing $\calG$ also contains all finite
intersections of elements of $\calG$ and is closed under arbitrary unions, hence contains
every union of elements of $\calB$. Conversely, the collection of all unions of elements
of $\calB$ is a topology containing $\calG$. Since $\tau$ is the smallest topology
containing $\calG$, it follows that every open set $U\in\tau$ can be written as a union
of elements of $\calB$.

The collection $\calB_0:=\{B\in\calB:\, \ov B\hbox{ is compact}\}$ covers $\Omega$. To see this, fix
$\omega\in\Omega$ and choose an open neighbourhood $U$ of
$\omega$ with compact closure. There exists $B\in\calB$ with
$\omega\in B\subseteq U$. Since $\ov B$ is a closed subset of
the compact set $\ov U$, it is compact, so $B\in\calB_0$.

Enumerating $\calB_0$ as $(B_n)_{n\ge1}$, we see that the sets $K_n:=\bigcup_{j=1}^n\ov{B_j}$ are compact, increase with $n$, and satisfy
$\bigcup_{n\ge1}K_n=\Omega$.
\end{proof}

\section*{Tychonoff's Theorem}\label{sec:Tichonov}

Let $I$ be a nonempty set and let $(X_i,\tau_i)_{i\in I}$ be a family of topological
spaces indexed by $I$. The {\em Cartesian product}\index{Cartesian product}
$X=\prod_{i\in I}X_i$ consists of all mappings $x:I\to\bigcup_{i\in I}X_i$
such that $x(i)\in X_i$ for every $i\in I$. For each $i\in I$ the
\emph{coordinate mapping} $p_i:X\to X_i$ is defined by $p_i(x):=x(i)$ for $x\in X$.
The \emph{product topology} on $X$ is the smallest topology for which all
coordinate mappings $p_i$ are continuous.
This topology is generated by the sets
$p_i^{-1}(U_i)$, where $U_i$ ranges over the open subsets of $X_i$.
If $I=\{i_1,\dots,i_k\}$ is finite, this topology coincides with the usual
product topology generated by sets of the form
$U_{i_1}\times\cdots\times U_{i_k}$. This immediately implies that a mapping $f:Y\to\prod_{i\in I}X_i$ is continuous if and only if every coordinate mapping $p_i\circ f$ is continuous.

\begin{theorem}[Tychonoff]\label{thm:Tychonov}\index{theorem!Tychonoff}
Endowed with the product topology, the product of any family of compact spaces is compact. If each of the spaces is Hausdorff, then so is the product.
\end{theorem}

\begin{proof}
Let $X=\prod_{i\in I}X_i$, where each $X_i$ is compact.
The assertion is immediate if $I=\emptyset$ or $X=\emptyset$,
so we assume otherwise.

Let $\mathscr C$ be a collection of closed subsets of $X$
with the finite intersection property. By Proposition
\ref{prop:FIP}, it suffices to prove that
$\bigcap_{C\in\mathscr C}C\ne\emptyset$.
Using Zorn's lemma, we extend $\mathscr C\cup\{X\}$ to a maximal
family $\mathscr M$ of subsets of $X$ with the finite
intersection property. Indeed, the union of any nonempty
totally ordered collection of such families again has the
finite intersection property and contains $\mathscr C\cup\{X\}$.

By maximality, $\mathscr M$ is closed under finite intersections, for
adjoining an intersection of finitely many members preserves
the finite intersection property. Consequently, if a subset
$A\subseteq X$ intersects every member of $\mathscr M$, then
$\mathscr M\cup\{A\}$ also has the finite intersection property,
and maximality gives $A\in\mathscr M$.

For each $i\in I$, the family $\{\ov{p_i(M)}:\, M\in\mathscr M\}$
has the finite intersection property. Indeed, for
$M_1,\dots,M_n\in\mathscr M$, the set
$p_i(M_1\cap\cdots\cap M_n)$ is nonempty and contained in
$\bigcap_{j=1}^n\ov{p_i(M_j)}$.
By compactness of $X_i$, we may therefore choose
$$
x(i)\in\bigcap_{M\in\mathscr M}\ov{p_i(M)},
\quad i\in I.
$$
These choices define a point $x\in X$.

If $U_i\subseteq X_i$ is open and contains $x(i)$, then
$U_i$ intersects $p_i(M)$ for every $M\in\mathscr M$.
Thus $p_i^{-1}(U_i)$ intersects every member of $\mathscr M$,
and hence belongs to $\mathscr M$.
Every basic open neighbourhood of $x$ is a finite intersection
of such sets, so it too belongs to $\mathscr M$ and therefore
intersects every $C\in\mathscr C$.
Since each $C$ is closed, this implies that $x\in C$.
We have proved that $x\in\bigcap_{C\in\mathscr C}C$,
and hence $X$ is compact.

If each $X_i$ is Hausdorff and $x\ne x'$, then
$x(i)\ne x'(i)$ for some $i\in I$.
Choose disjoint open neighbourhoods of $x(i)$ and $x'(i)$
in $X_i$. Their inverse images under $p_i$ are disjoint
open neighbourhoods of $x$ and $x'$ in $X$.
Thus $X$ is Hausdorff.
\end{proof}

\cleardoublepage  

\cleardoublepage  

\chapter{Metric Spaces}\label{sec:metric_spaces}

\blfootnote{\rm \small \begin{quote}This book has been published by Cambridge University Press in the series `Cambridge Studies in Advanced Mathematics'. The present version is free to view and download for personal use only. Not for re-distribution, re-sale or use in derivative works.\newline {\copyright} Jan van Neerven\end{quote}}

\noindent
We now introduce an important class of Hausdorff spaces, namely, the class of metric spaces.
All results of the previous appendix apply to metric spaces, but to make
the present appendix self-contained, some proofs are repeated. We also
discuss a number of additional topics.

\section*{Definition and General Properties}

\begin{definition}[Metric spaces]\label{def:metric-space}
A \emph{metric space}\index{metric space}
is a pair~$(X,d)$, where $X$ is a set and $d$ a \emph{metric}\index{metric} (or {\em distance function})\index{distance function} on $X$, that is, a function $d: X\times X\to [0,\infty)$ such that for all~$x$,
$y$,~$z$ in~$X$ the following conditions are satisfied:
\begin{enumerate}[label={\rm(\roman*)}, leftmargin=*]
\item $d(x,y)=0 \Leftrightarrow x=y$;
\item $d(x,y)=d(y,x)$;
\item $d(x,z)\le d(x,y)+d(y,z)$ \ (the \emph{triangle inequality}).\index{triangle inequality}\index{inequality!triangle}
\end{enumerate}
\end{definition}

In what follows we omit the distance function $d$ from the notation and write $X$ instead of the more cumbersome $(X,d)$ to denote a metric
space, except in situations in which confusion could arise.

Let $X$ be a metric space, let $x\in X$, and let $r>0$.
The set\index{$Bba$@$B(x;r)$}
$$
B(x;r):=\{y\in X:\,  d(x,y)<r\}
$$
is called the {\em open ball}\index{open!ball} with centre $x$ and radius $r$.
The set\index{$Bbb$@$\ov{B}(x;r)$}
$$
\ov{B}(x;r):=\{y\in X:\,  d(x,y)\le r\}
$$
is called the {\em closed ball}\index{closed!ball} with centre $x$ and radius $r$.

A subset $S$ of a metric space $X$ is called {\em open}\index{open!set} if for all $x\in S$ there exists an $r>0$ such that
$B(x;r)\subseteq S$. A subset $S$ of a metric space $X$ is called {\em closed}\index{closed!set} if its complement $\complement S = X\setminus S$ is open.
It is an easy consequence of the triangle inequality that
every open ball $B(x;r)$ is open and  every closed ball $\ov{B}(x;r)$ is closed; this
justifies the terminology `open ball' and `closed ball'.

The {\em interior}\index{interior} of a subset $S$ of a metric space $X$ is the union of all open subsets of $X$ contained in $S$
and is denoted by $S^\circ$. It is the largest open subset of $X$ contained in $S$.
The {\em closure}\index{closure} of a subset $S$ of a metric space $X$ is the intersection of all closed subsets of $X$ containing $S$
and is denoted by $\ov S$. It is the smallest closed subset of $X$ containing $S$.

The closure $\ov{B(x;r)}$ of an open ball is always contained in the closed ball $\ov{B}(x;r)$, but this inclusion may be strict.
For example, take $X = \Z$ with distance function $d(m,n)=|n-m|$. The open ball ${B(0;1)} = \{0\}$ is also closed, so its closure equals
$\ov{B(0;1)} = \{0\}$. On the other hand, $\ov{B}(0;1) = \{-1,0,1\}$.

For $x\in X$ and a nonempty subset $S\subseteq X$, let\index{$Dx$@$d(x,S)$} $$d(x,S):=\inf_{y\in S}d(x,y).$$ A routine argument shows that
$ |d(x,S) - d(x',S)| \le d(x,x'),$ and therefore $x\mapsto d(x,S)$ is continuous. If $S$ is nonempty and closed and $x\not\in S$, then $B(x;r)\cap S = \emptyset$ for some $r>0$, and therefore $d(x,S)\ge r >0$.

For any metric space $(X,d)$, the collection of its open sets is a Hausdorff topology. Indeed, if $x,x'\in X$ satisfy $x\not= x'$, then $r:= d(x,x')>0$ and the open balls $B(x;\frac12 r)$ and $B(x';\frac12 r)$ are disjoint. More is true:

\begin{proposition}\label{prop:metric-normal}
Every metric space is normal.
\end{proposition}

\begin{proof}
Let $F$ and $G$ be disjoint closed sets of a metric space $X$. We must find disjoint open sets $U$ and $V$ such that $F\subseteq U$ and $G\subseteq V$.
We may assume that $F$ and $G$ are both nonempty, since otherwise the result is trivial.

By the observations preceding the proposition, the function $$f(x):= \frac{d(x,F)}{d(x,F)+d(x,G)}$$
 is well defined, continuous, takes values in $[0,1]$, and satisfies $f \equiv 0$ on $F$ and $f\equiv 1$ on $G$.
 The sets $U := \{x\in X:\, f(x) <\frac12\}$ and
 $V := \{x\in X:\, f(x) >\frac12\}$ are open
 and have the desired properties.
\end{proof}

\section*{Convergence}

A sequence $(x_n)_{n\ge 1}$ in a metric space $X$
is called {\em convergent}\index{convergent} if there exists an $x\in X$ such that
for all $\varepsilon >0$ there is an index $N\ge 1$ with the property
$  d(x_n, x) < \varepsilon$ for all $n\ge N$.
We then write $$\lim_{n\to\infty} x_n = x$$
and call $x$ a {\em limit}\index{limit} of the sequence $(x_n)_{n\ge 1}$.
It is clear that
$ \lim_{n\to\infty} x_n = x$ if and only if $ \lim_{n\to\infty} d(x_n, x)=0.$
Limits are unique, for if $\lim_{n\to\infty} x_n = x$
and $\lim_{n\to\infty} x_n = y$, then for all indices $n$ the triangle inequality gives
$0\le  d(x,y) \le d(x,x_n) + d(x_n,y)\to 0+0=0$
and therefore $d(x,y)= 0$.

A subset $S$ of a metric space $X$ is called
{\em sequentially closed}\index{sequentially!closed}
if the limit of every sequence in $S$ that converges in $X$
belongs to $S$.

\begin{proposition}\label{prop:sequentially-closed}
For a subset $S$ of a metric space $X$, the following assertions are equivalent:
\begin{enumerate}[label={\rm(\arabic*)}, leftmargin=*]
\item\label{it:sequentially-closed1} $S$ is closed;
\item\label{it:sequentially-closed2} $S$ is sequentially closed.
\end{enumerate}
\end{proposition}

 \begin{proof}
 \ref{it:sequentially-closed1}$\Rightarrow$\ref{it:sequentially-closed2}: \ Let $(x_n)_{n\ge 1}$ be a sequence in $S$, convergent in $X$
 with limit $x$. We need to show that $x\in S$.
 Assume the contrary. Then $x\in \complement S$, and since
  $\complement S$ is open there is an  $\e>0$ such that
  $B(x;\eps)\subseteq \complement S$.
 On the other hand, since $(x_n)_{n\ge 1}$ converges to $x$ there is an  index $N\ge 1$ such that
 $d(x_n,x)<\e$ for all $n\ge N$, that is, $x_n \in B(x;\eps)$ (and hence $x_n\in \complement S$) for all $n\ge N$.
 This contradicts $x_n\in S$.

 \smallskip
 \ref{it:sequentially-closed2}$\Rightarrow$\ref{it:sequentially-closed1}: \ We need to show that $\complement S$ is open.
 Choose $x\in \complement S$ arbitrarily.
 We must show  that there is an $\e>0$ such that
 $B(x;\eps)\subseteq \complement S$.
 Suppose that such an $\e>0$ does not exist. Then for every $n\ge 1$
 we can find an  $x_n \in B(x;\frac1n)\cap S$. The resulting sequence $(x_n)_{n\ge 1}$ is contained in
 $S$ and satisfies $d(x_n,x) < \frac1n$ for all $n\ge 1$,
that is, we have $\lim_{n\to\infty} x_n = x$.
 Since $S$ is sequentially closed we conclude that $x\in S$,
 contradicting the assumption that $x\in\complement S$.
\end{proof}

As a corollary, we obtain the following useful criterion for determining which
elements belong to the closure of a given set:

\begin{proposition} For a subset $S$ of a metric space $X$ and a point $x\in X$, the following assertions are equivalent:
\begin{enumerate}[label={\rm(\arabic*)}, leftmargin=*]
\item\label{it:closed1} $x\in\overline S$;
\item\label{it:closed2} $S\cap B(x;\e) \not=\emptyset$ for all $\e>0$;
\item\label{it:closed3} there exists a sequence $(x_n)_{n\ge 1}$
in $S$ with $\lim_{n\to\infty}x_n=x$.
\end{enumerate}
In particular, a set $S$ is dense in the closed set $S'$ if and only
if every $s'\in S'$ is the limit of a sequence in $S$.
\end{proposition}

\begin{proof}
\ref{it:closed1}$\Rightarrow$\ref{it:closed2}: \ If $S\cap B(x;\e) = \emptyset$ for some $\e>0$,
then
$S\subseteq \complement  B(x;\e)$. Since $\complement  B(x;\e)$ is closed, this implies
$\ov{S}\subseteq  \complement  B(x;\e)$ and therefore $x\not\in\ov{S}$.

\smallskip
\ref{it:closed2}$\Rightarrow$\ref{it:closed3}: \
For every $n\ge 1$ we choose $x_n\in S\cap B(x;\frac1n)$.
In this way we obtain a sequence $(x_n)_{n\ge 1}$ in $S$ converging to $x$.

\smallskip
\ref{it:closed3}$\Rightarrow$\ref{it:closed1}: \ If $x\not\in\ov{S}$, then
$x\in \complement  \ov{S}$ and this set is open.
Hence there exists  an $\e>0$ such that $B(x;\e)\subseteq \complement \ov{S}$.
In particular, $d(x,y)\ge\e$ for all $y\in S$.
This implies that no sequence in $S$ can converge to $x$.
\end{proof}

\section*{Completeness}

A sequence $(x_n)_{n\ge 1}$ in a metric space $X$ is a
\emph{Cauchy sequence}\index{Cauchy sequence} if for all~$\e>0$ there is an index $N\ge 1 $
such that $ d(x_n,x_m)<\e$ for all  $n,m\ge N$.

\begin{proposition}\label{prop:Cauchy-sequence}
In a metric space $X$ the following assertions hold:
\begin{enumerate}[label={\rm(\arabic*)}, leftmargin=*]
\item\label{it:Cauchy-sequence1} every convergent sequence is a Cauchy sequence;
\item\label{it:Cauchy-sequence2} every Cauchy sequence with a convergent subsequence is convergent.
\end{enumerate}
\end{proposition}
\begin{proof}
\ref{it:Cauchy-sequence1}: \ Let $x$ be the limit  of the convergent sequence $(x_n)_{n\ge 1}$ and let $\e>0$ be arbitrary.
Choose $N\ge 1$ so large that
$d(x_k,x)< \e$ for all $k\ge N$. By the triangle inequality
it follows that, for all $n,m\ge N$, we have
$ d(x_n,x_m)\le d(x_n,x)+ d(x,x_m) <\e+ \e = 2\e.$

\smallskip
\ref{it:Cauchy-sequence2}: \  Let $(x_n)_{n\ge 1}$ be a Cauchy sequence with a subsequence $(x_{n_k})_{k\ge 1}$ convergent to $x$.
We check that $(x_n)_{n\ge 1}$ converges to $x$.
Choose $\e>0$ arbitrarily and let $N\ge 1$ be such that $d(x_n,x_m)<\e$ for all $m,n\ge N$.
Let $K\ge 1$ be such that for all $k\ge K$ we have both $n_k\ge N$ and $d(x_{n_k},x)<\e $. Now choose $k\ge K$ arbitrarily.
Then for all $n\ge N$ we have
$d(x_n,x)\le d(x_n,x_{n_k})+ d(x_{n_k},x) < \e + \e = 2\e.$
\end{proof}

A subset $S$ of a metric space $X$ is called {\em complete}\index{complete}
if every Cauchy sequence contained in $S$ converges to a limit in $S$.
In particular, $X$ is complete if every Cauchy sequence in $X$ is
convergent in $X$.
Every complete set $S$
is sequentially closed and hence closed. Conversely, if $X$ is complete, then every closed subset $S$ is complete.

\begin{theorem}[Completion]\label{thm:completion-metric}
If $X$ is a metric space, there exists a complete metric space $(\ov X, \ov d)$ and a mapping $i:X\to \ov X$
with the following properties:
\begin{enumerate}[label={\rm(\roman*)}, leftmargin=*]
 \item\label{it:completion-metric1}  $i$ is isometric, that is, $\ov d(ix,iy) = d(x,y)$ for all $x,y\in X$;
 \item\label{it:completion-metric2}  $i$ has dense range, that is, $i(X)$ is dense in $\ov X$.
\end{enumerate}
Moreover, if $(\ov{\ov{X}}, \ov{\ov{d}})$ is another complete metric space and $i': X \to \ov{\ov{X}}$ is a mapping satisfying {\rm\ref{it:completion-metric1}}
and {\rm\ref{it:completion-metric2}}, then the identity mapping on $X$ has a unique extension to an isometry from $\ov X$ onto $\ov{\ov{X}}$.
\end{theorem}

A more precise way of stating the last assertion is that there exists a unique
isometry $j$ from $\ov X$ onto $\ov{\ov{X}}$ which satisfies $j(ix) = i'x$ for all $x\in X$.

\begin{proof}
On the set of Cauchy sequences in $X$ we define an equivalence relation by declaring the Cauchy
 sequences $(x_n)_{n\ge 1}$ and  $(x_n')_{n\ge 1}$ equivalent if $\limn d(x_n,x_n')= 0$.
Let $\ov X$ be the set of all equivalence classes. The elements of $\ov X$ are denoted by $\ov x$.
On $\ov X$ we define a metric $\ov d$ by
$$ \ov d(\ov x,\ov y):= \limn d(x_n,y_n),$$
where $(x_n)_{n\ge 1}$ and $(y_n)_{n\ge 1}$ are Cauchy sequences representing $\ov x$ and $\ov y$, respectively.
Using the triangle inequality, one readily checks that
this limit indeed exists and is independent of the choice of the representing sequences.

If $x\in X$, the constant sequence $(x)_{n\ge 1}$ is Cauchy and therefore
defines an element of $\ov X$ which we shall denote by $[x]$. We thus obtain a mapping
$i:X\to \ov X$ by declaring $ix:= [x]$. From
$ \ov d(ix,iy) = \ov d([x],[y]) = \limn d(x,y) = d(x,y)$ we see that $i$ is isometric.

This gives property \ref{it:completion-metric1}. To prove property \ref{it:completion-metric2}, let $(x_n)_{n\ge 1}$ be a
Cauchy sequence in $X$ representing $\ov x\in \ov X$. For any $\e>0$ we may choose $N$ so large that
$d(x_n,x_m)<\e$ for all $m,n\ge N$. Then
$ \ov d(ix_N,\ov x) = \limn d(x_N,x_n) \le \e.$ This shows that
$i(X)$ is dense in $\ov X$.

To prove that $\ov X$ is complete, let $(\ov x_n)_{n\ge 1}$ be a Cauchy sequence in $\ov X$.
By the density of $i(X)$ in $\ov X$ we may pick elements $x_n\in X$ such that
$\ov d( \ov x_n, [x_n]) <\frac1n$.
From
\begin{align*}d(x_n, x_m) = \ov d([x_n],[x_m])  \le \ov d([x_n], \ov x_n) + \ov d(\ov x_n, \ov x_m) + \ov d(\ov x_m, [x_m]) \le \frac1n+ \ov d(\ov x_n,  \ov x_m) +\frac1m
\end{align*}
and the fact that the right-hand side tends to $0$ as $n,m\to\infty$,
we infer that $(x_n)_{n\ge 1}$ is a Cauchy sequence in $X$.
Let $\ov x\in \ov X$ be its equivalence class. As was shown in the proof of density, we
have $\limn \ov d([x_n], \ov x) = 0$.
Then
$$ \ov d( \ov x_n,  \ov x) \le   \ov d( \ov x_n, [x_n]) + \ov d( [x_n] ,\ov x) < \frac1n + \ov d( [x_n],  \ov x)$$
shows that $\limn \ov d( \ov x_n, \ov x) = 0$.
This proves the completeness of $\ov X$.

Let $(\ov{\ov{X}}, \ov{\ov{d}})$ and $i': X\to \ov{\ov{X}}$ satisfy \ref{it:completion-metric1} and \ref{it:completion-metric2}. We obtain a well-defined isometry $j$ from $\ov X$ onto $\ov{\ov{X}}$ by putting
$$ j\ov x :=
\limn i' x_n,$$
where $(x_n)_{n\ge 1}$ is a Cauchy sequence representing $\ov x$ and the limit on the right-hand side is taken in $\ov{\ov{X}}$;
this limit exists because the sequence $(i'x_n)_{n\ge 1}$ is Cauchy in the complete space $\ov{\ov{X}}$. This isometry is onto: its range is closed and contains the dense image of $X$.
If $k:\ov X\to \ov{\ov{X}}$  is any isometry satisfying $kix = i'x$ for all $x\in X$, then
$$k\ov x = \limn kix_n = \limn i' x_n = \limn jix_n = j\ov x$$ for all $\ov x\in \ov X$ represented by the Cauchy sequence $(x_n)_{n\ge 1}$, and therefore $j=k$.
\end{proof}

A complete metric space
$(\ov X,\ov d)$ is called a {\em completion}\index{completion} of $(X,d)$ if there
exists a mapping $i:X\to \ov X$ with properties \ref{it:completion-metric1} and \ref{it:completion-metric2}.
The second part of the theorem asserts that completions are `unique up to isometry'.
To avoid this minor ambiguity we agree
to call the metric space $(\ov X,\ov d)$ constructed in the above proof `the' completion of $(X,d)$.

\section*{Continuity}

Let $X$ and $Y$ be metric spaces and consider a mapping $f:X\to Y$.

\begin{definition}
 A mapping $f:X\to Y$ is {\em continuous at the point}\index{continuous!at a point} $x_0\in X$ if and only if for every~$\e>0$ there exists a  $\delta>0$ such that for all~$x\in X$ with $d_X(x,x_0)<\delta$ we have
$d_Y(f(x),f(x_0))<\e$, and $f$ is {\em continuous}\index{continuous} if $f$ is continuous at every point of $X$.
\end{definition}

This definition is consistent with the one in the previous appendix; this is clear from the fact that every open set $U$ in $X$ containing $x_0$ contains the open balls $B(x_0;\delta)$
for sufficiently small $\delta>0$ and every open set $V$ in $Y$ containing $f(x_0)$ contains the open balls $B(f(x_0);\eps)$
for sufficiently small $\eps>0$.

A mapping $f:X\to Y$ is called {\em sequentially continuous at the point $x_0\in X$}
if for every sequence $(x_n)_{n\ge 1}$ in $X$ with $\lim_{n\to\infty}x_n=x_0$
we have $\lim_{n\to\infty}f(x_n)=f(x_0)$.
We call $f$ {\em sequentially continuous}
\index{sequentially!continuous}
if $f$ is sequentially continuous at every point of~$X$.

\begin{proposition}\label{prop:sequentially-cont}
For a mapping $f:X\to Y$ between the metric spaces $X$ and $Y$, the following assertions are equivalent:
\begin{enumerate}[label={\rm(\arabic*)}, leftmargin=*]
\item\label{it:sequentially-cont1} $f$ is continuous at the point $x_0\in X$;
\item\label{it:sequentially-cont2} $f$ is sequentially continuous at the point $x_0\in X$.
\end{enumerate} In particular,  $f$ is continuous if and only if
$f$ is sequentially continuous.
\end{proposition}

 \begin{proof}
 \ref{it:sequentially-cont1}$\Rightarrow$\ref{it:sequentially-cont2}: \ Suppose that $\lim_{n\to\infty}x_n=x_0$
 in $X$.
 Choose $\e>0$ arbitrarily and choose, using the continuity of $f$ at $x_0$,
 a   $\delta>0$ such that $d_Y(f(x),f(x_0))<\e$
 for all $x\in X$ with $d_X(x,x_0)<\delta$. Since
 $\lim_{n\to\infty} d_X(x_n,x_0)=0$ we can find an index $N\ge 1$ such
 that $d_X(x_n,x_0)<\delta$ for all $n\ge N$.
 For all $n\ge N$ we then have $d_Y(f(x_n),f(x_0))<\e$.

 \smallskip
 \ref{it:sequentially-cont2}$\Rightarrow$\ref{it:sequentially-cont1}: \
 Suppose that there exists an $\e>0$ for which no $\delta>0$ can be found such that
 $d_Y(f(x),f(x_0))<\e$ for all $x\in X$ with $d_X(x,x_0)<\delta$.
 Then for every $n\ge 1$ we can find $x_n\in X$ with
 $d_X(x_n,x_0)<\frac1n$ and $d_Y(f(x_n),f(x_0)) \ge \e$.
 But this implies that $f$ is not sequentially continuous.
\end{proof}

A function $f:X\to Y$ is
 {\em uniformly continuous}\index{uniformly!continuous} if for every $\e>0$ there exists a $\delta>0$ such that $d_Y(f(x),f(y))<\e$ whenever $x,y\in X$ satisfy $d_X(x,y)<\delta$.
Every uniformly continuous function is continuous, but the converse is false even for bounded functions: the function
$f:(0,1)\to [-1,1]$, $f(x) = \sin(1/x)$, is continuous but not uniformly continuous.
In the next section we prove that if $X$ is compact, then every continuous function from
$X$ into another metric space is uniformly continuous.

Uniformly continuous functions have the following extension property:

\begin{proposition} Let $X$ and $Y$ be metric spaces, with $Y$ complete.
 If $f:X\to Y$ is uniformly continuous, there exists a unique uniformly continuous function
 $\ov f:\ov X\to Y$ extending $f$.
\end{proposition}
\begin{proof} We identify $X$ and $Y$ with their natural images in the completions $\ov X$ and $\ov Y$.
 Let $x\in \ov X$ and choose a sequence $x_n\to x$ with each $x_n$ in $X$. Then $(x_n)_{n\ge 1}$
 is Cauchy in $X$, and the uniform continuity of $f$ implies that $(f(x_n))_{n\ge 1}$ is
 Cauchy in $Y$. Since $Y$ is complete, this sequence converges to a limit, say $y$. We
 set $\ov f(x):= y$. We need to check that $\ov f$ is well defined (that is, $\ov f(x)$ does not
 depend on the choice of the approximating sequence) and is uniformly continuous.

 If $\wt x_n\to x$ with each $\wt x_n$ in $X$, then $\limn d_X(x_n,\wt x_n) = 0$.
 By the uniform continuity of $f$ it follows that $\limn d_Y(f(x_n),f(\wt x_n)) = 0$,
 and therefore $\limn f(\wt x_n) = \limn f(x_n)$. This proves that $\ov f$ is well defined.

 To prove that $\ov f$ is uniformly continuous, let $\e>0$. Using the definition
 of uniform continuity, choose $\delta>0$ so that
$d_Y(f(u),f(v))<\frac12\e$ whenever $d_X(u,v)<\delta$. If $\ov d_X(x,y)<\delta$ in $\ov X$, and if $(x_n)_{n\ge 1}$
 and $(y_n)_{n\ge 1}$ are approximating sequences in $X$, then
 $(f(x_n))_{n\ge 1}$
 and $(f(y_n))_{n\ge 1}$ are approximating sequences for $\ov f(x)$ and $\ov f(y)$ in $Y$,
and for large enough $n$ we have  $d_X(x_n,y_n)<\delta$ and $d_Y(f(x_n),f(y_n))<\frac12\e$. Hence
$$
d_Y(\ov f(x),\ov f(y))
=\limn d_Y(f(x_n),f(y_n))
\le \frac12\e<\e.
$$

It is clear from the construction that $\ov f$ extends $f$.

To prove uniqueness, suppose that $g:\overline X\to Y$ is another uniformly continuous extension of $f$. Then for any $x\in\overline X$ and any sequence $(x_n)\subseteq X$ with $x_n\to x$ we have $g(x)=\lim_{n} g(x_n)=\lim_{n} f(x_n)=\overline f(x)$. Thus $g=\overline f$.
\end{proof}

\section*{Compactness}

Let $(X,d)$ be a metric space. We recall that
a subset $S$ of a metric space $X$ is {\em compact}
if every open cover of $S$ has a finite subcover, and {\em relatively compact} if its closure $\ov S$ is compact.
In order to characterise compactness in terms of sequences, we introduce the following terminology.
A subset $S$ of a  metric space $X$ is called
{\em sequentially compact}\index{sequentially!compact} when every sequence in $S$ has a convergent subsequence with limit in $S$.

A subset
$S$ of a metric space $X$ is called {\em totally bound\-ed}\index{totally!bounded}\index{bounded!totally} if
for every $r>0$ there is a finite cover of $S$ with balls in $X$ of radius $r$.

\begin{theorem}[Compactness and total boundedness]\label{thm:totally-bounded}
For a  subset $S$ of a metric space $X$
the following assertions are equivalent:
\begin{enumerate}[label={\rm(\arabic*)}, leftmargin=*]
\item\label{it:totally-bounded1} $S$ is compact;
\item\label{it:totally-bounded2} $S$ is sequentially compact;
\item\label{it:totally-bounded3} $S$ is complete and totally bounded.
\end{enumerate}
\end{theorem}

\begin{proof}
\ref{it:totally-bounded1}$\Rightarrow$\ref{it:totally-bounded2}: \
Suppose that  $(x_n)_{n\ge 1}$ is a  sequence in $S$ not containing any subsequence
converging to an element of $S$. We will construct an open cover of $S$
without a finite subcover.

The assumption entails that for every $x\in S$ there is an $\e(x)>0$ with the property
that the open ball $B(x;\e(x))$ contains at most finitely many terms of the sequence $(x_n)_{n\ge 1}$.
Let  $\mathscr{U} = \{B(x;\e(x)):\,  x\in S\}$.
This is an open cover of $S$ without a finite subcover, as
every $U\in\mathscr{U}$ contains at most finitely many terms of
the sequence $(x_n)_{n\ge 1}$. But this sequence has infinitely many distinct
terms: otherwise we could immediately pick a convergent subsequence.

\smallskip
\ref{it:totally-bounded2}$\Rightarrow$\ref{it:totally-bounded3}: \ Suppose
that $S$ is not totally bounded.
We will construct a sequence in $S$ without any subsequence converging to an element of $S$.

By our assumption there exists an $\e>0$ such that
$S$ has no finite cover by $\e$-balls in $X$.
Choose $x_1\in S$ arbitrarily. The collection $\{B(x_1;\e)\}$ does not cover $S$, so
there is an $x_2\in S$ with $x_2\not\in
B(x_1;\e)$. Note that $d(x_1,x_2)\ge \e$.

The collection $\{B(x_1;\e), B(x_2;\e)\}$ is not
a cover of $S$, so there is an $x_3\in S$ with $x_3\not\in
B(x_1;\e)\cup B(x_2;\e)$. Note that $d(x_1,x_3)\ge \e$ and $d(x_2,x_3)\ge \e$.

Continuing this way we obtain a sequence $(x_n)_{n\ge 1}$ with the property that
$d(x_n,x_m)\ge \e$ for all choices of $n\not=m$. This sequence has no
Cauchy subsequence, and therefore no convergent subsequence.

Next we prove that $S$ is complete.
Suppose that $(x_n)_{n\ge 1}$ is a Cauchy sequence in $S$. Since $S$ is sequentially compact, this sequence has a  convergent subsequence with limit $x$ in $S$.
But then $(x_n)_{n\ge 1}$ itself converges to $x$. This shows that $S$ is complete.

\smallskip
\ref{it:totally-bounded3}$\Rightarrow$\ref{it:totally-bounded1}: \
Suppose, for a contradiction,  that $S$ is
complete and totally bounded but not compact.

Since $S$ is not compact, there is an open cover $\mathscr{U}$
of $S$ without a finite subcover.
Since $S$ is totally bounded, for every $n\ge 1$ we can find a
finite cover $\calB_n$ of $S$ consisting of $\frac1n$-balls.

There is a ball $B_1\in \calB_1$ such that $S\cap B_1$ cannot be covered by
finitely many open sets in $\mathscr{U}$.
Similarly, there is a ball $B_2\in\calB_2$ such that $S\cap B_1\cap B_2$ cannot be covered by finitely many open sets in $\mathscr{U}$.
Continuing in this way we find a sequence of balls $B_k\in\calB_k$ such that
$S\cap B_1\cap\cdots\cap B_k$ cannot be covered by
finitely many open sets in $\mathscr{U}$.

For each $k\ge 1$ select a point $y_k\in S\cap B_1\cap\cdots\cap B_k$.
The sequence $(y_k)_{k\ge 1}$ is a Cauchy sequence in $S$.
By completeness, the sequence $(y_k)_{k\ge 1}$ converges to a limit $x$ which belongs to  $S$.
Choose $U\in\mathscr{U}$ such that $x\in U$. Since $U$ is open, there exists
$r>0$ such that $B(x;r)\subseteq U$.

Write $B_N=B(x_N;\tfrac1N)$. Since $y_N\in B_N$, we have $d(x_N,y_N)<\tfrac1N$.
Choose $N\ge 1$ so large that $\tfrac1N<\tfrac14 r$ and $d(y_N,x)<\tfrac14 r$.
Then
$$
d(x_N,x)\le d(x_N,y_N)+d(y_N,x)<\tfrac12 r,
$$
and therefore
$$
B_1\cap\cdots\cap B_N \subseteq B_N = B(x_N;\tfrac1N)\subseteq B(x_N;\tfrac12 r)\subseteq B(x;r)\subseteq U.
$$
But this means that $S\cap B_1\cap\cdots\cap B_N$ is covered by the finite subcollection
$\{U\}$ of $\mathscr{U}$. This contradiction concludes the proof.
\end{proof}

The equivalence of \ref{it:totally-bounded1} and \ref{it:totally-bounded3}
implies that in a complete metric space, a subset is compact if and only if it is closed and totally bounded, and relatively compact if and only if it is totally bounded. The `only if' parts are trivial, and for the `if' parts we note that the closure of a totally bounded set is totally bounded; for if
$S$ can be covered by finitely many balls $B(x_n;\eps)$ with centres in $S$, the balls $B(x_n;2\eps)$ cover $\ov S$.
Now it remains to observe that a closed subset of a complete metric space is complete.

\begin{theorem}[Heine--Borel]\label{thm:Heine--Borel}\index{theorem!Heine--Borel}
A subset of $\R^d$ is compact if and only if it is
closed and bounded.
\end{theorem}
\begin{proof}
Compact subsets of a metric space are closed and bounded.
Conversely, let $S\subseteq\R^d$ be closed and bounded.
The completeness of $\R$ implies that of $\R^d$, since
Cauchy sequences in $\R^d$ converge coordinatewise.
Being closed in $\R^d$, $S$ is therefore complete.

To prove total boundedness, fix $\eps>0$ and choose $r>0$
such that $S\subseteq[-r,r]^d$.
Subdivide this cube into finitely many smaller cubes of side
length less than $\eps/\sqrt d$, and choose a point of $S$
in each small cube that intersects $S$.
Each small cube has diameter less than $\eps$, so the
$\eps$-balls centred at the chosen points cover $S$.
Thus $S$ is totally bounded, and Theorem
\ref{thm:totally-bounded} implies that $S$ is compact.
\end{proof}

\begin{theorem}\label{thm:compact-UC}
Let $(X,d_X)$ and $(Y,d_Y)$ be metric spaces with
$X$ compact. Every continuous mapping $f: X\to Y$ is uniformly continuous.
\end{theorem}

\begin{proof}
Let $\e>0$ be arbitrary. For every $x\in X$ we can find a  $\delta(x)>0$
such that for all $x'\in X$ with $d_X(x,x')<\delta(x)$ we have
$d_Y(f(x),f(x'))<\tfrac12\e$.
The collection
$$\mathscr{U} = \bigl\{B(x;\tfrac12\delta(x)):\,  x\in X\bigr\}$$
is an open cover of $X$, and therefore has a finite subcover, say
$$ \mathscr{U}' = \bigl\{B(x_j;\tfrac12\delta(x_j)):\,  j=1,\dots,n\bigr\}.$$
Let $\delta:=\min\{\frac12\delta(x_j): \, j=1,\dots,n\}$.

Suppose now that $x,x'\in X$ satisfy $d_X(x,x')<\delta$. We have
$x\in B(x_j;\tfrac12\delta(x_j))$ for some $1\le j\le n$
(since $\mathscr{U}'$ covers $X$).
Then $d_X(x,x_j) < \tfrac12\delta(x_j)$ and
$$d_X(x'\!,x_j) \le d_X(x'\!,x)+d_X(x,x_j) < \delta + \tfrac12\delta(x_j)\le \delta(x_j).$$
Consequently,
$ d_Y(f(x),f(x'))\le d_Y(f(x),f(x_j)) +  d_Y(f(x_j),f(x')) < \tfrac12\e +
\tfrac12\e = \e.$
\end{proof}

In many applications, the following special case of Tychonoff's theorem
suffices.

\begin{proposition}\label{prop:Tychonov-metric}
If $K_1,\dots,K_n$ are compact metric spaces, then their product $K:= K_1\times \cdots\times K_n$ is a compact metric space with respect to the {\em product metric}\index{product!metric}
$$ d(s,t):= \sum_{j=1}^n d_j(s_j,t_j), \quad s,t\in K.$$
\end{proposition}
\begin{proof} It is routine to verify that the above metric induces the product topology.  Given $\eps>0$, for each $j=1,\dots,n$ choose finitely many open $d_j$-balls of radius $\eps/n$ to cover $K_j$. Their products are open, contained in $d$-balls of radius $\eps$, and cover $K$.
Since $\eps>0$ was arbitrary, this shows that $K$ is totally bounded. Since the completeness of the spaces $K_j$ implies that $K$ is complete, this proves the result.
\end{proof}

\begin{definition}[Separability]\label{def:separable}\index{separable!metric space}
 A metric space is called {\em separable} if it contains a dense countable subset.
\end{definition}

\begin{proposition}
Every compact metric space is separable.
\end{proposition}
\begin{proof}
 For each $n=1,2,\dots$ we cover the metric space with finitely many open balls $B_1^{(n)}\!,\dots,B_{N_n}^{(n)}$ of radius $\frac1n$. We claim that the centres of all these balls form a dense subset. Indeed,
 if $U$ is a nonempty open set, then $U$ contains an open ball $B(x;r)$. Let $n\ge 1$ be such that $\frac1n<r$. For some $k\in \{1,\dots,N_n\}$ we have $x \in B_k^{(n)} =: B(c_k;\frac1n)$,
hence $d(x, c_k) < \frac1n < r$, and therefore $c_k\in B(x;r)\subseteq  U$.
\end{proof}

\section*{The Hausdorff--Alexandroff Theorem}

The aim of this section is to prove that every compact metric space is the continuous image of a compact Cantor-like subset of the unit interval $[0,1]$. This result will be needed in the proof of Theorem \ref{thm:Neumann-commute}.

\begin{theorem}[Hausdorff--Alexandroff]\label{thm:Hausdorff-Alexandroff}
\index{theorem!Hausdorff--Alexandroff}
Let $(X,d)$ be a nonempty compact metric space. Then there exist a compact set $L\subseteq[0,1]$
and a continuous surjection $\pi:L\to X$.
\end{theorem}

\begin{proof}
If $X$ consists of one point $x$, we take $L=\{0\}$ and $\pi(0):=x$. In the rest of the proof we may assume $X$ has at least two points.

\smallskip
{\em Step 1} --
For each $n\ge 1$ we will construct an integer $N_n\ge 1$ and a finite family of nonempty closed sets $F_{i_1,\dots ,i_n}\subseteq X$ with $i_k\in \{1,\dots, N_k\}$ for all $k=1,\dots,n$, such that for all $n\ge 1$ we have
\begin{enumerate}[label={\rm(\roman*)}, leftmargin=*]
\item\label{it:FS3} $X=\bigcup_{i_1=1}^{N_1}\cdots\bigcup_{i_n=1}^{N_n}F_{i_1,\dots, i_n}$;
\item\label{it:FS1} $\operatorname{diam}(F_{i_1,\dots, i_n})\le 2^{-n}$;
\item\label{it:FS2} $F_{i_1,\dots, i_{n}, i_{n+1}}\subseteq F_{i_1,\dots, i_n}$;
\item\label{it:FS4} $F_{i_1,\dots,i_n} = \bigcup_{i_{n+1}=1}^{N_{n+1}}
F_{i_1,\dots,i_n,i_{n+1}}$.
\end{enumerate}
For $n=1$, we may cover $X$ by $N_1\ge 2$ nonempty open sets of diameter $<\frac12$, and take their closures.

\begin{figure}[ht]
\centering
\begin{tikzpicture}[x=1cm,y=1cm,line cap=round,scale=0.75]

\def\dx{6.8}
\def\dy{-4.7}

\newcommand{\Xclip}{%
  plot[smooth cycle,tension=0.9] coordinates
  {(-0.4,0.7) (0.7,0.1) (2.3,0.2) (3.3,1.1) (3.1,2.6) (2.0,3.3) (0.6,3.1) (-0.6,2.2)}
}
\newcommand{\drawXdomain}{%
  \path[draw=black,fill=gray!10,line width=0.6pt] \Xclip;
}

\pgfmathsetmacro{\xL}{-0.8}
\pgfmathsetmacro{\xR}{3.6}
\pgfmathsetmacro{\yB}{0.0}
\pgfmathsetmacro{\yT}{3.6}
\pgfmathsetmacro{\xMid}{(\xL+\xR)/2}
\pgfmathsetmacro{\yMid}{(\yB+\yT)/2}
\pgfmathsetmacro{\Wsix}{(\xR-\xL)/6}
\pgfmathsetmacro{\Hthree}{(\yT-\yB)/3}
\pgfmathsetmacro{\xLMid}{(\xL+\xMid)/2}
\pgfmathsetmacro{\xRMid}{(\xMid+\xR)/2}

\tikzset{
  A/.style={fill=gray!25,opacity=0.75,draw=black,line width=0.5pt},
  B/.style={fill=gray!40,opacity=0.75,draw=black,line width=0.5pt},
  C/.style={fill=gray!55,opacity=0.75,draw=black,line width=0.5pt},
  D/.style={fill=gray!70,opacity=0.75,draw=black,line width=0.5pt}
}

\newcommand{\twoPieces}{%
  \begin{scope}
    \clip \Xclip;
    \path[A] (\xL,\yB) rectangle (\xMid,\yT);
    \path[C] (\xMid,\yB) rectangle (\xR,\yT);
    \draw[black,line width=0.6pt] (\xMid,\yB) -- (\xMid,\yT);
  \end{scope}
  \path[draw=black,fill=none,line width=0.6pt] \Xclip;
  \node at (0.35,2.5) {$F_1$};
  \node at (2.9,1.0) {$F_2$};
}

\newcommand{\sixPieces}{%
  \begin{scope}
    \clip \Xclip;

    \foreach \j/\sty in {0/A,1/B,2/C}{
      \pgfmathsetmacro{\ya}{\yB+\j*\Hthree}
      \pgfmathsetmacro{\yb}{\yB+(\j+1)*\Hthree}

      \path[\sty] (\xL,\ya) rectangle (\xMid,\yb);
      \path[\sty] (\xMid,\ya) rectangle (\xR,\yb);

      \draw[black,line width=0.45pt] (\xL,\ya) -- (\xMid,\ya);
      \draw[black,line width=0.45pt] (\xMid,\ya) -- (\xR,\ya);
    }

    \draw[black,line width=0.45pt] (\xL,\yT) -- (\xMid,\yT);
    \draw[black,line width=0.45pt] (\xMid,\yT) -- (\xR,\yT);

    \draw[black,line width=0.8pt] (\xMid,\yB) -- (\xMid,\yT);
  \end{scope}
  \path[draw=black,fill=none,line width=0.6pt] \Xclip;
}

\newcommand{\twelvePieces}{%
  \begin{scope}
    \clip \Xclip;

    \foreach \j in {0,1,2}{
      \pgfmathsetmacro{\ya}{\yB+\j*\Hthree}
      \pgfmathsetmacro{\yb}{\yB+(\j+1)*\Hthree}

      \path[B] (\xL,\ya) rectangle (\xLMid,\yb);
      \path[D] (\xLMid,\ya) rectangle (\xMid,\yb);

      \path[B] (\xMid,\ya) rectangle (\xRMid,\yb);
      \path[D] (\xRMid,\ya) rectangle (\xR,\yb);

      \draw[black,line width=0.45pt] (\xL,\ya) -- (\xMid,\ya);
      \draw[black,line width=0.45pt] (\xMid,\ya) -- (\xR,\ya);

      \draw[black,line width=0.45pt] (\xLMid,\ya) -- (\xLMid,\yb);
      \draw[black,line width=0.45pt] (\xRMid,\ya) -- (\xRMid,\yb);
    }

    \draw[black,line width=0.45pt] (\xL,\yT) -- (\xMid,\yT);
    \draw[black,line width=0.45pt] (\xMid,\yT) -- (\xR,\yT);

    \draw[black,line width=0.8pt] (\xMid,\yB) -- (\xMid,\yT);
  \end{scope}
  \path[draw=black,fill=none,line width=0.6pt] \Xclip;
}

\begin{scope}[shift={(0,0)}]
  \drawXdomain
  \node at (1.55,1.75) {$X$};
\end{scope}

\begin{scope}[shift={(\dx,0)}]
  \twoPieces
\end{scope}

\begin{scope}[shift={(0,\dy)}]
  \sixPieces
\end{scope}

\begin{scope}[shift={(\dx,\dy)}]
  \twelvePieces
\end{scope}

\end{tikzpicture}
\caption{Illustration of Step~1 for $N_1=2$, $N_2=3$, $N_3=2$.\label{fig:nested}}
\end{figure}

Suppose now that the families have been constructed up to level $n$.
For each multi-index $\alpha=(i_1,\dots,i_n)$, cover the compact metric
space $F_\alpha$, with its relative metric, by finitely many nonempty
relatively open sets of diameter $<2^{-(n+1)}$, and take their closures in
$F_\alpha$. These closures are nonempty closed subsets of $X$, are
contained in $F_\alpha$, have diameter $\le 2^{-(n+1)}$, and cover
$F_\alpha$.
Let these sets be $F_{\alpha,1},\dots,F_{\alpha,M_\alpha}$.
Since there are only finitely many multi-indices $\alpha$ at level $n$, we may put
$N_{n+1}:=\max_\alpha M_\alpha.$
Repeating sets if necessary, we may assume that for each $\alpha$ the cover
consists of exactly $N_{n+1}$ sets, denoted by
$F_{\alpha,1},\dots,F_{\alpha,N_{n+1}}$.
Then \ref{it:FS3}--\ref{it:FS4} hold at level $n+1$.
See Figure \ref{fig:nested}.

\smallskip
{\em Step 2} -- To each set $F_\alpha$ we associate a closed subinterval of $[0,1]$ as follows.

Set $I_\emptyset:=[0,1]$. Given an interval $I_\alpha=[a,b]\subseteq [0,1]$ at level $n-1$ (with $|\alpha|=n-1$), put $\delta:=(b-a)/(3N_n)$. For $i_n=1,\dots,N_n$ define
$$
I_{\alpha, i_n}:=\bigl[a+(3i_n-2)\delta,\ a+(3i_n-1)\delta\bigr].
$$
For each $n\ge 1$ the intervals $I_{\alpha, i_n}$ are pairwise disjoint closed subintervals of $I_\alpha$ of the same length $\delta_n$, where $\delta_1=1/(3N_1)$ and, inductively, $\delta_n=1/(3^n N_1\cdots N_n)$.
Note that any two distinct intervals $I_{\alpha, i_n}$ have distance at least
$\delta_n$.
See Figure \ref{fig:H-K}.

Define
$$
L:=\bigcap_{n\ge 1}\ \bigcup_{i_1=1}^{N_1}\cdots\bigcup_{i_n=1}^{N_n} I_{i_1,\dots, i_n}\subseteq[0,1].
$$
This set is a closed subset of $[0,1]$ and therefore compact. Moreover, $L\neq\emptyset$,
since the nested intervals $I_1\supseteq I_{1,1}\supseteq I_{1,1,1}\supseteq\cdots$ have nonempty intersection.

\begin{figure}[ht]
\centering
\begin{tikzpicture}[x=10cm,y=0.9cm,line cap=round]
  \def\None{2} 
  \def\Ntwo{3} 
  \def\Nthree{2} 

  \def\yzero{0}
  \def\yone{-1}
  \def\ytwo{-2}
  \def\ythree{-3}

  \newcommand{\baseline}[1]{\draw[gray!40,line width=0.4pt] (0,#1) -- (1,#1);}

  \baseline{\yzero}
  \baseline{\yone}
  \baseline{\ytwo}
  \baseline{\ythree}

  \node[left] at (0,\yzero) {$n=0$};
  \node[left] at (0,\yone)  {$n=1$};
  \node[left] at (0,\ytwo)  {$n=2$};
  \node[left] at (0,\ythree){$n=3$};

  \draw[line width=2.2pt] (0,\yzero) -- (1,\yzero);
  \node[above] at (0.5,\yzero) {$I_\emptyset=[0,1]$};

  \pgfmathsetmacro{\deltaone}{1/(3*\None)}            
  \pgfmathsetmacro{\deltatwo}{\deltaone/(3*\Ntwo)}    
  \pgfmathsetmacro{\deltathree}{\deltatwo/(3*\Nthree)}

  \foreach \i in {1,...,\None}{
    \pgfmathsetmacro{\a}{(3*\i-2)*\deltaone}
    \pgfmathsetmacro{\b}{(3*\i-1)*\deltaone}
    \draw[line width=2.2pt] (\a,\yone) -- (\b,\yone);
    \pgfmathsetmacro{\m}{(\a+\b)/2}
    \node[above] at (\m,\yone) {$I_{\i}$};
  }
  \node[above right] at (0,\yone) {\scriptsize $N_1=\None$};

  \foreach \i in {1,...,\None}{
    \pgfmathsetmacro{\leftone}{(3*\i-2)*\deltaone}
    \foreach \j in {1,...,\Ntwo}{
      \pgfmathsetmacro{\a}{\leftone + (3*\j-2)*\deltatwo}
      \pgfmathsetmacro{\b}{\leftone + (3*\j-1)*\deltatwo}
      \draw[line width=2.2pt] (\a,\ytwo) -- (\b,\ytwo);
    }
  }
  \node[above right] at (0,\ytwo) {\scriptsize $N_2=\Ntwo$};

  \foreach \i in {1,...,\None}{
    \pgfmathsetmacro{\leftone}{(3*\i-2)*\deltaone}
    \foreach \j in {1,...,\Ntwo}{
      \pgfmathsetmacro{\lefttwo}{\leftone + (3*\j-2)*\deltatwo}
      \foreach \k in {1,...,\Nthree}{
        \pgfmathsetmacro{\a}{\lefttwo + (3*\k-2)*\deltathree}
        \pgfmathsetmacro{\b}{\lefttwo + (3*\k-1)*\deltathree}
        \draw[line width=2.2pt] (\a,\ythree) -- (\b,\ythree);
      }
    }
  }
  \node[above right] at (0,\ythree) {\scriptsize $N_3=\Nthree$};

\end{tikzpicture}
\caption{First three steps in the construction of $L$ (here $N_1=2$, $N_2=3$, $N_3=2$).\label{fig:H-K}}
\end{figure}

\smallskip
{\em Step 3} --
For every $t\in L$ and each $n\ge1$ there is a unique multi-index $(i_1,\dots,i_n)$ with
$t\in I_{i_1,\dots, i_n}$. Denote it by $\alpha_n(t)$.
Set $K_n(t):=F_{\alpha_n(t)}$.
By \ref{it:FS2}, these sets are decreasing, and each $K_n(t)$ is nonempty and compact. Since $X$ is compact, it follows that $\bigcap_{n\ge 1} K_n(t)\neq\emptyset$.  On the other hand, by \ref{it:FS1}, if $x,y\in \bigcap_{n\ge 1}K_n(t)$, then $$ d(x,y)\le \operatorname{diam}(K_n(t))\le 2^{-n}\quad n\ge 1, $$ and therefore $d(x,y)=0$. Hence the intersection consists of exactly one point; define $\pi(t)$ to be this point.

To see that $\pi$ is continuous, fix $t\in L$ and $\varepsilon>0$. Choose $n\ge 1$ such that
$2^{-n}<\varepsilon$. If $s,t\in L$ satisfy $|s-t|< \delta_n$, the length of the intervals $I_{i_1,\dots,i_n}$, then $\alpha_n(s) = \alpha_n(t)$. It follows that
$\pi(s),\pi(t)\in F_{\alpha_n(t)}$, and therefore
$$d(\pi(s),\pi(t))\le \operatorname{diam}(F_{\alpha_n(t)})\le 2^{-n}<\varepsilon.$$

Finally, $\pi$ is surjective. Given $x\in X$, use \ref{it:FS4} inductively to choose indices
$i_1,i_2,\dots$ with $x\in F_{i_1,\dots, i_n}$ for all $n$. Then the nested intervals
$I_{i_1,\dots, i_n}$ have lengths tending to $0$, and hence their intersection consists of a single
point $t\in[0,1]$. By construction $t\in L$ and $x\in\bigcap_{n\ge1} F_{i_1,\dots, i_n}$, so
$\pi(t)=x$.
\end{proof}

Theorem \ref{thm:Hausdorff-Alexandroff} will be used together with the
following simple result.

\begin{proposition}[Lower semicontinuous right inverses]
\label{prop:min-section}
\index{right inverse}\index{compact!subset of $\R$}
Let $(X,d)$ be a compact metric space, let $L\subseteq [0,1]$ be compact, and let $\pi:L\to X$ be a
continuous surjection. Define $\sigma:X\to L$ by
$$
\sigma(x):=\min\{t\in L:\,  \pi(t)=x\},\quad x\in X.
$$
Then $\pi\circ\sigma=\mathrm{id}_X$, and $\sigma$ is lower semicontinuous.
\end{proposition}

\begin{proof}
For each $x\in X$, the set $\pi^{-1}(\{x\})$ is a nonempty
compact subset of $\R$, so it has a minimum.
Thus $\sigma$ is well defined and $\pi(\sigma(x))=x$.

For every $a\in\R$ we have
$$
\{x\in X:\, \sigma(x)\le a\}
=\pi\bigl(L\cap(-\infty,a]\bigr).
$$
The set on the right-hand side is compact, hence closed in $X$.
Therefore $\sigma$ is lower semicontinuous.
\end{proof}

\section*{The Brouwer Fixed Point Theorem}

For $n\ge1$ let $\Delta_n$ be the $n$-dimensional simplex in $\mathbb R^{n+1}$ with vertices $e_0,\dots,e_n$.
\[
\Delta_n
:=
\Bigl\{(x_0,\dots,x_n)\in[0,1]^{n+1}:
\sum_{j=0}^n x_j=1\Bigr\}.
\]
A finite triangulation of $\Delta_n$ is {\em Sperner-labelled} if each
of its vertices $v$ is assigned a label
$\ell(v)\in\{0,\dots,n\}$ in such a way that $v_j=0$ implies $\ell(v)\ne j$.

\begin{lemma}[Sperner]\label{lem:Sperner}\index{lemma!Sperner}
Every Sperner-labelled triangulation of $\Delta_n$ contains a simplex
whose vertices carry all the labels $0,\dots,n$.
\end{lemma}

\begin{proof}
We prove the stronger assertion that the number
of such simplices is odd. The case $n=1$ is immediate. For the induction step,
assume that the number of $(n-1)$-simplices in $F_n:= \{x\in\Delta_n:\,  x_n=0\}$ whose vertices
carry the labels $0,\dots,n-1$ is odd.

Count the pairs $(S,F)$ in which $S$ is an $n$-simplex of the
triangulation and $F$ is a face of $S$ whose vertices carry precisely
the labels $0,\dots,n-1$. If such a face lies in the interior, it belongs to
two $n$-simplices. If such a face lies in the boundary of $\Delta_n$,
it must lie in $F_n$, since on the face $x_j=0$ the label $j$ is
forbidden. Thus the number of pairs $(S,F)$ is odd.

An $n$-simplex contributes one pair if its vertices carry all labels
$0,\dots,n$, two pairs if they carry the labels $0,\dots,n-1$ with
one label repeated, and no pairs otherwise. It follows that the number
of simplices carrying all labels $0,\dots,n$ is odd.
\end{proof}

\begin{theorem}[Brouwer fixed point theorem]\label{thm:Brouwer}\index{theorem!Brouwer fixed point}
Every continuous mapping $F:\Delta_n\to\Delta_n$ has a fixed point.
\end{theorem}

\begin{proof}
Choose finite triangulations of $\Delta_n$ whose meshes tend to zero;
for instance, repeated barycentric subdivision has this property.

Consider one such triangulation. If one of its vertices is a fixed
point of $F$, there is nothing to prove. Otherwise, for each vertex
$v$ choose an index $\ell(v)$ such that $v_{\ell(v)}>F_{\ell(v)}(v).$
Such an index exists since the coordinates of both $v$ and $F(v)$
sum to $1$. Moreover, $v_j=0$ makes $v_j>F_j(v)$ impossible, so this
defines a Sperner labelling.

By Lemma \ref{lem:Sperner}, there is therefore a simplex with vertices
$v^{(0)},\dots,v^{(n)}$ such that
$v^{(j)}_j>F_j(v^{(j)})$ for all $j=0,\dots,n.$
Doing this for triangulations whose meshes tend to zero and using the compactness of
$\Delta_n$, we can pass to a subsequence such that all the corresponding vertices $v^{(j)}$
converge to a point $x\in\Delta_n$. By continuity we have $x_j\ge F_j(x)$ for all $j=0,\dots,n$.
Since the coordinates on both sides sum to $1$, equality holds for
every $j$. Hence $F(x)=x$.
\end{proof}

\cleardoublepage  

\cleardoublepage  

\chapter{Measure Spaces}

\blfootnote{\rm \small \begin{quote}This book has been published by Cambridge University Press in the series `Cambridge Studies in Advanced Mathematics'. The present version is free to view and download for personal use only. Not for re-distribution, re-sale or use in derivative works.\newline {\copyright} Jan van Neerven\end{quote}}

\noindent
This appendix reviews the basic elements of measure theory.

\section*{\texorpdfstring{$\sigma$}{sigma}-Algebras}

Let $\Omega$ be a set.

\begin{definition}[$\sigma$-Algebras]
A {\em $\sigma$-algebra}\index{sAlgebra@$\sigma$-algebra} in $\Omega$
is a collection $\F$ of subsets of $\Omega$ with the following
properties:

\begin{enumerate}[label={\rm(\roman*)}, leftmargin=*]
\item $\Omega\in \F$;
\item $F\in{\calF}$ implies $\complement F\in\F$;
\item $F_1, F_2,\hdots \in \F$ implies $\bigcup_{n\ge 1} F_n\in\F$\!.
\end{enumerate}
Here, $\complement F = \Omega\setminus F$ is the complement of $F$.

A {\em measurable space}\index{measurable space} is a pair $(\Om,\F)$, where $\Om$ is a set and
$\F$ is a $\sigma$-algebra in $\Om$. The sets in $\calF$ are often referred to as the {\em measurable subsets} of $\Om$.
\end{definition}

These properties express that $\F$ is nonempty, closed under taking
complements, and closed under taking countable unions.
From
$$\bigcap_{n\ge 1} F_n = \complement \Big( \bigcup_{n\ge 1} \complement F_n\Big)$$
it follows that $\F$ is closed under taking countable
intersections.  Clearly, $\F$ is closed under finite unions and intersections
as well.

For any set $\Omega$, its power set $2^\Omega$\index{power set}\index{$2^S$} (the set of all subsets of $\Omega$) is a $\sigma$-algebra. In general one is interested in smaller $\sigma$-algebras:

\begin{example}
When $\calC$ is any collection of subsets of $\Om$, the
{\em $\sigma$-algebra generated by $\calC$}
is defined as the intersection of all
$\sigma$-algebras in $\Omega$ containing $\calC$, and is denoted by $\sigma(\calC)$.\index{$S$@$\sigma(\calC)$}
It is the smallest $\sigma$-algebra containing $\calC$. Such $\sigma$-algebras arise in a variety of situations:

\begin{enumerate}[label={\rm(\arabic*)}, leftmargin=*]
\item When $(X,\tau)$ is a topological space, the $\sigma$-algebra $\calB(X)$\index{$Bz$@$\calB(X)$}
generated by $\tau$ is called the {\em Borel $\sigma$-algebra}\index{Borel!$\sigma$-algebra} of $(X,\tau)$.

\item Let $\Om$ and $\Om'$ be sets and let $f: \Omega \to \Om'$ be any function.
When $\calF'$ is a $\sigma$-algebra in $\Om'$, for $F'\in\calF'$ we define\index{$FinF$@$\{f\in F\}$}
$$\{f\in F'\} := \{\omega\in\Omega:\,  f(\omega)\in F'\}.$$
The collection
$$ \sigma(f) = \big\{\{f\in F'\}:\,  F'\in\calF'\big\}$$
is a $\sigma$-algebra in $\Om$, the {\em $\sigma$-algebra
generated by $f$}.
The $\sigma$-algebra generated by a family of functions is defined similarly.

\item Let $(\Om_1,\F_1)$, $(\Om_2,\F_2), \, \dots$ be a sequence of measurable spaces.
On the product $\Om:= \prod_{n\ge 1} \Om_n$, the  {\em product $\sigma$-algebra}\index{product!$\sigma$-algebra} $ \prod_{n\ge 1}\F_n$ is the $\sigma$-algebra generated by all sets of the form
$$ F_1\times\cdots\times F_N\times\Om_{N+1}\times\Om_{N+2}\times\cdots$$
with $N = 1,2,\dots$ and $F_n\in \F_n$ for $n=1,\dots,N$.

The product of finitely many measurable spaces $(\Om_1,\F_1),\,\dots\,(\Om_N,\F_N)$
is defined similarly; here one takes the $\sigma$-algebra in
$\Om := \Om_1\times\cdots\times\Om_N$
generated by all sets of the form
$ F_1\times\cdots\times F_N$ with $F_n\in \F_n$ for $n=1,\dots,N$.
This is precisely the $\sigma$-algebra generated by the coordinate projections $p_n: \Om\to \Om_n$, $n=1,\dots,N$.
For example, the reader may check that
$$ (\R^d\!,\calB(\R^d)) = \prod_{n=1}^d (\R,\calB(\R)).$$
For a proof, one may use the fact that every open set in $\R^d$ is a countable union of open rectangles
of the form $(a_1,b_1)\times\cdots \times (a_d,b_d)$.
\end{enumerate}
\end{example}

\section*{Measures}\label{sec:app-measures}

Let $(\Om,\F)$ be a measurable space.

\begin{definition}[Measures]\label{def:measure} A {\em measure}\index{measure} on $(\Om,\F)$ is a mapping
$ \mu: \F \to [0,\infty]$ with the following properties:

\begin{enumerate}[label={\rm(\roman*)}, leftmargin=*]
\item $\mu(\emptyset) = 0$;
\item for all disjoint sets $F_1, F_2,\dots$ in $\F$ we have
      $\mu(\bigcup_{n\ge 1} F_n) = \sum_{n\ge 1} \mu(F_n).$
\end{enumerate}

A triple $(\Om,\F\!,\mu)$, with $\mu$ a measure on a measurable space $(\Om,\F)$, is
called a {\em measure space}.\index{measure space}
\end{definition}

A measure space $(\Om,\F\!,\mu)$ is called {\em finite} if $\mu$ is a {\em finite}\index{finite!measure}
measure, that is, if $\mu(\Om)<\infty$. If $\mu(\Om)=1$, then $\mu$ is called
a {\em probability measure}\index{probability!measure} and $(\Om,\F\!,\mu)$ is  called a
{\em probability space}.\index{probability!space} In probability theory, it is customary to use the
symbol $\P$ for a probability measure. A measure space $(\Om,\F\!,\mu)$ is called {\em $\sigma$-finite}\index{S@$\sigma$-finite measure}
if there exist $F_1,F_2,\dots$ in $ \F$
such that $\bigcup_{n\ge 1} F_n = \Omega$ and $\mu(F_n)<\infty$ for all $n\ge 1$.
A {\em Borel measure} on a topological space $(X,\tau)$ is a measure $\mu:\calB(X)\to[0,\infty]$, where $\mathscr{B}(X)$ is the Borel $\sigma$-algebra of $X$.
\index{Borel!measure}

The following properties of measures are easily checked:
\begin{enumerate}[label={\rm(\roman*)}, leftmargin=*]
\item\label{it:property-measures1} if
$F_1\subseteq F_2$ in $\F$\!, then $\mu(F_1)\le \mu(F_2)$;
\item\label{it:property-measures2} if $ F_1, F_2, \dots$ in $\F$\!,
then $$\mu\Big(\bigcup_{n\ge 1} F_n\Big) \le \sumn \mu(F_n);$$
\item\label{it:property-measures3} if $ F_1\subseteq F_2 \subseteq \dots$ in $\F$\!,
then $$\mu\Big(\bigcup_{n\ge 1} F_n\Big) = \limn \mu(F_n);$$
\item\label{it:property-measures4} if $ F_1\supseteq F_2 \supseteq \dots$ in $\F$ and
$\mu(F_1)<\infty$,
then $$\mu\Big(\bigcap_{n\ge 1} F_n\Big) = \limn \mu(F_n).$$
\end{enumerate}
In \ref{it:property-measures3} and \ref{it:property-measures4}, the limits (in $[0,\infty]$) exist by monotonicity.

\section*{Dynkin's Lemma}

Let $(\Om,\F)$ be a measurable space. A subcollection $\calD\subseteq\F$ is a {\em Dynkin system}\index{Dynkin
system} if
\begin{itemize}
\item $\Omega\in\calD$;
\item if $D_1\subseteq D_2$ with $D_1,D_2\in \calD$, then $D_2\setminus D_1\in \calD$;
\item if $D_1\subseteq D_2\subseteq \dots $ with all $D_n\in \calD$, then
also $\bigcup_{n\ge 1} D_n\in \calD$.
\end{itemize}

\begin{lemma}[$\pi$--$\lambda$ theorem]\label{lem:pi-lambda}\index{theorem!$\pi-\lambda$}
Let $(\Omega,\F)$ be a measurable space and let $\mathscr A\subseteq\F$ be a subcollection closed under finite intersections and containing $\Om$.
Let $\mathscr D\subseteq\F$ be a Dynkin system containing $\mathscr A$.
Then $\sigma(\mathscr A)\subseteq \mathscr D$.
\end{lemma}

\begin{proof}
Let $\mathscr D_0$ denote the smallest Dynkin system in $\F$ containing $\mathscr A$.
It suffices to show that $\mathscr D_0$ is a $\sigma$-algebra, for then
$\sigma(\mathscr A)\subseteq \mathscr D_0\subseteq \mathscr D$.

Let $\calC = \{D_0\in\calD_0:\,  D_0\cap A\in \calD_0$ for all $A\in\calA\}$. This is a Dynkin system and
$\calA\subseteq \calC$ since $\calA$ is closed under taking finite intersections.
It follows that $\calD_0\subseteq \calC$, since $\calD_0$ is the smallest Dynkin system
containing $\calA$. But $\calC\subseteq \calD_0$, and therefore
$\calC = \calD_0$.

Now let $\calC' = \{D_0\in\calD_0:\,  D_0\cap D \in \calD_0$ for all $D\in\calD_0\}$.
This is a Dynkin system and the
fact that $\calC = \calD_0$ implies that $\calA\subseteq \calC'$\!.
Hence $\calD_0\subseteq \calC'$\!, since $\calD_0$ is the smallest Dynkin system
containing $\calA$.  But $\calC'\subseteq \calD_0$, and therefore
$\calC' = \calD_0$.

It follows that $\calD_0$ is closed under taking finite intersections.
We claim that $\calD_0$ is a $\sigma$-algebra.
Indeed, $\mathscr D_0$ is closed under taking complements. It is also closed under taking finite unions, for if $D$ and $E$ belong to $\mathscr D_0$, then so does
$$
D\cup E = \complement\bigl((\Omega\setminus D)\cap(\Omega\setminus E)\bigr).
$$
Now let $(D_n)_{n\ge1}\subseteq\mathscr D_0$ be any sequence and set
$S_n:=\bigcup_{k=1}^n D_k$. Note that $S_n\in\mathscr D_0$, $S_1\subseteq S_2\subseteq \dots$, and $\bigcup_{n\ge 1} S_n = \bigcup_{n\ge1}D_n$, and therefore $\bigcup_{n\ge1}D_n\in\mathscr D_0$. This proves the claim.

Now $\calA\subseteq\calD_0$ implies that also $\sigma(\calA)\subseteq \calD_0$.
\end{proof}

\begin{lemma}[Dynkin lemma]\label{lem:unique-app}\index{lemma!Dynkin}
Let $\mu_1$ and $\mu_2$ be finite measures on a measurable space $(\Omega,\F)$.

Let ${\mathscr{A}}\subseteq {\mathscr{F}}$
be a collection of sets with the following properties:
\begin{enumerate}[label={\rm(\roman*)}, leftmargin=*]
\item\label{it:Dynkini} $\Omega\in \mathscr{A}$;
\item\label{it:Dynkinii} ${\mathscr{A}}$ is closed under finite intersections;
\item\label{it:Dynkiniii} the $\sigma$-algebra generated by ${\mathscr{A}}$ equals  ${\calF}$.
\end{enumerate}
If $\mu_1(A)=\mu_2(A)$ for all $A\in\mathscr A$, then $\mu_1=\mu_2$ on $\F$.
\end{lemma}

\begin{proof}
Let
$$
\mathscr D:=\{D\in\F:\,  \mu_1(D)=\mu_2(D)\}.
$$
Then $\Omega\in\mathscr D$.
If $D_1\subseteq D_2$ with $D_1,D_2\in\mathscr D$, then
$\mu_j(D_2\setminus D_1)=\mu_j(D_2)-\mu_j(D_1)$ for $j=1,2$,
so $\mu_1(D_2\setminus D_1)=\mu_2(D_2\setminus D_1)$ and therefore $D_2\setminus D_1\in\mathscr D$.
If $D_1\subseteq D_2 \subseteq \dots $ and $D:= \bigcup_{n\ge 1}D_n$,  with all $D_n\in\mathscr D$, then
$\mu_j(D)=\lim_{n\to\infty}\mu_j(D_n)$ for $j=1,2$,
and therefore $D\in\mathscr D$. Thus $\mathscr D$ is a Dynkin system, and it contains $\mathscr A$ by assumption.
Lemma~\ref{lem:pi-lambda} therefore gives
$\F=\sigma(\mathscr A)\subseteq\mathscr D$, that is, $\mu_1=\mu_2$ on $\F$.
\end{proof}

\section*{Outer Measures}\label{sec:outer}

Let $S$ be a set and let $2^S$ be its power set, that is, the set of all subsets of $S$.

\begin{definition}[Outer measures]\label{def:outermeasure}
A mapping $\nu:2^S\to [0,\infty]$ is called an \emph{outer measure}\index{outer measure} if
\begin{enumerate}[label={\rm(\roman*)}, leftmargin=*]
\item\label{it:outermeasure1} $\nu(\emptyset) = 0$;
\item\label{it:outermeasure2} $A\subseteq B$ implies $\nu(A)\leq \nu(B)$;
\item\label{it:outermeasure3} for all $A_1,A_2,\hdots \in 2^S$ we have $$\nu\Big(\bigcup_{n\ge 1} A_n\Big) \leq \sum_{n\ge 1} \nu(A_n).$$
\end{enumerate}
\end{definition}

\begin{lemma}\label{lem:addoutermeas}
Let $\mathscr{C}\subseteq 2^S$ satisfy $\emptyset\in \mathscr{C}$ and suppose that $\mu:\mathscr{C}\to [0,\infty]$ satisfies $\mu(\emptyset) = 0$. For subsets $A\subseteq S$ define
\begin{equation}\label{eq:outermeasdef}
\mu^*(A) := \inf\Big\{\sum_{j\ge 1} \mu(C_j): \, A\subseteq \bigcup_{j\ge 1} C_j, \ \text{where} \ C_j\in \mathscr{C} \ \text{for all} \ j\geq 1\Big\}
\end{equation}
with the convention that $\mu^*(A) = \infty$ if the above set is empty. Then $\mu^*$ is an outer measure.
\end{lemma}
\begin{proof}
The mapping $\mu^*:2^S\to [0,\infty]$ clearly satisfies the conditions \ref{it:outermeasure1} and \ref{it:outermeasure2}
in Definition \ref{def:outermeasure}. To check condition \ref{it:outermeasure3}, let $A_1,A_2,\hdots$ be subsets of $S$ and let $\varepsilon>0$ be arbitrary. If $\mu^*(A_n) = \infty$ for some $n\geq1$, then \ref{it:outermeasure3}
trivially holds. We may therefore assume that $\mu^*(A_n) < \infty$ for all $n\geq1$. By the definition of $\mu^*$, for each fixed $n\geq 1$ we can find $C_{n,j}\in \mathscr{C}$ such that
$$A_n\subseteq \bigcup_{j\ge 1} C_{n,j}  \ \hbox{ and } \  \dps\sum_{j\ge 1} \mu(C_{n,j})\le  \mu^*(A_n) + 2^{-n}\varepsilon .$$
Then $\bigcup_{n\ge 1}  A_n\subseteq \bigcup_{n,j=1}^\infty C_{n,j}$, and, again by the definition of $\mu^*$,
$$\mu^*\Big(\bigcup_{n\ge 1} A_n\Big)
\leq \sum_{n,j=1}^\infty \mu(C_{n,j})
\leq \sum_{n\ge 1} (\mu^*(A_n) + 2^{-n}\varepsilon) = \varepsilon + \sum_{n\ge 1} \mu^*(A_n).
$$
Since $\varepsilon>0$ was arbitrary, this proves the required estimate.
\end{proof}

Let $\mu:2^S\to [0,\infty]$ be a mapping which satisfies $\mu(\emptyset) = 0$. A set $A\subseteq S$ is called {\em $\mu$-measurable} if
$$\mu(Q)  = \mu(Q\cap A) + \mu(Q\cap \complement A),\quad Q\in 2^S.$$
The collection of all $\mu$-measurable sets is denoted by $\M_{\mu}$.

\begin{theorem}[Measures from outer measures]\label{thm:outermeasismeas}
If $\nu:2^S\to [0,\infty]$ is an outer measure, then $\M_{\nu}$ is a $\sigma$-algebra and $\nu$ is a measure on $(S,\M_{\nu})$.
\end{theorem}

For the proof of the theorem, we need the following terminology.
A {\em ring}\index{ring} in $S$ is a subset $\mathscr{R}$ of $2^S$ with the following properties:
\begin{enumerate}[label={\rm(\roman*)}, leftmargin=*]
\item $\emptyset \in \mathscr{R}$;
\item $A,B\in \mathscr{R}$ implies $A\setminus B\in \mathscr{R}$;
\item\label{it:ring3} $A,B\in \mathscr{R}$ implies $A\cup B \in \mathscr{R}$.
\end{enumerate}
If $\mathscr{R}$ is a ring, the identity $A\cap B = A\setminus (A\setminus B)$ implies that if $A,B\in \mathscr{R}$, then $A\cap B \in \mathscr{R}$.

\begin{proof}[Proof of Theorem \ref{thm:outermeasismeas}]
We proceed in two steps.

\smallskip
{\em Step 1} --
We begin by checking that if $\mu:2^S\to [0,\infty]$ is any mapping which satisfies $\mu(\emptyset) = 0$, then $\M_{\mu}$ is a ring and $\mu$ is additive on $\M_{\mu}$.

It is clear that $\emptyset\in \M_{\mu}$.
To check that $\M_{\mu}$ is a ring, we verify the following:
\begin{enumerate}[\rm(a), leftmargin=*]
\item\label{it:prop-ring-alt1} $A\in \M_{\mu}$ implies $\complement A\in \M_{\mu}$;
\item\label{it:prop-ring-alt2} $A,B\in \M_{\mu}$ implies $A\cap B\in \M_{\mu}$.
\end{enumerate}
Given these properties it is straightforward to check that $\M_{\mu}$ is a ring. Indeed, this follows from the formulas $B\setminus A = B \cap \complement A$ and $A\cup B = \complement(\complement A \cap \complement B)$.

Property \ref{it:prop-ring-alt1} is clear. To check \ref{it:prop-ring-alt2}, let $A,B\in \M_{\mu}$ and set $C := A\cap B$. Let $Q\in 2^S$ be arbitrary. Observing that
$A \cap \complement B = \complement C\cap A$ and $\complement A = \complement C\cap \complement A$,
and making repeated use of the definition of $\M_{\mu}$, we have
\begin{align*}
\mu(Q) & = \mu(Q\cap A) +\mu(Q\cap \complement A)
\\ & = \mu(Q\cap A\cap B) + \mu(Q\cap A\cap \complement B) +\mu(Q\cap \complement A) &
\\ & = \mu(Q\cap C) + \mu(Q\cap \complement C\cap A) +\mu(Q\cap \complement C\cap \complement A)
\\ & = \mu(Q\cap C) + \mu(Q\cap \complement C).
\end{align*}
Therefore, $A \cap B = C\in \M_{\mu}$.

For later use, we prove a slightly stronger identity.
Let $A,B\in\M_\mu$ be disjoint and let $Q\subseteq S$ be
arbitrary. Applying the measurability condition for $A$
to the test set $Q\cap(A\cup B)$ gives
\begin{align*}
\mu(Q\cap(A\cup B))
&=\mu(Q\cap(A\cup B)\cap A)
  +\mu(Q\cap(A\cup B)\cap\complement A)\\
&=\mu(Q\cap A)+\mu(Q\cap B).
\end{align*}
Taking $Q=S$ shows that $\mu$ is additive on $\M_\mu$.

\smallskip
{\em Step 2} -- We now turn to the proof of the theorem.
From Step 1
we know that $\M_{\nu}$ is a ring and $\nu$ is additive on $\M_{\nu}$. In view of property \ref{it:prop-ring-alt1} it remains to check that for any disjoint sequence $(A_n)_{n\geq 1}$ in $\M_{\nu}$,
\begin{equation}\label{eq:sigmaalgcount}
A:=\bigcup_{n\ge 1} A_n\in \M_{\nu} \ \hbox{ and }  \  \displaystyle \nu(A) = \sum_{n\ge 1} \nu(A_n).
\end{equation}
Let $B_n = \bigcup_{j=1}^n A_j$ for each $n\geq 1$. Fix an arbitrary subset $Q$ of $S$.
By Step 1, for all $n\ge 1$ we have
$\complement A\subseteq \complement B_n$, $B_n\in \M_{\nu}$, and
\begin{align*}
\sum_{j=1}^n  \nu(Q\cap A_j) + \nu(Q\cap \complement A) & = \nu(Q\cap B_n) + \nu(Q\cap \complement A)
\\ & \leq \nu(Q\cap B_n) + \nu(Q\cap \complement B_n)
 = \nu(Q).
\end{align*}
Using the countable subadditivity of $\nu$ and then passing to the limit $n\to\infty$, we infer
\begin{equation}\label{eq:inequaequa}
\nu(Q\cap A) + \nu(Q\cap \complement A) \leq \sum_{j\ge 1} \nu(Q\cap A_j) + \nu(Q\cap \complement A)\leq \nu(Q).
\end{equation}
On the other hand, by subadditivity also the converse inequality $\nu(Q)\leq \nu(Q\cap A) + \nu(Q\cap \complement A)$ holds. This shows that $A\in \M_{\nu}$ and that the inequalities in \eqref{eq:inequaequa} are in fact equalities. Now \eqref{eq:sigmaalgcount} follows by taking $Q = A$ in \eqref{eq:inequaequa}.
\end{proof}

\section*{Carath\'eodory's Extension Theorem}

The following result holds for additive functions
$\mu:\mathscr R\to[0,\infty]$.

\begin{lemma}\label{lem:ringestimates}
Let $\mathscr{R}$ be a ring and $\mu:\mathscr{R}\to [0,\infty]$ be additive, that is,
$$\mu\Bigl(\bigcup_{j=1}^n A_j\Bigr) = \sum_{j=1}^n \mu(A_j)$$ holds for all disjoint sets $A_1,\dots,A_n\in \mathscr{R}$. The following assertions hold:
\begin{enumerate}[label={\rm(\arabic*)}, leftmargin=*]
\item\label{it:ringestimates1} if $A,B\in \mathscr{R}$ and $A\subseteq B$, then $\mu(A)\leq \mu(B)$;
\item\label{it:ringestimates2} if $A_1, A_2, \dots\in \mathscr{R}$ and $\bigcup_{j\ge 1} A_j \in \mathscr{R}$ and $\mu$ is countably additive on $\mathscr{R}$, then
 $$\mu\Big(\bigcup_{j\ge 1} A_j\Big) \leq \sum_{j\ge 1} \mu(A_j).$$
\end{enumerate}
\end{lemma}
\begin{proof} \ref{it:ringestimates1}: \  Writing $B = A\cup (B\setminus A)$, we see that
$$\mu(B) = \mu(A\cup (B\setminus A)) = \mu(A) + \mu(B\setminus A)\geq \mu(A).$$

\smallskip
\ref{it:ringestimates2}: \ The sets
$B_1 := A_1$, $B_2 := A_2\setminus A_1$, $B_3 := A_3\setminus (A_1\cup A_2),\ \hdots$
are disjoint and we have
$\bigcup_{j\ge 1} B_j  = \bigcup_{j\ge 1} A_j$. Therefore, by the countable additivity of $\mu$,
\[
\mu\Big(\bigcup_{j\ge 1} A_j\Big) = \mu\Big(\bigcup_{j\ge 1} B_j\Big) = \sum_{j\ge 1} \mu(B_j)\leq \sum_{j\ge 1} \mu(A_j).
\qedhere\]
\end{proof}

\begin{theorem}[Carath\'eodory's extension theorem]\label{thm:Caratheodory}
Let $\mathscr{R}$ be a ring in $S$ and suppose that $\mu:\mathscr{R}\to [0,\infty]$ is countably additive on $\mathscr{R}$ and satisfies $\mu(\emptyset) = 0$. Let $\mu^*$ be the associated outer measure. Then:
\begin{enumerate}[label={\rm(\arabic*)}, leftmargin=*]
 \item \label{it:Caratheodory1}
the outer measure $\mu^*$ restricts to a measure on $\sigma(\mathscr{R})$ extending $\mu$;
 \item \label{it:Caratheodory2}
 if $\mu$ is $\sigma$-finite on $\mathscr{R}$, and if $\mu_1$ and $\mu_2$ are $\sigma$-finite measures on $\sigma(\mathscr{R})$ that extend $\mu$, then $\mu_1=\mu_2$ on $\sigma(\mathscr{R})$.
\end{enumerate}
\end{theorem}

Explicitly, in \ref{it:Caratheodory2} we assume that there exist $R_1,R_2,\hdots\in\mathscr{R}$ such that
$S=\bigcup_{n\ge1}R_n$ and $\mu(R_n)<\infty$ for each $n\ge 1$. Replacing each $R_n$ by $\bigcup_{j=1}^n R_j$, we may assume that $R_1\subseteq R_2 \subseteq\hdots$ and $\bigcup_{n\ge 1} R_n = S$.

\begin{proof} By Theorem \ref{thm:outermeasismeas}, $\mu^*$
is a measure on the $\sigma$-algebra $\M_{\mu^*}$. We claim that it has the following properties:
\begin{enumerate}[label={\rm(\roman*)}, leftmargin=*]
\item\label{it:property-outer3} $\mathscr{R}\subseteq \M_{\mu^*}$;
\item\label{it:property-outer2} $\mu^*(A) = \mu(A)$ for all $A\in \mathscr{R}$.
\end{enumerate}
Clearly, part \ref{it:Caratheodory1} of the theorem follows from the claim, which actually shows that there is a further extension to the possibly larger $\sigma$-algebra $\M_{\mu^*}$.

\smallskip
{\em Step 1} -- In this step we prove \ref{it:property-outer3}. Let $A\in \mathscr{R}$ and $Q\subseteq S$ be given. The subadditivity of $\mu^*$ gives $ \mu^*(Q) \leq {\mu^*(Q\cap A) + \mu^*(Q\cap \complement A)}$.
The converse estimate  $\mu^*(Q\cap A) + \mu^*(Q\cap \complement A)\le \mu^*(Q)$ trivially holds if $\mu^*(Q)=\infty$. If $\mu^*(Q)<\infty$, choose $B_1, B_2, \dots \in \mathscr{R}$ such that $Q\subseteq \bigcup_{n\ge 1} B_n$. Then
$B_n\cap A$ and $B_n\cap \complement A = B_n\setminus A$ belong to $\mathscr{R}$ for all $n\geq 1$, and
$$Q\cap A\subseteq \bigcup_{n\ge 1} B_n\cap A \ \text{ and } \ Q\cap \complement A\subseteq \bigcup_{n\ge 1} B_n\cap \complement A.$$
Using first the definition of $\mu^*$ and then the additivity of $\mu$ on $\mathscr{R}$, we find
\begin{align*}
\mu^*(Q\cap A) + \mu^*(Q\cap \complement A) \leq \sum_{n\ge 1} \mu(B_n\cap A) + \sum_{n\ge 1} \mu(B_n\cap \complement A) = \sum_{n\ge 1} \mu(B_n).
\end{align*}
Taking the infimum over all admissible sequences $B_1, B_2, \dots$ as specified above,
we obtain ${\mu^*(Q\cap A) + \mu^*(Q\cap \complement A)\leq \mu^*(Q)}$. Combining both estimates, we conclude that $A\in \M_{\mu^*}$.

\smallskip
{\em Step 2} -- In this step we prove \ref{it:property-outer2}.
Let $A\in \mathscr{R}$. It is clear that $\mu^*(A)\leq \mu(A)$; this follows by taking $B_1 = A$ and $B_n = \emptyset$ for $n\geq 2$ in \eqref{eq:outermeasdef} (note that $\emptyset\in\mathscr{R}$). The converse estimate $\mu(A)\leq \mu^*(A)$ trivially holds if $\mu^*(A) = \infty$. If $\mu^*(A)<\infty$, choose $B_1, B_2, \dots \in \mathscr{R}$ such that $A\subseteq \bigcup_{n\ge 1} B_n$. Then, by Lemma \ref{lem:ringestimates}, where part \ref{it:ringestimates1} is applied to the inclusion $A\cap B_n\subseteq B_n$ and part \ref{it:ringestimates2} to the union $A = \bigcup_{n\ge 1} A\cap B_n$,
$$\mu(A)\leq \sum_{n\ge 1} \mu(A\cap B_n)\leq \sum_{n\ge 1} \mu(B_n).$$
Taking the infimum over all admissible sequences $B_1, B_2, \dots$ as specified above, we obtain $\mu(A) \leq \mu^*(A)$.

\smallskip
{\em Step 3} -- We now prove assertion \ref{it:Caratheodory2}.
For each $n\ge 1$, define finite measures on $\sigma(\mathscr{R})$ by
$$
\mu_i^{(n)}(A):=\mu_i(A\cap R_n),\quad A\in\sigma(\mathscr{R}),\ i=1,2.
$$
On the ring $\mathscr{R}$ we have, for every $R\in\mathscr{R}$,
$$
\mu_1^{(n)}(R)=\mu_1(R\cap R_n)=\mu(R\cap R_n)=\mu_2(R\cap R_n)=\mu_2^{(n)}(R),
$$
because $R\cap R_n\in\mathscr{R}$ and both $\mu_1,\mu_2$ extend $\mu$ there.
Since $\mathscr{R}$ is closed under finite intersections, Dynkin's lemma
(Lemma~\ref{lem:unique-app}, applied with generating family
$\mathscr{A} = \mathscr{R}\cup\{S\}$) implies that
$\mu_1^{(n)}=\mu_2^{(n)}$ on $\sigma(\mathscr{R})$.
Then, for any $A\in\sigma(\mathscr{R})$,
$$
\mu_1(A)=\lim_{n\to\infty}\mu_1(A\cap R_n)
=\lim_{n\to\infty}\mu_2(A\cap R_n)=\mu_2(A)
$$
by the general properties of measures discussed at the beginning of this appendix.
\end{proof}

To verify the countable additivity condition in Carath\'eodory's result, one may use the following sufficient condition.

\begin{proposition}\label{prop:suffsigma}
Let $\mathscr{R}$ be a ring in a set $S$ and let $\mu:\mathscr{R}\to [0,\infty]$ be an additive map with the property that $\mu(\emptyset) = 0$. If for each nonincreasing sequence $(A_n)_{n\geq 1}$ in $\mathscr{R}$ with $\bigcap_{n\ge 1}A_n =\emptyset$ we have $\limn\mu(A_n) = 0$, then $\mu$ is countably additive on $\mathscr{R}$.
\end{proposition}
\begin{proof}
Let $(B_j)_{j\ge1}$ be a pairwise disjoint sequence in
$\mathscr R$ with $B:=\bigcup_{j\ge1}B_j\in\mathscr R$.
For each $n\ge 1$ put
$A_n:=B\setminus\bigcup_{j=1}^{n-1}B_j$.
Then $A_n\in\mathscr R$ and $A_n\downarrow\emptyset$, so
$\mu(A_n)\to0$ by assumption.
Finite additivity gives
$$
\mu(B)=\sum_{j=1}^{n-1}\mu(B_j)+\mu(A_n),
\quad n\ge1.
$$
Letting $n\to\infty$, we obtain
$\mu(B)=\sum_{j\ge1}\mu(B_j)$
as required.
\end{proof}

We have the following variant of Carath\'eodory's theorem for open sets:

\begin{theorem}\label{thm:open-Caratheodory}
Let $X$ be a topological space and let $\mathscr U$ denote the collection
of its open subsets. Let $\tau:\mathscr U\to[0,\infty)$ satisfy
$\tau(\emptyset)=0$, and suppose that $\tau$ is monotone and countably
subadditive on $\mathscr U$.
Define, for $A\subseteq X$,
$$
        \tau^*(A):=\inf\{\tau(U):\,  A\subseteq U,\ U\in\mathscr U\}.
$$
Assume, moreover, that for all $U,V\in\mathscr U$ and all $\varepsilon>0$
there exists a closed set $F\subseteq U$ such that
$$
        \tau(U\cap V)+\tau(V\setminus F)\le \tau(V)+\varepsilon.
$$
Then $\tau^*$ is an outer measure, every open subset of $X$ is
$\tau^*$-measurable, and therefore the restriction of $\tau^*$ to
$\mathscr B(X)$ is a finite Borel measure. This measure extends $\tau$
on open sets and is outer regular:
$$
        \tau^*(B)=\inf\{\tau^*(U):\,  B\subseteq U,\ U\in\mathscr U\},
        \quad B\in\mathscr B(X).
$$
\end{theorem}

\begin{proof}
It is clear that $\tau^*(\emptyset)=0$ and that $\tau^*$ is monotone. To
prove countable subadditivity, let $A_1,A_2,\ldots\subseteq X$ and fix
$\varepsilon>0$. For each $n\ge1$ choose an open set $U_n\supseteq A_n$
such that
$$
        \tau(U_n)\le \tau^*(A_n)+2^{-n}\varepsilon .
$$
Then $\bigcup_{n\ge1}U_n$ is open and contains $\bigcup_{n\ge1}A_n$, whence
\begin{align*}
 \tau^*\Bigl(\bigcup_{n\ge1}A_n\Bigr)
 \le \tau\Bigl(\bigcup_{n\ge1}U_n\Bigr)
 \le \sum_{n\ge1}\tau(U_n)
 \le \sum_{n\ge1}\tau^*(A_n)+\varepsilon .
\end{align*}
Since $\varepsilon>0$ was arbitrary, $\tau^*$ is an outer measure.

Next we prove that every open set is $\tau^*$-measurable. Let $U\in\mathscr U$
and $Q\subseteq X$. By subadditivity of the outer measure $\tau^*$,
$$
        \tau^*(Q)\le \tau^*(Q\cap U)+\tau^*(Q\cap\complement U).
$$
It remains to prove the converse inequality. If $\tau^*(Q)=\infty$, there is
nothing to prove. Otherwise, let $\varepsilon>0$ and choose an open set
$V\supseteq Q$ such that
$\tau(V)\le \tau^*(Q)+\varepsilon$.
By assumption, there exists a closed set $F\subseteq U$ such that
$$
        \tau(U\cap V)+\tau(V\setminus F)\le \tau(V)+\varepsilon .
$$
Since $Q\cap U\subseteq U\cap V$ and $Q\cap\complement U\subseteq V\setminus F$,
we obtain
\begin{align*}
 \tau^*(Q\cap U)+\tau^*(Q\cap\complement U)
 \le \tau(U\cap V)+\tau(V\setminus F)
 \le \tau(V)+\varepsilon
 \le \tau^*(Q)+2\varepsilon .
\end{align*}
Letting $\varepsilon\downarrow0$ gives the desired inequality. Thus
$U\in\mathscr M_{\tau^*}$.

By Theorem~\ref{thm:outermeasismeas}, $\mathscr M_{\tau^*}$ is a
$\sigma$-algebra and $\tau^*$ restricts to a measure on it. Since all open
sets belong to $\mathscr M_{\tau^*}$, all Borel sets belong to
$\mathscr M_{\tau^*}$. If $U$ is open, then $\tau^*(U)\le\tau(U)$ by
definition, while the reverse inequality follows from the monotonicity of
$\tau$. Thus $\tau^*$ extends $\tau$ on open sets. Outer
regularity is immediate from the definition of $\tau^*$.
\end{proof}

\section*{Lebesgue Measure}

As a first application of Carath\'eodory's theorem
we construct the \emph{Lebesgue measure}\index{Lebesgue!measure}.

For $a = (a_1, \dots,a_d)$ and $b = (b_1, \dots, b_d)$ such that $a_j\leq b_j$ for $j=1,\dots,d$, we write $$(a,b] := \{x\in \R^d: \, a_j < x_j \le b_j, \ j=1,\dots,d\}.$$
The collection $\mathscr{I}^d$ of all finite unions of half-open rectangles is a ring and every set in $\mathscr{I}^d$ can be written as a finite union of {\em disjoint} half-open rectangles.

For $I = (a,b]$ let
$$|I| := \prod_{j=1}^d (b_j - a_j).$$
For $A\in \mathscr{I}^d$ of the form $A = I_1\cup \cdots \cup I_n$, with disjoint half-open rectangles $I_j\in \mathscr{I}^d$\!, define $\lambda_d:\mathscr{I}^d\to [0,\infty]$ by
$$\lambda_d(A) := \sum_{j=1}^n |I_j|.$$

We must check that this number is well defined. To this end, suppose that $A$ has two representations as unions of disjoint half-open rectangles, say $A = (a_1,b_1]\cup \cdots \cup(a_m, b_m] = (c_1,d_1]\cup \cdots \cup(c_n, d_n]$.
Set
$$
I_{ij}:=(a_i,b_i]\cap(c_j,d_j].
$$
For each fixed $i$, the nonempty rectangles $I_{ij}$ form a disjoint
partition of $(a_i,b_i]$, and similarly for each fixed $j$.
Refining in each coordinate by the finitely many endpoints occurring in
the two decompositions, the additivity of interval length and distributivity
of products give
$$
|(a_i,b_i]|=\sum_{j=1}^n|I_{ij}|,
\qquad
|(c_j,d_j]|=\sum_{i=1}^m|I_{ij}|.
$$
Therefore
$$
\sum_{i=1}^m |(a_i,b_i]|
=
\sum_{i=1}^m\sum_{j=1}^n|I_{ij}|
=
\sum_{j=1}^n |(c_j,d_j]|,
$$
which proves the asserted well-definedness.

When the dimension $d$ is fixed and there is no danger of confusion, we write $\lambda$ for $\lambda_d$.

\begin{lemma}\label{lem:Lebesguesigma}
The function $\lambda:\mathscr{I}^d\to [0,\infty]$ is $\sigma$-finite and countably additive on $\mathscr{I}^d$\!.
\end{lemma}

\begin{proof}
The $\sigma$-finiteness follows from writing
$\R^d=\bigcup_{n\ge1}(-n,n]^d$.
Finite additivity and $\lambda(\emptyset)=0$ follow directly
from the definition.
By Proposition \ref{prop:suffsigma}, it therefore suffices
to prove that $\lambda(A_n)\to0$ whenever $(A_n)_{n\ge1}$
is a nonincreasing sequence in $\mathscr I^d$ with
$\bigcap_{n\ge1}A_n=\emptyset$.

Fix such a sequence and let $\eps>0$.
Writing each $A_n$ as a finite union of disjoint half-open
rectangles and moving their lower endpoints slightly inwards,
we can choose $B_n\in\mathscr I^d$ such that
$\ov{B_n}\subseteq A_n$ and
$$
\lambda(A_n\setminus B_n)\le2^{-n}\eps.
$$
The closed sets $\ov{B_n}$ are contained in the compact set
$\ov{A_1}$ and have empty intersection.
Compactness therefore gives an $N\ge1$ such that
$\bigcap_{j=1}^N\ov{B_j}=\emptyset$.

For $n\ge N$, every point of $A_n$ belongs to each $A_j$,
$j\le N$, but lies outside at least one $B_j$.
Consequently,
$$
A_n\subseteq\bigcup_{j=1}^N(A_j\setminus B_j),
\quad n\ge N.
$$
Monotonicity and finite subadditivity now give
$$
\lambda(A_n)
\le\sum_{j=1}^N\lambda(A_j\setminus B_j)
\le\eps\sum_{j=1}^N2^{-j}
<\eps,
\quad n\ge N.
$$
Thus $\lambda(A_n)\to0$, as required.
\end{proof}

\begin{theorem}[Lebesgue measure]\label{thm:Lebesguemeas}
There exists a unique $\sigma$-finite Borel measure $\lambda$ on $\R^d$ such that $$ \lambda((a,b])=\prod_{j=1}^d (b_j-a_j) $$ for every half-open rectangle $(a,b]$.
Moreover, for all $h\in \R^d$ and $A\in \B(\R^d)$ we have $$\lambda(A+h) = \lambda(A),$$
where $A+h := \{x+h:x\in A\}$.
\end{theorem}

\begin{proof}
In Lemma \ref{lem:Lebesguesigma} we have shown that $\lambda$ is $\sigma$-finite and countably additive on the ring $\mathscr{I}^d$\!. Therefore, by Theorem \ref{thm:Caratheodory}, $\lambda$ admits an extension to a $\sigma$-finite measure on $\sigma(\mathscr{I}^d) = \B(\R^d)$.

If $\nu$ is any Borel measure on $\R^d$ such that
$\nu((a,b])=\prod_{j=1}^d (b_j-a_j)$ for every half-open rectangle $(a,b]$,
then $\nu$ is $\sigma$-finite on $\mathscr I^d$.
Since $\sigma(\mathscr I^d)=\calB(\R^d)$, it follows that
$\nu=\lambda$ on $\mathscr{B}(\R^d)$ by
the uniqueness assertion in Carath\'eodory's theorem.

To prove translation invariance, fix $h\in \R^d$\!.
We claim that for every $A\in \B(\R^d)$ the set $A+h$ belongs to $\B(\R^d)$. To see this, let $\mathscr{A}_h = \{A\in \B(\R^d):\, A+h\in \B(\R^d)\}$. This is a $\sigma$-algebra contained in $\B(\R^d)$.  For each open set $A$, the set $A+h$ is open and hence belongs to $\B(\R^d)$. It follows that $\mathscr{A}_h$ contains all open sets. Since $\B(\R^d)$ is the smallest $\sigma$-algebra containing all open sets, it follows that $\B(\R^d)\subseteq \mathscr{A}_h$. Since also $\mathscr{A}_h\subseteq\B(\R^d)$, we obtain equality $\mathscr{A}_h=\B(\R^d)$. This proves the claim.

For $A\in \B(\R^d)$ set $\mu_h(A) := \lambda(A+h)$. Then $\mu_h$ is a $\sigma$-finite Borel measure on $\B(\R^d)$ and for any half-open rectangle $I$, $\mu_h(I) = |I+h| = |I|= \lambda(I)$. By uniqueness, we find that $\mu_h(A) = \lambda(A)$ for all $A\in \B(\R^d)$.
\end{proof}

\section*{Product Measures}

As a second application of Carath\'eodory's theorem, we prove the existence of product measures.

\begin{theorem}[Product measures]\label{thm:product-measures}
Let $(\Om_j,\F_j,\mu_j)$, $j=1,\dots,n$, be $\sigma$-finite measure
spaces.
Then there exists a unique $\sigma$-finite measure $\mu =\prod_{j=1}^n \mu_j$
on the product $\sigma$-algebra $\F = \prod_{j=1}^n\F_j$ which satisfies
$$ \mu(F_1\times\cdots\times F_n) =
\prod_{j=1}^n \mu_j(F_j)$$
whenever the sets $F_j\in \F_j$ satisfy $\mu_j(F_j)<\infty$ for $j=1,\dots,n$.
\end{theorem}

The measure $\mu$ is called the {\em product}\index{product!measure}
of $\mu_1,\dots,\mu_n$.

\begin{proof}
Let $\mathscr{R}$ be the ring consisting of all finite unions of sets of the form $\prod_{j=1}^n F_j$ with $F_j\in \calF_j$ satisfying $\mu_j(F_j)<\infty$ for $j=1,\dots,n$.

For each $j=1,\dots,n$ choose sets $E_{j,m}\in\calF_j$ such that $\mu_j(E_{j,m})<\infty$ for all $m\ge1$ and $\Omega_j=\bigcup_{m\ge1}E_{j,m}$. Then every rectangle $E_{1,m_1}\times\cdots\times E_{n,m_n}$ belongs to $\mathscr R$, and $$ \Omega_1\times\cdots\times\Omega_n =\bigcup_{m_1,\dots,m_n\ge1} E_{1,m_1}\times\cdots\times E_{n,m_n}. $$ Next, $\sigma(\mathscr R)=\prod_{j=1}^n\calF_j$. Indeed, every finite-measure rectangle belongs to $\mathscr R$, and every measurable rectangle $F_1\times\cdots\times F_n$ can be written as a countable union of finite-measure rectangles by intersecting each $F_j$ with the sets $E_{j,m}$.

Since the intersection of
finitely many measurable rectangles of finite measure is a measurable rectangle of finite measure, every $R\in \mathscr{R}$ can be written
as a finite union of {\em disjoint} measurable rectangles of finite measure, say $R = R^{(1)}\cup\cdots\cup R^{(k)}$,
and we may define
$$ \mu(R):= \sum_{j=1}^k \mu(R^{(j)}),$$
where each $\mu(R^{(j)})$ is given by the product formula in the statement of the theorem.
The proof that $\mu(R)$ is well defined is similar to the corresponding proof for Lebesgue measure.
It is clear that $\mu$ is additive on $\mathscr{R}$,
and the above argument shows that $\mu$ is $\sigma$-finite on $\mathscr R$.

We claim that $\mu$ is countably additive on $\mathscr{R}$. Once we know this, the existence
of a unique $\sigma$-finite product measure follows from Carath\'eodory's theorem.

A quick proof of the claim is obtained by applying Proposition \ref{prop:suffsigma} in combination with the dominated convergence theorem. The reader may check that no circularity is
introduced by borrowing this result at this stage. Thus let
$(A_j)_{j\ge 1}$ be a nonincreasing sequence of sets in $\mathscr{R}$ satisfying $\bigcap_{j\ge 1}A_j =\emptyset$.
We must show that $\limj\mu(A_j) = 0$.
Since $A_j$ is a finite union of measurable rectangles of finite measure
and the identity holds for such rectangles by definition, we have
$$ \mu(A_j) = \int_{\Om_n}\cdots\int_{\Om_1} \one_{A_j}\ud \mu_1\cdots\ud\mu_n.$$
The asserted convergence now follows by applying dominated convergence
$n$ times.
\end{proof}

\begin{example} The Lebesgue measure on $(\R^d\!,\mathscr{B}(\R^d))$ is the product measure of $d$ copies
of the Lebesgue measure on $(\R,\mathscr{B}(\R))$.
\end{example}

\section*{Borel Measures on Metric Spaces}

For a Borel measure $\mu$ on a topological space $X$, the following
properties are easily seen to be equivalent by taking complements:
\begin{itemize}
\item for every Borel subset $B$ of $X$ and every $\e>0$
there is an open set $U$ in $X$ such that
$B\subseteq U$ and $\mu(U\setminus B)<\e$;
\item for every Borel subset $B$ of $X$ and every $\e>0$
there is a closed set $F$ in $X$ such that
$F\subseteq B$ and $\mu(B\setminus F)<\e$;
\item for every Borel subset $B$ of $X$ and every $\e>0$
there is an open set $U$ and a closed set $F$ in $X$ such that $F\subseteq B\subseteq U$ and $\mu(U\setminus F)<\e$.
\end{itemize}

\begin{definition}[Regular measures]\label{def:regular}
A Borel measure $\mu$ on a topological space $X$ is called {\em regular}\index{regular measure} if it satisfies the above equivalent conditions.
\end{definition}

\begin{proposition}[Regularity of Borel measures]\label{prop:regular} Every
finite Borel measure on a metric space is regular.
\end{proposition}
\begin{proof}
Let $\mu$ be a finite Borel measure on a metric space $X$.
Denote by $\calA(X)$ the collection of all Borel sets $A$ in $X$
which have the property that for all $\e>0$
there exist a closed set $F$ and an open set $U$
such that $F\subseteq A\subseteq U$ and $\mu(U\setminus F)<\e$.
We must prove that $\calA(X) = \calB(X)$, the Borel $\sigma$-algebra of $X$.

\smallskip
{\em Claim 1}: \ {\em $\calA(X)$ is a $\sigma$-algebra}. It is clear that
$\emptyset\in\calA(X)$ and that $\calA(X)$ is closed under taking complements.
To see that $\calA(X)$ is closed under taking countable unions, let
$(A_n)_{n\ge 1}$ be a sequence of sets in $\calA(X)$. Let $\e>0$ be given and let
$(F_n)_{n\ge 1}$ and $(U_n)_{n\ge 1}$ be sequences
of closed and open sets such that $F_n\subseteq A_n\subseteq U_n$ and
$\mu(U_n\setminus F_n) < \frac{\e}{2^n}$.
The set $U = \bigcup_{n\ge 1} U_n$ is open. In view of
$$ \mu \Big( U\setminus  \bigcup_{n\ge 1} F_n\Big) \le
\sumn\mu ( U_n\setminus F_n)  < \e$$
and the continuity property \ref{it:property-measures4} of finite measures mentioned at the beginning of this appendix, there exists an index $N\ge 1$ such that
$$ \mu \Big( U\setminus  \bigcup_{n=1}^N F_n\Big) < \e.$$
The set $F:= \bigcup_{n=1}^N F_n$ is closed, satisfies
$F \subseteq \bigcup_{n\ge 1} A_n\subseteq U$, and $\mu(U\setminus F)<\e$.

\smallskip{\em Claim 2}: \ {\em $\calA(X)$ contains all closed subsets of $X$}.
We already noted that $\emptyset \in \calA(X)$. If $F$ is a nonempty closed subset of $X$ define, for $k\ge 1$,
$U_k := \{x\in X:\,  d(x,F) < \frac1k\}$. Then each $U_k$ is open and we have
$\bigcap_{k\ge 1} U_k = F$\!. Hence, $\limk \mu(U_k)=\mu(F)$ and the claim
follows.

\smallskip
Combining the two claims, we see that $\calA(X)=\calB(X)$.
\end{proof}

To state the next proposition, we need the following terminology.
In this definition and the next, we limit ourselves to Hausdorff spaces in order to ensure that compact sets are Borel.

\begin{definition}[Tight measures] A finite Borel measure $\mu$ on a Hausdorff topological space $X$
is called {\em tight}\index{tight}
if for every $\eps>0$ there exists a compact set
$K$ in $X$ such that $\mu(X\setminus K)<\eps$.
\end{definition}

The following proposition gives a sufficient condition for tightness.

\begin{proposition}\label{prop:tight}
Every finite Borel measure $\mu$ on a separable complete metric space $X$ is tight.
\end{proposition}
\begin{proof}
Let $(x_n)_{n\ge 1}$ be a dense sequence in $X$ and fix $\eps>0$.
For each integer $k\ge 1$,
the closed balls $\ov B(x_n;\frac1k)$
cover $X$, and therefore there exists an index $N_k\ge 1$ such that
$$\mu\Bigl(\bigcup_{n=1}^{N_k} \ov B(x_n;\tfrac1{k})\Bigr) \ge
\mu(X)-\frac{\eps}{2^{k}}.$$
The set
$$ K := \bigcap_{k\ge 1} \bigcup_{n=1}^{N_k} \ov B(x_n;\tfrac1{k})$$
is closed and totally bounded.
Since $X$ is assumed to be complete, $K$ is compact by Theorem \ref{thm:totally-bounded}. Moreover,
$\mu(\complement K)\le \sum_{k\ge 1} \frac{\eps}{2^{k}} = \eps.$
\end{proof}

For separable complete metric spaces, this result implies the following strengthening of the regularity result provided by Proposition \ref{prop:regular}.

\begin{corollary}\label{cor:Radon} Let $\mu$ be a finite Borel measure on a separable complete metric space $X$.
Then for all Borel subsets $B$ in $X$ and all $\e>0$
there is a compact set $K$ in $X$ such that
$K\subseteq B$ and $\mu(B\setminus K)<\e$.
\end{corollary}
\begin{proof}
The measure $\mu$ is
regular by Proposition \ref{prop:regular}, so
for every Borel set $B$ there
is a closed set $F$ in $X$ such that $F\subseteq B$ and $\mu(B\setminus F)<\frac12\eps$.
By Proposition \ref{prop:tight} there is a compact set $G$ in $X$ such that $\mu(X\setminus G)<\frac12\eps$.
Then $K:=F\cap G$ is compact, contained in $B$, and satisfies $\mu(B\setminus K) < \eps$.
\end{proof}

\begin{definition}[Radon measures]\label{def:Radon}
A finite Borel measure $\mu$ on a Hausdorff topological space $X$ is called {\em Radon}\index{Radon!measure} if for every Borel subset $B$ of $X$ and all $\eps>0$
there is a compact set $K$ in $X$ and an open set $U$ in $X$ such that $K\subseteq B\subseteq U$ and $\mu(U\setminus K)<\eps$.
\end{definition}

\begin{proposition}\label{prop:Radon}
Every finite Borel measure $\mu$ on a separable complete metric space $X$ is Radon.
\end{proposition}
\begin{proof} Fix $\e>0$.
The measure $\mu$ is outer regular by Proposition \ref{prop:regular}, so
for every Borel set $B$ there
is an open set $U$ in $X$ such that $B\subseteq U$ and $\mu(U\setminus B)<\frac12\eps$.
By Corollary \ref{cor:Radon} there is a compact set $K$ in $X$ such that
$K\subseteq B$ and $\mu(B\setminus K)<\frac12\e$.
This proves the result.
\end{proof}

\cleardoublepage  

\cleardoublepage  

\chapter{Integration}

\blfootnote{\rm \small \begin{quote}This book has been published by Cambridge University Press in the series `Cambridge Studies in Advanced Mathematics'. The present version is free to view and download for personal use only. Not for re-distribution, re-sale or use in derivative works.\newline {\copyright} Jan van Neerven\end{quote}}

\noindent
In this appendix we review the Lebesgue integral.

\section*{Measurable Functions}

Let $(\Om_1,\F_1)$ and $(\Om_2,\F_2)$ be measurable spaces.
A function $f:\Om_1\to \Om_2$
is said to be {\em measurable}\index{measurable}
if $\{f\in F\}\in \F_1$ for all $F\in \F_2$.
Clearly, compositions of measurable functions are measurable.
If $\calC_2$ is a subset of $\F_2$ with the property that
$\sigma(\calC_2) = \F_2$, then a function $f:\Om_1\to\Om_2$ is measurable
if and only if $$\hbox{ $\{f\in C\}\in \F_1$ for all $C\in \calC_2$.}$$
Indeed, notice that $\{F\in \F_2:\,  \{f\in F\}\in \F_1\}$ is a sub-$\sigma$-algebra
of $\F_2$ containing $\calC_2$.

If $\mu_1$ is a measure on $(\Om_1,\F_1)$ and $f:\Om_1\to \Om_2$ is measurable, then
$$ f_\#(\mu_1)(F_2) := \mu_1\{f\in F_2\}, \quad F_2\in \F_2,$$
defines a measure $f_\#(\mu_1)$ on $(\Om_2,\F_2)$. This
measure is called the {\em image}\index{image!of a measure}
of $\mu_1$ under $f$. Measurable functions $f$ from a probability space $(\Om,\F\!,\P)$ to another measurable space are called {\em random variables}\index{random variable} and the image of the probability measure $\P$ under $f$ is
called the {\em distribution}\index{distribution} of $f$ under $\P$.

In most applications we are concerned with measurable functions $f$ from a measurable
space $(\Om,\F)$ to $(\K,\calB(\K))$. Such functions are said to be {\em Borel measurable}.\index{Borel!measurable}
Since open sets are Borel, continuous functions are Borel measurable.

In what follows we summarise some
measurability properties of Borel measurable functions. The adjective `Borel' will be
omitted except when confusion is likely to arise.
By the observation made earlier, a function $f:\Om\to \R$ is measurable if and only if
$\{f>a\}\in \F$ for all $a\in \R$, and a function $f:\Om\to \C$ is measurable if and only if
$\{\Re f>a,\, \Im f>b\}\in \F$ for all $a,b\in \R$.
From this it follows that
linear combinations, products, and quotients (if defined)
of measurable functions are measurable.
For example, if the real-valued functions $f$ and $g$ are measurable, then $f+g$ is measurable since
$$ \{f+g>a\} = \bigcup_{q\in\Q} \{f>q\}\cap\{g>a-q\}.$$
The measurability of the sum of two complex-valued measurable functions
is proved similarly.

If $f,g:\Om\to \R$ are measurable, then $f^2$ and $g^2$ are measurable. Indeed, if $a< 0$, then $\{f^2>a\} = \Om$, and if $a\ge  0$, then $\{f^2>a\} = \{-f>\sqrt a\} \cup \{f>\sqrt a\}$, and the same argument applies to $g$. Therefore
$ fg = \frac12 [(f+g)^2 - (f^2+g^2)]$ is measurable. If $f,g:\Om\to \C$ are measurable, the measurability of $fg$ follows by applying the preceding argument to the real and imaginary parts of $f$ and $g$.

If $g$ is measurable and zero-free, then $f/g = f \cdot 1/g$ is measurable, by the continuity of $x\mapsto 1/x$ on $\{x\in \K: \, x\not=0\}$.

If $f:\Om\to \C$ is measurable, then its complex conjugate $\ov f$ is measurable,
and therefore $\Re f = \frac12(f+\ov f)$ and $\Im f = \frac1{2i}(f-\ov f)$ are measurable.
Conversely, if $\Re f$ and $\Im f$ are measurable, so is $f = \Re f + i\Im f$.

If $f = \sup_{n\ge 1} f_n$ pointwise and each $f_n:\Om\to\R$ is measurable, then $f$ is
measurable since
$$\{ f > a\} = \bigcup_{n\ge 1} \{f_n>a\}.$$
It follows from $\inf_{n\ge 1} f_n = - \sup_{n\ge 1} (-f_n)$ that the pointwise
infimum of a sequence of measurable functions is measurable
as well.
It follows that the pointwise limit superior and limit inferior of
measurable functions are measurable, since
$$ \limsup_{n\to\infty} f_n = \lim_{n\to\infty}\Big(\sup_{k\ge n} f_k\Big)
= \inf_{n\ge 1}\Big(\sup_{k\ge n} f_k\Big)$$
and
$ \liminf_{n\to\infty} f_n = -\limsup_{n\to\infty} (-f_n).$
It follows that the pointwise limit $\limn f_n$ of a sequence of measurable
functions $f_n:\Om\to\R$ is measurable. By considering real and imaginary parts separately, the same conclusion holds for pointwise limits of functions $f_n:\Om\to\C$.

In the above considerations involving suprema, infima, and limits it is implicitly assumed that these suprema and infima exist and are finite pointwise.
This restriction can be lifted by considering functions $f:\Om\to
[-\infty,\infty]$.
Such  functions are said to be {\em Borel measurable}\index{Borel!measurable} if
the sets $\{f\in B\}$, $B\in \calB(\R)$, as well as $\{f=\infty\}$ and $\{f=-\infty\}$ are in
$\F$\!.

A {\em simple function}\index{simple function} is a
function $f:\Om\to\K$ that can be represented in the form
$$ f = \sum_{n=1}^N c_n \one_{F_n}$$
with coefficients $c_n\in\K$ and disjoint sets $F_n\in \F$ for all
$n=1,\dots,N$.

\begin{proposition}\label{prop:approx}
A function $f:\Om\to\K$ is measurable if and only if it is the pointwise limit of
a sequence of simple functions $f_n:\Om\to\K$. This sequence may be chosen to satisfy $0\le |f_n|\uparrow |f|$
pointwise.  If $f$ is bounded, it may in addition be arranged that the convergence is uniform. If $f$ is nonnegative, we may furthermore arrange that
$0\le f_n\uparrow f$ pointwise, and uniformly if $f$ is bounded.
\end{proposition}
\begin{proof} There is no loss of generality in taking $\K = \C$.

The `if' part is clear from the fact that measurability is preserved under taking pointwise limits. It remains to prove the `only if' part.

To prove the first assertion, for $j,k\in\Z$ and $n\in \N$ consider the rectangles $R_{jk}^{(n)} = [\frac{j}{2^n}, \frac{j+1}{2^n}) + i [\frac{k}{2^n}, \frac{k+1}{2^n})$ in the complex plane. Let $c_{jk}^{(n)}$ be the unique point in the closure of $R_{jk}^{(n)}$ with minimum modulus. Then the simple functions
$$f_n = \sum_{j,k=-2^{2n}}^{2^{2n}} c_{jk}^{(n)} \one_{\{f\in R_{jk}^{(n)}\}}$$ satisfy $f_n\to f$ and $0\le |f_n|\uparrow |f|$ pointwise as $n\to\infty$. If $f$ is bounded, the convergence is uniform.

To prove the second assertion, for $j\in\N$ and $n\in \N$ consider the intervals $I_{j}^{(n)} = [\frac{j}{2^n}, \frac{j+1}{2^n})$ in the nonnegative real line.
 Then the simple functions $$f_n = \sum_{j=0}^{2^{2n}}\frac{j}{2^n} \one_{\{f\in I_{j}^{(n)}\}}$$ satisfy $0\le f_n\uparrow f$ pointwise. If $f$ is bounded, the convergence is uniform.
\end{proof}

\section*{The Lebesgue Integral}\index{Lebesgue!integral}

The construction of the Lebesgue integral proceeds in two stages: first,
the integral of an arbitrary nonnegative measurable function is defined
(allowing the value $\infty$); second, integrability is introduced and
the integral of an integrable function is defined.

Let $(\Om,\F\!,\mu)$ be a measure space. For a nonnegative simple function $f
= \sum_{n=1}^N c_n \one_{F_n}$ we define
$$ \int_\Om f\ud \mu := \sum_{n=1}^N c_n \mu(F_n).$$
We allow $\mu(F_n)$ to be infinite; this causes no problems because the
coefficients $c_n$ are nonnegative (we use the convention $0\cdot\infty =
0$).
It is easy to check that this integral is well defined, in the sense that it
does not depend on the particular representation of $f$ as a simple function.
Also, the integral is linear with respect to addition and multiplication by
nonnegative scalars,
$$ \int_\Om af+bg\ud \mu = a \int_\Om f\ud \mu + b \int_\Om g\ud \mu,$$
and monotone in the sense that if $0\le f\le g$
pointwise, then
$$ \int_\Om f\ud \mu\le \int_\Om g\ud \mu.$$

In what follows, a {\em nonnegative function} is a function with values in
$[0,\infty]$. Recall that such a function is said to be {\em (Borel) measurable}\index{Borel!measurable} if
all sets $\{f\in B\}$ with $B\in \calB(\R)$, as well as the set $\{f=\infty\}$, are in
$\F$\!.

For a nonnegative measurable function $f$ we choose a sequence of simple
functions
 $0\le f_n\uparrow f$ (see Proposition \ref{prop:approx}; on $\{f=\infty\}$ set $f_n = n$) and define
 $$ \int_\Om f\ud \mu := \limn \int_\Om f_n\ud \mu.$$
The following lemma implies that this definition does not depend on the
approximating sequence.

\begin{lemma} For a nonnegative measurable function $f$ and nonnegative simple functions
$f_n$ and $g$
such that $0\le f_n\uparrow f$ and $ g\le f$ pointwise, we have
$$ \int_\Om g\ud \mu \le  \limn \int_\Om f_n\ud \mu.$$
\end{lemma}

\begin{proof}
Fix $0<\eps<1$ and put
$$
E_n:=\{f_n\ge(1-\eps)g\},\quad n\ge1.
$$
Since $f_n\uparrow f$ and $g\le f$, the sets $E_n$ increase
to $\Om$.
Moreover, $f_n\ge(1-\eps)g\one_{E_n}$, so monotonicity of
the integral for simple functions gives
$$
\int_\Om f_n\ud\mu
\ge(1-\eps)\int_\Om g\one_{E_n}\ud\mu.
$$
Write $g=\sum_{j=1}^m c_j\one_{F_j}$ with $c_j\ge0$ and
pairwise disjoint measurable sets $F_j$.
By continuity from below of $\mu$,
$$
\int_\Om g\one_{E_n}\ud\mu
=\sum_{j=1}^m c_j\mu(F_j\cap E_n)
\uparrow\sum_{j=1}^m c_j\mu(F_j)
=\int_\Om g\ud\mu.
$$
Consequently,
$$
\limn\int_\Om f_n\ud\mu
\ge(1-\eps)\int_\Om g\ud\mu.
$$
Letting $\eps\downarrow0$ proves the result.
\end{proof}

The integral is linear and monotone on the set of
nonnegative measurable functions. Indeed, if $f$ and $g$ are such functions and
 $0\le f_n\uparrow f$ and $0\le g_n\uparrow
g$, then for $a,b\ge 0$ we have $0\le af_n+bg_n\uparrow af+bg$ and therefore
\begin{align*} \int_\Om af+bg \ud \mu & = \limn \int_\Om af_n+bg_n\ud \mu
\\ & = a\limn \int_\Om f_n\ud \mu + b\limn \int_\Om g_n\ud \mu
 = a \int_\Om f\ud \mu+ b\int_\Om g\ud \mu.
\end{align*}
If such $f$ and $g$ satisfy $f\le g$ pointwise, then from $0\le f_n \le \max\{f_n,g_n\} \uparrow g$
we see that
$$\int_\Om f\ud \mu  = \limn \int_\Om f_n\ud \mu\le
\limn \int_\Om \max\{f_n,g_n\}\ud \mu \le \int_\Om g\ud \mu.$$

Let us now take a closer look at the role of null sets.
We begin with a simple observation.

\begin{proposition}\label{prop:null1}
If $f$ is a nonnegative measurable function, then:

\begin{enumerate}[label={\rm(\arabic*)}, leftmargin=*]
\item if $ \int_\Om f\ud \mu <\infty,$ then $\mu(\{f=\infty\}) = 0$;
\item if $ \int_\Om f\ud \mu =0,$ then $\mu(\{f>0\}) = 0$.
\end{enumerate}
\end{proposition}
\begin{proof}
For all $c>0$ we have $0\le c \one_{\{f=\infty\}}\le f$ and therefore
$$ 0\le c  \mu \{f=\infty\} \le  \int_\Om f\ud \mu.$$
The first result follows from this by letting $c\to\infty$.
For the second, note that for all $n\ge 1$ we have $ \frac1n \one_{\{f\ge
\frac1n\}}\le f$ and therefore
$$\frac1n \mu(\{f\ge \frac1n\})
\le \int_\Om f\ud \mu =0.$$ It follows that $\mu(\{f\ge \frac1n\}) = 0$.
Now note that $\{f>0\} = \bigcup_{n\ge 1} \{f\ge \frac1n\}$.
\end{proof}

\section*{The Monotone Convergence Theorem}

The next theorem is the cornerstone of integration theory.

\begin{theorem}[Monotone Convergence Theorem]\label{thm:MCT}\index{theorem!monotone convergence}
Let
$0\le f_1\le f_2\le \dots $ be a sequence of nonnegative measurable functions
converging pointwise to a function $f$. Then,
$$ \lim_{n\to\infty} \int_\Om f_n \ud \mu = \int_\Om f\ud \mu.$$
\end{theorem}
\begin{proof}
The function $f$ is measurable as a pointwise limit of measurable
functions. For $n\ge1$, put
$$
g_n:=2^{-n}\bigl\lfloor 2^n(f_n\wedge n)\bigr\rfloor,
$$
where $\lfloor t\rfloor$ denotes the greatest integer not exceeding
$t$. Each $g_n$ is a nonnegative simple function with $g_n\le f_n$.
Since $f_n\wedge n$ increases with $n$ and each dyadic subdivision refines
the preceding one, the sequence $(g_n)_{n\ge1}$ is increasing.
Moreover,
$$
0\le(f_n\wedge n)-g_n<2^{-n},
$$
so $g_n\uparrow f$ pointwise.

By the definition of the integral and its monotonicity,
\[
\int_\Om f\ud\mu
=\limn\int_\Om g_n\ud\mu
\le\limn\int_\Om f_n\ud\mu
\le\int_\Om f\ud\mu.
\qedhere\]
\end{proof}

\begin{example}\label{ex:subst} We
have the following {\em substitution formula}.\index{substitution formula}
For any measure $\mu$ on $(\Om_1,\F_1)$, measurable
$f:\Om_1\to\Om_2$, and nonnegative measurable
$\phi:\Om_2\to\R$,
$$ \int_{\Om_1} \phi\circ f\ud \mu = \int_{\Om_2} \phi\ud  f_\#(\mu).$$
To prove this, note that this is trivial for simple functions $\phi = \one_F$ with
$F\in \F_2$. By linearity, the identity extends to nonnegative simple functions
$\phi$, and by monotone convergence (using Proposition \ref{prop:approx}) to
nonnegative measurable functions $\phi$.
\end{example}

From the monotone convergence theorem, we deduce the following useful corollary.

\begin{theorem}[Fatou's Lemma]\label{lem:Fat}\index{lemma!Fatou}
Let $(f_n)_{n\ge 1}$ be a sequence of nonnegative
measurable functions on $(\Om,\F\!,\mu)$.
Then
$$ \int_\Om\liminf_{n\to\infty} f_n\ud \mu \le
\liminf_{n\to\infty}\int_\Om f_n\ud \mu.
$$
\end{theorem}
\begin{proof}
 From $\inf_{k\ge n} f_k \le f_m$, $m\ge n$, we infer
$$
\int_\Om \inf_{k\ge n} f_k \ud \mu \le \inf_{m\ge n}\int_\Om f_m \ud \mu.
$$
Hence, by the monotone convergence theorem,
\begin{align*}
 \int_\Om \liminf_{n\to\infty} f_n\ud \mu & = \int_\Om
\lim_{n\to\infty}\inf_{k\ge n} f_k\ud \mu
=  \lim_{n\to\infty}\int_\Om\inf_{k\ge n} f_k\ud \mu
\\ & \le \lim_{n\to\infty}\inf_{m\ge n}\int_\Om f_m \ud \mu
 = \liminf_{n\to\infty}\int_\Om f_n \ud \mu.
 \qedhere
\end{align*}
\end{proof}

\section*{The Dominated Convergence Theorem}

A measurable function $f:\Om\to \K$ is called {\em integrable}\index{integrable}
if $$\int_\Om |f|\ud \mu <\infty.$$
Clearly, if $f$ and $g$ are measurable and $|g|\le |f|$ pointwise, then
$g$ is integrable if $f$ is integrable. In particular, if $f$ is integrable,
then $\Re f$ and $\Im f$ are integrable. If a real-valued function $f$ is integrable, so are its positive part $f^+$ and its negative part $f^-$, and we define
$$ \int_\Om f\ud \mu := \int_\Om f^+\ud \mu - \int_\Om f^-\ud \mu.$$
For a complex-valued integrable function $f$ we define
$$ \int_\Om f\ud \mu := \int_\Om \Re f\ud \mu + i \int_\Om \Im f\ud \mu.$$

For a set $F\in\F$ we write
$$ \int_F f\ud \mu := \int_\Om \one_F f\ud \mu,$$
noting that $\one_F f$ is integrable.
The monotonicity and additivity properties of this integral
carry over to this more general situation, provided we assume that the
functions we integrate are integrable.

The next result is among the most useful in all of analysis.

\begin{theorem}[Dominated Convergence Theorem]\label{thm:DC}\index{theorem!dominated convergence}
Let $(f_n)_{n\ge 1}$ be a sequence of
integrable functions such that $\limn f_n = f$ pointwise.
If there exists an integrable function $g$
such that $|f_n|\le |g|$ for
all $n\ge 1$, then
$$ \limn \int_\Om |f_n-f|\ud \mu = 0.$$
In particular,
$$  \limn \int_\Om f_n\ud \mu  = \int_\Om f\ud \mu. $$
\end{theorem}
\begin{proof}
We make the preliminary observation that if $(h_n)_{n\ge 1}$ is
a sequence of nonnegative measurable functions such that
$\limn h_n = 0$ pointwise and
$h$ is a nonnegative
integrable function such that  $h_n\le h$ for
all $n\ge 1$, then by Fatou's lemma
\begin{align*}
 \int_\Om h\ud \mu & = \int_\Om  \liminf_{n\to\infty} (h-h_n) \ud \mu
\le \liminf_{n\to\infty} \int_\Om h- h_n\ud \mu =  \int_\Om h\ud \mu -
\limsup_{n\to\infty} \int_\Om h_n \ud \mu .
\end{align*}
Since $\int_\Om h\ud \mu$ is finite, it follows that $0\le \limsup_{n\to\infty} \int_\Om h_n \ud \mu\le 0$ and
therefore
$$\lim_{n\to\infty} \int_\Om h_n \ud \mu=0.$$
The theorem follows by applying this to $h_n = |f_n-f|$ and $h = 2|g|$.
\end{proof}

If $h:\Om\to[0,\infty]$ is measurable and $\mu(\{h\ne0\})=0$, then $\int_\Om h\ud\mu=0$.
Indeed, choose simple functions $0\le h_n\uparrow h$. Each $h_n$ has integral zero, and the assertion follows from the definition of the integral. Consequently,
nonnegative measurable functions that are equal almost everywhere have equal integrals. The corresponding assertion for integrable scalar-valued functions follows by considering
the positive and negative parts of their real and imaginary parts.

Consequently, in the main results of this appendix,
in particular in the monotone convergence theorem (Theorem \ref{thm:MCT}) and the
dominated convergence theorem (Theorem \ref{thm:DC}),
{\em we may replace pointwise convergence by pointwise convergence $\mu$-almost everywhere}, where the
latter means that there exists a set $N \in \F$ such that $\mu(N) = 0$ and pointwise convergence takes place on $\Om\setminus N$.
For instance, in the monotone convergence theorem it suffices to assume that $f$ is measurable and
the measurable functions $f_n$ satisfy
$0\le f_n \uparrow f$ pointwise $\mu$-almost everywhere; likewise, in the
dominated convergence theorem it suffices to assume that $f$ is measurable and the measurable functions $f_n$ satisfy $\limn f_n = f$ pointwise $\mu$-almost everywhere
and $|f_n|\le |g|$ pointwise $\mu$-almost everywhere for all $n$.

\section*{Fubini's Theorem}

\begin{proposition}\label{prop:Fub} Let $(\Om_1,\F_1)$ and $(\Om_2,\F_2)$ be measurable spaces
and let $f:\Om_1\times\Om_2\to [0,\infty]$ be measurable with respect to the product
$\sigma$-algebra $\calF_1\times \calF_2$. Then:
 \begin{enumerate}[label={\rm(\arabic*)}, leftmargin=*]
\item\label{it:Fub1} for all $\om_1\in\Om_1$ the
function $\om_2\mapsto f(\om_1,\om_2)$ is measurable;
\item\label{it:Fub2} for all $\om_2\in\Om_2$ the
function $\om_1\mapsto f(\om_1,\om_2)$ is measurable.

\end{enumerate}
\end{proposition}

\begin{proof}
Fix $\om_1\in\Om_1$ and define
$\iota_{\om_1}:\Om_2\to\Om_1\times\Om_2$ by
$\iota_{\om_1}(\om_2):=(\om_1,\om_2)$.
For every measurable rectangle $F_1\times F_2$,
$$
\iota_{\om_1}^{-1}(F_1\times F_2)
=
\begin{cases}
F_2, & \om_1\in F_1,\\
\emptyset, & \om_1\notin F_1.
\end{cases}
$$
Since these rectangles generate the product $\sigma$-algebra,
the mapping $\iota_{\om_1}$ is measurable.
The section $\om_2\mapsto f(\om_1,\om_2)$ is therefore measurable
as the composition $f\circ\iota_{\om_1}$.
This proves \ref{it:Fub1}; the proof of
\ref{it:Fub2} is entirely similar.
\end{proof}

\begin{theorem}[Tonelli--Fubini]\label{thm:Fubini1}\index{theorem!Tonelli--Fubini}
Let $(\Om_1,\F_1,\mu_1)$ and $(\Om_2,\F_2,\mu_2)$ be $\sigma$-finite measure
spaces. Let $f:\Om_1\times\Om_2\to [0,\infty]$ be nonnegative and measurable with respect to the product $\sigma$-algebra $\calF_1\times \calF_2$. Then:
\begin{enumerate}[label={\rm(\arabic*)}, leftmargin=*]
\item\label{it:Fubini1-1} the nonnegative function
$ \om_1\mapsto \int_{\Om_2} f(\om_1,\om_2)\ud \mu_2(\om_2)$
is measurable;
\item\label{it:Fubini1-2} the nonnegative function
$ \om_2\mapsto \int_{\Om_1} f(\om_1,\om_2)\ud \mu_1(\om_1)$
is measurable;
\item\label{it:Fubini1-3} we have
$$
 \int_{\Om_1}\int_{\Om_2} f\ud \mu_2\ud \mu_1 = \int_{\Om_2}\int_{\Om_1} f\ud \mu_1\ud \mu_2 = \int_{\Om_1\times\Om_2} f\ud (\mu_1\times\mu_2).
$$
\end{enumerate}
\end{theorem}

\begin{proof}
By symmetry, it suffices to prove \ref{it:Fubini1-1} and the
equality of the first and third terms in \ref{it:Fubini1-3}.

First suppose that $\mu_1$ and $\mu_2$ are finite.
For $D\in\calF_1\times\calF_2$, put
$$
J_D(\om_1):=\int_{\Om_2}\one_D(\om_1,\om_2)\ud\mu_2(\om_2),
\quad \om_1\in\Om_1.
$$
The sections are measurable by Proposition \ref{prop:Fub},
so these integrals are well defined.
Let $\calD$ be the collection of sets $D$ for which $J_D$ is
measurable and
$$
\int_{\Om_1}J_D\ud\mu_1=(\mu_1\times\mu_2)(D).
$$
Every measurable rectangle $E\times F$ belongs to $\calD$, since
$J_{E\times F}=\mu_2(F)\one_E$ and

\begin{align*}
\int_{\Om_1}J_{E\times F}\ud\mu_1
=\mu_1(E)\mu_2(F)
=(\mu_1\times\mu_2)(E\times F).
\end{align*}

We show that $\calD$ is a Dynkin system.
Clearly $\calD$ contains $\Om_1\times\Om_2$.
If $D_1\subseteq D_2$ with $D_1,D_2\in\calD$, then
$J_{D_2\setminus D_1}=J_{D_2}-J_{D_1}$ is measurable.
Since both measures are finite, subtracting the defining
integral identities gives $D_2\setminus D_1\in\calD$.
If $D_n\uparrow D$ with all $D_n\in\calD$, monotone convergence
in the section integrals gives $J_{D_n}\uparrow J_D$.
Thus $J_D$ is measurable, and another application of monotone
convergence gives
\begin{align*}
\int_{\Om_1}J_D\ud\mu_1
=\limn\int_{\Om_1}J_{D_n}\ud\mu_1
=\limn(\mu_1\times\mu_2)(D_n)
=(\mu_1\times\mu_2)(D).
\end{align*}
Hence $D\in\calD$.

The measurable rectangles form a collection closed under finite
intersections and generate $\calF_1\times\calF_2$.
By Lemma \ref{lem:pi-lambda} we therefore obtain
$\calD=\calF_1\times\calF_2$.
This proves the required assertions for indicator functions,
and linearity extends them to nonnegative simple functions.

For a general nonnegative measurable $f$, choose nonnegative
simple functions $s_n\uparrow f$.
Monotone convergence in the inner integral gives the required
measurability, and monotone convergence in the outer and product
integrals gives the desired identity.

Finally, suppose that $\mu_1$ and $\mu_2$ are $\sigma$-finite.
Choose measurable sets $E_n\uparrow\Om_1$ and $F_n\uparrow\Om_2$
with $\mu_1(E_n)<\infty$ and $\mu_2(F_n)<\infty$.
The product of the restricted measures on $E_n$ and $F_n$
is the restriction of $\mu_1\times\mu_2$ to $E_n\times F_n$,
by uniqueness of product measures.
Applying the finite-measure case on $E_n\times F_n$, and extending
the section integrals by zero outside $E_n$, proves the required
assertions for
$f_n:=\one_{E_n\times F_n}f$.
Since $f_n\uparrow f$, monotone convergence completes the proof.
\end{proof}

A variation of the Tonelli--Fubini theorem holds if we replace `nonnegative measurable' by
`integrable':

\begin{theorem}[Fubini]\label{it:Fubini2}\index{theorem!Fubini}
Let $(\Om_1,\F_1,\mu_1)$ and $(\Om_2,\F_2,\mu_2)$ be $\sigma$-finite measure spaces. If $f:\Om_1\times\Om_2\to \K$ is integrable with respect to the product measure $\mu_1\times\mu_2$, then the following assertions hold:

\begin{enumerate}[label={\rm(\arabic*)}, leftmargin=*]
\item\label{it:Fubini2-1} for almost all $\om_1\in \Om_1$, the function $\om_2\mapsto f(\om_1,\om_2)$ is integrable, and
the function
$ \om_1\mapsto \int_{\Om_2} f(\om_1,\om_2)\ud \mu_2(\om_2)$
is integrable with respect to $\mu_1$;
\item\label{it:Fubini2-2}  for almost all $\om_2\in \Om_2$, the function $\om_1\mapsto f(\om_1,\om_2)$ is integrable, and
the function
$ \om_2\mapsto \int_{\Om_1} f(\om_1,\om_2)\ud \mu_1(\om_1)$
is integrable with respect to $\mu_2$;
\item\label{it:Fubini2-3} we have
$$ \int_{\Om_1}\int_{\Om_2} f\ud \mu_2\ud \mu_1 =
\int_{\Om_2}\int_{\Om_1} f\ud \mu_1\ud \mu_2
=  \int_{\Om_1\times\Om_2} f\ud (\mu_1\times\mu_2).
$$
\end{enumerate}
\end{theorem}

It is understood that the integrals in \ref{it:Fubini2-1} and \ref{it:Fubini2-2} are defined to be $0$ on the exceptional null sets on which the corresponding sections fail to be integrable.
\begin{proof}
Considering real and imaginary parts separately, we may assume that $f$ is real-valued.
We apply Theorem \ref{thm:Fubini1} to the nonnegative measurable function $|f|$ to obtain the measurability of the function $\omega_1\mapsto\int_{\Om_2}|f(\omega_1,\cdot)|\ud\mu_2$ along with the equality $$\int_{\Om_1\times\Om_2}|f|\ud(\mu_1\times\mu_2)=\int_{\Om_1}\int_{\Om_2}|f|\ud\mu_2\ud\mu_1.$$
Similarly, the functions
$\omega_1\mapsto\int_{\Om_2}f^\pm(\omega_1,\cdot)\ud\mu_2$ are measurable and $$\int_{\Om_1\times\Om_2}f^\pm \ud(\mu_1\times\mu_2)=\int_{\Om_1}\int_{\Om_2}f^\pm\ud\mu_2\ud\mu_1.$$
Since $f=f^+-f^-=|f|-2f^-$ and all the integrals above are finite,
we conclude that \ref{it:Fubini2-1} holds and that the first and third
terms in \ref{it:Fubini2-3} are equal. The integrability assertion in
\ref{it:Fubini2-1} follows from the integrability of $\om_1\mapsto \int_{\Om_2}|f(\om_1,\cdot)|\ud\mu_2$. The proof of \ref{it:Fubini2-2} and of the equality of the second and
third terms in \ref{it:Fubini2-3} is entirely similar.
\end{proof}

\cleardoublepage  

\cleardoublepage  

\chapter{Notes}

\blfootnote{\rm \small \begin{quote}This book has been published by Cambridge University Press in the series `Cambridge Studies in Advanced Mathematics'. The present version is free to view and download for personal use only. Not for re-distribution, re-sale or use in derivative works.\newline {\copyright} Jan van Neerven\end{quote}}

\noindent
Historical\, perspectives\, on\, functional\, analysis\, are\, presented\, in\, \cite{Dieu, Monna, Pie}.
Among the many excellent textbooks on functional analysis, our favourites include
\cite{Bressan, Bre, Conway, DunSch1, EinWar, Lax, Rudin, Schechter, Werner, Yosida}.

\section*{Chapter \ref{ch:Banach}}

Exhaustive treatments of the theory of Banach spaces and the Bochner integral are given in \cite{AlbKal, Diestel, DunSch1, LiQue} and in \cite{DieUhl, DunSch1, HNVW1}, respectively.

The proof of the Mazur--Ulam theorem outlined in Problem \ref{prob:Ulam-Mazur} is taken from \cite{Nica}, where further references to its history are given.
The construction, in Problem \ref{prob:nonbdd-operator}, of a linear mapping from $\ell^2$ to $\K$ which fails to be bounded depends on
the existence of an algebraic basis of $\ell^2$.
This, in turn, is deduced with the help of Zorn's lemma. Since Zorn's
lemma is equivalent to the Axiom of Choice, this raises the question
whether a constructive example of an unbounded operator can be given.
Within Zermelo--Fraenkel set theory (ZF)  the answer is negative: it is consistent with ZF that every linear operator on a Banach space is bounded.
In fact, it is a theorem in ZF together with the so-called Axiom of
Determinacy \index{Axiom!of Determinacy} and the Countable Axiom of Choice that every linear operator on a Banach space is bounded \cite[Theorem 567H (c)]{Fremlin}.

\section*{Chapter \ref{ch:ClassicalBanach}}

The proofs of Propositions \ref{prop:Minkowski} and \ref{prop:Holder} are taken from \cite{Haase-Convexity}. The presentation in Section \ref{subsec:Lebesgue-diff}
follows \cite{HNVW1}, where more detailed information on this topic can be found. The classical reference is \cite{Stein}. The Fr\'echet--Kolmogorov compactness theorem is usually stated for bounded subsets of $L^p(\R^d)$. The fact that boundedness follows from the assumptions \ref{it:FK1} and \ref{it:FK2} was observed later by Sudakov; our argument follows \cite{HHM}.

Detailed expositions of Theorems \ref{thm:Jordan} and \ref{thm:RN} may be found in \cite{Bog1}.
For conditional expectations and their basic $L^p$-properties in Section
\ref{subsec:conditional-expectations}, see, for example, \cite{Schilling, Kal}.

Comprehensive treatments of vector lattices, Banach lattices, and positive operators are given in \cite{AliBur, LuxZaa, MeyNie, Schaefer, Zaa}.

For the Schauder fixed point theorem of Problem \ref{prob:Schauder}, see,
for example, \cite{Conway}.
The general question of when a bounded operator $T$ on $L^p(\Omega)$ admits a bounded extension $T\otimes I_X$ to $L^p(\Omega;X)$ is subtle; for a systematic discussion see \cite{HNVW1}.\index{extension!vector-valued} An interesting complement to Problem \ref{prob:LpX-pos}, which treats the case of positivity-preserving operators $T$, is the following result of Paley and Marcinkiewicz--Zygmund: If $T$ is bounded on $L^p(\Om)$ with $1\le p<\infty$ and $X$ is a Hilbert space, then $T \otimes I_X$ admits a unique extension to a bounded operator on $L^p(\Om;X)$ and its norm equals $\n T\n$.\index{extension!Hilbert space-valued} The proof uses properties of Gaussian random variables.

\section*{Chapter \ref{ch:Hilbert-spaces}}

Theorem \ref{thm:ONC} characterises Hilbert spaces up to isomorphism. More precisely, the following deep theorem has been
proved in \cite{LinTza}: A Banach space $X$ is isomorphic to a Hilbert space if and only if every closed subspace of $X$ is the range of a bounded linear projection on $X$.\index{Hilbert space!characterisation via complementation}

The proof of the Radon--Nikod\'ym theorem outlined in Problem \ref{prob:RN} is due to von Neumann and follows \cite{Rudin-RCA}.

It can be quite hard to decide whether a given subspace is dense in a given Hilbert space. The following example may illustrate this. Let $H$ denote the Hilbert space of all scalar sequences $c = (c_n)_{n\ge 1}$ for which the norm
$$ \n c \n_{H}^2:= \sum_{n\ge 1}\frac{1}{n^2}|c_n|^2$$
is finite.
For $x\in\R$ let $\{x\}$ denote its fractional part, that is, the unique real number
in $[0,1)$ such that $x = k+\{x\}$ for some integer $k\in \Z$. For $m = 1,2,3,\dots$ let $$c^{(m)} := (\{\frac{n}{m}\})_{n\ge 1}$$ and note that these sequences belong to $H$.
In a formulation due to \cite{Nyman1950} and \cite{Beurling1955}, refined by \cite{BaezDuarte2003}, the Riemann hypothesis is equivalent to a certain density statement in a Hilbert space of sequences. In the formulation of \cite{Bag}, the linear span of the sequences $(c^{(m)})_{m\ge1}$ defined above is dense in $H$ if and only if the Riemann hypothesis holds.
The Riemann hypothesis is one of the Clay Millennium Prize Problems and arguably the most famous unsolved problem in mathematics.

\section*{Chapter \ref{ch:duality}}

Our proof of Theorem \ref{thm:CK-dual} combines ideas of \cite{Fol} and \cite{RuzTur} with the open set version of Carath\'eodory's theorem proved in the Appendix.
One should be aware that different authors use slightly different definitions of Radon measures.

The norm-preserving real Hahn--Banach theorem was proved by
\cite{Hahn1927} and independently by Banach in 1929.
Banach also established the more general extension theorem
for linear functionals dominated by a sublinear functional
\cite{Banach1929}; see also \cite{Banach1932}.
The extension to complex scalars was added several years later by \cite{BohnenblustSobczyk1938}. Banach also proved the sequential version of the Banach--Alaoglu theorem; the general version is due to \cite{Alaoglu1940}.
A detailed survey of the Hahn--Banach theorem is given in \cite{Buskes}.

The weak and weak$\s$ topologies are examples of {\em locally convex} topologies.
For systematic introductions to this subject we recommend \cite{AliBur, Conway, Rudin}.
Theorem \ref{thm:locrefl} is a special case of the so-called {\em principle of local reflexivity}\index{principle of local reflexivity}. Its full formulation can be found, for example, in \cite{AlbKal}.

For Theorem \ref{thm:Prokhorov} and its ramifications we recommend \cite{Bog2}, where also more refined versions of the theorem can be found.

The marginal problem has a long history in probability theory and
is closely related to optimal transport. Kellerer's papers
\cite{Kellerer1964,Kellerer1984} develop existence and duality
theorems in much greater generality than the compact metric setting
used here. In this setting the criterion can be obtained directly
by extending a positive functional on sums of coordinate functions.
This also explains its close relationship with the Hahn--Banach
proof of Kantorovich duality in Subsection \ref{subsec:optimal-transport}.

The inequality in Corollary \ref{cor:epi-chsh-classical} is due to
\cite{ClauserHorneShimonyHolt1969}. For four probability measures on
$\{-1,1\}^2$ which agree on their common one-coordinate marginals,
Fine's theorem \cite{Fine1982} asserts that the eight CHSH
inequalities (corresponding to the eight possible inequivalent distributions of signs over the $E_{ij}$) are also sufficient for the existence of a common
joint distribution. The deduction of the CHSH inequality from
Kellerer's condition is implicit in \cite{KCBS}; our presentation follows \cite{WaaijerNeerven}.

The formulation of optimal transport in terms of couplings goes back
to \cite{Kantorovich1942}. Our proof of
Theorem \ref{thm:kantorovich-duality} adapts the Hahn--Banach argument
of \cite{Edwards2010} to compact metric spaces. For an account of the
Kantorovich--Rubinstein theorem and its extensions to more general
metric spaces, see \cite{Edwards2011}. Systematic introductions to
optimal transport, Wasserstein distances, and their applications are
given in \cite{AmbrosioGigliSavare2005, Santambrogio2015, Villani2009}.

The result of Problem \ref{prob:strictly-convex} is discussed in \cite{Phelps}. The converse also holds:
If every functional on a closed subspace of $X$ has a unique Hahn--Banach extension of the same norm, then $X\s$ is strictly convex; see \cite{Foguel, Taylor}.

\section*{Chapter \ref{ch:operators}}

General references on the theory of bounded operators include \cite{Beau, GGK1, Nik}.

Most treatments of the Fourier--Plancherel transform use the Schwartz space $\mathscr{S}(\R^d)$  of rapidly decreasing smooth functions instead of our $\calF^2(\R^d)$. Our treatment of the Fourier transform and the Hilbert transform is based on \cite{HNVW1}. The $L^p$-boundedness of the Hilbert transform is classical; we follow \cite{GraI}.

The theory of Fourier multiplier operators can be meaningfully extended to the $L^p$-setting,
where it becomes a powerful tool in the Calder\'on--Zygmund theory\index{Calder\'on--Zygmund theory} of singular integral operators.\index{singular integral} The prime example of such an operator is the Hilbert transform. Detailed treatments of the Hilbert transform and singular integral operators in the $L^p$-setting are given in \cite{GraI, Stein, Stein-HA}.
The theory of the Hilbert transform extends to higher dimensions where analogous statements hold for the {\em Riesz transforms}, whose Fourier multipliers are the functions $m_j\in L^\infty(\R^d)$ given by
\begin{align*}
m_j(\xi):=-i\frac{\xi_j}{|\xi|},\qquad \xi\ne0, \ j=1,\dots,d.
\end{align*}
An exhaustive treatment of these matters belongs to the realm of harmonic analysis; see \cite{GraI, Stein} and Chapter 5 of \cite{HNVW1}.

It was shown by Kwapie\'n that the Fourier--Plancherel transform extends to a bounded operator on $L^2(\R^d;X)$ if and only if $X$ is isomorphic to a Hilbert space (and is unitary if $X$ is a Hilbert space); by results of Bourgain and Burkholder, the Hilbert transform
(see Section \ref{sec:HT}) extends to a bounded operator on $L^p(\R;X)$ for some (equivalently, every) $1<p<\infty$ if and only if $X$ has the so-called {\em UMD property}\index{UMD property}; this abbreviation stands for `unconditionality of martingale differences'. Proofs of these results and their ramifications can be found in \cite{HNVW1, Pis-Mart}.
The UMD property also characterises the boundedness of the vector-valued extension of the It\^o stochastic integral of Problem \ref{prob:Ito}; see \cite{NVW-survey} and the references given therein.

The $\sigma$-finiteness assumption in the Riesz--Thorin theorem (Theorem \ref{thm:RieszThorin}) is actually redundant, but proving this introduces additional technicalities. The proof that $\n T_\C\n = \n T\n$ is taken from \cite{HNVW1}. In this reference, the proof of the Clarkson inequalities sketched in Problem \ref{prob:Clarkson} is attributed to J\"urgen Voigt.
It is a famous result of \cite{Beck75}\index{theorem!Beckner} that the constant $1$ in the Hausdorff--Young inequality
$$\n\calF f\n_{L^{q}(\R^d\!,m)} \le \n f\n_{L^p(\R^d\!,m)}$$
for the Fourier transform with respect to the normalised Lebesgue measure $m$, where $1\le p\le 2$ and $\frac1p+\frac1{q}=1$, can be improved to
$$ \n \calF f\n_{L^{q}(\R^d\!,m)}  \le (p^{1/p}/q^{1/q})^{d/2}\n f\n_{L^{p}(\R^d\!,m)}. $$
In the same paper, Beckner proved the corresponding improvement of Young's inequality and showed that both results are sharp. The proofs rely on (but go beyond) the techniques developed in Section \ref{sec:SQ}.
Counterexamples to \eqref{eq:TcT} in the range $q<p$ can be found in \cite{Riesz} and \cite{Maligr}.

The Marcinkiewicz interpolation theorem is of fundamental importance in the theory of singular integrals; see \cite{GraI, Stein, HNVW3}. Our treatment follows \cite{HNVW1}.
The $L^p$-bounded\-ness of the Hilbert transform, here derived as a consequence of the Riesz--Thorin theorem, can also be derived from the Marcinkiewicz interpolation theorem; the required weak $L^1$-bound is due to Kolmogorov; a full proof can be found in \cite{Duo}.

The result of Problem \ref{prob:Pettis} is due to \cite{Pet38}. It is no coincidence that the counter\-example for $p=1$ in part \ref{it:Pettis2} lives in the space $c_0$: part \ref{it:Pettis1} extends to $p=1$ for all Banach spaces not containing a closed subspace isomorphic to $c_0$. A proof of this fact is given in \cite{DieUhl}. Further results on the Pettis integral can be found in \cite{vDul, Mus,Tal}. For Problem \ref{prob:HT-weak11}, see \cite{Duo}.

\section*{Chapter \ref{ch:spectral}}

There are many excellent treatments of spectral theory, such as the monumental classic \cite{DunSch2} and the monographs \cite{Arv, Aup, Mul}.

The results of Section \ref{sec:spectrum} are standard.
A discussion of the conditions \ref{it:winding1} and \ref{it:winding2} for contours in Cauchy's theorem can be found in \cite{Rudin-RCA}.
The result of Problem \ref{prob:allan-ransford} is due to \cite{Gel}; the proof outlined here
is due to \cite{AllRan}.

\section*{Chapter \ref{ch:CompactOperators}}

The results of Sections \ref{sec:compact}--\ref{sec:spectrum-compact} are standard.
An excellent introduction to Fredholm theory is \cite{Schechter}. Our proof of Theorem \ref{thm:Gohberg-Krein} follows \cite{BleeBoo}. Theorem \ref{thm:Coburn} is from \cite{Cob1, Cob2}; see also \cite{Arv, Douglas}. Toeplitz operators are covered in more depth in
\cite{Arv, BotSil, Douglas, Nik}.

The result of Problem \ref{prob:Pitt-reverse} has the following interesting complement, due to Pitt: For all $1 \le p < q < \infty$, every bounded operator from $\ell^q$ to $\ell^p$ is compact.
An immediate consequence is that the spaces $\ell^p$ and $\ell^q$ are not isomorphic. The proof of Pitt's theorem requires some effort; see for instance \cite{AlbKal, Ryan}.
The result of Problem \ref{prob:Terzioglu} is due to \cite{Ter}.

By elementary $C^\star$-algebra techniques one can derive Theorem \ref{thm:HW}
as a simple corollary of Theorem \ref{thm:Coburn}. We begin by introducing some terminology. A {\em Banach algebra}\index{Banach!algebra} is a Banach space $\mathscr{A}$ endowed with
an associative bilinear multiplication mapping $(x,y)\mapsto xy$ such that $$\n x y\n \le \n x\n \n y\n$$ holds for all $x,y\in \mathscr{A}$.
A {\em unital Banach algebra}\index{unital} is a Banach algebra $\mathscr{A}$ with a unit element, that is, an element $e\in \mathscr{A}$ such that $$e x = x e = x, \quad x\in \mathscr{A}.$$
The {\em spectrum} $\sigma_{\mathscr{A}}(x)$ of an element $x$ of a unital Banach algebra $\mathscr{A}$ is the set of all $\la\in \C$ for which no $y\in \mathscr{A}$ can be found such that $(\la e-x)y = y(\la e-x) = e$, and most of the spectral theory contained in Chapter \ref{ch:spectral} can be routinely extended to this situation.

A {\em $C^\star$-algebra}\index{C@$C^\star$-algebra} is a Banach algebra $\mathscr{A}$ with an {\em involution}\index{involution}, that is, a mapping $x\mapsto x^\star$ on $\mathscr{A}$ satisfying $$(x+y)^\star = x^\star + y^\star, \qquad(cx)^\star  = \ov c x^\star, \qquad (x y)^\star = y^\star x^\star, \qquad (x^\star)^\star=x,$$
as well as
$$ \n x\n = \n x^\star\n, \qquad \n x^\star x\n = \n x\n \n x^\star\n$$
for all $x,y\in \mathscr{A}$ and $c\in\C$. According to the {\em Gelfand--Naimark theorem}\index{theorem!Gelfand--Naimark}, every $C^\star$-algebra
$\mathscr{A}$ is $\star$-isometric to a closed $\star$-subalgebra of $\calL(H)$ for a suitably chosen Hilbert space $H$ (a {\em $\star$-subalgebra} being a subalgebra closed under taking involutions and a {\em $\star$-isometric isomorphism} being an isometric isomorphism $ x\mapsto T_x $ with the additional properties that $T_{xy} = T_x \circ T_y$ and $T_{x^\star} = (T_x)^\star$ for all $x,y\in \mathscr{A}$). This theorem connects the abstract definition of $C^\star$-algebras given here with the concrete approach taken in Section \ref{sec:bicommutant}.

The following generalisation of Proposition \ref{prop:siAT} holds (see \cite[Proposition 1.23]{Fol-AHA} or \cite[Theorem 11.29]{Rudin}): If $\mathscr{A}$ is a closed unital $\star$-subalgebra of a unital $C^\star$-algebra $\mathscr{B}$, then
\begin{align}\label{eq:subalg-spectra}
\sigma_{\mathscr{A}}(x) = \sigma_{\mathscr{B}}(x), \quad  x\in \mathscr{A}.
\end{align} The proof follows that of Proposition \ref{prop:siAT}, except that one needs the
fact that if $x=x^\star$, then
$$
  \si_{\mathscr{B}}(x)\subseteq \R .
$$
This fact can be proved by combining the Gelfand--Naimark theorem and
Proposition \ref{prop:siAT}. An elementary proof can be given as follows. If
$u\in\mathscr B$ is unitary, that is, $uu^\star=u^\star u=e$, the
argument suggested in Problem \ref{prob:unitary-spectrum} proves that
$\si_{\mathscr{B}}(u)\subseteq\mathbb T$. Then, the argument suggested in
Problem \ref{prob:selfadjoint-spectrum} proves that if $x=x^\star$, then
$\si_{\mathscr{B}}(x)\subseteq\R$.

Using \eqref{eq:subalg-spectra}, let us now give a simple alternative proof of Theorem \ref{thm:HW} based on Theorem \ref{thm:Coburn} (cf. Section 4.3 in \cite{Arv}).
The results of Problem \ref{prob:Calkin} prove that for any Banach space $X$ the Calkin algebra
$\calL(X)/\mathscr{K}(X)$ is a unital Banach algebra. If $H$ is a Hilbert space, then
$\calL(H)/\mathscr{K}(H)$ is a unital $C^\star$-algebra.

In what follows we take $H:= H^2(\mathbb{D})$.
 Since $\mathscr{K}(H)$ is contained in the Toeplitz algebra $\mathscr{T}$,
it is meaningful to consider the quotient space $\mathscr{T}/\mathscr{K}(H)$. This space is a closed unital $\star$-subalgebra of $\calL(H)/\mathscr{K}(H)$. By Coburn's theorem, the mapping
$T_\phi+\mathscr K(H)\mapsto\phi$ sets up a $\star$-isometry from $ \mathscr{T}/\mathscr{K}(H)$ onto $C(\mathbb{T})$ and we have
the commuting diagram
\begin{figure}[ht]
 \begin{center}
 \begin{tikzcd}[scale cd=1.1]
  0  \arrow[r]
& \mathscr{K}(H) \arrow[r] \arrow[d, "="]
& \mathscr{T} \arrow[d, "\subseteq"]  \arrow[r, "\pi"]
& C(\mathbb{T}) \arrow[d, "j" ]  \arrow[r]
& 0
\\
   0  \arrow[r]
 & \mathscr{K}(H) \arrow[r]
 & \calL(H) \arrow[r,"\pi"]
 & \calL(H)/\mathscr{K}(H) \arrow[r]
 & 0
\end{tikzcd}
\end{center}
\end{figure}

\noindent
where each occurrence of $\pi$ denotes the relevant quotient mapping and $j$ is the composition of the $\star$-isometry from $C(\mathbb{T})$ onto $\mathscr{T}/\mathscr{K}(H)$ provided by Theorem \ref{thm:Coburn}
and the natural inclusion mapping from $\mathscr{T}/\mathscr{K}(H)$ into  $\calL(H)/\mathscr{K}(H)$.
As such, $j$ is injective.

Now suppose that $\phi\in C(\mathbb{T})$ is such that $T_\phi\in\mathscr{T}$ is Fredholm.
By Atkinson's theorem there exists an $S\in\calL(H^2(\mathbb D))$ such
that $I-T_\phi S$ and $I-ST_\phi$ are compact. This means that $T_\phi$
defines an invertible element in $\calL(H)/\mathscr K(H)$.
By \eqref{eq:subalg-spectra}, $T_\phi$ therefore defines an invertible
element in $\mathscr T/\mathscr K(H)$. Hence its inverse in this quotient
is represented by some element $R\in\mathscr T$; equivalently,
$$
  \pi(RT_\phi)=\pi(T_\phi R)=\pi(I).
$$
By Coburn's theorem, we may write $R=T_\psi+K$ with
$\psi\in C(\mathbb T)$ and $K\in\mathscr K(H)$. Since for
$\phi,\psi\in C(\mathbb T)$ one has
$T_\psi T_\phi-T_{\psi\phi}\in\mathscr K(H)$, it follows that
$$
  j(\psi\phi)=\pi(T_\psi T_\phi)=\pi(RT_\phi)=\pi(I)=j(\one).
$$
The injectivity of $j$ implies $\psi\phi=\one$. This is only possible if
$\phi$ is zero-free.

Problem \ref{prob:Cauchy-index} relates the Cauchy integral formula to Fredholm theory by interpreting the argument principle as an index formula. More generally,
Theorem \ref{thm:Gohberg-Krein} computes an analytic invariant, the Fredholm index, from a topological invariant, the winding number of the symbol. This illustrates the
same principle as the index theorem of \cite{AtiyahSinger1963}, which computes the index of an elliptic operator on a compact manifold without boundary from topological data associated with its principal symbol.

\section*{Chapter \ref{chap:Hilbert-operators}}

Our proof of Theorem \ref{thm:FPR} is taken from \cite{Rudin}.
The clever proof of Proposition \ref{prop:SMT-normal} is based on \cite{Whit}.
The proof of Theorem \ref{thm:Nagy} is taken from \cite{Davies-sgr}.
The proof of the Toeplitz--Hausdorff theorem proposed in Problem \ref{prob:Toep-Haus} is due to
\cite{Li94}. More about numerical ranges can be found in \cite{Gustafson}.

\section*{Chapter \ref{chap:spectral-theorem}}

Most treatments of the spectral theorem for normal operators proceed via the theory of $C^\star$-algebras; see, for example, \cite{Arv, Rudin}. This permits concise abstract proofs, but has the drawback that this theory depends on the existence of maximal ideals, a well-known consequence of Zorn's lemma.
Our approach avoids the use of Zorn's lemma.
Our treatment of the von Neumann bicommutant theorem is taken from \cite{Pedersen}.

Theorem \ref{thm:Neumann-commute} extends, by the same argument, to any finite family of pairwise commuting normal operators: the common operator $S$ may again be taken selfadjoint and the representing functions bounded and continuous. A related $k$-tuple formulation for commuting selfadjoint operators with bounded Borel representing functions may be found in \cite{Nagy}.

Alternative proofs of Theorem \ref{thm:vN-PVM} can be found in \cite{Conway-OpTh, DunSch3} and, for the selfadjoint case, \cite{Nagy}; the proof in the revised printing of the first edition used an idea of Marijn Waaijer (personal communication). \cite{Nagy} also contains an example that shows that the separability assumption cannot be omitted.

The presentation of Section \ref{sec:orthpol} follows \cite{Koelink}.

In \cite{Heuser}, a direct, albeit tricky, proof of the result of
Problem \ref{prob:UT} is given which uses only the continuous functional
calculus for selfadjoint operators.

\section*{Chapter \ref{ch:unbdd}}

References for this chapter include
\cite{AkhGla1, BirSol, DunSch2, EdmEva, Kato-PertTh, ReedSimon2, Rudin, Schm}.
Our proof of the spectral theorem combines elements of \cite{Rudin} and \cite{Schm}, and is elementary in that it avoids the use of $C^\star$-algebra techniques.
For selfadjoint operators a more direct construction of the measurable calculus can be given; see, for example, \cite{Rudin}.
Problem \ref{prob:Arhancet} is taken from \cite{Arhancet-curvature}.

\section*{Chapter \ref{chap:boundaryvalueproblems}}

The connections between functional analysis and the theory of partial differential equations are emphasised in \cite{Bressan, Bre, Jost}. The results of this chapter barely scratch the surface of what can be said in this context.

Sobolev spaces are treated in detail in \cite{Adams, Evans}. Some of our proofs are modelled on those presented in these references.
Our presentation of Propositions \ref{prop:fundvarcalculus} and \ref{prop:der-const}
follows \cite{HNVW1}.
The proof of Theorem \ref{thm:lower-order} follows an idea of \cite{Kry}.

Extension operators are treated in \cite{Adams, Evans}.
In this connection it is worth mentioning that if $D$ is a bounded $C^k$-domain, then every function $f\in C^k(\ov D)$ is the restriction of a function in $C^k(\R^d)$; the analogous result holds for functions
in $C^\infty(\ov D)$ when $D$ is a bounded $C^\infty$-domain.
In both cases, the extensions can be realised by a linear mapping. This result is due to \cite{Seeley}.
Theorem \ref{thm:trace-W1p} is the standard density result
for functions smooth up to the boundary, obtained by
Sobolev extension and mollification; see \cite{Adams, Evans}.
The related {\em Meyers--Serrin theorem}
\index{theorem!Meyers--Serrin} asserts that if $k\in \N$ and $1\le p<\infty$, then
$C^\infty(D)\cap W^{k,p}(D)$ is dense in $W^{k,p}(D)$
for every open set $D\subseteq\R^d$, without any regularity
assumption on its boundary. In this general result, the
approximating functions need not be smooth up to the boundary.

The proof of Theorem \ref{thm:H01-cont} follows \cite{AreUrb} and \cite{Bre}.
The $C^1$ assumptions in the second part of the theorem can be relaxed; see \cite{BieWar}.
If $D$ is bounded and has $C^1$-boundary $\partial D$, then for
$1\le p<\infty$ the mapping
$f\mapsto f|_{\partial D}$, initially defined on $C^\infty(\ov D)$,
extends uniquely to a bounded operator
$T:W^{1,p}(D)\to L^p(\partial D),$
called the {\em trace operator}.\index{trace operator}
Here, we think of $\partial D$ as being equipped with its surface measure. Moreover, for a function $f\in W^{1,p}(D)$ one has $f\in W_0^{1,p}(D)$ if and only if $Tf = 0$.
The details can be found in \cite{Adams, Bre, Evans}.

The construction of Example \ref{ex:noclass-Poisson} is based on \cite{AreUrb}; for a related example see Problem~4.9(a) in \cite{GilTru}.
Proofs of the $H^2$-regularity result mentioned in Remark \ref{rem:MR-d1} and its analogue for Neumann boundary conditions can be found in Chapter 6 of \cite{Evans}.

A comprehensive development of variational boundary value problems is given in \cite{LionsMagenes1968}.
A textbook treatment is given in \cite{Evans}.

Systematic treatments of the finite element method are presented in the monographs \cite{AH, Br}.

\section*{Chapter \ref{chap:forms}}

Excellent references for the theory of forms are \cite{ArendtVogtVoigt2026, Kato-PertTh, Ouh}; some of our proofs follow the latter reference. The identity
$\Dom(\aa)=\Dom(A^{1/2})$ for densely defined closed positive forms is commonly referred to as Kato's second representation theorem.

A generalisation of the classical form method was introduced by \cite{ArendtTerElst2012Forms, ArendtTerElst2012Sectorial}. Here the form domain $V$ is linked to the Hilbert space $H$ by a bounded linear mapping $j:V\to H$ with dense range which need not be injective. In this framework one can associate an $m$-sectorial operator, and hence a holomorphic semigroup, with suitable sectorial forms without assuming that the form is closable. This is useful, in particular, for degenerate differential operators and for boundary operators such as the Dirichlet-to-Neumann operator. The analogous theory for general accretive forms, including criteria for the associated operator to be $m$-accretive, is developed by \cite{ElstSauterVogt2015}.
A related sharp connection between numerical ranges and analyticity is obtained by \cite{HaaseKunstmannVogt2016}. For a symmetric $L^\infty$-contractive semigroup on $L^2$, they show that the numerical range of the negative of its generator on $L^p$, $1<p<\infty$, is contained in the optimal sector of angle $\arcsin|1-2/p|$,
and consequently the $L^p$-semigroup is contractive analytic on the optimal sector of angle
$\arccos|1-2/p|$.
This may be viewed as an $L^p$ counterpart of the relation between sectorial numerical-range estimates and contractive analyticity discussed here and in Chapter \ref{chap:semigroups}.

For the spectral theory of differential operators, see \cite{EdmEva-Elliptic, EdmEva} and the references therein; for variational methods, see \cite{Hen}.
More complete treatments of Dirichlet and Neumann Laplacians can be found in the lecture notes \cite{Are-ISEM} and the survey papers \cite{Are04, GreNgu}, where further references to the literature are given. A standard reference for the theory of elliptic second-order differential operators is \cite{GilTru}.

In some of the results about the Neumann Laplacian, the $C^1$ assumption on the boundary can be relaxed. For example, the Neumann Laplacian has compact resolvent, hence discrete spectrum, if $\partial D$ has the so-called segment property and $D$ lies
`on one side' of it.
The steps are as follows: if $\partial D$ has the segment property, then
$W^{1,2}(D)$ is compactly embedded in $L^2(D)$ by Theorems 5.4.4 and
5.4.17 of \cite{EdmEva}. Consequently, the Neumann Laplacian has compact
resolvent.
The $C^1$-boundary assumptions in Proposition \ref{prop:DeltaFredholm}\ref{it:DeltaFredholm2}, in the second assertion of Theorem \ref{thm:spectra-Laplacian}, and in the Neumann case of Theorem \ref{thm:CF} may be weakened accordingly.

The question ``Can one hear the shape of a drum?'' was asked in \cite{Kac} and answered in the negative in \cite{GWW}. For a concise treatment of Weyl's theorem, see \cite{Hig}. An example of a Jordan curve of positive area is given in \cite{Osgood}.

The inequality $\mu_n\le \la_n$ of Corollary \ref{cor:DirNeumEV} comparing the Dirichlet eigenvalues $\la_n$ and the Neumann eigenvalues $\mu_n$ is strengthened in \cite{Fri} to
 $ \mu_{n+1} \le \la_n$. After reducing to smooth domains, an important step in the proof is the {\em spectral flow identity}\index{spectral!flow}
$$ N_{\rm Neum}(\la) - N_{\rm Dir}(\la) =  n(\la),$$
for $\la>0$ satisfying $\la\not\in \sigma(-\Delta_{\rm Dir})\cup \sigma(-\Delta_{\rm Neum}),$
where
\begin{align*}
N_{\rm Dir}(\la)  & = \#\{\la_n\in  \sigma(-\Delta_{\rm Dir}):\, \la_n < \la\}, \\
N_{\rm Neum}(\la) & = \#\{\mu_n\in  \sigma(-\Delta_{\rm Neum}):\, \mu_n <\la\},
\end{align*}
and $n(\la)$ is the number of negative eigenvalues of the {\em Dirichlet-to-Neumann operator}\index{Dirichlet-to-Neumann operator} $R_\la$, counting multiplicities throughout. This is the operator with domain
$$
  \Dom(R_\la)
  :=
  \bigl\{
  f\in H^{1/2}(\partial D):
  \partial_\nu u|_{\partial D}\in L^2(\partial D)
  \bigr\},
$$
where $u\in H^1(D)$ is the unique weak solution of
$$
\begin{cases}
 -\Delta u  = \la u  &\hbox{on }\, D,\\
 \phantom{-\Delta}u = f &\hbox{on }\, \partial D .
\end{cases}
$$
For $f\in\Dom(R_\la)$ one sets
$$
  R_\la f:=\partial_\nu u|_{\partial D}.
$$
This $L^2(\partial D)$-realisation is selfadjoint, bounded below, and has
compact resolvent.

A simpler proof of Friedlander's theorem, based on a variant of the
Courant--Fischer theorem, was obtained by \cite{Fil}. \cite{LevWei} obtained the inequality $$\mu_{n+d}\le \la_n$$ for bounded convex domains $D$ in $\R^d$\!,
with strict inequality when $d\ge2$ and $\partial D$ is smooth.

Weyl's theorem has been extended to other types of boundary conditions, including Neumann boundary conditions, and positive selfadjoint elliptic operators.
For more details the reader is referred to \cite{SafVas}. Such extensions are nontrivial even for the Laplace operator because the domain monotonicity for Dirichlet eigenvalues
of Lemma \ref{lem:Weyl-monotone2} generally fails for boundary conditions other than Dirichlet. This is demonstrated by the following example, taken  from \cite{Funano}.
We use the notation $a\lesssim b$ to mean that $a\le Cb$ for a universal constant $C$.

For $1\le p\le 2$ let $B_{\ell_d^p}$ denote the open unit ball of
$\ell_d^p=\R^d$ endowed with the norm given by
$\n x\n_p^p=\sum_{j=1}^d|x_j|^p\!.$
If the positive real number $r_{d,p}$ is defined by the condition ${\rm vol}(r_{d,p}B_{\ell_d^p})=1$, then $r_{d,p}\eqsim d^{1/p}$.
The smallest positive Neumann eigenvalue for the Laplace operator on
$D':=r_{d,p}B_{\ell_d^p}$ can be shown to satisfy $$\mu_{2,D'}\gtrsim 1,$$
keeping in mind our convention that $\mu_{1,D'}=0$.

Approximating the segment in $D'$ connecting the origin and the point $(r_{d,p},0,0,\hdots,0)$ by a convex $C^1$-domain
$D\subseteq D'$, one can show that
$$ \mu_{2,D}\eqsim\frac1{r_{d,p}^2}\eqsim \frac1{d^{2/p}}.$$
On the positive side, the same reference proves the following
monotonicity result: if $D,D'\subseteq\R^d$ are bounded convex sets
with $C^1$ boundaries and if $D\subseteq D'$, then
for all $n\ge 1$ we have
\begin{equation*}
\mu_{n,D'}\lesssim d^2 \mu_{n,D},\quad n\ge 1.
\end{equation*}
Letting $p\downarrow1$ in the counterexample shows that the factor $d^2$ is essentially optimal.

The use of the landscape function to study localisation of
eigenfunctions was introduced by Filoche and Mayboroda
\cite{FilocheMayboroda2012}. For the Dirichlet Laplacian, this function
is the classical torsion function. More generally, for a bounded
nonnegative potential $V$ and Dirichlet boundary conditions, the
landscape of $-\Delta+V$ solves $(-\Delta+V)u=1$.
The interpretation of $1/u$ as an effective
potential, including energy identities and exponential decay estimates,
is developed in \cite{ArnoldEtAl2019EffectivePotential}.
The proof of Theorem \ref{thm:landscape-count} adapts the upper-bound
argument in Section 2 of \cite{DavidFilocheMayboroda2021} to the
Dirichlet Laplacian. That paper establishes substantially sharper
results for Schr\"odinger operators, including upper and lower bounds
for the eigenvalue counting function in terms of the landscape.
The effective Weyl approximation \eqref{eq:landscape-effective-Weyl}
and its numerical performance for disordered potentials are discussed
in \cite{ArnoldEtAl2019ComputingSpectra}. Its accuracy at low energies
is a separate issue from the leading asymptotic proved here.

Problems \ref{prob:form-elliptic1} and \ref{prob:form-elliptic3} are taken from \cite{Are-ISEM}.

\section*{Chapter \ref{chap:semigroups}}

Excellent introductions to the theory of $C_0$-semigroups include the
monographs  \cite{App, Davies-sgr, EngNag, Paz}.
For a discussion of the examples
in Section \ref{sec:examples-sg} we refer to these sources.
The monumental 1957 treatise \cite{HilPhi} is freely available
online.

Parts of Sections \ref{sec:semigroups}--\ref{sec:analytic-semigroups} and Figure \ref{fig:Gamma}
are taken from Appendix G of \cite{HNVW2}, which in turn is based on the corresponding material in the author's lecture notes for the 2006/07 Internet Seminar `Stochastic Evolution Equations', available on the author's webpage.

Theorem \ref{thm:Phillips} is due to \cite{Phi55}.
Theorem \ref{thm:Jorgensen} is a special case of a result of \cite{Jo}. The idea of using this result to prove Wiener's Tauberian theorem
is from \cite{Nee97}. The proof of Theorem \ref{thm:strong-stability} uses the classical
Fourier argument for Ingham's Tauberian theorem; see
\cite{ChillSeifert2016}.
Theorem \ref{thm:HilleYosida} was obtained independently by Hille and Yosida
near the end of the 1940s. An extension
to arbitrary $C_0$-semigroups, which is somewhat more technical to state, was found soon afterwards.
A detailed account of the theorem and its history
is given in \cite{EngNag}. The intimate connections between semigroups and the theory of Laplace transforms
are emphasised in \cite{ABHN}.

Fuller treatments of the abstract Cauchy problem are given in  \cite{Ama, Are04, Tan}.

Analytic semigroups are treated in detail in \cite{Lun}. A bounded analytic $C_0$-semigroup which is contractive on the positive real line may fail to be a contractive analytic $C_0$-semigroup; for an example see \cite{Fuhrman}, where Problem \ref{prob:GoldysvanNeerven} can be found.

Maximal regularity for bounded analytic $C_0$-semigroups on Hilbert spaces was first proved in \cite{Simon64}.
The result remains valid if $L^2$ is replaced by $L^p$ with $1<p<\infty$ throughout; this follows from rather deep extrapolation arguments for singular integral operators and falls outside our scope.
For a full treatment, as well as references to the extensive literature
on the subject, we refer the reader to \cite{HNVW3}, whose treatment we
follow. The method of applying maximal regularity to solve the time-dependent
problems in Section \ref{subsec:maxreg} goes back to \cite{CleLi} and has been extended to cover a wealth of other nonlinear problems.

Conversely, it is worth noting that every generator $-A$ of an analytic
$C_0$-semigroup on a complex Hilbert space $H$ which is contractive on some sector about the positive real axis can be
represented in divergence form, in the following precise sense:
there exists a Hilbert space $\HH$, a closed operator
$V: \Dom(V) \subseteq H\to \HH$ with dense domain and dense range,
and a bounded operator $B \in \calL(\HH)$ satisfying $\Re\iprod{Bx}{x}_\HH\ge \beta\n x\n_{\HH}^2$ for some $\beta>0$ and all $x\in\HH$, such that
$$A = V^\star B V.$$
More precisely, there exists a densely defined, closed, sectorial form $\aa$ in $H$ with domain $\Dom(\aa) = \Dom(V)$ such that $A$ is the operator associated with $\aa$ and
$$ \aa(g,h)= \iprod{ BV g}{ Vh}, \quad  g,h\in \Dom(V).$$
A proof of this result can be found in \cite{MaaNee}, where it is also shown that this representation is essentially unique.

Our proofs of Theorem \ref{thm:heat-pos} and Lemma \ref{lem:dominDirNeum} are taken from \cite{Are-ISEM}.

For the gradient-flow interpretation of the heat equation
discussed in Subsection \ref{subsec:heat-sgr}, see the
general theory of gradient flows of convex functionals
in Hilbert spaces developed in \cite{Brezis1973}. A systematic treatment
of gradient flows in metric spaces and spaces of probability
measures is given in \cite{AmbrosioGigliSavare2005}.
The dependence on the choice of ambient geometry is illustrated
by the work of Jordan, Kinderlehrer, and Otto
\cite{JordanKinderlehrerOtto1998}: the heat equation on $\R^d$
can also be interpreted as the gradient flow of the entropy
$\rho\mapsto\int_{\R^d}\rho\log\rho\ud x$ with respect to the
quadratic Wasserstein distance. This provides a different
variational interpretation of the same equation, alongside
its description as the $L^2$-gradient flow of the Dirichlet energy.

The Ornstein--Uhlenbeck semigroup has many interesting properties, for which the reader is referred to \cite{HNVW4,Jan,Nua}.
Probabilistically, up to a rescaling of time, it arises as the transition
semigroup associated
with the solution $(u_x(t))_{t\ge 0}$ of the stochastic differential equation
$$ {\rm d}u(t) = - \frac12 u(t)\ud t+ {\rm d}B_t, \quad t\ge 0,$$
with initial condition $u(0)=x$; the driving process $(B_t)_{t\ge 0}$ is a standard Brownian motion in $\R^d$. More precisely, for all $t\ge 0$ and $f\in L^2(\R^d\!,\gamma)$ one has
$$ {OU}(t/2)f(x) = \E (f(u_x(t)))$$ for almost all $x\in \R^d$.

The domain identification $\Dom(L) =  W^{2,p}(\R^d\!,\gamma)$ for the Ornstein--Uhlenbeck operator $L$ in $L^p(\R^d\!,\gamma)$ with $1<p<\infty$ is a special case of a more general result due to \cite{MPRS}. This paper also contains references to earlier papers on this subject, in particular regarding the special case $p=2$.
The Ornstein--Uhlenbeck semigroup extends to an analytic
$C_0$-contraction semigroup on $L^p(\R^d\!,\gamma)$ for $1<p<\infty$, with optimal angle
$\theta_p =  \arccos|1-2/p|.$
This result is due to \cite{Epp}, who also showed that the exact
domain of holomorphy is the set
\begin{align*}
E_p := \{z = x+iy\in\C:\,
\cos(\theta_p)|\sin(y)| < \sin(\theta_p)\sinh(x)\}.
\end{align*}
A simpler proof of this description of the domain of holomorphy was
given in \cite{NP18}.

The proof of Theorem \ref{thm:Strichartz} uses the standard
duality argument, often called the $TT^*$ argument; see
\cite[Section~2.3]{Tao2006}.
Lemma \ref{lem:Strichartz-fractional} is a special case of the
Hardy--Littlewood--Sobolev inequality. Its proof follows
Hedberg's argument based on the Hardy--Littlewood maximal
theorem \cite{Hedberg1972}.
The $L^p$-Strichartz estimate presented in the main text can be extended to an estimate in terms of the mixed norm $L^p(\mathbb R;L^q(\mathbb R^d))$ for suitable exponents $(p,q)$. Analogous estimates hold for wave equations. For a detailed treatment of these matters, the reader is referred to the lecture notes \cite{Schn}; this reference also contains the argument in the first part of Section \ref{subsec:wave}. An elementary introduction is \cite{Ste-Sha-4}.

The Gearhart--Pr\"uss theorem (see Problem \ref{prob:GP}) does not extend to general Banach spaces \cite{Are95}. A counterexample is obtained, after rescaling, from the $C_0$-semigroup given by $(S(t) f)(x) := f(xe^t)$ on \(L^p(1,\infty)\cap L^q(1,\infty)\), where \(1\le p<q<\infty\).

If $A$ is a densely defined operator acting in a Banach space $X$ with the property that $(-\infty,0)\subseteq \varrho(A)$ and $\n (\la+A)^{-1}\n \le M/\la$ for $\la>0$,
there is a unique densely defined closed operator $A^{1/2}$
such that \((A^{1/2})^2=A\) and, for every \(x\in\Dom(A)\),
$$ A^{1/2}x = \lim_{\eps\downarrow 0} (\eps+A)^{1/2}x = \frac1\pi\int_0^\infty t^{-1/2}(t+A)^{-1}Ax\ud t.$$
A proof of this result can be found in Section 3.8 of \cite{ABHN}.
In particular, it applies to $A = -B$ whenever $B$
is the generator of a uniformly bounded $C_0$-semigroup on $X$.
The result should be compared with Proposition \ref{prop:pos-sqrt}.

Let us apply this to the accretive divergence form operator $L_a:= -{\rm div}(a\nabla)$
associated with the densely defined form
$$
        \aa_a(f,g) :=
        \int_{\R^d} {a\nabla f}\cdot \ov{\nabla g}\ud x,
        \quad f,g\in H^{1}(\R^d),
$$
where the $d\times d$ matrix $a=(a_{ij})_{i,j=1}^d$ is assumed to have
bounded measurable coefficients satisfying the uniform
ellipticity condition stated in Example \ref{ex:exmaples-revisited-II}.
Since $-L_a$ generates an analytic $C_0$-semigroup of contractions,
the square root $L_a^{1/2}$ is well defined.
The {\em Kato square root problem}\index{Kato square root problem} is to
decide whether the domain equality
$$
        \Dom(L_a^{1/2})=\Dom(\aa_a)=H^1(\R^d)
$$
holds, with equivalence of homogeneous norms $\|L_a^{1/2}f\|\eqsim \|\nabla f\|$
for all functions $f$ in this common domain.
Beginning with the papers \cite{Kato-sqrt} and \cite{McI-Kato},
this problem has had a long and interesting history. It was finally resolved
in the affirmative, in the generality stated here, in \cite{AHLMT02}. This paper
also contains references to the various special cases that had been obtained before.
An alternative proof based on the theory of bisectorial operators was subsequently obtained in \cite{AKM06}.

\section*{Chapter \ref{chap-HilbertSchmidt-TraceClass}}

The results of Sections \ref{sec:Hilbert-Schmidt}--\ref{sec:trace-duality} are standard. Section \ref{subsec:partialtrace} follows \cite{Attal}.

The argument in Step 1 of Proposition \ref{prop:Schur} follows \cite[Lemma XI.6.21]{DunSch2}.

Our proof of Lidskii's theorem (Theorem \ref{thm:Lidskii}) follows \cite{Sim}. A survey of the connections between determinants and traces, containing a proof of MacMahon's formula as well as a treatment of Fredholm determinants, is \cite{Cartier}.
For much more on this topic the reader may consult \cite{Simon-TraceIdeals}.

For positive kernel operators, Theorem \ref{thm:Mercer} is due to \cite{Mercer}. The extension to general integral operators of trace class is taken from \cite{Bir}. In that paper it is also shown how to extend the result to general measure spaces provided the
corresponding $L^2$-space is separable. Further interesting results on this topic can be found in \cite{Bris}. It is worth noting that not every integral operator with continuous kernel is of trace class; a classical counterexample can be found in \cite{Carleman}.

Theorem \ref{thm:Fedosov} is taken from \cite{Murphy}. Theorem \ref{thm:Toeplitz-trace} is from \cite{HH}. The trace formula for the Dirichlet heat semigroup, Theorem \ref{thm:heat-trace},
is obtained here as an Abelian consequence of Weyl's eigenvalue asymptotics
from Theorem \ref{thm:Weyl}.
The derivation of Euler's formula from the trace of the inverse Dirichlet Laplacian on the interval is taken from \cite{Grieser}.

The proof outlined in Problem \ref{prob:HaaseL2Linfty} is taken from \cite{Are-ISEM}, where it is attributed to Markus Haase. Problem \ref{prob:ConnesConsani} is taken from \cite{ConCon}, Problem \ref{prob:Murphy1} is from \cite{Murphy}, and Problem \ref{prob:HeltonHowe} is from \cite{HH}.

\section*{Chapter \ref{chap:QM}}

Historical aspects of the interaction between functional analysis and quantum mechanics are well documented in \cite{Landsman-FA-QM}.
An excellent modern introduction to quantum mechanics from the mathematician's point of view is \cite{Hall}. More advanced treatments are offered in \cite{Landsman1998, Lan, Mac, Parth, Takh}.

The mathematical formulation of quantum mechanics in terms of
Hilbert spaces is due to \cite{Neu}. The question of which
geometric, logical, or operational assumptions lead to this
framework has since motivated numerous axiomatic approaches;
see \cite{Chiribella-DAP, Ludwig-AxiomaticI, Varad}.

One such approach uses orthomodularity, the requirement that
$Y+Y^\perp=H$ for every closed subspace $Y$ of an inner product
space $H$. Over the complex numbers this property characterises
Hilbert spaces \cite{Pir, AA}. Sol\`er's theorem \cite{Sol}
extends this result to orthomodular Hermitian spaces over
division rings with involution: the existence of an infinite orthonormal
sequence forces the space to be a real, complex, or quaternionic
Hilbert space. For a survey, see \cite{Hol}; for an application
to categorical characterisations of real and complex Hilbert
spaces, see \cite{HeuKor}.
\index{Hilbert space!characterisation via Hermitian forms}

A complementary approach describes classical and quantum systems
in the common language of $C^\star$-algebras. By the
Gelfand--Naimark theorem, every $C^\star$-algebra is
$\star$-isometrically isomorphic to a norm-closed
$\star$-subalgebra of $\calL(H)$ for a suitable Hilbert space $H$.
In the commutative case it is $\star$-isometrically isomorphic
to $C_0(X)$ for some locally compact Hausdorff space $X$;
if the algebra is unital, $X$ is compact.
Thus commutative algebras describe classical systems, while
noncommutative algebras accommodate quantum systems.
Comprehensive treatments of $C^\star$-algebras and their states
are offered in \cite{Black, BraRob, Pedersen, Take}.

For the normal states and countably additive observables
considered in this chapter, the more specific setting of
von Neumann algebras is particularly natural. Since every commutative von Neumann algebra
can be represented as a space $L^\infty(X,\mu)$, this embeds
the measure-theoretic description of classical systems into this
operator-algebraic framework.
In this setting an {\em observable} is described by a normal unital
completely positive linear mapping $\Psi:\mathscr A\to\mathscr M$,
where $\mathscr A$ is a commutative von Neumann algebra describing the
outcomes and $\mathscr M$ is the von Neumann algebra of the system.
Complete positivity means that, for every $n\ge1$, the entrywise
extension
$$
[a_{jk}]_{j,k=1}^n\mapsto[\Psi(a_{jk})]_{j,k=1}^n
$$
maps positive matrices over $\mathscr A$ to positive matrices over
$\mathscr M$. Since $\mathscr A$ is commutative, positivity already
implies complete positivity. Normality means that the restriction of $\Psi$ to the closed
unit ball is continuous for the weak operator topologies.
Sharp observables correspond to the mappings $\Psi$ which are
$\star$-homomorphisms.

To make the
connection with POVMs explicit, let $(\Omega,\calF,\nu)$
be a probability space and suppose that $\mathscr M\subseteq\calL(H)$.
A normal unital positive mapping
$\Psi:L^\infty(\Omega,\nu)\to\mathscr M$ defines a POVM by
$$
Q_F:=\Psi(\one_F),\quad F\in\calF.
$$
For positive functionals on $\calL(H)$, this characterisation
was proved in Theorem~\ref{thm:normal-continuity}.
Normality gives the countable additivity of $Q$, and we have
$$
\Psi(g)=\int_\Omega g\ud Q,\quad g\in L^\infty(\Omega,\nu).
$$
Conversely, every $\mathscr M$-valued POVM satisfying
$Q_F=0$ whenever $\nu(F)=0$ defines such a normal mapping.
Moreover, $Q$ is projection-valued if and only if $\Psi$ is a
$\star$-homomorphism.

The classical case is obtained by taking the target algebra to be
commutative as well. More precisely, let $(X,\calX,\mu)$ and
$(\Omega,\calF,\nu)$ be probability spaces, with
$(\Omega,\calF)$ standard Borel. Every normal unital positive mapping
$$
\Psi:L^\infty(\Omega,\nu)\to L^\infty(X,\mu)
$$
is represented by a Markov kernel $K:X\times\calF\to[0,1]$ through
$$
(\Psi g)(x)=\int_\Omega g(\omega)\,K(x,\ud\omega),
\quad \mu\hbox{-almost every }x,
$$
for each bounded measurable $g$. The kernel is unique up to changes
on a $\mu$-null set, and if $\nu(F)=0$, then $K(\cdot,F)=0$ almost everywhere.

Conversely, every Markov kernel with this property defines a normal
unital positive mapping between these $L^\infty$-spaces.
The mapping is a $\star$-homomorphism precisely when the kernel
is deterministic outside a $\mu$-null set: there is a measurable
$f:X\to\Omega$, with $f_\#\mu\ll\nu$, such that
$$
K(x,\cdot)=\delta_{f(x)}
\quad\mu\hbox{-almost everywhere}.
$$
In this case $\Psi g=g\circ f$.
In this representation, Markov kernels describe classical
unsharp observables, and measurable functions describe sharp
ones. The pointwise formulation in
Subsection~\ref{subsec:Markov-kernels-POVMs} uses bounded
measurable functions without identifying functions that agree
almost everywhere, and requires no reference measures.

Proofs of Gleason's theorem mentioned at the end of Section \ref{subsec:states} can be found in \cite{Lan, Parth}.
For effects and Busch's theorem (Theorem \ref{thm:Busch}), as well as their
role in the theory of unsharp observables, see \cite{BGL, BLPY}. The phase POVM of Section \ref{subsec:number-phase} goes back to \cite{GarWong}.

Physically, a spin-$\frac12$ degree of freedom carries the two-dimensional
irreducible unitary representation of $SU(2)$ and provides a physical
realisation of a qubit. More generally, the irreducible representation
on $\C^{n+1}$ has spin $n/2$. Thus, even $n$ gives integer spin and odd
$n$ half-integer spin; under the hypotheses of the relativistic
spin--statistics theorem these correspond to bosons and fermions,
respectively. More on this topic can be found in \cite{Sternberg,Woit}.

References on the theory of POVMs include
\cite{Berberian, DaviesOpen, Holevo, Landsman1998}.
Our proof of Naimark's dilation theorem
(Theorem~\ref{thm:Naimark}) is related to the one presented in
\cite{AkhGla2}. A derivation from Stinespring's dilation theorem
is sketched in \cite{HLLL}. In the
von Neumann algebra framework outlined above, the argument can be sketched as follows.
Let $Q$ be a POVM on $(\Omega,\calF)$.
 We may assume that $H\ne\{0\}$.
Since $H$ is separable, there exists a positive trace class operator
$T_0$ with $\tr(T_0)=1$ and $\ker T_0=\{0\}$. The probability
measure
$$
\nu(F):=\tr(Q_FT_0),\quad F\in\calF,
$$
satisfies $\nu(F)=0$ if and only if $Q_F=0$. Consequently, the
bounded functional calculus of Proposition~\ref{prop:POVM-calculus}
induces a normal unital completely positive mapping
$\Psi_Q:L^\infty(\Omega,\nu)\to\calL(H)$.
The normal version of Stinespring's theorem provides a Hilbert
space $\wt H$, an isometry $J:H\to\wt H$, and a normal unital
$\star$-homomorphism
$\Pi:L^\infty(\Omega,\nu)\to\calL(\wt H)$
such that
$$
\Psi_Q(g)=J^\star\Pi(g)J,\quad g\in L^\infty(\Omega,\nu).
$$
Since $\Pi$ is a normal $\star$-homomorphism,
$$
\wt P_F:=\Pi(\one_F),\quad F\in\calF,
$$
defines a projection-valued measure, and we have $Q_F=J^\star\wt P_FJ$ for all $F\in\calF$, that is,
$\wt P$ dilates $Q$.
Thus every unsharp observable is the compression of a sharp
observable on a larger Hilbert space. Taking the Stinespring
representation to be minimal gives the minimality condition in
Theorem~\ref{thm:Naimark}, since simple functions are dense
in $L^\infty(\Omega,\nu)$.

Theorem \ref{thm:commutative-POVM-kernel} is from \cite{Beneduci}.
Theorem \ref{thm:effective-model} is a measurable Markov kernel
version of the classical representation theorem
for statistical models due to \cite{Holevo}. In Holevo's formulation, any
separated finite-dimensional statistical model is obtained as a reduction of
a classical model whose phase space is the set of extreme states; the allowed
classical measurements are the transition probabilities obtained by restricting
the original affine measurements to these extreme states. This reference also derives the corresponding model for the qubit. For a systematic study of the equivalence relation defining the
Holevo space and its interpretation in terms of indiscernibility
of quantum states, see \cite{vanNeervenWaaijer2026}.

The characterisation of quantum correlation matrices in
Proposition \ref{prop:epi-quantum-correlation-duality} and the sharp
bound $2\sqrt2$ for the CHSH expression are due to
\cite{Tsirelson1980}, where the author's name is transliterated as
Cirel'son; see also \cite{Tsirelson1987}.
This characterisation provides a converse to
Lemma \ref{lem:epi-quantum-vectors}: every finite matrix of inner
products of real unit vectors can be realised as a quantum
correlation matrix, using finite-dimensional Hilbert spaces and
sharp binary observables. Together with the forward implication,
this identifies the real Grothendieck constant with the best
universal comparison constant for bipartite binary correlation
expressions.

Grothendieck's inequality goes back to
\cite{Grothendieck1953}. The proof given here follows Krivine's Gaussian and
tensor-power argument \cite{Krivine1979}; see also the exposition
in \cite{Briet2018Grothendieck}. The number
$\pi/(2\log{(}1+\sqrt2{)})$ is an upper bound for the real Grothendieck
constant, rather than its value. In fact, \cite{BravermanEtAl2013}
proved that the inequality between these two constants is strict.
More recently, Heilman, Jones and Malavolta
\cite{HeilmanJonesMalavolta2026} obtained the computer-assisted
bounds $$1.773\le K_G\le1.7799.$$
A comprehensive review of Grothendieck's inequality, equivalent formulations, and its ramifications is \cite{Pisier2012}.
The inequality also has a complex version, in which
the matrix entries and Hilbert space are complex and the scalar
signs in the definition of $\beta(M)$ are replaced by complex
numbers of modulus one. The optimal constant, called the
{\em complex Grothendieck constant} and denoted by $K_G^{\C}$,
is strictly smaller than $K_G$; its exact value is also unknown.
For further details, see \cite{Pisier2012}.

The result stated in Problem \ref{prob:phi-bicomm} is taken from \cite{Pedersen}.

\section*{Chapter \ref{ch:SymmQuant}}

The version of Wigner's theorem stated here (Theorem \ref{thm:Wigner}) is the projective, or pure-state, form of Wigner's original theorem. A complete proof is given in \cite{Lan}.
Bargmann's theorem mentioned in the text is in \cite{Bar}; a complete proof is also found in \cite{Parth}.

For introductions to LCA groups we recommend \cite{Fol-AHA, RudinLCA}. These references contain proofs of the basic facts about Haar measures and Fourier analysis on LCA groups that were assumed without proof in the main text. More complete treatments of covariance lead to the notion of {\em systems of imprimitivity} studied in \cite{Mac}. For in-depth discussions of covariance and the way in which it determines observables we recommend \cite{Holevo, Parth, Varad}. A discussion from the physicist's point of view is given in \cite{BGL}.
Theorem \ref{thm:Ugamma} is a special case of a generalisation of Stone's theorem (Theorem \ref{thm:Stone}) for arbitrary strongly continuous unitary representations of $G$; see Theorem 4.5 in \cite{Fol-AHA}.
Pauli's no-go theorem originates in his 1933 \emph{Handbuch} article \cite{Pauli1990-Reprint}; for modern operator-algebraic treatments of measurement and time observables see \cite{BLPY, Holevo}.

The presentation of the Stone--von Neumann theorem follows \cite{Foll-PhaseSpace} and \cite{Hall}. The theorem admits a generalisation to LCA groups, essentially due to Mackey. For modern references to the literature on this subject the reader is referred to the survey article \cite{Ros}.
The formula for the Orn\-stein--Uhlenbeck semigroup in Theorem \ref{thm:exptL} goes back at least to \cite{unter}; in its present form it is taken from \cite{NP18}.
From the viewpoint of analysis, the Weyl calculus is a theory of pseudodifferential operators, and many properties of this calculus are instances of the pseudodifferential calculus developed by H\"ormander; comprehensive treatments are given in \cite[Vol.~III]{Hormander1985} and \cite{Taylor1981}.

A {\em quantisation scheme} assigns to a classical observable $f$ on phase space an operator ${\rm Op}_h(f)$ on a Hilbert space in such a way that real-valued functions map to selfadjoint operators and Poisson brackets are recovered from commutators in the limit $h\downarrow 0$.
Theorem~\ref{thm:commutator-Poisson} goes back to Weyl's correspondence \cite{Weyl1927}. The phase space formulation of quantum mechanics was developed by \cite{Moyal1949}.
Modern treatments of Weyl calculus and the Moyal product can be found in
\cite{Foll-PhaseSpace, Zworski2012}.
In less formal terms, the {\em Groenewold--Van Hove theorem}
\index{theorem!Groenewold--Van Hove} asserts that no quantisation scheme can
map the Poisson bracket {\em exactly} to the commutator on the full classical
Poisson algebra of observables. In its standard polynomial formulation, one
considers the Poisson algebra of polynomial functions on $\R^{2d}$; under
the usual normalisation and irreducibility requirements, and with the
standard quantisations of the coordinate and momentum functions fixed, there
is no linear assignment $a\mapsto {\rm Op}_h(a)$ satisfying
$$
        ih\,{\rm Op}_h(\{a,b\})
        =
        [{\rm Op}_h(a),{\rm Op}_h(b)]
$$
for all polynomial observables $a,b$.

Treatments of second quantisation can be found in \cite{Jan, Par, Simon74}.
For a discussion from the physicist's point of view we recommend \cite{Tala}.
A different proof of Theorem \ref{thm:sec-quant-pos} is given in \cite{Simon74}.
Theorems \ref{thm:Segal} and \ref{thm:WT} are due to \cite{Seg}.

There is an interesting connection of second quantisation with number theory.
Heuristically, it says that, spectrally, the positive integers are
precisely the second quantised primes. The starting point for making this
rigorous is a theorem of \cite{BP}: if $T$ is a bounded operator on a
Hilbert space $H$, then the spectrum of its $n$-fold tensor product
$T^{\ot n}$ acting on the Hilbert space $H^{\ot n}$ equals
$$ \sigma(T^{\ot n}) = \{\la_1\cdots\la_n:\,  \la_j\in \sigma(T), \ j=1,\dots,n\}.$$
If $\n T\n<1$, taking direct sums gives the formula
$$ \sigma\Bigl(\bigoplus_{n\in\N}T^{\ot n}\Bigr) = \{1\}\cup\overline{\bigcup_{n\ge1}\bigl\{ \la_1\cdots\la_n: \, \la_j\in \sigma(T),\ 1\le j\le n\bigr\}}$$
with the contribution $1$ for the spectrum of $T^{\ot 0}:= I$.
Now let $\mathbb{P} = \{2,3,5,7,11,\dots\}$ be the set of primes and consider the Hilbert space $\ell^2(\mathbb{P})$. Denoting the standard basis vectors of this space by $e_2,e_3,e_5,\dots$, we consider the contraction
$$T: e_{p} \mapsto p^{-1} e_p, \quad p\in\mathbb{P}.$$
Then $\n T\n=\frac12$ and
$\sigma(T) = \{0\} \cup \{\frac1p:\, p\in \mathbb{P}\}$, and accordingly
$$ \sigma\Bigl(\bigoplus_{n\in\N}T^{\ot n}\Bigr) = \{0\}\cup\Bigl\{\frac1n:\, n\in \N, \, n\ge 1\Bigr\},$$
with each point $\frac1n$ being a pole. For symmetric tensor powers, uniqueness of prime factorisation implies that each pole $\frac1n$ occurs with multiplicity one.
This observation, as well as deeper connections, can be found in \cite{BostConnes, Connes}.

Most results in Section \ref{sec:SQ} generalise to infinite dimensions if one replaces the Gaussian measure $\gamma$ on $\R^d$ by a so-called $H$-isonormal process defined on a probability space $(\Om,\P)$, with a real Hilbert space $H$ replacing $\R^d$; see, for example, \cite{Nua}. Here the $L^2$-space is taken with respect to the
$\sigma$-algebra generated by the process.

\section*{Appendices} Most of the material is standard; some proofs are taken from
\cite{Fol, Kal}.
Zorn's lemma is equivalent to the Axiom of Choice. A proof of this fact and further equivalences can be found in \cite{Jech, Jech-set, Rubin}. The proof of Tychonoff's theorem is standard and is presented in classical topology texts such as \cite{Kelley-GeneralTopology}.

Theorem~\ref{thm:Hausdorff-Alexandroff} is a variant of a classical result, often attributed to
Hausdorff and Alexandroff, asserting that every nonempty compact metric space is a continuous
image of the Cantor set $\{0,1\}^{\N}$; see, for instance, \cite{Kelley-GeneralTopology}. Our proof follows a simple `nested covers' construction that produces a compact Cantor-like set $L\subseteq [0,1]$ together with a continuous surjection $L\to X$. Proposition \ref{prop:min-section} is from \cite{FanSong}. The proof of Brouwer's fixed point theorem based on Sperner's lemma follows the classical combinatorial route of
\cite{KnasterKuratowskiMazurkiewicz1929}; see also
\cite{Border}.

The treatment of Carath\'eodory's extension theorem is based on lecture notes by Mark Veraar.

%% file: literature.bib
@book {Adams,
    AUTHOR = {Adams, R. A. and Fournier, J. J. F.},
     TITLE = {Sobolev spaces},
    SERIES = {Pure and Applied Mathematics (Amsterdam)},
    VOLUME = {140},
   EDITION = {2nd},
 PUBLISHER = {Elsevier/Academic Press, Amsterdam},
      YEAR = {2003},
     PAGES = {xiv+305},
   MRCLASS = {46E35 (46-01 46-02 46B70 46Exx)},
  MRNUMBER = {2424078},
}

@book {AkhGla1,
    AUTHOR = {Akhiezer, N. I. and Glazman, I. M.},
     TITLE = {Theory of linear operators in {H}ilbert space. {V}ol. {I}},
    SERIES = {Monographs and Studies in Mathematics},
    VOLUME = {9},
 PUBLISHER = {Pitman, Boston, Mass.-London},
      YEAR = {1981},
     PAGES = {xxxii+312},
   MRCLASS = {47-01 (46-01)},
  MRNUMBER = {615736},
}

@book {AkhGla2,
    AUTHOR = {Akhiezer, N. I. and Glazman, I. M.},
     TITLE = {Theory of linear operators in {H}ilbert space. {V}ol. {II}},
    SERIES = {Monographs and Studies in Mathematics},
    VOLUME = {10},
 PUBLISHER = {Pitman, Boston, Mass.-London},
      YEAR = {1981},
     PAGES = {i--xii and 313--552 and xiii--xxxii},
   MRCLASS = {47-01 (46-01)},
  MRNUMBER = {615737},
}

@article{Alaoglu1940,
  author  = {Alaoglu, L.},
  title   = {Weak topologies of normed linear spaces},
  journal = {Ann. of Math. (2)},
  volume  = {41},
  number  = {1},
  pages   = {252--267},
  year    = {1940},
  publisher = {Ann. Math.},
}

@BOOK{AlbKal,
    Address        = {New York},
    Author         = {Albiac, F. and Kalton, N. J.},
    Pages          = {xii+373},
    Publisher      = {Springer},
    Series         = {Graduate Texts in Mathematics},
    Title          = {{T}opics in {B}anach space theory},
    Volume         = {233},
    Year           = {2006},
    Mrnumber       = {MR2192298 (2006h:46005)}
}

@book {AliBur,
    AUTHOR = {Aliprantis, C. D. and Burkinshaw, O.},
     TITLE = {Positive operators},
    SERIES = {Pure and Applied Mathematics},
    VOLUME = {119},
 PUBLISHER = {Academic Press, Inc., Orlando, FL},
      YEAR = {1985},
     PAGES = {xvi+367},
   MRCLASS = {47B55 (46-01 46A40 46B30 47-01)},
  MRNUMBER = {809372},
MRREVIEWER = {Yu. A. Abramovich},
}

@article {AllRan,
    AUTHOR = {Allan, G. R. and Ransford, T. J.},
     TITLE = {Power-dominated elements in a {B}anach algebra},
   JOURNAL = {Studia Math.},
  FJOURNAL = {Polska Akademia Nauk. Instytut Matematyczny. Studia
              Mathematica},
    VOLUME = {94},
      YEAR = {1989},
    NUMBER = {1},
     PAGES = {63--79},
   MRCLASS = {46H05},
  MRNUMBER = {1008239},
MRREVIEWER = {John Duncan},
}

@book {Ama,
    AUTHOR = {Amann, H.},
     TITLE = {Linear and quasilinear parabolic problems, Vol. {I}: {A}bstract linear theory},
    SERIES = {Monographs in Mathematics},
    VOLUME = {89},
 PUBLISHER = {Birkh\"{a}user Boston, Inc., Boston, MA},
      YEAR = {1995},
     PAGES = {xxxvi+335},
   MRCLASS = {34Gxx (35Kxx 46M35 46N20 47D06 47N20)},
  MRNUMBER = {1345385},
MRREVIEWER = {Paolo Acquistapace},
}

@book{AmbrosioGigliSavare2005,
  author    = {Ambrosio, L. and Gigli, N. and Savar{\'e}, G.},
  title     = {Gradient flows: {I}n metric spaces and in the
               space of probability measures},
  series    = {Lectures in Mathematics ETH Z{\"u}rich},
  publisher = {Birkh{\"a}user},
  address   = {Basel},
  year      = {2005},
  doi       = {10.1007/b137080}
}

@article {AA,
    AUTHOR = {Amemiya, I. and Araki, H.},
     TITLE = {A remark on {P}iron's paper},
   JOURNAL = {Publ. Res. Inst. Math. Sci. Ser. A},
  FJOURNAL = {Kyoto University. Research Institute for Mathematical
              Sciences. Publications},
    VOLUME = {2},
      YEAR = {1966/1967},
     PAGES = {423--427},
   MRCLASS = {06.40 (46.00)},
  MRNUMBER = {0213266},
MRREVIEWER = {P. A. Fillmore},
}

@book {App,
    AUTHOR = {Applebaum, D.},
     TITLE = {Semigroups of linear operators},
    SERIES = {London Mathematical Society Student Texts},
    VOLUME = {93},
 PUBLISHER = {Cambridge University Press, Cambridge},
      YEAR = {2019},
     PAGES = {x+223},
   MRCLASS = {47D03 (60J25)},
  MRNUMBER = {3967735},
MRREVIEWER = {Pierre Portal},
}

@book{Are-ISEM,
    AUTHOR = {Arendt, W.},
     TITLE = {Heat kernels},
      NOTE = {Lecture notes, 9th Internet Seminar 2005/06},
      YEAR = {2006},
}

@incollection {Are04,
    AUTHOR = {Arendt, W.},
 BOOKTITLE = {Evolutionary equations. {V}ol. {I}},
    SERIES = {Handb. Differ. Equ.},
     PAGES = {1--85},
     TITLE = {Semigroups and evolution equations: functional calculus,
              regularity and kernel estimates},
 PUBLISHER = {North-Holland, Amsterdam},
      YEAR = {2004},
   MRCLASS = {47D06 (34G10 34G20 35B65 35C15 35K90 47A60 47N20)},
  MRNUMBER = {2103696},
MRREVIEWER = {Atsushi Yagi},
}

@incollection {Are95,
    AUTHOR = {Arendt, W.},
 BOOKTITLE = {Evolution equations ({B}aton {R}ouge, {LA}, 1992)},
     pages = {21--28},
     TITLE = {Spectrum and growth of positive semigroups},
    SERIES = {Lecture Notes in Pure and Appl. Math.},
    volume = {168},
 publisher = {Dekker, New York},
      year = {1995},
}

@BOOK{ABHN,
    AUTHOR = {Arendt, W. and Batty, C. J. K. and Hieber, M.
              and Neubrander, F.},
     TITLE = {Vector-valued {L}aplace transforms and {C}auchy problems},
    SERIES = {Monographs in Mathematics},
    VOLUME = {96},
   EDITION = {2nd},
 PUBLISHER = {Birkh\"auser/Springer Basel AG, Basel},
      YEAR = {2011},
     PAGES = {xii+539},
   MRCLASS = {47D06 (35A22 35K15 35L15 35S10 44A10 47N20)},
  MRNUMBER = {2798103 (2012b:47109)},
}

@article{ArendtTerElst2012Sectorial,
  author  = {Arendt, W. and ter Elst, A. F. M.},
  title   = {Sectorial forms and degenerate differential operators},
  journal = {J. Operator Theory},
  volume  = {67},
  number  = {1},
  pages   = {33--72},
  year    = {2012}
}

@incollection{ArendtTerElst2012Forms,
  author    = {Arendt, W. and ter Elst, A. F. M.},
  title     = {From forms to semigroups},
  booktitle = {Spectral Theory, Mathematical System Theory,
               Evolution Equations, Differential and Difference Equations},
  series    = {Operator Theory: Advances and Applications},
  volume    = {221},
  pages     = {47--69},
  publisher = {Birkh{\"a}user Basel},
  year      = {2012},
}

@book {AreUrb,
    AUTHOR = {Arendt, W. and Urban, K.},
     TITLE = {Partial differential equations---an introduction to analytical
              and numerical methods},
    SERIES = {Graduate Texts in Mathematics},
    VOLUME = {294},
 PUBLISHER = {Springer, Cham},
      YEAR = {2023},
     PAGES = {xxiv+452},
   MRCLASS = {35-01},
  MRNUMBER = {4628013},
}

@book{ArendtVogtVoigt2026,
  author    = {Arendt, W. and Vogt, H. and Voigt, J.},
  title     = {Form methods for evolution equations},
  series    = {Operator Theory: Advances and Applications},
  volume    = {316},
  publisher = {Birkh{\"a}user},
  address   = {Cham},
  year      = {2026},
  note      = {Forthcoming}
}

@article {Arhancet-curvature,
   AUTHOR = {Arhancet, C.},
    TITLE = {Curvature, {D}olbeault--{D}irac operators, and an {$L^p$}-index theorem on compact {K}\"ahler manifolds},
  year          = {2026},
  journal       = {arXiv:2401.04203v9},
  primaryClass  = {math.FA},
}

@article{ArnoldEtAl2019EffectivePotential,
  author  = {Arnold, D. N. and David, G. and Filoche, M.
             and Jerison, D. and Mayboroda, S.},
  title   = {Localization of eigenfunctions via an effective potential},
  journal = {Comm. Partial Differential Equations},
  volume  = {44},
  number  = {11},
  pages   = {1186--1216},
  year    = {2019},
  doi     = {10.1080/03605302.2019.1626420}
}

@article{ArnoldEtAl2019ComputingSpectra,
  author  = {Arnold, D. N. and David, G. and Filoche, M.
             and Jerison, D. and Mayboroda, S.},
  title   = {Computing spectra without solving eigenvalue problems},
  journal = {SIAM J. Sci. Comput.},
  volume  = {41},
  number  = {1},
  pages   = {B69--B92},
  year    = {2019},
  doi     = {10.1137/17M1156721}
}

@book {Arv,
    AUTHOR = {Arveson, W.},
     TITLE = {A short course on spectral theory},
    SERIES = {Graduate Texts in Mathematics},
    VOLUME = {209},
 PUBLISHER = {Springer, New York},
      YEAR = {2002},
     PAGES = {x+135},
   MRCLASS = {47-01 (46-01 47A10)},
  MRNUMBER = {1865513},
MRREVIEWER = {Douglas R. Farenick},
}

@article{AtiyahSinger1963,
  author  = {Atiyah, M. F. and Singer, I. M.},
  title   = {The index of elliptic operators on compact manifolds},
  journal = {Bull. Amer. Math. Soc.},
  volume  = {69},
  number  = {3},
  pages   = {422--433},
  year    = {1963},
  doi     = {10.1090/S0002-9904-1963-10957-X}
}

@book {AH,
    AUTHOR = {Atkinson, K. and Han, W.},
     TITLE = {Theoretical numerical analysis. {A} functional analysis framework},
    SERIES = {Texts in Applied Mathematics},
    VOLUME = {39},
   EDITION = {3rd},
 PUBLISHER = {Springer, Dordrecht},
      YEAR = {2009},
     PAGES = {xvi+625},
   MRCLASS = {65-01},
  MRNUMBER = {2511061 (2010b:65001)},
}

@book{Attal,
author = {Attal, S.},
title  = {Tensor products and partial traces},
note   = {Lecture notes},
year   = {2013},
}

@book {Aup,
    AUTHOR = {Aupetit, B.},
     TITLE = {A primer on spectral theory},
    SERIES = {Universitext},
 PUBLISHER = {Springer, New York},
      YEAR = {1991},
     PAGES = {xii+193},
   MRCLASS = {46-01 (47-01 47A10)},
  MRNUMBER = {1083349},
MRREVIEWER = {Florian Horia Vasilescu},
}

@ARTICLE{AHLMT02,
    Author         = {Auscher, P. and Hofmann, S. and Lacey, M. and McIntosh, A. and Tchamitchian, Ph.},
    Journal        = {Ann. of Math. (2)},
    Number         = {2},
    Pages          = {633--654},
    Title          = {The solution of the {K}ato square root problem for second order elliptic operators on {${\mathbb{R}}\sp n$}},
    Volume         = {156},
    Year           = {2002},
    Mrnumber       = {MR1933726 (2004c:47096c)}
}

@ARTICLE{AKM06,
    Author         = {Axelsson, A. and Keith, S. and McIntosh, A.},
    Journal        = {Invent. Math.},
    Number         = {3},
    Pages          = {455--497},
    Title          = {Quadratic estimates and functional calculi of perturbed {D}irac operators},
    Volume         = {163},
    Year           = {2006},
    Mrnumber       = {MR2207232 (2007k:58029)}
}

@article{BaezDuarte2003,
  author  = {B{\'a}ez-Duarte, L.},
  title   = {A strengthening of the {N}yman-{B}eurling criterion for the {R}iemann hypothesis},
  journal = {Rend. Lincei Mat. Appl.},
  volume  = {14},
  number  = {1},
  pages   = {5--11},
  year    = {2003},
}

@article{Bag,
 author = {Bagchi, B.},
  title = {On {N}yman, {B}eurling and {B}aez-{D}uarte's {H}ilbert space reformulation of the {R}iemann hypothesis},
journal = {Proc. Indian Acad. Sci. Math. Sci.},
 volume = {116},
  number = {2},
   year = {2006},
  pages = {137--146},
}

@article{Banach1929,
  author  = {Banach, S.},
  title   = {Sur les fonctionnelles lin{\'e}aires {II}},
  journal = {Studia Math.},
  volume  = {1},
  number  = {1},
  pages   = {223--239},
  year    = {1929},
  doi     = {10.4064/sm-1-1-223-239}
}

@book{Banach1932,
  author = {Banach, S.},
  title = {Th{\'e}orie des op{\'e}rations lin{\'e}aires},
  series = {Monografje Matematyczne},
  volume = {1},
  address = {Warsaw},
  year = {1932}
}

@article {Bar,
    AUTHOR = {Bargmann, V.},
     TITLE = {On unitary ray representations of continuous groups},
   JOURNAL = {Ann. of Math. (2)},
  FJOURNAL = {Annals of Mathematics. Second Series},
    VOLUME = {59},
      YEAR = {1954},
     PAGES = {1--46},
   MRCLASS = {20.0X},
  MRNUMBER = {58601},
MRREVIEWER = {I. E. Segal},
}

@BOOK{Beau,
    AUTHOR = {Beauzamy, B.},
     TITLE = {Introduction to operator theory and invariant subspaces},
    SERIES = {North-Holland Mathematical Library},
 PUBLISHER = {North-Holland},
      YEAR = {1988},
}

@article {Beck75,
    AUTHOR = {Beckner, W.},
     TITLE = {Inequalities in {F}ourier analysis},
   JOURNAL = {Ann. of Math. (2)},
  FJOURNAL = {Annals of Mathematics. Second Series},
    VOLUME = {102},
      YEAR = {1975},
    NUMBER = {1},
     PAGES = {159--182},
   MRCLASS = {42A92 (42A18)},
  MRNUMBER = {MR0385456 (52 \#6317)},
MRREVIEWER = {I. I. Hirschman, Jr.},
}

@article{Beneduci,
  author  = {Beneduci, R.},
  title   = {Positive operator valued measures and {F}eller {M}arkov kernels},
  journal = {J. Math. Anal. Appl.},
  volume  = {442},
  number  = {1},
  year    = {2016},
  pages   = {50--71},
}

@book {Berberian,
    AUTHOR = {Berberian, S. K.},
     TITLE = {Notes on spectral theory},
    SERIES = {Van Nostrand Mathematical Studies, No. 5},
 PUBLISHER = {D. Van Nostrand Co., Inc., Princeton, NJ-Toronto,
              Ont.-London},
      YEAR = {1966},
     PAGES = {iii+121},
   MRCLASS = {47.30},
  MRNUMBER = {0190760},
MRREVIEWER = {J. L. B. Cooper},
      NOTE = {author's 2nd ed.},
}

@article{Beurling1955,
 AUTHOR = {Beurling, A.},
 TITLE = {On a closure problem related to the {R}iemann zeta-function},
 JOURNAL = {Proc. Natl. Acad. Sci. USA},
 VOLUME = {41},
 NUMBER = {5},
 YEAR = {1955},
 PAGES = {312--314},
}

@article {BieWar,
    AUTHOR = {Biegert, M. and Warma, M.},
     TITLE = {Regularity in capacity and the {D}irichlet {L}aplacian},
   JOURNAL = {Potential Anal.},
  FJOURNAL = {Potential Analysis. An International Journal Devoted to the
              Interactions between Potential Theory, Probability Theory,
              Geometry and Functional Analysis},
    VOLUME = {25},
      YEAR = {2006},
    NUMBER = {3},
     PAGES = {289--305},
   MRCLASS = {35J67 (31C15 35J05)},
  MRNUMBER = {2255350},
MRREVIEWER = {Jana Bj\"{o}rn},
}

@unpublished{Bir,
 AUTHOR = {Birman, M. Sh.},
  TITLE = {A proof of the {F}redholm trace formula as an application of a simple embedding
  for kernels of integral operators of trace class in {$L^2(\R^m)$}},
  NOTE = {Report LIT-MATK-R-89-30, Department of Mathematics, Link\"oping University},
  YEAR = {1989},
}

@book {BirSol,
    AUTHOR = {Birman, M. Sh. and Solomjak, M. Z.},
     TITLE = {Spectral theory of selfadjoint operators in {H}ilbert space},
    SERIES = {Mathematics and its Applications (Soviet Series)},
 PUBLISHER = {D. Reidel Publishing Co., Dordrecht},
      YEAR = {1987},
     PAGES = {xv+301},
   MRCLASS = {47-01 (47A10 47B25)},
  MRNUMBER = {1192782},
}

@book {Black,
    AUTHOR = {Blackadar, B.},
     TITLE = {Operator algebras: {T}heory of {$C^*$}-algebras and von {N}eumann algebras},
    SERIES = {Encyclopaedia of Mathematical Sciences},
    VOLUME = {122},
 PUBLISHER = {Springer, Berlin},
      YEAR = {2006},
     PAGES = {xx+517},
   MRCLASS = {46L05 (46L10 46L80)},
  MRNUMBER = {2188261},
MRREVIEWER = {Paul Jolissaint},
}

@book {BleeBoo,
    AUTHOR = {Bleecker, D. D. and Boo\ss -Bavnbek, B.},
     TITLE = {Index theory---with applications to mathematics and physics},
 PUBLISHER = {International Press, Somerville, MA},
      YEAR = {2013},
     PAGES = {xxii+769},
     MRCLASS = {58-02 (19L10 35S99 57R57 58J20 70S15)},
  MRNUMBER = {3113540},
MRREVIEWER = {Peter Haskell},
}

@book {Bog1,
    AUTHOR = {Bogachev, V. I.},
     TITLE = {Measure theory. {V}ol. {1}},
 PUBLISHER = {Springer, Berlin},
      YEAR = {2007},
   MRCLASS = {28-02 (28Axx 28Cxx 46G12 60G42 60G44)},
  MRNUMBER = {2267655},
MRREVIEWER = {Ren\'{e} L. Schilling},
}

@book {Bog2,
    AUTHOR = {Bogachev, V. I.},
     TITLE = {Measure theory. {V}ol. {2}},
 PUBLISHER = {Springer, Berlin},
      YEAR = {2007},
   MRCLASS = {28-02 (28Axx 28Cxx 46G12 60G42 60G44)},
  MRNUMBER = {2267655},
MRREVIEWER = {Ren\'{e} L. Schilling},
}

@article{BohnenblustSobczyk1938,
  author = {Bohnenblust, H. F. and Sobczyk, A.},
  title = {Extensions of functionals on complex linear spaces},
  journal = {Bull. Amer. Math. Soc.},
  volume = {44},
  number = {2},
  pages = {91--93},
  year = {1938}
}

@book{Border,
  author    = {Border, K. C.},
  title     = {Fixed point theorems with applications to economics and game theory},
  publisher = {Cambridge University Press},
  address   = {Cambridge},
  year      = {1985},
}

@article {BostConnes,
    AUTHOR = {Bost, J.-B. and Connes, A.},
     TITLE = {Hecke algebras, type {III} factors and phase transitions with
              spontaneous symmetry breaking in number theory},
   JOURNAL = {Selecta Math.},
  FJOURNAL = {Selecta Mathematica. New Series},
    VOLUME = {1},
      YEAR = {1995},
    NUMBER = {3},
     PAGES = {411--457},
   MRCLASS = {46L35 (11M06 11R18 11Z05 20G30 82B10 82B26)},
  MRNUMBER = {1366621},
MRREVIEWER = {Alain Valette},
}

@book {BotSil,
    AUTHOR = {B\"{o}ttcher, A. and Silbermann, B.},
     TITLE = {Analysis of {T}oeplitz operators},
    SERIES = {Springer Monographs in Mathematics},
   EDITION = {2nd},
 PUBLISHER = {Springer, Berlin},
      YEAR = {2006},
     PAGES = {xiv+665},
   MRCLASS = {47-02 (47B35 47G10 47L80)},
  MRNUMBER = {2223704},
MRREVIEWER = {Estelle L. Basor},
}

@book {BraRob,
    AUTHOR = {Bratteli, O. and Robinson, D. W.},
     TITLE = {Operator algebras and quantum statistical mechanics. {V}ol. 1},
    SERIES = {Texts and Monographs in Physics},
   EDITION = {2nd},
 PUBLISHER = {Springer, New York},
      YEAR = {1987},
     PAGES = {xiv+505},
   MRCLASS = {46Lxx (81-02 82-02)},
  MRNUMBER = {887100},
}

@article{BravermanEtAl2013,
  author  = {Braverman, M. and Makarychev, K. and Makarychev, Yu. and Naor, A.},
  title   = {The {G}rothendieck constant is strictly smaller than {K}rivine's bound},
  journal = {Forum Math. Pi},
  volume  = {1},
  pages   = {e4},
  year    = {2013},
  doi     = {10.1017/fmp.2013.4}
}

@book {Br,
    AUTHOR = {Brenner, S. C. and Scott, L. R.},
     TITLE = {The mathematical theory of finite element methods},
    SERIES = {Texts in Applied Mathematics},
    VOLUME = {15},
   EDITION = {3rd},
 PUBLISHER = {Springer, New York},
      YEAR = {2008},
   MRCLASS = {65-01 (65-02)},
  MRNUMBER = {2373954 (2008m:65001)},
}

@book {Bressan,
    AUTHOR = {Bressan, A.},
     TITLE = {Lecture notes on functional analysis},
    SERIES = {Graduate Studies in Mathematics},
    VOLUME = {143},
 PUBLISHER = {Amer. Math. Soc., Providence, RI},
      YEAR = {2013},
   MRCLASS = {46-01 (35-01)},
  MRNUMBER = {2987297},
MRREVIEWER = {Valentin Keyantuo},
}

@book {Bre,
    AUTHOR = {Brezis, H.},
     TITLE = {Functional analysis, {S}obolev spaces and partial differential
              equations},
    SERIES = {Universitext},
 PUBLISHER = {Springer, New York},
      YEAR = {2011},
   MRCLASS = {35-01 (46-01 46E35 46N20 47F05)},
  MRNUMBER = {2759829 (2012a:35002)},
MRREVIEWER = {Vicen{\c{t}}iu R{\u{a}}dulescu},
}

@book{Brezis1973,
  author    = {Br{\'e}zis, H.},
  title     = {Op{\'e}rateurs maximaux monotones et semi-groupes
               de contractions dans les espaces de {Hilbert}},
  series    = {North-Holland Mathematics Studies},
  volume    = {5},
  publisher = {North-Holland},
  address   = {Amsterdam},
  year      = {1973}
}

@misc{Briet2018Grothendieck,
  author       = {Bri{\"e}t, J.},
  title        = {Lecture notes on applications of
                  {Grothendieck}'s inequality: The inequality},
  year         = {2018},
  howpublished = {Lecture notes, NETWORKS Training Week 6,
                  Doorn, 9--13 April 2018},
  url          = {https://www.thenetworkcenter.nl/uploaded_files/inlineitem/JB_grothendieck_proof.pdf}
}

@article {Bris,
    AUTHOR = {Brislawn, C.},
     TITLE = {Kernels of trace class operators},
   JOURNAL = {Proc. Amer. Math. Soc.},
  FJOURNAL = {Proceedings of the American Mathematical Society},
    VOLUME = {104},
      YEAR = {1988},
    NUMBER = {4},
     PAGES = {1181--1190},
   MRCLASS = {47B38 (47B10)},
  MRNUMBER = {929421},
MRREVIEWER = {V. S. Sunder},
}

@article {BP,
    AUTHOR = {Brown, A. and Pearcy, C.},
     TITLE = {Spectra of tensor products of operators},
   JOURNAL = {Proc. Amer. Math. Soc.},
  FJOURNAL = {Proceedings of the American Mathematical Society},
    VOLUME = {17},
      YEAR = {1966},
     PAGES = {162--166},
   MRCLASS = {47.30},
  MRNUMBER = {188786},
MRREVIEWER = {H. H. Schaefer},
}

@book {BGL,
    AUTHOR = {Busch, P. and Grabowski, M. and Lahti, P. J.},
     TITLE = {Operational quantum physics},
    SERIES = {Lecture Notes in Physics. New Series m: Monographs},
    VOLUME = {31},
 PUBLISHER = {Springer, Berlin},
      YEAR = {1995},
   MRCLASS = {81P15 (47N50 81-02 81P05)},
  MRNUMBER = {1356220},
MRREVIEWER = {I. Gottlieb},
}

@book{BLPY,
  author    = {P. Busch and P. Lahti and J.-P. Pellonp{\"a}{\"a} and K. Ylinen},
  title     = {Quantum measurement},
  series    = {Theoretical and Mathematical Physics},
  publisher = {Springer},
  address   = {Cham},
  year      = {2016},
}

@article {Buskes,
    AUTHOR = {Buskes, G.},
     TITLE = {The {H}ahn-{B}anach theorem surveyed},
   JOURNAL = {Dissertationes Math.},
  FJOURNAL = {Dissertationes Mathematicae (Rozprawy Matematyczne)},
    VOLUME = {327},
      YEAR = {1993},
     PAGES = {49},
   MRCLASS = {46A22 (04A25 46-02 46M10 46S10 47B60)},
  MRNUMBER = {1240598},
MRREVIEWER = {W. A. J. Luxemburg},
}

@article {Carleman,
    AUTHOR = {Carleman, T.},
     TITLE = {\"{U}ber die {F}ourierkoeffizienten einer stetigen {F}unktion},
   JOURNAL = {Acta Math.},
  FJOURNAL = {Acta Mathematica},
    VOLUME = {41},
      YEAR = {1916},
    NUMBER = {1},
     PAGES = {377--384},
   MRCLASS = {DML},
  MRNUMBER = {1555157},
}

@incollection {Cartier,
    AUTHOR = {Cartier, P.},
 BOOKTITLE = {Conformal invariance and string theory ({P}oiana {B}ra\c{s}ov,
              1987)},
    SERIES = {Perspect. Phys.},
     PAGES = {443--557},
     TITLE = {A course on determinants},
 PUBLISHER = {Academic Press, Boston, MA},
      YEAR = {1989},
   MRCLASS = {58G26 (15A15 47B10 58C50 82B10)},
  MRNUMBER = {1010666},
MRREVIEWER = {Steven Rosenberg},
}

@article{ChillSeifert2016,
  author  = {Chill, R. and Seifert, D.},
  title   = {Quantified versions of {I}ngham's theorem},
  journal = {Bull. Lond. Math. Soc.},
  volume  = {48},
  number  = {3},
  pages   = {519--532},
  year    = {2016},
}

@article{Chiribella-DAP,
  author  = {Chiribella, G. and D'Ariano, G. M. and Perinotti, P.},
  title   = {Informational derivation of quantum theory},
  journal = {Phys. Rev. A},
  volume  = {84},
  year    = {2011},
  pages   = {012311},
  doi     = {10.1103/PhysRevA.84.012311},
}

@article{ClauserHorneShimonyHolt1969,
  author  = {Clauser, J. F. and Horne, M. A. and
             Shimony, A. and Holt, R. A.},
  title   = {Proposed experiment to test local hidden-variable theories},
  journal = {Phys. Rev. Lett.},
  volume  = {23},
  number  = {15},
  year    = {1969},
  pages   = {880--884},
  doi     = {10.1103/PhysRevLett.23.880}
}

@article {CleLi,
    AUTHOR = {Cl\'ement, Ph. and Li, S.},
     TITLE = {Abstract parabolic quasilinear equations and application to a
              groundwater flow problem},
   JOURNAL = {Adv. Math. Sci. Appl.},
  FJOURNAL = {Advances in Mathematical Sciences and Applications},
    VOLUME = {3\,},
      YEAR = {1993/94},
    NUMBER = {Special Issue},
     PAGES = {17--32},
   MRCLASS = {34G20 (35J60 35K55 47F05 47N20 76S05)},
  MRNUMBER = {1287921},
MRREVIEWER = {Paolo Acquistapace},
}

@article {Cob1,
    AUTHOR = {Coburn, L. A.},
     TITLE = {Weyl's theorem for nonnormal operators},
   JOURNAL = {Michigan Math. J.},
  FJOURNAL = {Michigan Mathematical Journal},
    VOLUME = {13},
      YEAR = {1966},
     PAGES = {285--288},
   MRCLASS = {47.10},
  MRNUMBER = {201969},
MRREVIEWER = {H. Widom},
}

@article {Cob2,
    AUTHOR = {Coburn, L. A.},
     TITLE = {The {$C^{\ast} $}-algebra generated by an isometry},
   JOURNAL = {Bull. Amer. Math. Soc.},
  FJOURNAL = {Bulletin of the American Mathematical Society},
    VOLUME = {73},
      YEAR = {1967},
     PAGES = {722--726},
   MRCLASS = {47.40 (46.00)},
  MRNUMBER = {213906},
MRREVIEWER = {Bernard Russo},
}

@book {Connes,
    AUTHOR = {Connes, A.},
     TITLE = {Noncommutative geometry},
 PUBLISHER = {Academic Press, Inc., San Diego, CA},
      YEAR = {1994},
     PAGES = {xiv+661},
   MRCLASS = {46Lxx (19K56 22D25 58B30 58G12 81T13 81V22 81V70)},
  MRNUMBER = {1303779},
MRREVIEWER = {John Roe},
}

@article {ConCon,
    AUTHOR = {Connes, A. and Consani, C.},
     TITLE = {The scaling {H}amiltonian},
   JOURNAL = {J. Operator Theory},
  FJOURNAL = {Journal of Operator Theory},
    VOLUME = {85},
      YEAR = {2021},
    NUMBER = {1},
     PAGES = {257--276},
   MRCLASS = {11M55 (46L87 58B34)},
  MRNUMBER = {4198972},
}

@book {Conway,
    AUTHOR = {Conway, J. B.},
     TITLE = {A course in functional analysis},
    SERIES = {Graduate Texts in Mathematics},
    VOLUME = {96},
   EDITION = {2nd},
 PUBLISHER = {Springer, New York},
      YEAR = {1990},
   MRCLASS = {46-01 (47-01)},
  MRNUMBER = {1070713},
}

@book {Conway-OpTh,
    AUTHOR = {Conway, J. B.},
     TITLE = {A course in operator theory},
    SERIES = {Graduate Studies in Mathematics},
    VOLUME = {21},
 PUBLISHER = {Amer. Math. Soc., Providence, RI},
      YEAR = {2000},
     PAGES = {xvi+372},
   MRCLASS = {47-01 (46-01)},
  MRNUMBER = {1721402},
MRREVIEWER = {Victor\ S.\ Shul\cprime man},
}

@article{DavidFilocheMayboroda2021,
  author  = {David, G. and Filoche, M. and Mayboroda, S.},
  title   = {The landscape law for the integrated density of states},
  journal = {Adv. Math.},
  volume  = {390},
  pages   = {107946},
  year    = {2021},
  doi     = {10.1016/j.aim.2021.107946}
}

@book {Davies-sgr,
    AUTHOR = {Davies, E. B.},
     TITLE = {One-parameter semigroups},
    SERIES = {London Mathematical Society Monographs},
    VOLUME = {15},
 PUBLISHER = {Academic Press, Inc.,
              London-New York},
      YEAR = {1980},
   MRCLASS = {47D05},
  MRNUMBER = {591851},
MRREVIEWER = {S. Kurepa},
}

@book {DaviesOpen,
    AUTHOR = {Davies, E. B.},
     TITLE = {Quantum theory of open systems},
 PUBLISHER = {Academic Press, London-New York},
      YEAR = {1976},
   MRCLASS = {81.00 (46N05)},
  MRNUMBER = {0489429},
MRREVIEWER = {H. Araki},
}

@book {Diestel,
    AUTHOR = {Diestel, J.},
     TITLE = {Sequences and series in {B}anach spaces},
    SERIES = {Graduate Texts in Mathematics},
    VOLUME = {92},
 PUBLISHER = {Springer, New York},
      YEAR = {1984},
   MRCLASS = {46B15 (04A15 04A20 46-01 46B20)},
  MRNUMBER = {737004},
MRREVIEWER = {J. Bourgain},
}

@BOOK{DieUhl,
    Address        = {Providence, RI},
    Author         = {Diestel, J. and Uhl, Jr., J. J.},
    Pages          = {xiii+322},
    Publisher      = {Amer. Math. Soc.},
    Title          = {{V}ector measures},
    Year           = {1977},
    Mrnumber       = {MR0453964 (56 \#12216)}
}

@book {Dieu,
    AUTHOR = {Dieudonn{\'e}, J.},
     TITLE = {History of functional analysis},
    SERIES = {North-Holland Mathematics Studies},
    VOLUME = {49},
 PUBLISHER = {North-Holland Publishing Co., Amsterdam-New York},
      YEAR = {1981},
   MRCLASS = {46-03 (01A60)},
  MRNUMBER = {605488 (83d:46001)},
MRREVIEWER = {R. G. Bartle},
}

@book {Douglas,
    AUTHOR = {Douglas, R. G.},
     TITLE = {Banach algebra techniques in operator theory},
    SERIES = {Graduate Texts in Mathematics},
    VOLUME = {179},
   EDITION = {2nd},
 PUBLISHER = {Springer, New York},
      YEAR = {1998},
   MRCLASS = {47-01 (46-01)},
  MRNUMBER = {1634900},
}

@book {vDul,
    AUTHOR = {{\VAN{Dulst}{van}{van}} Dulst, D.},
     TITLE = {Characterizations of {B}anach spaces not containing {$\ell{^1}$}},
    SERIES = {CWI Tract},
    VOLUME = {59},
 PUBLISHER = {CWI},
   ADDRESS = {Amsterdam},
      YEAR = {1989},
     PAGES = {iv+163},
   MRCLASS = {46B25 (46B20 46G10)},
  MRNUMBER = {1002733 (90h:46037)},
MRREVIEWER = {Giovanni Emmanuele},
}

@book {DunSch1,
    AUTHOR = {Dunford, N. and Schwartz, J. T.},
     TITLE = {Linear operators. {P}art {I}: General theory},
    SERIES = {Wiley Classics Library},
      NOTE = {Reprint of the 1958 original},
 PUBLISHER = {John Wiley \& Sons, Inc., New York},
      YEAR = {1988},
   MRCLASS = {47-01 (46-01)},
  MRNUMBER = {1009162},
}

@book {DunSch2,
    AUTHOR = {Dunford, N. and Schwartz, J. T.},
     TITLE = {Linear operators. {P}art {II}: Spectral theory, selfadjoint operators in {H}ilbert space},
    SERIES = {Wiley Classics Library},
      NOTE = {Reprint of the 1963 original},
 PUBLISHER = {John Wiley \& Sons, Inc., New York},
      YEAR = {1988},
   MRCLASS = {47-01 (46-01)},
  MRNUMBER = {1009163},
}

@book {DunSch3,
    AUTHOR = {Dunford, N. and Schwartz, J. T.},
     TITLE = {Linear operators. {P}art {III}: Spectral operators},
    SERIES = {Wiley Classics Library},
      NOTE = {Reprint of the 1971 original},
 PUBLISHER = {John Wiley \& Sons, Inc., New York},
      YEAR = {1988},
}

@book {Duo,
    AUTHOR = {Duoandikoetxea, J.},
     TITLE = {Fourier analysis},
    SERIES = {Graduate Studies in Mathematics},
    VOLUME = {29},
      NOTE = {Translated and revised from the 1995 Spanish original},
 PUBLISHER = {Amer. Math. Soc., Providence, RI},
      YEAR = {2001},
     PAGES = {xviii+222},
   MRCLASS = {42-01},
  MRNUMBER = {1800316},
MRREVIEWER = {Loukas Grafakos},
}

@book {EdmEva-Elliptic,
    AUTHOR = {Edmunds, D. E. and Evans, W. D.},
     TITLE = {Elliptic differential operators and spectral analysis},
    SERIES = {Springer Monographs in Mathematics},
 PUBLISHER = {Springer, Cham},
      YEAR = {2018},
     PAGES = {xiii+322},
   MRCLASS = {35-02 (35Jxx 35Pxx 37A10 47B25)},
  MRNUMBER = {3931747},
MRREVIEWER = {G\"{u}nter Berger},
}

@book {EdmEva,
    AUTHOR = {Edmunds, D. E. and Evans, W. D.},
     TITLE = {Spectral theory and differential operators},
    SERIES = {Oxford Mathematical Monographs},
   EDITION = {2nd},
 PUBLISHER = {Oxford University Press, Oxford},
      YEAR = {2018},
     PAGES = {xviii+589},
   MRCLASS = {47-02 (34Lxx 35Pxx 47Axx 47Bxx 47F05)},
  MRNUMBER = {3823299},
}

@article{Edwards2010,
  author  = {Edwards, D. A.},
  title   = {A simple proof in {Monge--Kantorovich} duality theory},
  journal = {Studia Math.},
  volume  = {200},
  number  = {1},
  pages   = {67--77},
  year    = {2010},
  doi     = {10.4064/sm200-1-4}
}

@article{Edwards2011,
  author  = {Edwards, D. A.},
  title   = {On the {Kantorovich--Rubinstein} theorem},
  journal = {Expo. Math.},
  volume  = {29},
  number  = {4},
  pages   = {387--398},
  year    = {2011},
  doi     = {10.1016/j.exmath.2011.06.005}
}

@book {EinWar,
    AUTHOR = {Einsiedler, M. and Ward, Th.},
     TITLE = {Functional analysis, spectral theory, and applications},
    SERIES = {Graduate Texts in Mathematics},
    VOLUME = {276},
 PUBLISHER = {Springer, Cham},
      YEAR = {2017},
     PAGES = {xiv+614},
   MRCLASS = {46-01 (35J25 35P10 35P20 37A99 47-01 47A60)},
  MRNUMBER = {3729416},
MRREVIEWER = {Michael M. Neumann},
}

@article{ElstSauterVogt2015,
  author  = {{\VAN{Elst}{ter}{ter}} Elst, A. F. M. and Sauter, M. and Vogt, H.},
  title   = {A generalisation of the form method for accretive forms
             and operators},
  journal = {J. Funct. Anal.},
  volume  = {269},
  number  = {3},
  pages   = {705--744},
  year    = {2015},
}

@BOOK{EngNag,
    Address        = {New York},
    Author         = {Engel, K.-J. and Nagel, R.},
    Pages          = {xxii+586},
    Publisher      = {Springer},
    Series         = {Graduate Texts in Mathematics},
    Title          = {{O}ne-parameter semigroups for linear evolution equations},
    Volume         = {194},
    Year           = {2000},
    Mrnumber       = {MR1721989 (2000i:47075)}
}

@article {Epp,
    AUTHOR = {Epperson, J. B.},
     TITLE = {The hypercontractive approach to exactly bounding an operator
              with complex {G}aussian kernel},
   JOURNAL = {J. Funct. Anal.},
  FJOURNAL = {Journal of Functional Analysis},
    VOLUME = {87},
      YEAR = {1989},
    NUMBER = {1},
     PAGES = {1--30},
   MRCLASS = {42B20 (47G10)},
  MRNUMBER = {1025881},
}

@book {Evans,
    AUTHOR = {Evans, L. C.},
     TITLE = {Partial differential equations},
    SERIES = {Graduate Studies in Mathematics},
    VOLUME = {19},
   EDITION = {2nd},
 PUBLISHER = {Amer. Math. Soc., Providence, RI},
      YEAR = {2010},
   MRCLASS = {35-01},
  MRNUMBER = {2597943 (2011c:35002)},
MRREVIEWER = {Diego M. Maldonado},
}

@article{FanSong,
  author        = {Fan, L. and Song, S.},
  title         = {A new proof of the {W}eyl--von {N}eumann--{B}erg theorem},
  year          = {2025},
  journal       = {arXiv:2502.19028},
  primaryClass  = {math.FA},
}

@article{FilocheMayboroda2012,
  author  = {Filoche, M. and Mayboroda, S.},
  title   = {Universal mechanism for {Anderson} and weak localization},
  journal = {Proc. Natl. Acad. Sci. USA},
  volume  = {109},
  number  = {37},
  pages   = {14761--14766},
  year    = {2012},
  doi     = {10.1073/pnas.1120432109}
}

@article{Fil,
  title = {On an inequality between {D}irichlet and {N}eumann eigenvalues for the {L}aplace operator},
 author = {Filonov, N.},
journal = {St. Petersburg Math. J.},
 volume = {16},
 number = {2},
  pages = {413--416},
   year = {2005}
}

@article{Fine1982,
  author  = {Fine, A.},
  title   = {Hidden variables, joint probability, and the {Bell} inequalities},
  journal = {Phys. Rev. Lett.},
  volume  = {48},
  number  = {5},
  year    = {1982},
  pages   = {291--295},
  doi     = {10.1103/PhysRevLett.48.291}
}

@article {Foguel,
    AUTHOR = {Foguel, S. R.},
     TITLE = {On a theorem by {A}. {E}. {T}aylor},
   JOURNAL = {Proc. Amer. Math. Soc.},
  FJOURNAL = {Proceedings of the American Mathematical Society},
    VOLUME = {9},
      YEAR = {1958},
     PAGES = {325},
   MRCLASS = {46.00},
  MRNUMBER = {93696},
MRREVIEWER = {A. E. Taylor},
}

@book{Fol,
    AUTHOR = {Folland, G. B.},
     TITLE = {Real analysis},
    SERIES = {Pure and Applied Mathematics (New York)},
   EDITION = {2nd},
 PUBLISHER = {John Wiley \& Sons Inc.},
   ADDRESS = {New York},
      YEAR = {1999},
   MRCLASS = {00A05 (26-01 28-01 46-01)},
  MRNUMBER = {1681462 (2000c:00001)},
}

@book {Fol-AHA,
    AUTHOR = {Folland, G. B.},
     TITLE = {A course in abstract harmonic analysis},
    SERIES = {Textbooks in Mathematics},
   EDITION = {2nd},
 PUBLISHER = {CRC Press, Boca Raton, FL},
      YEAR = {2016},
   MRCLASS = {43-01 (22-01 42-01 46-01)},
}

@book {Foll-PhaseSpace,
    AUTHOR = {Folland, G. B.},
     TITLE = {Harmonic analysis in phase space},
    SERIES = {Annals of Mathematics Studies},
    VOLUME = {122},
 PUBLISHER = {Princeton University Press, Princeton, NJ},
      YEAR = {1989},
   MRCLASS = {22E30 (43A80 58G15 81S30)},
  MRNUMBER = {983366},
}

@book {Fremlin,
    AUTHOR = {Fremlin, D. H.},
     TITLE = {Measure theory. {V}ol. 5. {S}et-theoretic measure theory.
              {P}art {I}},
      NOTE = {Corrected reprint of the 2008 original},
 PUBLISHER = {Torres Fremlin, Colchester},
      YEAR = {2015},
   MRCLASS = {03-02 (03-01 03E15 03E70 28-01 28-02 28A05 28C15)},
  MRNUMBER = {3723040},
MRREVIEWER = {Klaas Pieter Hart},
}

@article{Fri,
  title = {Some inequalities between {D}irichlet and {N}eumann eigenvalues},
 author = {Friedlander, L.},
journal = {Arch. Ration. Mech. Anal.},
 volume = {116},
 number = {2},
  pages = {153--160},
   year = {1991},
}

@article {Fuhrman,
    AUTHOR = {Fuhrman, M.},
     TITLE = {Analyticity of transition semigroups and closability of
              bilinear forms in {H}ilbert spaces},
   JOURNAL = {Studia Math.},
  FJOURNAL = {Studia Mathematica},
    VOLUME = {115},
      YEAR = {1995},
    NUMBER = {1},
     PAGES = {53--71},
   MRCLASS = {47D07 (47B99 47N30 60J35)},
  MRNUMBER = {1347432},
MRREVIEWER = {J.\ A.\ van Casteren},
}

@article{Funano,
  author = {Funano, K.},
   title = {A note on domain monotonicity for the {N}eumann eigenvalues of the {L}aplacian},
 journal = {Illinois J. Math.},
  volume = {67},
  number = {4},
  pages = {677--686},
   year = {2023},
}

@article {GarWong,
    AUTHOR = {Garrison, J. C. and Wong, J.},
     TITLE = {Canonically conjugate pairs, uncertainty relations, and phase
              operators},
   JOURNAL = {J. Math. Phys.},
  FJOURNAL = {Journal of Mathematical Physics},
    VOLUME = {11},
      YEAR = {1970},
     PAGES = {2242--2249},
   MRCLASS = {81.46},
  MRNUMBER = {266532},
MRREVIEWER = {H. Araki},
}

@article {Gel,
    AUTHOR = {Gelfand, I. M.},
     TITLE = {Zur {T}heorie der {C}haraktere der {A}belschen topologischen
              {G}ruppen},
   JOURNAL = {Mat. Sb. (N.S.)},
    VOLUME = {51},
     number = {9},
      YEAR = {1941},
     PAGES = {49--50},
   MRCLASS = {20.0X},
  MRNUMBER = {0004635},
MRREVIEWER = {S. Bochner},
}

@book {GilTru,
    AUTHOR = {Gilbarg, D. and Trudinger, N. S.},
     TITLE = {Elliptic partial differential equations of second order},
    SERIES = {Classics in Mathematics},
      NOTE = {Reprint of the 1998 edition},
 PUBLISHER = {Springer, Berlin},
      YEAR = {2001},
     PAGES = {xiv+517},
   MRCLASS = {35-02 (35Jxx)},
  MRNUMBER = {1814364},
}

@book{GGK1,
  title={Basic classes of linear operators},
  author={Goh\-berg, I. and Goldberg, S. and Kaashoek, M. {A.}},
  year={2003},
  publisher={Birkh{\"a}user Verlag, Basel}
}

@article {GWW,
    AUTHOR = {Gordon, C. and Webb, D. L. and Wolpert, S.},
     TITLE = {One cannot hear the shape of a drum},
   JOURNAL = {Bull. Amer. Math. Soc. (N.S.)},
  FJOURNAL = {American Mathematical Society. Bulletin. New Series},
    VOLUME = {27},
      YEAR = {1992},
    NUMBER = {1},
     PAGES = {134--138},
   MRCLASS = {58G25 (35R30)},
  MRNUMBER = {1136137},
MRREVIEWER = {Robert Brooks},
}

@book {GraI,
    AUTHOR = {Grafakos, L.},
     TITLE = {Classical {F}ourier analysis},
    SERIES = {Graduate Texts in Mathematics},
    VOLUME = {249},
   EDITION = {2nd},
 PUBLISHER = {Springer},
   ADDRESS = {New York},
      YEAR = {2008},
   MRCLASS = {42-01 (42Bxx)},
  MRNUMBER = {2445437 (2011c:42001)},
MRREVIEWER = {Andreas Seeger},
}

@article{GreNgu,
  title = {Geometrical structure of {L}aplacian eigenfunctions},
 author = {Grebenkov, D. S. and Nguyen, B.-T.},
journal = {SIAM Rev.},
 volume = {55},
 number = {4},
  pages = {601--667},
   year = {2013},
}

@article {Grieser,
    AUTHOR = {Grieser, D.},
     TITLE = {\"{U}ber {E}igenwerte, {I}ntegrale und {$\frac{\pi^2}{6}$}: die
              {I}dee der {S}purformel},
   JOURNAL = {Math. Semesterber.},
  FJOURNAL = {Mathematische Semesterberichte},
    VOLUME = {54},
      YEAR = {2007},
    NUMBER = {2},
     PAGES = {199--217},
   MRCLASS = {58J50 (11M36 35P10 47G10)},
  MRNUMBER = {2342260},
MRREVIEWER = {Daniel Belti\c{t}\u{a}},
}

@article{Grothendieck1953,
  author  = {Grothendieck, A.},
  title   = {R{\'e}sum{\'e} de la th{\'e}orie m{\'e}trique des produits
             tensoriels topologiques},
  journal = {Bol. Soc. Mat. S{\~a}o Paulo},
  volume  = {8},
  year    = {1953},
  pages   = {1--79}
}

@book {Gustafson,
    AUTHOR = {Gustafson, K. E. and Rao, D. K. M.},
     TITLE = {Numerical range},
    SERIES = {Universitext},
 PUBLISHER = {Springer, New York},
      YEAR = {1997},
     PAGES = {xiv+189},
   MRCLASS = {47A12 (15A60 47A10 65J15)},
  MRNUMBER = {1417493},
}

@article {Haase-Convexity,
    AUTHOR = {Haase, M. H. A.},
     TITLE = {Convexity inequalities for positive operators},
   JOURNAL = {Positivity},
    VOLUME = {11},
      YEAR = {2007},
    NUMBER = {1},
     PAGES = {57--68},
   MRCLASS = {39B62 (26D10 39B72 47A63 47B65)},
  MRNUMBER = {2297322 (2008d:39034)},
MRREVIEWER = {B. Choczewski},
}

@article{HaaseKunstmannVogt2016,
  author  = {Haase, M. H. A. and Kunstmann, P. C. and Vogt, H.},
  title   = {On the numerical range of generators of symmetric
             {$L_\infty$}-contractive semigroups},
  journal = {Arch. Math. (Basel)},
  volume  = {107},
  number  = {5},
  pages   = {553--559},
  year    = {2016},
}

@article{Hahn1927,
  author = {Hahn, H.},
  title = {{\"U}ber lineare {G}leichungssysteme in linearen {R}{\"a}umen},
  journal = {J. Reine Angew. Math.},
  volume = {157},
  pages = {214--229},
  year = {1927}
}

@book {Hall,
    AUTHOR = {Hall, B. C.},
     TITLE = {Quantum theory for mathematicians},
    SERIES = {Graduate Texts in Mathematics},
    VOLUME = {267},
 PUBLISHER = {Springer, New York},
      YEAR = {2013},
     PAGES = {xvi+554},
   MRCLASS = {81-02 (46N50)},
  MRNUMBER = {3112817},
  }

@article {HLLL,
    AUTHOR = {Han, D. and Larson, D. R. and Liu, B. and Liu, R.},
     TITLE = {Operator-valued measures, dilations, and the theory of frames},
   JOURNAL = {Mem. Amer. Math. Soc.},
  FJOURNAL = {Memoirs of the American Mathematical Society},
    VOLUME = {229},
      YEAR = {2014},
    NUMBER = {1075},
   MRCLASS = {46G10 (42C15 47A20 47C15)},
  MRNUMBER = {3186831},
MRREVIEWER = {Xunxiang Guo},
}

@article{HHM,
  title={An improvement of the {K}olmogorov--{R}iesz compactness theorem},
  author={Hanche-Olsen, H. and Holden, H. and Malinnikova, E.},
  journal={Expo. Math.},
  volume={37},
  number={1},
  pages={84--91},
  year={2019},
  publisher={Elsevier}
}

@article{Hedberg1972,
  author  = {Hedberg, L. I.},
  title   = {On certain convolution inequalities},
  journal = {Proc. Amer. Math. Soc.},
  volume  = {36},
  number  = {2},
  pages   = {505--510},
  year    = {1972},
  doi     = {10.1090/S0002-9939-1972-0312232-4}
}

@article{HeilmanJonesMalavolta2026,
  author = {Heilman, S. and Jones, C. and Malavolta, G.},
  title  = {Primal/dual method for the {Grothendieck} constant},
  year   = {2026},
 journal = {arXiv:2610.10477},
  url    = {https://arxiv.org/abs/2610.10477}
}

@inproceedings {HH,
    AUTHOR = {Helton, J. W. and Howe, R. E.},
    SERIES = {Lecture Notes in Math., Vol. 345},
     TITLE = {Integral operators: commutators, traces, index and homology},
 BOOKTITLE = {Proceedings of a {C}onference {O}perator {T}heory ({D}alhousie
              {U}niv., {H}alifax, {N}.{S}., 1973)},
     PAGES = {141--209},
      YEAR = {1973},
   MRCLASS = {47B47 (47G05)},
  MRNUMBER = {0390829},
MRREVIEWER = {James A. Deddens},
}

@book {Hen,
    AUTHOR = {Henrot, A.},
     TITLE = {Extremum problems for eigenvalues of elliptic operators},
    SERIES = {Frontiers in Mathematics},
 PUBLISHER = {Birkh\"{a}user Verlag, Basel},
      YEAR = {2006},
     PAGES = {x+202},
   MRCLASS = {35P05 (35J05 35J15 35P15 47F05 49R50)},
  MRNUMBER = {2251558},
MRREVIEWER = {G\"{u}nter Berger},
}

@article{HeuKor,
    title = {Axioms for the category of {H}ilbert spaces},
   author = {Heunen, C. and Kornell, A.},
  journal = {Proc. Natl. Acad. Sci. USA},
   volume = {119},
   number = {9},
    pages = {e2117024119},
     year = {2022},
publisher = {Nat. Acad. Sci.}
}

@book {Heuser,
    AUTHOR = {Heuser, H.},
     TITLE = {Funktionalanalysis},
    SERIES = {Math. Leitf\"aden},
   EDITION = {4th},
 PUBLISHER = {B. G. Teubner, Stuttgart},
      YEAR = {2006},
     PAGES = {696},
   MRCLASS = {46-01 (47-01)},
  MRNUMBER = {2380292},
}

@incollection {Hig,
    AUTHOR = {Higson, N.},
 BOOKTITLE = {Contemporary developments in algebraic {$K$}-theory},
    SERIES = {ICTP Lect. Notes, XV},
     PAGES = {443--536},
     TITLE = {The local index formula in noncommutative geometry},
 PUBLISHER = {Abdus Salam Int. Cent. Theoret. Phys., Trieste},
      YEAR = {2004},
   MRCLASS = {58J42 (19K56 58J22)},
  MRNUMBER = {2175637},
MRREVIEWER = {Evgeniy V. Troitski\u{\i}},
}

@BOOK{HilPhi,
    Address        = {Providence, RI},
    Author         = {Hille, E. and Phillips, R. S.},
    Note           = {rev. ed},
    Pages          = {xii+808},
    Publisher      = {Amer. Math. Soc.},
    Series         = {Amer. Math. Soc. Colloq. Publ., vol. 31},
    Title          = {{F}unctional analysis and semi-groups},
    Year           = {1957},
    Mrnumber       = {MR0089373 (19,664d)}
}

@book {Holevo,
    AUTHOR = {Holevo, A.},
     TITLE = {Probabilistic and statistical aspects of quantum theory},
    SERIES = {Quaderni/ Monographs},
    VOLUME = {1},
   EDITION = {2nd},
 PUBLISHER = {Edizioni della Normale, Pisa},
      YEAR = {2011},
     PAGES = {xvi+323},
   MRCLASS = {81-02 (46N50 60K40 62P35 81S25)},
  MRNUMBER = {2797301},
}

@article {Hol,
    AUTHOR = {Holland, Jr., S. S.},
     TITLE = {Orthomodularity in infinite dimensions; a theorem of {M}.
              {S}ol\`er},
   JOURNAL = {Bull. Amer. Math. Soc. (N.S.)},
  FJOURNAL = {American Mathematical Society. Bulletin. New Series},
    VOLUME = {32},
      YEAR = {1995},
    NUMBER = {2},
     PAGES = {205--234},
   MRCLASS = {06C15 (11E39 16W99 46C15 51D99 81P10)},
  MRNUMBER = {1307904},
MRREVIEWER = {Gudrun Kalmbach},
}

@book{Hormander1985,
  author    = {H{\"o}rmander, L.},
  title     = {The analysis of linear partial differential operators {III}: {P}seudo-differential operators},
  series    = {Grundlehren der mathematischen Wissenschaften},
  volume    = {274},
  publisher = {Springer},
  year      = {1985},
}

@book {HNVW1,
    AUTHOR = {Hyt\"onen, T. P. and van Neerven, J. M. A. M. and Veraar, M. C. and
              Weis, L. W.},
     TITLE = {Analysis in {B}anach spaces. {V}ol. {I}: {M}artingales and
              {L}ittlewood-{P}aley theory},
    SERIES = {Ergebnisse der Mathematik und ihrer Grenzgebiete. 3. Folge},
    VOLUME = {63},
 PUBLISHER = {Springer, Cham},
      YEAR = {2016},
     PAGES = {xvi+614},
   MRCLASS = {46-02 (42B35 46E30)},
  MRNUMBER = {3617205},
MRREVIEWER = {Adam Os\polhk ekowski},
}

@book {HNVW2,
    AUTHOR = {Hyt\"{o}nen, T. P. and van Neerven, J. M. A. M. and Veraar, M. C. and Weis, L. W.},
     TITLE = {Analysis in {B}anach spaces. {V}ol. {II}: {P}robabilistic methods and operator theory},
    SERIES = {Ergebnisse der Mathematik und ihrer Grenzgebiete. 3. Folge},
    VOLUME = {67},
 PUBLISHER = {Springer, Cham},
      YEAR = {2017},
     PAGES = {xxi+616},
   MRCLASS = {46-02 (42B35 46E30 47-02 60B11 60H30)},
  MRNUMBER = {3752640},
MRREVIEWER = {Adam Os\polhk ekowski},
}

@book {HNVW3,
    AUTHOR = {Hyt\"onen, T. P. and van Neerven, J. M. A. M. and Veraar, M. C. and
              Weis, L. W.},
     TITLE = {Analysis in {B}anach spaces. {V}ol. {III}: {H}armonic analysis and operator theory},
      SERIES = {Ergebnisse der Mathematik und ihrer Grenzgebiete. 3. Folge},
    VOLUME = {76},
 PUBLISHER = {Springer, Cham},
      YEAR = {2023},
}

@book {HNVW4,
    AUTHOR = {Hyt\"onen, T. P. and van Neerven, J. M. A. M. and Veraar, M. C. and
              Weis, L. W.},
     TITLE = {Analysis in {B}anach spaces. {V}ol. {IV}: {S}tochastic analysis},
     NOTE   = {In preparation},
    YEAR = {2026+},
}

@BOOK{Jan,
    Address        = {Cambridge},
    Author         = {Janson, S.},
    Pages          = {x+340},
    Publisher      = {Cambridge University Press},
    Series         = {Cambridge Tracts in Mathematics},
    Title          = {{G}aussian {H}ilbert spaces},
    Volume         = {129},
    Year           = {1997},
    Mrnumber       = {MR1474726 (99f:60082)}
}

@book {Jech,
    AUTHOR = {Jech, Th. J.},
     TITLE = {The axiom of choice},
    SERIES = {Studies in Logic and the Foundations of Mathematics},
    VOLUME = {75},
 PUBLISHER = {North-Holland, Amsterdam-London},
      YEAR = {1973},
     PAGES = {xi+202},
   MRCLASS = {04A25 (02K05)},
  MRNUMBER = {0396271},
MRREVIEWER = {D. Pincus},
}

@book {Jech-set,
    AUTHOR = {Jech, Th. J.},
     TITLE = {Set theory},
    SERIES = {Springer Monogr. Math.},
    EDITION = {3rd},
 PUBLISHER = {Springer, Berlin},
      YEAR = {2003},
     PAGES = {xiv+769},
   MRCLASS = {03Exx (03-01 03-02)},
  MRNUMBER = {1940513},
MRREVIEWER = {Eva Coplakova},
}

@article{JordanKinderlehrerOtto1998,
  author  = {Jordan, R. and Kinderlehrer, D. and Otto, F.},
  title   = {The variational formulation of the
             {Fokker--Planck} equation},
  journal = {SIAM J. Math. Anal.},
  volume  = {29},
  number  = {1},
  pages   = {1--17},
  year    = {1998},
  doi     = {10.1137/S0036141096303359}
}

@article{Jo,
   title = {Spectral theory for infinitesimal generators of one-parameter groups of isometries: {T}he min-max principle and compact
perturbations},
  author = {Jorgensen, P.},
 journal = {J. Math. Anal. Appl.},
  volume = {90},
   pages = {343--370},
    year = {1982}
}

@book {Jost,
    AUTHOR = {Jost, J.},
     TITLE = {Partial differential equations},
    SERIES = {Graduate Texts in Mathematics},
    VOLUME = {214},
   EDITION = {3rd},
 PUBLISHER = {Springer, New York},
      YEAR = {2013},
     PAGES = {xiv+410},
   MRCLASS = {35-01},
  MRNUMBER = {3012036},
}

@article {Kac,
    AUTHOR = {Kac, M.},
     TITLE = {Can one hear the shape of a drum?},
   JOURNAL = {Amer. Math. Monthly},
  FJOURNAL = {American Mathematical Monthly},
    VOLUME = {73},
      YEAR = {1966},
    NUMBER = {4, part II},
     PAGES = {1--23},
   MRCLASS = {57.50 (00.00)},
  MRNUMBER = {201237},
MRREVIEWER = {I. Stakgold},
}

@article{Kantorovich1942,
  author  = {Kantorovich, L. V.},
  title   = {On the translocation of masses},
  journal = {Dokl. Akad. Nauk SSSR},
  volume  = {37},
  number  = {7--8},
  pages   = {227--229},
  year    = {1942},
  note    = {In Russian.}
}

@book{Kal,
  author    = {Kallenberg, O.},
  title     = {Foundations of modern probability},
  edition   = {3rd},
  series    = {Probability Theory and Stochastic Modelling},
  volume    = {99},
  publisher = {Springer},
  address   = {Cham},
  year      = {2021},
}

@book {Kato-PertTh,
    AUTHOR = {Kato, T.},
     TITLE = {Perturbation theory for linear operators},
    SERIES = {Classics in Mathematics},
      NOTE = {Reprint of the 1980 edition},
 PUBLISHER = {Springer, Berlin},
      YEAR = {1995},
     PAGES = {xxii+619},
   MRCLASS = {47A55 (46-00 47-00)},
  MRNUMBER = {1335452 (96a:47025)},
}

@Article{Kato-sqrt,
    Author = {Kato, T.},
    Title = {Fractional powers of dissipative operators},
    FJournal = {Journal of the Mathematical Society of Japan},
    Journal = {J. Math. Soc. Japan},
    Volume = {13},
    Pages = {246--274},
    Year = {1961},
    Publisher = {Mathematical Society of Japan, Tokyo},
    Language = {English},
    Zbl = {0113.10005}
}

@article{Kellerer1964,
  author  = {Kellerer, H. G.},
  title   = {Ma{\ss}theoretische {M}arginalprobleme},
  journal = {Math. Ann.},
  volume  = {153},
  year    = {1964},
  pages   = {168--198},
  doi     = {10.1007/BF01360315}
}

@article{Kellerer1984,
  author  = {Kellerer, H. G.},
  title   = {Duality theorems for marginal problems},
  journal = {Z. Wahrscheinlichkeitstheorie verw. Gebiete},
  volume  = {67},
  year    = {1984},
  pages   = {399--432},
  doi     = {10.1007/BF00532047}
}

@book{Kelley-GeneralTopology,
  author    = {Kelley, J. L.},
  title     = {General topology},
  publisher = {Springer},
  year      = {1975},
  series    = {Graduate Texts in Mathematics},
  volume    = {27},
  note      = {Reprint of the corrected 1961 edition}
}

@article{KCBS,
  author = {Klyachko, A. A. and Can, M. A. and Binicioğlu, S. and Shumovsky, A. S.},
   title = {Simple test for hidden variables in spin-1 systems},
 journal = {Phys. Rev. Lett.},
  volume = {101},
    year = {2008},
   pages = {020403},
}

@article{KnasterKuratowskiMazurkiewicz1929,
  author  = {Knaster, B. and Kuratowski, K. and Mazurkiewicz, S.},
  title   = {Ein {B}eweis des {F}ixpunktsatzes f{\"u}r $n$-dimensionale {S}implexe},
  journal = {Fundamenta Mathematicae},
  volume  = {14},
  number  = {1},
  year    = {1929},
  pages   = {132--137},
}

@article{Koelink,
  title={8 Lectures on quantum groups and $q$-special functions},
  author={Koelink, H. T.},
  journal={arXiv:q-alg/9608018},
  year={1996}
}

@article{Krivine1979,
  author  = {Krivine, J.-L.},
  title   = {Constantes de {Grothendieck} et fonctions de type positif
             sur les sph{\`e}res},
  journal = {Adv. in Math.},
  volume  = {31},
  number  = {1},
  year    = {1979},
  pages   = {16--30},
  doi     = {10.1016/0001-8708(79)90017-3}
}

@book {Kry,
    AUTHOR = {Krylov, N. V.},
     TITLE = {Lectures on elliptic and parabolic equations in {S}obolev
              spaces},
    SERIES = {Graduate Studies in Mathematics},
    VOLUME = {96},
 PUBLISHER = {Amer. Math. Soc., Providence, RI},
      YEAR = {2008},
     PAGES = {xviii+357},
   MRCLASS = {35-01 (35J15 35K10 35S05 46E35 46N20)},
  MRNUMBER = {2435520},
MRREVIEWER = {Vicen\c{t}iu D. R\u{a}dulescu}
}

@article{Landsman-FA-QM,
  title={Quantum theory and functional analysis},
  author={Landsman, N. P.},
  journal={arXiv:1911.06630},
  year={2019}
}

@book {Landsman1998,
    AUTHOR = {Landsman, N. P.},
     TITLE = {Mathematical topics between classical and quantum mechanics},
    SERIES = {Springer Monographs in Mathematics},
 PUBLISHER = {Springer, New York},
      YEAR = {1998},
     PAGES = {xx+529},
   MRCLASS = {81R15 (00A05 46Lxx 70Hxx 81-02 81S10)},
  MRNUMBER = {1662141},
MRREVIEWER = {Markus J. Pflaum},
}

@book {Lan,
    AUTHOR = {Landsman, N. P.},
     TITLE = {Foundations of quantum theory},
    SERIES = {Fundamental Theories of Physics},
    VOLUME = {188},
 PUBLISHER = {Springer, Cham},
      YEAR = {2017},
     PAGES = {xv+881},
   MRCLASS = {81-01 (81P10 81R05 81R10)},
  MRNUMBER = {3643288},
MRREVIEWER = {Michael M. Dediu},
}

@book {Lax,
    AUTHOR = {Lax, P. D.},
     TITLE = {Functional analysis},
 PUBLISHER = {Wiley},
      YEAR = {2002},
}

@article{LevWei,
title = {Inequalities between {D}irichlet and {N}eumann eigenvalues},
 author = {Levine, H. A. and Weinberger, H. F.},
journal = {Arch. Ration. Mech. Anal.},
 volume = {94},
 number = {3},
  pages = {193--208},
   year = {1986},
}

@article {Li94,
    AUTHOR = {Li, C.-K.},
     TITLE = {{$C$}-numerical ranges and {$C$}-numerical radii},
   JOURNAL = {Linear Multilinear Algebra},
  FJOURNAL = {Linear and Multilinear Algebra},
    VOLUME = {37},
      YEAR = {1994},
    NUMBER = {1--3},
     PAGES = {51--82},
   MRCLASS = {15A60 (15-02 47A12)},
  MRNUMBER = {1313758},
MRREVIEWER = {Natalia Bebiano},
}

@book {LiQue,
    AUTHOR = {Li, D. and Queff{\'e}lec, H.},
     TITLE = {Introduction \`a l'\'etude des espaces de {B}anach},
    SERIES = {Cours Sp\'ecialis\'es},
    VOLUME = {12},
 PUBLISHER = {Soci\'et\'e Math\'ematique de France, Paris},
      YEAR = {2004},
     PAGES = {xxiv+627},
   MRCLASS = {46B09 (46Bxx 60B12 60G15)},
  MRNUMBER = {2124356 (2006c:46012)},
MRREVIEWER = {Anatolij M. Plichko},
}

@article{LinTza,
    AUTHOR = {Lindenstrauss, J. and Tzafriri, L.},
     TITLE = {On the complemented subspaces problem},
   JOURNAL = {Israel J. Math.},
  FJOURNAL = {Israel Journal of Mathematics},
    VOLUME = {9},
      YEAR = {1971},
     PAGES = {263--269},
   MRCLASS = {46.10},
  MRNUMBER = {0276734 (43 \#2474)},
MRREVIEWER = {W. J. Davis},
}

@book{LionsMagenes1968,
  author    = {Lions, J.-L. and Magenes, E.},
  title     = {Probl\`emes aux limites non homog\`enes et applications. {V}ols. 1, 2},
  publisher = {Dunod},
  address   = {Paris},
  year      = {1968}
}

@book{Ludwig-AxiomaticI,
  author    = {G. Ludwig},
  title     = {An axiomatic basis for quantum mechanics. {V}ol. 1: {D}erivation of {H}ilbert space structure},
  publisher = {Springer},
  address   = {Berlin},
  year      = {1985},
  doi       = {10.1007/978-3-642-70029-3},
}

@book {Lun,
    AUTHOR = {Lunardi, A.},
     TITLE = {Analytic semigroups and optimal regularity in parabolic
              problems},
    SERIES = {Progress in Nonlinear Differential Equations and their
              Applications},
    VOLUME = {16},
 PUBLISHER = {Birkh\"{a}user, Basel},
      YEAR = {1995},
     PAGES = {xviii+424},
   MRCLASS = {47D06 (34G20 35Kxx 46M35 46N20 47-02 47N20 58D25)},
  MRNUMBER = {1329547},
MRREVIEWER = {Paolo Acquistapace},
}

@book{LuxZaa,
  title = {Riesz spaces},
  author = {Luxemburg, W. A. J. and Zaanen, A. C.},
  year = {1971},
  publisher = {North-Holland}
}

@article {MaaNee,
    AUTHOR = {Maas, J. and van Neerven, J. M. A. M.},
     TITLE = {Boundedness of {R}iesz transforms for elliptic operators on
              abstract {W}iener spaces},
   JOURNAL = {J. Funct. Anal.},
  FJOURNAL = {Journal of Functional Analysis},
    VOLUME = {257},
      YEAR = {2009},
    NUMBER = {8},
     PAGES = {2410--2475},
     CODEN = {JFUAAW},
   MRCLASS = {60H07 (35J05 42B25 46E15 47A60 47D06 47F05)},
  MRNUMBER = {2555015 (2011c:60177)},
}

@book {Mac,
    AUTHOR = {Mackey, G. W.},
     TITLE = {Induced representations of groups and quantum mechanics},
 PUBLISHER = {W. A. Benjamin, Inc., New York-Amsterdam; Boringhieri,
              Turin},
      YEAR = {1968},
     PAGES = {viii+167},
   MRCLASS = {22-01},
  MRNUMBER = {0507212},
}

@inproceedings{Maligr,
  title={On the norms of operators in the real and the complex case},
  author={Maligranda, L.},
  booktitle={Seminar on {B}anach spaces and related topics},
  pages={67--71},
  year={1997}
}

@incollection{McI-Kato,
    AUTHOR = {McIntosh, A.},
 BOOKTITLE = {Nonlinear partial differential equations and their
              applications. {C}oll\`ege de {F}rance {S}eminar, {V}ol. {III} ({P}aris, 1980/1981)},
    SERIES = {Res. Notes in Math.},
    VOLUME = {70},
     PAGES = {252--267},
     TITLE = {On representing closed accretive sesquilinear forms as
              {$(A^{1/2}u,\,A^{\ast 1/2}v)$}},
 PUBLISHER = {Pitman, Boston, Mass.-London},
      YEAR = {1982},
   MRCLASS = {47B44},
  MRNUMBER = {670278},
}

@article{Mercer,
  title={Functions of positive and negative type, and their connection to the theory of integral equations},
  author={Mercer, J.},
  journal={Philos. Trans. R. Soc. Lond. A},
  volume={209},
  number={441--458},
  pages={415--446},
  year={1909},
  publisher={The Royal Society London}
}

@article{MPRS,
  title = {The domain of the {O}rnstein--{U}hlenbeck operator on an {$L^p$}-space with invariant measure},
 author = {Metafune, G. and Pr{\"u}ss, J. and Rhandi, A. and Schnaubelt, R.},
journal = {Ann. Sc. Norm. Super. Pisa Cl. Sci. (5)},
 volume = {1},
 number = {2},
  pages = {471--485},
   year = {2002}
}

@book {MeyNie,
    AUTHOR = {Meyer-Nieberg, P.},
     TITLE = {Banach lattices},
    SERIES = {Universitext},
 PUBLISHER = {Springer},
   ADDRESS = {Berlin},
      YEAR = {1991},
     PAGES = {xvi+395},
   MRCLASS = {46B42 (46A40 47B60)},
  MRNUMBER = {1128093 (93f:46025)},
MRREVIEWER = {Yu. A. Abramovich},
}

@book {Monna,
    AUTHOR = {Monna, A. F.},
     TITLE = {Functional analysis in historical perspective},
 PUBLISHER = {John Wiley \& Sons, New York-Toronto, Ont.},
      YEAR = {1973},
     PAGES = {viii+167},
   MRCLASS = {46-03 (01A55 01A60)},
  MRNUMBER = {0482022 (58 \#2112)},
MRREVIEWER = {Thomas Hawkins},
}

@article{Moyal1949,
  author  = {Moyal, J. E.},
  title   = {Quantum mechanics as a statistical theory},
  journal = {Proc. Cambridge Philos. Soc.},
  volume  = {45},
  year    = {1949},
  pages   = {99--124},
}

@book {Mul,
    AUTHOR = {M\"{u}ller, V.},
     TITLE = {Spectral theory of linear operators and spectral systems in
              {B}anach algebras},
    SERIES = {Operator Theory: Advances and Applications},
    VOLUME = {139},
   EDITION = {2nd},
 PUBLISHER = {Birkh\"{a}user, Basel},
      YEAR = {2007},
     PAGES = {x+439},
   MRCLASS = {47A10 (46H05 46H30 47A13)},
  MRNUMBER = {2355630},
}

@article {Murphy,
    AUTHOR = {Murphy, G. J.},
     TITLE = {Fredholm index theory and the trace},
   JOURNAL = {Proc. Roy. Irish Acad. Sect. A},
  FJOURNAL = {Proceedings of the Royal Irish Academy. Section A.
              Mathematical and Physical Sciences},
    VOLUME = {94},
      YEAR = {1994},
    NUMBER = {2},
     PAGES = {161--166},
   MRCLASS = {47A53 (46L99 47D25)},
  MRNUMBER = {1369029},
MRREVIEWER = {Steve Wright},
}

@INCOLLECTION{Mus,
    Address        = {Amsterdam},
    Author         = {Musia{\l}, K.},
    Booktitle      = {{H}andbook of measure theory, {V}ol. {I}, {II}},
    Pages          = {531--586},
    Publisher      = {North-Holland},
    Title          = {Pettis integral},
    Year           = {2002},
    Mrnumber       = {MR1954622 (2004d:28026)}
}

@book {Nagy,
    AUTHOR = {Sz.-Nagy, B.},
     TITLE = {Spektraldarstellung linearer {T}ransformationen des
              {H}ilbertschen {R}aumes},
    SERIES = {Ergebnisse der Mathematik und ihrer Grenzgebiete, vol. 39},
 PUBLISHER = {Springer, Berlin-New York},
      YEAR = {1967},
     PAGES = {vi+81},
   MRCLASS = {47.00},
  MRNUMBER = {0213890},
}

@article{Nee97,
  title = {Elementary operator-theoretic proof of {W}iener's {T}auberian theorem},
 author = {{\VAN{Neerven}{van}{van}} Neerven, J. M. A. M.},
   year = {1997},
journal = {Rend. Istit. Mat. Univ. Trieste, Suppl. Vol. XXVIII},
 pages = {281--286},
}

@incollection{NVW-survey,
       title = {Stochastic integration in {B}anach spaces -- a survey},
      author = {{\VAN{Neerven}{van}{van}} Neerven, J. M. A. M. and Veraar, M. C.  and Weis, L. W.},
   booktitle = {Stochastic analysis: A series of lectures},
      series = {Progress in Probability},
   publisher = {Birkh\"auser Verlag},
      volume = {68},
        year = {2015},
}

@article {NP18,
    AUTHOR = {{\VAN{Neerven}{van}{van}} Neerven, J. M. A. M. and Portal, P.},
     TITLE = {The {W}eyl calculus with respect to the {G}aussian measure and
              restricted {$L^p$}-{$L^q$} boundedness of the
              {O}rnstein-{U}hlenbeck semigroup in complex time},
   JOURNAL = {Bull. Soc. Math. France},
  FJOURNAL = {Bulletin de la Soci\'{e}t\'{e} Math\'{e}matique de France},
    VOLUME = {146},
      YEAR = {2018},
    NUMBER = {4},
     PAGES = {691--712},
   MRCLASS = {47A60 (60G15 81S05)},
  MRNUMBER = {3936540},
MRREVIEWER = {Eric Stachura},
}

@article{vanNeervenWaaijer2026,
  author  = {{\VAN{Neerven}{van}{van}} Neerven, J. M. A. M. and Waaijer, M.},
  title   = {Indiscernibility of quantum states},
  journal = {J. Math. Phys.},
  volume  = {67},
  number  = {5},
  pages   = {052102},
  year    = {2026},
  doi     = {10.1063/5.0302377}
}

@book {Neu,
    AUTHOR = {{\VAN{Neumann}{von}{von}} Neumann, J.},
     TITLE = {Mathematische {G}rundlagen der {Q}uantenmechanik},
    SERIES = {Grundlehren der mathematischen Wissenschaften},
    VOLUME = {38},
      NOTE = {Reprint of the first edition, 1932},
 PUBLISHER = {Springer, Berlin-New York},
      YEAR = {1968},
     PAGES = {v+262},
   MRCLASS = {81.00},
  MRNUMBER = {0223138},
}

@article {Nica,
    AUTHOR = {Nica, B.},
     TITLE = {The {M}azur--{U}lam theorem},
   JOURNAL = {Expo. Math.},
  FJOURNAL = {Expositiones Mathematicae},
    VOLUME = {30},
      YEAR = {2012},
    NUMBER = {4},
     PAGES = {397--398},
   MRCLASS = {46B04},
  MRNUMBER = {2997831},
}

@book {Nik,
   author = {Nikolski, N. K.},
    title = {Operators, functions, and systems: an easy reading, {V}ol. 1},
   series = {Mathematical Surveys and Monographs},
   volume = {92},
publisher = {Amer. Math. Soc., Providence, RI},
     year = {2002}
}

@book {Nua,
    AUTHOR = {Nualart, D.},
     TITLE = {The {M}alliavin calculus and related topics},
    SERIES = {Probability and its Applications (New York)},
   EDITION = {2nd},
 PUBLISHER = {Springer, Berlin},
      YEAR = {2006},
     PAGES = {xiv+382},
   MRCLASS = {60-02 (60H07 60H30)},
  MRNUMBER = {2200233 (2006j:60004)},
MRREVIEWER = {Daniel Ocone},
}

@phdthesis{Nyman1950,
  author      = {Nyman, B.},
  title       = {On some groups and semigroups of translations},
  school      = {University of Uppsala},
  year        = {1950},
}

@article {Osgood,
    AUTHOR = {Osgood, W. F.},
     TITLE = {A {J}ordan curve of positive area},
   JOURNAL = {Trans. Amer. Math. Soc.},
  FJOURNAL = {Transactions of the American Mathematical Society},
    VOLUME = {4},
      YEAR = {1903},
    NUMBER = {1},
     PAGES = {107--112},
   MRCLASS = {28A75},
  MRNUMBER = {1500628},
}

@book {Ouh,
    AUTHOR = {Ouhabaz, E. M.},
     TITLE = {Analysis of heat equations on domains},
    SERIES = {London Mathematical Society Monographs Series},
    VOLUME = {31},
 PUBLISHER = {Princeton University Press, Princeton, NJ},
      YEAR = {2005},
     PAGES = {xiv+284},
   MRCLASS = {35-02 (35J70 35K20 47D06 47D07 47F05)},
  MRNUMBER = {2124040},
MRREVIEWER = {Sergey G. Pyatkov},
}

@book {Parth,
    AUTHOR = {Parthasarathy, K. R.},
     TITLE = {Mathematical foundations of quantum mechanics},
    SERIES = {Texts and Readings in Mathematics},
    VOLUME = {35},
 PUBLISHER = {Hindustan Book Agency},
   ADDRESS = {New Delhi},
      YEAR = {2005},
     PAGES = {viii+170},
   MRCLASS = {81-01 (46N50 47N50)},
  MRNUMBER = {2188339 (2007g:81001)},
MRREVIEWER = {Pavel B. Kurasov},
}

@book {Par,
    AUTHOR = {Parthasarathy, K. R.},
     TITLE = {An introduction to quantum stochastic calculus},
    SERIES = {Modern Birkh\"{a}user Classics},
      NOTE = {2012 reprint of the 1992 original},
 PUBLISHER = {Birkh\"{a}user/Springer Basel AG, Basel},
      YEAR = {2012},
     PAGES = {xii+290},
   MRCLASS = {01A75},
  MRNUMBER = {3012668},
}

@book{Pauli1990-Reprint,
  author    = {W. Pauli},
  title     = {Die allgemeinen {P}rinzipien der {W}ellenmechanik},
  note      = {Annotated reprint of the original 1933 article},
  publisher = {Springer},
  address   = {Berlin},
  year      = {1990},
}

@BOOK{Paz,
    Address        = {New York},
    Author         = {Pazy, A.},
    Pages          = {viii+279},
    Publisher      = {Springer},
    Series         = {Applied Math. Sciences},
    Title          = {Semigroups of linear operators and applications to partial differential equations},
    Volume         = {44},
    Year           = {1983},
    Mrnumber       = {MR710486 (85g:47061)}
}

@book{Pedersen,
    AUTHOR = {Pedersen, G. K.},
     TITLE = {{$C^*$}-algebras and their automorphism groups},
    SERIES = {Pure and Applied Mathematics (Amsterdam)},
 PUBLISHER = {Academic Press, London},
      YEAR = {2018},
     PAGES = {xviii+520},
   MRCLASS = {46Lxx (46-02)},
  MRNUMBER = {3839621},
}

@ARTICLE{Pet38,
    Author         = {Pettis, B. J.},
    Journal        = {Trans. Amer. Math. Soc.},
    Number         = {2},
    Pages          = {277--304},
    Title          = {On integration in vector spaces},
    Volume         = {44},
    Year           = {1938},
    Mrnumber       = {MR1501970}
}

@article {Phelps,
    AUTHOR = {Phelps, R. R.},
     TITLE = {Uniqueness of {H}ahn-{B}anach extensions and unique best
              approximation},
   JOURNAL = {Trans. Amer. Math. Soc.},
  FJOURNAL = {Transactions of the American Mathematical Society},
    VOLUME = {95},
      YEAR = {1960},
     PAGES = {238--255},
   MRCLASS = {46.00},
  MRNUMBER = {113125},
MRREVIEWER = {V. Pt\'{a}k},
}

@article{Phi55,
Author         = {Phillips, R. S.},
Title          = {The adjoint semi-group},
Journal        = {Pacific J. Math.},
Volume         = {5},
Pages          = {269--283},
Year           = {1955}
}

@book {Pie,
    AUTHOR = {Pietsch, A.},
     TITLE = {History of {B}anach spaces and linear operators},
 PUBLISHER = {Birkh\"auser Boston Inc.},
   ADDRESS = {Boston, MA},
      YEAR = {2007},
     PAGES = {xxiv+855},
   MRCLASS = {46-03 (01A05 01A60 46Bxx 47-03)},
  MRNUMBER = {2300779 (2008i:46002)},
MRREVIEWER = {Joram Lindenstrauss},
}

@article {Pir,
    AUTHOR = {Piron, C.},
     TITLE = {Axiomatique quantique},
   JOURNAL = {Helv. Phys. Acta},
  FJOURNAL = {Helvetica Physica Acta. Physica Theoretica. Societatis
              Physicae Helveticae Commentaria Publica},
    VOLUME = {37},
      YEAR = {1964},
     PAGES = {439--468},
   MRCLASS = {81.02},
  MRNUMBER = {204048},
MRREVIEWER = {J. C. T. Pool},
}

@book{Pis-Mart,
  title={Martingales in {B}anach spaces},
  author={Pisier, G.},
  series={Cambridge Studies in Advanced Mathematics},
  volume={155},
  year={2016},
  publisher={Cambridge University Press}
}

@article{Pisier2012,
  author  = {Pisier, G.},
  title   = {{Grothendieck}'s theorem, past and present},
  journal = {Bull. Amer. Math. Soc. (N.S.)},
  volume  = {49},
  number  = {2},
  year    = {2012},
  pages   = {237--323},
  doi     = {10.1090/S0273-0979-2011-01348-9}
}

@book {ReedSimon2,
    AUTHOR = {Reed, M. and Simon, B.},
     TITLE = {Methods of modern mathematical physics. {V}ol. {II}: {F}ourier analysis, self-adjointness},
 PUBLISHER = {Academic Press},
   ADDRESS = {New York},
      YEAR = {1975},
}

@article {Riesz,
    AUTHOR = {Riesz, M.},
     TITLE = {Sur les maxima des formes bilin\'{e}aires et sur les
              fonctionnelles lin\'{e}aires},
   JOURNAL = {Acta Math.},
  FJOURNAL = {Acta Mathematica},
    VOLUME = {49},
      YEAR = {1926},
    NUMBER = {3--4},
     PAGES = {465--497},
   MRCLASS = {DML},
  MRNUMBER = {1555250},
}

@incollection {Ros,
    AUTHOR = {Rosenberg, J.},
 BOOKTITLE = {Operator algebras, quantization, and noncommutative geometry},
    SERIES = {Contemp. Math.},
    VOLUME = {365},
     PAGES = {331--353},
     TITLE = {A selective history of the {S}tone-von {N}eumann theorem},
 PUBLISHER = {Amer. Math. Soc., Providence, RI},
      YEAR = {2004},
   MRCLASS = {46-03 (01A60 22D10 22D25 46L60 81Q60 81S05)},
  MRNUMBER = {2106827},
}

@book {Rubin,
    AUTHOR = {Rubin, H. and Rubin, J. E.},
     TITLE = {Equivalents of the axiom of choice},
    SERIES = {Studies in Logic and the Foundations of Mathematics},
 PUBLISHER = {North-Holland Publishing Co., Amsterdam-London},
      YEAR = {1970},
     PAGES = {xxiii+134},
   MRCLASS = {02K20},
  MRNUMBER = {0434812},
}

@book{RudinLCA,
  author    = {Rudin, W.},
  title     = {Fourier analysis on groups},
  publisher = {Interscience Publishers},
  address   = {New York},
  year      = {1962},
  series    = {Interscience Tracts in Pure and Applied Mathematics},
  volume    = {12}
}

@book{Rudin-RCA,
    AUTHOR = {Rudin, W.},
     TITLE = {Real and complex analysis},
   EDITION = {3rd},
 PUBLISHER = {McGraw-Hill Book Co., New York},
      YEAR = {1987},
     PAGES = {xiv+416},
   MRCLASS = {00A05 (26-01 30-01 46-01)},
  MRNUMBER = {924157 (88k:00002)},
}

@BOOK{Rudin,
    AUTHOR = {Rudin, W.},
     TITLE = {Functional analysis},
    SERIES = {International Series in Pure and Applied Mathematics},
   EDITION = {2nd},
 PUBLISHER = {McGraw-Hill, Inc., New York},
      YEAR = {1991},
     PAGES = {xviii+424},
   MRCLASS = {46-01 (47-01)},
  MRNUMBER = {1157815 (92k:46001)},
}

@book {RuzTur,
    AUTHOR = {Ruzhansky, M. and Turunen, V.},
     TITLE = {Pseudo-differential operators and symmetries},
    SERIES = {Pseudo-Differential Operators. Theory and Applications},
    VOLUME = {2},
 PUBLISHER = {Birkh\"{a}user Verlag, Basel},
      YEAR = {2010},
     PAGES = {xiv+709},
   MRCLASS = {35-02 (35S05 43A77 43A80 43A85 47G30 58J40)},
  MRNUMBER = {2567604},
MRREVIEWER = {Fabio Nicola},
}

@book {Ryan,
    AUTHOR = {Ryan, R. A.},
     TITLE = {Introduction to tensor products of {B}anach spaces},
    SERIES = {Springer Monographs in Mathematics},
 PUBLISHER = {Springer London, Ltd., London},
      YEAR = {2002},
     PAGES = {xiv+225},
   MRCLASS = {46B28 (46-01 47A07 47B10 47L20)},
  MRNUMBER = {1888309},
MRREVIEWER = {Andreas Defant},
}

@book{SafVas,
  title = {The asymptotic distribution of eigenvalues of partial differential operators},
  author = {Safarov, Yu. and Vassilev, D.},
  volume = {155},
  year = {1997},
  series = {Transl. Math. Monogr.},
  publisher={Amer. Math. Soc., Providence, RI.}
}

@book{Santambrogio2015,
  author    = {Santambrogio, F.},
  title     = {Optimal transport for applied mathematicians:
               {C}alculus of variations, {PDEs}, and modeling},
  series    = {Progress in Nonlinear Differential Equations
               and Their Applications},
  volume    = {87},
  publisher = {Birkh{\"a}user},
  address   = {Cham},
  year      = {2015},
  doi       = {10.1007/978-3-319-20828-2}
}

@book {Schaefer,
    AUTHOR = {Schaefer, H. H.},
     TITLE = {Banach lattices and positive operators},
    SERIES = {Grundlehren der Mathematischen Wissenschaften},
    VOLUME = {215},
 PUBLISHER = {Springer},
   ADDRESS = {New York},
      YEAR = {1974},
     PAGES = {xi+376},
   MRCLASS = {46A40 (47B55 47D20)},
  MRNUMBER = {MR0423039 (54 \#11023)},
MRREVIEWER = {A. C. Zaanen},
}

@book {Schechter,
    AUTHOR = {Schechter, M.},
     TITLE = {Principles of functional analysis},
    SERIES = {Graduate Studies in Mathematics},
    VOLUME = {36},
   EDITION = {2nd},
 PUBLISHER = {Amer. Math. Soc., Providence, RI},
      YEAR = {2002},
     PAGES = {xxii+425},
   MRCLASS = {46-01 (47-01)},
  MRNUMBER = {1861991},
MRREVIEWER = {R. G. Bartle},
}

@book{Schilling,
  author    = {Schilling, R. L.},
  title     = {Measures, integrals and martingales},
  edition   = {2nd},
  publisher = {Cambridge University Press},
  address   = {Cambridge},
  year      = {2017},
}

@book {Schm,
    AUTHOR = {Schm\"udgen, K.},
     TITLE = {Unbounded self-adjoint operators on {H}ilbert space},
    SERIES = {Graduate Texts in Mathematics},
    VOLUME = {265},
 PUBLISHER = {Springer, Dordrecht},
      YEAR = {2012},
     PAGES = {xx+432},
   MRCLASS = {47-01 (47B25 47E05)},
  MRNUMBER = {2953553},
MRREVIEWER = {G. V. Rozenblum},
}

@book{Schn,
    AUTHOR = {Hundertmark, D. and Machinek, L. and Meyries, M. and Schnaubelt, R.},
     TITLE = {Operator semigroups and dispersive equations},
    NOTE = {Lecture notes, 16th Internet Seminar 2012/13},
    year = {2013},
}

@article {Seeley,
    AUTHOR = {Seeley, R. T.},
     TITLE = {Extension of {$C^{\infty }$} functions defined in a half
              space},
   JOURNAL = {Proc. Amer. Math. Soc.},
  FJOURNAL = {Proceedings of the American Mathematical Society},
    VOLUME = {15},
      YEAR = {1964},
     PAGES = {625--626},
   MRCLASS = {46.38},
  MRNUMBER = {165392},
MRREVIEWER = {J. L. Lions},
}

@article {Seg,
    AUTHOR = {Segal, I. E.},
     TITLE = {Tensor algebras over {H}ilbert spaces. {I}},
   JOURNAL = {Trans. Amer. Math. Soc.},
  FJOURNAL = {Transactions of the American Mathematical Society},
    VOLUME = {81},
      YEAR = {1956},
     PAGES = {106--134},
   MRCLASS = {46.2X},
  MRNUMBER = {MR0076317 (17,880d)},
MRREVIEWER = {F. H. Brownell},
}

@article {Simon64,
    AUTHOR = {De Simon, L.},
     TITLE = {Un'applicazione della teoria degli integrali singolari allo
              studio delle equa\-zioni differenziali lineari astratte del
              primo ordine},
   JOURNAL = {Rend. Sem. Mat. Univ. Padova},
  FJOURNAL = {Rendiconti del Seminario Matematico della Universit\`a di
              Padova. The Mathematical Journal of the University of Padova},
    VOLUME = {34},
      YEAR = {1964},
     PAGES = {205--223},
   MRCLASS = {34.95},
  MRNUMBER = {0176192},
MRREVIEWER = {J. Horv\'ath},
}

@article {Sim,
    AUTHOR = {Simon, B.},
     TITLE = {Notes on infinite determinants of {H}ilbert space operators},
   JOURNAL = {Adv. Math.},
  FJOURNAL = {Advances in Mathematics},
    VOLUME = {24},
      YEAR = {1977},
    NUMBER = {3},
     PAGES = {244--273},
   MRCLASS = {47B10},
  MRNUMBER = {482328},
}

@book{Simon74,
    AUTHOR = {Simon, B.},
     TITLE = {The {$P(\phi )_{2}$} {E}uclidean (quantum) field theory},
    SERIES = {Princeton Series in Physics},
 PUBLISHER = {Princeton University Press, Princeton, NJ},
      YEAR = {1974},
     PAGES = {xx+392},
   MRCLASS = {81.60},
  MRNUMBER = {0489552},
MRREVIEWER = {Raymond Streater},
}

@book{Simon-TraceIdeals,
    title = {Trace ideals and their applications},
   author = {Simon, B.},
   series = {Mathematical Surveys and Monographs},
   volume = {120},
   edition = {2nd},
     year = {2005},
publisher = {American Mathematical Soc.}
}

@article {Sol,
    AUTHOR = {Sol\`er, M. P.},
     TITLE = {Characterization of {H}ilbert spaces by orthomodular spaces},
   JOURNAL = {Comm. Algebra},
  FJOURNAL = {Communications in Algebra},
    VOLUME = {23},
      YEAR = {1995},
    NUMBER = {1},
     PAGES = {219--243},
   MRCLASS = {46C15 (06C15 11E39)},
  MRNUMBER = {1311786},
MRREVIEWER = {Hans Arwed Keller},
}

@book{Stein-HA,
  title = {Harmonic analysis: real-variable methods, orthogonality, and oscillatory integrals},
  author = {Stein, E. M.},
  year = {1993},
  publisher = {Princeton University Press}
}

@book {Stein,
    AUTHOR = {Stein, E. M.},
     TITLE = {Singular integrals and differentiability properties of
              functions},
    SERIES = {Princeton Mathematical Series, No. 30},
 PUBLISHER = {Princeton University Press, Princeton, NJ},
      YEAR = {1970},
     PAGES = {xiv+290},
   MRCLASS = {46.38 (26.00)},
  MRNUMBER = {0290095},
MRREVIEWER = {R. E. Edwards},
}

@book{Ste-Sha-4,
    title = {Functional analysis},
     note = {Volume 4 of Princeton Lectures in Analysis},
   author = {Stein, E. M. and Shakarchi, R.},
     year = {2011},
publisher = {Princeton University Press, Princeton, NJ}
}

@book {Sternberg,
    AUTHOR = {Sternberg, S.},
     TITLE = {Group theory and physics},
 PUBLISHER = {Cambridge University Press, Cambridge},
      YEAR = {1994},
     PAGES = {xiv+429},
   MRCLASS = {20-01 (20C30 20C35 22E46)},
  MRNUMBER = {1287387},
MRREVIEWER = {J. L\~{o}hmus},
}

@book {Take,
    AUTHOR = {Takesaki, M.},
     TITLE = {Theory of operator algebras. {I}},
    SERIES = {Encyclopaedia of Mathematical Sciences},
    VOLUME = {124},
      NOTE = {Reprint of the first (1979) edition},
 PUBLISHER = {Springer, Berlin},
      YEAR = {2002},
     PAGES = {xx+415},
   MRCLASS = {46Lxx (46-01)},
  MRNUMBER = {1873025},
}

@book {Takh,
    AUTHOR = {Takhtajan, L. A.},
     TITLE = {Quantum mechanics for mathematicians},
    SERIES = {Graduate Studies in Mathematics},
    VOLUME = {95},
 PUBLISHER = {Amer. Math. Soc., Providence, RI},
      YEAR = {2008},
     PAGES = {xvi+387},
   MRCLASS = {81-01 (37J05 47N50 70H05)},
  MRNUMBER = {2433906},
MRREVIEWER = {H. Hogreve},
}

@book{Tala,
   author = {Talagrand, M.},
    title = {What is a quantum field theory?},
publisher = {Cambridge University Press},
     year = {2022}
}

@ARTICLE{Tal,
    Author         = {Talagrand, M.},
    Journal        = {Mem. Amer. Math. Soc.},
    Number         = {307},
    Title          = {Pettis integral and measure theory},
    Volume         = {51},
    Year           = {1984},
    Mrnumber       = {MR756174 (86j:46042)}
}

@book {Tan,
    AUTHOR = {Tanabe, H.},
     TITLE = {Equations of evolution},
    SERIES = {Monographs and Studies in Mathematics},
    VOLUME = {6},
 PUBLISHER = {Pitman, Boston, Mass.-London},
      YEAR = {1979},
   MRCLASS = {47D05 (34G10 35K22 47H05)},
  MRNUMBER = {533824},
MRREVIEWER = {A. Pazy},
}

@book{Tao2006,
  author    = {Tao, T.},
  title     = {Nonlinear dispersive equations:
               {L}ocal and global analysis},
  series    = {CBMS Regional Conference Series in Mathematics},
  volume    = {106},
  publisher = {Amer. Math. Soc.},
  address   = {Providence, RI},
  year      = {2006},
  doi       = {10.1090/cbms/106}
}

@article {Taylor,
    AUTHOR = {Taylor, A. E.},
     TITLE = {The extension of linear functionals},
   JOURNAL = {Duke Math. J.},
  FJOURNAL = {Duke Mathematical Journal},
    VOLUME = {5},
      YEAR = {1939},
     PAGES = {538--547},
   MRCLASS = {46.3X},
  MRNUMBER = {345},
MRREVIEWER = {R. P. Agnew},
}

@book{Taylor1981,
  author    = {Taylor, M. E.},
  title     = {Pseudodifferential operators},
  publisher = {Princeton University Press},
  year      = {1981},
}

@article {Ter,
    AUTHOR = {Terzio\u{g}lu, T.},
     TITLE = {A characterization of compact linear mappings},
   JOURNAL = {Arch. Math. (Basel)},
  FJOURNAL = {Archiv der Mathematik},
    VOLUME = {22},
      YEAR = {1971},
     PAGES = {76--78},
   MRCLASS = {47B05},
  MRNUMBER = {0291865},
MRREVIEWER = {A. Pietsch},
}

@article{Tsirelson1980,
  author  = {Tsirel'son, B. S.},
  title   = {Quantum generalizations of {Bell}'s inequality},
  journal = {Lett. Math. Phys.},
  volume  = {4},
  year    = {1980},
  pages   = {93--100},
  doi     = {10.1007/BF00417500}
}

@article{Tsirelson1987,
  author  = {Tsirel'son, B. S.},
  title   = {Quantum analogues of the {Bell} inequalities.
             {The} case of two spatially separated domains},
  journal = {J. Soviet Math.},
  volume  = {36},
  number  = {4},
  pages   = {557--570},
  year    = {1987},
  doi     = {10.1007/BF01663472}
}

@article {unter,
    AUTHOR = {Unterberger, A.},
     TITLE = {Oscillateur harmonique et op\'{e}rateurs pseudo-diff\'{e}rentiels},
   JOURNAL = {Ann. Inst. Fourier (Grenoble)},
  FJOURNAL = {Universit\'{e} de Grenoble. Annales de l'Institut Fourier},
    VOLUME = {29},
      YEAR = {1979},
    NUMBER = {3},
     PAGES = {xi, 201--221},
   MRCLASS = {58G15 (47G05 81C99)},
  MRNUMBER = {552965},
MRREVIEWER = {M. S. Agranovich},
}

@book {Varad,
    AUTHOR = {Varadarajan, V. S.},
     TITLE = {Geometry of quantum theory},
   EDITION = {2nd},
 PUBLISHER = {Springer},
   ADDRESS = {New York},
      YEAR = {1985},
     PAGES = {xviii+412},
   MRCLASS = {81Bxx (81-02)},
  MRNUMBER = {805158 (87a:81009)},
MRREVIEWER = {S. Gudder},
}

@book{Villani2009,
  author    = {Villani, C.},
  title     = {Optimal transport: {O}ld and new},
  series    = {Grundlehren der mathematischen Wissenschaften},
  volume    = {338},
  publisher = {Springer},
  address   = {Berlin},
  year      = {2009},
}

@article{WaaijerNeerven,
  author       = {Waaijer, M. and van Neerven, J. M. A. M. },
  title        = {Holevo classicalisation and context integration: {G}lobal versus protocolwise  integration},
  year         = {2026},
  journal      = {arXiv:2509.08487v2},
  primaryClass = {quant-ph},
}

@book{Werner,
    AUTHOR = {Werner, D.},
     TITLE = {Funktionalanalysis},
   EDITION = {3rd},
 PUBLISHER = {Springer, Berlin},
      YEAR = {2000},
     PAGES = {xii+501},
   MRCLASS = {46-01 (47-01)},
  MRNUMBER = {1787146 (2001i:46001)},
MRREVIEWER = {Manfred Wolff},
}

@article{Weyl1927,
  author  = {Weyl, H.},
  title   = {Quantenmechanik und {G}ruppentheorie},
  journal = {Z. Phys.},
  volume  = {46},
  year    = {1927},
  pages   = {1--46},
}

@article {Whit,
    AUTHOR = {Whitley, R.},
     TITLE = {The spectral theorem for a normal operator},
   JOURNAL = {Amer. Math. Monthly},
  FJOURNAL = {American Mathematical Monthly},
    VOLUME = {75},
      YEAR = {1968},
     PAGES = {856--861},
   MRCLASS = {47.30},
  MRNUMBER = {236746},
MRREVIEWER = {J. B. Deeds},
}

@book {Woit,
    AUTHOR = {Woit, P.},
     TITLE = {Quantum theory, groups and representations},
 PUBLISHER = {Springer, Cham},
      YEAR = {2017},
     PAGES = {xxii+668},
   MRCLASS = {81-01 (81Rxx)},
  MRNUMBER = {3726869},
MRREVIEWER = {Michael M. Dediu},
}

@book {Yosida,
    AUTHOR = {Yosida, K.},
     TITLE = {Functional analysis},
    SERIES = {Grundlehren der Mathematischen Wissenschaften},
    VOLUME = {123},
   EDITION = {6th},
 PUBLISHER = {Springer, Berlin-New York},
      YEAR = {1980},
     PAGES = {xii+501},
   MRCLASS = {46-01 (47-01)},
  MRNUMBER = {617913},
}

@book {Zaa,
    AUTHOR = {Zaanen, A. C.},
     TITLE = {Introduction to operator theory in {R}iesz spaces},
 PUBLISHER = {Springer, Berlin},
      YEAR = {1997},
     PAGES = {xii+312},
   MRCLASS = {47B60 (46A40 47B65)},
  MRNUMBER = {1631533},
MRREVIEWER = {Yu. A. Abramovich},
}

@book{Zworski2012,
  author    = {Zworski, M.},
  title     = {Semiclassical analysis},
  series    = {Graduate Studies in Mathematics},
  volume    = {138},
  publisher = {Amer. Math. Soc.},
  year      = {2012},
}
